\documentclass[twoside,11pt,dvipsnames, openleft]{book}
\usepackage{jmlr2e}
\usepackage{bm}
\usepackage{cancel}

\usepackage{pythonhighlight}
\usepackage{extarrows}

\renewcommand{\arraystretch}{1.15}

\usepackage[english]{babel}

\usepackage{imakeidx}
\makeindex[columns=2, title=Alphabetical Index]

\usepackage{stackengine}

\usepackage{dsfont}

\usepackage{extarrows}

\usepackage[framemethod=tikz]{mdframed}

\usepackage{amsmath}
\usepackage{arydshln}

\numberwithin{equation}{chapter}

\usepackage{blkarray}

\usepackage[absolute,overlay]{textpos}

\usepackage[final]{pdfpages}

\usepackage[noframe]{showframe}
\usepackage{framed}
\usepackage{lipsum}

{\endMakeFramed}

\usepackage{xcolor}

\usepackage{tcolorbox}
\usepackage{graphicx}             
\usepackage[hang]{subfigure}

\usepackage{algorithm}
\usepackage{algpseudocode}
\usepackage{graphics}
\usepackage{epsfig}

\usepackage{blkarray}
\usepackage{listings}

\usepackage[Bjornstrup]{fncychap}

\newcommand{\colortitlechap}{\color[RGB]{60,113,183}} 
\newcommand{\colornumberchap}{\color[RGB]{60,113,183}} 
\newcommand{\colorbackchap}{\colorbox[RGB]{200,200,200}} 

\makeatletter
\renewcommand{\DOCH}{%
		\settowidth{\py}{\CNoV\thechapter}
		\addtolength{\py}{-10pt}%
		\fboxsep=0pt%
		\colorbackchap{\rule{0pt}{40pt}\parbox[b]{\textwidth}{\hfill}}%
		\kern-\py\raise20pt%
		\hbox{\colornumberchap\CNoV\thechapter}\\%
	}

\renewcommand{\DOTI}[1]{%
		\nointerlineskip\raggedright%
		\fboxsep=\myhi%
		\vskip-1ex%
		\colorbackchap{\parbox[t]{\mylen}{\CTV\FmTi{\colortitlechap#1}}}\par\nobreak%
		\vskip 40\p@%
	}

\renewcommand{\DOTIS}[1]{%
		\fboxsep=0pt%
		\colorbackchap{\rule{0pt}{40pt}\parbox[b]{\textwidth}{\hfill}}\\%
		\nointerlineskip\raggedright%
		\fboxsep=\myhi%
		\colorbackchap{\parbox[t]{\mylen}{\CTV\FmTi{\colortitlechap#1}}}\par\nobreak%
		\vskip 40\p@%
	}
\makeatother

\usepackage{fancyhdr, blindtext}

\newcommand{\changefonts}{%
	\fontsize{9}{11}\selectfont
}

\fancypagestyle{plain}{%

	\fancyhf{}  
}

\usepackage[colorlinks]{hyperref}
\usepackage{hyperref}
\hypersetup{
    colorlinks=true, 
    linktoc=all,     
    linkcolor=mydarkblue,  
    anchorcolor=blue,
    citecolor=mydarkgreen,
}
\usepackage[hyperpageref]{backref} 

\usepackage{booktabs}  

\usepackage{bookmark}
\bookmarksetup{
	numbered,
	open
}

\usepackage{tikz}
\usetikzlibrary{mindmap,trees,backgrounds}
\usetikzlibrary{arrows.meta}
\usetikzlibrary{patterns}

\usetikzlibrary{decorations.text}

\usetikzlibrary{calc,backgrounds}

\newcommand*\circled[1]{\tikz[baseline=(char.base)]{
		\node[shape=circle,draw,inner sep=2pt] (char) {#1};}}
	
\usepackage{changepage}                 

\usepackage{hhline}

\usetikzlibrary{positioning,shapes,shadows,arrows}
\tikzstyle{condition}=[rectangle, draw=black, rounded corners, fill=colorqr, drop shadow,
text centered, anchor=north, text=black, text width=3cm]
\tikzstyle{abstract}=[rectangle, draw=black, rounded corners, fill=blue!30, drop shadow,
text centered, anchor=north, text=black, text width=3cm]
\tikzstyle{comment}=[rectangle, draw=black, rounded corners, fill=color1, drop shadow,
text centered, anchor=north, text=black, text width=3cm]
\tikzstyle{myarrow}=[->, >=open triangle 90, thick]
\tikzstyle{line}=[-, thick]

\usepackage{tikz-3dplot}
\tikzset{>=latex}
\usetikzlibrary{matrix}

\usepackage{sidecap}
\sidecaptionvpos{figure}{t}

\usepackage{verbatimbox}

\usepackage[mathscr]{euscript}

\usepackage{stmaryrd}

\newlength{\offsetpage}
\newenvironment{widepage}{\begin{adjustwidth}{-\offsetpage/2}{-\offsetpage/2}
		\addtolength{\textwidth}{\offsetpage}}%
	{\end{adjustwidth}}
	{\end{adjustwidth}}

\usepackage{minitoc}   
\let\cleardoublepage\clearpage
\usepackage[titletoc]{appendix}

\usepackage[labelfont=bf]{caption}
\newcommand\myhrulefill[1]{\leavevmode\leaders\hrule height#1\hfill\kern0pt}
\DeclareCaptionFormat{myformat}{
	{\color[RGB]{0,0,0}\myhrulefill{0.12em}}
	\\#1#2#3
}
\usepackage{afterpage}

\makeatletter \renewcommand*\cleardoublepage{
\clearpage
\if@twoside   
	\ifodd\c@page 
		\hbox{}\newpage
		\if@twocolumn\hbox{}   
			\newpage
		\fi
	\fi
\fi
} \makeatother
\let\originalpart=\part
\def\part#1{\cleardoublepage\clearpage \pagecolor{\partcolor} \originalpart{#1}\nopagecolor }

\usepackage{yfonts,color,lettrine}
\definecolor{caligraphcolor}{HTML}{74AECB}

\usepackage{setspace}

\AtBeginDocument{\setlength{\DefaultFindent}{0.5em}}
\def\algoalign#1{\parbox[t]{\dimexpr\linewidth-\algorithmicindent}{#1}}

\firstpageno{1}

\newcommand{\projectS}{\mathcal{P}_{\mathbb{S}}}
\newcommand{\projectT}{\mathcal{P}_{\mathbb{T}}}
\newcommand{\project}{\mathcal{P}}
\newcommand{\bprox}{\mathcal{B}\mathrm{rox}}

\newcommand{\bproxfphi}{\mathcal{B}\mathrm{rox}_{f,\phi}}
\newcommand{\bproxtphi}{\mathcal{B}\mathrm{rox}_{\widetilde{f},\phi}}
\newcommand{\bproxhphi}{\mathcal{B}\mathrm{rox}_{\widehat{f},\phi}}

\newcommand{\prox}{\mathcal{P}\mathrm{rox}}
\newcommand{\proxf}{\mathcal{P}\mathrm{rox}_f}

\newcommand{\indicatorS}{\delta_{\mathbb{S}}}
\newcommand{\indicatorG}{\delta}
\newcommand{\topone}{{(1)}}
\newcommand{\topzero}{{(0)}}
\newcommand{\toptminus}{{(t-1)}}
\newcommand{\toptminusTOP}{{(t-1)\top}}
\newcommand{\toptzero}{{(t)}}
\newcommand{\toptzeroTOP}{{(t)\top}}

\newcommand{\toptone}{{(t+1)}}
\newcommand{\toptoneTOP}{{(t+1)\top}}
\newcommand{\fbest}{f_{\text{best}}}
\newcommand{\Fbest}{F_{\text{best}}}
\newcommand{\extreal}{\widehat{\mathbb{R}}}
\newcommand{\dom}{\mathrm{dom}}
\newcommand{\epi}{\mathrm{epi}}
\newcommand{\lev}{\mathrm{Lev}}
\newcommand{\closure}{\text{cl}}
\newcommand{\interior}{\text{int}}
\newcommand{\relint}{\text{relint}}
\newcommand{\aff}{\text{aff}}

\newcommand{\proxgradL}{\mathcal{T}_{L}}

\newcommand{\gradmapL}{\mathcal{G}_{L}}

\newcommand{\fejer}{{Fej\'er}}

\newcommand{\supp}{\mathrm{supp}}
\newcommand\comple[1]{\cancel{#1}}
\newcommand{\lfd}{\sF_{\text{LFD}}}
\newcommand{\fd}{\sF_{\text{FD}}}
\newcommand{\sfd}{\sF_{\text{SFD}}}

\newcommand{\holders}{\text{H{\"o}lder's} }
\newcommand{\cond}{\text{cond} }

\newcommand{\gapforall}{\,\,}
\newcommand{\gap}{\,\,\,\,\,\,\,\,}  
\mathchardef\mhyphen="2D

\newcommand{\exampbar}{\hfill $\square$\par}
\newcommand{\argmax}{\operatorname*{\text{arg max}}}
\newcommand{\argmin}{\operatorname*{\text{arg min}}}

\newcommand\mathopmin[1]{\mathop{\min}_{#1}}
\newcommand\mathoplim[1]{\mathop{\lim}_{#1}}

\newcommand{\trans}[1]{\ensuremath{#1^{ \top}}}

\newcommand{\rms}{\text{RMS}}

\def\ceil#1{\lceil #1 \rceil}    
\def\floor#1{\lfloor #1 \rfloor}

\newcommand{\indicator}{\mathds{1}}
\def\1{\bm{1}}

\newcommand{\conv}{\text{conv}}
\newcommand{\cone}{\text{cone}}

\newcommand\inner[2]{\left\langle#1, #2\right\rangle}
\newcommand\abs[1]{\left\lvert#1\right\rvert}
\newcommand\absbig[1]{\big\lvert#1\big\rvert}
\newcommand\norm[1]{\left\lVert#1\right\rVert}
\newcommand\normbig[1]{\big\lVert#1\big\rVert}
\newcommand\normzero[1]{\left\lVert#1\right\rVert_0}
\newcommand\normone[1]{\left\lVert#1\right\rVert_1}
\newcommand\normonebig[1]{\big\lVert#1\big\rVert_1}

\newcommand\normtwo[1]{\left\lVert#1\right\rVert_2}
\newcommand\normtwobig[1]{\big\lVert#1\big\rVert_2}

\newcommand\normf[1]{\left\lVert#1\right\rVert_F}
\newcommand\normfbig[1]{\big\lVert#1\big\rVert_F}
\newcommand\norminf[1]{\left\lVert#1\right\rVert_{\infty}}

\newcommand\normA[1]{\left\lVert#1\right\rVert_{\bA}}
\newcommand\normB[1]{\left\lVert#1\right\rVert_{\bB}}
\newcommand\normC[1]{\left\lVert#1\right\rVert_{\bC}}

\newcommand\normH[1]{\left\lVert#1\right\rVert_{\bH}}
\newcommand\normZ[1]{\left\lVert#1\right\rVert_{\bZ}}
\newcommand\normAbig[1]{\big\lVert#1\big\rVert_{\bA}}

\newcommand\normCbig[1]{\big\lVert#1\big\rVert_{\bC}}

\newcommand\normZbig[1]{\big\lVert#1\big\rVert_{\bZ}}
\newcommand\innerproduct[1]{\left\langle#1\right\rangle}

\newcommand{\longeqf}{\scalebox{4}[1]{=}}

\newcommand\longeqtextfour[1]{\overset{#1}{\longeqf}}

\newcommand{\R}{\mathbb{R}}
\newcommand{\real}{\mathbb{R}}

\newcommand{\symmetric}{\mathbb{S}}
\newcommand{\psd}{\mathbb{S}_{+}}
\newcommand{\pd}{\mathbb{S}_{++}}
\newcommand{\nsd}{\mathbb{S}_{-}}
\newcommand{\nd}{\mathbb{S}_{--}}
\newcommand{\integer}{\mathbb{Z}}
\newcommand{\naturalset}{\mathbb{N}}

\newcommand{\KL}{D_{\mathrm{KL}}}

\definecolor{titlepagecolor}{cmyk}{75,68,67,90}
\definecolor{titlepagecolor2}{rgb}{1.0, 0.08, 0.58}
\definecolor{emerald}{rgb}{0.31, 0.78, 0.47}
\definecolor{deeppink}{HTML}{D14064}
\definecolor{lowpink}{HTML}{ffe6ec}
\newcommand{\partcolor}{gray!65} 
\definecolor{lowblue}{HTML}{E1EBFE}

\let\oldforall\forall
\renewcommand{\forall}{\oldforall\, }

\usepackage{color}   

\definecolor{mylightbluetitle}{RGB}{60,113,183}
\definecolor{mylightbluetext}{rgb}{0,0.08,0.45}
\definecolor{structurecolorblue}{RGB}{60,113,183}
\definecolor{structurecolorgreen}{RGB}{63,145,182}
\colorlet{structurecolor}{structurecolorblue}
\definecolor{structurecolorelegant}{RGB}{60,113,183}
\definecolor{structurecolorlt}{RGB}{31,119,185}

\definecolor{structurecolorHighTheoremBlue}{RGB}{220,227,248}
\definecolor{structurecolorHighTheoremGreen}{RGB}{188,222,231}
\colorlet{structurecolorHighTheorem}{structurecolorHighTheoremBlue}
\newcommand{\mdframecolorTheorem}{gray!35}  

\definecolor{winestain}{rgb}{0.5,0,0}
\definecolor{mydarkblue}{rgb}{0,0.08,0.45}
\definecolor{mydeepblue}{rgb}{0,0.08,0.65}
\definecolor{mydarkred}{rgb}{0.70,0.00,0.00}
\definecolor{mydarkgreen}{rgb}{0.00,0.30,0.00}
\definecolor{mydarkyellow}{RGB}{197,151,13}
\definecolor{mydarkpurple}{RGB}{90,35,140}
\definecolor{mydarkgray}{RGB}{64,64,64}

\definecolor{color0}  {RGB}{174,225,254} 
\definecolor{color1}  {RGB}{220,227,248} 
\definecolor{color2}  {RGB}{28,130,185} 
\definecolor{color3}  {RGB}{255,253,250} 
\definecolor{colormiddleright}  {RGB}{245,253,250} 
\definecolor{colorbottomleft}  {RGB}{255,243,250} 
\definecolor{coloruppermiddle}  {RGB}{255,253,230} 
\definecolor{colormiddleleft}  {RGB}{255,244,237}
\definecolor{colorcr}  {RGB}{249,253,232} 
\definecolor{colorreduction}  {RGB}{255,235,254} 
\definecolor{colorqr}  {RGB}{254,221,199} 
\definecolor{colorbiconjugate}  {RGB}{251,149,161} 
\definecolor{colorsvd}  {RGB}{215,247,235} 
\definecolor{colorupperright}  {RGB}{239,246,251} 
\definecolor{colorspectral}  {RGB}{206,226,243} 
\definecolor{colorbottomright}  {RGB}{220,224,236} 
\definecolor{coloreigenvalue}  {RGB}{197,203,224} 
\definecolor{colorcp} {RGB}{217, 234, 186} 
\definecolor{colorcpborder} {RGB}{233, 243, 216} 
\definecolor{colorupperleft}  {RGB}{235,243,240} 
\definecolor{colorsemidefinite}  {RGB}{217,232,226} 
\definecolor{colormiddle} {RGB}{235, 240,255}
\definecolor{colorlu}  {RGB}{220,227,255} 
\definecolor{colorals}  {RGB}{240,230,255} 
\definecolor{coloralsbkg}  {RGB}{248,243,255} 
\definecolor{canaryyellow}{rgb}{1.0, 0.75, 0.0}
\definecolor{bluepigment}{rgb}{0.0, 0.0, 1.0}
\definecolor{canarypurple}{RGB}{208, 13, 241}
\definecolor{colorGreenOcre}{RGB}{51,102,0} 
\definecolor{colorBlue2}{RGB}{200,207,248}
\definecolor{shadecolor}{gray}{0.75}

\newcommand{\hadaprod}{\circledast}

\newcommand{\Cov}{\mathbb{C}\mathrm{ov}}

\newcommand{\Exp}{\mathbb{E}}
\newcommand{\Var}{\mathbb{V}\mathrm{ar}}

\newcommand{\cspace}{\mathcal{C}}
\newcommand{\nspace}{\mathcal{N}}
\newcommand{\bzero}{\mathbf{0}}
\newcommand{\bone}{\mathbf{1}}
\newcommand{\diag}{\mathrm{diag}}
\newcommand{\defect}{\mathrm{def}}

\newcommand{\rank}{\mathrm{rank}}

\newcommand{\trace}{\mathrm{tr}}

\newcommand{\spn}{\mathrm{span}}

\newcommand{\domain}{\mathrm{dom}}

\DeclareMathOperator{\sign}{sign}   

\newcommand{\bgamma}{\boldsymbol\gamma}
\newcommand{\bepsilon}{\boldsymbol\epsilon}
\newcommand{\bLambda}{\boldsymbol\Lambda}

\newcommand{\bzeta}{\boldsymbol\zeta}
\newcommand{\blambda}{\boldsymbol\lambda}
\newcommand{\bxi}{\boldsymbol\xi}

\newcommand{\bmu}{\boldsymbol\mu}

\newcommand{\bSigma}{\boldsymbol\Sigma}

\newcommand{\widebarbx}{\overline{\bm{x}}}
\newcommand{\widebarby}{\overline{\bm{y}}}

\newcommand{\widehatf}{\widehat{f}}

\newcommand{\widehatx}{\widehat{x}}

\newcommand{\widehatbd}{\widehat{\bm{d}}}
\newcommand{\widehatbe}{\widehat{\bm{e}}}

\newcommand{\widehatbx}{\widehat{\bm{x}}}
\newcommand{\widehatby}{\widehat{\bm{y}}}

\newcommand{\mathcalC}{\mathcal{C}}
\newcommand{\mathcalD}{\mathcal{D}}
\newcommand{\mathcalE}{\mathcal{E}}

\newcommand{\mathcalG}{\mathcal{G}}

\newcommand{\mathcalI}{\mathcal{I}}

\newcommand{\mathcalL}{\mathcal{L}}

\newcommand{\mathcalO}{\mathcal{O}}
\newcommand{\mathcalP}{\mathcal{P}}

\newcommand{\mathcalT}{\mathcal{T}}

\newcommand{\mathcalV}{\mathcal{V}}
\newcommand{\mathcalW}{\mathcal{W}}
\newcommand{\mathcalX}{\mathcal{X}}

\newcommand{\widetildeblambda}{\widetilde{\bm{\lambda}}}
\newcommand{\widetildebA}{\widetilde{\bm{A}}}
\newcommand{\widetildebB}{\widetilde{\bm{B}}}

\newcommand{\widetildebS}{\widetilde{\bm{S}}}

\newcommand{\widetildebZ}{\widetilde{\bm{Z}}}

\newcommand{\widetildebb}{\widetilde{\bm{b}}}

\newcommand{\widetildebd}{\widetilde{\bm{d}}}

\newcommand{\widetildebg}{\widetilde{\bm{g}}}

\newcommand{\widetildebl}{\widetilde{\bm{l}}}

\newcommand{\widetildebu}{\widetilde{\bm{u}}}

\newcommand{\widetildebx}{\widetilde{\bm{x}}}
\newcommand{\widetildeby}{\widetilde{\bm{y}}}

\newcommand{\widetildeB}{\widetilde{B}}

\newcommand{\widetildeR}{\widetilde{R}}

\newcommand{\widetildeb}{\widetilde{b}}
\newcommand{\widetildec}{\widetilde{c}}

\newcommand{\widetildef}{\widetilde{f}}

\newcommand{\widetildet}{\widetilde{t}}

\newcommand{\widetildex}{\widetilde{x}}

\newcommand{\widetildez}{\widetilde{z}}

\newcommand{\ba}{\bm{a}}
\newcommand{\bA}{\bm{A}}
\newcommand{\bb}{\bm{b}}
\newcommand{\bB}{\bm{B}}
\newcommand{\bc}{\bm{c}}
\newcommand{\bC}{\bm{C}}  
\newcommand{\bd}{\bm{d}}
\newcommand{\bD}{\bm{D}}
\newcommand{\be}{\bm{e}}
\newcommand{\bE}{\bm{E}}
\newcommand{\bff}{\bm{f}}
\newcommand{\bF}{\bm{F}}
\newcommand{\bg}{\bm{g}}
\newcommand{\bG}{\bm{G}}
\newcommand{\bh}{\bm{h}}
\newcommand{\bH}{\bm{H}}

\newcommand{\bI}{\bm{I}}

\newcommand{\bJ}{\bm{J}}

\newcommand{\bl}{\bm{l}}
\newcommand{\bL}{\bm{L}}
\newcommand{\bmm}{\bm{m}}
\newcommand{\bM}{\bm{M}}
\newcommand{\bn}{\bm{n}}
\newcommand{\bN}{\bm{N}}

\newcommand{\bp}{\bm{p}}
\newcommand{\bP}{\bm{P}}
\newcommand{\bq}{\bm{q}}
\newcommand{\bQ}{\bm{Q}}
\newcommand{\br}{\bm{r}}
\newcommand{\bR}{\bm{R}}
\newcommand{\bs}{\bm{s}}
\newcommand{\bS}{\bm{S}}
\newcommand{\bt}{\bm{t}}
\newcommand{\bT}{\bm{T}}
\newcommand{\bu}{\bm{u}}
\newcommand{\bU}{\bm{U}}
\newcommand{\bv}{\bm{v}}
\newcommand{\bV}{\bm{V}}
\newcommand{\bw}{\bm{w}}
\newcommand{\bW}{\bm{W}}
\newcommand{\bx}{\bm{x}}
\newcommand{\bX}{\bm{X}}
\newcommand{\by}{\bm{y}}
\newcommand{\bY}{\bm{Y}}
\newcommand{\bz}{\bm{z}}
\newcommand{\bZ}{\bm{Z}}

\newcommand{\whbx}{\widehat{\bx}}
\newcommand{\whbd}{\widehat{\bd}}
\newcommand{\whbg}{\widehat{\bg}}

\newcommand{\whf}{\widehat{f}}

\def\vmu{{\bm{\mu}}}

\def\va{{\bm{a}}}

\def\vx{{\bm{x}}}

\def\mA{{\bm{A}}}
\def\mB{{\bm{B}}}

\def\mH{{\bm{H}}}

\def\mJ{{\bm{J}}}

\def\mX{{\bm{X}}}

\def\mSigma{{\bm{\Sigma}}}

\DeclareMathAlphabet{\mathsfit}{\encodingdefault}{\sfdefault}{m}{sl}
\SetMathAlphabet{\mathsfit}{bold}{\encodingdefault}{\sfdefault}{bx}{n}

\def\sA{{\mathbb{A}}}
\def\sB{{\mathbb{B}}}
\def\sC{{\mathbb{C}}}
\def\sD{{\mathbb{D}}}
\def\sF{{\mathbb{F}}}
\def\sG{{\mathbb{G}}}
\def\sH{{\mathbb{H}}}
\def\sI{{\mathbb{I}}}
\def\sJ{{\mathbb{J}}}
\def\sK{{\mathbb{K}}}
\def\sL{{\mathbb{L}}}

\def\sS{{\mathbb{S}}}
\def\sT{{\mathbb{T}}}

\def\sV{{\mathbb{V}}}

\def\sX{{\mathbb{X}}}

\def\ra{{\textnormal{a}}}
\def\rb{{\textnormal{b}}}
\def\rc{{\textnormal{c}}}

\def\rx{{\textnormal{x}}}

\def\rX{{\textnormal{X}}}
\def\rY{{\textnormal{Y}}}

\def\rva{{\mathbf{a}}}

\def\rvx{{\mathbf{x}}}

\def\rmA{{\mathbf{A}}}

\def\eva{{a}}

\makeatletter
\newcommand{\paragrapharrow}{%
	\@startsection{paragraph}{4}{\z@}%
	{3.25ex \@plus1ex \@minus.2ex}%
	{-1em}%
	{\normalfont\normalsize\bfseries$\blacktriangleright$\ }}
\makeatother

\usepackage[noframe]{showframe}
\usepackage{framed}
\usepackage{lipsum}

\newcommand{\roundcornertheorem}{0pt}
\newcommand{\linewidththeorem}{0.1pt}
\newcommand{\frametitlerulewidththeorem}{0.1pt}
\newcommand{\innerbottommargintheorem}{2pt}
\newcommand{\innerleftmargintheorem}{2pt}
\newcommand{\innerrightmargintheorem}{2pt}
\newcommand{\innertopmargintheorem}{2pt}
\newcommand{\outerlinewidththeorem}{1pt}

\newcommand{\textremark}{3pt}

\newenvironment{theoremHigh}[1][]{%
\refstepcounter{theo}%
\ifstrempty{#1}%
{\newcommand{\theoName}{}}
{\newcommand{\theoName}{:~(#1)}}
\mdfsetup{
	backgroundcolor=structurecolorHighTheorem,
	linecolor=structurecolorHighTheorem,
	frametitlerulewidth=\frametitlerulewidththeorem,
	roundcorner=\roundcornertheorem,
	linewidth=\linewidththeorem,
	innerbottommargin=\innerbottommargintheorem,
	innerleftmargin=\innerleftmargintheorem,
	innerrightmargin=\innerrightmargintheorem,
	innertopmargin=\innertopmargintheorem,
	outerlinewidth=\outerlinewidththeorem,
	topline=false,
	innertopmargin=-5pt,
	innerbottommargin=1pt,
	linewidth=0,
	startinnercode=\paragraph{{\strut Theorem~\thetheo\theoName}}
}
\begin{mdframed}[]\relax%
}{\end{mdframed}}

\newenvironment{theorem}[1][]{%
	\refstepcounter{theo}%
	\ifstrempty{#1}%
	{\newcommand{\theoName}{}}
	{\newcommand{\theoName}{:~(#1)}}
	\mdfsetup{
		backgroundcolor=\mdframecolorTheorem,
		linecolor=\mdframecolorTheorem,
		frametitlerulewidth=\frametitlerulewidththeorem,
		roundcorner=\roundcornertheorem,
		linewidth=\linewidththeorem,
		innerbottommargin=\innerbottommargintheorem,
		innerleftmargin=\innerleftmargintheorem,
		innerrightmargin=\innerrightmargintheorem,
		innertopmargin=\innertopmargintheorem,
		outerlinewidth=\outerlinewidththeorem,
		topline=false,
		innertopmargin=-5pt,
		innerbottommargin=1pt,
		linewidth=0,
		startinnercode=\paragraph{{\strut Theorem~\thetheo\theoName}}
	}
	\begin{mdframed}[]\relax%
	}{\end{mdframed}}
\newenvironment{corollary}[1][]{%
	\refstepcounter{theo}%
\ifstrempty{#1}%
{\newcommand{\theoName}{}}
{\newcommand{\theoName}{:~(#1)}}
\mdfsetup{
	backgroundcolor=\mdframecolorTheorem,
	linecolor=\mdframecolorTheorem,
	frametitlerulewidth=\frametitlerulewidththeorem,
	roundcorner=\roundcornertheorem,
	linewidth=\linewidththeorem,
	innerbottommargin=\innerbottommargintheorem,
	innerleftmargin=\innerleftmargintheorem,
	innerrightmargin=\innerrightmargintheorem,
	innertopmargin=\innertopmargintheorem,
	outerlinewidth=\outerlinewidththeorem,
	topline=false,
	innertopmargin=-5pt,
	innerbottommargin=1pt,
	linewidth=0,
	startinnercode=\paragraph{{\strut Corollary~\thetheo\theoName}}
}
	\begin{mdframed}[]\relax%
	}{\end{mdframed}}
\newenvironment{lemma}[1][]{%
	\refstepcounter{theo}%
\ifstrempty{#1}%
{\newcommand{\theoName}{}}
{\newcommand{\theoName}{:~(#1)}}
\mdfsetup{
	backgroundcolor=\mdframecolorTheorem,
	linecolor=\mdframecolorTheorem,
	frametitlerulewidth=\frametitlerulewidththeorem,
	roundcorner=\roundcornertheorem,
	linewidth=\linewidththeorem,
	innerbottommargin=\innerbottommargintheorem,
	innerleftmargin=\innerleftmargintheorem,
	innerrightmargin=\innerrightmargintheorem,
	innertopmargin=\innertopmargintheorem,
	outerlinewidth=\outerlinewidththeorem,
	topline=false,
	innertopmargin=-5pt,
	innerbottommargin=1pt,
	linewidth=0,
	startinnercode=\paragraph{{\strut Lemma~\thetheo\theoName}}
}
\begin{mdframed}[]\relax%
}{\end{mdframed}}

\newenvironment{proposition}[1][]{%
\refstepcounter{theo}%
\ifstrempty{#1}%
{\newcommand{\theoName}{}}
{\newcommand{\theoName}{:~(#1)}}
\mdfsetup{
	backgroundcolor=\mdframecolorTheorem,
	linecolor=\mdframecolorTheorem,
	frametitlerulewidth=\frametitlerulewidththeorem,
	roundcorner=\roundcornertheorem,
	linewidth=\linewidththeorem,
	innerbottommargin=\innerbottommargintheorem,
	innerleftmargin=\innerleftmargintheorem,
	innerrightmargin=\innerrightmargintheorem,
	innertopmargin=\innertopmargintheorem,
	outerlinewidth=\outerlinewidththeorem,
	topline=false,
	innertopmargin=-5pt,
	innerbottommargin=1pt,
	linewidth=0,
	startinnercode=\paragraph{{\strut Proposition~\thetheo\theoName}}
}
\begin{mdframed}[]\relax%
}{\end{mdframed}}

\usepackage{amsthm}
\newtheoremstyle{normalfontstyle} 
{3pt}                           
{3pt}                           
{\normalfont}                   
{}                              
{\bfseries}                     
{}                             
{ }                             
{}                              

\usepackage{thmtools}
\declaretheoremstyle[
spaceabove=3pt,
spacebelow=3pt,
headfont=\bfseries,
notefont=\bfseries, 
notebraces={(}{)}, 
bodyfont=\normalfont,
postheadspace=1em, 
]{normalfontboldhead}

\theoremstyle{normalfontstyle}
\newcommand{\BlackBox}{\rule{1.5ex}{1.5ex}}  
\renewenvironment{proof}{\par\noindent{\bf Proof\ }}{\hfill\BlackBox\\[2mm]}

\declaretheorem[style=normalfontboldhead, name=Definition, numberlike=theo]{definitionT}
\newmdenv[skipabove=7pt,
skipbelow=7pt,
rightline=false,
leftline=true,
topline=false,
bottomline=false,
linecolor=mydarkblue,
innerleftmargin=5pt,
innerrightmargin=5pt,
innertopmargin=0pt,
leftmargin=2cm,
rightmargin=0cm,
linewidth=4pt,
innerbottommargin=0pt]{dBox}
\newenvironment{definition}{\begin{dBox}\begin{definitionT}}{\end{definitionT}\end{dBox}}

\declaretheorem[style=normalfontboldhead, name=Exercise, numberlike=theo]{exerciseC}
\newmdenv[skipabove=7pt,
skipbelow=7pt,
rightline=false,
leftline=true,
topline=false,
bottomline=false,
linecolor=mydarkgreen,
innerleftmargin=5pt,
innerrightmargin=5pt,
innertopmargin=0pt,
leftmargin=2cm,
rightmargin=0cm,
linewidth=4pt,
innerbottommargin=0pt]{eBox}
\newenvironment{exercise}{\begin{eBox}\begin{exerciseC}}{\end{exerciseC}\end{eBox}}

\declaretheorem[style=normalfontboldhead, name=Remark, numberlike=theo]{remarekC}
\newmdenv[skipabove=7pt,
skipbelow=7pt,
rightline=false,
leftline=true,
topline=false,
bottomline=false,
linecolor=mydarkpurple,
innerleftmargin=5pt,
innerrightmargin=5pt,
innertopmargin=0pt,
leftmargin=2cm,
rightmargin=0cm,
linewidth=4pt,
innerbottommargin=0pt]{rBox}
\newenvironment{remark}{\begin{rBox}\begin{remarekC}}{\end{remarekC}\end{rBox}}

\declaretheorem[style=normalfontboldhead, name=Assumption, numberlike=theo]{assumptionC}
\newmdenv[skipabove=7pt,
skipbelow=7pt,
rightline=false,
leftline=true,
topline=false,
bottomline=false,
linecolor=mydarkpurple,
innerleftmargin=5pt,
innerrightmargin=5pt,
innertopmargin=0pt,
leftmargin=2cm,
rightmargin=0cm,
linewidth=4pt,
innerbottommargin=0pt]{asBox}
\newenvironment{assumption}{\begin{asBox}\begin{assumptionC}}{\end{assumptionC}\end{asBox}}

\declaretheorem[style=normalfontboldhead, name=Example, numberlike=theo]{exampleC}
\newmdenv[skipabove=7pt,
skipbelow=7pt,
rightline=false,
leftline=false,
topline=false,
bottomline=false,
linecolor=mydarkgreen,
innerleftmargin=1pt,
innerrightmargin=5pt,
innertopmargin=0pt,
leftmargin=2cm,
rightmargin=0cm,
linewidth=4pt,
innerbottommargin=0pt]{xBox}
\newenvironment{example}{\begin{xBox}\begin{exampleC}}{\exampbar\end{exampleC}\end{xBox}}

\usepackage{adforn}  
\newcommand{\xchaptertitle}{Chapter~\thechapter~}
\newcommand{\problemname}{Problems}
\newenvironment{problemset}[1][\xchaptertitle~\problemname]{
\vspace*{10pt}
\begin{center}
\phantomsection\addcontentsline{toc}{section}{\texorpdfstring{\xchaptertitle~\problemname}{\problemname}}
\markright{#1}
\textcolor{structurecolor}{\Large\bfseries\adftripleflourishleft~#1~\adftripleflourishright}
\end{center}
\begin{enumerate}[ref=\thechapter.\theenumi]}{
\end{enumerate}}

\usepackage[shortlabels]{enumitem}  
\usepackage{enumitem}
\newcommand*{\eitemi}{\tikz \draw [baseline, ball color=structurecolor,draw=none] circle (2pt);}
\newcommand*{\eitemii}{\tikz \draw [baseline, fill=structurecolor,draw=none,circular drop shadow] circle (2pt);}
\newcommand*{\eitemiii}{\tikz \draw [baseline, fill=structurecolor,draw=none] circle (2pt);}
\setlist[enumerate,1]{label=\color{black}\arabic*.,itemsep=0pt,partopsep=0pt,parsep=\parskip,topsep=3pt}
\setlist[enumerate,2]{label=\color{black}(\alph*).,itemsep=0pt,partopsep=0pt,parsep=\parskip,topsep=3pt}
\setlist[enumerate,3]{label=\color{black}\Roman*.,itemsep=0pt,partopsep=0pt,parsep=\parskip,topsep=3pt}
\setlist[enumerate,4]{label=\color{black}\Alph*.,itemsep=0pt,partopsep=0pt,parsep=\parskip,topsep=3pt}
\setlist[itemize,1]{label={\eitemi},itemsep=0pt,partopsep=0pt,parsep=\parskip,topsep=3pt}
\setlist[itemize,2]{label={\eitemii},itemsep=0pt,partopsep=0pt,parsep=\parskip,topsep=3pt}
\setlist[itemize,3]{label={\eitemiii},itemsep=0pt,partopsep=0pt,parsep=\parskip,topsep=3pt}

\begin{document}

\frontmatter
\pagecolor{lowblue}
\newpage
\thispagestyle{empty}
\begingroup
\hypersetup{
	colorlinks=true, 
	linktoc=all,     
	linkcolor=mydarkblue,  
	anchorcolor=blue,
	citecolor=mydarkgreen,
	urlcolor=mydarkblue,
}
%
%
%
\endgroup
\clearpage

\clearpage
\nopagecolor 

\thispagestyle{empty}  
\title{Practical Topics in Optimization}

\author{
\begin{center}
\name Jun Lu \\ 
\email jun.lu.locky@gmail.com \\
\copyright \,\, 2025$\sim$ \,\, Jun Lu\\
\end{center}
}

\maketitle

\chapter*{\centering \begin{normalsize}Preface\end{normalsize}}

In an era where data-driven decision-making and computational efficiency are paramount, optimization plays a foundational role in advancing fields such as mathematics, computer science, operations research, machine learning, and beyond. From refining machine learning models to improving resource allocation and designing efficient algorithms, optimization techniques serve as essential tools for tackling complex problems. This book aims to provide both an introductory guide and a comprehensive reference, equipping readers with the necessary knowledge to understand and apply optimization methods within their respective fields.

Our primary goal is to demystify the inner workings of optimization algorithms, including black-box and stochastic optimizers, by offering both formal and intuitive explanations. Starting from fundamental mathematical principles, we derive key results to ensure that readers not only learn how these techniques work but also understand when and why to apply them effectively. By striking a careful balance between theoretical depth and practical application, this book serves a broad audience, from students and researchers to practitioners seeking robust optimization strategies.

A significant focus is placed on gradient descent, one of the most widely used optimization algorithms in machine learning. The text examines basic gradient descent, momentum-based methods, conjugate gradient techniques, and stochastic optimization strategies, highlighting their historical development and modern applications. Special attention is given to stochastic gradient descent (SGD) and its pivotal role in optimizing deep neural network or transformer structures, particularly in overcoming computational challenges and escaping saddle points. The book traces the evolution of these methods from the pioneering work of Robbins and Monro in 1951 to their current status as core technologies in deep learning or large language models.

Beyond gradient-based techniques, the book explores modern optimization methods such as proximal algorithms, augmented Lagrangian methods, the alternating direction method of multipliers (ADMM), and trust region approaches. It also discusses least squares problems, sparse optimization, and learning rate adaptation strategies, all of which play critical roles in large-scale and high-dimensional optimization settings.



\paragraph{Keywords.}
Optimality conditions, First-order methods, Second-order methods, Newton's method, Gradient descent, Stochastic gradient descent, Steepest descent, Greedy search, Conjugate descent and conjugate gradient, Learning rate annealing, Adaptive learning rate.

\vspace{5em}
\noindent

\newpage
\begingroup
\hypersetup{
linkcolor=structurecolor,
linktoc=page,  
}
\dominitoc
\pdfbookmark{\contentsname}{toc} 
\tableofcontents 
\endgroup

\chapter*{Notation}\label{notation}
\index{Notation}


This section provides a concise reference describing notation used throughout this
book.
If you are unfamiliar with any of the corresponding mathematical concepts,
the book describes most of these ideas in Chapter~\ref{chapter_introduction}.

\noindent
\begin{minipage}{\textwidth}
\centerline{\bf General Notations}
\bgroup
\def\arraystretch{1.5}
\begin{tabular}{cp{4.25in}}
$\displaystyle \triangleq$ & equals by definition \\ 
$\displaystyle :=, \leftarrow$ & equals by assignment \\ 
$\displaystyle \pi$ &	 3.141592$\ldots$ \\
$\displaystyle e, \exp$ &	 2.71828$\ldots$
\end{tabular}
\egroup
\end{minipage}

\noindent
\begin{minipage}{\textwidth}
\centerline{\bf Numbers and Arrays}
\bgroup
\def\arraystretch{1.5}
\begin{tabular}{cp{4.25in}}
$\displaystyle a$   & A scalar (integer or real)\\
$\displaystyle \ba$ & A vector\\
$\displaystyle \bone$ & A all-ones vector\\
$\displaystyle \bA$ & A matrix\\
$\displaystyle \bI_n$ & Identity matrix with $n$ rows and $n$ columns\\
$\displaystyle \bI$   & Identity matrix with dimensionality implied by context\\
$\displaystyle \be_i$ & Standard basis/canonical orthonormal  vector $[0,\dots,0,1,0,\dots,0]$ with a 1 at position $i$\\
$\displaystyle \text{diag}(\va)$ & A square, diagonal matrix with diagonal entries given by $\ba$\\
$\displaystyle \text{diag}(\bA, \bB,\ldots)$&  a block-diagonal matrix=($\bA\oplus \bB\oplus \ldots$), the direct sum notation\\
$\displaystyle \ra$   & A scalar random variable\\
$\displaystyle \rva$  & A vector-valued random variable\\
$\displaystyle \rmA$  & A matrix-valued random variable\\
\end{tabular}
\egroup
\index{Scalar}
\index{Vector}
\index{Matrix}
\index{Tensor}
\end{minipage}

\index{Sets}
\vspace{0.2in}
\noindent
\begin{minipage}{\textwidth}
\centerline{\bf Sets}
\bgroup
\def\arraystretch{1.5}
\begin{tabular}{cp{4.25in}}
$\displaystyle \sS, \abs{\sS}, \comple{\sS}$ & A set $\sS$, its cardinality, and its complement\\
$\displaystyle \varnothing$ & The null set \\
$\displaystyle \real, \real_+, \real_{++},\extreal $ & The set of real/nonnegative/positive/extended real numbers \\
$\displaystyle \naturalset, \integer$ & The set of natural/integer numbers \\
$\sB_p(\bc, r), \sB_p[\bc, r]$ & Open/close ball: $ \{\bx\in\real^n: \norm{\bx-\bc}_p <r \text{ or } \norm{\bx-\bc}_p \leq r \}$ \\
$\displaystyle \sF$ & The set of either real or complex numbers \\
$\displaystyle \{0, 1\}$ & The set containing 0 and 1 \\
$\displaystyle \{0, 1, \ldots, n \}$ & The set of all integers between $0$ and $n$\\
$\displaystyle [a, b]$ & The real interval including $a$ and $b$\\
$\displaystyle (a, b]$ & The real interval excluding $a$ but including $b$\\
$\displaystyle \sA \backslash \sB$ & Set subtraction, i.e., the set containing the elements of $\sA$ that are not in $\sB$\\

$\displaystyle \lev[f, \alpha]$  & Level set of a function: $\{ \bx \in \real^n \mid f(\bx) \leq \alpha \}$\\
\end{tabular}
\egroup
\index{Scalar}
\index{Vector}
\index{Matrix}
\index{Tensor}

\vspace{0.4in}
\noindent
\begin{minipage}{\textwidth}
\centerline{\bf Calculus}
\bgroup
\def\arraystretch{1.5}
\begin{tabular}{cp{4.25in}}
$\displaystyle\frac{d y} {d x}$ & Derivative of $y$ with respect to $x$\\ [2ex]
$\displaystyle \frac{\partial y} {\partial x} $ & Partial derivative of $y$ with respect to $x$ \\
$\displaystyle  f^\prime(\bx; \bd)$ or $D_{\bd}f(\bx)$& Directional derivative of $f$\\
$\displaystyle \frac{\partial f}{\partial x_i} (\bx)$, $D_{\be_i}f(\bx)$, or $\partial_i f(\bx)$ & $i$-th partial derivative of $f$\\
$\displaystyle \nabla_{\bx} f(\bx), \nabla f(\bx), \text{ or } \bg$ & Gradient of $f$ with respect to $\bx$ \\
$\displaystyle \nabla_{\bX} y $ & Matrix derivatives of $y$ with respect to $\bX$ \\
$\displaystyle \frac{\partial f}{\partial \vx}, \nabla f(\bx), \bJ_f(\bx), \text{ or } \bJ(\bx) $ & Jacobian matrix $\mJ \in \R^{m\times n}$ of $f: \R^n \rightarrow \R^m$\\
$\displaystyle \nabla_\vx^2 f(\vx), \nabla^2 f(\vx), \bH, \text{ or }\mH( f)(\vx)$ & The Hessian matrix of $f$ at input point $\vx$\\
$\displaystyle \int f(\vx) d\vx $ & Definite integral over the entire domain of $\vx$ \\
$\displaystyle \int_\sS f(\vx) d\vx$ & Definite integral with respect to $\vx$ over the set $\sS$ \\
\end{tabular}
\egroup
\index{Derivative}
\index{Integral}
\index{Jacobian matrix}
\index{Hessian matrix}
\end{minipage}

\index{Set}
\end{minipage}

\index{Matrix indexing}
\vspace{0.2in}
\noindent
\begin{minipage}{\textwidth}
\centerline{\bf Indexing}
\bgroup
\def\arraystretch{1.5}
\begin{tabular}{cp{4.25in}}
	$\displaystyle \eva_i$ & Element $i$ of vector $\va$, with indexing starting at 1 \\
	$\displaystyle \ba_{-i}$ & All elements of vector $\va$ except for element $i$ \\
	$\displaystyle  \bA[i,j]=a_{ij}$ & Element $(i, j)$ of matrix $\mA$ \\
	$\displaystyle \mA_{i, :}=\bA[i,:]$ & Row $i$ of matrix $\mA$ \\
	$\displaystyle \mA_{:, i}=\bA[:,i]=\ba_i$ & Column $i$ of matrix $\mA$ \\
\end{tabular}
\egroup
\end{minipage}

\vspace{0.4in}
\noindent
\begin{minipage}{\textwidth}
\centerline{\bf Probability and Information Theory}
\bgroup
\def\arraystretch{1.5}
\begin{tabular}{cp{4.25in}}
$\displaystyle \ra \bot \rb$ & The random variables $\ra$ and $\rb$ are independent\\
$\displaystyle \ra \bot \rb \mid \rc $ & They are conditionally independent given $\rc$\\
$\displaystyle P(\ra)$ & A probability distribution over a discrete variable\\
$\displaystyle p(\ra)$ & A probability distribution over a continuous variable, or over
a variable whose type has not been specified\\
$\displaystyle \ra \sim P$ & Random variable $\ra$ has distribution $P$\\
$\displaystyle  \Exp_{\rx\sim P} [ f(x) ]\text{ or } \Exp [f(x)]$ & Expectation of $f(x)$ with respect to $P(\rx)$ \\
$\displaystyle \Var[f(x)] $ &  Variance of $f(x)$ under $P(\rx)$ \\
$\displaystyle \Cov[f(x),g(x)] $ & Covariance of $f(x)$ and $g(x)$ under $P(\rx)$\\
$\displaystyle H(\rx) $ & Shannon entropy of the random variable $\rx$\\
$\displaystyle \KL ( P \Vert Q ) $ & Kullback-Leibler divergence of P and Q \\
$\displaystyle \mathcal{N} ( \vx ; \vmu , \mSigma)$ & Gaussian distribution %
over $\vx$ with mean $\vmu$ and covariance $\mSigma$ \\
\end{tabular}
\egroup
\index{Independence}
\index{Conditional independence}
\index{Variance}
\index{Covariance}
\index{Kullback-Leibler divergence}
\index{Shannon entropy}
\end{minipage}

\vspace{0.4in}
\noindent
\begin{minipage}{\textwidth}
\centerline{\bf Functions}
\bgroup
\def\arraystretch{1.5}
\begin{tabular}{cp{3.25in}}
$\displaystyle f: \sA \rightarrow \sB$ & The function $f$ with domain $\sA$ and range $\sB$\\
$\displaystyle f^\prime(\bx; \bd), \nabla f(\bx)$ & Directional derivative, gradient \\
$\displaystyle f'(\bx), \partial f(\bx)$ & Subgradient, subdifferential \\
$\displaystyle \indicatorS(\bx)$ & Indicator function, $\indicatorS(\bx)=0$ if $\bx\in\sS$, and $\indicatorS(\bx)=\infty$ otherwise\\
$\displaystyle \indicator\{\mathrm{condition}\}$ & is 1 if the condition is true, 0 otherwise\\
$\projectS(\bx), \proxf(\bx), \bproxfphi(\bx)$ & Projection/proximal/Bregman-proximal operator\\
$\proxgradL^{f,g}(\bx)\triangleq\prox_{\frac{1}{L} g}\left(\bx - \frac{1}{L}\nabla f(\bx)\right)$ & Prox-grad operator\\
$\gradmapL^{f,g}(\bx)\triangleq L\big(\bx-\proxgradL^{f,g}(\bx)\big)$ & Gradient mapping\\
$\displaystyle f \circ g $ & Composition of the functions $f$ and $g$ \\
$\displaystyle \log(x), \ln(x)$ & Natural logarithm of $x$ \\
$\displaystyle \sigma(x),\, \text{Sigmoid}(x)$ & Logistic sigmoid, i.e., $\displaystyle \frac{1} {1 + \exp\{-x\}}$ \\
$\displaystyle \phi(\eta)$ & $\phi(\eta)\triangleq f(\bx^\toptzero - \eta \bg^\toptzero)$ for line search algorithms\\
$\displaystyle \psi_t(\bd)$ &  $\psi_t(\bd) \triangleq f(\bx^\toptzero) + \bd^\top \bg^\toptzero + \frac{1}{2} \bd^\top \nabla^2 f(\bx^\toptzero) \bd$ for trust region methods\\
$\displaystyle \norm{\bx}_p $ & $\ell_p$ norm of $\vx$ \\
$\displaystyle \norm{\bx}_1, \norm{\bx}_2, \norm{\bx}_\infty$ & $\ell_1$ norm, $\ell_2$ norm, $\ell_\infty$ norm\\
$\displaystyle \normf{\bA}, \normtwo{\bA}$ & Frobenius and spectral norms \\
$\displaystyle  [x]_+$ & Positive part of $x$, i.e., $\max(0,x)$\\
\end{tabular}
\egroup
\index{Sigmoid}
\index{Softplus}
\index{Norm}
\end{minipage}
Sometimes we use a function $f$ whose argument is a scalar but apply
it to a vector or a matrix: $f(\vx)$ or $f(\mX)$.
This denotes the application of $f$ to the
array element-wise. For example, if $\by = [\bx]_+$, then $y_i = [x_i]_+$ for all valid values of $i$.

\vspace{0.2in}
\noindent
\begin{minipage}{\textwidth}
\centerline{\bf Linear Algebra Operations}
\bgroup
\def\arraystretch{1.5}
\begin{tabular}{cp{4.25in}}
$\displaystyle \abs{\bA}$ & Componentwise absolute matrix $\mA$ \\
$\displaystyle \frac{[\bA]}{[\bB]}$ &  Componentwise division between two matrices\\
$\displaystyle \bA^\top$ & Transpose of matrix $\mA$ \\
$\displaystyle \bA^{-1}$ & Inverse of $\mA$\\
$\displaystyle \bA^+$ & Moore-Penrose pseudo-inverse of $\mA$\\
$\displaystyle \bA \circledast \bB $ & Element-wise (Hadamard) product of $\mA$ and $\mB$ \\
$\displaystyle \det(\bA)$ & Determinant of $\bA$ \\
$\displaystyle \trace(\bA)$ & Trace of $\bA$ \\
$\displaystyle \mathrm{rref}(\bA)$ & Reduced row echelon form of $\bA$ \\
$\displaystyle \cspace(\bA)$ & Column space of $\bA$ \\
$\displaystyle \nspace(\bA)$ & Null space of $\bA$ \\
$\displaystyle \mathcalV, \mathcalW$ & A general subspace \\
$\displaystyle \dim(\mathcalV)$ & Dimension of the space $\mathcalV$ \\
$\displaystyle \defect(\bA)$ & Defect or nullity of $\bA$ \\
$\displaystyle \rank(\bA)$ & Rank of $\bA$ \\
$\displaystyle \trace(\bA)$ & Trace of $\bA$ \\
$\displaystyle \lambda_{\min}(\bA),\lambda_{\max}(\bA)$ & Smallest/largest eigenvalues of $\bA$\\
$\displaystyle \sigma_{\min}(\bA),\sigma_{\max}(\bA)$ & Smallest/largest singuar values of $\bA$\\
$\displaystyle d_{\min}(\bA),d_{\max}(\bA)$ & Smallest/largest diagonal values of $\bA$\\
$\displaystyle \bA\geq \bzero, \bA>\bzero$ & Nonnegative/positive matrix $\bA$\\
$\displaystyle \bA\succeq \bzero\, (\bA\in\psd^n), \bA\succ \bzero\, (\bA\in\pd^n)$ & Positive semidefinite/definite matrix $\bA$\\
subdiagonal/superdiagonal & Entries below/above the main diagonal
\end{tabular}
\egroup
\index{Transpose}
\index{Element-wise product|see {Hadamard product}}
\index{Hadamard product}
\index{Determinant}
\end{minipage}


\vspace{0.4in}
\noindent
\begin{minipage}{\textwidth}
\centerline{\bf Abbreviations}
\bgroup
\def\arraystretch{1.5}
\begin{tabular}{cp{4.25in}}
PD & Positive definite  \\
PSD & Positive semidefinite \\
SC, SS & Strongl convexity, Strong smoothness\\
RSC, RSS & Restricted strong convexity, Restricted strong smoothness\\
KKT & Karush-Kuhn-Tucker conditions\\
SVD &  Singular value decomposition \\
PCA & Principal component analysis \\
OLS & Ordinary least squares\\
LS & Least squares \\
GD & Gradient descent  method\\
PGD & Projected (sub)gradient descent  method\\
CG & Conditional gradient method \\
GCG & Generalized conditional gradient method \\
FISTA & Fast Proximal Gradient method\\
SGD & Stochastic gradient descent \\
LM & Levenberg-Marquardt method\\
ADMM & Alternating direction methods of multipliers\\
ADPMM & Alternating direction proximal methods of multipliers\\
IHT & Iterative hard-thresholding method \\
LASSO & Least absolute selection and shrinkage operator\\

\end{tabular}
\egroup
\end{minipage}

\clearpage


\mainmatter
\newpage
\chapter{Introduction and Mathematical Tools}\label{chapter_introduction}
\begingroup
\hypersetup{
linkcolor=structurecolor,
linktoc=page,  
}
\minitoc \newpage
\endgroup
\section*{Introduction and Background}
\addcontentsline{toc}{section}{Introduction and Background}

In today's era, where data-driven decisions and computational efficiency are paramount, optimization techniques have become an indispensable cornerstone across various fields including mathematics, computer science, operations research, and machine learning. These methods play a pivotal role in solving complex problems, whether it involves enhancing machine learning models, improving resource allocation, or designing effective algorithms. Aimed at providing a comprehensive yet accessible reference, this book is designed to equip readers with the essential knowledge for understanding and applying optimization methods within their respective domains.

Our primary goal is to unveil the inner workings of optimization algorithms, ranging from black-box optimizers to stochastic ones, through both formal and intuitive explanations. Starting from fundamental mathematical principles, we derive key results to ensure that readers not only learn how these techniques work but also understand when and where to apply them effectively. We strike a careful balance between theoretical depth and practical application to cater to a broad audience, seeking robust optimization strategies.

The book begins by laying down the necessary mathematical foundations, such as linear algebra, inequalities, norms, and matrix decompositions, which form the backbone for discussing optimality conditions, convexity, and various optimization techniques. As the discussion progresses, we delve deeper into first-order and second-order methods, constrained optimization strategies, and stochastic optimization approaches. Each topic is treated with clarity and rigor, supplemented by detailed explanations, examples, and formal proofs to enhance comprehension.

Special emphasis is placed on gradient descent---a quintessential optimization algorithm widely used in machine learning. The text meticulously examines basic gradient descent, momentum-based methods, conjugate gradient techniques, and stochastic optimization strategies, highlighting their historical development and modern applications. Particular attention is given to stochastic gradient descent (SGD) and its crucial role in optimizing deep neural network or transformer structures, especially concerning overcoming computational challenges and escaping saddle points.

At the heart of the book also lies the exploration of constrained optimization, covering methods for handling equality and inequality constraints, augmented Lagrangian methods, alternating direction method of multipliers (ADMM), and trust region methods. Additionally, second-order methods such as Newton's method, quasi-Newton methods, and conjugate gradient methods are discussed, along with insights into least squares problems (both linear and nonlinear) and sparse optimization problems. The final section delves into stochastic optimization, encompassing topics like SGD, variance reduction techniques, and learning rate annealing and warm-up strategies.

Throughout the book, compactness, clarity, and mathematical rigor are prioritized, with detailed explanations, examples, and proofs provided to ensure a profound understanding of optimization principles and their applications. 
The primary aim of this book is to provide a self-contained introduction to the concepts
and mathematical tools used in optimization analysis and complexity, and based on them
to present major optimization algorithms and their applications to the reader.
We clearly realize our inability to cover all the useful and interesting topics concerning
optimization methods. Therefore, we recommend additional literature for more detailed introductions to related areas. Some excellent references include \citet{beck2017first, boyd2004convex, sun2006optimization, bertsekas2015convex}.

In the following sections of this chapter, we provide a concise overview of key concepts in linear algebra and calculus. Further important ideas will be introduced as needed to enhance understanding. It is important to note that this chapter does not aim for an exhaustive exploration of these subjects. Those interested in deeper study should refer to advanced texts on linear algebra and calculus.

\section{Linear Algebra}
\index{Vector space}
\index{Matrix space}
The vector space $\real^n$ consists of all $n$-dimensional column vectors with real components. Our main focus in this book will be problems within the $\real^n$ vector space, though we will occasionally examine other vector spaces, such as the nonnegative vector space. Similarly, the matrix space $\real^{m\times n}$ comprises all real-valued matrices of dimensions $m\times n$.

Scalars will be represented using non-bold font, possibly with subscripts (e.g., $a$, $\alpha$, $\alpha_i$). Vectors will be denoted by \textbf{bold} lowercase letters, potentially with subscripts (e.g., $\bmu$, $\bx$, $\bx_n$, $\bz$), while matrices will be indicated by \textbf{bold} uppercase letters, also potentially with subscripts (e.g., $\bA$, $\bL_j$). 
The $i$-th element of a vector $\bz$ will be written as $z_i$ in  non-bold font.
While a matrix $\bA\in\real^{m\times n}$, with $(i,j)$-th element being $a_{ij}$, can be denoted as $\bA=\{a_{ij}\}$ or $\bA=\{a_{ij}\}_{i,j=1}^{m,n}$.

When we fix a subset of indices, we form subarrays. Specifically, the entry at the $i$-th row and $j$-th column value of matrix $\bA$ (entry ($i,j$) of $\bA$) is denoted by  $a_{ij}$. 
We also adopt \textbf{Matlab-style notation}: the submatrix of   $\bA$ spanning from the $i$-th row to the $j$-th row and from the $k$-th column to the $m$-th column is written as $\bA_{i:j,k:m}=\bA[i:j,k:m]$. 
A colon indicates all elements along a dimension; for example,  $\bA_{:,k:m}=\bA[:,k:m]$ represents columns column $k$ through $m$  of  $\bA$, and $\bA_{:,k}$ denotes the $k$-th column of $\bA$. Alternatively, the $k$-th column of  $\bA$ can be more compactly written as $\ba_k$.

\index{Matlab-style notation}

When the indices are non-continuous, given ordered subindex sets $\sI$ and $\sJ$, $\bA[\sI, \sJ]$ denotes the submatrix of $\bA$ formed by selecting the rows and columns of $\bA$ indexed by $\sI$ and $\sJ$, respectively.
Similarly, $\bA[:, \sJ]$ indicates the submatrix of $\bA$ obtained by extracting all rows from $\bA$ and only the columns specified by $\sJ$, where again the colon operator signifies that all indices along that dimension are included. 

\begin{definition}[Matlab Notation]\label{definition:matlabnotation}
Let $\bA\in \real^{m\times n}$, and $\sI=\{i_1, i_2, \ldots, i_k\}$ and $\sJ=\{j_1, j_2, \ldots, j_l\}$ be two index sets. Then, $\bA[\sI,\sJ]$ represents   the $k\times l$ submatrix
$$
\bA[\sI,\sJ]=
\begin{bmatrix}
a_{i_1,j_1} & a_{i_1,j_2} &\ldots & a_{i_1,j_l}\\
a_{i_2,j_1} & a_{i_2,j_2} &\ldots & a_{i_2,j_l}\\
\vdots & \vdots&\ddots & \vdots\\
a_{i_k,j_1} & a_{i_k,j_2} &\ldots & a_{i_k,j_l}\\
\end{bmatrix}.
$$
Similarly, $\bA[\sI,:]$ denotes a $k\times n$ submatrix, and $\bA[:,\sJ]$ denotes a $m\times l$ submatrix.
For a vector $\bv\in\real^m$, $\bv_{\sI}\triangleq\bv[\sI]$ denotes a $\real^k$ subvector, and $\bv_{\comple{\sI}}\triangleq\bv[\comple{\sI}]$ denotes a $\real^{m-k}$ subvector, where $\comple{\sI}$ is the complement of  set $\sI$: $\comple{\sI} \triangleq\{1,2,\ldots,m\}\backslash\sI$.

Note that it does not matter whether the index sets $\sI$ and $\sJ$ are row vectors or column vectors; what's crucial is which axis they index  (either rows  or columns of $\bA$).
Additionally, the ranges of the indices are given as follows:
$$
\left\{
\begin{aligned}
0&\leq \min(\sI) \leq \max(\sI)\leq m;\\
0&\leq \min(\sJ) \leq \max(\sJ)\leq n.
\end{aligned}
\right.
$$
\end{definition}

In all cases, vectors are presented in column form rather than as rows. A row vector will be denoted by the transpose of a column vector, such as $\ba^\top$. 
A column vector with specific values is separated by the semicolon symbol $``;"$, for instance, $\bx=[1;2;3]$ is a column vector in $\real^3$. 
Similarly, a  row vector with specific values is split by the comma symbol $``,"$, e.g., $\by=[1,2,3]$ is a row vector with 3 values. 
Additionally, a column vector can also be represented as the transpose of a row vector, e.g., $\by=[1,2,3]^\top$ is also a column vector.

The transpose of a matrix $\bA$ will be denoted as $\bA^\top$, and its inverse as $\bA^{-1}$. 
We  denote the $p \times p$ identity matrix by $\bI_p$ (or simply by $\bI$ when the size is clear from context). A vector or matrix of all zeros will be denoted by a {boldface} zero, $\bzero$, with the size inferred from context, or we denote $\bzero_p$ to be the vector of all zeros with $p$ entries.
Similarly, a vector or matrix of all ones will be denoted by a {boldface} one, $\bone$, whose size is inferred  from  context, or we denote $\bone_p$ to be the vector of all ones with $p$ entries.
We  frequently omit the subscripts of these matrices when the dimensions are clear from  context.

We   use $\be_1, \be_2, \ldots, \be_n$ to represent the standard (unit) basis of $\real^n$, where $\be_i$ is the vector whose $i$-th component is one while all  others are zero.

\begin{definition}[Nonnegative Orthant, Positive Orthant, and Unit-Simplex]\label{definition:simplex}
The \textit{nonnegative orthant} is a subset of $\real^n$ that consists of all vectors in $\real^n$ with nonnegative components and is denoted by $\real_+^n$: 
$$
\real_+^n=\{[x_1, x_2, \ldots, x_n]^\top\in\real^n: x_1, x_2, \ldots, x_n\geq 0\}.
$$
Similarly, the \textit{positive orthant} comprises all vectors in $\real^n$ with strictly positive components and is denoted by $\real_{++}^n$:
$$
\real_{++}^n=\{[x_1, x_2, \ldots, x_n]^\top\in\real^n: x_1, x_2, \ldots, x_n> 0\}.
$$
The \textit{unit-simplex} (or simply simplex) is a subset of $\real^n$ comprising all nonnegative vectors whose components sum to one:
$$
\Delta_n=\{\bx=[x_1, x_2, \ldots, x_n]^\top\in\real^n:  x_1, x_2, \ldots, x_n\geq  0,
\,\, \sum_{i=1}^{n}x_i=1 \}.
$$
\end{definition}

\index{Eigenvalue}
\index{Eigenvector}
\begin{definition}[Eigenvalue, Eigenvector]
Given any vector space $\sF$ and a linear map $\bA: \sF \rightarrow \sF$ (or simply a real matrix $\bA\in\real^{n\times n}$), a scalar $\lambda \in \sK$ is called an \textit{eigenvalue, or proper value, or characteristic value} of $\bA$, if there exists some nonzero vector $\bu \in \sF$ such that
\begin{equation*}
\bA \bu = \lambda \bu.
\end{equation*}
The vector $\bu$ is then called an \textit{eigenvector} of $\bA$ associated with $\lambda$.
\end{definition}
The pair $(\lambda, \bu)$ is termed an \textit{eigenpair}. Intuitively, these definitions imply that multiplying matrix $\bA$ by the vector $\bu$ results in a new vector that is in the same direction as $\bu$, but its length is scaled by $\lambda$ (an eigenvector $\bu$ of a matrix $\bA$ represents a direction that remains unchanged when transformed into the coordinate system defined by the columns of $\bA$.). Any eigenvector $\bu$ can be scaled by a scalar $s$ such that $s\bu$ remains an eigenvector of $\bA$. To avoid ambiguity, it is common practice to normalize eigenvectors to have unit length and ensure the first entry is positive (or negative), since both $\bu$ and $-\bu$ are valid eigenvectors.
Note that real-valued matrices can have complex eigenvalues, but all eigenvalues of symmetric matrices are real (see Theorem~\ref{theorem:spectral_theorem}).

In linear algebra, every vector space has a basis, and any vector within the space can be expressed as a linear combination of the basis vectors. We define the span and dimension of a subspace using the basis.

\index{Subspace}
\begin{definition}[Subspace]
A nonempty subset $\mathcalV$ of $\real^n$ is called a \textit{subspace} if $x\ba+y\ba\in \mathcalV$ for every $\ba,\bb\in \mathcalV$ and every $x,y\in \real$.
\end{definition}

\index{Span}
\begin{definition}[Span]
If every vector $\bv$ in subspace $\mathcalV$ can be expressed as a linear combination of $\{\ba_1, \ba_2, \ldots,$ $\ba_m\}$, then $\{\ba_1, \ba_2, \ldots, \ba_m\}$ is said to \textit{span} $\mathcalV$.
\end{definition}

\index{Linearly independent}
In this context, we frequently use the concept of linear independence for a set of vectors. Two equivalent definitions are given below.
\begin{definition}[Linearly Independent]
A set of vectors $\{\ba_1, \ba_2, \ldots, \ba_m\}$ is said to be \textit{linearly independent} if there is no combination that  can yield $x_1\ba_1+x_2\ba_2+\ldots+x_m\ba_m=0$ unless all $x_i$'s are equal to  zero. An equivalent definition states that $\ba_1\neq \bzero$, and for every $k>1$, the vector $\ba_k$ does not belong to the span of $\{\ba_1, \ba_2, \ldots, \ba_{k-1}\}$.
\end{definition}

\index{Basis}
\index{Dimension}
\begin{definition}[Basis and Dimension]
A set of vectors $\{\ba_1, \ba_2, \ldots, \ba_m\}$ forms a \textit{basis} of a subspace $\mathcalV$ if they are linearly independent and span $\mathcalV$. Every basis of a given subspace contains the same number of vectors, and this number is called the \textit{dimension} of the subspace $\mathcalV$. By convention, the subspace $\{\bzero\}$ is defined to have a dimension of zero. Additionally, every subspace of nonzero dimension has an orthogonal basis, meaning that a basis for the subspace can always be chosen to be orthogonal.
\end{definition}

\index{Column space}
\begin{definition}[Column Space (Range)]
If $\bA$ is an $m \times n$ real matrix, its \textit{column space (or range)} is the set of all linear combinations of its columns, formally defined as:
\begin{equation*}
\mathcalC (\bA) = \{ \by\in \real^m: \exists\, \bx \in \real^n, \, \by = \bA \bx \}.
\end{equation*}
Similarly, the \textit{row space} of $\bA$ is the set spanned by its rows, which is equivalent to the column space of $\bA^\top$:
\begin{equation*}
\mathcalC (\bA^\top) = \{ \bx\in \real^n: \exists\, \by \in \real^m, \, \bx = \bA^\top \by \}.
\end{equation*}
\end{definition}

\index{Null space}
\begin{definition}[Null Space (Nullspace, Kernel)]\label{definition:nullspace}
If $\bA$ is an $m \times n$ real matrix, its \textit{null space} (also called the \textit{kernel} or \textit{nullspace}) is the set of all vectors that are mapped to zero by $\bA$:
\begin{equation*}
\nspace (\bA) = \{\by \in \real^n:  \, \bA \by = \bzero \}.
\end{equation*}
Similarly, the null space of $\bA^\top$ is defined as 	
\begin{equation*}
\nspace (\bA^\top) = \{\bx \in \real^m:  \, \bA^\top \bx = \bzero \}.
\end{equation*}
\end{definition}

Both the column space of $\bA$ and the null space of $\bA^\top$ are subspaces of $\real^m$. In fact, every vector in $\nspace(\bA^\top)$ is orthogonal  to every vector in $\cspace(\bA)$, and vice versa. Similarly, every vector in $\nspace(\bA)$ is also perpendicular to every vector in $\cspace(\bA^\top)$, and vice versa.

\index{Rank}
\begin{definition}[Rank]
The \textit{rank} of a matrix $\bA\in \real^{m\times n}$ is the dimension of its column space. 
In other words, the rank of $\bA$ is the maximum number of linearly independent columns of $\bA$, which is also equal to the maximum number of linearly independent rows of $\bA$. 
A fundamental property of matrices is that $\bA$ and its transpose $\bA^\top$ always have the same rank. 
A matrix $\bA$ is considered to have \textit{full rank} if its rank is equal to $min\{m,n\}$.

A specific example of rank computation arises when considering the outer product of two vectors. Given a vector $\bu \in \real^m$ and a vector $\bv \in \real^n$, then the $m\times n$ matrix $\bu\bv^\top$ formed by their outer product always has rank 1. In short, the rank of a matrix is equal to:
\begin{itemize}
\item The number of linearly independent columns.
\item The number of linearly independent rows.
\end{itemize}
These two quantities are always equal (see Theorem~\ref{theorem:fundamental-linear-algebra}).
\end{definition}

\index{Orthogonal complement}
\begin{definition}[Orthogonal Complement in General]
The \textit{orthogonal complement} of a subspace $\mathcalV\subseteq\real^m$, denoted as $\mathcalV^\perp\subseteq \real^m$, consists of all vectors in $\real^m$ that are perpendicular to every vector in $\mathcalV$. Formally,
$$
\mathcalV^\perp = \{\bv\in \real^m : \bv^\top\bu=0, \,\,\, \forall \bu\in \mathcalV  \}.
$$
These two subspaces are mutually exclusive yet collectively span the entire space.
The dimensions of $\mathcalV$ and $\mathcalV^\perp$ sum up to the dimension of the entire space. Furthermore, it holds that $(\mathcalV^\perp)^\perp=\mathcalV$.
\end{definition}

\index{Orthogonal complement}
\begin{definition}[Orthogonal Complement of Column Space]
For an  $m \times n$ real matrix $\bA$, the orthogonal complement of its column space $\mathcalC(\bA)$, denoted as $\mathcalC^{\bot}(\bA)$, is the subspace defined as:
\begin{equation*}
\begin{aligned}
\mathcalC^{\bot}(\bA) &= \{\by\in \real^m: \, \by^\top \bA \bx=\bzero, \, \forall \bx \in \real^n \} \\
&=\{\by\in \real^m: \, \by^\top \bv = \bzero, \, \forall \bv \in \mathcalC(\bA) \}.
\end{aligned}
\end{equation*}
\end{definition}
We can then identify the four fundamental subspaces associated with any matrix $\bA\in \real^{m\times n}$ of rank $r$:
\begin{itemize}
\item $\cspace(\bA)$: Column space of $\bA$, i.e., linear combinations of columns, with dimension $r$.

\item $\nspace(\bA)$: Null space of $\bA$, i.e., all $\bx$ satisfying $\bA\bx=\bzero$, with dimension $n-r$.

\item  $\cspace(\bA^\top)$: Row space of $\bA$, i.e., linear combinations of rows, with dimension $r$.

\item  $\nspace(\bA^\top)$: Left null space of $\bA$, i.e., all $\by$ satisfying $\bA^\top \by=\bzero$, with dimension $m-r$.
\end{itemize}
Furthermore, $\nspace(\bA)$ is the orthogonal complement of $\cspace(\bA^\top)$, and $\cspace(\bA)$ is the orthogonal complement of $\nspace(\bA^\top)$; see Theorem~\ref{theorem:fundamental-linear-algebra}.

\index{Fundamental subspace}

\index{Orthogonal matrix}
\begin{definition}[Orthogonal Matrix, Semi-Orthogonal Matrix]
A real square matrix $\bQ\in\real^{n\times n}$ is called an \textit{orthogonal matrix} if its inverse is equal to its  transpose, that is, $\bQ^{-1}=\bQ^\top$ and $\bQ\bQ^\top = \bQ^\top\bQ = \bI$. 
In other words, suppose $\bQ=[\bq_1, \bq_2, \ldots, \bq_n]$, where $\bq_i \in \real^n$ for all $i \in \{1, 2, \ldots, n\}$, then $\bq_i^\top \bq_j = \delta(i,j)$, where $\delta(i,j)$ is the Kronecker delta function. 
One important property of an orthogonal matrix is that it preserves vector norms: $\norm{\bQ\bx}= \norm{\bx}$ for all $\bx\in\real^n$.

If $\bQ$ contains only $\gamma$ of these columns with $\gamma<n$, then $\bQ^\top\bQ = \bI_\gamma$ stills holds, where $\bI_\gamma$ is the $\gamma\times \gamma$ identity matrix. 
But $\bQ\bQ^\top=\bI$ will not hold. 
In this case, $\bQ$ is called \textit{semi-orthogonal}.
\end{definition}

From an introductory course on linear algebra, we have the following remark regarding  the equivalent claims of nonsingular matrices.
\begin{remark}[List of Equivalence of Nonsingularity for a Matrix]
For a square matrix $\bA\in \real^{n\times n}$, the following claims are equivalent:
\begin{itemize}
\item $\bA$ is nonsingular;~\footnote{The source of the name is a result of the singular value decomposition (SVD); see, for example, \citet{lu2021numerical}.}
\item $\bA$ is invertible, i.e., $\bA^{-1}$ exists; 
\item $\bA\bx=\bb$ has a unique solution $\bx = \bA^{-1}\bb$;
\item $\bA\bx = \bzero$ has a unique, trivial solution: $\bx=\bzero$;
\item Columns of $\bA$ are linearly independent;
\item Rows of $\bA$ are linearly independent;
\item Determinant $\det(\bA) \neq 0$; 
\item Dimension $\dim(\nspace(\bA))=0$;
\item $\nspace(\bA) = \{\bzero\}$, i.e., the null space is trivial;
\item $\cspace(\bA)=\cspace(\bA^\top) = \real^n$, i.e., the column space or row space spans the entire $\real^n$;
\item $\bA$ has full rank $r=n$;
\item The reduced row echelon form is $\bR=\bI$;
\item $\bA^\top\bA$ is symmetric positive definite;
\item $\bA$ has $n$ nonzero (positive) singular values;
\item All eigenvalues are nonzero.
\end{itemize}
\end{remark}

These equivalences are fundamental in linear algebra and will be useful in various contexts. Conversely, the following remark establishes the equivalent conditions for a singular matrix.
\begin{remark}[List of Equivalence of Singularity for a Matrix]
For a square matrix $\bA\in \real^{n\times n}$ with eigenpair $(\lambda, \bu)$, the following claims are equivalent:
\begin{itemize}
\item $(\bA-\lambda\bI)$ is singular; 
\item $(\bA-\lambda\bI)$ is not invertible;
\item $(\bA-\lambda\bI)\bx = \bzero$ has  nonzero  solutions $\bx\neq \bzero$, and $\bx=\bu$ is one of such solutions;
\item $(\bA-\lambda\bI)$ has linearly dependent columns;
\item Determinant $\det(\bA-\lambda\bI) = 0$; 
\item Dimension $\dim(\nspace(\bA-\lambda\bI))>0$;
\item Null space of $(\bA-\lambda\bI)$ is nontrivial;
\item Columns of $(\bA-\lambda\bI)$ are linearly dependent;
\item Rows of $(\bA-\lambda\bI)$ are linearly dependent;
\item $(\bA-\lambda\bI)$ has rank $r<n$;
\item Dimension of column space = dimension of row space = $r<n$;
\item $(\bA-\lambda\bI)^\top(\bA-\lambda\bI)$ is symmetric semidefinite;
\item $(\bA-\lambda\bI)$ has $r<n$ nonzero (positive) singular values;
\item Zero is an eigenvalue of $(\bA-\lambda\bI)$.
\end{itemize}
\end{remark}

\section{Well-Known Inequalities}\label{section:inequalities}
In this section, we introduce several well-known inequalities that will be frequently used throughout the discussion.

The AM-GM inequality is a fundamental tool in competitive mathematics, particularly useful for determining the maximum or minimum values of multivariable functions or expressions. It establishes a relationship between  the \textit{arithmetic mean (AM)} and the \textit{geometric mean (GM)}. \index{AM-GM inequality}
\begin{theorem}[AM-GM Inequality]
For  any nonnegative real numbers $x_1,x_2,\ldots,x_n$,  the following inequality holds:
$$
\frac{\sum_{i=1}^{n} x_i}{n} \geq \sqrt[n]{\prod_{i=1}^{n}x_i}.
$$ 
That is, the \textit{geometric mean} of a set of nonnegative numbers does not exceed their \textit{arithmetic mean}.
Equality holds if and only if all the numbers are equal.
\end{theorem}

When $n=2$, the AM-GM inequality can be expressed as:
\begin{equation}\label{equation:amgm_ineq}
a^2+b^2\geq 2\sqrt{a^2b^2} = 2\abs{ab}\geq ab.
\end{equation}

\begin{proposition}[Weighted AM-GM Inequality]\label{proposition:weighted_amgm}
For any nonnegative real numbers $x_1,x_2,\ldots,x_n$ and nonnegative weights $w_1,w_2,\ldots, w_n$,  the following inequality holds:
$$
\frac{\sum_{i=1}^{n} w_i x_i}{\sum_{i=1}^{n}w_i} \geq
\sqrt[(\sum w_i)]{\prod_{i=1}^{n}x_i^{w_i}}.
$$ 
When $w_1=w_2=\ldots=w_n=1$, this reduces to  the standard AM-GM inequality. Alternatively, when $\bw\in \Delta_n$, i.e., $\bw$ is a unit-simplex with $\sum_{i=1}^{n}w_i=1$, it follows that 
$$
\sum_{i=1}^{n} w_i x_i \geq
\prod_{i=1}^{n}x_i^{w_i},
\quad\text{where}\quad \bw\in \Delta_n.
$$
\end{proposition}
\begin{proof}[of Proposition~\ref{proposition:weighted_amgm}]
For simplicity, we will only prove the second part, as the first part can be established similarly.
Applying  Jensen's inequality (Theorem~\ref{theorem:jensens_ineq}) to the convex function $f(x)=-\ln (x)$ with $x_1, x_2,  \ldots, x_n>0$ and $\bw \in \Delta_n$, we have 
$
-\ln\left( \sum_{i=1}^{n} w_i x_i\right) \leq -\sum_{i=1}^{n} w_i \ln(\bx_i).
$
Taking the exponent of both sides, we obtain 
$
\sum_{i=1}^{n}w_i x_i \geq \exp^{\sum_{i=1}^{n} w_i \ln(x_i)}
=\prod_{i=1}^{n}x_i^{w_i}.
$
This completes the proof.
\end{proof}

\subsection{Cauchy-Schwarz Inequality}
\index{Cauchy–Schwarz inequality}
The \textit{Cauchy–Schwarz inequality} is one of the most important and widely used inequalities in mathematics.
\begin{proposition}[Cauchy-Schwarz Matrix (Vector) Inequality]\label{proposition:cauchy-schwarz-inequ}
For any $m\times n$ matrices $\bX$ and $\bY$, we have 
$$
\normbig{\bX^\top \bY} \leq \Vert\bX\Vert \cdot \Vert\bY\Vert.
$$
This is a special form of the Cauchy-Schwarz inequality, where the inner product is defined as $\langle\bX,\bY\rangle = \norm{\bX^\top \bY}$.
Similarly, for any vectors $\bu$ and $\bv$, we have 
\begin{equation}\label{equation:vector_form_cauchyschwarz}
\absbig{\bu^\top \bv} \leq \normtwo{\bu} \cdot  \normtwo{\bv},
\end{equation}
where the equality holds if and only if $\bu$ and $\bv$ are linearly dependent.
In the two-dimensional case, it becomes
$$
(ac+bd)^2 \leq (a^2 +b^2)(c^2+d^2).
$$
\end{proposition}
The vector form of the Cauchy-Schwarz inequality plays an important role in various branches of modern mathematics, including Hilbert space theory and numerical analysis \citep{wu2009various}. For simplicity, we provide only the proof for the vector form of the Cauchy-Schwarz inequality. 
To see this, given two vectors $\bu,\bv\in \real^n$, we have 
$$
\begin{aligned}
&0\leq \sum_{i=1}^{n}\sum_{j=1}^{n} (u_i v_j - u_j v_i)^2 = 
\sum_{i=1}^{n}\sum_{j=1}^{n} u_i^2v_j^2 + \sum_{i=1}^{n}\sum_{j=1}^{n} v_i^2 u_j^2 - 2\sum_{i=1}^{n}\sum_{j=1}^{n} u_iu_j v_iv_j\\
&=\left(\sum_{i=1}^{n} u_i^2\right) \left(\sum_{j=1}^{n} v_j^2\right) +
\left(\sum_{i=1}^{n} v_i^2\right) \left(\sum_{j=1}^{n} u_j^2\right) - 
2\left(\sum_{i=1}^{n} u_iv_i \right)^2
=2 \normtwo{\bu}^2 \cdot \normtwo{\bv}^2 -2 \abs{\bu^\top\bv}^2,
\end{aligned}
$$
from which the result follows.
The equality holds if and only if $\bu = k\bv$ for some constant $k\in \real$, i.e., $\bu$ and $\bv$ are linearly dependent.

\paragraph{Angle between two vectors.} From Equation~\eqref{equation:vector_form_cauchyschwarz}, given two vectors $\bx,\by$, we note that 
$$
-1\leq 
\frac{\bx^\top\by}{\Vert\bx\Vert_2 \Vert\by\Vert_2} 
\leq 1.
$$
This two-side inequality illustrates the concept of the angle between two vectors.
\begin{definition}[Angle Between Vectors]\label{definition:angle_bet_vec_ineq}
The angle between two vectors $\bx$ and $\by$ is the number $\theta\in [0,\pi]$ such that 
$$
\cos (\theta) = \frac{\bx^\top\by}{\Vert\bx\Vert_2 \Vert\by\Vert_2} .
$$
\end{definition}
The definition of the angle between vectors will be useful in the discussion of line search strategies (Section~\ref{section:gd_conv_line_search}).

Starting from the vector form of the Cauchy-Schwarz inequality, let $x_i\triangleq u_i^2$, $y_i\triangleq v_i^2$ for all $i=\{1,2,\ldots,n\}$, where $\bu,\bv\in\real^n$ in Equation~\eqref{equation:vector_form_cauchyschwarz}. We then obtain:
\begin{equation}\label{equation:gen_ca_sc_1}
\left( \sum_{i=1}^{n} x_i \right)^{1/2} \left( \sum_{i=1}^{n} y_i \right)^{1/2}
\geq  \sum_{i=1}^{n} (x_iy_i)^{1/2} .
\end{equation}
More generally, we have the generalized Cauchy-Schwarz inequality.
\begin{theorem}[Generalized Cauchy-Schwarz Inequality]
Given a set of vectors $\ba,\bb,\ldots, \bz\in \real^n$, and weights $\lambda_1, \lambda_2, \ldots, \lambda_z$ with $\lambda_1+\lambda_2+\ldots+\lambda_z=1$, it follows that 
$$
\small
\begin{aligned}
& \left( \sum_{i=1}^{n} a_i \right)^{\lambda_a}
\left( \sum_{i=1}^{n} b_i \right)^{\lambda_b}
\ldots
\left( \sum_{i=1}^{n} z_i \right)^{\lambda_z}
\geq 
(a_1^{\lambda_a} b_1^{\lambda_b} \ldots z_1^{\lambda_z})+
(a_2^{\lambda_a} b_2^{\lambda_b} \ldots z_2^{\lambda_z})+
\ldots +
(a_n^{\lambda_a} b_n^{\lambda_b} \ldots z_n^{\lambda_z}).
\end{aligned}
$$
The equality holds if $a_i=b_i=\ldots=z_i$ for all $i\in\{1,2,\ldots,n\}$.
\end{theorem}


\subsection{Young's Inequality}
Young's inequality is a special case of the weighted AM-GM inequality (Proposition~\ref{proposition:weighted_amgm}) and has a wide range of applications.
\begin{theorem}[Young's Inequality]\label{theorem:holder-inequality1}
For nonnegative numbers $x,y \geq 0$, and positive real numbers $p,q>1$ with $\frac{1}{p}+\frac{1}{q} = 1$, it follows that
\begin{equation}\label{equation:holder-v1}
xy \leq  \frac{1}{p} x^p +\frac{1}{q} y^q.
\end{equation}
\end{theorem}

\begin{figure}[h]
\centering
\vspace{-0.35cm}
\subfigtopskip=2pt
\subfigbottomskip=2pt
\subfigcapskip=-5pt
\subfigure[Case 1: $x< y^{1/(p-1)}$ and $p\geq 2$.]{\label{fig:holder1}
\includegraphics[width=0.38\linewidth]{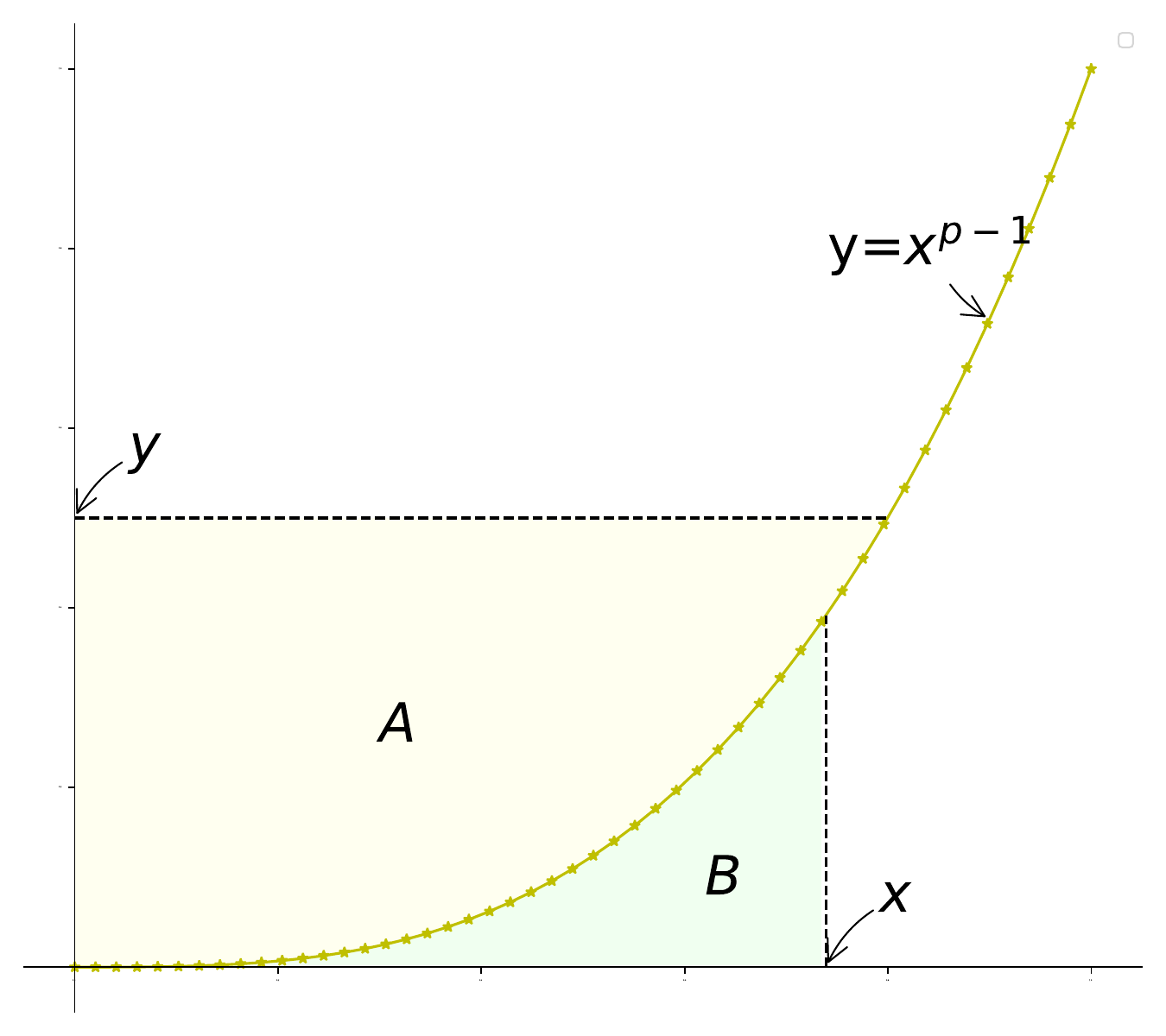}}
\quad 
\subfigure[Case 2: $x\geq  y^{1/(p-1)}$ and $p\geq 2$.]{\label{fig:holder2}
\includegraphics[width=0.38\linewidth]{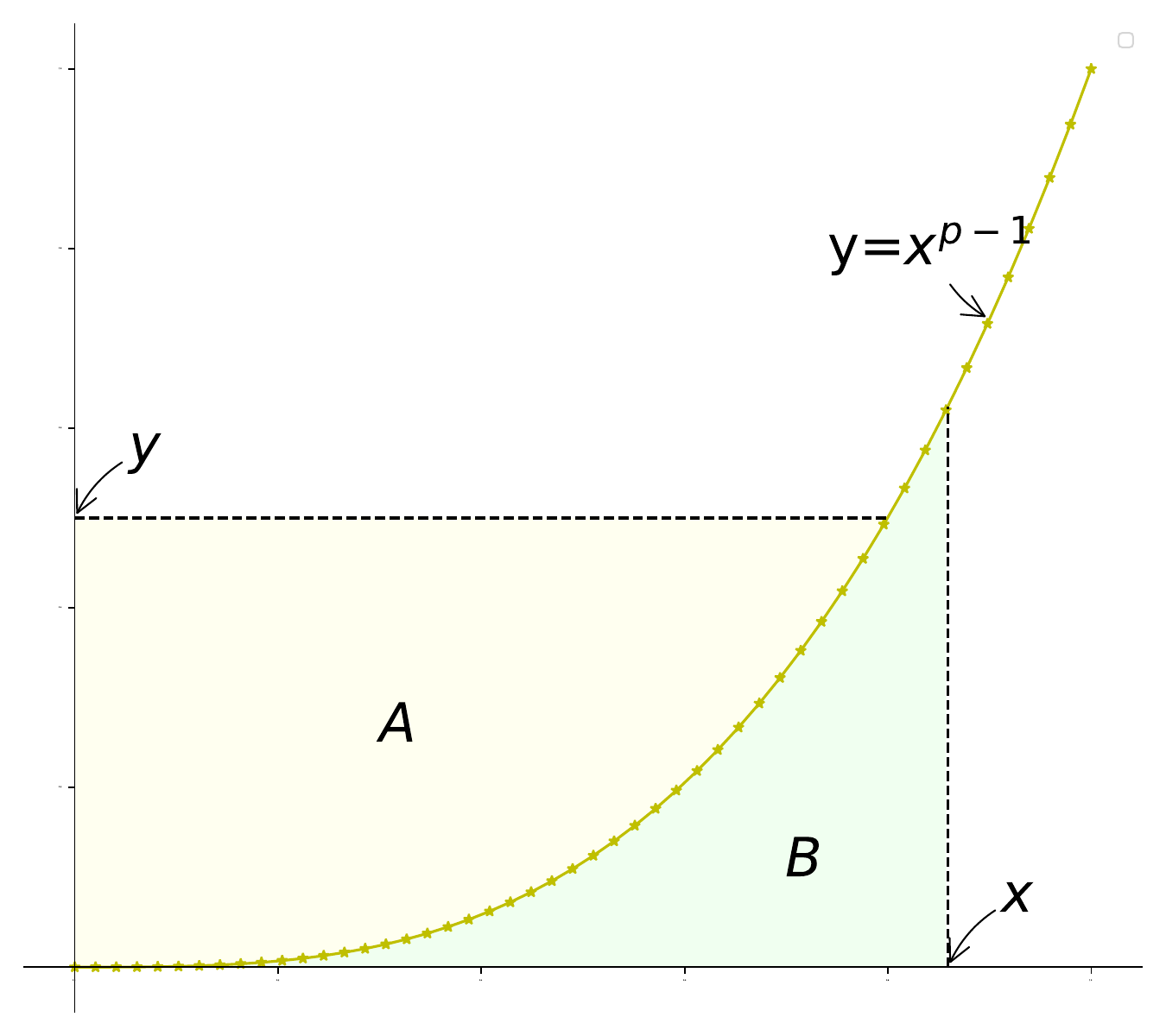}}
\subfigure[Case 3: $x< y^{1/(p-1)}$ and $1<p< 2$.]{\label{fig:holder3}
\includegraphics[width=0.38\linewidth]{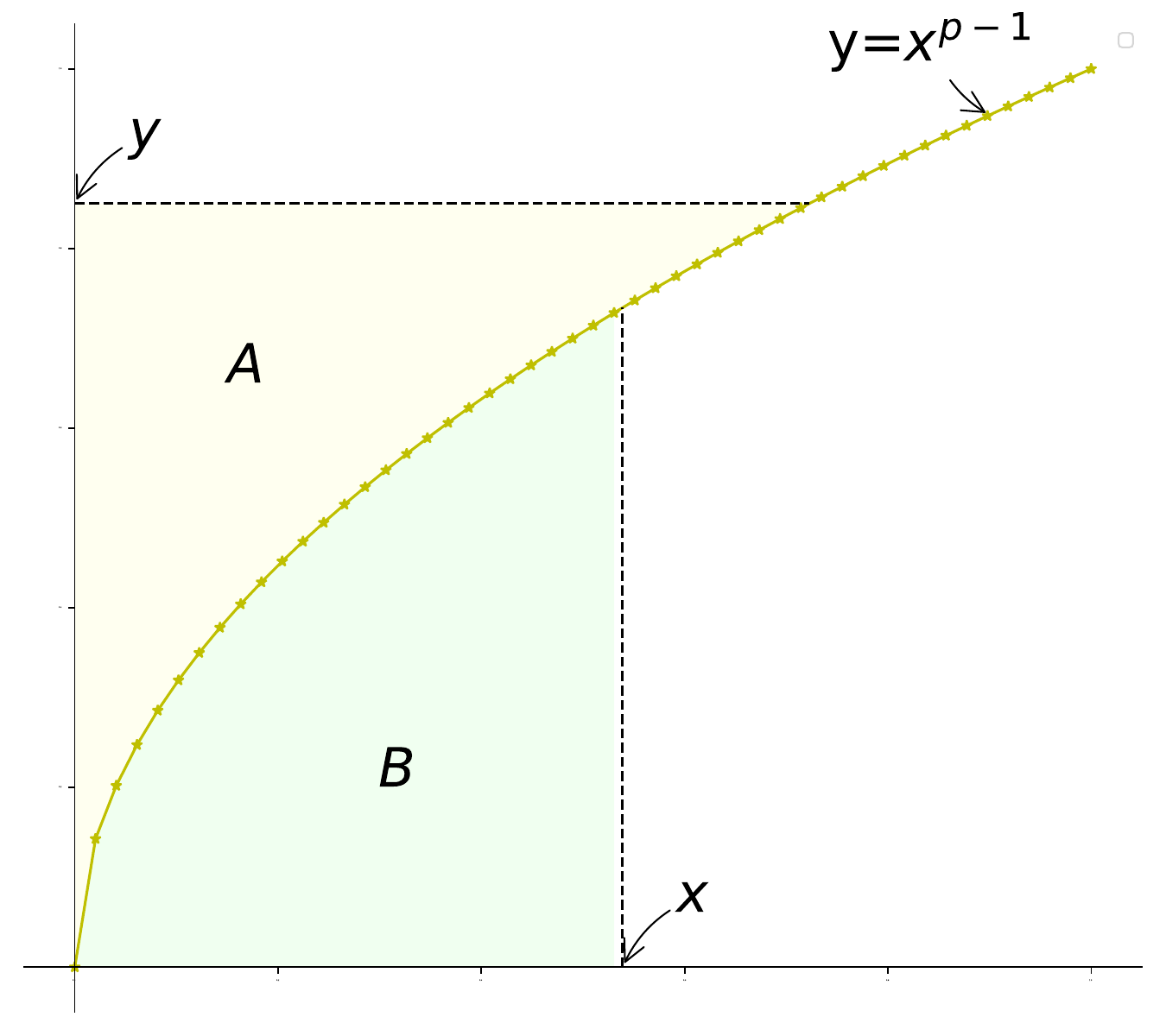}}
\quad 
\subfigure[Case 4: $x\geq  y^{1/(p-1)}$ and $1<p< 2$.]{\label{fig:holder4}
\includegraphics[width=0.38\linewidth]{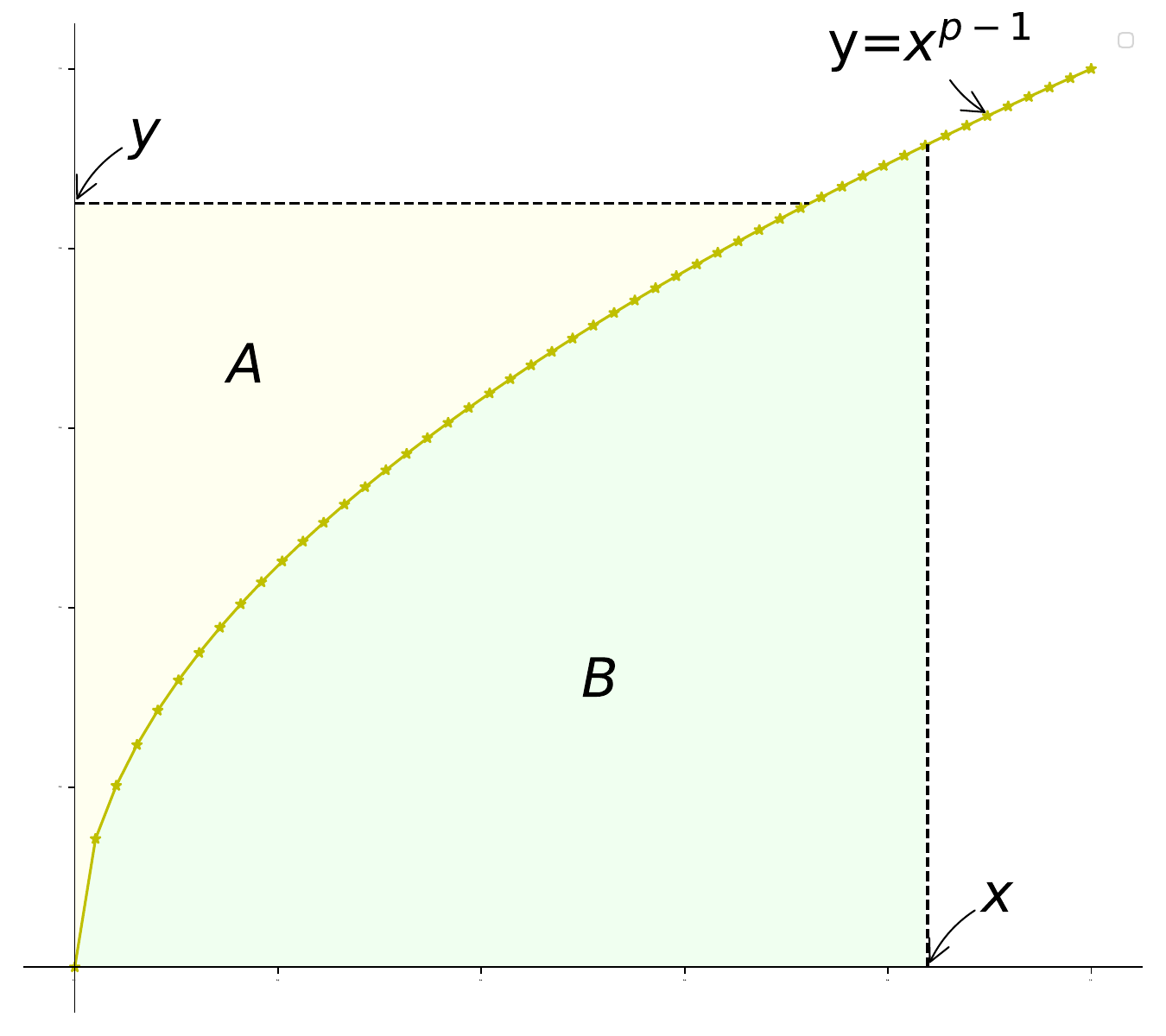}}
\caption{Demonstration of Young's inequality for different cases.}
\label{fig:holderineqd}
\end{figure}
\begin{proof}[of Theorem~\ref{theorem:holder-inequality1}]
The area $xy$ is bounded above by the sum of the areas of the two trapezoids with curved edges (the shaded regions $A$ and $B$) as shown in Figure~\ref{fig:holderineqd}:
$$
\mathrm{area } \,\,\,B =\int_0^x x^{p-1} dx, \qquad \mathrm{area } \,\,\,A=\int_{0}^y y^{1/(p-1)}dy.
$$
That is, 
$$
\small
\begin{aligned}
xy &\leq \int_0^x x^{p-1} dx + \int_{0}^y y^{1/(p-1)}dy 
= \frac{1}{p} x^p + \int_{0}^y y^{q/p}dy 
= \frac{1}{p} x^p +  (\frac{q}{p}+1)^{-1} y^{q/p+1} 
= \frac{1}{p} x^p +\frac{1}{q} y^q.
\end{aligned}
$$
This completes the proof.
\end{proof}

\begin{SCfigure}
\centering
\includegraphics[width=0.48\textwidth]{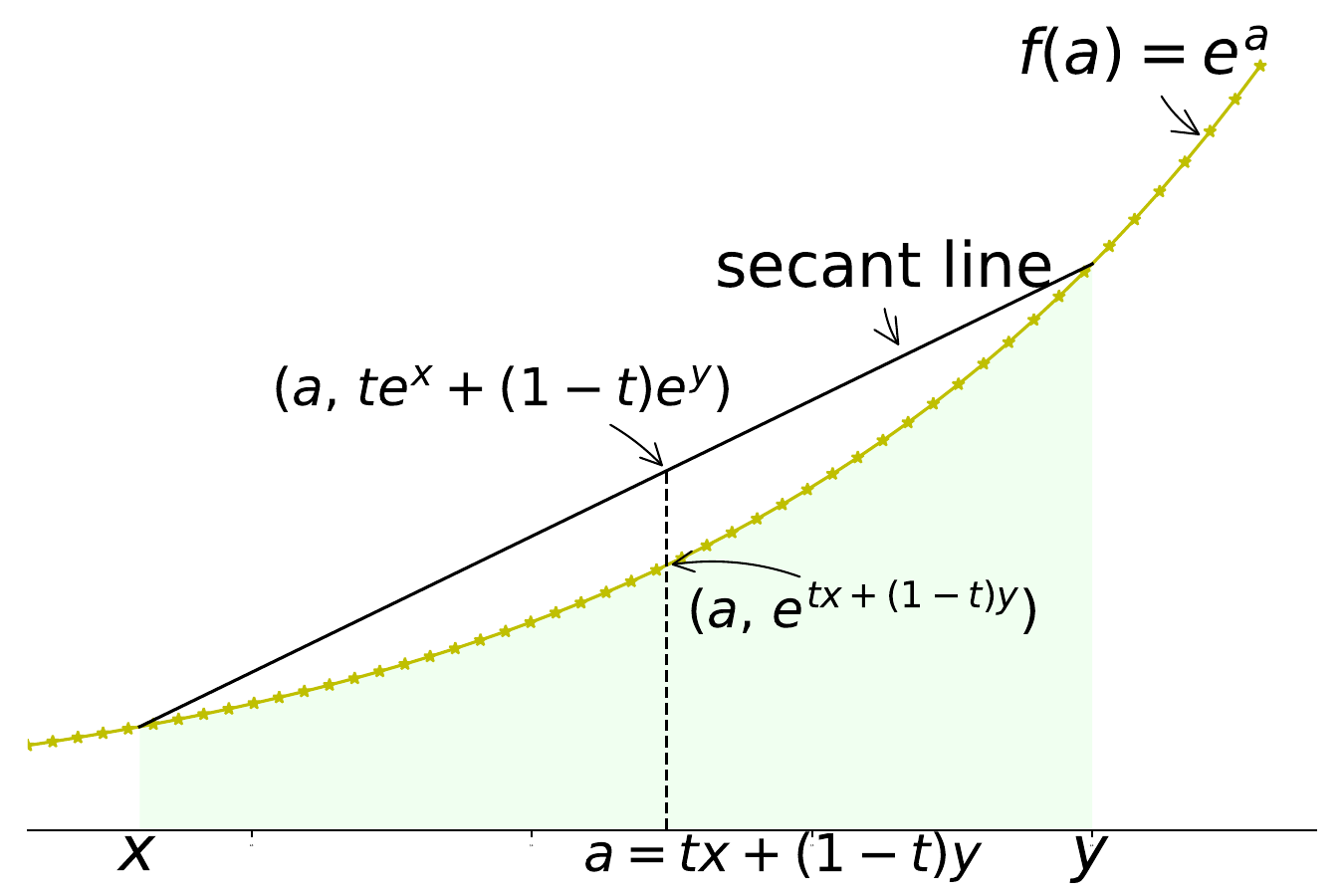}
\caption{Demonstration of interpolation inequality for $e^x$.}
\label{fig:older_interpolation}
\end{SCfigure}

There are several conceptually different ways to prove Young's inequality. As an alternative, we can use the \textit{interpolation inequality} as follows.
\begin{lemma}[Interpolation Inequality for Exponential Functions]\label{lemma:interpolation_inequality}
For $t\in [0,1]$, we have 
$$
e^{tx+(1-t)y} \leq t e^x + (1-t)e^y. 
$$
\end{lemma}
\begin{proof}[of Lemma~\ref{lemma:interpolation_inequality}]
The function of the scant line through the points $(x,e^x)$ and $(y, e^y)$ on the graph of $f(a)=e^a$ is (see Figure~\ref{fig:older_interpolation}):
$$
f(t) = (tx+(1-t)y, te^x+(1-t)e^y).
$$
Since the function $f(a)=e^a$ is convex, we have 
$$
e^{tx+(1-t)y}\leq t e^x + (1-t)e^y. 
$$
This completes the proof.
\end{proof}


Using the interpolation inequality for the exponential function, we provide an alternative proof of Theorem~\ref{theorem:holder-inequality1}.
\begin{proof}[of Theorem~\ref{theorem:holder-inequality1}, an alternative way]
We observe that
$$
\begin{aligned}
xy &= \exp\{\log x+\log y\} 
= \exp\left \{ \frac{q-1}{q} \frac{q}{q-1} \log x + \frac{1}{q} q \log y \right\}.
\end{aligned}
$$
Applying Lemma~\ref{lemma:interpolation_inequality}, we obtain
$$
\begin{aligned}
xy&\leq \frac{q-1}{q} e^{\frac{q}{q-1} \log x} + \frac{1}{q} e^{q\log y} 
= \frac{q-1}{q} x^{\frac{q}{q-1}} + \frac{1}{q} y^q
= \frac{1}{p} x^p +\frac{1}{q}y^q,
\end{aligned}
$$
from which the result follows.
\end{proof}


\subsection{H{\"o}lder's Inequality}

\holders inequality, named after \textit{Otto  H{\"o}lder}, is widely used in optimization, machine learning, and many other fields. 
It is also a generalization of the (vector) Cauchy-Schwarz inequality.

\begin{theorem}[\holders  Inequality]\label{theorem:holder-inequality}
Suppose $p,q>1$ such that $\frac{1}{p}+\frac{1}{q} = 1$. Then,  for any vector $\bx,\by\in \real^n$, we have
$$
\sum_{i=1}^{n}x_i y_i
\leq 
\abs{\sum_{i=1}^{n}x_i y_i}
\leq \sum_{i=1}^{n}\abs{x_i} \abs{y_i} \leq \left(\sum_{i=1}^{n}  \abs{x_i}^p\right)^{1/p}  \left(\sum_{i=1}^{n} \abs{y_i}^q\right)^{1/q}=\norm{\bx}_p\norm{\by}_q,
$$
where $\norm{\bx}_p = \left(\sum_{i=1}^{n}  \abs{x_i}^p\right)^{1/p} $ is known as the \textbf{$\ell_p$ norm} or \textbf{$p$ norm} of the vector $\bx$ (see Section~\ref{section:vector-norm}). The equality holds if the two sequences $\{\abs{x_i}^p\}$ and $\{\abs{y_i}^q\}$ are linearly dependent~\footnote{To be more concrete, the equality attains if and only if $\abs{\bx^\top\by}=\abs{\bx}^\top\abs{\by}$ and 
$$
\left\{
\begin{aligned}
&\abs{\bx}\hadaprod\abs{\by}=\norminf{\by}\abs{\bx}, &\text{if }& p=1; \\ 
&\abs{\bx}\hadaprod\abs{\by}=\norminf{\bx}\abs{\by}, &\text{if }& p=\infty;\\ 
&\norm{\by}_q^{1/p} \abs{\bx}^{\hadaprod1/q} = \norm{\bx}_p^{1/q}\abs{\by}^{\hadaprod1/p}, &\text{if }&1<p<\infty, \\ 
\end{aligned}
\right.
$$
where $\hadaprod$ denotes Hadamard product/power; see \citet{bernstein2008matrix} and the references therein.
}.
When $p=q=2$, this reduces to the vector Cauchy-Schwarz inequality \eqref{equation:vector_form_cauchyschwarz}.
\end{theorem}
\begin{proof}[of Theorem~\ref{theorem:holder-inequality}]
Let $ u \triangleq \frac{\abs{x_i}}{\Vert\bx\Vert_p}$ and $ v \triangleq \frac{\abs{y_i}}{\Vert\by\Vert_q}$. From Equation~\eqref{equation:holder-v1}, it follows that
$$
uv = \frac{\abs{x_i}\abs{y_i}}{\norm{\bx}_p\norm{\by}_q} \leq \frac{1}{p} \frac{\abs{x_i}^p}{\norm{\bx}_p^p} + \frac{1}{q}\frac{\abs{y_i}^q}{\norm{\by}_q^q}, \qquad \forall i \in \{1,2,\ldots,n\}.
$$
Therefore,
$$
\sum_{i=1}^{n} \frac{\abs{x_i}\abs{y_i}}{\norm{\bx}_p\norm{\by}_q} \leq \frac{1}{p\norm{\bx}_p^p} \sum_{i=1}^{n}\abs{x_i}^p  +
\frac{1}{q\norm{\by}_q^q}   \sum_{i=1}^{n}\abs{y_i}^q = \frac{1}{p}+\frac{1}{q} = 1.
$$
That is, 
$
\sum_{i=1}^{n} \abs{x_i}\abs{y_i}  \leq  \norm{\bx}_p\norm{\by}_q.
$
It is trivial that $\sum_{i=1}^{n}x_i y_i\leq \sum_{i=1}^{n}\abs{x_i} \abs{y_i} $,
from which the result follows.
\end{proof}


Minkowski's inequality follows immediately from the \holders inequality.
\begin{theorem}[Minkowski's Inequality]\label{theorem:minkoski}
Given nonnegative reals $x_1, x_2, \ldots, x_n\geq 0$ and $y_1, y_2, \ldots, y_n\geq 0$, and $p\geq 1$, we have 
$$
\left( \sum_{i=1}^{n} (x_i+y_i)^p  \right)^{1/p} \leq  
\left( \sum_{i=1}^{n}x_i^p \right)^{1/p}+
\left( \sum_{i=1}^{n}y_i^p \right)^{1/p}.
$$
The equality holds if and only if the sequence $x_1, x_2, \ldots, x_n$ and $y_1, y_2, \ldots, y_n$ are proportional. If $p<1$, the inequality sign is reversed.
\end{theorem}
Minkowski's inequality is essentially  the triangle inequality of $\ell_p$ norms (Section~\ref{section:vector-norm}).
\begin{proof}[of Theorem~\ref{theorem:minkoski}]
If $p=1$, the inequality is immediate. We then assume $p>1$.
Observe that
$$
\sum_{i=1}^{n} (x_i+y_i)^p = 
\sum_{i=1}^{n} x_i (x_i+y_i)^{p-1}
+ 
\sum_{i=1}^{n} y_i(x_i+y_i)^{p-1}.
$$
Using \holders inequality, given $p,q$ such that $p>1$ and $\frac{1}{p}+\frac{1}{q}=1$, we have 
$$
\footnotesize
\begin{aligned}
\sum_{i=1}^{n} x_i (x_i+y_i)^{p-1}
&\leq 
\bigg( \sum_{i=1}^{n} x_i^p \bigg)^{1/p}
\bigg( \sum_{i=1}^{n} (x_i+y_i)^{(p-1)q} \bigg)^{1/q}
=\bigg( \sum_{i=1}^{n} x_i^p \bigg)^{1/p}
\bigg( \sum_{i=1}^{n} (x_i+y_i)^{p} \bigg)^{1/q},\\
\sum_{i=1}^{n} y_i (x_i+y_i)^{p-1}
&\leq 
\bigg( \sum_{i=1}^{n} y_i^p \bigg)^{1/p}
\bigg( \sum_{i=1}^{n} (x_i+y_i)^{(p-1)q} \bigg)^{1/q}
=\bigg( \sum_{i=1}^{n} y_i^p \bigg)^{1/p}
\bigg( \sum_{i=1}^{n} (x_i+y_i)^{p} \bigg)^{1/q}.
\end{aligned}
$$
Combining these results, we obtain
$$
\sum_{i=1}^{n} (x_i+y_i)^p
\leq 
\bigg(
\bigg( \sum_{i=1}^{n} x_i^p \bigg)^{1/p}
+
\bigg( \sum_{i=1}^{n} y_i^p \bigg)^{1/p}
\bigg)
\bigg( \sum_{i=1}^{n} (x_i+y_i)^{p} \bigg)^{1/q}.
$$
This completes the proof.
\end{proof}

\section{Norm}\label{appendix:matrix-norm}
The concept of a norm is essential for evaluating the magnitude of vectors and, consequently, allows for the definition of certain metrics on linear spaces equipped with norms. Norms provide a measure of the magnitude of a vector or matrix, which is useful in many applications, such as determining the length of a vector in Euclidean space or the size of a matrix in a multidimensional setting.
Additionally, norms enable us to define distances between vectors or matrices. 
The distance between two vectors $\bu$ and $\bv$ can be computed using the norm of their difference $\norm{\bu-\bv}$. This is critical for tasks involving proximity measures, such as clustering algorithms in machine learning.~\footnote{We only discuss the norms (and inner products) for real vector or matrix spaces. Most results can be applied directly to complex cases.}
\begin{definition}[Vector Norm and Matrix Nrom\index{Matrix norm}\index{Vector norm}]\label{definition:matrix-norm}
Given a norm $\norm{\cdot}: \real^n\rightarrow \real$ on vectors or a norm $\norm{\cdot}: \real^{m\times n}\rightarrow \real$ on matrices, for any vector $\bx \in \real^{n}$ and  any matrix $\bA \in \real^{m\times n}$, we have~\footnote{When $\norm{\bx}=\bzero$ for some nonzero vector $\bx$, the norm is called a \textit{semi-norm}.}~\footnote{When the vector has a single element, the norm can be understood as the absolute value operation.}~\footnote{When $\norm{\bI}=1$ for a matrix norm, the matrix norm is said to be \textit{normalized}.}
\begin{enumerate}
\item \textit{Positive homogeneity}. $\norm{\lambda \bA} = \abs{\lambda} \cdot \norm{\bA}$ or $\norm{\lambda \bx} = \abs{\lambda} \cdot \norm{\bx}$ for any $\lambda \in \real$.
\item \textit{Triangle inequality, a.k.a., subadditivity}. $\norm{\bA+\bB} \leq \norm{\bA}+\norm{\bB}$, or $\norm{\bx+\by} \leq \norm{\bx}+\norm{\by}$ for any matrices $\bA, \bB\in \real^{m\times n}$ or vectors $\bx,\by\in \real^n$.
\begin{enumerate}
\item The triangle inequality also indicates~\footnote{
$\big| \, \norm{\bx} - \norm{\by}\, \big| \leq \norm{\bx-\by}$ since $$
\left.
\begin{aligned}
\norm{\bx} &= \norm{\bx-\by+\by} \leq \norm{\bx-\by} + \norm{\by};\\
\norm{\by} &= \norm{\by-\bx+\bx} \leq \norm{\by-\bx} + \norm{\bx},\\
\end{aligned}
\right\}
\implies 
\begin{aligned}
\norm{\bx}-\norm{\by}  &\leq \norm{\bx-\by}; \\
\norm{\by} -\norm{\bx} &\leq  \norm{\by-\bx},
\end{aligned}
$$
}
\begin{equation}
\big| \, \norm{\bA} - \norm{\bB}\, \big| \leq \norm{\bA-\bB}
\quad\text{or}\quad
\big| \, \norm{\bx} - \norm{\by}\, \big| \leq \norm{\bx-\by}.
\end{equation}

\end{enumerate}
\item \textit{Nonnegativity}. $\norm{\bA} \geq 0$ or $\norm{\bx}\geq 0$.
\begin{enumerate}
\item  the equality holds if and only if $\bA=\bzero $ or $\bx=\bzero$. 
\end{enumerate}
\end{enumerate}
\end{definition}
\index{Vector norm}
\index{Matrix norm}
The vector space $\real^n$, together with a given norm $\norm{\cdot}$, is called a \textit{normed vector space}.
On the other hand, one way to define norms for matrices is by viewing a matrix $\bA\in\real^{m\times n}$ as a vector in $\real^{mn}$ through its vectorization.
What distinguishes a matrix norm is a property called \textit{submultiplicativity}: $\norm{\bA\bB}\leq \norm{\bA}\norm{\bB}$ if $\norm{\cdot}$ is a submultiplicative matrix norm (see discussions below). Almost all of the matrix norms we discuss are submultiplicative (Frobenius norm in Proposition~\ref{propo:submul_frob} and spectral norm in Proposition~\ref{propo:submul_spec}, both of which are special cases of the \textit{Schatten $p$-norm}, as a consequence of the singular value decomposition \citep{lu2021numerical}).

\index{Orthogonally invariant}
\index{Initarily invariant}
\begin{definition}[Orthogonally Invariant Norms]\label{definition:unitarily_invaria}
A matrix norm on $\real^{m\times n}$ is \textit{orthogonally invariant} if $\norm{\bU\bA\bV} =\norm{\bA}$ for all orthogonal $\bU\in\real^{m\times m}$ and $\bV\in\real^{n\times n}$ and for all $\bA\in\real^{m \times n}$; and  it is \textit{weakly orthogonally invariant} if  $\norm{\bU\bA\bU^\top} =\norm{\bA}$ for all orthogonal $\bU\in\real^{m\times m}$ and $\bA\in\real^{m\times m}$ is square.
Similarly, a vector norm on $\real^{n}$ is \textit{orthogonally invariant} if $\norm{\bQ\bx}=\norm{\bx}$ for all orthogonal $\bQ\in\real^{n\times n}$ and for all $\bx\in\real^n$.~\footnote{The term \textit{unitarily invariant} is used more frequently in the literature for complex matrices.}
\end{definition}


\index{Semi-norm}
\index{Semi-inner product}
\index{Inner product}
\begin{definition}[Inner Product]\label{definition:inner_prod}
In most cases, the norm can be derived from the vector \textit{inner product} $\langle\cdot, \cdot\rangle: \real^n\times \real^n\rightarrow \real$ (the inner product of vectors $\bx,\by\in\real^n$ is given by $\langle\bx,\by\rangle$), which satisfies the following three axioms:~\footnote{When $\inner{\bx}{\bx}=0$ for some nonzero $\bx$, the inner product is called a \textit{semi-inner product}.}
\begin{enumerate}
\item \textit{Commutativity}. $\langle\bx,\by\rangle = \langle\by,\bx\rangle$ for any $\bx,\by\in \real^n$. 
\item \textit{Linearity}. $\langle\lambda_1 \bx_1+\lambda_2\bx_2, \by\rangle = \lambda_1\langle\bx,\by\rangle+\lambda_2\langle\bx_2,\by\rangle$ for any $\lambda_1,\lambda_2 \in \real$ and $\bx,\by\in \real^n$.
\item \textit{Nonnegativity}. $\langle \bx,\bx \rangle \geq 0$ for any $\bx \in \real^n$.
\begin{enumerate}
\item $\langle \bx,\bx \rangle = 0$ if and only if $\bx=\bzero$.
\end{enumerate}
\end{enumerate}
Similarly, an inner product for matrices can be defined as a function $\langle\cdot, \cdot\rangle: \real^{m\times n}\times \real^{m\times n}\rightarrow \real$.
\end{definition}

For example, the \textit{Euclidean inner product}, defined as $\inner{\bu}{\bv}=\bu^\top\bv, \forall\bu,\bv\in\real^n$, is an inner product on vectors; the \textit{Frobenius inner product}, defined as $\inner{\bA}{\bB} = \trace(\bA^\top\bB)$, is an inner product on matrices.
Unless stated otherwise, we will use the Euclidean inner product throughout this book.

\begin{exercise}[Cauchy-Schwarz Inequality for Inner Product]
Given an inner product $\langle\cdot, \cdot\rangle: \real^n\times \real^n\rightarrow \real$, show that 
$$
\abs{\langle \bu, \bv\rangle}^2 \leq \inner{\bu}{\bu}\inner{\bv}{\bv} 
\gap 
\text{for all }
\gap 
\bu,\bv\in\real^n.
$$
\textit{Hint: Consider the vector $\bx \triangleq \inner{\bv}{\bv}\bu - \inner{\bu}{\bv}\bv$ and analyze $\inner{\bx}{\bx}\geq 0$.}
\end{exercise}

\begin{exercise}[Norm Derived from Inner Product, and its Property]\label{exercise:norm_innr_pro}
Let $\langle\cdot, \cdot \rangle$ be an inner product on $\real^n\times \real^n$. Then the function $\norm{\cdot}=\langle\cdot, \cdot\rangle^{1/2}:  \real^n\rightarrow [0, \infty)$ is a norm on $\real^n$.
Show that the function satisfies 
\begin{equation}\label{equation:norm_inn_pro}
\textbf{(Parallelogram identity):}\,\, \norm{\bu+\bv}^2+\norm{\bu-\bv}^2 = 2(\norm{\bu}^2+\norm{\bv}^2), \,\, \forall \bu,\bv\in\real^n.
\end{equation} 
(This is known as the \textit{parallelogram identity}: consider a parallelogram in a vector space where the sides are represented by the vectors $\bu$ and $\bv$, the diagonals of this parallelogram are represented by the vectors $\bu+\bv$ and $\bu-\bv$.)
Prove the following \textit{polarization identity}:
\begin{equation}\label{equation:norm_pol_pro}
\textbf{(Polarization identity):} \gap \inner{\bu}{\bv} = \frac{1}{4}(\norm{\bu+\bv}^2-\norm{\bu-\bv}^2).
\end{equation}
\end{exercise}
Therefore, the vector space $\real^n$, when equipped with an inner product $\inner{\cdot}{\cdot}$, is called an \textit{inner product space}, which is also a normed vector space with the derived norm.

\index{$\ell_p$ norm}
\index{$\ell_2$ norm}
\index{$\ell_1$ norm}
\index{$\ell_\infty$ norm}
\subsection{Vector Norm}\label{section:vector-norm}
The vector norm is derived from the definition of the inner product. In most cases, the inner product in $\real^n$ is implicitly  referred to as  the  \textit{dot product (a.k.a., Euclidean inner product)}, defined by 
$
\langle\bx,\by\rangle = \sum_{i=1}^{n} x_i y_i.
$
The $\ell_2$ \textit{norm} (also called the $\ell_2$ vector norm) is induced by the dot product and is given by
$
\normtwo{\bx} = \sqrt{\langle\bx,\bx\rangle} = \sqrt{\sum_{i=1}^{n}x_i^2}.
$
More generally, for any $p\geq 1$ \footnote{When $p<1$, the $\ell_p$ function does not satisfy the third axiom of Definition~\ref{definition:matrix-norm}, meaning it is not a valid norm.}, the $\ell_p$ \textit{norm} (a.k.a., the $\ell_p$ vector norm or \holders norm) is given by 
\begin{equation}\label{equation:l_p_norm}
\textbf{($\ell_p$ {norm})}:\qquad \norm{\bx}_p = \sqrt[p]{ \sum_{i=1}^{n}\abs{x_i}^p  },
\quad 
p\geq 1.
\end{equation}
This norm satisfies positive homogeneity and nonnegativity by definition, while the triangle inequality follows from Minkowski's inequality (Theorem~\ref{theorem:minkoski}; alternatively, we can prove the triangle inequality for the $\ell_p$ norm using H{\"o}lder's inequality; see Exercise~\ref{exercise:triangle_ineq}).
From the general definition of the $\ell_p$ norm,  the $\ell_1$ and $\ell_\infty$ \textit{norms}  can be obtained by, respectively, 
$$
\norm{\bx}_1 = \sum_{i=1}^{n} \abs{x_i}
\gap \text{and}\gap
\norm{\bx}_\infty = \mathop{\max}_{i=1,2,\ldots,n} \abs{x_i} .
$$

Following the definition of the $\ell_p$ norm, we can obtain the famous H{\"o}lder's inequality in Theorem~\ref{theorem:holder-inequality}. Thus, the $\ell_p$ norm is sometimes referred to as  \textit{H{\"o}lder's norm}. 
Conversely, H{\"o}lder's inequality can be used to prove the validity of the  $\ell_p$ norm. This proof is left as an exercise; see, for example, \citet{lu2021numerical}.
\begin{exercise}[Triangle Inequality]\label{exercise:triangle_ineq}
Prove the triangle inequality for the $\ell_p$ norm using H{\"o}lder's inequality.
For the special case $p=2$, use the Cauchy-Schwartz inequlaity (Proposition~\ref{proposition:cauchy-schwarz-inequ}) to prove the triangle inequality.
\end{exercise}

\begin{exercise}[Orthogonal Invariance of $\ell_2$]\label{exercise:orthogo_ell2}
Show that the $\ell_2$ norm is orthogonally invariant: $\norm{\bQ\bx}_2=\normtwo{\bx}$ for all $\bx\in\real^n$ if $\bQ$ is orthogonal.
Investigate whether this property holds for other  $\ell_p$ norms.
\end{exercise}

\section*{Properties}
For any vector norm, we also have the following properties.
\begin{proposition}[$\ell_p$ Norm Inequalities]\label{prop:lp_norm_ineqs}
Given any $\bx\in\real^n$, the function $f(p)=\norm{\bx}_p$ is nonincreasing for $p=[1,\infty)$. Therefore, 
\begin{equation}
\norm{\bx}_\infty\leq \norm{\bx}_q \leq \norm{\bx}_p\leq \normtwo{\bx} \leq \normone{\bx},
\quad \text{if }  \infty \geq q \geq p\geq 2\geq 1.
\end{equation}
If $\bx$ has at least two nonzero elements, the inequalities can be strict.
This behavior can also be observed in the figures of unit balls (Figure~\ref{fig:p-norm-comparison-3d}).
Moreover, for any $\infty\geq p \geq 1$, we also have the bound
\begin{equation}
\norm{\bx}_\infty \leq \norm{\bx}_p \leq  n^{1/p} \norm{\bx}_\infty,
\quad \text{if }  \infty \geq  p \geq 1.
\end{equation}
\end{proposition}
\begin{proof}[of Proposition~\ref{prop:lp_norm_ineqs}]
Let $S(p) \triangleq \norm{\bx}_p^p$. Taking the logarithm of $f(p)$ gives
$\ln f(p) = \frac{1}{p} \ln S(p). $
Differentiate both sides with respect to $p$:
$$\frac{f'(p)}{f(p)} = \frac{d}{dp} \left( \frac{1}{p} \ln S(p) \right)
= -\frac{1}{p^2} \ln S(p) + \frac{1}{p} \cdot \frac{S'(p)}{S(p)}
. 
$$
This shows
$$f'(p) = f(p) \left( -\frac{1}{p^2} \ln S(p) + \frac{1}{p} \cdot \frac{S'(p)}{S(p)} \right)
=
\frac{1}{p}\norm{\bx}_p^{1-p}\sum_{i=1}^{n} \alpha_i,
$$
where 
$$
S'(p)=\sum_{i=1}^{n}\abs{x_i}^p \ln\abs{x_i}
\quad \text{and} \quad 
\alpha_i =
\left\{
\begin{aligned}
&\abs{x_i}^p (\ln \abs{x_i} - \ln \norm{\bx}_p), &\text{if } x_i\neq 0;\\
&0, &\text{if } x_i=0.
\end{aligned}
\right.
$$
This shows $f'(p)\leq 0$ since $\ln \abs{x_i} - \ln \norm{\bx}_p\leq 0$.
If $\bx$ has at least two nonzero elements, $f'(p)< 0$ and $f(p)$ is decreasing.

For the second part, following the definition of the $\ell_p$ norm, we notice that 
$$
(\norm{\bx}_\infty )^p = \left(\mathop{\max}_{i=1,2,\ldots,n} \abs{x_i}\right)^p \leq \sum_{i=1}^{n} \abs{x_i}^p \leq n \mathop{\max}_{i=1,2,\ldots,n} \abs{x_i}^p = n (\norm{\bx}_\infty )^p,
$$
from which the result follows.
\end{proof}

\begin{proposition}[Orthogonally Invariant Vector Norm]\label{proposition:orinv_vecnorm}
Let $\norm{\cdot}$ be an orthogonally invariant norm on $\real^n$, and let $\norm{\cdot}_2$ be the $\ell_2$ norm. Then, for any $\bx\in\real^n$, we have $\norm{\bx}=\normtwo{\bx}\norm{\be_1}$.
\end{proposition}
\begin{proof}[of Proposition~\ref{proposition:orinv_vecnorm}]
Define $\gamma\triangleq\normtwo{\bx}$ and let $\widehat{\bx}\triangleq\bx/\gamma$, so that $\bx=\gamma \widehat{\bx}$. 
Consider  an orthogonal matrix $\bU$  whose first column is $\widehat{\bx}$.
By the definition of a norm and the orthogonally invariant property, we have $\norm{\bx}=\gamma\norm{\widehat{\bx}}=\gamma\norm{\bU\be_1}=\gamma\norm{\be_1}$. 
\end{proof}
The result also implies that the $\ell_2$ norm is the only orthogonally invariant norm satisfying $\norm{\be_1}=1$.

We conclude this section by introducing an important property of vector norms that  is frequently useful. 
The following theorem on the equivalence of vector norms states that if a vector is small in one norm, it is also small in another norm, and vice versa.
\begin{theorem}[Equivalence of Vector Norms]\label{theorem:equivalence-vector-norm}
Let $\norm{\cdot}_a$ and $\norm{\cdot}_b$ be two different vector norms, where $\norm{\cdot}_a$, $\norm{\cdot}_b$: $\real^n\rightarrow \real$. Then, there exist positive scalars $\alpha$ and $\beta$ such that for all $\bx \in \real^n$, the following inequality holds:
$$
\alpha\norm{\bx}_a \leq \norm{\bx}_b \leq \beta\norm{\bx}_a.
$$
This implies 
$$
\frac{1}{\beta}\norm{\bx}_b  \leq \norm{\bx}_a \leq \frac{1}{\alpha}\norm{\bx}_b.
$$
This justifies the term ``equivalence" in the theorem.
\end{theorem}
\begin{proof}[of Theorem~\ref{theorem:equivalence-vector-norm}]
In advanced calculus, it is stated that $\mathop{\sup}_{\bx\in \sS} f(\bx)$ is attained for some vector $\bx\in \sS$ as long as $f(\cdot)$ is continuous and $\sS$ is a compact set (closed and bounded); see Theorem~\ref{theorem:weierstrass_them}. When the supremum is an element in the set $\sS$, this supremum is known as the maximum such that $\mathop{\sup}_{\bx\in\sS} f(\bx) =\mathop{\max}_{\bx\in\sS} f(\bx)$.
Without loss of generality, we assume $\bx\neq \bzero$. Then we have:
$$
\begin{aligned}
\norm{\bx}_b &= \frac{\norm{\bx}_b}{\norm{\bx}_a} \norm{\bx}_a 
\leq \mathop{\sup}_{\bz\neq 0} \frac{\norm{\bz}_b}{\norm{\bz}_a} \norm{\bx}_a
= \mathop{\sup}_{\norm{\by}_a=1} \norm{\by}_b \norm{\bx}_a
=\norm{\bx}_a \mathop{\max}_{\norm{\by}_a=1} \norm{\by}_b.
\end{aligned}
$$
The last equality holds since $\{\by: \norm{\by}_a=1\}$ is a compact set.
By setting $\beta\triangleq \mathop{\max}_{\norm{\by}_a=1} \norm{\by}_b  $, we have $\norm{\bx}_b \leq \beta\norm{\bx}_a$.
From the above argument, there exists a $\tau$ such that 
$
\norm{\bx}_a \leq \tau \norm{\bx}_b.
$
Letting $\alpha \triangleq \frac{1}{\tau}$, we obtain $\alpha\norm{\bx}_a \leq \norm{\bx}_b$, from which the result follows.
\end{proof}

\begin{example}[Equivalence of Vector Norms]
The following inequalities hold for all $\bx\in\real^n$:
$$
\begin{aligned}
\norm{\bx}_\infty &\leq \norm{\bx}_1 \leq n\norm{\bx}_\infty; \\
\norm{\bx}_\infty &\leq \normtwo{\bx} \leq \sqrt{n}\norm{\bx}_\infty; \\
\normtwo{\bx} &\leq \norm{\bx}_1 \leq \sqrt{n}\normtwo{\bx}. \\
\end{aligned}
$$
This demonstrates the equivalence of the $\ell_1, \ell_2$, and $\ell_\infty$ vector norms.
\end{example}

\section*{Dual Norm}
Consider the $\ell_p$ vector norm. From \holders inequality (Theorem~\ref{theorem:holder-inequality}), we have 
$
\bx^\top\by \leq \norm{\bx}_p \norm{\by}_q,
$
where $p,q>1$ satisfy $\frac{1}{p}+\frac{1}{q}=1$, and $\bx,\by\in \real^n$. Equality holds if the two sequences $\{\abs{x_i}^p\}$ and $\{\abs{y_i}^q\}$ are linearly dependent. This implies
\begin{equation}\label{equation:dual_norm_equa}
\mathop{\max}_{\norm{\by}_q=1} \bx^\top\by = \norm{\bx}_p.
\end{equation}
For this reason, $\norm{\cdot}_q$ is called the \textit{dual norm} of $\norm{\cdot}_p$.
On the other hand, for each $\bx\in \real^n$ with $\norm{\bx}_p=1$, there exists a vector $\by\in \real^n$ such that $\norm{\by}_q=1$ and $\bx^\top\by=1$.
Notably, the $\ell_2$ norm is self-dual, while the $\ell_1$ and $\ell_\infty$ norms are dual to each other.

\begin{definition}[Set of Primal Counterparts]\label{definition:set_primal}
Let $ \norm{\cdot} $ be any norm on $ \real^n $. Then the \textit{set of primal counterparts of $\ba$} is defined as 
\begin{equation}
\Lambda_{\ba} = \argmax_{\bu \in \real^n} \left\{ \innerproduct{\ba, \bu} \mid  \norm{\bu} \leq 1 \right\}.
\end{equation}
That is, $\innerproduct{\ba, \ba^\dagger} = \norm{\ba}_*$ for any $\ba^\dagger\in \Lambda_{\ba}$, where $\norm{\cdot}_{*}$ denotes the dual norm.
It follows that 
\begin{enumerate}[(i)]
\item If $\ba \neq \bzero$, then $\norm{\ba^\dagger} = 1$ for any $\ba^\dagger \in \Lambda_{\ba}$.
\item If $\ba = \bzero$, then $\Lambda_{\ba} = \{\bx\in\real^n\mid  \norm{\bx} \leq 1\}$.
\end{enumerate}
\end{definition}

\begin{example}[Set of Primal Counterparts]\label{example:set_primal_count}
A few examples for the sets of primal counterparts are shown below:
\begin{itemize}
\item If the norm is the  $\ell_2$ norm, then for any $\ba \neq \bzero$,
$
\Lambda_{\ba} = \left\{ {\ba}/{\normtwo{\ba}} \right\}.
$

\item If the norm is the  $\ell_1$ norm, then for any $\ba \neq \bzero$,
$$
\Lambda_{\ba} = \left\{ \sum_{i \in \sI(\ba)} \lambda_i \sign(a_i) \be_i \mid \sum_{i \in \sI(\ba)} \lambda_i = 1, \lambda_j \geq 0, j \in \sI(\ba) \right\},
$$
where $\sI(\ba) \triangleq \argmax_{i=1,2,\ldots,n} |a_i|$.

\item If the norm is the  $\ell_\infty$ norm, then for any $\ba \neq \bzero$,
$$
\Lambda_{\ba} = \left\{ \bx \in \real^n \mid x_i = \sign(a_i), i \in \sI_{\neq}(\ba), |x_j| \leq 1, j \in \sI_0(\ba) \right\},
$$
where
$
\sI_{\neq}(\ba) \triangleq \left\{ i \in \{1, 2, \ldots, n\} \mid a_i \neq 0 \right\} $ and $ \sI_0(\ba) \triangleq \left\{ i \in \{1, 2, \ldots, n\} \mid a_i = 0 \right\}.
$
\end{itemize}
These examples play a crucial role in the development of non-Euclidean gradient descent methods, which will be discussed in Sections~\ref{section:als-gradie-descent-taylor} and \ref{section:noneucli_gd}.
\end{example}

\index{Unit ball}
\section*{Unit Ball}
\begin{figure}[h]
\centering  
\vspace{-0.25cm} 
\subfigbottomskip=2pt 
\subfigcapskip=-5pt 
\subfigure[$p=\infty$.]{\label{fig:p-norm-3d1}
\includegraphics[width=0.18\linewidth]{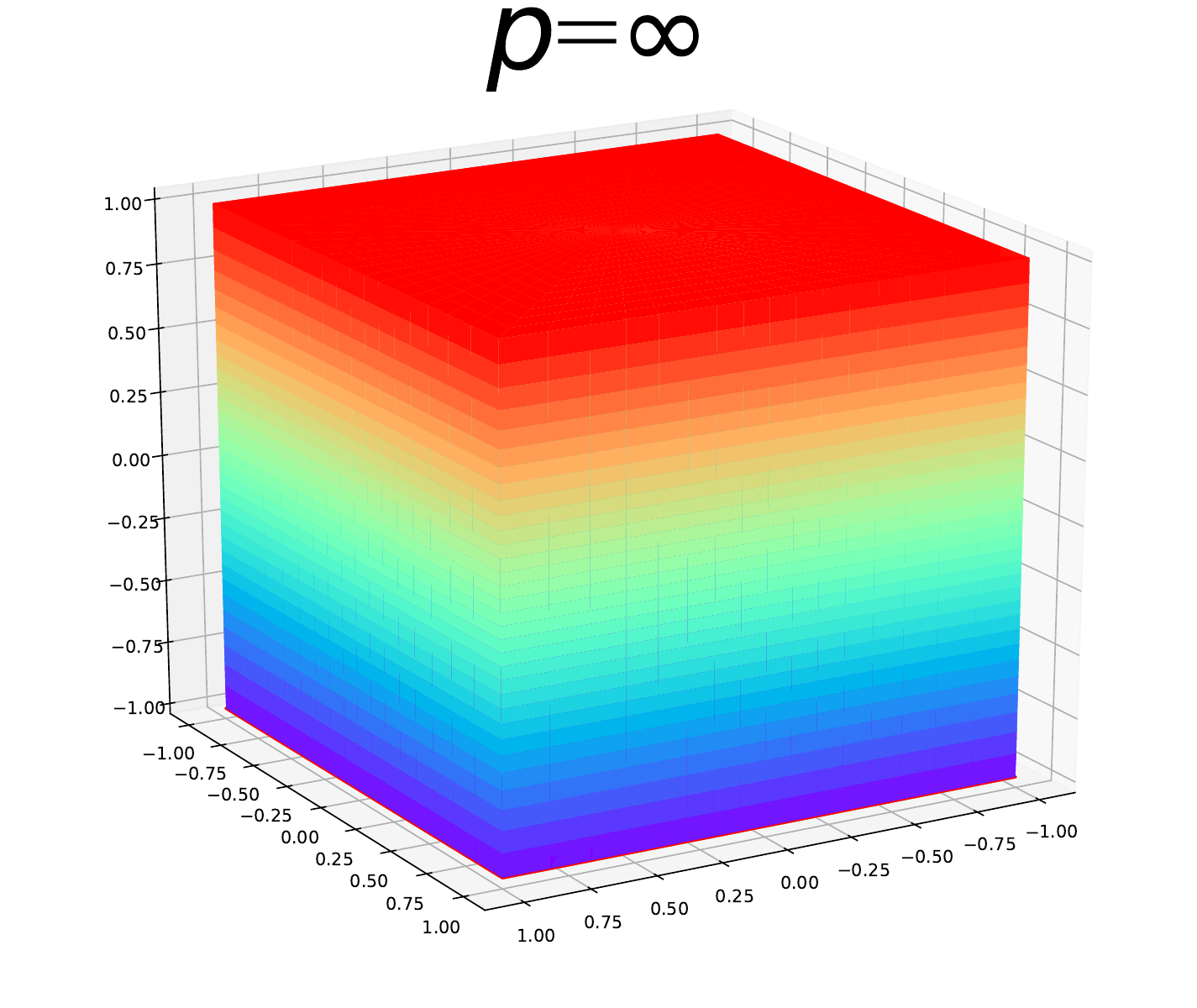}}
\subfigure[$p=2$.]{\label{fig:p-norm-3d2}
\includegraphics[width=0.18\linewidth]{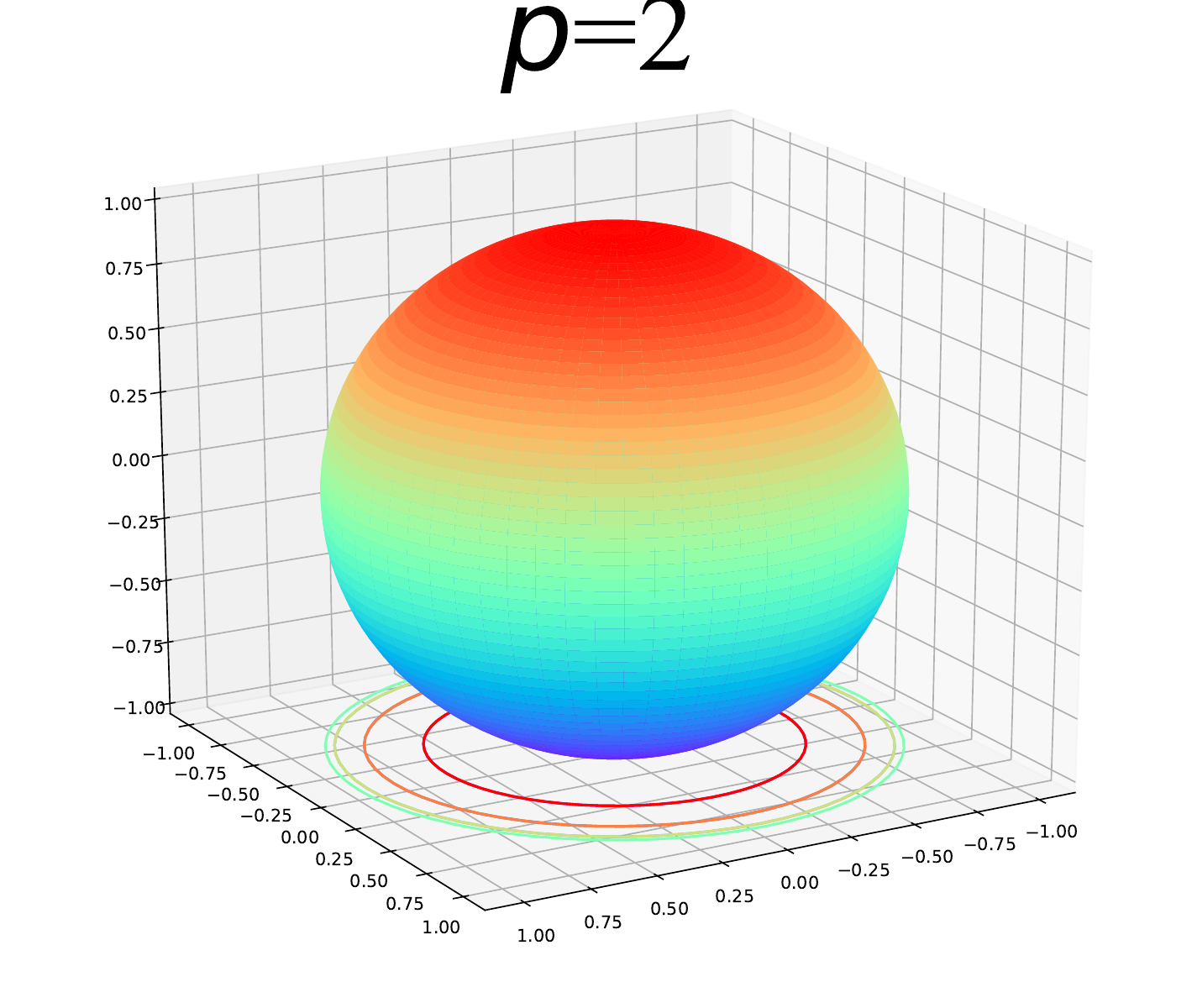}}
\subfigure[$p=1$.]{\label{fig:p-norm-3d3}
\includegraphics[width=0.18\linewidth]{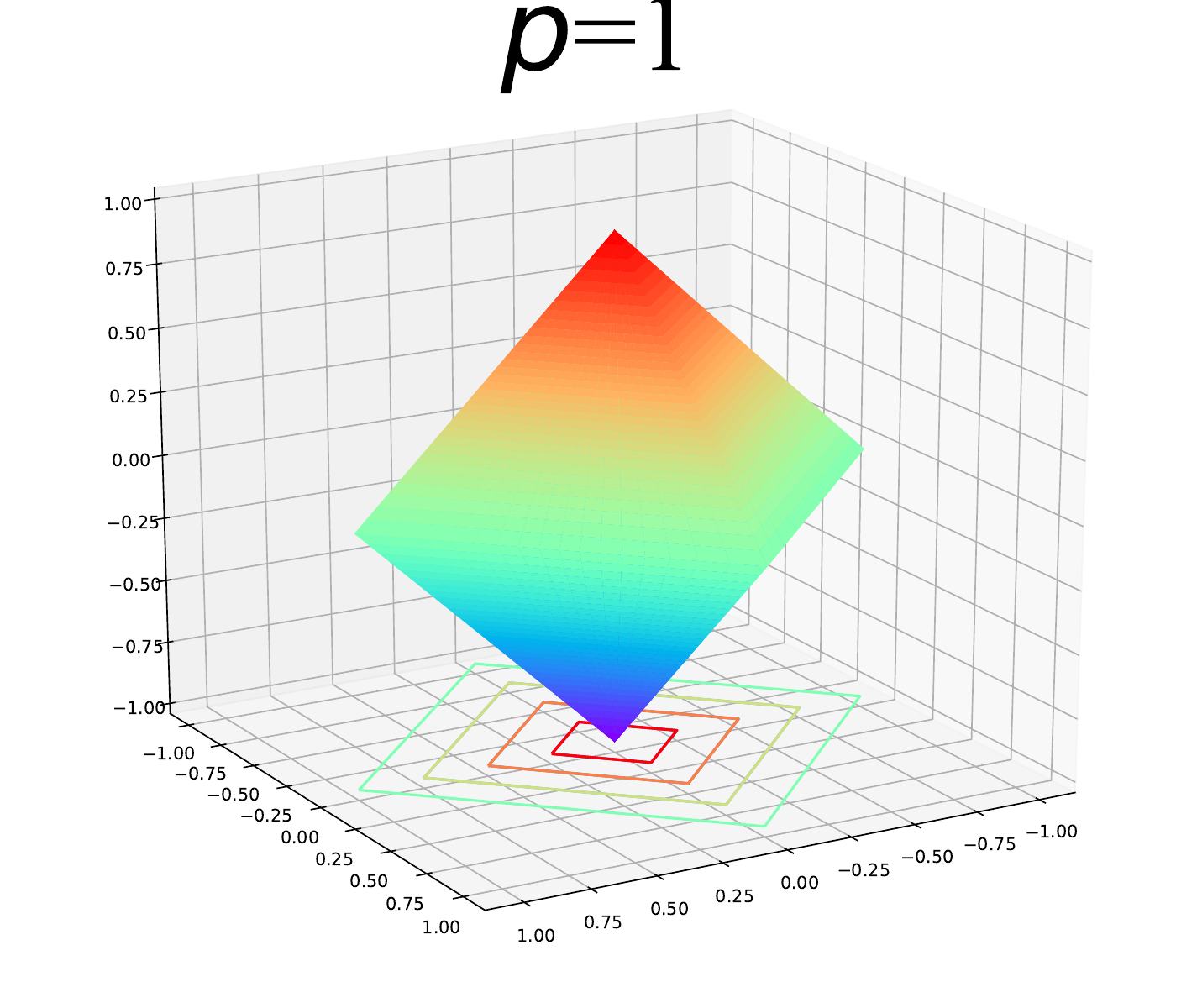}}
\subfigure[$p=0.5$.]{\label{fig:p-norm-3d4}
\includegraphics[width=0.18\linewidth]{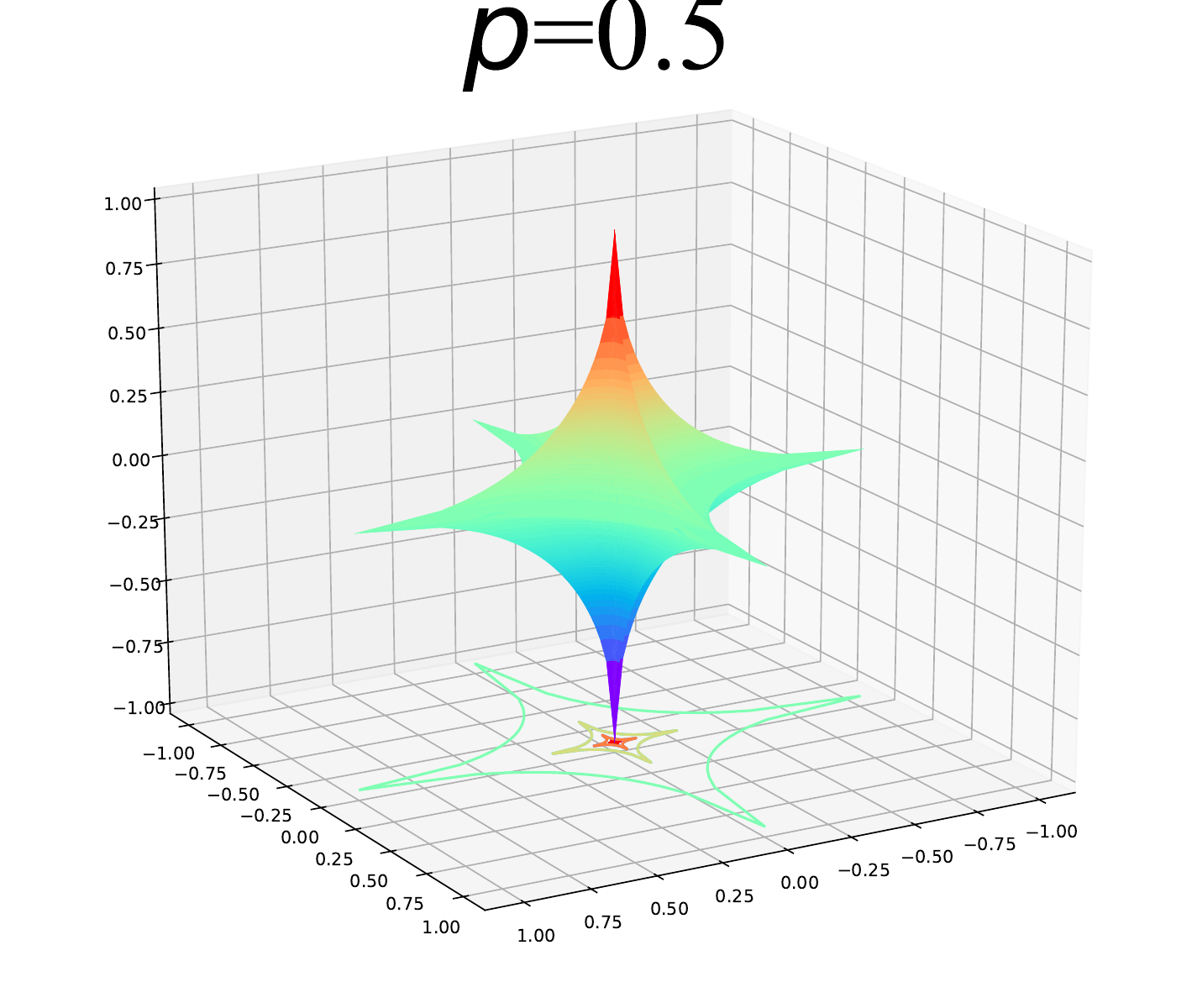}}
\subfigure[$p=0$.]{\label{fig:p-norm-3d5}
\includegraphics[width=0.18\linewidth]{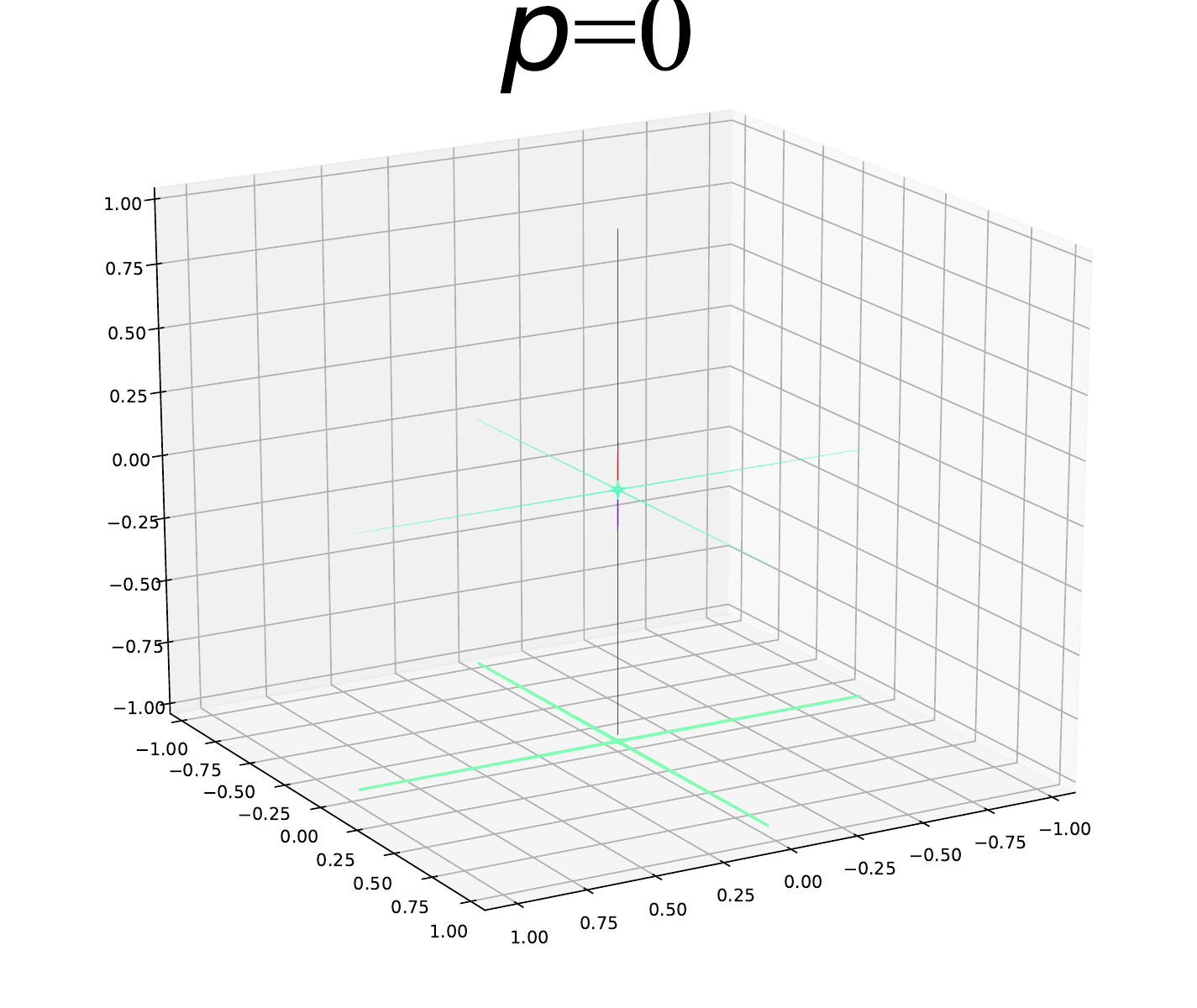}}
\caption{Unit ball of $\ell_p$ norms in three-dimensional space. When $p<1$, the metric does not qualify as a norm since it does not satisfy the third axiom of the norm in Definition~\ref{definition:matrix-norm}.}
\label{fig:p-norm-comparison-3d}
\end{figure}

The \textit{unit ball} of a norm is the set of all points whose distance from the origin (i.e., the zero vector) equals 1.
If the distance is defined by the $\ell_p$ norm,  the unit ball is the collection of 
$$
\sB_p = \{\bx: \norm{\bx}_p=1\}.
$$
The comparison of the  $\ell_p$ norm in three-dimensional space with different values of $p$ is depicted in  Figure~\ref{fig:p-norm-comparison-3d}.

Vector norms are fundamental in machine learning. In Section~\ref{section:pre_ls}, we will discuss how least squares aim to minimize the squared distance between the observed value $\bb$ and the predicted value $\bA\bx$: $\normtwo{\bA\bx-\bb}$, which is the $\ell_2$ norm of $\bA\bx-\bb$.
Alternatively, minimizing the $\ell_1$ norm of the difference between the observed and predicted values can yield a more robust estimation of $\bx$, particularly in the presence of outliers  \citep{zoubir2012robust}.

\section*{$\bQ$-Inner Product, $\bQ$-Norm, $\bS$-Norm}
The standard dot product is not the only possible inner product that can be defined on $\real^n$. Given a positive definite $n\times n$ matrix $\bQ$, a $\bQ$\textit{-dot product} can be defined as~\footnote{When $\bQ$ is positive semidefinite, the $\bQ$-dot product is a semi-inner product.}
$$
\langle \bx,\by \rangle_{\bQ} = \bx^\top\bQ\by.
$$
One can verify that the $\bQ$-dot product defined above satisfies the three axioms of an inner product, as discussed at the beginning of this section. When $\bQ=\bI$, this reduces to the standard Euclidean inner product (dot product). From the $\bQ$-dot product, we can define the corresponding \textit{$\bQ$-$\ell_2$-norm} as 
$
\norm{\bx}_{\bQ,2} = \sqrt{\bx^\top\bQ\bx}.
$
More generally, given a positive definite matrix $\bQ\in\real^{n\times n}$ and  a general norm $\norm{\cdot}$ on $\real^n$, we define the \textit{$\bQ$-norm} as:
\begin{equation}\label{equation:q_norm}
\textbf{$\bQ$-norm:} \qquad\norm{\bx}_{\bQ} = \norm{\sqrt{\bQ}\bx}
\end{equation}	
Similarly, given a matrix $\bS\in\real^{m\times n}$ with full column rank and a general norm $\norm{\cdot}$ on $\real^m$,  the \textit{$\bS$-norm} is defined as:
\begin{equation}\label{equation:s_norm}
\textbf{$\bS$-norm:} \qquad \norm{\bx}_{\bS}=\norm{\bS\bx},
\end{equation}	
which is a norm on $\real^n$.
To establish this rigorously, we introduce the following lemma.
\begin{lemma}[Construct Norms from Other Norms]\label{lemma:construct_norm}
Let $\norm{\cdot}$ be a norm on $\real^m$. Given a matrix $\bS\in \real^{m\times n}$, then $\norm{\bS(\cdot)}$ defines a semi-norm~\footnote{It may be zero for nonzero vectors in $\real^n$.} on $\real^n$. 
If $\bS$ has full column rank $n$, then $\norm{\bS(\cdot)}$ is a norm on $\real^m$.
\end{lemma}
\begin{proof}[of Lemma~\ref{lemma:construct_norm}]
For any vectors $\bx, \by\in \real^m$, following from Definition~\ref{definition:matrix-norm}, we have 
$$
\norm{\bx}\geq 0, \gap \norm{\lambda \bx} = \abs{\lambda}\cdot  \norm{\bx}, \gap \norm{\bx+\by} \leq \norm{\bx} + \norm{\by}.
$$
Now, for any $\bc, \bd \in \real^n$, it follows that:
$$
\norm{\bS\bc}\geq 0, \gap \norm{\lambda \bS\bc} = \abs{\lambda}\cdot  \norm{\bS\bc}, \gap \norm{\bS\bc+\bS\bd} \leq \norm{\bS\bc} + \norm{\bS\bb}.
$$
Therefore, $\norm{\bS(\cdot)}$ is a semi-norm.
If $\bS$ has full column rank, $\bS\bc=\bzero$ only if $\bc=\bzero$, ensuring that $\norm{\bS(\cdot)}$ satisfies the norm definiteness property. Hence, it is a norm. This completes the proof.
\end{proof}


\begin{exercise}[Weighted $\ell_p$ Norm]
Let $w_1, w_2, \ldots,w_p$ be positive real numbers, and let $p\geq 1$. Show that  the weighted $\ell_p$ norm $\norm{\bx}=\sqrt[p]{ \sum_{i=1}^{n}w_i\abs{x_i}^p  }$  is a valid norm on $\real^n$. \textit{Hint: Show that it is a norm of the form $\norm{\bS\bx}_p$.} 
\end{exercise}


\subsection{Matrix Norm}\label{appendix:matrix-norm-sect2}

\index{Submultiplicativity}
\index{Submultiplicative matrix norms}
\subsection*{Submultiplicativity of Matrix Norms}
In some texts, a matrix norm that is not \textit{submultiplicative} is referred to as a \textit{vector norm on matrices} or a \textit{generalized matrix norm}.
The submultiplicativity of a matrix norm plays a crucial role in the analysis of square matrices. However, it's important to note that the definition of a matrix norm applies to both square and rectangular matrices.
For a submultiplicative matrix norm $\norm{\cdot}$ satisfying $\norm{\bA\bB}\leq \norm{\bA}\norm{\bB}$, considering $\bA\in\real^{n\times n}$, it follows that 
\begin{equation}\label{equation:power_subm}
\norm{\bA^2} \leq \norm{\bA}^2
\quad\implies\quad
\norm{\bA^k} \leq \norm{\bA}^k, \forall k\in\{1,2,\ldots,\}.
\end{equation}
Therefore, if the matrix is idempotent, i.e., $\bA^2=\bA$, we have $\norm{\bA}\geq 1$. This also indicates:
\begin{equation}\label{equation:power_subm2}
\norm{\bI}\geq 1,
\quad\text{if}\quad
\norm{\cdot} \text{ is submultiplicative}.
\end{equation}
On the other hand, if $\bA$ is nonsingular, then for submultiplicative norms, we have the inequality:
$$
1\leq \norm{\bI}=\norm{\bA\bA^{-1}}\leq \norm{\bA}\norm{\bA^{-1}}.
$$  
This means that for a submultiplicative norm, $\norm{\bI}\geq 1$, and it is considered \textit{normalized} if and only if $\norm{\bI}=1$.

\paragraph{Bounds on the spectral radius.}
The property of submultiplicativity can be utilized to establish bounds on the \textit{spectral radius} of a matrix  (i.e., the largest magnitude of the eigenvalues of a matrix).
Given an eigenpair $(\lambda,\bx)$ of a matrix $\bA\in\real^{n\times n}$, and consider the matrix $\bX=\bx\bone^\top=[\bx,\bx,\ldots,\bx]$. 
For a submultiplicative matrix norm $\norm{\cdot}$ on $\real^{n\times n}$, it holds that 
$
\abs{\lambda}\norm{\bX}
=
\norm{\lambda\bX}
=
\norm{\bA\bX}
\leq \norm{\bA}\norm{\bX}.
$
Thus, the spectral radius is bounded by the matrix norm $\norm{\bA}$.
Similarly, we can prove that $\abs{\lambda^{-1}}\leq \norm{\bA^{-1}}$.
Consequently, we obtain lower and upper bounds on the spectral radius $\rho(\bA)$:
\begin{equation}\label{equation:bounds_submulnorm}
\textbf{Bounds on spectral radius:} \qquad \frac{1}{\norm{\bA^{-1}}}\leq \rho(\bA)\leq \norm{\bA}.
\end{equation}
Additionally, the submultiplicative property aids in understanding bounds on the spectral radius of the product of matrices. 
For instance, given matrices $\bA$ and $\bB$ (provided  the matrix product $\bA\bB$ is defined), we have
\begin{equation}
\rho(\bA\bB)\leq \norm{\bA\bB} \leq \norm{\bA}\norm{\bB}.
\end{equation}

\index{Convergent matrices}
\paragraph{Power of square matrices.}
When $\norm{\bA}\leq 1$ for a square matrix $\bA\in\real^{n\times n}$, \eqref{equation:power_subm} shows that $\norm{\bA^k} \leq \norm{\bA}^k\stackrel{k\rightarrow \infty}{=}0$.
This holds true regardless of the specific norm used, provided it is submultiplicative. Thus, if $\bA^k$ tends to the zero matrix when $k\rightarrow \infty$ if $\norm{\bA}\leq 1$.
A matrix with this property is called \textit{convergent}.
The convergence of a matrix can also be characterized by its spectral radius.
\begin{exercise}[Convergence Matrices]\label{exer:conv_mat}
Let $\bA\in\real^{n\times n}$. Show that $\mathoplim{k\rightarrow \infty}\bA^k=\bzero$ if and only if the spectral radius $\rho(\bA)<1$.
\textit{Hint: Examine $\bA^k\bx=\lambda^k\bx$.}
\end{exercise}

\begin{exercise}[Convergence Matrices]
Let $\bA\in\real^{n\times n}$ and  $\epsilon>0$. Show that there is a constant $C=C(\bA,\epsilon)$ such that $\abs{(\bA^k)_{ij}} \leq C(\rho(\bA)+\epsilon)^k$ for $k\in\{1,2,\ldots\}, i,j\in\{1,2,\ldots,n\}$.
\textit{Hint: Use Exercise~\ref{exer:conv_mat} and examine $\bB=\frac{1}{\rho(\bA)+\epsilon}\bA$, whose spectral radius is strictly less than 1.}
\end{exercise}

\index{Gelfand formula}
The Gelfand formula offers a method to estimate the spectral radius of a matrix using submultiplicative norms.
\begin{exercise}[Gelfand Formula]\label{exercise:gelfand_formula}
Let $\bA\in\real^{n\times n}$ and let $\norm{\cdot}$ be a submultiplicative matrix norm on $\real^{n\times n}$. 
Show that $\rho(\bA)=\mathoplim{k\rightarrow \infty}\norm{\bA^k}^{1/k}$.
\textit{Hint: Consider $\rho(\bA^k)=\rho(\bA)^k\leq \norm{\bA}$ and examine $\bB=\frac{1}{\rho(\bA)+\epsilon}\bA$, whose spectral radius is strictly less than 1.}
\end{exercise}

\index{Frobenius}
\subsection*{Frobenius Norm}
The norm of a matrix serves a similar purpose to the norm of a vector.
One important matrix norm is the Frobenius norm, which can be considered the matrix equivalent of the  $\ell_2$ vector norm.
\begin{definition}[Frobenius Norm]\label{definition:frobenius}
The Frobenius norm of a matrix $\bA\in \real^{m\times n}$ is defined as 
$$
\norm{\bA}_F = \sqrt{\sum_{i=1,j=1}^{m,n} (a_{ij})^2}=\sqrt{\trace(\bA\bA^\top)}=\sqrt{\trace(\bA^\top\bA)} = \sqrt{\sigma_1^2+\sigma_2^2+\ldots+\sigma_r^2},
$$
where the values of $\sigma_i$ are the singular values of $\bA$, and $r$ is the rank of $\bA$. 
This represents the square root of the sum of the squares of all elements in $\bA$. 
\end{definition}
The equivalence of $\sqrt{\sum_{i=1,j=1}^{m,n} (a_{ij})^2}$, $\sqrt{\trace(\bA\bA^\top)}$, and $\sqrt{\trace(\bA^\top\bA)}$ is straightforward. 
The equivalence between $\sqrt{\trace(\bA\bA^\top)}$ and $\sqrt{\sigma_1^2+\sigma_2^2+\ldots+\sigma_r^2}$ can be shown using the singular value decomposition (SVD). Suppose $\bA$ admits the SVD $\bA = \bU\bSigma\bV^\top$, then:
$$
\sqrt{\trace(\bA\bA^\top)} = \sqrt{\trace(\bU\bSigma\bV^\top \bV\bSigma\bU^\top)} = \sqrt{\trace(\bSigma^2)}=\sqrt{\sigma_1^2+\sigma_2^2+\ldots+\sigma_r^2}.
$$
Apparently, the Frobenius norm can  also be defined using the  vector $\ell_2$ norm such that $\norm{\bA}_F = \sqrt{\sum_{i=1}^{n} \norm{\ba_i}^2}$, where $\ba_i$ for all $i \in \{1,2,\ldots, n\}$ are the columns of $\bA$.

\index{Submultiplicativity}
\begin{proposition}[Submultiplicativity of Frobenius]\label{propo:submul_frob}
The Frobnenius norm is submultiplicative. That is, $\norm{\bA\bB}_F\leq \norm{\bA}_F\norm{\bB}_F$.
\end{proposition}
\begin{proof}[of Proposition~\ref{propo:submul_frob}]
Suppose $\bA\in\real^{m\times n}$ and $\bB\in\real^{n\times p}$. We have 
$$
\norm{\bA\bB}_F = 
\bigg(\sum_{i,j=1}^{m,p} \big(\sum_{k=1}^{n}a_{ik}b_{kj} \big)^2 \bigg)^{1/2}
\leq 
\bigg(\sum_{i,j=1}^{m,p} \big(\sum_{k=1}^{n}a_{ik}^2 \big)\big(\sum_{k=1}^{n}b_{kj}^2 \big) \bigg)^{1/2}
=
\norm{\bA}_F\norm{\bB}_F.
$$
This completes the proof.
\end{proof}

\index{Orthogonally invariance}
\begin{proposition}[Orthogonally Invariance of Frobenius]\label{proposition:frobenius-orthogonal-equi}
Let $\bA\in \real^{m\times n}$ be given, and let $\bU\in \real^{m\times m}$ and $\bV\in \real^{n\times n}$ be orthogonal matrices. Then,
$
\norm{\bA }_F =  \norm{\bU\bA\bV }_F.
$
\end{proposition}
\begin{proof}[of Proposition~\ref{proposition:frobenius-orthogonal-equi}]
We observe that 
$$
\begin{aligned}
\norm{\bU\bA\bV}_F  &=\sqrt{\trace((\bU\bA\bV)(\bU\bA\bV)^\top)} = \sqrt{\trace(\bU\bA\bA^\top\bU^\top)}\\
&= \sqrt{\trace(\bA\bA^\top\bU^\top\bU)} = \sqrt{\trace(\bA\bA^\top)}=
\norm{\bA}_F ,
\end{aligned}
$$
where the third equality holds due to the cyclic property of the trace.
\end{proof}

\index{Schur inequality}
The Frobenius norm is defined as the square root of the sum of the squares of the matrix elements. Additionally, we have the well-known \textit{Schur inequality}, which follows from the definition of the Frobenius norm.
\begin{theorem}[Schur Inequality]\label{theorem:schur_inequality}
Let $\lambda_1, \lambda_2, \ldots, \lambda_n$ be  real eigenvalues of the matrix $\bA\in \real^{n\times n}$. Then,
$
\sum_{i=1}^{n} \abs{\lambda_i}^2 \leq \sum_{i=1,j=1}^{n,n} \abs{a_{ij}}^2 = \norm{\bA}_F^2,
$
which means the sum of the squared absolute values of the eigenvalues is bounded by the Frobenius norm of the matrix.
\end{theorem}
\begin{proof}[of Theorem~\ref{theorem:schur_inequality}]
Suppose the Schur decomposition of $\bA$ is given by $\bA=\bQ\bU\bQ^\top$ (see, for example, \citet{lu2021numerical}), where the diagonal of $\bU$ contains the eigenvalues of $\bA$. By the orthogonal invariance Proposition~\ref{proposition:frobenius-orthogonal-equi}, we have
$
\norm{\bU}_F = \norm{\bQ\bU\bQ^\top}_F = \norm{\bA}_F. 
$
Therefore,
$$
\sum_{i=1}^{n} \abs{\lambda_i}^2 = \sum_{i=1}^{n} u_{ii}^2 \leq \sum_{i=1}^{n} u_{ii}^2 + \sum_{ i\neq j} u_{ij}^2 = \norm{\bU}_F^2
\quad\implies \quad
\sum_{i=1}^{n} \abs{\lambda_i}^2  \leq \norm{\bA}_F^2.
$$
This completes the proof.
\end{proof}
\begin{exercise}[Schur Inequality]
Prove the Schur inequality for general matrices that do not necessarily have real eigenvalues. Under what conditions does equality hold?
\end{exercise}

\subsection*{Spectral Norm}
Another important matrix norm that is extensively used is the \textit{spectral norm}.
\begin{definition}[Spectral Norm]\label{definition:spectral_norm_app}
The spectral norm of a matrix $\bA\in \real^{m\times n}$ is defined as 
$$
\normtwo{\bA} = \mathop{\max}_{\bx\neq\bzero} \frac{\normtwo{\bA\bx}}{\normtwo{\bx}}  =\mathop{\max}_{\bu\in \real^n: \norm{\bu}_2=1}  \norm{\bA\bu}_2 ,
$$
which corresponds to the largest singular value of $\bA$, i.e., $\normtwo{\bA} = \sigma_{\max}(\bA)$~\footnote{When $\bA$ is an $n\times n$  positive semidefinite matrix, $\normtwo{\bA}=\sqrt{\lambda_{\max}(\bA^2)}=\lambda_{\max}(\bA)$, i.e., the largest eigenvalue of $\bA$.}. The second equality holds because scaling $\bx$ by a nonzero scalar does not affect the ratio:
$
\frac{\norm{\lambda \cdot \bA\by}_2 }{\norm{\lambda \cdot \by}_2 }  =  \norm{\bA\by}_2.
$
The definition also indicates the matrix-vector inequality:
$$
\text{ $\normtwo{\bA\bx} \leq \normtwo{\bA}\normtwo{\bx}$ for all vectors $\bx\in\real^n$.}
$$
\end{definition}
To see why the spectral norm of a matrix is equal to its largest singular value, consider the singular value decomposition $\bA=\bU\bSigma\bV^\top$. We have:
$$
\begin{aligned}
\norm{\bU\bSigma\bV^\top}_2 &= \max_{\bx \neq \bzero} \frac{\norm{\bU\bSigma\bV^\top \bx}_2}{\normtwo{\bx}} \stackrel{*}{=} \max_{\bV^\top\bx \neq \bzero} \frac{\norm{\bSigma\bV^\top \bx}_2}{\normtwo{\bx}}\\
&= \max_{\bV^\top\bx \neq \bzero} \frac{\norm{\bSigma\bV^\top \bx}_2}{\norm{\bV^\top\bx}_2} \frac{\norm{\bV^\top\bx}_2}{\normtwo{\bx}} 
=\max_{\by \neq \bzero} \frac{\norm{\bSigma \by}_2}{\norm{\by}_2}
\leq \sigma_{\max}(\bA),
\end{aligned}
$$
where the equality ($*$) holds since $\bV$ is orthogonal and the $\ell_2$ vector norm is orthogonally invariant. The inequality holds because the largest singular value of $\bA$ is the maximum value of $\frac{\norm{\bSigma \by}_2}{\norm{\by}_2}$ over all nonzero vectors $\by$.
Alternatively, this can be shown by noting that: $\normtwo{\bA\bx}^2=\abs{\bx^\top\bA^\top\bA\bx}\leq \sigma_{\max}^2(\bA)$.

\index{Submultiplicativity}
\begin{proposition}[Submultiplicativity of Spectral]\label{propo:submul_spec}
The spectral norm is submultiplicative. That is, $\norm{\bA\bB}_2\leq \normtwo{\bA}\norm{\bB}_2$.
\end{proposition}
\begin{proof}[of Proposition~\ref{propo:submul_spec}]
The definition of the spectral norm shows that
$$
\begin{aligned}
\norm{\bA\bB}_2 &= \max_{\bx \neq \bzero} \frac{\norm{\bA\bB \bx}_2}{\normtwo{\bx}} \stackrel{*}{=} \max_{\bB\bx \neq \bzero} \frac{\norm{\bA\bB \bx}_2}{\normtwo{\bx}} 
= \max_{\bB\bx \neq \bzero} \frac{\norm{\bA\bB \bx}_2}{\norm{\bB\bx}_2} \frac{\norm{\bB\bx}_2}{\normtwo{\bx}} \\
&\leq \max_{\by \neq \bzero} \frac{\norm{\bA \by}_2}{\norm{\by}_2} \max_{\bx \neq \bzero} \frac{\norm{\bB\bx}_2}{\normtwo{\bx}}
=\normtwo{\bA}\norm{\bB}_2,
\end{aligned}
$$
where the equality ($*$) holds because if the maximum is obtained when $\bx\neq \bzero$ and $\bB\bx=\bzero$, the norm is $\norm{\bA\bB}=\bzero$. This holds only when $\bA\bB=\bzero$, and the submultiplicativity holds obviously.

Alternatively, we have 
$
\norm{\bA\bB}_2 = \max_{\norm{\bu}=1}\norm{\bA\bB \bu}_2,
$
where $\norm{\bA\bB \bu}_2\leq \normtwo{\bA}\norm{\bB\bu}_2\leq\normtwo{\bA}\norm{\bB}_2 $. This again indicates $\norm{\bA\bB}_2\leq \normtwo{\bA}\norm{\bB}_2 $ implicitly.
\end{proof}

\index{Orthogonally invariance}
\begin{proposition}[Orthogonally Invariance of Spectral]\label{proposition:submul_spec}
Let $\bA\in \real^{m\times n}$ be given, and let $\bU\in \real^{m\times m}$ and $\bV\in \real^{n\times n}$ be orthogonal matrices. Then,
$
\norm{\bA }_2 =  \norm{\bU\bA\bV }_2.
$
\end{proposition}
\begin{proof}[of Proposition~\ref{proposition:submul_spec}]
We notice that 
$$
\begin{aligned}
\normtwo{\bU\bA\bV}  &=\mathop{\max}_{\bx\neq\bzero} \frac{\norm{\bU\bA\bV\bx}_2}{\normtwo{\bx}}
=\mathop{\max}_{\bx\neq\bzero} \frac{\norm{\bA\bV\bx}_2}{\normtwo{\bx}}
=\mathop{\max}_{\bx\neq\bzero} \frac{\norm{\bA\bV\bx}_2}{\norm{\bV\bx}_2}
=\mathop{\max}_{\by\neq\bzero} \frac{\norm{\bA\by}_2}{\norm{\by}_2}
=\normtwo{\bA} ,
\end{aligned}
$$
This completes the proof.
\end{proof}

We conclude that the spectral norm is a normalized matrix norm.
\begin{proposition}[Normalized Spectral]\label{proposition:norma_spec}
The spectral norm is normalized such that:
\begin{enumerate}
\item \textit{Normalization}. When $\bA\neq \bzero$, we have $\norm{\frac{1}{\normtwo{\bA}} \bA}_2=1$.
\item \textit{Normalized}. $\norm{\bI}_2=1$.
\end{enumerate}
\end{proposition}

\section{Symmetry, Definiteness, and Quadratic Models}
In this section, we will introduce the fundamental concepts of symmetric matrices, positive definite matrices, and quadratic functions (models), which are essential in the modeling of many optimization problems.
\subsection{Symmetric Matrices}\label{section:symmetric_mats}
We introduce four important properties of symmetric matrices.
The following proposition states  that symmetric matrices have only real eigenvalues.
\begin{proposition}[Symmetric Property-I: Real Eigenvalues]\label{proposition:real-eigenvalues-spectral}
The eigenvalues of any symmetric matrix are all real. 
\end{proposition}
\begin{proof}[of Proposition~\ref{proposition:real-eigenvalues-spectral}]
Suppose  $\lambda$ is an eigenvalue of a symmetric matrix $\bA$, and let it be a complex number $\lambda=a+ib$, where $a,b$ are real. Its complex conjugate is $\bar{\lambda}=a-ib$. 
Similarly, let the corresponding complex eigenvector be $\bx = \bc+i\bd$ with its conjugate $\bar{\bx}=\bc-i\bd$, where $\bc, \bd$ are real vectors. Then, we have 
$$
\bA \bx = \lambda \bx\qquad   \underrightarrow{\text{ leads to }}\qquad  \bA \bar{\bx} = \bar{\lambda} \bar{\bx}\qquad   \underrightarrow{\text{ transpose to }}\qquad  \bar{\bx}^\top \bA =\bar{\lambda} \bar{\bx}^\top.
$$
We take the dot product of the first equation with $\bar{\bx}$ and the last equation with $\bx$:
$$
\bar{\bx}^\top \bA \bx = \lambda \bar{\bx}^\top \bx \qquad \text{and } \qquad \bar{\bx}^\top \bA \bx = \bar{\lambda}\bar{\bx}^\top \bx.
$$
Then we have the equality $\lambda\bar{\bx}^\top \bx = \bar{\lambda} \bar{\bx}^\top\bx$. Since $\bar{\bx}^\top\bx = (\bc-i\bd)^\top(\bc+i\bd) = \bc^\top\bc+\bd^\top\bd$ is a real number, we deduce that the imaginary part of $\lambda$ must be zero. Thus, $\lambda$ is real.
\end{proof}

\begin{proposition}[Symmetric Property-II: Orthogonal Eigenvectors]\label{proposition:orthogonal-eigenvectors}
The eigenvectors  corresponding to distinct eigenvalues of a symmetric matrix are orthogonal. 
Consequently, these eigenvectors can be normalized to form an orthonormal basis, since 
$$
\bA\bx = \lambda \bx 
\qquad \implies \qquad 
\bA\frac{\bx}{\normtwo{\bx}} = \lambda \frac{\bx}{\normtwo{\bx}},
$$ 
which preserves the eigenvalue.
\end{proposition}
\begin{proof}[of Proposition~\ref{proposition:orthogonal-eigenvectors}]
Let $\lambda_1, \lambda_2$ be distinct eigenvalues of $\bA$, with corresponding eigenvectors $\bx_1, \bx_2$, satisfying $\bA\bx_1=\lambda \bx_1$ and $\bA\bx_2 = \lambda_2\bx_2$. We have the following equalities:
$$
\bA\bx_1=\lambda_1 \bx_1 
\qquad \implies \qquad 
\bx_1^\top \bA =\lambda_1 \bx_1^\top 
\qquad \implies \qquad 
\bx_1^\top \bA \bx_2 =\lambda_1 \bx_1^\top\bx_2,
$$
and 
$$
\bA\bx_2 = \lambda_2\bx_2 
\qquad \implies \qquad 
\bx_1^\top\bA\bx_2 = \lambda_2\bx_1^\top\bx_2,
$$
implying $\lambda_1 \bx_1^\top\bx_2=\lambda_2\bx_1^\top\bx_2$. Since  $\lambda_1\neq \lambda_2$, the eigenvectors are orthogonal.
\end{proof}

\begin{proposition}[Symmetric Property-III: Orthonormal Eigenvectors for Duplicate Eigenvalue]\label{proposition:eigen-multiplicity}
Let $\bA\in\real^{n\times n}$ be symmetric.
If $\bA$ has a duplicate eigenvalue $\lambda_i$ with multiplicity  $k\geq 2$, then there exist $k$ orthonormal eigenvectors corresponding to $\lambda_i$.
\end{proposition}
\begin{proof}[of Proposition~\ref{proposition:eigen-multiplicity}]
We start by noting that there exists at least one unit-length eigenvector $\bx_{i1}$ corresponding to $\lambda_i$. 
Furthermore, for such an eigenvector $\bx_{i1}$, we can consistently find  $n-1$ additional orthonormal vectors $\by_2, \by_3, \ldots, \by_n$ such that $\{\bx_{i1}, \by_2, \by_3, \ldots, \by_n\}$ constitutes an orthonormal basis of $\real^n$. 
Define the matrices $\bY_1$ and $\bP_1$ as follows:
$$
\bY_1\triangleq[\by_2, \by_3, \ldots, \by_n] \qquad \text{and} \qquad \bP_1\triangleq[\bx_{i1}, \bY_1].
$$
Since $\bA$ is symmetric, we then have
$
\bP_1^\top\bA\bP_1 = \footnotesize\begin{bmatrix}
\lambda_i &\bzero \\
\bzero & \bY_1^\top \bA\bY_1
\end{bmatrix}.
$
Since $\bP_1$ is nonsingular and orthogonal, $\bA$ and $\bP_1^\top\bA\bP_1$ are similar matrices such that they have the same eigenvalues  (see, for example, \citet{lu2021numerical}).
Using the determinant of block matrices~\footnote{If matrix $\bM$ has a block formulation: $\bM=\scriptsize\begin{bmatrix}
		\bA & \bB \\
		\bC & \bD 
	\end{bmatrix}$, then $\det(\bM) = \det(\bA)\det(\bD-\bC\bA^{-1}\bB)$.}, we get: 
$$
\det(\bP_1^\top\bA\bP_1 - \lambda\bI_n) =
(\lambda_i - \lambda )\det(\bY_1^\top \bA\bY_1 - \lambda\bI_{n-1}).
$$
If $\lambda_i$ has a multiplicity of $k\geq 2$, then the term $(\lambda_i-\lambda)$ appears $k$ times in the characteristic  polynomial resulting from the determinant $\det(\bP_1^\top\bA\bP_1 - \lambda\bI_n)$, i.e., this term appears  $k-1$ times in the characteristic  polynomial from $\det(\bY_1^\top \bA\bY_1 - \lambda\bI_{n-1})$. In other words, $\det(\bY_1^\top \bA\bY_1 - \lambda_i\bI_{n-1})=0$, and $\lambda_i$ is an eigenvalue of $\bY_1^\top \bA\bY_1$ with multiplicity $k-1$.

Let $\bB\triangleq\bY_1^\top \bA\bY_1$. Since $\det(\bB-\lambda_i\bI_{n-1})=0$, the null space of $\bB-\lambda_i\bI_{n-1}$ is nonempty. Suppose $(\bB-\lambda_i\bI_{n-1})\bn = \bzero$, i.e., $\bB\bn=\lambda_i\bn$, where $\bn$ is an eigenvector of $\bB$.

From $
\bP_1^\top\bA\bP_1 = 
\footnotesize
\begin{bmatrix}
	\lambda_i &\bzero \\
	\bzero & \bB
\end{bmatrix},
$
we have $
\bA\bP_1 
\footnotesize
\begin{bmatrix}
	z \\
	\bn 
\end{bmatrix} 
= 
\bP_1
\footnotesize
\begin{bmatrix}
	\lambda_i &\bzero \\
	\bzero & \bB
\end{bmatrix}
\begin{bmatrix}
	z \\
	\bn 
\end{bmatrix}$, where $z$ is any scalar. 
From the left side of this equation:
\begin{equation}\label{equation:spectral-pro4-right}
	\begin{aligned}
		\bA\bP_1 
		\begin{bmatrix}
			z \\
			\bn 
		\end{bmatrix} 
		&=
		\begin{bmatrix}
			\lambda_i\bx_{i1}, \bA\bY_1
		\end{bmatrix}
		\begin{bmatrix}
			z \\
			\bn 
		\end{bmatrix} 
		=\lambda_iz\bx_{i1} + \bA\bY_1\bn.
	\end{aligned}
\end{equation}
From the right side of the equation:
\begin{equation}\label{equation:spectral-pro4-left}
	\begin{aligned}
		\bP_1
		\begin{bmatrix}
			\lambda_i &\bzero \\
			\bzero & \bB
		\end{bmatrix}
		\begin{bmatrix}
			z \\
			\bn 
		\end{bmatrix}
		&=
		\begin{bmatrix}
			\bx_{i1} & \bY_1
		\end{bmatrix}
		\begin{bmatrix}
			\lambda_i &\bzero \\
			\bzero & \bB
		\end{bmatrix}
		\begin{bmatrix}
			z \\
			\bn 
		\end{bmatrix}
		=
		\begin{bmatrix}
			\lambda_i\bx_{i1} & \bY_1\bB 
		\end{bmatrix}
		\begin{bmatrix}
			z \\
			\bn 
		\end{bmatrix}\\
		&= \lambda_i z \bx_{i1} + \bY_1\bB \bn 
		=\lambda_i z \bx_{i1} + \lambda_i \bY_1 \bn,
	\end{aligned}
\end{equation}
where the last equality is due to  $\bB \bn=\lambda_i\bn$.
Combining \eqref{equation:spectral-pro4-left} and \eqref{equation:spectral-pro4-right}, we obtain 
$$
\bA\bY_1\bn = \lambda_i\bY_1 \bn,
$$
which means $\bY_1\bn$ is an eigenvector of $\bA$ corresponding to the eigenvalue $\lambda_i$ (the same eigenvalue corresponding to $\bx_{i1}$). Since $\bY_1\bn$ is a linear combination of $\by_2, \by_3, \ldots, \by_n$, which are orthonormal to $\bx_{i1}$, it can be chosen to be orthonormal to $\bx_{i1}$.

To conclude, if there exists an  eigenvector, $\bx_{i1}$, corresponding to the eigenvalue  $\lambda_i$ with a multiplicity  $k\geq 2$, we can construct a second eigenvector by choosing a vector from the null space of $(\bB-\lambda_i\bI_{n-1})$, as constructed above. 

Now, suppose  we have constructed the second eigenvector $\bx_{i2}$, which is orthonormal to $\bx_{i1}$.  
For such eigenvectors $\bx_{i1}$ and $\bx_{i2}$, we can always find  $n-2$ additional orthonormal vectors $\by_3, \by_4, \ldots, \by_n$ such  that $\{\bx_{i1},\bx_{i2}, \by_3, \by_4, \ldots, \by_n\}$ forms an orthonormal basis for $\real^n$. 
Place these vectors  $\by_3, \by_4, \ldots, \by_n$ into matrix $\bY_2$ and $\{\bx_{i1},\bx_{i2},  \by_3, \by_4, \ldots, \by_n\}$ into matrix $\bP_2$:
$$
\bY_2\triangleq[\by_3, \by_4, \ldots, \by_n] \qquad \text{and} \qquad \bP_2\triangleq[\bx_{i1}, \bx_{i2},\bY_1].
$$
Since $\bA$ is symmetric, we then have
$$
\bP_2^\top\bA\bP_2 = 
\small
\begin{bmatrix}
	\lambda_i & 0 &\bzero \\
	0& \lambda_i &\bzero \\
	\bzero & \bzero & \bY_2^\top \bA\bY_2
\end{bmatrix}
\triangleq
\begin{bmatrix}
	\lambda_i & 0 &\bzero \\
	0& \lambda_i &\bzero \\
	\bzero & \bzero & \bC
\end{bmatrix},
$$
where $\bC\triangleq\bY_2^\top \bA\bY_2$ such that $\det(\bP_2^\top\bA\bP_2 - \lambda\bI_n) = (\lambda_i-\lambda)^2 \det(\bC - \lambda\bI_{n-2})$. If the multiplicity of $\lambda_i$ is $k\geq 3$, then $\det(\bC - \lambda_i\bI_{n-2})=0$, and the null space of $\bC - \lambda_i\bI_{n-2}$ is not empty. Thus, we can still find a vector $\bn$ from the null space of $\bC - \lambda_i\bI_{n-2}$ such that  $\bC\bn = \lambda_i \bn$. Now we can construct a vector $\footnotesize\begin{bmatrix}
	z_1 \\
	z_2\\
	\bn
\end{bmatrix}\in \real^n $, where $z_1$ and $ z_2$ are any scalar values, such that 
$$
\bA\bP_2\small\begin{bmatrix}
	z_1 \\
	z_2\\
	\bn
\end{bmatrix} = \normalsize\bP_2 
\small
\begin{bmatrix}
	\lambda_i & 0 &\bzero \\
	0& \lambda_i &\bzero \\
	\bzero & \bzero & \bC
\end{bmatrix}
\begin{bmatrix}
	z_1 \\
	z_2\\
	\bn
\end{bmatrix}.
$$
Similarly, from the left side of the above equation, we will get $\lambda_iz_1\bx_{i1} +\lambda_iz_2\bx_{i2}+\bA\bY_2\bn$. From the right side of the above equation, we will get $\lambda_iz_1\bx_{i1} +\lambda_i z_2\bx_{i2}+\lambda_i\bY_2\bn$. As a result, 
$$
\bA\bY_2\bn = \lambda_i\bY_2\bn,
$$
where $\bY_2\bn$ is an eigenvector of $\bA$, orthogonal to both $\bx_{i1}$ and $\bx_{i2}$. This eigenvector can also be normalized to ensure orthonormality with the first two eigenvectors.

The process can continue, ultimately yielding a set of $k$ orthonormal eigenvectors corresponding to $\lambda_i$.

In fact, the dimension of the null space of $\bP_1^\top\bA\bP_1 -\lambda_i\bI_n$ is equal to the multiplicity $k$. It also follows that if the multiplicity of $\lambda_i$ is $k$, there cannot be more than $k$ orthogonal eigenvectors corresponding to $\lambda_i$. 
If there were more than $k$, it would lead to the conclusion that there are more than $n$  orthogonal eigenvectors in $\real^n$, which is a contradiction.
\end{proof}

For any matrix multiplication, the rank of the resulting matrix is at most the rank of the input matrices.
\begin{lemma}[Rank of $\bA\bB$]\label{lemma:rankAB}
Let $\bA\in \real^{m\times n}$ and  $\bB\in \real^{n\times k}$ be any matrices. Then, the rank of the product  $\bA\bB\in \real^{m\times k}$ satisfies $\rank$($\bA\bB$)$\leq$min($\rank$($\bA$), $\rank$($\bB$)).
\end{lemma}
\begin{proof}[of Lemma~\ref{lemma:rankAB}]
For the matrix product $\bA\bB$, we observe the following: 
\begin{itemize}
\item Every row of $\bA\bB$ is a linear combination of the rows of $\bB$. Hence, the row space of $\bA\bB$ is contained within the row space of $\bB$, implying that $\rank$($\bA\bB$)$\leq$$\rank$($\bB$).

\item Similarly, every column of $\bA\bB$ is a linear combination of the columns of $\bA$. Therefore, the column space of $\bA\bB$ is a subspace of the column space of $\bA$, which gives $\rank$($\bA\bB$)$\leq$$\rank$($\bA$).
\end{itemize}
Combining these two results, we conclude that $\rank$($\bA\bB$)$\leq$min($\rank$($\bA$), $\rank$($\bB$)).
\end{proof}
This lemma establishes that the rank of a symmetric matrix is equal to the number of its nonzero eigenvalues.
\begin{proposition}[Symmetric Property-IV: Rank of Symmetric Matrices]\label{proposition:rank-of-symmetric}
Let $\bA$ be an $n\times n$ real symmetric matrix. Then, 
$$
\text{rank($\bA$) =
	the total number of nonzero eigenvalues of $\bA$. }
$$
In particular, $\bA$ has full rank if and only if it is nonsingular. Moreover, the column space of $\bA$, denoted $\cspace(\bA)$, is spanned by the eigenvectors corresponding to its nonzero eigenvalues. 
\end{proposition}
\begin{proof}[of Proposition~\ref{proposition:rank-of-symmetric}]
For any symmetric matrix $\bA$, we have $\bA$, in spectral form, as $\bA = \bQ \bLambda\bQ^\top$ and also $\bLambda = \bQ^\top\bA\bQ$; see \eqref{equation:spec_decom}, which is a result of Symmetric Properties-I$\sim$III (Propositions~\ref{proposition:real-eigenvalues-spectral}$\sim$\ref{proposition:eigen-multiplicity}). 
By Lemma~\ref{lemma:rankAB}, we know that for any matrix multiplication, $\rank$($\bA\bB$)$\leq$min($\rank$($\bA$), $\rank$($\bB$)).
Applying this result:
\begin{itemize}
\item From $\bA = \bQ \bLambda\bQ^\top$, we have $\rank(\bA) \leq \rank(\bQ \bLambda) \leq \rank(\bLambda)$.
	
\item From $\bLambda = \bQ^\top\bA\bQ$, we have $\rank(\bLambda) \leq \rank(\bQ^\top\bA) \leq \rank(\bA)$.
\end{itemize}
Since both inequalities hold in opposite directions, we conclude that $\rank(\bA) = \rank(\bLambda)$, which is the total number of nonzero eigenvalues.

Furthermore, $\bA$ is nonsingular if and only if all  its eigenvalues are nonzero, which establishes that $\bA$ has full rank if and only if it is nonsingular.
\end{proof}

\subsection{Positive Definiteness}
A symmetric matrix can be  categorized into positive definite, positive semidefinite, negative definite, negative semidefinite, and indefinite types as follows.

\begin{definition}[Positive Definite and Positive Semidefinite\index{Positive definite}\index{Positive semidefinite}]\label{definition:psd-pd-defini}
A symmetric matrix $\bA\in \real^{n\times n}$ ($\bA\in\symmetric^n$) is considered \textit{positive definite (PD)} if $\bx^\top\bA\bx>0$ for all nonzero $\bx\in \real^n$, denoted by $\bA\succ \bzero$ or $\bA\in\pd^n$.
Similarly, a symmetric matrix $\bA\in \real^{n\times n}$ is called \textit{positive semidefinite (PSD)} if $\bx^\top\bA\bx \geq 0$ for all $\bx\in \real^n$,  denoted by $\bA\succeq\bzero$ or $\bA\in\psd^n$. 
\footnote{A symmetric matrix $\bA\in\real^{n\times n}$ is called \textit{negative definite (ND) } if $\bx^\top\bA\bx<0$ for all nonzero $\bx\in\real^n$, denoted by $\bA\prec\bzero$ or $\bA\in\nd^n$; 
a symmetric matrix $\bA\in\real^{n\times n}$ is called \textit{negative semidefinite (NSD)}  if $\bx^\top\bA\bx\leq 0$ for all $\bx\in\real^n$, denoted by $\bA\preceq\bzero$ or $\bA\in\nsd^n$;
and a symmetric matrix $\bA\in\real^{n\times n}$ is called \textit{indefinite (ID)}  if there exist vectors $\bx$ and $\by\in\real^n$ such that $\bx^\top\bA\bx<0$ and $\by^\top\bA\by>0$.
}
\end{definition}

Given a negative definite matrix $\bA$, then $-\bA$ is a positive definite matrix; if $\bA$ is negative semidefintie matrix, then $-\bA$ is positive semidefinite.
The following property demonstrates that PD and PSD matrices have special forms of eigenvalues.
\begin{theorem}[PD Property-I: Eigenvalue Characterization Theorem]\label{theorem:eigen_charac}
A symmetric matrix $\bA$ is positive definite if and only if all its eigenvalues are positive. Similarly, a matrix $\bA$ is positive semidefinite if and only if all its eigenvalues are nonnegative.
Additionally,  a matrix is indefinite if and only if it possesses at least one positive eigenvalue and at least one negative eigenvalue.
Furthermore, we have the following implications:
\begin{itemize}
\item $\bA-\mu\bI\succeq \bzero$ if and only if $\lambda_{\min}(\bA) \geq \mu$;
\item  $\bA-\mu\bI\succ \bzero$ if and only if $\lambda_{\min}(\bA) > \mu$;
\item $\bA-\mu\bI\preceq \bzero$ if and only if $\lambda_{\max}(\bA) \leq \mu$;
\item $\bA-\mu\bI\prec \bzero$ if and only if $\lambda_{\max}(\bA) < \mu$;
\item $\lambda_{\min}(\bA)\bI\preceq \bA \preceq \lambda_{\max}(\bA)\bI$,
\end{itemize}
\noindent where $\lambda_{\min}(\bA)$ and $\lambda_{\max}(\bA)$ denote the minimum and maximum eigenvalues of $\bA$, respectively. 
The notation $\bB \prec \bC$ means that $\bC-\bB$ is PSD.
\end{theorem}
Given the eigenpair $(\lambda, \bx)$ of $\bA$, the forward implication can be shown that $\bx^\top\bA\bx=\lambda\bx^\top\bx>0$ such that $ \lambda=\frac{\bx^\top\bA\bx}{\bx^\top\bx}>0$ (resp. $\geq 0$) if $\bA$ is PD (resp. PSD).
The complete proof of this equivalence can be derived using the spectral theorem (Theorem~\ref{theorem:spectral_theorem}); the details can be found in \citet{lu2021numerical}.
This theorem provides an alternative definition of positive definiteness and positive semidefiniteness in terms of the eigenvalues of the matrix, which is fundamental for the Cholesky decomposition (Theorems~\ref{theorem:cholesky-factor-exist}).

Although not all components of a positive definite matrix are necessarily positive, the diagonal components of such a matrix are guaranteed to be positive, as stated in the following result.
\begin{theorem}[PD Property-II: Positive Diagonals of PD Matrices]\label{theorem:positive-in-pd}
The diagonal elements of a positive definite matrix $\bA\in\real^{n\times n}$ are all \textit{positive}. Similarly, the diagonal elements of a positive semidefinite matrix $\bB\in\real^{n\times n}$ are all \textit{nonnegative}.
Additionally,  the diagonal elements of an indefinite matrix $\bC\in\real^{n\times n}$ contains at least one positive diagonal and at least one negative diagonal.
\end{theorem}
\begin{proof}[of Theorem~\ref{theorem:positive-in-pd}]
By the definition of positive definite matrices, we have $\bx^\top\bA \bx >0$ for all nonzero vectors $\bx$. Consider the specific case where $\bx=\be_i$, with $\be_i$ being the $i$-th unit basis vector having the $i$-th entry equal to 1 and all other entries equal to 0. Then,
$$
\be_i^\top\bA \be_i = a_{ii}>0, \gapforall \forall i \in \{1, 2, \ldots, n\},
$$	
where $a_{ii}$ is the $i$-th diagonal component. The proofs for the second and the third parts follow a similar argument. This completes the proof.
\end{proof}

\subsection{Quadratic Functions and Quadratic Models}
\begin{figure}[h]
\centering  
\vspace{-0.35cm} 
\subfigtopskip=2pt 
\subfigbottomskip=2pt 
\subfigcapskip=-5pt 
\subfigure[Positive definite matrix: $\bA =\scriptsize \begin{bmatrix}
200 & 0 \\ 0 & 200
\end{bmatrix}$.]{\label{fig:quadratic_PD1}
\includegraphics[width=0.485\linewidth]{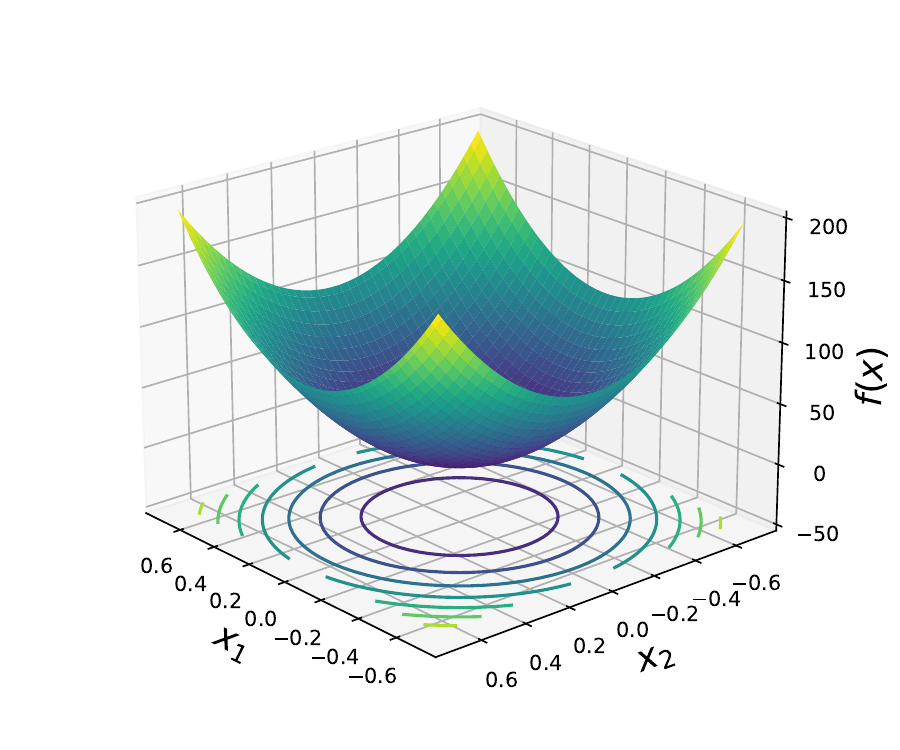}}
\subfigure[Negative definite matrix: $\bA = \scriptsize\begin{bmatrix}
-200 & 0 \\ 0 & -200
\end{bmatrix}$.]{\label{fig:quadratic_ND1}
\includegraphics[width=0.485\linewidth]{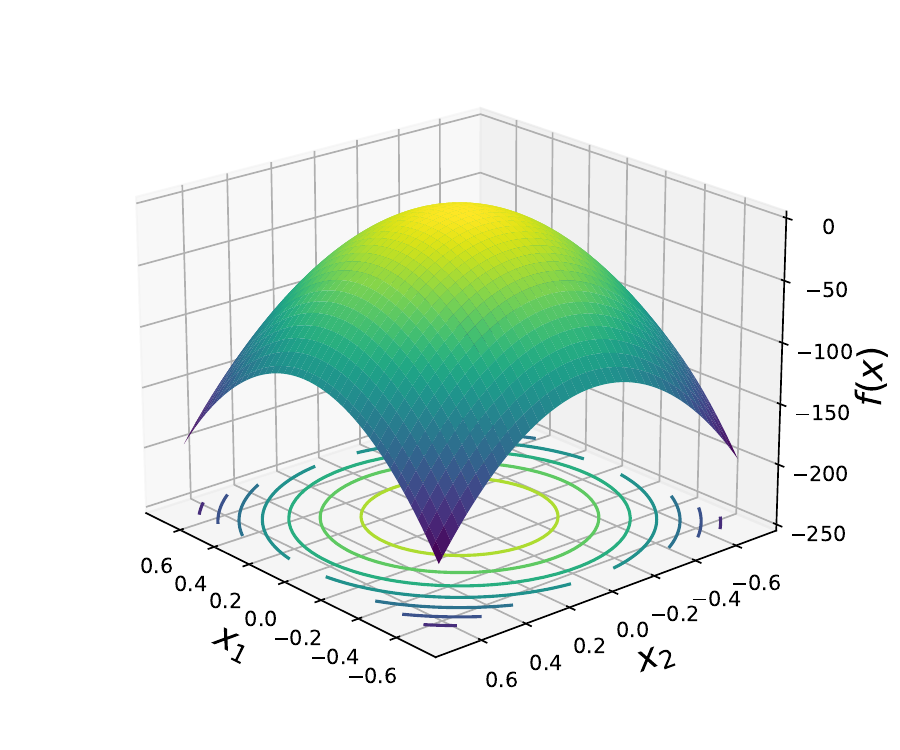}}
\subfigure[Semidefinite matrix: $\bA = \scriptsize\begin{bmatrix}
200 & 0 \\ 0 & 0
\end{bmatrix}$. A line runs through the bottom of the valley is the set of solutions.]{\label{fig:quadratic_singular1}
\includegraphics[width=0.485\linewidth]{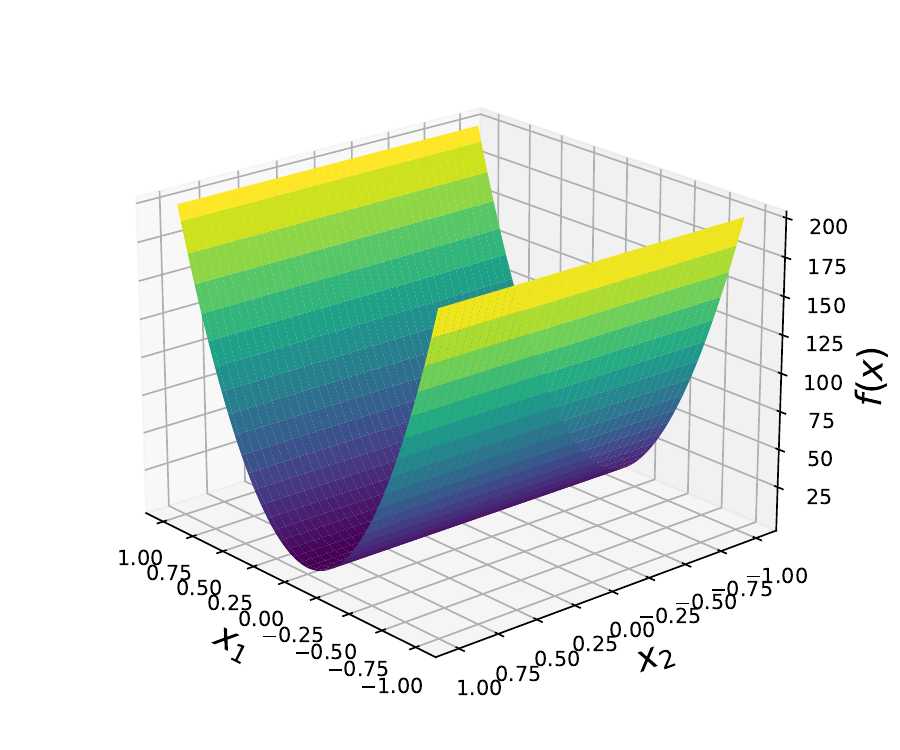}}
\subfigure[Indefinte matrix: $\bA = \scriptsize\begin{bmatrix}
200 & 0 \\ 0 & -200
\end{bmatrix}$.]{\label{fig:quadratic_saddle1}
\includegraphics[width=0.485\linewidth]{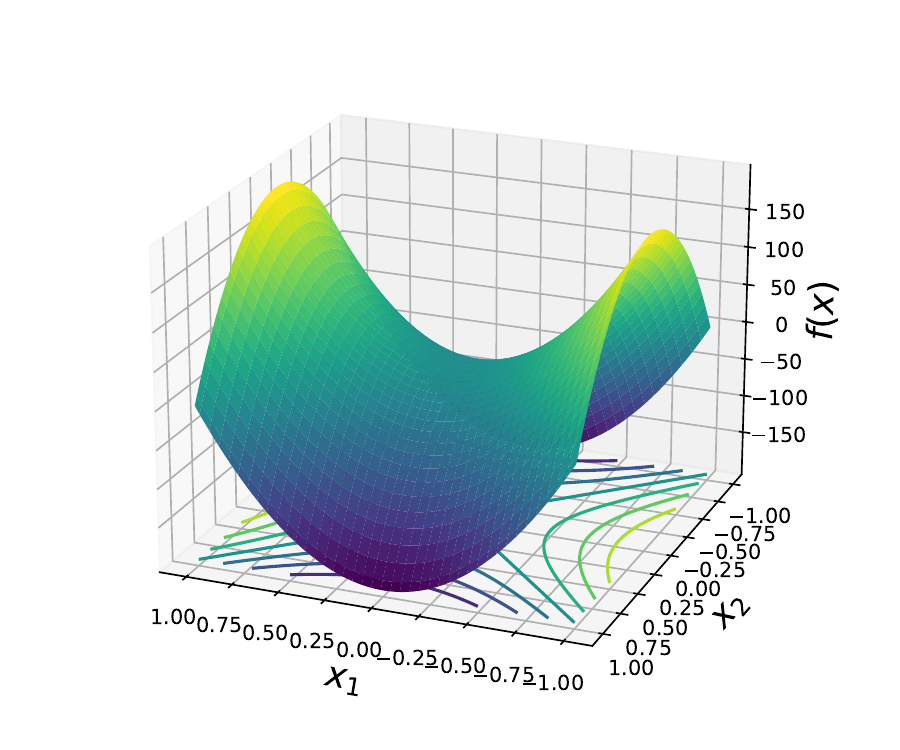}}
\caption{Loss surfaces for different quadratic forms, providing the surface plots and contour plots (\textcolor{mylightbluetext}{blue}=low,
\textcolor{mydarkyellow}{yellow}=high), where the upper graphs are the surface plots, and the lower ones are their projection (i.e., contours).}
\label{fig:different_quadratics1}
\end{figure}

We  further  discuss  linear systems with different types of matrices, the  quadratic form:
\begin{equation}\label{equation:quadratic-form-general-form1}
f(\bx) = \frac{1}{2} \bx^\top \bA \bx - \bb^\top \bx + c, \gap \bx\in \real^n,
\end{equation}
where $\bA\in \real^{n\times n}$, $\bb \in \real^n$, and $c$ is a scalar constant. Though the quadratic form in Equation~\eqref{equation:quadratic-form-general-form1} is an extremely simple model, it is rich enough to approximate many other functions, e.g., the Fisher information matrix \citep{amari1998natural}, and it captures key features of pathological curvature. The gradient of $f(\bx)$ at point $\bx$ is given by 
$$
\nabla f(\bx) = \frac{1}{2} (\bA^\top +\bA) \bx - \bb.
$$
Any optimum point of the function is the solution to the linear system $\frac{1}{2} (\bA^\top +\bA) \bx=  \bb $:
$$
\bx^* = 2(\bA^\top +\bA)^{-1}\bb.
$$
If $\bA$ is symmetric, the equation reduces to 
$
\nabla f(\bx) = \bA \bx - \bb.
$
Then the optimum point of the function is the solution of the linear system $\bA\bx=\bb$, where $\bA$ and $\bb$ are known matrix or vector, and $\bx$ is an unknown vector; and the optimum point of $\bx$ is thus given by 
$$
\bx^*  = \bA^{-1}\bb, \quad\text{if $\bA$ is nonsingular.}
$$
For different types of matrix $\bA$, the loss surface of $f(\bx)$ will vary as shown in Figure~\ref{fig:different_quadratics1}. When $\bA$ is positive definite, the surface is a \textit{convex bowl}; when $\bA$ is negative definite, on the contrary, the surface is a \textit{concave bowl}. $\bA$ also could be singular, in which case $\bA\bx-\bb=\bzero$ has more than one solution, and the set of solutions is a line (in the two-dimensional case) or a hyperplane (in the high-dimensional case).
This situation is similar to the case of a semidefinite quadratic form, as shown in Figure~\ref{fig:quadratic_singular1}.
Moreover, $\bA$ could be none of the above, then there exists a \textit{saddle point} (see Figure~\ref{fig:quadratic_saddle1}), where the gradient descent may fail (see discussions in the following chapters). In such senses, other methods, e.g., perturbed gradient descent \citep{jin2017escape}, can be applied to escape  saddle points. 

\paragraph{Diagonally dominant matrices.}
A specific form of diagonally dominant matrices constitutes a significant subset of positive semidefinite matrices.
\begin{definition}[Diagonally Dominant Matrices]
Given a symmetric matrix $\bA\in\real^{n\times n}$,  $\bA$ is called diagonally dominant if 
$$
\abs{a_{ii}} \geq \sum_{j\neq i} \abs{a_{ij}}, \qquad \forall i\in\{1,2,\ldots, n\};
$$
and $\bA$ is called strictly diagonally dominant if 
$$
\abs{a_{ii}} > \sum_{j\neq i} \abs{a_{ij}}, \qquad \forall i\in\{1,2,\ldots, n\}.
$$
\end{definition}

We  now show that \textit{diagonally dominant matrices} with nonnegative diagonal elements are positive semidefinite and that \textit{strictly diagonally dominant matrices} with positive diagonal elements are positive definite. 

\begin{theorem}[Positive Definiteness of Diagonally Dominant Matrices]\label{theorem:pd_diag_domi}
Given a symmetric matrix  $\bA\in\real^{n\times n}$, 
\begin{enumerate}[(i)]
\item If $\bA$ is diagonally dominant with nonnegative diagonals, then $\bA$ is positive semidefinite; 
\item If $\bA$ is strictly diagonally dominant with positive diagonals, then $\bA$ is positive definite.
\end{enumerate}

\end{theorem}
\begin{proof}[of Theorem~\ref{theorem:pd_diag_domi}]
\textbf{(i).} Suppose $\bA$ is not PSD and has a negative eigenvalue $\lambda$ associated with an eigenvector $\bv$ such that $\bA\bv=\lambda\bv$. 
Let $v_i$ be the element of $\bv$ with largest magnitude.
Consider the $i$-th element of $(\bA-\lambda\bI)\bv=\bzero$, we have
$$
\abs{a_{ii} - \lambda}\cdot \abs{v_i}
=
\abs{\sum_{j\neq i} a_{ij} v_j}
\leq 
\left( \sum_{ j\neq i} \abs{a_{ij}} \right) \abs{v_i}
\leq \abs{a_{ij}} \cdot \abs{v_i}.
$$
This implies $\abs{a_{ii} - \lambda} \leq \abs{a_{ij}}$ and leads to a contradiction.

\paragraph{(ii).}  From part (i), we known that $\bA$ is positive semidefinite. Suppose $\bA$ is not positive definite and has a zero eigenvalue $0$ associated with a nonzero eigenvector $\bv$ such that $\bA\bv=\bzero$.   Similarly, we have 
$$
\abs{a_{ii}}\cdot \abs{v_i}
=
\abs{\sum_{j\neq i} a_{ij} v_j}
\leq 
\left( \sum_{ j\neq i} \abs{a_{ij}} \right) \abs{v_i}
< \abs{a_{ij}} \cdot \abs{v_i},
$$
which is impossible and the result follows.
\end{proof}

\begin{exercise}[Quadratic Form]\label{exercise:quad_form_psd}
Consider the quadratic form $f(\bx) = \frac{1}{2} \bx^\top \bA \bx - \bb^\top \bx + c$, where $\bA\in\real^{n\times n}$ is symmetric, $ \bx\in \real^n$, and $c\in\real$. Show that the following two claims are equivalent~\footnote{This result finds its application in non-convex quadratic constrained quadratic problems (QCQPs) \citep{beck2014introduction}.}:
\begin{itemize}
\item $f(\bx) \geq 0$ for all $\bx\in\real^n$.
\item $\scriptsize\begin{bmatrix}
\bA & -\bb\\
-\bb^\top & 2c
\end{bmatrix}\succeq \bzero$.
\end{itemize} 
\textit{Hint: Apply the eigenvalue characterization theorem on $\bA$.}
\end{exercise}

\section{Differentiability and Differential Calculus}\label{section:differ_calc}
\index{Continuously differentiability}
\index{Second-order partial derivative}
\begin{definition}[Directional Derivative, Partial Derivative]\label{definition:partial_deri}
Given a function $f$ defined over a set $\sS\subseteq \real^n$ and a nonzero vector $\bd\in\real^n$, the \textit{directional derivative} of $f$ at $\bx$ with respect to the direction $\bd$ is given by, if the limit exists, 
$$
\mathop{\lim}_{\mu\rightarrow 0}
\frac{f(\bx+\mu\bd) - f(\bx)}{\mu}.
$$
And it is denoted by $f^\prime(\bx; \bd)$ or $D_{\bd}f(\bx)$. 
The directional derivative is sometimes referred to as the  \textit{G\^ateaux derivative}.

For any $i\in\{1,2,\ldots,n\}$, the directional derivative at $\bx$ with respect to the direction of the $i$-th standard basis vector  $\be_i$ (if it exists) is called the $i$-th \textit{partial derivative} and is denoted by $\frac{\partial f}{\partial x_i} (\bx)$, $D_{\be_i}f(\bx)$, or $\partial_i f(\bx)$.
\end{definition}

\begin{figure*}[h]
\centering  
\subfigtopskip=2pt 
\subfigbottomskip=9pt 
\subfigcapskip=-5pt 
\includegraphics[width=0.75\textwidth]{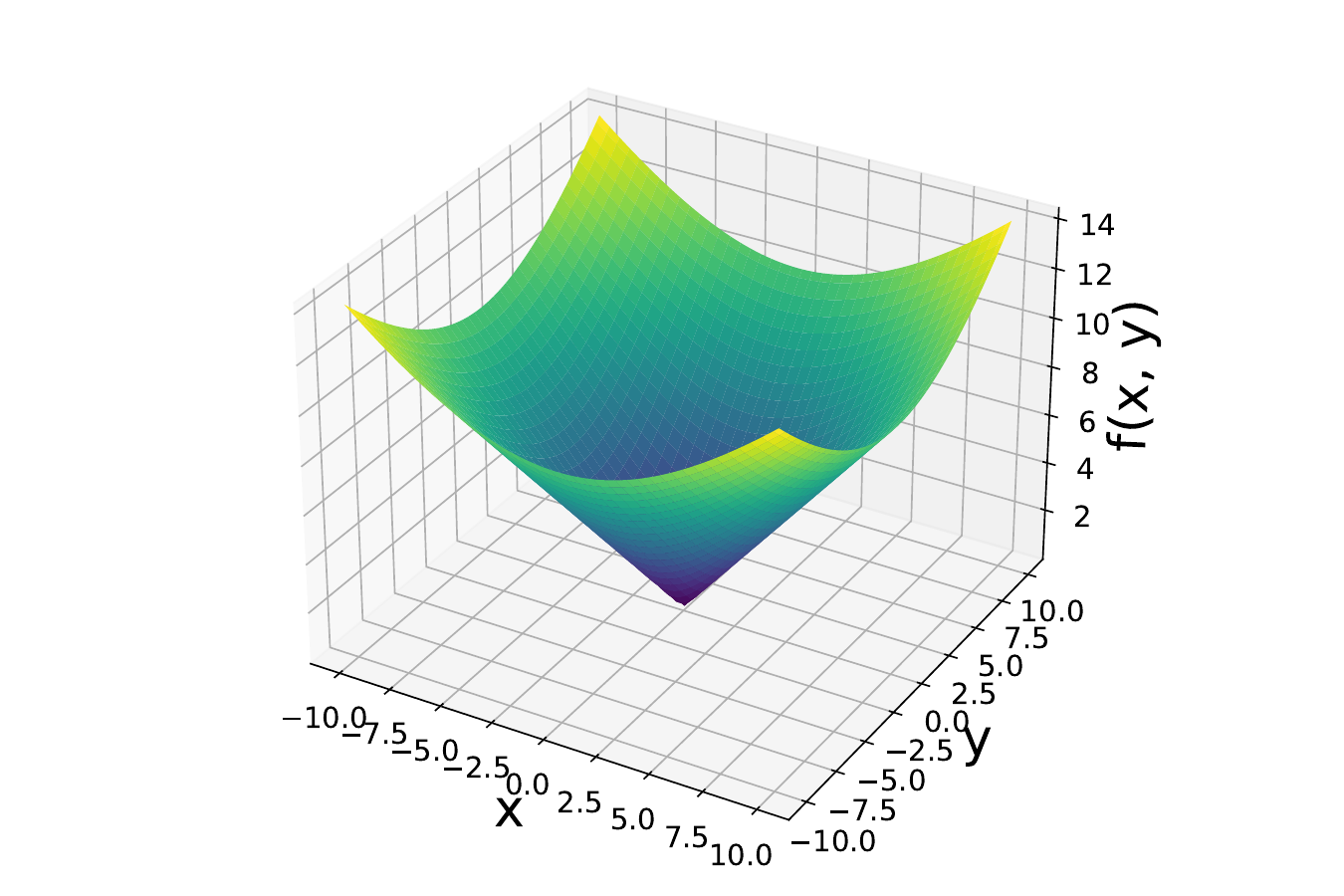}
\caption{Plot for the function $f(x, y) = \sqrt{x^2+y^2}$, in which case any directional derivative for the direction $\bd=[a,b]^\top$ with $a\neq 0$ and $b\neq 0$ at point $[0,0]^\top$ exists. However, the partial derivatives at this point do not exist.}
\label{fig:direc_deri_no_partial}
\end{figure*}

It's important to note that even if a function can have a directional derivative in every direction $\bd$ at some points, its partial derivatives may not exist.
For example, given the direction $\bd=[a,b]^\top$ with $a\neq 0$ and $b\neq 0$, the directional derivation of the function $f(x,y)=\sqrt{x^2+y^2}$ (see Figure~\ref{fig:direc_deri_no_partial}) at $[0,0]^\top$ can be obtained by 
$$
f^\prime(0,0; \bd) = \mathop{\lim}_{\mu\rightarrow 0}
\frac{\sqrt{(0+\mu a)^2+(0+\mu b)^2} -\sqrt{0^2+0^2} }{\mu} = \sqrt{a^2+b^2}.
$$
When $\bd$ is a unit vector, the directional derivative is 1.
However, it can be shown that the partial derivatives are 
$$
\begin{aligned}
\frac{\partial f}{\partial x}(0,0) &= \mathop{\lim}_{h\rightarrow 0} \frac{\sqrt{(0+h)^2+0^2}}{h} = \mathop{\lim}_{h\rightarrow 0} \frac{\abs{h}}{h},\\
\frac{\partial f}{\partial y}(0,0) &= \mathop{\lim}_{h\rightarrow 0} \frac{\sqrt{0^2+(0+h)^2}}{h} = \mathop{\lim}_{h\rightarrow 0} \frac{\abs{h}}{h}.
\end{aligned}
$$
When $h>0$, the partial derivatives are 1; when $h<0$, the partial derivaties are $-1$. Therefore, the partial derivatives do not exist.

If all the partial derivatives of a function $f$ exist at a point $\bx\in\real^n$, then the \textit{gradient} of $f$ at $\bx$ is defined as the column vector containing all the partial derivatives:
$$
\nabla f(\bx)=
\begin{bmatrix}
\frac{\partial f}{\partial x_1} (\bx),
\frac{\partial f}{\partial x_2} (\bx),
\ldots,
\frac{\partial f}{\partial x_n} (\bx)	
\end{bmatrix}^\top
\in \real^n.
$$

\begin{exercise}
Show that a function can have partial derivatives but is not necessarily continuous (Definition~\ref{definition:conti_funs}). 
Additionally, demonstrate that a function can be continuous but does not necessarily have partial derivatives.
\end{exercise}

A function $f$ defined over an open set $\sS\subseteq \real^n$ (Definition~\ref{definition:open_close_sets}) is called \textit{differentiable}
if all its partial derivatives exist (i.e., the derivatives exist in univariate cases, and the gradients exist in multivariate cases). This is actually the definition of \textit{Fr\'echet differentiability}.

\index{(Fr\'echet) differentiability}
\begin{definition}[(Fr\'echet) Differentiability]
Given a function $f:\real^n\rightarrow \real$, the function $f$ is said to be differentiable at $\bx$ if there exists a vector $\bg\in\real^n$ such that 
$$
\mathop{\lim}_{\bd\rightarrow\bzero} \frac{f(\bx+\bd)-f(\bx)-\bg^\top\bd}{\normtwo{\bd}}=0.
$$
The unique vector $\bg$ is equal to the gradient $\nabla f(\bx)$.
\end{definition}

\index{Continuously differentiable}
Moreover, a function $f$ defined over an open set $\sS\subseteq \real^n$ is called \textit{continuously differentiable} over $\sS$ if all its partial derivatives exist and are also continuous on $\sS$.

A differentiable function may not be continuously differentiable. For example, consider  the following function $f$:
$$
f(x) = 
\begin{cases}
x^2 \sin(\frac{1}{x}), & \text{if } x \neq 0, \\
0, & \text{if } x = 0.
\end{cases}
$$
This function is differentiable everywhere, including at $x=0$. To see this, for $x \neq 0$, $f(x)$ is a product of two differentiable functions ($x^2$ and $\sin(1/x)$), hence it is differentiable.
At $x = 0$, we can calculate the derivative using the limit definition:
$$
f^\prime(0) = \mathop{\lim}_{\mu\rightarrow 0}
\frac{f(\mu) - f(0) }{\mu}
=\mathop{\lim}_{\mu\rightarrow 0} \frac{\mu^2 \sin(\frac{1}{\mu}) - 0 }{\mu}
=\mu \sin(\frac{1}{\mu}) 
$$
Since $\abs{\sin(1/\mu)}\leq 1$ for all $\mu \neq 0$, the limit exists and is equal to 0. Thus, $f(x)$ is differentiable at $x = 0$ with $f^\prime(0) = 0$.
However,  $f(x)$ is not continuously differentiable at $x = 0$. The derivative of $f(x)$ when $x\neq 0$ is 
$$
f^\prime(x) = 2x \sin(\frac{1}{x}) - \cos(\frac{1}{x}).
$$
The limit $\mathop{\lim}_{x\rightarrow 0}f^\prime(x)$ does not exist because the sine and cosine  functions oscillate as $x$ approaches 0.

In the setting of \textit{differentiability}, the directional derivative and gradient have the following relationship:
\begin{equation}\label{equation:direc_contdiff}
f^\prime(\bx; \bd) = \nabla f(\bx)^\top \bd, \gap \text{for all }\bx\in\sS \text{ and }\bd\in\real^n.
\end{equation} 
\begin{proof}
The formula is obviously correct for $\bd = \bzero$. We then assume that $\bd \neq \bzero$. The differentiability of $f$ implies that
$$
0 = \lim_{\mu \to 0^+} \frac{f(\bx + \mu \bd) - f(\bx) - \langle \nabla f(\bx), \mu \bd \rangle}{\normtwo{\mu \bd}} 
= \lim_{\mu \to 0^+} \left[ \frac{f(\bx + \mu \bd) - f(\bx)}{\mu \normtwo{\bd}} - \frac{\langle \nabla f(\bx), \bd \rangle}{\normtwo{\bd}} \right].
$$
Therefore,
$$
\begin{aligned}
&f'(\bx; \bd) = \lim_{\mu \to 0^+} \frac{f(\bx + \mu \bd) - f(\bx)}{\mu}\\
&= \lim_{\mu \to 0^+} \left\{ \normtwo{\bd} \left[ \frac{f(\bx + \mu \bd) - f(\bx)}{\mu \normtwo{\bd}} - \frac{\langle \nabla f(\bx), \bd \rangle}{\normtwo{\bd}} \right] + \langle \nabla f(\bx), \bd \rangle \right\}
&= \langle \nabla f(\bx), \bd \rangle.
\end{aligned}
$$
This proves \eqref{equation:direc_contdiff}.
\end{proof}

Recalling the definition of differentiability, we also have:
\begin{equation}
\mathop{\lim}_{\bd\rightarrow \bzero}
\frac{f(\bx+\bd) - f(\bx) - \nabla f(\bx)^\top \bd}{\normtwo{\bd}} = 0\gap 
\text{for all }\bx\in\sS,
\end{equation}
or 
\begin{equation}
f(\by) = f(\bx)+\nabla f(\bx)^\top (\by-\bx) + o(\normtwo{\by-\bx}),
\end{equation}
where the \textit{small-oh} function $o(\cdot): \real_+\rightarrow \real$ is a one-dimensional function satisfying $\frac{o(\mu)}{\mu}\rightarrow 0$ as $\mu\rightarrow 0^+$.~\footnote{Note that we also use the standard \textit{big-Oh} notation to describe the asymptotic behavior of functions.
Specifically, the notation $g(\bd) = \mathcalO(\normtwo{\bd}^p)$ means that there are positive numbers
$C_1$ and $\delta$ such that $\abs{g(\bd)} \leq  C_1 \normtwo{\bd}^p$ for all $\normtwo{\bd}\leq\delta$. In practice it is
often equivalent to $\abs{g(\bd)} \approx C_2\normtwo{\bd}^p$ for  sufficiently small $\bd$, where $C_2$ is  another positive constant.
The \textit{soft-Oh} notation is employed to hide poly-logarithmic factors i.e., $f = \widetilde{\mathcalO}(g)$ will
imply $f = \mathcalO(g \log^c(g))$ for some absolute constant $c$.}
Therefore, any differentiable function ensures that $\mathop{\lim}_{\bx\rightarrow \ba} f(\bx) = f(\ba)$. Hence, any differentiable function is continuous. 

On the other hand, in the setting of differentiability, although the partial derivatives may not be continuous, the partial derivatives exist for all differentiable points.

The partial derivative $\frac{\partial f}{\partial x_i} (\bx)$ is also a real-valued function of $\bx\in\sS$ that can be partially differentiated. The $j$-th partial derivative of $\frac{\partial f}{\partial x_i} (\bx)$ is defined as 
$$
\frac{\partial^2 f}{\partial x_j\partial x_i} (\bx)=
\frac{\partial \left(\frac{\partial f}{\partial x_i} (\bx)\right)}{\partial x_j} (\bx).
$$
This is called the ($j,i$)-th \textit{second-order partial derivative} of function $f$.
A function $f$ defined over an open set $\sS\subseteq$ is called \textit{twice continuously differentiable} over $\sS$ if all the second-order partial derivatives exist and are continuous over $\sS$. In the setting of twice continuously differentiability, the second-order partial derivative are symmetric:
$$
\frac{\partial^2 f}{\partial x_j\partial x_i} (\bx)=
\frac{\partial^2 f}{\partial x_i\partial x_j} (\bx).
$$
The \textit{Hessian} of the function $f$ at a point $\bx\in\sS$ is defined as the symmetric $n\times n$ matrix 
$$
\nabla^2f(\bx)=
\begin{bmatrix}
\frac{\partial^2 f}{\partial x_1^2} (\bx) & 
\frac{\partial^2 f}{\partial x_1\partial x_2} (\bx) & \ldots &
\frac{\partial^2 f}{\partial x_1\partial x_n} (\bx)\\
\frac{\partial^2 f}{\partial x_2\partial x_1} (\bx) & 
\frac{\partial^2 f}{\partial x_2\partial x_2} (\bx) & \ldots &
\frac{\partial^2 f}{\partial x_2\partial x_n} (\bx)\\
\vdots & 
\vdots & \ddots &
\vdots\\
\frac{\partial^2 f}{\partial x_n\partial x_1} (\bx) & 
\frac{\partial^2 f}{\partial x_n\partial x_2} (\bx) & \ldots &
\frac{\partial^2 f}{\partial x_n^2} (\bx)
\end{bmatrix}.
$$

We then provide a simple proof of Taylor's expansion   for one-dimensional functions. 
\begin{theorem}[Taylor's Expansion with Lagrange Remainder]
Let $f(x): \real\rightarrow \real$ be $k$-times continuously differentiable on the closed interval $I$ with endpoints $x$ and $y$, for some $k\geq 0$. If $f^{(k+1)}$ exists on the interval $I$, then there exists a $\xi \in (x,y)$ such that 
$$
\begin{aligned}
f(x)& = f(y) + f^\prime(y)(x-y) +\ldots + \frac{f^{(k)}(y)}{k!}(x-y)^k
+ \frac{f^{(k+1)}(\xi)}{(k+1)!}(x-y)^{k+1}\\
&=\sum_{i=0}^{k} \frac{f^{(i)}(y)}{i!} (x-y)^i + \frac{f^{(k+1)}(\xi)}{(k+1)!}(x-y)^{k+1}.
\end{aligned}
$$ 
Taylor's expansion can be extended to a function of vector $f(\bx):\real^n\rightarrow \real$ or a function of matrix $f(\bX): \real^{m\times n}\rightarrow \real$.
\end{theorem}
Taylor's expansion, or also known as  \textit{Taylor's series}, approximates the function $f(x)$ around a value  $y$ using a polynomial in a single variable $x$. To understand the origin of this series, we recall from the elementary calculus course that the approximated function  of $\cos (\theta)$ around $\theta=0$ is given by 
$
\cos (\theta) \approx 1-\frac{\theta^2}{2}.
$
This means that $\cos(\theta)$ can be approximated by a second-degree polynomial. 
If we want to approximate $\cos(\theta)$ more generally with a second-degree polynomial $ f(\theta) = c_1+c_2 \theta+ c_3 \theta^2$, an intuitive approach is to match the function and its derivatives at $\theta=0$. That is,
$$\left\{
\begin{aligned}
\cos(0) &= f(0); \\
\cos^\prime(0) &= f^\prime(0);\\
\cos^{\prime\prime}(0) &= f^{\prime\prime}(0);\\
\end{aligned}
\right.
\quad\implies\quad 
\left\{
\begin{aligned}
1 &= c_1; \\
-\sin(0) &=0= c_2;\\
-\cos(0) &=-1= 2c_3.\\
\end{aligned}
\right.
$$
Solving these equations yields $f(\theta) = c_1+c_2 \theta+ c_3 \theta^2 = 1-\frac{\theta^2}{2}$, which matches our initial approximation $\cos (\theta) \approx 1-\frac{\theta^2}{2}$ around  $\theta=0$. 

For high-dimensional functions, we have the following approximation results.
\begin{theorem}[Mean Value Theorem]\label{theorem:mean_approx}
Let $f(\bx):\sS\rightarrow \real$ be a  continuously differentiable function over an open set $\sS\subseteq\real^n$, and given two points $\bx, \by\in\sS$. Then, there exists a point $\bxi\in[\bx,\by]$ such that 
$$
f(\by) = f(\bx)+ \nabla f(\bxi)^\top (\by-\bx).
$$

\end{theorem}

\begin{theorem}[Linear Approximation Theorem]\label{theorem:linear_approx}
Let $f(\bx):\sS\rightarrow \real$ be a twice continuously differentiable function over an open set $\sS\subseteq\real^n$, and let $\bx, \by\in\sS$. Then, there exists a point $\bxi\in[\bx,\by]$ such that 
$$
f(\by) = f(\bx)+ \nabla f(\bx)^\top (\by-\bx) + \frac{1}{2} (\by-\bx)^\top \nabla^2 f(\bxi) (\by-\bx),
$$ 
or 
$$
f(\by) = f(\bx)+\nabla f(\bx)^\top (\by-\bx) + o(\normtwo{\by-\bx}),
$$
or
$$
f(\by) = f(\bx)+\nabla f(\bx)^\top (\by-\bx) + \mathcalO(\normtwo{\by-\bx}^2).
$$
\end{theorem}

\begin{theorem}[Quadratic Approximation Theorem]\label{theorem:quad_app_theo}
Let $f(\bx):\sS\rightarrow \real$ be a twice continuously differentiable function over an open set $\sS\subseteq\real^n$, and let $\bx, \by\in\sS$. Then it follows that 
$$
f(\by) = f(\bx)+ \nabla f(\bx)^\top (\by-\bx) + \frac{1}{2} (\by-\bx)^\top \nabla^2 f(\bx) (\by-\bx)
+
o(\normtwo{\by-\bx}^2),
$$
or 
$$
f(\by) = f(\bx)+ \nabla f(\bx)^\top (\by-\bx) + \frac{1}{2} (\by-\bx)^\top \nabla^2 f(\bx) (\by-\bx)
+
\mathcalO(\normtwo{\by-\bx}^3).
$$
\end{theorem}

\section{Fundamental Theorems}
This section introduces fundamental theorems from optimization, linear algebra, and calculus that are frequently utilized in subsequent discussions.
\begin{theorem}[Fundamental Theorem of Optimization]\label{theorem:funda_opt}
Let $\ba,\bb\in\real^n$ be two vectors. The inner product of these two vectors can be expressed as a sum of norms:
$$
\ba^\top\bb = \frac{1}{2} \left( \normtwo{\ba}^2 + \normtwo{\bb}^2 - \normtwo{\ba-\bb}^2 \right).
$$
\end{theorem}

The following theorem provides an elementary proof that the row rank and column rank of any matrix $\bA\in \real^{m\times n}$ are equal, which also highlights the fundamental theorem of linear algebra. \index{Matrix rank}
\begin{theorem}[Row Rank Equals Column Rank\index{Matrix rank}]\label{theorem:equal-dimension-rank}
The dimension of the column space of a matrix $\bA\in \real^{m\times n}$ is equal to the dimension of its row space. 
In other words, the row rank and the column rank of a matrix $\bA$ are identical.
\end{theorem}

\begin{proof}[of Theorem~\ref{theorem:equal-dimension-rank}]
Firstly, observe that the null space of $\bA$ is orthogonal complementary to the row space of $\bA$: $\nspace(\bA) \bot \cspace(\bA^\top)$ (where the row space of $\bA$ corresponds to the column space of $\bA^\top$). 
This means that any vector in the null space of $\bA$ is orthogonal to every vector in the row space of $\bA$. 
To illustrate this, assume $\bA$ has rows $\ba_1^\top, \ba_2^\top, \ldots, \ba_m^\top$ and let $\bA=[\ba_1^\top; \ba_2^\top; \ldots; \ba_m^\top]$. For any vector $\bx\in \nspace(\bA)$, we have $\bA\bx = \bzero$, implying $[\ba_1^\top\bx; \ba_2^\top\bx; \ldots; \ba_m^\top\bx]=\bzero$. And since the row space of $\bA$ is spanned by $\ba_1^\top, \ba_2^\top, \ldots, \ba_m^\top$, it follows that $\bx$ is perpendicular to all vectors in $\cspace(\bA^\top)$, confirming $\nspace(\bA) \bot \cspace(\bA^\top)$.

Next, suppose the dimension of the row space of $\bA$ is $r$. \textcolor{mylightbluetext}{Let $\br_1, \br_2, \ldots, \br_r\in\real^n$ be a basis for the row space}. 
Consequently, the $r$ vectors $\bA\br_1, \bA\br_2, \ldots, \bA\br_r$ are in the column space of $\bA$ and are linearly independent. To see this, suppose we have a linear combination of the $r$ vectors: $x_1\bA\br_1 + x_2\bA\br_2+ \ldots+ x_r\bA\br_r=\bzero$, that is, $\bA(x_1\br_1 + x_2\br_2+ \ldots+ x_r\br_r)=\bzero$, and the vector $\bv\triangleq x_1\br_1 + x_2\br_2+ \ldots+ x_r\br_r$ is in null space of $\bA$. But since $\{\br_1, \br_2, \ldots, \br_r\}$ is a basis for the row space of $\bA$, $\bv$ must also lie in the row space of $\bA$. We have shown that vectors from the null space of $\bA$ are perpendicular to vectors from the row space of $\bA$; thus, it holds that $\bv^\top\bv=0$ and $x_1=x_2=\ldots=x_r=0$. Therefore, \textcolor{mylightbluetext}{$\bA\br_1, \bA\br_2, \ldots, \bA\br_r$ are in the column space of $\bA$, and they are linearly independent}. This implies that the dimension of the column space of $\bA$ is larger than $r$. Therefore, \textbf{the row rank of $\bA\leq $ column rank of $\bA$}. 

By applying the same argument to $\bA^\top$, we deduce that \textbf{column rank of $\bA\leq $ row rank of $\bA$}, completing the proof that the row and column ranks of $\bA$ are equal.
\end{proof}

From this proof, it follows that if $\{\br_1, \br_2, \ldots, \br_r\}$ is a set of vectors in $\real^n$ forming a basis for the row space, then \textcolor{black}{$\{\bA\br_1, \bA\br_2, \ldots, \bA\br_r\}$ constitutes a basis for the column space of $\bA$}. 
This observation is formalized  in the following lemma.

\begin{lemma}[Column Basis from Row Basis]\label{lemma:column-basis-from-row-basis}
Let $\bA\in \real^{m\times n}$ be a matrix, and suppose that $\{\br_1, \br_2, \ldots, \br_r\}$ is a set of vectors in $\real^n$, which forms a basis for the row space. Then, $\{\bA\br_1, \bA\br_2, \ldots, \bA\br_r\}$ constitutes a basis for the column space of $\bA$.
\end{lemma}

\index{Fundamental theorem}
For any matrix $\bA\in \real^{m\times n}$, it is straightforward to verify that any vector in the row space of $\bA$ is orthogonal to any vector in the null space of $\bA$. Suppose $\bx_n \in \nspace(\bA)$, then $\bA\bx_n = \bzero$, meaning that $\bx_n$ is orthogonal to every row of $\bA$, thereby confirming this property.

Similarly, we can also demonstrate that any vector in the column space of $\bA$ is perpendicular to any vector in the null space of $\bA^\top$. 
Furthermore, the column space of $\bA$ together with the null space of $\bA^\top$ span the entire space $\real^m$, which is a key aspect of the fundamental theorem of linear algebra.

The fundamental theorem consists of two aspects: the dimensions of the subspaces and their orthogonality. 
The orthogonality can be easily verified as shown above. 
Additionally, when the row space has dimension $r$, the null space has dimension $n-r$. 
This is not immediately obvious and is proven in the following theorem.

\begin{figure}[h!]
\centering
\includegraphics[width=0.98\textwidth]{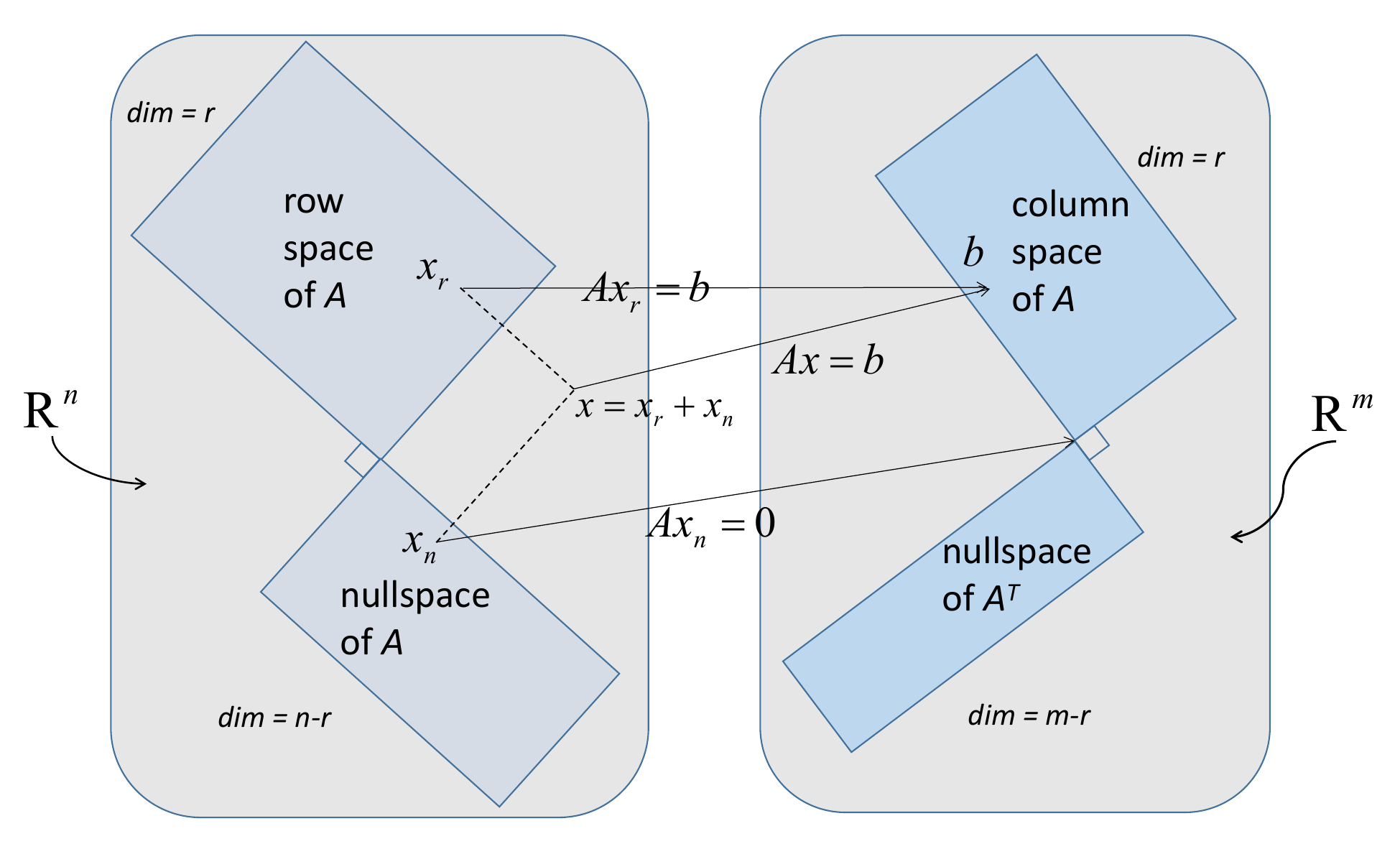}
\caption{Two pairs of orthogonal subspaces in $\real^n$ and $\real^m$. $\dim(\cspace(\bA^\top)) + \dim(\nspace(\bA))=n$ and $\dim(\nspace(\bA^\top)) + \dim(\cspace(\bA))=m$. The null space component maps to zero as $\bA\bx_n = \bzero \in \real^m$. The row space component maps to the column space as $\bA\bx_r = \bA(\bx_r+\bx_n)=\bb \in \cspace(\bA)$.}
\label{fig:lafundamental}
\end{figure}
\begin{theorem}[Fundamental Theorem of Linear Algebra]\label{theorem:fundamental-linear-algebra}
Orthogonal Complement and Rank-Nullity Theorem: Let $\bA\in \real^{m\times n}$ be a matrix. Then, 
\begin{itemize}
\item $\nspace(\bA)$ is the orthogonal complement to the row space $\cspace(\bA^\top)$ in $\real^n$: $\dim(\nspace(\bA))+\dim(\cspace(\bA^\top))=n$;
\item $\nspace(\bA^\top)$ is the orthogonal complement to the column space $\cspace(\bA)$ in $\real^m$: $\dim(\nspace(\bA^\top))+\dim(\cspace(\bA))=m$;
\item If $\bA$ has rank $r$, then $\dim(\cspace(\bA^\top)) = \dim(\cspace(\bA)) = r$, implying that $\dim(\nspace(\bA)) = n-r$ and $\dim(\nspace(\bA^\top))=m-r$.
\end{itemize}
\end{theorem}

\begin{proof}[of Theorem~\ref{theorem:fundamental-linear-algebra}]
Building upon the proof of Theorem~\ref{theorem:equal-dimension-rank}, consider a set of vectors $\br_1, \br_2, \ldots, \br_r$  in $\real^n$ forming a basis for the row space.
Consequently, \textcolor{black}{$\{\bA\br_1, \bA\br_2, \ldots, \bA\br_r\}$ constitutes a basis for the column space of $\bA$}. 
Let $\bn_1, \bn_2, \ldots, \bn_k \in \real^n$ form a basis for the null space of $\bA$. Following again from the proof of Theorem~\ref{theorem:equal-dimension-rank}, it follows that $\nspace(\bA) \bot \cspace(\bA^\top)$,  indicating the orthogonality between $\br_1, \br_2, \ldots, \br_r$ and $\bn_1, \bn_2, \ldots, \bn_k$. Consequently, the set   $\{\br_1, \br_2, \ldots, \br_r, \bn_1, \bn_2, \ldots, \bn_k\}$ is linearly independent in $\real^n$.

For any vector $\bx\in \real^n $, $\bA\bx$ is in the column space of $\bA$. Then it can be expressed as a linear combination of $\bA\br_1, \bA\br_2, \ldots, \bA\br_r$: $\bA\bx \triangleq \sum_{i=1}^{r}a_i\bA\br_i$, which states that $\bA(\bx-\sum_{i=1}^{r}a_i\br_i) = \bzero$ and $\bx-\sum_{i=1}^{r}a_i\br_i$ is thus in $\nspace(\bA)$. Since $\{\bn_1, \bn_2, \ldots, \bn_k\}$ is a basis for the null space of $\bA$, $\bx-\sum_{i=1}^{r}a_i\br_i$ can be expressed as a linear combination of $\bn_1, \bn_2, \ldots, \bn_k$: $\bx-\sum_{i=1}^{r}a_i\br_i = \sum_{j=1}^{k}b_j \bn_j$, i.e., $\bx=\sum_{i=1}^{r}a_i\br_i + \sum_{j=1}^{k}b_j \bn_j$. That is, any vector $\bx\in \real^n$ can be expressed as $\{\br_1, \br_2, \ldots, \br_r, \bn_1, \bn_2, \ldots, \bn_k\}$ and the set forms a basis for $\real^n$. Thus, the dimensions satisfy: $r+k=n$, i.e., $\dim(\nspace(\bA))+\dim(\cspace(\bA^\top))=n$. Similarly, we can prove that $\dim(\nspace(\bA^\top))+\dim(\cspace(\bA))=m$.
\end{proof}

Figure~\ref{fig:lafundamental} visually represents the orthogonality of these subspaces and illustrates how $\bA$ maps a vector $\bx$ into the column space. The row space and null space dimensions add up to $n$, while the dimensions of the column space of $\bA$ and the null space of $\bA^\top$ sum to $m$. The null space component maps to zero as $\bA\bx_{\bn} = \bzero \in \real^m$, which is the intersection of the column space of $\bA$ and the null space of $\bA^\top$. 
Conversely, the row space component transforms into the column space as $\bA\bx_{\br} = \bA(\bx_{\br} + \bx_{\bn})=\bb\in \real^m$.

We provide the fundamental theorem of calculus below, which  plays a pivotal role in linking differential and integral calculus, offering profound insights into the behavior of functions. 
\begin{theorem}[Fundamental Theorem of Calculus]\label{theorem:fund_theo_calculu}
Let $f(\cdot):\real^n\rightarrow \real$. The fundamental theorem of calculus describes the difference between function values:
\begin{subequations}
\begin{equation}\label{equation:fund_theo_calculu3}
f(\by) - f(\bx) = \int_{0}^{1} \langle \nabla f(\bx+\mu(\by-\bx)), \by-\bx  \rangle d\mu, 
\end{equation}
or the difference between gradients:
\begin{align}
\nabla f(\by) -\nabla f(\bx) &= \left( \int_{0}^{1} \nabla^2f(\bx+\mu(\by-\bx))d\mu \right) \cdot (\by-\bx); \label{equation:fund_theo_calculu1}\\
\nabla f(\bx+\alpha\bd) - \nabla f(\bx) &= \int_{0}^{\alpha} \nabla^2 f(\bx+\mu\bd)\bd d\mu.
\end{align}
\end{subequations}
where $ \nabla f(\by) $ denotes the gradient of $ f $ evaluated at $ \by $ (see Section~\ref{section:differ_calc}), and $ \innerproduct{\cdot, \cdot} $ represents the inner product.
Furthermore, use directional derivative for twice continuously differentiable functions, we also have
\begin{equation}
f(\by) = f(\bx) + \innerproduct{\nabla f(\bx),  (\by-\bx)} + \int_{0}^{1} (1-\mu) \frac{\partial^2 f(\bx+\mu(\by-\bx))}{\partial \mu^2} d\mu.
\end{equation}
\end{theorem}
The first equality \eqref{equation:fund_theo_calculu3} can be derived by letting $g(t)\triangleq f(\bx+t(\by-\bx))$. Then we can obtain $f(\by)$ as 
$
f(\by) = g(1) = g(0) + \int_{0}^{1} g'(t)dt.
$
The other equalities  can be demonstrated similarly. 

\index{Implicit function theorem}
We present the implicit function theorem without a proof. Further  details can be found in \citet{krantz2002implicit}.
\begin{theorem}[Implicit Function Theorem]\label{theorem:implic_func_theorem}
Let $ \bh(\bx)=[h_i(\bx)]_{i\in\{1,2,\ldots,p\}} : \mathbb{R}^n \rightarrow \mathbb{R}^p $ ($p<n$) and $ \widetilde{\bx} = [\widetilde{y}_1, \ldots, \widetilde{y}_p, \widetilde{z}_1, \ldots, \widetilde{z}_{n-p}] = [\widetilde{\by}, \widetilde{\bz}] \in\real^n $ satisfy:
\begin{enumerate}[(i)]
\item $ \bh(\widetilde{\bx}) = \bzero $.
\item $ \bh(\bx) $ is continuously differentiable in a neighborhood of $ \widetilde{\bx} $.
\item The $ p \times p $ Jacobian matrix
$$
\small
\begin{bmatrix}
\frac{\partial h_1(\widetilde{\bx})}{\partial y_1} & \cdots & \frac{\partial h_1(\widetilde{\bx})}{\partial y_p} \\
\vdots & \ddots & \vdots \\
\frac{\partial h_p(\widetilde{\bx})}{\partial y_1} & \cdots & \frac{\partial h_p(\widetilde{\bx})}{\partial y_p}
\end{bmatrix}
$$
is nonsingular.
\end{enumerate}
Then, there exists $ \varepsilon > 0 $ along with functions $ \bs(\bz) = [s_1(\bz), \ldots, s_p(\bz)]^\top\in\real^p $ such that for all $ \bz \in \sB(\widetilde{\bz}, \varepsilon) $, $ \bh(\bs(\bz), \bz) = \bzero $. Moreover, each $ s_k(\bz) $ is continuously differentiable for $ k \in \{1,2, \ldots, p\} $. Furthermore, for all $ i \in \{1,2, \ldots, p\} $ and $ j \in \{1,2, \ldots, n-p\} $, we have:
$$
\sum_{k=1}^p \frac{\partial h_i(\by, \bz)}{\partial y_k} \frac{\partial s_k(\bz)}{\partial z_j} + \frac{\partial h_i(\by, \bz)}{\partial z_j} = 0.
$$
\end{theorem}

\index{Inverse function theorem}
\begin{theorem}[Inverse Function Theorem]
Let $ F: \mathbb{R}^n \to \mathbb{R}^n $ be a continuously differentiable function, and let $ \bx^* $ be a point in $ \mathbb{R}^n $. If the derivative of $ F $ at $ \bx^* $, denoted as $ \nabla F(\bx^*) $, is an invertible linear transformation (meaning the Jacobian matrix $ \bJ_F(\bx^*) $ is nonsingular), then there exists an open neighborhood $ \sS $ of $ \bx^* $ such that:
\begin{enumerate}[(i)]
\item $ F $ maps $ \sS $ bijectively onto an open set $ \sV = F(\sS) $,
\item The inverse function $ F^{-1}: \sV \to \sS $ is also continuously differentiable,
\item For each $ \by \in \sV $, the derivative of the inverse function at $ \by $ is the inverse of the derivative of $ F $ at $ \bx $, where $ F(\bx) = \by $; i.e.,
$$
\nabla (F^{-1})(\by) = [\nabla F(\bx)]^{-1}.
$$
\end{enumerate}
\end{theorem}
In simpler terms, if a function locally behaves like a linear transformation (i.e., its derivative is nonzero and invertible), then near that point, an inverse function exists that ``undoes" the original function. This inverse function is also smooth (continuously differentiable).
The inverse function theorem is a powerful tool that has wide-ranging applications in mathematics, physics, economics, and engineering. For example, the theorem can be used to solve systems of nonlinear equations by transforming them into a system of linear equations.

\begin{theorem}[Von-Neumann Lemma]\label{theorem:von_neumannlem}
Let $\norm{\cdot}$ be a submultiplicative matrix norm with $\norm{\bI} = 1$, and let $\bE \in \real^{n \times n}$. If $\norm{\bE} < 1$, then $\bI - \bE$ is nonsingular, and
\begin{equation}\label{equation:von_neumannlem1}
(\bI - \bE)^{-1} = \sum_{t=0}^{\infty} \bE^t,
\end{equation}
\begin{equation}\label{equation:von_neumannlem2}
\norm{(\bI - \bE)^{-1}} \leq \frac{1}{1 - \norm{\bE}}.
\end{equation}
If $\bA \in \real^{n \times n}$ is nonsingular and $\norm{\bA^{-1}(\bB - \bA)} < 1$, then $\bB$ is nonsingular and satisfies
\begin{equation}\label{equation:von_neumannlem3}
\bB^{-1} = \sum_{t=0}^{\infty} (\bI - \bA^{-1}\bB)^t \bA^{-1},
\end{equation}
and
\begin{equation}\label{equation:von_neumannlem4}
\norm{\bB^{-1}} \leq \frac{\norm{\bA^{-1}}}{1 - \norm{\bA^{-1}(\bB - \bA)}}.
\end{equation}
\end{theorem}
\begin{proof}[of Theorem~\ref{theorem:von_neumannlem}]
Since $\norm{\bE} < 1$, then
$ \bS_t \triangleq \bI + \bE + \bE^2 + \cdots + \bE^t $
defines a Cauchy sequence and thus $\bS_t$ converges. 
Taking the limit, we obtain
$ \sum_{t=0}^{\infty} \bE^t = \lim_{t \to \infty} \bS_t = (\bI - \bE)^{-1} $
which proves \eqref{equation:von_neumannlem1} and \eqref{equation:von_neumannlem2}.
For the second part, since $\bA$ is nonsingular and $\norm{\bA^{-1}(\bB - \bA)} = \norm{-(\bI - \bA^{-1}\bB)} < 1$, we set $\bE \triangleq \bI - \bA^{-1}\bB$ and apply \eqref{equation:von_neumannlem1} and \eqref{equation:von_neumannlem2}, immediately yielding \eqref{equation:von_neumannlem3} and \eqref{equation:von_neumannlem4}. 
\end{proof}

\begin{theorem}[Inverse Function Theorem: a Variant]\label{theorem:invfunctheo_var}
Let $ F : \sS \subset \mathbb{R}^n \to \mathbb{R}^n $ be continuously differentiable in a neighborhood of $ \bx^* \in \sS $, and suppose that $ \nabla F(\bx^*) $ is nonsingular. Then there exist $ \delta > 0 $, $ \xi > 0 $, and $ \epsilon > 0 $, such that when $ \normtwo{\by - \bx^*} < \delta $ and $ \by \in \sS $, $ \nabla F(\by) $ is nonsingular and
\begin{equation}\label{equation:invfunctheo_var1}
\normtwo{\nabla F(\by)^{-1}} \leq \xi,
\end{equation}
where $\normtwo{\cdot}$ for a matrix denotes the spectral norm.
Additionally, the inverse $ \nabla F(\by)^{-1} $ is continuous at $ \bx^* $, meaning
$$
\normtwo{\nabla F(\by)^{-1} - \nabla F(\bx^*)^{-1}} < \epsilon.
$$
\end{theorem}
\begin{proof}
Define $\alpha \triangleq \normtwo{\nabla F(\bx^*)^{-1}}$. Since $F$ is continuously differentiable, for a given $\beta < \alpha^{-1}$, choose $\delta$ such that when $\normtwo{\by - \bx^*} < \delta$ with $\by \in \sS$,
$$
\normtwo{\nabla F(\bx^*) - \nabla F(\by)} \leq \beta.
$$
It follows from Von-Neumann lemma (Theorem~\ref{theorem:von_neumannlem}) that $\nabla F(\by)$ is invertible, and \eqref{equation:von_neumannlem4} shows that 
$$
\normtwo{\nabla F(\by)^{-1}} \leq  \frac{\normtwo{\nabla F(\bx^*)^{-1}}}{1 - \normtwo{\nabla F(\bx^*)^{-1}(\nabla F(\by) - \nabla F(\bx^*))}}
\leq \alpha / (1 - \beta \alpha).
$$
Therefore, \eqref{equation:invfunctheo_var1} holds with $\xi \triangleq \alpha / (1 - \beta \alpha)$. Thus,
$$
\begin{aligned}
\normtwo{\nabla F(\bx^*)^{-1} - \nabla F(\by)^{-1}}
&= \normtwo{\nabla F(\bx^*)^{-1}(\nabla F(\by) - \nabla F(\bx^*))\nabla F(\by)^{-1}}
\leq \alpha \beta \xi 
\triangleq \epsilon,
\end{aligned}
$$
which shows that the continuity of $\nabla F$ guarantees the continuity of $(\nabla F)^{-1}$.
\end{proof}

\section{Matrix Decomposition}
This section provide several matrix decomposition approaches, which can be useful for the proofs in the optimization methods.
\subsection{Cholesky Decomposition}\label{section:choleskydecomp}
Positive definiteness or positive semidefiniteness is one of the most desirable properties a matrix can have. In this section, we introduce decomposition techniques for positive definite matrices, with a focus on the well-known \textit{Cholesky decomposition}.
The Cholesky decomposition is named after a French military officer and mathematician, Andr\'{e}-Louis Cholesky (1875{\textendash}1918), who developed this method in his surveying work. It is primarily used to solve linear systems involving positive definite matrices.

Here, we establish the existence of the Cholesky decomposition using an inductive approach. Alternative proofs also exist, such as those derived from the LU decomposition  \citep{lu2022matrix}.

\index{Cholesky decomposition}
\begin{theorem}[Cholesky Decomposition]\label{theorem:cholesky-factor-exist}
Every positive definite (PD) matrix $\bA\in \real^{n\times n}$ can be factored as 
$$
\bA = \bR^\top\bR,
$$
where $\bR \in \real^{n\times n}$ is an upper triangular matrix \textbf{with positive diagonal elements}. This decomposition is called the \textit{Cholesky decomposition}  of $\bA$, and $\bR$ is known as the \textit{Cholesky factor} or \textit{Cholesky triangle} of $\bA$.
Specifically, the Cholesky decomposition is unique (Corollary~\ref{corollary:unique-cholesky-main}).
\end{theorem}
Alternatively, $\bA$ can be factored as $\bA=\bL\bL^\top$, where $\bL=\bR^\top$ is a lower triangular matrix \textit{with positive diagonals}.
\begin{proof}[of Theorem~\ref{theorem:cholesky-factor-exist}]
We will prove by induction that every $n\times n$ positive definite matrix $\bA$ has a decomposition $\bA=\bR^\top\bR$. The $1\times 1$ case is trivial by setting $R\triangleq\sqrt{A}$, so that $A=R^2$. 

Suppose any $k\times k$ PD matrix $\bA_k$ has a Cholesky decomposition. We must show that any $(k+1)\times(k+1)$ PD matrix $\bA_{k+1}$ can also be factored as this Cholesky decomposition, then we complete the proof.

For any $(k+1)\times(k+1)$ PD matrix $\bA_{k+1}$, write $\bA_{k+1}$ as
$
\bA_{k+1} \triangleq \footnotesize\begin{bmatrix}
\bA_k & \bb \\
\bb^\top & d
\end{bmatrix}.
$
Since $\bA_k$ is PD, by the inductive hypothesis, it admits a Cholesky decomposition  $\bA_k = \bR_k^\top\bR_k$. Define the upper triangular matrix
$
\bR_{k+1}\triangleq\footnotesize\begin{bmatrix}
\bR_k & \br\\
0 & s
\end{bmatrix}.
$
Then,
$$
\bR_{k+1}^\top\bR_{k+1} = 
\begin{bmatrix}
\bR_k^\top\bR_k & \bR_k^\top \br\\
\br^\top \bR_k & \br^\top\br+s^2
\end{bmatrix}.
$$
Therefore, if we can prove $\bR_{k+1}^\top \bR_{k+1} = \bA_{k+1}$ is the Cholesky decomposition of $\bA_{k+1}$ (which requires the value $s$ to be positive), then we complete the proof. That is, we need to prove
$$
\begin{aligned}
\bb &= \bR_k^\top \br, \\
d &= \br^\top\br+s^2.
\end{aligned}
$$
Since $\bR_k$ is nonsingular, we can solve uniquely for $\br$ and $s$: 
$$
\begin{aligned}
\br &= \bR_k^{-\top}\bb, \\
s &= \sqrt{d - \br^\top\br} = \sqrt{d - \bb^\top\bA_k^{-1}\bb},
\end{aligned}
$$
where we assume $s$ is nonnegative. However, we need to further prove that $s$ is not only nonnegative, but also positive. Since $\bA_k$ is PD, from Sylvester's criterion (see \citet{lu2022matrix}), and the fact that if matrix $\bM$ has a block formulation: $\bM=\footnotesize\begin{bmatrix}
\bA & \bB \\
\bC & \bD 
\end{bmatrix}$, then $\det(\bM) = \det(\bA)\det(\bD-\bC\bA^{-1}\bB)$, we have
$$
\det(\bA_{k+1}) = \det(\bA_k)\det(d- \bb^\top\bA_k^{-1}\bb) =  \det(\bA_k)(d- \bb^\top\bA_k^{-1}\bb)>0.
$$
Because $ \det(\bA_k)>0$, we then obtain that $(d- \bb^\top\bA_k^{-1}\bb)>0$, and this implies $s>0$.
This completes the proof.
\end{proof}

\index{Uniqueness}
\begin{corollary}[Uniqueness of Cholesky Decomposition\index{Uniqueness}]\label{corollary:unique-cholesky-main}
The Cholesky decomposition $\bA=\bR^\top\bR$ for any positive definite matrix $\bA\in \real^{n\times n}$ is unique.
\end{corollary}
\begin{proof}[of Corollary~\ref{corollary:unique-cholesky-main}]
Suppose the Cholesky decomposition is not unique. Then, there exist two distinct decompositions such that $\bA=\bR_1^\top\bR_1 = \bR_2^\top\bR_2$.  Rearranging, we obtain
$$
\bR_1\bR_2^{-1}= \bR_1^{-\top} \bR_2^\top.
$$
Since the inverse of an upper triangular matrix is also upper triangular, and the product of two upper triangular matrices is upper triangular, \footnote{Similarly, the inverse of a lower triangular matrix is lower triangular, and the product of two lower triangular matrices is also lower triangular.} we conclude that the left-hand side of the above equation is an upper triangular matrix, while the right-hand side is a lower triangular matrix. Consequently, $\bR_1\bR_2^{-1}= \bR_1^{-\top} \bR_2^\top$ must be a diagonal matrix, and $\bR_1^{-\top} \bR_2^\top= (\bR_1^{-\top} \bR_2^\top)^\top = \bR_2\bR_1^{-1}$.
Let $\bLambda \triangleq \bR_1\bR_2^{-1}= \bR_2\bR_1^{-1}$ be the diagonal matrix. We notice that the diagonal value of $\bLambda$ is the product of the corresponding diagonal values of $\bR_1$ and $\bR_2^{-1}$ (or $\bR_2$ and $\bR_1^{-1}$). Explicitly, writing the matrices as
$$
\bR_1=
\small
\begin{bmatrix}
r_{11} & r_{12} & \ldots & r_{1n} \\
0 & r_{22} & \ldots & r_{2n}\\
\vdots & \vdots & \ddots & \vdots\\
0 & 0 & \ldots & r_{nn}
\end{bmatrix},
\qquad 
\normalsize
\bR_2=
\small
\begin{bmatrix}
s_{11} & s_{12} & \ldots & s_{1n} \\
0 & s_{22} & \ldots & s_{2n}\\
\vdots & \vdots & \ddots & \vdots\\
0 & 0 & \ldots & s_{nn}
\end{bmatrix},
$$
we find that
$$
\begin{aligned}
\bR_1\bR_2^{-1}=
\small
\begin{bmatrix}
\frac{r_{11}}{s_{11}} & 0 & \ldots & 0 \\
0 & \frac{r_{22}}{s_{22}} & \ldots & 0\\
\vdots & \vdots & \ddots & \vdots\\
0 & 0 & \ldots & \frac{r_{nn}}{s_{nn}}
\end{bmatrix}
=
\small
\begin{bmatrix}
\frac{s_{11}}{r_{11}} & 0 & \ldots & 0 \\
0 & \frac{s_{22}}{r_{22}} & \ldots & 0\\
\vdots & \vdots & \ddots & \vdots\\
0 & 0 & \ldots & \frac{s_{nn}}{r_{nn}}
\end{bmatrix}
\normalsize
=\bR_2\bR_1^{-1}.
\end{aligned}
$$ 
Since both $\bR_1$ and $\bR_2$ have positive diagonals, this implies $r_{11}=s_{11}, r_{22}=s_{22}, \ldots, r_{nn}=s_{nn}$. 
Thus, we conclude that $\bLambda = \bR_1\bR_2^{-1}= \bR_2\bR_1^{-1}  =\bI$, which implies $\bR_1=\bR_2$, contradicting our initial assumption that the decomposition was not unique. Therefore, the Cholesky decomposition is unique.
\end{proof}		

\paragraph{Computing Cholesky decomposition element-wise.}

It is common to compute the Cholesky decomposition using element-wise equations derived directly from solving the matrix equation $\bA=\bR^\top\bR$. Observing that the $(i,j)$-th entry of $\bA$ is given by $a_{ij} = \bR_{:,i}^\top \bR_{:,j} = \sum_{k=1}^{i} r_{ki}r_{kj}$ if $i<j$. This further implies the  following recurrence relation: if $i<j$, we have
$$
\begin{aligned}
a_{ij} &= \bR_{:,i}^\top \bR_{:,j} = \sum_{k=1}^{i} r_{ki}r_{kj} 
= \sum_{k=1}^{i-1} r_{ki}r_{kj} + r_{ii}r_{ij}
\implies
r_{ij} = (a_{ij} - \sum_{k=1}^{i-1} r_{ki}r_{kj})/r_{ii},
\gap 
\text{if }i<j.
\end{aligned}
$$
For the diagonal entries ($i=j$), we have:
$$
\begin{aligned}
a_{jj} &= \sum_{k=1}^{j} r_{kj}^2=\sum_{k=1}^{j-1} r_{kj}^2 + r_{jj}^2
&\implies
r_{jj} = \sqrt{a_{jj} - \sum_{k=1}^{j-1} r_{kj}^2}.
\end{aligned}
$$
If we equate the elements of $\bR$ by taking a column at a time and start with $r_{11} = \sqrt{a_{11}}$, the element-level algorithm is formulated in Algorithm~\ref{alg:compute-choklesky-element-level}.

\begin{algorithm}[H] 
\caption{Cholesky Decomposition Element-Wise: $\bA=\bR^\top\bR$} 
\label{alg:compute-choklesky-element-level} 
\begin{algorithmic}[1] 
\Require 
Positive definite matrix $\bA$ with size $n\times n$; 
\State Calculate first element of $\bR$ by $r_{11} \leftarrow \sqrt{a_{11}}$; 
\For{$j=1$ to $n$} \Comment{Compute the $j$-th column of $\bR$}
\For{$i=1$ to $j-1$} 
\State $r_{ij} \leftarrow (a_{ij} - \sum_{k=1}^{i-1} r_{ki}r_{kj})/r_{ii}$, since $i<j$;
\EndFor
\State $r_{jj} \leftarrow \sqrt{a_{jj}- \sum_{k=1}^{j-1}r_{kj}^2}$;
\EndFor
\State Output $\bA=\bR^\top\bR$.
\end{algorithmic} 
\end{algorithm}
On the other hand, Algorithm~\ref{alg:compute-choklesky-element-level}  can be modified to compute the Cholesky decomposition in the form $\bA=\bL\bD\bL^\top$, where $\bL$ is unit lower triangular and $\bD$ is diagonal, as outlined in Algorithm~\ref{alg:compute-choklesky-_ldl}, whose Step 3 and Step 5 are derived from (since $l_{ii}=1, \forall i\in\{1,2,\ldots,n\}$):
$$
\begin{aligned}
a_{jj}&=\sum_{k=1}^{j-1}d_{kk} l_{jk}^2 + d_{jj};\\
a_{ij}&= d_{jj} l_{ij}+ \sum_{k=1}^{j-1} d_{kk} l_{ik}l_{jk}, \gap \text{if }i>j.
\end{aligned}
$$
\begin{exercise}
Derive the complexity of Algorithm~\ref{alg:compute-choklesky-_ldl}.
\end{exercise}
This form of Cholesky decomposition is useful for determining the condition number of a PD matrix \citep{lu2021numerical}. 
In essence, the condition number of a function measures the sensitivity of the output value to small changes in the input; a smaller condition number indicates better numerical stability. 
For positive definite linear systems, the condition number is defined as the ratio of the largest eigenvalue to the smallest eigenvalue.
The condition number of a positive definite matrix is lower bounded by the diagonal matrix in the Cholesky decomposition (see Problem~\ref{problem:cond_pd}):
\begin{equation}\label{equation:cond_pd_ineq}
\cond(\bA) \geq \cond(\bD).
\end{equation}
This can be proven by showing that $\lambda_{\max}\geq d_{\max}$ and $\lambda_{\min}\leq d_{\min}$, where $\lambda_{\max}$ and $\lambda_{\min}$ are the largest and smallest eigenvalue of $\bA$, and $d_{\max}$ and $d_{\min}$ are the largest and smallest diagonals of $\bD$.
Therefore, this form of the Cholesky decomposition can be utilized  to modify Newton's method; see \S~\ref{section:modified_damp_new}.

\begin{algorithm}[h] 
\caption{Cholesky Decomposition Element-Wise: $\bA=\bL\bD\bL^\top$}  
\label{alg:compute-choklesky-_ldl} 
\begin{algorithmic}[1] 
\Require 
Positive definite matrix $\bA$ with size $n\times n$; 
\For{$j=1$ to $n$} \Comment{Compute the $j$-th column of $\bL$}
\State $l_{jj}\leftarrow1$;
\State $c_{jj}\leftarrow a_{jj}-\sum_{k=1}^{j-1}d_{kk} l_{jk}^2$;
\State $d_{jj}\leftarrow c_{jj}$
\For{$i=j+1$ to $n$} 
\State $c_{ij}\leftarrow a_{ij}-\sum_{k=1}^{j-1}d_{kk} l_{ik}l_{jk}$, since $i>j$;
\State $l_{ij}\leftarrow \frac{c_{ij}}{d_{jj}}$;
\EndFor
\EndFor
\State Output $\bA=\bL\bD\bL^\top$, where $\bD=\diag(d_{11}, d_{22},\ldots,d_{nn})$.
\end{algorithmic} 
\end{algorithm}

\subsection{Eigenvalue and Spectral Decomposition}\label{section:eigendecomp}\label{section:spectraldecomp}\label{section:existence-of-spectral}
\textit{Eigenvalue decomposition} is also known as to diagonalize the matrix. If a matrix $\bA\in\real^{n\times n}$ has distinct eigenvalues, its corresponding eigenvectors are guaranteed to be linearly independent, allowing $\bA$ to be diagonalized. It is important to note that without $n$ linearly independent eigenvectors, diagonalization is not possible.
\index{Eigenvalue decomposition}
\begin{theorem}[Eigenvalue Decomposition]\label{theorem:eigenvalue-decomposition}
Any square matrix $\bA\in \real^{n\times n}$ with linearly independent eigenvectors can be factored as 
$$
\bA = \bX\bLambda\bX^{-1},
$$
where $\bX\in\real^{n\times n}$ is a matrix whose columns are the eigenvectors of $\bA$, and $\bLambda\in\real^{n\times n}$ is a diagonal matrix given by $\diag(\lambda_1, $ $\lambda_2, \ldots, \lambda_n)$ with $\lambda_1, \lambda_2, \ldots, \lambda_n$ denoting the eigenvalues of $\bA$.
\end{theorem}
\begin{proof}[of Theorem~\ref{theorem:eigenvalue-decomposition}]
Let $\bX=[\bx_1, \bx_2, \ldots, \bx_n]$ be the linearly independent eigenvectors of $\bA$. Clearly, we have
$$
\bA\bx_1=\lambda_1\bx_1,\qquad \bA\bx_2=\lambda_2\bx_2, \qquad \ldots, \qquad\bA\bx_n=\lambda_n\bx_n.
$$
Expressing this in matrix form, we obtain
$$
\bA\bX = [\bA\bx_1, \bA\bx_2, \ldots, \bA\bx_n] = [\lambda_1\bx_1, \lambda_2\bx_2, \ldots, \lambda_n\bx_n] = \bX\bLambda.
$$
Since the eigenvectors are linearly independent, $\bX$ is invertible, leading to
$
\bA = \bX\bLambda \bX^{-1}.
$
This completes the proof.
\end{proof}

An important advantage of the decomposition $\bA =\bX\bLambda\bX^{-1}$ is that it allows efficient computation of matrix powers.
\index{$m$-th Power}
\begin{remark}[$m$-th Power]\label{remark:power-eigenvalue-decom}
If $\bA$ admits an eigenvalue decomposition $\bA =\bX\bLambda\bX^{-1}$, then its $m$-th power  can be computed as  $\bA^m = \bX\bLambda^m\bX^{-1}$.
\end{remark}

The existence of eigenvalue decomposition depends on the linear independence of the eigenvectors of $\bA$. This condition is automatically satisfied in specific cases.
\begin{lemma}[Different Eigenvalues]\label{lemma:diff-eigenvec-decompo}
If the eigenvalues $\lambda_1, \lambda_2, \ldots, \lambda_n$ of $\bA\in \real^{n\times n}$ are all distinct, then the associated eigenvectors are necessarily linearly independent. Consequently, any square matrix with unique eigenvalues is diagonalizable. 
\end{lemma}
\begin{proof}[of Lemma~\ref{lemma:diff-eigenvec-decompo}]
Assume that $\bA$ has distinct eigenvalues $\lambda_1, \lambda_2, \ldots, \lambda_n$ and that the eigenvectors $\bx_1,\bx_2, \ldots, \bx_n$ are linearly dependent. Then, there exists a nonzero vector $\bc = [c_1,c_2,\ldots,c_{n-1}]^\top$ such that
$
\bx_n = \sum_{i=1}^{n-1} c_i\bx_{i}. 
$
Applying $\bA$ to both sides, we get
$$
\begin{aligned}
\bA \bx_n &= \bA (\sum_{i=1}^{n-1} c_i\bx_{i}) 
=c_1\lambda_1 \bx_1 + c_2\lambda_2 \bx_2 + \ldots + c_{n-1}\lambda_{n-1}\bx_{n-1}.
\end{aligned}
$$
and 
$$
\begin{aligned}
\bA \bx_n &= \lambda_n\bx_n
=\lambda_n (c_1\bx_1 +c_2\bx_2+\ldots +c_{n-1} \bx_{n-1}).
\end{aligned}
$$
Equating both expressions, we obtain
$
\sum_{i=1}^{n-1} (\lambda_n - \lambda_i)c_i \bx_i = \bzero .
$
This leads to a contradiction since $\lambda_n \neq \lambda_i$ for all $i\in \{1,2,\ldots,n-1\}$, from which the result follows.
\end{proof}

The \textit{spectral theorem}, also referred to as  the \textit{spectral decomposition} for symmetric matrices, states that every symmetric matrix has real eigenvalues and can be diagonalized using a (real) orthonormal basis \footnote{Note that the spectral decomposition for \textit{Hermitian matrices} states that Hermitian matrices also have real eigenvalues but are diagonalizable using a complex orthonormal basis.}. 
\index{Spectral decomposition}
\begin{theorem}[Spectral Decomposition]\label{theorem:spectral_theorem}
A real matrix $\bA \in \real^{n\times n}$ is symmetric if and only if there exists an orthogonal matrix $\bQ$ and a diagonal matrix $\bLambda$ such that
\begin{equation}\label{equation:spec_decom}
\bA = \bQ \bLambda \bQ^\top,
\end{equation}
where the columns of $\bQ = [\bq_1, \bq_2, \ldots, \bq_n]$ are eigenvectors of $\bA$ and are mutually orthonormal, and the entries of $\bLambda=\diag(\lambda_1, \lambda_2, \ldots, \lambda_n)$ are the corresponding \textbf{real} eigenvalues of $\bA$. 
Moreover, the rank of $\bA$ is equal to the number of nonzero eigenvalues. 
This result is known as the \textit{spectral decomposition} or the \textit{spectral theorem} for real symmetric matrices. Specifically, the following properties hold:
\begin{enumerate}[(i)]
\item A symmetric matrix has only \textbf{real eigenvalues}.
\item The eigenvectors are orthogonal and can be chosen to be \textbf{orthonormal} by normalization.
\item The rank of $\bA$ is equal to the number of nonzero eigenvalues (Proposition~\ref{proposition:rank-of-symmetric}).
\item If the eigenvalues are distinct, the eigenvectors are necessarily linearly independent.
\end{enumerate}
\end{theorem}

The existence of the spectral decomposition follows directly from Symmetric Properties-I$\sim$III (Propositions~\ref{proposition:real-eigenvalues-spectral}$\sim$\ref{proposition:eigen-multiplicity}). Note that the Symmetric Property-IV (Proposition~\ref{proposition:rank-of-symmetric}) is proved using the main result in the spectral decomposition.
Similar to the eigenvalue decomposition, spectral decomposition allows for efficient computation of the $m$-th power of a matrix.
\index{$m$-th Power}
\begin{remark}[$m$-th Power]\label{remark:power-spectral}
If a matrix $\bA$ admits a spectral decomposition $\bA=\bQ\bLambda\bQ^\top$, then its $m$-th power can be computed as  $\bA^m = \bQ\bLambda^m\bQ^\top$.
\end{remark}

\subsection{QR Decomposition}
In many applications, we are interested in the column space of a matrix $\bA=[\ba_1, \ba_2, \ldots, \ba_n] \in \real^{m\times n}$. The successive subspaces spanned by the columns $\ba_1, \ba_2, \ldots$ of $\bA$ plays a crucial role in understanding the structure and properties of the matrix:
$$
\cspace([\ba_1])\,\,\,\, \subseteq\,\,\,\, \cspace([\ba_1, \ba_2]) \,\,\,\,\subseteq\,\,\,\, \cspace([\ba_1, \ba_2, \ba_3])\,\,\,\, \subseteq\,\,\,\, \ldots,
$$
where $\cspace([\ldots])$ denotes the subspace spanned by the enclosed vectors, alternatively expressed as $\cspace([\ldots])=\spn\{\ldots\}$.
The concept of orthogonal or orthonormal bases within the column space is fundamental to many algorithms, enabling efficient computations and interpretations. \textit{QR factorization} is a widely used technique for analyzing and decomposing matrices in a way that explicitly reveals their column space structure.

QR decomposition constructs a sequence of orthonormal vectors $\bq_1, \bq_2, \ldots$ that span the same successive subspaces:
$$
\big\{\cspace([\bq_1])=\cspace([\ba_1])\big\}\subseteq 
\big\{\cspace([\bq_1, \bq_2])=\cspace([\ba_1, \ba_2])\big\}\subseteq
\big\{\cspace([\bq_1, \bq_2, \bq_3])=\cspace([\ba_1, \ba_2, \ba_3])\big\} 
\subseteq \ldots.
$$
\begin{theorem}[QR Decomposition]\label{theorem:qr-decomposition}
Let $\bA=[\ba_1, \ba_2, \ldots, \ba_n] \in\real^{m\times n}$ be any matrix (regardless of whether its columns are linearly independent or dependent) with $m\geq n$. Then, it can be factored as 
$$
\bA = \bQ\bR,~\footnote{If $\bA$ is complex, then $\bQ$ is unitary (or semi-unitary with orthonormal columns) and $\bR$ is complex upper triangular.}
$$
where 
\begin{enumerate}
\item  \textbf{Reduced}: $\bQ$ is $m\times n$ with orthonormal columns, and $\bR$ is an $n\times n$ upper triangular matrix, known as the \textbf{reduced QR decomposition};
\item \textbf{Full}: $\bQ$ is $m\times m$ with orthonormal columns, and $\bR$ is an $m\times n$ upper triangular matrix, known as the \textbf{full QR decomposition}. If we further restrict the upper triangular matrix to be a square matrix, the full QR decomposition takes the form
$$
\bA = \bQ\begin{bmatrix}
\bR_0\\
\bzero
\end{bmatrix},
$$
where $\bR_0$ is an $n\times n$ upper triangular matrix. 
\end{enumerate}

Specifically, when $\bA$ has full rank, i.e., $\bA$  has linearly independent columns, then $\bR$ also exhibits linearly independent columns, and $\bR$ is nonsingular in the \textit{reduced} case, meaning its diagonal elements are nonzero.
Under this condition, when we further restrict elements on the diagonal of $\bR$ to be  \textbf{positive} and $\rank(\bA)=n$, the \textit{reduced} QR decomposition is \textbf{unique} \citep{lu2022matrix}. However, the \textit{full} QR decomposition is generally not unique since the rightmost $(m-n)$ columns in $\bQ$ can be arranged in any order.
\end{theorem}

%
%

\paragraph{Project vector $\ba$ onto vector $\bb$.}
The projection of a vector $\ba$ onto another vector $\bb$ finds the closest vector to $\ba$ that lies along the line spanned by $\bb$.
The projected vector, denoted by $\widehat{\ba}$, is a scalar multiple of $\bb$. 
Let $\widehat{\ba} \triangleq \widehat{x} \bb$. By construction, $\ba-\widehat{\ba}$ is perpendicular to $\bb$, leading to the following result:
$$
\text{$\ba^\perp \triangleq(\ba-\widehat{\ba}) \perp \bb 
\quad\implies\quad 
(\ba-\widehat{x}\bb)^\top\bb=0
\quad\implies\quad 
\widehat{x}$ = $\frac{\ba^\top\bb}{\bb^\top\bb}$ and $\widehat{\ba} = \frac{\ba^\top\bb}{\bb^\top\bb}\bb = \frac{\bb\bb^\top}{\bb^\top\bb}\ba$}.
$$
The above discussion leads to the following lemma.
\begin{lemma}[Finding an Orthogonal Vector]
Given two unit vectors $\bu$ and $\bv$ (i.e., $\normtwo{\bu}=\normtwo{\bv}=1$), then $\bw\triangleq\bv-\bu\bu^\top\bv$ is orthogonal to $\bu$.
If $\bu$ and $\bv$ are not unit vectors, they can be normalized $\bu:=\frac{\bu}{\normtwo{\bu}}$ and $\bv:=\frac{\bv}{\normtwo{\bv}}$ to achieve the same result.
\end{lemma}

Given three linearly independent vectors $\{\ba_1, \ba_2, \ba_3\}$ and the space spanned by the three linearly independent vectors $\cspace{([\ba_1, \ba_2, \ba_3])}$, i.e., the column space of the matrix $[\ba_1, \ba_2, \ba_3]$, we aim to construct three orthogonal vectors $\{\bb_1, \bb_2, \bb_3\}$ such that $\cspace{([\bb_1, \bb_2, \bb_3])}$ = $\cspace{([\ba_1, \ba_2, \ba_3])}$, i.e., the column space remains unchanged. 
We then normalize these orthogonal vectors to obtain an orthonormal set. This process produces three mutually orthonormal vectors $\bq_1 \triangleq \frac{\bb_1}{\normtwo{\bb_1}}$, $\bq_2 \triangleq \frac{\bb_2}{\normtwo{\bb_2}}$, and $\bq_2 \triangleq \frac{\bb_2}{\normtwo{\bb_2}}$.

\begin{figure}[H]
\centering  
\vspace{-0.35cm} 
\subfigtopskip=2pt 
\subfigbottomskip=2pt 
\subfigcapskip=-5pt 
\subfigure[Project $\ba_2$ onto the space perpendicular to $\bb_1$.]{\label{fig:gram-schmidt1}
\includegraphics[width=0.4\linewidth]{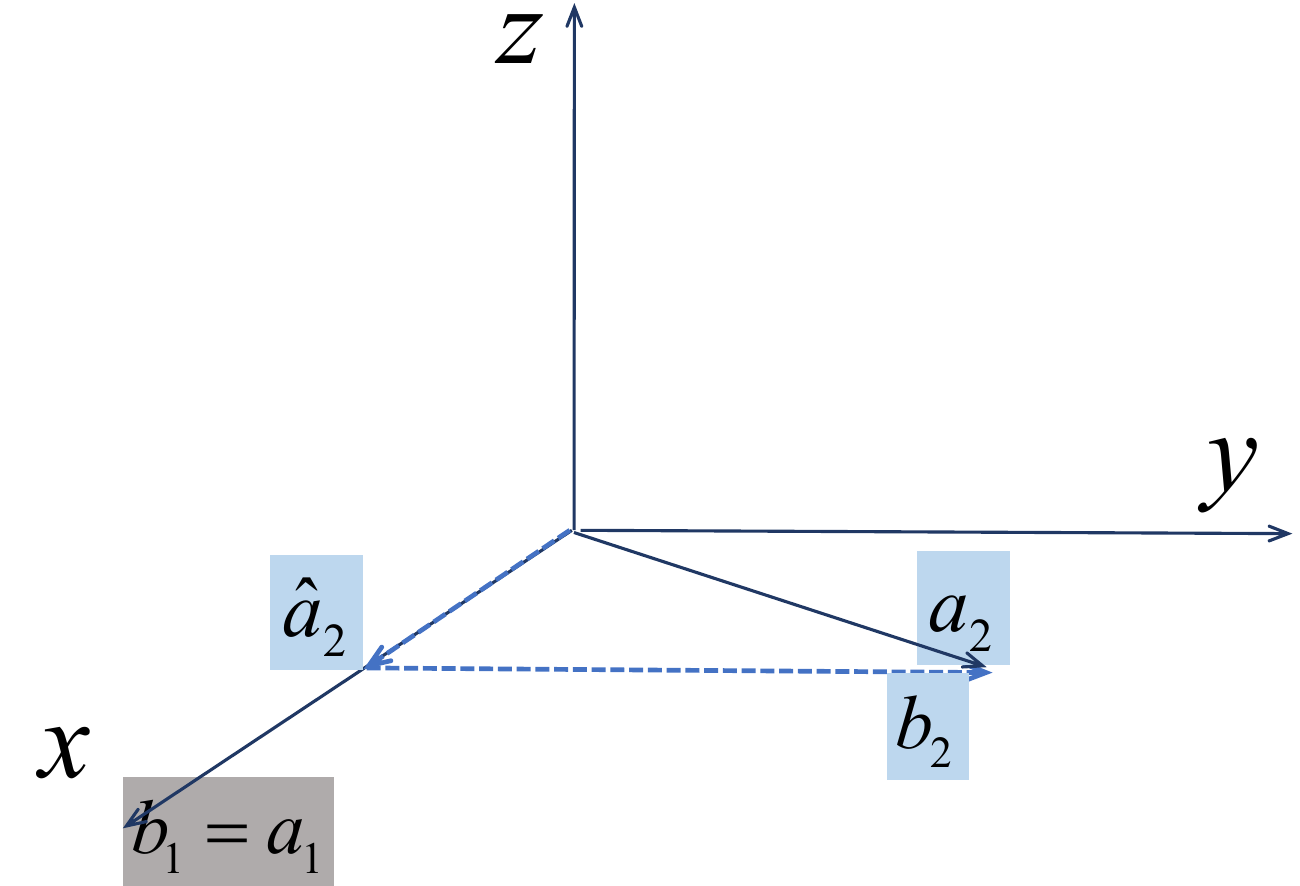}}
\quad 
\subfigure[Project $\ba_3$ onto the space perpendicular to $\bb_1, \bb_2$.]{\label{fig:gram-schmidt2}
\includegraphics[width=0.4\linewidth]{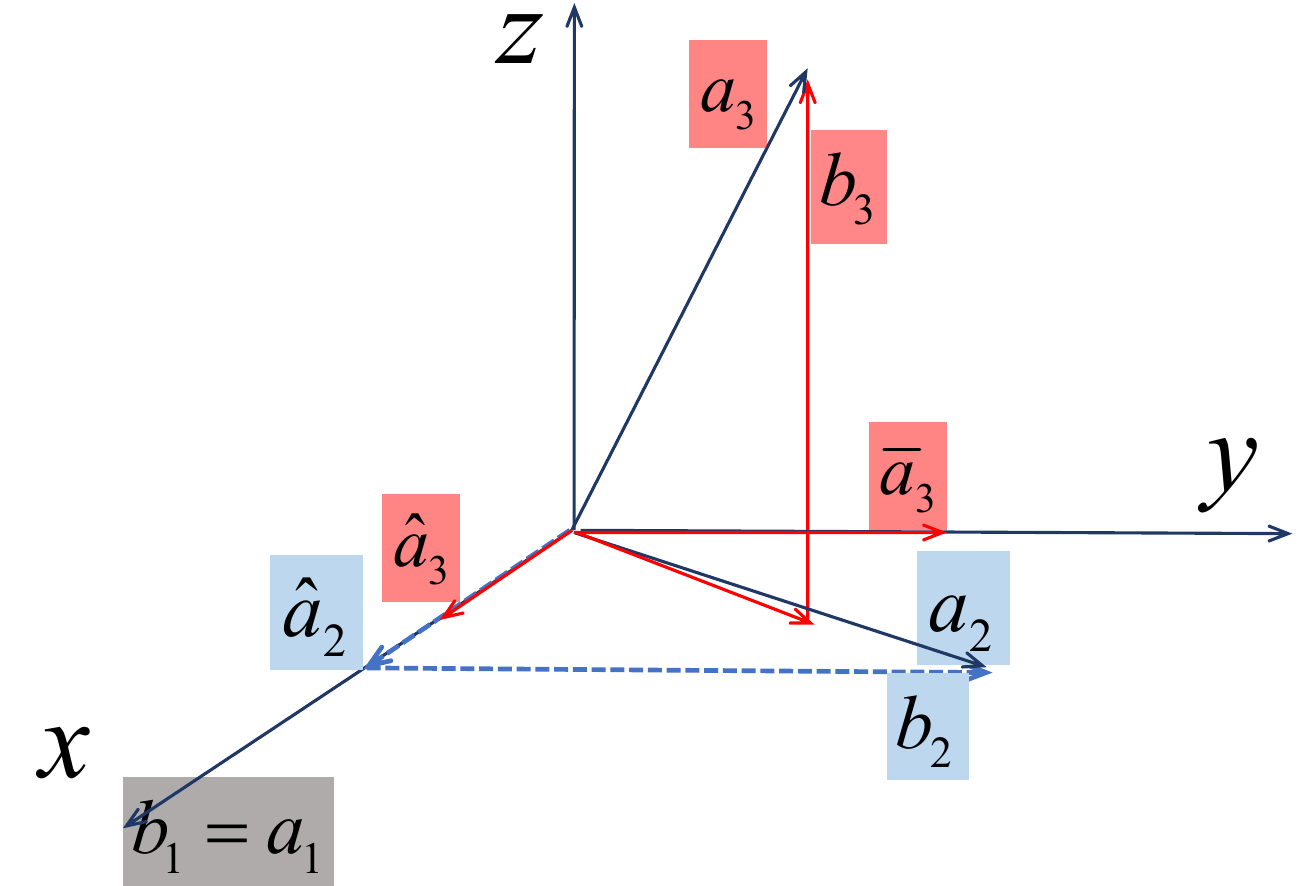}}
\caption{The Gram-Schmidt process.}
\label{fig:gram-schmidt-12}
\end{figure}

For the first vector, we simply set $\bb_1 \triangleq \ba_1$. The second vector $\bb_2$ must be perpendicular to the first one. This is achieved by considering the vector $\ba_2$ and subtracting its projection along $\bb_1$:
\begin{equation}
\begin{aligned}
\bb_2 &
= \left(\bI- \frac{\bb_1 \bb_1^\top}{\bb_1^\top\bb_1} \right)\ba_2  
= \ba_2-  \underbrace{\frac{ \bb_1^\top \ba_2}{\bb_1^\top\bb_1} \bb_1}_{\triangleq\widehat{\ba}_2},
\end{aligned}
\end{equation}
where the first equation shows that $\bb_2$ results from the multiplication of the matrix $\big(\bI- \frac{\bb_1 \bb_1^\top}{\bb_1^\top\bb_1} \big)$ and the vector $\ba_2$, signifying the projection of $\ba_2$ onto the orthogonal complement space of $\cspace{([\bb_1])}$. The second equality in the above equation shows that $\ba_2$ is a linear combination of $\bb_1$ and $\bb_2$.
Clearly, the space spanned by $\bb_1, \bb_2$ coincides with the space spanned by $\ba_1, \ba_2$. 
Figure~\ref{fig:gram-schmidt1} illustrates this process, where we designate \textbf{the direction of $\bb_1$ as the $x$-axis in the Cartesian coordinate system}, and $\widehat{\ba}_2$ represents the projection of $\ba_2$ onto the line defined by $\bb_1$. 
The figure clearly demonstrates that  the component of $\ba_2$ that is perpendicular to $\bb_1$ is $\bb_2 = \ba_2 - \widehat{\ba}_2$, as derived from the Pythagorean theorem.

For the third vector $\bb_3$, it must be perpendicular to both  $\bb_1$ and $\bb_2$, which corresponds to the vector $\ba_3$ subtracting its projection along the plane spanned by $\bb_1$ and $\bb_2$:
\begin{equation}\label{equation:gram-schdt-eq2}
\begin{aligned}
\bb_3 &
= \left(\bI- \frac{\bb_1 \bb_1^\top}{\bb_1^\top\bb_1}  - \frac{\bb_2 \bb_2^\top}{\bb_2^\top\bb_2} \right)\ba_3  
= \ba_3- \underbrace{\frac{ \bb_1^\top\ba_3}{\bb_1^\top\bb_1} \bb_1}_{\triangleq\widehat{\ba}_3} - \underbrace{\frac{ \bb_2^\top\ba_3}{\bb_2^\top\bb_2}  \bb_2}_{\triangleq\bar{\ba}_3},   
\end{aligned}
\end{equation}
where the first equation shows that $\bb_3$ is a multiplication of the matrix $\left(\bI- \frac{\bb_1 \bb_1^\top}{\bb_1^\top\bb_1}  - \frac{\bb_2 \bb_2^\top}{\bb_2^\top\bb_2} \right)$ and the vector $\ba_3$, signifying the projection of $\ba_3$ onto the orthogonal complement space of $\cspace{([\bb_1, \bb_2])}$. The second equality in the above equation shows that $\ba_3$ is a linear combination of $\bb_1, \bb_2, $ and $\bb_3$. This property is essential in the idea of the QR decomposition.
Once again, it can be shown that the space spanned by $\bb_1, \bb_2, \bb_3$ is identical to the space spanned by $\ba_1, \ba_2, \ba_3$. 
Figure~\ref{fig:gram-schmidt2} illustrates this step, where we designate \textbf{the direction of $\bb_2$ as the $y$-axis in the Cartesian coordinate system}. Here, $\widehat{\ba}_3$ is the projection of $\ba_3$ onto line $\bb_1$, while $\bar{\ba}_3$ represents the projection of $\ba_3$ onto line $\bb_2$. 
It can be shown that the component of $\ba_3$ perpendicular to both $\bb_1$ and $\bb_2$ is $\bb_3=\ba_3-\widehat{\ba}_3-\bar{\ba}_3$ from the figure.

Finally, we normalize each vector by dividing its  length, resulting in three orthonormal vectors $\bq_1 \triangleq \frac{\bb_1}{\normtwo{\bb_1}}$, $\bq_2 \triangleq \frac{\bb_2}{\normtwo{\bb_2}}$, and $\bq_2 \triangleq \frac{\bb_2}{\normtwo{\bb_2}}$.

This procedure can be extended beyond three vectors to an arbitrary set and is known as the \textit{Gram-Schmidt process}. After applying this process, the matrix $\bA$ is transformed into an upper triangular form. The method is named after \textit{Jørgen Pedersen Gram} and \textit{Erhard Schmidt}, though  it appeared earlier in the work of \textit{Pierre-Simon Laplace} within the context of Lie group decomposition theory.

As previously mentioned, the key idea behind QR decomposition is to construct a sequence of orthonormal vectors $\bq_1, \bq_2, \ldots$, which successively span the same subspaces as the original vectors:
$$
\big\{\cspace([\bq_1])=\cspace([\ba_1])\big\} \subseteq
\big\{\cspace([\bq_1, \bq_2])=\cspace([\ba_1, \ba_2])\big\} \subseteq
\big\{\cspace([\bq_1, \bq_2, \bq_3])=\cspace([\ba_1, \ba_2, \ba_3])\big\} 
\subseteq \ldots.
$$
This implies that  any vector $\ba_k$ resides in the space spanned by $\cspace([\bq_1, \bq_2, \ldots, \bq_k])$.~\footnote{Moreover, any vector $\bq_k$ resides in the space spanned by $\cspace([\ba_1, \ba_2, \ldots, \ba_k])$.} Once these orthonormal vectors are determined, the reconstruction of $\ba_i$'s from the orthogonal matrix $\bQ=[\bq_1, \bq_2, \ldots, \bq_n]$ requires an upper triangular matrix $\bR$,  such that $\bA = \bQ\bR$.

\index{Gram–Schmidt}
While the Gram-Schmidt process is a widely used algorithm for QR decomposition, it is not the only approach. Alternative methods, such as \textit{Householder reflections} and \textit{Givens rotations}, offer more numerical stability in the presence of round-off errors and may modify the order in which the columns of $\bA$ are processed \citep{lu2021numerical}.

\subsection{Singular Value Decomposition}
In eigenvalue decomposition, we factor the matrix into a diagonal matrix. However, this is not always possible. 
When $\bA$ lacks linearly independent eigenvectors, such diagonalization cannot be achieved.  
The \textit{singular value decomposition (SVD)} fills this gap. 
Instead of using an eigenvector matrix, SVD decomposes a matrix into two orthogonal matrices. Below, we state the theorem for SVD.

\index{Karhunen-Loewe expansion}
\index{Decomposition: SVD}
\begin{theoremHigh}[Reduced SVD for Rectangular Matrices]\label{theorem:reduced_svd_rectangular}
Let $\bA\in\real^{m\times n}$ be any real $m\times n$ matrix with rank $r$. Then, it can be factored as
$$
\bA = \bU \bSigma \bV^\top,
$$ 
where $\bSigma\in \real^{r\times r}$ is a real diagonal matrix, $\bSigma=\diag(\sigma_1, \sigma_2 \ldots, \sigma_r)$ with $\sigma_1 \geq \sigma_2 \geq \ldots \geq \sigma_r>0$ and 
\footnote{Note when $\bA$ is complex,  both $\bU$ and $\bV$ become semi-unitary. However, $\bSigma$ remains real and diagonal because the eigenvalues of Hermitian matrices $\bA\bA^*$ or $\bA^*\bA$ are real.}
\begin{itemize}
\item The $\sigma_i$'s are the nonzero \textit{singular values} of $\bA$, in the meantime, which are the positive square roots of the nonzero \textit{eigenvalues} of $\trans{\bA} \bA$ and $ \bA \trans{\bA}$.

\item The columns of $\bU\in \real^{m\times r}$ contain the $r$ eigenvectors of $\bA\bA^\top$ corresponding to the $r$ nonzero eigenvalues of $\bA\bA^\top$ (semi-orthogonal). 

\item The columns of $\bV\in \real^{n\times r}$ contain the $r$ eigenvectors of $\bA^\top\bA$ corresponding to the $r$ nonzero eigenvalues of $\bA^\top\bA$ (semi-orthogonal). 


\item The columns of $\bU$ and $\bV$ are called the \textit{left and right singular vectors} of $\bA$, respectively. 

\item Furthermore, the columns of both $\bU$ and $\bV$ are mutually orthonormal (by spectral theorem~\ref{theorem:spectral_theorem}). 
\end{itemize}

In particular, the matrix decomposition can be expressed as the sum of outer products of vectors:
$$
\bA = \bU \bSigma \bV^\top = \sum_{i=1}^r \sigma_i \bu_i \bv_i^\top,
$$ which is a sum of $r$ rank-one matrices.

\end{theoremHigh}

By appending $m-r$ additional orthonormal columns to the  $r$  eigenvectors of   $\bA\bA^\top$, we obtain an orthogonal matrix $\bU\in \real^{m\times m}$. 
The same process applies to the columns of $\bV$.
We then present the full SVD  in the following theorem. 
Differences between the reduced and full SVD are highlighted in \textcolor{mylightbluetext}{blue} text.
\begin{theoremHigh}[Full SVD for Rectangular Matrices]\label{theorem:full_svd_rectangular}
Let $\bA\in\real^{m\times n}$ be any real $m\times n$ matrix with rank $r$. Then, it can be factored as
$$
\bA = \bU \bSigma \bV^\top,
$$ 
where the  top-left part of $\bSigma\in $\textcolor{mylightbluetext}{$\real^{m\times n}$} is a real diagonal matrix, that is $\bSigma=\footnotesize\begin{bmatrix}
\bSigma_1 & \bzero \\
\bzero & \bzero
\end{bmatrix}$, where $\bSigma_1=\diag(\sigma_1, \sigma_2 \ldots, \sigma_r)\in \real^{r\times r}$ with $\sigma_1 \geq \sigma_2 \geq \ldots \geq \sigma_r>0$ and 
\begin{itemize}
\item The $\sigma_i$'s are the nonzero \textit{singular values} of matrix $\bA$, which are the positive square roots of the nonzero \textit{eigenvalues} of $\trans{\bA} \bA$ and $ \bA \trans{\bA}$. 

\item  The columns of $\bU\in \textcolor{mylightbluetext}{\real^{m\times m}}$ contain the $r$ eigenvectors of $\bA\bA^\top$ corresponding to the $r$ nonzero eigenvalues of $\bA\bA^\top$ \textcolor{mylightbluetext}{and $m-r$ extra orthonormal vectors from $\nspace(\bA^\top)$}. 

\item  The columns of $\bV\in \textcolor{mylightbluetext}{\real^{n\times n}}$ contain the $r$ eigenvectors of $\bA^\top\bA$ corresponding to the $r$ nonzero eigenvalues of $\bA^\top\bA$ \textcolor{mylightbluetext}{and $n-r$ extra orthonormal vectors from $\nspace(\bA)$}.

\item  The columns of  $\bU$ and $\bV$ are called the \textit{left and right singular vectors} of $\bA$, respectively. 

\item Furthermore, the columns of $\bU$ and $\bV$ are mutually orthonormal (by spectral theorem~\ref{theorem:spectral_theorem}), and \textcolor{mylightbluetext}{$\bU$ and $\bV$ are orthogonal matrices}. 
\end{itemize}

In particular, the matrix decomposition can be expressed as the sum of outer products of vectors: $ \bA = \bU \bSigma \bV^\top = \sum_{i=1}^r \sigma_i \bu_i \bv_i^\top$, which is a sum of $r$ rank-one matrices.
\end{theoremHigh}

\begin{figure}[H]
\centering  
\vspace{-0.35cm} 
\subfigtopskip=2pt 
\subfigbottomskip=2pt 
\subfigcapskip=-5pt 
\subfigure[Reduced SVD decomposition.]{\label{fig:svdhalf}
\includegraphics[width=0.47\linewidth]{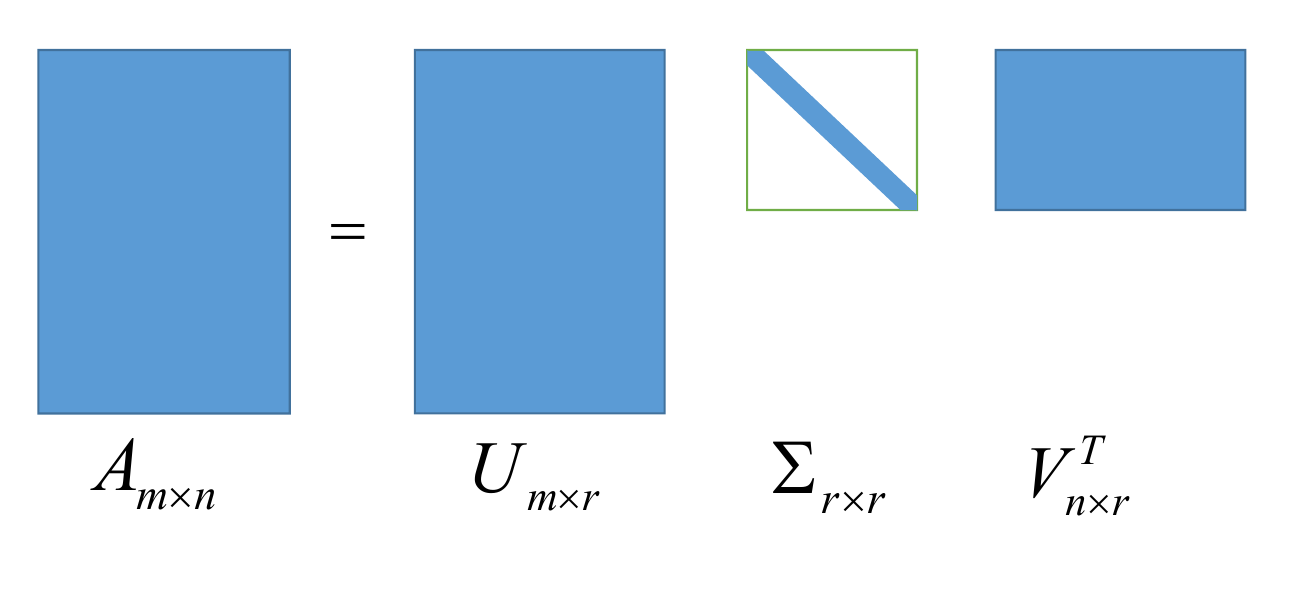}}
\quad 
\subfigure[Full SVD decomposition.]{\label{fig:svdall}
\includegraphics[width=0.47\linewidth]{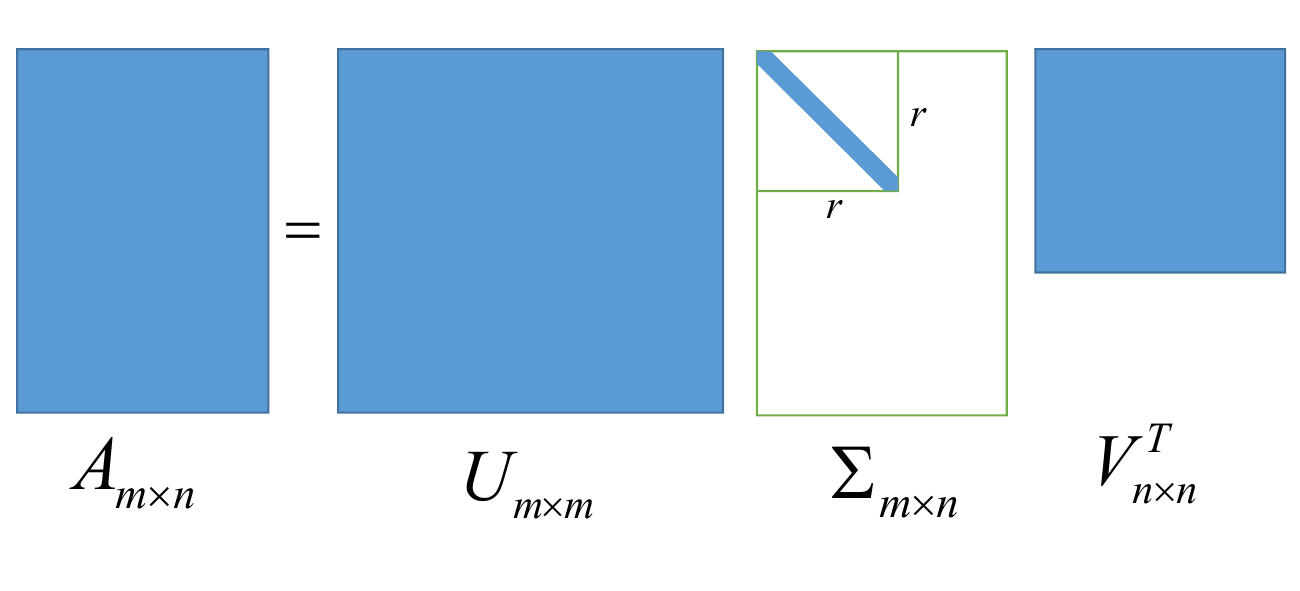}}
\caption{Comparison between the reduced and full SVD.}
\label{fig:svd-comparison}
\end{figure}

In image processing,  SVD is also known as the \textit{Karhunen-Loewe expansion}.
The comparison between the reduced and the full SVD is shown in Figure~\ref{fig:svd-comparison}, where white entries are zero, and \textcolor{mylightbluetext}{blue} entries are not necessarily zero.

\section{Classes of Rate of Convergence}
Before we delve into specific  algorithms for solving optimization problems, it is essential to establish criteria for evaluating the speed at which these algorithms converge, as most of them are iterative methods. We define the convergence of a sequence as follows. Note that the $t$-th element in a sequence is denoted by a superscript in parentheses. 
For instance, $\bA^{(t)}$ denotes the $t$-th matrix in a sequence, and $\ba^{(t)}$ represents the $t$-th vector in a sequence. \index{Rate of convergence}
\begin{definition}[Convergence of a Sequence]
Let $\mu^{(1)}, \mu^{(2)}, \ldots \in \real$ be an infinite sequence of scalars. The sequence $\{\mu^{(t)}\}_{t>0}$ is said to converge to $\mu^*$ if 
$$
\mathop{\lim}_{t\rightarrow \infty}  \abs{\mu^{(t)} - \mu^*} = 0.
$$
Similarly, let $\bmu^{(1)}, \bmu^{(2)}, \ldots \in \real^n$ be an infinite sequence of vectors. The sequence $\{\bmu^{(t)}\}_{t>0}$ is said to converge to $\bmu^\star$ if 
$$
\mathop{\lim}_{t\rightarrow \infty} \norm{\bmu^{(t)} - \bmu^\star} = 0.
$$
\end{definition}
The convergence of a sequence of vectors or matrices depends on the norm used. It's important to note that, according to the equivalence of vector or matrix norms (Theorem~\ref{theorem:equivalence-vector-norm}), if a sequence of vectors or matrices converges in one norm, it will also converge in all other norms. 

\begin{definition}[Linear Convergence]\label{definition:linear-convergence}
A sequence $\{\mu^{(t)}\}_{t>0}$ with limit $\mu^*$ is said to converge \textit{linearly} if there exists a constant \textcolor{mylightbluetext}{$c\in (0,1)$} such that 
$$
\abs{\mu^{(t+1)} - \mu^*} \leq c \abs{\mu^{(t)} - \mu^*}.
$$
In other words, the \textit{linearly convergent sequence} has the following property:
$$
\mathop{\lim}_{t\rightarrow \infty}  \frac{\abs{\mu^{(t+1)} - \mu^*}}{\abs{\mu^{(t)} - \mu^*}} = c \in (0,1).
$$
\end{definition}
For example, the sequence $\mu^{(t)} = 4+ (1/4)^t$ converges linearly to $\mu^* = 4$ since 
$$
\mathop{\lim}_{t\rightarrow \infty}  \frac{\abs{\mu^{(t+1)} - \mu^*}}{\abs{\mu^{(t)} - \mu^*}} = \frac{1}{4} \in (0,1).
$$
\begin{definition}[Superlinear Convergence]\label{definition:superlinear_convergence}
A sequence $\{\mu^{(t)}\}_{t>0}$ with limit $\mu^*$ is said to converge \textit{superlinearly} if there exists a constant \textcolor{mylightbluetext}{$c_t >0$ with $c_t \rightarrow 0$}  such that 
$$
\abs{\mu^{(t+1)} - \mu^*} \leq c_t \abs{\mu^{(t)} - \mu^*}.
$$
In other words, the \textit{superlinearly convergent sequence} has the following property:
$$
\mathop{\lim}_{t\rightarrow \infty}  \frac{\abs{\mu^{(t+1)} - \mu^*}}{\abs{\mu^{(t)} - \mu^*}} =0.
$$
\end{definition}
For example, the sequence $\mu^{(t)} = 4+\left(\frac{1}{t+4}\right)^{t+3}$ converges superlinearly to $\mu^* = 4$ since 
$$
\mathop{\lim}_{t\rightarrow \infty}  \frac{\abs{\mu^{(t+1)} - \mu^*}}{\abs{\mu^{(t)} - \mu^*}}= 
\left(\frac{k+4}{k+5}\right)^{k+3}  \frac{1}{k+5}
= 0.
$$

\begin{definition}[Quadratic Convergence]\label{definition:quadratic-convergence}
A sequence $\{\mu^{(t)}\}_{t>0}$ with limit $\mu^*$ is said to converge \textit{quadratically} if there exists a constant \textcolor{mylightbluetext}{$c>0$} such that 
$$
\abs{\mu^{(t+1)} - \mu^*} \leq c \abs{\mu^{(t)} - \mu^*}^2.
$$
In other words, the \textit{quadratically convergent sequence} has the following property:
$$
\mathop{\lim}_{t\rightarrow \infty}  \frac{\abs{\mu^{(t+1)} - \mu^*}}{\abs{\mu^{(t)} - \mu^*}^2} = c .
$$
\end{definition}
For example, the sequence $\mu^{(t)} = 4+ (1/4)^{2^t}$ converges quaeratically to $\mu^* = 4$ since 
$$
\mathop{\lim}_{t\rightarrow \infty}  \frac{\abs{\mu^{(t+1)} - \mu^*}}{\abs{\mu^{(t)} - \mu^*}^2} = 1.
$$

\begin{definition}[Cubic Convergence]
A sequence $\{\mu^{(t)}\}_{t>0}$ with limit $\mu^*$ is said to converge \textit{cubically} if there exists a constant \textcolor{mylightbluetext}{$c>0$} such that 
$$
\abs{\mu^{(t+1)} - \mu^*} \leq c \abs{\mu^{(t)} - \mu^*}^3.
$$
In other words, the \textit{cubically convergent sequence} has the following property:
$$
\mathop{\lim}_{t\rightarrow \infty}  \frac{\abs{\mu^{(t+1)} - \mu^*}}{\abs{\mu^{(t)} - \mu^*}^3} = c .
$$
\end{definition}
For example, the sequence $\mu^{(t)} = 4+ (1/4)^{3^t}$ converges cubically to $\mu^* = 4$ since 
$$
\mathop{\lim}_{t\rightarrow \infty}  \frac{\abs{\mu^{(t+1)} - \mu^*}}{\abs{\mu^{(t)} - \mu^*}^3} = 1.
$$

\begin{problemset}
\item Prove that the sum or maximum  of two norms is also a norm.
\item Prove that the $\ell_1$ and $\ell_\infty$ norms cannot be derived from an inner product.
\textit{Hint: Consider the standard basis vectors $\be_1$ and $\be_2$, and examine the parallelogram identity in \eqref{equation:norm_inn_pro} or the polarization identity in \eqref{equation:norm_pol_pro}.}
\item \textbf{$k$-Norm.}
Show that the $k$-norm $\norm{\bx}^{[k]}$ on $\bx\in\real^m$, defined as the sum of the largest $k$ absolute elements of a vector, is a valid norm. 
Specifically, verify that $\norm{\bx}_\infty = \norm{\bx}^{[1]}$ and $\norm{\bx}^{[n]}=\norm{\bx}_1$. 
Additionally, prove the following sequence of inequalities:
\begin{equation}
\norm{\bx}_\infty = \norm{\bx}^{[1]} 
\leq 
\norm{\bx}^{[2]}
\leq\ldots 
\leq 
\norm{\bx}^{[n]}=\norm{\bx}_1.
\end{equation}

\item \label{problem:cond_pd} Given the Cholesky decomposition of a PD matrix: $\bA=\bL\bD\bL^\top$, show that $\cond(\bA)\geq \cond(\bD)$.

\item In Proposition~\ref{proposition:orthogonal-eigenvectors}, we demonstrated that  eigenvectors corresponding to distinct eigenvalues  of a symmetric matrix are orthogonal. More generally, prove that the eigenvectors corresponding to distinct eigenvalues of any matrix are linearly independent.

\item Show that if $\bA$ is positive semidefinite and nonsingular, then $\bA^{-1}$ is positive definite.

\item \textbf{Quadratic form.} Consider   the quadratic form $L(\bx) = \frac{1}{2} \bx^\top \bA \bx - \bb^\top \bx + c$ in Equation~\eqref{equation:quadratic-form-general-form1}, where $\bA\in\real^{n\times n}$, $ \bx\in \real^n$, and $c\in\real$. Suppose $\bA$ is positive semidefinite. 
Show that $L(\bx)$ is bounded below over $\real^n$ if and only if $\bb$ is in the column space of $\bA$.

\item \textbf{Quadratic form.} Consider  again  the quadratic form $L(\bx) = \frac{1}{2} \bx^\top \bA \bx - \bb^\top \bx + c$ in Equation~\eqref{equation:quadratic-form-general-form1}. Show that $L(\bx)$ is coercive~\footnote{A function $f(\bx):\real^n\rightarrow \real$ is called coercive if $\mathop{\lim}_{\bx\rightarrow \infty} f(x)=\infty$.} if and only if $\bA$ is PD.

\item  \label{problem:symm_square} \textbf{Square of symmetric matrices.} Let $\bA\in\real^{n\times n}$ be symmetric. Show that the eigenvalues of $\bA^2$ (can only be nonnegative) are the squares  of the eigenvalues of $\bA\in\real^{n\times n}$. What is the relationship between the eigenvectors?

\item \label{problem:part_ortho} Consider the partition of an orthogonal matrix
$
\bQ = 
\scriptsize
\begin{bmatrix}
	\underset{p\times p}{\bA} &\underset{p\times q}{\bB} \\
	\underset{q\times p}{\bC} & \underset{q\times q}{\bD}
\end{bmatrix}\in\real^{n\times n}.
$
Show that $\rank(\bB)=\rank(\bC)$ and $\rank(\bD)=n+\rank(\bA)-2p$.
\end{problemset}

\newpage
\chapter{Optimization and Optimality Conditions}\label{chapter:opt_cond}
\begingroup
\hypersetup{
linkcolor=structurecolor,
linktoc=page,  
}
\minitoc \newpage
\endgroup

\section{Classes of Sets and Functions}

We  briefly introduce different notions of sets and functions.

\subsection{Interior Points and Closed Sets}

Given a specific norm definition, we introduce the concepts of an open ball and a closed ball as follows:
\begin{definition}[Open Ball, Closed Ball]\label{definition:open_closed_ball}
Let $\norm{\cdot}_p: \real^n\rightarrow \real_+$ be the $\ell_p$ norm function. The \textit{open ball} centered at $\bc\in\real^n$ with radius $r$  is defined as 
$$
\sB_p(\bc, r) \triangleq \{\bx\in\real^n\mid  \norm{\bx-\bc}_p <r\}.
$$
Similarly, the \textit{closed ball} centered at $\bc\in\real^n$ with radius $r$  is defined as 
$$
\sB_p[\bc,r] \triangleq \{\bx\in\real^n\mid \norm{\bx-\bc}_p \leq r\}.
$$
For example, $\sB_2[\bzero,1]$ represents  the unit closed ball w.r.t. to  the $\ell_2$ norm.
To simplify notation, we omit the subscript 2 for  $\ell_2$ norms and  $\bzero$ for balls centered at zero, e.g., $ \sB[1] \triangleq \sB_2[\bzero,1]  $.
As a special case, the notation $\sB_0[k] \triangleq\sB_0[\bzero, k]$ denotes the set of \textit{$k$-sparse vectors}, i.e., containing vectors that have only $k$ (or less) nonzero elements.
More generally, let $\norm{\cdot}$ be any norm, the induced open and closed balls are denoted as 
$$
\sB_{\norm{\cdot}} (\bc, r)
\qquad \text{and}\qquad 
\sB_{\norm{\cdot}} [\bc, r].
$$
\end{definition}

Subsequently, we define the concepts of interior and relative interior points. These concepts are foundational in topology and convex analysis and are based on the definitions of balls.
\index{Interior points}
\index{Relative interior points}
\begin{definition}[Interior and Relative Interior Points]
The \textit{interior} of a set $ \sS $ in a topological space is the set of all points in $ \sS $ that have a neighborhood completely contained within $ \sS $. Formally, the interior of $ \sS $, denoted by $ \interior(\sS) $, is given by:
$$
\interior(\sS) = \{ \bx \in \sS \mid \exists \epsilon > 0 \text{ such that } \sB(\bx, \epsilon) \subseteq \sS \}.
$$
where $ \sB(\bx, \epsilon) $ is an open ball centered at $ \bx $ with radius $ \epsilon $.

The \textit{relative interior} of a set $ \sS $ in a topological space, typically within the context of convex analysis, is the interior of $ \sS $ relative to its affine hull. The \textit{affine hull} of $ \sS $, denoted by $ \aff(\sS) $, is the smallest affine space containing $ \sS $ (Definition~\ref{definition:coms_hulls}). The relative interior of $ \sS $, denoted by $ \relint(\sS) $, is given by:
$$
\relint(\sS) = \{ \bx \in \sS \mid \exists \epsilon > 0 \text{ such that } \sB_{\aff(\sS)}(\bx, \epsilon) \subseteq \sS \}.
$$
where $ \sB_{\aff(\sS)}(\bx, \epsilon) $ is the open ball centered at $ \bx $ with radius $ \epsilon $ in the subspace topology of $ \aff(\sS) $.

In other words, the \textit{interior} of a set is concerned with the open balls in the ambient space; the \textit{relative interior} is concerned with open balls in the affine hull of the set, which is particularly important in higher dimensions and in convex analysis where the set may lie in a lower-dimensional subspace.
\end{definition}

\begin{example}[Interior and Relative Interior Points]
\textbf{Interval in $\real$.} Let $ \sS = [0, 1] $. Then $ \interior(\sS) = (0, 1) $ and $ \relint(\sS) = (0, 1) $. Here, the interior and relative interior are the same because the set $ [0, 1] $ is a subset of $\real$ and $\aff(\sS) = \real$.

\paragraph{Line Segment in $\real^2$.} Let $ \sS = \{ (x, y) \in \real^2 \mid 0 \leq x \leq 1, y = 0 \} $ with
$\aff(\sS)=\{(x,0)\in\real^2 \mid x\in\real\}$. Then  $ \interior(\sS) = \varnothing $ and $ \relint(\sS) = \{ (x, 0) \mid 0 < x < 1 \} $.
The interior is empty because there are no open balls in $\real^2$ entirely contained within $ \sS $. However, the relative interior is the open interval (0, 1) in the context of the line segment (the affine hull is the line $ y = 0 $).

\paragraph{Triangle in $\real^2$.} Let $ \sS $ be a triangle with vertices at $(0,0)$, $(1,0)$, and $(0,1)$. Then the interior consists of all points inside the triangle, excluding the boundary; while the  relative interior is the same as the interior in this case because the affine hull of the triangle is the entire plane $ \real^2 $.
\end{example}

\index{Open sets}
\index{Closed sets}
\index{Compact sets}
\index{Level sets}
\begin{definition}[Open, Closed, Compact, and Level Sets]\label{definition:open_close_sets}
A set $\sS_1$ is said to be \textit{open} if it consists of only interior points. That is, for every $\bx\in\sS$, there exists a scalar $s>0$ such that $\sB(\bx, s)\subseteq\sS_1$.

A set $\sS_2$ is said to be \textit{closed} if its complement $\comple{\sS_2}$ is open.
Alternatively, $\sS_2$ is closed if it contains all the limit points of convergent sequence of points in $\sS_2$; that is, for every sequence of points $\{\bx_i\}_{i\geq 1}\subseteq \sS_2$ satisfying $\bx_i\rightarrow \bx^\star$ as $i\rightarrow \infty$, it follows that $\bx^\star\in\sS_2$.
In the meantime, the \textit{closure} of a set $\sS$, denoted by $\closure(\sS)$, is defined as the smallest closed set containing $\sS$:
$$
\closure(\sS) = \cap\left\{ \sT\mid  \sS \subseteq \sT, \sT \text{ is closed} \right\}.
$$

A set $\sS \subseteq \real^n$ is called \textit{bounded} if there exists $M > 0$ for which $\sS \subseteq \sB(\bzero, M)$, where $\sB(\bzero, M)$ is the ball of finite radius $M$ centered at the origin. 
A set $\sS \subseteq \real^n$ is called \textit{compact} if it is closed and bounded.

Given a function $ f: \sS \rightarrow \real $ where $ \sS \subseteq \real^n $, the \textit{level set} of $ f $ at a particular value $ c \in \real $ is defined as:
$$
\lev[f, c] = \{ \bx \in \sS \mid f(\bx) \leq c \}.~\footnote{Note that we use square brackets instead of parentheses to indicate that the equality can be obtained; same as the definition of closed balls in Definition~\ref{definition:open_closed_ball}.}
$$
In other words, the level set $ \lev[f, c] $ consists of all points $ \bx $ in the domain $ \sS $ for which the function $ f $ evaluates to smaller than or equal to  the constant $ c $.
\end{definition} 

\subsection{Continuous Functions}

A continuous function is one that does not exhibit abrupt changes in value---referred to as \textit{discontinuities}. More formally, it means that small changes in the input of the function lead to correspondingly small changes in its output. Continuous functions are foundational in calculus and analysis, boasting several important properties.
\begin{definition}[Continuous Functions]\label{definition:conti_funs}
Let $f:\sS \subseteq \real^n\rightarrow \real$. The function $f$ is said to be \textit{continuous} at a point $\by \in \sS$ if for every $\epsilon > 0$, there exists a $\delta > 0$ such that for all $\bx$ in $\sS$:
\begin{equation}
\forall \epsilon > 0, \exists \delta > 0 :\quad \forall \bx \in \sS,
\normtwo{\bx - \by} < \delta 
\quad\implies \quad
\abs{f(\bx) - f(\by)} < \epsilon.
\end{equation}
In simpler terms, this definition indicates that we can make the values of $f(\bx)$ arbitrarily close to $f(\by)$ by selecting  $\bx$ sufficiently close to $\by$.

If $f$ is continuous at every point in its domain $\sS$, then $f$ is said to be \textit{continuous} on $\sS$ or simply continuous.
\end{definition}

\begin{remark}[Properties of Continuous Functions]
For a continuous function, we have the following properties:
\begin{enumerate}
\item \textit{Intermediate value theorem.} If $f$ is continuous on a closed interval $[a, b]$ and $y$ is any number between $f(a)$ and $f(b)$, then there exists a $c \in [a, b]$ such that $f(c) = y$.
\item \textit{Extreme value theorem.} If $f$ is continuous on a closed interval $[a, b]$, then $f$ attains both a maximum and minimum value on that interval (Theorem~\ref{theorem:weierstrass_them}).
\item \textit{Composition of continuous functions.} If $g$ is continuous at $c$ and $f$ is continuous at $g(c)$, then the composition $f \circ g$ is continuous at $c$.
\item \textit{Arithmetic operations.} If $f$ and $g$ are continuous at $c$, then so are $f + g$, $f - g$, $f g$, and $\frac{f}{g}$ (provided $g(c) \neq 0$).
\item \textit{Continuity and limits.} A function $f$ is continuous at $c$ if and only if $\lim_{x \to c} f(x) = f(c)$.
\end{enumerate}
\end{remark}

\index{Continuous}
\begin{example}[Continuous Functions]
Following are some continuous functions:
\begin{itemize}
\item Polynomial functions are continuous everywhere.
\item Rational functions are continuous wherever they are defined (i.e., where the denominator is not zero).
\item Trigonometric functions like sine and cosine are continuous everywhere.
\item Exponential and logarithmic functions are continuous on their domains.
\item Absolute value function $\abs{x}$ is continuous everywhere.
\end{itemize}
There are different types of discontinuities:
\begin{itemize}
\item \textit{Removable discontinuity.} The function has a hole at a point, but it can be ``filled" to become continuous.
\item \textit{Jump discontinuity.} The left-hand limit and right-hand limit exist but are not equal.
\item \textit{Infinite discontinuity.} The function approaches  positive or negative infinity at a point.
\item \textit{Oscillating discontinuity.} The function does not approach a single value as the input approaches a point.
\end{itemize}
Continuous functions are easier to analyze and manipulate, while discontinuities can lead to complex behaviors that require special handling. Recognizing the type of discontinuity helps in determining appropriate strategies for dealing with these functions in various applications.
\end{example}

\begin{definition}[Uniform Continuity]\label{definition:uniform_cont}
A function $ f:\sS \subseteq \real^n\rightarrow \real $ is said to be \textit{uniformly continuous} on a set $ \sS $ if for every $ \epsilon > 0 $, there exists a $ \delta > 0 $ such that for all $ \bx $ and $ \by $ in $ \sS $, whenever the distance between $ \bx $ and $ \by $ is less than $ \delta $ (that is, $ \abs{\bx - \by} < \delta $), it follows that the distance between $ f(\bx) $ and $ f(\by) $ is less than $ \epsilon $ (that is, $ \abs{f(\bx) - f(\by)} < \epsilon $).
Formally, $ f: \sS \rightarrow \real $ is uniformly continuous on $ \sS $ if:
$$
\forall \epsilon > 0, \exists \delta > 0 :\quad \forall \bx, \textcolor{mylightbluetext}{\by} \in \sS, \normtwo{\bx - \by} < \delta 
\quad\implies\quad
\abs{f(\bx) - f(\by)} < \epsilon.
$$

The key difference between uniform continuity and ordinary (pointwise) continuity is the choice of $ \delta $. For pointwise continuity, $ \delta $ can depend on both $ \epsilon $ and the point $ \by $ in the domain. However, for uniform continuity, $ \delta $ must work for all points in the set $ \sS $ simultaneously, given an $ \epsilon $.
\end{definition}
Uniform continuity is a stronger form of continuity. Every uniformly continuous function is continuous, but not every continuous function is uniformly continuous. Uniform continuity is especially important in analysis because it guarantees certain desirable  properties, like the preservation of Cauchy sequences and the ability to interchange limits and function evaluations under suitable conditions.

\index{Uniformly continuous}
\begin{example}[Continuity vs Uniform Continuity]
Consider $ f(x) = x^2 $. For any given $ y \in \real $, we can always find a $\delta$ depending on $ y $ and $\epsilon$ to satisfy the definition of continuity. Hence, $ f(x) = x^2 $ is continuous at every point $ y \in \real $.

\paragraph{Not uniformly continuous.} To show that $ f(x) = x^2 $ is not uniformly continuous on $\real$, consider two sequences $ x_t = t $ and $ y_t = t + \frac{1}{t} $ where $ t $ is an integer. As $ t \to \infty $, the distance between $ x_t $ and $ y_t $ approaches  0 ($ \abs{x_t - y_t} = \frac{1}{t} \to 0 $), but the difference in function values does not approach  0:
$$
\abs{f(x_t) - f(y_t)} = \abs{t^2 - \left(t + \frac{1}{t}\right)^2} = \abs{t^2 - \left(t^2 + 2 + \frac{1}{t^2}\right)} = \abs{-2 - \frac{1}{t^2} } \to 2
$$
Since the difference in function values does not become arbitrarily small even when the inputs are arbitrarily close, $ f(x) = x^2 $ is not uniformly continuous on $\real$.
This shows that while $ f(x) = x^2 $ is continuous everywhere on the real line, it fails to be uniformly continuous because the required $\delta$ in the definition of uniform continuity cannot be chosen independently of $ x $.
\end{example}
\begin{remark}[Consequence of Uniform Continuity]\label{remark:conse_uniform_cont}
Let $\nabla f$ be uniformly continuous. Then,
$$
\innerproduct{\nabla f(\bx_1), \bd} = \innerproduct{\nabla f(\bx_2), \bd} + o(\normtwo{\bd}), \quad \text{where }\bd\triangleq\bx_1-\bx_2.
$$
To see this,
by the definition of uniform continuity, for any $\epsilon > 0$, there exists a $\delta > 0$ such that:
$$
\normtwo{\bd} < \delta \quad\implies\quad\normtwo{\nabla f(\bx_1) - \nabla f(\bx_2)} < \epsilon.
$$
This implies that as $\normtwo{\bd} \to 0$, we can choose $\delta$ arbitrarily small, and the gradient difference satisfies:
$
\normtwo{\nabla f(\bx_1) - \nabla f(\bx_2)} = o(1).
$
Using the linearity of the inner product, we can write:
$$
\innerproduct{\nabla f(\bx_1), \bd} = \innerproduct{\nabla f(\bx_2), \bd} + \innerproduct{\nabla f(\bx_1) - \nabla f(\bx_2), \bd},
$$
where
$
\innerproduct{\nabla f(\bx_1) - \nabla f(\bx_2),  \bd} = \normtwo{\nabla f(\bx_1) - \nabla f(\bx_2)} \cdot \normtwo{\bd} \cdot \cos(\theta),
$
and $\theta$ is the angle between $\nabla f(\bx_1) - \nabla f(\bx_2)$ and $\bd$.
Since $\normtwo{\nabla f(\bx_1) - \nabla f(\bx_2)} = o(1)$, we have:
$
\innerproduct{\nabla f(\bx_1) - \nabla f(\bx_2), \bd} = o(\normtwo{\bd}).
$
Therefore, combining these results, we get the desired result.
\end{remark}

\index{Lower semicontinuity}
\index{Lower semicontinuous}
\begin{definition}[Lower Semicontinuity]
A function $ f: \real^n \rightarrow \real\cup \{-\infty, \infty\} $ is called \textit{lower semicontinuous at $ \bx^* \in \real^n $} if
$$
f(\bx) \leq \liminf_{t \rightarrow \infty} f(\bx^\toptzero)
$$
for any sequence $ \{\bx^\toptzero\}_{t \geq 1} \subseteq \real^n $ for which $ \bx^\toptzero \rightarrow \bx^* $ as $ t \rightarrow \infty $.
In other words, the limit inferior of the function values at points approaching $\bx^*$ must be greater than or equal to the function value at $\bx^*$ . This means that the function does not have any sudden drops as you approach $\bx^*$.
A function $ f: \real^n \rightarrow  \real \cup \{-\infty, \infty \} $ is called \textit{lower semicontinuous} if it is {lower semicontinuous} at each point in $ \real^n $.
\end{definition}

Note that including $\{-\infty, \infty\}$ in the codomain of a function when defining lower semicontinuityis not strictly necessary for all functions, but it allows for greater generality and flexibility in mathematical analysis. 
Some functions may have points where they tend towards infinity, either positively or negatively. By including $\{-\infty, \infty\}$, we can still discuss the continuity properties of such functions at those points.

\begin{figure}[h!]
\centering                      
\vspace{-0.35cm}                
\subfigtopskip=2pt            
\subfigbottomskip=2pt           
\subfigcapskip=-5pt           
\subfigure[A lower semicontinuous funtion.]{\label{fig:aaa-1}
	\includegraphics[width=0.31\linewidth]{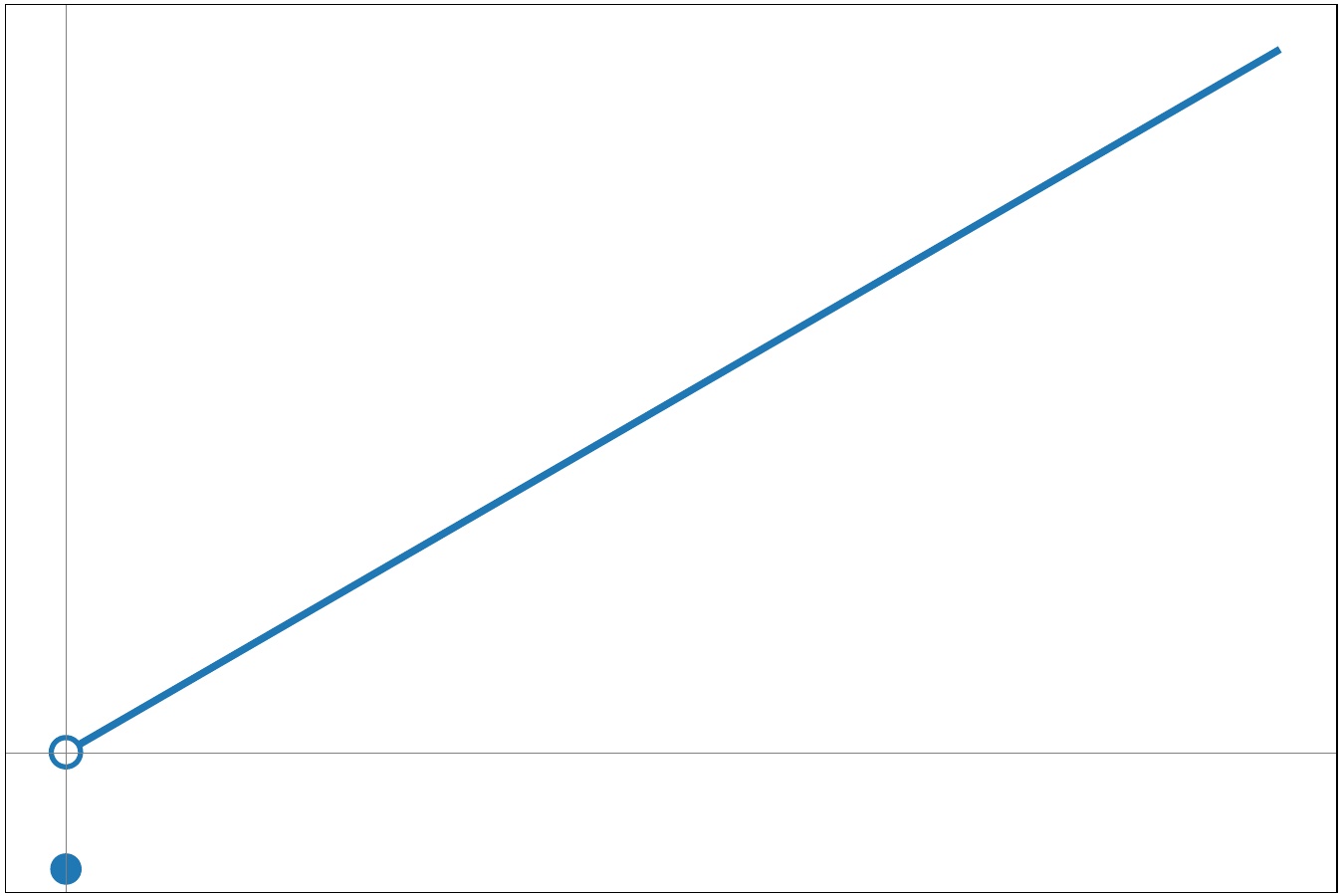}}
\subfigure[A non-lower semicontinuous funtion]{\label{fig:aaa-2}
	\includegraphics[width=0.31\linewidth]{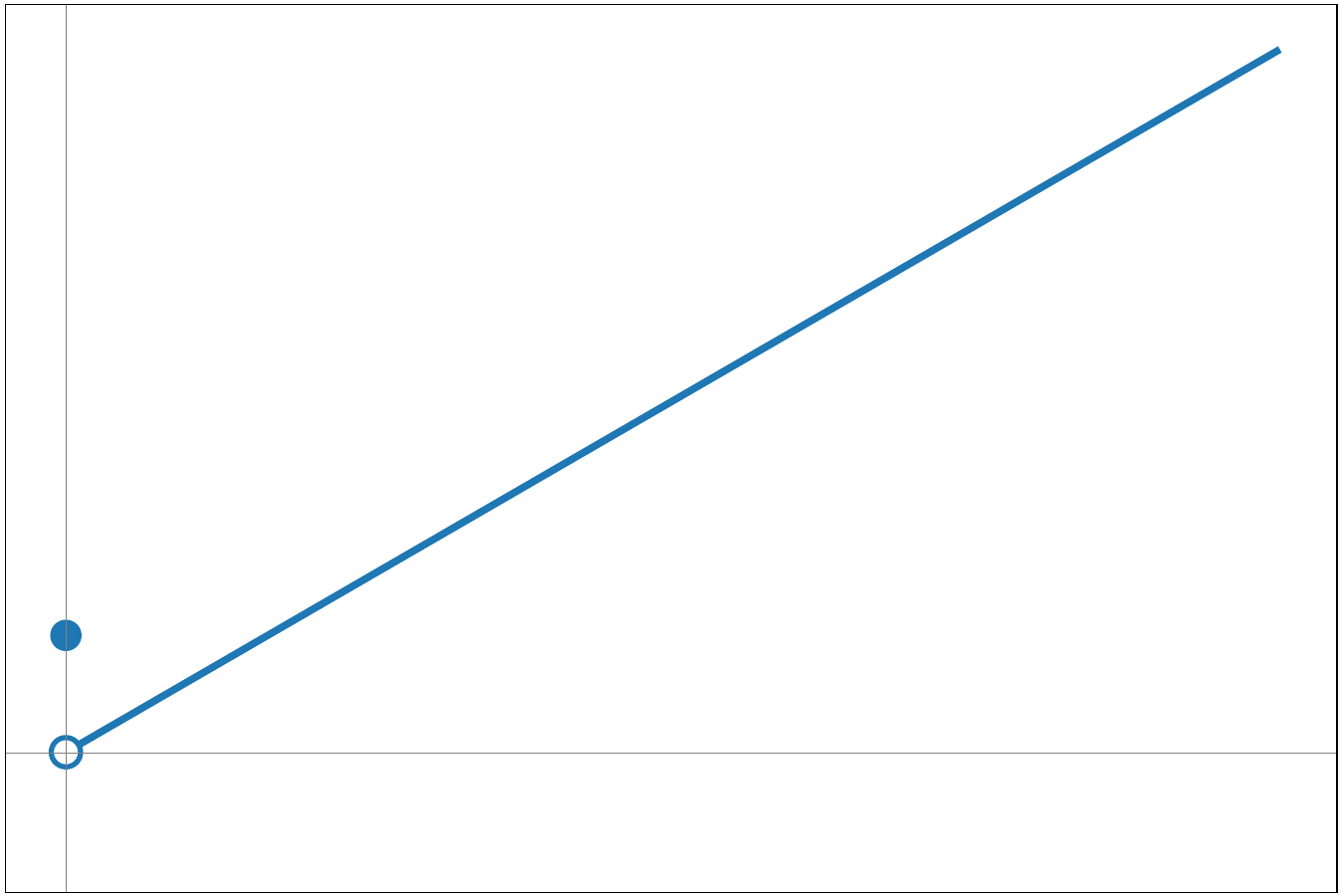}}
\subfigure[A lower semicontinuous funtion.]{\label{fig:aaa-3}
	\includegraphics[width=0.31\linewidth]{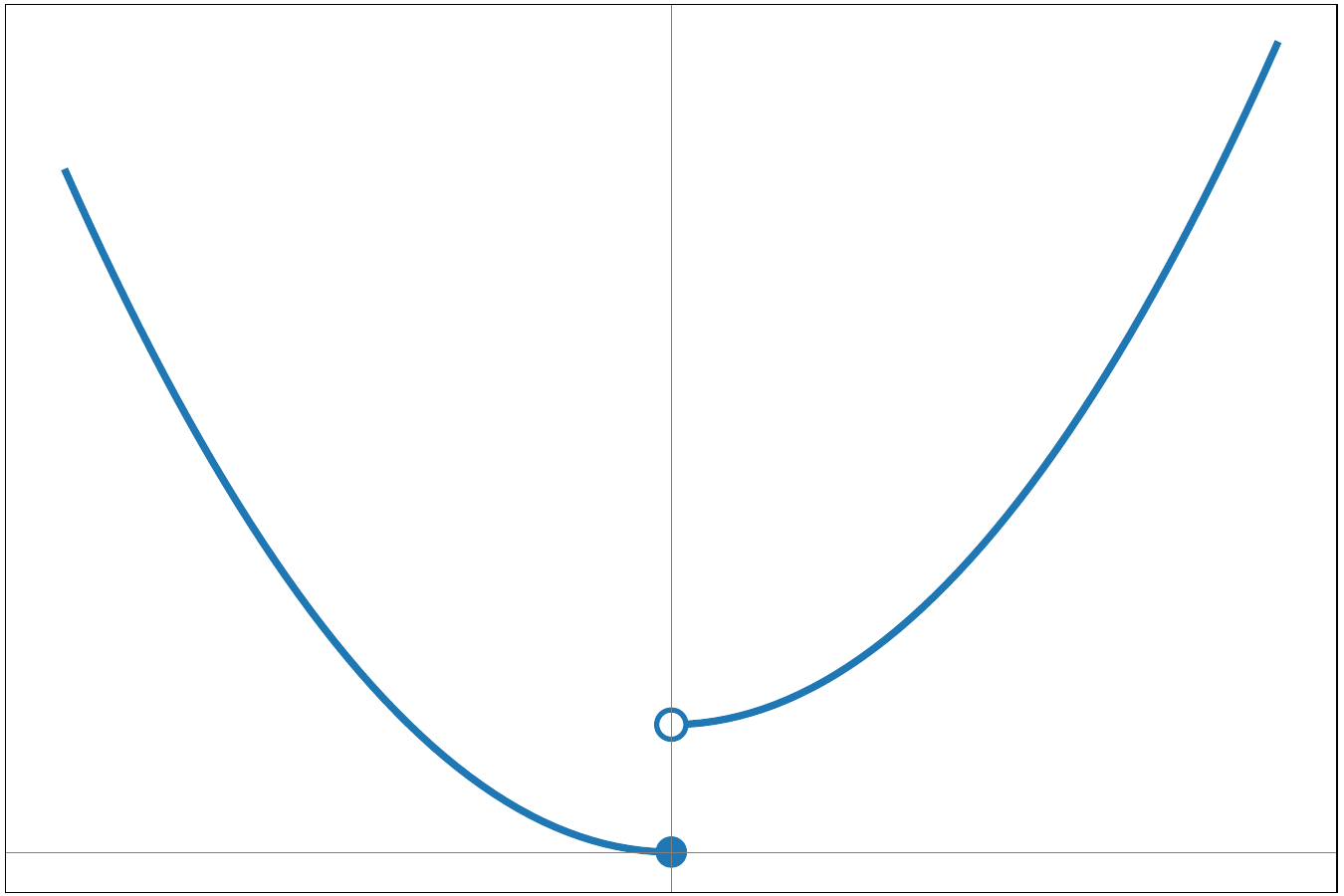}}
\caption{Lower semicontinuous and non-lower semicontinous funtions. In numerical methods and iterative algorithms, lower semicontinuity can help in proving convergence properties. For instance, if a sequence generated by an algorithm is minimizing a lower semicontinuous function, it can be shown that the limit points of the sequence are stationary points.}
\label{fig:low_nonlowersemis}
\end{figure}

Figure~\ref{fig:low_nonlowersemis} depicts examples of lower semicontinuous and non-lower semicontinuous functions.
In optimization problems, lower semicontinuity plays a key role because it guarantees that if a sequence of feasible points approaches a minimum, then the limit of this sequence will also be a feasible point and will achieve the minimum value. This property ensures that optimal solutions exist under certain conditions.
In many applications, especially in calculus of variations and optimal control, lower semicontinuity is used to prove the existence of minimizers of functionals. If a functional is lower semicontinuous and coercive (it tends to infinity as the norm of its argument grows), then under some compactness assumptions, one can show that there exists a minimizer; see \ref{weier2_prop_close} in Theorem~\ref{theorem:weierstrass_them}.
In numerical methods and iterative algorithms, lower semicontinuity can help in proving convergence properties. For instance, if a sequence generated by an algorithm is minimizing a lower semicontinuous function, it can be shown that the limit points of the sequence are stationary points; for example, the lower semicontinuity can be applied to prove the uniqueness of the optimizer for a closed and strongly convex function (Theorem~\ref{theorem:exi_close_sc}).
For this reason, in most of our discussions, we consider lower semicontinuous (closed) functions for evaluation.

\subsection{Extended Real-Valued Functions}
Extended real-valued functions are those that can take values from the extended real number line, which includes all real numbers along with two additional elements: positive infinity ($+\infty$) and negative infinity ($-\infty$). These functions are widely  used in mathematical analysis, optimization theory, and measure theory, providing a convenient framework to handle limits and boundary conditions.

\index{Extended real-valued functions}
\index{Proper functions}
\index{Closed functions}
\index{Epigraph}
\begin{definition}[Extended Real-Valued, Proper, Closed Functions]\label{definition:extrel_pro_clo_funcs}
An \textit{extended real-valued function} $f$ is a function defined on a set $\sS$ (often a subset of $\real^n$ or another topological space) that maps each point $\bx$ in $\sS$ to a value in the \textit{extended real number line} $\extreal = \real \cup \{-\infty, \infty \}$. Formally, an extended real-valued function is written as:
$$
f : \sS \rightarrow \extreal
\qquad \text{or}\qquad 
f : \sS \rightarrow \real \cup \{-\infty, \infty \}.
$$
The \textit{effective domain} or simply the \textit{domain} and the \textit{epigraph} of an extended real-valued function are the sets:
$$
\begin{aligned}
\textbf{(domain)}:\qquad &\dom(f) &\triangleq& \{\bx\in\sS: f(\bx)<\infty\};\\
\textbf{(epigraph)}:\qquad & \epi(f) &\triangleq& \{(\bx, r) \in \sS \times \real : f(\bx) \leq r\}.
\end{aligned}
$$
That is, the epigraph of an extended real-valued function $f$ is the set of points lying on or above its graph.

In the meantime,  a function $f: \sS \rightarrow \extreal$ is called \textit{proper} if it satisfies the following two conditions:
\begin{enumerate}
\item \textit{Non-identically infinite.} The function does not attain the value $+\infty$ everywhere in its domain. In other words, there exists at least one point $\bx \in \sS$ such that $f(\bx) < +\infty$.
\item \textit{No negative infinity values.} The function never takes on the value $-\infty$. That is, for all $\bx \in \sS, f(\bx) > -\infty$.
\end{enumerate}
\noindent These conditions ensure that the function has some finite values and avoids the undefined behavior associated with $-\infty$.
For a proper function, this domain is nonempty because of the first condition above.

A function $f$ is said to be \textit{closed} if its epigraph is a closed set in the product topology of $\sS \times \real$. In other words, for every sequence $(\bx^{(t)}, r_t)$ in $\text{epi}(f)$ that converges to a limit $(\bx^*, r^*)$, we have $f(\bx^*) \leq r^*$, meaning that the limit point $(\bx^*, r^*)$ also belongs to the epigraph.
\end{definition}

\begin{definition}[Coerciveness]\label{definition:coerciveness}
Let $f: \real^n \rightarrow \real$ be a continuous function defined over $\real^n$. The function $f$ is called \textit{coercive} if
$$
\lim_{\normtwo{\bx} \rightarrow \infty} f(\bx) = \infty.
$$
\end{definition}

\begin{exercise}[Coerciveness of Quadratic Functions]\label{exercise:coerci_quad}
	Let $f(\bx)=\frac{1}{2}\bx^\top\bA\bx+\bb^\top\bx+c$, where $\bA\in\real^{n\times n}$ is symmetric, $\bb\in\real^n$, and $c\in\real$. Show that $f(\bx)$ is coercive if $\bA$ is positive definite.
\end{exercise}

Proper functions are particularly significant  in convex analysis and optimization theory because they allow for a well-defined notion of minimization. When dealing with optimization problems, one often seeks to minimize a function over a given set. A proper function ensures that there are feasible points to consider (points where the function value is finite), which is necessary for the problem to be meaningful. Therefore, unless otherwise stated, we only consider proper functions in this book.
Moreover, many results in convex analysis require the objective function to be proper. For example, the Fenchel-Moreau theorem, which relates a function to its biconjugate, requires the function to be proper, convex, and lower semi-continuous.

\begin{example}[Proper Functions]\label{example:proper_funcs}
	The following functions are proper:
	\begin{itemize}
		\item A linear function $f(x) = ax + b$ on $\real$ is proper since it never reaches $\pm\infty$.
		\item The indicator function $\indicatorS(\bx)$, which is $0$ if $\bx$ belongs to a set $\sS$ and $+\infty$ otherwise, is proper if $\sS$ is nonempty.
		\item The function $f(x) = \frac{1}{x}$ on $\real$ is not proper since it tends to $+\infty$ as $x$ approaches $0$ from the positive side and to $-\infty$ from the negative side.
	\end{itemize}
Consider why the non-proper function $f(x) = \frac{1}{x}$ is not well-defined for a minimization problem.
\end{example}

\index{Indicator function}
\begin{exercise}[Closedness of Indicator]\label{exercise_closed_indica}
Show that the indicator function $\indicatorS$, which is $0$ if $\bx$ belongs to a set $\sS$ and $+\infty$ otherwise,  is closed if and only if $\sS$ is a closed set.
\end{exercise}

\begin{theorem}[Equivalence of Closedness, Lower Semicontinuity, and Closedness of Level Sets]\label{theorem:equiv_close_clkos_semicon}
Let $ f: \real^n\rightarrow  \real \cup \{-\infty, \infty \} $. Then the following three statements  are equivalent:

\begin{enumerate}
\item[(i)] $ f $ is lower semicontinuous.
\item[(ii)] $ f $ is closed.
\item[(iii)] For any $ \alpha \in \real $, the level set
$
\lev[f, \alpha] = \{ \bx \in \real^n \mid f(\bx) \leq \alpha \}
$
is closed.
\end{enumerate}
\end{theorem}
\begin{proof}[of Theorem~\ref{theorem:equiv_close_clkos_semicon}]
\textbf{(i) $\implies$ (ii).} Assume  that $ f $ is lower semicontinuous. We will demonstrate  that $ \epi(f) $ is closed. 
Let  $ \{(\bx^\toptzero, y^\toptzero)\}_{t \geq 1} \subseteq \epi(f) $ be a sequence such that $ (\bx^\toptzero, y^\toptzero) \rightarrow (\bx^*, y^*) $ as $ t \rightarrow \infty $. By definition, for any $ t \geq 1 $,
$$
f(\bx^\toptzero) \leq y^\toptzero.
$$
Due to the lower semicontinuity of $ f $ at $ \bx^* $, we have
$$
f(\bx^*) \leq \liminf_{t \rightarrow \infty} f(\bx^\toptzero) \leq \liminf_{t \rightarrow \infty} y^\toptzero = y^*,
$$
which shows that $ (\bx^*, y^*) \in \epi(f) $, proving that $f$ is closed.

\paragraph{(ii) $\implies$ (iii).} Suppose that $ f $ is closed, i.e., the epigraph $ \epi(f) $ is closed. 
For any $ \alpha \in \real $, we aim to show that  $ \lev[f, \alpha] $ is closed. If $ \lev[f, \alpha] = \varnothing $, the claim holds trivially. 
Otherwise, take a sequence $ \{\bx^\toptzero\}_{t \geq 1} \subseteq \lev[f, \alpha] $ that converges to $ \widehatbx $. Obviously, $ (\bx^\toptzero, \alpha) \in \epi(f) $ for any $ t $ and $ (\bx^\toptzero, \alpha) \rightarrow (\widehatbx, \alpha) $ as $ t \rightarrow \infty $. By the closedness of $ \epi(f) $, it follows that $ (\widehatbx, \alpha) \in \epi(f) $, establishing the fact that $ \widehatbx \in \lev[f, \alpha] $.

\paragraph{(iii) $\implies$ (i).} Assuming all level sets of $ f $ are closed, we prove that it is lower semicontinuous. Suppose for contradiction that $ f $ is not lower semicontinuous, meaning that there exists $ \bx^* \in \real^n $ and $ \{\bx^\toptzero\}_{t \geq 1} \subseteq \real^n $ such that $ \bx^\toptzero \rightarrow \bx^* $ and $ \liminf_{t \rightarrow \infty} f(\bx^\toptzero) < f(\bx^*) $. 
Choose  $ \alpha $ satisfying
\begin{equation}\label{equation:equiv_lowsemi}
\liminf_{t \rightarrow \infty} f(\bx^\toptzero) < \alpha < f(\bx^*). 
\end{equation}
Then there exists a subsequence $ \{\bx_{n_t}\}_{t \geq 1} $ such that $ f(\bx_{n_t}) \leq \alpha $ for all $ t \geq 1 $. By the closedness of the level set $ \lev[f, \alpha] $ and the fact that $ \bx_{n_t} \rightarrow \bx^* $ as $ t\rightarrow \infty $, it follows that $ f(\bx^*) \leq \alpha $, which is a contradiction to \eqref{equation:equiv_lowsemi}, establishing that (iii) implies (i).
\end{proof}

\subsection{Convex Functions}
This subsection introduces the concepts of affine, convex, and conic combinations and hulls within the context of vector sets in $\real^n$, providing definitions and mathematical formulations for these fundamental geometric constructs.
\begin{definition}[Affine/Convex/Conic Combinations and  Hulls]\label{definition:coms_hulls}
Given a set of vectors $\bx_1,\bx_2,\ldots,\bx_p\in\real^n$, 
an \textit{affine combination} of these $p$ vectors is a vector of the form  $\sum_{i=1}^{p} \lambda_i\bx_i$, where $\sum_{i=1}^{p} \lambda_i=1$, i.e., a linear combination of these points where the coefficients (not necessarily nonnegative) sum to 1; 
a \textit{convex combination} of these $p$ vectors is a vector of the form $\sum_{i=1}^{p} \lambda_i\bx_i$, where $\lambda_i\geq 0$ for $i\in\{1,2,\ldots,p\}$ satisfying $\sum_{i=1}^{p}\lambda_i=1$ (i.e., $\{\lambda_i\}$ belongs to the unit-simplex in $\real^p$); 
when only requiring $\lambda_i\geq 0$ without the constraint on their sum, the combination is referred to as a \textit{conic combination} of these vectors. 

The \textit{affine hull} of a set $ \sS $ is the smallest affine space (a translation of a vector subspace) that contains $ \sS $.
Formally, given a set $ \sS \subseteq \real^n $, the affine hull of $ \sS $, denoted as $ \aff(\sS) $, is defined as the set of all affine combinations of points from $ \sS $:
$$
\aff(\sS) \triangleq \left\{ \sum_{i=1}^p \lambda_i \bx_i \mid   \bx_1,\bx_2,\ldots,\bx_p\in\sS, \lambda_i \in \real, \bone^\top\blambda=1,p\in\naturalset \right\}.
$$
The \textit{convex hull} of $\sS$, denoted as $\conv(\sS)$, is the set comprising all the convex combinations of vectors from $\sS$: 
$$
\conv(\sS) \triangleq \left\{ \sum_{i=1}^{p} \lambda_i\bx_i\mid \bx_1,\bx_2,\ldots,\bx_p\in\sS, \blambda\in\Delta_p, p\in\naturalset \right\}.
$$ 
Similarly, the \textit{conic hull} of $\sS$, denoted as $\cone(\sS)$, is the set comprising all the conic combinations of vectors from $\sS$:
$$
\cone(\sS) \triangleq \left\{\sum_{i=1}^{p} \lambda_i\bx_i\mid \bx_1,\bx_2,\ldots,\bx_p\in\sS, \blambda\in\real^n_+, p\in\naturalset  \right\}.
$$
\end{definition}

\begin{exercise}[Closedness of Conic Hull of a Finite Set]\label{exercise:closed_fini_cone}
Let $\bx_1, \bx_2, \ldots,\bx_p\in\real^n$. Show that $\cone(\{\bx_1, \bx_2, \ldots,\bx_p\})$ is closed.
\end{exercise}

\begin{exercise}[Closure]
Show that 
$$
\begin{aligned}
\closure(\real_{++}^n) &= \real_+^n; \\
\closure(\sB(\bx,r)) &= \sB[\bx,r],\quad \text{where $\bx\in\real^n$, $r\in\real_{++}$};\\
\closure((\bx,\by)) &= [\bx,\by], \quad \text{where $\bx,\by\in\real^n$ and $\by\geq \bx$}.
\end{aligned}
$$
\end{exercise}

\begin{exercise}[Interior and Closure]\label{exercise:int_clos}
Let $\sS\subseteq\real^n$ be a convex set with a nonempty interior. Show that 
$$
\closure(\interior(\sS)) = \closure(\sS)
\qquad\text{and}\qquad 
\interior(\closure(\sS)) =\interior(\sS).
$$
\end{exercise}

\begin{figure}[h]
\centering       
\vspace{-0.25cm}                 
\subfigtopskip=2pt               
\subfigbottomskip=-2pt         
\subfigcapskip=-10pt      
\includegraphics[width=0.98\textwidth]{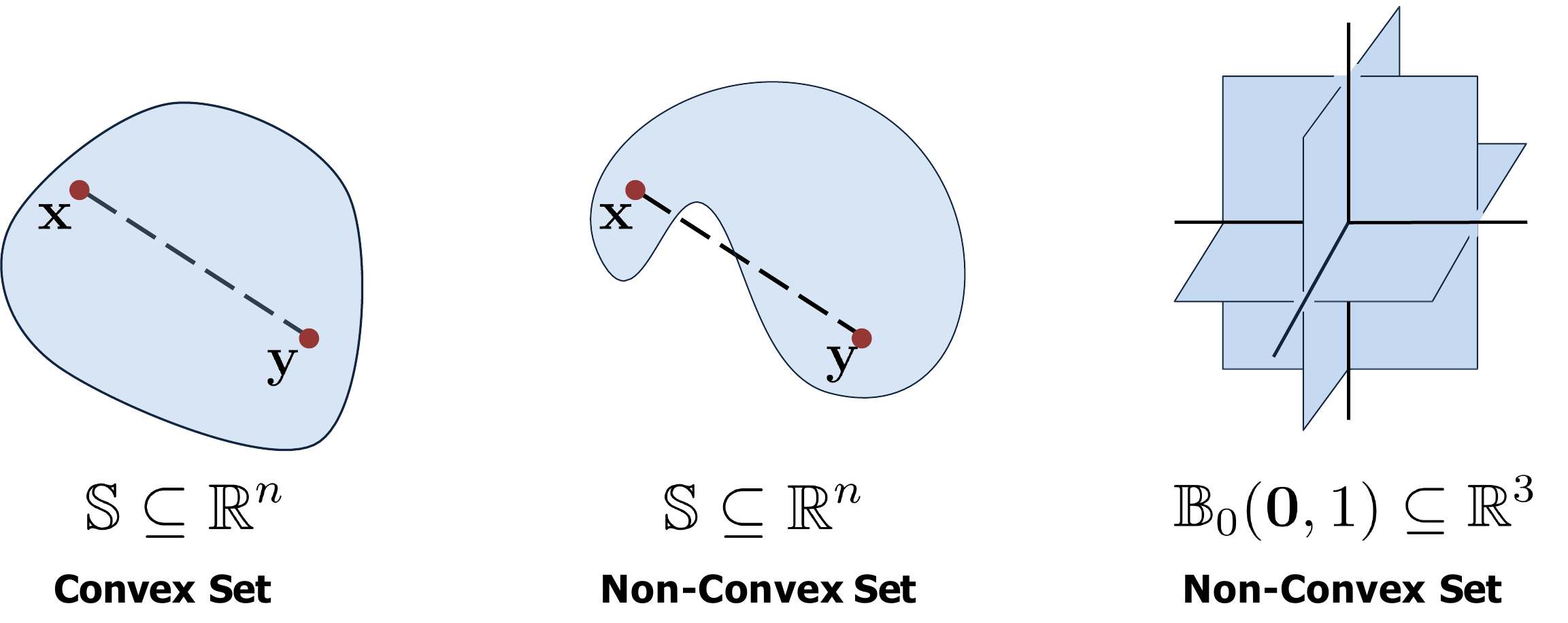}
\caption{A set is considered	 convex if it includes all convex combinations of its points. If there exists even one convex combination that lies outside the set, then by definition, the set is not convex. Therefore, a convex set must have a shape without any inward ``dents" or ``bulges". It's worth noting that the collection of sparse vectors does not satisfy this criterion and thus forms a non-convex set. The illustration provided is adapted from \citet{jain2017non}.}
\label{fig:cvxset}
\end{figure}

For vectors $\bx$ and $\by$,  the point $\lambda\bx+(1-\lambda)\by$, where $\lambda\in[0,1]$, is known as  a convex combination of the two vectors.
A set that is closed under arbitrary convex combinations is referred to as a convex set.
The standard definition follows:
\begin{definition}[Convex Set]\label{definition:convexset}
A set $\sS\subseteq \real^n$ is called \textit{convex} if, for any $\bx,\by\in\sS$ and $\lambda\in[0,1]$, the point $\lambda\bx+(1-\lambda)\by$ also belongs to $\sS$.
\end{definition}
Geometrically, convex sets  contain all line segments
that join two points within the set (Figure~\ref{fig:cvxset}). Consequently, these sets do not feature any concave indentations.

\begin{definition}[Cone and Convex Cone]\label{definition:convex_cone}
A set $ \sS \subseteq \real^n $ is a \textit{cone} if for every $ \bx \in \sS $, $ \alpha \bx \in \sS $ for any $ \alpha \geq 0 $.
A set $ \sS $ is a \textit{convex cone}  if it is both a cone and a convex set.
\end{definition}

A related concept is that of convex functions, which exhibit specific behavior under convex combinations. We now recall the definition:
\begin{definition}[Convex Functions]\label{definition:convexfuncs}
A function $f:\sS\rightarrow \real$ defined over a convex set $\sS\subseteq \real^n$ is called \textit{convex} if 
$$
f(\lambda\bx+(1-\lambda)\by)
\leq \lambda f(\bx) + (1-\lambda) f(\by), \text{ for any }\bx,\by\in\sS, \lambda\in[0,1].
$$
Moreover,   $f$ is called \textit{strictly convex} if 
$$
f(\lambda\bx+(1-\lambda)\by)
< \lambda f(\bx) + (1-\lambda) f(\by),\text{ for any }\bx\neq \by\in\sS, \lambda\in(0,1).
$$
\end{definition}
A well-known inequality derived from the concept of convex functions is provided below without a proof.
\begin{theorem}[Jensen's Inequality]\label{theorem:jensens_ineq}
Let $f: \sS \rightarrow \real$ be a convex function defined on a convex subset $\sS \subseteq \real^{n}$. For any finite sequence of points $\bx_{1}, \bx_{2}, \ldots, \bx_{m} \in \sS$ and any sequence of nonnegative weights $\lambda_{1}, \lambda_{2}, \lambda_{3}, \ldots, \lambda_{m}$ such that $\sum_{i=1}^{m} \lambda_{i} = 1$, Jensen's inequality states:
$$
f\left(\sum_{i=1}^{m} \lambda_{i} \bx_{i}\right) \leq \sum_{i=1}^{m} \lambda_{i} f(\bx_{i}).
$$
If $f$ is concave, the inequality is reversed:
$$
f\left(\sum_{i=1}^{m} \lambda_{i} \bx_{i}\right) \geq \sum_{i=1}^{m} \lambda_{i} f(\bx_{i}).
$$
In the context of probability theory, if $\rvx$ is a random vector with values in $\sS$ and $f$ is a convex function, Jensen's inequality can be stated as follows:
$$
f(\Exp[\rvx]) \leq \Exp[f(\rvx)],
$$
where $\Exp[\cdot]$ denotes the expectation operator over the random vector $\rvx$. For a concave function, the inequality is again reversed.
\end{theorem}

Similarly, we provide the definition of quasi-convexity and strict quasi-convexity.
\begin{definition}[Quasi-Convex Function]\label{definition:quasi_convex}
Let  $f:\sS\subseteq\real^n \rightarrow \real$ be a function defined over the convex set $\sS$. 
Then $f$ is called \textit{quasi-convex} if for any $\mu\in\real$ the level set $\lev[f, \mu]$ is convex. Formally, it satisfies that 
$$
f(\lambda \bx + (1-\lambda)\by) \leq  \max\{f(\bx), f(\by)\}, \text{ for any }\bx,\by\in\sS, \lambda\in[0,1].
$$
Similarlly, $f$ is called \textit{strictly quasi-convex} if 
$$
f(\lambda \bx + (1-\lambda)\by) <  \max\{f(\bx), f(\by)\},\text{ for any }\bx\neq \by\in\sS, \lambda\in(0,1).
$$ 
\end{definition}

The geometric interpretation of a  quasi-convex function is that the value of the function on any line segment within its domain does not exceed the maximum value of the function at the endpoints of the segment.
Quasi-convexity is a broader property than convexity since all convex functions are also quasi-convex, but not all quasi-convex functions are convex. For example, $f(x)=\sqrt{\abs{x}}$ is a quasi-convex function, but it is not convex.
Although a  quasi-convex function is not necessarily a convex function,  its level sets are convex, and it can include some well-behaved non-convex functions (Theorem~\ref{theorem:uni_qua_conv}).

In scenarios involving multiple variables, the input vector  $\bx$ can be partitioned into $\bx=(\ba, \bb)\in\real^n$. 
When discussing properties of functions in this context, convexity can sometimes be referred to as \textit{joint convexity}, meaning that the function is convex with respect to all components of the input simultaneously.
Additionally, there are instances where a function exhibits \textit{partial convexity}, often termed \textit{marginal convexity}, indicating that the function is convex when considered with respect to one or more subsets of its input variables while keeping others constant.

\begin{definition}[Marginally Convex]\label{definition:marginal_convexfuncs}
A function $f(\ba,\bb):\sS\subseteq\real^n\rightarrow \real$ (i.e., $(\ba,\bb)\in\real^n$), defined over a convex set $\sS$, is called \textit{marginally convex} if 
$$
f(\lambda\bx+(1-\lambda)\by, \bb)
\leq \lambda f(\bx, \bb) + (1-\lambda) f(\by, \bb), \text{ for any }(\bx,\bb),(\by,\bb)\in\sS, \lambda\in[0,1].
$$
\end{definition}

\begin{exercise}
Show that the norm functions $\normtwo{\bx}^2$  and $\normone{\bx}$ are convex.
\end{exercise}

\begin{exercise}
Suppose the continuously differentiable function $f$ is convex. Show that the directional derivative (Definition~\ref{definition:partial_deri}) $g(\bd)\triangleq f^\prime(\bx;\bd)$ is also convex.
\end{exercise}

\begin{exercise}[Convexity of Quadratic Functions]\label{exercise:conv_quad}
Let $f(\bx)=\frac{1}{2}\bx^\top\bA\bx+\bb^\top\bx+c$, where $\bA\in\real^{n\times n}$ is symmetric, $\bb\in\real^n$, and $c\in\real$. Show that $f(\bx)$ is convex (resp. strict convex) if and only if $\bA\succeq \bzero$ (resp. $\bA\succ\bzero$)
\end{exercise}

\begin{exercise}[Closure of Convex]\label{exercise:clo_conv}
Let $\sS\subseteq \real^n$ be a convex set. Show that $\closure(\sS)$ is also convex.
\end{exercise}

\begin{exercise}[Convexity of Indicator]\label{exercise_convex_indica}
	Show that the indicator function $\indicatorS$, which is $0$ if $\bx$ belongs to a set $\sS$ and $+\infty$ otherwise,  is convex if and only if $\sS$ is a convex set.
\end{exercise}

\paragraph{First-order characterizations of convex functions.}
Convex functions do not necessarily have to be differentiable. However, when they are  differentiable, such functions can be characterized by the gradient inequality.
\begin{theorem}[Gradient Inequality]\label{theorem:conv_gradient_ineq}
Let $f:\sS\rightarrow \real$ be a {continuously differentiable} function defined on a convex set $\sS\subseteq\real^n$. Then, $f$ is convex over $\sS$ if and only if 
\begin{equation}\label{equation:conv_gradient_ineq1}
f(\bx) +\nabla f(\bx)^\top (\by-\bx)\leq f(\by), \text{ for any $\bx,\by\in\sS$}.
\end{equation}
Similarly, the function is strictly convex over $\sS$ if and only if 
\begin{equation}\label{equation:conv_gradient_ineq2}
f(\bx) +\nabla f(\bx)^\top (\by-\bx)< f(\by), \text{ for any $\bx\neq \by\in\sS$}.
\end{equation}
This indicates that the graph of a convex function lies above its tangent plane at any point. For concave or strictly concave functions, the inequality signs are reversed.
\end{theorem}
\begin{proof}[of Theorem~\ref{theorem:conv_gradient_ineq}]
For brevity, we only prove \eqref{equation:conv_gradient_ineq1}, and \eqref{equation:conv_gradient_ineq2} can be proved analogously.
Suppose first that $ f $ is convex. Let $ \bx, \by \in \sS $ and $ \lambda \in [0, 1] $. If $ \bx = \by $, then \eqref{equation:conv_gradient_ineq1} trivially holds. We will therefore assume that $ \bx \neq \by $. Then, the definition of convexity shows that 
$$
f(\lambda \by + (1-\lambda) \bx) \leq \lambda f(\by) + (1-\lambda) f(\bx)
\;\implies\;
\frac{f(\bx + \lambda (\by - \bx)) - f(\bx)}{\lambda} \leq f(\by) - f(\bx).
$$
Taking $ \lambda \to 0^+ $, the left-hand side converges to the directional derivative of $ f $ at $ \bx $ in the direction $ \by - \bx $ (Definition~\ref{definition:partial_deri}), whence we have 
$
f'(\bx; \by - \bx) \leq f(\by) - f(\bx).
$
Since $ f $ is continuously differentiable, it follows that $ f'(\bx; \by - \bx) = \nabla f(\bx)^\top (\by - \bx) $ by \eqref{equation:direc_contdiff}, and hence \eqref{equation:conv_gradient_ineq1} follows.

Conversely, assume that the gradient inequality holds. Let $ \bx, \by \in \sS $, and let $ \lambda \in [0, 1] $. We will show that $ f(\lambda \bx + (1 - \lambda) \by) \leq \lambda f(\bx) + (1 - \lambda) f(\by) $. The cases for $\lambda=0$ or $\lambda=1$ holds trivially. We will therefore assume that $\lambda\in(0,1)$.
Let $ \ba \triangleq \lambda \bx + (1 - \lambda) \by \in \sS $. Then
$
 -\frac{\lambda}{1 - \lambda}  (\bx - \ba)  = (\by - \ba).
$
Combining with the gradient inequality implies that 
$$
f(\ba) + \nabla f(\ba)^\top (\bx - \ba) \leq f(\bx)
\qquad\text{and}\qquad 
f(\ba) - \frac{\lambda}{1 - \lambda} \nabla f(\ba)^\top (\bx - \ba) \leq f(\by).
$$
Multiplying the first inequality by $ \frac{\lambda}{1 - \lambda} $ and adding it to the second one, we obtain
$
\frac{1}{1 - \lambda} f(\ba) \leq \frac{\lambda}{1 - \lambda} f(\bx) + f(\by),
$
establishing the desired result.
\end{proof}

\begin{theorem}[Monotonicity of Gradient of Convex Functions]\label{theorem:monoton_convgrad}
Let $f:\sS\rightarrow \real$ be a {continuously differentiable} function defined over a convex set $\sS\subseteq\real^n$. Then, given any $\bx,\by\in\sS$, $f$ is convex over $\sS$ if and only if 
\begin{equation}\label{equation:monoton_convgrad}
\big(\nabla f(\bx) - \nabla f(\by)\big)^\top (\bx-\by)\geq 0.
\end{equation}
Similarly, given any $\bx\neq \by\in\sS$, $f$ is strictly convex over $\sS$ if and only if 
\begin{equation}\label{equation:monoton_convgrad2}
	\big(\nabla f(\bx) - \nabla f(\by)\big)^\top (\bx-\by)> 0.
\end{equation}
\end{theorem}
\begin{proof}[of Theorem~\ref{theorem:monoton_convgrad}]
For brevity, we only prove \eqref{equation:monoton_convgrad}, and \eqref{equation:monoton_convgrad2} can be proved analogously.
Assume  that $ f $ is convex over $ \sS $. Then by the gradient inequality, for any $ \bx, \by \in \sS $,  we have 
$$
f(\bx) \geq f(\by) + \nabla f(\by)^\top (\bx - \by)
\qquad\text{and}\qquad 
f(\by) \geq f(\bx) + \nabla f(\bx)^\top (\by - \bx).
$$
Summing the two inequalities yields the desired inequality in \eqref{equation:monoton_convgrad}.

Conversely, assume that \eqref{equation:monoton_convgrad} holds, and let $ \bx, \by \in \sS $. Define $g(\mu) \triangleq f(\bx + \mu (\by - \bx)),  \mu \in [0, 1]$ as a one-dimensional function.
By the fundamental theorem of calculus (Theorem~\ref{theorem:fund_theo_calculu}), we have
$$
\begin{aligned}
f(\by) 
&= g(1) = g(0) + \int_0^1 g'(\mu) \, d\mu
= f(\bx) + \int_0^1 (\by - \bx)^\top \nabla f\big(\bx + \mu (\by - \bx)\big) \, d\mu\\
&= f(\bx) + \nabla f(\bx)^\top (\by - \bx) + \int_0^1 (\by - \bx)^\top \big(\nabla f(\bx + \mu (\by - \bx)) - \nabla f(\bx)\big) \, d\mu\\
&\geq f(\bx) + \nabla f(\bx)^\top (\by - \bx),
\end{aligned}
$$
where the last inequality follows from the monotonicity of the gradient. This shows that $f$ is convex by Definition~\ref{definition:convexfuncs}.
\end{proof}

\paragraph{Second-order characterizations of convex functions.}

When the function is further assumed to be twice continuously differentiable,
its convexity can be characterized by the positive semidefiniteness (see Definition~\ref{definition:psd-pd-defini}) of the Hessian matrix
\begin{theorem}[PSD Hessian of Convex Functions]\label{theorem:psd_hess_conv}
Let $f:\sS\subseteq\real^n\rightarrow \real$ be a twice continuously differentiable function defined on an \textbf{open} convex set $\sS$, then the function is convex \textbf{if and only if} $\nabla^2 f(\bx)\succeq \bzero $ for any $\bx\in\sS$.~\footnote{Note that the openness of $\sS$ can be relaxed to prove the convexity of $\nabla^2 f(\bx)\succeq \bzero $ for any $\bx\in\interior(\sS)$.}
Moreover, \textbf{if} $\nabla^2 f(\bx)\succ \bzero $ for any $\bx\in\sS$, then the function is strictly convex over $\sS$~\footnote{Notice that the former condition is both sufficient and necessary, the latter condition is merely sufficient but not necessary. For example, $f(x)=x^6$ is strictly convex, but $f^{\prime\prime}(x)=30x^4$ is equal to zero at $x=0$.}. 
\end{theorem}
\begin{proof}[of Theorem~\ref{theorem:psd_hess_conv}]
Assume that $\nabla^2 f(\bx) \succeq \bzero$ for all $\bx \in \sS$, and let $\bx, \by \in \sS$. 
By the linear approximation theorem (Theorem~\ref{theorem:linear_approx}), there exists $\bxi \in [\bx, \by]$ (and hence $\bxi \in \sS$) such that
\begin{equation}\label{equation:psd_hess_conv1}
f(\by) = f(\bx) + \nabla f(\bx)^\top (\by - \bx) + \frac{1}{2} (\by - \bx)^\top \nabla^2 f(\bxi) (\by - \bx), \quad \bxi \in [\bx, \by].
\end{equation}
Given that  $\nabla^2 f(\bxi) \succeq \bzero$, it follows that $(\by - \bx)^\top \nabla^2 f(\bxi) (\by - \bx) \geq 0$, whence by \eqref{equation:psd_hess_conv1} we have  $f(\by) \geq f(\bx) + \nabla f(\bx)^\top (\by - \bx)$, establishing convexity by Theorem~\ref{theorem:conv_gradient_ineq}.

Conversely, assume that $f$ is convex over $\sS$. Let $\bx \in \sS$ and let $\by \in \real^n$. Since $\sS$ is open, there exists a scalar $\varepsilon>0$ such that  $\bx + \lambda \by \in \sS$ for any  $\lambda \in(0, \varepsilon)$. 
The gradient inequality (Theorem~\ref{theorem:conv_gradient_ineq}) shows that 
\begin{equation}\label{equation:psd_hess_conv2}
f(\bx + \lambda \by) \geq f(\bx) + \lambda \nabla f(\bx)^\top \by.
\end{equation}
Additionally, by the quadratic approximation theorem (Theorem~\ref{theorem:quad_app_theo}), we have
$$
f(\bx + \lambda \by) = f(\bx) + \lambda \nabla f(\bx)^\top \by + \frac{\lambda^2}{2} \by^\top \nabla^2 f(\bx) \by + o(\lambda^2 \normtwo{\by}^2),
$$
which combined with \eqref{equation:psd_hess_conv2} yields the inequality
$
\frac{\lambda^2}{2} \by^\top \nabla^2 f(\bx) \by + o(\lambda^2 \normtwo{\by}^2) \geq 0
$
for any $\lambda \in (0, \varepsilon)$. Dividing the latter inequality by $\lambda^2$ yields that 
$
\frac{1}{2} \by^\top \nabla^2 f(\bx) \by + \frac{o(\lambda^2 \normtwo{\by}^2)}{\lambda^2} \geq 0.
$
Taking $\lambda \to 0^+$, we conclude that
$
\by^\top \nabla^2 f(\bx) \by \geq 0
$
for any $\by \in \real^n$, implying that $\nabla^2 f(\bx) \succeq \bzero$ for any $\bx \in \sS$.
The strict convexity can be proved analogously.
\end{proof}

\begin{exercise}[Preservation of Closedness and Convexity]\label{exercise:pres_conv_clos}
Prove the following results:
\begin{enumerate}[(i)]
\item Let $\bA\in\real^{n\times n}$, and let $f: \real^n \rightarrow (-\infty, \infty]$ be an extended real-valued closed (resp. convex) function. Then the function 
$
g(\bx) = f(\bA\bx+\bb)
$
is closed (resp. convex).

\item Let $f_1, f_2, \ldots, f_m: \real^n \rightarrow (-\infty, \infty]$ be extended real-valued closed (resp. convex) functions, and let $\sigma_1, \sigma_2, \ldots, \sigma_m \in \real_+$. Then the function $f = \sum_{i=1}^{m} \sigma_i f_i$ is closed (resp. convex).

\item Let $f_i: \real^n \rightarrow (-\infty, \infty], i \in \sI$ be extended real-valued closed (resp. convex) functions, where $\sI$ is a given index set. Then the function
$
f(\bx) = \max_{i \in \sI} f_i(\bx)
$
is closed (resp. convex).
\end{enumerate}
\textit{Hint: the epigraph of $\max_{i \in \sI} f_i(\bx)$ is the intersection of the epigraphs of $f_i$'s.}
\end{exercise}

\subsection{Subgradient and Conjugate Functions}\label{section:sub_conjug}

The gradient inequality for convex functions applies specifically to continuously differentiable functions. However, this concept can be extended through the notion of a \textit{subgradient}, which plays a crucial role in optimization, especially for non-differentiable or non-smooth functions.

\begin{definition}[Subgradient and Subdifferential]\label{definition:subgrad}
Let $f: \sS\subseteq\real^n \rightarrow \real$. A vector $\bg\in\real^n$ is called a \textit{subgradient} of $f$ at $\bx\in\sS$ if
$$
f(\by) \geq  f(\bx) + \innerproduct{\bg, \by-\bx} \quad \text{for all $\by\in\sS$}.
$$
The inequality is called a \textit{subgradient inequality}.
The set of all subgradients of $f$ at $\bx$, denoted by $\partial f(\bx)$, is called the \textit{subdifferential} of $f$ at $\bx$:
$$
\partial f(\bx) \triangleq\left\{\bg\in\real^n\mid f(\by) \geq  f(\bx) + \innerproduct{\bg, \by-\bx}\text{ for all $\by\in\sS$} \right\}.
$$
Any subgradient of $f$ at $\bx$ can be denoted as
$$
f'(\bx) \in \partial f(\bx).
$$
\end{definition}
It's worth noting that the notion of a subgradient is applicable not only to convex functions but also extends to non-convex settings. For convex functions, however, it is assured that the subdifferential at any point within the domain is nonempty (Theorem~\ref{theorem:nonemp_relint_conv}).

\begin{exercise}[Subdifferential of Norms]\label{exercise:sub_norms}
Let $f(\bx)=\normtwo{\bx}$. Show that 
$$
\partial f(\bx)=
\begin{cases}
\left\{\frac{\bx}{\normtwo{\bx}}\right\}, & \bx\neq \bzero;\\
\sB_2[\bzero, 1], & \bx= \bzero.
\end{cases}
$$
Additionally, let $g(\bx)=\normone{\bx}$. Show that $\partial g(\bx) = \sum_{i=1}^{n} \partial g_i(\bx)$, where $g_i(\bx) = \abs{x_i}$ and 
$$
\partial  g_i(\bx)=
\begin{cases} 
\{\sign(x_i) \be_i\}, & x_i \neq 0, \\
[-\be_i, \be_i], & x_i = 0.
\end{cases}
$$
This indicates that 
$$
\sign(\bx)\in\partial g(\bx),
$$ 
where $\sign(\bx)$ returns the sign for each component.
\end{exercise}

An established result in convex analysis asserts that the relative interior of a convex set is always nonempty. Moreover, for a proper convex function, it is guaranteed to be subdifferentiable at every relative interior point of its effective domain. The following theorem encapsulates these findings without proof.
\begin{theorem}[Nonemptiness of Relative Interior \citep{rockafellar2015convex, beck2017first}]\label{theorem:nonemp_relint_conv}
Let $\sS\subseteq \real^n$ be a nonempty convex set. Then, the relative interior of $\sS$, denoted as $\relint(\sS)$, is nonempty.
Moreover, let $f : \real^n\rightarrow  (-\infty, \infty]$ be a proper convex function, 	and let $\bx \in \relint(\dom(f))$. 
Consequently, the subdifferential of $f$ at $\bx$,  $\partial f(\bx)$, is nonempty.
\end{theorem}

We then provide the definition and properties of conjugate functions.
\begin{definition}[Conjugate Functions]\label{definition:conjug_func}
Let $f : \real^n \rightarrow [-\infty,\infty]$ be an extended real-valued function. The conjugate function  $f^* : \real^n \rightarrow  \real \cup \{-\infty, \infty \}$ of $f$ is defined by
$$
f^*(\by) = \max_{\bx \in \real^n} \{ \innerproduct{\by, \bx } - f(\bx) \}, \quad \by \in \real^n.
$$
\end{definition}
Fenchel's inequality follows immediately from the definition of conjugate functions.
\index{Fenchel's inequality}
\begin{theorem}[Fenchel's Inequality]\label{theorem:fenchel_ineq}
Let $f: \real^n \rightarrow (-\infty, \infty]$ be a proper and  extended real-valued  function. Then for any $\bx,\by \in \real^n$,
$$
f(\bx) + f^*(\by) \geq \innerproduct{\by, \bx}.
$$
\end{theorem}

\begin{exercise}\label{exercise:biconjugate}
Let $f:\real^n\rightarrow [-\infty,\infty]$ be an extended real-valued function. Show that $f(\bx)\geq f^{**}(\bx)$ for any $\bx\in\real^n$.
Additionally, let $g:\real^n\rightarrow (-\infty,\infty]$ be a proper closed and convex function. Show that $g(\bx)= g^{**}(\bx)$
\end{exercise}

\begin{exercise}[Closedness and Convexity of Conjugate]\label{exercise:closedconv_conj}
Let $f:\real^n\rightarrow (-\infty,\infty]$ be a proper  function. Show that $f^*$ is  closed and convex.
\end{exercise}
\begin{exercise}[Properness of Conjugate]\label{exercise:proper_conj}
Let $f:\real^n\rightarrow (-\infty,\infty]$ be a proper convex function. Show that $f^*$ is also proper.
\end{exercise}

The main result concerning the subdifferential of a conjugate function is the so-called \textit{conjugate subgradient theorem}.

\index{Conjugate subgradient theorem}
\begin{theorem}[Conjugate Subgradient Theorem]\label{theorem:conju_subgra}
Let $ f : \real^n \to (-\infty, \infty] $ be a  proper and \textbf{convex} function. The following two claims are equivalent for any $ \bx, \by \in \real^n $:
\begin{enumerate}[(i)]
\item $ \innerproduct{\bx, \by} = f(\bx) + f^*(\by) $. 
\item $ \by \in \partial f(\bx) $.
\end{enumerate}
\noindent If  $ f $ is additionally  \textbf{closed}, then (i) and (ii) are equivalent to
\begin{itemize}[(iii)]
\item $ \bx \in \partial f^*(\by) $.
\end{itemize}
\noindent Alternatively, let $ f : \real^n \to (-\infty, \infty] $ be a proper \textbf{closed convex} function. Then for any $\bx, \by \in \real^n$,
\begin{subequations}
\begin{align}
\partial f(\bx) &= \argmax_{\bu \in \real^n} \left\{ \innerproduct{\bx, \bu} - f^*(\bu) \right\};\\
\partial f^*(\by) &= \argmax_{\bv \in \real^n} \left\{ \innerproduct{\by, \bv} - f(\bv) \right\}.
\end{align}
\end{subequations}
\end{theorem}
\begin{proof}[of Theorem~\ref{theorem:conju_subgra}]
By the subgradient inequality of $ \by \in \partial f(\bx) $, it follows that 
$$
f(\bu) \geq f(\bx) + \innerproduct{\by, \bu - \bx} 
\quad\implies\quad 
\innerproduct{\by, \bx } - f(\bx) \geq \innerproduct{ \by, \bu} - f(\bu)
\text{ for all } \bu \in \real^n.
$$
Taking the maximum over $ \bu $, the above inequality is equivalent to 
$$
\innerproduct{\by, \bx} - f(\bx) \geq f^*(\by),
$$
which by Fenchel's inequality (Theorem~\ref{theorem:fenchel_ineq}) is the same as the equality $\innerproduct{\bx, \by} = f(\bx) + f^*(\by)$, establishing the equivalence between (i) and (ii). 
When $f$ is proper closed and convex, by Exercise~\ref{exercise:biconjugate}, $f^{**} = f$, which  implies that (i) is equivalent to
$
\innerproduct{\bx, \by}= f^*(\by) + f^{**}(\bx)
$. By the above argument, we conclude that (i) is equivalent to $\bx \in \partial f^*(\by)$.
\end{proof}

\index{Fenchel's duality theorem}
An important result regarding conjugate functions is known as  \textit{Fenchel's duality theorem}. We briefly discuss the result here; and details can be found, for example, in \citet{bertsekas2009convex, rockafellar2015convex, beck2017first}.
Consider the following (primal) problem and its equivalent statement with an auxiliary variable:
$$
\begin{aligned}
\text{(P)}:\qquad  \quad& \min_{\bx \in \real^n} f(\bx) + g(\bx)\\
=&\min_{\bx, \by \in \real^n} \{f(\bx) + g(\by) : \bx = \by\}.
\end{aligned}
$$
Construct the Lagrangian function for (P) (see Section~\ref{section:gen_kkt_cond} for more details):
$$
L(\bx, \by, \blambda) = f(\bx) + g(\by) + \innerproduct{\blambda, \by - \bx} = -\left\{\innerproduct{\blambda, \bx} - f(\bx)\right\} - \left\{(-\blambda, \by) - g(\by)\right\}.
$$
The dual objective function is computed by minimizing the Lagrangian w.r.t. the primal variables $\bx, \by$:
$$
q(\blambda) = \min_{\bx, \by\in \real^n} L(\bx, \by, \blambda) = -f^*(\blambda) - g^*(-\blambda).
$$
Therefore, \textit{Fenchel's dual} can be expressed using the conjugate functions of the primal functions w.r.t. the Lagrangian multipliers $\blambda$:
$$
\text{(D)}:\qquad \min_{\blambda \in \real^n} \{f^*(\blambda) + g^*(-\blambda)\}.
$$
Fenchel's duality theorem then shows that strong duality (equivalence of (P) and (D)) holds for convex functions.
\begin{theorem}[Fenchel's Duality Theorem \citet{bertsekas2009convex, rockafellar2015convex, beck2017first}]\label{theorem:fen_dual_theo}
Let $f, g : \real^n \to (-\infty, \infty]$ be proper convex functions. If $\relint(\dom(f)) \cap \relint(\dom(g)) \neq \emptyset$~\footnote{By Theorem~\ref{theorem:nonemp_relint_conv}, the subdifferential in this set is nonempty.}, then
$$
\min_{\bx \in \real^n} \{f(\bx) + g(\bx)\} = \min_{\blambda \in \real^n} \{f^*(\blambda) + g^*(-\blambda)\}.
$$
\end{theorem}

\subsection{Lipschitz, Strongly Convex, Strongly Smooth Functions}
This section explores the properties of Lipschitz, strongly convex, and strongly smooth functions, defining these concepts and examining their implications on function behavior, including the boundedness of subdifferential sets for convex functions.
\index{Lipschitzness}
\begin{definition}[Lipschitz Functions]\label{definition:lipschi_funs}
Let $L\geq 0$, and let $f:\sS\rightarrow \real$ be a function defined over $\sS\subseteq  \real^n$. Then, the function $f$ is called \textit{L-Lipschitz} if, for every $\bx,\by\in\sS$, it follows that
$$
\abs{f(\bx)-f(\by)} \leq L\cdot \normtwo{\bx-\by}.
$$
If the function is continuously differentiable, then its gradient is called \textit{Lipschitz continuous} over $\real^n$ (or the function is called \textit{Lipschitz continuously differentiable}) if, for every $\bx,\by\in\sS$, it follows that
$$
\normtwo{\nabla f(\bx) - \nabla f(\by)} \leq L\cdot \normtwo{\bx-\by}.~\footnote{The concept of gradient Lipschitzness can be defined using dual norms: $\norm{\nabla f(\bx) - \nabla f(\by)}_* \leq L\cdot \norm{\bx-\by}$, where $\norm{\cdot}_*$ is the dual norm of $\norm{\cdot}$. Notably, the $\ell_2$ vector norm is self-dual. In our case, we only consider the $\ell_2$ vector norm.}
$$
The class of functions with a Lipschitz gradient and constant $L$ is denoted by $C_L^{1,1}(\real^n)$ or simply by $C_L^{1,1}$.

Similarly, the function is called \textit{twice Lipschitz continuously differentiable} with constant $L$, denoted by $C_L^{2,2}(\real^n)$, if, for every $\bx,\by\in\sS$, it follows that
$$
\normtwo{\nabla^2 f(\bx) - \nabla^2 f(\by)} \leq L\cdot \normtwo{\bx-\by}.
$$

More generally, let $\sS\subseteq\real^n$. We denote by $C_L^{k,p}(\sS)$ the class of functions possessing  the
following properties:
\begin{itemize}
\item Any $f\in C_L^{k,p}(\sS)$ is $k$ times continuously differentiable on $\sS$.
\item  Its $p$-th derivative is Lipschitz continuous on $\sS$ with constant $L$:
$\normtwo{\nabla^p f(\bx) - \nabla^p f(\by)} \leq L\cdot \normtwo{\bx-\by}.$
for all $\bx,\by\in\sS$. 
\end{itemize}
In this book, we usually work with $p = 0$, $p=1$, and $p = 2$. Clearly, we always have $p\leq k$. And it follows that $C_L^{t,p}(\sS) \subseteq  C_L^{k,p}(\sS)$ if $t\geq k$.
For simplicity, we also denotes $C^k(\sS)$ as the set of $k$ times continuously differentiable functions on $\sS$.
\end{definition}

\begin{example}\label{example:lipschitz_spar}
Consider the following examples:\begin{itemize}
\item Since $\absbig{\normone{\bx} - \normone{\by}} \leq \normone{\bx-\by} \leq \sqrt{n}\normtwo{\bx-\by}$ for any $\bx,\by\in\real^n$, $f_1(\bx)=\normone{\bx}$ is thus $\sqrt{n}$-Lipschitz.

\item Since $\absbig{\normtwo{\bb-\bA\bx} - \normtwo{\bb-\bA\by}} \leq \normtwo{(\bb-\bA\bx) - (\bb-\bA\by)} \leq \normtwo{\bA}\normtwo{\bx-\by}$, $f_2(\bx) = \normtwo{\bb-\bA\bx}$ is $\normtwo{\bA}$-Lipschitz, where $\normtwo{\bA}$ represents the spectral norm of $\bA$.
\end{itemize}
The function $f_3(\bx) \triangleq f_2^2(\bx) = \normtwo{\bb-\bA\bx}^2$ is not Lipschitz. However, it can be shown that $f_3(\bx)$ is Lipschitz continuously differentiable with constant $2\normtwo{\bA^\top\bA}$.
\end{example}

\begin{exercise}
Let $ f_1 \in C^{k,p}_{L_1}(\sS) $, $ f_2 \in C^{k,p}_{L_2}(\sS) $, $ \alpha_1, \alpha_2 \in \real $, and 
$ L_3 = \abs{\alpha_1} \cdot L_1 + \abs{\alpha_2} \cdot L_2 $.
Show that $ \alpha_1 f_1 + \alpha_2 f_2 \in C^{k,p}_{L_3}(\sS) $.
\end{exercise}

\begin{exercise}\label{exercise:cl2111_hess_bound}
Show that a function $ f(\cdot) $ belongs to the class $ C^{2,1}_L(\real^n) \subset C^{1,1}_L(\real^n) $ if and only if for all $ \bx \in \real^n $ we have
$
\normtwo{\nabla^2 f(\bx)} \leq L.
$~\footnote{Since the spectral norm of $\nabla^2 f(\bx)$ is the largest singular value of $\nabla^2 f(\bx)$, where $\nabla^2 f(\bx)$ is symmetric. This also indicates that $-L\bI\preceq \nabla^2 f(\bx) \preceq L\bI$.}
\textit{Hint: use the fundamental theorem of calculus (Theorem~\ref{theorem:fund_theo_calculu}).}
\end{exercise}

\begin{theorem}[Lipschitz and Boundedness of Subdifferential Sets]\label{theorem:lipsc_equiv}
Let $f:\sS \rightarrow \real$ be a  \textbf{convex} function, and let $\sX \subseteq \interior(\sS)$. Consider the following two claims:
\begin{enumerate}[(i)]
\item $\abs{f(\bx) - f(\by)} \leq L\normtwo{\bx - \by}$ for any $\bx, \by \in \sX$.
\item $\normtwo{\bg} \leq L$ for any $\bg \in \partial f(\bx), \bx \in \sX$.
\end{enumerate}
Then,
\begin{enumerate}[(a)]
\item The implication (ii) $\implies$ (i) holds,
\item If $\sX$ is open, then statement  (i) holds if and only if  statement  (ii) holds.
\end{enumerate}
\end{theorem}
\begin{proof}[of Theorem~\ref{theorem:lipsc_equiv}]
\textbf{(a).} Assume that (ii) is satisfied, and let $\bx, \by \in \sX$. Let $\bg_x \in \partial f(\bx)$ and $\bg_y \in \partial f(\by)$.
By the definitions of $\bg_x, \bg_y$ and the  Cauchy-Schwarz inequality (Proposition~\ref{proposition:cauchy-schwarz-inequ}),
\begin{align*}
f(\bx) - f(\by) &\leq \innerproduct{\bg_x, \bx - \by} \leq \normtwo{\bg_x} \normtwo{\bx - \by} \leq L\normtwo{\bx - \by}; \\
f(\by) - f(\bx) &\leq \innerproduct{\bg_y, \by - \bx} \leq \normtwo{\bg_y} \normtwo{\bx - \by} \leq L\normtwo{\bx - \by},
\end{align*}
which establishes the validity of (i).

\paragraph{(b).} The implication (ii) $\implies$ (i) was shown in (a). Conversely, assume that (i) is satisfied. Take $\bx \in \sX$ and $\bg \in \partial f(\bx)$. Let $\widetildebg\triangleq\frac{\bg}{\normtwo{\bg}}$ such that $\normtwo{\widetildebg} = 1 $ and $ \innerproduct{\widetildebg, \bg} = \normtwo{\bg}$. Since $\sX$ is open, we can take $\epsilon > 0$ small enough such that $\bx + \epsilon \widetildebg \in \sX$. 
Given the subgradient inequality, it follows that
$
f(\bx + \epsilon \widetildebg) \geq f(\bx) + \epsilon \innerproduct{\bg, \widetildebg},
$
whence we have
$$
\epsilon \normtwo{\bg}  \leq f(\bx + \epsilon \widetildebg) - f(\bx) \leq L \normtwo{\bx + \epsilon \widetildebg - \bx} = L\epsilon
\quad\implies\quad \normtwo{\bg} \leq L.
$$
This completes the proof.
\end{proof}

\begin{exercise}[Gradient Lipschitzness of Quadratic Functions]
Let $f(\bx) = \frac{1}{2} \bx^\top\bA\bx+\bb^\top\bx+c$, where $\bA$ is symmetric, $\bb\in\real^n$, and $c\in\real$. Show that the gradient of $f(\bx)$ is $\normtwo{\bA}$-Lipschitz. This also indicates that the gradient function is 0-Lipschitz when $\bA=\bzero$ (an affine function). 
\end{exercise}

Below, we provide rigorous definitions of \textit{strong convexity (SC)} and \textit{strong smoothness (SS, or simply smoothness)}. Although both terms can be defined for non-differentiable functions using subgradients in convex cases, we will focus on differentiable functions in this book.
\begin{figure}[htp]
\centering       
\vspace{-0.25cm}                 
\subfigtopskip=2pt               
\subfigbottomskip=-2pt         
\subfigcapskip=-10pt      
\includegraphics[width=0.98\textwidth]{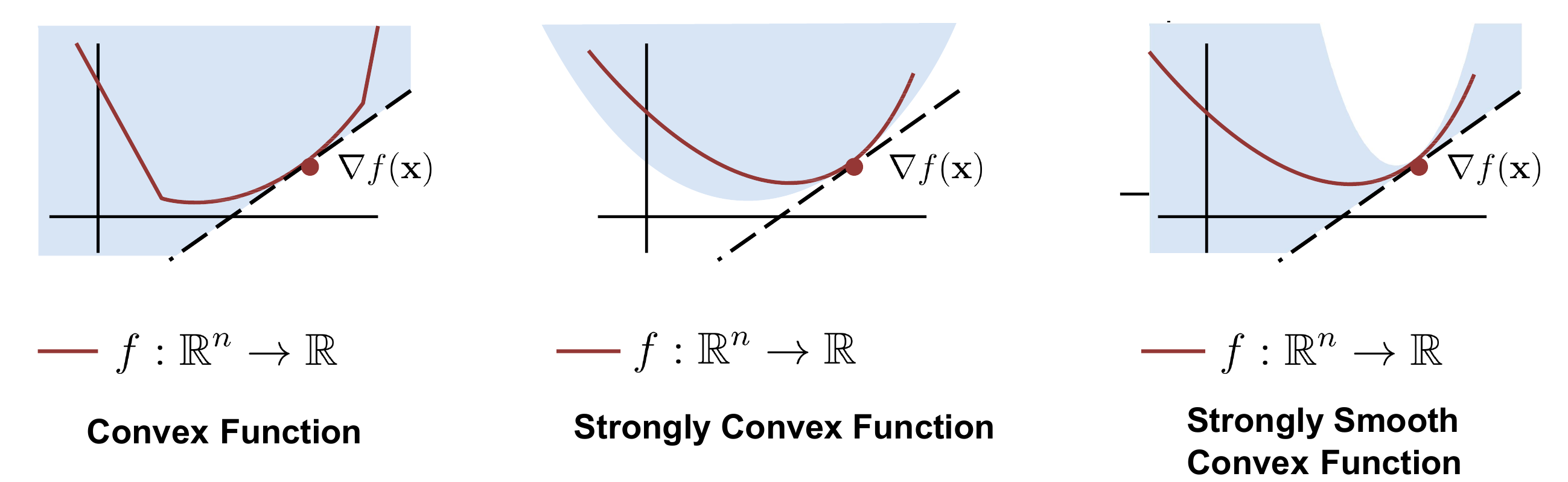}
\caption{A convex function is always bounded below by its tangent at any point. In contrast, strongly convex functions are not only bounded below but also have a quadratic lower bound that limits how slowly they can grow. Similarly, smooth functions have an upper bound on their growth rate, which prevents them from increasing too rapidly. These bounds ensure that strongly convex and smooth functions cannot grow either too slowly or too quickly compared to quadratic functions. In each figure, the shaded area represents the region where the function's curve is allowed to lie. The figure is adapted from \citet{jain2017non}.}
\label{fig:cvxfunc}
\end{figure}
\index{Strongly convex}
\index{Strongly smooth}
\index{Smooth}
\begin{definition}[Strongly Convex/Smooth Functions]\label{definition:scss_func}
Let $f:\sS\rightarrow \real$ be a  differentiable function defined over a \textbf{convex} $\sS\subseteq  \real^n$~\footnote{Note the notions of strong convexity and smoothness are primarily meaningful for differentiable functions (at least mostly for differentiable) defined over a convex set. For example, $\bS$ can be set to $\real^n$. The reason will be clear in the sequel.}. Then, $f$ is called \textit{$\alpha$-strongly convex (SC, or simply $\alpha$-convex)} and \textit{$\beta$-strongly smooth (SS, or simply $\beta$-smooth)} if, for every $\bx,\by\in\sS$, it follows that 
\begin{equation}\label{equation:scss_func1}
\frac{\alpha}{2} \normtwo{\bx-\by}^2
\leq f(\by)-f(\bx)-\nabla f(\bx)^\top (\by-\bx)
\leq \frac{\beta}{2} \normtwo{\bx-\by}^2.~\footnote{In many texts, the second inequality is called the \textit{descent lemma} for $\beta$-strongly smooth functions since $f(\by) \leq f(\bx)+\nabla f(\bx)^\top (\by-\bx)+\frac{\beta}{2} \normtwo{\bx-\by}^2$, i.e., an update of the function value  from $f(\bx)$ to $f(\by)$ is upper-bounded.}
\end{equation}
When the function is twice differentiable, Theorem~\ref{theorem:psd_hess_conv} also indicates  that 
\begin{equation}\label{equation:scss_func2}
\alpha\bI \preceq \nabla^2 f(\bx) \preceq \beta\bI, \quad\text{for all } \bx\in\interior(\sS).~\footnote{We relax to $\bx\in\interior(\sS)$ as opposed to claiming an open set in Theorem~\ref{theorem:psd_hess_conv}.}
\end{equation}
That is, the smallest eigenvalue of $\nabla^2 f(\bx)$ is at least $\alpha$, and the largest eigenvalue of $\nabla^2 f(\bx)$ is at most $\beta$ (Theorem~\ref{theorem:eigen_charac}):
\begin{equation}\label{equation:scss_func3_ra}
\alpha\leq \frac{\bv^\top \nabla^2 f(\bx) \bv}{\normtwo{\bv}^2} \leq \beta, \quad\text{for all } \bx\in\sS, \text{and nonzero } \bv\in\real^n.~\footnote{$\frac{\bv^\top\bA\bv}{\normtwo{\bv}^2}$ is called the \textit{Rayleigh quotient} of vector $\bv$ associated with the matrix $\bA$, which lies between the largest and smallest eigenvalues of $\bA$; see, for example, \citet{lu2021numerical} for more details.}
\end{equation}
It is evident that if the function is $\beta_1$-strongly smooth, it is also $\beta_2$-strongly smooth for any $\beta_2\geq \beta_1$.
A $\alpha$-strongly convex function is necessarily a convex function (when $\alpha=0$, it reduces to a standard convex function). 
While if the function is $\alpha_1$-strongly convex, it is also $\alpha_2$-strongly convex for any $\alpha_2\in(0, \alpha_1]$.
\end{definition}
Strong convexity and smoothness play crucial roles in analyzing gradient descent algorithms (see Chapter~\ref{chapter:gd_convg}). Definition~\ref{definition:res_scss_func} extends these notions to non-convex sets.
Strong convexity ensures the function curves upwards sharply, providing a lower bound on its curvature and preventing gradients from vanishing too quickly. Conversely, strong smoothness imposes a quadratic upper bound on the curvature, ensuring that gradients do not change too rapidly, thereby controlling the magnitude of gradient descent steps and avoiding overly aggressive updates (see Figure~\ref{fig:cvxfunc}).

\begin{example}[Convex, SC]
Let $f(\bx) = \bx^\top \bA \bx + \bb^\top \bx + c$.
Then the function is 
\begin{itemize}
\item Convex if and only if $\bA \succeq \bzero$; strictly convex if and only if $\bA \succ \bzero$.
\item Concave if and only if $\bA \preceq \bzero$; strictly concave if and only if $\bA \prec \bzero$.
\end{itemize}
The proofs are easy if we use the second-order characterization of convexity (Theorem~\ref{theorem:psd_hess_conv}).
\end{example}

\begin{exercise}[Sum of SC and Convex, SS]\label{exercise:sum_sc_conv}
Prove the following two results:
\begin{itemize}
\item Let $f$ be a $\alpha$-SC function and $g$ be a convex function. Show that $f+g$ is $\alpha$-SC. 
\item Let $f$ be a $\beta_1$-SS function and $g$ be a $\beta_2$-SS function. Show that $f+g$ is $(\beta_1+\beta_2)$-SS. 
\end{itemize}
\end{exercise}
\begin{exercise}[Level Sets of SC and SS]
Let $f:\real^n\rightarrow \real$ be $\alpha$-SC and $\beta$-SS, and define the level set $\lev[f, \gamma] \triangleq \{\bx \mid f(\bx)\leq\gamma\}$. Show that 
$$
\sB_2\Big[\bx^*, \sqrt{\frac{2}{\beta}(\gamma-f(\bx^*))}\Big]
\subseteq
\lev[f, \gamma]\subseteq 
\sB_2\Big[\bx^*, \sqrt{\frac{2}{\alpha}(\gamma-f(\bx^*))}\Big],
$$
where $\bx^*$ is the global minimum point of $f$ (such that $f(\bx^*)\leq f(\bx)$ for all $\bx$; Definition~\ref{definition:local_global_opt}).
That is, the level set lies between the two balls.
\end{exercise}

\begin{theorem}[SS Property-O: Equivalence between Gradient Lipschitzness and Smoothness]\label{theorem:equi_gradsch_smoo}
Let $f:\sS\rightarrow \real$ be a  differentiable function defined over a \textbf{convex} set $\sS\subseteq  \real^n$.
If $f$ is gradient  Lipschitz continuous with constant $\beta$, i.e., $f\in C_{\beta}^{1,1}(\sS)$, then it is also $\beta$-smooth.
\end{theorem}
\begin{proof}[of Theorem~\ref{theorem:equi_gradsch_smoo}]
By the fundamental theorem of calculus (Equation~\eqref{equation:fund_theo_calculu3}), it follows that 
$
f(\by) - f(\bx) = \innerproduct{\nabla f(\bx), \by - \bx}+ \int_0^1 \innerproduct{\nabla f(\bx + \mu(\by - \bx)) - \nabla f(\bx), \by - \bx} d\mu.
$
Since $\sS$ is a convex set, the function $f(\bx + \mu(\by - \bx))$ is well-defined for $\mu\in[0,1]$.
This results in 
$$
\begin{aligned}
&\gap \abs{f(\by) - f(\bx) - \innerproduct{\nabla f(\bx), \by - \bx}} =  \abs{\int_0^1 \innerproduct{\nabla f(\bx + \mu(\by - \bx)) - \nabla f(\bx), \by - \bx} d\mu } \\
&\overset{\dag}{\leq} \int_0^1 \normtwo{\nabla f(\bx + \mu(\by - \bx)) - \nabla f(\bx)}\cdot  \normtwo{\by - \bx} d\mu
\leq   \int_0^1 \beta \mu \normtwo{\by - \bx}^2 d\mu 
= \frac{\beta}{2} \normtwo{\by - \bx}^2,
\end{aligned}
$$
where the inequality $(\dag)$ follows from the  Cauchy-Schwarz inequality (Proposition~\ref{proposition:cauchy-schwarz-inequ}).
\end{proof}

If $ f $ is $\beta$-smooth and has a global minimum point $ \bx^* $, an important corollary is that we can use the quadratic upper bound to estimate the magnitude of $f(\bx) - f(\bx^*)$, where $ \bx $ can be any point in the domain.
\begin{theorem}[SS\&SC Property-I: Bound of $f(\bx) - f(\bx^*)$]\label{theorem:smoo_prop2_bound}
Let $ f:\real^n\rightarrow \real $ be  a differentiable function  defined on $ \real^n $  and have a global minimum point $ \bx^* $. If $ f(\bx) $ is  $ \beta $-smooth and $\alpha$-strongly convex, then
$$
\frac{1}{2\beta} \normtwo{\nabla f(\bx)}^2 \leq f(\bx) - f(\bx^*)
\leq \frac{1}{2\alpha}  \normtwo{\nabla f(\bx)}^2,\quad  \text{ for any } \bx\in\real^n.
$$
\end{theorem}
\begin{proof}[of Theorem~\ref{theorem:smoo_prop2_bound}]
Since $ \bx^* $ is a global minimum point, for any $\by\in\real^n$, applying the quadratic upper bound from the smoothness gives
$$
f(\bx^*) \leq f(\by) \leq f(\bx) + \nabla f(\bx)^\top (\by - \bx) + \frac{\beta}{2} \normtwo{\by - \bx}^2.
$$
Fixing $ \bx $, note that the above inequality holds for any $ \by $, thus we can take the infimum on the right side of the inequality:
$$
\begin{aligned}
f(\bx^*) &\leq \inf_{\by \in \real^n} \left\{ f(\bx) + \nabla f(\bx)^\top (\by - \bx) + \frac{\beta}{2} \normtwo{\by - \bx}^2 \right\} 
= f(\bx) - \frac{1}{2\beta} \normtwo{\nabla f(\bx)}^2,
\end{aligned}
$$
where the minimum is attained at $\by = \bx-\frac{1}{\beta}\nabla f(\bx)$ (since the domain is $\real^n$).

For the second part, since $f$ is strongly convex, we have 
$$
\begin{aligned}
f(\by)&\geq f(\bx) +\nabla f(\bx)^\top (\by-\bx) + \frac{\alpha}{2} \normtwo{\by-\bx}^2\\
&\geq f(\bx) +\nabla f(\bx)^\top (\widetildeby-\bx) + \frac{\alpha}{2} \normtwo{\widetildeby-\bx}^2\\
&=f(\bx)-\frac{1}{2\alpha} \normtwo{\nabla f(\bx)}^2,
\end{aligned}
$$
where $\widetildeby \triangleq \bx-(1/\alpha) \nabla f(\bx)$ minimizes the right-hand side of the above inequality. 
Letting $\by\triangleq\bx^*$ completes the proof.
\end{proof}

The second result in the above theorem shows that 
\begin{equation}
\normtwo{\nabla f(\bx)} \leq (2\alpha\varepsilon)^{1/2}
\quad\implies\quad 
f(\bx) - f(\bx^*)\leq \varepsilon,
\end{equation}
where $\varepsilon$ represents a positive  \textit{tolerance}.
This generalizes the \textit{optimality condition for convex functions} stated  in Theorem~\ref{theorem:fetmat_opt}.
Additionally, the theorem indicates that a strongly smooth function exhibits a specific growth characteristic near the minimizer.

The next theorem states that a proper closed and strongly convex function has a unique minimizer and  satisfies a certain growth property around this minimizer.
\begin{theorem}[SC Property-II: Existence and Uniqueness of a Minimizer of Closed SC Functions]\label{theorem:exi_close_sc}
Let $f :  \real^n \rightarrow (-\infty, \infty]$ be a differentiable~\footnote{Note that this theorem also holds for non-differentiable functions, and subgradient can be used in the proof by Theorem~\ref{theorem:nonemp_relint_conv}.}, proper closed, and $\alpha$-strongly convex function ($\alpha > 0$). Then,
\begin{enumerate}[(i)]
\item  $f$ has a unique minimizer.
\item  $f(\bx) - f(\bx^*) \geq \frac{\alpha}{2} \normtwo{\bx - \bx^*}^2$ for all $\bx \in \real^n$, where $\bx^*$ is the unique minimizer of $f$.
\end{enumerate}
\end{theorem}
\begin{proof}[of Theorem~\ref{theorem:exi_close_sc}]
\textbf{(i).} Let $\bx_0\in\real^n$. By  strong convexity, it follows that
$$
\begin{aligned}
&f(\bx) \geq f(\bx_0) + \innerproduct{\nabla f(\bx_0), \bx - \bx_0} + \frac{\alpha}{2} \normtwo{\bx - \bx_0}^2 \\
&\quad\implies 
f(\bx) \geq f(\bx_0) - \frac{1}{2\alpha} \normtwo{\nabla f(\bx_0)}^2 + \frac{\alpha}{2} \normtwo{\bx - \left(\bx_0 - \frac{1}{\alpha} \nabla f(\bx_0)\right)}^2 \quad \text{for any } \bx \in \real^n.
\end{aligned}
$$
This indicates that
$$
\lev[f, f(\bx_0)] \subseteq \sB_{2}\left[\bx_0 - \frac{1}{\alpha} \nabla f(\bx_0), \frac{1}{\alpha} \normtwo{\nabla f(\bx_0)}\right].
$$
Since $f$ is closed, the above level set is closed (Theorem~\ref{theorem:equiv_close_clkos_semicon}); and since it is contained in a ball, it is also bounded. Therefore, $\lev[f, f(\bx_0)]$ is compact. We can thus deduce that the optimal set for minimizing $f$ over $\dom(f)$ is identical to the optimal set for minimizing $f$ over the nonempty compact set $\lev[f, f(\bx_0)]$. Invoking Weierstrass theorem  for proper closed functions (\ref{weier2_prop_close}.\ref{weier2_prop_close_v1} in Theorem~\ref{theorem:weierstrass_them}), it follows that a minimizer exists. To show the uniqueness, assume that $\widetildebx$ and $\bx^*$ are minimizers of $f$. Then $f(\widetildebx)  = f(\bx^*)$, where $f(\bx^*)$ is the minimal value of $f$. Then by property (iv) of Theorem~\ref{theorem:charac_stronconv},
$$
f(\bx^*) \leq f\left(\frac{1}{2}\widetildebx + \frac{1}{2}\bx^*\right) \leq \frac{1}{2} f(\widetildebx) + \frac{1}{2} f(\bx^*) - \frac{\alpha}{8} \normtwo{\widetildebx - \bx^*}^2 = f(\bx^*) - \frac{\alpha}{8} \normtwo{\widetildebx - \bx^*}^2,
$$
implying that $\widetildebx = \bx^*$, thus proving the uniqueness of the minimizer of $f$.

\paragraph{(ii).} Let $\bx^*$ be the unique minimizer of $f$. Then by optimality condition for convex functions (Theorem~\ref{theorem:fetmat_opt}), $ \nabla f(\bx^*) = \bzero$. Thus, by the definition of strong convexity,
$$
f(\bx) - f(\bx^*) \geq \langle \bzero, \bx - \bx^* \rangle + \frac{\alpha}{2} \normtwo{\bx - \bx^*}^2 = \frac{\alpha}{2} \normtwo{\bx - \bx^*}^2 ,
\;\text{ for any $\bx \in \real^n$},
$$
establishing claim (ii).
\end{proof}

The definition of smoothness does not require the function to be convex. However, when it indeed is, the function can be equivalently characterized as follows. 
\begin{theorem}[Characterization Theorem of SS and Convexity \citep{beck2017first}]\label{theorem:charac_smoo}
Let $f: \sS\subseteq\real^n \rightarrow \real$ be a differentiable \textbf{convex} function defined on a convex set $\sS$, and let $\beta > 0$. Then the following claims are equivalent:
\begin{enumerate}[(i)]
\item $f$ is $\beta$-smooth.
\item The function $g(\bx)\triangleq \frac{\beta}{2}\bx^\top\bx - f(\bx)$ is convex.
\item $f(\by) \leq f(\bx) + \innerproduct{\nabla f(\bx), \by - \bx} + \frac{\beta}{2} \normtwo{\bx - \by}^2$ for all $\bx, \by \in \sS$.
\item $ f(\by) \geq f(\bx) + \innerproduct{\nabla f(\bx), \by - \bx} + \frac{1}{2\beta} \normtwo{\nabla f(\bx) - \nabla f(\by)}^2 $ for all $\bx, \by \in \sS$.
\item $f(\lambda \bx + (1 - \lambda)\by) \geq \lambda f(\bx) + (1 - \lambda)f(\by) - \frac{\beta}{2} \lambda (1 - \lambda) \normtwo{\bx - \by}^2$ for any $\bx, \by \in \sS$ and $\lambda \in [0, 1]$.
\end{enumerate}
Reversing $\bx$ and $\by$ in (iii) and (iv), we can also obtain 
\begin{equation}\label{equ:charac_smoo_3} 
\frac{1}{\beta} \normtwo{\nabla f(\bx) - \nabla f(\by)}^2 
\leq  \innerproduct{\nabla f(\bx) - \nabla f(\by), \bx - \by} 
\leq  
\beta\normtwo{\bx-\by}^2,  
\text{ for all  $\bx, \by \in \sS$.}
\end{equation}
\end{theorem}
\begin{proof}[of Theorem~\ref{theorem:charac_smoo}]
\textbf{(i)$\implies$(ii).} For any $\bx,\by\in\sS$, it follows that 
$$
\begin{aligned}
\big(\nabla g(\bx) - \nabla g(\by)\big)^\top (\bx - \by) 
&= \beta \normtwo{\bx - \by}^2 - \big(\nabla f(\bx) - \nabla f(\by)\big)^\top (\bx - \by)\\
&\geq \beta \normtwo{\bx - \by}^2 - \normtwo{\bx - \by} \normtwo{\nabla f(\bx) - \nabla f(\by)} \geq 0.
\end{aligned}
$$
where the last inequality follows from Theorem~\ref{theorem:equi_gradsch_smoo} and the definition of smoothness.
Thus, $ g(\bx) $ is a convex function by Theorem~\ref{theorem:monoton_convgrad}.

\paragraph{(ii)$\implies$(iii).}
Since $ g(\bx) \triangleq \frac{\beta}{2} \normtwo{\bx}^2 - f(\bx) $ is convex, this implies
$$
g(\by) \geq g(\bx) + \nabla g(\bx)^\top (\by - \bx), \quad \forall \bx, \by \in \sS,
$$
where
$
\nabla g(\bx) = \beta \bx - \nabla f(\bx).
$
Substituting the gradient into the convexity condition yields that
\begin{align}
\frac{\beta}{2} \normtwo{\by}^2 - f(\by) \geq \frac{\beta}{2} \normtwo{\bx}^2 - f(\bx) + \big(\beta \bx - \nabla f(\bx)\big)^\top (\by - \bx) \nonumber\\
\quad\implies\quad f(\by) - f(\bx) \leq \nabla f(\bx)^\top (\by - \bx) + \frac{\beta}{2} \normtwo{\by - \bx}^2. \label{equation:charac_smoo11}
\end{align}

\paragraph{(iii)$\implies$(iv).} We can assume that $\nabla f(\bx) \neq \nabla f(\by)$ since otherwise the inequality (iv) is trivial by the convexity of $f$. For a fixed $\bx \in \sS$, consider the function
$$ 
g_{\bx}(\by) \triangleq f(\by) - f(\bx) - \innerproduct{\nabla f(\bx), \by - \bx}, \quad \by \in \sS, 
$$
which is convex w.r.t. $\by$.
Note that $\nabla g_{\bx}(\by) = \nabla f(\by)-\nabla f(\bx)$ for any $\by\in\sS$.
By the condition of (iii), for any $\by, \bu \in \sS$, we have 
\begin{equation}\label{equation:char_ss_eqq1}
\begin{aligned}
g_{\bx}(\bu) &= f(\bu) - f(\bx) - \innerproduct{\nabla f(\bx), \bu - \bx} \\
&\leq \big\{f(\by) + \innerproduct{\nabla f(\by), \bu - \by} + \frac{\beta}{2} \normtwo{\bu - \by}^2\big\} - f(\bx) - \innerproduct{\nabla f(\bx), \bu - \bx} \\
&= f(\by) - f(\bx) - \innerproduct{\nabla f(\bx), \by - \bx} + \innerproduct{\nabla f(\by) - \nabla f(\bx), \bu - \by} + \frac{\beta}{2} \normtwo{\bu - \by}^2 \\
&= g_{\bx}(\by) + \innerproduct{\nabla g_{\bx}(\by), \bu - \by} + \frac{\beta}{2} \normtwo{\bu - \by }^2,
\end{aligned}
\end{equation}
Since $\nabla g_{\bx}(\bx) = \bzero$, which by the convexity of $g_{\bx}$ implies that $\bx$ is a global minimizer of $g_{\bx}$ such that
$
0=g_{\bx}(\bx) \leq g_{\bx}(\bu),\, \text{for all } \bu \in \sS.  
$
Let $\by \in \sS$, and let $\widetildebg \triangleq \frac{\nabla g_{\bx}(\by)}{\normtwo{\nabla g_{\bx}(\by)}}$ such that  $\normtwo{\widetildebg} = 1$ and $\innerproduct{\nabla g_{\bx}(\by), \widetildebg} = \normtwo{\nabla g_{\bx}(\by)}$. Defining
$ \bu \triangleq \by - \frac{\normtwo{\nabla g_{\bx}(\by)}}{\beta} \widetildebg  $, the optimum of $g_{\bx}(\bx)$ shows that 
$$
0 = g_{\bx}(\bx) \leq g_{\bx} \left( \by - \frac{\normtwo{\nabla g_{\bx}(\by)}}{\beta} \widetildebg \right). 
$$
Combining the preceding inequality with \eqref{equation:char_ss_eqq1}, we obtain
\begin{align*}
0 &= g_{\bx}(\bx) 
\leq g_{\bx}(\by) - \frac{\normtwo{\nabla g_{\bx}(\by)}}{\beta} \innerproduct{\nabla g_{\bx}(\by), \widetildebg} + \frac{\normtwo{\nabla g_{\bx}(\by)}^2 \cdot \normtwo{\widetildebg}^2 }{2\beta} 
= g_{\bx}(\by) - \frac{\normtwo{\nabla g_{\bx}(\by)}^2}{2\beta}  \\
&= f(\by) - f(\bx) - \innerproduct{\nabla f(\bx), \by - \bx} - \frac{1}{2\beta} \normtwo{\nabla f(\bx) - \nabla f(\by)}^2,
\end{align*}
establishing (iv).

\paragraph{(iii)$\implies$(v).}
Let $\bx, \by \in \sS$ and $\lambda \in [0, 1]$, and denote $\bxi \triangleq \lambda \bx + (1 - \lambda) \by$. Since $\sS$ is convex, we also have $\bxi\in\sS$.
Invoking (iii) with $\bx,\bxi$ and $\by,\bxi$, and plugging in the expression of $\bxi$ yields that 
$$
\begin{aligned}
	f(\bx) &\leq f(\bxi) + (1 - \lambda) \innerproduct{\nabla f(\bxi), \bx - \by} + \frac{\beta (1 - \lambda)^2}{2} \normtwo{\bx - \by}^2, \\
	f(\by) &\leq f(\bxi) + \lambda \innerproduct{\nabla f(\bxi), \by - \bx} + \frac{\beta \lambda^2}{2} \normtwo{\bx - \by}^2.
\end{aligned}
$$
Multiplying the first inequality by $\lambda$ and the second by $1 - \lambda$ and adding them yields the inequality (v).

\paragraph{(v)$\implies$(iii).} Rearranging terms in the inequality (v), we obtain that it is equivalent to
$$
f(\by) \leq f(\bx) + \frac{f\big(\bx + (1 - \lambda)(\by - \bx)\big) - f(\bx)}{1 - \lambda} + \frac{\beta}{2} \lambda \normtwo{\bx - \by}^2. 
$$
Taking $\lambda \to 1^-$, the preceding inequality becomes
$ f(\by) \leq f(\bx) + f'(\bx; \by - \bx) + \frac{\beta}{2} \normtwo{\bx - \by}^2$.
Since $f$ is differentiable, \eqref{equation:direc_contdiff} shows that $f'(\bx; \by - \bx) = \innerproduct{\nabla f(\bx), \by - \bx}$, establishing the desired result. 
\end{proof}

\begin{theorem}[Characterization Theorem of SC \citep{beck2017first}]\label{theorem:charac_stronconv}
Let $ f:\sS\subseteq\real^n\rightarrow \real $ be a convex and differentiable function defined on a convex set $\sS$, and let $ \alpha > 0 $. Then the following  claims are equivalent:
\begin{enumerate}[(i)]
\item $ f $ is $ \alpha $-strongly convex.
\item $
f(\by) \geq f(\bx) + \innerproduct{\nabla f(\bx) , \by - \bx} + \frac{\alpha}{2} \normtwo{\by - \bx}^2
$
for any $ \bx , \by \in \sS$.
\item $
\innerproduct{\nabla f(\bx)  - \nabla f(\by) , \bx - \by} \geq \alpha \normtwo{\bx - \by}^2
$
for any $ \bx, \by \in \sS $.
\item $f(\lambda \bx+(1-\lambda)\by )\leq \lambda f(\bx)+(1-\lambda)f(\by) -\frac{\alpha}{2}\lambda(1-\lambda)\normtwo{\bx-\by}^2$ for any $\lambda\in[0,1]$ and $\bx,\by\in\sS$.
\end{enumerate}
\end{theorem}
The proof is similar to that of Theorem~\ref{theorem:charac_smoo} and is left as an exercise.
Notably, combining Theorem~\ref{theorem:charac_smoo} and Theorem~\ref{theorem:charac_stronconv}, if $f:\sS\subseteq\real^n\rightarrow \real$ is differentiable, $\alpha$-SC, and $\beta$-SS defined on a convex set $\sS$ with $0<\alpha<\beta$, then for any $\bx,\by\in\sS$,
\begin{subequations}
\begin{align}
\frac{\alpha}{2} \normtwo{\bx-\by}^2
\;&\leq f(\by) -f(\bx) - \innerproduct{\nabla f(\bx), \by - \bx}
&\leq&  \frac{\beta}{2} \normtwo{\bx - \by}^2; \label{equation:sssc_eq1} \\
\alpha \normtwo{\bx - \by}^2 
\;&\leq \innerproduct{\nabla f(\bx)  - \nabla f(\by) , \bx - \by} 
&\leq& \beta\normtwo{\bx-\by}^2;\label{equation:sssc_eq2} \\
\frac{\alpha\zeta}{2}\normtwo{\bx-\by}^2 
\;&\leq \big\{\lambda f(\bx)+(1-\lambda)f(\by)\} - f\big(\lambda \bx+(1-\lambda)\by \big) 
&\leq& \frac{\beta\zeta}{2}  \normtwo{\bx - \by}^2, \label{equation:sssc_eq3}
\end{align}
\end{subequations}
where $\zeta\triangleq \lambda(1-\lambda)$ for any $\lambda\in[0,1]$, and \eqref{equation:sssc_eq1} is equivalent to Definition~\ref{definition:scss_func}.
In some texts, the definitions of SC and SS are defined using \eqref{equation:sssc_eq3} instead since it does not require the function to be differentiable.
Again, we focus on differentiable functions in this book. 
When the function is not differentiable, the gradient in Theorem~\ref{theorem:charac_smoo} or  Theorem~\ref{theorem:charac_stronconv} can be replaced with subgradients due to the nonemptyness of subdifferential by Theorem~\ref{theorem:nonemp_relint_conv}.
Additionally, we have the following characterization theorem for SC and SS functions.
\begin{theorem}[Characterization Theorem of SC and SS]\label{theorem:charac_smoo_n_stronconv}
Let $f:\sS\subseteq\real^n\rightarrow \real$ be a differentiable, $\alpha$-strongly convex, and $\beta$-smooth function defined on a convex set $\sS$ with $\alpha\leq \beta$. Then, for any $\bx,\by\in\sS$, it follows that 
\begin{equation}\label{equation:baillon_hadded}
\innerproduct{\nabla f(\bx) - \nabla f(\by), \bx -\by}
\geq
\frac{\alpha\beta}{\alpha+\beta} \normtwo{\bx-\by}^2 + \frac{1}{\alpha+\beta} \normtwo{\nabla f(\bx) - \nabla f(\by)}^2.
\end{equation}
When $\alpha=0$, this reduces to the implication in \eqref{equ:charac_smoo_3} of Theorem~\ref{theorem:charac_smoo}.
\end{theorem}
\begin{proof}[of Theorem~\ref{theorem:charac_smoo_n_stronconv}]
Let $g(\bx) \triangleq f(\bx) -\frac{1}{2} \alpha\normtwo{\bx}^2$. Then $\nabla g(\bx) = \nabla f(\bx) -\alpha\bx$.
Using (iii) in Theorem~\ref{theorem:charac_stronconv} and \eqref{equ:charac_smoo_3} in Theorem~\ref{theorem:charac_smoo}, this implies $g(\bx)$ is $(\beta-\alpha)$-smooth and 
$$
\frac{1}{\beta-\alpha} \normtwo{\nabla g(\bx) - \nabla g(\by)}^2 
\leq  \innerproduct{\nabla g(\bx) - \nabla g(\by), \bx - \by},
$$
which is equivalent to the result.
\end{proof}

\begin{figure}
\centering       
\vspace{-0.25cm}                 
\subfigtopskip=0pt               
\subfigbottomskip=0pt         
\subfigcapskip=0pt      
\includegraphics[width=0.98\textwidth]{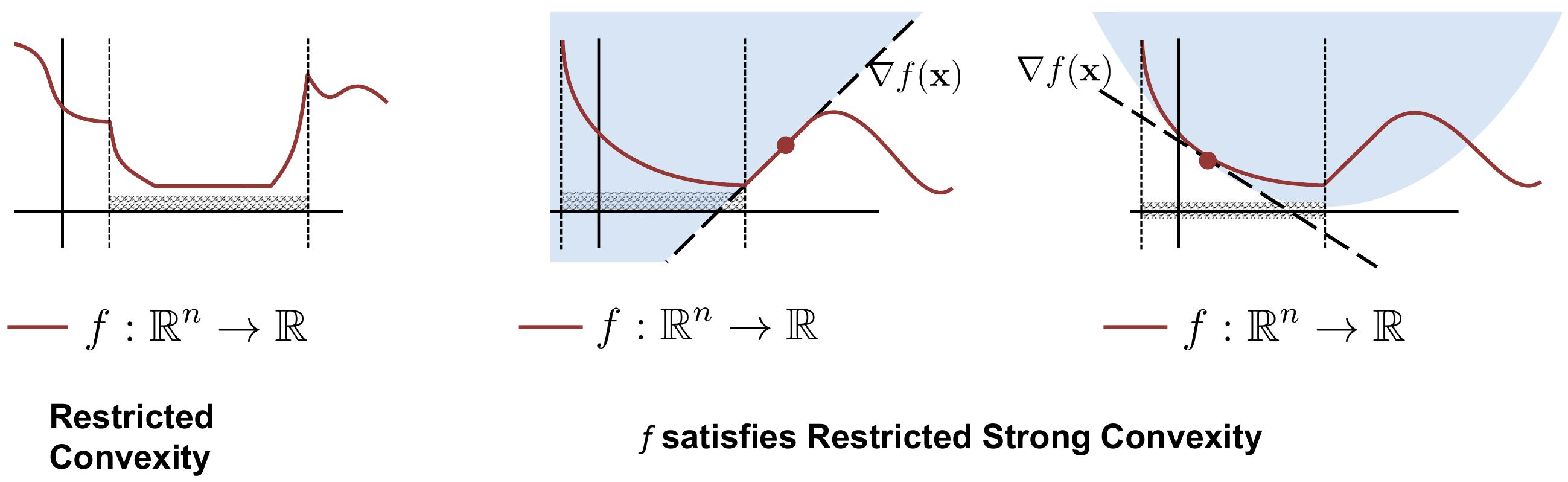}
\caption{
An illustration of restricted convexity properties. The first function $f$ exhibits non-convex behavior over the entire real line, yet it shows convex characteristics within the shaded region delineated by the dotted vertical lines. The second and third functions are also non-convex overall but demonstrate restricted strong convexity. Outside the shaded area (once more, marked by the dotted vertical lines), these functions fail to maintain convexity, as evidenced by their curves falling below their tangents. However, within the defined region, they exhibit strong convexity. Figure is adapted from \citet{jain2017non}.}
\label{fig:rssfunc}
\end{figure}

In non-convex optimization problems, we may also be interested  in the concept of restricted convexity/smoothness for functions.
\begin{definition}[Restricted Strong Convexity/Smoothness]\label{definition:res_scss_func}
Let $f:\sS\subseteq \real^n\rightarrow \real$ be a  differentiable function. 
Then $f$ is said to satisfy \textit{restricted convexity} over a (possibly non-convex) set $\sS\subseteq  \real^n$ if for every $\bx,\by\in\sS$, it follows that 
\begin{equation}\label{equation:res_scss_func0}
f(\by) \geq f(\bx) +\nabla f(\bx)^\top (\by-\bx).
\end{equation}
Additionally, $f$ is said to satisfy \textit{$\alpha$-restricted strong convexity (RSC)} and \textit{$\beta$-restricted strong smoothness (RSS)} over a (possibly non-convex) set $\sS\subseteq  \real^n$ if for every $\bx,\by\in\sS$, it follows that 
\begin{equation}\label{equation:res_scss_func1}
\frac{\alpha}{2} \normtwo{\bx-\by}^2
\leq f(\by)-f(\bx)-\nabla f(\bx)^\top (\by-\bx)
\leq \frac{\beta}{2} \normtwo{\bx-\by}^2.
\end{equation}
When the function is twice differentiable, Theorem~\ref{theorem:psd_hess_conv} also indicates  that 
\begin{equation}\label{equation:res_scss_func2}
\alpha\bI \preceq \nabla^2 f(\bx) \preceq \beta\bI,
\end{equation}
for any $\bx\in\interior(\sS)$ and there exists a vector $\by\in\real^n$ such that $\bx+\lambda\by\in\sS$ for some  sufficiently small $\lambda$ (see proof of Theorem~\ref{theorem:psd_hess_conv}).
That is, 
\begin{equation}\label{equation:res_scss_func3_ra}
\alpha\leq \frac{\bv^\top \nabla^2 f(\bx) \bv}{\normtwo{\bv}^2} \leq \beta,
\end{equation}
for any $\bx\in\interior(\sS)$ and $\bv\in\real^n$ such that $\bx+\lambda\bv\in\sS$ for some sufficiently small $\lambda$.
\end{definition}
Figure~\ref{fig:rssfunc} depicts the situation of restricted convexity, RSC, RSC. 
We note that when $f(\bx) = \frac{1}{2m}\normtwo{\bA\bx-\bb}^2$, where $\bA\in\real^{m\times n}$ and $\bb\in\real^m$, \eqref{equation:res_scss_func3_ra} becomes $\alpha\leq \frac{\bv^\top (\frac{1}{m}\bA^\top\bA) \bv}{\normtwo{\bv}^2} \leq \beta$. 
This yields the RSC and RSS properties for sparse problems (Definition~\ref{definition:res_scss_mat}).

\section{Projection, Proximal, Bregman-Proximal, and Separation}\label{section:proj_prox_sep}

In this section, we explore the fundamental concepts of projection and proximal operators, which are essential tools in optimization theory and convex analysis. These concepts are particularly important for projected/proximal gradient descent methods discussed in Chapter~\ref{chapter:gd_convg}, as well as for proving various theorems.

To begin, let us define a divergence measure based on a convex and differentiable function.
\begin{definition}[Bregman Distance]\label{definition:breg_dist}
The \textit{Bregman distance (or Bregman divergence)} generalizes the concept of squared Euclidean distance and measures the difference between two points in a space.
It is defined for a  convex and differentiable function $ \phi: \sS \rightarrow \real $, where $ \sS \subseteq \real^n $ is a convex set. The Bregman distance from point $ \by $ to point $ \bx $ with respect to $ \phi $ is given by:
$$
\mathcalD_{\phi}(\bx,\by) = \phi(\bx) - \phi(\by) - \innerproduct{\nabla \phi(\by), \bx - \by }.
$$
\end{definition}

For brevity, we summarize the key properties of the Bregman distance in the following remark.
\begin{remark}[Properties of Bregman Distance]\label{remark:bregnan_dist}
The definition of the Bregman distance has the following properties:
\begin{enumerate}
\item \textit{Nonnegativity.}
$
\mathcalD_{\phi}(\bx, \by) \geq 0, \text{for all} \quad \bx, \by \in \sS.
$
When $\phi$ is strictly convex, the quality holds if and only if $ \bx = \by $. 

\item \textit{Convexity in the first argument.}
For a fixed $ \by $, $ \mathcalD_{\phi}(\cdot, \by) $ is convex in the first argument. This property is important in optimization algorithms as it ensures that minimizing $ \mathcalD_{\phi}(\bx, \by) $ over $ \bx $ leads to a well-defined problem.

\item \textit{Not symmetric:}
In general, $ \mathcalD_{\phi}(\bx, \by) \neq \mathcalD_{\phi}(\by, \bx) $. Thus, Bregman distances/divergences do not satisfy the symmetry property required for a distance metric. The term \textit{Bregman distance} is used for historical reasons.

\item \textit{Linearity in the gradient.}
If $ \phi $ is a quadratic function, then the Bregman distance simplifies to a scaled version of the squared Euclidean distance. More generally, if $ \phi $ is linear, the Bregman distance becomes zero.

\item \textit{Pythagorean property.}
Suppose $ \sT $ is a nonempty, closed, and convex subset of $ \sS $, and let $ \projectT^{\phi}(\by) $ denote the projection of $ \by $ onto $ \sT $ with respect to $ \mathcalD_{\phi} $ (Definition~\ref{definition:projec_prox_opt}). Then for any $ \bx \in \sT $,
$$
\mathcalD_{\phi}(\bx, \by) = \mathcalD_{\phi}(\bx, \projectT^{\phi}(\by)) + \mathcalD_{\phi}(\projectT^{\phi}(\by), \by)
$$
This property is analogous to the Pythagorean theorem in Euclidean geometry and is useful in proving convergence results for iterative algorithms.

\item \textit{Three-point property.}
Given three points $ \bx, \by, \bz \in \sS $, the following identity holds:
\begin{equation}\label{equation:three_point_breg}
\mathcalD_{\phi}(\bx, \by) + \mathcalD_{\phi}(\by, \bz) = \mathcalD_{\phi}(\bx, \bz) + \langle \nabla \phi(\bz) - \nabla \phi(\by), \bx - \by \rangle
\end{equation}
This can be used to establish bounds on the Bregman distance.

\item \textit{Local behavior.}
When $\bx$ is near $ \by $, $ \mathcalD_{\phi}(\bx, \by) $ behaves like the second-order Taylor expansion of $ \phi $ around $ \by $ (by the quadratic approximation theorem in Theorem~\ref{theorem:quad_app_theo}). Specifically, if $ \phi $ is twice continuously differentiable, then for small $ \normtwo{\bx - \by} $,
$$
\mathcalD_{\phi}(\bx, \by) \approx \frac{1}{2} (\bx - \by)^\top \nabla^2 \phi(\by) (\bx - \by)
$$
where $ \nabla^2 \phi(\by) $ is the Hessian matrix of $ \phi $ at $ \by $.
\end{enumerate}
\end{remark}

The above remark shows that when $\phi(\bx)$ is strictly convex, the equality holds if and only if $\bx=\by$.
Two commonly used examples of strongly convex functions are provided below.
\begin{example}[Bregman with Euclidean and Negative Entropy]\label{example:breg_examp}
Note that when $\phi(\bx)\triangleq\frac{1}{2}\normtwo{\bx}^2: \real^n\rightarrow (-\infty,\infty]$, $\mathcalD_{\phi}(\bx,\by)$ becomes 
$$
\mathcalD_{\phi}(\bx,\by)=\frac{1}{2}\normtwo{\bx-\by}^2,
$$
which is equivalent to the Euclidean distance.
When $\phi(\bx)\triangleq
\begin{cases}
\sum_{i=1}^{n}x_i\ln (x_i)& \bx\in\real^n_+;\\
\infty&\text{otherwise},
\end{cases} 
$, $\mathcalD_{\phi}(\bx,\by)$ becomes 
$$
\begin{aligned}
\mathcalD_{\phi}(\bx,\by)
&=
\sum_{i=1}^{n} x_i \ln(x_i) - \sum_{i=1}^{n} y_i \ln(y_i)
-
\sum_{i=1}^{n} (\ln(y_i+1))(x_i-y_i)\\
&=\sum_{i=1}^{n} x_i\ln(x_i/y_i),
\end{aligned}
$$
which is also known as the \textit{Kullback-Leibler divergence} between nonnegative vectors $\bx$ and $\by$.
\end{example}

Bregman distances find applications in various fields including optimization, machine learning, statistics, and information theory. They are particularly useful in designing algorithms for convex optimization problems because they provide a way to measure progress towards an optimal solution while respecting the geometry induced by the underlying convex function $\phi$. For instance, Bregman distances play a central role in mirror descent methods, which generalize gradient descent to non-Euclidean geometries (Section~\ref{section:mirror}).

Next, we provide  definitions of the projection and proximal operators, and a variant called \textit{Bregman-proximal} operators.
\begin{definition}[Projection, Proximal, Bregman-Proximal Operators]\label{definition:projec_prox_opt}
Given a \textbf{closed} set $\sS\subseteq \real^n$ and a point $\by\in\real^n$, the \textit{(orthogonal) projection} of $\by$, denoted by $\projectS(\by):\real^n\rightarrow \sS$, is defined as 
\begin{equation}\label{equation:projec_def}
\textbf{(Projection)}:\quad 
\widetildeby\triangleq\projectS(\by)
\triangleq
\mathop{\argmin}_{\bx\in\sS} \normtwo{\bx-\by}
=
\mathop{\argmin}_{\bx\in\sS} \normtwo{\bx-\by}^2.
\end{equation}
Given a \textbf{proper} function $f:\real^n\rightarrow (-\infty,\infty]$ and a point $\by\in\real^n$, the \textit{proximal mapping} of $f$, denoted by $\proxf(\by)$, is defined as
\begin{equation}\label{equation:prox_def}
\textbf{(Proximal)}:\quad 
\widehatby \triangleq \proxf(\by) \triangleq 
\mathop{\argmin}_{\bx\in\real^n} \left( f(\bx) + \frac{1}{2}\normtwo{\bx-\by}^2 \right).
\end{equation}
Equipped with the definition of the indicator function $\indicatorS$  (Example~\ref{example:proper_funcs}),  if $\sS$ is nonempty, the definitions of projection and proximal operators shows that 
\begin{equation}\label{equation:equi_proj_prox}
\prox_{\indicatorS}(\by) = \projectS(\by)\quad \text{for all }\by\in\sS.
\end{equation}
If we use the Bregman distance (Definition~\ref{definition:breg_dist}) instead of the Euclidean distance in \eqref{equation:prox_def}, we have the \textit{Bregman-proximal mapping}, denoted by $\bprox_f(\by)$, as 
\begin{equation}\label{equation:bregprox_def}
\textbf{(Bregman)}:\qquad 
\widebarby \triangleq \bproxfphi(\by) \triangleq 
\mathop{\argmin}_{\bx\in\real^n} \big( f(\bx) + \mathcalD_\phi(\bx, \by) \big),
\end{equation}
where $\phi:\real^n\rightarrow (-\infty,\infty]$ is a proper \textbf{convex} function due to the definition of the Bregman distance (Definition~\ref{definition:breg_dist}). By Example~\ref{example:breg_examp}, 
$$
\bproxfphi(\by) \equiv \proxf(\by)
\quad 
\text{if $\phi(\bx)=\frac{1}{2}\normtwo{\bx}^2$}.
$$
\end{definition}

Unless otherwise stated, we only consider closed sets for projection operators and proper functions for proximal operators.
The following example provides an informal explanation for this restriction.
\begin{example}[Counterexample for Non-Closed $\sS$ and Non-Proper Functions]\label{example:non_clo_proj}
Consider the open set $\sS = \{x \in \real : x > 0\}$ and let $y = 0$. The projection problem then becomes:
$
\widehat{y} = \arg\min_{x > 0} x^2.
$
\begin{itemize}
\item The infimum of $x^2$ is $0$, but $x = 0$ is not in $\sS$. Thus, there is no point in $\sS$ achieving the infimum.
\item Any sequence $\{x_t\}_{t>0} \subset \sS$ with $x_t \rightarrow 0$ minimizes $x^2$, but the ``solution" is not unique since it depends on the sequence chosen.
\end{itemize}
Therefore, projection operators are not well-defined for non-closed sets. Similarly, for non-proper functions, the proximal operator is also not well-defined.
\end{example}

\subsection{Properties of Projection Operators}
We introduce some key properties of projection   operators.
\begin{lemma}[Projection Property-O]\label{lemma:proj_prop0}
Let $\sS\subseteq \real^n$ be \textbf{any closed set}, and let $\by\in\real^n$ such that $\widetildeby\triangleq\projectS(\by)$ is the projection of $\by$ onto the set $\sS$. Then, for all $\bx\in\sS$, we have $\normtwo{\widetildeby - \by}\leq \normtwo{\bx-\by}$.

\end{lemma}
\begin{proof}[of Lemma~\ref{lemma:proj_prop0}]
This proof follows by simply observing that the projection step solves the the optimization
problem $\projectS(\by)=\mathop{\argmin}_{\bx\in\sS} \normtwo{\bx-\by}^2$. 
\end{proof}

Note that the inequality in Projection Property-O holds for all closed sets, whether convex or not.
However, the following three properties necessarily hold only for convex sets.
\begin{lemma}[Projection Property-I]\label{lemma:proj_prop1}
Let $\sS\subseteq \real^n$ be a closed \textbf{convex set}, and let $\by\in\real^n$ such that $\widetildeby\triangleq\projectS(\by)$. Then, for all $\bx\in\sS$, we have $\innerproduct{\bx-\widetildeby, \by-\widetildeby}\leq 0$, i.e., the angle between the two vectors is greater than or equal to 90\textdegree.
\end{lemma}
\begin{proof}[of Lemma~\ref{lemma:proj_prop1}]
Assume for contradiction that $\innerproduct{\bx-\widetildeby, \by-\widetildeby}> 0$ for some $\bx\in\sS$. 
Since $\sS$ is convex and $\widetildeby, \bx \in \sS$, for any $\lambda \in [0, 1]$, we have
$
\bxi_{\lambda} \triangleq \lambda \cdot \bx + (1 - \lambda) \cdot \widetildeby \in \sS.
$
We will now show that  there exists a value of  $\lambda \in [0, 1]$ such that $\normtwo{\by - \bxi_{\lambda}}^2 < \normtwo{\by - \widetildeby}^2$, which contradicts the fact that $\widetildeby$ is the closest point in the convex set to $\by$ and prove the lemma. 
To see this, let $\lambda_1\triangleq\frac{2\innerproduct{\bx - \widetildeby, \by - \widetildeby}}{\normtwo{\bx - \widetildeby}^2}$ for the moment. We have 
$$
\normtwo{\by-\bxi_{\lambda_1}}^2
= 
\normtwo{\by-\widetildeby}^2 - 2\lambda_1 \innerproduct{\bx-\widetildeby, \by-\widetildeby} + \lambda_1^2 \normtwo{\bx-\widetildeby}^2
=\normtwo{\by-\widetildeby}^2.
$$
Let $g(\lambda) \triangleq \normtwo{\by-\widetildeby}^2 - 2\lambda \innerproduct{\bx-\widetildeby, \by-\widetildeby} + \lambda^2 \normtwo{\bx-\widetildeby}^2$. It can be shown that $\lambda=\lambda_1$ does not obtain the minimum of $g(\lambda)$ by the property of quadratic function.
Since we assumed $\langle \bx - \widetildeby, \by - \widetildeby \rangle > 0$,
letting $0 < \lambda < \min\left\{1, \frac{2\innerproduct{\bx - \widetildeby, \by - \widetildeby}}{\normtwo{\bx - \widetildeby}^2}\right\}$, we have $\normtwo{\by - \bxi_{\lambda}}^2 < \normtwo{\by - \widetildeby}^2$. This contradict the property of projection operators, and completes the proof.
\end{proof}

\begin{lemma}[Projection Property-II]\label{lemma:proj_prop2}
Let $\sS\subseteq \real^n$ be a closed  \textbf{convex set}, and let $\by\in\real^n$ such that $\widetildeby\triangleq\projectS(\by)$. Then, for all $\bx\in\sS$, we have $\normtwo{\widetildeby - \bx}^2\leq \normtwo{\by-\bx}^2 - \normtwo{\by-\widetildeby}^2$, which also implies $\normtwo{\widetildeby - \bx} \leq \normtwo{\by-\bx}$ (the former property is related to the Pythagorean theorem).
\end{lemma}
\begin{proof}[of Lemma~\ref{lemma:proj_prop2}]
We have the following elementary inequalities
$$
\begin{aligned}
\normtwo{\by-\bx}^{2} &= \normtwo{(\widetildeby-\bx)-(\widetildeby-\by)}^{2} 
= \normtwo{\widetildeby-\bx}^{2} + \normtwo{\widetildeby-\by}^{2} - 2\langle\widetildeby-\bx, \widetildeby-\by\rangle \\
&\stackrel{\dag}{\geq} \normtwo{\widetildeby-\bx}^{2} + \normtwo{\widetildeby-\by}^{2}
\geq \normtwo{\widetildeby-\bx}^{2},
\end{aligned}
$$
where inequality ($\dag$) follows from the Projection Property-I.
\end{proof}

\begin{lemma}[Projection Property-III]\label{lemma:proj_prop3}
Let $\sS\subseteq \real^n$ be a closed  \textbf{convex set}, and let $\by\in\real^n$. Then, $\projectS(\by)
\triangleq
\mathop{\argmin}_{\bx\in\sS} \normtwo{\bx-\by}^2$ has a \textbf{unique} optimal solution.
\end{lemma}
\begin{proof}[Lemma~\ref{lemma:proj_prop3}]
Since  the  objective function in $\projectS(\by)
\triangleq
\mathop{\argmin}_{\bx\in\sS} \normtwo{\bx-\by}^2$ is a quadratic function associated with a positive definite
matrix, it follows by Exercise~\ref{exercise:coerci_quad} that the objective function is coercive and hence, by Theorem~\ref{theorem:att_coer}, that the problem has at least one optimal solution. In addition, since the
objective function is strictly convex (again, since the objective function is quadratic associated with a
positive definite matrix), it follows by Theorem~\ref{theorem:stric_op_str_conv} that there exists only one optimal
solution.
\end{proof}

\begin{figure}[h]
\centering       
\vspace{-0.25cm}                 
\subfigtopskip=2pt               
\subfigbottomskip=-2pt         
\subfigcapskip=-10pt      
\includegraphics[width=0.98\textwidth]{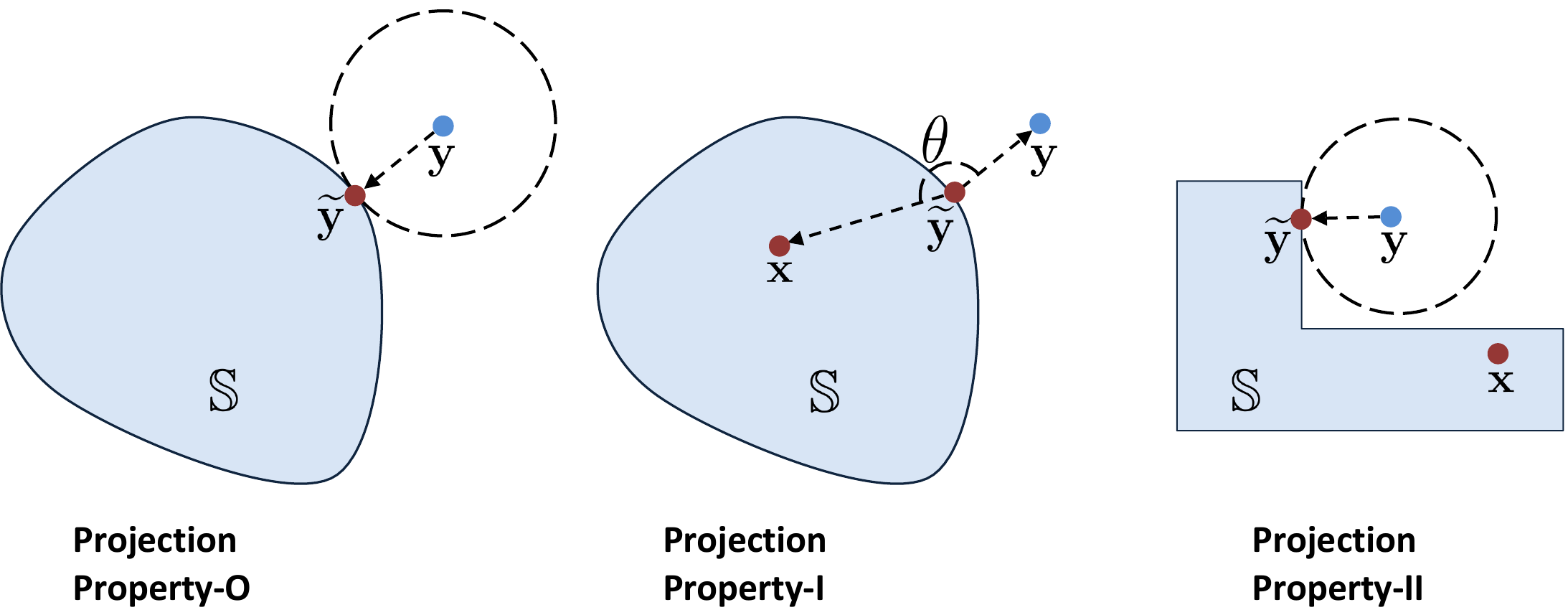}
\caption{
Projection operators identify the closest point within a set to a given point outside the set. For convex sets, Projection Property-I guarantees that the angle $\theta$ formed between the line segment connecting the original point to its projection and any line segment from the projection to another point in the set is never acute. Sets that adhere to Projection Property-I also comply with Projection Property-II, which asserts that projecting onto these sets will not increase the distance between the projected point $\by$ and any other point within the set.
However, non-convex sets may violate Projection Property-II. 
}
\label{fig:proj_prop}
\end{figure}

\paragraph{Convex set.} Note that Projection Properties-I, II, and III are also called first-order properties and can be violated if the underlying set is non-convex. However, Projection Property-O, often called a zeroth-order property, always holds regardless of whether the underlying set is convex or not; see Figure~\ref{fig:proj_prop}.

\begin{example}[Projection onto Subsets of $\real^n$]
Below are examples of nonempty closed and convex sets along with their corresponding orthogonal projections:
\begin{enumerate}[(i)]
\item \textbf{Nonnegative orthant $\sS_1 = \real^n_+$}: $\project_{\sS_1}(\by) = \{\max\{y_i, 0\}\}_{i=1}^n$.
\item \textbf{Box $\sS_2 = \text{Box}[\bl, \bu]\triangleq\{\bx\in\real^n\mid \bl\leq \bx\leq \bu\}$}: $\project_{\sS_2}(\by) = \left\{\min(\max(y_i, l_i), u_i)\right\}_{i=1}^n$.
\item \textbf{Closed $\ell_2$ Ball $\sS_3 = \sB_2[\ba, \alpha]$}: $\project_{\sS_3}(\by)=\ba + \frac{\alpha}{\max\{\normtwo{\by - \ba}, \alpha\}}(\by - \ba)$.
\item \textbf{Half-space $\sS_4 = \{\bx\mid \bc^\top \bx \leq \beta\}$}: $\project_{\sS_4}(\by)=\by - \frac{\max\{\bc^\top \by - \beta, 0\}}{\normtwo{\bc}^2}\bc$.
\item \textbf{Affine set $\sS_5 = \{\bx \in \real^n \mid \bA\bx = \bb\}$}: $\project_{\sS_5}(\by) = \by - \bA^\top(\bA\bA^\top)^{-1}(\bA\by - \bb)$.
\end{enumerate}
where $\bl \in  [-\infty, \infty)^n, \bu \in (-\infty, \infty]^n$ satisfying $\bl \leq \bu$, $\ba \in \real^n$, $\alpha > 0$, $\bc \in \real^n \setminus \{\bzero\}$,  $\beta \in \real$, $\bA \in \real^{m \times n}$ has full row rank such that $\bA\bA^\top$ is nonsingular, and $\bb \in \real^m$.
The results in (i)--(iv) holds trivially. 
For (v),
it's equivalent to solving  the following optimization problem:
$$
\min_{\bx \in \real^n} \quad \frac{1}{2} \normtwo{ \bx - \by}^2 \quad\text{s.t.}\quad\bA\bx = \bb,
$$
where $\bA \in \real^{m \times n}$, $\bb \in \real^m$ and $\by \in \real^n$ are given matrices and vectors. 
We will use the result of KKT conditions to prove this (Section~\ref{section:gen_kkt_cond}); feel free to skip for a first reading.
We assume that the matrix $\bA$ has full row rank (for matrices that do not have full row rank, we can remove linearly dependent rows to obtain a reduced constraint condition). 
For the equality constraint, we introduce the Lagrange multiplier $\blambda \in \real^m$ and construct the Lagrangian function
$$
L(\bx, \blambda) = \frac{1}{2} \normtwo{\bx - \by}^2 + \blambda^\top(\bA\bx - \bb).
$$
Since the quadratic function is a convex function and there is only an equality constraint, the Slater's condition is satisfied. Theorem~\ref{theorem:kktslat1_slater} shows that $\bx^*$ is a global optimum if and only if there exists $\blambda^* \in \real^m$ such that
$
\begin{cases}
\bx^* - \by + \bA^\top\blambda = \bzero, \\
\bA\bx^* = \bb.
\end{cases}
$
This yields that 
$$
\bA \bx^* - \bA \by + \bA\bA^\top\blambda = \bzero
\qquad \implies\qquad 
\blambda = (\bA\bA^\top)^{-1}(\bA\by - \bb).
$$
Substituting $\blambda$ back into the first KKT condition obtains the desired result.
\end{example}

\subsection{Properties of Proximal Operators}
Analogously, for the proximal operator, if we require the function to be closed, it can be shown under some coerciveness assumptions (Definition~\ref{definition:coerciveness}) that $\proxf(\by)$ is never an empty set.
\begin{lemma}[Proximal Property-O: Nonemptiness of  Prox under Closedness and Coerciveness]\label{lemma:prox_prop0}
Let $f: \real^n \rightarrow (-\infty, \infty]$ be a proper \textbf{closed} function, and assume that the following condition holds:
$$
\text{the function } \bx \mapsto f(\bx) + \frac{1}{2} \normtwo{\bx - \by}^2 \text{ is coercive for any } \by \in \real^n. 
$$
Then $\proxf(\by)$ is \textbf{nonempty} for any $\by \in \real^n$.
\end{lemma}
\begin{proof}[of Lemma~\ref{lemma:prox_prop0}]
For any $\bx \in \real$, the proper function $g(\bx) \triangleq f(\bx) + \frac{1}{2} \normtwo{\bx - \by}^2$ is closed as a sum of two closed functions. 
Since by assumption it is also coercive, Theorem~\ref{theorem:att_coer} (which requires the function to be proper closed and coercive) ensures that $\proxf(\by)$, which consists of the minimizers of $g$, is nonempty.
\end{proof}

\begin{lemma}[Proximal Property-I]\label{lemma:prox_prop1}
Let $ f: \real^n \rightarrow (-\infty, \infty] $ be a proper \textbf{closed and convex} function. Then for any $ \by, \widehatby \in \real^n $, the following three claims are equivalent:
\begin{enumerate}[(i)]
\item $ \widehatby = \proxf(\by) $.
\item $ \by - \widehatby \in \partial f(\widehatby) $.
\item $ \innerproduct{\by - \widehatby, \bz - \widehatby} \leq f(\bz) - f(\widehatby) $ for any $ \bz \in \real^n $.
\end{enumerate}
\end{lemma}
\begin{proof}[of Lemma~\ref{lemma:prox_prop1}]
By definition, $ \widehatby = \proxf(\by) $ if and only if $ \widehatby $ is the minimizer of the problem
$
\mathopmin{\bv} \left\{ f(\bv) + \frac{1}{2} \normtwo{\bv - \by}^2 \right\},
$
which, by the optimality condition (Theorem~\ref{theorem:fetmat_opt}) and the sum rule of subdifferential calculus, is equivalent to the relation
$
\bzero \in \partial f(\widehatby) + \widehatby - \by. 
$
This shows the equivalence between claims (i) and (ii). Finally, by the definition of the subgradient (Definition~\ref{definition:subgrad}), the membership relation of claim (ii) is equivalent to (iii).
\end{proof}

When $ f = \indicatorS $, with $ \sS $ being a nonempty closed and convex set, the equivalence between claims (i) and (iii) in the Proximal Property-I amounts to the Projection Property-I (Lemma~\ref{lemma:proj_prop1}).

To maintain consistency, we also introduce the Proximal Property-II, though it's not so that interesting in practice.
\begin{lemma}[Proximal Property-II]\label{lemma:prox_prop2}
Let $ f: \real^n \rightarrow (-\infty, \infty] $ be a proper \textbf{closed and convex} function. Then for any $ \widehatby \triangleq \proxf(\by) $ and $\bz\in\real^n$, it follows that
$
\normtwo{\by-\bz} \geq \normtwo{\widehatby-\bz}^2 + \normtwo{\widehatby-\by}^2 + 2\left( f(\widehatby) - f(\bz) \right).
$
\end{lemma}

\begin{lemma}[Proximal Property-III]\label{lemma:prox_prop3}
Let $f: \real^n \rightarrow (-\infty, \infty]$ be a proper \textbf{closed and convex} function. Then $\proxf(\by)$ is a \textbf{singleton} for any $\by \in \real^n$.
\end{lemma}

\begin{proof}[of Lemma~\ref{lemma:prox_prop3}]
For any $\bx \in \real^n$,
\begin{equation}\label{equation:prox_prop0_e1}
\proxf(\by) = \arg\min_{\bx \in \real^n} g(\by;\bx) ,
\end{equation}
where $g(\by;\bx) = f(\bx) + \frac{1}{2}\normtwo{\bx - \by}^2$. The function $g(\by;\bx)$ is a closed and strongly convex function as a sum of the closed and strongly convex function $\frac{1}{2}\normtwo{\bx - \by}^2$ and the closed and convex function $f$ (Exercise~\ref{exercise:sum_sc_conv} and Exercise~\ref{exercise:pres_conv_clos}). The properness of $g(\by;\bx)$ immediately follows from the properness of $f$. Therefore, by the SC Property-II (Theorem~\ref{theorem:exi_close_sc}), there exists a unique minimizer to the problem in \eqref{equation:prox_prop0_e1}.
\end{proof}
The Proximal Property-III also proves the Projection Property-III (Lemma~\ref{lemma:proj_prop3}) by noting that $\prox_{\indicatorS}(\by) = \projectS(\by)$ for any $\by\in\sS$ if $\sS$ is convex and closed (the indicator function is closed convex due to Exercise~\ref{exercise_closed_indica}).

When $f$ is proper closed and convex, the preceding lemma shows that $\proxf(\by)$ is a singleton for any $\by \in \real^n$. In these cases, which will constitute the vast majority of cases that will be discussed in this book, we will treat $\proxf$ as a single-valued mapping from $\real^n$ to $\real^n$, meaning that we will write $\proxf(\by) = \widehatby$ instead of $\proxf(\by) = \{\widehatby\}$.

\paragraph{Convexity and Closedness.} Note again  that Proximal Properties-I, II, and III are also called first-order properties and can be violated if the underlying function is non-convex. However, Projection Property-O, often called a zeroth-order property, always holds, whether the underlying function is convex or not.
In all cases, we require the function to be closed to ensure the attainment due to the Weierstrass theorem (Theorem~\ref{theorem:weierstrass_them}) or Theorem~\ref{theorem:att_coer}.

\begin{theorem}[Projection and Proximal Property-IV: Nonexpansiveness]\label{theorem:proj_nonexpan}
Let $\sS\subseteq \real^n$ be a closed \textbf{convex set}, and let $f:\real^n\rightarrow (-\infty, \infty]$ be a proper \textbf{closed and convex} function. Then, 
\begin{enumerate}
\item \textit{Firm nonexpansiveness.} For any $\bx, \by\in\real^n$, it follows that 
$$
\begin{aligned}
\normtwo{\projectS(\bx)-\projectS(\by)}^2 &\leq \innerproduct{\projectS(\bx)-\projectS(\by), \bx-\by};\\
\normtwo{\proxf(\bx)-\proxf(\by)}^2 &\leq \innerproduct{ \proxf(\bx)-\proxf(\by), \bx-\by}.
\end{aligned}
$$
\item \textit{Nonexpansiveness.} For any $\bx, \by\in\real^n$, it follows that 
$$
\begin{aligned}
\normtwo{\projectS(\bx)-\projectS(\by)} &\leq \normtwo{\bx-\by};\\
\normtwo{\proxf(\bx)-\proxf(\by)}&\leq \normtwo{\bx-\by}.
\end{aligned}
$$
\end{enumerate} 
\end{theorem}
\begin{proof}[of Theorem~\ref{theorem:proj_nonexpan}]
\textbf{Projection operator.} By the Projection Property-I (Lemma~\ref{lemma:proj_prop1}),  for any $\bu \in \real^n$ and $\bv \in \sS$, we have that
$
\innerproduct{ \bv - \projectS(\bu), \bu - \projectS(\bu)} \leq 0. 
$
Invoking this inequality with $\bu \triangleq \bx, \bv \triangleq \projectS(\by)$ and $\bu \triangleq \by, \bv \triangleq \projectS(\bx)$:
$$
\begin{aligned}
\innerproduct{\projectS(\by) - \projectS(\bx), \bx - \projectS(\bx)} &\leq 0 
\qquad\text{and}\qquad
\innerproduct{\projectS(\bx) - \projectS(\by), \by - \projectS(\by)} &\leq 0. 
\end{aligned}
$$
Adding the two inequalities yields that
$
\innerproduct{\projectS(\by) - \projectS(\bx),  \bx - \by + \projectS(\by) - \projectS(\bx)} \leq 0,
$
showing the firm nonexpansiveness of projection operators.

To prove the nonexpansiveness, note that if $\projectS(\bx) = \projectS(\by)$, the inequality  holds trivially. We will therefore assume that $\projectS(\bx) \neq \projectS(\by)$. By the Cauchy-Schwarz inequality we have
$$
\innerproduct{\projectS(\bx) - \projectS(\by),  \bx - \by} \leq \normtwo{\projectS(\bx) - \projectS(\by)} \cdot \normtwo{\bx - \by},
$$
which combined with the firm nonexpansivess yields the desired result.

\paragraph{Proximal operator.} Denoting $ \bu \triangleq \proxf(\bx) $, $ \bv \triangleq \proxf(\by) $, by the equivalence of (i) and (ii) in the Proximal Property-I (Lemma~\ref{lemma:prox_prop1}), it follows that
$$
\bx - \bu \in \partial f(\bu)
\qquad \text{and}\qquad  
\by - \bv \in \partial f(\bv).
$$
By the subgradient inequality, we have
$$
f(\bv) \geq f(\bu) + \langle \bx - \bu, \bv - \bu \rangle
\qquad \text{and}\qquad  
f(\bu) \geq f(\bv) + \langle \by - \bv, \bu - \bv \rangle.
$$
Adding the two inequalities yields that
$
0 \geq \innerproduct{\by - \bx + \bu - \bv, \bu - \bv},
$
showing the firm nonexpansiveness of proximal operators.

To prove the nonexpansiveness, if $ \proxf(\bx) = \proxf(\by) $, then the inequality is trivial. Therefore, we assume that $ \proxf(\bx) \neq \proxf(\by) $. By the Cauchy-Schwarz inequality, it follows that
$$
\innerproduct{\proxf(\bx) - \proxf(\by), \bx - \by} \leq \normtwo{\proxf(\bx) - \proxf(\by)} \cdot \normtwo{\bx - \by}.
$$
which combined with the firm nonexpansivess yields the desired result.
\end{proof}
Note that the (firm) nonexpansiveness of projection operators can be  directly inferred from the properties of proximal operators by setting $f\triangleq\indicatorS$.

\begin{SCfigure}
	\centering
	\includegraphics[width=0.55\textwidth]{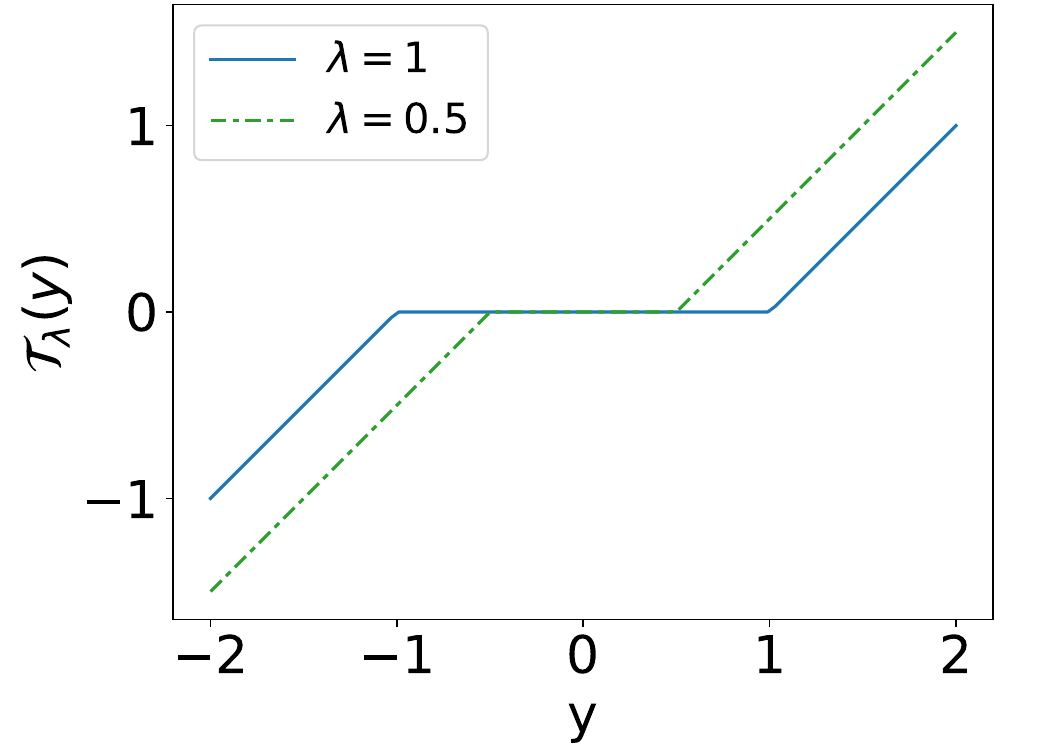} 
	\caption{Illustration of the soft thresholding function $\mathcalT_{\lambda}(y)$ for different values of $\lambda$.}
	\label{fig:soft_threshold}
\end{SCfigure}

\begin{example}[Soft Thresholding, Proximal of $\ell_1$ Norms]\label{example:soft_thres}
Let  $f(\bx) = \lambda \normone{\bx} =\sum_{i=1}^{n} \psi(x_i)$, where $\bx\in\real^n$, $\lambda > 0$, and $\psi(u) = \lambda |u|$. 
Then $
\proxf(\bx) = \big\{\prox_{\psi}(x_i)\big\}_{i=1}^n,
$
where $\prox_{\psi}(y) = \mathcalT_{\lambda}(y)$, and  $\mathcalT_{\lambda}$ is known as the \textit{soft thresholding function} and is defined as
$$
\mathcalT_{\lambda}(y) \triangleq
[\abs{y} - \lambda]_+ \cdot \sign(y)=
\begin{cases}
	y - \lambda, & y \geq \lambda, \\
	0, & |y| < \lambda, \\
	y + \lambda, & y \leq -\lambda.
\end{cases}
$$
where $[u]_+ = \max\{u, 0\}$ for any $u\in\real$. See Figure~\ref{fig:soft_threshold} for an illustration of this soft thresholding function.
More compactly, the proximal operation of $f$ can be denoted as 
$$
\prox_f(\bx) = \mathcalT_\lambda(\bx) = \left\{\mathcalT_\lambda(x_i)\right\}_{i=1}^n 
= [\abs{\bx} - \lambda \bone]_+ \hadaprod \sign(\bx),
$$
where $\hadaprod$ denotes the Hadamard product and $\bone$ denotes a vector of all ones.
\end{example}
\subsection{Properties of Bregman-Proximal Operators}

The zeroth-order property of the Bregman-proximal operator  indicates the unique attainment of a sum of two functions, which is essential for other Bregman-proximal properties.
\begin{lemma}[Bregman-Proximal Property-O: Attainment]\label{lemma:breg_prox_propo}
Let
\begin{itemize}
\item $\phi : \real^n \to (-\infty, \infty]$ be a proper \textbf{closed and convex} function differentiable over $\dom(\partial \phi)$.
\item $f : \real^n \to (-\infty, \infty]$ be a proper \textbf{closed and convex} function satisfying $\dom(f) \subseteq \dom(\phi)$.
\item $\phi + \delta_{\dom(f)}$ be $\sigma$-strongly convex ($\sigma > 0$).
\end{itemize}
\noindent
Then the minimizer of the problem
\begin{equation}\label{equation:breg_prox_propoe1}
\min_{\bx \in \real^n} \{f(\bx) + \phi(\bx)\}
\end{equation}
is \textbf{uniquely} attained at a point in $\dom(f) \cap \dom(\partial \phi)$.
\end{lemma}
\begin{proof}[of Lemma~\ref{lemma:breg_prox_propo}]
Problem \eqref{equation:breg_prox_propoe1} is the same as
\begin{equation}\label{equation:breg_prox_propoe2}
\min_{\bx \in \real^n} \omega(\bx), \quad \text{where $\omega \triangleq f + \phi$}.
\end{equation}
The function $\omega$ is closed since both $f$ and $\phi$ are closed (Exercise~\ref{exercise:pres_conv_clos}); it is proper by the fact that $\dom(\omega) = \dom(f) \neq \emptyset$. Since $\phi + \delta_{\dom(f)}$ is $\sigma$-strongly convex and $f$ is convex, their sum $f + \phi + \delta_{\dom(f)} = f + \phi = \omega$ is $\sigma$-strongly convex. To conclude, $\omega$ is proper closed and $\sigma$-strongly convex, and hence, by the SC Property-II (Theorem~\ref{theorem:exi_close_sc}), problem \eqref{equation:breg_prox_propoe2} has a unique minimizer $\bx^*$ in $\dom(\omega) = \dom(f)$. To show that $\bx^* \in \dom(\partial \phi)$, note that by first-order condition of convex optimization problem  (Theorem~\ref{theorem:fetmat_opt}), $\bzero \in \partial \omega(\bx^*)$, and in particular $\partial \omega(\bx^*) \neq \emptyset$. Therefore,  by the sum rule of subdifferential calculus, $\partial \omega(\bx^*) = \partial f(\bx^*) + \partial \phi(\bx^*)$, it follows in particular that $\partial \phi(\bx^*) \neq \emptyset$, meaning that $\bx^* \in \dom(\partial \phi)$.
\end{proof}

For the Bregman-proximal operator, we have a similar property as Lemma~\ref{lemma:prox_prop1} for the proximal operator.
\begin{lemma}[Bregman-Proximal Property-I]\label{lemma:breg_prox_prop1}
Let 
\begin{itemize}
\item $\phi : \real^n \to (-\infty, \infty]$ be a proper \textbf{closed and convex} function differentiable over $\dom(\partial \phi)$;
\item $f : \real^n \to (-\infty, \infty]$ be a proper \textbf{closed and convex} function satisfying $\dom(f) \subseteq \dom(\phi)$;
\item $\phi + \delta_{\dom(f)}$ be $\sigma$-strongly convex ($\sigma > 0$).
\end{itemize}
\noindent Then for any $ \by, \widehatby \in \dom(\partial \phi) $, the following two claims are equivalent:
\begin{enumerate}[(i)]
\item $ \widebarby = \bproxfphi(\by) \triangleq \mathop{\argmin}_{\bx \in \real^n} \left\{ f(\bx) + \mathcalD_\phi(\bx, \by) \right\}$ .
\item $ \innerproduct{\nabla \phi(\by) - \nabla \phi(\widebarby), \bz - \widebarby} \leq f(\bz) - f(\widebarby)$ for any $ \bz \in \dom(f) $.
\end{enumerate}
\end{lemma}
\begin{proof}[of Lemma~\ref{lemma:breg_prox_prop1}]
By the definition of the Bregman distance (Definition~\ref{definition:breg_dist}), $\widebarby = \bproxfphi(\by)$ can be rewritten as
\begin{equation}
\widebarby = \mathop{\argmin}_{\bx \in \real^n} \left\{ f(\bx) - \innerproduct{\nabla \phi(\by), \bx} + \phi(\bx) \right\}.
\end{equation}
The fact that $\widebarby \in \dom(\partial \phi)$ follows by invoking Lemma~\ref{lemma:breg_prox_propo} with $f(\bx) - \innerproduct{\nabla \phi(\by), \bx}$ taking the role of $f(\bx)$. Using the first-order condition of convex optimization problems (Theorem~\ref{theorem:fetmat_opt}), it follows by $\widebarby = \bproxfphi(\by)$ that there exists $f'(\widebarby) \in \partial f(\widebarby)$ for which
$
f'(\widebarby) + \nabla \phi(\widebarby) - \nabla \phi(\by) = \bzero.
$
Hence, by the subgradient inequality (Definition~\ref{definition:subgrad}), for any $\bz \in \dom(f)$,
$$
\innerproduct{\nabla \phi(\by) - \nabla \phi(\widebarby), \bz - \widebarby} = \innerproduct{f'(\widebarby), \bz - \widebarby} \leq f(\bz) - f(\widebarby),
$$
which completes the proof.
\end{proof}

\subsection{Separation Theorems}
The following strict separation theorem, derived from the properties of projection operators, is important for proving  Farkas' lemma (Lemma~\ref{lemma:farka_lemm}), which in turn is crucial for proving KKT conditions.
\begin{theorem}[Strict Separation Theorem]\label{theorem:stric_sep_theo}
Let $\sS \subseteq \real^n$ be a closed \textbf{convex} set, and let $\by \notin \sS$. Then there exist $\bw \in \real^n \setminus \{\bzero\}$ and $\alpha \in \real$ such that
$$
\bw^\top \by > \alpha
\qquad\text{and}\qquad
\bw^\top \bx \leq \alpha, \quad \text{for all } \bx \in \sS.
$$
\end{theorem}
\begin{proof}[of Theorem~\ref{theorem:stric_sep_theo}]
Let $\widetildeby \triangleq \projectS(\by) \in \sS$. By the Projection Property-I (Lemma~\ref{lemma:proj_prop1}), 
$$
\innerproduct{\bx - \widetildeby, \by - \widetildeby} \leq 0
\quad\implies\quad (\by - \widetildeby)^\top \bx \leq (\by - \widetildeby)^\top \widetildeby
,\quad  \text{for all } \bx \in \sS.
$$
Denote $\bw \triangleq \by - \widetildeby \neq \bzero$ (since $\by \notin \sS$) and $\alpha \triangleq (\by - \widetildeby)^\top \widetildeby$. Then we have that $\bw^\top \bx \leq \alpha$ for all $\bx \in \sS$. On the other hand,
$$
\bw^\top \by = (\by - \widetildeby)^\top \by = (\by - \widetildeby)^\top (\by - \widetildeby) + (\by - \widetildeby)^\top \widetildeby =\normtwo{\by - \widetildeby}^2 + \alpha > \alpha,
$$
establishing the desired result.
\end{proof}

More generally, we have the following separation theorem.
\begin{theorem}[General Separation Theorem]\label{theorem:stric_sep_theo2}
Let $\sS, \sT \subseteq \real^n$ be \textbf{convex} sets with empty intersection $\sS \cap \sT = \varnothing$. Then there exists a point $\ba \in \real^n$ and a number $\beta \in \real$ such that
\begin{enumerate}[(i)]
\item For all $\bx \in \sS$, we have $\innerproduct{\ba, \bx} \geq \beta$.
\item For all $\bx \in \sT$, we have $\innerproduct{\ba, \bx} \leq \beta$.
\end{enumerate}
If $\sS$ and $\sT$ are closed and at least one of them is bounded, then we can replace the inequalities by strict inequalities.
\end{theorem}

The case we're most concerned with is when both sets are compact (i.e., closed and bounded). We highlight its proof here.
\begin{proof}[of Theorem~\ref{theorem:stric_sep_theo2} for Compact Sets]
In this case, the Cartesian product $\sS \times \sT$ is also compact. Therefore, the distance function $\normtwo{\bx - \by}$ attains its minimum over $\sS \times \sT$ (Theorem~\ref{theorem:weierstrass_them}). Taking $\bu, \bv$ to be two points that achieve the minimum. A separating hyperplane is given by the hyperplane perpendicular to $\bv - \bu$ that passes through the midpoint between $\bu$ and $\bv$. That is, $\ba = \bv - \bu$ and $\beta =\frac{-\innerproduct{\ba, \bv+\bu}}{2}$ (use the result of Problem~\ref{prob:dist_hyper} to prove this). 
For the sake of contradiction, suppose there is a point $\bz$ on this hyperplane contained in one of the two sets, say, $\sS$. Then the line segment from $\bu$ to $\bz$ is also contained in $\sS$ by convexity. We can then find a point along the line segment that is closer to $\bv$ than $\bu$ is, thus contradicting our assumption.
\end{proof}

\begin{theorem}[Supporting Hyperplane Theorem]\label{theorem:supp_hyp_theo}
Let $\sS \subseteq \real^n$ be a \textbf{convex} set (not necessarily closed) and let $\by \notin \sS$. Then there exists $\bzero \neq \bw \in \real^n$ such that
$$
\bw^\top \bx \leq \bw^\top \by, \quad  \text{ for any } \bx \in \sS.
$$
\end{theorem}
\begin{proof}[of Theorem~\ref{theorem:supp_hyp_theo}]
Since the set is not closed, we cannot apply the properties of projection operators directly; however, we can apply these properties to the closure $\closure(\sS)$.
Since $\by \notin \interior(\sS)$, it follows that $\by \notin \interior(\closure(\sS))$ (by Exercise~\ref{exercise:int_clos}, $\interior(\sS) = \interior(\closure(\sS))$). Therefore, there exists a sequence $\{\by^\toptzero\}_{t \geq 1}$ satisfying $\by^\toptzero \notin \closure(\sS)$ such that $\by^\toptzero \rightarrow \by$. Since $\closure(\sS)$ is convex by Exercise~\ref{exercise:clo_conv} and closed by its definition, it follows by the strict separation theorem (Theorem~\ref{theorem:stric_sep_theo}) that there exists $\bzero \neq \bw^\toptzero \in \real^n$ such that
$$
(\bw^\toptzero)^\top \bx < (\bw^\toptzero)^\top \by^\toptzero,\quad \text{for all } \bx \in \closure(\sS), \text{ where } \by^\toptzero \notin \closure(\sS).
$$
Dividing the  inequality by $\normtwo{\bw^\toptzero} \neq \bzero$, we obtain
\begin{equation}\label{equation:supp_eq1}
\frac{(\bw^\toptzero)^\top}{\normtwo{\bw^\toptzero}}(\bx - \by^\toptzero) < 0,\quad \text{ for any } \bx \in \closure(\sS), \text{ where } \by^\toptzero \notin \closure(\sS).
\end{equation}
Since the sequence $\left\{\frac{\bw^\toptzero}{\normtwo{\bw^\toptzero}}\right\}_{t \geq 1}$ is bounded, it follows that there exists a subsequence $\left\{\frac{\bw^\toptzero}{\normtwo{\bw^\toptzero}}\right\}_{t \in \sT}$ indexed by $\sT$ such that $\frac{\bw^\toptzero}{\normtwo{\bw^\toptzero}} \rightarrow \bw$ as $t \rightarrow \infty$ for some $\bw \in \real^n$. Obviously, $\normtwo{\bw} = 1$ and hence in particular $\bw \neq \bzero$. Taking the limit as $t \rightarrow \infty$ in inequality \eqref{equation:supp_eq1}, we obtain that
$$
\bw^\top(\bx - \by) \leq 0, \quad  \text{ for any } \bx \in \closure(\sS),
$$
which completes the proof since $\sS \subseteq \closure(\sS)$.
\end{proof}

\section{Optimization and Optimality Conditions}

Optimization holds a pivotal role across various scientific disciplines and practical applications, including economics, operations research, network analysis, and the optimal design of mechanical or electrical systems, among others.
In this book, we will explore numerical methods designed to solve continuous optimization problems. In these problems, the objective is formulated as a real-valued function of multiple variables, with the goal of identifying the set of input values that yields the minimum possible function value. Such optimization challenges often emerge from parameter estimation tasks, including data fitting as a specific application.
\begin{definition}[Optimization Problem]\label{definition:opt_probs_all}
Let $f: \real^n \rightarrow \real$, the \textit{unconstrained optimization problem} is defined as
$$
\textbf{(P1)}:\qquad \text{Find}\quad \bx^* = \mathop{\argmin}_{\bx} f(\bx).
$$
For finding the minimum over a set $\sS$, the \textit{constrained optimization problem} becomes
$$
\textbf{(P2)}:\qquad \text{Find}\quad \bx^* = \mathop{\argmin}_{\bx} f(\bx)\quad \text{s.t.}\quad\bx\in\sS.
$$
Similarly, we may also want to find the minimum point of a composite function:
$$
\textbf{(P3)}:\qquad \text{Find}\quad \bx^* = \mathop{\argmin}_{\bx} \left\{F(\bx)\triangleq f(\bx)+g(\bx) \right\},
$$
where $f$ possesses some differentiability properties (e.g., smoothness) but is not assumed to be convex, while $g$ is convex but may lack special differentiability properties.
When the $\sS$ in the problem (P2) has the following form, we deal with the constrained optimization problem with equality and inequality constraints:
$$
\begin{aligned}
\textbf{(P4)}:\qquad \text{Find}\quad 
&\bx^* &=& \mathop{\argmin}_{\bx} f(\bx) \\
&\text{s.t.}  &c_i&(\bx)=0, \,\,i\in\mathcalE, \,\,\text{with $\abs{\mathcalE}=p$};\\
&  &c_i&(\bx)\leq 0, \,\,i\in\mathcalI, \,\,\text{with $\abs{\mathcalI}=q$}.\\
\end{aligned}
$$
Since we focus on minimization, we assume in all cases, unless otherwise stated, that the underlying  function is proper (Definition~\ref{definition:extrel_pro_clo_funcs}).
\end{definition}
In some cases we want a maximizer of a function. This is easily determined if we find a minimizer of the function with opposite sign.
In this sense, problem (P1) is called a \textit{unconstrained optimization}, (P2) is called a \textit{constrained optimization}, and (P3) is called a non-convex composite problem. 
The function $f$ or $F$ is called the \textit{objective function} or \textit{cost function}, and $\bx^*$ is the \textit{minimum point (or simply minimizer)}.
Different forms of $g$  lead to different optimization types:
\begin{itemize}
\item When $g(\bx)=0$,  problem (P3) reduces to an unconstrained problem (P1).
\item When $g(\bx)=\indicatorS(\bx)$,  problem (P3)  transforms into a constrained problem (P2).
\item When $g(\bx)=\lambda\normone{\bx}$,  problem (P3) becomes an $\ell_1$ regularized problem. 
\end{itemize}

Regarding (P3), it is assumed that $f$ can be accessed through a first-order oracle and that $g$ is known and ``simple," where simplicity will become clearer upon examining the algorithm description in Chapter~\ref{chapter:gd_convg}.

The ideal scenario in optimization occurs when the objective function has a single, unique minimizer known as the global minimizer. However, some functions may have multiple minimizers, or even an infinite number, e.g., a cosine function.  In these cases, identifying any one of these minimizers can be sufficient for solving the problem at hand.

In many practical applications, objective functions feature both a global minimizer and numerous local minimizers. Developing algorithms that reliably locate the global minimizer in such scenarios presents a significant challenge.
Throughout this book, we will focus on methods designed to identify a local minimizer of the objective function. Once a minimizer is found, it remains uncertain whether it is the global minimizer or merely one of several local minimizers. Additionally, there is no guarantee that the algorithm will find the local minimizer closest to the initial starting point. To explore various local minimizers, one approach is to conduct multiple runs with different starting points. Ideally, one could also analyze intermediate results from algorithms designed for global optimization.

A \textit{local minimizer (or local minimum point)} for a function $f$ is an argument vector giving the smallest function value inside a certain region, defined by $\epsilon$; while a \textit{global minimizer (or global minimum point)} is an argument vector given the smallest function value over the feasible set. 
In general, local minimizer and global minimizer are collectively referred to as \textit{minimizers or optimal points}.
Rigorously, we have the following definition.
\begin{definition}[Local/Global Minimum and Maximum]\label{definition:local_global_opt}
Let $f : \sS\subseteq \real^n \rightarrow \real$ be defined on a set $\sS \subseteq \real^n$. Then,
\begin{enumerate}
\item $\bx^* \in \sS$ is called a \textit{local minimum point} of $f$ over $\sS$ if there exists $\epsilon > 0$ for which $f(\bx^*) \leq f(\bx)$ for any $\bx \in \sS \cap \sB(\bx^*, \epsilon)$,
\item $\bx^* \in \sS$ is called a \textit{strict local minimum point} of $f$ over $\sS$ if there exists $\epsilon > 0$ for which $f(\bx^*) < f(\bx)$ for any $\bx \neq \bx^* \in \sS \cap \sB(\bx^*, \epsilon)$,
\item $\bx^* \in \sS$ is called a \textit{local maximum point} of $f$ over $\sS$ if there exists $\epsilon > 0$ for which $f(\bx^*) \geq f(\bx)$ for any $\bx \in \sS \cap \sB(\bx^*, \epsilon)$,
\item $\bx^* \in \sS$ is called a \textit{strict local maximum point} of $f$ over $\sS$ if there exists $\epsilon > 0$ for which $f(\bx^*) > f(\bx)$ for any $\bx \neq \bx^* \in \sS \cap \sB(\bx^*, \epsilon)$.
\end{enumerate}
Similarly,
\begin{enumerate}
\item $\bx^* \in \sS$ is called a \textit{global minimum point} of $f$ over $\sS$ if $f(\bx) \geq f(\bx^*)$ for any $\bx \in \sS$,
\item $\bx^* \in \sS$ is called a \textit{strict global minimum point} of $f$ over $\sS$ if $f(\bx) > f(\bx^*)$ for any $\bx \neq \bx^* \in \sS$,
\item $\bx^* \in \sS$ is called a \textit{global maximum point} of $f$ over $\sS$ if $f(\bx) \leq f(\bx^*)$ for any $\bx \in \sS$,
\item $\bx^* \in \sS$ is called a \textit{strict global maximum point} of $f$ over $\sS$ if $f(\bx) < f(\bx^*)$ for any $\bx \neq \bx^* \in \sS$.
\end{enumerate}
\end{definition}

\subsection*{Fundamental Theorems for Optimality Conditions}
For obtaining a \textit{global optimal point} (a.k.a. a global extremum point or a global optimizer, which can be either a global minimum or maximum), the Weierstrass theorem provides conditions under which such optimizers can be attained.
We  present  the Weierstrass theorem and its variants for optimization problems involving continuous or proper closed functions.
\begin{theorem}[Weierstrass Theorem and Variants]\label{theorem:weierstrass_them}
We consider the Weierstrass theorem and its variants for different types of functions and sets:
\begin{enumerate}[(i)]
\item \label{weier1_continus} Let $f:\sS\rightarrow (-\infty, \infty]$ be a \textbf{continuous function} defined over a \textit{nonempty and compact (closed and bounded) set} $\sS\subseteq\real^n$.
Then, there exists a global minimum point  and a global maximum point of $f$ over $\sS$.

\item \label{weier1_continus_lev}  Let $f:\sS\rightarrow (-\infty, \infty]$ be a \textbf{continuous function} defined over a \textit{nonempty and closed set} $\sS\subseteq\real^n$. Suppose that all the level sets $\lev[f, \alpha]=\{\bx\in\sS: f(\bx)\leq \alpha\}$ are bounded. Then, $f$ has a global minimum point over $\sS$.

\item  \label{weier2_continus_coerc} Let $f:\sS\rightarrow  (-\infty, \infty]$ be a \textbf{continuous and coercive function} and let $\sS\subseteq\real^n$ be a \textit{nonempty and closed set}. Then, $f$ has a global minimum point over $\sS$. (\text{The coerciveness ensures the function is bounded over a subset.})

\item \label{weier2_prop_close} Let $f:\sS\rightarrow (-\infty, \infty]$ be a \textbf{proper closed function}, and one of the following is satisfied:
\begin{enumerate}[(a)]
\item \label{weier2_prop_close_v1}  $f$ is defined over a \textit{nonempty and compact set} $\sS$. 
\item \label{weier2_prop_close_v2}  There exists a \textit{nonempty and bounded level set} $\lev[f, \alpha]=\{\bx\in\sS: f(\bx)\leq \alpha\}$. 
\item \label{weier2_prop_close_v3} $f$ is coercive, i.e., $\mathop{\lim}_{\normtwo{\bx}\rightarrow\infty}=\infty$.
\end{enumerate}
Then, there exists a global minimum point of $f$ over $\sS$. And the set of minimizers $\{\bx\in\sS:  f(\bx)\leq  f(\by), \forall \by\in \sS\}$  of $\mathopmin{\bx\in\sS} f(\bx)$ is nonempty and compact.
\end{enumerate}

\end{theorem}
Note the coerciveness in \ref{weier2_continus_coerc} ensures the function is bounded over a subset. The three conditions of \ref{weier2_prop_close}.\ref{weier2_prop_close_v1}$\sim$\ref{weier2_prop_close}.\ref{weier2_prop_close_v3} essentially ensure that the minimum value of $f(\bx)$ cannot be attained at infinity.
The results in \ref{weier2_prop_close} are used more extensively in optimization analysis. To provide a counterexample, let $f=\exp(-x)$ defined on $\real$, which is a proper closed function.
However, the domain or some level sets of $f$ are not bounded, and $f$ is not coercive. Therefore, the attainment of the global optimal point is not guaranteed.

\index{Fermat's theorem}
\textit{Fermat's theorem}, also known as \textit{Fermat's theorem on stationary points}, is a fundamental result in calculus and mathematical optimization. It provides a necessary condition for a function to have a local optimum (either a local maximum or a local minimum) at a point inside the domain of the function.
For a univariate function, Fermat's theorem states the optimality condition for a optimal point that lies in the interior of a set, i.e., a one-dimensional constrained optimization problem.
\begin{theorem}[Fermat's Theorem: Necessary Condition for Univariate Functions]\label{theorem:fermat_theorem}
Let $f: (a,b)\rightarrow \real$ be a one-dimensional differentiable function defined over an interval ($a, b$). 
If a point $x^*\in(a,b)$ (i.e., $x^*\in\interior([a,b])$) is a local maximum or minimum, then $f^\prime(x^*)=0$. 
\end{theorem}
In other words, if a function $f$ has a local maximum or minimum at a point $x^*$ where the function is differentiable, then the slope of the tangent line at that point is zero; the derivative of the function at that point must be zero.
It's important to note that this condition is necessary but \textbf{not} sufficient for $x^*$ to be an optimum. There are cases where the derivative is zero but the point is neither a maximum nor a minimum, such as in the case of an inflection point. Additionally, Fermat's theorem does \textbf{not} apply to boundary points of the domain of $f$ or to points where $f$ is not differentiable.

\subsection*{Descent Directions}
The first-order optimality conditions utilize gradient (first-order) information to evaluate  the optimality of a given point. Here, we first consider the case where the objective function is differentiable and provide the definition of descent directions for constrained and unconstrained optimization problems.
\begin{definition}[Unconstrained and Constrained Descent Directions]\label{definition:uncons_des_direct}
Consider the unconstrained optimization problem (P1) in Definition~\ref{definition:opt_probs_all},
where $ f $ is continuously differentiable over  $\real^n $. Then a vector $ \bd \neq \bzero $ is called an \textit{unconstrained descent direction} at $ \bx$ if $ f'(\bx; \bd)=\nabla f(\bx)^\top \bd < 0 $.

Alternatively, consider the constrained optimization problem (P2) in Definition~\ref{definition:opt_probs_all}, where $f$ is defined on $\sS\subseteq\real^n$. Then a vector $ \bd \neq \bzero $ is called a \textit{constrained descent direction or a feasible descent directions} at $\bx$ if $ f'(\bx; \bd)=\nabla f(\bx)^\top \bd < 0 $, and there exists $ \varepsilon > 0 $ such that $ \bx + \mu\bd \in \sS $ for all $ \mu \in [0, \varepsilon] $.
When the setting is clear from the context, we will simply refer to \textit{unconstrained descent directions} or \textit{constrained descent directions} as \textit{descent directions} for brevity.
\end{definition}

\begin{lemma}[Unconstrained and Constrained Descent Property]\label{lemma:descent_property}
Let $ f $ be a continuously differentiable function over $\sS\subseteq\real^n$ (resp. $\real^n$), and let $\bx \in \sS$ (resp. $\bx\in\real^n$). Suppose that $\bd$ is a descent direction of $ f $ at $\bx$. Then, there exists $\varepsilon > 0$ such that 
$$
f(\bx + \mu \bd) < f(\bx),\quad \text{for any $ \mu \in (0, \varepsilon] $.}
$$
\end{lemma}
\begin{proof}[of Lemma~\ref{lemma:descent_property}]
Since $ f'(\bx; \bd)=\nabla f(\bx)^\top \bd < 0 $, it follows from the definition of the directional derivative (Definition~\ref{definition:partial_deri}) that
$$
\lim_{\mu \to 0^+} \frac{f(\bx + \mu \bd) - f(\bx)}{\mu} = f'(\bx; \bd) < 0.
$$
Therefore, there exists an $\varepsilon > 0$ such that
$
\frac{f(\bx + \mu \bd) - f(\bx)}{\mu} < 0
$
for any $ \mu \in (0, \varepsilon] $, which completes the proof.
\end{proof}

\index{Sufficient decrease condition}
The definition and property of descent directions formalize the concept of moving from a current point in a direction that ensures a reduction in the objective function's value, with specific conditions for both unconstrained and constrained settings.
Building on this concept, we can leverage the properties of descent directions to establish criteria that guarantee a sufficient decrease in the function value when taking a step in such a direction.
\begin{lemma}[Sufficient Decrease Condition from Descent Directions]\label{lemma:valid_suff_des_armi}
Let $ f $ be a continuously differentiable function over $\real^n$, and let $\bx \in \real^n$. Suppose that $\bzero \neq \bd \in \real^n$ is a descent direction of $f$ at $\bx$ and let $\alpha \in (0,1)$. Then there exists $\varepsilon > 0$ such that the inequality
$$
f(\bx) - f(\bx + \mu \bd) \geq -\alpha \mu \nabla f(\bx)^\top \bd
$$
holds for all $\mu \in [0, \varepsilon]$.
\end{lemma}
\begin{proof}[of Lemma~\ref{lemma:valid_suff_des_armi}]
Since $f$ is continuously differentiable, by the linear approximation theorem (Theorem~\ref{theorem:linear_approx})
\begin{align}
&f(\bx + \mu \bd) = f(\bx) + \mu \nabla f(\bx)^\top \bd + o(\mu \normtwo{\bd}) \nonumber\\
&\quad\implies 
f(\bx) - f(\bx + \mu \bd) = -\alpha \mu \nabla f(\bx)^\top \bd - (1 - \alpha) \mu \nabla f(\bx)^\top \bd - o(\mu \normtwo{\bd}). \label{equation:valid_suff_des_armi2}
\end{align}
Since $\bd$ is a descent direction of $f$ at $\bx$, we have
$$
\lim_{\mu \to 0^+} \frac{(1 - \alpha) \mu \nabla f(\bx)^\top \bd + o(\mu \normtwo{\bd})}{\mu} = (1 - \alpha) \nabla f(\bx)^\top \bd < 0.
$$
Therefore, there exists an $\varepsilon > 0$ such that the inequality $
(1 - \alpha) \mu \nabla f(\bx)^\top \bd + o(\mu \normtwo{\bd}) < 0
$ holds for all $\mu \in (0, \varepsilon]$, which combined with \eqref{equation:valid_suff_des_armi2} implies the desired result.
\end{proof}
The sufficient decrease condition guarantees the existence of the Armijo condition, the Goldstein condition, and the Wolfe condition in descent methods with line search (Section~\ref{section:line_search}).

On the other hand, using gradient information, we can construct arbitrary descent directions.
\begin{theorem}[Unconstrained Descent Directions]\label{theorem:uncons_des_dir}
Let $f:\real^n\rightarrow \real$ be a differentiable function.
If $ \nabla f(\bx) \neq \bzero $ and $ \bB $ is a positive definite matrix, then
$$
\bd_1 = -\bB\nabla f(\bx) \qquad \text{and} \qquad \bd_2 = -\bB^{-1}\nabla f(\bx)
$$
are descent directions.
\end{theorem}
\begin{proof}[of Theorem~\ref{theorem:uncons_des_dir}]
A positive definite matrix $ \bB \in \real^{n \times n} $ satisfies
$
\bu^\top \bB \bu > 0, \text{ for all }  \bzero\neq \bu \in \real^n.
$
If we take $ \bu \triangleq \bd_1 $ and exploit the symmetry of $ \bB $, we get
$$
\bd_1^\top \nabla f(\bx) = -\nabla f(\bx)^\top \bB^\top \nabla f(\bx) = -\nabla f(\bx)^\top \bB \nabla f(\bx) < 0.
$$
Similarly, if we let $ \bu \triangleq \bd_2 $, we get
$$
\bd_2^\top \nabla f(\bx) = \bd_2^\top (-\bB \bd_2) = -\bd_2^\top \bB \bd_2 < 0.
$$
Thus, the condition in Definition~\ref{definition:uncons_des_direct} is satisfied in both cases.
\end{proof}
Note that when $\bB\triangleq\bI$, the descent direction becomes the negative gradient direction, which is also known as the \textit{steepest descent direction} and will be explored further in Chapter~\ref{chapter:gradient-descent}.

\paragrapharrow{Optimality conditions and stationary conditions.}
In the remainder of this section, we will introduce  \textbf{optimality conditions} for problems (P1), (P2), and (P3) as defined in Definition~\ref{definition:opt_probs_all}. The optimality conditions for problem (P4) are discussed in Section~\ref{section:gen_kkt_cond}.
In a nutshell, these optimality conditions are outlined in Table~\ref{tab:optimality_conditions}.
\begin{table}[H]
\centering
\caption{Optimality conditions for problems (P1), (P2), and (P3) as defined in Definition~\ref{definition:opt_probs_all}, where ``Nec." is short for ``Necessary" and ``Suff." is short for ``Sufficient".}
\label{tab:optimality_conditions}
\begin{tabular}{|c|c|c|}
\hline
\textbf{Problem} & \textbf{First-Order Condition} & \textbf{Second-Order Condition} \\ \hline
Differentiable (P1) & $\nabla f(\bx^*) = \bzero$ (Nec.) & 
\begin{tabular}{@{}l@{}}
$\nabla^2 f(\bx^*) \succeq \bzero$ (Nec.) \\
$\nabla^2 f(\bx^*) \succ \bzero$ (Suff.)
\end{tabular} \\ \hline
Nonsmooth (P1) & $\bzero \in \partial f(\bx^*)$ (Nec.) & --- \\ \hline
Convex (P1) & 
\begin{tabular}{@{}l@{}}
$\bzero \in \partial f(\bx^*)$ (Nec./Suff.) \\
$\bx^* = \proxf(\bx^*)$  (Nec./Suff.)
\end{tabular} 
& --- \\ \hline
Convex (P2) & $\nabla f(\bx^*)^\top (\bx-\bx^*) \geq 0$ (Nec./Suff.) & ---\\ \hline 
Composite (P3) & $-\nabla f(\bx^*) \in \partial g(\bx^*)$ (Nec.) & --- \\ \hline
\end{tabular}
\end{table}
Additionally, \textbf{stationary conditions} of problems (P1), (P2), and (P3) are defined in Definition~\ref{definition:stat_point}, Definition~\ref{definition:stat_point_uncons_convset}, and Definition~\ref{definition:stat_opt_cond_p3}, respectively:
\begin{subequations}\label{equation:stationar_p123}
\begin{align}
&\textbf{(P1)}:\qquad  \nabla f(\bx^*) &&=\bzero;\\
&\textbf{(P2)}:\qquad  \nabla f(\bx^*)^\top (\bx - \bx^*)& &\geq 0 \text{ for any $\bx\in\sS$};\\
&\textbf{(P3)}:\qquad  -\nabla f(\bx^*) &&\in \partial g(\bx^*).
\end{align}
\end{subequations}

\subsection{Unconstrained Optimization}
Unconstrained differentiable optimization problems are typically expressed as (P1) in Definition~\ref{definition:opt_probs_all}:
\begin{equation}
\text{(P1)}:\qquad \min_{\bx \in \real^n} f(\bx),
\end{equation}
where it is assumed that $ f $ is a continuously differentiable function. Definition~\ref{definition:local_global_opt} introduced the definitions of local and global minimizers. Given a point $\widetildebx$, we want to determine whether this point is a local or global minimizer of the function $ f $. If we proceed from the definition, we need to judge all points in its neighborhood, which is impractical. Therefore, we need a simpler way to verify whether a point is an extremum point. We call these \textit{optimality conditions}, which mainly include first-order optimality conditions and second-order optimality conditions.

\index{Optimality condition}

\subsection*{Local Optimality}
Most objective functions, especially those with several local minimizers,
contain local maximizers and other points which satisfy a necessary condition for a local minimizer. The following theorems help us identify  such points and distinguish  local minimizers from  irrelevant points.

Definition~\ref{definition:uncons_des_direct} and Lemma~\ref{lemma:descent_property} show that as long as there is a descent direction, the function does not reach  a local minimum. Only when $\nabla f(\bx^*)=\bzero$ at $\bx^*$ there are no descent directions for the unconstrained optimization problem (P1). Rigorously, we have the following first-order necessary condition for a local minimum point.
\begin{theoremHigh}[First-Order Necessary Condition for a Local Minimum]\label{theorem:fermat_fist_opt}
Let $f: \real^n \rightarrow \real$ be a  differentiable function. If $\bx^*$ is a local minimizer for $f$, then
$$
\nabla f(\bx^*) = \bzero.
$$
\end{theoremHigh}
\begin{proof}[of Theorem~\ref{theorem:fermat_fist_opt}]
Let $i \in \{1, 2, \ldots, n\}$ and consider the one-dimensional function $g(\mu) = f(\bx^* + \mu \be_i)$. Note that $g$ is differentiable at $\mu = 0$ and that $g'(0) = \frac{\partial f}{\partial x_i}(\bx^*)$. Since $\bx^*$ is a local minimum point of $f$, it follows that $\mu = 0$ is a local minimum of $g$, which immediately implies that $g'(0) = 0$ by Theorem~\ref{theorem:fermat_theorem}. This equality is exactly the same as $\frac{\partial f}{\partial x_i}(\bx^*) = 0$. Since this holds for any $i \in \{1, 2, \ldots, n\}$, we obtain $\nabla f(\bx^*) = \bzero$.
\end{proof}

\begin{figure}[h]
\centering       
\vspace{-0.25cm}                 
\subfigtopskip=2pt               
\subfigbottomskip=-2pt         
\subfigcapskip=-10pt      
\includegraphics[width=0.98\textwidth]{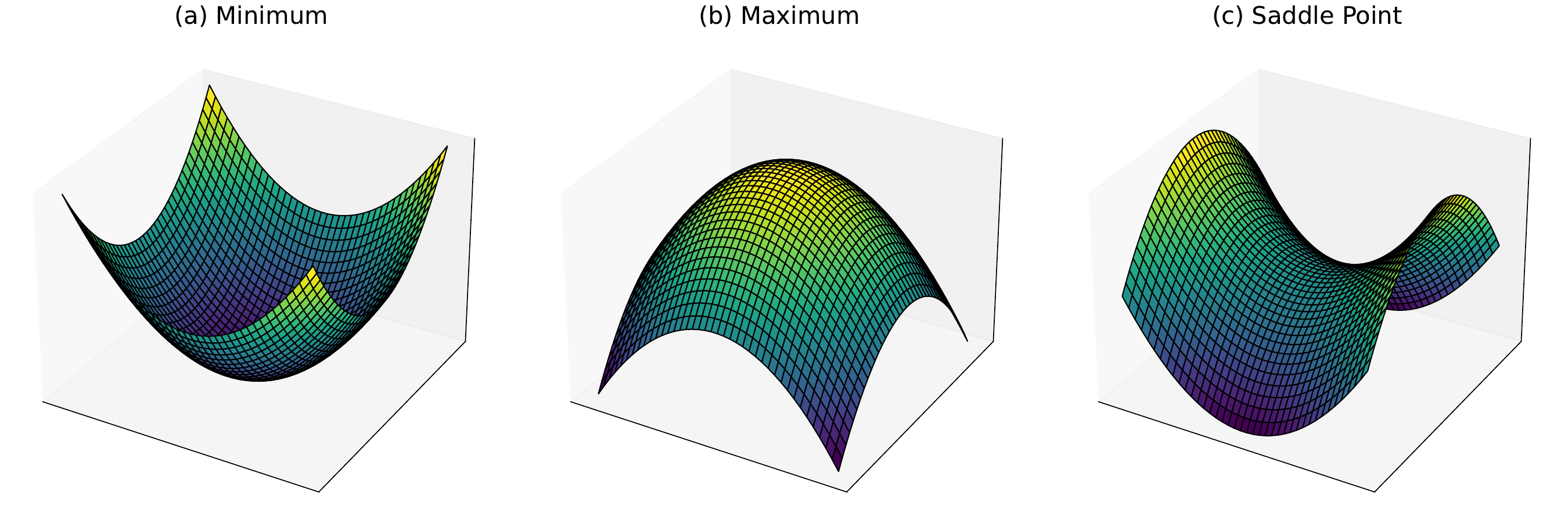}
\caption{Illustration of stationary points, showing minimum, maximum, and saddle points.
}
\label{fig:saddle_point_all}
\end{figure}

\index{Stationary point}
The local minimizers are among the points with $\nabla f(\bx) = \bzero$. They have a special name.
\begin{definition}[Stationary Point, Saddle Point for (P1)]\label{definition:stat_point}
Let $f: \real^n \rightarrow \real$ be a differentiable function. If $\nabla f(\widetildebx) = \bzero$, then $\widetildebx$ is said to be a \textit{stationary point} for $f$.
In the meantime, a stationary point $\widetildebx$ is called a \textit{saddle point} of $f$ if it is neither a local minimum point nor a local maximum point of $f$ (where $\nabla^2 (\widetildebx)$ has both positive and negative eigenvalues); see Figure~\ref{fig:saddle_point_all}.
\end{definition}

Stationary points include local maximizers,  local minimizers, and ``the rest" (which we call \textit{saddle points}). 
Without additional assumptions, even if the first-order necessary conditions are satisfied, we cannot determine whether a stationary point is a local minimum. 
To distinguish between them, we need one extra term in the Taylor expansion. Provided that $f$ is twice continuously differentiable, by the quadratic approximation theorem (Theorem~\ref{theorem:quad_app_theo}),
\begin{equation}\label{equation:necess_cond_loca2}
f(\bx + \bd) = f(\bx) + \bd^\top \nabla f(\bx) + \frac{1}{2} \bd^\top \nabla^2 f(\bx) \bd + \mathcalO(\normtwo{\bd}^3).
\end{equation}
For a stationary point $\widetildebx$, \eqref{equation:necess_cond_loca2} takes the form
$$
f(\widetildebx + \bd) = f(\widetildebx) + \frac{1}{2} \bd^\top \nabla^2 f(\widetildebx) \bd + \mathcalO(\normtwo{\bd}^3).
$$
If the second term is positive for all $\bd$, we say that the matrix $\nabla^2 f(\widetildebx)$ is positive definite (Definition~\ref{definition:psd-pd-defini}). Further, we can take $\normtwo{\bd}$ so small that the remainder term is negligible, and it follows that $\widetildebx$ is a local minimizer. Rigorously, we have the following results.
\begin{theoremHigh}[Second-Order Necessary Condition for a Local Minimum]\label{theorem:second_nec_loca}
Let $f: \real^n \rightarrow \real$ be a twice continuously differentiable function.  If $\bx^*$ is a local minimizer (resp. local maximizer), then $\nabla^2 f(\bx^*)$ is positive semidefinite (resp. negative semidefinite).
\end{theoremHigh}
\begin{proof}[of Theorem~\ref{theorem:second_nec_loca}]
We will prove the result for a  local minimum; the result for a  local maximum follows by considering the function $-f$.
Since $\bx^*$ is a local minimum point, there exists an open ball neighboring $\bx^*$ with a radius $\tau$: $\sB(\bx^*, \tau)$ such that
$
f(\bx) \geq f(\bx^*)
$
for all $\bx \in \sB(\bx^*, \tau)$.
Let $\bd \in \real^n$ be a nonzero vector. For any $0 < \mu < \frac{\tau}{\normtwo{\bd}}$, we have $\bx_\mu = \bx^* + \mu \bd \in \sB(\bx^*, \tau)$, and thus  for any such $\mu$
\begin{equation}\label{equation:second_nec_loca1}
f(\bx_\mu) \geq f(\bx^*).  
\end{equation}
On the other hand, by the linear approximation theorem (Theorem~\ref{theorem:linear_approx}), it follows that there exists a vector $\bu_\mu \in [\bx^*, \bx_\mu]$ such that
\begin{equation}\label{equation:second_nec_loca2}
\begin{aligned}
f(\bx_\mu) - f(\bx^*) 
&= \nabla f(\bx^*)^\top (\bx_\mu - \bx^*) + \frac{\mu^2}{2} \bd^\top \nabla^2 f(\bu_\mu) \bd
= \frac{\mu^2}{2} \bd^\top \nabla^2 f(\bu_\mu) \bd,
\end{aligned}
\end{equation}
where the last equality follows from the first-order necessary condition in Theorem~\ref{theorem:fermat_fist_opt}.

Combining \eqref{equation:second_nec_loca1} and \eqref{equation:second_nec_loca2}, it follows that for any $\mu \in (0, \frac{\tau}{\normtwo{\bd}})$, the inequality $\bd^\top \nabla^2 f(\bu_\mu) \bd \geq 0$ holds. Finally, using the fact that $\bu_\mu \rightarrow \bx^*$ as $\mu \rightarrow 0^+$, and the continuity of the Hessian ($f$ is twice continuously differentiable), we obtain that $\bd^\top \nabla^2 f(\bx^*) \bd \geq 0$. Since this inequality holds for any $\bd \in \real^n$, we conclude the desired result.
\end{proof}

Theorem~\ref{theorem:second_nec_loca} provides a  necessary condition for a local optimum $\bx^*$.
The Taylor expansion \eqref{equation:necess_cond_loca2}  also forms the basis for the proof of the following theorem.

\begin{theoremHigh}[Second-Order Sufficient Condition for a Local Minimum]\label{theorem:second_nec_nonstrict_loca}
Let $f:  \real^n \rightarrow \real$ be a  twice continuously differentiable function, and suppose that $\bx^*$ is a stationary point. 
Further, assume that $\nabla^2 f(\bx)$ is positive semidefinite  in a neighbourhood of $\widetildebx$; that is, there exist a $\tau>0$ such that 
$$
\nabla^2 f(\bx)\succeq \bzero, \text{ for any }\bx\in\sB(\widetildebx, \tau).
$$
Then $\widetildebx$ is a local minimizer.
\end{theoremHigh}
\begin{proof}[of Theorem~\ref{theorem:second_nec_nonstrict_loca}]
By the linear approximation theorem (Theorem~\ref{theorem:linear_approx}), it follows that for any $\bx \in \sB(\widetildebx, \tau)$, there exists a vector $\bxi \in [\bx^*, \bx]$ (hence $\bxi\in\sB(\widetildebx, \tau)$) for which
$$
f(\bx) - f(\bx^*) = \frac{1}{2} (\bx - \bx^*)^\top \nabla^2 f(\bxi) (\bx - \bx^*), \quad\text{for any $\bx \in \sB(\widetildebx, \tau)$}.
$$
Since $\nabla^2 f(\bxi) \succeq \bzero$, we have that $f(\bx) \geq f(\bx^*)$, establishing the result that $\bx^*$ is a local minimum point of $f$ over $\sB(\widetildebx, \tau)$.
\end{proof}

\begin{theoremHigh}[Second-Order Sufficient Condition for a Strict Local Minimum]\label{theorem:second_suff_loca}
Let $f:  \real^n \rightarrow \real$ be a  twice continuously differentiable function, and suppose that $\bx^*$ is a stationary point. 
If $\nabla^2 f(\bx^*) \succ \bzero$ (resp. $\nabla^2 f(\bx^*) \prec \bzero$), then $\bx^*$ is a strict local minimum point (resp. strict local maximum point) of $f$.
\end{theoremHigh}
\begin{proof}[of Theorem~\ref{theorem:second_suff_loca}]
We will prove the result for a strict local minimum; the result for a strict local maximum follows by considering the function $-f$. Suppose  that $\bx^*$ is a stationary point satisfying $\nabla^2 f(\bx^*) \succ \bzero$. Since the Hessian is continuous ($f$ is assumed to be twice continuously differentiable), it follows that there exists an open ball $\sB(\bx^*, \tau)$ for which $\nabla^2 f(\bx) \succ \bzero$ for any $\bx \in \sB(\bx^*, \tau)$. By the linear approximation theorem (Theorem~\ref{theorem:linear_approx}), it follows that for any $\bx \in \sB(\bx^*, \tau)$, there exists a vector $\bxi \in [\bx^*, \bx]$ (and hence $\bxi \in \sB(\bx^*, \tau)$) for which
\begin{equation}\label{equation:second_suff_loc1}
f(\bx) - f(\bx^*) = \frac{1}{2} (\bx - \bx^*)^\top \nabla^2 f(\bxi) (\bx - \bx^*),\quad \text{ for any $\bx \in \sB(\bx^*, \tau)$}.
\end{equation}
Since $\nabla^2 f(\bxi) \succ \bzero$, it follows by \eqref{equation:second_suff_loc1} that for any $\bx \in \sB(\bx^*, \tau)$ with $\bx \neq \bx^*$, the inequality $f(\bx) > f(\bx^*)$ holds, implying that $\bx^*$ is a strict local minimum point of $f$ over $\sS$.
\end{proof}

\begin{exercise}[Sufficient Condition for a Saddle Soint]\label{exercise:second_suff_locasadd}
Let $f: \real^n \rightarrow \real$ be a twice continuously differentiable function, and suppose that $\bx^*$ is a stationary point. If $\nabla^2 f(\bx^*)$ is an indefinite matrix, then $\bx^*$ is a saddle point of $f$ over $\sS$.
\end{exercise}

%

In conclusion, the local maximizers and the \textit{saddle points}, can be characterized by the following remark.
\begin{remark}[Characterization Theorem of Stationary Points]\label{remark:charac_statpoint}
Assume that $\widetildebx$ is a stationary point and that $\nabla^2 f(\widetildebx) \neq \bzero$. Then,
\begin{enumerate}
\item If $\nabla^2 f(\widetildebx)$ is PD: $\widetildebx$ is a strict local minimizer (Theorem~\ref{theorem:second_suff_loca}).
\item If $\nabla^2 f(\widetildebx)$ is PSD: $\widetildebx$ is a local minimizer or a saddle point (Theorem~\ref{theorem:second_nec_loca}).
\item If $\nabla^2 f(\widetildebx)$ is ND: $\widetildebx$ is a strict local maximizer (Theorem~\ref{theorem:second_suff_loca}).
\item If $\nabla^2 f(\widetildebx)$ is NSD: $\widetildebx$ is a local maximizer or a saddle point (Theorem~\ref{theorem:second_nec_loca}).
\item If $\nabla^2 f(\widetildebx)$ is neither definite nor semidefinite: $\widetildebx$ is a saddle point (Exercise~\ref{exercise:second_suff_locasadd}).
\item If $\nabla^2 f(\widetildebx) = \bzero$, then we need higher order terms in the Taylor expansion in order to find the local minimizers among the stationary points.
\end{enumerate}
\end{remark}

\subsection*{Global Optimality}
The conditions outlined in the previous section can, at best, ensure local optimality of stationary points, as they rely solely on local information---specifically, the values of the gradient and the Hessian at a given point (or in the neighborhood of the stationary points as in Theorem~\ref{theorem:second_nec_nonstrict_loca}). To guarantee global optimality, conditions must incorporate global information. For instance, if the Hessian of the function is always positive semidefinite, then all stationary points are also global minima. This property, which we  refer to as convexity (Theorem~\ref{theorem:psd_hess_conv}), ensures that any stationary point found is not just locally but globally optimal.
\begin{theoremHigh}[Global Optimality]\label{theorem:global_optima}
Let $f$ be a twice continuously differentiable function defined over $\real^n$. Suppose that $\nabla^2 f(\bx) \succeq \bzero$ for any $\bx \in \real^n$. Let $\bx^* \in \real^n$ be a stationary point of $f$. Then $\bx^*$ is a global minimum point of $f$.
\end{theoremHigh}
\begin{proof}[of Theorem~\ref{theorem:global_optima}]
By the linear approximation theorem (Theorem~\ref{theorem:linear_approx}), it follows that for any $\bx \in \real^n$, there exists a vector $\bxi \in [\bx^*, \bx]$ for which
$$
f(\bx) - f(\bx^*) = \frac{1}{2} (\bx - \bx^*)^\top \nabla^2 f(\bxi) (\bx - \bx^*).
$$
Since $\nabla^2 f(\bxi) \succeq \bzero$, we have that $f(\bx) \geq f(\bx^*)$, establishing the result that $\bx^*$ is a global minimum point of $f$.
\end{proof}

\subsection{Constrained Optimization}\label{section:constrain_opt}

In this book, we will also consider the constrained optimization problem (P2) given by
$$
\begin{aligned}
\text{(P2)} \qquad & \min f(\bx)  \quad  \text{s.t.}\quad \bx \in \sS\subseteq\real^n,
\end{aligned}
$$
where $f$ is a continuously differentiable function and $\sS$ can be either a convex or non-convex set. 
While stationarity is a necessary condition for an unconstrained local optimal point, the situation becomes more complex when dealing with constrained problems of the form (P2). Instead of examining \textit{stationary points of a function}, we need to consider the notion of \textit{stationary points of a problem}. This subsection discusses the optimality conditions for general constrained optimization problems; for optimality conditions specific to equality and inequality constrained optimization problems, refer to Section~\ref{section:gen_kkt_cond}, which covers KKT conditions.

As previously mentioned, a well-known result for a one-dimensional function $f$ defined and differentiable over an interval $ (a, b) $ is that if a point $ x^* \in (a, b) $ is a local maximum or minimum point, then $ f'(x^*) = 0 $.
This is also known as \textit{Fermat's theorem} (Theorem~\ref{theorem:fermat_theorem}). 
The first-order necessary condition for an unconstrained optimization problem (P1) is also based on this result (Theorem~\ref{theorem:fermat_fist_opt}).
On the other hand, the multidimensional extension of this result states that the gradient is also zero at local optimal points. We refer to such an optimality condition as a \textit{first-order optimality condition for constrained optimization problems}, as it is expressed in terms of  first-order derivatives. In what follows, we will also discuss \textit{second order optimality conditions for constrained optimization problems} that use in addition information on the second order (partial) derivatives.

\index{Fermat's theorem}
\begin{theoremHigh}[First-Order Necessary Condition for an Interior Local Minimum]\label{theorem:fermat_fist_opt_int}
Let $ f : \sS \to \real $ be a function defined on a  set $ \sS \subseteq \real^n $ (not necessarily a convex set). Suppose that $ \bx^* \in \interior(\sS) $~\footnote{The condition can be modified to such that $\bx^*\in\sS$ where $\sS$ is an open set. Note the difference between Theorem~\ref{theorem:fermat_fist_opt_int} and Theorem~\ref{theorem:fermat_fist_opt}.} is a local minimum point and that all the partial derivatives of $ f $ exist at $ \bx^* $. Then $ \nabla f(\bx^*) = \bzero $, i.e., the gradient vanishes at all local minimum points.
\end{theoremHigh}
\begin{proof}[of Theorem~\ref{theorem:fermat_fist_opt_int}]
Let $ i \in \{1, 2, \ldots, n\} $ and consider the one-dimensional function $ g(\mu) = f(\bx^* + \mu \be_i) $. Note that $ g $ is differentiable at $ \mu = 0 $ and that $ g'(0) = \frac{\partial f}{\partial x_i}(\bx^*) $. Since $ \bx^* $ is a local minimum point of $ f $, it follows that $ \mu = 0 $ is a local minimum of $ g $, which immediately implies that $ g'(0) = 0 $ by Fermat's theorem (Theorem~\ref{theorem:fermat_theorem}). The latter equality is exactly the same as $ \frac{\partial f}{\partial x_i}(\bx^*) = 0 $. Since this is true for any $ i \in \{1, 2, \ldots, n\} $, it holds that $ \nabla f(\bx^*) = \bzero $.
\end{proof}
Note again that  this optimality condition is a necessary condition; however, there could be vanished points which are not local maximum or minimum point. For example, the derivative of the function $f(x)=x^5$ is zero at $x=0$, but this point is neither a local minimum point nor a local maximum point.

\begin{theoremHigh}[Second-Order Necessary Condition for an Interior Local Minimum]\label{theorem:second_nec_loca_int}
Let $f: \sS\subseteq\real^n \rightarrow \real$ be a twice continuously differentiable function defined on $\sS$.  If $\bx^*\in\interior(\sS)$~\footnote{The condition can be modified to such that $\bx^*\in\sS$ where $\sS$ is an open set. Note the difference between Theorem~\ref{theorem:second_nec_loca_int} and Theorem~\ref{theorem:second_nec_loca}.} is a local minimizer (resp. local maximizer), then $\nabla^2 f(\bx^*)$ is positive semidefinite (resp. negative semidefinite).
\end{theoremHigh}
\begin{proof}[of Theorem~\ref{theorem:second_nec_loca_int}]
Since $\bx^*\in\interior(\sS)$ is a local minimum point, there exists an open ball neighboring $\bx^*$ with a radius $\tau$: $\sB(\bx^*, \tau) \subseteq \interior(\sS)$ such that
$
f(\bx) \geq f(\bx^*)
$
for all $\bx \in \sB(\bx^*, \tau)$.
The proof follows from the same argument as the proof of Theorem~\ref{theorem:second_nec_loca} by using the first-order necessary condition in Theorem~\ref{theorem:fermat_fist_opt_int} instead of Theorem~\ref{theorem:fermat_fist_opt}.
\end{proof}


%

\begin{theoremHigh}[Second-Order Sufficient Condition for a Strict $\&$ Interior Local Minimum]\label{theorem:second_suff_loca_int}
Let $f: \sS \subseteq \real^n \rightarrow \real$ be a twice continuously differentiable defined on $\sS$, and suppose that $\bx^*\in\interior(\sS)$~\footnote{Again, the condition can be modified to such that $\bx^*\in\sS$ where $\sS$ is an open set. Note the difference between Theorem~\ref{theorem:second_suff_loca_int} and Theorem~\ref{theorem:second_suff_loca}.} is a stationary point. 
If $\nabla^2 f(\bx^*) \succ \bzero$ (resp. $\nabla^2 f(\bx^*) \prec \bzero$), then $\bx^*$ is a strict local minimum point (resp. strict local maximum point) of $f$ over $\sS$.
\end{theoremHigh}
\begin{proof}[of Theorem~\ref{theorem:second_suff_loca_int}]
Suppose  that $\bx^*$ is a stationary point satisfying $\nabla^2 f(\bx^*) \succ \bzero$. Since the Hessian is continuous ($f$ is assumed twice continuously differentiable) and $\bx^*\in\interior(\sS)$, it follows that there exists an open ball $\sB(\bx^*, \tau) \subseteq \sS$ for which $\nabla^2 f(\bx) \succ \bzero$ for any $\bx \in \sB(\bx^*, \tau)$.
The proof then follows from the same argument as the proof of Theorem~\ref{theorem:second_suff_loca}.
\end{proof}

\begin{exercise}[Sufficient Condition for a Saddle Point]\label{exercise:second_suff_locasadd_int}
Let $f: \sS\subseteq \real^n \rightarrow \real$ be a function defined on an open set $\sS \subseteq \real^n$. Suppose that $f$ is twice continuously differentiable over $\sS$ and that $\bx^*$ is a stationary point. If $\nabla^2 f(\bx^*)$ is an indefinite matrix, then $\bx^*$ is a saddle point of $f$ over $\sS$.
\end{exercise}

Note that we will introduce stationary points of constrained optimization problems over a convex set in  Definition~\ref{definition:stat_point_uncons_convset} and Theorem~\ref{theorem:stat_point_uncons_convset}. 
For non-convex constrained sets, the descent direction (Definition~\ref{definition:uncons_des_direct}) can serve as a first-order necessary condition.
\begin{theoremHigh}[First-Order Necessary Condition for (P2)]\label{theorem:_nonconv_fea_loca_optim}
Consider the constrained optimization problem (P2) in Definition~\ref{definition:opt_probs_all}, where $\sS$ can be either convex or non-convex.
If $ \bx^* $ is a local optimal solution of (P2), then there are \textbf{no} constrained/feasible descent directions at $ \bx^* $ (Definition~\ref{definition:uncons_des_direct}).

\end{theoremHigh}
\begin{proof}[of Theorem~\ref{theorem:_nonconv_fea_loca_optim}]
The proof is by contradiction. If there is a constrained/feasible descent direction at $\bx^*$, that is, a vector $ \bd $ and $ \varepsilon_1 > 0 $ such that $ \bx^* + \mu\bd \in \sS $ for all $ \mu \in [0, \varepsilon_1] $ and $ \nabla f(\bx^*)^\top \bd < 0 $ (Definition~\ref{definition:uncons_des_direct}), then by  Lemma~\ref{lemma:descent_property}, there is an $ \varepsilon_2 < \varepsilon_1 $ such that $ f(\bx^* + \mu\bd) < f(\bx^*) $ for all $ \mu \in [0, \varepsilon_2] $, which leads to a contradiction to the local optimality of $ \bx^* $.
\end{proof}

Having established the criteria for a stationary point to be a global minimum in Theorem~\ref{theorem:global_optima}, we now turn our attention to the conditions under which a global minimum is guaranteed to exist.
\begin{theorem}[Attainment Under Coerciveness]\label{theorem:att_coer}
Let $f: \real^n\rightarrow \real$ be a \textcolor{black}{proper closed  (or continuous)} and coercive function, and let $\sS\subseteq\real^n$ be a nonempty closed set. Then $f$ has a global minimum
point over $\sS$.
\end{theorem}
\begin{proof}[of Theorem~\ref{theorem:att_coer}]
Let $\bx_0 \in \sS$ be an arbitrary point in $\sS$. Since the function is coercive (Definition~\ref{definition:coerciveness}), it follows that there exists a $B > 0$ such that
\begin{equation}\label{equation:att_coer}
	f(\bx) > f(\bx_0) \text{ for any } \bx \text{ such that } \normtwo{\bx} > B. 
\end{equation}
Since any global minimizer $\bx^*$ of $f$ over $S$ satisfies $f(\bx^*) \leq f(\bx_0)$, it follows from \eqref{equation:att_coer} that the set of global minimizers of $f$ over $\sS$ is the same as the set of global minimizers of $f$ over $\sS \cap \sB[0, B]$. The set $\sS \cap \sB[0, B]$ is compact and nonempty, and thus by the Weierstrass theorem (\ref{weier2_prop_close} in Theorem~\ref{theorem:weierstrass_them}, which requires the function to be proper closed and coercive; while \ref{weier1_continus} requires the function to be continuous), there exists a global minimizer of $f$ over $\sS \cap \sB[0, B]$ and hence also over $\sS$.
\end{proof}

\subsection{Constrained Optimization over Convex Sets}\label{section:constr_convset}

We restrict the problem to be an optimization problem over a \textbf{convex} set $\sS$. 
In this setting, if the function $f$ is also convex, the problem becomes a \textit{convex optimization problem} (\S~\ref{section:convex_optimiz}).
 Previously, in Definition~\ref{definition:stat_point}, we introduced stationary points for unconstrained problems. In contrast, constrained problems necessitate a different definition.
\begin{definition}[Stationary Points of Constrained Problems on a Convex Set]\label{definition:stat_point_uncons_convset}
	Consider the constrained optimization problem (P2) in Definition~\ref{definition:opt_probs_all}.
	Let $f$ be a continuously differentiable function over a closed \textbf{convex} set $\sS$. Then $\bx^* \in \sS$ is called a stationary point of (P2) if $\nabla f(\bx^*)^\top (\bx - \bx^*) \geq 0$ for any $\bx \in \sS$.
\end{definition}

Note that Definition~\ref{definition:stat_point} defines  stationary points for functions; while Definition~\ref{definition:stat_point_uncons_convset} pertains to the definition of stationary points for optimization problems.


\begin{theoremHigh}[First-Order Necessary Condition for (P2) on a Convex Set]\label{theorem:stat_point_uncons_convset}
Consider the constrained optimization problem (P2) in Definition~\ref{definition:opt_probs_all}.
Let $f$ be a continuously differentiable function over a closed \textbf{convex} set $\sS$, and let $\bx^*$ be a local minimum of (P2). Then $\bx^*$ is a stationary point of (P2).
\end{theoremHigh}
\begin{proof}[of Theorem~\ref{theorem:stat_point_uncons_convset}]
Let $\bx^*$ be a local minimum of (P2), and assume in contradiction that $\bx^*$ is not a stationary point of (P2). Then there exists $\bx \in \sS$ such that $\nabla f(\bx^*)^\top (\bx - \bx^*) < 0$. By Lemma~\ref{lemma:descent_property}, it follows that there exists $\varepsilon \in (0, 1)$ such that $f(\bx^* + \mu \bd) < f(\bx^*)$ for all $\mu \in (0, \varepsilon)$, where $\bd\triangleq \bx-\bx^*$. Since $\sS$ is convex we have that $\bx^* + \mu \bd = (1 - \mu)\bx^* + \mu\bx \in \sS$, leading to the conclusion that $\bx^*$ is not a local minimum point of (P2), which results in a contradiction to the assumption that $\bx^*$ is a local minimum point of (P2).
\end{proof}

When the projection onto the convex set $\sS$ can be expressed in a closed-form, the necessity and sufficiency of the first-order condition has the following result.
\begin{theoremHigh}[Necessity/Sufficiency for (P2) on a Convex Set under Projection]\label{theorem:stat_point_uncons_convset_proj}
Consider the constrained optimization problem (P2) in Definition~\ref{definition:opt_probs_all}.
Let $f$ be a continuously differentiable function over a closed \textbf{convex} set $\sS$, and let $\eta > 0$ be an arbitrary positive scalar. Then $\bx^*$ is a stationary point (P2)
if and only if
$$
\bx^* = \projectS(\bx^* - \eta \nabla f(\bx^*)).
$$
\end{theoremHigh}
\begin{proof}[of Theorem~\ref{theorem:stat_point_uncons_convset_proj}]
By the Projection Property-I (Lemma~\ref{lemma:proj_prop1}), it follows that $\bx^* = \projectS(\bx^* - \eta \nabla f(\bx^*))$ if and only if
$$
\begin{aligned}
&\innerproduct{\bx^* - \eta \nabla f(\bx^*) - \bx^*, \bx - \bx^*} \leq 0 \\
&\quad \implies\quad \nabla f(\bx^*)^\top (\bx - \bx^*) \geq 0, \text{ for any } \bx \in \sS.
\end{aligned}
$$
This completes the proof.
\end{proof}

\subsection{Convex Optimization: Unconstrained}\label{section:convex_opt_uncons}

Subdifferential sets are extremely useful in characterizing minimum points. 
One of the most fundamental optimality conditions states that a point is a global minimum of a proper extended real-valued convex function if and only if $\bzero$ belongs to the subdifferential set at the point. In essence, this is a generalization of Fermat's optimality condition at points of differentiability: $\nabla f(\bx^*) = \bzero$. 
\begin{theoremHigh}[Necessity/Sufficiency  of Unconstrained Convex]\label{theorem:fetmat_opt}
	Let $f: \real^n \rightarrow (-\infty, \infty]$ be a proper convex function. Then,
	$$
	\bx^* \in \mathop{\argmin}_{\bx \in \real^n}f(\bx) 
	$$
	if and only if $\bzero \in \partial f(\bx^*)$.
\end{theoremHigh}
\begin{proof}[of Theorem~\ref{theorem:fetmat_opt}]
	By the definition of the subgradient (Definition~\ref{definition:subgrad}), it holds that $\bx^* \in \mathop{\argmin}_{\bx \in \real^n}f(\bx) $ if and only if 
	$$
	f(\bx) \geq f(\bx^*) + \innerproduct{\bzero, \bx - \bx^*} \quad \text{for any }  \bx \in \domain(f),
	$$
	which is the same as the inclusion $\bzero \in \partial f(\bx^*)$.
\end{proof} 

A direct consequence of Theorem~\ref{theorem:fetmat_opt} and the Proximal Property-I (Lemma~\ref{lemma:prox_prop1}) is that for a proper closed and convex function, $ \bx = \proxf(\bx) $ if and only if $ \bx $ is a minimizer of $ f $.
\begin{theoremHigh}[Necessity/Sufficiency  of Unconstrained Convex Optimization Under Proximal]\label{theorem:opt_cond_prox}
Let $ f:  \real^n \rightarrow (-\infty, \infty]  $ be a proper closed and convex function. Then, $ \bx^* $ is a minimizer of $ f $ if and only if $ \bx^* = \proxf(\bx^*) $.
\end{theoremHigh}
\begin{proof}[of Theorem~\ref{theorem:opt_cond_prox}]
By Theorem~\ref{theorem:fetmat_opt}, $ \bx^* $ is a minimizer of $ f $ if and only if $ \bzero \in \partial f(\bx^*) $, that is, if and only if $\bx^* - \bx^* \in \partial f(\bx^*) $, which, by the Proximal Property-I (Lemma~\ref{lemma:prox_prop1}),  is the same as $ \bx^* = \proxf(\bx^*) $.
\end{proof}

\subsection{Convex Optimization: Constrained}\label{section:convex_optimiz}

A \textit{convex optimization problem} (or simply a convex problem) involves minimizing a convex function over a convex set:
\begin{equation}\label{equation:convex_optim1} 
\textbf{(Convex Optimization)}:\qquad
\begin{aligned}
&\min \quad f(\bx) \quad &&\text{(convex function)}\\
&\text{s.t.}\quad \bx\in\sS \quad &&\text{(convex set)}.
\end{aligned}
\end{equation}

\begin{theorem}[Local is Global in Convex Optimization]\label{theorem:local_glob_conv}
Let $f: \sS \rightarrow \real$ be a convex function (resp. strictly convex function) defined over the convex set $\sS$. 
If $\bx^* \in \sS$ be a local minimum of $f$ over $\sS$, then $\bx^*$ is the global minimum (resp. strict global minimum) of $f$ over $\sS$.
\end{theorem}
\begin{proof}[of Theorem~\ref{theorem:local_glob_conv}]
Since $\bx^*$ is a local minimum of $f$ over $\sS$, there exists a scalar $\tau > 0$ such that $f(\bx) \geq f(\bx^*)$ for any $\bx \in \sS$ satisfying $\bx \in \sB[\bx^*, \tau]$. Now let $\by \in \sS$ satisfy $\by \neq \bx^*$. It suffices to show that $f(\by) \geq f(\bx^*)$. Let $\lambda \in (0, 1]$ be such that $\bx^* + \lambda(\by - \bx^*) \in \sB[\bx^*, \tau]$. An example of such $\lambda$ is $\lambda = \frac{\tau}{\normtwo{\bx^* - \by}}$. Since $\bx^* + \lambda(\by - \bx^*) \in \sB[\bx^*, \tau] \cap \sS$, it follows that
$
f(\bx^*) \leq f(\bx^* + \lambda(\by - \bx^*)),
$
and hence by Jensen's inequality
$
f(\bx^*) \leq f(\bx^* + \lambda(\by - \bx^*)) \leq (1 - \lambda) f(\bx^*) + \lambda f(\by).
$
Therefore, we obtain $f(\bx^*) \leq f(\by)$.

A slight modification of the above argument shows that any local minimum of a strictly convex function over a convex set is indeed  a strict global minimum of the function over the set.
\end{proof}

The optimal set of the convex problem \eqref{equation:convex_optim1}  is the set of all minimizers, that is, $\sX^*=\argmin\{f(\bx) : \bx \in \sS\}$. This definition of an optimal set is also valid for general problems. 
A notable property of convex problems is that their optimal sets are also convex.

\begin{theorem}[Convexity of the Optimal Set in Convex Optimization]\label{theorem:stric_op_str_conv}
Let $f: \sS \rightarrow \real$ be a convex function defined over the convex set $\sS \subseteq \real^n$. Then the set of optimal solutions of the problem,
$
\sX^*=\argmin\{f(\bx) : \bx \in \sS\},
$ 
is convex. If, in addition, $f$ is strictly convex over $\sS$, then there exists \textbf{at most one} optimal solution.
\end{theorem}
\begin{proof}[of Theorem~\ref{theorem:stric_op_str_conv}]
If $\sX^* = \varnothing$, the result follows trivially. 
We then assume that $\sX^* \neq \varnothing$ and denote the optimal value by $f^*$. Let $\bx, \by \in \sX^*$ and $\lambda \in [0, 1]$. Then, by Jensen's inequality $f(\lambda \bx + (1 - \lambda) \by) \leq \lambda f^* + (1 - \lambda) f^* = f^*$, and hence $\lambda \bx + (1 - \lambda) \by$ is also optimal, i.e., belongs to $\sX^*$, establishing the convexity of $\sX^*$. Suppose now that $f$ is strictly convex and $\sX^*$ is nonempty; to show that $\sX^*$ is a singleton, suppose in contradiction that there exist $\bx, \by \in \sX^*$ such that $\bx \neq \by$. Then $\frac{1}{2}\bx + \frac{1}{2}\by \in \sS$, and by the strict convexity of $f$ we have
$$
f\left(\frac{1}{2}\bx + \frac{1}{2}\by\right) < \frac{1}{2}f(\bx) + \frac{1}{2}f(\by) = \frac{1}{2}f^* + \frac{1}{2}f^* = f^*,
$$
which leads to a contradiction to the fact that $f^*$ is the optimal value.
\end{proof}

For strictly quasi-convex functions (Definition~\ref{definition:quasi_convex}), we have the following theorem concerning the existence and uniqueness of solutions:
\begin{theorem}[Unique Miminizer of Closed Strictly Quasi-Convex Functions]\label{theorem:uni_qua_conv}
Let $\sS$ be a nonempty, closed, and convex subset of $\real^n$, and let $f: \sS \rightarrow (-\infty, +\infty]$ be a proper, closed, and strictly quasi-convex function. Then there exists a unique $\bx^*$ such that
$$
f(\bx^*) < f(\bx), \quad \forall \bx \in \sS \setminus \{\bx^*\}.
$$
\end{theorem}
\begin{proof}[of Theorem~\ref{theorem:uni_qua_conv}]
By the Weierstrass theorem (Theorem~\ref{theorem:weierstrass_them}), $f$ has at least one global minimizer $\bx^*$. Suppose there is another global minimizer $\by^*$; then $f(\bx^*) = f(\by^*)$. According to the definition of a strictly quasi-convex function (Definition~\ref{definition:quasi_convex}), for any $\lambda\in(0,1)$, we have
$$
f(\lambda \bx^* + (1 - \lambda)\by^*) < \max\{f(\bx^*), f(\by^*)\} = f(\bx^*),
$$
which contradicts the global optimality of $\bx^*$.
\end{proof}
From the definition of a strictly quasi-convex function, any strictly convex function is also strictly quasi-convex, but a convex function is not necessarily strictly quasi-convex. Using the above conclusion, for any closed strictly convex function defined on a bounded convex set (such as $f(x) = x^2$), its optimal solution is unique. However, for a general convex function, the optimal solution may not be unique. For example, the function $f(x) = \max\{x, 0\}$ has any $x \leq 0$ as an optimal solution.

Stationarity is a \textbf{necessary optimality condition} for local optimality. However, when the objective function is additionally assumed to be convex, stationarity is a \textbf{necessary and sufficient condition} for optimality.

\begin{theoremHigh}[Necessity/Sufficiency  of Constrained Convex]\label{thm:conv_stationary_optimality}
Let $ f:\sS\subseteq\real^n\rightarrow\real $ be a continuously differentiable convex function over a closed and convex set $ \sS$. Then $ \bx^* $ is a stationary point  (Definition~\ref{definition:stat_point_uncons_convset}) of 
$$
\text{(P2+Convex)} \qquad \min_{\bx \in \sS} f(\bx)
$$
if and only if $ \bx^* $ is an optimal solution of (P2+Convex).
\end{theoremHigh}
\begin{proof}[of Theorem~\ref{thm:conv_stationary_optimality}]
If $ \bx^* $ is an optimal solution of (P2+Convex), then by Theorem~\ref{theorem:stat_point_uncons_convset}, it follows that $ \bx^* $ is a stationary point of (P2+Convex). To prove the sufficiency of the stationarity condition, assume that $ \bx^* $ is a stationary point (Definition~\ref{definition:stat_point_uncons_convset}) of (P2+Convex). For any $ \bx \in \sS $, we have:
$$
f(\bx) \geq f(\bx^*) + \nabla f(\bx^*)^\top (\bx - \bx^*) \geq f(\bx^*),
$$
where the first inequality follows from the gradient inequality for convex functions (Theorem~\ref{theorem:conv_gradient_ineq}), and the second inequality follows from the definition of a stationary point (Definition~\ref{definition:stat_point_uncons_convset}). This shows that $ \bx^* $ is indeed the global minimum point of (P2+Convex), completing the proof.
\end{proof}

\begin{theoremHigh}[Sufficiency of Stationarity of Constrained Convex]\label{theorem:suff_sta_conv}
Let $ f:\sS\subseteq\real^n\rightarrow \real $ be a continuously differentiable and convex function defined on a convex set $ \sS $. Suppose that $ \nabla f(\bx^*) = \bzero $ for some $ \bx^* \in \sS $. Then $ \bx^* $ is a global minimizer of $ f $ over $ \sS $.
\end{theoremHigh}
\begin{proof}[of Theorem~\ref{theorem:suff_sta_conv}]
Let $ \bz \in \sS $. Plugging $ \bx = \bx^* $ and $ \by = \bz $ in the gradient inequality (Theorem~\ref{theorem:conv_gradient_ineq}), we obtain that
$$
f(\bz) \geq f(\bx^*) + \nabla f(\bx^*)^\top (\bz - \bx^*),
$$
which by the fact that $ \nabla f(\bx^*) = \bzero $ implies that $ f(\bz) \geq f(\bx^*) $, thus establishing that $ x^* $ is the global minimizer of $ f $ over $ \sS $. 
\end{proof}
Note that Theorem~\ref{theorem:suff_sta_conv} establishes only the sufficiency of the stationarity condition (of the function, not the problem; see Definition~\ref{definition:stat_point} and Definition~\ref{definition:stat_point_uncons_convset}) $\nabla f(\bx^*) = \bzero$ for guaranteeing that $\bx^*$ is a global optimal solution. When $\sS$ is not the entire space, this condition is not necessary . However, when $\sS=\real^n$, then by Theorem~\ref{theorem:fermat_fist_opt} (necessary condition for local optima) and  Theorem~\ref{theorem:local_glob_conv} (local=optimal for convex cases),  this becomes both a necessary and sufficient condition; the same argument in Theorem~\ref{theorem:fetmat_opt}.

\begin{example}[Equality and Inequality Constraints]
Consider the following problem:
\begin{equation}\label{equation:conv_eqineq_const}
\begin{aligned}
\min & \quad f(\bx) \\
\text{s.t.} & \quad g_i(\bx) \leq 0, \quad i = \{1, 2, \ldots, m\},\\
&\quad h_j(\bx) = 0, \quad j = \{1,2, \ldots, p\},
\end{aligned}
\end{equation}
where $f, g_1, g_2, \ldots,g_m:\real^n\rightarrow \real$ are convex functions, and $h_1, h_2, \ldots,h_p$ are affine functions.
Since the objective function is convex and the feasible set is a convex set, this problem is a convex optimization problem.
The feasible set can be written as
$$
\sS = \bigg( \bigcap_{i=1}^m \lev[g_i, 0] \bigg) \cap \bigg( \bigcap_{j=1}^p \{ \bx \mid  h_j(\bx) = 0 \} \bigg),
$$
which is a convex set since $\sS$ is  an intersection of level sets of convex functions and hyperplanes, both of which are also convex sets.
\end{example}

\begin{example}[Linear Programming (LP)]\label{example:linear_program}
A \textit{linear programming (LP) problem}  involves minimizing a linear objective function subject to linear equalities and inequalities:
$$
\text{(LP)}\qquad 
\begin{aligned}
\min & \quad \bc^\top \bx +d \\
\text{s.t.}
&\quad \bG\bx \leq \bh, \\
&\quad \bA\bx = \bb,
\end{aligned}
$$
where $ \bG \in \real^{m \times n} $, $ \bh \in \real^m $, $ \bA \in \real^{p \times n} $, $ \bb \in \real^p $, and $ \bc \in \real^n $. 
This constitutes a convex optimization problem since affine functions are inherently convex. The constant $d$ in the objective function can be omitted without affecting the optimal set.
Consider the following LP problem:
$$
\text{(Standard LP)}\qquad 
\begin{aligned}
\min & \quad \bc^\top \bx \\
\text{s.t.} & \quad \bA\bx = \bb, \\
& \quad \bx \geq \bzero.
\end{aligned}
$$
This formulation is often referred to as the ``standard form" of LP in literature.
To convert an LP to its standard form, we can  introduce nonnegative slack variables $\bs$  such that $\bG\bx +\bs=  \bh$ and   express $\bx$ as the difference of two nonnegative variables $\bx\triangleq\bx^+-\bx^-$, where $\bx^+, \bx^-\geq \bzero$:
$$
\text{(LP$'$)}\qquad 
\begin{aligned}
\min & \quad \bc^\top \bx^+ - \bc^\top\bx^- +d \\
\text{s.t.}
&\quad \bG\bx^+ - \bG\bx^-  +\bs=  \bh,   \\
&\quad \bA\bx^+ - \bA\bx^-  = \bb, \\
&\quad \bs \geq \bzero, \bx^+ \geq \bzero, \bx^-\geq \bzero, 
\end{aligned}
\quad\iff \quad
\begin{aligned}
\min & \quad [\bc^\top , - \bc^\top, \bzero] \widetildebx +d \\
\text{s.t.}
&\quad  [\bG, -\bG, \bI] \widetildebx =\bh,   \\
&\quad [\bA, - \bA, \bzero ] \widetildebx  = \bb, \\
&\quad \widetildebx \geq \bzero, 
\end{aligned}
$$
where 
$
\widetildebx\triangleq
\scriptsize
\begin{bmatrix}
\bx^+ \\  \bx^- \\\bs 
\end{bmatrix}
$. This is in a standard form with variables $\widetildebx$. For each feasible solution $\bx^*$ of (LP), we can set $\bs\triangleq\bh-\bG\bx^*$, $x_i^+\triangleq \max\{0, x_i\}$, and $x_i^-\triangleq \max\{0, -x_i\}$ for all $i\in\{1,2,\ldots, n\}$, which is a feasible solution of (LP$'$). This shows the optimal value of (LP$'$) is less than or equal to the optimal value of (LP). 
Conversely, suppose $\bx^+, \bx^-$, and $\bs$ are feasible solutions of (LP$'$). Then, $\bx\triangleq\bx^+-\bx^-$ is a feasible solution of (LP), demonstrating that the optimal value of (LP) is less than or equal to the optimal value of (LP$'$).  Therefore, the optimization problems (LP) and (LP$'$) are equivalent.
\end{example}

\begin{example}[Quadratic Programming]
The convex optimization problem \eqref{equation:conv_eqineq_const} is called a \textit{quadratic program (QP)} if the objective function is quadratic with positive semidefinite (convex), and the constraint functions  are affine. A general form can be written as:
$$
\begin{aligned}
\text{min}& \quad \bx^\top \bP \bx + 2 \bq^\top \bx + r \\
\text{s.t.} & \quad \bG\bx \leq \bh, \\
&\quad \bA\bx=\bb,
\end{aligned}
$$
where $\bP \in \real^{n \times n}$ is positive semidefinite, $\bq \in \real^n$, $\bG\in\real^{m\times n}$, $\bh\in\real^m$, $\bA \in \real^{p \times n}$, and $\bb \in \real^p$. 
If both the objective function and inequality constraints are convex quadratic:
$$
\begin{aligned}
\text{min}& \quad \bx^\top \bP \bx + 2 \bq^\top \bx + r \\
\text{s.t.} & \quad \bx^\top\bG_i\bx + 2\bh_i^\top\bx+ c_i \leq 0,\quad i=\{1,2,\ldots,m\},\\
&\quad \bA\bx=\bb,
\end{aligned}
$$
where $\bP, \bG_1, \bG_2, \ldots,\bG_m$ are positive semidefinite, the problem is called a \textit{quadratically constrained quadratic program (QCQP)}.
In a QCQP, we minimize a convex quadratic function over a feasible region that is the intersection of ellipsoids when $\bG_i\succ \bzero$.
\end{example}

\subsection{Non-Convex Composite Optimization}
We finally consider optimality conditions for a composite problem (P3) as defined in Definition~\ref{definition:opt_probs_all}:
$$
\textbf{(P3)}:\qquad \text{Find}\quad \bx^* = \mathop{\argmin}_{\bx} \left\{F(\bx)\triangleq f(\bx)+g(\bx) \right\}.
$$
\begin{theoremHigh}[Optimality Conditions for the Composite Problem (P3)]\label{theorem:opt_cond_p3}
Let $f: \real^n \rightarrow (-\infty, \infty]$ be a \textcolor{black}{proper} function, and let $g: \real^n \rightarrow (-\infty, \infty]$ be a \textcolor{black}{proper} \textbf{convex} function such that $\dom(g) \subseteq \interior(\dom(f))$. Then, problem (P3) has the following properties: 
\begin{enumerate}[(i)]
\item  \textit{Necessary condition.} If $\bx^* \in \dom(g)$ is a local optimal solution of (P3) and $f$ is differentiable at $\bx^*$, then
\begin{equation}\label{equation:opt_cond_p3_e1}
-\nabla f(\bx^*) \in \partial g(\bx^*).
\end{equation}

\item \textit{Necessary and sufficient condition for convex problems.} Suppose additionally that $f$ is convex. If $f$ is differentiable at $\bx^* \in \dom(g)$, then $\bx^*$ is a global optimal solution of (P3) if and only if \eqref{equation:opt_cond_p3_e1} is satisfied.
\end{enumerate}
\end{theoremHigh}
\begin{proof}[of Theorem~\ref{theorem:opt_cond_p3}]
\textbf{(i).} Let $\by \in \dom(g)$. By the convexity of $\dom(g)$, for any $\lambda \in (0, 1)$, the point $\bx_\lambda = (1 - \lambda)\bx^* + \lambda \by$ is in $\dom(g)$. 
By the local optimality of $\bx^*$, it follows that, for sufficiently small $\lambda$,
$$
\begin{aligned}
&f(\bx_\lambda) + g(\bx_\lambda) \geq f(\bx^*) + g(\bx^*)\\
&\implies\quad f\big((1 - \lambda)\bx^* + \lambda \by\big) + g\big((1 - \lambda)\bx^* + \lambda \by\big) \geq f(\bx^*) + g(\bx^*).
\end{aligned}
$$
Using the convexity of $g$, it follows that
$$
\begin{aligned}
&f\big((1 - \lambda)\bx^* + \lambda \by\big) + (1 - \lambda)g(\bx^*) + \lambda g(\by) \geq f(\bx^*) + g(\bx^*)\\
&\implies\quad \frac{f\big((1 - \lambda)\bx^* + \lambda \by\big) - f(\bx^*)}{\lambda} \geq g(\bx^*) - g(\by).
\end{aligned}
$$
Taking $\lambda \rightarrow 0^+$ in the above inequality and by \eqref{equation:direc_contdiff}: $f^\prime(\bx^*; \by - \bx^*) = \innerproduct{\nabla f(\bx^*), \by - \bx^*}$), we have
$$
f^\prime(\bx^*; \by - \bx^*) =\innerproduct{\nabla f(\bx^*), \by - \bx^*}\geq g(\bx^*) - g(\by),
$$
where we used the fact that since $f$ is differentiable at $\bx^*$, its directional derivatives exist. 
Therefore, for any $\by \in \dom(g)$,
$$
g(\by) \geq g(\bx^*) + \innerproduct{-\nabla f(\bx^*), \by - \bx^*},
$$
showing that indeed $-\nabla f(\bx^*) \in \partial g(\bx^*)$.

\paragraph{(ii).} Suppose further that $f$ is convex. If $\bx^*$ is an optimal solution of (P3), then we have already shown in part (i) that \eqref{equation:opt_cond_p3_e1} is satisfied. Now assume that \eqref{equation:opt_cond_p3_e1} is satisfied. Then, for any $\by \in \dom(g)$,
\begin{equation}\label{equation:opt_cond_p3_ii1}
g(\by) \geq g(\bx^*) + \innerproduct{-\nabla f(\bx^*), \by - \bx^*}.
\end{equation}
By the convexity of $f$, for any $\by \in \dom(g)$,
\begin{equation}\label{equation:opt_cond_p3_ii2}
f(\by) \geq f(\bx^*) + \innerproduct{\nabla f(\bx^*), \by - \bx^*}.
\end{equation}
Adding \eqref{equation:opt_cond_p3_ii1} and \eqref{equation:opt_cond_p3_ii2}, we find
$$
f(\by) + g(\by) \geq f(\bx^*) + g(\bx^*), \text{ for any $\by \in \dom(g)$}.
$$
This proves that $\bx^*$ is an optimal solution of (P3).
\end{proof}

Note that Definition~\ref{definition:stat_point} defines  stationary points for functions; while Definition~\ref{definition:stat_point_uncons_convset} pertains to the definition of stationary points for optimization problems.
The condition \eqref{equation:opt_cond_p3_e1} is an important optimality condition for the composite problem (P3), and we will refer
to it as the ``stationarity" condition of (P3). Their comparison is shown in \eqref{equation:stationar_p123}.
\begin{definition}[Stationarity of (P3)]\label{definition:stat_opt_cond_p3}
Consider the problem (P3) in Definition~\ref{definition:opt_probs_all}. Let $f: \real^n \rightarrow (-\infty, \infty]$ be \textcolor{black}{proper}, and let $g: \real^n \rightarrow (-\infty, \infty]$ be a \textcolor{black}{proper} \textbf{convex} function such that $\dom(g) \subseteq \interior(\dom(f))$. 
A point $\bx^*$ in which $f$ is differentiable is called a \textit{stationary point} of (P3) if
\begin{equation}
-\nabla f(\bx^*) \in \partial g(\bx^*).
\end{equation}
\end{definition}

\section{Constrained Optimization: Equality/Inequality Constraints}\label{section:gen_kkt_cond}
We have already observed that finding global solutions can be challenging in the absence of constraints (Remark~\ref{remark:charac_statpoint}). However, introducing constraints can improve the situation. The feasible set may exclude many local minima, making it easier to identify the global minimum from the remaining candidates.

In Section~\ref{section:constrain_opt}, we introduced optimality conditions for constrained optimization problems that do not rely on the algebraic specification of the set $\sS$, but rather on its geometric properties. This section addresses the optimality conditions that depend on the algebraic specifications, specifically equality and inequality constraints, as seen in problem (P4) of Definition~\ref{definition:opt_probs_all}. For simplicity, we use different letters for the equality and inequality constraints:
\begin{equation}\label{equation:P4_in_kkt}
(\text{P4}')\qquad 
\begin{aligned}
\min & \quad f(\bx) \\
\text{s.t.} & \quad g_i(\bx) \leq 0, \quad i = \{1, 2, \ldots, m\} \triangleq\mathcalI,\\
&\quad h_j(\bx) = 0, \quad j = \{1,2, \ldots, p\}\triangleq\mathcalE.
\end{aligned}
\end{equation}
That is, the \textit{feasible set} is $\sS = \{\bx \mid g_i(\bx) \leq 0, i\in\mathcalI;    h_j(\bx) = 0, j\in \mathcalE\}$.
For the sake of convenience,  any constraint (either equality or inequality) can be denoted by 
$$
c_i, \quad i \in \mathcalE \cup \mathcalI,
$$
where $c_i(\bx) = g_i(\bx)$ if $i\in\mathcalI$ and  $c_i(\bx) = h_i(\bx)$ if $i\in\mathcalE$.

\subsection{Mathematical Tools}
Before delving into the optimality conditions of $(\text{P4}')$, we provide some mathematical background.
Let $\bg(\bx) \triangleq \left\{g_{i}(\bx)\right\}_{i\in\{1,2,\ldots,m\}}:\real^n\rightarrow \real^{m}$ and $\bh(\bx) \triangleq \left\{h_{i}(\bx)\right\}_{i\in\{1,2,\ldots,p\}}:\real^n\rightarrow \real^{p}$. Then
\begin{equation}\label{equation:jacob_inkkt}
\nabla \bg(\bx) = 
\begin{bmatrix}
\nabla g_1(\bx)^\top \\
\nabla g_2(\bx)^\top \\
\vdots \\
\nabla g_m(\bx)^\top \\
\end{bmatrix}
\in\real^{m\times n}
\qquad \text{and}\qquad 
\nabla \bh(\bx) = 
\begin{bmatrix}
\nabla h_1(\bx)^\top \\
\nabla h_2(\bx)^\top \\
\vdots \\
\nabla h_p(\bx)^\top \\
\end{bmatrix}
\in\real^{p\times n}
\end{equation}
are Jacobian matrices, where each row corresponds to the transpose of the gradient of the respective function.
\begin{definition}[Sets of (P4$'$)]\label{definition:setsp4prime}
Given a  point $\bx$ in $(\text{P4}')$, we define the following sets: 
\begin{itemize}
	\item $\sD(\bx) \triangleq \{\bd \mid \nabla f(\bx)^\top \bd < 0\}$, the cone (Definition~\ref{definition:convex_cone}) of ``improving" directions of $f(\bx)$  at   $\bx$  (not including the origin $\bzero$ in the cone), i.e., the set of descent directions (Definition~\ref{definition:uncons_des_direct}) at $\bx$.
	\item $\sI(\bx)\triangleq \{i \mid g_i(\bx) = 0\}$, the set of indices of the binding inequality constraints at $\bx$, i.e., the set of active constraints.
	\item $\sG(\bx)\triangleq \{\bd \mid \nabla g_i(\bx)^\top \bd < 0 \text{ for all } i \in \sI(\bx) \}$, the cone of ``inward" pointing directions for the binding constraints at $\bx$, (not including the origin $\bzero$ in the cone).
	\item $\sH(\bx)\triangleq \{\bd \mid \nabla h_i(\bx)^\top \bd = 0 \text{ for all } i = 1,2, \ldots, p\}$, the set of tangent directions for the equality constraints at $\bx$.
	\item $\sA(\bx)\triangleq \sI(\bx) \cup \mathcalE$, the set of union of active constraints and equality constraints, i.e., the \textit{general active set}.
\end{itemize}
\end{definition}

From an introductory course in optimization, there is an important relationship between $ \bg^* $, the gradient of the cost function, and $ \ba_i^* $ ($i\in \mathcalI\cup \mathcalE$), the gradients of the constraint functions, all evaluated at a local minimizer. This leads to the introduction of the \textit{Lagrangian function}:
\begin{definition}[Lagrange's Function]\label{definition:lag_func}
Given the differentiable objective function $ f $ and the constraints $ c_i, i=1, 2, \ldots, q $, where $q \triangleq m+p$ is the total number of constraints for problem (P4$'$) in \eqref{equation:P4_in_kkt}. \textit{Lagrange's function (or Lagrangian function)} is defined as:
$$
L(\bx, \blambda) = f(\bx) - \sum_{i=1}^{q} \lambda_i c_i(\bx).
$$
The scalars $ \left\{ \lambda_i \right\} $ are called \textit{Lagrangian multipliers}.
The gradient of $ L $ with respect to $ \bx $ is denoted by $ \nabla_{\bx} L \triangleq \nabla L $, and we see that
$$
\nabla_{\bx} L(\bx, \blambda) = \nabla f(\bx) - \sum_{i=1}^{q} \lambda_i \nabla c_i(\bx).
$$
\end{definition}

When characterizing the directions in which we can step away from $\bx$ while remaining feasible, 
a related concept is called \textit{tangent}, which is a limiting direction of a feasible sequence.

\begin{definition}[Tangent Vector, Tangent Cone]
The vector $\bd$ is said to be a \textit{tangent (or tangent vector)} to the set $\sS$ at a point $\bx\in\sS$ if there are a feasible sequence $\{\bu^\toptzero\}_{t>0}$ approaching $\bx$ and a sequence of positive scalars $\{c_t\}_{t>0}$ with $c_t \to 0$ such that
\begin{equation}
	\lim_{t \to \infty} \frac{\bu^\toptzero - \bx}{c_t} = \bd.
\end{equation}
The set of all tangents to $\sS$ at $\bx$ is called the \textit{tangent cone} and is denoted by $\mathcalT_\sS(\bx)$.
\end{definition}

It is easy to see that the tangent cone is indeed a cone (Definition~\ref{definition:convex_cone}). If $\bd$ is a tangent vector with corresponding sequences $\{\bu^\toptzero\}$ and $\{c_t\}$, then by replacing each $c_t$ by $\alpha^{-1} c_t$, for any $\alpha > 0$, we find that $\alpha \bd \in \mathcalT_\sS(\bx^*)$. We obtain that $\bzero \in \mathcalT_\sS(\bx)$ by setting $\bu^\toptzero \triangleq \bx$ in the definition of feasible sequence.

We turn now to the \textit{linearized feasible direction set}, which we define as follows.

\begin{definition}[(Linearized, Sequential) Feasible Directions]\label{definition:lfd_sfd}
Let $\bx\in\sS$ be a feasible point. If there exists $\varepsilon>0$ such that $\bx+\mu\bd\in\sS$ for all $\mu\in[0, \varepsilon]$, then $\bd$ is called a \textit{feasible direction} of $\sS$ at $\bx$. The set of all feasible directions of $\sS$ at $\bx$ is denoted as 
$$
\fd(\bx) \triangleq \left\{\bd\neq\bzero \mid \bx + \mu\bd \in \sS, \forall \mu\in[0, \varepsilon]\right\}.
$$
Given further the general active set $\sA(\bx)\triangleq\sI(\bx) \cup \mathcalE$ of Definition~\ref{definition:setsp4prime}, the set of \textit{linearized feasible directions}, denoted as $\lfd(\bx)$, is
$$
\lfd(\bx) \triangleq \left\{ \bd \; \middle| \;
\begin{array}{ll}
	\bd^\top \nabla c_i(\bx) = 0, & \text{for all } i \in \mathcalE, \\
	\bd^\top \nabla c_i(\bx) \leq  0, & \text{for all } i \in \sI(\bx)
\end{array}
\right\}.
$$
As with the tangent cone, it is easy to verify that $\lfd(\bx)$ is a cone.
If there exist sequences $\{\bd^\toptzero\}_{t>0}$ and $\{\mu_t>0\}_{t>0}$ such that $\bx^\toptzero\triangleq\bx+\mu_t\bd^\toptzero \in\sS$ for all $t>0$ and $\bd^\toptzero\rightarrow\bd, \mu_t\rightarrow 0$, then the limiting direction $\bd$ is called the \textit{sequential feasible direction} of $\sS$ at $\bx$. The set of all sequential feasible directions of $\sS$ at $\bx$ is 
$$
\sfd(\bx) \triangleq \left\{ \bd \; \middle| \;
\begin{array}{ll}
\bx+\mu_t\bd^\toptzero\in\sS, & \forall t, \\
\bd^\toptzero\rightarrow \bd, \mu_t\rightarrow 0 &
\end{array}
\right\}.
$$
Note that 
\begin{itemize}
\item $\bx^\toptzero\neq \bx$ for all $t>0$.
\item $\lim_{t\to \infty}\bx^\toptzero = \bx$.
\item $\bx^\toptzero\in\sS$ for all $t$ sufficiently large.
\item If we set $\mu_t\triangleq\normtwo{\bx^\toptzero - \bx}$, then $\bd^\toptzero = \frac{\bx^\toptzero - \bx}{\normtwo{\bx^\toptzero - \bx}} \rightarrow \bd$. Therefore, 
$$
\mathcalT_\sS(\bx) = \sfd(\bx) \cup\{\bzero\}.
$$
\end{itemize}
\end{definition}
Note the difference between $\lfd(\bx)$ and $\sD(\bx)\triangleq \{\bd \mid \nabla f(\bx)^\top \bd < 0\}$ we have defined previously.
Intuitively speaking, vectors in the linearized feasible direction set $\lfd(\bx)$ should ensure that 
\begin{itemize}
\item The gradient of equality constraint functions remains perpendicular to these directions, thus preserving the value of $c_i(\bx)$, $i \in \mathcal{E}$ (or $h_j(\bx),  j = \{1,2, \ldots, p\}$).
\item  For the active constraints $\sI(\bx)=\sA(\bx) \cap \mathcal{I}$, the value of inequality constraint functions $g_i(\bx)$ should not increase along these directions. Therefore, the linearized feasible direction for $c_i(\bx)$, $i \in \sA(\bx) \cap \mathcal{I}$ (or $g_i(\bx),  i = \{1,2, \ldots, m\}$), can be a descent direction.
\end{itemize}
However, for constraints in the non-active set, there is no requirement for any specific direction in the linearized feasible direction set.

It is important to note that the definition of tangent cone does not rely on the algebraic specification of the set $\sS$, only on its geometry. The linearized feasible direction set does, however, depend on the definition of the constraint functions $c_i$, $i \in \mathcalE \cup \mathcalI$.

Using the definition of general active set $ \sA(\bx) $, we define the notion of linear independence constraint qualification (LICQ), the importance of which will be clear in the sequel.
\index{Linear independence constraint qualification}
\index{LICQ}
\begin{definition}[Linear Independence Constraint Qualification (LICQ)]\label{definition:licq}
Consider the problem in \eqref{equation:P4_in_kkt}.
Given a feasible point $ \bx $ and the corresponding general active set $ \sA(\bx) $. If the gradients of the constraint functions corresponding to the active set, i.e., $ \nabla c_i(\bx), i \in \sA(\bx) $ ($c_i$ can be either the equality constraint or the inequality constraint function), are linearly independent, then the \textit{linear independence constraint qualification (LICQ)} holds at the point $ \bx $.
\end{definition}

We end up this subsection by providing Farkas' lemma and its variants, which is a necessary tool for proving the KKT conditions of $(\text{P4}')$ in the sequel.
\begin{lemma}[Farkas' Lemma and Variants]\label{lemma:farka_lemm}
Let $\bc \in \real^n$ and $\bA \in \real^{m \times n}$. Then \textbf{exactly one} of the following systems has a solution:
\begin{enumerate}[(i)]
\item $\bA\bx \leq \bzero, \bc^\top \bx > 0.$
\item $\bA^\top \by = \bc, \by \geq \bzero.$  (i.e., $\bc$ is a conic combination of rows of $\bA$.)
\end{enumerate}
\noindent Equivalently, this implies  the following two claims are equivalent:
\begin{enumerate}[(a)]
\item The implication $\bA\bx \leq \bzero \implies \bc^\top \bx \leq 0$ holds true.
\item There exists $\by\in \real_{+}^{m}$ such that $\bA^\top \by = \bc$.
\end{enumerate}

The above Farka's lemma also indicates the following result.
Given matrices $ \bD, \bE, $ and $ \bF $ of appropriate dimensions, \textbf{exactly one} of the following two systems has a solution:
\begin{enumerate}[(1)]
\item $ \bD\bx < \bzero, \bE \bx \leq \bzero, \bF \bx = \bzero, $
\item $ \bD^\top \bu + \bE^\top \bw + \bF^\top \bv = \bzero, $ $ \bu \geq \bzero, \bw \geq \bzero, \bone^\top \bu = 1. $
\end{enumerate}
\end{lemma}

\begin{proof}[of Lemma~\ref{lemma:farka_lemm}]
\textbf{(a)$\iff$(b).}
Since the system (i, ii) is equivalent to the system (a, b), we only prove the latter here.
Suppose that system (b) is feasible, meaning that there exists $\by\in \real_{+}^{m}$ such that $\bA^\top \by = \bc$. To see that the implication (a) holds, suppose that $\bA\bx \leq \bzero$ for some $\bx\in \real^{n}$. Since $\by$ is nonnegative, premultiplying the inequality $\bA\bx \leq \bzero$  by $\by^\top$ yields
$
\bc^\top\bx = \by^\top \bA\bx \leq \bzero,
$
where $\bc^\top = \by^\top \bA$ by system (b).

Conversely, assume that the implication (a) holds.  Suppose in contradiction that system (b) is infeasible, and consider the following closed and convex set:
$$
\sS = \left\{ \bu \in \real^n \mid  \bu = \bA^\top \by \text{ for some } \by \in \real_{+}^{m} \right\}.
$$
The closedness of $\sS$ follows from Exercise~\ref{exercise:closed_fini_cone}. The infeasibility of (b) means that $\bc \notin \sS$. By Theorem~\ref{theorem:stric_sep_theo}, it follows that there exists a vector $\bw \in \real^{m} \setminus \{\bzero\}$ and $\alpha \in \real$ such that 
\begin{equation}\label{equation:fakar_eq1}
\bw^\top \bc > \alpha
\qquad \text{and}\qquad 
\bw^\top \bu \leq \alpha \text{ for all } \bu \in \sS. 
\end{equation}
Since $\bzero \in \sS$, we can conclude that $\alpha \geq 0$, and hence also that $\bw^\top \bc > 0$. In addition, the second equality in \eqref{equation:fakar_eq1} is equivalent to
\begin{equation}\label{equation:fakar_eq2}
(\bA\bw)^\top \by=\bw^\top \bA^\top\by \leq \alpha, \text{ for all } \by \in \real_{+}^{m}.
\end{equation}
This implies that $\bA\bw \leq 0$. Indeed, if an index $i \in \{1,2,\ldots, m\}$ existed such that $(\bA\bw)_{i} > 0$, then choosing $\by \triangleq \gamma \be_{i}$ would result in $(\bA\bw)^\top \by = \gamma (\bA\bw)_{i}$, which is an expression that approaches $\infty$ as $\gamma \rightarrow \infty$. Taking a large enough $\gamma$ will contradict \eqref{equation:fakar_eq2}, thereby contradicting our assumption that implication (a) holds, leading us to conclude that system (b) must be feasible.

\paragraph{(1)$\iff$(2).}
Suppose (1) does not have a solution. Then the following system  has no solution:
$$
\begin{cases}
\bD\bx + \bone\gamma \leq \bzero, & \gamma > \bzero; \\
\bE\bx \leq \bzero; \\
\bF\bx \leq \bzero; \\
-\bF\bx \leq \bzero
\end{cases}
\quad\iff\quad
\begin{bmatrix}
\bD & \bone \\
\bE & \bzero \\
\bF & \bzero \\
-\bF & \bzero
\end{bmatrix}
\begin{bmatrix}
\bx \\
\gamma
\end{bmatrix}
\leq \bzero, \quad [0, \ldots, 0, 1] \cdot
\begin{bmatrix}
\bx \\
\gamma
\end{bmatrix}
> 0.
$$
From Farkas' lemma, there exists a vector $ [\bu; \bw; \bv^1; \bv^2] \geq 0 $ such that:
$$
\begin{bmatrix}
\bD & \bone \\
\bE & 0 \\
\bF & 0 \\
-\bF & 0
\end{bmatrix}^\top
\begin{bmatrix}
\bu \\
\bw \\
\bv^1 \\
\bv^2
\end{bmatrix}
=
\begin{bmatrix}
0 \\
\vdots \\
0 \\
1
\end{bmatrix}.
$$
This can be rewritten as:
$$
\bD^\top \bu + \bE^\top \bw + \bF^\top (\bv^1 - \bv^2) = 0 \,, \quad \bone^\top \bu = 1 \,.
$$
Setting $\bv \triangleq\bv^1 -\bv^2$ completes the proof.
\end{proof}

\subsection{KKT Conditions for Linearly Constrained Problems}
To begin, we suppose the equality and inequality constraints of (\text{P4}$'$) in \eqref{equation:P4_in_kkt} are linear functions.
Below, we outline the\textit{Karush-Kuhn-Tucker (KKT) conditions} for problems with linear constraints.

\begin{theoremHigh}[KKT Conditions for Linearly Constrained Problems, Theorem 10.7 in \citet{beck2014introduction}]\label{theorem:kktcond}
Consider the minimization problem
$$
(\text{P})\qquad
\begin{aligned}
\min & \quad f(\bx), \quad\\
\text{s.t.}&\quad g_i(\bx) =\ba_i^\top \bx - b_i\leq 0, \quad i\in\{1,2,\ldots, m\},\\
&\quad h_j(\bx)=\bc_j^\top\bx-d_j=0, \quad j\in\{1,2,\ldots, p\},
\end{aligned}
$$
where $f(\cdot)$ is  continuously differentiable over $\real^n$, $\ba_1,\ba_2,\ldots, \ba_m,\bc_1,\bc_2,\ldots, \bc_p\in\real^n$, and $b_1,b_2,\ldots,b_m, d_1,d_2,\ldots,d_p\in\real$.
Then, 
\begin{enumerate}[(i)]
\item \textbf{Necessity of the KKT conditions.}
Let $\bx^*$ be a local minimum point of (P). Then, there exist $\lambda_1, \lambda_2, \ldots,\lambda_m\geq 0$ and $\mu_1, \mu_2, \ldots, \mu_p\in\real$ such that 
\begin{subequations}
\begin{align}
	&\nabla f(\bx^*) + \sum_{i=1}^{m} \lambda_i\ba_i +\sum_{j=1}^{p}\mu_j\bc_j=\bzero ,\label{equa:kkt1}\\
	&\lambda_i (\ba_i^\top\bx^*-b_i)=0, \quad i=\{1,2,\ldots, m\}. \label{equa:kkt2}
\end{align}
\end{subequations}
where \eqref{equa:kkt2} is called the \textit{complementary slackness} condition.

\item \textbf{Sufficiency in the convex case.} Suppose further that the function $f$ is convex over $\real^n$, and $\bx^*$ is a feasible solution of (P) for which there exist $\lambda_1, \lambda_2, \ldots,\lambda_m\geq 0$ and $\mu_1, \mu_2, \ldots, \mu_p\in\real$ such that \eqref{equa:kkt1} and \eqref{equa:kkt2} are satisfied. Then, $\bx^*$ is an optimal solution of (P).
\end{enumerate}

Moreover, we can define the \textit{Lagrangian function} (Definition~\ref{definition:lag_func}):
$$
L(\bx, \blambda, \bmu) = f(\bx) + \sum_{i=1}^{m} \lambda_ig_i(\bx) + \sum_{j=1}^{p} \mu_jh_j(\bx).
$$
The KKT condition \eqref{equa:kkt1} can be derived  from setting the gradient of the Lagrangian function with respect to $\bx$ to zero:
$$
\nabla_{\bx} L(\bx, \blambda, \bmu) 
= \nabla f(\bx) + \sum_{i=1}^{m} \lambda_i \nabla g_i(\bx) + \sum_{j=1}^{p} \mu_j \nabla h_j(\bx)
=\bzero.
$$
\end{theoremHigh}
To prove the desired result, we first need to establish the KKT conditions for linear inequality constraints as presented in \eqref{equation:kkt_ineq}. For brevity, these conditions are encapsulated within a single theorem.
\begin{proof}[of Theorem~\ref{theorem:kktcond}]
\textbf{(i).} Consider the equivalent problem
$$
(\text{P}^\prime)\qquad
\begin{aligned}
\min\;  f(\bx) \quad 
\text{s.t.}\quad &\ba_i^\top \bx - b_i\leq 0, \quad i\in\{1,2,\ldots, m\},\\
&\bc_j^\top\bx-d_j\leq 0, \quad j\in\{1,2,\ldots, p\}, \\
-&\bc_j^\top\bx+d_j\leq0, \quad j\in\{1,2,\ldots, p\},
\end{aligned}
$$
Since $ \bx^* $ is a local minimum solution of (P), it is also a local minimum solution of ($\text{P}^\prime$).

\paragraph{(i). KKT for linear inequality constraints.} To validate our claim, consider the following more general inequality-constrained optimization problem:
$$
\text{(Q)}\qquad
\begin{aligned}
\min\;  f(\bx) \quad 
\text{s.t.}\quad &\bg_i^\top \bx - f_i\leq 0, \quad i\in\{1,2,\ldots, l\}.
\end{aligned}
$$
Assume $ \widehatbx $ is a local minimum point of (Q). We aim to show that there exist multipliers $\kappa_1, \kappa_2, \ldots, \kappa_l\geq 0$ such that
\begin{equation}\label{equation:kkt_ineq}
\textbf{(KKT for Ineq.)}: \qquad \nabla f(\widehatbx) + \sum_{i=1}^{l} \kappa_i \bg_i = \bzero 
\qquad\text{and}\qquad
\kappa_i (\bg_i^\top\widehatbx -f_i) = 0, \,\forall i.
\end{equation}
This is known as the \textit{KKT conditions for linear inequality constraints}. The result in Theorem~\ref{theorem:kktcond} provides a more general result.
To see this,
it follows by Theorem~\ref{theorem:stat_point_uncons_convset} that $ \widehatbx $ is a stationary point, meaning that $ \nabla f(\widehatbx)^\top (\bx - \widehatbx) \geq 0 $ for every $ \bx \in \real^n $ satisfying $ \bg_i^\top \bx \leq f_i $ for any $ i = \{1, 2, \ldots, l\} $.
Making the change of variables $ \by \triangleq \bx - \widehatbx $, we obtain that $ \nabla f(\widehatbx)^\top \by \geq 0 $ for any $ \by \in \real^n $ satisfying $ \bg_i^\top (\by + \widehatbx) \leq f_i $ for any $ i = \{1, 2, \ldots, l\} $, that is, for any $ \by \in \real^n $ satisfying
$$
\bg_i^\top \by \leq f_i - \bg_i^\top \widehatbx, \,\,\forall i\quad\implies\quad 
\qquad \begin{cases}
\bg_i^\top \by \leq 0, & i \in \sI(\widehatbx), \\
\bg_i^\top \by \leq f_i - \bg_i^\top \widehatbx, & i \notin \sI(\widehatbx),
\end{cases}
$$
where $\sI(\widehatbx) \triangleq \{i \mid  \bg_i^\top \widehatbx = f_i\}$ denotes  the \textit{set of active constraints} (Definition~\ref{definition:setsp4prime}).
We will show that in fact the second set of inequalities in the latter system can be removed, that is, that the following implication is valid:
\begin{equation}\label{equation:kkt_pre1}
\bg_i^\top \by \leq 0 \text{ for all } i \in \sI(\widehatbx) \quad\implies\quad  \nabla f(\widehatbx)^\top \by \geq 0.
\end{equation}
Suppose then that $ \by $ satisfies $ \bg_i^\top \by \leq 0 $ for all $ i \in \sI(\widehatbx) $. Since $ f_i - \bg_i^\top \widehatbx > 0 $ for all $ i \notin \sI(\widehatbx) $, it follows that there exists a small enough $ \gamma > 0 $ for which $ \bg_i^\top (\gamma \by) \leq f_i - \bg_i^\top \widehatbx $ for  $i \notin \sI(\widehatbx)$. Thus, since in addition $ \bg_i^\top (\gamma \by) \leq 0 $ for any $ i \in \sI(\widehatbx) $, it follows by the stationarity condition that $ \nabla f(\widehatbx)^\top (\gamma \by) \geq 0 $, and hence that $ \nabla f(\widehatbx)^\top \by \geq 0 $. This shows the implication in \eqref{equation:kkt_pre1}.
Thus, by Farkas' lemma (Lemma~\ref{lemma:farka_lemm}), it follows that there exist $ \kappa_i \geq 0, i \in \sI(\widehatbx) $, such that
$
-\nabla f(\widehatbx) = \sum_{i \in \sI(\widehatbx)} \kappa_i \bg_i.
$
Defining $ \kappa_i \triangleq 0 $ for all $ i \notin \sI(\widehatbx) $, we get that $ \kappa_i (\bg_i^\top \widehatbx - f_i) = 0 $ for all $ i = \{1, 2, \ldots, l\} $ and that $\nabla f(\widehatbx) + \sum_{i=1}^{l} \kappa_i \bg_i = 0$, confirming \eqref{equation:kkt_ineq}.

\paragraph{(i). Main result.}
Returning to problem ($\text{P}^\prime$), the above derivation for problem (Q) shows that  that there exist multipliers $ \lambda_1, \lambda_2, \ldots, \lambda_m \geq 0 $ and $ \mu_1^+, \mu_1^-, \mu_2^+, \mu_2^-, \ldots, \mu_p^+, \mu_p^- \geq 0 $ such that
\begin{equation}\label{equation:kkt_tm1}
\nabla f(\bx^*) + \sum_{i=1}^{m} \lambda_i \ba_i + \sum_{j=1}^{p} \mu_j^+ \bc_j - \sum_{j=1}^{p} \mu_j^- \bc_j = \bzero
\end{equation}
and
\begin{equation}\label{equation:kkt_tm2}
\begin{aligned}
\lambda_i (\ba_i^\top \bx - b_i) &= 0, \quad i = \{1, 2, \ldots, m\}, \\
\mu_j^+ (\bc_j^\top \bx - d_j) &= 0, \quad j = \{1, 2, \ldots, p\}, \\
\mu_j^- (-\bc_j^\top \bx + d_j) &= 0, \quad j = \{1, 2, \ldots, p\}.
\end{aligned}
\end{equation}
We thus obtain that \eqref{equa:kkt1} and \eqref{equa:kkt1} are satisfied with $ \mu_j \triangleq \mu_j^+ - \mu_j^-, \,\, j = \{1, 2, \ldots, p\} $.
Note that $\mu_j$ are not necessarily nonnegative from this deduction.

\paragraph{(ii).} To prove the second part, we still need to prove the sufficiency of \eqref{equation:kkt_ineq} for the linearly inequality constrained problem (Q). 
\paragraph{(ii). Sufficiency of KKT for linear inequality constraints.}
Suppose that $ \widehatbx $ is a feasible solution of (Q) satisfying \eqref{equation:kkt_ineq}. Let $ \bx $ be any feasible solution of (Q). Define the function
$$
s(\bx) \triangleq f(\bx) + \sum_{i=1}^{l} \kappa_i (\bg_i^\top \bx - f_i).
$$
Then by \eqref{equation:kkt_ineq}, it follows that $ \nabla s(\widehatbx) = \bzero $. Since $ s $ is convex, it follows by Theorem~\ref{theorem:suff_sta_conv} that $ \widehatbx $ is a minimizer of $ s $ over $ \real^n $, which combined with \eqref{equation:kkt_ineq} implies that
$$
f(\bx) = f(\widehatbx) + \sum_{i=1}^{l} \kappa_i (\bg_i^\top \widehatbx - f_i) = s(\widehatbx) \leq s(\bx) = f(\bx) + \sum_{i=1}^{m} \kappa_i (\bg_i^\top \bx - f_i) \leq f(\bx),
$$
where the last inequality follows from the fact that $ \kappa_i \geq 0 $ and $ \bg_i^\top \bx - f_i \leq 0 $ for $ i = \{1, 2, \ldots, l\} $. This proves that $ \widehatbx $ is a global optimal solution of (Q).

\paragraph{(ii). Main result.}
Returning to the sufficiency of \eqref{equa:kkt1} and \eqref{equa:kkt2} for  equality and inequality constrained optimization,
suppose that $ \bx^* $ satisfies \eqref{equa:kkt1} and \eqref{equa:kkt1}. Then it also satisfies \eqref{equation:kkt_tm1} and \eqref{equation:kkt_tm2} with
$$
\mu_j^+ = [\mu_j]_+\triangleq \max\{\mu_j, 0\} \qquad \text{and} \qquad \mu_j^- = [\mu_j]_- = -\min\{\mu_j, 0\},
$$
which by the above deduction implies that $ \bx^* $ is an optimal solution of $ (\text{P}^\prime) $, and thus is also an optimal solution of (P).
\end{proof}

\subsection{KKT Conditions for Nonlinearly Constrained Problems}

Equipped with various feasible directions as defined in Definition~\ref{definition:lfd_sfd}, we provide some optimality conditions for the constrained optimization problem $(\text{P4}')$.
Recalling that Theorem~\ref{theorem:_nonconv_fea_loca_optim} establishes the first-order necessary condition for a general constrained optimization (P2), which states that  there are no feasible descent directions (Definition~\ref{definition:uncons_des_direct}). When we have algebraic specifications  of the constrained set in $(\text{P4}')$, we have the following result.
\begin{theoremHigh}[Geometric First-Order Necessary Condition-I: First-Order Necessary Condition for $(\text{P4}')$:]\label{theorem:_nonconv_fea_loca_optim_p4p}
Let $f:\sS\rightarrow \real$ be continuously differentiable where $\sS$ is given by $(\text{P4}')$,  and consider the constrained optimization problem $(\text{P4}')$. 
If $ \bx^* $ is a local optimal solution of $(\text{P4}')$ and $f, c_i$ for $i\in\mathcalI\cup\mathcalE$ are differentiable at $\bx^*$, then 
$$
\bd^\top \nabla f(\bx^*)\geq 0, \quad \text{for all $\bd\in\mathcalT_\sS(\bx^*)$},
$$
which is also equivalent to
$$
\mathcalT_\sS(\bx^*) \cap \sD(\bx^*) = \varnothing, \quad  \text{ where  $\sD(\bx^*) \triangleq \{\bd \mid \nabla f(\bx^*)^\top \bd < 0\}$}.
$$
\end{theoremHigh}
\begin{proof}[of Theorem~\ref{theorem:_nonconv_fea_loca_optim_p4p}]
The proof is by contradiction. Assume by contradiction that at point $\bx^*$ with $\mathcalT_\sS(\bx^*) \cap \{\bd \mid \big(\nabla f(\bx^*)\big)^\top \bd < 0\} \neq \varnothing$, let $\bd \in \mathcalT_\sS(\bx^*) \cap \{\bd \mid \big(\nabla f(\bx^*)\big)^\top \bd < 0\}$. According to the definition of tangent vectors, there exist $\{c_t\}_{t>0}$ and $\{\bd^\toptzero\}_{t>0}$ such that $\bx^* + c_t \bd^\toptzero \in \sS$, where $c_t \to 0$ and $\bd^\toptzero \to \bd$. Since $\big(\nabla f(\bx^*)\big)^\top \bd < 0$, for sufficiently large $t$, we have by the linear approximation theorem (Theorem~\ref{theorem:linear_approx}) that 
$$
\begin{aligned}
&f(\bx^* + c_t \bd^\toptzero) = f(\bx^*) + c_t \big(\nabla f(\bx^*)\big)^\top \bd^\toptzero +  o(c_t)\\
&= f(\bx^*) + c_t \big(\nabla f(\bx^*)\big)^\top \bd + c_t \big(\nabla f(\bx^*)\big)^\top (\bd^\toptzero - \bd) +  o(c_t)\\
&= f(\bx^*) + c_t \big(\nabla f(\bx^*)\big)^\top \bd + o(c_t)
< f(\bx^*),
\end{aligned}
$$
where the last inequality follows from the fact that $\frac{o(c_t)}{c_t}\rightarrow 0$ as $c_t\rightarrow 0^+$ and $\big(\nabla f(\bx^*)\big)^\top \bd<0$.
This contradicts the local minimality of $\bx^*$ and completes the proof.
\end{proof}

\begin{theoremHigh}[Geometric First-Order Necessary Condition-II]\label{theorem:geo_for_nec}
Let $ \bx^* $ be a local minimum of the problem $(\text{P4}')$
\begin{align*}
\min & \quad f(\bx) \\
\text{s.t.} & \quad g_i(\bx) \leq 0, \quad i = \{1, 2, \ldots, m\},\\
&\quad h_j(\bx) = 0, \quad j = \{1,2, \ldots, p\},
\end{align*}
where $f, g_1, g_2, \ldots, g_m$ are continuously differentiable functions over $ \real^n $. 
Suppose either:
\begin{enumerate}[(i)]
\item $ \bh(\bx) $ is a linear function, i.e., $ \bh(\bx) = \bA\bx - \bb $ for $ \bA \in \real^{p \times n} $, or
\item $ \nabla h_i(\bx^*), i = \{1,2, \ldots, p\} $, are linearly independent,
\end{enumerate}
\noindent Then
$$
\sD(\bx^*) \cap \sG(\bx^*) \cap \sH(\bx^*) = \varnothing. 
$$
\end{theoremHigh}
\begin{proof}[of Theorem~\ref{theorem:geo_for_nec}]
\textbf{(i).}
Assume $ \bh(\bx) $ is the linear function $ \bh(\bx) = \bA\bx - \bb $, whereby $ \nabla h_i(\bx^*) = \bA[i,:] $ for $ i \in \{1,2, \ldots, p\} $ and $ \sH(\bx^*) = \{\bd \mid \bA\bd = \bzero\} $. Suppose $ \bd \in \sD(\bx^*) \cap \sG(\bx^*) \cap \sH(\bx^*) $.
Since because $ \bd \in \sG(\bx^*) $, then there exist $\varepsilon_1>0$ such that for all $\eta\in(0, \varepsilon_1)$,  $ g_i(\bx^* + \eta \bd) < g_i(\bx^*) = 0 $ for $ i \in \sI(\bx^*) $ (Theorem~\ref{theorem:_nonconv_fea_loca_optim}). 
Since $g_i$ are continuous and we have  the fact that $ g_i(\bx^*) < 0 $ for $ i \notin \sI(\bx^*) $, there exists $\varepsilon_2$ such that $ g_i(\bx^* + \eta \bd) < 0 $ for all $\eta\in(0, \varepsilon_2)$. Let $\eta\in(0, \min\{\varepsilon_1, \varepsilon_2\}]$, we have $ g_i(\bx^*) < 0 $ for all $i\in\{1,2,\ldots,p\}$.
Furthermore $ \bh(\bx^* + \eta \bd) = (\bA\bx^* - \bb) + \eta \bA\bd = \bzero $.

Therefore, $ \bx^* + \eta \bd \in \sS $ for all $\eta\in(0, \min\{\varepsilon_1, \varepsilon_2\})$, where $\sS\triangleq\{\bx\mid \bg(\bx)\leq \bzero, \bh(\bx)=\bzero\}$. On the other hand, for all $\eta\in(0, \min\{\varepsilon_1, \varepsilon_2\})$, it holds that $ f(\bx^* + \eta \bd) < f(\bx^*) $. This contradicts the assumption that $ \bx^* $ is a local minimum of $(\text{P4}')$.

\paragraph{(ii).}
The proof when $ \bh(\bx) $ is nonlinear is a bit involved, which relies on the implicit function theorem (Theorem~\ref{theorem:implic_func_theorem}).
Let $ \bA \triangleq \nabla \bh(\bx^*) \in \real^{p \times n} $. Then $ \bA $ has full row rank, and its columns (along with corresponding elements of $ \bx^* $) can be re-arranged so that $ \bA = [\bB, \bN] $ and $ \bx^* = [\by^*; \bz^*] $, where $ \bB\in\real^{p\times p} $ is nonsingular and $\bN\in\real^{p\times (n-p)}$. Let $ \bz $ lie in a small neighborhood of $ \bz^* $: $\bz\in \sB(\bz^*, \varepsilon)$. Then, from the implicit function theorem (Theorem~\ref{theorem:implic_func_theorem}), there exists a function $ \bs(\bz): \real^{n-p}\rightarrow \real^p $ such that $ \bh(\bs(\bz), \bz) = \bzero $.

Suppose in contradiction that $ \bd \in \sD(\bx^*) \cap \sG(\bx^*) \cap \sH(\bx^*)  $, and denote $\bd$ as $ \bd \triangleq [\bq; \bp] $, where $\bq\in\real^p$ and $\bp\in\real^{n-p}$. Then $ \bzero = \bA\bd = \bB\bq + \bN\bp $, whereby $ \bq = -\bB^{-1}\bN\bp $.
Let 
$$
\bz(\theta) \triangleq \bz^* + \theta \bp, \qquad
\by(\theta) \triangleq \bs(\bz(\theta)) = \bs(\bz^* + \theta \bp), 
\qquad\text{and}\qquad
\bx(\theta) \triangleq [\by(\theta); \bz(\theta)], 
$$
where $\bz(\theta)  \in \sB(\bz^*, \varepsilon)$ is in the small neighborhood of $\bz^*$.
We will derive a contradiction by showing that $ \bd $ is an improving feasible direction, i.e., for small $ \theta > 0 $, $ \bx(\theta) $ is feasible and $ f(\bx(\theta)) < f(\bx^*) $.

To show feasibility of $ \bx(\theta) $, note that for $ \theta > 0 $ sufficiently small such that $\bz(\theta)  \in \sB(\bz^*, \varepsilon)$, it follows from the implicit function theorem (Theorem~\ref{theorem:implic_func_theorem}) that:
$$
\bh(\bx(\theta)) = \bh(\bs(\bz(\theta)), \bz(\theta)) = \bzero,
$$
and for $ i = \{1,2, \ldots, p\} $, we have:
$$
\frac{\partial h_i(\bx(\theta))}{\partial \theta} = \sum_{k=1}^p \frac{\partial h_i(\bs(\bz(\theta)), \bz(\theta))}{\partial y_k} \frac{\partial s_k(\bz(\theta))}{\partial \theta} + \sum_{k=1}^{n-p} \frac{\partial h_i(\bs(\bz(\theta)), \bz(\theta))}{\partial z_k} \frac{\partial z_k(\theta)}{\partial \theta} = 0. 
$$
Let $ r_k \triangleq \frac{\partial s_k(\bz(\theta))}{\partial \theta} $, and recall that $ \frac{\partial z_k(\theta)}{\partial \theta} = p_k $. At $ \theta = 0 $, the above equation system can then be re-written as $ \bzero = \bB\br + \bN\bp $, or $ \br = -\bB^{-1}\bN\bp = \bq $. Therefore,
$$
\frac{\partial x_k(\theta)}{\partial \theta} = d_k, \quad\text{for}\quad k = \{1,2, \ldots, n\}. 
$$
For $ i \in \sI(\bx^*) $,
$$
\begin{aligned}
g_i(\bx(\theta)) &= g_i(\bx^*) + \theta \frac{\partial g_i(\bx(\theta))}{\partial \theta} \bigg|_{\theta=0} + \abs{\theta} \alpha_i(\theta)
= \theta \sum_{k=1}^n \frac{\partial g_i(\bx(\theta))}{\partial x_k} \frac{\partial x_k(\theta)}{\partial \theta} \bigg|_{\theta=0}\\
&= \theta \nabla g_i(\bx^*)^\top \bd + \abs{\theta} \alpha_i(\theta),
\end{aligned}
$$
where $ \alpha_i(\theta) \rightarrow 0 $ as $ \theta \rightarrow 0 $. Hence $ g_i(\bx(\theta)) < 0 $ for all $ i = 1,2, \ldots, m $ for $ \theta > 0 $
sufficiently small, and therefore, $ \bx(\theta) $ is feasible for any $ \theta > 0 $ sufficiently small.

On the other hand,
$$ f(\bx(\theta)) = f(\bx^*) + \theta \nabla f(\bx^*)^\top \bd + \abs{\theta} \alpha(\theta) < f(\bx^*) $$
for $ \theta > 0 $ sufficiently small, where $ \alpha(\theta) \rightarrow 0 $ as $ \theta \rightarrow 0 $. But this contradicts the local optimality of $ \bx^* $. Therefore no such $ \bd $ can exist, and the theorem is proved.
\end{proof}
Note that Theorem~\ref{theorem:geo_for_nec} essentially states that if a point $ \bx^* $ is (locally) optimal, there is no direction $ \bd $ which is both a feasible direction (i.e., such that $ g(\bx^* + \lambda \bd) \leq 0 $ and $ h(\bx^* + \lambda \bd) \approx 0 $ for small $ \lambda > 0 $) and is an improving direction (i.e., such that $ f(\bx^* + \lambda \bd) < f(\bx^*) $ for small $ \lambda > 0 $), which makes sense intuitively.
Using Theorem~\ref{theorem:geo_for_nec}, we prove the following Fritz-John conditions.
\index{Fritz-John conditions}
\begin{theoremHigh}[Fritz-John Necessary Conditions for Equality/Inequality Constraints]\label{theorem:frijohn_eqineq}
Let $ \bx^* $ be a local minimum point  of the problem $(\text{P4}')$
\begin{align*}
\min & \quad f(\bx) \\
\text{s.t.} & \quad g_i(\bx) \leq 0, \quad i = \{1, 2, \ldots, m\},\\
&\quad h_j(\bx) = 0, \quad j = \{1,2, \ldots, p\},
\end{align*}
where $ f, g_1, \ldots, g_m $ are continuously differentiable functions over $ \real^n $. Then, there exist multipliers $ \lambda_0, \lambda_1, \ldots, \lambda_m \geq 0, \mu_1, \mu_2, \ldots,\mu_p\in\real $, which are not all zeros (i.e., nontriviality: $[\lambda_0, \blambda, \bmu]^\top\neq \bzero$), such that
\begin{subequations}
\begin{align}
	\lambda_0 \nabla f(\bx^*) &+ \sum_{i=1}^m \lambda_i \nabla g_i(\bx^*) +\sum_{j=1}^{p} \mu_j\nabla h_j(\bx^*) = \bzero,\label{equation:frijob_full1} \\
	&\lambda_i g_i(\bx^*) = 0, \quad i = \{1, 2, \ldots, m\}.
\end{align}
\end{subequations}
Note that \eqref{equation:frijob_full1} can be compactly written as $\lambda_0 \nabla f(\bx^*) = \nabla \bg(\bx)^\top \blambda + \nabla \bh(\bx)^\top \bmu$ using Jacobian matrices (see \eqref{equation:jacob_inkkt}).

\end{theoremHigh}
\begin{proof}[of Theorem~\ref{theorem:frijohn_eqineq}]
If the vectors $ \{\nabla h_i(\bx^*)\} $ are linearly dependent, then there exists $ \bmu \neq \bzero $ such that $ \sum_{j=1}^p \mu_j \nabla h_i(\bx^*) = \bzero $. Setting $ [\lambda_0, \blambda] \triangleq \bzero $ establishes the result.

Suppose now that the vectors $ \{\nabla h_i(\bx^*)\} $ are linearly independent. Then we can apply Theorem~\ref{theorem:geo_for_nec} and assert that $ \sD(\bx^*) \cap \sG(\bx^*) \cap \sH(\bx^*) = \varnothing $. Without loss of generality, assume for simplicity that $ \sI(\bx^*) = \{1,2, \ldots, l\} $ with $l\leq m$. Let
$$
\bD \triangleq \begin{bmatrix}
\nabla f(\bx^*)^\top \\
\nabla g_1(\bx^*)^\top \\
\vdots \\
\nabla g_l(\bx^*)^\top
\end{bmatrix}
\quad\text{and}\quad
\bF \triangleq \begin{bmatrix}
\nabla h_1(\bx^*)^\top \\
\vdots \\
\nabla h_p(\bx^*)^\top
\end{bmatrix}.
$$
Then there is no $ \bd $ that satisfies $ \bD\bd < \bzero $ and  $ \bF\bd = \bzero $. From the system (1 $\&$ 2) of Lemma~\ref{lemma:farka_lemm}, there exists $ (\lambda_0, \lambda_1, \ldots, \lambda_l) $ and $ (\mu_1,\mu_2 \ldots, \mu_p) $ such that
$$
\lambda_0 \nabla f(\bx^*) + \sum_{i=1}^l \lambda_i \nabla g_i(\bx^*) + \sum_{j=1}^p \mu_j \nabla h_j(\bx^*) = \bzero,
$$
with $ \lambda_0 + \lambda_1 + \cdots + \lambda_l = 1 $ and $ (\lambda_0, \lambda_1, \ldots, \lambda_l) \geq 0 $. Define $ \lambda_{l+1}, \ldots, \lambda_m \triangleq 0 $. Then, $ (\lambda_0, \blambda) \geq \bzero $ and $ (\lambda_0, \blambda) \neq \bzero $, and for any $ i $, either $ g_i(\bx^*) = 0 $, or $ \lambda_i = 0 $. Furthermore,
$$
\lambda_0 \nabla f(\bx^*) + \sum_{i=1}^m \lambda_i \nabla g_i(\bx^*) + \sum_{i=1}^p \mu_i \nabla h_i(\bx^*) = \bzero,
$$
which obtains the desired result.
\end{proof}

\index{KKT conditions}
\index{Linear independence constraint qualification}
\index{LICQ}

A major drawback of the Fritz-John conditions is in the fact that they allow $ \lambda_0 $ to be zero. The case $ \lambda_0 = 0 $ is not particularly informative since condition \eqref{equation:frijob_full1} then becomes
$$ \sum_{i=1}^m \lambda_i \nabla g_i(\bx^*) +\sum_{j=1}^{p}\mu_j\nabla h_j(\bx^*) =\bzero , $$
which indicates that the gradients of the active constraints $ \{\nabla g_i(\bx^*)\}_{i \in \sI(\bx^*)} $ and the equality constraints $\{\nabla h_j(\bx^*)\}_{j\in\{1,2,\ldots,p\}}$ are linearly dependent.
This situation provides no information about the objective function, implying many points satisfying the Fritz-John conditions may not be local minima.

If we add an assumption that the gradients of the active/equality constraints are linearly independent at $\bx^*$, then we can establish the KKT conditions, which are the same
as the Fritz-John conditions with $\lambda_0 = 1$.
The condition that the gradients of the active/quality constraints are linearly independent is one of many types of assumptions that are referred to in the literature as ``\textit{constraint qualifications}." \citep{bazaraa2006nonlinear, beck2014introduction}.
We summarize  this \textit{linear independence constraint qualification (LICQ)} result in the following theorem, which follows directly from the Fritz-John conditions stated in  Theorem~\ref{theorem:frijohn_eqineq}.
\begin{theoremHigh}[KKT  Conditions for Equality/Inequality under LICQ]\label{theorem:kkt_licq}
Consider the optimization problem $(\text{P4}')$
\begin{align*}
\min & \quad f(\bx) \\
\text{s.t.} & \quad g_i(\bx) \leq 0, \quad i = \{1, 2, \ldots, m\},\\
&\quad h_j(\bx) = 0, \quad j = \{1,2, \ldots, p\},
\end{align*}
where $ f, g_1, g_2, \ldots, g_m, h_1, h_2, \ldots, h_p $ are continuously differentiable functions over $ \real^n $. 
\begin{enumerate}[(i)]
\item \textbf{Necessity of the KKT conditions.} Let $ \bx^* $ be a local minimum of $(\text{P4}')$. Suppose that the gradients of the active constraints and the equality constraints
$$ \{\nabla g_i(\bx^*) : i \in \sI(\bx^*)\} \cup \{\nabla h_j(\bx^*) : j = 1, 2, \ldots, p\} $$
are linearly independent (where $ \sI(\bx^*) \triangleq \{i : g_i(\bx^*) = 0\} $, i.e., LICQ in Definition~\ref{definition:licq}). Then there exist multipliers $ \lambda_1, \lambda_2, \ldots, \lambda_m \geq 0 $ and $ \mu_1, \mu_2, \ldots, \mu_p \in \real $ such that
\begin{subequations}
\begin{align}
	\nabla f(\bx^*) &+ \sum_{i=1}^m \lambda_i \nabla g_i(\bx^*) + \sum_{j=1}^p \mu_j \nabla h_j(\bx^*) = \bzero, \label{equation:kkt_licq1}\\
	&\lambda_i g_i(\bx^*) = 0, \quad i = \{1, 2, \ldots, m\}. \label{equation:kkt_licq2}
\end{align}
\end{subequations}
\item \label{suff:kkt_licq} \textbf{Sufficiency in the convex case.}
Suppose further that the function $f, g_1, g_2\ldots, g_m$ are \textbf{convex} over $\real^n$, $h_1, h_2, \ldots, h_p$ are \textbf{affine} functions, and $\bx^*$ is a feasible solution of $(\text{P4}')$ for which there exist $\lambda_1, \lambda_2, \ldots,\lambda_m\geq 0$ and $\mu_1, \mu_2, \ldots, \mu_p\in\real$ such that \eqref{equation:kkt_licq1} and \eqref{equation:kkt_licq2} are satisfied. Then, $\bx^*$ is an optimal solution of $(\text{P4}')$.
\end{enumerate}
\end{theoremHigh}
\begin{proof}[of Theorem~\ref{theorem:kkt_licq}]
The first part is a direct consequence of the Fritz-John necessary conditions. 
For the sufficiency under the convexity, 
let $ \bx $ be a feasible solution of $(\text{P4}')$. We will show that $ f(\bx) \geq f(\bx^*) $. Note that the function
$$ s(\bx) = f(\bx) + \sum_{i=1}^m \lambda_i g_i(\bx) + \sum_{j=1}^p \mu_j h_j(\bx) $$
is convex, and since $ \nabla s(\bx^*) = \nabla f(\bx^*) + \sum_{i=1}^m \lambda_i \nabla g_i(\bx^*) + \sum_{j=1}^p \mu_j \nabla h_j(\bx^*) = \bzero $, it follows by Theorem~\ref{theorem:suff_sta_conv} that $ \bx^* $ is a minimizer of $ s(\cdot) $ over $ \real^n $, and in particular $ s(\bx^*) \leq s(\bx) $ for any feasible point $\bx$.
Therefore, the KKT conditions and the feasibility of $\bx^*$ imply that 
\begin{align*}
f(\bx^*) &= f(\bx^*) + \sum_{i=1}^m \lambda_i g_i(\bx^*) + \sum_{j=1}^p \mu_j h_j(\bx^*)
= s(\bx^*)  \\
&\leq s(\bx) 
= f(\bx) + \sum_{i=1}^m \lambda_i g_i(\bx) + \sum_{j=1}^p \mu_j h_j(\bx)
\leq f(\bx).
\end{align*}
This completes the proof.
\end{proof}

\index{KKT points}
The above theorem introduces two concepts: KKT points and regularity, which we formulate in the following definitions.
\begin{definition}[KKT Points]\label{definition:kkt_point}
Consider the optimization problem $(\text{P4}')$ in \eqref{equation:P4_in_kkt}, where $ f, g_1, \ldots, g_m, h_1, h_2, \ldots, h_p $ are continuously differentiable functions over $ \real^n $. A feasible point $ \bx^* $ is called a \textit{KKT point} if there exist $ \lambda_1, \lambda_2, \ldots, \lambda_m \geq 0 $ and $ \mu_1, \mu_2, \ldots, \mu_p \in \real $ such that
$$
\begin{aligned}
\nabla f(\bx^*) &+ \sum_{i=1}^m \lambda_i \nabla g_i(\bx^*) +\sum_{j=1}^{p} \mu_j\nabla h_j(\bx^*) = \bzero, \\
&\lambda_i g_i(\bx^*) = 0, \quad i = \{1, 2, \ldots, m\}.
\end{aligned}
$$
\end{definition}

\index{Regularity}
\begin{definition}[Regularity]
Consider again the optimization problem $(\text{P4}')$ in \eqref{equation:P4_in_kkt}, where $ f, g_1, \ldots, g_m, h_1, h_2, \ldots, h_p $ are continuously differentiable functions over $ \real^n $. A feasible point $ \bx^* $ is called \textit{regular} if the gradients of the active constraints among the inequality constraints and of the equality constraints
$$ \big\{\nabla g_i(\bx^*) : i \in \sI(\bx^*)\big\} \cup \big\{\nabla h_j(\bx^*) : j = \{1, 2, \ldots, p\}\big\} $$
are linearly independent.
\end{definition}

\subsection{KKT Conditions for Convex Optimization Problems}

Specifically, for convex optimization problems (Section~\ref{section:convex_optimiz}), we can use another set of conditions to state the KKT conditions, i.e., \textit{Slater's condition}.
\begin{theoremHigh}[KKT Conditions under  Slater's Condition]\label{theorem:kktslat1_slater}
Consider the minimization problem
\begin{equation}\label{equation:kktslat1_slater1}
\begin{aligned}
	\min &\quad  f(\bx) \\
	\text{s.t.}
	&\quad g_i(\bx) \leq 0, \quad i\in\{1,2,\ldots, m\}, \quad \text{(convex)}\\
	&\quad h_j(\bx)= 0, \quad j\in\{1,2,\ldots, p\}, \quad \text{(affine)}\\
\end{aligned}
\end{equation}
where $f$ is  continuously differentiable over $\real^n$, $g_1, \ldots, g_m$ are \textbf{convex} and continuously differentiable over $\real^n$, $h_1, h_2, \ldots h_p$  are \textbf{affine}. 
Then, we have the following results:
\begin{enumerate}[(i)]
\item \textbf{Necessity of the KKT conditions.}
Let $\bx^*$ be a local minimum point of \eqref{equation:kktslat1_slater1}. 
Suppose that $\{\nabla h_j(\bx^*)\}_{j\in\{1,2,\ldots,p\}}$ are linearly independent~\footnote{When the number of affine constraint is zero, this condition can be relaxed.}, and there exists $\widehatbx\in\real^n$ such that the following  Slater's condition holds:
$$
\textbf{(Slater's condition)}: \quad
\left\{
\begin{aligned}
	g_i(\widehatbx) &\textcolor{mylightbluetext}{<}0, \quad i\in\{1,2,\ldots, m\}, \\
	h_j(\widehatbx)&= 0, \quad j\in\{1,2,\ldots, p\},\\
\end{aligned}
\right.
$$
Then, there exist $\lambda_1, \lambda_2, \ldots,\lambda_m\geq 0$ and $\mu_1, \mu_2, \ldots, \mu_p\in\real$ such that 
\begin{subequations}
\begin{align}
	&\nabla f(\bx^*) + \sum_{i=1}^{m} \lambda_i\nabla g_{i}(\bx^*) +\sum_{j=1}^{p}\mu_j \nabla h_j(\bx^*)  =\bzero,\label{equa:kktslat1_gen1}\\
	&\lambda_i g_i(\bx^*)=0, \quad i=\{1,2,\ldots, m\}.\label{equa:kktslat1_gen2}
\end{align}
\end{subequations}

\item \textbf{Sufficiency in the convex case.} Suppose further that the function $f$ is \textbf{convex} over $\real^n$, and $\bx^*$ is a feasible solution of \eqref{equation:kktslat1_slater1} for which there exist $\lambda_1, \lambda_2, \ldots,\lambda_m\geq 0$ and $\mu_1, \mu_2, \ldots, \mu_p\in\real$ such that \eqref{equa:kktslat1_gen1} and \eqref{equa:kktslat1_gen2} are satisfied. Then, $\bx^*$ is an optimal solution of \eqref{equation:kktslat1_slater1}. 
This is equivalent to the claim in Theorem~\ref{theorem:kkt_licq}\ref{suff:kkt_licq}.
\end{enumerate}

\end{theoremHigh}
\begin{proof}[of Theorem~\ref{theorem:kktslat1_slater}]
As in the proof of the necessity of the KKT conditions under LICQ (Theorem~\ref{theorem:kkt_licq}), since $ \bx^* $ is an optimal solution of \eqref{equation:kktslat1_slater1}, then the Fritz-John necessary conditions (Theorem~\ref{theorem:frijohn_eqineq}) are satisfied. That is, there exist nonnegative scalars $ \widetilde{\lambda}_0, \widetilde{\lambda}_1, \widetilde{\lambda}_2, \ldots,\widetilde{\lambda}_m\geq 0  $ and real scalars $\widetilde{\mu}_1, \widetilde{\mu}_2, \ldots, \widetilde{\mu}_p\in\real$, which are not all zeros, such that
\begin{equation}\label{equation:kktslat1_slater_prov1}
\begin{aligned}
	&\widetilde{\lambda}_0\nabla f(\bx^*) + \sum_{i=1}^{m} \widetilde{\lambda}_i\nabla g_{i}(\bx^*) +\sum_{j=1}^{p}\widetilde{\mu}_j \nabla h_j(\bx^*)  =\bzero,\\
	&\widetilde{\lambda}_i g_i(\bx^*)=0, \quad i=\{1,2,\ldots, m\}.
\end{aligned}
\end{equation}
To complete the proof, we need to show that $ \widetilde{\lambda}_0 > 0 $, and then the conditions \eqref{equa:kktslat1_gen1} and \eqref{equa:kktslat1_gen2} will be satisfied with $ \lambda_i = \frac{\widetilde{\lambda}_i}{\widetilde{\lambda}_0}, i = \{1, 2, \ldots, m\}$ and $ \mu_j = \frac{\widetilde{\mu}_j}{\widetilde{\lambda}_0}, j = \{1, 2, \ldots, p\}$. To prove that $ \widetilde{\lambda}_0 > 0 $, assume, for contradiction, that it is zero then
\begin{equation}\label{equation:kktslat1_slater_prov2}
\sum_{i=1}^{m} \widetilde{\lambda}_i\nabla g_{i}(\bx^*) +\sum_{j=1}^{p}\widetilde{\mu}_j \nabla h_j(\bx^*) =\bzero.
\end{equation}
By the convexity of $g_i$, for all $ i = \{1, 2, \ldots, m\} $, we have
$$
\begin{aligned}
&0 > g_i(\widehatbx) \geq g_i(\bx^*) + \nabla g_i(\bx^*)^\top (\widehatbx - \bx^*).\\
\end{aligned}
$$
Multiplying the $ i $-th inequality by $ \widetilde{\lambda}_i $ and summing over $ i = \{1, 2, \ldots, m\} $, we obtain
$$
\begin{aligned}
0 > &\sum_{i=1}^m \widetilde{\lambda}_i g_i(\bx^*) + \left( \sum_{i=1}^m \widetilde{\lambda}_i \nabla g_i(\bx^*) \right)^\top (\widehatbx - \bx^*)\\
\end{aligned}
$$
where the inequality is strict since not all the $ \widetilde{\lambda}_i $ are zero (use the fact that  $\{\nabla h_j(\bx^*)\}_{j\in\{1,2,\ldots,p\}}$ are linearly independent such that $\widetilde{\lambda}_i$'s are nonnegative and sum to 1; see the proof of the Fritz-John conditions in Theorem~\ref{theorem:frijohn_eqineq}). Plugging the identities \eqref{equation:kktslat1_slater_prov1} and \eqref{equation:kktslat1_slater_prov2} into the above inequality, we  obtain the impossible statement that $ 0 > 0 $, thus establishing the result. 
\end{proof}

\begin{problemset}
\item Using the proof of Theorem~\ref{theorem:psd_hess_conv}, prove the equivalence among \eqref{equation:scss_func1}, \eqref{equation:scss_func2}, and \eqref{equation:scss_func3_ra} for the definition of strongly convexity and strongly smoothness.

\item Prove the properties of Bregman distances as  outlined in Remark~\ref{remark:bregnan_dist}.

\item Let $ f \in C^{2,2}_L(\real^n) $. Show that  for all $ \bx, \by \in \real^n $,  the following inequalities hold:
$$
\begin{aligned}
&\normtwo{\nabla f(\by) - \nabla f(\bx) - \nabla^2 f(\bx)(\by - \bx)} \leq \frac{L}{2} \normtwo{\by - \bx}^2;\\
&\abs{f(\by) - f(\bx) - \innerproduct{\nabla f(\bx), \by - \bx} - \frac{1}{2} \innerproduct{\nabla^2 f(\bx)(\by - \bx), \by - \bx}} \leq \frac{L}{6} \normtwo{ \by - \bx}^3.
\end{aligned}
$$
Additionally, if $\normtwo{\bx-\by}\leq M$, demonstrate that:
$$
\nabla^2 f(\bx) -L M \bI \preceq \nabla^2 f(\by) \preceq \nabla^2 f(\bx)+ LM  \bI.
$$
\textit{Hint: Refer to Theorem~\ref{theorem:equi_gradsch_smoo} and Exercise~\ref{exercise:cl2111_hess_bound}.}

\item \label{prob:dist_hyper} \textbf{Distance between a vector and a hyperplane.} Given a nonzero vector $\bzero\neq \ba\in\real^n$ and a scalar $\beta$,  we define  the hyperplane $H(\ba, \beta) \triangleq \{\bx\in\real^n:\ba^\top\bx+\beta=0\}$. For any vector $\by\in\real^n$, show that the distance between the vector $\by$ and the hyperplane $H(\ba, \beta)$ is given by 
\begin{equation}
	d(\by, H(\ba, \beta)) = \frac{\abs{\ba^\top\by + \beta}}{\normtwo{\ba}}.
\end{equation}
\textit{Hint: Select  two random points on the plane and first show that $\ba$ is orthogonal to the plane.}

\end{problemset}

\newpage
\chapter{Descent Methods}\label{chapter:gradient-descent} 
\begingroup
\hypersetup{linkcolor=structurecolor,
linktoc=page,  
}
\minitoc \newpage
\endgroup

\index{Gradient descent}
\section{Descent Methods and Gradient Descent Method}\label{section:gradient-descent-all} 
The \textit{gradient descent} (GD) method is a specific type of descent method used to find the (local or global) minimum of a differentiable function, whether convex or non-convex. This function is commonly referred to as the \textit{cost function} (also known as the \textit{loss function} or \textit{objective function}).
It stands out as one of the most popular algorithms to perform optimization and by far the
most common way to optimize machine learning, deep learning, and various optimization problems. 
This is especially true for optimizing neural networks and transformer networks \citep{lecun2015deep, goodfellow2016deep, vaswani2017attention}. 
In the context of machine learning, the cost function measures the difference between a model’s predicted output and the actual output. Neural networks, transformer networks, and machine learning models in general seek to find a set of parameters $\bx\in \real^n$ (also known as weights) that optimize an objective function $f(\bx)$. This is expressed as the unconstrained optimization problem (P1) in Definition~\ref{definition:opt_probs_all}:
$$
\textbf{(P1)}:\qquad \text{Find}\quad \bx^* = \mathop{\argmin}_{\bx} f(\bx).
$$
In some special cases, (P1) can be solved analytically by finding the stationary points of the function using the first-order optimality condition in Theorem~\ref{theorem:fermat_fist_opt}: $\nabla f(\bx)=\bzero$. However, in most cases, the problem must be solved using iterative methods.
All the optimization methods discussed in this book, including gradient descent, fall under the category of \textit{iterative methods}. Denoting $t=1,2,\ldots$ as the iteration number, these methods generate a sequence of vectors:
\begin{equation}
\bx^{(1)}, \bx^{(2)}, \ldots, \bx^{(T)} \in \dom(f)
\end{equation}
\footnote{Some texts denote the starting point as $\bx^{(0)}$, but in this book, we use $\bx^{(1)}$.}
such that as $T\rightarrow \infty$, the sequence converges to the optimal solution $\bx^*$, and the objective function value $f(\bx^{(T)})$ approaches the optimal minimum $f(\bx^*)$, under certain mild conditions.
At each iteration $t$, an \textit{update step (or a descent step)} $\bh^\toptzero$ is applied to update the parameters. Denoting the parameters at the $t$-th iteration as $\bx^\toptzero$,  the update rule  is given by:
\begin{equation}\label{equation:sg_update_rule}
(\textbf{GD update}):\qquad 
\bx^\toptone \leftarrow \bx^\toptzero + \bh^\toptzero.
\end{equation}

\index{Learning rate}
\index{Strict SGD}
\index{Mini-batch SGD}
\index{Local minima}
\paragrapharrow{Gradient descent.}
The most basic form of gradient descent is the \textit{vanilla update}, where the parameters move in the opposite direction of the gradient. This follows the \textit{steepest descent direction} since gradients are orthogonal to level curves (also known as level surfaces, see Lemma~\ref{lemm:direction-gradients}):
\begin{equation}\label{equation:gd-equaa-gene}
\bh^\toptzero = -\eta_t \bg^\toptzero\triangleq -\eta_t \nabla f(\bx^\toptzero),
\end{equation}
where the positive value $\eta_t$ denotes the \textit{learning rate  (or stepsize, step length, step size)} that depends on specific problems.
The term $\bg^\toptzero\triangleq\nabla f(\bx^\toptzero) \in \real^n$ represents the gradient of the parameters.
The {learning rate} $\eta_t$ controls how large of a step to take in the direction of negative gradient so that we can reach a (local) minimum.
The method that follows the negative gradient direction (i.e., $\bd^\toptzero \triangleq -\nabla f(\bx^\toptzero)$ in Algorithm~\ref{alg:struc_gd_gen}) is called the \textit{steepest descent method (or gradient method)}. 
The choice of descent direction is ``the best" (locally; see \eqref{equation:steep_des}) and we could combine it with an exact line search to determine the learning rate  (Section~\ref{section:exact_line_search}). A method like this converges, but the final convergence is linear and often very slow. 

Examples in \citet{madsen2010and, boyd2004convex} show how the gradient descent method with exact line search and finite computer precision can fail to find the minimizer of a second degree polynomial. However, for many problems, it performs well in the early stages of the iterative process.
Considerations like this has lead to the so-called \textit{hybrid methods}, which---as the name suggests---are based on two different methods. One which is good in the initial stage, like the \textit{gradient method}, and another method which is good in the final stage, like \textit{Newton's method} (Section~\ref{section:new_methods}). A key challenge with hybrid methods is designing an effective mechanism to switch between the two approaches at the appropriate time.

The descent property of such a step in \eqref{equation:gd-equaa-gene} is guaranteed by Definition~\ref{definition:uncons_des_direct} and Lemma~\ref{lemma:descent_property} or  (by setting $\bB\triangleq\bI$ in Theorem~\ref{theorem:uncons_des_dir}).
In \eqref{equation:sg_update_rule}, $\bh^\toptzero$ is referred to as a \textit{descent step}. While a direction $\bd^\toptzero$ in Definition~\ref{definition:uncons_des_direct} satisfying $\innerproduct{\bg^\toptzero,\bd^\toptzero}<0$ is called a \textit{descent direction}. In most cases, the relationship between the descent step and descent direction follows a scale by the learning rate:
\begin{equation}\label{equation:desce_step_direc}
\begin{aligned}
	&\textbf{(Descent direction)}: \qquad&&\bd^\toptzero = -\bg^\toptzero;\\
	&\textbf{(Descent step)}:\qquad &&\bh^\toptzero=\eta_t\bd^\toptzero.
\end{aligned}
\end{equation}
In many cases, when the learning rate is equal to 1, the above two terms are used \textbf{interchangeably},  then the descent direction and the descent step are the same; for example, Newton's method in Chapter~\ref{chapter:second_order}.

\paragrapharrow{Gradient descent by calculus}
An intuitive analogy to understand gradient descent is to imagine the path of a river starting from a mountain peak and flowing downhill to reach the lowest point at its base.
Similarly, the goal of gradient descent is to find the lowest point in the landscape defined by the objective function $f(\bx)$, where $\bx$ represents  a $n$-dimensional input variable. Our task is to use an algorithm that guides us to a (local) minimum of $f(\bx)$.
To better understand this process, consider moving a ball a small distance $d_1$ along the $x_1$ axis, a small amount $d_2$ along the $x_2$ axis, and so on up to $d_n$ along the $x_n$ axis. 
Calculus informs us of the variation in the objective function  $f(\bx)$ as follows:
$$
\Delta f(\bx) \approx \frac{\partial f}{\partial x_1}d_1 + \frac{\partial f}{\partial x_2}d_2 + \ldots + \frac{\partial f}{\partial x_n}d_n.
$$
Our challenge is to choose $d_1, d_2, \ldots, d_n$ such that they cause $\Delta f(\bx)$ to be negative, thereby decreasing the objective function towards minimization.
Let $\bd=[d_1,d_2, \ldots , d_n]^\top$ denote  the vector of  changes in $\bx$, and let $\nabla f(\bx)=[\frac{\partial f}{\partial x_1},\frac{\partial f}{\partial x_2}, \ldots , \frac{\partial f}{\partial x_n}]^\top$ denotes the gradient vector of $f(\bx)$ \footnote{Note the difference between $\Delta f(\bx)$ and $\nabla f(\bx)$.}. Then it follows that
$$
\Delta f(\bx) \approx \frac{\partial f}{\partial x_1}d_1 + \frac{\partial f}{\partial x_2}d_2 +\ldots  + \frac{\partial f}{\partial x_n}d_n = \innerproduct{\nabla f(\bx), \bd}.
$$ 
In the context of descending the function, our aim is to ensure that $\Delta f(\bx)$ is negative. 
This ensures that moving from $\bx^\toptzero$ to $\bx^\toptone = \bx^\toptzero+\bd^\toptzero$  (from $t$-th iteration to $(t+1)$-th iteration) results in a reduction of the loss function $f(\bx^\toptone) = f(\bx^\toptzero) + \Delta f(\bx^\toptzero)$, given that $\Delta f(\bx^\toptzero) \leq 0$.
It can be demonstrated that if the update step is defined as $\bd^\toptzero=-\eta_t \nabla f(\bx^\toptzero)$, where $\eta_t$ is the learning rate, the following relationship holds:
$$
\Delta  f(\bx^\toptzero) \approx -\eta_t \nabla f(\bx^\toptzero)^\top\nabla f(\bx^\toptzero) = -\eta_t\normtwobig{\nabla f(\bx^\toptzero)}^2 \leq 0. 
$$
Specifically,  $\Delta  f(\bx^\toptzero) < 0$ unless we are already at the optimal point with zero gradients.
This analysis validates the approach of gradient descent. 


\index{Convex functions}
\paragrapharrow{Gradient descent for convex functions.}
We further explore the application of gradient descent in (unconstrained) convex problems (see Section~\ref{section:convex_opt_uncons} for optimality conditions of this problem).
If the objective function $f(\bx)$ is (continuously differentiable) convex, then the relationship $\innerproduct{\nabla f(\bx^\toptzero), (\bx^\toptone-\bx^\toptzero)}\geq 0$ implies $f(\bx^\toptone) \geq f(\bx^\toptzero)$. This can be derived from the gradient inequality of a continuously differentiable convex function, i.e., $f(\bx^\toptone)- f(\bx^\toptzero)\geq \innerproduct{\nabla f(\bx^\toptzero), (\bx^\toptone-\bx^\toptzero)}$; see Theorem~\ref{theorem:conv_gradient_ineq}. 

In this sense, to ensure a reduction in the objective function from the point $\bx^\toptzero$ to $\bx^\toptone$, it is imperative to ensure  $\innerproduct{\nabla f(\bx^\toptzero), (\bx^\toptone-\bx^\toptzero)}\leq 0$. 
In the context of  gradient descent, the choice of $\eta_t\bd^\toptzero = \bx^\toptone-\bx^\toptzero$ aligns with the negative gradient $-\nabla f(\bx^\toptzero)$. However, there are many other descent methods, such as \textit{(non-Euclidean) greedy descent}, \textit{normalized steepest descent}, \textit{Newton step}, and so on. The core principle of these methods is to ensure that $\innerproduct{\nabla f(\bx^\toptzero), (\bx^\toptone-\bx^\toptzero)}= \innerproduct{\nabla f(\bx^\toptzero), \bd^\toptzero} \leq 0$, provided the objective function is convex.

\paragrapharrow{Gradient descent with momentum.}
\textit{Gradient descent with momentum} is an improvement over basic gradient descent, frequently used in machine learning and deep learning to minimize the loss function and update model parameters (Section~\ref{section:sgd_momentum}). While standard gradient descent updates parameters solely based on the current gradient, momentum-based gradient descent introduces a \textit{momentum} term to accelerate convergence and smooth the optimization path.

In this approach, the momentum term enables the algorithm to build velocity in directions with a steady but small gradient, helping it overcome local minima and saddle points. By incorporating a fraction of the previous update into the current one, this technique mimics inertia, allowing the algorithm to continue moving in the same direction despite minor fluctuations in the gradient. Consequently, this method not only speeds up convergence but also reduces oscillations, particularly in regions where the surface curvature varies significantly across different dimensions.
At each iteration $t$, the process involves two key steps:
\begin{subequations}\label{equation:gd_with_momentum}
\begin{align}
\textbf{(Velocity update)}:\qquad \bd^\toptzero &\leftarrow \rho \bd^\toptminus - \eta_t \nabla f(\bx^\toptzero);\\
\textbf{(Parameter update)}:\qquad \bx^\toptone &\leftarrow \bx^\toptzero + \bd^\toptzero.
\end{align}
\end{subequations}
By incorporating past gradients into the update rule, gradient descent with momentum enables more efficient traversal across the error surface, particularly in complex landscapes, leading to faster convergence and improved performance.
In summary, The gradient descent with momentum approach is advantageous for the following reasons:
\begin{itemize}
\item At saddle points, the gradient of the cost function becomes nearly zero or entirely negligible. This results in minimal or no updates to the weights, causing the network's learning process to stagnate and effectively halt.
\item The trajectory taken by the gradient descent method tends to be quite erratic, even when employing mini-batch processing. This jittery path can impede efficient convergence towards the minimum (Section~\ref{section:sgd_momentum}).
\end{itemize}
We will analyze the behaviors of standard gradient descent and momentum-based gradient descent for quadratic models in Sections~\ref{section:quadratic_vanilla_GD}, \ref{section:quadratic-in-steepestdescent}, and \ref{section:quadratic-in-momentum}.

\paragrapharrow{Steepest descent.}
The linear approximation theorem (Theorem~\ref{theorem:linear_approx}) states that 
\begin{equation}\label{equation:gd_gree_tay11}
	f(\bx^\toptzero + \eta \bd) = f(\bx^\toptzero)+ \eta \bd^\top \nabla  f(\bx^\toptzero )+ \mathcalO(\normtwo{\eta \bd}^2).
\end{equation}
From \eqref{equation:gd_gree_tay11} and by the definition of directional derivative, we observe that when taking  a step $\eta \bd$ with a positive stepsize $\eta$,  the relative reduction in function value satisfies
$$
\lim_{\eta \rightarrow 0} \frac{f(\bx^\toptzero) - f(\bx^\toptzero + \eta \bd)}{\eta \normtwo{\bd}} = -\frac{1}{\normtwo{\bd}} \bd^{\top} \nabla f(\bx^\toptzero) = \normtwobig{\nabla f(\bx^\toptzero)} \cos (\theta),
$$
where $\theta$ is the angle between the vectors $\bd$ and $-\nabla f(\bx^\toptzero)$. This equation indicates that we get the greatest gain rate if $\theta = 0$, meaning the optimal descent direction is the steepest descent direction $\bd_{\text{sd}}^\toptzero$, given by
\begin{equation}\label{equation:steep_des}
\bd_{\text{sd}}^\toptzero = -\nabla f(\bx^\toptzero).
\end{equation}
That is, the steepest descent method coincides with the gradient descent method.

\paragrapharrow{Stochastic gradient descent.} 
In many cases, the function $f(\bx)$ is defined over a datasets $\mathcalD=\{\bs_1, \bs_2, \ldots, \bs_D\}$ such that $f(\bx)$ and its gradient $\nabla f(\bx)$ can be expressed as 
\begin{equation}
f(\bx)\triangleq f(\mathcalD; \bx) = \frac{1}{D} \sum_{d=1}^{D} f(\bs_d; \bx)
\qquad \text{and}\qquad 
\nabla f(\bx) \triangleq \frac{1}{D}  \sum_{d=1}^{D}\nabla f(\bs_d; \bx),
\end{equation} 
respectively.
While if we follow the negative gradient of a single sample or a batch of samples iteratively, the local estimate of the direction can be obtained and is known as the \textit{stochastic gradient descent} (SGD) \citep{robbins1951stochastic}.
The SGD method can be categorized  into two types:
\begin{itemize}
\item \textbf{The strict SGD:} Computes the gradient using only one randomly selected data point per iteration: $\nabla f(\bx^\toptzero) \approx \nabla f(\bs_d; \bx^\toptzero)$.
\item \textbf{The mini-batch SGD:} A compromise between full gradient descent and strict SGD, where a small subset (mini-batch) of the dataset is used to compute an estimate of the gradient: $\nabla f(\bx^\toptzero) \approx \frac{1}{\abs{\sS}}  \sum_{d\in\sS}\nabla f(\bs_d; \bx^\toptzero)$.
\end{itemize}
The SGD method is particular useful when the number of \textit{training entries} (i.e., the data used for updating/training the model, while the data used for final evaluation is called the \textit{test entries or test data}) are substantial, as computing the full gradient can be computationally expensive or even resulting in that the gradients from different input samples may cancel out and the final update is small.
However, since the gradient is estimated using only a subset of the data, the updates can be noisy.
In the SGD framework, the objective function is stochastic, composed of a sum of subfunctions evaluated at different subsamples of the data (see Chapter~\ref{chapter:stochastic_opt}). However, a drawback of the vanilla update (both GD and SGD) lies in its susceptibility to getting trapped in local minima \citep{rutishauser1959theory}.

\paragrapharrow{Choice of stepsize.}
For a small stepsize, gradient descent ensures a monotonic improvement at every iteration, guaranteeing convergence, albeit to a local minimum. However, the speed of the vanilla gradient descent method is generally slow, and it can exhibit a linear rate in case of poor curvature conditions. 
While choosing a stepsize larger than an optimal threshold may cause divergence in terms of the objective function.
Determining an optimal learning rate (whether global or per-dimension) becomes more of an art than science for many problems. 
Previous work has attempted to alleviate the need for manually selecting a global learning rate \citep{zeiler2012adadelta}, though such methods remain sensitive to other hyperparameters. We will introduce \textit{(exact or inexact) line search strategies} to determine the stepsize more systematically in Section~\ref{section:line_search}.

\paragrapharrow{Descending Property.}
Most (if not all) optimization methods incorporate mechanisms to enforce the descending property:
\begin{equation}\label{equation:des_prob1}
f(\bx^\toptone) < f(\bx^\toptzero). 
\end{equation}
This prevents convergence to a maximizer and also makes it less probable that we get convergence to a saddle point (Definition~\ref{definition:stat_point}). 
If the objective function has several minimizers, the final solution depends on the starting point $\bx^{(1)}$. We do not know which of the minimizers that will be found; the specific minimizer found is not necessarily the one closest to $\bx^{(1)}$.

As mentioned previously, in many cases the method produces vectors which converge towards the minimizer in two clearly different stages: the ``global stage" where $\bx^{(1)}$ is far from the solution and we want the method to produce iterates which move steadily towards the optimizer $\bx^*$, and the ``final stage" where $\bx^\toptzero$ is close to $\bx^*$ and seek faster convergence.

The global convergence properties of a method describe its behavior when initialization occurs at a point $\bx^{(1)}$, which is not close to a (local) minimizer $\bx^{*}$. 
Ideally, the iterates should move steadily toward a neighborhood of  $\bx^{*}$. For instance, there are methods for which it is possible to prove that any accumulation point (i.e., limit of a subseries) of $\{\bx^\toptzero\}_{t>0}$ is a stationary point (Definition~\ref{definition:stat_point}), meaning the gradient vanishes:
$$
\nabla f(\bx^\toptzero) \rightarrow \bzero \qquad \text{ for } \qquad t \rightarrow \infty. 
$$
While this does not eliminate the possibility of convergence to a saddle point or maximizer, the descending property \eqref{equation:des_prob1} typically prevents such cases in practice.
In this ``global phase", our primary concern is ensuring that losses do not increase (except for possibly the initial steps).
To analyze convergence in terms of iterates rather than function values, a natural potential function is 
$$
\text{$e_t\triangleq\normtwobig{\be^\toptzero},\quad$ where $\be^\toptzero\triangleq \bx^\toptzero - \bx^*$.}
$$
Let $\{\be^\toptzero\}_{t>0}$ denote the error sequence. The requirement for progress is:
$$
\normtwobig{\be^\toptone} < \normtwobig{\be^\toptzero}\qquad  \text{ for }\qquad t > t'.
$$

In the final stages of the iteration where the $\bx^\toptzero$ are close to $\bx^{*}$, we expect faster convergence. 
Local convergence analysis describes how quickly the iterates approach $\bx^{*}$ to a desired accuracy. Some methods exhibit  linear convergence (Definition~\ref{definition:linear-convergence}):
$$
\normtwobig{\be^\toptone}\leq c_{1}\normtwobig{\be^\toptzero}, \quad \text{with } 0 < c_{1} < 1 \text{ and } \bx^\toptzero \text{ close to } \bx^{*}.  
$$
However, higher-order convergence is preferable. For instance, quadratic convergence (Definition~\ref{definition:quadratic-convergence}) satisfies:
$$
\normtwobig{\be^\toptone}\leq c_{2}\normtwobig{\be^\toptzero}^{2}, \quad \text{with } c_{2} > 0 \text{ and } \bx^\toptzero \text{ close to } \bx^{*}. 
$$
Few practical methods achieve quadratic convergence, but superlinear convergence (Definition~\ref{definition:superlinear_convergence}) is a common goal:
$$
\frac{\normtwobig{\be^\toptone}}{\normtwobig{\be^\toptzero}} \rightarrow 0 \quad\text{ for } t \rightarrow \infty.
$$
Superlinear convergence is faster than linear convergence, though typically not as rapid as quadratic convergence.

\paragrapharrow{Framework of a descent method.}

\begin{algorithm}[H] 
	\caption{Structure of  Descent Methods}
	\label{alg:struc_gd_gen}
	\begin{algorithmic}[1] 
		\Require A function $f(\bx)$; 
		\State {\bfseries Input:}  Initialize $\bx^{(1)}$;
		\For{$t=1,2,\ldots$}
		\State Find a descent direction $\bd^\toptzero$ such that $\innerproduct{\bd^\toptzero, \bg^\toptzero}<0$;
		\State Pick a stepsize $\eta_t$;
		\State $\bx^{(t+1)} \leftarrow \bx^\toptzero + \eta_t \bd^\toptzero$;
		\EndFor
		\State (Output Option 1) Output  $\bx_{\text{final}}\leftarrow \bx^{(T)}$;
		\State (Output Option 2) Output  $\bx_{\text{avg}}\leftarrow \frac{1}{T}(\sum_{t=1}^{t}\bx^\toptzero)$ or $\sum_{t=1}^{T} \frac{2t}{T(T+1)} \bx^\toptzero$;
		\State (Output Option 3) Output  $\bx_{\text{best}}\leftarrow \argmin_{t\in\{1,2,\ldots,T\}} f(\bx^\toptzero)$;
	\end{algorithmic} 
\end{algorithm}

The methods presented in this book are descent methods, meaning they satisfy the descending condition \eqref{equation:des_prob1} at each iteration. Each iteration consists of:
\begin{itemize}
\item Finding a descent direction $\bd^\toptzero$ at the $t$-th iteration.
\item Determining a stepsize $\eta_t$ giving a good decrease in the function value.
\end{itemize}
This sequence of operations forms the foundation of descent algorithms, see Algorithm~\ref{alg:struc_gd_gen}.
The search direction $ \bd^\toptzero $ at each iteration must be a descent direction (Definition~\ref{definition:uncons_des_direct}, Theorem~\ref{theorem:uncons_des_dir}). 
This ensures that we can reduce $ f(\bx) $ by choosing an appropriate walking distance, and thus we can satisfy the descending condition \eqref{equation:des_prob1}. 

\paragrapharrow{Stopping criteria.}
Ideally, a stopping criterion should indicate when the current error is sufficiently small:
$$
\textbf{(ST1)}:\qquad \normtwobig{\be^{\toptzero}} < \delta_{1}.
$$
Another ideal condition would be when the current function value is close enough to the minimum:
$$
\textbf{(ST2)}:\qquad f(\bx^\toptzero ) - f(\bx^{*}) < \delta_{2}.
$$
Both conditions reflect the convergence $ \bx^\toptzero  \rightarrow \bx^{*} $. 
However, they are impractical because $ \bx^{*} $ and $ f(\bx^{*}) $ are (in most cases) unknown~\footnote{In some cases, the value $f(\bx^*)$ is known. For example, in the convex feasibility problem, we seek feasible points within a convex set, where $f(\bx^*)$ is zero.}. 
Instead, we rely on approximations:
\begin{equation}\label{equation:des_stopcri1}
\textbf{(ST3)}:\qquad \normtwobig{\bx^\toptone - \bx^\toptzero} < \varepsilon_{1} \qquad \text{or} \qquad f(\bx^\toptzero ) - f(\bx^\toptone) < \varepsilon_{2}.
\end{equation}
We must emphasize that even if \eqref{equation:des_stopcri1} is fulfilled with small $ \varepsilon_{1} $ and $ \varepsilon_{2} $, it does not guarantee that  $ \normtwobig{\be^{\toptzero}} $ or $ f(\bx^\toptzero ) - f(\bx^{*}) $ are small.

Another form of convergence, mentioned earlier in this chapter, is $ \nabla f (\bx^\toptzero ) \rightarrow \bzero $ for $ t \rightarrow \infty $ (Theorem~\ref{theorem:fermat_fist_opt}). This leads to another commonly used stopping criterion:
\begin{equation}\label{equation:des_stopcri2}
\textbf{(ST4)}:\qquad \normtwobig{\nabla f(\bx^\toptzero )} < \varepsilon_{3},
\end{equation}
which is included in many implementations of descent methods.

Another useful approach involves leveraging the property of converging function values. 
The quadratic approximation (Theorem~\ref{theorem:quad_app_theo}) of $ f $ at $ \bx^{*} $ is
$$
f(\bx^\toptzero ) \approx f(\bx^{*}) + (\bx^\toptzero  - \bx^{*})^{\top} \nabla f(\bx^{*}) + \frac{1}{2}(\bx^\toptzero  - \bx^{*})^{\top} \nabla^2 f(\bx^{*})(\bx^\toptzero - \bx^{*}).
$$
Since $ \bx^{*} $ is a local minimizer, we have $ \nabla f(\bx^{*}) = \bzero $ and $ \bH^{*} \triangleq \nabla^2 f(\bx^{*}) $ is positive semidefinite (Theorem~\ref{theorem:second_nec_loca}). This simplifies to:
$
f(\bx^\toptzero ) - f(\bx^{*}) \approx \frac{1}{2}(\bx^\toptzero  - \bx^{*})^{\top} \bH^{*}(\bx^\toptzero  - \bx^{*}).
$
Thus, another stopping criterion can be defined as:
$$
\textbf{(ST5)}:\qquad \frac{1}{2}(\bx^\toptone - \bx^\toptzero )^{\top} \bH^\toptzero(\bx^\toptone - \bx^\toptzero ) < \varepsilon_{4} \quad \text{with} \quad \bx^\toptzero  \approx \bx^{*}. 
$$
Here, $ \bx^\toptzero  - \bx^{*} $ is approximated by $ \bx^\toptone - \bx^\toptzero  $ and $ \bH^{*} $ is approximated by
$
\bH^\toptzero \triangleq \nabla^2 f(\bx^\toptzero ).
$

In the following sections, we delve into a detailed exploration of the gradient descent method, examining its variations and adaptations from different perspectives. This comprehensive analysis aims to provide a deeper understanding of the algorithm, its formulations, challenges, and practical applications.


\index{Greedy search}
\index{Non-Euclidean gradient descent}
\index{Dual norm}
\section{Gradient Descent by Greedy Search and Variants}\label{section:als-gradie-descent-taylor}

We now consider  the \textit{greedy search} method such that $\bx^\toptone    \leftarrow \mathop{\arg \min}_{\bx^\toptzero} f(\bx^\toptzero)$ under some mild conditions. 
The linear approximation theorem (Theorem~\ref{theorem:linear_approx}) shows that 
\begin{equation}\label{equation:gd_gree_tay}
f(\bx^\toptzero + \eta \bd) = f(\bx^\toptzero)+ \eta \bd^\top \nabla  f(\bx^\toptzero )+ \mathcalO(\normtwo{\eta \bd}^2).
\end{equation}
For small values of  $\eta$, the term $\mathcalO(\normtwo{\eta \bd}^2)$ becomes negligible compared to the middle term. Therefore, we can approximate $f(\bx^\toptzero  + \eta \bd)$ as
\begin{equation}\label{equation:gd_gree_approx}
f(\bx^\toptzero + \eta \bd) \approx f(\bx^\toptzero ) + \eta \bd^\top \nabla  f(\bx^\toptzero ),
\end{equation}
when $\eta$ is sufficiently small. 
The second term on the right-hand side, $ \bd^\top \nabla  f(\bx^\toptzero ) $, is the \textit{directional derivative} of $ f $ at $ \bx^\toptzero $ in the direction $ \bd $. 
To reiterate, it indicates the approximate change in $ f $ for a small step $ \bd $. The step $ \bd $ is a descent direction if the directional derivative is negative.

To address how to choose $\bd$ to make the directional derivative as negative as possible, note that since directional derivative $ \bd^\top \nabla  f(\bx^\toptzero ) $ is linear in $ \bd $, it can be made arbitrarily negative by increasing $ \bd $ (provided $ \bd $ is a descent direction, i.e., $ \bd^\top \nabla  f(\bx^\toptzero ) < 0 $). To make this question meaningful, we must limit the size of $ \bd $, or normalize by its length.

Let $ \norm{\cdot} $ be any norm on $ \real^n $. We define a \textit{normalized greedy descent direction} (with respect to the norm $ \norm{\cdot} $) as
\begin{equation}\label{equation:norma_greedes}
\bd_{\text{ngd}}^\toptzero \in \argmin_{\bd} \left\{ \bd^\top \nabla  f(\bx^\toptzero ) \text{ s.t. } \norm{\bd} = 1 \right\}.
\end{equation}
(Note that there may be multiple minimizers.) A normalized greedy descent direction $ \bd_{\text{ngd}}^\toptzero $ is a step of unit norm that provides the largest decrease in the linear approximation of $ f $.
By the definition of the dual norm \eqref{equation:dual_norm_equa}, it follows that $\innerproduct{\bd_{\text{ngd}}^\toptzero, \nabla  f(\bx^\toptzero )} =-\norm{\nabla  f(\bx^\toptzero )}_{*}$, where $\norm{\cdot}_{*}$ denotes the dual norm.

Since the problem in \eqref{equation:norma_greedes} can  equivalently be stated using the constraint $\norm{\bd}\leq 1$, $\bd_{\text{ngd}}^\toptzero$ also lies within the set of primal counterparts of $\nabla  f(\bx^\toptzero ) $, whose existence is shown in Definition~\ref{definition:set_primal}.

It is also convenient to consider an unnormalized greedy descent step $\bd_{\text{ugd}}^\toptzero$ by scaling the normalized greedy descent direction in a particular way:
\begin{equation}\label{equation:unnorma_greedes}
\bd_{\text{ugd}}^\toptzero \triangleq \norm{\nabla f(\bx^\toptzero)}_{*} \bd_{\text{ngd}}^\toptzero.
\end{equation}
The reason for this particular unnormalization is that it aligns with the negative gradient (the steepest descent direction) when the underlying norm is the $\ell_2$ norm.
 Note that for the greedy descent step, we have
$$
\nabla f(\bx^\toptzero)^\top \bd_{\text{ugd}}^\toptzero = \norm{\nabla f(\bx^\toptzero)}_{*} \nabla f(\bx^\toptzero)^\top \bd_{\text{ngd}}^\toptzero = -\norm{\nabla f(\bx^\toptzero)}_{*}^2.
$$
When exact line search is used (see Section~\ref{section:line_search}),  scale factors in the descent direction do not affect the outcome, so either the normalized or unnormalized direction can be used.

\paragrapharrow{Greedy search for $\ell_2$ norm.}
If we take the norm $\norm{\cdot}$ to be the $\ell_2$ norm, we find that the greedy descent direction in \eqref{equation:norma_greedes} is simply the negative gradient, i.e., 
\begin{equation}\label{equation:greedy_l2_ngdugd}
\bd_{\text{ngd}}^\toptzero = - \frac{\nabla f(\bx^\toptzero )}{\normtwobig{\nabla f(\bx^\toptzero )}}
\qquad \text{and}\qquad 
\bd_{\text{ugd}}^\toptzero = -\nabla f(\bx^\toptzero ).
\end{equation}
The greedy descent method for the $\ell_2$ norm coincides with the gradient descent method (or the steepest descent method).
The above equality also shows that the unnormalized greedy search direction  corresponds to the negative gradient direction or the steepest descent direction.

\paragrapharrow{Greedy search for $\ell_1$ norm.}
As another example,  consider the greedy descent method for the $\ell_1$ norm. A normalized greedy descent direction can be characterized as
$$
\bd_{\text{ngd}}^\toptzero \in \argmin_{\bd} \left\{ \bd^\top \nabla  f(\bx^\toptzero ) \text{ s.t. } \normone{\bd} \leq 1 \right\}.
$$
We use `$\in$' since the solution of the problem may not be unique.
Let $i$ be any index for which $\norminf{\nabla f(\bx^\toptzero)} = \abs{(\nabla f(\bx^\toptzero))_i}$. 
By the definition of the $\ell_\infty$ norm and the dual norm (\eqref{equation:l_p_norm} and \eqref{equation:dual_norm_equa}), a normalized greedy descent direction $\bd_{\text{ngd}}^\toptzero$ for the $\ell_1$ norm is given by
$$
\bd_{\text{ngd}}^\toptzero = -\sign\left(\frac{\partial f}{\partial x_i}(\bx^\toptzero)\right) \be_i,
$$
where $\be_i$ is the $i$-th  unit basis vector (see Example~\ref{example:set_primal_count}). An unnormalized greedy descent step is then
$$
\bd_{\text{ugd}}^\toptzero = \bd_{\text{ngd}}^\toptzero \norminf{\nabla f(\bx^\toptzero)} = -\frac{\partial f}{\partial x_i} (\bx^\toptzero) \be_i.
$$
which is a descent direction since $\innerproduct{\bd_{\text{ugd}}^\toptzero, \nabla f(\bx^\toptzero)}<0$ (assuming $\nabla f(\bx^\toptzero) \neq\bzero$).
Thus, the normalized greedy descent step in the $\ell_1$ norm can always be chosen to be a (positive or negative) standard basis vector, representing the coordinate axis direction along which the approximate decrease in $f$ is greatest.
Note that the index for which $\norminf{\nabla f(\bx^\toptzero)} = \abs{(\nabla f(\bx^\toptzero))_i}$ may not be unique (see Example~\ref{example:set_primal_count}). In such cases, a convex combination of these descent directions can be used as the final descent direction.

The greedy descent algorithm in the $\ell_1$ norm has a  natural interpretation: At each iteration we select a component of $\nabla f(\bx^\toptzero)$ with maximum absolute value (though the component may not be unique), and then decrease or increase the corresponding component of $\bx^\toptzero$, according to the sign of $(\nabla f(\bx^\toptzero))_i$. The algorithm is sometimes called a \textit{coordinate-descent algorithm} because only one component of the variable $\bx$ is updated at each iteration, potentially simplifying or even trivializing the line search.

\paragrapharrow{Greedy search for $\bQ$-norm.}
We further consider the $\bQ$-norm in \eqref{equation:q_norm}: 
$
\norm{\bx}_{\bQ} = (\bx^\top \bQ \bx)^{1/2} = \normtwo{\bQ^{1/2} \bx}
$
for any $\bx\in\real^n$,
where $\bQ$ is positive definite. The normalized greedy descent direction is given by
$$
\begin{aligned}
\bd_{\text{ngd}}^\toptzero 
&=\argmin_{\bd} \left\{ \bd^\top \nabla  f(\bx^\toptzero ) \text{ s.t. } \norm{\bd}_{\bQ} \leq 1 \right\}\\
&=-\normtwo{\bQ^{-1/2} \nabla f(\bx^\toptzero)}^{-1/2}\bQ^{-1} \nabla f(\bx^\toptzero). 
\end{aligned}
$$
This can be solved using the KKT conditions or the definition of the dual norm.
The dual norm is given by $\norm{\bx}_{*} = \normtwo{\bQ^{-1/2} \bx} =\norm{\bx}_{\bQ^{-1}}$ for any $\bx\in\real^n$, so the greedy descent step with respect to $\norm{\cdot}_{\bQ}$ is given by
\begin{equation}\label{equation:qnorm_ugd}
	\bd_{\text{ugd}}^\toptzero = -\bQ^{-1} \nabla f(\bx^\toptzero).
\end{equation}
This is a descent direction by Theorem~\ref{theorem:uncons_des_dir}.

\paragrapharrow{Change of variables in $\bQ$-norm.}

An interesting alternative interpretation of the greedy descent direction $\bd_{\text{ugd}}^\toptzero $ is as the gradient search direction after applying a change of coordinates to the problem. Let $\widetildebx \triangleq \bQ^{1/2} \bx$; thus, $\norm{\bx}_{\bQ} = \normtwo{\widetildebx}$. Using this change of coordinates, we can solve the original problem of minimizing $f$ by solving the equivalent problem of minimizing the function $\widetildef : \real^n \rightarrow \real$, given by
$$
\widetildef(\widetildebx) \triangleq f(\bQ^{-1/2} \widetildebx) = f(\bx).
$$
If we apply the gradient method to $\widetildef$, the search direction at a point $\widetildebx^\toptzero$ (which corresponds to the point $\bx^\toptzero = \bQ^{-1/2} \widetildebx^\toptzero$ for the original problem) is
$$
\widetildebd^\toptzero = -\nabla \widetildef(\widetildebx^\toptzero) = -\bQ^{-1/2} \nabla f(\bQ^{-1/2} \widetildebx^\toptzero) = -\bQ^{-1/2} \nabla f(\bx^\toptzero).
$$
Since $\widetildebx = \bQ^{1/2} \bx$ by definition, the search direction in the original space is obtained by mapping $\widetildebd^\toptzero$ back using $\bQ^{-1/2}$: 
$$
\bd^\toptzero = \bQ^{-1/2} \widetildebd^\toptzero = -\bQ^{-1} \nabla f(\bx^\toptzero)
$$
which corresponds to the unnormalized greedy search direction in \eqref{equation:qnorm_ugd}. In other words, the greedy descent method in the $\bQ$-norm $\norm{\cdot}_{\bQ}$ can be thought of as the gradient method applied to the problem after the change of variables $\widetildebx^\toptzero = \bQ^{1/2} \bx^\toptzero$ for each iteration $t$.

\section{Geometrical Interpretation of Gradient Descent} 
\begin{figure}[h]
\centering   
\vspace{-0.25cm}  
\subfigtopskip=2pt  
\subfigbottomskip=2pt  
\subfigcapskip=-5pt  
\subfigure[A two-dimensional convex function $f(\bx)$.]{\label{fig:alsgd1}
\includegraphics[width=0.47\linewidth]{./imgs/quadratic_singular.pdf}}
\subfigure[$f(\bx)=c$ is a constant.]{\label{fig:alsgd2}
\includegraphics[width=0.44\linewidth]{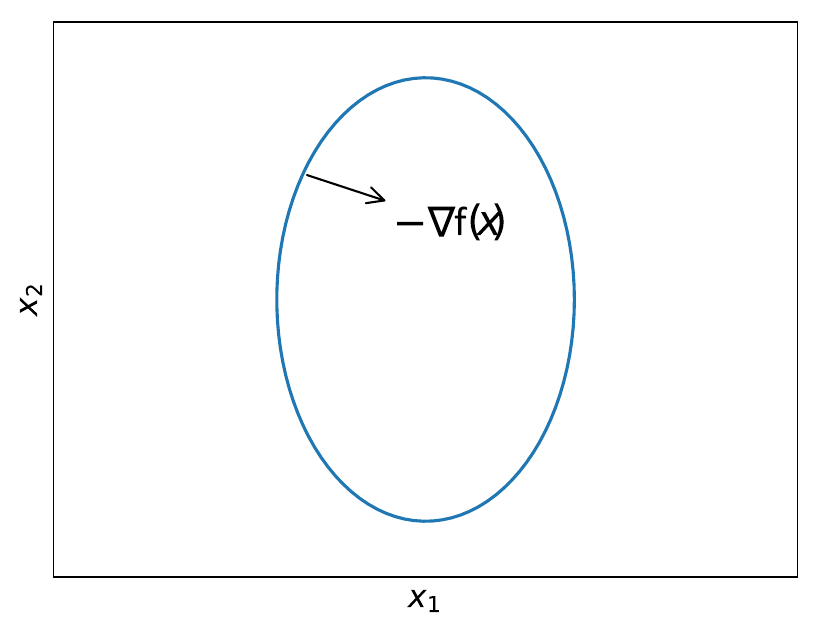}}
\caption{Figure~\ref{fig:alsgd1} shows a convex function surface plot and its contour plot (\textcolor{mylightbluetext}{blue}=low, \textcolor{mydarkyellow}{yellow}=high), where the upper graph represents  the surface plot, and the lower one is its projection (i.e., contour). Figure~\ref{fig:alsgd2}: $-\nabla f(\bx)$ directs the reduction in loss for the convex function $f(\bx)$.}
\label{fig:alsgd-geometrical}
\end{figure}

\begin{lemma}[Direction of Gradients]\label{lemm:direction-gradients}
An important fact is that gradients are orthogonal to level curves (also known as level surfaces).
\end{lemma}
\begin{proof}[of Lemma~\ref{lemm:direction-gradients}: Informal]
To prove this, we need to show that the gradient is orthogonal to the tangent of the level curve. Let's start with the two-dimensional case. Suppose the level curve has the form $f(x,y)=c$. 
This equation implicitly defines a relationship  between $x$ and $y$, such that $y=y(x)$, where $y$ can be considered as a function of $x$. Therefore, the level curve can be expressed as:
$$
f(x, y(x)) = c.
$$
Applying the chain rule gives us:
$$
\frac{\partial f}{\partial x} \underbrace{\frac{dx}{dx}}_{=1} + \frac{\partial f}{\partial y} \frac{dy}{dx}=0.
$$
This implies that the gradient is perpendicular to the tangent vector:
$$
\left\langle \frac{\partial f}{\partial x}, \frac{\partial f}{\partial y}\right\rangle
\cdot 
\left\langle \frac{dx}{dx}, \frac{dy}{dx}\right\rangle=0.
$$
Now, let's generalize this to higher dimensions. Consider a level set defined by a vector $\bx\in \real^n$: $f(\bx) = f(x_1, x_2, \ldots, x_n)=c$. Each variable $x_i$ can be regarded as a function of a parameter $t$ along the level set $f(\bx)=c$: $f(x_1(t), x_2(t), \ldots, x_n(t))=c$. Differentiating both sides with respect to $t$ using the chain rule yields:
$$
\frac{\partial f}{\partial x_1} \frac{dx_1}{dt} + \frac{\partial f}{\partial x_2} \frac{dx_2}{dt}
+\ldots + \frac{\partial f}{\partial x_n} \frac{dx_n}{dt}
=0.
$$
Therefore, the gradients is perpendicular to the tangent in the $n$-dimensional case:
$$
\left\langle \frac{\partial f}{\partial x_1}, \frac{\partial f}{\partial x_2}, \ldots, \frac{\partial f}{\partial x_n}\right\rangle
\cdot 
\left\langle \frac{dx_1}{dt}, \frac{dx_2}{dt}, \ldots, \frac{dx_n}{dt}\right\rangle=0.
$$
This completes the proof.
\end{proof}
The lemma above provides a profound geometric interpretation of gradient descent. In the process of minimizing a (convex) function $f(\bx)$, gradient descent strategically moves in the direction opposite to the gradient, which reduces the loss. Figure~\ref{fig:alsgd-geometrical} illustrates a two-dimensional scenario where  $-\nabla f(\bx)$ guides the decrease in loss for a (convex) function $f(\bx)$. 

\index{Regularization}
\section{Geometrical Interpretation of Regularization}\label{section:geom_int_regu}
\begin{figure}[h]
\centering
\includegraphics[width=0.95\textwidth]{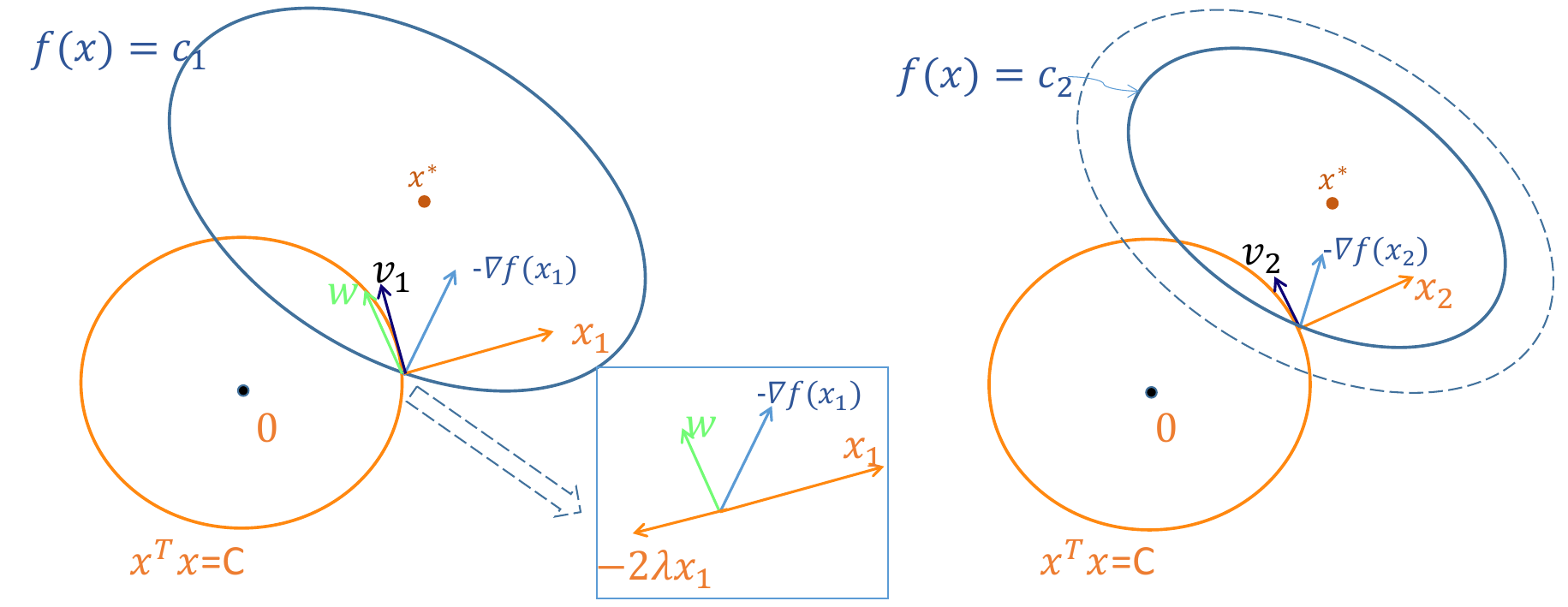}
\caption{Constrained gradient descent with $\bx^\top\bx\leq C$. The \textcolor{mydarkgreen}{green} vector $\bw$ is the projection of $\bv_1$ into $\bx^\top\bx\leq C$, where $\bv_1$ is the component of $-\nabla f(\bx)$ perpendicular to $\bx_1$. The right picture shows the next step after the update in the left picture. $\bx^*$ denotes the optimal solution of \{$\min f(\bx)$\}.}
\label{fig:alsgd3}
\end{figure}
\textit{Regularization} is a machine learning technique employed to prevent overfitting and improve model generalization. Overfitting occurs when a model is overly complex and fits the training data too closely, resulting in poor performance on  unseen data. 
To mitigate this issue, regularization introduces a constraint or a penalty term into the loss function used for model optimization, discouraging the development of overly complex models. 
This creates  a trade-off between having a simple, generalizable model and fitting the training data well. 
Common types of regularization include $\ell_1$ regularization, $\ell_2$ regularization (Tikhonov regularization), and elastic net regularization (a combination of $\ell_1$ and $\ell_2$ regularizations). 
Regularization finds extensive applications in machine learning algorithms such as linear regression, logistic regression, and neural networks.

Gradient descent also reveals the geometric significance of regularization. To avoid confusion, we denote the loss function without regularization by $f(\bz)$ and the loss with the $\ell_2$ regularization by $F(\bx) \triangleq f(\bx)+\lambda \normtwo{\bx}^2$, where $f(\bx): \real^n \rightarrow \real$ (this notation is exclusive to this subsection). When minimizing $f(\bx)$, the descent method searches for a solution in $\real^n$. 
However, in machine learning, an exhaustive search across the entire space may lead to overfitting. A partial remedy involves searching within a subset of the vector space, such as searching in $\bx^\top\bx < C$ for some constant $C$. That is,
$$
\argmin_{\bx} \, \big\{f(\bx) \gap  \text{s.t.} \gap \bx^\top\bx\leq C\big\}.
$$
This constrained search helps prevent overfitting by introducing regularization through the addition of a penalty term in the optimization process.
In the previous discussion, a basic gradient descent approach proceeds in the direction of $-\nabla f(\bx)$,  updating $\bx$ by $\bx^+\leftarrow \bx-\eta \nabla f(\bx)$ for a small stepsize $\eta$. 
When the level curve is $f(\bx)=c_1$ and the descent approach is situated at $\bx=\bx_1$, where $\bx_1$ is the intersection of $\bx^\top\bx=C$ and $f(\bx)=c_1$, the descent direction $-\nabla f(\bx_1)$ will be perpendicular to the level curve of $f(\bx_1)=c_1$, as shown in the left picture of Figure~\ref{fig:alsgd3}. 
However, if we further restrict that the optimal value can only be in the subspace $\bx^\top\bx\leq C$, the trivial descent direction $-\nabla f(\bx_1)$ will lead $\bx_2=\bx_1-\eta \nabla f(\bx_1)$ outside of $\bx^\top\bx\leq C$. 

To address this,  the step $-\nabla f(\bx_1)$ is decomposed  into 
$$
-\nabla f(\bx_1) = a\bx_1 + \bv_1,
$$ 
where $a\bx_1$ is the component perpendicular to the curve of $\bx^\top\bx=C$, and $\bv_1$ is the component parallel to the curve of $\bx^\top\bx=C$. Keeping only the step $\bv_1$, then the update 
$$
\bx_2 = \text{project}(\bx_1+\eta \bv_1) = \text{project}\left(\bx_1 + \eta
\underbrace{(-\nabla f(\bx_1) -a\bx_1)}_{\bv_1}\right)\footnote{where the project($\bx$) operator
will project the vector $\bx$ to the closest point inside $\bx^\top\bx\leq C$. Notice here the direct update $\bx_2 = \bx_1+\eta \bv_1$ can still make $\bx_2$ outside the curve of $\bx^\top\bx\leq C$.}
$$ 
will lead to a smaller loss from $f(\bx_1)$ to $f(\bx_2)$ while still satisfying the constraint $\bx^\top\bx\leq C$. 
This technique is known as the \textit{projected gradient descent} (Section~\ref{section:pgd}). It is not hard to see that the update $\bx_2 = \text{project}(\bx_1+\eta \bv_1)$ is equivalent to finding a vector $\bw$ (depicted by the \textcolor{mydarkgreen}{green} vector in the left panel of Figure~\ref{fig:alsgd3}) such that $\bx_2=\bx_1+\bw$ lies inside the curve of $\bx^\top\bx\leq C$. Mathematically, the vector $\bw$ can be obtained as $-\nabla f(\bx_1) -2\lambda \bx_1$ for some $\lambda$, as shown in the middle panel of Figure~\ref{fig:alsgd3}. This aligns with the negative gradient of $F(\bx)=f(\bx)+\lambda\normtwo{\bx}^2$ such that 
$$
-\nabla F(\bx_1) = -\nabla f(\bx_1) - 2\lambda \bx_1,
$$
and 
$$
\begin{aligned}
\bw &= -\nabla F(\bx_1) 
\qquad \implies \qquad 
\bx_2 = \bx_1+ \bw =\bx_1 -  \nabla F(\bx_1).
\end{aligned}
$$
In practice, a small stepsize $\eta$ can be applied to prevent crossing  the curve boundary of $\bx^\top\bx\leq C$:
$$
\bx_2  =\bx_1 -  \eta\nabla F(\bx_1).
$$

\index{Line search}
\index{Steepest descent}
\index{Exact line search method} 
\index{Inexact line search method} 
\index{Soft line search method}
\index{Greedy search}
\section{Descent Methods with Line Search}\label{section:line_search}

For the unconstrained optimization problem (P1) in Definition~\ref{definition:opt_probs_all}, we will use the process of finding the minimum value of  $f(\bx)$ as an analogy to going downhill. Imagine a person standing at some point $\bx$, where $f(\bx)$ represents the height at that location. To find the lowest point, two things need to be determined at point $\bx$: first, which direction to take for the next step; second, how far to walk along that direction before choosing the next downhill direction. Once these factors are determined, the process can be repeated until reaching the minimum value of $f(\bx)$.

The second phase is then solved using line search strategies.
The mathematical description of the line search method is as follows: Given the current iteration point $\bx^\toptzero$, first determine the direction vector $\bd^\toptzero$ using some method, then determine the positive scalar $\eta_t$. The next iteration point can be written as:
$ \bx^\toptone \leftarrow \bx^\toptzero + \eta_t \bd^\toptzero.  $
We call $\bd^\toptzero$ the \textit{search direction} or \textit{descent direction} at the iteration point $\bx^\toptzero$, and $\eta_t$ the corresponding stepsize. 
As mention in the previous sections, we require that $\bd^\toptzero$ is a descent direction, i.e., $\innerproduct{\bd^\toptzero, \nabla f(\bx^\toptzero)} < 0$ (Definition~\ref{definition:uncons_des_direct}, assume $\nabla f(\bx^\toptzero)\neq \bzero$). This descent property ensures that the function value decreases along this direction. The key to the line search method is determining a good direction $\bd^\toptzero \in \real^n$ and an appropriate stepsize $\eta_t$.

Gradient descent, steepest descent, greedy search methods are some examples of how to select the search direction $\bd^\toptzero$.
In this section, we will focus on selecting $\eta_t$.  
While methods for selecting $\bd^\toptzero$ vary widely, techniques for choosing $\eta_t$ are similar across different algorithms. First, construct the auxiliary function
$ \phi(\eta) = f(\bx^\toptzero + \eta \bd^\toptzero), $
where $\bd^\toptzero$ is the given descent direction, and $\eta > 0$ is the independent variable of this auxiliary function. The geometric meaning of $\phi(\eta)$ is straightforward: it is the restriction of the objective function $f(\bx)$ to the ray $\{\bx^\toptzero + \eta \bd^\toptzero \mid \eta > 0\}$. Note that $\phi(\eta)$ is a univariate function, making it relatively easy to study.

The goal of line search is to select an appropriate $\eta_t$ such that $\phi(\eta_t)$ is minimized as much as possible. However, this task is still challenging: $\eta_t$ should ensure sufficient decrease in $f$, while also minimizing computational effort. 
Balancing these aspects is crucial. A natural approach is to find $\eta_t$ such that
$$ 
\eta_t = \arg\min_{\eta > 0} \phi(\eta), 
$$
i.e., $\eta_t$ is the optimal stepsize. This line search algorithm is called the \textit{exact line search method}. Although the exact line search method often leads to solutions under most circumstances, it typically requires significant computation to determine $\eta_t$, making it impractical for many applications.
The \textit{inexact line search method} or \textit{soft line search method}, on the other hand, does not require $\eta_t$ to be the minimum point of $\phi(\eta)$, but only requires $\phi(\eta_t)$ to satisfy certain inequality conditions. Due to its simpler structure, the inexact line search method is more commonly used in practice. We will now introduce the details of these algorithms.

\subsection{Exact Line Search Conditions}\label{section:exact_line_search}
In the preceding sections, we derived the gradient descent, where the update step at step $t$ is $ -\eta_t\bg^\toptzero\triangleq-\eta_t \nabla  f(\bx^\toptzero)$, and the learning rate $\eta_t$ controls how large of a step to take in the direction of negative gradient. 
Exact line search is a method that directly determines the optimal learning rate to achieve the most significant improvement in the objective function. Formally, at the $t$-th step of gradient descent, exact line search solves the following problem:
$$
\eta_t = \underset{\eta}{\arg\min}\,\,   f(\bx^\toptzero - \eta \bg^\toptzero).
$$
After performing the gradient update $\bx^\toptone \leftarrow \bx^\toptzero - \eta_t \bg^\toptzero$, the gradient is computed at $\bx^\toptone$ for the next step $t+1$. More generally, let $\bd^\toptzero$ be the descent direction; then, the descent method with line search can be described by:
$$
\eta_t = \underset{\eta}{\arg\min}\,\,   f(\bx^\toptzero + \eta \bd^\toptzero).
$$
In line search methods, the loss function at iteration $t$ can be  expressed in terms of $\eta$ as
$
\phi(\eta) \triangleq f(\bx^\toptzero + \eta \bd^\toptzero).
$
Consequently, the problem can be formulated as finding
$$
\eta_t = \underset{\eta}{\arg\min} \,\, f(\bx^\toptzero + \eta \bd^\toptzero)=\underset{\eta}{\arg\min} \,\,  \phi(\eta).
$$
A necessary condition on $\eta_t$ is $\phi^\prime(\eta_t)=0$, as stated by Theorem~\ref{theorem:fermat_fist_opt}.
Since the derivative of $\phi$ is  $\phi^\prime(\eta) = \innerproduct{\nabla f(\bx^\toptzero + \eta \bd^\toptzero), \bd^\toptzero}$, this shows that either $\nabla f(\bx^\toptzero + \eta_t \bd^\toptzero) = \bzero$, , indicating that we have found a stationary point for $f$, or if $\nabla f(\bx^\toptzero + \eta_t \bd^\toptzero) \neq \bzero$, then  $\phi^\prime(\eta_t)=0$ indicates 
$$
\innerproduct{\nabla f(\bx^\toptone),  \bd^\toptzero} = 0.
$$
The two cases mean that the exact line search will stop at a point where the local gradient is orthogonal to the search direction. We formalize this finding in the following lemma, which will be crucial in the development of conjugate gradient methods (Section~\ref{section:conjugate-descent}).
\begin{lemma}[Orthogonality in Exact Line Search]\label{lemm:linear-search-orghonal}
The gradient of optimal point $\bx^\toptone=\bx^\toptzero + \eta_t \bd^\toptzero $ of an exact line search is orthogonal to the current update direction $\bd^\toptzero$:
$$
\innerproduct{\nabla f(\bx^\toptone), \bd^\toptzero} = 0.
$$
\end{lemma}

When $\eta=0$, we have
\begin{equation}\label{equation:linesearc-eta0}
	\phi^\prime (0) = \innerproduct{\bd^\toptzero,  \bg^\toptzero} \leq 0.
\end{equation}
A crucial property in typical line search settings is that the loss function
$\phi(\eta)$, when expressed in terms of $\eta$, is often a unimodal function. 
If a value $\eta_{\max}$ is identified such that $\phi^\prime(\eta_{\max}) > 0$, the optimal learning rate is then in the range of $[0, \eta_{\max}]$. Line search methods are used to find the optimal $\eta_t$ within this range, satisfying the optimal condition $\phi^\prime(\eta_t)=0$.
Following this rule, we then introduce conditions for soft line search algorithms and we will introduce some prominent line search approaches: polynomial interpolation line search, bisection line search, golden-section line search, and Armijo condition search in the following sections.

\subsection{Soft Line Search Conditions}

In the early days of optimization, exact line search was dominant. Now, soft line search is used more frequently, and new methods rarely require exact line search.
One advantage of soft line search over exact line search is its speed. If the initial guess for the step length is a rough approximation to the minimizer in the given direction, the line search can terminate immediately if certain mild criteria are met. While the result of exact line search is typically a good approximation, it often requires additional function evaluations, which can make descent methods with exact line search less efficient overall. Despite this, exact line search can sometimes find local minima in fewer iterations compared to soft line search.

At the start of an iteration with a descent method, where $\bx^\toptzero$ is far
from the solution $\bx^*$, the imprecision of soft line search results is less critical. This further supports the preference for soft line search in practice. 

In the soft line search method, the selection of $\eta_t$ must meet certain requirements known as line search conditions. The appropriateness of these conditions directly affects the convergence of the algorithm. Choosing inappropriate line search conditions can lead to failure in convergence. We illustrate this with an example.

\begin{example}[Divergence of Descent Methods]
Consider a one-dimensional unconstrained optimization problem
$$
\min_{x} \quad f(x) = x^2,
$$
with the initial point $ x^{(1)} = 1 $. Since the problem is one-dimensional, the descent direction only has two possibilities: $\{-1, 1\}$. We choose $ d^\toptzero = -\sign(x^\toptzero) $ and require that the stepsize satisfies the condition of monotonic decrease of the function value at the iteration point, i.e., $ f(x^\toptzero + \eta_t d^\toptzero) < f(x^\toptzero) $. Consider the following two stepsizes:
$$
\eta_{t,1} = \frac{1}{3^{t+1}}
\qquad \text{and}\qquad
\eta_{t,2} = 1 + \frac{2}{3^{t+1}},
$$
Performing some algebra yields the sequences:
$$
x_1^\toptzero = \frac{1}{2} \left(1 + \frac{1}{3^t}\right)
\qquad \text{and}\qquad
\qquad x_2^\toptzero = \frac{(-1)^t}{2} \left(1 + \frac{1}{3^t}\right).
$$
Clearly, both sequences $\{f(x_1^\toptzero)\}$ and $\{f(x_2^\toptzero)\}$ are  monotonically decreasing. However, the sequence $\{x_1^\toptzero\}$ converges to a point that is not a minimum point, while the sequence $\{x_2^\toptzero\}$ oscillates around the origin and does not converge.
\end{example}

The reason for the divergence in the above example is insufficient decrease in the function value $ f(x^\toptzero) $ during the iteration process, preventing the algorithm from converging to the minimum point. To avoid this issue, more reasonable soft line search conditions must be introduced to ensure convergence.

\subsection*{Armijo Condition and Backtracking Methods}

We first introduce the \textit{Armijo condition}, which is a commonly used soft line search condition.
The purpose of the Armijo condition is to ensure that each iteration sufficiently decreases the function value.
\begin{definition}[Armijo Condition/Sufficient Decrease Condition \citep{armijo1966minimization}]\label{definition:arminjo_condition}
Let $f:\real^n\rightarrow\real$ be a differentiable function, and let $ \bd^\toptzero $ denote  a  descent direction at the point $ \bx^\toptzero $. 
Then the stepsize $\eta=\eta_t$ satisfies the \textit{Armijo condition (a.k.a., the sufficient decrease condition)} if
\begin{equation}\label{equation:arminjo_condition}
f(\bx^\toptzero + \eta \bd^\toptzero) \leq f(\bx^\toptzero) + c_1 \eta \innerproduct{\nabla f(\bx^\toptzero),  \bd^\toptzero}, 
\end{equation}
where $c_1 \in (0, 1)$ is a fixed constant.
The existence of such a constant $c_1$ is guaranteed by Lemma~\ref{lemma:valid_suff_des_armi}.
\end{definition}

\begin{figure}[h!]
\centering  
\vspace{-0.35cm} 
\subfigtopskip=2pt 
\subfigbottomskip=2pt 
\subfigcapskip=-5pt 
\subfigure[Armijo condition.]{\label{fig:descent_armijo}
\includegraphics[width=0.48\linewidth]{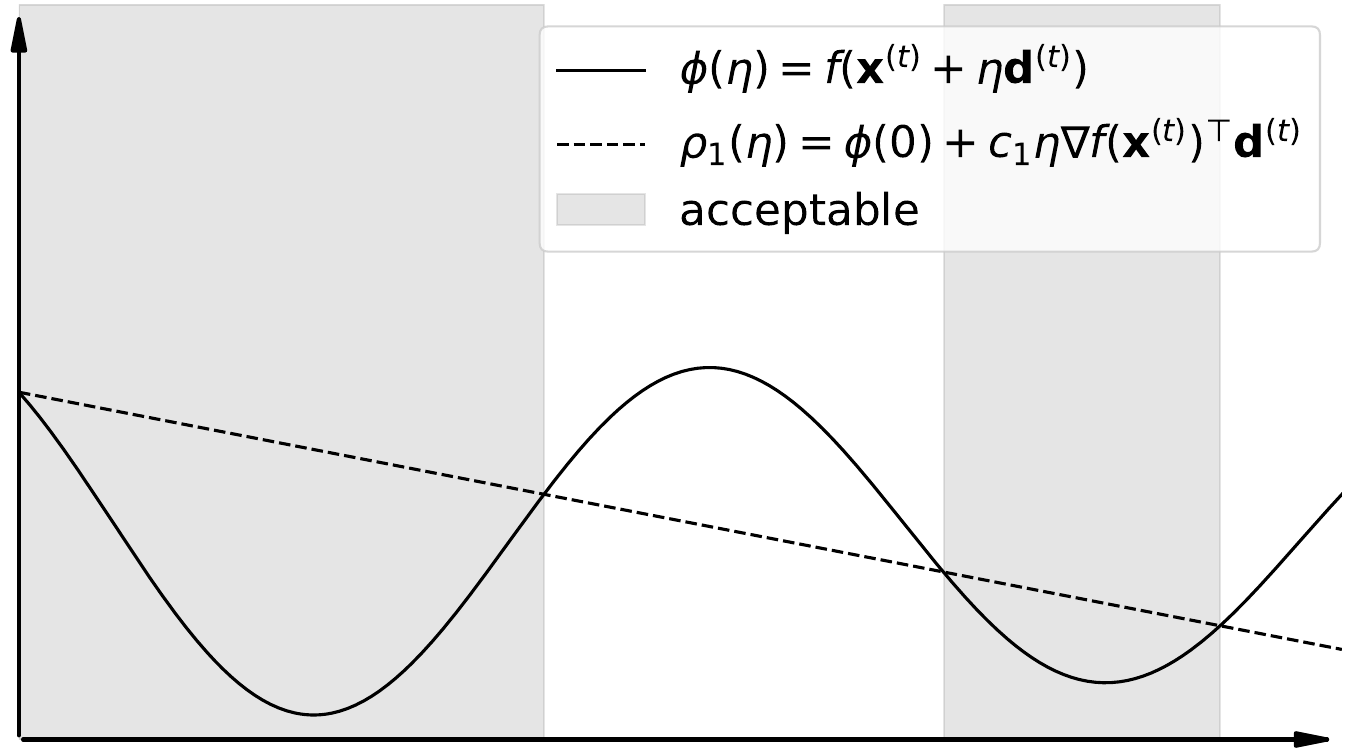}}
\subfigure[Goldstein condition. Same as Figure~\ref{fig:descent_goldstein2}.]{\label{fig:descent_goldstein}
\includegraphics[width=0.48\linewidth]{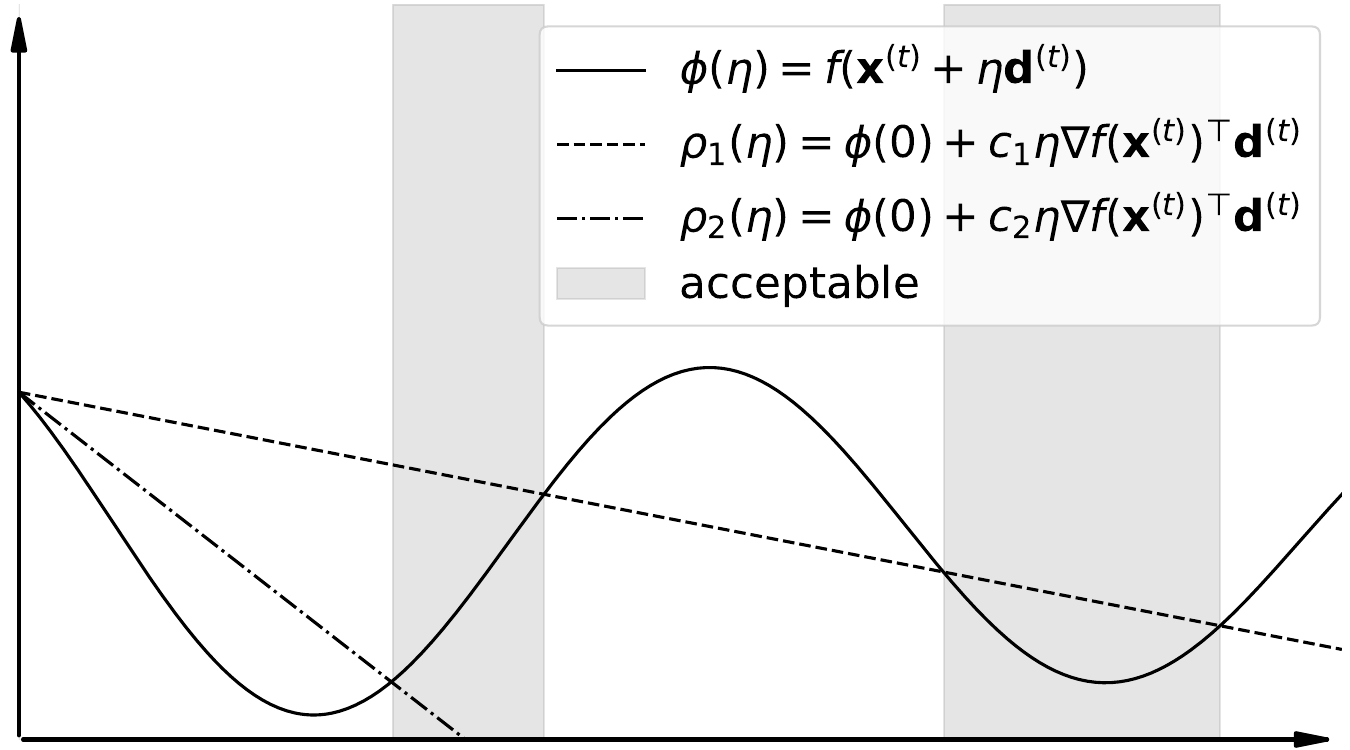}}
\caption{The Armijo condition ensures that any point $(\eta, \phi(\eta))$ must lie below the line $\rho_1(\eta)$, while the Goldstein condition ensures that any point must lie between the lines $\rho_1(\eta)$ and $\rho_2(\eta)$.}
\label{fig:descent_armijo_goldstein}
\end{figure}

The Armijo condition \eqref{equation:arminjo_condition} has a very intuitive geometric meaning, indicating that the point $(\eta, \phi(\eta))$ must lie below the line
$$
\rho_1(\eta) \triangleq \phi(0) + c_1 \eta \innerproduct{\nabla f(\bx^\toptzero),  \bd^\toptzero}, 
\quad
\text{where }\phi(\eta) \triangleq f(\bx^\toptzero + \eta \bd^\toptzero).
$$
As shown in Figure~\ref{fig:descent_armijo}, all points in the shaded area satisfy the Armijo condition. Note that $\bd^\toptzero$ is the descent direction, which means the slope of $\rho_1(\eta)$ is negative. Choosing $\eta$ that satisfies condition \eqref{equation:arminjo_condition} indeed ensures the decrease of the function value. In practical applications, the parameter $c_1$ is usually chosen to be a very small positive number, such as $c_1 = 10^{-4}$, making the Armijo condition  easy to satisfy. However, using only the Armijo condition alone cannot guarantee the convergence of the iteration. This is because $\eta = 0$ trivially satisfies condition \eqref{equation:arminjo_condition}, which means the sequence of iterations does not change. Studying such a stepsize is meaningless. Therefore, the Armijo condition needs to be used in conjunction with other conditions.

In the implementation of optimization algorithms, finding a stepsize that satisfies the Armijo condition is relatively easy. A commonly used algorithm is the \textit{backtracking method}. Given an initial guess $\widehat{\eta}$, the backtracking method continuously reduces the trial stepsize by a factor until it finds the first point that satisfies the Armijo condition  \eqref{equation:arminjo_condition}. Specifically, the backtracking method selects
$$
\eta_t = \gamma^{j} \widehat{\eta},
$$
where
$
j = \min \left\{ i = 0, 1, \cdots \mid f(\bx^\toptzero + \gamma^i \widehat{\eta} \bd^\toptzero) \leq f(\bx^\toptzero) + c_1 \gamma^i \widehat{\eta} \innerproduct{\nabla f(\bx^\toptzero),  \bd^\toptzero} \right\},
$
and $\gamma \in (0, 1)$ is a given constant. The basic process of the backtracking method is shown in Algorithm~\ref{alg:gd_line_search}.

\begin{algorithm}[h] 
\caption{Descent Method with Backtracking Line Search at Iteration $t$}
\label{alg:gd_line_search}
\begin{algorithmic}[1] 
\Require A function $f(\bx)$ and $t$-th iteration $\bx^\toptzero$; 
\State {\bfseries Input:}  Choose an initial stepsize $\widehat{\eta}$, parameters $\gamma, c_1 \in (0, 1)$;
\State {\bfseries Input:}  Initialize $\eta \leftarrow \widehat{\eta}$;
\While{$f(\bx^\toptzero + \eta \bd^\toptzero) > f(\bx^\toptzero) + c_1 \eta \innerproduct{\nabla f(\bx^\toptzero),  \bd^\toptzero}$}
\State Let $\eta \leftarrow \gamma \eta$;
\EndWhile
\State {\bfseries Return:}  $\eta_t \leftarrow \eta$;
\end{algorithmic} 
\end{algorithm}

This algorithm is called the backtracking method because the trial values of $\eta$ are reduced from large to small, ensuring that the output $\eta_t$ is as large as possible while still satisfying the Armijo condition. Additionally, Algorithm~\ref{alg:gd_line_search} will not run indefinitely because $\bd^\toptzero$ is a descent direction.

\subsection*{Goldstein Condition}

When $\eta$ is sufficiently small, the Armijo condition always holds. In practical applications,
to overcome the shortcomings of the Armijo condition, we usually set a lower bound for $\eta$ to prevent the stepsize from being too small. 
Since the Armijo condition only requires that the point $(\eta, \phi(\eta))$ must lie below a certain line, we can also require that this point lies above another line. This is the   \textit{Goldstein condition}.

\begin{definition}[Goldstein Condition \citep{goldstein1965steepest}]\label{definition:goldstein_cond}
Let $f:\real^n\rightarrow\real$ be a differentiable function, and let $ \bd^\toptzero $ denote  a  descent direction at the point $ \bx^\toptzero $. 
Then the stepsize $\eta=\eta_t$ satisfies the \textit{Goldstein condition} if
\begin{subequations}\label{equation:goldsteins}
\begin{align}
f(\bx^\toptzero + \eta \bd^\toptzero) &\leq f(\bx^\toptzero) + c_1 \eta \innerproduct{\nabla f(\bx^\toptzero),  \bd^\toptzero}, \label{equation:goldstein1} \\
f(\bx^\toptzero + \eta \bd^\toptzero) &\geq f(\bx^\toptzero) + (1 - c_1) \eta \innerproduct{\nabla f(\bx^\toptzero),  \bd^\toptzero},\label{equation:goldstein2}
\end{align}
\end{subequations}
where $c_1 \in \left(0, \frac{1}{2}\right)$  is a fixed constant. Let 
$
\phi(\eta) \triangleq f(\bx^\toptzero + \eta \bd^\toptzero)$. 
Then the Goldstein condition can  be equivalently stated as
\begin{subequations}\label{equation:goldsteins_all2}
\begin{align}
\phi(\eta) &\leq \phi(0) + c_1 \eta \phi^\prime(0)\triangleq\rho_1(\eta), \label{equation:goldsteinv2_1} \\
\phi(\eta) &\geq \phi(0) + (1 - c_1) \eta \phi^\prime(0)\triangleq\rho_2(\eta).\label{equation:goldsteinv2_2}
	\end{align}
\end{subequations}
\end{definition}

Similarly, the Goldstein condition \eqref{equation:goldsteins} also has an  intuitive geometric interpretation, indicating that the point $(\eta, \phi(\eta))$ must lie between two lines:
$$
\begin{aligned}
\rho_1(\eta) &\triangleq \phi(0) + c_1 \eta \innerproduct{\nabla f(\bx^\toptzero),  \bd^\toptzero},\\
\rho_2(\eta) &\triangleq \phi(0) + (1 - c_1) \eta \innerproduct{\nabla f(\bx^\toptzero),  \bd^\toptzero}.
\end{aligned}
$$
As shown in Figure~\ref{fig:descent_goldstein}, all points in the shaded area satisfy the Goldstein condition. We also note that the Goldstein condition effectively prevents excessively small stepsizes $\eta$.
Note that the restriction $c_1 < \frac{1}{2}$ is necessary so that $\rho_2$ lies below $\rho_1$. 

\begin{figure}[h!]
\centering  
\vspace{-0.35cm} 
\subfigtopskip=2pt 
\subfigbottomskip=2pt 
\subfigcapskip=-5pt 
\subfigure[Goldstein condition, where $c_1<1$ and $c_2=1-c_1$; same as Figure~\ref{fig:descent_goldstein}.]{\label{fig:descent_goldstein2}
	\includegraphics[width=0.48\linewidth]{./imgs/descent_goldstein.pdf}}
\subfigure[Wolfe condition.]{\label{fig:descent_wolfe}
	\includegraphics[width=0.48\linewidth]{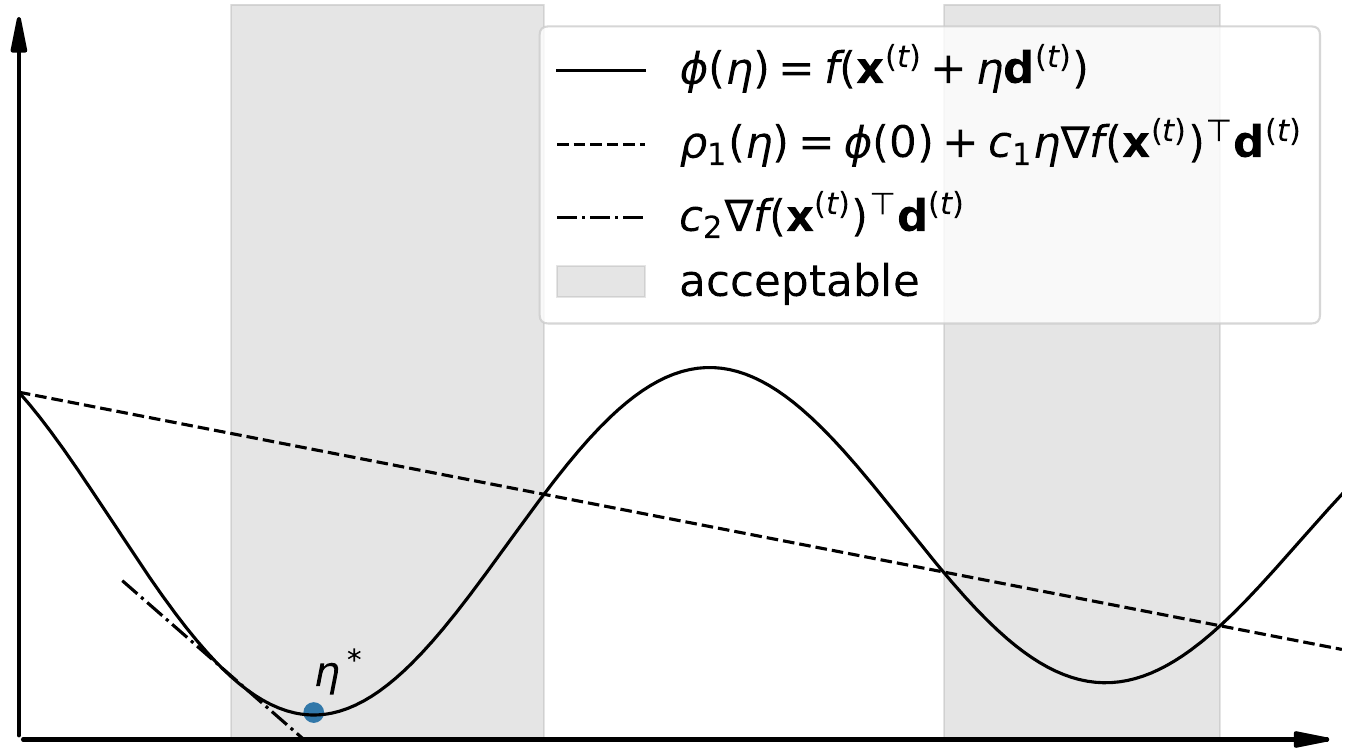}}
\caption{The Goldstein condition ensures that any point must lie between the lines $\rho_1(\eta)$ and $\rho_2(\eta)$, while the Wolfe condition ensures that any point must lie below the line $\rho_1(\eta)$ and the slope must be greater than $c_2$ times the initial slope $\phi^\prime(0)$.}
\label{fig:descent_goldstein_wolfe}
\end{figure}

\subsection*{Wolfe Condition}

The Goldstein condition ensures that the function value decreases sufficiently, but it may exclude the optimal function value. As shown in Figure~\ref{fig:descent_goldstein}, the minimum point of the one-dimensional function $\phi(\eta)$ does not necessarily lie within the shaded area that satisfies the Goldstein condition. Therefore, we introduce the \textit{Wolfe condition}.
\begin{definition}[Wolfe Condition]\label{definition:wolfe_cond}
Let $f:\real^n\rightarrow\real$ be a differentiable function, and let $ \bd^\toptzero $ denote  a  descent direction at the point $ \bx^\toptzero $. 
Then the stepsize $\eta=\eta_t$ satisfies the \textit{Wolfe condition} if
\begin{subequations}\label{equation:wolfe_all}
\begin{align}
f(\bx^\toptzero + \eta \bd^\toptzero) &\leq f(\bx^\toptzero) + c_1 \eta \innerproduct{\nabla f(\bx^\toptzero),  \bd^\toptzero}, \label{equation:wolfe_1} \\
\innerproduct{\nabla f(\bx^\toptzero + \eta \bd^\toptzero),  \bd^\toptzero} &\geq c_2 \innerproduct{\nabla f(\bx^\toptzero),  \bd^\toptzero}, \label{equation:wolfe_2}
\end{align}
\end{subequations}
where $0<c_1 < c_2<1$ are given constants.
Alternatively, let 
$\phi(\eta) \triangleq f(\bx^\toptzero + \eta \bd^\toptzero)$.
Then the Wolfe condition \eqref{equation:wolfe_all} can be equivalently stated as 
\begin{subequations}\label{equation:wolfe_var}
\begin{align}
\phi(\eta) &\leq \phi(0)+ c_1 \eta \phi^\prime(0)  \triangleq\rho_1(\eta), \label{equation:wolfevar_1} \\
\phi^\prime(\eta) &\geq c_2 \phi^\prime(0). \label{equation:wolfevar_2}
\end{align}
\end{subequations}
\end{definition}

In the Wolfe condition,
\begin{itemize}
\item The first inequality \eqref{equation:wolfe_1} is the Armijo condition, which requires that $\eta_t$ must be sufficiently small to provide a useful decrease in the objective function.
\item The second inequality \eqref{equation:wolfe_2} is the essential requirement of the Wolfe condition. Note that $\innerproduct{\nabla f(\bx^\toptzero + \eta \bd^\toptzero),  \bd^\toptzero}$ is  the derivative of $\phi(\eta)$. The Wolfe condition requires that the slope of the tangent line to $\phi(\eta)$ at $\eta=\eta_t$ cannot be less than $c_2$ times the slope at $\eta = 0$, i.e., $\eta_t$ also needs to be sufficiently large so that we have moved away from the starting tangent of the curve $y=\phi(\eta)$ for $\eta\geq 0$. 
\end{itemize}
As shown in Figure~\ref{fig:descent_wolfe}, all points in the interval shaded area satisfy the Wolfe condition. Note that at the minimum point $\eta^*$ of $\phi(\eta)$, $\phi^\prime(\eta^*) = \innerproduct{\nabla f(\bx^\toptzero + \eta^* \bd^\toptzero),  \bd^\toptzero} = 0$, so $\eta^*$ always satisfies condition \eqref{equation:wolfe_2}. Choosing a smaller $c_1$ can make $\eta^*$ also satisfy condition \eqref{equation:wolfe_1}, i.e., the Wolfe condition includes the exact solution of the line search problem in most cases. In practical applications, the parameter $c_2$ is usually set to 0.9.

\subsection*{Strong Wolfe Condition}
In addition, a stepsize may satisfy the Wolfe condition without being close to a minimizer of $\phi(\eta)$ (see Figure~\ref{fig:descent_wolfe}). 
The \textit{strong Wolfe condition} modifies the curvature condition to force $\eta_t$ to lie in at least a broad neighborhood of a local minimizer or stationary point of $\phi$. 
\begin{definition}[Strong Wolfe Condition]\label{definition:strong_wolfe_cond}
Let $f:\real^n\rightarrow\real$ be a differentiable function, and let $ \bd^\toptzero $ denote  a  descent direction at the point $ \bx^\toptzero $. 
Then the stepsize $\eta=\eta_t$ satisfies the \textit{strong Wolfe condition} if
\begin{subequations}\label{equation:strong_wolfe_all}
\begin{align}
f(\bx^\toptzero + \eta \bd^\toptzero) &\leq f(\bx^\toptzero) + c_1 \eta \innerproduct{\nabla f(\bx^\toptzero),  \bd^\toptzero}, \label{equation:strong_wolfe_1} \\
\abs{\innerproduct{\nabla f(\bx^\toptzero + \eta \bd^\toptzero),  \bd^\toptzero}} &\leq c_2 \abs{\innerproduct{\nabla f(\bx^\toptzero),  \bd^\toptzero}}, \label{equation:strong_wolfe_2}
\end{align}
\end{subequations}
where $0<c_1 < c_2<1$ are given constants.
Alternatively, let 
$\phi(\eta) \triangleq f(\bx^\toptzero + \eta \bd^\toptzero)$.
Then the strong Wolfe condition \eqref{equation:strong_wolfe_all} can be equivalently stated as 
\begin{subequations}\label{equation:strong_wolfe_var}
\begin{align}
\phi(\eta) &\leq \phi(0)+ c_1 \eta \phi^\prime(0)  \triangleq\rho_1(\eta), \label{equation:strong_wolfevar_1} \\
\abs{\phi^\prime(\eta)} &\geq c_2 \abs{\phi^\prime(0)}. \label{equation:strong_wolfevar_2}
\end{align}
\end{subequations}
\end{definition}

\subsection*{Nonmonotone Line Search Conditions}

The three conditions introduced above have a common feature: they generate a monotonic iteration sequence. In practical applications, nonmonotone algorithms sometimes perform better. This requires the use of nonmonotone line search conditions.

\begin{definition}[Grippo Condition \citep{grippo1986nonmonotone}]\label{definition:grippo_condition}
Let $f:\real^n\rightarrow\real$ be a differentiable function, and let $ \bd^\toptzero $ denote  a  descent direction at the point $ \bx^\toptzero $. 
Given a positive integer $ k > 0 $, then the stepsize $\eta=\eta_t$ satisfies the \textit{Grippo condition} if
\begin{equation}
f(\bx^\toptzero + \eta \bd^\toptzero) \leq \max_{0 \leq j \leq \min(t, k)} f(\bx^{(t-j)}) + c_1 \eta \innerproduct{\nabla f(\bx^\toptzero),  \bd^\toptzero},
\end{equation}
where $ c_1 \in (0, 1) $ and $k\geq 0$ are fixed constants. When $k=0$, this reduces to the Armijo condition.
\end{definition}
The Grippo condition is similar to the Armijo condition, but with a key difference. The Armijo condition requires the function value $ f(\bx^\toptone) $ of the next iteration to decrease sufficiently compared to the current iteration's function value $ f(\bx^\toptzero) $. 
In contrast, the Grippo condition only requires the function value of the next step to be lower than the maximum of the previous $ k $ steps' function values. This makes the Grippo condition less restrictive than the Armijo condition and does not require the sequence  $ \{f(\bx^\toptzero)\}_{t>0} $ to be monotonic.

\subsection{Line Search Algorithm}

This section introduces line search algorithms used in practical applications. Previous discussions have introduced the backtracking method (Algorithm~\ref{alg:gd_line_search}) and noted that it can be used to find a stepsize that satisfies the Armijo condition \eqref{equation:arminjo_condition}.
In fact, by simply modifying the termination condition of the algorithm, the backtracking method can be applied to other line search conditions, such as the  nonmonotone line search conditions in Definitions~\ref{definition:grippo_condition}. The implementation of the backtracking method is simple and straightforward, making it one of the most commonly used line search algorithms. 
However, it has some drawbacks: first, it cannot guarantee finding a stepsize that satisfies the Wolfe condition \eqref{equation:wolfe_2}, which is essential for some optimization algorithms; 
second, the backtracking method reduces the stepsize exponentially, making it sensitive to the initial guess $ \widehat{\eta} $ and parameters $ \gamma $ and $ c_1 $ in Algorithm~\ref{alg:gd_line_search}.

The selection of parameter $\gamma$ is relatively sensitive. 
If $\gamma$ is too large, each reduction in step length during the trial is very small, resulting in low backtracking efficiency. 
Conversely, if $\gamma$ is too small, the backtracking is overly aggressive, leading to a final step length that is too small, potentially missing opportunities to select a larger stepsize. 
Below is a brief introduction to other types of line search algorithms. See also \citet{madsen2010and, aggarwal2020linear,  lu2022gradient} for a reference.

\subsection*{Polynomial Interpolation Line Search}
To improve the efficiency of the backtracking method, we use a \textit{polynomial interpolation-based line search} algorithm. Assuming an initial step length $\widehat{\eta}_1$ is given, if this step length does not satisfy the Armijo condition, we reduce the trial step length using polynomial interpolation rather than reducing $\widehat{\eta}_1$ by a constant multiple. 
Specifically, we construct a quadratic interpolation function $q(\eta)$ using $\phi(0)$, $\phi^\prime(0)$, $\phi(\widehat{\eta}_1)$ to, i.e., we find a quadratic function $q(\eta)$ satisfying:
$$
q(0) = \phi(0), \qquad q^\prime(0) = \phi^\prime(0), \qquad q(\widehat{\eta}_1) = \phi(\widehat{\eta}_1).
$$
Since a quadratic function has only three parameters, these three conditions uniquely determine $q(\eta)$. 
It is easy to verify that the minimum value point of $q(\eta)$ lies within $(0, \widehat{\eta}_1)$. 
We then take this minimum value point, $\widehat{\eta}_2$,  as the next trial point and repeat the process until the Armijo condition is satisfied.

\begin{figure}[ht]
\centering
\begin{minipage}{1\textwidth}
\begin{algorithm}[H]
\caption{Soft Line Search with Wolfe Condition  at Iteration $t$}
\label{alg:soft_line_search_wolfe}
\begin{algorithmic}[1]
\Require A function $f(\bx)$ and $t$-th iteration $\bx^\toptzero$; 
\State {\bfseries Input:}  Choose the maximum of step length $\eta_{\max}$ and the maximum number of trials $k_{\max}$;
\If{$\phi^\prime(0) = 0$} \Comment{(i)}
\State $\eta \leftarrow 0$;
\Else
\State $k \leftarrow 0; \quad \gamma \leftarrow c_2  \phi^\prime(0);$
\State $[a, b] \leftarrow \big[0, \min\{1, \eta_{\max}\}\big]$;  \Comment{(ii)}
\State $//$ $b$ is sufficiently small but not large enough $\implies$ increase the range;
\While{$(\phi(b) \leq \rho_1(b))$ and $(\phi^\prime(b) \leq \gamma)$ and $(b < \eta_{\max})$ and $(k < k_{\max})$} \Comment{(iii)}
\State $k \leftarrow k + 1$; 
\State $[a, b] \leftarrow \big[b, \min\{2b, \eta_{\max}\}\big]$;
\EndWhile
\State $\eta \leftarrow b$;
\State $//$ $b$ is too large, and $a$ is either 0 or too small $\implies$ $[a,b]$ contains acceptable points;
\While{$((\phi(\eta) > \rho_1(\eta))$ or $(\phi^\prime(\eta) < \gamma))$ and $(k < k_{\max})$}
\State $k \leftarrow k + 1$
\State Refine $\eta$ and $[a, b]$ by Algorithm~\ref{alg:soft_line_search_wolfe_refine}; \Comment{(iv)}
\EndWhile
\If{$\phi(\eta) \geq \phi(0)$} \Comment{(v)}
\State $\eta \leftarrow 0$;
\EndIf
\EndIf
\State {\bfseries Return:}  Final $\eta$;
\end{algorithmic}
\end{algorithm}
\end{minipage}\hfill
\begin{minipage}{1\textwidth}
\begin{algorithm}[H]
\caption{Find the Stepsize by Interpolation}
\label{alg:soft_line_search_wolfe_refine}
\begin{algorithmic}[1]
\Require A function $f(\bx)$, the interval $[a,b]$ where $a$ is small enough and $b$ is large enough; 
\State $D \leftarrow b - a; \quad c \leftarrow \left( \phi(b) - \phi(a) - D \cdot \phi^\prime(a) \right) / D^2$ 
\If{$c > 0$}  \Comment{Minimization of the polynomial $q(\zeta)$ in \eqref{equation:wolfe_poly}}
\State $\eta \leftarrow a - \phi^\prime(a) / (2c)$
\State $\eta \leftarrow \min\left\{\max\{\eta, a + 0.1D\}, b - 0.1D\right\}$ \Comment{Make  $\eta$  in the middle 80\% of $[a,b]$}
\Else  \Comment{Take $\eta$ as the midpoint of $[a, b]$, i.e., bisection}
\State $\eta \leftarrow (a + b) / 2$
\If{$\phi(\eta) < \rho_1(\eta)$} \Comment{$\eta$ is sufficiently small}
\State $[a,b] \leftarrow [a, \eta]$;
\Else
\State $[a,b]\leftarrow [\eta, b]$;
\EndIf
\EndIf
\State {\bfseries Return:}  Final candidate $\eta$ and refined interval $[a,b]$;
\end{algorithmic}
\end{algorithm}
\end{minipage}
\end{figure}

The interpolation-based line search algorithm effectively reduces the number of trials but still cannot guarantee that the stepsize satisfies the Wolfe condition. To address this, Algorithm~\ref{alg:soft_line_search_wolfe} can be used to find the stepsize iteratively \citep{fletcher2000practical, frandsen1999unconstrained, madsen2010and}.
Algorithm~\ref{alg:soft_line_search_wolfe} maintains an interval $[a,b]$ until $b$ satisfies the upper bound of the Wolfe condition by \eqref{equation:wolfevar_1}:
\begin{enumerate}[(i)]
\item The step shows $\bx^\toptzero$ is a stationary point ($\nabla f(\bx^\toptzero) = \bzero \Rightarrow \phi^\prime(0) = 0$) or $\bd^\toptzero$ is not downhill, then we do nothing.

\item The initial choice $b=1$ is used because in many optimization methods (e.g., Newton's method in Chapter~\ref{chapter:second_order}) $\eta=1$ is a very good guess in the final steps of the iteration. The upper bound $\eta_{\max}$ must be provided by the user to prevent  infinite loop if $f$ is unbounded.

\item The condition shows that $b$ is sufficiently small, but not large enough. So we increase the range from $[a, b]$ to $\big[b, \min\{2b, \eta_{\max}\}\big]$

\item The interval $[a,b]$ contains acceptable points, and use Algorithm~\ref{alg:soft_line_search_wolfe_refine} to find the optimal stepsize $\eta$ or decrease the range of the interval.

\item The algorithm may have stopped abnormally, e.g., by exceeding the permitted number $k_{\max}$ of function evaluations. If the current value of $\eta$ does not decrease the objective function, then we return $\eta = 0$.
\end{enumerate}
Sub-Algorithm~\ref{alg:soft_line_search_wolfe_refine} receives an interval $[a,b]$ where $a$ is the lower bound and $b$ is the upper bound.
We first construct a  second-order polynomial
\begin{equation}\label{equation:wolfe_poly}
q(\zeta) = \phi(a) + \phi^\prime(a) \cdot (\zeta - a) + c \cdot (\zeta - a)^2,
\end{equation}
where  $D \triangleq b - a$ and $c \triangleq \left( \phi(b) - \phi(a) - D \cdot \phi^\prime(a) \right) / D^2$ 
satisfying 
$$
q(a) = \phi(a), \quad q^\prime(a) = \phi^\prime(a), \quad \text{and}\quad q(b) = \phi(b).
$$
If $c > 0$, then the polynomial $q(\cdot)$ has a minimum, and we let $\eta$ be the minimizer. Otherwise we take $\eta$ as the midpoint of $[a, b]$.
In this case, if $\phi(\eta)$ is sufficiently small, then the right-hand part of $[a, b]$ contains points that satisfy both  constraints of the Wolfe condition \eqref{equation:wolfe_var}. Otherwise, $[\eta, b]$ is sure to contain acceptable points.

\index{Bisection line search}
\subsection*{Bisection Line Search}
As mentioned previously, the loss function $\phi(\eta)$, when expressed in terms of $\eta$, is often a unimodal function. 
In the \textit{bisection line search} method, we start by setting the interval $[a,b]$ as $[\eta_{\min}, \eta_{\max}]$, where $\eta_{\min}$ and $ \eta_{\max}$ serve as the lower and upper bounds, respectively, for the learning rate $\eta$ ($\eta_{\min}$ can be set to 0 as specified by \eqref{equation:linesearc-eta0}). The bisection line search involves evaluating the loss function $\phi(\eta)$ at the midpoint $\frac{a+b}{2}$. Given that $\phi^\prime(a)<0$ and $\phi^\prime(b)>0$, the bisection line search follows that
$$
\left\{
\begin{aligned}
\text{set } a &\triangleq \frac{a+b}{2} \text{, \gap if $\phi^\prime\left(\frac{a+b}{2}\right)<0$}; \\
\text{set } b &\triangleq \frac{a+b}{2} \text{, \gap if $\phi^\prime\left(\frac{a+b}{2} \right)>0$}.
\end{aligned}
\right.
$$
The procedure is repeated until the interval between $a$ and $b$ becomes sufficiently small.

The bisection line search is also known as the \textit{binary line search}. And in some cases, the derivative of $\phi(\eta)$ cannot be easily obtained; then the interval is narrowed by evaluating the objective function at two closely spaced
points around $\frac{a+b}{2} $. To be more concrete, assume $\phi(\eta)$ is convex (since we are in the descent setting), we evaluate the loss function at $\frac{a+b}{2}$ and $\frac{a+b}{2}+\epsilon$, where $\epsilon$ is a  small numerical value, e.g.,  $\epsilon=1e-8$.  
This allows us to evaluate whether the function is increasing or decreasing at $\frac{a+b}{2}$ by determining which of the two evaluations is larger. If the function is increasing
at $\frac{a+b}{2}$, the interval is narrowed to $[a,\frac{a+b}{2}+\epsilon]$; otherwise, it is narrowed to
$[\frac{a+b}{2}, b]$.
$$
\left\{
\begin{aligned}
\text{set } b &= \frac{a+b}{2}+\epsilon,  &\text{ \gap if increasing at $\frac{a+b}{2}$};\\
\text{set } a &= \frac{a+b}{2}, &\text{ \gap otherwise}. \\
\end{aligned}
\right.
$$
This iterative process continues until the range is sufficiently small or the required level of accuracy is achieved in the interval.

\index{Golden-section line search}
\subsection*{Golden-Section Line Search}
\begin{figure}[h]
\centering  
\vspace{-0.35cm} 
\subfigtopskip=2pt 
\subfigbottomskip=2pt 
\subfigcapskip=-5pt 
\subfigure[$\eta=a$ yields the minimum.]{\label{fig:golden_1}
\includegraphics[width=0.23\linewidth]{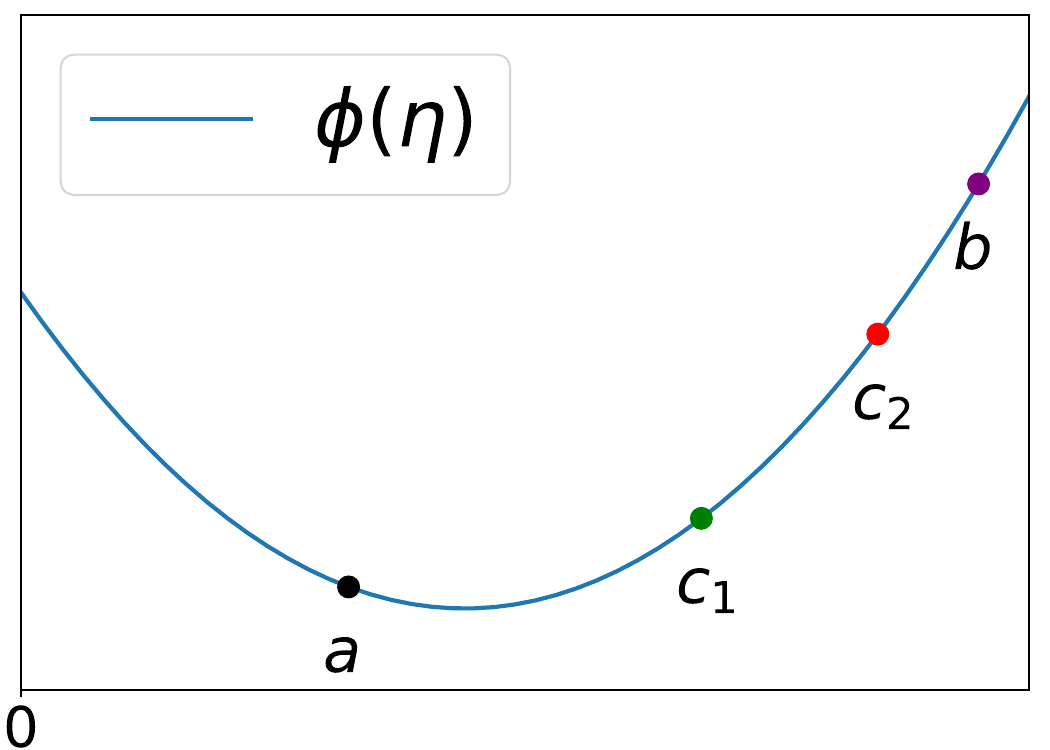}}
\subfigure[$\eta=c_1$ yields the minimum.]{\label{fig:golden_2}
\includegraphics[width=0.23\linewidth]{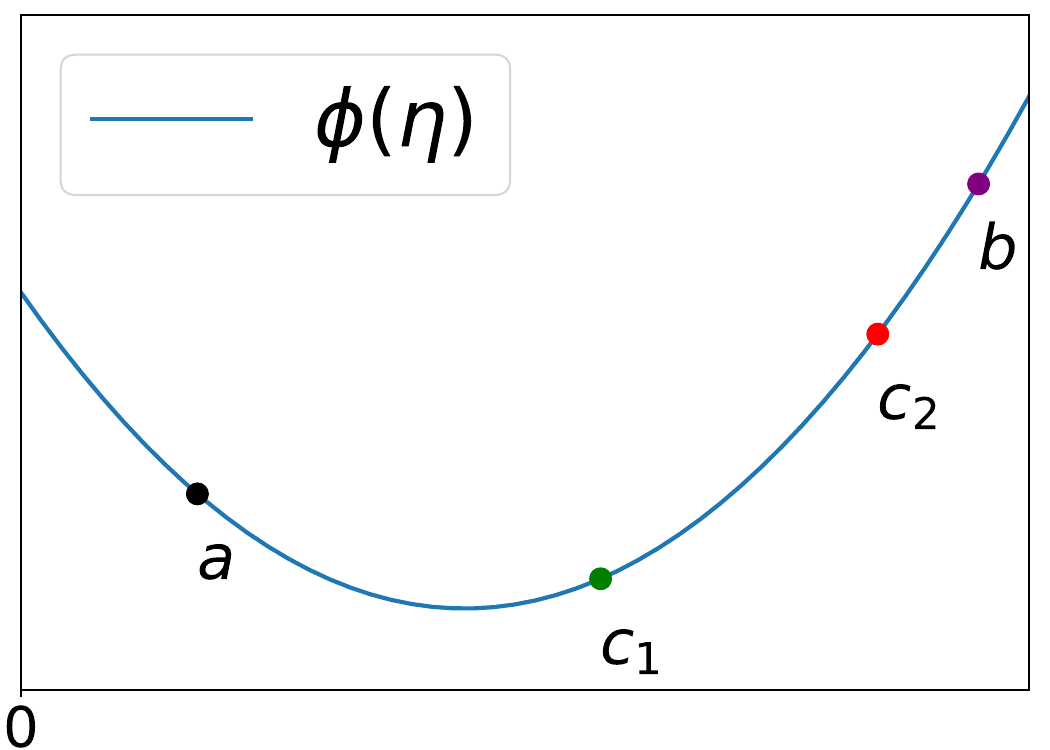}}
\subfigure[$\eta=c_2$ yields the minimum.]{\label{fig:golden_3}
\includegraphics[width=0.23\linewidth]{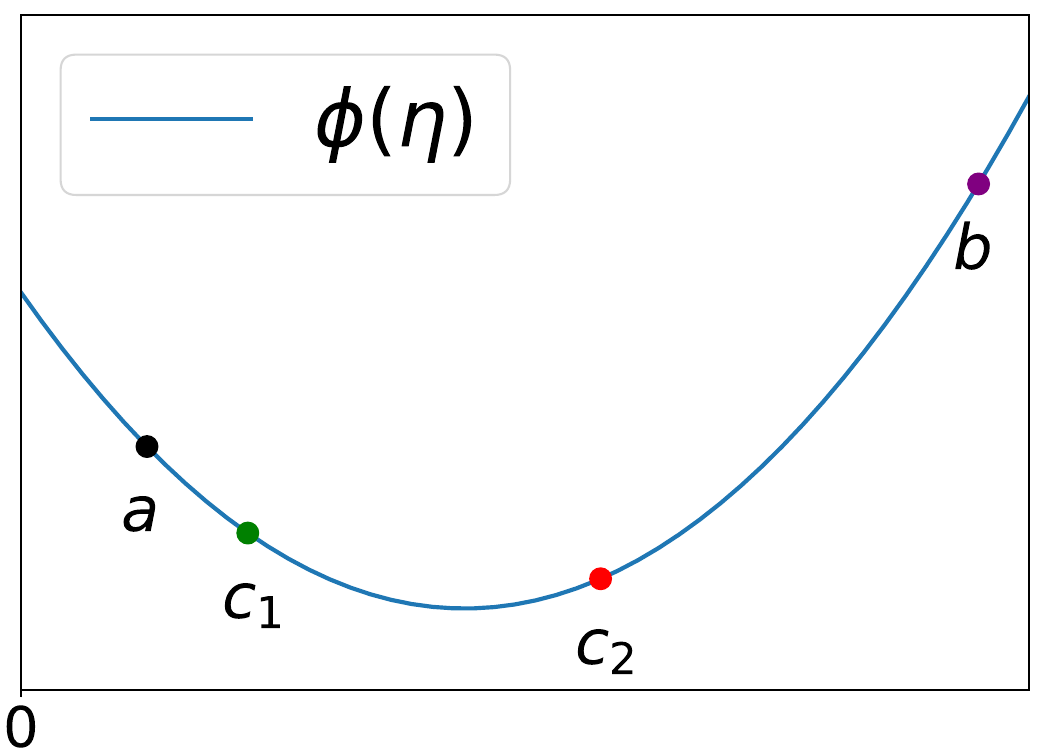}}
\subfigure[$\eta=b$ yields the minimum.]{\label{fig:golden_4}
\includegraphics[width=0.23\linewidth]{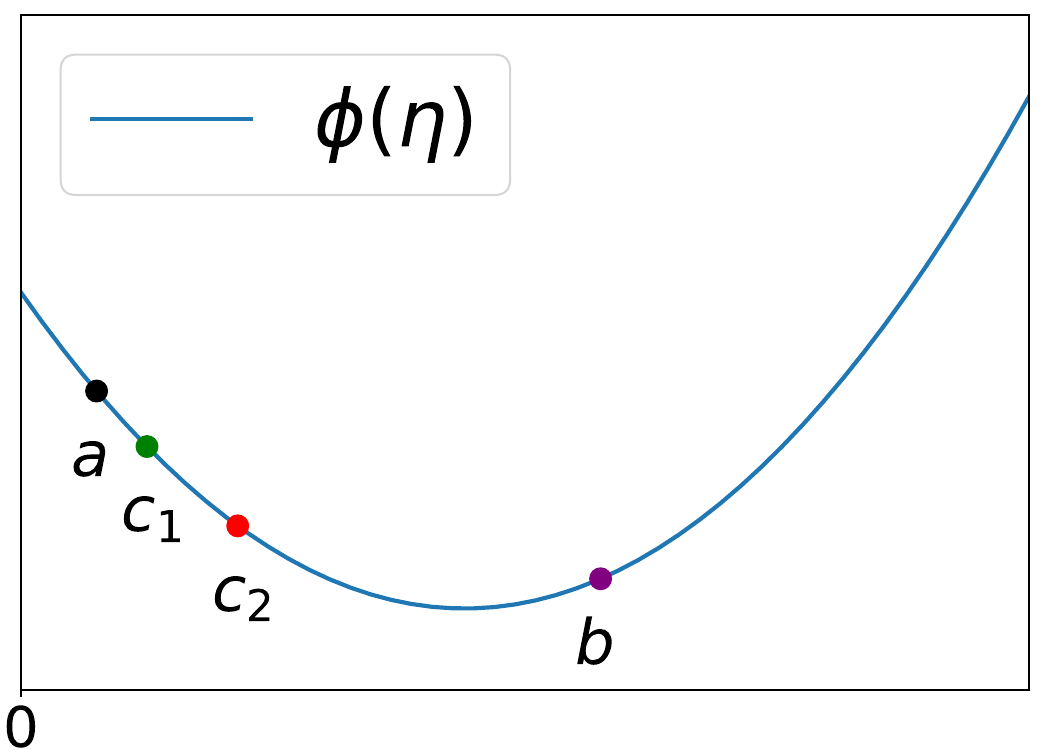}}
\caption{Demonstration of four different update ways in golden-section line search.}
\label{fig:conjguatecy-golden_1234}
\end{figure}

Similar to the bisection line search, the \textit{golden-section line search} also identifies the best learning rate $\eta$ for a unimodal function $\phi(\eta)$. 
It starts with the interval $[a,b]$ as $[0, \eta_{\max}]$. However, instead of selecting a midpoint, the golden-section search designates two points of $c_1$ and $c_2$ such that $a<c_1<c_2<b$. 
The procedure is as follows:  
\begin{itemize}
\item If $\eta=a$ results in the minimum value for $\phi(\eta)$ (among the four values $\phi(a), \phi(c_1), \phi(c_2)$, and $\phi(b)$), we can exclude the interval $(c_1, b]$. 
\item If $\eta=c_1$ yields the minimum value, we can exclude the interval $(c_2, b]$. 
\item If $\eta=c_2$ yields the minimum value, we can exclude the interval $[a, c_1)$. 
\item If $\eta=b$ yields the minimum value, we can exclude the interval $[a, c_2)$.
\end{itemize}
These four situations are illustrated in Figure~\ref{fig:conjguatecy-golden_1234}.
In other words, at least one of the intervals $[a,c_1]$ and $[c_2, b]$ can be discarded  in the golden-section search
method.
By excluding one of the four intervals, the new bounds $[a, b]$ are adjusted accordingly, and the process iterates until the range is sufficiently small.

\subsection{Convergence Analysis}\label{section:gd_conv_line_search}

This subsection discusses the convergence of the algorithm derived using different line search criteria. While the conclusions are relatively weak because they are established in the framework of general line search algorithms, they provide essential insights into the requirements for the convergence of these algorithms.
\begin{theorem}[Descent Lemma and Zoutendijk  under SS and Wolfe]\label{theorem:zoutendijk_cond}
Let $f:\real^n\rightarrow\real$ be a differentiable and $\beta$-smooth function that is bounded below. 	
Consider the general update rule of  descent methods in Algorithm~\ref{alg:struc_gd_gen}, where $\bd^\toptzero$ is a descent direction at the point $ \bx^\toptzero $ (not necessarily the negative gradient), $\eta_t$ is the step length, and  the Wolfe condition (Definition~\ref{definition:wolfe_cond}) is satisfied in each iteration. 
Then, the descent lemma holds for each iteration $t$:
$$
f(\bx^\toptone) -f(\bx^\toptzero) \leq  c_1 \frac{c_2 - 1}{\beta} \cos^2 (\theta_t) \normtwobig{\nabla f(\bx^\toptzero)}^2,
$$
and the following \textit{Zoutendijk condition} holds:
\begin{equation}\label{equation:zoutendijk_cond}
\sum_{t=1}^{\infty} \cos^2 (\theta_t) \normtwobig{\nabla f(\bx^\toptzero)}^2 < +\infty.~\footnote{This implies that $\cos^2 (\theta_t) \normtwobig{\nabla f(\bx^\toptzero)}^2\rightarrow 0$ as $t\rightarrow \infty$.}
\end{equation}
where $\cos (\theta_t)$ is the cosine of the angle between the negative gradient $-\nabla f(\bx^\toptzero)$ and the descent direction $\bd^\toptzero$, i.e.,
$
\cos (\theta_t) = \frac{-\innerproduct{\nabla f(\bx^\toptzero),  \bd^\toptzero}}{\normtwobig{\nabla f(\bx^\toptzero)} \normtwobig{\bd^\toptzero}}.
$
\end{theorem}
\begin{proof}[of Theorem~\ref{theorem:zoutendijk_cond}]
Since $\bx^\toptone = \bx^\toptzero +\eta_t\bd^\toptzero$, by the Wolfe condition \eqref{equation:wolfe_2},
$$
\innerproduct{ \nabla f(\bx^\toptone) - \nabla f(\bx^\toptzero) ,  \bd^\toptzero} \geq (c_2 - 1) \innerproduct{\nabla f(\bx^\toptzero),  \bd^\toptzero}.
$$
By the Cauchy-Schwarz inequality (Proposition~\ref{proposition:cauchy-schwarz-inequ}) and the  $\beta$-smoothness,
$$
\innerproduct{ \nabla f(\bx^\toptone) - \nabla f(\bx^\toptzero) ,  \bd^\toptzero}
\leq \normtwo{\nabla f(\bx^\toptone) - \nabla f(\bx^\toptzero)} \normtwobig{\bd^\toptzero} \leq \eta_t \beta \normtwobig{\bd^\toptzero}^2.
$$
Combining these  two inequalities yields that 
$$
\eta_t \geq \frac{c_2 - 1}{\beta} \frac{\innerproduct{\nabla f(\bx^\toptzero),  \bd^\toptzero}}{\normtwobig{\bd^\toptzero}^2}.
$$
Since $\bd^\toptzero$ is a descent direction, implying that $\innerproduct{\nabla f(\bx^\toptzero),  \bd^\toptzero} < 0$, substituting the above into the Wolfe condition \eqref{equation:wolfe_1}, whence we have the following \textit{descent inequality}~\footnote{Such inequalities are called descent lemma for descent methods and play a core role in proving the convergence results for gradient descent methods and their variants in Chapter~\ref{chapter:gd_convg}.}:
$$
f(\bx^\toptone) -f(\bx^\toptzero) \leq  c_1 \frac{c_2 - 1}{\beta} \frac{\big(\innerproduct{\nabla f(\bx^\toptzero),  \bd^\toptzero}\big)^2}{\normtwobig{\bd^\toptzero}^2}
= c_1 \frac{c_2 - 1}{\beta} \cos^2 (\theta_t) \normtwobig{\nabla f(\bx^\toptzero)}^2,
$$
where the last equality follows from the definition of $\theta_t$.
Performing telescopic cancellations over $t\in\{1,2,\ldots,T\}$, we have
$$
f(\bx^{(T+1)}) \leq f(\bx^{(1)}) + c_1 \frac{c_2-1}{\beta} \sum_{t=1}^T \cos^2 (\theta_t) \normtwobig{\nabla f(\bx^\toptzero)}^2.
$$
Since the function $f$ is bounded below and $0 < c_1 < c_2 < 1$, it follows that $c_1 (c_2-1) < 0$. Therefore, as $t \to \infty$,
$
\sum_{t=1}^\infty \cos^2 (\theta_t) \normtwobig{\nabla f(\bx^\toptzero)}^2 < +\infty,
$
establishing the desired result.
\end{proof}

Theorem~\ref{theorem:zoutendijk_cond} indicates that as long as the iterations satisfy the Wolfe condition, for a smooth and bounded below function, the Zoutendijk condition always holds. In fact, using the Goldstein condition can also lead to similar conclusions \citep{nocedal1999numerical, sun2006optimization}. The Zoutendijk condition characterizes the nature of line search conditions, and together with the selection of the descent direction $\bd^\toptzero$, we can obtain convergence results for descent methods.

In the aforementioned theorem, given a descent direction $\bd$ at point $\bx$, we introduce the angle between $\bd$ and $-\nabla f(\bx)$:
$$
\cos (\theta) = \frac{-\bd^\top \nabla f(\bx)}{\normtwo{\bd}\normtwo{\nabla f(\bx)}}.
$$
The existence of the angle is guaranteed by the Cauchy-Schwartz inequality in Definition~\ref{definition:angle_bet_vec_ineq} and provide a type of descent method known as the \textit{absolute descent method}.

\begin{definition}[Absolute Descent Method]\label{definition:absolute_descent_direction}
Let $f:\real^n\rightarrow\real$ be a differentiable function.
Consider the general update rule of  descent methods in Algorithm~\ref{alg:struc_gd_gen}, where $\bd^\toptzero$ is a descent direction at the point $ \bx^\toptzero $ (not necessarily the negative gradient). 
Let $\theta_t$ be the angle between the negative gradient $-\nabla f(\bx^\toptzero)$ and the descent direction $\bd^\toptzero$: 
$\cos(\theta_t) = \frac{-\innerproduct{\bd^\toptzero \nabla f(\bx^\toptzero)}}{\normtwo{\bd^\toptzero}\normtwo{\nabla f(\bx^\toptzero)}}$, and assume for any $t>0$, there exists a constant $\zeta > 0$ such that
\begin{equation}\label{equation:abs_des_dire2}
\theta_t < \frac{\pi}{2} - \zeta,\quad \forall t>0, \quad\text{ with } \zeta>0 \text{ independent of $ t $},
\end{equation}
ensuring $\theta_t\in(0, \frac{\pi}{2})$.
That is, the angle $\theta_t$ between $\bd^\toptzero$ and $-\nabla f(\bx^\toptzero)$ is \textit{uniformly bounded away} from $90^{\circ}$. Note that if $\theta=90^{\circ}$, then $\bd^\toptzero$ is almost not a descent direction.
\end{definition}
The restriction that $\zeta$ must be constant in all the steps is necessary for the global convergence result for iterative methods.
Then we have the following convergence results for descent methods.
\begin{theoremHigh}[Convergence of Descent Method under Zoutendijk]\label{theorem:conv_line_search}
Under the same conditions as Theorem~\ref{theorem:zoutendijk_cond}, let $\bd^\toptzero$ satisfy the absolute descent direction in \eqref{equation:abs_des_dire2}.
Then the sequence $\{\bx^\toptzero\}_{t>0}$ converges to a stationary point of $f$.
That is, if  the Zoutendijk condition $
\sum_{t=1}^{\infty} \cos^2 (\theta_t) \normtwobig{\nabla f(\bx^\toptzero)}^2 < +\infty
$ is satisfied~\footnote{We should also note that the Zoutendijk condition is derived from the Wolfe condition; More generally, if  exact line search is applicable, the Zoutendijk condition is also satisfied since  exact line search represents a special case of the Wolfe condition.}, it follows that
$$
\lim_{t \to \infty} \nabla f(\bx^\toptzero) = \bzero.
$$
\end{theoremHigh}
\begin{proof}[of Theorem~\ref{theorem:conv_line_search}]
Assume, for the sake of contradiction, that the conclusion is false. That is, there exists a subsequence $\{t_k\}_{k>0}$ and a positive constant $\delta > 0$ such that
$$
\normtwo{\nabla f(\bx^{(t_k)})} \geq \delta, \quad k = 1, 2, \ldots.
$$
Based on the assumption of $\theta_t$, for any $t>0$,
$
\cos (\theta_t) > \cos \left( \frac{\pi}{2} - \zeta \right) = \sin (\zeta) > 0.
$
Focusing only on the terms corresponding to$t_k$ in Equation~\eqref{equation:zoutendijk_cond}, and we have
$$
\sum_{t=1}^\infty \cos^2 (\theta_t) \normtwobig{\nabla f(\bx^\toptzero)}^2 \geq \sum_{k=1}^\infty \cos^2 (\theta_{t_k}) \normtwo{\nabla f(\bx^{(t_k)})}^2
> \sum_{l=1}^\infty\sin^2 (\zeta) \cdot \delta^2 \to +\infty.
$$
This leads to a contradiction to Theorem~\ref{theorem:zoutendijk_cond}. Therefore, it must hold that 
$
\lim_{t \to \infty} \nabla f(\bx^\toptzero) = \bzero,
$
which completes the proof.
\end{proof}

Theorem~\ref{theorem:conv_line_search} is derived from the Zoutendijk condition, which essentially imposes a constraint on the absolute descent direction. Specifically, it ensures that each descent direction $\bd^\toptzero$ does not become asymptotically orthogonal to the negative gradient direction.  
This condition is geometrically intuitive: when the descent direction $\bd^\toptzero$ and the gradient are nearly orthogonal, according to  Taylor's expansion, the objective function value $f(\bx^\toptzero)$ would experience little to no change. 
Consequently, we require $\bd^\toptzero$ to maintain a consistent lower bound on the angle relative to the gradient's orthogonal direction. 
Apparently, the steepest descent direction and the greedy search directions we introduced in the previous sections satisfy this condition.

Under the Goldstein condition, we can also obtain a convergence result under mild assumptions.
\begin{theoremHigh}[Convergence of Descent Method under Goldstein or Wolfe]
Let $f:\real^n\rightarrow\real$ be a differentiable function. 	
Consider the general update rule of  descent methods in Algorithm~\ref{alg:struc_gd_gen}, where $\bd^\toptzero$ is a descent direction at the point $ \bx^\toptzero $ (not necessarily the negative gradient), and let $\eta_t$ in Algorithm~\ref{alg:struc_gd_gen} be defined by the Goldstein condition (Definition~\ref{definition:goldstein_cond}) or Wolfe condition (Definition~\ref{definition:wolfe_cond}).
Suppose that  the descent direction $\bd^\toptzero$ satisfy  the absolute descent direction in \eqref{equation:abs_des_dire2} for each iteration. If $\nabla f$ exists and is uniformly continuous on the level set $\sL\triangleq\lev[f, f(\bx^{(1)})]\triangleq\{\bx \mid f(\bx) \leq f(\bx^{(1)})\}$, then either 
$$
\text{$\nabla f(\bx^\toptzero) = \bzero$ for some $t$, $\quad$ or $f(\bx^\toptzero) \to -\infty$,$\quad$ or $\nabla f(\bx^\toptzero) \to \bzero$.}
$$
\end{theoremHigh}
\begin{proof}
Let $\eta_t$ be defined  by the Goldstein condition (Definition~\ref{definition:goldstein_cond}). 
Suppose that, for all $t$, $\bg^\toptzero \triangleq \nabla f(\bx^\toptzero) \neq \bzero$ (whence $\bh^\toptzero \triangleq \eta_t \bd^\toptzero \neq \bzero$) and $f(\bx^\toptzero)$ is bounded below, it follows that $f(\bx^\toptzero) - f(\bx^\toptone) \to 0$, hence $-\bg^\toptzeroTOP \bh^\toptzero \to 0$ from the first inequality of the Goldstein condition (or the first inequality of the  Wolfe condition) in \eqref{equation:goldstein1}.

Now, assume that $\bg^\toptzero \to \bzero$ does not hold. Then there exist a constant $\varepsilon > 0$ and a subsequence $\sT$ such that $\normtwobig{\bg^\toptzero} \geq \varepsilon$  and $\normtwobig{\bh^\toptzero} \to 0$ for $t\in\sT$. Based on the assumption of $\theta_t$, for any $t>0$,
$
\cos (\theta_t) > \cos \left( \frac{\pi}{2} - \zeta \right) = \sin (\zeta) > 0,
$
whence we have
$$
-\bg^\toptzeroTOP \bh^\toptzero \geq \sin (\zeta) \normtwobig{\bg^\toptzero} \normtwobig{\bh^\toptzero} \geq \varepsilon \sin (\zeta) \normtwobig{\bh^\toptzero}.
$$
By the mean value theorem (Theorem~\ref{theorem:mean_approx}), there exist a vector $\bxi^\toptzero$  on the line segment $(\bx^\toptzero, \bx^\toptone)$ such that 
$$
f(\bx^\toptone) = f(\bx^\toptzero) + \nabla f(\bxi^\toptzero)^\top \bh^\toptzero,
$$
By the uniform continuity of $\nabla f$, we have $\nabla f(\bxi^\toptzero) \to \bg^\toptzero$ when $\bh^\toptzero \to 0$. So
$$
f(\bx^\toptone) = f(\bx^\toptzero) + \bg^\toptzeroTOP \bh^\toptzero + o(\normtwobig{\bh^\toptzero}).
$$
Therefore,  it follows that 
$$
\frac{f(\bx^\toptone) - f(\bx^\toptzero) }{\bg^\toptzeroTOP \bh^\toptzero} \to 1 \quad \text{ as }\quad \normtwobig{\bh^\toptzero} \to 0,
$$
which contradicts the second inequality of the Goldstein condition in \eqref{equation:goldstein2}. Hence, we conclude that $ \bg^\toptzero \to \bzero $, completing the proof.

Similarly, for the Wolfe condition, we rearrange its second inequality in \eqref{equation:wolfe_2} as follows:
$$
-\innerproduct{\bg^\toptzero, \bd^\toptzero} 
\leq  
\frac{\innerproduct{(\bg^\toptone - \bg^\toptzero), \bd^\toptzero}}{1-c_2}
\leq 
\frac{\normtwo{\bg^\toptone - \bg^\toptzero} \normtwo{\bd^\toptzero}}{1-c_2}
$$
where $c_2\in(0,1)$.
Since the uniform continuity $\nabla f(\bx) $ implies that $\normtwo{\bg^\toptone - \bg^\toptzero} = o(1)$ for $t\in\sT$ (Remark~\ref{remark:conse_uniform_cont}) and given the assumption that $\normtwobig{\bg^\toptzero} \geq \varepsilon$ for $t\in \sT$, the above inequality also implies that 
$$
\cos(\theta_t) =\frac{-\innerproduct{\bg^\toptzero, \bd^\toptzero} }{\normtwo{\bg^\toptzero}\normtwo{\bd^\toptzero}}
\leq \frac{o(1)}{(1-c_2)\normtwo{\bg^\toptzero}}
\leq o(1),
\quad \text{for }t\in \sT,
$$
which contradicts the assumption of the absolute descent direction in \eqref{equation:abs_des_dire2}. Thus, $\bg^\toptzero \rightarrow \bzero$.
This completes the proof.
\end{proof}

\section{Descent Methods with Trust Region: An Overview}\label{section:trs_intro}
Equation~\eqref{equation:gd_gree_approx}, based on linear approximation theorem (Theorem~\ref{theorem:linear_approx}), provides us with a linear approximation to $ f $ near a given point $ \bx^\toptzero $:
\begin{subequations}\label{equation:trs_first_all}
\begin{align}
	&f(\bx^\toptzero+\bd) \approx \psi_t(\bd)   \label{equation:trs_first1} \\
	&\text{with } \psi_t(\bd) \triangleq f(\bx^\toptzero) + \bd^\top \nabla f(\bx^\toptzero).  \label{equation:trs_first2}
\end{align}
\end{subequations}
Similarly, we obtain a quadratic approximation (Theorem~\ref{theorem:quad_app_theo}) to $ f $
\begin{subequations}\label{equation:trs_second_all}
\begin{align}
	& f(\bx^\toptzero+\bd) \approx \psi_t(\bd) \label{equation:trs_second1} \\
	& \text{with } \psi_t(\bd) \triangleq f(\bx^\toptzero) + \bd^\top \nabla f(\bx^\toptzero) + \frac{1}{2} \bd^\top \nabla^2 f(\bx^\toptzero) \bd. \label{equation:trs_second2}
\end{align}
\end{subequations}
In both cases, $ \psi_t(\bd) $ provides a good approximation to $ f(\bx^\toptzero+\bd) $ only when $ \bd $ is sufficiently small. 
These approximations motivate the following iterative step, which is determined by solving the following model problem:
\begin{equation}\label{equation:trs_subpro_gd}
\begin{aligned}
	& \bd_{\text{tr}}^\toptzero = \mathop{\argmin}_{\bd \in \sS}\,\psi_t(\bd) , \\
	& \text{where } \sS \triangleq \{ \bd \mid \normtwo{\bd} \leq \Delta \}, \, \Delta > 0.
\end{aligned}
\end{equation}
The region $ \sS $ is referred to as the \textit{trust region} and $ \psi_t(\bd) $ is given by \eqref{equation:trs_first2} or \eqref{equation:trs_second2}.
Since \eqref{equation:trs_second2} is used extensively in the literature \citep{nocedal1999numerical}, we defer its analysis to Section~\ref{section:des_trust_reg}, where we introduce second-order methods.

We use $ \bd = \bd_{\text{tr}}^\toptzero $ as a candidate to our next step, and reject it if $ f(\bx^\toptzero+\bd) \geq f(\bx^\toptzero) $. 
The reduction in the objective function value determines the size of the trust region for the next step. 
Specifically, the actual reduction is compared to the predicted reduction from the approximation function, leading to the definition of the \textit{gain factor}:
\begin{equation}\label{equation:gain_fact_gd}
\textbf{(Gain factor)}: \qquad \nu_t \triangleq \frac{f(\bx^\toptzero) - f(\bx^\toptzero+\bd_{\text{tr}}^\toptzero)}{\psi_t(\bzero) - \psi_t(\bd_{\text{tr}}^\toptzero)}. 
\end{equation}
When $\nu_t$ is small,  the approximation poorly reflects $f$, whereas a large $\nu_t$ indicates good agreement. 
Consequently, the gain factor regulates the trust region's size for the next step (or, if
$\nu_t\leq  0$,  $\bd_{\text{tr}}^\toptzero$ is rejected and it triggers a re-evaluation of $\bd_{\text{tr}}^\toptzero$).
These ideas are summarized in Algorithm~\ref{alg:trust_region0}.
The numerical values in the algorithm, 0.75, 2, 0.25, and 1/3, are chosen based on practical experience; while $\gamma$ is selected for convergence issues (Theorem~\ref{theorem:global_conv_trs_eta0} and Theorem~\ref{theorem:glo_conv_trs_othe}). Although minor variations in these values do not significantly impact performance, the parameters $p_1$ and $p_2$ in the expressions $\Delta_{t+1} \leftarrow p_1 \Delta_t$ and $\Delta_{t+1} \leftarrow \Delta_t / p_2$ must be chosen carefully to prevent oscillations in the trust region radius.

\begin{algorithm}
\caption{Descent Method with Trust Region}
\label{alg:trust_region0}
\begin{algorithmic}[1]
\Require A first  or twice  differentiable function $f(\bx)$ for \eqref{equation:trs_first_all} or \eqref{equation:trs_second_all}, respectively; 
\State {\bfseries Input:}   Set the maximum radius $ \Delta_{\max} $, initial radius $ \Delta_1 $, initial point $ \bx^{(1)} $, accept radius $\gamma\in[0,\frac{1}{4})$;
\For{$t=1,2,\ldots$}
\State $\bd_{\text{tr}}^{\toptzero} \leftarrow \text{Solution of trust region subproblem \eqref{equation:trs_subpro_gd}}$;
\State $\nu_t \leftarrow \text{gain factor \eqref{equation:gain_fact_gd}}$;
\If{$\nu_t > 0.75$ and $\normtwobig{\bd^\toptzero} =\Delta_t$} \Comment{very good step, and the step is at the border}
\State $\Delta_{t+1} \leftarrow \min\{2 \Delta_t, \Delta_{\max}\}$; \Comment{larger trust region}
\EndIf
\If{$\nu_t < 0.25$} \Comment{poor step}
\State $\Delta_{t+1} \leftarrow \Delta_t / 3$; \Comment{smaller trust region}
\EndIf
\If{$\nu_t > \gamma$} \Comment{reject step if $\nu_t \leq \gamma$}
\State $\bx^\toptone \leftarrow \bx^\toptzero + \bd_{\text{tr}}^{\toptzero}$;
\EndIf
\EndFor
\State {\bfseries Return:} final $\bx\leftarrow \bx^\toptzero$;
\end{algorithmic}
\end{algorithm}

\index{Quadratic form}
\index{Fisher information matrix}
\index{Positive definite}
\index{Positive semidefinite}
\index{Symmetry}
\section{Quadratic Model using Gradient Descent with Fixed Stepsize}\label{section:quadratic_vanilla_GD}

We further explore gradient descent with a fixed stepsize applied to the simplest model: the  quadratic function,
\begin{equation}\label{equation:quadratic-form-general-form}
f(\bx) = \frac{1}{2} \bx^\top \bA \bx - \bb^\top \bx + c, \gap \bx\in \real^n,
\end{equation}
where $\bA\in \real^{n\times n}$, $\bb \in \real^n$, and $c$ is a scalar constant.
The function is convex (resp. strictly convex) if $\bA$ is positive semidefinite (resp. positive definite); see Exercise~\ref{exercise:conv_quad}. Although the quadratic form in \eqref{equation:quadratic-form-general-form} is a highly simplified model, it is versatile enough to approximate many other functions, such as the Fisher information matrix \citep{amari1998natural}, while also capturing key characteristics of pathological curvature. The gradient of $f(\bx)$ at a given point $\bx$ is given by 
\begin{equation}\label{equation:unsymmetric_gd_gradient}
\nabla f(\bx) = \frac{1}{2} (\bA^\top +\bA) \bx - \bb.
\end{equation}
The stationary point (if exists) of the function is the solution of the linear system $\frac{1}{2} (\bA^\top +\bA) \bx=  \bb $:
\begin{equation}\label{equation:gd_solution_unsymmetric}
\bx^* = 2(\bA^\top +\bA)^{-1}\bb.
\end{equation}
If $\bA$ is symmetric (for most of our discussions, we will restrict to be symmetric $\bA$ or even {positive definite}), the equation reduces to 
\begin{equation}\label{equation:symmetric_gd_gradient}
\nabla f(\bx) = \bA \bx - \bb.
\end{equation}
Then the unique minimum of the function is the solution of the linear system $\bA\bx=\bb$ (Theorem~\ref{theorem:fermat_fist_opt}); and the optimal point of $\bx$ is thus given by 
$$
\bx^* = \bA^{-1}\bb
$$
if $\bA$ is nonsingular.

\begin{figure}[h!]
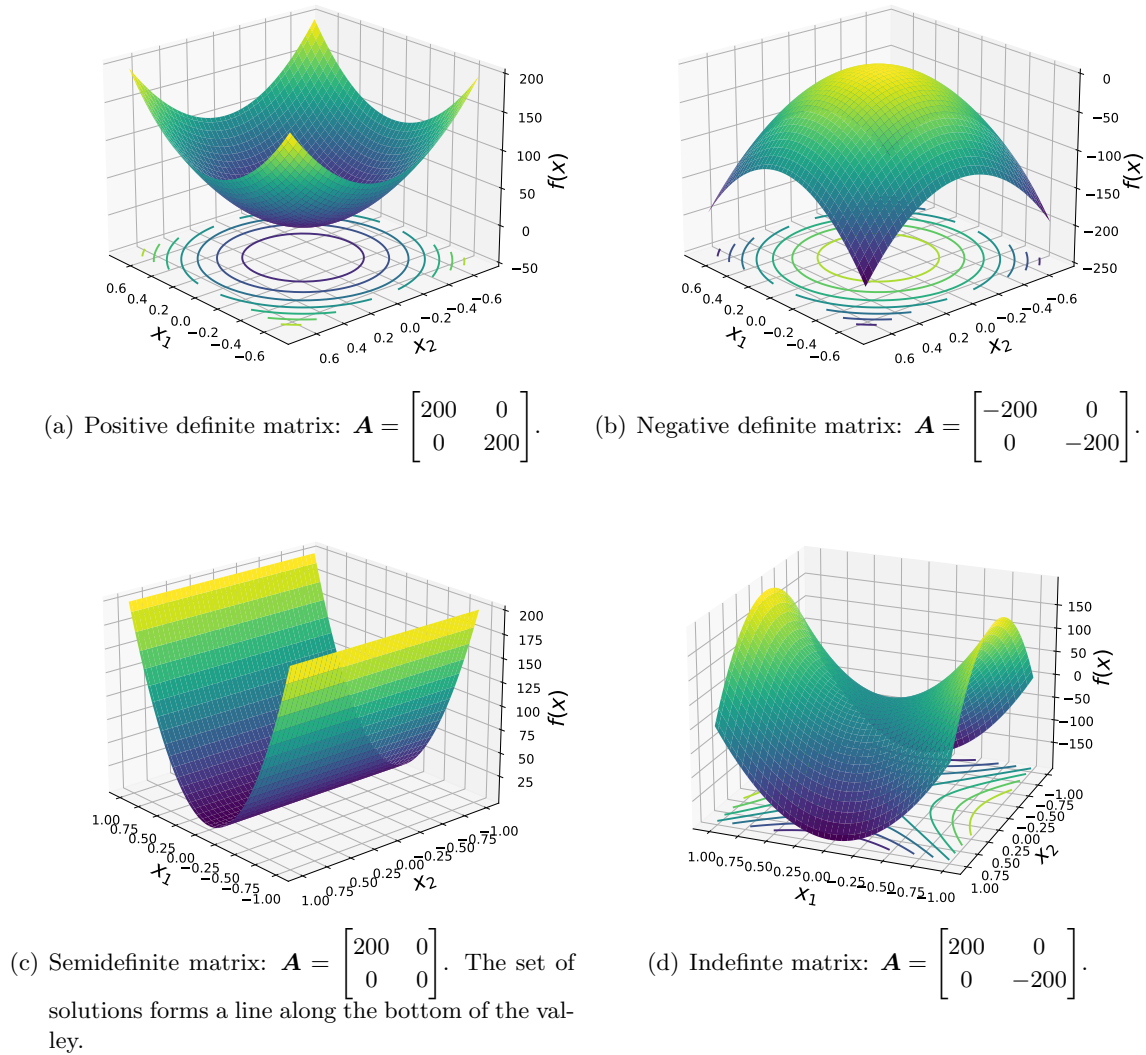

\centering  
\vspace{-0.35cm} 
\subfigtopskip=2pt 
\subfigbottomskip=2pt 
\subfigcapskip=-5pt 
\subfigure[Positive definite matrix: $\bA = \begin{bmatrix}
	200 & 0 \\ 0 & 200
\end{bmatrix}$.]{\label{fig:quadratic_PD}
	\includegraphics[width=0.485\linewidth]{imgs/quadratic_PD.pdf}}
\subfigure[Negative definite matrix: $\bA = \begin{bmatrix}
	-200 & 0 \\ 0 & -200
\end{bmatrix}$.]{\label{fig:quadratic_ND}
	\includegraphics[width=0.485\linewidth]{imgs/quadratic_ND.pdf}}
\subfigure[Semidefinite matrix: $\bA = \begin{bmatrix}
	200 & 0 \\ 0 & 0
\end{bmatrix}$. The set of solutions forms a line along the bottom of the valley.]{\label{fig:quadratic_singular}
	\includegraphics[width=0.485\linewidth]{imgs/quadratic_singular.pdf}}
\subfigure[Indefinte matrix: $\bA = \begin{bmatrix}
	200 & 0 \\ 0 & -200
\end{bmatrix}$.]{\label{fig:quadratic_saddle}
	\includegraphics[width=0.485\linewidth]{imgs/quadratic_saddle.pdf}}
\caption{Loss surfaces for different quadratic forms.}
\label{fig:different_quadratics}
\end{figure}
For different types of matrix $\bA$, the loss surface of $f(\bx)$ will be different, as illustrated in Figure~\ref{fig:different_quadratics}. When $\bA$ is positive definite, the surface forms a convex bowl; when $\bA$ is negative definite, on the contrary, the surface becomes a concave bowl. $\bA$ also could be singular, in which case $\bA\bx-\bb=\bzero$ has more than one solution, and the set of solutions is a line (in the two-dimensional case ) or a hyperplane (in the high-dimensional case). This behavior is similar to that of a semidefinite quadratic form, as shown in Figure~\ref{fig:quadratic_singular}. If $\bA$ does not fall into any of these categories, a saddle point emerges (see Figure~\ref{fig:quadratic_saddle}), posing  challenges for gradient descent.  
In such cases, alternative methods, e.g., perturbed GD \citep{jin2017escape, du2017gradient}, can be employed to navigate away from saddle points.

\begin{figure}[htp]
\centering  
\vspace{-0.35cm} 
\subfigtopskip=2pt 
\subfigbottomskip=2pt 
\subfigcapskip=-5pt 
\subfigure[Contour and the descent direction. The red dot is the optimal point.]{\label{fig:quadratic_vanillegd_contour}
	\includegraphics[width=0.31\linewidth]{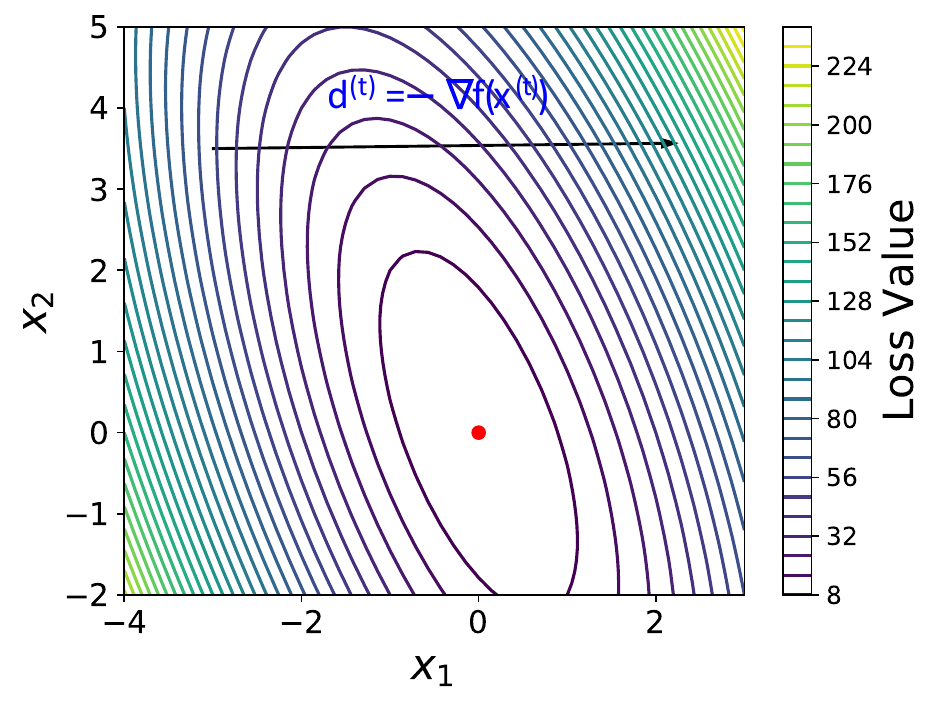}}
\subfigure[GD, $\eta=0.02$.]{\label{fig:quadratic_vanillegd_contour2}
	\includegraphics[width=0.31\linewidth]{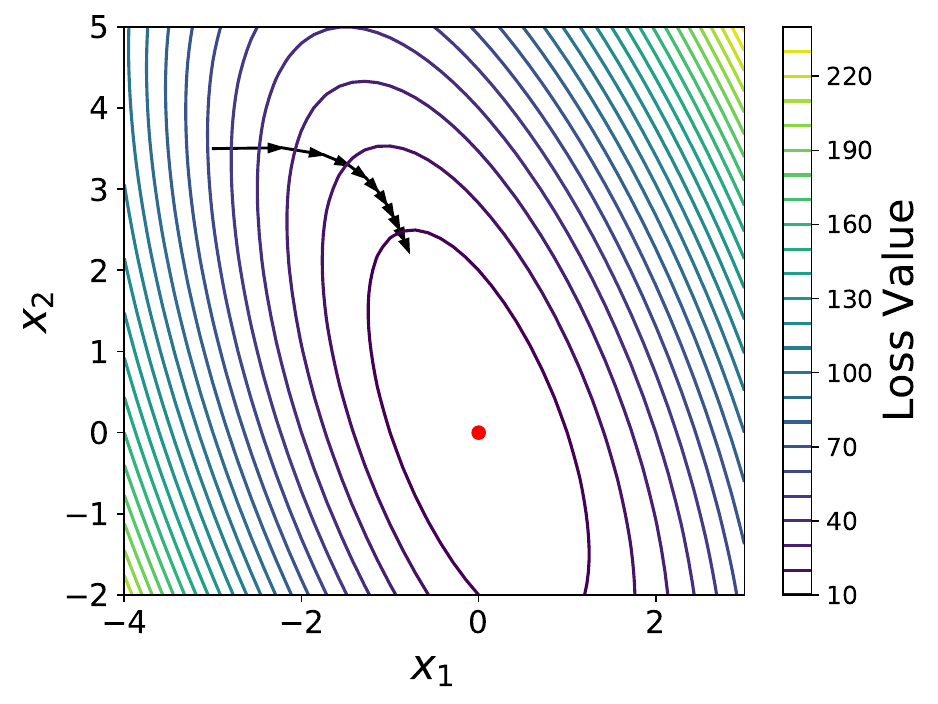}}
\subfigure[GD, $\eta=0.08$.]{\label{fig:quadratic_vanillegd_contour8}
	\includegraphics[width=0.31\linewidth]{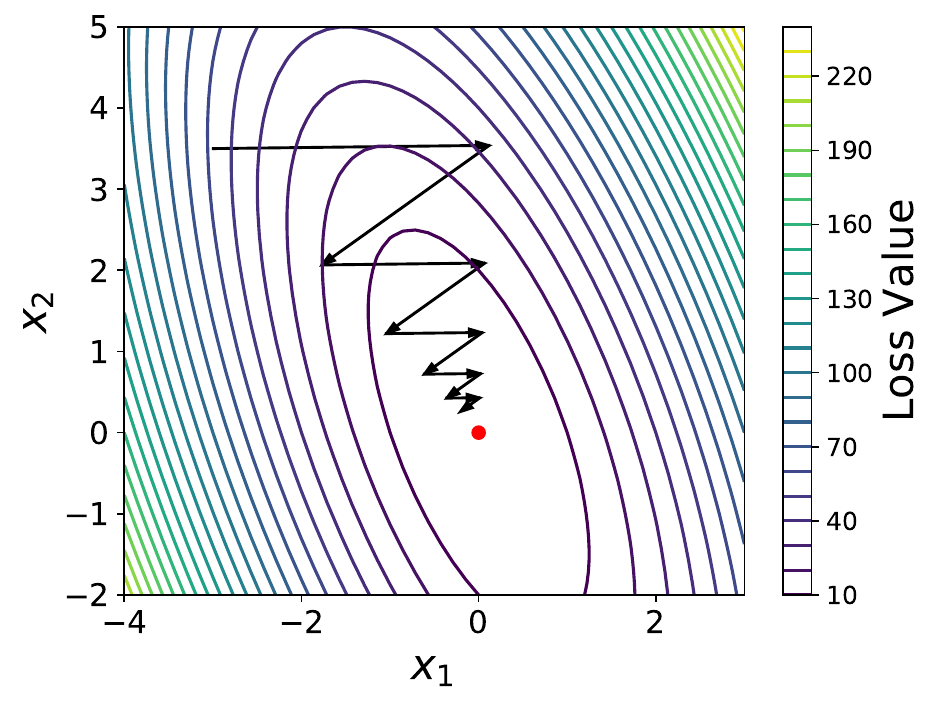}}
\caption{Illustration of  gradient descent with fixed stepsize applied to a quadratic form with $\bA=\scriptsize\begin{bmatrix}
		20 & 7 \\ 5 & 5
	\end{bmatrix}$, $\bb=\bzero$, and $c=0$. The procedure is at $\bx^\toptzero=[-3,3.5]^\top$ for the $t$-th iteration.}
\label{fig:quadratic_vanillegd}
\end{figure}

\index{Quadratic form}
\index{Saddle point}

In this context, gradient descent is not strictly necessary, as the minimum can be computed directly if the data matrix $\bA$ is nonsingular and its inverse can be determined efficiently. However, our focus is on the iterative updates of the convex quadratic function.
Suppose we initialize the process at $\bx^{(1)}\in \real^n$. At each time step $t$, the simplest update rule involves fixing the learning rate $\eta$ and selecting a descent direction $\bd^\toptzero$, leading to the update:
$$
\bx^\toptone \leftarrow  \bx^\toptzero + \eta \bd^\toptzero. 
$$
This results in a monotonically decreasing sequence of function values $\{f(\bx^\toptzero)\}$. 
Specifically, when the descent direction is chosen to be the negative gradient $\bd^\toptzero = -\bg^\toptzero=-(\bA\bx^\toptzero-\bb)$ ($\bA$ is positive definite), the update becomes 
\begin{equation}\label{equation:vanilla-gd-update}
\text{(GD with fixed stepsize)}:\quad \bx^\toptone = \bx^\toptzero - \eta (\bA\bx^\toptzero-\bb).
\end{equation}
A concrete example is given in Figure~\ref{fig:quadratic_vanillegd}, where $\bA=\scriptsize\begin{bmatrix}
20 & 7 \\ 5 & 5
\end{bmatrix}$, $\bb=\bzero$, and $c=0$. Suppose at $t$-th iteration, $\bx^\toptzero=[-3,3.5]^\top$. Figure~\ref{fig:quadratic_vanillegd_contour} shows the descent direction given by the negative gradient at  $\bx^\toptzero$; Figure~\ref{fig:quadratic_vanillegd_contour2} and Figure~\ref{fig:quadratic_vanillegd_contour8} illustrate  10 iterations of gradient descent with learning rates $\eta=0.02$ and $\eta=0.08$, respectively.

\index{Spectral decomposition}
\paragrapharrow{Closed-form for GD.}
When $\bA$ is positive definite, it admits a spectral decomposition (Theorem~\ref{theorem:spectral_theorem}):
$$
\bA=\bQ\bLambda\bQ^\top \in \real^{n\times n} 
\quad\implies\quad 
\bA^{-1} = \bQ\bLambda^{-1}\bQ^\top,
$$ 
where $\bQ = [\bq_1, \bq_2, \ldots , \bq_n]$ consists of mutually orthonormal eigenvectors of $\bA$, and $\bLambda = \diag(\lambda_1, \lambda_2, \ldots , \lambda_n)$ contains the corresponding real eigenvalues.  
Since $\bA$ is positive definite,  the eigenvalues are all positive (Theorem~\ref{theorem:eigen_charac}). By convention, we order the eigenvalues such that $\lambda_1\geq \lambda_2\geq \ldots \geq \lambda_n$. Define the following iterate vector at iteration $t$ as
\begin{equation}\label{equation:vanilla-yt}
\by^\toptzero \triangleq \bQ^\top(\bx^\toptzero - \bx^*),
\end{equation}
where $\bx^* = \bA^{-1}\bb$.
By the update rule $\bx^\toptone \leftarrow \bx^\toptzero-\eta\nabla f(\bx^\toptzero)$,
it then follows that
$$
\small
\begin{aligned}
\by^\toptone
&\triangleq \bQ^\top(\bx^\toptone - \bx^*) = \bQ^\top(\bx^\toptzero - \eta(\bA\bx^\toptzero-\bb) - \bx^*)\\
&= \by^\toptzero - \eta \bQ^\top (\bQ\bLambda\bQ^\top\bx^\toptzero-\bb) 
= \by^\toptzero - \eta  (\bLambda\bQ^\top\bx^\toptzero-\bQ^\top\bb) \\
&= \by^\toptzero - \eta \bLambda\bQ^\top (\bx^\toptzero-\bx^*) = \by^\toptzero - \eta \bLambda \by^\toptzero \\
&= (\bI - \eta \bLambda)\by^\toptzero  = (\bI - \eta \bLambda)^t\by^{(1)}. \\
\end{aligned}
$$
From this, the error at each iteration follows:
\begin{equation}\label{equation:vanilla-gd-closedform}
\normtwo{\bx^\toptone - \bx^*}^2 = \normtwo{\bQ\by^\toptone}^2 = \normtwo{\bQ(\bI - \eta \bLambda)^t\by^{(1)}}^2 = \normtwo{\sum_{i=1}^{n} y_{i}^{(1)} \cdot (1-\eta \lambda_i)^t \bq_i}^2,
\end{equation}
where $\by^{(1)}$ depends on the initial parameter $\bx^{(1)}$, and $y_{i}^{(1)}$ is the $i$-th element of $\by^{(1)}$. 
Intuitively, $\by^\toptone$ represents  the error in the $\bQ$-basis at iteration $t+1$. 
Equation~\eqref{equation:vanilla-gd-closedform} highlights that the learning rate should satisfy:
\begin{equation}\label{equation:vanillagd-quandr-rate-chgoices}
\abs{1-\eta\lambda_i} < 1, \gap \forall \,\, i\in \{1,2,\dots, n\}.
\end{equation}
The error consists of $n$ terms, each with its own convergence rate, governed by $\abs{1-\eta\lambda_i}$. 
The closer this value is to 1, the slower the convergence in that dimension \citep{shewchuk1994introduction, o2015adaptive, goh2017momentum}.

\begin{figure}[h]
\centering
\includegraphics[width=0.5\textwidth]{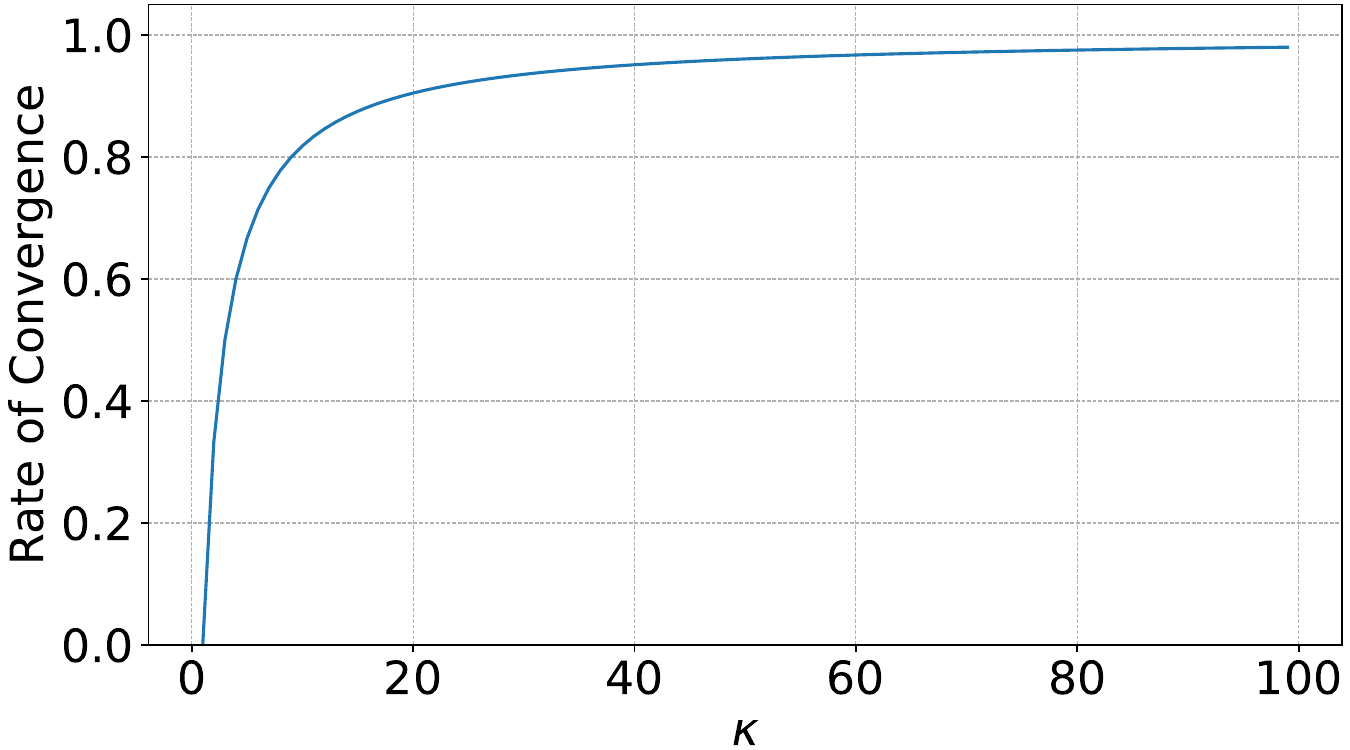}
\caption{Rate of convergence (per iteration) in  GD with fixed stepsize. The $y$-axis is $\frac{\kappa-1}{\kappa+1}$.}
\label{fig:rate_convergen_vanillaGD}
\end{figure}
\index{Rate of convergence}

To ensure convergence, the learning rate must satisfy that $\abs{1-\eta\lambda_i} < 1$, implying $0<\eta\lambda_i <2$ for $i$ in $\{1,2,\ldots, n\}$. Therefore, the overall rate of convergence is determined by the slowest component:
$$
\text{rate}(\eta) = \max \{\abs{1-\eta\lambda_1},  \abs{1-\eta\lambda_n}\},
$$
since $\lambda_1\geq \lambda_2\geq \ldots \geq \lambda_n$. The optimal learning rate is obtained when the first and last eigenvectors converge at the same rate, i.e., $\eta\lambda_1-1 =1- \eta\lambda_n$:
\begin{equation}\label{equation:eta-vanilla-gd}
\text{optimal } \eta = \underset{\eta}{\arg\min} \text{ rate}(\eta) = \frac{2}{\lambda_1+\lambda_n},
\end{equation}
and 
\begin{equation}\label{equation:vanialla-gd-rate}
\text{optimal rate}  = \underset{\eta}{\min} \text{ rate}(\eta) =
\frac{\lambda_1/\lambda_n - 1}{\lambda_1/\lambda_n + 1}
\triangleq\frac{\kappa - 1}{\kappa + 1},
\end{equation}
where $\kappa \triangleq \frac{\lambda_1}{\lambda_n} \geq 1$ is known as the \textit{condition number} (see, for example, \citet{lu2021numerical} for more details). When $\kappa=1$, convergence occurs in a single step. However, as the condition number increases, gradient descent slows down. Figure~\ref{fig:rate_convergen_vanillaGD} illustrates this effect: poorly conditioned matrices $(\kappa\gg 1)$ significantly hinder the convergence of GD.

\index{Quadratic form}
\section{Quadratic Model using Gradient Descent with Line Search}\label{section:quadratic-in-steepestdescent}
Following the discussion of the quadratic form in  gradient descent with a fixed stepsize, we now consider  the quadratic form in gradient descent with exact line search. 
Assuming that  $\bA$ is positive definite, we define $\phi(\eta)$ as follows:
$$
\begin{aligned}
\phi(\eta) \triangleq f(\bx^\toptzero+\eta \bd^\toptzero) &= \frac{1}{2} (\bx^\toptzero+\eta \bd^\toptzero)^\top \bA (\bx^\toptzero+\eta \bd^\toptzero) - \bb^\top (\bx^\toptzero+\eta \bd^\toptzero) +c\\
&= f(\bx^\toptzero) + \eta \bd^\toptzeroTOP \underbrace{\left(\bA\bx^\toptzero - \bb \right)}_{=\nabla f(\bx^\toptzero)=\bg^\toptzero} +\frac{1}{2} \eta^2 \bd^\toptzeroTOP \bA\bd^\toptzero,
\end{aligned}
$$
which is a quadratic function of $\eta$.
A closed-form solution for the line search exists (since $\bd^\toptzeroTOP \bA\bd^\toptzero>0$ for $\bd^\toptzero \neq 0 $ due to positive definiteness):
\begin{equation}\label{equation:eta-gd-steepest}
\eta_t = - \frac{\bd^\toptzeroTOP \bg^\toptzero}{ \bd^\toptzeroTOP \bA\bd^\toptzero }.
\end{equation}
When the search direction is the negative gradient $\bd^\toptzero=-\bg^\toptzero$, the descent update becomes
\begin{equation}\label{equation:steepest-quadratic}
\begin{aligned}
\text{(GD with line search)}: \quad \bx^\toptone 
&=  \bx^\toptzero + \eta_t \bd^\toptzero
= \bx^\toptzero  -  \frac{\bd^\toptzeroTOP \bg^\toptzero}{ \bd^\toptzeroTOP \bA\bd^\toptzero } \bd^\toptzero \\
& =\bx^\toptzero   - \frac{\bg^\toptzeroTOP \bg^\toptzero}{ \bg^\toptzeroTOP \bA\bg^\toptzero } \bg^\toptzero.
\end{aligned}
\end{equation}

\begin{figure}[h]
\centering  
\vspace{-0.35cm} 
\subfigtopskip=2pt 
\subfigbottomskip=2pt 
\subfigcapskip=-5pt 
\subfigure[Contour plot and the descent direction. The red dot represents the optimal point.]{\label{fig:quadratic_steepest_contour_tilt}
\includegraphics[width=0.485\linewidth]{imgs/quadratic_steepest_contour_tilt.pdf}}
\subfigure[Intersection of the loss surface and the vertical plane through the descent direction.]{\label{fig:quadratic_steepest_surface_tilt}
\includegraphics[width=0.485\linewidth]{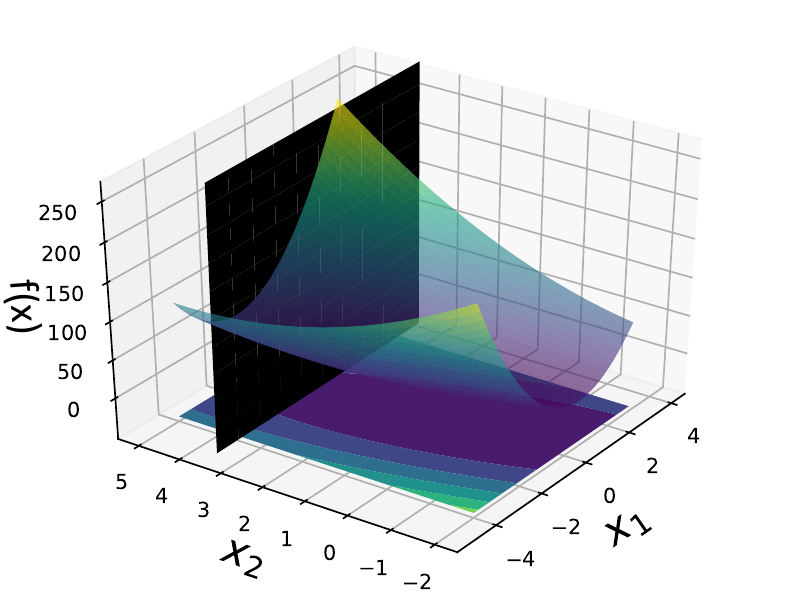}}
\subfigure[Intersection of the loss surface and the vertical plane through the descent direction in two-dimensional space.]{\label{fig:quadratic_steepest_tilt_intersection}
\includegraphics[width=0.485\linewidth]{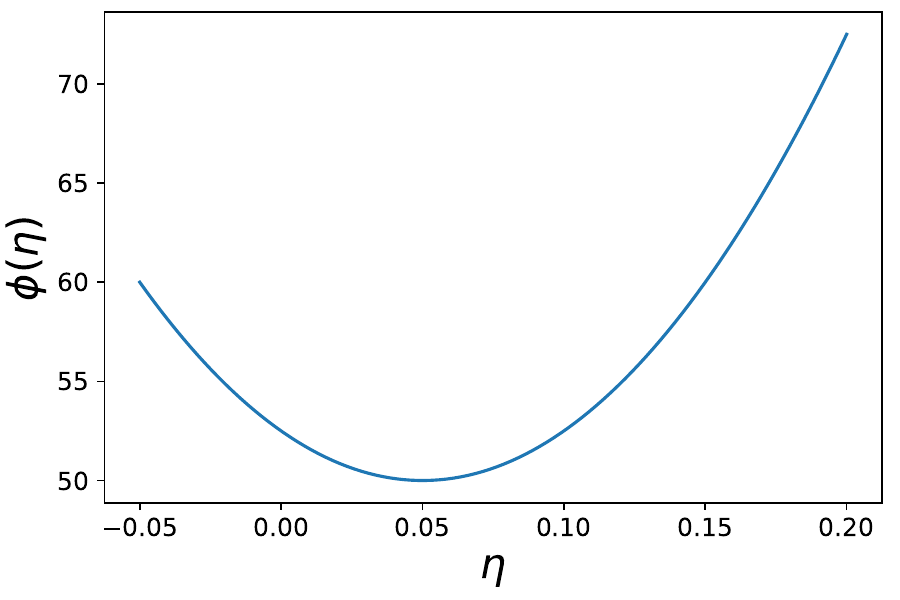}}
\subfigure[Various gradients along the line through the descent direction, where the gradient at the center point is orthogonal to the gradient of the previous step.]{\label{fig:quadratic_steepest_tilt_gradient_direction}
\includegraphics[width=0.485\linewidth]{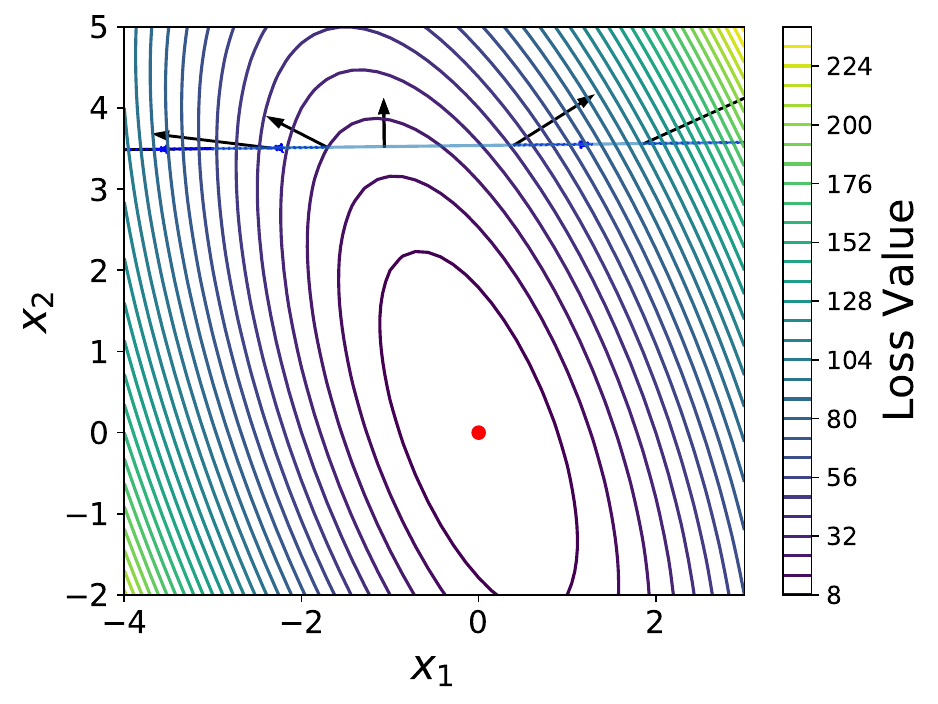}}
\caption{Illustration of GD with line search applied to a quadratic form with $\bA=\scriptsize\begin{bmatrix}
20 & 7 \\ 5 & 5
\end{bmatrix}$, $\bb=\bzero$, and $c=0$. The procedure is at $\bx^\toptzero=[-3,3.5]^\top$ for the $t$-th iteration.
}
\label{fig:quadratic_steepest_tilt}
\end{figure}

\begin{figure}[h]
\centering  
\vspace{-0.35cm} 
\subfigtopskip=2pt 
\subfigbottomskip=2pt 
\subfigcapskip=-5pt 
\subfigure[GD, $\eta=0.02$.]{\label{fig:momentum_gd_conjugate2}
\includegraphics[width=0.31\linewidth]{./imgs/steepest_gd_mom-0_lrate-2.pdf}}
\subfigure[GD, $\eta=0.08$.]{\label{fig:momentum_gd_conjugate8}
\includegraphics[width=0.31\linewidth]{./imgs/steepest_gd_mom-0_lrate-8.pdf}}
\subfigure[GD with line search.]{\label{fig:conjguatecy_zigzag2}
\includegraphics[width=0.31\linewidth]{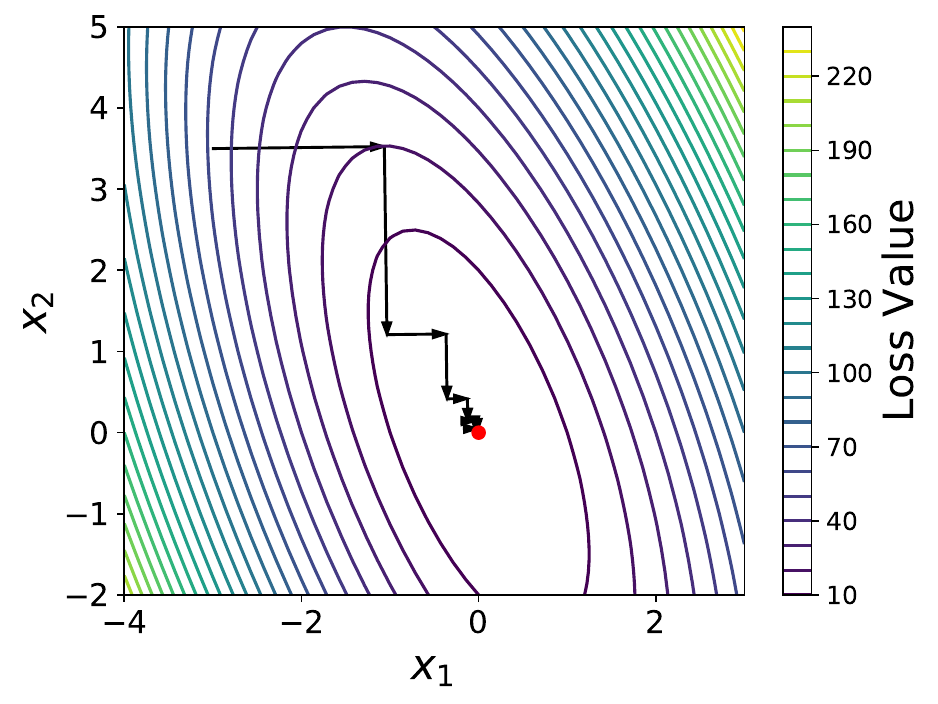}}
\caption{Illustration of GD with line search applied to a quadratic form with $\bA=\scriptsize\begin{bmatrix}
20 & 7 \\ 5 & 5
\end{bmatrix}$, $\bb=\bzero$, and $c=0$. The procedure is at $\bx^\toptzero=[-3,3.5]^\top$ for the $t$-th iteration. The example is executed for 10 iterations by  GD with $\eta=0.02$,  GD with $\eta=0.08$, and GD with line search, respectively. We observe the cumbersome choice of learning rates in vanilla GD and the zigzag trajectory in GD with line search due to the orthogonality between each gradient and the previous gradient (Lemma~\ref{lemm:linear-search-orghonal}).
}
\label{fig:quadratic_steepest_tilt22}
\end{figure}

A concrete example is presented in Figure~\ref{fig:quadratic_steepest_tilt}, where $\bA=\scriptsize\begin{bmatrix}
20 & 7 \\ 5 & 5
\end{bmatrix}$, $\bb=\bzero$, and $c=0$. 
Suppose at the $t$-th iteration, the parameter is positioned at $\bx^\toptzero=[-3,3.5]^\top$. 
Additionally, Figure~\ref{fig:quadratic_steepest_contour_tilt} shows  the descent direction using the negative gradient. 
The GD with line search involves selecting the learning rate $\eta_t$ by minimizing $\phi(\eta) = f(\bx^\toptzero + \eta\bd^\toptzero)$, which is equivalent to determining the point on the intersection of the vertical plane through the descent direction and the paraboloid defined by the loss function $f(\bx)$, as shown in Figure~\ref{fig:quadratic_steepest_surface_tilt}. Figure~\ref{fig:quadratic_steepest_tilt_intersection} further illustrates the parabola resulting from the intersection of the two surfaces. 
In Figure~\ref{fig:quadratic_steepest_tilt_gradient_direction}, various gradients along the line through the descent direction are displayed, where the gradient at the center point is orthogonal to the gradient of the previous step $\nabla f(\bx^\toptone)^\top \bd^\toptzero$, as proved in Lemma~\ref{lemm:linear-search-orghonal}. The black arrows represent the gradients, and the blue arrows are the projection of these gradients along $\bd^\toptzero = -\nabla f(\bx^\toptzero)$. 
An intuitive explanation for this orthogonality at the minimum is as follows: the slope of the parabola (Figure~\ref{fig:quadratic_steepest_tilt_intersection}) at any point is equal to the magnitude of the projection of the gradients onto the search direction (Figure~\ref{fig:quadratic_steepest_tilt_gradient_direction}) \citep{shewchuk1994introduction}. These projections represent the rate of increase of the loss function $f(\bx)$ as the point traverses the search line; and the minimum of $f(\bx)$ occurs where the projection is zero, which corresponds to the point where the gradient is orthogonal to the search line.

The example is executed for 10 iterations in Figure~\ref{fig:momentum_gd_conjugate2}, Figure~\ref{fig:momentum_gd_conjugate8}, and Figure~\ref{fig:conjguatecy_zigzag2} using  GD with $\eta=0.02$,  GD with $\eta=0.08$, and GD with (exact) line search, respectively. It is evident that  GD involves cumbersome selection of learning rates; and the zigzag trajectory is observed in the  GD with line search, resulting from the orthogonality between each gradient and the previous gradient (Lemma~\ref{lemm:linear-search-orghonal}). However, this limitation will be addressed in conjugate descent (Section~\ref{section:conjugate-descent}).

\index{Spectral decomposition}
\index{Quadratic form}
\subsubsection{Special Case:  Error Vector as an Eigenvector}
To analyze the convergence behavior of GD with line search, we examine specific scenarios. Following \citet{shewchuk1994introduction}, we first introduce two key definitions: the \textit{error vector}, which represents the difference between the parameter $\bx^\toptzero$ at the $t$-th iteration and the optimal parameter $\bx^*$, and the \textit{residual vector}, which measures the difference between the target vector $\bb$ and the predicted value $\bA\bx^\toptzero$ at iteration $t$.

\begin{definition}[Error and Residual Vector]\label{definition:error-gd-}
At iteration $t$, the \textit{error vector} is defined as $\be^\toptzero \triangleq \bx^\toptzero - \bx^*$, which quantifies the distance of the iterate from the optimal solution, where $\bx^* =\bA^{-1}\bb$ when $\bA$ is positive definite. 
Substituting this into Equation~\eqref{equation:steepest-quadratic}, the update for the error vector become
\begin{equation}\label{equation:steepest-quadratic-error}
\be^\toptone = \be^\toptzero - \frac{\bg^\toptzeroTOP \bg^\toptzero}{ \bg^\toptzeroTOP \bA\bg^\toptzero } \bg^\toptzero.
\end{equation}
Additionally, the \textit{residual vector} is defined as $\br^\toptzero \triangleq \bb - \bA\bx^\toptzero$, which indicates how far the iterate is from the correct value of $\bb$. 
Since $\bA$ is positive definite, the residual is equal to the negative gradient and the descent direction, i.e., $\br^\toptzero = \bd^\toptzero = -\bg^\toptzero$ (we may use $-\bg^\toptzero$ and $\br^\toptzero$ interchangeably when $\bA$ is positive definite). 
\end{definition}

\begin{figure}[h]
\centering  
\vspace{-0.35cm} 
\subfigtopskip=2pt 
\subfigbottomskip=2pt 
\subfigcapskip=-5pt 
\subfigure[GD with line search, $\bA\be^\toptzero=\lambda_t\be^\toptzero$.]{\label{fig:steepest_gd_bisection_eigenvector}
\includegraphics[width=0.481\linewidth]{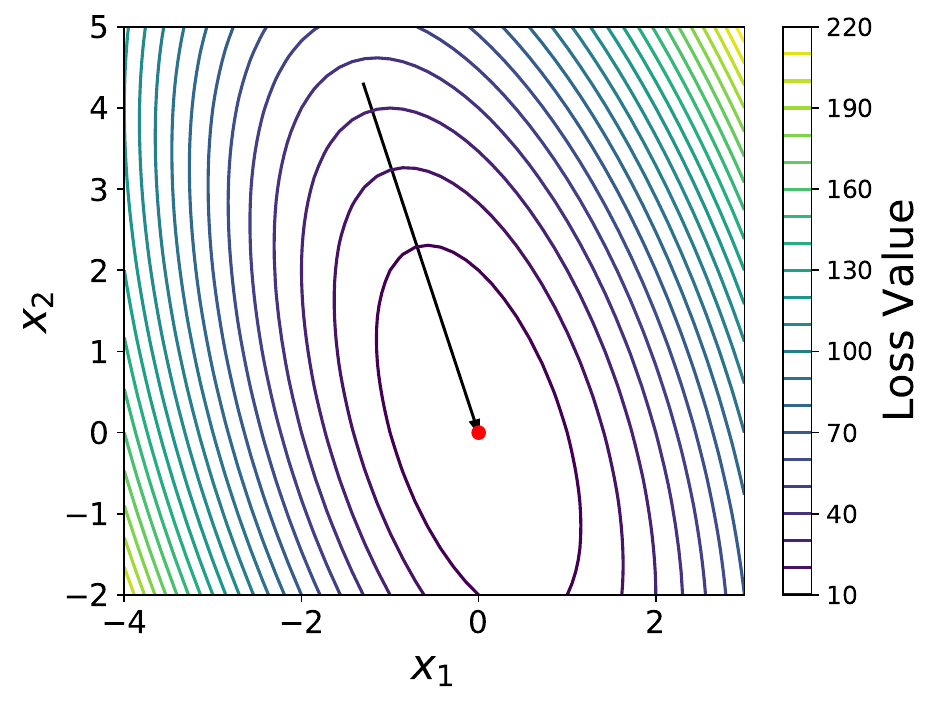}}
\subfigure[GD with line search, $\lambda_1=\lambda_2$.]{\label{fig:steepest_gd_bisection_eigenvector_sameeigenvalue}
\includegraphics[width=0.481\linewidth]{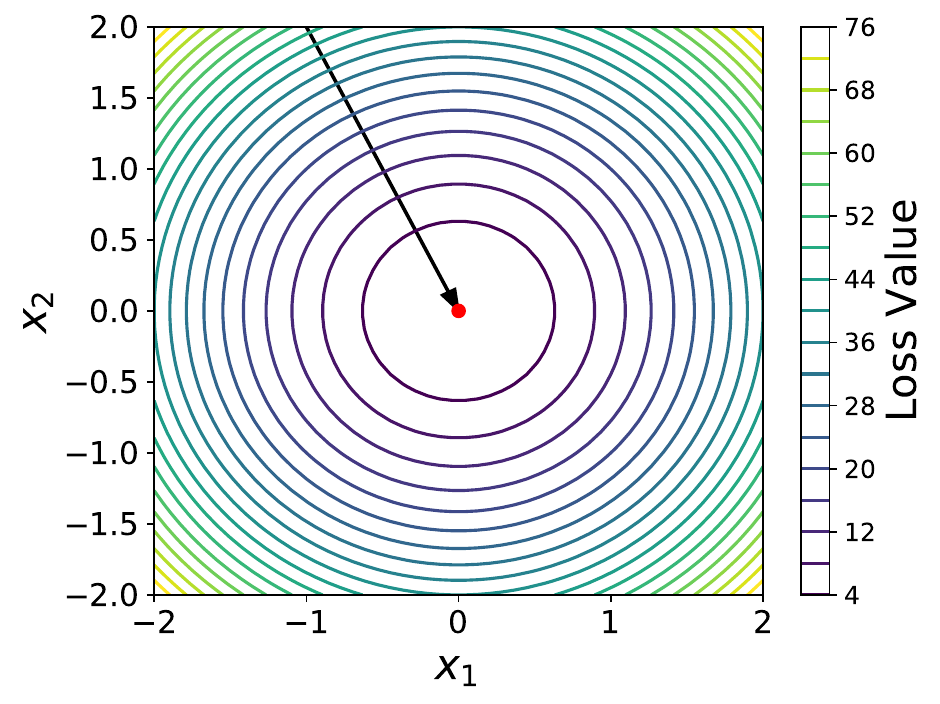}}
\caption{Illustration of special cases for GD with line search applied to quadratic forms. $\bA=\scriptsize\begin{bmatrix}
20 & 5 \\ 5 & 5
\end{bmatrix}$, $\bb=\bzero$, $c=0$, and starting point to descent is $\bx^\toptzero=[-1.3, 4.3]^\top$ for Figure~\ref{fig:steepest_gd_bisection_eigenvector}. $\bA=\scriptsize\begin{bmatrix}
20 & 0 \\ 0 & 20
\end{bmatrix}$, $\bb=\bzero$, $c=0$, and starting point to descent is $\bx^\toptzero=[-1, 2]^\top$ for Figure~\ref{fig:steepest_gd_bisection_eigenvector_sameeigenvalue}.}
\label{fig:steepest_specialcases}
\end{figure}

We first consider the case where the error vector $\be^\toptzero$ at iteration $t$ is an eigenvector of $\bA$ corresponding to eigenvalue $\lambda_t$, i.e., $\bA\be^\toptzero = \lambda_t\be^\toptzero $. Then, the gradient vector (for positive definite $\bA$ by \eqref{equation:symmetric_gd_gradient}) becomes
$$
\bg^\toptzero = \bA\bx^\toptzero-\bb = \bA\bigg(\bx^\toptzero-\underbrace{\bA^{-1}\bb}_{=\bx^*}\bigg) =\bA\be^\toptzero = \lambda_t\be^\toptzero,
$$ 
which is also an eigenvector of $\bA$ corresponding to the eigenvalue  $\lambda_t$, i.e., $\bA\bg^\toptzero = \lambda_t \bg^\toptzero$.
By \eqref{equation:steepest-quadratic-error}, the update for the $(t+1)$-th iteration is 
$$
\begin{aligned}
\be^\toptone &= \be^\toptzero   - \frac{\bg^\toptzeroTOP \bg^\toptzero}{ \bg^\toptzeroTOP \bA\bg^\toptzero } \bg^\toptzero
=\be^\toptzero - \frac{\bg^\toptzeroTOP \bg^\toptzero}{ \lambda_t \bg^\toptzeroTOP \bg^\toptzero } (\lambda_t\be^\toptzero) = \bzero .
\end{aligned}
$$
Thus, when $\be^\toptzero$ is an eigenvector of $\bA$, the GD with line search converges to the solution in just one additional step. 
A concrete example is shown in Figure~\ref{fig:steepest_gd_bisection_eigenvector}, where $\bA=\scriptsize\begin{bmatrix}
20 & 5 \\ 5 & 5
\end{bmatrix}$, $\bb=\bzero$, and $c=0$.

\subsubsection{Special Case: Spherical}
When $\bA$ is positive definite, it admits a spectral decomposition (Theorem~\ref{theorem:spectral_theorem}):
$$
\bA=\bQ\bLambda\bQ^\top \in \real^{n\times n} 
\qquad\implies\qquad
\bA^{-1} = \bQ\bLambda^{-1}\bQ^\top,
$$ 
where the columns of $\bQ = [\bq_1, \bq_2, \ldots , \bq_n]$ are mutually orthonormal eigenvectors of $\bA$, and the entries of $\bLambda = \diag(\lambda_1, \lambda_2, \ldots , \lambda_n)$ are the  corresponding eigenvalues of $\bA$, which are positive real and ordered as $\lambda_1\geq \lambda_2\geq \ldots \geq \lambda_n$. Since the eigenvectors are chosen to be mutually orthonormal, they satisfy:
$$
\bq_i^\top \bq_j =
\left\{
\footnotesize
\begin{aligned}
&1, \gap i=j;\\
&0, \gap i\neq j.
\end{aligned}
\right.
$$
Moreover, the eigenvectors form a complete basis for $\real^n$, allowing any error vector $\be^\toptzero  \in \real^n$ to be expressed as a linear combination of these eigenvectors:
\begin{equation}\label{equation:steepest-et-eigen-decom}
\be^\toptzero = \sum_{i=1}^{n} \alpha_i \bq_i,
\end{equation}
where $\alpha_i$ represents  the component of $\be^\toptzero$ along the direction of $\bq_i$.
Then the gradient vector  can be obtained by  
\begin{equation}\label{equation:steepest-eigen-decom-part}
\begin{aligned}
\bg^\toptzero = \bA\bx^\toptzero-\bb = \bA\bigg(\bx^\toptzero-\underbrace{\bA^{-1}\bb}_{=\bx^*}\bigg) =\bA\be^\toptzero
 = 
\bA  \sum_{i=1}^{n} \alpha_i \bq_i
= 
\sum_{i=1}^{n} \alpha_i \lambda_i\bq_i,
\end{aligned}
\end{equation} 
which is also a linear combination of eigenvectors, with each component scaled by the corresponding eigenvalue $\lambda_i$.
Again, by \eqref{equation:steepest-quadratic-error}, the update for the $(t+1)$-th iteration is  computed as follows:
$$
\begin{aligned}
\be^\toptone &= \be^\toptzero   - \frac{\bg^\toptzeroTOP \bg^\toptzero}{ \bg^\toptzeroTOP \bA\bg^\toptzero } \bg^\toptzero
=\be^\toptzero - \frac{\sum_{i=1}^{n}\alpha_i^2\lambda_i^2}{ \sum_{i=1}^{n}\alpha_i^2\lambda_i^3 } \sum_{i=1}^{n} \alpha_i \lambda_i\bq_i 
 .
\end{aligned}
$$
The above equation indicates  that when only one component of $\alpha_i$'s is nonzero (i.e., $\be^\toptzero = \alpha_i\bq_i$), the convergence can be  achieved in only one additional step as well, as illustrated in Figure~\ref{fig:steepest_gd_bisection_eigenvector}.

More specially, when $\lambda_1=\lambda_2=\ldots =\lambda_n=\lambda$, i.e., all the eigenvalues are the same, it then follows that 
$$
\begin{aligned}
\be^\toptone &= \be^\toptzero   - \frac{\bg^\toptzeroTOP \bg^\toptzero}{ \bg^\toptzeroTOP \bA\bg^\toptzero } \bg^\toptzero
=\be^\toptzero - \frac{\sum_{i=1}^{n}\alpha_i^2}{ \sum_{i=1}^{n}\alpha_i^2 }  \be^\toptzero 
= \bzero .
\end{aligned}
$$
Thus, for any arbitrary initial error vector $\be^\toptzero$ (which is always an eigenvector of $\bA$), the algorithm converges to the solution in a single step as well. A specific example is shown in Figure~\ref{fig:steepest_gd_bisection_eigenvector_sameeigenvalue}, where $\bA=\scriptsize\begin{bmatrix}
20 & 0 \\ 0 & 20
\end{bmatrix}$, $\bb=\bzero$, and $c=0$.

\subsubsection{General Convergence Analysis for Symmetric PD Quadratic}\label{section:general-converg-steepest}
To explore general convergence results, we introduce the \textit{energy norm}~\footnote{Also known as the $\bQ$-norm in \eqref{equation:q_norm}.} for the error vector, defined as $\norm{\be}_{\bA} = (\be^\top\bA\be)^{1/2}$. It can be shown that minimizing $\normbig{\be^\toptzero}_{\bA}$ is equivalent to minimizing $f(\bx^\toptzero)$ due to the following relation:
\begin{equation}\label{equation:energy-norm-equivalent}
\norm{\be}_{\bA}^2  = 2f(\bx^\toptzero) \underbrace{- 2f(\bx^*) -2\bb^\top\bx^*}_{\text{constant}}.
\end{equation}
Using  the definition of the energy norm, Equation~\eqref{equation:steepest-quadratic-error}, and the  positive definiteness of $\bA$, we express the update of the energy norm sequence as follows:
\begin{equation}\label{equation:energy-norm-leq}
\small
\begin{aligned}
\normbig{\be^\toptone}_{\bA}^2 
&= (\be^\toptone)^\top \bA\be^\toptone 
= \bigg(\be^\toptzero - \frac{\bg^\toptzeroTOP \bg^\toptzero}{ \bg^\toptzeroTOP \bA\bg^\toptzero } \bg^\toptzero\bigg)^\top \bA \bigg(\be^\toptzero - \frac{\bg^\toptzeroTOP \bg^\toptzero}{ \bg^\toptzeroTOP \bA\bg^\toptzero } \bg^\toptzero\bigg)\\
&=\normbig{\be^\toptzero}_{\bA}^2 + \left(\frac{\bg^\toptzeroTOP \bg^\toptzero}{ \bg^\toptzeroTOP \bA\bg^\toptzero }\right)^2 \bg^\toptzeroTOP \bA\bg^\toptzero 
- 2 \frac{\bg^\toptzeroTOP \bg^\toptzero}{ \bg^\toptzeroTOP \bA\bg^\toptzero } \bg^\toptzeroTOP \bA \be^\toptzero   \\
&\stackrel{\dag}{=}\normbig{\be^\toptzero}_{\bA}^2 - \frac{(\bg^\toptzeroTOP \bg^\toptzero)^2}{ \bg^\toptzeroTOP \bA\bg^\toptzero }   
=\normbig{\be^\toptzero}_{\bA}^2 \cdot\bigg(1- \frac{(\bg^\toptzeroTOP \bg^\toptzero)^2}{ \bg^\toptzeroTOP \bA\bg^\toptzero \cdot \be^\toptzeroTOP \bA\be^\toptzero}\bigg)\\
&\stackrel{\ddag}{=}\normbig{\be^\toptzero}_{\bA}^2  \cdot\left(1- \frac{(\sum_{i=1}^{n}\alpha_i^2\lambda_i^2)^2}{(\sum_{i=1}^{n}\alpha_i^2\lambda_i^3) \cdot (\sum_{i=1}^{n}\alpha_i^2\lambda_i)}\right)
\triangleq\normbig{\be^\toptzero}_{\bA}^2 \cdot r^2, 
\end{aligned}
\end{equation}
where the equality ($\dag$) follows from $\bA\be^\toptzero = \bg^\toptzero$, and   the equality ($\ddag$) follows from \eqref{equation:steepest-et-eigen-decom} and \eqref{equation:steepest-eigen-decom-part}.
The value $r^2 \triangleq\left(1- \frac{(\sum_{i=1}^{n}\alpha_i^2\lambda_i^2)^2}{(\sum_{i=1}^{n}\alpha_i^2\lambda_i^3) \cdot (\sum_{i=1}^{n}\alpha_i^2\lambda_i)}\right)$
determines the rate of convergence. 
By convention, we assume the eigenvalues of $\bA$ are ordered in descending magnitude and they are all positive, i.e.,  $\lambda_1 \geq \lambda_2\geq \ldots \geq \lambda_n>0$, since $\bA$ is positive definite. 
Then the condition number is defined as $\kappa \triangleq \frac{\lambda_1}{\lambda_n}$. 
Additionally, let $\kappa_i \triangleq \frac{\lambda_i}{\lambda_n}$ and $\sigma_i \triangleq \frac{\alpha_i}{\alpha_1}$. 
Then it follows that:
$$
r^2 = 
1- 
\frac{ (\kappa^2+ \sum_{i=\textcolor{mylightbluetext}{2}}^{n}\sigma_i^2\kappa_i^2)^2}
{ (\kappa^3 + \sum_{i=\textcolor{mylightbluetext}{2}}^{n}\sigma_i^2\kappa_i^3) \cdot 
( \kappa +\sum_{i=\textcolor{mylightbluetext}{2}}^{n}\sigma_i^2\kappa_i)}
.
$$
Thus, the rate of convergence is influenced by $\kappa$, $\sigma_i$'s, and $\kappa_i$'s, where $|\kappa_i|\geq  1$ for $i\in \{2,3,\ldots,n\}$. 

\begin{figure}[h!]
\centering                       
\vspace{-0.15cm}                 
\subfigtopskip=2pt               
\subfigbottomskip=2pt             
\subfigcapskip=-5pt              
\subfigure[Demonstration of the rate of convergence in GD with line search in the two-dimensional case. When $\sigma_2=0, \infty$, or $\kappa=1$, the rate of convergence is 0, resulting in immediate convergence within a single step. 
These two cases correspond to $\be^\toptzero$ being an eigenvector of $\bA$ or the eigenvalues being identical, respectively, as illustrated  in Figure~\ref{fig:steepest_specialcases}.]{\label{fig:converge_contour_steepest}
	\includegraphics[width=0.48\linewidth]{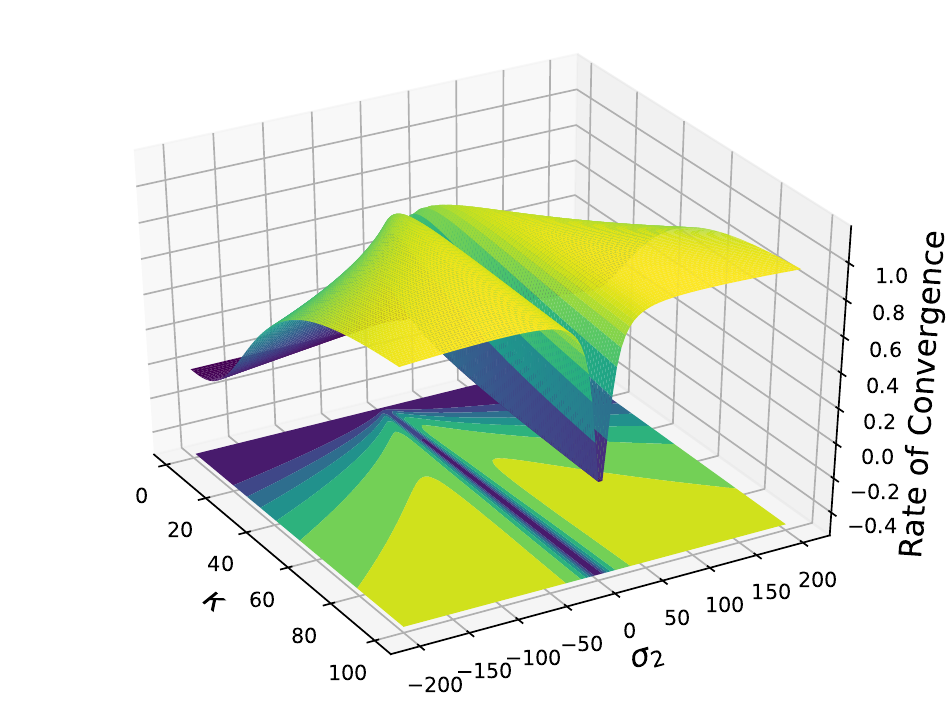}}
\subfigure[Upper bound on the rate of convergence (per iteration) for GD with line search in the two-dimensional case. The $y$-axis represents $\frac{\kappa-1}{\kappa+1}$.]{\label{fig:rate_convergen_steepest}
	\includegraphics[width=0.48\linewidth]{./imgs/rate_convergen_steepest.pdf}}
\caption{Rate of convergence for GD with line search, the lower the better.}
\label{fig:caaa}
\end{figure}

\paragrapharrow{Two-dimensional case.} Specifically, when $n=2$, we have 
\begin{equation}\label{equation:2d-rate-steepest}
r^2=
1- 
\frac{ (\kappa^2+ \sigma_2^2)^2}
{ (\kappa^3 + \sigma_2^2) \cdot 
( \kappa +\sigma_2^2)} .
\end{equation}
Figure~\ref{fig:converge_contour_steepest} illustrates  the value $r^2$ as a function of $\kappa$ and $\sigma_2$. When $n=2$, from \eqref{equation:steepest-et-eigen-decom}, we have 
\begin{equation}\label{equation:steepest-et-eigen-decom-2d}
\be^\toptzero = \alpha_1 \bq_1+\alpha_2\bq_2. 
\end{equation}
This confirms the two special examples shown in Figure~\ref{fig:steepest_specialcases}: when $\be^\toptzero$ is an eigenvector of $\bA$, it follows that:
$$
\begin{aligned}
\text{case 1: }	&\alpha_2=0 \qquad\implies\qquad  \sigma_2 = \alpha_2/\alpha_1 \rightarrow 0;\\
\text{case 2: }	&\alpha_1=0 \qquad\implies\qquad  \sigma_2 = \alpha_2/\alpha_1 \rightarrow \infty.
\end{aligned}
$$
That is, the slope of $\sigma_2$ is either zero or infinite, the rate of convergence approaches  zero and it converges instantly in just one step (example in Figure~\ref{fig:steepest_gd_bisection_eigenvector}). Similarly, when the eigenvalues are identical ($\kappa=1$),  the rate of convergence is also zero (example in Figure~\ref{fig:steepest_gd_bisection_eigenvector_sameeigenvalue}).

\index{Rate of convergence}

\paragrapharrow{Worst case.} We recall that $\sigma_2 = \alpha_2/\alpha_1$ determines the error vector $\be^\toptzero$ (Equation~\eqref{equation:steepest-et-eigen-decom} or Equation~\eqref{equation:steepest-et-eigen-decom-2d}), which, in turn, determines the point $\bx^\toptzero$ in the two-dimensional case. 
A natural question is: what is the worst possible starting point for convergence?
For a fixed $\kappa$ (i.e., $\bA$ and the loss function $f(\bx)$ are fixed), assuming $t=1$,  we seek the initial point $\bx^{(1)}$ that results in the slowest descent. 
It can be shown that the rate of convergence in \eqref{equation:2d-rate-steepest} is maximized when $\sigma_2=\pm \kappa$:
$$
\begin{aligned}
r^2 &\leq 1-\frac{4 \kappa^2}{\kappa^5+2\kappa^4+\kappa^3}= \frac{(\kappa-1)^2}{(\kappa+1)^2}.
\end{aligned}
$$
Substituting into \eqref{equation:energy-norm-leq}, we obtain 
$$
\begin{aligned}
\normbig{\be^\toptone}_{\bA}^2 &\leq \normbig{\be^\toptzero}_{\bA}^2 \cdot \frac{(\kappa-1)^2}{(\kappa+1)^2}
\quad\implies\quad
\normbig{\be^\toptone}_{\bA} \leq \normbig{\be^{(1)}}_{\bA} \cdot\left(\frac{\kappa-1}{\kappa+1}\right)^t.
\end{aligned}
$$
The upper bound on the convergence rate per iteration is shown in Figure~\ref{fig:rate_convergen_steepest}. 
As expected, the more \textit{ill-conditioned} the matrix, the slower the convergence of gradient descent with line search. 
Notably, the upper bound on the convergence rate matches that of  GD with a fixed stepsize; see \eqref{equation:vanialla-gd-rate}. 
However, a key difference is that the rate of GD with a fixed stepsize in \eqref{equation:vanialla-gd-rate} is obtained by selecting a specific learning rate, as shown in \eqref{equation:eta-vanilla-gd}. 
This makes the rate of GD with a fixed stepsize more of  a tight bound. In practice, GD with a fixed stepsize converges slower than GD with line search, as evident in the examples shown in Figure~\ref{fig:momentum_gd_conjugate2}, Figure~\ref{fig:momentum_gd_conjugate8}, and Figure~\ref{fig:conjguatecy_zigzag2}.

\index{Spectral decomposition}
\index{Eigenvalue decomposition}
\index{Quadratic form}
\section{Quadratic Model using Momentum}\label{section:quadratic-in-momentum}
We have discussed the gradient descent method with momentum in \eqref{equation:gd_with_momentum}.
Following the discussions of the quadratic form in GD with a fixed stepsize (Section~\ref{section:quadratic_vanilla_GD})  and GD with line search (Section~\ref{section:quadratic-in-steepestdescent}), we now examine  the quadratic form in GD with momentum. The general update is:
$$
\begin{aligned}
\bd^\toptzero &= \rho\bd^\toptminus - \eta_t\nabla f(\bx^\toptzero); \\
\bx^\toptone &= \bx^\toptzero + \bd^\toptzero,
\end{aligned}
$$
where $\nabla f(\bx^\toptzero) = \bA\bx^\toptzero - \bb$ under the assumption that  $\bA$ is positive definite for the quadratic form. 
Fixing the learning rate $\eta_t = \eta$, the update simplifies to:
$$
\begin{aligned}
\bd^\toptzero &= \rho\bd^\toptminus - \eta(\bA\bx^\toptzero - \bb); \\
\bx^\toptone &= \bx^\toptzero + \bd^\toptzero.
\end{aligned}
$$
Again, define the iterate vectors as follows:
$$
\left\{
\begin{aligned}
\by^\toptzero &\triangleq \bQ^\top(\bx^\toptzero - \bx^*);\\
\bz^\toptzero &\triangleq \bQ^\top  \bd^\toptzero,
\end{aligned}
\right.
$$
where $\bx^* = \bA^{-1}\bb$ under the assumption that $\bA$ is positive definite, and $\bA=\bQ\bLambda\bQ^\top$ is the spectral decomposition of matrix $\bA$. This construction leads to the following update rule:
$$
\begin{aligned}
\bz^\toptzero &= \rho \bz^\toptminus - \eta\bLambda \by^\toptzero; \\
\by^\toptone &= \by^\toptzero + \bz^\toptzero,
\end{aligned}
\quad\implies\quad
\begin{bmatrix}
\bz^\toptzero \\
\by^\toptone
\end{bmatrix}
= 
\begin{bmatrix}
\rho \bI & -\eta \bLambda \\
\rho \bI & -\eta \bLambda + \bI
\end{bmatrix}
\begin{bmatrix}
\bz^\toptminus \\
\by^\toptzero
\end{bmatrix}.
$$
Alternatively, this results in the per-dimension update:
\begin{equation}\label{equation:momentum-quadra-generalformula}
\begin{bmatrix}
z_{i}^\toptzero \\
y_{i}^\toptone
\end{bmatrix}
= 
\begin{bmatrix}
\rho  & -\eta\lambda_i \\
\rho  & 1-\eta \lambda_i 
\end{bmatrix}^t
\begin{bmatrix}
z_{i}^{(0)} \\
y_{i}^{(1)}
\end{bmatrix}
\triangleq
\bB^t
\begin{bmatrix}
z_{i}^{(0)} \\
y_{i}^{(1)}
\end{bmatrix},
\end{equation}
where $z_{i}^\toptzero$ and $y_{i}^\toptzero$ are $i$-th element of $\bz^\toptzero$ and $\by^\toptzero$, respectively, and 
$\bB\triangleq
\footnotesize
\begin{bmatrix}
\rho  & -\eta\lambda_i \\
\rho  & 1-\eta \lambda_i 
\end{bmatrix}$.
Note here, $\bz^{(0)}$ is initialized as a zero vector, and $\by^{(1)}$ is initialized as $\bQ^\top (\bx^{(1)}-\bx^*)$, where $\bx^{(1)}$ represents the initial parameter.
Suppose the eigenvalue decomposition (Theorem~\ref{theorem:eigenvalue-decomposition}) of $\bB$ admits
$$
\bB = \bC\bD \bC^{-1},
$$
where the columns of $\bC$ contain eigenvectors of $\bB$, and $\bD=\diag(\alpha,\beta)$ is a diagonal matrix  containing the eigenvalues of $\bB$. Then $\bB^t = \bC\bD^t\bC^{-1}$ (Remark~\ref{remark:power-eigenvalue-decom}). 
Alternatively, the eigenvalues of $\bB$ can be computed by solving $\det(\bB-\alpha\bI)=0$:
$$
\alpha, \beta = \frac{(\rho+1-\eta\lambda_i) \pm \sqrt{(\rho+1-\eta\lambda_i)^2 -4\rho}}{2}.
$$
We then have by \citet{williams1992n} that
$$
\bB^t=
\left\{
\begin{aligned}
&\alpha^t \frac{\bB-\beta\bI}{\alpha-\beta} - \beta^t \frac{\bB-\alpha\bI}{\alpha-\beta}, \gap &\text{if $\alpha\neq \beta$};\\
&\alpha^{t-1}(t\bB - (t-1)\alpha\bI), \gap &\text{if $\alpha=\beta$}.
\end{aligned}
\right.
$$
Substituting this into \eqref{equation:momentum-quadra-generalformula} yields the following expression:
$$
\begin{bmatrix}
z_{i}^\toptzero \\
y_{i}^\toptone
\end{bmatrix}
=
\bB^{t}
\begin{bmatrix}
z_{i}^{(0)} \\
y_{i}^{(1)}
\end{bmatrix},
$$
where the rate of convergence is controlled by the slower one, $\max\{\abs{\alpha}, \abs{\beta}\}$. When we have $\max\{\abs{\alpha}, \abs{\beta}\}<1$, the GD with momentum is guaranteed to converge. 
In the case of  $\rho=0$, the momentum reduces to the basic GD, the condition for convergence becomes
$$
\max\{\abs{\alpha}, \abs{\beta}\} = \abs{1-\eta\lambda_i} <1, \gap \forall \,\, i\in \{1,2,\ldots, n\},
$$
which aligns with  that in \eqref{equation:vanillagd-quandr-rate-chgoices}.

\begin{figure}[h]
\centering   
\vspace{-0.35cm}  
\subfigtopskip=2pt  
\subfigbottomskip=2pt  
\subfigcapskip=-5pt  
\subfigure[Momentum $\rho=0.2$, convergence rate$\approx0.79$.]{\label{fig:momentum_rho2}
\includegraphics[width=0.31\linewidth]{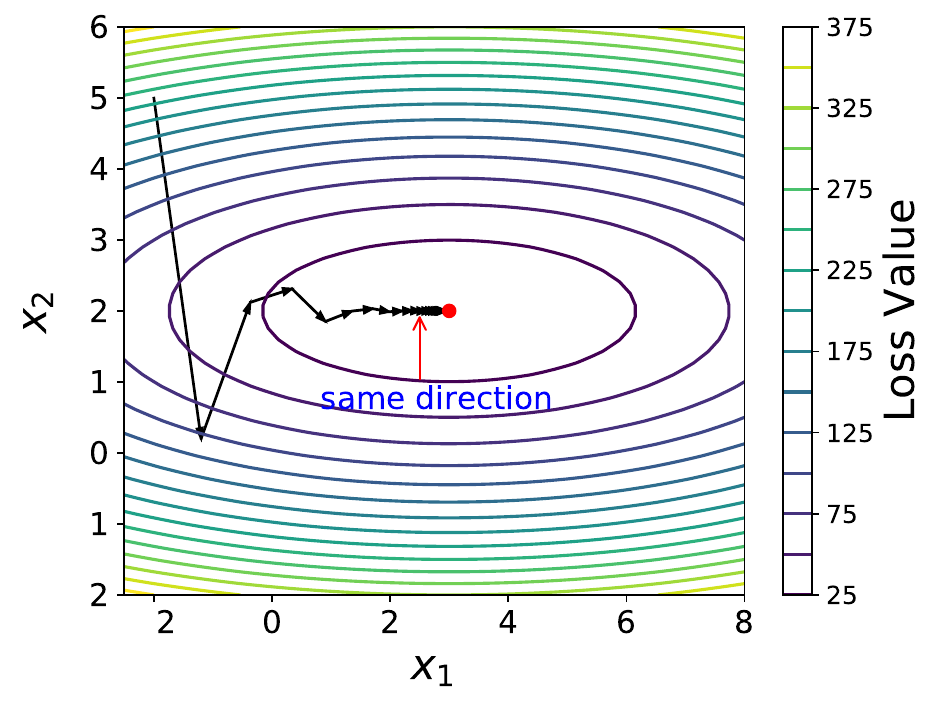}}
\subfigure[Momentum $\rho=0.8$, convergence rate$\approx0.89$.]{\label{fig:momentum_rho8}
\includegraphics[width=0.31\linewidth]{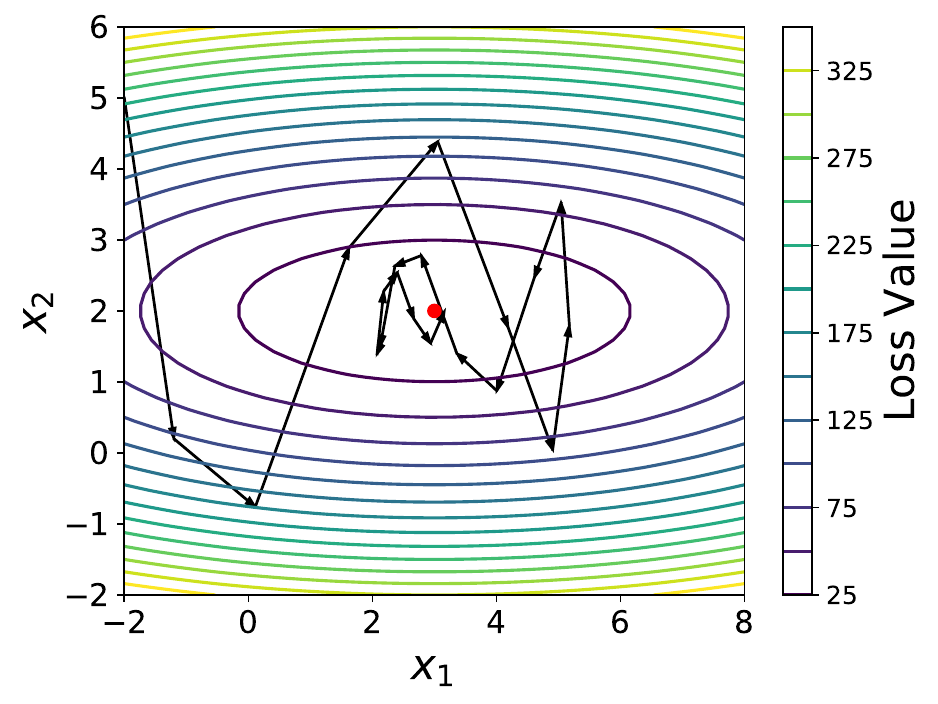}}
\subfigure[Momentum $\rho=1$, convergence rate=1.]{\label{fig:momentum_rho10}
\includegraphics[width=0.31\linewidth]{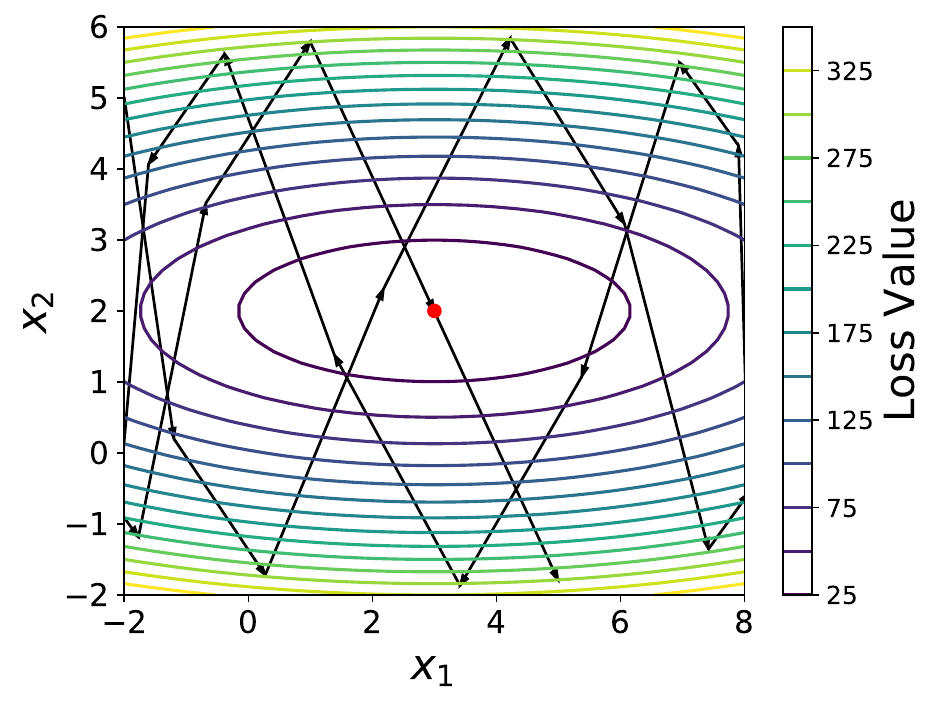}}
\caption{Momentum creates its own oscillations. The learning rates $\eta$ are set to be 0.04 for all scenarios.}
\label{fig:momentum_own_oscilation}
\end{figure}

As an example, consider the quadratic form with $\bA=\footnotesize\begin{bmatrix}
4 & 0\\ 0 & 40
\end{bmatrix}$, which has eigenvalues  $\lambda_1=4$ and $\lambda_2=40$, respectively. Then, the corresponding matrix $\bB$ in \eqref{equation:momentum-quadra-generalformula} is 
$$
\bB_1 = 
\begin{bmatrix}
\rho & -4 \eta\\
\rho & 1-4\eta
\end{bmatrix}
\gap \text{or}\gap 
\bB_2 = 
\begin{bmatrix}
\rho & -40\eta \\
\rho & 1-40\eta
\end{bmatrix}.
$$
It can be shown that when $\eta=0.04, \rho=0.8$, the rate of convergence is approximately  $0.89$; Figure~\ref{fig:momentum_rho8} displays the updates for 20 iterations. 
Although the trajectory exhibits a zigzag pattern, the method still converges. 
However, when $\eta=0.04, \rho=1$, the rate of convergence is equal to 1, indicating no actual progress toward convergence; Figure~\ref{fig:momentum_rho10} shows the updates for 20 iterations, where the trajectory diverges, even though it intermittently passes through the optimal point.

\begin{problemset}
\item \textbf{Constructing symmetric form.} Let $P(\bA)\triangleq\frac{1}{2}(\bA+\bA^\top)$ if $\bA\in\real^{n\times n}$. Show that 
\begin{itemize}
	\item \textit{Null space.} $\nspace(\bA)\subset \nspace(P(\bA))$ and $\nspace(\bA^\top)\subset \nspace(P(\bA))$ such that $\rank(P(\bA))\leq \rank(\bA)$.
	\item When $\rank(P(\bA))= \rank(\bA)$, it follows that $\bA$, $\bA^\top$, and $P(\bA)$ have the same null space.
\end{itemize}
\textit{Hint: Consider the quadratic form $\bx^\top\bA\bx$ and $\bx^\top\bA^\top\bx$.}

\item Let $\bA,\bB\in\real^{n\times n}$ be  positive semidefinite matrices. Show that $\bA+\bB$ is also positive semidefinite.
\item Let $\bA\in\real^{n\times n}$ and $\bB\in\real^{m\times m}$ be  two symmetric matrices. Show that the following two claims are equivalent:
\begin{enumerate}
	\item $\bA$ and $\bB$ are positive semidefinite.
	\item $\scriptsize\begin{bmatrix}
		\bA& \bzero \\
		\bzero & \bB
	\end{bmatrix}$ is positive semidefinite.
\end{enumerate}
\item 
Given a matrix $\bB\in\real^{n\times k}$ and let $\bA=\bB\bB^\top$. Show that 
\begin{enumerate}
	\item $\bA$ is positive semidefinite.
	\item $\bA$ is positive definite if and only if $\bB$ has full row rank.
\end{enumerate}
\item Show that if $\bA$ is positive semidefinite and nonsingular, then $\bA^{-1}$ is positive definite.
\end{problemset}

\newpage
\chapter{First-Order Methods}\label{chapter:gd_convg}
\begingroup
\hypersetup{linkcolor=structurecolor,
	linktoc=page,  
}
\minitoc \newpage
\endgroup

\section{Background and Mathematical Tools}

First-order methods are essential in large-scale optimization, particularly in the fields of machine learning and deep learning. This chapter explores advanced gradient-based techniques, including projected gradient descent, proximal gradient methods, and mirror descent. Additionally, it delves into acceleration strategies such as Nesterov’s accelerated gradient and adaptive approaches like the conditional gradient method. These methodologies form the backbone of modern optimization techniques.

Our primary focus will be on the convergence results of algorithms for problems defined in Definition~\ref{definition:opt_probs_all}. Specifically, we will examine:
\begin{subequations}\label{equation:p1p2}
\begin{align}
\text{(P1)}:\qquad &\min \{f(\bx)\};\\
\text{(P2)}:\qquad &\min \{f(\bx)\} \quad\text{s.t.}\quad \bx\in\sS;\\
\text{(P3)}:\qquad &\min \left\{F(\bx)\triangleq f(\bx)+g(\bx) \right\}.
\end{align}
\end{subequations}
To be more specific, we will consider the \textit{gradient descent (GD) method} for (P1); the \textit{projected gradient descent (PGD)}, \textit{mirror descent}, and \textit{conditional gradient (CG) descent} methods for (P2); the \textit{proximal gradient}, the \textit{G-mirror descent}, and the \textit{generalized conditional gradient (GCG)} methods   for (P3):

When seeking a maximizer of a function, one can simply find a minimizer of the function with opposite sign. In this context, problem (P1) is referred to \textit{unconstrained optimization}, (P2) as \textit{constrained optimization}, and (P3) as a non-convex composite problem. 
The function $f$ or $F$ is called the \textit{objective function} (a.k.a.,  \textit{cost function} or \textit{loss function}), while we use  $\bx^*$ to represent  the \textit{minimum point (or simply minimizer)}.
For (P3), different forms of $g$ lead to various types of optimization problems:
\begin{itemize}
	\item When $g(\bx)=0$,  problem (P3) reduces to an unconstrained problem (P1).
	\item When $g(\bx)=\indicatorS(\bx)$,  problem (P3) becomes a constrained problem (P2).
	\item When $g(\bx)=\lambda\normone{\bx}$,  problem (P3) turns into an $\ell_1$ regularized problem. 
\end{itemize}

\begin{figure}[H]
\centering
\begin{widepage}
\centering
\resizebox{0.6\textwidth}{!}{%
\begin{tikzpicture}[>=latex]

\tikzstyle{state} = [draw, very thick, fill=white, rectangle, minimum height=3em, minimum width=7.8em, node distance=3em, font={\sffamily\bfseries}]
\tikzstyle{stateEdgePortion} = [black,thick];
\tikzstyle{stateEdge} = [stateEdgePortion,->];
\tikzstyle{stateEdge2} = [stateEdgePortion,<->];
\tikzstyle{edgeLabel} = [pos=0.5, text centered, font={\sffamily\small}];

\node[state, name=cg, node distance=7em, xshift=-9em, yshift=-15em, fill={colorlu}] {CG (P2)};
\node[state, name=gcg, fill={colorlu}, right=of cg, xshift=18em] {GCG (P3)};
\node[state, name=gd, fill={colorlu}, above=of cg, yshift=0.1em, xshift=0em] {\parbox{5.8em}{GD $\;\;$(P1)\\Line Search}};
\node[state, name=pgd, fill={colorlu}, above=of cg, yshift=0.1em, xshift=14em] {PGD (P2)};
\node[state, name=prox, fill={colorlu}, right=of gd, xshift=18em] {\parbox{6.6em}{Proximal $\;\;$(P3)\\FISTA/Nesterov}};
\node[state, name=mirror, draw, fill={colorlu}, above=of gd, yshift=0.2em, xshift=0em] {Mirror (P2)};
\node[state, name=gmirror, draw, fill={colorlu}, right=of mirror, xshift=18em] {G-Mirror (P3)};

\draw (mirror.east) 
edge[stateEdge] node[edgeLabel, xshift=0.5em, yshift=-0.8em]{Composite} 
($(gmirror.west) + (0em,0em)$);

\draw (pgd.east)
edge[stateEdge] node[edgeLabel, xshift=-0.1em, yshift=0.5em,]{Composite} 
(prox.west) ;

\draw ($(cg.east)$)
edge[stateEdge] node[edgeLabel, xshift=0.5em, yshift=0.5em]{Composite} 
(gcg.west) ;

\draw ($(pgd.west) + (0em,0em)$)
edge[stateEdge, bend left=0] node[edgeLabel, xshift=-0.0em,yshift=0.5em]{Unconstrained} 
(gd.east) ;
\draw ($(pgd.north) + (0em,0em)$)
edge[stateEdge, bend left=-12.5] node[edgeLabel, xshift=1.5em,yshift=-0em]{Non-Euclidean Distance} 
($(mirror.south) + (2.5em,0em)$) ;
\draw ($(pgd.south)$)
edge[stateEdge, bend left=+12.5] node[edgeLabel, xshift=1em, yshift=0em]{Convex Combination Update} 
($(cg.north) + (2.5em,0em)$) ;

\draw ($(prox.north) + (0em,0em)$)
edge[stateEdge] node[edgeLabel, xshift=-1.5em,yshift=0em]{Non-Euclidean Distance} 
(gmirror.south) ;
\draw ($(prox.south)$)
edge[stateEdge] node[edgeLabel, xshift=-1em, yshift=0em]{Convex Combination Update} 
(gcg.north) ;

\begin{pgfonlayer}{background}
\draw [join=round,cyan,dotted,fill={colormiddle}] ($(cg.south west -| mirror.west) + (-0.5em, -0.5em)$) rectangle ($(mirror.north east -| prox.north east) + (0.6em, 0.5em)$);		
\end{pgfonlayer}

\end{tikzpicture}
}
\end{widepage}
\caption{
Relationship and classification of first-order optimization methods can be understood in the context of different problem types, specifically (P1), (P2), and (P3). The diagram depicts an evolution from fundamental techniques, such as PGD, which are suitable for simpler problem formulations, to more sophisticated approaches like proximal gradient methods and GCG, designed for handling composite and constrained problems.
At the core of these methodologies lies PGD, which serves as a foundational algorithm. Advanced methods build upon it by incorporating concepts such as convex combination updates and non-Euclidean metrics, which enhance their ability to address complex optimization challenges.
}
\label{fig:first_order_opt_map}
\end{figure}
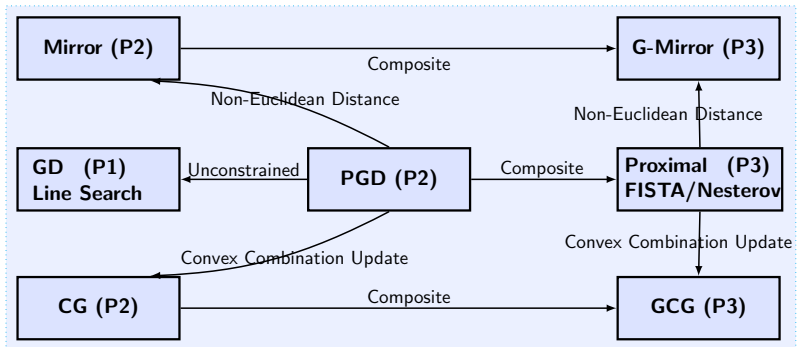

Part of the relationship among these methods is illustrated in Figure~\ref{fig:first_order_opt_map}. The diagram shows an evolution from fundamental techniques, such as PGD, which are suitable for simpler problem formulations, to more sophisticated approaches like proximal gradient methods and GCG, designed for handling composite and constrained problems.
At the core of these methodologies lies PGD, serving as a foundational algorithm. Advanced methods build upon PGD by incorporating concepts such as convex combination updates and non-Euclidean metrics, enhancing their ability to address complex optimization challenges. For this reason, we begin our analysis with PGD and its unconstrained counterpart, GD.

\section*{Descent Lemmas}
We start by proving the descent lemmas for GD and PGD.
Previously, we introduced the descent inequality under the Zoutendijk condition in the proof of Theorem~\ref{theorem:zoutendijk_cond}: $f(\bx^\toptone) -f(\bx^\toptzero) \leq  c_1 \frac{c_2 - 1}{\beta} \cos^2 (\theta_t) \normtwobig{\nabla f(\bx^\toptzero)}^2$, where $0 < c_1 < c_2 < 1$.
This kind of descent inequality is crucial for proving the convergence or the convergence rate of the underlying algorithm.
Analogously, definitions of smoothness and strong convexity lead to descent lemmas for GD or PGD updates: 
$
\bx^+\leftarrow \bx-\eta\nabla f(\bx)
 \text{ or } 
\bx^+ \leftarrow \projectS\left(\bx - \eta \nabla f(\bx)\right),
$
where $\sS$ denotes a set and $\projectS(\cdot)$ denotes the projection operation (Definition~\ref{definition:projec_prox_opt}).
\noindent
\begin{lemma}[Descent Lemmas for Problems (P1) and (P2)]\label{lemma:gdupd_sm_sconv}
For a differentiable function $f$, we  consider the gradient descent update rule $\bx^+\leftarrow \bx-\eta\nabla f(\bx)$ or the projected gradient descent update rule $\bx^+ \leftarrow \projectS\left(\bx - \eta \nabla f(\bx)\right)$ with a constant stepsize $\eta$. 
Then,
\begin{enumerate}[(i)]
\item \label{descent_lem_gdss} \textit{GD for SS.}  Let $f: \real^n \rightarrow \real$ be a differentiable  and $\beta$-smooth function. Then the gradient descent update $\bx^+\leftarrow \bx-\frac{1}{\beta}\nabla f(\bx)$  ensures that:
$$
\textbf{(SS)}:\qquad  f(\bx^+)-f(\bx) \leq -\frac{1}{2\beta} \normtwo{\nabla f(\bx)}^2.
$$
This indicates that the gradient update decreases the function value by an amount proportional to the squared norm of the gradient.~\footnote{See Theorem~\ref{theorem:pgd_smooth} for more insights.}
\item \label{descent_lem_gdscss} \textit{GD for SC and SS.} Let $f: \real^n \rightarrow \real$ be a differentiable, $\alpha$-strongly convex, and $\beta$-smooth function. Then for any $\bx, \by \in \real^n$ and an update of the form $\bx^+ \leftarrow \bx - \frac{1}{\beta} \nabla f(\bx)$:
$$
\textbf{(SC+SS)}:\qquad f(\bx^+) - f(\by) \leq \nabla f(\bx)^\top (\bx - \by) - \frac{1}{2\beta} \normtwo{\nabla f(\bx)}^2 - \frac{\alpha}{2} \normtwo{\bx - \by}^2.~\footnote{In practice, $\by$ can be set to the minimizer $\bx^*$ of the function; see Theorem~\ref{theorem:gd_sc_ss} for more insights.}
$$

\item \label{descent_lem_pgdss} \textit{PGD for  SS.} Let $f: \sS\rightarrow \real$ be a differentiable and $\beta$-smooth function, where $\sS\subseteq\real^n$ is a convex set. Then for any $ \bx, \by \in \sS$, an update $\bx^+ \in \sS$ of the form  $\bx^+ \leftarrow \projectS\left(\bx - \frac{1}{\beta} \nabla f(\bx)\right)$~\footnote{See Algorithm~\ref{alg:pgd_gen}.}, and the function $g: \sS \rightarrow \real^n$ defined as $g(\bx) \triangleq \beta(\bx - \bx^+)$ (i.e., $\bx^+ =\bx-\frac{1}{\beta}g(\bx)$):
$$
\textbf{(PGD+SS)}:\qquad f(\bx^+) - f(\by) \leq g(\bx)^\top (\bx - \by) - \frac{1}{2\beta} \normtwo{g(\bx)}^2.
$$

\item \label{descent_lem_pgdscss} \textit{PGD for SC and SS.} Let $f: \sS\rightarrow \real$ be a differentiable, $\alpha$-strongly convex, and $\beta$-smooth function, where $\sS\subseteq\real^n$ is a convex set. Then for any $ \bx, \by \in \sS$, an update $\bx^+ \in \sS$ of the form $\bx^+ \leftarrow \projectS\left(\bx - \frac{1}{\beta} \nabla f(\bx)\right)$, and the function $g: \sS \rightarrow \real^n$ defined as $g(\bx) \triangleq \beta(\bx - \bx^+)$ (i.e., $\bx^+ =\bx-\frac{1}{\beta}g(\bx)$):
$$
\textbf{(PGD+SC+SS)}:\qquad f(\bx^+) - f(\by) \leq g(\bx)^\top (\bx - \by) - \frac{1}{2\beta} \normtwo{g(\bx)}^2 - \frac{\alpha}{2} \normtwo{\bx - \by}^2.
$$
When $\sS=\real^n$, we have $g(\bx)=\nabla f(\bx)$ and this reduces to (ii); when $\alpha=0$, it reduces to (iii); when $\sS=\real^n$, $\by=\bx$, and $\alpha=0$, this reduces to (i).
\end{enumerate}
\end{lemma}
\begin{proof}[of Lemma~\ref{lemma:gdupd_sm_sconv}]
Note that the first part follows directly from the definition of smoothness. 
Since the first three parts are special cases of \ref{descent_lem_pgdscss}, we only prove \ref{descent_lem_pgdscss} here.
By Projection Property-I (Lemma~\ref{lemma:proj_prop1}), it follows that $\innerproduct{\bx^+-(\bx-\frac{1}{\beta}\nabla f(\bx)), \bx^+-\by}\leq 0$. This implies
\begin{equation}\label{equation:gdupd_sm_sconv1}
\nabla f(\bx)^\top (\bx^+ - \by) \leq g(\bx)^\top (\bx^+ - \by).
\end{equation}
Thus, by the smoothness and strong convexity (Definition~\ref{definition:scss_func}), 
$$
\begin{aligned}
f(\bx^+)  - f(\by)  
&= \left(f(\bx^+) - f(\bx)\right) + \left(f(\bx) - f(\by)\right) \\
&\leq \left(\nabla f(\bx)^\top (\bx^+ - \bx) + \frac{\beta}{2} \normtwo{\bx^+ - \bx}^2\right)  + 
\left(\nabla f(\bx)^\top (\bx - \by) - \frac{\alpha}{2} \normtwo{\bx - \by}^2\right) \\
&\leq   g(\bx)^\top (\bx^+ - \by) + \frac{1}{2\beta} \normtwo{g(\bx)}^2 - \frac{\alpha}{2} \normtwo{\bx - \by}^2 \\
&=  g(\bx)^\top (\bx - \by) - \frac{1}{2\beta} \normtwo{ g(\bx)}^2 - \frac{\alpha}{2} \normtwo{\bx - \by}^2.
\end{aligned}
$$
This completes the proof.
\end{proof}
In the above proof, since $\bx\in\sS$,  by Projection Property-II (Lemma~\ref{lemma:proj_prop2}) and letting $\bz\triangleq \bx-\frac{1}{\beta}\nabla f(\bx)$, we can also show that  that
\begin{equation}\label{equation:gd_des_cor}
\normtwo{\bx^+ - \bx} \leq \normtwo{\bz-\bx}
\quad\implies\quad 
\normtwo{g(\bx)} \leq \normtwo{\nabla f(\bx)}.
\end{equation}

\begin{algorithm}[h] 
\caption{Gradient Descent Method}
\label{alg:gd_gen}
\begin{algorithmic}[1] 
\Require A function $f(\bx)$; 
\State {\bfseries Input:}  Initialize $\bx^{(1)}$;
\For{$t=1,2,\ldots$}
\State Pick a stepsize $\eta_t$;
\State $\bx^{(t+1)} \leftarrow \bx^{(t)} - \eta_t \nabla f(\bx^{(t)})$ or $\bx^{(t+1)} \leftarrow \bx^{(t)} - \eta_t f^\prime(\bx^{(t)})$, where $ f^\prime(\bx^{(t)})\in \partial f(\bx^{(t)})$;
\EndFor
\State (Output Option 1) Output  $\bx_{\text{final}}\leftarrow \bx^{(T)}$;
\State (Output Option 2) Output  $\bx_{\text{avg}}\leftarrow \frac{1}{T}(\sum_{t=1}^{t}\bx^{(t)})$ or $\sum_{t=1}^{T} \frac{2t}{T(T+1)} \bx^{(t)}$;
\State (Output Option 3) Output  $\bx_{\text{best}}\leftarrow \argmin_{t\in\{1,2,\ldots,T\}} f(\bx^{(t)})$;
\end{algorithmic} 
\end{algorithm}

\section{Gradient Descent}\label{section:gd_basic}
For the unconstrained problem (P1) given in \eqref{equation:p1p2}, we consider the gradient descent method as described in Algorithm~\ref{alg:gd_gen}~\footnote{Although convergence results are also applicable to non-differentiable functions using subgradient descent, for brevity, we will focus solely on differentiable cases in this section.}, i.e., using negative gradient direction (steepest descent direction) for the descent method in Algorithm~\ref{alg:struc_gd_gen}.
For convex and smooth functions, the gradient descent method achieves a convergence rate of $\mathcalO(1/T)$.
\begin{theoremHigh}[GD for Convex and SS Functions: $\mathcalO(1/T)$]\label{theorem:pgd_smooth}
Let $f:\real^n\rightarrow \real$ be a  differentiable, convex, and  $\beta$-smooth function defined on $\real^n$. 
Suppose $\{\bx^\toptzero\}_{t > 0}$ is the sequence generated by the gradient descent method (Algorithm~\ref{alg:gd_gen}) for solving  problem (P1)
with a constant stepsize $\eta = \frac{1}{\beta}$, and let $\bx^*$ be an optimizer of $f$. 
Then, for any $T>0$,
$$
f(\bx^{(T)}) - f(\bx^*) \leq \frac{2\beta \normtwo{\bx^{(1)} -\bx^*}^2}{T-1}.
$$
\end{theoremHigh}
\begin{proof}[of Theorem~\ref{theorem:pgd_smooth}]
Let $\delta_t \triangleq f(\bx^{(t)}) - f(\bx^*)$. 
By convexity (Theorem~\ref{theorem:conv_gradient_ineq}), it follows that
$$
\delta_t \leq \nabla f(\bx^{(t)})^\top (\bx^{(t)} -\bx^*) \leq \normtwobig{\bx^{(t)} -\bx^*} \cdot \normtwobig{\nabla f(\bx^{(t)})}.
$$
By the update rule and the definition of smoothness (Lemma~\ref{lemma:gdupd_sm_sconv}), we have
$$
f(\bx^{(t+1)}) - f(\bx^{(t)}) \leq -\frac{1}{2\beta} \normtwobig{\nabla f(\bx^{(t)})}^2
\quad\implies\quad 
\delta_{t+1} \leq \delta_t - \frac{1}{2\beta} \normtwobig{\nabla f(\bx^{(t)})}^2.
$$
Combining the above two inequalities yields 
\begin{equation}\label{equation:pgd_smooth1}
\delta_{t+1} \leq \delta_t - \frac{\delta_t^2}{2\beta\normtwobig{\bx^{(t)} -\bx^*}^2}.
\end{equation}
Since the gradient at stationary points vanishes: $\nabla f(\bx^*) = \bzero$, we have
$$
\begin{aligned}
\normtwobig{\bx^{(t+1)} -\bx^*}^2 
&= \normtwo{\bx^{(t)} - \frac{1}{\beta} \nabla f(\bx^{(t)}) -\bx^*}^2 \\
&= \normtwobig{\bx^{(t)} -\bx^*}^2 - \frac{2}{\beta} \nabla f(\bx^{(t)})^\top (\bx^{(t)} -\bx^*) + \frac{1}{\beta^2} \normtwobig{\nabla f(\bx^{(t)})}^2 \\
&\stackrel{\dag}{\leq} \normtwobig{\bx^{(t)} -\bx^*}^2 - \frac{1}{\beta^2} \normtwobig{\nabla f(\bx^{(t)})}^2 
\leq \normtwobig{\bx^{(t)} -\bx^*}^2.
\end{aligned}
$$
where the inequality ($\dag$) follows from \eqref{equ:charac_smoo_3} in Theorem~\ref{theorem:charac_smoo}.
This shows $\normtwobig{\bx^{(t)} -\bx^*}$ is decreasing with $t$, which with \eqref{equation:pgd_smooth1} implies
$
\delta_{t+1} \leq \delta_t - \frac{\delta_t^2}{2\beta \normtwo{\bx^{(1)} - \bx^*}^2} .
$
Let $z \triangleq \frac{1}{2\beta \normtwo{\bx^{(1)} -\bx^*}^2}$. Then,
$$
z\delta_t^2 + \delta_{t+1} \leq \delta_t 
\iff 
z\frac{\delta_t}{\delta_{t+1}} + \frac{1}{\delta_t} \leq \frac{1}{\delta_{t+1}} 
\implies 
z  \leq \frac{1}{\delta_{t+1}} - \frac{1}{\delta_t} 
\implies 
z(T-1) \leq \frac{1}{\delta_T},
$$
where the last inequality follows by performing telescopic cancellations.
\end{proof}

\begin{theoremHigh}[GD for SC and SS Functions]\label{theorem:gd_sc_ss}
Let $f: \real^n \rightarrow \real$ be a differentiable, $\alpha$-strongly convex, and $\beta$-smooth function. 
Suppose $\{\bx^\toptzero\}_{t > 0}$ is the sequence generated by the gradient descent method (Algorithm~\ref{alg:gd_gen}) for solving  problem (P1),
and let $\bx^*$ be an optimizer of $f$. If the constant stepsize is $\eta=\frac{1}{\beta}$, then,
$$
\normtwo{\bx^{(t+1)} - \bx^*}^2 \leq \exp\left(-t\frac{\alpha}{\beta}\right) \normtwo{\bx^{(1)} - \bx^*}^2.
$$
For the function value, if the constant stepsize is  $\eta = \frac{2}{\alpha + \beta}$, then
$$
\begin{aligned}
\normtwo{\bx^{(t+1)} - \bx^*}^2 &\leq\exp\left(-\frac{4t}{\alpha + 1}\right) \normtwo{\bx^{(1)} - \bx^*}^2;\\
f(\bx^{(t+1)}) - f(\bx^*)& \leq \frac{\beta}{2} \exp\left(-\frac{4t}{\alpha + 1}\right) \normtwo{\bx^{(1)} - \bx^*}^2.
\end{aligned}
$$
\end{theoremHigh}
\begin{proof}[of Theorem~\ref{theorem:gd_sc_ss}]
\textbf{When $\eta=\frac{1}{\beta}$.}
Using the update rule,
$$
\small
\begin{aligned}
&\normtwo{\bx^{(t+1)} - \bx^*}^2 
= \normtwo{\bx^{(t)} - \frac{1}{\beta} \nabla f(\bx^{(t)}) - \bx^*}^2\\
&= \normtwo{\bx^{(t)} - \bx^*}^2 - \frac{2}{\beta} \nabla f(\bx^{(t)})^\top (\bx^{(t)} - \bx^*) + \frac{1}{\beta^2} \normtwobig{\nabla f(\bx^{(t)})}^2 \\
&\stackrel{\dag}{\leq} \left(1 - \frac{\alpha}{\beta}\right) \normtwo{\bx^{(t)} - \bx^*}^2 
\leq \left(1 - \frac{\alpha}{\beta}\right)^t \normtwo{\bx^{(1)} - \bx^*}^2 
\leq \exp\left(-t\frac{\alpha}{\beta}\right) \normtwo{\bx^{(1)} - \bx^*}^2,
\end{aligned}
$$
where inequality ($\dag$) follows from  the descent lemma (\ref{descent_lem_gdscss} in Lemma~\ref{lemma:gdupd_sm_sconv}), and we have used the fact that $1 - r \leq \exp(-r)$ for all $r \in \real$.

\paragraph{When $\eta = \frac{2}{\alpha + \beta}$.}
By the smoothness and $\nabla f(\bx^*) = \bzero$, we have 
$
f(\bx^{(t)}) - f(\bx^*) \leq \frac{\beta}{2} \normtwo{\bx^{(t)} - \bx^*}^2
$
(Theorem~\ref{theorem:charac_smoo}).
Thus, by the characterization theorem of SC and SS (Theorem~\ref{theorem:charac_smoo_n_stronconv}),
$$
\small
\begin{aligned}
\normtwo{\bx^{(t+1)} - \bx^*}^2 &= \normtwo{\bx^{(t)} - \eta \nabla f(\bx^{(t)}) - \bx^*}^2 \\
&= \normtwobig{\bx^{(t)} -\bx^*}^2 - 2\eta \innerproduct{\nabla f(\bx^{(t)}), (\bx^{(t)} - \bx^*)} + \eta^2 \normtwobig{\nabla f(\bx^{(t)})}^2 \\
&\leq \left(1 - 2\frac{\eta \alpha \beta}{\beta + \alpha}\right) \normtwo{\bx^{(t)} - \bx^*}^2 + \left(\eta^2 - 2\frac{\eta}{\beta + \alpha}\right) \normtwobig{\nabla f(\bx^{(t)})}^2 \\
&= \left(\frac{\alpha - 1}{\alpha + 1}\right)^2 \normtwo{\bx^{(t)} - \bx^*}^2 
\leq \exp\left(-\frac{4t}{\alpha + 1}\right) \normtwo{\bx^{(1)} - \bx^*}^2.
\end{aligned}
$$
This completes the proof.
\end{proof}

\index{Condition number}
\paragrapharrow{Condition number.}
We observe that the convergence rate of the gradient descent algorithm for an $\alpha$-SC and $\beta$-SS function is of the form 
$$
\normtwo{\bx^{(t+1)} - \bx^*}^2 \leq \exp\left(-t\frac{\alpha}{\beta}\right) \normtwo{\bx^{(1)} - \bx^*}^2.
$$
The number $\kappa \triangleq \frac{\beta}{\alpha}$ is known as the \textit{condition number} 
of the optimization problem.
The condition number of the objective function significantly affects the convergence rate of algorithms.
Indeed, if $\kappa = \frac{\beta}{\alpha}$ is small, then $\exp\left(-\frac{\alpha}{\beta}\right) = \exp\left(-\frac{1}{\kappa}\right)$ will be small, ensuring rapid convergence. However, if $\kappa \gg 1$, then $\exp\left(-\frac{1}{\kappa}\right) \approx 1$, potentially leading to slow convergence.

The concept of the condition number is central to numerical optimization. Below, we provide an informal definition of this concept. In subsequent sections, we will see how the condition number appears repeatedly in the context of the convergence of various optimization algorithms for both convex and non-convex problems. 
The exact numerical form of the condition number (for instance, here it is $\beta/\alpha$) will vary depending on the application~\footnote{In a linear system, the condition number can be defined as the ratio of the largest eigenvalue to the smallest eigenvalue of the underlying matrix; see, for example, \citet{lu2021numerical} for more details.}. 
However, in general, all these definitions of condition number will satisfy the following property.

\begin{definition}[Condition Number---Informal]
The condition number of a function $f: \sS \rightarrow \real$ is a scalar $\kappa \in \real$ that bounds how much the function value can change relative to a perturbation of the input.
\end{definition}

Functions with a small condition number are stable, meaning changes to their input do not significantly affect the function output values. Conversely, functions with a large condition number may exhibit abrupt changes in output values even with slight changes to the input. To better understand this concept, consider a differentiable function $f$ that is $\alpha$-SC and $\beta$-SS. 
Suppose $\bx^*$  is a stationary point of $f$, i.e., $\nabla f(\bx^*) = \bzero$. 
For a general function, such a point could be a local optimum or a saddle point. However, since $f$ is strongly convex, $\bx^*$ represents the unique global minimum of $f$. Therefore, for any other point  $\by\neq \bx^*$,
$$
\frac{\alpha}{2} \normtwo{\bx^* - \by}^2 \leq f(\by) - f(\bx^*) \leq \frac{\beta}{2} \normtwo{\bx^* - \by}^2.
$$
Dividing throughout by $\frac{\alpha}{2} \normtwo{\bx^* - \by}^{2}$ yields
$$
\frac{f(\by) - f(\bx^*)}{\frac{\alpha}{2} \normtwo{\bx^* - \by}^2} \in [1, \frac{\beta}{\alpha}] \triangleq [1, \kappa].
$$
Thus, upon perturbing the input from the global minimum $\bx^*$ to a point $\by$ at a distance $\normtwo{\bx^* - \by} = \epsilon$  away, the function value changes moderately---it increases by at least $\frac{\alpha \epsilon^2}{2}$ but no more than $\kappa \cdot \frac{\alpha \epsilon^2}{2}$. Such predictable behavior in response to perturbations makes it easier for optimization algorithms to achieve fast convergence.

\paragrapharrow{Change of variables and choice of $\bQ$-norm for greedy descent.}
In the next section, we will establish convergence results for non-Euclidean gradient descent methods. For now, let's focus on the $\bQ$-norm.
The choice of  $\bQ$ used to define the greedy descent direction can dramatically affect the convergence rate.
In Section~\ref{section:als-gradie-descent-taylor}, we showed that the greedy descent method with  $\bQ$-norm is the same as the gradient method applied to the problem after the change of variables $\widetildebx^\toptzero = \bQ^{1/2} \bx^\toptzero$. We know that the gradient method works well when the condition numbers $\frac{\beta}{\alpha}$ are moderate, and works poorly when the condition numbers are large. 
It follows that if the condition number is moderate after the change of variables $\widetildebx^\toptzero = \bQ^{1/2} \bx^\toptzero$, the greedy descent method will perform well according to Theorem~\ref{theorem:gd_sc_ss}.

This observation suggests a strategy for choosing $\bQ$: It should be selected so that the condition number of the underlying function, transformed by $\bQ^{-1/2}$, is well-conditioned. 
For example, if an approximation $\widetildebB$ of the Hessian at the optimal point $\nabla^2 f(\bx^*)$ were known, a very good choice for  $\bQ$ would be $\bQ = \widetildebB$, because  the Hessian of $f$ at the optimum would then satisfy
\begin{equation}\label{equation:gree_cond_opt}
\widetildebB^{-1/2} \nabla^2 f(\bx^*) \widetildebB^{-1/2} \approx \bI,
\end{equation}
indicating a likely low condition number.

\begin{algorithm}[H] 
\caption{Non-Euclidean Gradient  Descent Method}
\label{alg:non_euclidean_gd_gen}
\begin{algorithmic}[1] 
\Require A function $f(\bx)$; 
\State {\bfseries Input:}  Initialize $\bx^{(1)}$;
\For{$t=1,2,\ldots$}
\State Pick a stepsize $\eta_t$;
\State Pick $\bd_{\text{ngd}}^\toptzero \in \Lambda_{\nabla f(\bx^\toptzero)}$ and let $\bd_{\text{ugd}}^\toptzero \triangleq \norm{\nabla f(\bx^\toptzero)}_{*} \bd_{\text{ngd}}^\toptzero$;
\State Update $\bx^\toptone \leftarrow \bx^\toptzero - \eta_t\bd_{\text{ugd}}^\toptzero$;
\EndFor
\State {\bfseries Return:}   $\bx_{\text{final}}\leftarrow \bx^{(T)}$;
\end{algorithmic} 
\end{algorithm}

\section{Non-Euclidean Gradient Descent}\label{section:noneucli_gd}
For the unconstrained problem (P1) in \eqref{equation:p1p2}, we consider  the \textit{non-Euclidean gradient descent} method in Algorithm~\ref{alg:non_euclidean_gd_gen}, which is derived  from greedy searching using different norms (Section~\ref{section:als-gradie-descent-taylor}).
For any gradient $\nabla f(\bx^\toptzero)$ at $\bx^\toptzero$, Equations~\eqref{equation:norma_greedes} and \eqref{equation:unnorma_greedes} provide normalized and unnormalized greedy search directions, respectively, where we use the latter for illustration since the unnormalized greedy search directions in the $\ell_2$ Euclidean case is equivalent to the negative gradient direction, establishing guidance for the selection of  stepsizes; see Theorem~\ref{theorem:noneu_conv_ss}.

Letting $ \norm{\cdot} $ be any norm on $ \real^n $, by the conjugate subgradient theorem (Theorem~\ref{theorem:conju_subgra}), the set of primal counterparts (Definition~\ref{definition:set_primal}) can be equivalently stated as 
$$
\Lambda_{\ba} = \partial h(\ba), \text{ where } h(\cdot) \triangleq \norm{\cdot}_*,  \text{ for any }\ba\in\real^n,
$$
where  $\norm{\cdot}_*$ denotes the dual norm.
Using these definitions, the update rule for the non-Euclidean gradient descent method at the $t$-th iteration can be stated as
\begin{enumerate}
\item Pick $\bd_{\text{ngd}}^\toptzero \in \Lambda_{\nabla f(\bx^\toptzero)}$ and let $\bd_{\text{ugd}}^\toptzero \triangleq \norm{\nabla f(\bx^\toptzero)}_{*} \bd_{\text{ngd}}^\toptzero$;
\item Update $\bx^\toptone \leftarrow  \bx^\toptzero - \eta_t\bd_{\text{ugd}}^\toptzero$;
\end{enumerate}

Note that smoothness was defined under the $\ell_2$ norm in Definition~\ref{definition:scss_func}. In this section, letting $ \norm{\cdot} $ be any norm on $ \real^n $, we assume that $f$ is $\beta$-smooth w.r.t. the underlying norm:
\begin{subequations}
\begin{equation}\label{equation:smooth_noneul}
f(\by)-f(\bx)-\nabla f(\bx)^\top (\by-\bx)
\leq \frac{\beta}{2} \norm{\bx-\by}^2, \text{ for any }\bx,\by\in\real^n.
\end{equation}
Similar to the SS Property-O (Theorem~\ref{theorem:equi_gradsch_smoo}), this also shows 
\begin{equation}\label{equation:smooth_noneu2}
\norm{\nabla f(\bx) - \nabla f(\by)}_* \leq \beta  \norm{\bx-\by}.
\end{equation}
\end{subequations}
We are now ready to present the convergence results of the non-Euclidean gradient method.
\begin{theoremHigh}[Non-Euclidean Method for SS]\label{theorem:noneu_conv_ss}
Let $f:\real^n\rightarrow \real$ be a $\beta$-smooth function w.r.t. to the norm $\norm{\cdot}$. Let $\{\bx^\toptzero\}_{t > 0}$ be the sequence generated by the non-Euclidean gradient method (Algorithm~\ref{alg:non_euclidean_gd_gen}) for solving  problem (P1)
with  a constant stepsize corresponding to $\eta_t\triangleq\eta \in\left(0, \frac{1}{\beta}\right)$. Then,
\begin{enumerate}[(i)]
\item It satisfies that $  f(\bx^\toptzero) - f(\bx^\toptone) \geq \eta(1-\frac{\beta\eta}{2})\norm{\nabla f(\bx^\toptzero)}_*^2$.
\item The sequence $\{ f(\bx^\toptzero) \}_{t > 0}$ is nonincreasing; in addition, $f(\bx^\toptone) < f(\bx^\toptzero)$ if and only if $\nabla f(\bx^\toptzero) \neq \bzero$.
\item If the sequence $\{ f(\bx^\toptzero) \}_{t > 0}$ is bounded below, then $\nabla f(\bx^\toptzero) \to \bzero$ as $t \to \infty$;
\item If the optimal value of (P1) is finite and let $\bx^*$ be any optimal point. Then, for any $T>0$,
\begin{equation}
\min_{t=1,2,\ldots,T} \norm{\nabla f(\bx^\toptzero)}_* \leq 
\frac{\sqrt{f(\bx^{(1)}) - f(\bx^*)}}{\sqrt{\eta(1-\frac{\beta\eta}{2})T}}.
\end{equation}
\item All limit points of the sequence $\{ \bx^\toptzero \}_{t > 0}$ are stationary points of problem (P1).
\end{enumerate}
\end{theoremHigh}
\begin{proof}[of Theorem~\ref{theorem:noneu_conv_ss}]
\textbf{(i, ii).}
By the smoothness \eqref{equation:smooth_noneul} w.r.t. the norm $\norm{\cdot}$, the definition of the stepsize $\eta\in\left(0, \frac{1}{\beta}\right)$, and the update rule of Algorithm~\ref{alg:non_euclidean_gd_gen}, we have
\begin{equation}\label{equation:noneu_conv_ss1}
\small
\begin{aligned}
f(\bx^\toptone) &\leq f(\bx^\toptzero) + \innerproduct{\nabla f(\bx^\toptzero), \bx^\toptone - \bx^\toptzero} 
+ \frac{\beta}{2} \norm{\bx^\toptone - \bx^\toptzero}^2 \\
&= f(\bx^\toptzero) - \eta{\norm{\nabla f(\bx^\toptzero)}_*} \innerproduct{\nabla f(\bx^\toptzero), \bd_{\text{ngd}}^\toptzero} 
+\frac{\beta\eta^2}{2} \norm{\nabla f(\bx^\toptzero)}_*^2 \\
&\stackrel{\dag}{=} f(\bx^\toptzero) - \eta{\norm{\nabla f(\bx^\toptzero)}_*^2} +\frac{\beta\eta^2}{2} \norm{\nabla f(\bx^\toptzero)}_*^2 \\
&= f(\bx^\toptzero) - \eta(1-\frac{\beta\eta}{2})\norm{\nabla f(\bx^\toptzero)}_*^2,
\end{aligned}
\end{equation}
where the equality ($\dag$) follows from the definition of the dual norm. This inequality readily implies that $f(\bx^\toptzero) \geq f(\bx^\toptone)$ since $\eta\in\left(0, \frac{1}{\beta}\right)$ and that if $\nabla f(\bx^\toptzero) \neq \bzero$, then $f(\bx^\toptone) < f(\bx^\toptzero)$.

\paragraph{(iii).} Since the sequence $\{ f(\bx^\toptzero) \}_{t > 0}$ is nonincreasing and bounded below, it converges. Thus, in particular $f(\bx^\toptzero) - f(\bx^\toptone) \to 0$ as $t \to \infty$, which, combined with \eqref{equation:noneu_conv_ss1}, implies that $\nabla f(\bx^\toptzero) \to \bzero$ as $t \to \infty$.

\paragraph{(iv).} Let $B\triangleq\eta(1-\frac{\beta\eta}{2})$. By (i), for any $t > 0$,
$ f(\bx^\toptzero) - f(\bx^\toptone) \geq B \norm{\nabla f(\bx^\toptzero)}_*^2. $
Performing telescopic sum over $t = 1,2, \ldots, T$, we obtain
$$ 
f(\bx^{(1)}) - f(\bx^{(T)}) \geq B \sum_{t=1}^T \norm{\nabla f(\bx^\toptzero)}_*^2 \geq T B \min_{t=1,2\ldots,T} \norm{\nabla f(\bx^\toptzero)}_*^2. 
$$
Since $f(\bx^{(T)}) \geq f(\bx^*)$, the result follows.

\paragraph{(v).} Let $\widetildebx$ be a limit point of $\{ \bx^\toptzero \}_{t > 0}$. Then there exists a subsequence $\{ \bx^{(t_j)} \}_{j \geq 0}$ converging to $\widetildebx$. By the triangle inequality and the smoothness in \eqref{equation:smooth_noneu2},  for any $j \geq 0$
$$
\norm{\nabla f(\widetildebx)}_* \leq \norm{\nabla f(\bx^{(t_j)}) - \nabla f(\widetildebx)}_* + \norm{\nabla f(\bx^{(t_j)}) }_* \leq \beta \norm{\bx^{(t_j)} - \widetildebx} + \norm{\nabla f(\bx^{(t_j)})}_*. 
$$
Since the right-hand side of the above inequality goes to $0$ as $j \to \infty$, it follows that $\nabla f(\widetildebx) = \bzero$.
\end{proof}

\begin{remark}[Smoothness Parameter]
Note that to prove the result, we require the function to be smooth. In many cases, the smoothness parameter $\beta$ might not be known in advance. In such scenarios, similar results can be obtained using backtracking (Algorithm~\ref{alg:gd_line_search}) or exact line search as long as these methods satisfy:
$$
f(\bx^\toptone) \leq f(\bx^\toptzero) + \innerproduct{\nabla f(\bx^\toptzero), \bx^\toptone - \bx^\toptzero} + \frac{1}{2\eta_t} \norm{\bx^\toptone - \bx^\toptzero}^2
$$ 
for each iteration.
This discussion applies to many results in this section; thus, we will not repeat the details for them subsequently.
\end{remark}

Similar to Theorem~\ref{theorem:pgd_smooth}, we can establish an $\mathcalO(1/T)$ rate of convergence for convex functions using non-Euclidean gradient descent methods.
To demonstrate this, we require an additional boundedness assumption:
\begin{itemize}
\item For the initial point $\bx^{(1)}$ and any optimal point $\bx^*$, there exists $R > 0$ such that
$
\max_{\bx, \bx^*} \{ \norm{\bx^* - \bx} \mid f(\bx) \leq f(\bx^{(1)})  \} \leq R
$.
\end{itemize}
The proof of the convergence rate is based on the following lemma.
\begin{lemma}\label{lemma:noneu_cvx_sequence}
Let $\{a_t\}_{t > 0}$ be a sequence of nonnegative real numbers satisfying for any $t > 0$:
$$ a_t - a_{t+1} \geq \frac{1}{\gamma} a_t^2 $$
for some $\gamma > 0$. Then, for any $T \geq 1$,
\begin{equation}
a_T \leq \frac{\gamma}{T-1}.
\end{equation}
\end{lemma}
\begin{proof}
Let $t$ be a positive integer. If $a_t = 0$, then the result holds trivially. We then assume that $a_t > 0$. By the monotonicity of $\{a_t\}_{t > 0}$, we have that $a_1, a_2, \ldots, a_T > 0$. For any $t = 1, 2, \ldots, T$,
$$
\frac{1}{a_{t+1}} - \frac{1}{a_{t}} = \frac{a_{t} - a_{t+1}}{a_{t} a_{t+1}} \geq \frac{1}{\gamma} \frac{a_{t}^2}{a_{t} a_{t+1}} = \frac{1}{\gamma} \frac{a_{t}}{a_{t+1}} \geq \frac{1}{\gamma}.
$$
where the last inequality follows from the monotonicity of the sequence. Telescoping the sum over $t = 1, 2, \ldots, T$ yields that 
$
\frac{1}{a_T} \geq \frac{1}{a_1} + \frac{T-1}{\gamma} \geq \frac{T-1}{\gamma}
$.
\end{proof}

\begin{theoremHigh}[Non-Euclidean Method for Convex and SS: $\mathcalO(1/T)$]\label{theorem:noneu_conv_sscvx}
Let $f:\real^n\rightarrow \real$ be a $\beta$-smooth and \textbf{convex} function w.r.t. to the norm $\norm{\cdot}$. Let $\{\bx^\toptzero\}_{t > 0}$ be the sequence generated by the non-Euclidean gradient method (Algorithm~\ref{alg:non_euclidean_gd_gen}) for solving  problem (P1)
with  a constant stepsize corresponding to $\eta_t\triangleq\eta \in\left(0, \frac{1}{\beta}\right)$~\footnote{Similar results can be obtained using backtracking or exact line search as long as it satisfies $f(\bx^\toptone) \leq f(\bx^\toptzero) + \innerproduct{\nabla f(\bx^\toptzero), \bx^\toptone - \bx^\toptzero} + \frac{1}{2\eta_t} \norm{\bx^\toptone - \bx^\toptzero}^2$ for each iteration.}.
Suppose further that for the initial point $\bx^{(1)}$ and any optimal point $\bx^*$, there exists $R > 0$ such that
$
\max_{\bx, \bx^*} \{ \norm{\bx^* - \bx} \mid f(\bx) \leq f(\bx^{(1)}) \} \leq R
$.
Then,
\begin{itemize}
\item It follows that 
$
f(\bx^\toptzero) - f(\bx^\toptone) \geq \frac{1}{C} \big(f(\bx^\toptzero) - f(\bx^*)\big)^2, 
$
where $C\triangleq\frac{R^2}{\eta(1-\frac{\beta\eta}{2})}$ is a positive constant. 
\item For any $T \geq 1$,
$
f(\bx^{(T)}) - f(\bx^*) \leq \frac{C}{T-1}
$.
\end{itemize}
\end{theoremHigh}
\begin{proof}[of Theorem~\ref{theorem:noneu_conv_sscvx}]
By Theorem~\ref{theorem:noneu_conv_ss}(ii), $\{f(\bx^\toptzero)\}_{t > 0}$ is nonincreasing, and in particular for any $t > 0$ it holds that $f(\bx^\toptzero) \leq f(\bx^{(1)})$. Therefore, for any $\bx^* $ and $t > 0$,
$
\norm{\bx^\toptzero - \bx^*} \leq R
$. 
By Theorem~\ref{theorem:noneu_conv_ss}(i), we have 
\begin{equation}
f(\bx^\toptzero) - f(\bx^\toptone) \geq \eta\big(1-\frac{\beta\eta}{2}\big) \norm{\nabla f(\bx^\toptzero)}_*^2.
\end{equation}
This establishes the first claim.
By the convexity of $f$ (Theorem~\ref{theorem:conv_gradient_ineq}) and the \holders inequality (Theorem~\ref{theorem:holder-inequality}), for any optimal point $\bx^*$,
\begin{equation}
\small
\begin{aligned}
f(\bx^\toptzero) - f(\bx^*) 
\leq \innerproduct{\nabla f(\bx^\toptzero), \bx^\toptzero - \bx^*} 
\leq \norm{\nabla f(\bx^\toptzero)}_* \norm{\bx^\toptzero - \bx^*}
\leq R \norm{\nabla f(\bx^\toptzero)}_*. 
\end{aligned}
\end{equation}
Combining the above two inequalities yields that 
$$ 
f(\bx^\toptzero) - f(\bx^\toptone) \geq \eta(1-\frac{\beta\eta}{2}) \norm{\nabla f(\bx^\toptzero)}_*^2 \geq \frac{\eta(1-\frac{\beta\eta}{2})}{R^2} \big(f(\bx^\toptzero) - f(\bx^*)\big)^2. 
$$
Applying Lemma~\ref{lemma:noneu_cvx_sequence} proves the second claim.
\end{proof}

\begin{algorithm}[h] 
\caption{Projected (Sub)Gradient Descent Method}
\label{alg:pgd_gen}
\begin{algorithmic}[1] 
\Require A function $f(\bx)$ and a set $\sS$; 
\State {\bfseries Input:}  Initialize $\bx^{(1)}$;
\For{$t=1,2,\ldots$}
\State Pick a stepsize $\eta_t$;
\State $\by^{(t+1)} \leftarrow \bx^{(t)} - \eta_t \nabla f(\bx^{(t)})$ or $\by^{(t+1)} \leftarrow \bx^{(t)} - \eta_t f^\prime(\bx^{(t)})$, where $ f^\prime(\bx^{(t)})\in \partial f(\bx^{(t)})$;
\State $\bx^{(t+1)} \in \mathcalP_{\sS}(\by^{(t+1)})$;
\EndFor
\State (Output Option 1) Output  $\bx_{\text{final}}\leftarrow \bx^{(T)}$;
\State (Output Option 2) Output  $\bx_{\text{avg}}\leftarrow \frac{1}{T}(\sum_{t=1}^{t}\bx^{(t)})$ or $\sum_{t=1}^{T} \frac{2t}{T(T+1)} \bx^{(t)}$;
\State (Output Option 3) Output  $\bx_{\text{best}}\leftarrow \argmin_{t\in\{1,2,\ldots,T\}} f(\bx^{(t)})$;
\end{algorithmic} 
\end{algorithm}


\section{Projected Gradient Descent}\label{section:pgd}
Standard gradient descent is typically used to solve unconstrained optimization problems. Given a set $\sS$, the \textit{projected gradient descent (PGD)} method  (Algorithm~\ref{alg:pgd_gen}) addresses constrained optimization problems of the form:
\begin{equation}\label{equation:p2_pgd}
\text{(P2)}: \qquad \min_{\bx\in\real^n } f(\bx) \quad \text{s.t.}\quad \bx\in\sS.
\end{equation}
Sections~\ref{section:constrain_opt} and \ref{section:constr_convset} discuss optimality conditions for this type of constrained optimization problem.
The PGD method represents the simplest approach to tackling the constrained optimization problem (P2).
As previously mentioned, the PGD method forms the cornerstone of first-order optimization algorithms, serving as a fundamental technique; refer to Figure~\ref{fig:first_order_opt_map}.
Numerous algorithms, such as the proximal gradient method, are derived from this basic algorithm.
This algorithm can be viewed as a fixed-point iteration for solving the stationarity condition
$$
\bx^* =\mathcalP_{\sS}(\bx^* - \eta\nabla f(\bx^*))
$$
if the set $\sS$ is closed convex (Theorem~\ref{theorem:stat_point_uncons_convset_proj}, assuming the stepsize is constant).

To illustrate further, in PGD, the update rule from the $t$-th to the $(t+1)$-th iteration can be equivalently expressed as:
\begin{tcolorbox}[colback=white,colframe=black]
\begin{minipage}{1\textwidth}
\begin{equation}\label{equation:pgd_decom_raw}
\small
\textbf{(PGD)}:\quad 
\begin{aligned}
\bx^\toptone 
&\leftarrow \mathcalP_{\sS}(\bx^{(t)} - \eta_t\nabla f(\bx^{(t)})) \\
&= \mathop{\argmin}_{\bx\in\sS} \normtwo{\bx - \left(\bx^{(t)} - \eta_t \nabla f(\bx^{(t)})\right)}^2\\
&=\mathop{\argmin}_{\bx\in\sS} 
\left\{
\frac{1}{2\eta_t} \normtwo{\bx-\bx^\toptzero}^2 + f(\bx^\toptzero)+ \innerproduct{\nabla f(\bx^\toptzero), \bx-\bx^\toptzero} 
\right\},
\end{aligned}
\end{equation}
\end{minipage}
\end{tcolorbox}
\noindent
where we add a constant $f(\bx^\toptzero)$ term in the third form since the optimization is over $\bx$.
Therefore, the update rule in PGD can be interpreted  as a minimization of the sum of a linearization of the smooth component  around the current iterate plus a quadratic term. The optimization of the linearization ensures that we can make a significant progress, while the quadratic term acts as a regularization mechanism, ensuring the update remains within a neighborhood of the current iterate.
Understanding this decomposition and the intuition behind the PGD updates is crucial for developing other methods, such as the proximal gradient method, conditional gradient method, and mirror descent method.

We start  by proving an $\mathcalO(1/\sqrt{T})$ rate of convergence for convex and Lipschitz functions (Definition~\ref{definition:lipschi_funs}) with a constant stepsize. 
This result depends on an upper bound on the distance between the initial guess $\bx^{(1)}$ and any optimal point $\bx^*$.
Additionally, we assume the total number of iterations $T$ is known in advance.
Note that Theorem~\ref{theorem:pgd_lipschitz_dyna} provides an alternative convergence result  when dynamic stepsizes are employed.
\begin{theoremHigh}[PGD for Convex and Lipschitz Functions: $\mathcalO(1/\sqrt{T})$]\label{theorem:pgd_lipschitz}\footnote{Theorem~\ref{theorem:pgd_lipschitz_dyna} provides an alternative convergence result with dynamic stepsizes.}
Let $f:\sS\rightarrow \real$ be a proper convex,  differentiable, and $L$-Lipschitz function defined over the closed and convex domain $\sS\subseteq\real^n$. 
Let $\bx^{(1)}, \bx^{(2)}, \ldots, \bx^{(T)}$ be the sequence of $T$ steps generated by the projected gradient descent method (Algorithm~\ref{alg:pgd_gen}) for solving  problem (P2)
with  a constant stepsize corresponding to  $\eta = \frac{R}{L\sqrt{T}}$, where $R$ is the upper bound on the distance $\normtwo{\bx^{(1)}-\bx^{*}}$ from the initial point $\bx^{(1)}$ to an optimal point $\bx^{*} \in \arg\min_{\bx\in\sS} f(\bx)$.~\footnote{Under the assumption that  $ \sS $ is \textbf{compact}.}
Then, for any $T>0$,
$$
\small
\begin{aligned}
f\left(\frac{1}{T}\sum_{t=1}^{T}\bx^{(t)}\right) - f\left(\bx^{*}\right) \leqslant \frac{RL}{\sqrt{T}}.
\end{aligned}
$$
This implies that for any integer $T$ satisfying $T\geq \frac{R^2L^2}{\epsilon^2}$, it follows that 
$$
\small
\begin{aligned}
f\left(\frac{1}{T}\sum_{t=1}^{T}\bx^{(t)}\right) - f\left(\bx^{*}\right) \leq \epsilon.~\footnote{If $R$ is not known, similar result can be obtained that $f\left(\frac{1}{T}\sum_{t=1}^{T}\bx^{(t)}\right) - f\left(\bx^{*}\right) \leq \epsilon$ for $T=\mathcalO(\frac{1}{\epsilon^2})$ and $\eta = \frac{1}{\sqrt{T}}$.}
\end{aligned}
$$
\end{theoremHigh}
\begin{proof}[of Theorem~\ref{theorem:pgd_lipschitz}]
The proof begins by  bounding the difference in function values
$
f(\bx^{(t)}) - f(\bx^*)
$ (called the \textit{sub-optimality} at the $t$-th iterate):
\begin{equation}\label{equation:pgd_lipschitz1}
\small
\begin{aligned}
f(\bx^{(t)}) - f(\bx^*) &\stackrel{\dag}{\leq} \innerproduct{\nabla f(\bx^{(t)}),  (\bx^{(t)} - \bx^*)}
\stackrel{*}{=} \frac{1}{\eta} \innerproduct{(\bx^{(t)} - \by^{(t+1)}), (\bx^{(t)} - \bx^*)}\\
&\stackrel{\ddag}{=} \frac{1}{2\eta} \left(\normtwo{\bx^{(t)} - \by^{(t+1)}}^2 + \normtwo{\bx^{(t)} - \bx^*}^2  - \normtwo{\by^{(t+1)} - \bx^*}^2 \right) \\
&\stackrel{*}{=} \frac{1}{2\eta} \left( \normtwo{\bx^{(t)} - \bx^*}^2 - \normtwo{\by^{(t+1)} - \bx^*}^2 \right) + \frac{\eta}{2} \normtwobig{\nabla f(\bx^{(t)})}^2\\
&\stackrel{+}{\leq} \frac{1}{2\eta} \left( \normtwo{\bx^{(t)} - \bx^*}^2 - \normtwo{\bx^{(t+1)} - \bx^*}^2 \right) + \frac{\eta L^2}{2},
\end{aligned}
\end{equation}
where the inequality ($\dag$) follows from convexity (Theorem~\ref{theorem:conv_gradient_ineq}),  the equalities ($*$) follow from the update rule of PGD, 
the equality ($\ddag$) follows from Theorem~\ref{theorem:funda_opt}, and  the inequality ($+$) follows from the Lipschitz condition (Theorem~\ref{theorem:lipsc_equiv}) and Projection Property-II (Lemma~\ref{lemma:proj_prop2}).

Now, sum these differences from $t=1$ to $t=T$:
$$
\small
\begin{aligned}
\sum_{t=1}^{T} \left(f(\bx^{(t)}) - f(\bx^*)\right) 
&\leq \frac{1}{2\eta} \sum_{t=1}^{T} \left( \normtwo{\bx^{(t)} - \bx^*}^2 - \normtwo{\bx^{(t+1)} - \bx^*}^2 \right) + \frac{\eta L^2 T}{2}\\
&\leq \frac{1}{2\eta} \normtwo{\bx^{(1)} - \bx^*}^2 + \frac{\eta L^2 T}{2}
\leq \frac{R^2}{2\eta} + \frac{\eta L^2 T}{2}.
\end{aligned}
$$~\footnote{If $R$ is not known, we can set $\bx^{(1)}=\bzero$ and $\eta=\frac{1}{\sqrt{T}}$, obtaining a similar complexity.}
Finally, by convexity and $\eta = \frac{R}{L\sqrt{T}} $,
$$
\small
\begin{aligned}
f\left(\frac{1}{T} \sum_{t=1}^T \bx^{(t)}\right) - f(\bx^*) \leq \left(\frac{1}{T} \sum_{t=1}^T f(\bx^{(t)})\right) - f(\bx^*) 
\leq \frac{R^2}{2\eta T} + \frac{\eta L^2}{2}
= \frac{R L}{\sqrt{T}}.
\end{aligned}
$$
This completes the proof.
\end{proof}

Let $e_t \triangleq f(\bx^{(t)}) - f(\bx^*)$. 
The equality (+) in \eqref{equation:pgd_lipschitz1} also tells us that the sub-optimality $e_t$ is small  if the consecutive iterates $\bx^\toptzero$ and $\bx^\toptone$ are close to each other, since $\frac{\eta L^2}{2}$ is small due to $\eta = \mathcalO(\frac{1}{\sqrt{T}})$.
This observation is quite useful since it tells us that once PGD stops making a lot of progress, it actually converges to the optimum.

\begin{figure}[h]
\centering       
\vspace{-0.25cm}                 
\subfigtopskip=2pt               
\subfigbottomskip=-2pt         
\subfigure[$f(a) + f(a+1) + \cdots + f(b) \leq \int_{a-1}^b f(t) \, dt$. ]{\label{fig:pgd_sumlem1}
\includegraphics[width=0.48\linewidth]{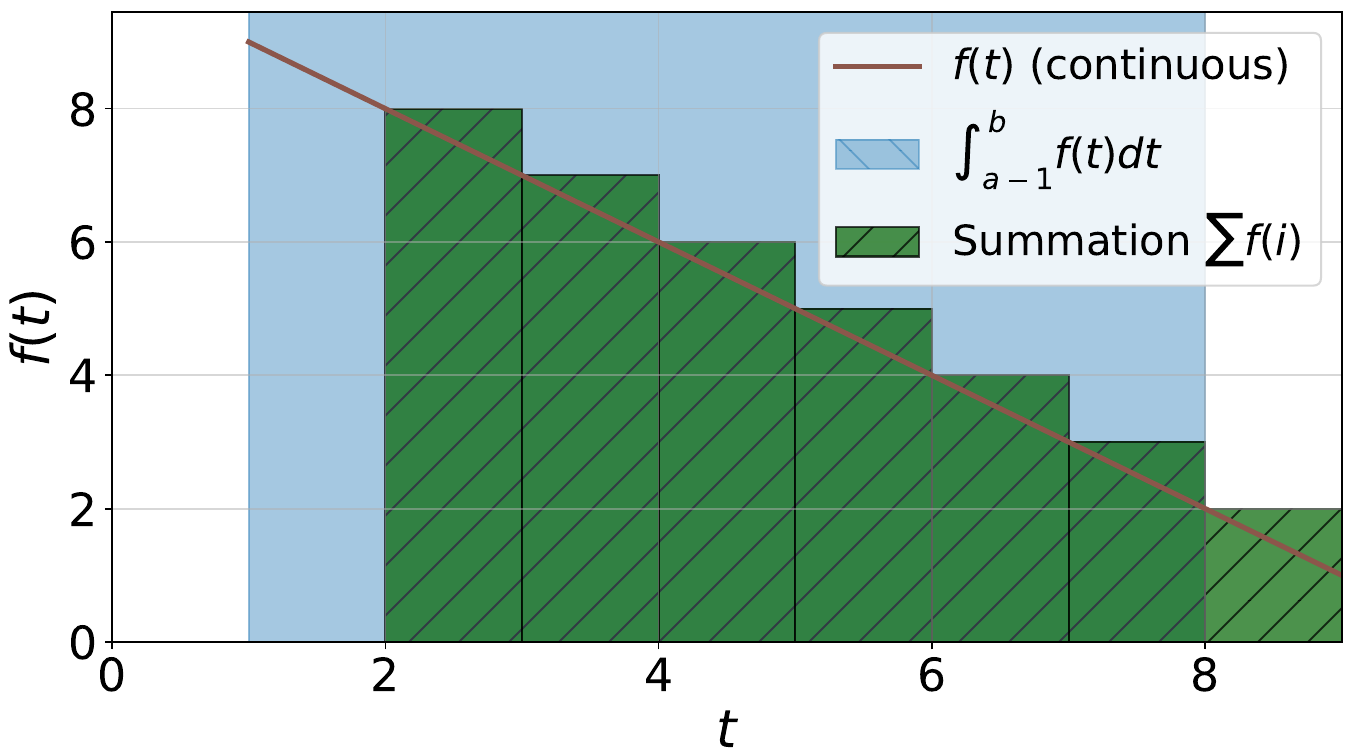}}
\subfigure[$\int_a^{b+1} f(t) \, dt \leq f(a) + f(a+1) + \cdots + f(b)$.]{\label{fig:pgd_sumlem2}
\includegraphics[width=0.48\linewidth]{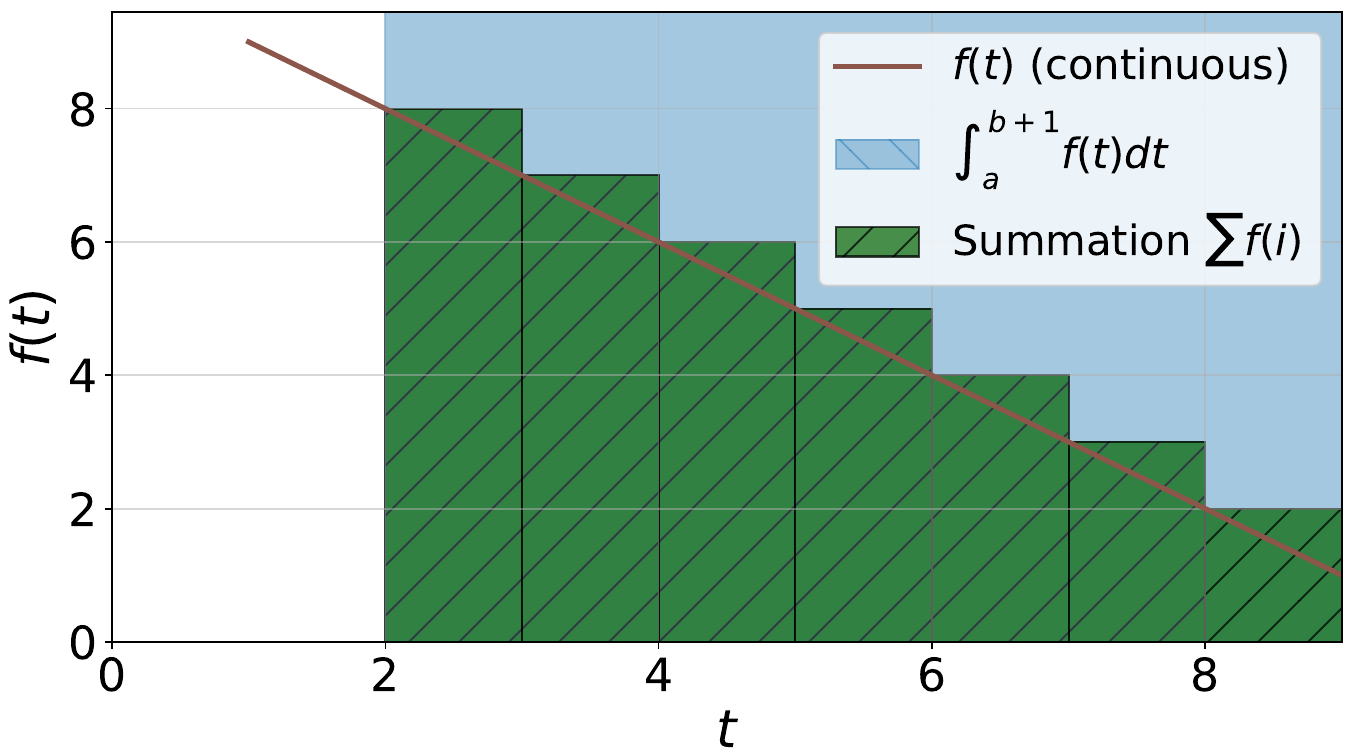}}
\caption{An example illustrating Lemma~\ref{lemma:pgd_lem1}, where we set $a\triangleq 1$ and $b\triangleq 8$ in the lemma.}
\label{fig:pgd_sumlem}
\end{figure}

Theorem~\ref{theorem:pgd_lipschitz} provides the convergence result for convex and Lipschitz function using predefined stepsizes for each iteration.
Alternative, we provide convergence results using dynamic stepsizes. To see this, we need the following lemmas.
\begin{lemma}\label{lemma:pgd_lem1}
Let $ f: [a-1, b+1] \rightarrow \real $ be a continuous nonincreasing function over $[a-1, b+1]$, where $a$ and $b$ are integer numbers satisfying $a \leq b$. Then,
$$
\int_a^{b+1} f(t) \, dt \leq f(a) + f(a+1) + \cdots + f(b) \leq \int_{a-1}^b f(t) \, dt.
$$
See Figure~\ref{fig:pgd_sumlem} for an illustrating example.
\end{lemma}

\begin{lemma}\label{lemma:pgd_lem2}
Let $ \Phi \in \real $. Then,
\begin{enumerate}[(i)]
\item \label{pgd_lem2_v1} For any $ T \geq 5 $,
$
\frac{\Phi + \sum_{t=1}^T \frac{1}{t+1}}{\sum_{t=1}^T \frac{1}{\sqrt{t+1}}} \leq \frac{\Phi  + \ln(T+1)}{\sqrt{T+1}}.
$

\item \label{pgd_lem2_v2} For any $ T \geq 2 $,
$
\frac{\Phi + \sum_{t=\lceil T/2 \rceil}^T \frac{1}{t+1}}{\sum_{t=\lceil T/2 \rceil}^T \frac{1}{\sqrt{t+1}}} \leq \frac{4(\Phi + \ln(3))}{\sqrt{T+2}}.
$
\end{enumerate}
\end{lemma}
\begin{proof}
\textbf{(i).} Using Lemma~\ref{lemma:pgd_lem1}, we obtain the following inequalities:
$$
\begin{aligned}
\sum_{t=1}^T \frac{1}{t+1} &\leq  \int_0^T \frac{1}{x+1} \, dx =  \ln(T+1);\\
\sum_{t=1}^T \frac{1}{\sqrt{t+1}} &\geq \int_1^{T+1} \frac{1}{\sqrt{x+1}} \, dx = 2\sqrt{T+2} - 2\sqrt{2} \geq \sqrt{T+1},
\end{aligned}
$$
where the last inequality holds for all $ T \geq 5 $. This yields the desired result in (i).
\paragraph{(ii).} Using Lemma~\ref{lemma:pgd_lem1}, we can also obtain the following inequalities for any $ T \geq 2 $:
$$
\begin{aligned}
\sum_{t=\lceil T/2 \rceil}^T \frac{1}{t+1} &\leq \int_{\lceil T/2 \rceil - 1}^T \frac{1}{x+1}\, dx  
= \ln\left(\frac{T+1}{\lceil  T/2 \rceil}\right) 
\leq \ln\left(2 + \frac{2}{T}\right) 
\leq \ln(3)
\end{aligned}
$$
and
$$
\begin{aligned}
\sum_{t=\lceil T/2 \rceil}^T \frac{1}{\sqrt{t+1}} &\geq \int_{\lceil T/2 \rceil}^{T+1} \frac{1}{\sqrt{x+1}}\,dx = 2\sqrt{T+2} - 2\sqrt{\lceil T/2 \rceil + 1} \\
&\geq 2\sqrt{T+2} - 2\sqrt{T/2 + 2} = \frac{4(T+2) - 4(T/2 + 2)}{2\sqrt{T+2} + 2\sqrt{T/2 + 2}} \\
&= \frac{T}{\sqrt{T+2} + \sqrt{T/2 + 2}} \geq \frac{T}{2\sqrt{T+2}}
\geq \frac{1}{4}\sqrt{T+2},
\end{aligned}
$$
where the last inequality holds since $ T \geq 2 $. This yields the desired result in (ii).
\end{proof}

The descent lemma of GD and PGD (Lemma~\ref{lemma:gdupd_sm_sconv}) demonstrates the function value inequality for GD or PGD updates. 
The following lemma illustrates the relationship between iterations of PGD.
\begin{lemma}[Iterate Inequality for PGD]\label{lemma:iterate_ineq_pgd}
Let $f: \sS\rightarrow \real$ be a  proper closed and  \textbf{convex} function, where $\sS\subseteq\real^n$ is a closed and convex set such that $\sS\subseteq \interior(\dom(f))$~\footnote{This ensures the nonemptness of the subdifferential by Theorem~\ref{theorem:nonemp_relint_conv}.}. Define $\bx^\toptone  \in \sS$ as $\bx^\toptone  \leftarrow \projectS\left(\bx^\toptzero - \eta_t f^\prime(\bx^\toptzero)\right)$, where $f^\prime(\bx^\toptzero)\in \partial f(\bx^\toptzero)$ denotes any subgradient. 
Let $\bx^*$ be any optimizer. Then,
\begin{equation}\label{equation:iterate_ineq_pgd1}
\normtwo{\bx^\toptone  - \bx^*}^2 \leq \normtwo{\bx^\toptzero - \bx^*}^2 - 2 \eta_t  \left(f(\bx^\toptzero) - f(\bx^*)\right) + \eta_t^2 \normtwo{f^\prime(\bx^\toptzero)}^2. 
\end{equation}
Let $\{\bx^\toptzero\}_{t> 0}$ be the sequence generated by this update rule with positive stepsizes $\{\eta_t\}_{t > 0}$. Then, for any   nonnegative integer $T>0$, performing telescopic cancellations using \eqref{equation:iterate_ineq_pgd1} shows that
\begin{equation}\label{equation:iterate_ineq_pgd2}
\sum_{t=1}^{T} \eta_t \big(f(\bx^\toptzero) - f(\bx^*)\big) \leq \frac{1}{2} \normtwo{\bx^{(1)} - \bx^*}^2 + \frac{1}{2} \sum_{t=1}^{T} \eta_t^2 \normtwo{f^\prime(\bx^\toptzero)}^2. 
\end{equation}
\end{lemma}
\begin{proof}[of Lemma~\ref{lemma:iterate_ineq_pgd}]
In this case, a minimum subgradient should be
$$
\begin{aligned}
\normtwo{\bx^\toptone - \bx^*}^2 &= \normtwo{\projectS\big(\bx^\toptzero - \eta_t f^\prime(\bx^\toptzero)\big) - \projectS(\bx^*)}^2 
\stackrel{\dag}{\leq} \normtwo{\bx^\toptzero - \eta_t f^\prime(\bx^\toptzero) - \bx^*}^2 \\
&= \normtwo{\bx^\toptzero - \bx^*}^2 - 2 \eta_t \innerproduct{f^\prime(\bx^\toptzero), \bx^\toptzero - \bx^*} + \eta_t^2 \normtwo{f^\prime(\bx^\toptzero)}^2 \\
&\stackrel{\ddag}{\leq} \normtwo{\bx^\toptzero - \bx^*}^2 - 2 \eta_t \left(f(\bx^\toptzero) - f(\bx^*)\right) + \eta_t^2 \normtwo{f^\prime(\bx^\toptzero)}^2,
\end{aligned}
$$
where the inequality ($\dag$) follows from the nonexpansiveness of the  projection operator (Theorem~\ref{theorem:proj_nonexpan}), and the inequality ($\ddag$) follows from the subgradient inequality (Definition~\ref{definition:subgrad}, or from the gradient inequality if $f$ is differentiable; Theorem~\ref{theorem:conv_gradient_ineq}).
\end{proof}
We are now ready to prove the convergence result.

\begin{theoremHigh}[PGD for Convex and Lipschitz Functions with Dynamic Stepsizes: $ \mathcalO(\ln(T)/\sqrt{T}) $]\label{theorem:pgd_lipschitz_dyna}
Let $f: \sS \rightarrow \real$ be a proper closed,  convex, and  $L$-Lipschitz function, where $\sS\subseteq\real^n$ is a closed and  convex set. 
Let $\{\bx^\toptzero\}_{t > 0}$ be the sequence generated by the projected subgradient method (Algorithm~\ref{alg:pgd_gen}) for solving  problem (P2).
Suppose further that  $\left\{\fbest^\toptzero \triangleq \min\left\{f(\bx^{(1)}), f(\bx^{(2)}), \ldots, f(\bx^{(t)})\right\}\right\}_{t>0}$ is the sequence of best achieved values.
Then,
\begin{enumerate}[(i)]
\item If $\frac{\sum_{t=1}^{T} \eta_t^2}{\sum_{t=1}^{T} \eta_t} \rightarrow 0  $ as $T \rightarrow \infty$, then $\fbest^{(T)} - f(\bx^*) \rightarrow 0$ as $T \rightarrow \infty$.
\item If the stepsize for each $t$-th iteration is $\eta_t = \frac{1}{\normtwo{f^\prime(\bx^\toptzero)}\sqrt{t+1}}$ if $f^\prime(\bx^\toptzero) \neq \bzero$ and $\eta_t = \frac{1}{L}$ otherwise. Then, for any $ T \geq 5 $,
$$
\max\left\{ \fbest^{(T)} - f(\bx^*) , f(\widetildebx^{(T)}) - f(\bx^*) \right\} \leq \frac{L}{2} \frac{\normtwo{\bx^{(1)} - \bx^*}^2 + \ln(T+1)}{\sqrt{T+1}},
$$
where $
\widetildebx^{(T)} = \frac{1}{\sum_{t=1}^{T} \eta_t} \sum_{t=1}^{T} \eta_t \bx^\toptzero
$,~\footnote{Let $\Phi_T \triangleq \sum_{t=1}^{T} \eta_t $. We note that $\widetildebx^{(T+1)} = \frac{\Phi_T}{\Phi_{T+1}} \widetildebx^{(T)} + \frac{\eta_{T+1}}{\Phi_{T+1}} \bx^{(T+1)}$, such that the sequence $\{\widetildebx^{(t)}\}_{t>0}$ can be computed recursively.}
i.e., a weighted average of the iterates.
\end{enumerate}

\end{theoremHigh}
\begin{proof}[of Theorem~\ref{theorem:pgd_lipschitz_dyna}]
\textbf{(i).}  Using \eqref{equation:iterate_ineq_pgd2} along with the inequality $f(\bx^\toptzero) \geq \fbest^{(T)}$ for any $t=\{1,2,\ldots,T\}$, we obtain
\begin{equation}\label{equation:pgd_lipschitz_dyna_1}
\fbest^{(T)} - f(\bx^*) \leq \frac{1}{2} \frac{\normtwo{\bx^{(1)} - \bx^*}^2 + \sum_{t=1}^{T} \eta_t^2 \normtwo{f^\prime(\bx^\toptzero)}^2}{\sum_{t=1}^{T} \eta_t}. 
\end{equation}
Since
$\sum_{t=1}^{T} \eta_t \rightarrow \infty$ as  $T \rightarrow \infty$ and $\normtwo{f^\prime\left(\bx^\toptzero\right)} \leq L$, the desired result in (i) follows directly from \eqref{equation:pgd_lipschitz_dyna_1}.

\paragraph{(ii).}
By Jensen's inequality (Theorem~\ref{theorem:jensens_ineq})
$
f(\widetildebx^{(T)}) \leq \frac{1}{\sum_{t=1}^{T} \eta_t} \sum_{t=1}^{T} \eta_t f(\bx^\toptzero)
$
and \eqref{equation:iterate_ineq_pgd2}, we have
$$
f(\widetildebx^{(T)}) - f(\bx^*) \leq \frac{1}{2} \frac{\normtwo{\bx^{(1)} - \bx^*}^{2} + \sum_{t=1}^{T} \eta_t^{2} \normtwo{f^\prime\left(\bx^\toptzero\right)}^{2}}{\sum_{t=1}^{T} \eta_t}. 
$$
The above two inequalities indicate
$$
\max\left\{\fbest^{(T)} - f(\bx^*), f(\widetildebx^{(T)}) - f(\bx^*)\right\} \leq \frac{1}{2} \frac{\normtwo{\bx^{(1)} - \bx^*}^{2} + \sum_{t=1}^{T} \eta_t^{2} \normtwo{f^\prime\left(\bx^\toptzero\right)}^{2}}{\sum_{t=1}^{T} \eta_t}.
$$
By the definition of $ \eta_t $, $ \eta_t^{2} \normtwo{f^\prime\left(\bx^\toptzero\right)}^{2} \leq \frac{1}{t+1} $ (satisfied as an equality when $ f^\prime\left(\bx^\toptzero\right) \neq \bzero $ and as a strict inequality when $ f^\prime\left(\bx^\toptzero\right) = \bzero $); in addition, since $ \normtwo{f^\prime\left(\bx^\toptzero\right)} \leq L $, we have $ \eta_t \geq \frac{1}{L \sqrt{t+1}} $. Therefore,
$$
\max\left\{\fbest^{(T)} - f(\bx^*), f(\widetildebx^{(T)}) - f(\bx^*)\right\} \leq \frac{L}{2} \frac{\normtwo{\bx^{(1)} - \bx^*}^{2} + \sum_{t=1}^{T} \frac{1}{t+1}}{\sum_{t=1}^{T} \frac{1}{\sqrt{t+1}}}.
$$
Invoking Lemma~\ref{lemma:pgd_lem2}\ref{pgd_lem2_v1} with $ \Phi \triangleq \normtwo{\bx^{(1)} - \bx^*}^{2} $ obtains the desired result.
\end{proof}

\begin{figure}[h]
\centering       
\vspace{-0.25cm}                 
\subfigtopskip=2pt               
\subfigbottomskip=-2pt         
\subfigure[$\eta_t=\frac{1}{\sqrt{t+1}}$.]{\label{fig:gd_conv_stepsizes_1}
\includegraphics[width=0.48\linewidth]{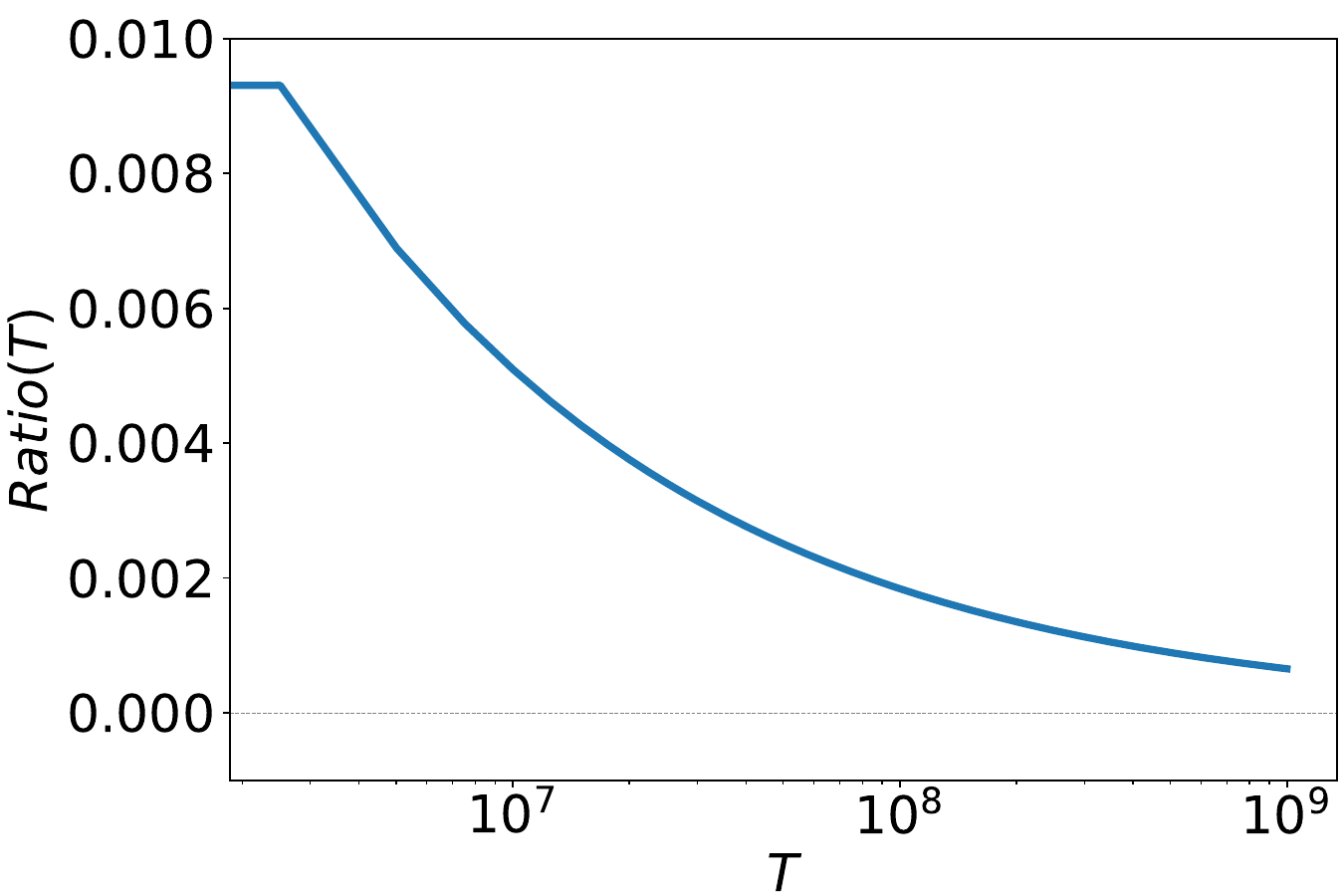}}
\subfigure[$\eta_t=\frac{1}{t}$.]{\label{fig:gd_conv_stepsizes_2}
\includegraphics[width=0.48\linewidth]{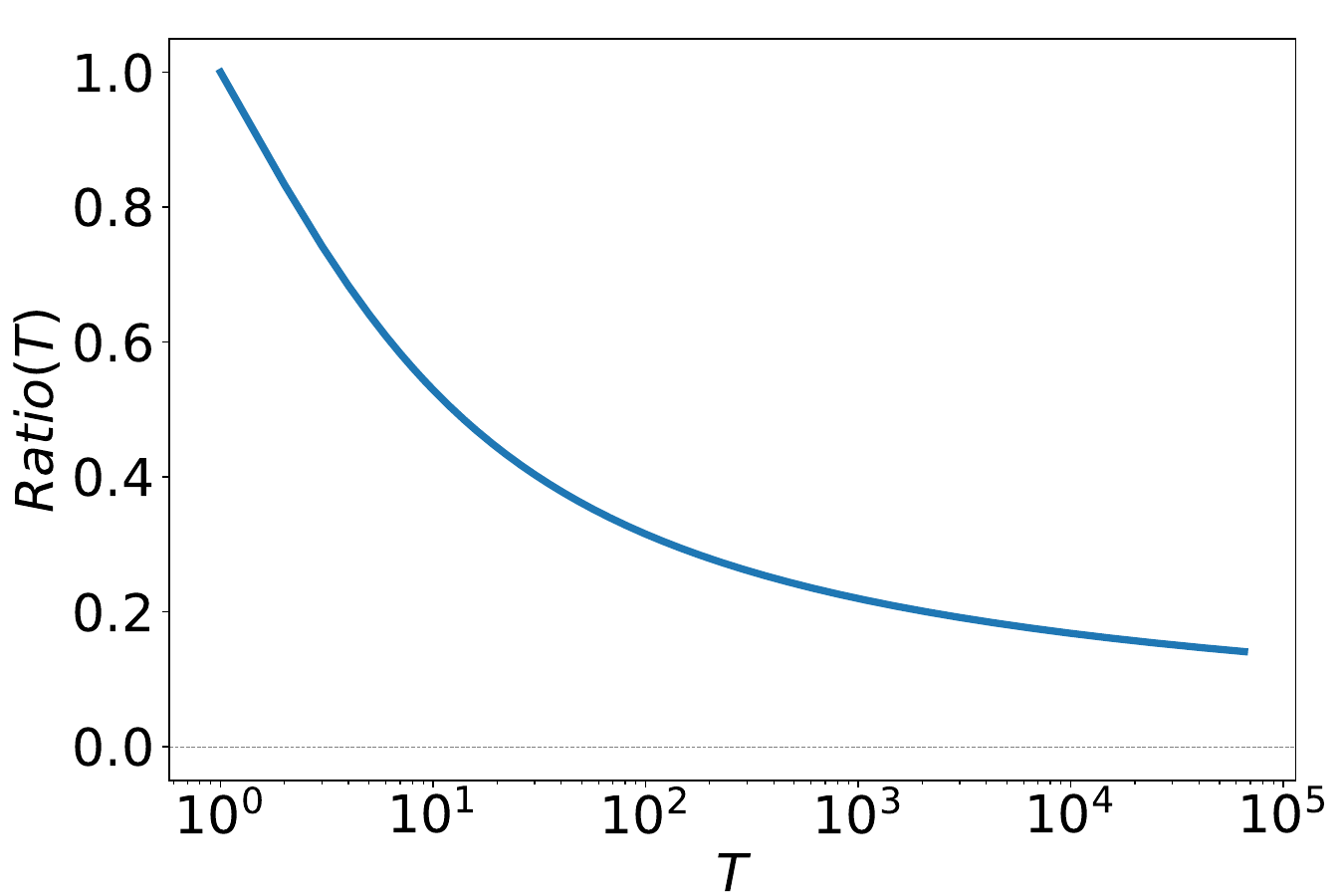}}
\caption{Illustration of the ratio  $\text{Ratio}(T)={(\sum_{t=1}^{T} \eta_t^2)}/{(\sum_{t=1}^{T} \eta_t)}$ as a function of the total number of iterations $T$ for $\eta_t\triangleq\frac{1}{\sqrt{t+1}}$ and $\eta_t\triangleq\frac{1}{t}$.}
\label{fig:gd_conv_stepsizes}
\end{figure}

According to Theorem~\ref{theorem:pgd_lipschitz_dyna}, we can choose, for example, the stepsizes as
$\eta_t \triangleq \frac{1}{\sqrt{t+1}}$ for each iteration.
This ensures convergence of function values to $ f(\bx^*) $,  since  $\sum_{t=1}^{T} \frac{1}{\sqrt{t+1}}$ is of the order of $\sqrt{T}$, whereas  $\sum_{t=1}^{T} \frac{1}{t+1}$ grows logarithmically with $T$.
Figure~\ref{fig:gd_conv_stepsizes_1} illustrates the ratio of $\frac{\sum_{t=1}^{T} \eta_t^2}{\sum_{t=1}^{T} \eta_t}$ as a function of the total number of iterations $T$ (note that Lemma~\ref{lemma:pgd_lem2}\ref{pgd_lem2_v1} demonstrates 
$
\frac{ \sum_{t=1}^T \frac{1}{t+1}}{\sum_{t=1}^T \frac{1}{\sqrt{t+1}}} \leq \frac{ \ln(T+1)}{\sqrt{T+1}}
$ in this scenario).
On the other hand, when  $ \eta_t \triangleq \frac{1}{t} $, the numerator becomes
$
\sum_{t=1}^{T} \eta_t^2  = \sum_{t=1}^{T} \frac{1}{t^2}
$.
The series $ \sum_{t=1}^{\infty} \frac{1}{t^2} $ converges (it is a $p$-series with $ p = 2 > 1 $), and its partial sums grow logarithmically.
The denominator is
$
\sum_{t=1}^{T} \eta_t = \sum_{t=1}^{T} \frac{1}{t}
$.
The harmonic series $ \sum_{t=1}^{\infty} \frac{1}{t} $ diverges, but its partial sums grow like $ \ln(T) $.
Therefore, the ratio also converges to 0 as $T\rightarrow\infty$ (Figure~\ref{fig:gd_conv_stepsizes_2}).

The $ \mathcalO(\ln(T)/\sqrt{T}) $ rate of convergence proved in Theorem~\ref{theorem:pgd_lipschitz_dyna} is less favorable compared to the $ \mathcalO(1/\sqrt{T}) $ rate established in Theorem~\ref{theorem:pgd_lipschitz} for the projected gradient descent method using a constant stepsize, where the stepsize $\eta = \frac{R}{L\sqrt{T}}$ depends on the total number of steps $T$. 
However, under the  assumption that the feasible set $ \sS $ is compact, it is also possible to achieve an $ \mathcalO(1/\sqrt{T}) $ convergence rate with dynamic stepsizes.

\begin{theoremHigh}[PGD for Convex and Lipschitz Functions with Dynamic Stepsizes: $ (\mathcalO(1/\sqrt{T}) $]\label{theorem:pggd_conv_lip_dyn}
Let $f: \sS \rightarrow \real$ be a proper closed,  convex, and  $L$-Lipschitz function, where $\sS\subseteq\real^n$ is a closed and  convex set, and assume that $ \sS $ is \textbf{compact}. Let $ \Omega $ be an upper bound on the half-squared diameter of $ \sS $:
$
\Omega \geq \max_{\bx,\by \in \sS} \frac{1}{2} \normtwo{\bx - \by}^2.
$
Let $\{\bx^\toptzero\}_{t > 0}$ sequence generated by the projected subgradient method (Algorithm~\ref{alg:pgd_gen}) for solving  problem (P2) with stepsizes chosen as either
$$
\eta_t = \frac{\sqrt{2\Omega}}{L \sqrt{t+1}} 
\qquad\text{or}\qquad
\eta_t = \begin{cases} 
\frac{\sqrt{2\Omega}}{\normtwo{f^\prime(\bx^\toptzero)}\sqrt{t+1}}, & f^\prime(\bx^\toptzero) \neq \bzero, \\
\frac{\sqrt{2\Omega}}{L \sqrt{t+1}}, & f^\prime(\bx^\toptzero) = \bzero.
\end{cases}
$$
Then, for all $T \geq 2$, it holds that 
$$
\fbest^{(T)} - f(\bx^*) \leq \frac{2(1 + \ln(3)) L \sqrt{2\Omega}}{\sqrt{T+2}},
$$
where $\left\{\fbest^\toptzero \triangleq \min\left\{f(\bx^{(1)}), f(\bx^{(2)}), \ldots, f(\bx^{(t)})\right\}\right\}_{t>0}$  is the sequence of best achieved values.
\end{theoremHigh}
\begin{proof}[of Theorem~\ref{theorem:pggd_conv_lip_dyn}]
By the iterate inequality in \eqref{equation:iterate_ineq_pgd1},  for any $t>0$,
$$
\frac{1}{2} \normtwo{\bx^\toptone - \bx^*}^2 \leq \frac{1}{2} \normtwo{\bx^\toptzero - \bx^*}^2 - \eta_t \big(f(\bx^\toptzero) - f(\bx^*)\big) + \frac{\eta_t^2}{2} \normtwo{f^\prime(\bx^\toptzero)}^2.
$$
Summing this inequality from $t = \lceil T/2 \rceil, \lceil T/2 \rceil + 1, \ldots, T$, we obtain
\begin{align*}
\sum_{t=\lceil T/2 \rceil}^{T} \eta_t \big(f(\bx^\toptzero) - f(\bx^*)\big) &\leq \frac{1}{2} \normtwo{\bx^{(\lceil T/2 \rceil)} - \bx^*}^2 - \frac{1}{2} \normtwo{\bx^{(T+1)} - \bx^*}^2 + \sum_{t=\lceil T/2 \rceil}^{T} \frac{\eta_t^2}{2} \normtwo{f^\prime(\bx^\toptzero)}^2 \\
&\leq \Omega + \sum_{t=\lceil T/2 \rceil}^{T} \frac{\eta_t^2}{2} \normtwo{f^\prime(\bx^\toptzero)}^2 
\leq \Omega + \Omega \sum_{t=\lceil T/2 \rceil}^{T} \frac{1}{t+1},
\end{align*}
where the last inequality follows from the definition of the stepsizes and the boundedness of the subgradients.
Since $\eta_t \geq \frac{\sqrt{2\Omega}}{L \sqrt{t+1}}$ and $f(\bx^\toptzero) \geq \fbest^{(T)}$ for all $t \leq T$, it follows that
$$
\sum_{t=\lceil T/2 \rceil}^{T} \eta_t \big(f(\bx^\toptzero) - f(\bx^*)\big) \geq \left( \sum_{t=\lceil T/2 \rceil}^{T} \frac{\sqrt{2\Omega}}{L \sqrt{t+1}} \right) \big(\fbest^{(T)} - f(\bx^*)\big). \
$$
Therefore, combining the above two inequalities and  invoking Lemma~\ref{lemma:pgd_lem2}\ref{pgd_lem2_v2} yields
$$
\fbest^{(T)} - f(\bx^*) \leq  \frac{L \sqrt{\Omega}}{\sqrt{2}} \frac{ \left(1 + \sum_{t=\lceil T/2 \rceil}^{T} \frac{1}{t+1}\right)}{\sum_{t=\lceil T/2 \rceil}^{T} \frac{1}{\sqrt{t+1}}}
\leq  \frac{L \sqrt{\Omega}}{\sqrt{2}} \frac{4(1 + \ln(3))}{\sqrt{T+2}}.
$$
This completes the proof.
\end{proof}

Based on   from Theorem~\ref{theorem:pggd_conv_lip_dyn}, assume further that $f$ is strongly convex.
The following theorem demonstrates that the convergence rate of the projected gradient descent for $\alpha$-strongly convex functions is similar to that of $\beta$-smooth functions, achieving an error bound of $\epsilon \leq \mathcalO(1/T)$, which is faster compared to the case involving only Lipschitz continuity.
\begin{theoremHigh}[PGD for SC and Lipschitz Functions: $\mathcalO(1/T)$]\label{theorem:pgd_strongconvex}
Let $f: \sS \rightarrow \real$ be a proper, differentiable, $\alpha$-strongly convex, and  $L$-Lipschitz function, where $\sS\subseteq\real^n$ is a closed and  convex set.
Let $\{\bx^\toptzero\}_{t > 0}$ be the sequence generated by the projected gradient descent method (Algorithm~\ref{alg:pgd_gen}) for solving  problem (P2) with an adaptive stepsize $\eta_t = \frac{2}{\alpha(t+1)}$ for each iteration $t$.
Assume  the function has a minimizer  $\bx^*$. Then,~\footnote{Note the result also holds for projected subgradient descent with non-differentiable functions.}
\begin{enumerate}[(i)]
\item  Define the sequence $\{\widetildebx^\toptzero \triangleq \sum_{i=1}^{t} \frac{2i}{t(t+1)} \bx^{(i)}\}$ (a convex combination of iterates). Then, for any $T>0$,
$$
\small
\begin{aligned}
f\left(\widetildebx^{(T)}\right) - f(\bx^*) \leq \frac{2L^2}{\alpha(T+1)}.
\end{aligned}
$$
\item Define  $\left\{\fbest^\toptzero \triangleq \min\left\{f(\bx^{(1)}), f(\bx^{(2)}), \ldots, f(\bx^{(t)})\right\}\right\}$ as the sequence of best achieved values. Then, for any $T>0$,
$$
\small
\begin{aligned}
\fbest^{(T)} -f(\bx^*) &\leq \frac{2L^2}{\alpha(T+1)}.
\end{aligned}
$$
It also follows that 
$$
\small
\begin{aligned}
\normtwo{\bx^{(i_T)} - \bx^{*}} \leq \frac{2 L}{\alpha\sqrt{T+1}},
\end{aligned}
$$
where $i_T\in\mathop{\argmin}_{t=1,2,\ldots,T} f(\bx^\toptzero)$.
\end{enumerate}
\noindent
This implies for any integer $T$ satisfying $T\geq \frac{2L^2}{\alpha\epsilon} -1$, it follows that 
$$
f\left(\widetildebx^{(T)}\right) - f(\bx^*)  \leq \epsilon 
\qquad\text{and}\qquad 
\fbest^{(T)} -f(\bx^*) \leq \epsilon.
$$
\end{theoremHigh}
\begin{proof}[of Theorem~\ref{theorem:pgd_strongconvex}]
\textbf{(i).}
To prove this part, we first establish bounds on the difference in function values
$
f(\bx^{(t)}) - f(\bx^*)
$  using the strong convexity (Definition~\ref{definition:scss_func}):
$$
\small
\begin{aligned}
f(\bx^{(t)}) &- f(\bx^*) \leq \nabla f(\bx^{(t)})^\top (\bx^{(t)} -\bx^*) - \frac{\alpha}{2} \normtwobig{\bx^{(t)} -\bx^*}^2 \\
&\stackrel{\dag}{=} \frac{1}{2\eta_t} \left( \normtwobig{\bx^{(t)} -\bx^*}^2 + \normtwo{\bx^{(t)} - \by^{(t+1)}}^2 - \normtwo{\by^{(t+1)} -\bx^*}^2 \right) - \frac{\alpha}{2} \normtwobig{\bx^{(t)} -\bx^*}^2 \\
&= \frac{1}{2\eta_t} \left( \normtwobig{\bx^{(t)} -\bx^*}^2 - \normtwo{\by^{(t+1)} -\bx^*}^2 \right) + \frac{\eta_t}{2} \normtwobig{\nabla f(\bx^{(t)})}^2 - \frac{\alpha}{2} \normtwobig{\bx^{(t)} -\bx^*}^2 \\
&\stackrel{+}{\leq} \frac{1}{2\eta_t} \left( \normtwobig{\bx^{(t)} -\bx^*}^2 - \normtwobig{\bx^{(t+1)} -\bx^*}^2 \right) + \frac{\eta_t}{2} \normtwobig{\nabla f(\bx^{(t)})}^2 - \frac{\alpha}{2} \normtwobig{\bx^{(t)} -\bx^*}^2 \\
&\leq \left( \frac{1}{2\eta_t} - \frac{\alpha}{2} \right) \normtwobig{\bx^{(t)} -\bx^*}^2 - \frac{1}{2\eta_t} \normtwobig{\bx^{(t+1)} -\bx^*}^2 + \frac{\eta_t L^2}{2},
\end{aligned}
$$
where the equality ($\dag$) follows from Theorem~\ref{theorem:funda_opt} and the update of PGD, the inequality ($+$) follows from Projection Property-II (Lemma~\ref{lemma:proj_prop2}), and the last inequality follows from Theorem~\ref{theorem:lipsc_equiv}.
By multiplying $t$ on both sides and substituting the stepsize $\eta_t$ by $\frac{2}{\alpha(t+1)}$, we get
$$
\small
\begin{aligned}
t\big(f(\bx^{(t)}) - f(\bx^*)\big) \leq \frac{L^2}{\alpha} + \frac{\alpha}{4} \left( t(t-1) \normtwobig{\bx^{(t)} -\bx^*}^2 - t(t+1) \normtwobig{\bx^{(t+1)} -\bx^*}^2 \right).
\end{aligned}
$$
Finally, we can find the upper bound of the function value by telescoping the sum of  $T$ steps:
$$
\small
\begin{aligned}
f&\Big( \sum_{t=1}^T \frac{2t}{T(T+1)} \bx^{(t)} \Big) 
\leq \sum_{t=1}^T \frac{2t}{T(T+1)} f(\bx^{(t)}) \\
&\leq \frac{2}{T(T+1)} \sum_{t=1}^T \left( t f(\bx^*) + \frac{L^2}{\alpha} + \frac{\alpha}{4} \left( t(t-1) \normtwobig{\bx^{(t)} -\bx^*}^2 - t(t+1) \normtwobig{\bx^{(t+1)} -\bx^*}^2 \right) \right) \\
&= \frac{2}{T(T+1)} \sum_{t=1}^T t f(\bx^*) + \frac{2L^2}{\alpha(T+1)} - \frac{\alpha}{2} \normtwo{\bx^{(T+1)} -\bx^*}^2 
\leq f(\bx^*) + \frac{2L^2}{\alpha(T+1)}.
\end{aligned}
$$

\paragraph{(ii).} Using the update rule,
$$
\small
\begin{aligned}
\normtwo{\bx^\toptone - \bx^*}^2 &= \normtwo{\projectS\left(\bx^\toptzero - \eta_t \nabla f(\bx^\toptzero)\right) - \projectS\left(\bx^*\right)}^2 
\stackrel{\dag}{\leq} \normtwo{\bx^\toptzero - \eta_t \nabla f(\bx^\toptzero) - \bx^*}^2 \\
&= \normtwo{\bx^\toptzero - \bx^*}^2 - 2\eta_t \innerproduct{\nabla f(\bx^\toptzero), \bx^\toptzero - \bx^*} + \eta_t^2 \normtwobig{\nabla f(\bx^\toptzero)}^2\\
&\stackrel{\ddag}{\leq} \left(1 - \alpha \eta_t\right) \normtwo{\bx^\toptzero - \bx^*}^2 - 2 \eta_t \left(f(\bx^\toptzero) - f(\bx^*)\right) + \eta_t^2 L^2.
\end{aligned}
$$
where the inequality ($\dag$) follows from the nonexpansiveness of the  projection operator (Theorem~\ref{theorem:proj_nonexpan}), and the inequality ($\ddag$) follows from the definition of Lipschitzness and  strong convexity  (Theorem~\ref{theorem:lipsc_equiv}, Definition~\ref{definition:scss_func}).
Plugging $\eta_t = \frac{2}{\alpha(t+1)}$ into the inequality and multiplying the above by $t$, this implies
$$
t\left(f(\bx^\toptzero) - f(\bx^*)\right) \leq \frac{\alpha t(t-1)}{4} \normtwo{\bx^\toptzero - \bx^*}^2 - \frac{\alpha t(t+1)}{4} \normtwo{\bx^\toptone - \bx^*}^2 + \frac{t}{\alpha(t+1)} L^2.
$$
Since $f(\bx^\toptzero) \geq \fbest^{(T)}$ for all $t = \{1,2,\ldots,T\}$, telescoping the sum yields
$$
\left.
\begin{aligned}
&\left(\sum_{t=1}^T t\right)\left(\fbest^{(t)} - f(\bx^*)\right) 
\leq \sum_{t=1}^T t\left(f(\bx^\toptzero) - f(\bx^*)\right) \\
&\leq - \frac{\alpha}{4}T(T+1)\normtwo{\bx^{T+1} - \bx^*}^2 + \frac{L^2}{\alpha} \sum_{t=1}^T \frac{t}{t+1} \leq \frac{L^2 T}{\alpha}
\end{aligned}
\right\}
\,\,\implies\,\, 
\fbest^{(T)} - f(\bx^*) \leq \frac{2L^2}{\alpha(T+1)}.
$$
Note that $\fbest^{(T)} = f(\bx^{(i_T)})$. 
Since $f$ is strongly convex, SC Property-II (Theorem~\ref{theorem:exi_close_sc}) implies that 
$$
\frac{\alpha}{2} \normtwo{\bx^{(i_T)} - \bx^{*}}^{2} \leq \fbest^{(T)} - f(\bx^*) \leq \frac{2 L^{2}}{\alpha(T+1)}
\quad\implies\quad 
\normtwo{\bx^{(i_T)} - \bx^*} \leq \frac{2 L}{\alpha\sqrt{T+1}}.
$$
This completes the proof.
\end{proof}

Similar to Theorem~\ref{theorem:gd_sc_ss}, we can also demonstrate the same result for the constrained scenario using projected gradient descent for strongly convex and smooth functions.
\begin{theoremHigh}[PGD for SC and SS Functions]\label{theorem:pgd_sc_ss}
Let $f: \sS \rightarrow \real$ be a proper differentiable, $\alpha$-strongly convex, and $\beta$-smooth function, where $\sS\subseteq\real^n$ is a closed and convex set. 
Let $\{\bx^\toptzero\}_{t > 0}$ be the sequence generated by the projected gradient descent method (Algorithm~\ref{alg:pgd_gen}) for solving  problem (P2) with a constant stepsize $\eta=\frac{1}{\beta}$.
Assume  the function has a minimizer  $\bx^*$. Then,
$$
\normtwo{\bx^{(T+1)} - \bx^*}^2 \leq \exp\left(-T\frac{\alpha}{\beta}\right) \normtwo{\bx^{(1)} - \bx^*}^2.
$$
This also shows  that $\normtwo{\bx^{(T)} -\bx^*}^2 \leq \epsilon$ after at most $T = \mathcalO\left(\frac{\beta}{\alpha} \ln \frac{1}{\epsilon}\right)$ steps.
\end{theoremHigh}
The proof of Theorem~\ref{theorem:pgd_sc_ss} follows exactly as in Theorem~\ref{theorem:gd_sc_ss}, substituting the appropriate projected gradient descent update for the standard gradient descent update using the descent lemma (Lemma~\ref{lemma:gdupd_sm_sconv}).
\begin{proof}[of Theorem~\ref{theorem:pgd_sc_ss}]
Let $g(\bx) \triangleq \beta\left(\bx - \projectS\left(\bx - \frac{1}{\beta} \nabla f(\bx)\right)\right)$. Then,
$$
\small
\begin{aligned}
\normtwo{\bx^{(t+1)} - \bx^*}^2 
&= \normtwo{\bx^{(t)} -\frac{1}{\beta}g(\bx^{(t)}) - \bx^*}^2 
= \normtwo{\bx^{(t)} - \bx^*}^2 - \frac{2}{\beta} g(\bx^{(t)})^\top (\bx^{(t)} - \bx^*) + \frac{1}{\beta^2} \normtwo{g(\bx^{(t)})}^2 \\
&\stackrel{\dag}{\leq} \left(1 - \frac{\alpha}{\beta}\right) \normtwo{\bx^{(t)} - \bx^*}^2 
\leq \left(1 - \frac{\alpha}{\beta}\right)^t \normtwo{\bx^{(1)} - \bx^*}^2 
\leq \exp\left(-t\frac{\alpha}{\beta}\right) \normtwo{\bx^{(1)} - \bx^*}^2,
\end{aligned}
$$
where inequality ($\dag$)  follows from  the descent lemma (\ref{descent_lem_pgdscss} in Lemma~\ref{lemma:gdupd_sm_sconv}).
\end{proof}

Following \citet{jain2017non}, we provide an alternative proof for the convergence result of PGD under SC and SS. 
This proof also highlights the roles of SC and SS in optimization algorithms.
\begin{proof}[of Theorem~\ref{theorem:pgd_sc_ss}: Alternative Way] Alternatively, we may prove the theorem in the following steps.
\paragraph{Strong smoothness: the gradient does not vanish too fast.} We use it to show that PGD always makes significant progress in each iteration:
$$
\begin{aligned}
& f(\bx^\toptone) - f(\bx^\toptzero) \leq \innerproduct{\nabla f(\bx^\toptzero), \bx^\toptone - \bx^\toptzero} + \frac{\beta}{2} \normtwo{\bx^\toptzero - \bx^\toptone}^{2} \\
&= \innerproduct{\nabla f(\bx^\toptzero), \bx^\toptone - \bx^{*}} + \innerproduct{\nabla f(\bx^\toptzero), \bx^{*} - \bx^\toptzero} + \frac{\beta}{2} \normtwo{\bx^\toptzero - \bx^\toptone}^{2} \\
&= \frac{1}{\eta} \innerproduct{\bx^\toptzero - \by^\toptone, \bx^\toptone - \bx^*} + \innerproduct{\nabla f(\bx^\toptzero), \bx^{*} - \bx^\toptzero} + \frac{\beta}{2} \normtwo{\bx^\toptzero - \bx^\toptone}^{2},
\end{aligned}
$$
where the last equality follows from the update $\by^{(t+1)} \leftarrow \bx^{(t)} - \eta \nabla f(\bx^{(t)})$.

\paragraph{Projection step.}To simplify the expression, which contains the term $\by^\toptone$, we  apply Projection Property-I (Lemma~\ref{lemma:proj_prop1}) to eliminate it, yielding $\innerproduct{\bx^\toptone-\bx^*, \bx^\toptone-\by^\toptone}\leq 0$:
$$
\begin{aligned}
\innerproduct{\bx^\toptzero - \by^\toptone, \bx^\toptone -\bx^*} &\leq \innerproduct{\bx^\toptzero - \bx^\toptone, \bx^\toptone -\bx^*} 
=-\innerproduct{\bx^\toptone-\bx^\toptzero, \bx^\toptone -\bx^*}\\
&= \frac{1}{2} \left(\normtwo{\bx^\toptzero -\bx^*}^2 - \normtwo{\bx^\toptzero - \bx^\toptone}^2 - \normtwo{\bx^\toptone -\bx^*}^2\right).
\end{aligned}
$$
where the last equality follows from the fundamental theorem of optimization (Theorem~\ref{theorem:funda_opt}).
Using $\eta = 1/\beta$ and combining the above results, we get
$$
\begin{aligned}
& f(\bx^\toptone) - f(\bx^\toptzero) \leq \innerproduct{\nabla f(\bx^\toptzero), \bx^{*} - \bx^\toptzero} + \frac{\beta}{2} \left(\normtwo{\bx^\toptzero - \bx^{*}}^{2} - \normtwo{\bx^\toptone - \bx^{*}}^{2}\right).
\end{aligned}
$$

\paragraph{Strong convexity.} The above expression is suitable  for a telescoping sum, except for the inner product term. Fortunately, this can be eliminated using strong convexity.
$$
\innerproduct{\nabla f(\bx^\toptzero), \bx^{*} - \bx^\toptzero} \leq f\left(\bx^{*}\right) - f(\bx^\toptzero) - \frac{\alpha}{2} \normtwo{\bx^\toptzero - \bx^{*}}^{2}
$$
Combining this with the previous results gives us
$$
\begin{aligned}
&0\leq   f(\bx^\toptone) - f(\bx^*) \leq \frac{\beta - \alpha}{2} \normtwo{\bx^\toptzero -\bx^*}^2 - \frac{\beta}{2} \normtwo{\bx^\toptone -\bx^*}^2\\
&\implies \ \normtwo{\bx^\toptone -\bx^*}^2 \leq (1-\frac{\alpha}{\beta}) \normtwo{\bx^\toptzero -\bx^*}^2 
\leq 
\exp\left(-t\frac{\alpha}{\beta}\right) \normtwo{\bx^{(1)} -\bx^*}^2.
\end{aligned}
$$
This completes the proof.
\end{proof}

All the above convergence results for PGD are established under the assumption that $f$  is convex and the set $\sS$ is also convex. Analysis of PGD for non-convex functions is more restricted. Several approaches address this limitation, the simplest being to transform the constraint set into a convex one, potentially by taking its convex hull, which is a common approach in relaxation methods. However, a less drastic alternative widely used in non-convex optimization literature exists. Note that the convergence results for the PGD algorithm in Theorem~\ref{theorem:pgd_sc_ss} do not require the objective function to be convex (or strongly convex/strongly smooth) over all of $\real^n$, but only over a subset of the space. When the function is restricted strongly convex and smooth (Definition~\ref{definition:res_scss_func}), convergence is also guaranteed.

\begin{theoremHigh}[PGD for RSC and RSS Functions (Non-Convex)]\label{theorem:pgd_rsc_rss}
Let $ f:\sS\subseteq\real^n\rightarrow \real $ be a proper closed (possibly non-convex) function satisfying the $\alpha$-RSC and $\beta$-RSS properties (Definition~\ref{definition:res_scss_func}) over a (possibly non-convex) constraint set $\sS$ with $\beta / \alpha < 2$.
Suppose $\bx^*$ is an optimizer such that $\nabla f(\bx^*)=\bzero$.
Let Algorithm~\ref{alg:pgd_gen} be executed with a constant stepsize $\eta = \frac{1}{\beta}$. 
Then,
$$
f(\bx^{(T+1)}) - f(\bx^*)\leq \exp\left(-T\frac{2\alpha-\beta}{\alpha}\right) \big(f(\bx^{(1)}) - f(\bx^*)\big)
$$
This also shows that $ f(\bx^{(T+1)}) -f(\bx^*)  \leq  \epsilon $ after at most  $ T = \mathcalO\left( \frac{\alpha}{2\alpha - \beta} \ln \frac{1}{\epsilon} \right) $ steps.
\end{theoremHigh}

This result holds even when the step length is chosen within a range that is large enough but still smaller than $1/\beta$. However, setting $\eta = \frac{1}{\beta}$  simplifies the proof and focuses attention on key concepts.

\begin{proof}[of Theorem~\ref{theorem:pgd_rsc_rss}]
We will use RSS to track the global convergence of the algorithm and RSC to locally assess the progress made by the algorithm in each iteration.

\paragraph{Restricted strong smoothness.}
Since both $\bx^\toptzero, \bx^\toptone \in \sS$ due to the projection steps and $\eta = 1/\beta$, we apply the $\beta$-RSS property (Definition~\ref{definition:res_scss_func}):
$$
\begin{aligned}
f(\bx^\toptone) - f(\bx^\toptzero) 
&\leq \innerproduct{\nabla f(\bx^\toptzero), \bx^\toptone - \bx^\toptzero} + \frac{\beta}{2} \normtwo{ \bx^\toptzero - \bx^\toptone}^2\\
&= \frac{1}{\eta} \innerproduct{\bx^\toptzero - \by^\toptone, \bx^\toptone - \bx^\toptzero} + \frac{\beta}{2} \normtwo{\bx^\toptzero - \bx^\toptone }^2\\
&= \frac{\beta}{2} \left( \normtwo{\bx^\toptone - \by^\toptone}^2 - \normtwo{\bx^\toptzero - \by^\toptone}^2 \right),
\end{aligned}
$$
where the last equality follows from the fundamental theorem of optimization (Theorem~\ref{theorem:funda_opt}) with $\ba\triangleq \bx^\toptzero - \by^\toptone$ and $\bb\triangleq\bx^\toptzero-\bx^\toptone$.

\paragraph{Projection update rule.}
To handle the term $\by^\toptone$, unlike before, we cannot apply Projection Properties I, II, or III as non-convex projections may not satisfy them. Instead, we use Projection Property-O (Lemma~\ref{lemma:proj_prop0}), which applies to all projections:
\begin{equation}\label{equation:pgd_rsc_rss1}
\begin{aligned}
	f(\bx^\toptone) - f(\bx^\toptzero) 
	&\leq \frac{\beta}{2} \left( \normtwo{\bx^* - \by^\toptone}^2 - \normtwo{\bx^\toptzero - \by^\toptone}^2 \right)\\
	&\stackrel{\dag}{=} \frac{\beta}{2} \left( \normtwo{\bx^* - \bx^\toptzero}^2 + 2 \innerproduct{ \bx^* - \bx^\toptzero, \bx^\toptzero - \by^\toptone} \right)\\
	&= \frac{\beta}{2} \normtwo{\bx^* - \bx^\toptzero}^2 + \innerproduct{\bx^* - \bx^\toptzero, \nabla f(\bx^\toptzero)},
\end{aligned}
\end{equation}
where the  equality ($\dag$) follows from the fundamental theorem of optimization (Theorem~\ref{theorem:funda_opt}) with $\ba\triangleq \bx^* - \bx^\toptzero$ and $\bb\triangleq\by^\toptone - \bx^\toptzero$.

\paragraph{Restricted strong convexity.}
Since both $\bx^\toptzero, \bx^* \in \sS$, we apply the $\alpha$-RSC property (Definition~\ref{definition:res_scss_func}) to them. However, we do so in two ways:
$$
\begin{aligned}
f(\bx^*) - f(\bx^\toptzero) 
&\geq \innerproduct{\nabla f(\bx^\toptzero), \bx^* - \bx^\toptzero } + \frac{\alpha}{2} \normtwo{\bx^\toptzero - \bx^*}^2;\\
f(\bx^\toptzero) - f(\bx^*) 
&\geq \innerproduct{\nabla f(\bx^*), \bx^\toptzero - \bx^*} + \frac{\alpha}{2} \normtwo{\bx^\toptzero - \bx^* }^2 = \frac{\alpha}{2} \normtwo{\bx^\toptzero - \bx^*}^2,
\end{aligned}
$$
where in the second inequality, we use $\nabla f(\bx^*) = \bzero$. 
By manipulating these equations, we get:
$$
\innerproduct{\nabla f(\bx^\toptzero), \bx^* - \bx^\toptzero} + \frac{\beta}{2} \normtwo{\bx^* - \bx^\toptzero }^2 \leq \left( 2 - \frac{\beta}{\alpha} \right) \left( f(\bx^*) - f(\bx^\toptzero) \right).
$$
Plugging this into \eqref{equation:pgd_rsc_rss1}:
\begin{equation}\label{equation:pgd_rsc_rss2}
f(\bx^\toptone) - f(\bx^\toptzero) \leq \left( 2 - \frac{\beta}{\alpha} \right) \left( f(\bx^*) - f(\bx^\toptzero) \right).
\end{equation}
Equation~\eqref{equation:pgd_rsc_rss2} indicates that the larger the gap between $f(\bx^*)$ and $f(\bx^\toptzero)$, the greater the decrease in the objective value when moving from $\bx^\toptzero$ to $\bx^\toptone$. The form of the result is also quite fortunate as it assures us that we will cover a constant fraction $\left( 2 - \frac{\beta}{\alpha} \right)$ of the remaining ``distance'' to $\bx^*$ at each step. Rearranging this gives
$$
f(\bx^\toptone) - f(\bx^*)\leq (\kappa - 1) \big(f(\bx^\toptzero) - f(\bx^*)\big),
$$
where $\kappa \triangleq \beta / \alpha$. Note that we always have $\kappa \geq 1$ and by assumption $\kappa = \beta / \alpha < 2$, so that we always have $\kappa - 1 \in [0, 1)$. This proves the result using the fact that $1 - r \leq \exp(-r)$ for all $r \in \real$.
\end{proof}

We observe that the condition number once again plays a crucial role in determining the convergence rate of the algorithm, this time for a non-convex problem. However, it is noted that the condition number is defined differently here, utilizing the RSC and RSS  constants instead of the standard SC  and SS  constants used in Section~\ref{section:gd_basic}.

Readers will notice that while there was no restriction on the condition number $\kappa$ in the analysis of the GD or PGD algorithms (refer to Theorem~\ref{theorem:gd_sc_ss} or Theorem~\ref{theorem:pgd_sc_ss}), the analysis for PGD with RSC and RSS functions does require $\kappa < 2$. It is worth noting that this restriction can be relaxed for specific problems; however, doing so significantly complicates the analysis. Addressing this issue in a general context is beyond the scope of this book; see, for example, \citet{jain2017non} for more details.

\begin{algorithm}[h] 
\caption{Proximal Gradient Method}
\label{alg:prox_gd_gen}
\begin{algorithmic}[1] 
\Require A function $f(\bx)$ and a closed convex function $g$ (usually non-smooth) satisfying (A1) and (A2) in Assumption~\ref{assumption:proximal_grad}; 
\State {\bfseries Input:}  Initialize $\bx^{(1)}$;
\For{$t=1,2,\ldots$}
\State Pick a stepsize $\eta_t$;
\State $\by^{(t+1)} \leftarrow \bx^{(t)} - \eta_t \nabla f(\bx^{(t)})$;
\State $\bx^{(t+1)} \leftarrow \prox_{\eta_t g}(\by^{(t+1)}) \triangleq \mathcalT_{L_t}^{f,g}(\bx^\toptzero)$;
\EndFor
\State (Output Option 1) Output  $\bx_{\text{final}}\leftarrow \bx^{(T)}$;
\State (Output Option 2) Output  $\bx_{\text{avg}}\leftarrow \frac{1}{T}(\sum_{t=1}^{t}\bx^{(t)})$ or $\sum_{t=1}^{T} \frac{2t}{T(T+1)} \bx^{(t)}$;
\State (Output Option 3) Output  $\bx_{\text{best}}\leftarrow \argmin_{t\in\{1,2,\ldots,T\}} f(\bx^{(t)})$;
\end{algorithmic} 
\end{algorithm}
\section{Proximal Gradient Method}\label{section:prox_gd}
We have introduced similarity between  projection and proximal operators in Section~\ref{section:proj_prox_sep}.
Noting the equivalence $\prox_{\indicatorS}(\by) = \projectS(\by)$, we can replace the projection step in Algorithm~\ref{alg:pgd_gen} with 
$$
\left\{\bx^{(t+1)} \leftarrow \mathcalP_{\sS}(\by^{(t+1)})\right\}
\qquad\leadsto \qquad
\left\{\bx^{(t+1)} \leftarrow \prox_{\indicatorS}(\by^{(t+1)})\right\}.
$$
When discussing the convergence results for the PGD approach, we require the underlying set $\sS$ to be closed convex~\footnote{Again, we should note that the notion of projection operator is only interesting for closed sets. Unless otherwise stated, when we discuss the projection operator, we assume inherently that the set is closed. Similarly, we assume the function is proper when we discuss the proximal operator.} such that the projection operator satisfies Projection Properties-O, I, II, III (Lemmas~\ref{lemma:proj_prop0}$\sim$\ref{lemma:proj_prop3}), which are essential for deriving convergence results for PGD.
The \textit{proximal gradient descent method} (Algorithm~\ref{alg:prox_gd_gen}) extends the projection step by 
\begin{equation}\label{equation:prox_obj1}
\left\{\bx^{(t+1)} \leftarrow \prox_{\indicatorS}(\by^{(t+1)})\right\}
\qquad\leadsto \qquad
\left\{\bx^{(t+1)} \leftarrow \prox_{\textcolor{mylightbluetext}{\eta_t g}}(\by^{(t+1)})\right\},
\end{equation}
where $\eta_t$ is the stepsize at $t$-th iteration (we will see its role shortly).
Similarly, we require the function $g$ to be closed and convex such that the proximal steps satisfy  Proximal Properties-O, I, II, III, IV (Lemmas~\ref{lemma:prox_prop0}$\sim$\ref{lemma:prox_prop3}).
Apparently, when $g=\indicatorS$, the method reduces to a PGD approach.

In PGD, the update rule from $t$-th to $(t+1)$-th iteration is 
\begin{center}
\framebox{
\begin{minipage}{0.95\textwidth}
\begin{equation}\label{equation:pgd_decom}
\small
\textbf{(PGD)}:\quad 
\begin{aligned}
\bx^\toptone 
&\leftarrow \mathop{\argmin}_{\bx\in\sS} \normtwo{\bx - \left(\bx^{(t)} - \eta_t \nabla f(\bx^{(t)})\right)}^2\\
&=\mathop{\argmin}_{\bx\in\sS} 
\left\{
\frac{1}{2\eta_t} \normtwo{\bx-\bx^\toptzero}^2 + f(\bx^\toptzero)+ \innerproduct{\nabla f(\bx^\toptzero), \bx-\bx^\toptzero} 
\right\},
\end{aligned}
\end{equation}
\end{minipage}
}
\end{center}
where we add a constant $f(\bx^\toptzero)$ term in the second form since the optimization is over $\bx$.
Therefore, the update rule in PGD can be seen as a minimization of the sum of the linearization of the smooth part around the current iterate plus a quadratic term. The optimization of the linearization ensures significant progress, while the quadratic term acts as regularization, ensuring the update stays within a neighborhood of the current iterate $\bx^\toptzero$.
The proximal gradient method add a non-smooth term $g$ into the objective function in \eqref{equation:pgd_decom} while consider a unconstrained problem:
\begin{subequations}\label{equation:prox_decom}
\begin{tcolorbox}[colback=white,colframe=black]
\begin{minipage}{1\textwidth}
\small
\begin{align}
\bx^\toptone 
&\leftarrow\mathop{\argmin}_{\bx\in\textcolor{mylightbluetext}{\real^n}} 
\left\{
\frac{1}{2\eta_t} \normtwo{\bx-\bx^\toptzero}^2 + f(\bx^\toptzero)+ \innerproduct{\nabla f(\bx^\toptzero), \bx-\bx^\toptzero }
+\textcolor{mylightbluetext}{g(\bx)}
\right\}\\
\textbf{(Prox)}:\qquad\,\,\,&=\mathop{\argmin}_{\bx\in\textcolor{mylightbluetext}{\real^n}} \eta_t g(\bx) +\frac{1}{2}\normtwo{\bx-\left( \bx^\toptzero - \eta_t\nabla f(\bx^\toptzero) \right)}\\
&=\prox_{\eta_t g}\left(\bx^\toptzero - \eta_t\nabla f(\bx^\toptzero)\right)
\triangleq \mathcalT_{L_t}^{f,g}(\bx^\toptzero), \label{equation:prox_decom3}
\end{align}
\end{minipage}
\end{tcolorbox}
\noindent where $L_t \triangleq\frac{1}{\eta_t}$, and $\mathcalT_{L}^{f,g}(\bx)\triangleq\prox_{\frac{1}{L} g}\left(\bx - \frac{1}{L}\nabla f(\bx)\right)$ denotes the \textit{prox-grad operator}. 
Define 
$$
\mathcalG_{L}^{f,g}(\bx)\triangleq L\big(\bx-\mathcalT_{L}^{f,g}(\bx)\big).
$$
Then the update in \eqref{equation:prox_decom3} can be equivalently denoted as 
\begin{equation}
\bx^\toptone 
\leftarrow 
\prox_{\eta_t g}\left(\bx^\toptzero - \eta_t\nabla f(\bx^\toptzero)\right)
\triangleq \mathcalT_{L_t}^{f,g}(\bx^\toptzero)
\triangleq \bx^\toptzero - \frac{1}{L_t} \mathcalG_{L_t}^{f,g}(\bx^\toptzero).
\end{equation}
\end{subequations}
The last form resembles a standard gradient descent update. 
Therefore, the term $\mathcalG_{L_t}^{f,g}(\bx^\toptzero)$ is called a \textit{gradient mapping}. 
When the identities of $f$ and $g$ are clear from the context, we will omit the superscripts and write $\mathcalT_{L_t}(\cdot)$ or $\mathcalG_{L_t}(\cdot)$ instead of $\mathcalT_{L_t}^{f,g}(\cdot)$ or $\mathcalG_{L_t}^{f,g}(\cdot)$.

\index{Gradient mapping}\index{Prox-grad operator}

Equation~\eqref{equation:prox_decom} explains the scaling term in \eqref{equation:prox_obj1}.
This demonstrates that the proximal gradient method aims to optimize the composite function:
$$
\text{(P3)}:\qquad \min \{F(\bx) \triangleq f(\bx)+g(\bx)\},
$$
where $f$ and $g$ satisfy (A1) and (A2) in Assumption~\ref{assumption:proximal_grad}.

\begin{assumption}[Proximal Gradient Method]\label{assumption:proximal_grad}
In problem (P3), we often assume 
\begin{itemize}
\item[(A1)] $ g: \real^n \rightarrow (-\infty, \infty] $ is proper, closed, and convex.~\footnote{A few examples are discussed in the paragraph below Definition~\ref{definition:opt_probs_all}.}
\item[(A2)] $ f: \real^n \rightarrow (-\infty, \infty] $ is proper and closed, $ \dom(f) $ is convex, $ \dom(g) \subseteq \textcolor{mylightbluetext}{\interior}(\dom(f)) $, and $ f $ is $ \beta $-smooth over $ \textcolor{mylightbluetext}{\interior}(\dom(f)) $.
\end{itemize}
The iterative methods introduced will use the function $g$ to generate a ``candidate" update. Therefore, we assume the function $g$ is proper closed and convex.
We do not assume the function $g$ has a compact domain; since the proximal operator has a quadratic term such that the function is strongly convex. This ensures the update has a minimum (Lemma~\ref{lemma:prox_prop3}, Theorem~\ref{theorem:exi_close_sc}).
\end{assumption}

Assumption (A1) ensures that $ g(\bx) $ is proper, closed, and convex, which guarantees that the proximal step has a unique solution (Lemma~\ref{lemma:prox_prop3}).
The domain compatibility  $ \dom(g) \subseteq \interior(\dom(f)) $ in (A2) ensures that the proximal operator is applied within the region where $ f(\bx) $ is differentiable, avoiding issues at the boundary of $ \dom(f) $. The closedness of $f$ and $g$ ensures that $F\triangleq f+g$ is also closed such that iterative algorithms can find stationary points (see discussion in Figure~\ref{fig:low_nonlowersemis}).

Several stopping criteria for gradient descent were discussed in Section~\ref{section:gradient-descent-all}. The descent lemma for the proximal gradient method shows that the gradient mapping can serve as a stopping criterion. The following theorem establishes that the stationary point of (P3) corresponds to the condition $\mathcalG_{L}^{f,g}(\bx)=\bzero$ (Theorem~\ref{theorem:opt_cond_p3}).

\begin{theoremHigh}[Optimality Condition using Gradient Mapping]\label{theorem:opt_gdmap}
Let $ f $ be a proper \textcolor{black}{closed} and $\beta$-smooth function,  and let $ g $ be a proper closed and convex function satisfying assumption (A1) and (A2) in Assumption~\ref{assumption:proximal_grad}.
Let $L\triangleq\beta$.
Then,
\begin{enumerate}[(i)]
\item When  $ g_0(\bx) = 0 $ for all $\bx\in\real^n$, $ \mathcalG_L^{f,g_0}(\bx) = \nabla f(\bx) $ for any $ \bx \in \interior(\dom(f)) $.
\item \label{opt_gdmap_item2} For $ \bx^* \in \interior(\dom(f)) $, it holds that $ \mathcalG_L^{f,g}(\bx^*) = \bzero $ if and only if $ \bx^* $ is a stationary point of problem (P3).
\item  Suppose additionally that $ f $ is convex. Then, for $ \bx^* \in \dom(g) $, $ \mathcalG_L^{f,g}(\bx^*) = \bzero $ if and only if $ \bx^* $ is an optimal solution of problem (P3).
\end{enumerate}
\end{theoremHigh}
\begin{proof}[of Theorem~\ref{theorem:opt_gdmap}]
\textbf{(i).} Since $ \prox_{\frac{1}{L}g_0}(\by) = \by $ for all $ \by \in \real^n $ if $g_0(\bx)=0$ for all $\bx\in\real^n$, it follows that
$
\mathcalG_L^{f,g_0}(\bx) = L(\bx - \mathcalT_L^{f,g_0}(\bx)) 
= L\left(\bx - \left(\bx - \frac{1}{L}\nabla f(\bx)\right)\right) = \nabla f(\bx).
$
\paragraph{(ii).} $ \mathcalG_L^{f,g}(\bx^*) = \bzero $ if and only if $ \bx^* = \prox_{\frac{1}{L}g}\left(\bx^* - \frac{1}{L}\nabla f(\bx^*)\right) $. By  Proximal Property-I (Lemma~\ref{lemma:prox_prop1}), the latter relation holds if and only if
$$
\bx^* - \frac{1}{L}\nabla f(\bx^*) - \bx^* \in \frac{1}{L}\partial g(\bx^*)
\qquad \iff\qquad 
-\nabla f(\bx^*) \in \partial g(\bx^*),
$$
which is exactly the condition for stationarity (Theorem~\ref{theorem:opt_cond_p3}).
\paragraph{(iii).} When $f$ is convex,  Theorem~\ref{theorem:opt_cond_p3} shows $\bx^*$ is an optimal solution of (P3).
\end{proof}

Previously, we observed that descent lemmas (e.g., Lemma~\ref{lemma:gdupd_sm_sconv}) are crucial for proving convergence results. For the proximal gradient descent method, there is a similar form of the descent lemma.
\begin{lemma}[Descent Lemma for Proximal Gradient Method]\label{lemma:des_lem_prox}
Let $ f $ be a proper \textcolor{black}{closed} and $\beta$-smooth function, and let  $ g $ be a proper closed and convex function satisfying assumption (A1) and (A2) in Assumption~\ref{assumption:proximal_grad}. Consider the composite function $ F \triangleq f + g $ and the prox-grad operator $ \mathcalT_{L} \triangleq \mathcalT_{L}^{f,g} $. 
Then, for any $ \bx \in \interior(\dom(f)) $ and $ L \in \left(\frac{\beta}{2}, \infty\right) $, the following inequality holds:
\begin{equation}\label{equation:des_lem_prox_res1}
F(\mathcalT_{L}(\bx)) -F(\bx) \leq \frac{\frac{\beta}{2} -L}{L^2} \normtwo{\mathcalG_{L}^{f,g}(\bx)}^2 \leq  0, 
\end{equation}
where $ \mathcalG_{L}^{f,g} : \interior(\dom(f)) \rightarrow \real^n $ is the operator defined by $\mathcalG_{L}^{f,g}(\bx)\triangleq L\big(\bx-\mathcalT_{L}(\bx)\big) $ for all $ \bx \in \interior(\dom(f)) $.

In addition, for any $\by\in\real^n$, $\bx \in \interior(\dom(f))$, consider the \textit{proximal gap} $\mathcalD_f(\by, \bx) \triangleq f(\by) - f(\bx) - \innerproduct{\nabla f(\bx), \by - \bx }$~\footnote{When $f$ is convex, this is also known as the Bregman distance (Definition~\ref{definition:breg_dist}).} and  $L\triangleq\beta$. Then it holds that 
\begin{equation}\label{equation:des_lem_prox_res2}
F(\mathcalT_{L}(\bx)) - F(\by)  \leq   \frac{L}{2} \normtwo{\by-\bx}^2 -\frac{L}{2}\normtwo{\by - \mathcalT_{L}(\bx)}^2 - \mathcalD_f(\by, \bx).~\footnote{Similar to \ref{descent_lem_gdscss} and \ref{descent_lem_pgdscss} in the descent lemma for GD and PGD (Lemma~\ref{lemma:gdupd_sm_sconv}), we can set $\by$ to the optimal point $\by\triangleq\bx^*$ to obtain convergence results.}
\end{equation}
When $\by=\bx$, \eqref{equation:des_lem_prox_res2} reduces to \eqref{equation:des_lem_prox_res1} by letting $L\triangleq\beta$.
\end{lemma}
\begin{proof}[of Lemma~\ref{lemma:des_lem_prox}]
\textbf{Equation~\eqref{equation:des_lem_prox_res1}.}
For simplicity, denote $ \bx^+ \triangleq \mathcalT_{L}(\bx) $. By the smoothness of $f$ (Definition~\ref{definition:scss_func}),
$$
f(\bx^+) \leq f(\bx) + \innerproduct{\nabla f(\bx), \bx^+ - \bx} + \frac{\beta}{2} \normtwo{\bx - \bx^+}^2. 
$$
By the Proximal Property-I (Lemma~\ref{lemma:prox_prop1}), since $ \bx^+ = \prox_{\frac{1}{L}g}(\bx - \frac{1}{L}\nabla f(\bx)) $, we have
$$
\begin{aligned}
&\innerproduct{\bx - \frac{1}{L}\nabla f(\bx) - \bx^+, \bx - \bx^+} \leq \frac{1}{L} g(\bx) - \frac{1}{L} g(\bx^+)\\
&\quad\implies\quad  g(\bx^+)    \leq  g(\bx)-L \normtwo{\bx^+ - \bx}^2  -\innerproduct{\nabla f(\bx), \bx^+ - \bx} .
\end{aligned}
$$
Combining the two inequalities yields
$$
f(\bx^+) + g(\bx^+) \leq f(\bx) + g(\bx) + \left(-L + \frac{\beta}{2}\right) \normtwo{\bx^+ - \bx}^2.
$$
Hence, taking into account the definitions of $ \bx^+ $, $ \mathcalG_{L}^{f,g}(\bx) $, and the identities $ F(\bx) = f(\bx) + g(\bx) $, $ F(\bx^+) = f(\bx^+) + g(\bx^+) $, the desired result follows.

\paragraph{Equation~\eqref{equation:des_lem_prox_res2}.} For the second part, consider the function
$$
\varphi(\bu) \triangleq f(\bx) + \innerproduct{\nabla f(\bx), \bu - \bx} + g(\bu) + \frac{L}{2} \normtwo{\bu - \bx}^2. 
$$
which is the first objective function in \eqref{equation:prox_decom} for each update of the proximal gradient method.
Since $\varphi$ is an $L$-strongly convex function and $ \mathcalT_{L}(\bx) = \arg\min_{\bu \in \real^n} \varphi(\bu) $, it follows by Theorem~\ref{theorem:exi_close_sc}(ii) that
\begin{equation}\label{equation:des_lem_prox1}
\varphi(\by) - \varphi(\mathcalT_L(\bx)) \geq \frac{L}{2} \normtwo{\by - \mathcalT_L(\bx)}^2. 
\end{equation}
This implies, by the smoothness of $f$, that
\begin{equation}\label{equation:des_lem_prox2}
\begin{aligned}
\varphi(\mathcalT_L(\bx)) &= f(\bx) + \innerproduct{\nabla f(\bx), \mathcalT_L(\bx) - \bx} + \frac{L}{2} \normtwo{\mathcalT_L(\bx) - \bx}^2 + g(\mathcalT_L(\bx)) \\
&\geq f(\mathcalT_L(\bx)) + g(\mathcalT_L(\bx)) = F(\mathcalT_L(\bx)).
\end{aligned}
\end{equation}
Thus, combining \eqref{equation:des_lem_prox1} and \eqref{equation:des_lem_prox2} yields that 
$$
\begin{aligned}
\underbrace{f(\bx) + \innerproduct{\nabla f(\bx), \by - \bx} + g(\by) + \frac{L}{2} \normtwo{\by - \bx}^2}_{=\varphi(\by)} - F(\mathcalT_L(\bx)) &\geq \frac{L}{2} \normtwo{\by - \mathcalT_L(\bx)}^2, \quad \forall \text{$\by \in \real^n $}.
\end{aligned}
$$
This proves the desired result.
\end{proof}

A direct consequence of Lemma~\ref{lemma:des_lem_prox} is the monotonicity of  function values and  iterates in proximal gradient methods.
\begin{lemma}[Monotonicity and Sufficient Decrease of Proximal Gradient]\label{lemma:mont_prox}
Let $ f $ be a proper \textcolor{black}{closed} and $\beta$-smooth function, and let $ g $ be a proper closed and convex function satisfying assumption (A1) and (A2) in Assumption~\ref{assumption:proximal_grad}.
Let $ \{\bx^\toptzero\}_{t > 0} $ be the sequence generated by the proximal gradient method (Algorithm~\ref{alg:prox_gd_gen}) for solving problem (P3). Then, for any $L_t\triangleq L \in\left(\frac{\beta}{2}, \infty\right)\,\forall t>0$, \eqref{equation:des_lem_prox_res1} shows that
\begin{equation}\label{equation:mont_prox_1}
\left(
\begin{aligned}
	\textbf{Monotonicity and }\\
	\textbf{Sufficient Decrease}
\end{aligned}
\right): \qquad
F(\bx^\toptone) - F(\bx^\toptzero) \leq\frac{\frac{\beta}{2} -L}{L^2} \normtwo{\mathcalG_{L}^{f,g}(\bx^\toptzero)}^2 \leq 0.
\end{equation}
Since $\mathcalG_{L}^{f,g}(\bx)\neq  \bzero$ by Theorem~\ref{theorem:opt_gdmap}\ref{opt_gdmap_item2} and equals to $\bzero$ only for stationary points, this also shows that 
\begin{itemize}
	\item $F(\bx^\toptzero) - F(\bx^\toptone) \geq 0$.
	\item $F(\bx^\toptzero) > F(\bx^\toptone)$  if $\bx^\toptzero$  is not a stationary point of problem (P3).
\end{itemize}

Suppose additionally  that $f$ is convex and $L_t \triangleq L \triangleq \beta$. Then, for any optimal point $\bx^*$ of $F=f+g$, it holds that 
$$
\begin{aligned}
\textbf{(\textbf{\fejer} \text{monotonicity})}:\qquad &\normtwo{\bx^\toptone - \bx^*} \leq \normtwo{\bx^\toptzero - \bx^*};\\
\textbf{(Gradient mapping)}:\qquad & \normtwo{\mathcalG_{L}(\bx^\toptone)} \leq \normtwo{\mathcalG_{L}(\bx^\toptzero)},
\end{aligned}
$$
where $ \mathcalG_{L} \triangleq \mathcalG_{L}^{f, g} $ and $ \mathcalT_{L} \triangleq \mathcalT_{L}^{f, g} $.
\end{lemma}
\begin{proof}[of Lemma~\ref{lemma:mont_prox}]
$\textbf{\textbf{\fejer} \text{monotonicity}}.$
Since $f$ is convex, $\mathcalD_f$ qualifies  as a Bregman distance such that $\mathcalD_f(\bx^*, \bx^\toptzero)\geq 0$ (Remark~\ref{remark:bregnan_dist}).
For any $ t > 0 $, using the descent lemma in \eqref{equation:des_lem_prox_res2}, we obtain
$$
\begin{aligned}
\frac{2}{\beta} \left(  F(\bx^\toptone) - F(\bx^*) \right) &\leq   \normtwo{\bx^* - \bx^\toptzero}^2 -\normtwo{\bx^* - \bx^\toptone}^2 - \frac{2}{\beta} \mathcalD_f(\bx^*, \bx^\toptzero) \\
&\leq  \normtwo{\bx^* - \bx^\toptzero}^2 - \normtwo{\bx^* - \bx^\toptone}^2.
\end{aligned}
$$

\paragraph{Gradient mapping.}
Let $ \ba \triangleq \nabla f(\bx^\toptone) - \nabla f(\bx^\toptzero) $ and $ \bb \triangleq \bx^\toptone - \bx^\toptzero $. By the convexity and smoothness of $f$, Theorem~\ref{theorem:charac_smoo} shows that
$$
\normtwo{\ba}^2 \leq L \innerproduct{\ba, \bb}
\quad\implies \quad
\normtwo{\ba - \frac{L}{2} \bb}^2 \leq \frac{L^2}{4} \normtwo{\bb}^2
\quad\implies \quad
\normtwo{\frac{1}{L} \ba - \frac{1}{2} \bb} \leq \frac{1}{2} \normtwo{\bb}.
$$
Using the triangle inequality,
$$
\normtwo{\frac{1}{L} \ba - \bb } \leq \normtwo{ \frac{1}{L} \ba - \bb + \frac{1}{2} \bb} + \frac{1}{2} \normtwo{\bb} \leq \normtwo{\bb}.
$$
Plugging the expressions for $ \ba $ and $ \bb $ into the above inequality, we obtain that
$$
\normtwo{\bx^\toptzero - \frac{1}{L} \nabla f(\bx^\toptzero) - \bx^\toptone + \frac{1}{L} \nabla f(\bx^\toptone)} \leq \normtwo{\bx^\toptone - \bx^\toptzero}.
$$
Since $g$ is convex, using Proximal Property-IV (Lemma~\ref{theorem:proj_nonexpan}),
$$
\begin{aligned}
&\normtwo{ \mathcalG_{L}(\bx^\toptone)} =  L \normtwo{\bx^\toptone - \mathcalT_{L}(\bx^\toptone)} 
\\
&= L \normtwo{\prox_{\frac{1}{L} g} \left( \bx^\toptzero - \frac{1}{L} \nabla f(\bx^\toptzero) \right) - \prox_{\frac{1}{L} g} \left( \bx^\toptone - \frac{1}{L} \nabla f(\bx^\toptone) \right)} \\
&\leq L \normtwo{\bx^\toptzero - \frac{1}{L} \nabla f(\bx^\toptzero) - \bx^\toptone + \frac{1}{L} \nabla f(\bx^\toptone)} 
\leq L \normtwo{\bx^\toptone - \bx^\toptzero}  = \normtwo{\mathcalG_{L}(\bx^\toptzero)}.
\end{aligned}
$$
This completes the proof.
\end{proof}

We then prove the convergence results for the proximal gradient methods with $F$ being smooth, convex and smooth, and strongly convex and smooth, respectively.
\begin{theoremHigh}[Proximal Gradient for SS $f$]\label{theorem:prox_conv_ss}
Let $ f $ be a proper \textcolor{black}{closed} and $\beta$-smooth function, and let  $ g $ be a proper closed and convex function satisfying assumption (A1) and (A2) in Assumption~\ref{assumption:proximal_grad}.
Let $ \{\bx^\toptzero\}_{t > 0} $ be the sequence generated by the proximal gradient method (Algorithm~\ref{alg:prox_gd_gen}) for solving problem (P3)  with a constant stepsize defined by $ L_t \triangleq L \in \left(\frac{\beta}{2}, \infty\right) $. Then,
\begin{enumerate}[(i)]
\item $ \mathcalG_L(\bx^\toptzero) \rightarrow \bzero $ as $ t \rightarrow \infty $.
\item 
$
\min_{t=1,\ldots,T} \normtwo{\mathcalG_L(\bx^\toptzero)} \leq \frac{\sqrt{F(\bx^{(1)}) - F(\bx^*)}}{\sqrt{CT}},
$
where $ C\triangleq -\frac{\frac{\beta}{2} -L}{L^2}$;
\item All limit points of the sequence $ \{\bx^\toptzero\}_{t > 0} $ are stationary points of problem (P3).
\end{enumerate}
\end{theoremHigh}

\begin{proof}[of Theorem~\ref{theorem:prox_conv_ss}]
\textbf{(i).} Since the sequence $ \{F(\bx^\toptzero)\}_{t > 0} $ is nonincreasing and bounded below, and is equal to zero only if $\bx^\toptzero$ is a stationary point (Lemma~\ref{lemma:mont_prox}), it converges. Thus, in particular, $ F(\bx^\toptzero) - F(\bx^\toptone) \rightarrow 0 $ as $ t \rightarrow \infty $, which, combined with \eqref{equation:mont_prox_1}, implies that $ \normtwo{\mathcalG_L(\bx^\toptzero)} \rightarrow \bzero $ as $ t \rightarrow \infty $.

\paragraph{(ii).} Performing telescopic cancellations on \eqref{equation:mont_prox_1} over $t=\{1,2,\ldots,T\}$, we have
$$
F(\bx^{(1)})-F(\bx^*) 
\geq 
F(\bx^{(1)}) - F(\bx^{(T+1)}) \geq C \sum_{t=1}^{T} \normtwo{\mathcalG_L(\bx^\toptzero)}^2 \geq CT \min_{t\in\{1,2,\ldots,T\}} \normtwo{\mathcalG_L(\bx^\toptzero)}^2.
$$
This concludes the result.

\paragraph{(iii).} Suppose that $ \widetildebx $ is a limit point of $ \{\bx^\toptzero\}_{t > 0} $. Then there exists a subsequence $ \{\bx^{(t_j)}\}_{j > 0} $ converging to $ \widetildebx $. For any $ j > 0 $,
$$
\normtwo{\mathcalG_L(\widetildebx)} \leq \normtwo{\mathcalG_L(\bx^{(t_j)}) - \mathcalG_L(\widetildebx)} + \normtwo{\mathcalG_L(\bx^{(t_j)})} \leq (2L + \beta) \normtwo{\bx^{(t_j)} - \widetildebx} + \normtwo{\mathcalG_L(\bx^{(t_j)})},
$$
where the second inequality follows from the Lipschitzness of $\mathcalG_L(\cdot)$ (Problem~\ref{prob:grad_map_lipschitz}). Since $\mathcalG_L(\bx^{(t_j)})$ tends to $ \bzero $ as $ j \rightarrow \infty $, it follows that $ \mathcalG_L(\widetildebx) = \bzero $, which by Theorem~\ref{theorem:opt_gdmap} implies that $ \widetildebx $ is a stationary point of problem (P3).
\end{proof}

\begin{theoremHigh}[Proximal Gradient for Convex and SS $f$: $\mathcalO(1/T)$]\label{theorem:prox_conv_ss_cvx}
Let $ f $ be a proper \textcolor{black}{closed} and $\beta$-smooth function, and let   $ g $ be a proper closed and convex function satisfying assumption (A1) and (A2) in Assumption~\ref{assumption:proximal_grad}.
Additionally, assume  that $ f $ is \textbf{convex}. 
Let $ \{\bx^\toptzero\}_{t > 0} $ be the sequence generated by the proximal gradient method (Algorithm~\ref{alg:prox_gd_gen}) for solving problem (P3)  with a constant stepsize rule in which $ L_t  \triangleq \beta $ for all $ t > 0 $. 
Then, for any optimizer $ \bx^*$ of $F\triangleq f+g$ and $ T > 1 $, it follows that
$$
F(\bx^{(T)}) - F(\bx^*) \leq \frac{ \beta \normtwo{\bx^{(1)} - \bx^*}^2}{2(T-1)}.
$$
\end{theoremHigh}
\begin{proof}[of Theorem~\ref{theorem:prox_conv_ss_cvx}]
Since $f$ is convex, $\mathcalD_f$ qualifies a Bregman distance such that $\mathcalD_f(\bx^*, \bx^\toptzero)\geq 0$ (Remark~\ref{remark:bregnan_dist}).
For any $ t > 0 $, using the descent lemma in \eqref{equation:des_lem_prox_res2}, we obtain
$$
\begin{aligned}
\frac{2}{\beta} \left(  F(\bx^\toptone) - F(\bx^*) \right) &\leq   \normtwo{\bx^* - \bx^\toptzero}^2 -\normtwo{\bx^* - \bx^\toptone}^2 - \frac{2}{\beta} \mathcalD_f(\bx^*, \bx^\toptzero) \\
&\leq  \normtwo{\bx^* - \bx^\toptzero}^2 - \normtwo{\bx^* - \bx^\toptone}^2.
\end{aligned}
$$
Performing telescopic cancellations over $ t = \{1, 2, \ldots, T-1\} $, we obtain
$$
\sum_{t=1}^{T-1} \left(  F(\bx^\toptone) - F(\bx^*) \right) \leq   \frac{ \beta}{2}\normtwo{\bx^* - \bx^{(1)}}^2 - \frac{ \beta}{2}\normtwo{\bx^* - \bx^{(T)}}^2
\leq  
\frac{ \beta}{2}\normtwo{\bx^* - \bx^{(1)}}^2.
$$
By the monotonicity of $ \{ F(\bx^\toptzero) \}_{t > 0} $ (Equation~\eqref{equation:des_lem_prox_res1} and Lemma~\ref{lemma:mont_prox}), we can conclude that
$$
\begin{aligned}
(T-1) \left( F(\bx^{(T)}) - F(\bx^*) \right) &\leq \sum_{t=1}^{T-1} \left( F(\bx^\toptone) - F(\bx^*) \right) \leq \frac{ \beta}{2} \normtwo{\bx^* - \bx^{(1)}}^2,
\end{aligned}
$$
which implies the desired result.
\end{proof}

\begin{theoremHigh}[Proximal Gradient for SC and SS $f$]\label{theorem:prox_conv_ss_sc}
Let $ f $ be a proper \textcolor{black}{closed} and $\beta$-smooth function, and let    $ g $ be a proper closed and convex function satisfying assumption (A1) and (A2) in Assumption~\ref{assumption:proximal_grad}. 
Additionally, assume that $ f $ is \textbf{$ \alpha $-strongly convex}.  
Let $ \{\bx^\toptzero\}_{t > 0} $ be the sequence generated by the proximal gradient method (Algorithm~\ref{alg:prox_gd_gen}) for solving problem (P3)  with a constant stepsize rule in which $ L_t  \triangleq \beta $ for all $ t > 0 $.
Then, for any optimizer $ \bx^*$ of $F\triangleq f+g$ and $ t >1 $,  it follows that
$$
\begin{aligned}
\normtwo{\bx^\toptone - \bx^*}^2 &\leq \exp\left(-t\frac{\alpha}{\beta}\right)  \normtwo{\bx^{(1)} - \bx^*}^2;\\
F(\bx^\toptone) - F(\bx^*) &\leq \frac{ \beta}{2} \exp\left(-t\frac{\alpha}{\beta}\right)  \normtwo{\bx^{(1)} - \bx^*}^2.
\end{aligned}
$$
This also shows  that $\normtwo{\bx^{(T)} -\bx^*}^2 \leq \epsilon$ after at most $T = \mathcalO\left(\frac{\beta}{\alpha} \ln \frac{1}{\epsilon}\right)$ steps.
\end{theoremHigh}

\begin{proof}[of Theorem~\ref{theorem:prox_conv_ss_sc}]
Since $f$ is $\alpha$-strongly convex, $\mathcalD_f$ qualifies a Bregman distance such that $\mathcalD_f(\bx^*, \bx^\toptzero)\geq 0$ (Remark~\ref{remark:bregnan_dist}), and  it follows by the definition of strong convexity (Definition~\ref{definition:scss_func}) that
$$
\mathcalD_f(\bx^*, \bx^\toptzero) = f(\bx^*) - f(\bx^\toptzero) - \innerproduct{\nabla f(\bx^\toptzero), \bx^* - \bx^\toptzero} \geq \frac{\alpha}{2} \normtwo{\bx^\toptzero - \bx^*}^2.
$$
For any $ t > 0 $, using the descent lemma in \eqref{equation:des_lem_prox_res2}, we obtain
\begin{equation}\label{equation:prox_conv_ss_sc}
\begin{aligned}
	F(\bx^\toptone) - F(\bx^*) &\leq  \frac{\beta}{2} \normtwo{\bx^* - \bx^\toptzero}^2 -\frac{\beta}{2}\normtwo{\bx^* - \bx^\toptone}^2 - \mathcalD_f(\bx^*, \bx^\toptzero)\\
&\leq   \frac{\beta-\alpha}{2} \normtwo{\bx^* - \bx^\toptzero}^2 -\frac{\beta}{2}\normtwo{\bx^* - \bx^\toptone}^2.
\end{aligned}
\end{equation}
Since $ \bx^* $ is a minimizer of $ F $, $  F(\bx^\toptone) -F(\bx^*) \geq 0 $, and hence, the above inequality implies that
$$
\normtwo{\bx^\toptone - \bx^*}^2  \leq \left(1 - \frac{\alpha}{ \beta}\right) \normtwo{\bx^\toptzero - \bx^*}^2 
\leq \left(1 - \frac{\alpha}{ \beta}\right)^t \normtwo{\bx^{(1)} - \bx^*}^2.
$$
On the other hand, \eqref{equation:prox_conv_ss_sc} also shows 
$$
\begin{aligned}
	F(\bx^\toptone) - F(\bx^*) &\leq \frac{\beta - \alpha}{2} \normtwo{\bx^\toptzero - \bx^*}^2 - \frac{\beta}{2} \normtwo{\bx^\toptone - \bx^*}^2 
	\leq \frac{ \beta - \alpha}{2} \normtwo{\bx^\toptzero - \bx^*}^2 \\
	&= \frac{ \beta}{2} \left(1 - \frac{\alpha}{ \beta}\right) \normtwo{\bx^\toptzero - \bx^*}^2 
	\leq \frac{ \beta}{2} \left(1 - \frac{\alpha}{ \beta}\right)^{t} \normtwo{\bx^{(1)} - \bx^*}^2,
\end{aligned}
$$
which completes the proof.
\end{proof}

Once again, we observe that the condition number $\kappa\triangleq\frac{\beta}{\alpha}$ plays a crucial role in determining the convergence rate of the algorithm. The rate of convergence for the proximal gradient method is identical to those described in Theorems \ref{theorem:gd_sc_ss} and \ref{theorem:pgd_sc_ss}. However, in this case, the algorithm is applied to solve the composite problem (P3).

\begin{algorithm}[h] 
\caption{Proximal Point Method}
\label{alg:prox_point_gen}
\begin{algorithmic}[1] 
\Require A proper closed convex function $g$ ; 
\State {\bfseries Input:}  Initialize $\bx^{(1)}$;
\For{$t=1,2,\ldots$}
\State Pick a stepsize $\eta_t$;
\State $\bx^{(t+1)} \leftarrow \prox_{\eta_t g}(\bx^{(t)})$;
\EndFor
\State Output  $\bx_{\text{final}}\leftarrow \bx^{(T)}$;
\end{algorithmic} 
\end{algorithm}

\begin{figure}[h]
\centering  
\vspace{-0.15cm} 
\subfigbottomskip=2pt 
\subfigcapskip=-5pt 
\includegraphics[width=0.78\textwidth]{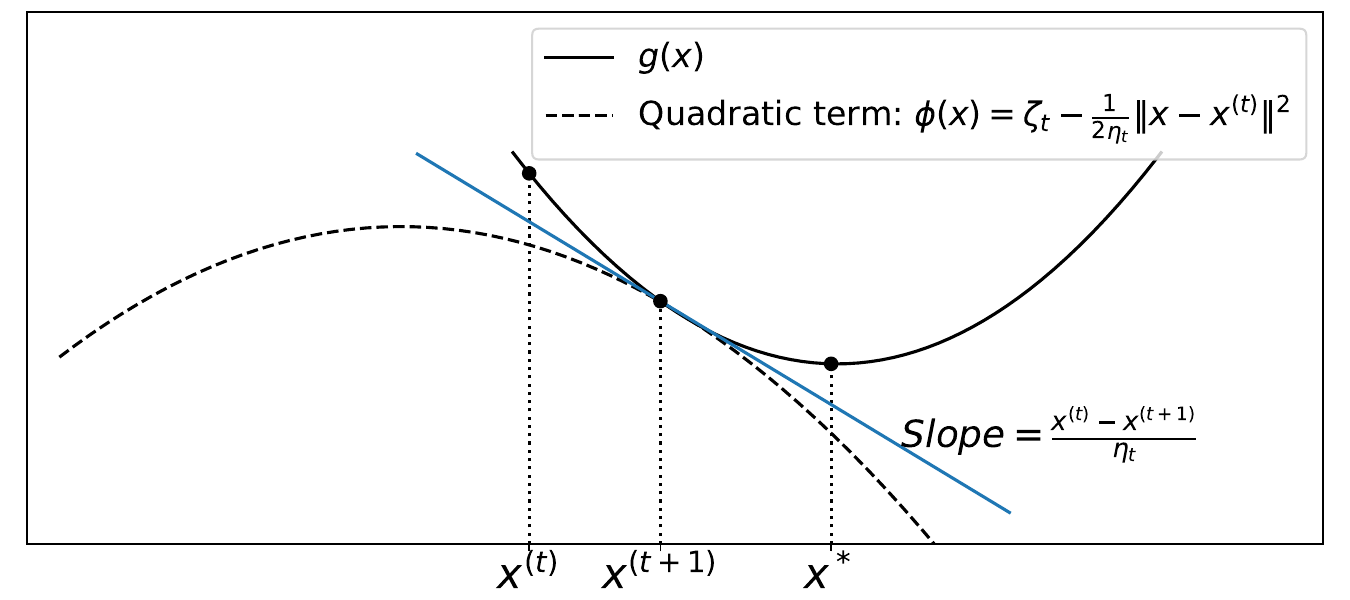}
\caption{Geometric view of the proximal point method. The minimum of the function 
$
g(\bx) + \frac{1}{2\eta_t} \normtwo{\bx - \bx^\toptzero}^2
$
is attained at a unique point $ \bx^\toptone $. The value $ \zeta_t $ in the quadratic term represents the scalar adjustment needed to raise the graph of
$
-\frac{1}{2\eta_t} \normtwo{\bx - \bx^\toptzero}^2
$
so that it just touches the graph of $ g $. The slope shown in the figure,
$
\frac{\bx^\toptzero - \bx^\toptone}{\eta_t},
$
is the common gradient of $ g(\bx) $ and 
$
-\frac{1}{2\eta_t} \normtwo{\bx - \bx^\toptzero}^2
$
at the minimizing point $ \bx^\toptone $.
}
\label{fig:prox_point}
\end{figure}
\section{Proximal Point Method}\label{section:proximal_point}

Consider the problem
\begin{equation}\label{equation:prob_prox_point}
\min_{\bx \in \real^n} g(\bx),
\end{equation}
where $ g : \real^n \to (-\infty, \infty] $ is a proper closed and convex function. Problem~\eqref{equation:prob_prox_point} is actually a special case of the composite problem (P3) with $ f \triangleq 0 $. The update step of the proximal gradient method in this scenario  takes the form
$$
\bx^\toptone = \prox_{\eta_t g}(\bx^\toptzero).
$$
Taking $ \eta_t = \eta$ for some $ \eta > 0 $, we obtain the \textit{proximal point method} (Algorithm~\ref{alg:prox_point_gen}).
See Figure~\ref{fig:prox_point} for an illustration of the optimization procedure at each step.

The proximal point method is generally not practical because its general step requires minimizing the function  $ g(\bx) + \frac{1}{2\eta} \normtwo{ \bx - \bx^\toptzero}^2 $ at the $t$-th iteration, which is typically as challenging as solving the original problem of minimizing $ g $. 
However, since the proximal point method is a special case of the proximal gradient method, we can derive its main convergence results from the corresponding results on the proximal gradient method. Specifically, given that the smooth part $ f = 0 $ is 0-smooth, we can use any constant stepsize to ensure convergence. Consequently, Theorem~\ref{theorem:prox_conv_ss_cvx} implies the following result.

\begin{theoremHigh}[Proximal Point Method for Convex]\label{theorem:prox_point}
Let $ g : \real^n \to (-\infty, \infty] $ be a proper closed and convex function. Assume that problem
$$
\min_{\bx \in \real^n} g(\bx)
$$
has a nonempty optimal set. Let $\bx^*$ be an optimal point of $g$, and let $ \{ \bx^\toptzero \}_{t>0} $ be the sequence generated by the proximal point method (Algorithm~\ref{alg:prox_point_gen}) with parameter $ \eta > 0 $. Then,
\begin{enumerate}[(a)]
\item $ g(\bx^\toptzero) - g(\bx^*) \leq \frac{\normtwo{\bx^\topone - \bx^*}^2}{2\eta t} $  for $ t>0 $;
\item The sequence $ \{ \bx^\toptzero \}_{t>0} $ converges to some optimal point $\bx^*$ of $g$.
\end{enumerate}
\end{theoremHigh}


\noindent
\begin{minipage}[t]{0.505\linewidth}
\begin{algorithm}[H]
\caption{FISTA-V1}
\label{alg:fistav1}
\begin{algorithmic}[1]
\State $\by^{(1)} = \bx^{(1)} \in \real^n$, $\gamma_1=\frac{1}{\zeta_1}=1$;
\For{$t=1,2,\ldots$}
\State Choose $\eta_{t}$ and $\zeta_t=\frac{1}{\gamma_t}$;
\State $\bx^{(t+1)} \leftarrow \prox_{\eta_{t} g}(\by^{(t)} - \eta_{t} \nabla f(\by^{(t)}))$;
\State $\small\begin{aligned}
&\by^{(t+1)} \leftarrow \bx^{(t+1)} + \frac{\zeta_t-1}{\zeta_{t+1}} (\bx^{(t+1)} - \bx^{(t)});\\
&= \bx^{(t+1)} + (\frac{\gamma_{t+1}}{\gamma_t}-\gamma_{t+1}) (\bx^{(t+1)} - \bx^{(t)});
\end{aligned}$
\EndFor
\State Output  $\bx_{\text{final}}\leftarrow \bx^{(T)}$;
\end{algorithmic}
\end{algorithm}
\end{minipage}%
\hfil 
\begin{minipage}[t]{0.495\linewidth}
\begin{algorithm}[H]
\caption{FISTA-V2}
\label{alg:fistav2}
\begin{algorithmic}[1]
\State $\bu^{(1)} =\by^{(1)} = \bx^{(1)} \in \real^n$, $\gamma_1=1$;
\For{$t=1,2,\ldots$}
\State Choose $\eta_{t}$ and $\gamma_t$;
\State $\bx^{(t+1)} \leftarrow \prox_{\eta_{t} g}(\by^{(t)} - \eta_{t} \nabla f(\by^{(t)}))$;
\State $\bu^{(t+1)} \leftarrow \bx^{(t)} + \frac{1}{\gamma_{t}}(\bx^{(t+1)} - \bx^{(t)})$;
\State $\by^{(t+1)} \leftarrow (1 - \gamma_{t+1})\bx^{(t+1)} + \gamma_{t+1} \bu^{(t+1)}$;
\EndFor
\State Output  $\bx_{\text{final}}\leftarrow \bx^{(T)}$;
\end{algorithmic}
\end{algorithm}
\end{minipage}
\section{Fast Proximal Gradient Method---FISTA}\label{section:fista}
The \textit{fast proximal gradient} method, also known as the \textit{fast iterative shrinkage-thresholding algorithm (FISTA)} (Algorithm~\ref{alg:fistav2}) is a variant of the proximal gradient method \citep{beck2009fast}. 
Recall that the proximal gradient method achieves a convergence rate of  $\mathcalO\left(\frac{1}{T}\right)$ for a convex and smooth function (Theorem~\ref{theorem:prox_conv_ss_cvx}). Naturally, we aim to accelerate this method, leading us to the introduction of FISTA in this section.

FISTA comprises two steps: first, it calculates a new point based on the direction from the previous two iterations; second, it performs a proximal gradient iteration at this new point. Specifically,
\begin{subequations}\label{equation:fista_intui}
\begin{align}
\by^{(t)}   &\leftarrow \bx^{(t)} + \frac{t-2}{t+1}(\bx^{(t)} - \bx^{(t-1)});\\
\bx^{(t+1)} &\leftarrow \prox_{\eta_t g}(\by^{(t)} - \eta_t \nabla f(\by^{(t)})).
\end{align}
\end{subequations}
The key difference from the standard proximal gradient method is the inclusion of a  \textit{momentum term}, $\frac{t-2}{t+1}(\bx^{(t)} - \bx^{(t-1)})$. 
This momentum term accelerates convergence by leveraging the difference between the current and previous iterates. 
Importantly, this approach has minimal impact on the computational cost per iteration but significantly enhances performance. 
If $\eta_t$ is chosen as a fixed stepsize  less than or equal to $\frac{1}{\beta}$ ($f$ is assumed to be $\beta$-smooth), the convergence rate improves to $\mathcalO\left(\frac{1}{T^2}\right)$. We will detail the derivation in the convergence analysis section.
With this intuition, the complete FISTA algorithm is shown in Algorithm~\ref{alg:fistav1}, where we let $\zeta_t=\frac{1}{\gamma_t}$ for all $t>0$. When $\frac{1}{\gamma_t}=\zeta_t \triangleq \frac{t+1}{2}$, Algorithm~\ref{alg:fistav1} is equivalent to the update in \eqref{equation:fista_intui}.

To facilitate convergence analysis, an equivalent formulation of the FISTA algorithm can be presented. This reformulation splits the second step of the original algorithm into two separate operations, as shown in Algorithm~\ref{alg:fistav2}. While Algorithm~\ref{alg:fistav1} highlights the momentum term explicitly, making it more explanatory, Algorithm~\ref{alg:fistav2} simplifies the convergence analysis (see Theorem~\ref{theorem:fista_conv_ssf}).

\paragrapharrow{Convergence requirement.}
For this algorithm framework, we need to determine how to choose the stepsizes $\eta_t$ and $\gamma_t = \frac{1}{\zeta_t}$, which influence the convergence rate of the algorithm. 
First, we outline the conditions under which Algorithm~\ref{alg:fistav2} achieves a convergence rate of $\mathcalO\left(\frac{1}{T^2}\right)$ (the detailed proof will follow):
\begin{subequations}\label{equation:fista_require}
\small
\begin{align}
&f(\bx^\toptzero) \leq f(\by^\toptzero) + \innerproduct{\nabla f(\by^\toptzero), \bx^\toptzero - \by^\toptzero} + \frac{1}{2\eta_t} \normtwo{\bx^\toptzero - \by^\toptzero}^2, \text{ i.e., $\eta_t \leq \frac{1}{\beta}$};  \label{equation:fista_require1}\\
&\gamma_1 = 1, \quad \frac{(1 - \gamma_{t+1})\eta_{t+1}}{\gamma_{t+1}^2} \leq \frac{\eta_{t}}{\gamma_{t}^2}, \quad t > 1; \label{equation:fista_require2}\\
&\frac{\gamma_t^2}{\eta_t} = \mathcalO\left(\frac{1}{t^2}\right). \label{equation:fista_require3}
\end{align}
\end{subequations}
\paragrapharrow{Choice of $\gamma_t = \frac{1}{\zeta_t}$.}
It can be observed  that when $\eta_t = \frac{1}{\beta}$ and $\frac{1}{\zeta_t}=\gamma_t = \frac{2}{t+1}$ for $t>0$, the above conditions are satisfied. Moreover, the choice of $\gamma_t$ is not unique; for example, we can take
$$
\small
\begin{aligned}
\gamma_1 = 1, \quad \frac{1}{\gamma_{t+1}} = \frac{1}{2} \left( 1 + \sqrt{1 + \frac{4}{\gamma_{t}^2}} \right)
\iff 
\zeta_{t+1} = \frac{1 + \sqrt{1 + 4\zeta_t^2}}{2},
\quad t>0.
\end{aligned}
$$
to obtain the sequence $\{\gamma_t\} = \{\frac{1}{\zeta_t}\}$, and the derived algorithm still converges at a rate of $\mathcalO\left(\frac{1}{T^2}\right)$.
Using the following lemma, it can also be shown that the above choice of $\{\gamma_t\} = \{\frac{1}{\zeta_t}\}$ also satisfies the requirements in specified  \eqref{equation:fista_require}.

\begin{lemma}\label{lemma:seq_gamma_zeta}
Let $\{\frac{1}{\gamma_t}\}_{t>0}=\{\zeta_t\}_{t> 0}$ be the sequence (see Algorithm~\ref{alg:fistav1}) defined by
$ \zeta_1 = 1, \; \zeta_{t+1} = \frac{1 + \sqrt{1 + 4\zeta_t^2}}{2},  t \geq 0. $
Then, $\zeta_t \geq \frac{t+1}{2}$ for all $t > 0$.
\end{lemma}
\begin{proof}[of Lemma~\ref{lemma:seq_gamma_zeta}]
The proof is by induction. Obviously, for $t = 1$, $\zeta_1 = 1 \geq \frac{1+1}{2}$. Suppose that the claim holds for $t$ such that $\zeta_t \geq \frac{t+1}{2}$. We will prove that $\zeta_{t+1} \geq \frac{t+2}{2}$. By the recursive relation defining the sequence and the induction assumption,
$ \zeta_{t+1} = \frac{1 + \sqrt{1 + 4\zeta_t^2}}{2} \geq \frac{1 + \sqrt{1 + (t+1)^2}}{2} \geq \frac{1 + \sqrt{(t+1)^2}}{2} = \frac{t+2}{2}$, which completes the proof.
\end{proof}

\paragrapharrow{Monotone issue.}
The original FISTA algorithm is not a descent algorithm. Here, we present a descent variant of FISTA that only requires modifying a single line of Algorithm~\ref{alg:fistav1}: $\bx^\toptone \in \arg\min\{F(\bx) : \bx = \bx^\toptzero, \bz^\toptzero\}$ to ensure the condition $F(\bx^\toptone) \leq \min\{F(\bz^\toptzero), F(\bx^\toptzero)\}$ (Algorithm~\ref{alg:fistav_monotone}). After calculating the proximal operator, instead of immediately updating the iterate point, we check if the function value has decreased. Only if there is a reduction in the function value do we proceed to update the iterate point.
More compactly,  Steps 6 to 10 in Algorithm~\ref{alg:fistav_monotone} can be equivalently stated as
$$
\by^\toptone = \bx^\toptone + \frac{\zeta_t}{\zeta_{t+1}} (\bz^\toptzero - \bx^\toptone) + \frac{\zeta_t - 1}{\zeta_{t+1}} (\bx^\toptone - \bx^\toptzero).
$$

\noindent
\begin{minipage}[h]{1\linewidth}
\begin{algorithm}[H]
\caption{Monotine FISTA \citep{beck2009fastgradient}, Compare to Algorithm~\ref{alg:fistav1}}
\label{alg:fistav_monotone}
\begin{algorithmic}[1]
\Require A function $f(\bx)$ and a closed convex function $g$ (usually non-smooth) satisfying (A1) and (A2) in Assumption~\ref{assumption:proximal_grad}; 
\State {\bfseries Input:} Initialize $\by^{(1)} = \bx^{(1)} \in \real^n$, $\gamma_1=\frac{1}{\zeta_1}=1$;
\For{$t=1,2,\ldots$}
\State Choose $\eta_{t}$ and $\zeta_t=\frac{1}{\gamma_t}$;
\State $\bz^\toptzero \leftarrow \prox_{\eta_{t} g}(\by^{(t)} - \eta_{t} \nabla f(\by^{(t)}))$;
\State choose $\bx^\toptone$ such that $F(\bx^\toptone) \leq \min\{F(\bz^\toptzero), F(\bx^\toptzero)\}$;
\If{$\bx^\toptone$ is $\bz^\toptzero$}
\State $\small\begin{aligned}
\by^{(t+1)} \leftarrow \bx^{(t+1)} + \frac{\zeta_t-1}{\zeta_{t+1}} (\bz^{(t)} - \bx^{(t)})
\end{aligned}$ \Comment{Same as FISTA}
\ElsIf{$\bx^\toptone$ is $\bx^\toptzero$}
\State $\small\begin{aligned}
\by^{(t+1)} \leftarrow \bx^\toptzero + \frac{\zeta_t}{\zeta_{t+1}} (\bz^\toptzero - \bx^\toptzero);
\end{aligned}$
\EndIf
\EndFor
\State {\bfseries Return:}   $\bx_{\text{final}}\leftarrow \bx^{(T)}$;
\end{algorithmic}
\end{algorithm}
\end{minipage}%
\hfil 
\begin{minipage}[t]{1\linewidth}
\begin{algorithm}[H]
\caption{Line Search Algorithm at $t$-th Iteration}
\label{alg:fistav_line}
\begin{algorithmic}[1]
\State {\bfseries Input:}  $\eta_t = \eta_{t-1} > 0$, $\rho < 1$. Reference point $\by^\toptzero$ and its gradient $\nabla f(\by^\toptzero)$;
\State Calculate candidate update $\bx^\toptone \leftarrow \prox_{\eta_t g}(\by^\toptzero - \eta_t \nabla f(\by^\toptzero))$;
\While {condition \eqref{equation:fista_require1} is not satisfied for $\bx^\toptone, \by^\toptzero$}
\State Reducing stepsize $\eta_t \gets \rho \eta_t$;
\State Recalculate $\bx^\toptone \leftarrow \prox_{\eta_t g}(\by^\toptzero - \eta_t \nabla f(\by^\toptzero))$;
\EndWhile
\State {\bfseries Return:}  update $\bx^\toptone$, stepsize $\eta_t$;
\end{algorithmic}
\end{algorithm}
\end{minipage}

\paragrapharrow{Line search methods.}
In Algorithms~\ref{alg:fistav1} and \ref{alg:fistav2}, the stepsize must satisfy $\eta_t \leq \frac{1}{\beta}$, under which the condition \eqref{equation:fista_require1} is satisfied. However, for most problems, we do not know the Lipschitz constant of the function $\nabla f$. In this case, condition \eqref{equation:fista_require1} can still be satisfied by employing a line search to determine an appropriate $\eta_t$, while choosing $\gamma_t$ such that conditions \eqref{equation:fista_require2} and \eqref{equation:fista_require3} are simultaneously satisfied, making the algorithm achieve a convergence rate of $\mathcalO\left(\frac{1}{T^2}\right)$.
The line search Algorithm~\ref{alg:fistav_line} is one such method used to determine the stepsize $\eta_t$ dynamically. 
For each iteration $t$, the initial stepsize $\eta_t$ is taken as the previous stepsize $\eta_{t-1}$. 
By progressively  reducing the stepsize $\eta_t$, condition \eqref{equation:fista_require1} is satisfied. Note that when $\eta_t$ is sufficiently small, this condition will always be met, thus preventing the line search from failing to terminate. 
It is straightforward  to verify that the other two conditions \eqref{equation:fista_require2} and \eqref{equation:fista_require3} are also satisfied during the iterative  process.

The following theorem provides the convergence rate of the FISTA algorithm with a fixed stepsize. The convergence with a dynamic stepsize using line search algorithms can be proved analogously and we will not repeat the details.
\begin{theoremHigh}[FISTA for  Convex and SS $f$, $\mathcalO(1/T^2)$:  Compare to Theorem~\ref{theorem:prox_conv_ss_cvx}]\label{theorem:fista_conv_ssf}
Let $ f $ be a proper \textcolor{black}{closed} and $\beta$-smooth function  and  $ g $ be a proper closed and convex function satisfying assumption (A1) and (A2) in Assumption~\ref{assumption:proximal_grad}. 
Additionally, assume  that $f$ is \textbf{convex}. 
Let $ \{\bx^\toptzero\}_{t > 0} $ be the sequence generated by the FISTA method (Algorithm~\ref{alg:fistav2}) for solving problem (P3)  with a constant stepsize rule in which $\eta_t \triangleq \frac{1}{\beta}$ for all $ t > 0 $. Then, for any optimizer $ \bx^*$ of $F\triangleq f+g$ and $ T > 1 $, it follows that
$$
F(\bx^{(T)}) - F(\bx^*) \leq \frac{2\beta}{T^2} \normtwo{\bx^{(1)} - \bx^*}^2.~\footnote{If the Lipschitz constant $\beta$ is not known beforehand, line search Algorithm~\ref{alg:fistav_line} can be applied to achieve a same rate of convergence.}
$$
\end{theoremHigh}
\begin{proof}[of Theorem~\ref{theorem:fista_conv_ssf}]
Since $\bx^\toptone = \prox_{\eta_t g}(\by^\toptzero - \eta_t \nabla f(\by^\toptzero))$, by the Proximal Property-I (Lemma~\ref{lemma:prox_prop1}), we have
$$
-\bx^\toptone + \by^\toptzero - \eta_t \nabla f(\by^\toptzero) \in \eta_t \partial g(\bx^\toptone).
$$
By the subgradient inequality (Definition~\ref{definition:subgrad}), for any $\bx\in\real^n$, we have
$$
\eta_t g(\bx) \geq \eta_t g(\bx^\toptone) + \innerproduct{-\bx^\toptone + \by^\toptzero - \eta_t \nabla f(\by^\toptzero), \bx - \bx^\toptone}.
$$
By the $\beta$-smoothness of $f$ (Definition~\ref{definition:scss_func}) and $\eta_t = \frac{1}{\beta}$, we obtain
$$
f(\bx^\toptone) \leq f(\by^\toptzero) + \innerproduct{\nabla f(\by^\toptzero), \bx^\toptone - \by^\toptzero} + \frac{1}{2\eta_t} \normtwobig{\bx^\toptone - \by^\toptzero}^2.
$$
Combining the above two inequalities, for any $\bx$, we have
\begin{equation}\label{equation:fista_conv_ssf1}
\small
\begin{aligned}
&F(\bx^\toptone) = f(\bx^\toptone) + g(\bx^\toptone)\\
&\leq g(\bx) + f(\by^\toptzero) + \innerproduct{\nabla f(\by^\toptzero), \bx - \by^\toptzero} + \frac{1}{\eta_t} \innerproduct{\bx^\toptone - \by^\toptzero, \bx - \bx^\toptone} + \frac{1}{2\eta_t} \normtwobig{\bx^\toptone - \by^\toptzero}^2\\
&\stackrel{\dag}{\leq} g(\bx) + f(\bx) + \frac{1}{\eta_t} \innerproduct{\bx^\toptone - \by^\toptzero, \bx - \bx^\toptone} + \frac{1}{2\eta_t} \normtwobig{\bx^\toptone - \by^\toptzero}^2\\
&= F(\bx) + \frac{1}{\eta_t} \innerproduct{\bx^\toptone - \by^\toptzero, \bx - \bx^\toptone} + \frac{1}{2\eta_t} \normtwobig{\bx^\toptone - \by^\toptzero}^2,
\end{aligned}
\end{equation}
where the inequality ($\dag$) follows from the convexity of $f$.
Invoking \eqref{equation:fista_conv_ssf1} with $\bx \triangleq \bx^{(t)}$ and $\bx \triangleq \bx^*$, for which we multiply by $1 - \gamma_t$ and $\gamma_t$, respectively, and add them up to get
\begin{equation}\label{equation:fista_conv_ssf2}
\small
\begin{aligned}
&F(\bx^\toptone) - F(\bx^*) - (1 - \gamma_t)\big(F(\bx^{(t)}) - F(\bx^*)\big)\\
&\leq \frac{1}{\eta_t} \innerproduct{\bx^\toptone - \by^\toptzero, (1 - \gamma_t)\bx^{(t)} + \gamma_t \bx^* - \bx^\toptone} + \frac{1}{2\eta_t} \normtwobig{\bx^\toptone - \by^\toptzero}^2\\
&\stackrel{\dag}{=} \frac{1}{2\eta_t} \left\{\normtwo{\by^\toptzero - (1 - \gamma_t)\bx^{(t)} - \gamma_t \bx^*}^2 - \normtwo{\bx^\toptone - (1 - \gamma_t)\bx^{(t)} - \gamma_t \bx^*}^2\right\}\\
&= \frac{\gamma_t^2}{2\eta_t} \left\{\normtwo{\bu^{(t)} - \bx^*}^2 - \normtwo{\bu^{(t+1)} - \bx^*}^2\right\},
\end{aligned}
\end{equation}
where the equality ($\dag$) follows from the fact that $2\innerproduct{\ba-\bb,\bc-\ba} +\normtwo{\ba-\bb}^2 = \normtwo{\bb-\bc}^2-\normtwo{\ba-\bc}^2$ for any $\ba,\bb,\bc\in\real^n$, and the last equality follows from the update rules:
$$
\begin{aligned}
\bu^{(t+1)} &= \bx^{(t)} + \frac{1}{\gamma_{t}}(\bx^{(t+1)} - \bx^{(t)})
\qquad\text{and}\qquad
\by^{(t)} = (1 - \gamma_t)\bx^{(t)} + \gamma_t \bu^{(t)}.
\end{aligned}
$$
Recalling that $\eta_t$ and $\gamma_t$ are chosen such that
$
\frac{1 - \gamma_t}{\gamma_t^2} \eta_t \leq \frac{1}{\gamma_{t-1}^2} \eta_{t-1},
$
we obtain an iterate inequality for consecutive iterations
$$
\small
\begin{aligned}
&\frac{\eta_t}{\gamma_t^2} \big(F(\bx^\toptone) - F(\bx^*)\big) + \frac{1}{2} \normtwo{\bu^{(t+1)} - \bx^*}^2 
\leq (1-\gamma_t)\frac{\eta_t}{\gamma_t^2} \big(F(\bx^{(t)}) - F(\bx^*)\big) + \frac{1}{2} \normtwo{\bu^{(t)} - \bx^*}^2\\
&\leq \frac{\eta_{t-1}}{\gamma_{t-1}^2} \big(F(\bx^{(t)}) - F(\bx^*)\big) + \frac{1}{2} \normtwo{\bu^{(t)} - \bx^*}^2
\leq \frac{\eta_1}{\gamma_1^2} \big(F(\bx^{(2)}) - F(\bx^*)\big) + \frac{1}{2} \normtwo{\bu^{(2)} - \bx^*}^2\\
&\leq \frac{(1 - \gamma_1)\eta_1}{\gamma_1^2} \big(F(\bx^{(1)}) - F(\bx^*)\big) + \frac{1}{2} \normtwo{\bu^{(1)} - \bx^*}^2
=\frac{1}{2} \normtwo{\bx^{(1)} - \bx^*}^2,
\end{aligned}
$$
where the last equality follows from the fact that  $\gamma_1 = 1$ and $\bu^{(1)} = \bx^{(1)}$.
Using Lemma~\ref{lemma:seq_gamma_zeta} concludes the result.
\end{proof}

\noindent
\begin{minipage}[t]{0.495\linewidth}
\begin{algorithm}[H]
\caption{FISTA-V2, Same as Algo~\ref{alg:fistav2}}
\label{alg:fistav2_comp}
\begin{algorithmic}[1]
	\State $\bu^{(1)} =\by^{(1)} = \bx^{(1)} \in \real^n$, $\gamma_1=1$;
	\For{$t=1,2,\ldots$}
	\State Choose $\eta_{t}$ and $\gamma_t$;
	\State $\bx^{(t+1)} \leftarrow \prox_{\eta_{t} g}(\by^{(t)} - \eta_{t} \nabla f(\by^{(t)}))$;
	\State $\bu^{(t+1)} \leftarrow \bx^{(t)} + \frac{1}{\gamma_{t}}(\bx^{(t+1)} - \bx^{(t)})$;
	\State $\by^{(t+1)} \leftarrow (1 - \gamma_{t+1})\bx^{(t+1)} + \gamma_{t+1} \bu^{(t+1)}$;
	\EndFor
	\State Output  $\bx_{\text{final}}\leftarrow \bx^{(T)}$;
\end{algorithmic}
\end{algorithm}
\end{minipage}%
\hfil 
\begin{minipage}[t]{0.505\linewidth}
\begin{algorithm}[H]
\caption{Nesterov Accelerated Method}
\label{alg:nesterov}
\begin{algorithmic}[1]
\State $\bx^{(1)} =\by^{(1)} = \widetildebx^{(1)} \in \real^n$, $\gamma_1=1$;
\For{$t=1,2,\ldots$}
\State Choose $\eta_{t}$ and $\gamma_{t+1}$;
\State  $\widetildebx^\toptone \leftarrow \prox_{\frac{\eta_t}{\gamma_t}g}\left(\widetildebx^{(t)} - \frac{\eta_t}{\gamma_t} \nabla f(\by^{(t)})\right)$.
\State  $\bx^\toptone \leftarrow (1 - \gamma_t)\bx^{(t)} + \gamma_t \widetildebx^\toptone$.
\State  $\by^\toptone \leftarrow (1 - \gamma_{t+1})\bx^{(t+1)} + \gamma_{t+1} \widetildebx^{(t+1)}$.
\EndFor
\State Output  $\bx_{\text{final}}\leftarrow \bx^{(T)}$;
\end{algorithmic}
\end{algorithm}
\end{minipage}

\section{Nesterov Accelerated Method}

This section will introduce  another acceleration algorithm besides the FISTA algorithm, which is an extension of Nesterov's algorithms \citep{nesterov1988approach, nesterov2005smooth}. 
The Nesterov accelerated method (Algorithm~\ref{alg:nesterov}) differs from the FISTA algorithm in that the sequences $\{\widetildebx^\toptzero\}$, $\{\by^\toptzero\}$, and $\{\bx^\toptzero\}$ are guaranteed to remain within the domain of definition (since all the updates are convex combinations of old iterates, which bears a resemblance to the conditional gradient method in Section~\ref{section:cond_gd}).
In contrast, the sequence $\{\by^\toptzero\}$ in the FISTA algorithm (Algorithm~\ref{alg:fistav2_comp}) does not necessarily stay within the domain.
Note that  the FISTA method modifies both the base iterate and the descent direction in the proximal operation; while the Nesterov accelerated method only updates the descent direction: 
$$
\bx^{(t+1)} \leftarrow \prox_{\eta_{t} g}(\by^{(t)} - \eta_{t} \nabla f(\by^{(t)}))
\qquad\leadsto \qquad
\widetildebx^\toptone \leftarrow \prox_{\frac{\eta_t}{\gamma_t}g}\left(\widetildebx^{(t)} - \frac{\eta_t}{\gamma_t} \nabla f(\by^{(t)})\right).
$$
Same as the FISTA method,  the Nesterov accelerated method uses $\gamma_t = \frac{2}{t+1}$ (or the recursive sequence $\{\gamma_t\}$ from Lemma~\ref{lemma:seq_gamma_zeta}) and 
$\eta_t = \frac{1}{\beta}$ to achieve a convergence rate of $\mathcalO\left(\frac{1}{T^2}\right)$.
The convergence analysis of the  Nesterov accelerated method can employ techniques analogous to those used for the FISTA method.

\begin{theoremHigh}[Nesterov for  Convex and SS $f$, Compare to Theorem~\ref{theorem:fista_conv_ssf}]\label{theorem:nesterov_acc}
Let $ f $ be a proper \textcolor{black}{closed} and $\beta$-smooth function  and  $ g $ be a proper closed and convex function satisfying assumption (A1) and (A2) in Assumption~\ref{assumption:proximal_grad}. 
Additionally, assume that $f$ is \textbf{convex}. 
Let $ \{\bx^\toptzero\}_{t > 0} $ be the sequence generated by the Nesterov accelerated method (Algorithm~\ref{alg:nesterov}) for solving problem (P3)  with a constant stepsize rule in which $\eta_t \triangleq \frac{1}{\beta}$ and $\gamma_t = \frac{2}{t+1}$ (or the recursive sequence $\{\gamma_t\}$ from Lemma~\ref{lemma:seq_gamma_zeta}) for all $ t > 0 $. Then, for any optimizer $ \bx^*$ of $F\triangleq f+g$ and $ T > 1 $, it follows that
$$
F(\bx^\toptzero) - F(\bx^*) \leq \frac{2\beta}{T^2} \normtwo{\bx^{(1)} - \bx^*}^2.
$$
\end{theoremHigh}
\begin{proof}[of Theorem~\ref{theorem:nesterov_acc}]
Since $\widetildebx^\toptone = \prox_{\frac{\eta_t}{\gamma_t}g}\left(\widetildebx^{(t)} - \frac{\eta_t}{\gamma_t} \nabla f(\by^{(t)})\right)$, by the Proximal Property-I (Lemma~\ref{lemma:prox_prop1}), we have
$$
\widetildebx^{(t)} - \frac{\eta_t}{\gamma_t} \nabla f(\by^{(t)}) - \widetildebx^\toptone 
\in \frac{\eta_t}{\gamma_t} \partial g(\widetildebx^\toptone),
$$
By the subgradient inequality (Definition~\ref{definition:subgrad}), for any $\bx\in\real^n$, we have
$$
g(\bx) \geq  g(\widetildebx^\toptone) + 
\frac{\gamma_t}{\eta_t} \innerproduct{\widetildebx^{(t)} - \frac{\eta_t}{\gamma_t} \nabla f(\by^{(t)}) - \widetildebx^\toptone , \bx - \widetildebx^\toptone}.
$$
By the convexity of $g$ and the update rule $\bx^\toptone = (1 - \gamma_t)\bx^{(t)} + \gamma_t \widetildebx^\toptone$, 
$$
g(\bx^\toptone) \leq (1 - \gamma_t) g(\bx^{(t)} ) + \gamma_t g(\widetildebx^\toptone).
$$
Substituting  the expression for $g(\widetildebx^\toptone)$ into the preceding inequality obtains
$$
g(\bx^\toptone) \leq (1 - \gamma_t) g(\bx^{(t)})
+ \gamma_t \left[ g(\bx) - \innerproduct{\frac{\gamma_t}{\eta_t} (\widetildebx^{(t)} - \widetildebx^{(t+1)}) - \nabla f(\by^\toptzero), \bx - \widetildebx^{(t+1)}} \right].
$$
Given that $\eta_t = \frac{1}{\beta}$, by the $\beta$-smoothness of $f$ (Definition~\ref{definition:scss_func}), we obtain
$$
f(\bx^\toptone) \leq 
f(\by^\toptzero) + \innerproduct{\nabla f(\by^\toptzero), \bx^\toptone - \by^\toptzero} + \frac{1}{2\eta_t} \normtwo{\bx^\toptone - \by^\toptzero}^2.
$$
Subtracting the two updare rules:
$\bx^\toptone = (1 - \gamma_t)\bx^{(t)} + \gamma_t \widetildebx^\toptone$ and
$\by^\toptzero = (1 - \gamma_t)\bx^{(t)} + \gamma_t \widetildebx^{(t)}$, 
we have 
$\bx^{(t+1)} - \by^\toptzero = \gamma_t (\widetildebx^{(t+1)} - \widetildebx^{(t)})$. 
Substituting this equality into the above inequality, using the update rule for $\bx^\toptone$ again, and letting $\be^{(t)} \triangleq \widetildebx^{(t+1)} - \widetildebx^{(t)}$, 
we get
$$
\footnotesize
\begin{aligned}
&f(\bx^\toptone) \leq 
f(\by^\toptzero) + \innerproduct{\nabla f(\by^\toptzero), (1 - \gamma_t)\bx^{(t)} + \gamma_t \widetildebx^\toptone - \by^\toptzero} + \frac{\gamma_t^2}{2\eta_t} \normtwobig{\be^{(t)}}^2\\
&= (1 - \gamma_t) \Big\{f(\by^\toptzero) + \innerproduct{\nabla f(\by^\toptzero), \bx^{(t)} - \by^\toptzero}\Big\}
+ \gamma_t \Big\{f(\by^\toptzero) + \innerproduct{\nabla f(\by^\toptzero), \widetildebx^\toptone - \by^\toptzero }\Big\} 
+ \frac{\gamma_t^2}{2\eta_t} \normtwobig{\be^{(t)}}^2\\
&\leq (1 - \gamma_t) f(\bx^{(t)}) + \gamma_t \left\{f(\by^\toptzero) + \innerproduct{\nabla f(\by^\toptzero), \widetildebx^\toptone - \by^\toptzero}\right\}
+ \frac{\gamma_t^2}{2\eta_t} \normtwobig{\be^{(t)}}^2,
\end{aligned}
$$
where the last inequality follows from the convexity of $f$ (Theorem~\ref{theorem:conv_gradient_ineq}).

Combining the inequalities for $f(\bx^\toptone)$ and $g(\bx^\toptone)$ with $\bx\triangleq\bx^*$, along with the convexity of $f$ such that $
f(\bx^*) \geq f(\by^\toptzero) + \innerproduct{\nabla f(\by^\toptzero), \bx^* - \by^\toptzero}
$, yields that 
$$
\small
\begin{aligned}
&F(\bx^\toptone) - F(\bx^*) - (1 - \gamma_t)\big(F(\bx^{(t)}) - F(\bx^*)\big)\\
&\leq 
\frac{\gamma_t^2}{\eta_t}\innerproduct{\widetildebx^\toptone -\widetildebx^\toptzero, \bx^*-\widetildebx^\toptone}  + \frac{\gamma_t^2}{2\eta_t} \normtwobig{\widetildebx^{(t+1)} - \widetildebx^{(t)}}^2\\
&\stackrel{\dag}{=} \frac{\gamma_t^2}{2\eta_t} \left\{ \normtwo{\widetildebx^\toptzero-\bx^*}^2 - \normtwo{\widetildebx^\toptone-\bx^*}^2 \right\}
\end{aligned}
$$
where the equality ($\dag$) follows from the fact that $2\innerproduct{\ba-\bb,\bc-\ba} +\normtwo{\ba-\bb}^2 = \normtwo{\bb-\bc}^2-\normtwo{\ba-\bc}^2$ for any $\ba,\bb,\bc\in\real^n$.
Recalling that $\eta_t$ and $\gamma_t$ are chosen such that
$
\frac{1 - \gamma_t}{\gamma_t^2} \eta_t \leq \frac{1}{\gamma_{t-1}^2} \eta_{t-1}
$ (same as the FISTA method),
we obtain an iterate inequality for consecutive iterations
$$
\small
\begin{aligned}
	&\frac{\eta_t}{\gamma_t^2} \big(F(\bx^\toptone) - F(\bx^*)\big) + \frac{1}{2} \normtwo{\widetildebx^{(t+1)} - \bx^*}^2 
	\leq (1-\gamma_t)\frac{\eta_t}{\gamma_t^2} \big(F(\bx^{(t)}) - F(\bx^*)\big) + \frac{1}{2} \normtwo{\widetildebx^{(t)} - \bx^*}^2\\
	&\leq \frac{\eta_{t-1}}{\gamma_{t-1}^2} \big(F(\bx^{(t)}) - F(\bx^*)\big) + \frac{1}{2} \normtwo{\widetildebx^{(t)} - \bx^*}^2
	\leq \frac{\eta_1}{\gamma_1^2} \big(F(\bx^{(2)}) - F(\bx^*)\big) + \frac{1}{2} \normtwo{\widetildebx^{(2)} - \bx^*}^2\\
	&\stackrel{\dag}{\leq} \frac{(1 - \gamma_1)\eta_1}{\gamma_1^2} \big(F(\bx^{(1)}) - F(\bx^*)\big) + \frac{1}{2} \normtwo{\widetildebx^{(1)} - \bx^*}^2
	=\frac{1}{2} \normtwo{\bx^{(1)} - \bx^*}^2,
\end{aligned}
$$
where the last equality follows from fact that  $\gamma_1 = 1$ and $\bx^{(1)} = \widetildebx^{(1)}$.
Using Lemma~\ref{lemma:seq_gamma_zeta} concludes the result.
\end{proof}

Similarly, note that the key steps in deriving Theorem~\ref{theorem:nesterov_acc} rely on conditions in  \eqref{equation:fista_require}. Therefore, if the smoothness constant is not known beforehand, the line search Algorithm~\ref{alg:fistav_line} can guarantee the same convergence results.

\begin{algorithm}[h] 
\caption{Conditional Gradient  Method}
\label{alg:cond_gen}
\begin{algorithmic}[1] 
\Require A function $f(\bx)$ and a set $\sS$; 
\State {\bfseries Input:}  Initialize $\bx^{(0)}$;
\For{$t=0,1,2,\ldots$}
\State Pick a stepsize $\eta_t\in[0,1]$;
\State $\widetildebx^{(t)} \leftarrow \mathop{\argmin}_{\bx\in \sS} \innerproduct{\nabla f(\bx^{(t)}), \bx}$; \Comment{linear optimization}
\State  $\bx^{(t+1)} \leftarrow \bx^{(t)} + \eta_t (\widetildebx^{(t)} - \bx^{(t)}) = (1-\eta_t)\bx^{(t)} +  \eta_t \widetildebx^{(t)} $;  \Comment{``convex" update}
\EndFor
\State {\bfseries Return:}   $\bx_{\text{final}} \leftarrow \bx^\toptzero$;
\end{algorithmic} 
\end{algorithm}
\section{Conditional Gradient  (Frank-Wolfe) Method}\label{section:cond_gd}
We investigate  the \textit{conditional gradient method}, commonly referred to as the \textit{Frank-Wolfe (FW) algorithm}. 
This method provides an appealing alternative to constrained optimization problems, particularly by addressing potential computational inefficiencies associated with the projection step in projected gradient descent methods.
Consider the following optimization problem:
\begin{equation}\label{equation:cgd_prob}
\text{(P2)}:\qquad \mathop{\min}_{\bx} f(\bx)\gap \text{s.t.}\gap \bx\in\sS.
\end{equation}
In the conditional gradient Method, the next step is computed as a convex combination of the current iterate and a minimizer of a linearized form of the objective function over the feasible set $\sS$ (Algorithm~\ref{alg:cond_gen}).

Recall that in the PGD method, the update rule from the $t$-th to the $(t+1)$-th iteration can be expressed as follows:
\begin{center}
\framebox{
\begin{minipage}{0.95\textwidth}
\begin{equation}\label{equation:pgd_decom_cg}
\small
\textbf{(PGD)}:\qquad 
\begin{aligned}
\bx^\toptone 
&\leftarrow \mathop{\argmin}_{\bx\in\sS} \normtwo{\bx - \left(\bx^{(t)} - \eta_t \nabla f(\bx^{(t)})\right)}^2\\
&=\mathop{\argmin}_{\bx\in\sS} 
\left\{
\frac{1}{2} \normtwobig{\bx-\bx^\toptzero}^2+ \eta_t\innerproduct{\nabla f(\bx^\toptzero), \bx} 
\right\},
\end{aligned}
\end{equation}
\end{minipage}
}
\end{center}
where the quadratic term acts as a regularization ensuring the update remains within a neighborhood of the current iterate $\bx^\toptzero$.
In contrast, the CG method employs a convex combination of the new update found by linear search  $\mathop{\argmin}_{\bx\in\sS}\innerproduct{\nabla f(\bx^{(t)}), \bx}$ and the current iterate $\bx^\toptzero$,  avoiding the use of quadratic regularization. This sum can also be viewed as a form of regularization, ensuring the new iterate does not stray too far from the current iterate:
\begin{tcolorbox}[colback=white,colframe=black]
\begin{minipage}{1\textwidth}
\begin{equation}\label{equation:cgmethod}
\small
\textbf{(CG)}:\qquad 
\begin{aligned}
\bx^\toptone 
&\leftarrow (1-\eta_t)\bx^{(t)} + \eta_t\mathop{\argmin}_{\bx\in\sS}\innerproduct{\nabla f(\bx^{(t)}), \bx} .
\end{aligned}
\end{equation}
\end{minipage}
\end{tcolorbox}
Similarly, CG can achieve a rate of $\mathcalO(T)$ for convex and smooth functions.
\begin{theoremHigh}[FW for Convex and SS: $\mathcalO(T)$]\label{theorem:fw_under_smoo}
Let $f: \sS \rightarrow \real$ be a proper, differentiable, convex, and $\beta$-smooth function defined over the convex domain $\sS\subseteq\real^n$. 
Suppose that $f$ attains its global minimum at a point $\bx^* \in \sS$, and that $D$ is the diameter of $\sS$, defined as $D = \max_{\bx, \by \in \sS} \normtwo{\bx - \by}$.
Let $ \{\bx^\toptzero\}_{t \geq 0} $~\footnote{For simplicity, note that we suppose the iterate starts by index 0 for the CG method.} be the sequence generated by the Frank-Wolfe method (Algorithm~\ref{alg:cond_gen}) for solving problem (P2)  with a dynamic stepsize rule in which $\eta_t = \frac{2}{t+2}$ at the $t$-th iteration. 
Then, for any $T\geq0$,
$$
f(\bx^{(T)}) - f(\bx^*) \leq \frac{2\beta D^2}{T+2}.
$$
\end{theoremHigh}

\begin{proof}[of Theorem~\ref{theorem:fw_under_smoo}]
By smoothness, convexity, and the progress rule of conditional gradient method, we have
$$
\begin{aligned}
f(\bx^{(t+1)}) 
&\leq f(\bx^{(t)}) + \nabla f(\bx^{(t)})^\top (\bx^{(t+1)} - \bx^{(t)}) + \frac{\beta}{2} \normtwo{ \bx^{(t+1)} - \bx^{(t)}}^2\quad &&(\text{smoothness})\\
&=f(\bx^{(t)}) + \eta_t \nabla f(\bx^{(t)})^\top (\widetildebx^{(t)} - \bx^{(t)}) + \frac{\eta_t^2 \beta}{2} \normtwo{ \widetildebx^{(t)} - \bx^{(t)} }^2\quad &&\text{(update rule)}\\
&\leq f(\bx^{(t)}) + \eta_t \nabla f(\bx^{(t)})^\top (\bx^* - \bx^{(t)}) + \frac{\eta_t^2 \beta D^2}{2}. &&\text{(linear update rule)}
\end{aligned}
$$
Due to convexity, we also have:
$
\nabla f(\bx^{(t)})^\top (\bx^* - \bx^{(t)}) \leq f(\bx^*) - f(\bx^{(t)}).
$
Combining the two inequalities yields
\begin{equation}\label{equation:fw_under_smoo1}
f(\bx^{(t+1)}) - f(\bx^*) \leq (1 - \eta_t)(f(\bx^{(t)}) - f(\bx^*)) + \frac{\eta_t^2 \beta D^2}{2}.
\end{equation}
We use induction in order to prove $f(\bx^{(t)}) - f(\bx^*) \leq \frac{2\beta D^2}{t+2}$ based on the above equation.

\paragraph{Base case $t = 0$.} When $t = 0$, we have $\eta_0 = \frac{2}{2} = 1$, and \eqref{equation:fw_under_smoo1} reduces to
$$
\begin{aligned}
f(\bx^{(1)}) - f(\bx^*) 
&\leq (1 - \eta_0)(f(\bx_0) - f(\bx^*)) + \frac{\beta}{2} \normtwo{\bx^{(1)} -\bx^*}^2 
=  \frac{\beta D^2}{2}
\leq \frac{2\beta D^2}{2}.
\end{aligned}
$$
\paragraph{Inductive step.} Proceeding by induction, we assume that 
$
f(\bx^{(t)}) - f(\bx^*) \leq \frac{2\beta D^2}{t+2}
$
holds for all integers up to $t$ and we show the claim for $t+1$. By \eqref{equation:fw_under_smoo1},
\begin{align*}
f(\bx^{(t+1)}) - f(\bx^*) &\leq \left(1 - \frac{2}{t+2}\right) \left(f(\bx^{(t)}) - f(\bx^*)\right) + \frac{2\beta D^2}{(t+2)^2} \\
&\leq \left(1 - \frac{2}{t+2}\right) \frac{2\beta D^2}{t+2} + \frac{2\beta D^2}{(t+2)^2} 
= 2\beta D^2 \frac{t+1}{(t+2)^2} \\
&= 2\beta D^2 \frac{t+1}{t+2} \cdot \frac{1}{t+2} 
\leq 2\beta D^2 \frac{t+2}{t+3} \cdot \frac{1}{t+2} 
= 2\beta D^2 \frac{1}{t+3}.
\end{align*}
Thus, the inequality also holds for the $t+1$ case.
\end{proof}

The CG method is useful for understanding the \textit{power method} for computing the eigenpair of a matrix.
\begin{example}[Power Method]
Consider the specific problem 
$$
\mathop{\min}_{\bx} f(\bx)=-\frac{1}{2} \bx^\top\bA\bx,\gap \text{s.t.}\gap \sS=\{\norm{\bx}_2\leq 1\},
$$
where $\bA$ is positive semidefinite. 
Then, $\bt^{(t)} = \argmin_{\bt\in\sS} \innerproduct{\bt, \nabla f(\bx^{(t)})}$ is given by $\frac{\bA\bx^{(t)}}{\norm{\bA\bx^{(t)}}_2}$. The conditional gradient update is:
$$
\bx^{(t+1)} = (1-s_t)\bx^{(t)} + s_t \frac{\bA\bx^{(t)}}{\norm{\bA\bx^{(t)}}_2}.
$$
If we set the stepsizes using an exact line search strategy, then 
$$
s_t =\mathop{\argmin}_{{s\in[0,1]}} f(\bx^{(t)}+s_t(\bt^{(t)} - \bx^{(t)})).
$$
When $\bx^{(t)}$ is not the optimal solution of the problem, and since $f(\bx)$ is concave, we can choose $s_t=1$. Therefore, the update becomes 
$$
\bx^{(t+1)} = \frac{\bA\bx^{(t)}}{\norm{\bA\bx^{(t)}}_2}.
$$
This is equivalent to the power method for finding the eigenvector corresponding to the maximal eigenvalue. The sequence for the eigenvalue is
$$
\lambda^{{(t+1)}} = \bx^{(t+1)\top}\bA\bx^{(t+1)}.
$$
See, for example, \citet{lu2021numerical} for more details.
\end{example}

\begin{algorithm}[h] 
\caption{Generalized Conditional Gradient  Method}
\label{alg:gen_cond_gen}
\begin{algorithmic}[1] 
\Require A function $f(\bx)$ and a convex function $g(\bx)$; 
\State {\bfseries Input:}  Initialize $\bx^{(1)}$;
\For{$t=1,2,\ldots$}
\State Pick a stepsize $\eta_t\in[0,1]$;
\State $\widetildebx^{(t)} \leftarrow \mathop{\argmin}_{\bx\in \textcolor{mylightbluetext}{\real^n}} \innerproduct{\nabla f(\bx^{(t)}), \bx} + g(\bx)$; \Comment{``linear" optimization}
\State  $\bx^{(t+1)} \leftarrow \bx^{(t)} + \eta_t (\widetildebx^{(t)} - \bx^{(t)}) = (1-\eta_t)\bx^{(t)} +  \eta_t \widetildebx^{(t)} $;  \Comment{``convex" update}
\EndFor
\State {\bfseries Return:}   $\bx_{\text{final}} \leftarrow\bx^\toptzero$;
\end{algorithmic} 
\end{algorithm}
\section{Generalized Conditional Gradient Method}

We have shown that the conditional gradient method addresses the same constrained problem, \eqref{equation:cgd_prob}, as the projected gradient descent approach. However, the \textit{generalized conditional gradient (GCG)} method (Algorithm~\ref{alg:gen_cond_gen}) focuses on solving an unconstrained composite problem:
\begin{equation}\label{equation:gene_cgd_prob}
\text{(P3)}:\qquad \min \{F(\bx) \triangleq f(\bx)+g(\bx)\}.
\end{equation}
In this context, we assume that $f$ and $g$ satisfy conditions (B1) and (B2) in Assumption~\ref{assumption:gcg}:
\begin{assumption}[GCG]\label{assumption:gcg}
For problem (P3), we assume
\begin{itemize}
	\item[(B1)] $ g: \real^n \rightarrow (-\infty, \infty] $ is proper, closed, and convex and \textcolor{mylightbluetext}{$\dom(g)$ is compact}.
	\item[(B2)] $ f: \real^n \rightarrow (-\infty, \infty] $ is proper, closed, and  $ \beta $-smooth over $ \dom(f) $ ($ \beta > 0 $), which is assumed to be an open and convex set satisfying $ \dom(g) \subseteq \dom(f) $.
\end{itemize}
\end{assumption}
The rationale behind these assumptions for the GCG method includes:
\begin{itemize}
\item \textit{Compactness of $\dom(g)$.}  (B1) ensures that the subproblem
$
\underset{\by \in \dom(g)}{\arg\min} \nabla f(\bx)^\top \by  + g(\bx)
$
has a solution. The compactness of $ \dom(g) $ ensures boundedness and avoids unbounded solutions in the subproblem.

\item \textit{Open and convex $\dom(f)$.} The openness of $ \dom(f) $ in (B2) allows the linear subproblem to remain well-defined and ensures that iterates stay in the feasible region.

\item \textit{Domain compatibility.}  The inclusion $ \dom(g) \subseteq \dom(f) $ ensures compatibility, i.e., the linear subproblem remains within the domain of $ f $, and $ f(\bx) $ can be evaluated and differentiated.
\end{itemize}
The  generalized conditional gradient method shares similarities with the standard conditional gradient method, but differs by computing a minimizer of the sum of the linearized smooth part of $f$ and the convex function $g$. Since the objective consists of a sum of convex and affine functions, the problem remains relatively straightforward to solve \citep{beck2017first}.

A natural optimality measure in proximal gradient methods is the gradient mapping. However, the analysis of the conditional gradient method relies on a different optimality measure, which we will refer to as the \textit{conditional gradient norm}.

\begin{definition}[Conditional Gradient Norm]\label{definition:cond_grad_norm}
	Suppose that $f$ is $\beta$-smooth and $g$ is \textcolor{black}{closed} convex. Then, the conditional gradient norm is the function $\mathcalW : \dom(f) \rightarrow \real$ defined by
	$$
	\begin{aligned}
		\mathcalW(\bx) &\triangleq \innerproduct{\nabla f(\bx), \bx - \widetildebx} + g(\bx) - g(\widetildebx)\\
		&= \max_{\bv \in \real^n} \{\innerproduct{\nabla f(\bx), \bx - \bv} + g(\bx) - g(\bv)\}\\
		&= \innerproduct{\nabla f(\bx), \bx} + g(\bx) + g^*(-\nabla f(\bx)).
	\end{aligned}
	$$
where $\widetildebx\in\mathop{\argmin}_{\bv\in\real^n} \innerproduct{\nabla f(\bx), \bv} + g(\bv)$, and $g^*$ represents the conjugate function of $g$ (Definition~\ref{definition:conjug_func}).
\end{definition}

\begin{theorem}[Conditional Gradient Norm as an Optimality Measure]\label{theorem:cgn_optmea}
Suppose that $f$ is $\beta$-smooth and $g$ is \textcolor{black}{closed} convex, satisfying assumption (B1) and (B2) in Assumption~\ref{assumption:gcg}. Then,
\begin{enumerate}[(i)]
\item $\mathcalW(\bx) \geq 0$ for any $\bx \in \dom(f)$;
\item $\mathcalW(\bx^*) = 0$ if and only if $- \nabla f(\bx^*) \in \partial g(\bx^*)$, i.e., if and only if $\bx^*$ is a stationary point of problem (P3); see Theorem~\ref{theorem:opt_cond_p3} and Definition~\ref{definition:stat_opt_cond_p3}.
\end{enumerate}
\end{theorem}

\begin{proof}
\textbf{(i).} This follows from the third expression of the conditional gradient norm in Definition~\ref{definition:cond_grad_norm}  and Fenchel's inequality (Theorem~\ref{theorem:fenchel_ineq}).

\paragraph{(ii).} By part (i), it follows that $\mathcalW(\bx^*) = 0$ if and only if $\mathcalW(\bx^*) \leq 0$, which is the same as the relation (using the second expression of conditional gradient norm in Definition~\ref{definition:cond_grad_norm} for $\mathcalW(\bx^*)$)
$$
\begin{aligned}
&\innerproduct{ \nabla f(\bx^*), \bx^* - \bv} + g(\bx^*) - g(\bv) \leq 0 \quad \text{for all } \bv \in \real^n\\
&\qquad\qquad\implies\qquad  g(\bv) \geq g(\bx^*) + \innerproduct{-\nabla f(\bx^*), \bv - \bx^*},
\end{aligned}
$$
which is equivalent to the relation $-\nabla f(\bx^*) \in \partial g(\bx^*)$, namely, to stationarity (Definition~\ref{definition:stat_opt_cond_p3}).
\end{proof}

Similar to descent lemma for GD and PGD (Lemma~\ref{lemma:gdupd_sm_sconv}),
the basic inequality that will be used in the analysis of the GCG method is the following recursive inequality.

\begin{lemma}[Descent Lemma for GCG]\label{lemma:desc_gcg}
Suppose that $f$ is $\beta$-smooth and $g$ is \textcolor{black}{closed} convex, satisfying assumption (B1) and (B2) in Assumption~\ref{assumption:gcg}. Let $\bx \in \dom(g)$ and $\eta \in [0, 1]$. 
Consider problem (P3).
Then,
\begin{equation}
F(\bx + \eta(\widetildebx - \bx)) \leq F(\bx) - \eta \mathcalW(\bx) + \frac{\eta^2 \beta}{2} \normtwo{\widetildebx - \bx}^2. 
\end{equation}
\end{lemma}
\begin{proof}[of Lemma~\ref{lemma:desc_gcg}]
By the smoothness of $f$ and  the convexity of $g$, we can write the following:
$$
\small
\begin{aligned}
F(\bx + &\eta(\widetildebx - \bx)) 
= f(\bx + \eta(\widetildebx - \bx)) + g(\bx + \eta(\widetildebx - \bx)) \\
&\leq f(\bx) + \eta \innerproduct{\nabla f(\bx),  \widetildebx -\bx } + \frac{\eta^2 \beta}{2} \normtwo{\widetildebx - \bx}^2 + (1 - \eta)g(\bx) + \eta g(\widetildebx) \\
&= F(\bx) + \eta\big(\innerproduct{\nabla f(\bx), \widetildebx -\bx } - g(\bx) + g(\widetildebx)\big) + \frac{\eta^2 \beta}{2} \normtwo{\widetildebx - \bx}^2 \\
&= F(\bx) - \eta \mathcalW(\bx) + \frac{\eta^2 \beta}{2} \normtwo{\widetildebx - \bx}^2.
\end{aligned}
$$
This completes the proof.
\end{proof}

\begin{lemma}[Monotonicity and Sufficient Decrease for GCG]\label{lemma:mono_suff_gcg}
Suppose that $f$ is $\beta$-smooth and $g$ is \textcolor{black}{closed} convex, satisfying assumption (B1) and (B2) in Assumption~\ref{assumption:gcg}. 
Let $ \{\bx^\toptzero\}_{t > 0} $ be the sequence generated by the generalized conditional gradient method (Algorithm~\ref{alg:gen_cond_gen}) for solving problem (P3) with adaptive stepsizes $\eta_t \triangleq \left\{ 1, \frac{\mathcalW(\bx^\toptzero)}{\beta \normtwo{\widetildebx^\toptzero -\bx^\toptzero}^2} \right\}$  for each $t$-th iteration. Then, for any $ t > 0 $,
$$
F(\bx^\toptzero) - F(\bx^\toptone) \geq \frac{1}{2} \min \left\{ \mathcalW(\bx^\toptzero), \frac{\mathcalW^2(\bx^\toptzero)}{\beta D^2} \right\}
$$
where $ D $ is an upper bound on the diameter of $ \dom(g) $:
$
D \geq \max_{\bx,\by \in \dom(g)} \normtwo{\bx - \by}.
$
Since $\mathcalW(\bx^\toptzero)\geq 0$ by Theorem~\ref{theorem:cgn_optmea} and equals to 0 only for stationary points, this also shows that 
\begin{itemize}
\item $F(\bx^\toptzero) - F(\bx^\toptone) \geq 0$.
\item $F(\bx^\toptzero) > F(\bx^\toptone)$  if $\bx^\toptzero$  is not a stationary point of problem (P3).
\end{itemize}
\end{lemma}
\begin{proof}[of Lemma~\ref{lemma:mono_suff_gcg}]
For any iteration $ t > 0 $ and the update rule in Algorithm~\ref{alg:gen_cond_gen}, it holds that $ \bx^\toptone = \bx^\toptzero + \eta_t (\widetildebx^\toptzero - \bx^\toptzero) $, where
$
\eta_t = \min \left\{ 1, \frac{\mathcalW(\bx^\toptzero)}{\beta \normtwo{\widetildebx^\toptzero - \bx^\toptzero}^2} \right\}.
$
By the descent lemma for GCG (Lemma~\ref{lemma:desc_gcg}), we have
\begin{equation}\label{equation:mono_suff_gcg}
F(\bx^\toptzero) - F(\bx^\toptone) \geq \eta_t \mathcalW(\bx^\toptzero) - \frac{\eta_t^2 \beta}{2} \normtwo{\widetildebx^\toptzero - \bx^\toptzero}^2.
\end{equation}
Then,
\begin{itemize}
\item \textit{Case 1: $\frac{\mathcalW(\bx^\toptzero)}{\beta \normtwo{\widetildebx^\toptzero - \bx^\toptzero}^2} < 1$.}
We have 
$
F(\bx^\toptzero) - F(\bx^\toptone) \geq \frac{\mathcalW^2(\bx^\toptzero)}{2\beta \normtwo{\widetildebx^\toptzero - \bx^\toptzero}^2} \geq \frac{\mathcalW^2(\bx^\toptzero)}{2\beta D^2}.
$

\item \textit{Case 2: $\frac{\mathcalW(\bx^\toptzero)}{\beta \normtwo{\widetildebx^\toptzero - \bx^\toptzero}^2} \geq 1$, $\eta_t=1$.}
We have
$
F(\bx^\toptzero) - F(\bx^\toptone) \geq \mathcalW(\bx^\toptzero) - \frac{\beta}{2} \normtwo{\widetildebx^\toptzero - \bx^\toptzero}^2 \geq \frac{1}{2} \mathcalW(\bx^\toptzero).
$
\end{itemize}
Combining the two cases obtains the desired result.
\end{proof}

When $ f $ is convex, the conditional gradient norm is lower bounded by the distance to optimality in terms of function values.
\begin{lemma}[Conditional Gradient Norm for Convex Functions]\label{lemma:con_gcg_cgnorm}
Suppose that $f$ is $\beta$-smooth and $g$ is \textcolor{black}{closed} convex, satisfying assumption (B1) and (B2) in Assumption~\ref{assumption:gcg}, and suppose additionally that $ f $ is convex. Then, for any $ \bx \in \dom(g) $ and any optimizer $\bx^*$ of $F$,
$$
\mathcalW(\bx) \geq F(\bx) - F(\bx^*).
$$
\end{lemma}
\begin{proof}[of Lemma~\ref{lemma:con_gcg_cgnorm}]
For any $ \bx \in \dom(g) $,
$$
\small
\begin{aligned}
\mathcalW(\bx) 
&= \innerproduct{\nabla f(\bx), \bx - \widetildebx} + g(\bx) - g(\widetildebx) 
= \innerproduct{\nabla f(\bx), \bx} + g(\bx) - \big(\innerproduct{\nabla f(\bx), \widetildebx} + g(\widetildebx)\big) \\
&\geq \innerproduct{\nabla f(\bx), \bx} + g(\bx) - \big(\innerproduct{\nabla f(\bx), \bx^*} + g(\bx^*)\big) 
= \innerproduct{\nabla f(\bx), \bx - \bx^*} + g(\bx) - g(\bx^*) \\
&\geq f(\bx) - f(\bx^*) + g(\bx) - g(\bx^*)  
= F(\bx) - F(\bx^*).
\end{aligned}
$$
where the last inequality follows from the convexity of $f$. This completes the proof.
\end{proof}

Using these lemmas, we then prove the convergence of GCG for smooth functions.

\begin{theoremHigh}[GCG for SS $f$]\label{theorem:gcg_ss}
Suppose that $f$ is $\beta$-smooth and $g$ is \textcolor{black}{closed} convex, satisfying assumption (B1) and (B2) in Assumption~\ref{assumption:gcg}. 
Let $ \{\bx^\toptzero\}_{t > 0} $ be the sequence generated by the generalized conditional gradient method (Algorithm~\ref{alg:gen_cond_gen}) for solving problem (P3) with adaptive stepsizes $\eta_t \triangleq \left\{ 1, \frac{\mathcalW(\bx^\toptzero)}{\beta \normtwo{\widetildebx^\toptzero -\bx^\toptzero}^2} \right\}$  for each $t$-th iteration.
Then,
\begin{enumerate}[(i)]
\item $\mathcalW(\bx^\toptzero) \rightarrow 0$ { as } $t \rightarrow \infty$.
\item For any $t > 0$, $\underset{t=\{1,\ldots,T\}}{\min}\mathcalW(\bx^\toptzero) \leq \max \left\{ \frac{2(F(\bx^{(1)}) - F(\bx^*))}{T}, \frac{\sqrt{2\beta D^2 (F(\bx^{(1)}) - F(\bx^*))}}{\sqrt{T}} \right\}$, where $ D $ is an upper bound on the diameter of $ \dom(g) $.
\item All limit points of the sequence $ \{\bx^\toptzero\}_{t > 0} $ are stationary points of problem (P3).
\end{enumerate}
\end{theoremHigh}
\begin{proof}[of Theorem~\ref{theorem:gcg_ss}]
\textbf{(i).} Since $ \{F(\bx^\toptzero)\}_{t>0} $ is nonincreasing and bounded below by $ F(\bx^*) $ (Lemma~\ref{lemma:mono_suff_gcg}), it follows that it is convergent, and in particular, $ F(\bx^\toptzero) - F(\bx^\toptone) \rightarrow 0 $ as $ t \rightarrow \infty $. Therefore, by the sufficient decrease in Lemma~\ref{lemma:mono_suff_gcg}, it follows that $ \min \left\{ \mathcalW(\bx^\toptzero), \frac{\mathcalW^2(\bx^\toptzero)}{\beta D^2} \right\} \rightarrow 0 $ as $ t \rightarrow \infty $, implying that $ \mathcalW(\bx^\toptzero) \rightarrow 0 $ as $ t \rightarrow \infty $.

\paragraph{(ii).} By the sufficient decrease in Lemma~\ref{lemma:mono_suff_gcg}, for all $t  \geq 0 $,
$$
F(\bx^\toptzero) - F(\bx^\toptone) \geq \frac{1}{2} \min \left\{ \mathcalW(\bx^\toptzero), \frac{\mathcalW^2(\bx^\toptzero)}{\beta D^2} \right\}.
$$
Telescoping the sum over $ t = 1,2,\ldots,T  $,
\begin{equation*}
\begin{aligned}
F(\bx^{(1)}) - F(\bx^*)  &\geq F(\bx^{(1)}) - F(\bx^{(T+1)}) \geq \frac{1}{2} \sum_{t=1}^{T} \min \left\{ \mathcalW(\bx^\toptzero), \frac{\mathcalW^2(\bx^\toptzero)}{\beta D^2} \right\}\\
&\geq
\frac{T}{2} \min_{t \in \{1, 2, \ldots, T\}} \left[ \min \left\{ \mathcalW(\bx^\toptzero), \frac{\mathcalW^2(\bx^\toptzero)}{\beta D^2} \right\} \right].
\end{aligned}
\end{equation*}
This implies  that there exists an $ t \in \{1, 2, \ldots, T\} $ for which
$$
\min \left\{ \mathcalW(\bx^\toptzero), \frac{\mathcalW^2(\bx^\toptzero)}{\beta D^2} \right\} \leq \frac{2\big(F(\bx^{(1)}) - F(\bx^*)\big)}{T}.
$$
That is,
$$
\mathcalW(\bx^\toptzero) \leq \max \left\{ \frac{2(F(\bx^{(1)}) - F(\bx^*))}{T}, \frac{\sqrt{2\beta D^2 (F(\bx^{(1)}) - F(\bx^*))}}{\sqrt{T}} \right\}.
$$

\paragraph{(iii).} Suppose that $ \bx^+ $ is a limit point of $ \{\bx^\toptzero\}_{t>0} $. Then there exists a subsequence $ \{\bx^{(t_j)}\}_{j > 0} $ that converges to $ \bx^+ $. By the definition of the conditional gradient norm $ \mathcalW(\cdot) $, it follows that for any $ \bv\in\real^n $,
$$
\mathcalW(\bx^{(t_j)}) \geq \innerproduct{\nabla f(\bx^{(t_j)}), \bx^{(t_j)} - \bv} + g(\bx^{(t_j)}) - g(\bv).
$$
Taking the limit as $ j \rightarrow \infty $ and using the fact that $ \mathcalW(\bx^{(t_j)}) \rightarrow 0 $ as $ j \rightarrow \infty $, along with the continuity of $ \nabla f $ and the closedness of $ g $, we obtain that
$$
0 \geq \innerproduct{\nabla f(\bx^+), \bx^+ - \bv} + g(\bx^+) - g(\bv) \text{ for any } \bv \in \real^n,
$$
which is the same as the relation $ -\nabla f(\bx^+) \in \partial g(\bx^+) $,  indicating stationarity.
\end{proof}

We will demonstrate that the iterates of function values in the GCG approach for convex and smooth functions follow a quadratic inequality form. The following lemma summarizes the convergence results for such a form.
\begin{lemma}\label{lemma:gcg_conv_ss_usg}
Let $B$ be a positive integer, and let $\{a_t\}_{t > 0}$ and $\{b_t\}_{t > 0}$ be nonnegative sequences satisfying for any $t > 0$
\begin{equation}\label{equation:ineq_gcg_conv_ss_usg}
a_{t+1} \leq a_t - \gamma_t b_t + \frac{A}{2} \gamma_t^2,
\end{equation}
where $\gamma_t = \frac{2}{t + 2B}$ and $A$ is a positive number. Suppose that $a_t \leq b_t$ for all $t>0$. Then,
\begin{enumerate}[(i)]
\item For any $T \geq 1$, $a_T \leq \frac{2\max\{A, a_1(B-1)\}}{T + 2B - 2}$.
\item For any $T \geq 3$,
$
\min_{t=\lfloor T/2 \rfloor + 2, \ldots, T} b_t \leq \frac{8\max\{A, a_1(B-1)\}}{T - 2}.
$
\end{enumerate}
\end{lemma}
\begin{proof}[of Lemma~\ref{lemma:gcg_conv_ss_usg}]
\textbf{(i).} 
Since $a_{t+1} \leq a_t - \gamma_t b_t + \frac{A}{2} \gamma_t^2$ and $a_t \leq b_t$ for all $t>0$, it follows that
$a_{t+1} \leq (1 - \gamma_t) a_t + \frac{A}{2} \gamma_t^2$.
Therefore,
$$
\begin{aligned}
a_2 &\leq (1 - \gamma_1) a_1 + \frac{A}{2} \gamma_1^2; \\
a_3 &\leq (1 - \gamma_2) a_2 + \frac{A}{2} \gamma_2^2 = (1 - \gamma_2)(1 - \gamma_1) a_1 + (1 - \gamma_2) \frac{A}{2} \gamma_1^2 + \frac{A}{2} \gamma_2^2; \\
a_4 &\leq (1 - \gamma_3) a_3 + \frac{A}{2} \gamma_3^2 = a_1 \prod_{t=1}^{3} (1 - \gamma_t)
+ \frac{A}{2} \big((1 - \gamma_3)(1 - \gamma_2)  \gamma_1^2 + (1 - \gamma_3)  \gamma_2^2 +  \gamma_3^2\big);\\
\vdots &\leq  \vdots.\\
\end{aligned}
$$
That is, for any $T>1$, the general term formula for $a_T$ is 
\begin{equation}\label{equation:gcg_conv_ss_usg1}
a_T \leq a_1 \prod_{t=1}^{T-1} (1 - \gamma_t) + \frac{A}{2} \sum_{i=1}^{T-1} \left[ \prod_{t=i+1}^{T-1} (1 - \gamma_t) \right] \gamma_i^2.
\end{equation}
For the first term of \eqref{equation:gcg_conv_ss_usg1}, we have
$$
a_1 \prod_{t=1}^{T-1} (1 - \gamma_t) = a_1 \prod_{t=1}^{T-1} \frac{t + 2B - 2}{t + 2B} = a_1 \frac{(2B - 1)2B }{(T + 2B - 2)(T + 2B - 1)}. 
$$
For the second term of \eqref{equation:gcg_conv_ss_usg1}, since $\gamma_t = \frac{2}{t + 2B}$, it follows that
$$
\begin{aligned}
& \sum_{i=1}^{T-1} \left[ \prod_{t=i+1}^{T-1} (1 - \gamma_t) \right] \gamma_i^2 
=  \sum_{i=1}^{T-1} \left[ \prod_{t=i+1}^{T-1} \frac{t + 2B - 2}{t + 2B} \right] \frac{4}{(i + 2B)^2} \\
&=  \sum_{i=1}^{T-1} \frac{(i + 2B - 1)(i + 2B)}{(T + 2B - 2)(T + 2B - 1)}  \frac{4}{(i + 2B)^2} 
= \sum_{i=1}^{T-1} \frac{4(i + 2B - 1)}{(T + 2B - 2)(T + 2B - 1)(i + 2B)}   \\
&\leq \frac{4(T-1)}{(T + 2B - 2)(T + 2B - 1)}. 
\end{aligned}
$$
Substituting these into \eqref{equation:gcg_conv_ss_usg1} yields
$$
\begin{aligned}
a_T &\leq   \frac{\big(a_1(2B - 2)2B\big) + \big(2A(T-1)\big)}{(T + 2B - 2)(T + 2B - 1)} 
=2\frac{a_1(B - 1)2B + A(T-1)}{(T + 2B - 2)(T + 2B - 1)} 
\\
&\leq \frac{2 \max\{A, a_1(B - 1)\}(T + 2B - 1)}{(T + 2B - 2)(T + 2B - 1)} 
= \frac{2 \max\{A, a_1(B - 1)\}}{T + 2B - 2}.
\end{aligned}
$$
\paragraph{(ii).} 
Summing the inequality~\eqref{equation:ineq_gcg_conv_ss_usg} over $t = j, j+1, \ldots, T$, we obtain that
$a_{T+1} \leq a_j - \sum_{t=j}^{T} \gamma_t b_t + \frac{A}{2} \sum_{t=j}^{T} \gamma_t^2. $
Let $\widetildeb\triangleq\min_{t \in \{j, j+1, \ldots, T\}} b_t$, where $j\geq 1$. Then, using the result of part (i),
$$
\begin{aligned}
\widetildeb\bigg( \sum_{t=j}^{T} \gamma_t \bigg)  
&\leq a_j + \frac{A}{2} \sum_{t=j}^{T} \gamma_t^2 
\leq \frac{2 \max\{A, a_1(B-1)\}}{j+2B-2 } + 2A \sum_{t=j}^{T} \frac{1}{(t+2B)^2}\\
&\leq \frac{2 \max\{A, a_1(B-1)\}}{j+2B-2 } + 2A \sum_{t=j}^{T} \left(\frac{1}{t+2B-1} - \frac{1}{t+2B}\right)\\
&\leq \frac{2 \max\{A, a_1(B-1)\}}{j+2B-2 } + 2A  \left(\frac{1}{j+2B-1} - \frac{1}{T+2B}\right)
\leq \frac{4 \max\{A, a_1(B-1)\}}{j+2B-2 }.
\end{aligned}
$$
Since  $\sum_{t=j}^{T} \gamma_t =  \sum_{t=j}^{T} \frac{2}{t+2B} \geq 2 \frac{T-j+1}{T+2B}$,
which, combined with the above inequality and letting $j \triangleq \lfloor T/2 \rfloor + 2$,  concludes that for any $T \geq 3$:
$$ 
\begin{aligned}
\widetildeb=\min_{t \in \{\lfloor T/2 \rfloor + 2, \ldots, T\}} b_t 
&\leq \frac{2 \max\{A, a_1(B-1)\}(T+2B)}{(\lfloor T/2 \rfloor + 2B)(T - \lfloor T/2 \rfloor - 1)}
\leq \frac{2 \max\{A, a_1(B-1)\}(T+2B)}{( T/2  + 2B-0.5)(T - \lfloor T/2 \rfloor - 1)}\\
&=\frac{4 \max\{A, a_1(B-1)\}(T+2B)}{( T  + 4B-1)(T - \lfloor T/2 \rfloor - 1)}
\leq \frac{4 \max\{A, a_1(B-1)\}}{T/2 - 1},
\end{aligned} 
$$
which yields the desired result.
\end{proof}

Similar to the proximal gradient method for solving a composite problem (P3) with $f$ being convex and smooth (Theorem~\ref{theorem:prox_conv_ss_cvx}), using the above lemma, the GCG can also achieves a rate of $\mathcalO(1/T)$.

\begin{theoremHigh}[GCG for Convex and SS $f$: $\mathcalO(1/T)$]\label{theorem:gcg_conv_ss}
Suppose that $f$ is $\beta$-smooth and $g$ is \textcolor{black}{closed} convex, satisfying assumption (B1) and (B2) in Assumption~\ref{assumption:gcg}, and suppose additionally that $ f $ is convex.
Let $ \{\bx^\toptzero\}_{t > 0} $ be the sequence generated by the generalized conditional gradient method (Algorithm~\ref{alg:gen_cond_gen}) for solving problem (P3) with either a predefined stepsize $ \eta_t \triangleq \gamma_t\triangleq \frac{2}{t+2} $ or an adaptive stepsize $\eta_t\triangleq \theta_t \triangleq \min\left\{ 1, \frac{\mathcalW(\bx^\toptzero)}{\beta \normtwo{\widetildebx^\toptzero -\bx^\toptzero}^2} \right\}$  for each $t$-th iteration. Let $ D $ be an upper bound on the diameter of $ \dom(g) $: $ D \geq \max_{\bx,\by \in \dom(g)} \normtwo{\bx - \by}.$
Then,
\begin{enumerate}[(i)]
\item For  any $T\geq 1 $,  $F(\bx^{(T)}) - F(\bx^*) \leq \frac{2\beta D^2}{T}$.
\item For any $T \geq 3$, $\min_{t=\lfloor T/2 \rfloor + 2, \ldots, T} \mathcalW(\bx^\toptzero) \leq \frac{8\beta D^2}{T-2}$.
\end{enumerate}
\end{theoremHigh}
\begin{proof}[of Theorem~\ref{theorem:gcg_conv_ss}]
By the descent lemma for GCG (Lemma~\ref{lemma:desc_gcg}),
$$
F\big(\bx^\toptzero + \eta_t (\widetildebx^\toptzero - \bx^\toptzero)\big) - F(\bx^*) \leq F(\bx^\toptzero) - F(\bx^*) - \eta_t \mathcalW(\bx^\toptzero) + \frac{\eta_t^2 \beta}{2} \normtwo{\widetildebx^\toptzero - \bx^\toptzero}^2,
$$
where $ \widetildebx^\toptzero =\mathop{\argmin}_{\bx\in {\real^n}} \innerproduct{\nabla f(\bx^{(t)}), \bx} + g(\bx) $. Specifically, if a predefined stepsize $ \gamma_t \triangleq \frac{2}{t+2} $ is used, then
$$
F\big(\bx^\toptzero + \gamma_t (\widetildebx^\toptzero - \bx^\toptzero)\big) - F(\bx^*) \leq F(\bx^\toptzero) - F(\bx^*) - \gamma_t \mathcalW(\bx^\toptzero) + \frac{\gamma_t^2 \beta}{2} \normtwo{\widetildebx^\toptzero - \bx^\toptzero}^2.
$$
In the adaptive stepsize strategy, $  \theta_t \triangleq \min \left\{ 1, \frac{\mathcalW(\bx^\toptzero)}{\beta \normtwo{\widetildebx^\toptzero - \bx^\toptzero}^2} \right\} $, where $ \theta_t $ satisfies
$
\theta_t = \mathop{\argmin}_{\theta \in [0,1]} \left\{ -\theta \mathcalW(\bx^\toptzero) + \frac{\theta^2 \beta}{2} \normtwo{\widetildebx^\toptzero - \bx^\toptzero}^2 \right\}. 
$
This implies
\begin{align*}
F\big(\bx^\toptzero + \theta_t (\widetildebx^\toptzero - \bx^\toptzero)\big) - F(\bx^*) &\leq F(\bx^\toptzero) - F(\bx^*) - \theta_t \mathcalW(\bx^\toptzero) + \frac{\theta_t^2 \beta}{2} \normtwo{\widetildebx^\toptzero - \bx^\toptzero}^2 \\
&\leq F(\bx^\toptzero) - F(\bx^*) - \gamma_t \mathcalW(\bx^\toptzero) + \frac{\gamma_t^2 \beta}{2} \normtwo{\widetildebx^\toptzero - \bx^\toptzero}^2,
\end{align*}
Thus, both stepsize strategies satisfy
$$
F(\bx^\toptone) - F(\bx^*) \leq F(\bx^\toptzero) - F(\bx^*) - \gamma_t \mathcalW(\bx^\toptzero) + \frac{\gamma_t^2 \beta D^2}{2}.
$$
Invoking Lemma~\ref{lemma:gcg_conv_ss_usg} with $ a_t \triangleq F(\bx^\toptzero) - F(\bx^*), b_t \triangleq \mathcalW(\bx^\toptzero), A \triangleq \beta D^2 $, and $ B \triangleq 1 $, and noting that $ a_t \leq b_t $ by Lemma~\ref{lemma:con_gcg_cgnorm}, both parts (i) and (ii) follow.
\end{proof}

\section{Mirror Descent Method}\label{section:mirror}
\begin{algorithm}[h] 
\caption{Mirror Descent Method}
\label{alg:mirror_des}
\begin{algorithmic}[1] 
\Require A  function $f(\bx)$, a convex function $\phi(\bx)$, and a set $\sS$; 
\State {\bfseries Input:}  Initialize $\bx^{(1)}$;
\For{$t=1,2,\ldots$}
\State Pick a stepsize $\eta_t$;
\State $\bx^{(t+1)} \leftarrow \mathop{\argmin}_{\bx\in\sS} 
\left\{  \innerproduct{\eta_t f'(\bx^\toptzero)-\nabla \phi(\bx^\toptzero), \bx}    + \phi(\bx) \right\}$; \Comment{$f'(\bx^\toptzero)\in \partial f(\bx^\toptzero)$}
\EndFor
\State (Output Option 1) Output  $\bx_{\text{final}}\leftarrow \bx^{(T)}$;
\State (Output Option 2) Output  $\bx_{\text{best}}\leftarrow \argmin_{t\in\{1,2,\ldots,T\}} f(\bx^{(t)})$;
\end{algorithmic} 
\end{algorithm}

We showed in \eqref{equation:pgd_decom_raw} that the PGD method is equivalent to a minimization of the sum of the linearization of the smooth part around the current iterate plus a quadratic term, over which we derived the proximal gradient descent for the unconstrained composite optimization problem (P3); see \eqref{equation:prox_decom}. We reiterate the equivalence  as follows:
\begin{center}
\framebox{
\begin{minipage}{0.95\textwidth}
\begin{equation}\label{equation:pgd_decom_mirr}
\small
\textbf{(PGD)}:\quad 
\begin{aligned}
\bx^\toptone 
&\leftarrow\mathop{\argmin}_{\bx\in\sS} \normtwo{\bx - \left(\bx^{(t)} - \eta_t  f'(\bx^{(t)})\right)}^2\\
&=\mathop{\argmin}_{\bx\in\sS} 
\left\{
\frac{1}{2\eta_t} \normtwo{\bx-\bx^\toptzero}^2 + f(\bx^\toptzero)+ \innerproduct{f'(\bx^\toptzero), \bx-\bx^\toptzero} 
\right\}.
\end{aligned}
\end{equation}
\end{minipage}
}
\end{center}
This expression represents the sum of a linearization term and an Euclidean distance term between $\bx$ and the current iterate $\bx^\toptzero$. Here, $f'(\bx^\toptzero)\in \partial f(\bx^\toptzero)$ denotes any subgradient of $f$ at $\bx^\toptzero$; when $f$ is differentiable, this reduces to the gradient.
The \textit{mirror descent method} extends this idea from using Euclidean distance to employing a non-Euclidean distance defined by the Bregman distance (Definition~\ref{definition:breg_dist}):
\begin{tcolorbox}[colback=white,colframe=black]
\begin{minipage}{1\textwidth}
\begin{subequations}\label{equation:mirr_rule}
\small
\begin{align}
\bx^\toptone 
&\leftarrow\mathop{\argmin}_{\bx\in\sS} 
\left\{
\frac{1}{\eta_t} \mathcalD_{\phi}(\bx,\bx^\toptzero) + \cancel{f(\bx^\toptzero)}+ \innerproduct{f'(\bx^\toptzero), \bx-\cancel{\bx^\toptzero}} 
\right\}             \\
\textbf{(Mirror)}:\qquad\qquad  &=\mathop{\argmin}_{\bx\in\sS} 
\left\{  \innerproduct{\eta_t f'(\bx^\toptzero)-\nabla \phi(\bx^\toptzero), \bx}    + \phi(\bx) \right\}   \\
&=\mathop{\argmin}_{\bx\in\textcolor{mylightbluetext}{\real^n}} 
\left\{  \innerproduct{\eta_t f'(\bx^\toptzero)-\nabla \phi(\bx^\toptzero), \bx}    + \widetilde{\phi}(\bx) \right\}   \\
&=\mathop{\argmin}_{\bx\in\textcolor{mylightbluetext}{\real^n}} 
\left\{  \widetildef(\bx)    + \mathcalD_{\phi}(\bx,\bx^\toptzero) \right\} \triangleq \bproxtphi(\bx), \label{equation:mirr_rule_v4}
\end{align}
\end{subequations}
\end{minipage}
\end{tcolorbox}
\noindent
where we replace $\frac{1}{2}\normtwo{\bx-\bx^\toptzero}^2$ term in PGD with a Bregman distance $\big( \mathcalD_{\phi}(\bx,\bx^\toptzero) \triangleq \phi(\bx) - \phi(\bx^\toptzero) - \innerproduct{\nabla \phi(\bx^\toptzero), \bx - \bx^\toptzero }\big)$ and remove constant terms. The function $\phi(\bx)$ needs to be a convex function by Definition~\ref{definition:breg_dist}, and:
$$
\begin{aligned}
\widetildef(\bx) &\triangleq \innerproduct{\eta_t f'(\bx^\toptzero), \bx} + \indicatorS(\bx)
\quad &\implies&\quad \text{closed convex if $\sS$ is closed convex};\\
\widetilde{\phi}(\bx) &\triangleq \phi(\bx)+\indicatorS(\bx) 
\quad &\implies&\quad \text{$\widetilde{\phi}(\bx) = \infty$ if $\bx\notin\sS$},
\end{aligned}
$$
where $\indicatorS(\bx)$ is closed convex if and only if $\sS$ is closed convex (Exercise~\ref{exercise_closed_indica} and Exercise~\ref{exercise_convex_indica}).
The operator $\bproxtphi(\cdot)$ in \eqref{equation:mirr_rule_v4} denotes the Bregman-proximal operator (Definition~\ref{definition:projec_prox_opt}), which is an extension of the proximal operator.
When $\phi(\bx) = \frac{1}{2}\normtwo{\bx}^2$, $\bproxtphi(\bx) = \prox_{\widetildef}(\bx)$ for any $\bx\in\real^n$; then the mirror descent is equivalent  to PGD (Section~\ref{section:pgd}).

\begin{assumption}[Mirror Descent]\label{assumption:mirror}
For the analysis of the mirror descent method, we assume that 
\begin{itemize}
	\item  $f: \sS\rightarrow \real$ is a  proper closed and  \textbf{convex} function, where $\sS\subseteq\real^n$ is a closed and convex set such that $\sS\subseteq \interior(\dom(f))$~\footnote{This ensures the nonemptness of the subdifferential by Theorem~\ref{theorem:nonemp_relint_conv}.}.
	\item  $\phi$ is a differentiable, proper closed and \textbf{convex} function such that $\sS\subseteq \dom(\phi)$ and $\phi+\indicatorS$ is $\alpha$-strongly convex ($\alpha>0$).  
\end{itemize}
\noindent Under this assumption, the update of the mirror descent method in \eqref{equation:mirr_rule_v4} has a \textbf{unique} minimizer in  $\sS \cap \dom(\partial \phi)$ by invoking the Bregman-Proximal Property-O (Lemma~\ref{lemma:breg_prox_propo}) with $f(\bx) \triangleq \widetildef(\bx)$.
\end{assumption}

We will thus discuss the convergence of the mirror descent method under this assumption.
Using the Bregman-Proximal Property-I (Lemma~\ref{lemma:breg_prox_prop1}, under the convexity and closeness of $\widetildef(\bx)$) and the three-point property of the Bregman distance (Remark~\ref{remark:bregnan_dist}, under the convexity of $\phi(\bx)$), we can now establish an iterate inequality satisfied by the sequence generated by the mirror descent method. The inequality can be seen as a generalization of Lemma~\ref{lemma:iterate_ineq_pgd} for PGD approaches.
\begin{lemma}[Iterate Inequality for Mirror Descent]\label{lemma:itera_ineq_mirr}
Assume that Assumption~\ref{assumption:mirror} holds.
Let $\{\bx^\toptzero\}_{t > 0}$ is the sequence generated by the mirror descent method (Algorithm~\ref{alg:mirror_des}) with positive stepsizes $\{\eta_t\}_{t > 0}$ (Algorithm~\ref{alg:mirror_des}). Then, for   any optimizer $\bx^*$ of $f$ and any $t > 0$,
\begin{equation}\label{equation:itera_ineq_mirre1}
\mathcalD_\phi(\bx^*, \bx^\toptone)  \leq \mathcalD_\phi(\bx^*, \bx^\toptzero) - \eta_t \big(f(\bx^\toptzero) - f(\bx^*)\big) + \frac{\eta_t^2}{2\alpha} \normtwo{f'(\bx^\toptzero)}^2,~
\footnote{Note that when $\phi(\bx)=\frac{1}{2}\normtwo{\bx}^2$, $\mathcalD_{\phi}(\bx,\by)=\frac{1}{2}\normtwo{\bx-\by}^2$ by Example~\ref{example:breg_examp} and $\phi$ is $1$-SC; the inequality reduces to \eqref{equation:iterate_ineq_pgd1} for PGD.}
\end{equation}
where $f^\prime(\bx^\toptzero)\in \partial f(\bx^\toptzero)$ denotes any subgradient. 
For any nonnegative integer $T>0$, performing telescopic cancellations using \eqref{equation:itera_ineq_mirre1} shows that 
\begin{equation}\label{equation:itera_ineq_mirre2}
\sum_{t=1}^{T} \eta_t \big(f(\bx^\toptzero) - f(\bx^*)\big) \leq \mathcalD_\phi(\bx^*, \bx^{(1)})  + \frac{1}{2\alpha} \sum_{t=1}^{T} \eta_t^2 \normtwo{f^\prime(\bx^\toptzero)}^2. 
\end{equation}
\end{lemma}
\begin{proof}[of Lemma~\ref{lemma:itera_ineq_mirr}]
By the update formula of the mirror descent method in \eqref{equation:mirr_rule} for $\bx^\toptone$ and the Bregman-Proximal Property-I (Lemma~\ref{lemma:breg_prox_prop1}) invoked with $\by \triangleq \bx^\toptzero$ and the underlying function $\widetildef(\bx) \triangleq \eta_t \innerproduct{f'(\bx^\toptzero), \bx} + \indicatorS(\bx)$ (and hence $\widebarby \triangleq \bx^\toptone$), we have for any $\bz \in \sS$,
$$
\innerproduct{\nabla \phi(\bx^\toptzero) - \nabla \phi(\bx^\toptone), \bz - \bx^\toptone} 
\leq \eta_t \innerproduct{f'(\bx^\toptzero), \bz - \bx^\toptone}.
$$
By the three-point property (Remark~\ref{remark:bregnan_dist}, with $\bx \triangleq \bz$, $\by \triangleq \bx^\toptone$, and $\bz \triangleq \bx^\toptzero$),
$$
\mathcalD_\phi(\bz, \bx^\toptone) + \mathcalD_\phi(\bx^\toptone, \bx^\toptzero) - \mathcalD_\phi(\bz, \bx^\toptzero)
=
\innerproduct{\nabla \phi(\bx^\toptzero) - \nabla \phi(\bx^\toptone), \bz - \bx^\toptone}.
$$
Combining the preceding two quantities, using the subgradient inequality (Definition~\ref{definition:subgrad}), and letting $\bz\triangleq\bx^*$,  we have 
$$
\small
\begin{aligned}
&\eta_t (f(\bx^\toptzero) - f(\bx^*))\stackrel{*}{\leq} \eta_t \innerproduct{f'(\bx^\toptzero), \bx^\toptzero - \bx^*} \\
&\quad\leq \mathcalD_\phi(\bx^*, \bx^\toptzero) - \mathcalD_\phi(\bx^*, \bx^\toptone) - \mathcalD_\phi(\bx^\toptone, \bx^\toptzero) + \eta_t \innerproduct{ f'(\bx^\toptzero), \bx^\toptzero - \bx^\toptone }\\
&\quad \stackrel{\dag}{\leq} \mathcalD_\phi(\bx^*, \bx^\toptzero) - \mathcalD_\phi(\bx^*, \bx^\toptone) - \frac{\alpha}{2} \normtwobig{\bx^\toptone - \bx^\toptzero}^2 + \innerproduct{\frac{\eta_t}{\sqrt{\alpha}} f'(\bx^\toptzero), \sqrt{\alpha} (\bx^\toptzero - \bx^\toptone)} \\
&\quad \stackrel{\ddag}{\leq} \mathcalD_\phi(\bx^*, \bx^\toptzero) - \mathcalD_\phi(\bx^*, \bx^\toptone) - \frac{\alpha}{2} \normtwobig{\bx^\toptone - \bx^\toptzero}^2 + \frac{\eta_t^2}{2\alpha} \normtwobig{f'(\bx^\toptzero)}^2 + \frac{\alpha}{2} \normtwobig{\bx^\toptone - \bx^\toptzero}^2 \\
&\quad = \mathcalD_\phi(\bx^*, \bx^\toptzero) - \mathcalD_\phi(\bx^*, \bx^\toptone) + \frac{\eta_t^2}{2\alpha} \normtwobig{f'(\bx^\toptzero)}^2,
\end{aligned}
$$
where the inequality ($*$) follows from the subgradient inequality of $f$, the inequality $(\dag)$ follows from the $\alpha$-SC of $\phi+\indicatorS$, and the inequality $(\ddag)$ follows from the fundamental theorem of optimization (Theorem~\ref{theorem:funda_opt}). 
This completes the proof.
\end{proof}

Under a boundedness assumption on $\mathcalD_\phi(\bx, \bx^{(1)})$ over $\sS$, we can deduce a useful bound on the sequence of best achieved function values defined by
\begin{equation}
\fbest^\toptzero \triangleq \min\left\{f(\bx^{(1)}), f(\bx^{(2)}), \ldots, f(\bx^{(t)})\right\}.
\end{equation}

Similar to the convergence of PGD for convex functions in Theorem~\ref{theorem:pgd_lipschitz_dyna}, we then obtain the convergence results for mirror descent methods.
\begin{theoremHigh}[Mirror Descent for Convex and Lipschitz, $\mathcalO(\ln (T) / \sqrt{T})$]\label{theorem:mirr_best_bound}
Considering the same setting as Lemma~\ref{lemma:itera_ineq_mirr},
assume that Assumption~\ref{assumption:mirror} holds.
Let $\{\bx^\toptzero\}_{t > 0}$ be the sequence generated by the mirror descent method (Algorithm~\ref{alg:mirror_des}) with positive stepsizes $\{\eta_t\}_{t > 0}$ and $\bx^*$ be any optimal point. 
Suppose that $\normtwo{f'(\bx)} \leq L$ for all $\bx \in \sS$, where $L > 0$. 
Then, for any optimal point $ \bx^*$,
\begin{enumerate}[(i)]
\item  If $\mathcalD_\phi(\bx, \bx^{(1)})$ is bounded over $\sS$ such that 
$
\Omega(\bx^{(1)}) \geq \max_{\bx \in \sS} \mathcalD_\phi(\bx, \bx^{(1)})
$, then for any $T > 0$,
$$
\fbest^{(T)} - f(\bx^*) \leq \frac{\Omega(\bx^{(1)}) + \frac{L^2}{2\alpha} \sum_{t=1}^T \eta_t^2}{\sum_{t=1}^T \eta_t}.
$$

\item If $ \frac{\sum_{t=1}^{T} \eta_t^2}{\sum_{t=1}^{T} \eta_t} \to 0 $ as $ T \to \infty $, then $ \fbest^{(T)} \to f(\bx^*) $ as $ T \to \infty $.

\item If $
\eta_t = \frac{\sqrt{2\alpha}}{L \sqrt{t+1}}
$ or 
$\eta_t = \scriptsize\begin{cases} 
	\sqrt{2\alpha}/{(\normtwo{f'(\bx^\toptzero)} \sqrt{t+1})}, & f'(\bx^\toptzero) \neq \bzero; \\
	{\sqrt{2\alpha}}/{(L \sqrt{t+1})}, & f'(\bx^\toptzero) = \bzero,
\end{cases}$
then for all $ T\geq 5 $,
$$
\fbest^{(T)} - f(\bx^*) \leq \frac{L}{\sqrt{2\alpha}} \frac{\mathcalD_{\phi}(\bx^*, \bx^{(1)}) + \ln(T+1)}{\sqrt{T+1}}.
$$
\end{enumerate}
\end{theoremHigh}
\begin{proof}[of Theorem~\ref{theorem:mirr_best_bound}]
\textbf{(i).}
By Lemma~\ref{lemma:itera_ineq_mirr},  it follows that for any $T > 0$,
\begin{equation}
\begin{aligned}
\sum_{t=1}^{T} \eta_t \big(\fbest^{(T)} - f(\bx^*)\big) 
&\leq 	
\sum_{t=1}^{T} \eta_t \big(f(\bx^\toptzero) - f(\bx^*)\big) \\
&\leq \mathcalD_\phi(\bx^*, \bx^{(1)})  + \frac{1}{2\alpha} \sum_{t=1}^{T} \eta_t^2 \normtwo{f^\prime(\bx^\toptzero)}^2
\leq \Omega(\bx^{(1)}) + \frac{L^2}{2\alpha} \sum_{t=1}^T \eta_t^2,
\end{aligned}
\end{equation}
which obtains the desired result in (i).
\paragraph{(ii).}
By (i), we have 
$$
\fbest^{(T)} - f(\bx^*)
\leq 
\frac{\mathcalD_\phi(\bx^*, \bx^{(1)})  + \frac{L}{2\alpha} \sum_{t=1}^{T} \eta_t^2}{\sum_{t=1}^{T} \eta_t}.
$$
The claim in (ii) follows immediately.
\paragraph{(iii).}
Note that in either case, it follows that  $ \eta_t^2 \normtwo{f'(\bx^\toptzero)}^2 \leq \frac{2\alpha}{t+1} $ and $ \eta_t \geq \frac{\sqrt{2\alpha}}{L \sqrt{t+1}} $, whence we have
$$
\fbest^{(T)} - f(\bx^*) \leq \frac{L}{\sqrt{2\alpha}} \frac{\mathcalD_{\phi}(\bx^*, \bx^{(1)}) + \sum_{t=1}^{T} \frac{1}{{t+1}}}{\sum_{t=1}^{T} \frac{1}{\sqrt{t+1}}}.
$$
Using Lemma~\ref{lemma:pgd_lem2} concludes the result in (iii).
\end{proof}

\begin{algorithm}[h] 
\caption{Generalized Mirror Descent Method}
\label{alg:gen_mirror_des}
\begin{algorithmic}[1] 
\Require Functions $f(\bx), g(\bx)$, a convex function $\phi(\bx)$; 
\State {\bfseries Input:}  Initialize $\bx^{(1)}$;
\For{$t=1,2,\ldots$}
\State Pick a stepsize $\eta_t$;
\State $\bx^{(t+1)} \leftarrow \mathop{\argmin}_{\bx\in\real^n} 
\left\{  \innerproduct{\eta_t f'(\bx^\toptzero)-\nabla \phi(\bx^\toptzero), \bx}+ \textcolor{mylightbluetext}{\eta_t g(\bx)}  + \phi(\bx) \right\}$;
\EndFor
\State (Output Option 1) Output  $\bx_{\text{final}}\leftarrow \bx^{(T)}$;
\State (Output Option 2) Output  $\bx_{\text{best}}\leftarrow \argmin_{t\in\{1,2,\ldots,T\}} f(\bx^{(t)})$;
\end{algorithmic} 
\end{algorithm}
\section{Generalized Mirror Descent Method}\label{section:gmirror}

We have shown that the mirror descent method addresses the same constrained problem \eqref{equation:p2_pgd} as the PGD approach. 
The \textit{generalized mirror (G-mirror) descent method (a.k.a., Bregman proximal gradient method)} (Algorithm~\ref{alg:gen_mirror_des}) tackles the unconstrained composite problem:
\begin{equation}\label{equation:gene_mirror_prob}
\text{(P3)}:\qquad \min \{F(\bx) \triangleq f(\bx)+g(\bx)\}.
\end{equation}

Similar to proximal gradient and generalized conditional gradient methods, by incorporating an  objective function  $\eta_t g(\bx)$ into the update rule  of the mirror descent, we obtain:
\begin{center}
\framebox{
\begin{minipage}{0.95\textwidth}
\begin{subequations}\label{equation:gmirr_rule}
\small
\begin{align}
\textbf{(G-Mirror)}:\qquad\bx^\toptone 
&\leftarrow\mathop{\argmin}_{\bx\in\textcolor{mylightbluetext}{\real^n}} 
\left\{  \innerproduct{\eta_t f'(\bx^\toptzero)-\nabla \phi(\bx^\toptzero), \bx} + \textcolor{mylightbluetext}{\eta_t g(\bx)}   + \phi(\bx) \right\}   \\
&=\mathop{\argmin}_{\bx\in\textcolor{black}{\real^n}} 
\left\{  \widehatf(\bx)    + \mathcalD_{\phi}(\bx,\bx^\toptzero) \right\} \triangleq \bproxhphi(\bx), \label{equation:gmirr_rule_v4}
\end{align}
\end{subequations}
\end{minipage}
}
\end{center}
where we remove the constraint over $\sS$ and add a $\eta_t g(\bx)$ term to the mirror descent update in \eqref{equation:mirr_rule}. The function $\phi(\bx)$ must be a convex function by Definition~\ref{definition:breg_dist}, and:
$$
\begin{aligned}
\widehatf(\bx) &\triangleq \innerproduct{\eta_t f'(\bx^\toptzero), \bx} + \eta_t g(\bx)
\quad &\implies&\quad \text{closed convex if $g(\bx)$ is closed convex}.\\
\end{aligned}
$$
Again, the operator $\bproxhphi(\cdot)$ in \eqref{equation:gmirr_rule_v4} denotes the Bregman-proximal operator (Definition~\ref{definition:projec_prox_opt}), which is an extension of the proximal operator.
When $\phi(\bx) = \frac{1}{2}\normtwo{\bx}^2$, $\bproxhphi(\bx) = \prox_{\widehatf}(\bx)$ for any $\bx\in\real^n$; then the G-mirror descent method is equivalent  to the proximal gradient method (Section~\ref{section:prox_gd}), and the update rule is: 
$$
\bx^\toptone \leftarrow \prox_{\eta_t g}(\bx^\toptzero - \eta_t f'(\bx^\toptzero)).
$$

\begin{assumption}[G-Mirror Descent]\label{assumption:gen_mirror}
For the analysis of the G-mirror descent method, we assume that 
\begin{itemize}
\item  $f, g: \real^n\rightarrow \real$ are proper closed and \textbf{convex} functions, where $\dom(g)\subseteq\interior(\dom(f))$~\footnote{This ensures the nonemptness of the subdifferential by Theorem~\ref{theorem:nonemp_relint_conv}.}. And $\normtwo{f'(\bx)}\leq L$ for any $\bx\in\dom(g)$, where $L>0$ and $f'(\bx)\in\partial f(\bx)$ denotes a subgradient.
\item  $\phi$ is a differentiable, proper closed and \textbf{convex} function such that $\dom(g)\subseteq \dom(\phi)$ and $\phi+\indicatorG_{\dom(g)}$ is $\alpha$-strongly convex ($\alpha>0$).  
\end{itemize}
\noindent Under this assumption, the update of the G-mirror descent method in Algorithm~\ref{alg:gen_mirror_des} has a \textbf{unique} minimizer in  $\dom(g) \cap \dom(\partial \phi)$ by invoking the Bregman-Proximal Property-O (Lemma~\ref{lemma:breg_prox_propo}) with $f(\bx) \triangleq\widehatf(\bx)$.
\end{assumption}

The analysis of the G-mirror method is based on arguments similar to those used in Section~\ref{section:mirror} for analyzing the mirror descent method. 
We start  by proving a technical lemma establishing an inequality similar to the one derived in Lemma~\ref{lemma:itera_ineq_mirr}. 
Note that in addition to our basic assumptions, we assume that $ g $ is a nonnegative function and that the stepsizes are nonincreasing.
We again define best achieved function value at each iteration as: 
\begin{equation}
\Fbest^\toptzero \triangleq  \min\left\{F(\bx^{(1)}), F(\bx^{(2)}), \ldots, F(\bx^{(t)})\right\}.
\end{equation}
\begin{lemma}[Iterate Inequality for G-Mirror Descent]\label{lemma:itera_ineq_gmirr}
Assume that Assumption~\ref{assumption:gen_mirror} holds, and suppose  that $ g $ is a nonnegative function. Let $ \{\bx^\toptzero\}_{t > 0} $ be the sequence generated by the G-mirror method (Algorithm~\ref{alg:gen_mirror_des}) with positive nonincreasing stepsizes $ \{\eta_t\}_{t > 0} $. Then, for any optimal point $ \bx^*$ and $ T > 0 $,
\begin{equation}\label{equation:itera_ineq_gmirre1}
\small
\begin{aligned}
\mathcalD_\phi(\bx^*, \bx^\toptone)  \leq \mathcalD_\phi(\bx^*, \bx^\toptzero) 
- \eta_t \big(f(\bx^\toptzero) +  g(\bx^\toptone)- F(\bx^*)\big) + \frac{\eta_t^2}{2\alpha} \normtwo{f'(\bx^\toptzero)}^2,
\end{aligned}
\end{equation}
where $f^\prime(\bx^\toptzero)\in \partial f(\bx^\toptzero)$ denotes any subgradient. 
For any nonnegative integer $T>0$, performing telescopic cancellations using \eqref{equation:itera_ineq_gmirre1} shows that 
\begin{equation}\label{equation:itera_ineq_gmirre2}
\Fbest^{(T)} - F(\bx^*) \leq \frac{\eta_1 g(\bx^{(1)}) + \mathcalD_\phi(\bx^*, \bx^{(1)}) + \frac{1}{2\alpha} \sum_{t=1}^T \eta_t^2 \normtwo{f'(\bx^\toptzero)}^2}{\sum_{t=1}^T \eta_t}.
\end{equation}
\end{lemma}
\begin{proof}[of Lemma~\ref{lemma:itera_ineq_gmirr}]
By the update formula of the G-mirror descent method in \eqref{equation:gmirr_rule} for $\bx^\toptone$ and the Bregman-Proximal Property-I (Lemma~\ref{lemma:breg_prox_prop1}) invoked with $\by \triangleq \bx^\toptzero$ and the underlying function $\widehatf(\bx) \triangleq \innerproduct{\eta_t f'(\bx^\toptzero), \bx} + \eta_t g(\bx)$ (and hence $\widebarby \triangleq \bx^\toptone$), we have for any $\bz \in \real^n$,
$$
\innerproduct{\nabla \phi(\bx^\toptzero) - \nabla \phi(\bx^\toptone), \bz - \bx^\toptone} 
\leq \eta_t \innerproduct{f'(\bx^\toptzero), \bz - \bx^\toptone} + \eta_t g(\bz) - \eta_t g(\bx^\toptone).
$$
By the three-point property (Remark~\ref{remark:bregnan_dist}, with $\bx \triangleq \bz$, $\by \triangleq \bx^\toptone$, and $\bz \triangleq \bx^\toptzero$),
$$
\mathcalD_\phi(\bz, \bx^\toptone) + \mathcalD_\phi(\bx^\toptone, \bx^\toptzero) - \mathcalD_\phi(\bz, \bx^\toptzero)
=\innerproduct{\nabla \phi(\bx^\toptzero) - \nabla \phi(\bx^\toptone), \bz - \bx^\toptone}.
$$
Combining the preceding two quantities, using the subgradient inequality (Definition~\ref{definition:subgrad}), and letting $\bz\triangleq\bx^*$,  we have 
$$
\begin{aligned}
 &\eta_t \left[ f(\bx^\toptzero) + g(\bx^\toptone) - F(\bx^*) \right] \stackrel{*}{\leq}\eta_t \innerproduct{f'(\bx^\toptzero), \bx^\toptzero - \bx^*} + \eta_t g(\bx^\toptone) - \eta_t g(\bx^*)\\
&\leq \mathcalD_\phi(\bx^*, \bx^\toptzero) - \mathcalD_\phi(\bx^*, \bx^\toptone) - \mathcalD_\phi(\bx^\toptone, \bx^\toptzero) + \eta_t \innerproduct{ f'(\bx^\toptzero), \bx^\toptzero - \bx^\toptone}\\
&\stackrel{\dag}{\leq} \mathcalD_\phi(\bx^*, \bx^\toptzero) - \mathcalD_\phi(\bx^*, \bx^\toptone) - \frac{\alpha}{2} \normtwo{\bx^\toptone - \bx^\toptzero}^2 + \innerproduct{ \frac{\eta_t}{\sqrt{\alpha}} f'(\bx^\toptzero), \sqrt{\alpha} (\bx^\toptzero - \bx^\toptone)}\\
&\stackrel{\ddag}{\leq} \mathcalD_\phi(\bx^*, \bx^\toptzero) - \mathcalD_\phi(\bx^*, \bx^\toptone) - \frac{\alpha}{2} \normtwo{\bx^\toptone - \bx^\toptzero}^2 + \frac{\eta_t^2}{2\alpha} \normtwo{f'(\bx^\toptzero)}^2 + \frac{\alpha}{2} \normtwo{\bx^\toptone - \bx^\toptzero}^2\\
&= \mathcalD_\phi(\bx^*, \bx^\toptzero) - \mathcalD_\phi(\bx^*, \bx^\toptone) + \frac{\eta_t^2}{2\alpha} \normtwo{f'(\bx^\toptzero)}^2.
\end{aligned}
$$
where the inequality ($*$) follows from the subgradient inequality of $f$, the inequality $(\dag)$ follows from the $\alpha$-SC of $\phi+\indicatorG_{\dom(g)}$, and the inequality $(\ddag)$ follows from the fundamental theorem of optimization (Theorem~\ref{theorem:funda_opt}). 
This establishes \eqref{equation:itera_ineq_gmirre1}.
Performing telescopic cancellations using \eqref{equation:itera_ineq_gmirre1} shows that 
$$
\sum_{t=1}^T \eta_t \left[ f(\bx^\toptzero) + g(\bx^\toptone) - F(\bx^*) \right] 
\leq \mathcalD_\phi(\bx^*, \bx^{(1)}) - \mathcalD_\phi(\bx^*, \bx^{(T+1)}) + \frac{1}{2\alpha} \sum_{t=1}^T \eta_t^2 \normtwo{f'(\bx^\toptzero)}^2.
$$
Since $\{\eta_t\}$ is nonincreasing and $g$ is nonnegative, adding the term $ \eta_1 g(\bx^{(1)}) - \eta_T g(\bx^{(T+1)}) $ to both sides and using the nonnegativity of the Bregman distance, we get
$$
\small
\begin{aligned}
\bigg( \sum_{t=1}^T \eta_t \bigg) \big( \Fbest^{(T)} - F(\bx^*) \big)
&\leq \sum_{t=1}^T \eta_t \left[ F(\bx^\toptzero) - F(\bx^*) \right] \\
&\leq\eta_1 (F(\bx^{(1)}) - F(\bx^*)) + \sum_{t=2}^T \left[ \eta_t f(\bx^\toptzero) + \eta_{t-1} g(\bx^\toptzero) - \eta_t F(\bx^*) \right]\\
&\leq \eta_1 g(\bx^{(1)}) - \eta_T g(\bx^{(T+1)}) + \mathcalD_\phi(\bx^*, \bx^{(1)}) + \frac{1}{2\alpha} \sum_{t=1}^T \eta_t^2 \normtwo{f'(\bx^\toptzero)}^2\\
&\leq\eta_1 g(\bx^{(1)}) + \mathcalD_\phi(\bx^*, \bx^{(1)}) + \frac{1}{2\alpha} \sum_{t=1}^T \eta_t^2 \normtwo{f'(\bx^\toptzero)}^2,
\end{aligned}
$$
which concludes the result.
\end{proof}

Using Lemma~\ref{lemma:itera_ineq_gmirr}, we can also establish the rate of convergence of the G-mirror method with a dynamic stepsize rule.
\begin{theoremHigh}[G-Mirror Descent for Convex: $\mathcalO(\ln T / \sqrt{T})$]\label{theorem:gmirror_dyn}
Assume that Assumption~\ref{assumption:gen_mirror} holds, and suppose  that $ g $ is a nonnegative function. 
Let $\{\bx^\toptzero\}_{t > 0}$ be the sequence generated by the G-mirror method (Algorithm~\ref{alg:gen_mirror_des}) with positive nonincreasing stepsizes $ \{\eta_t\}_{t > 0} $.
Suppose that $\normtwo{f'(\bx)} \leq L$ for all $\bx \in \real^n$, where $L > 0$. 
Then, for any optimal point $ \bx^*$,
\begin{enumerate}[(i)]
\item  Let
$
\Omega(\bx^{(1)}) \geq \max_{\bx \in \sS} \mathcalD_\phi(\bx, \bx^{(1)})
$. Then, for any $T > 0$,
$$
\Fbest^{(T)} - F(\bx^*) \leq \frac{\eta_1 g(\bx^{(1)}) + \Omega(\bx^{(1)}) + \frac{L^2}{2\alpha} \sum_{t=1}^T \eta_t^2}{\sum_{t=1}^T \eta_t}.
$$
Note that we can initialize $\bx^{(1)}$ such that $g(\bx^{(1)})=0$ to simplify the result.
\item If $ \frac{\sum_{t=1}^{T} \eta_t^2}{\sum_{t=1}^{T} \eta_t} \to 0 $ as $ T \to \infty $, then $ \Fbest^{(T)} \to F(\bx^*) $ as $ T \to \infty $.

\item If $
\eta_t = \frac{\sqrt{2\alpha}}{L \sqrt{t+1}}
$,~\footnote{Note that the second choice of the stepsize in (iii) of Theorem~\ref{theorem:mirr_best_bound} can not be applied since it might not be a nonincreasing sequence.}
then for all $ T\geq 5 $,
\begin{equation}\label{eq:dynamic_rate}
\Fbest^{(T)} - F(\bx^*) 
\leq 
\frac{L}{\sqrt{2\alpha}} \frac{\frac{\sqrt{\alpha}}{L} g(\bx^{(1)}) + \mathcalD_\phi(\bx^*, \bx^{(1)}) + \ln(T+1)}{\sqrt{T+1}}.
\end{equation}
\end{enumerate}
\end{theoremHigh}
\begin{proof}[of Theorem~\ref{theorem:gmirror_dyn}]
\textbf{(i, ii).}
The result in (i) follows immediately by Lemma~\ref{lemma:itera_ineq_gmirr}, which in turn concludes the claim in (ii).
\paragraph{(iii).}
By Lemma~\ref{lemma:itera_ineq_gmirr} and the choice of the stepsizes, we have 
$$
\Fbest^{(T)} - F(\bx^*) 
\leq 
\frac{ \frac{\sqrt{\alpha}}{L} g(\bx^{(1)}) + \mathcalD_\phi(\bx^*, \bx^{(1)}) + \frac{1}{2\alpha} \sum_{t=1}^T \eta_t^2 \normtwo{f'(\bx^\toptzero)}^2}{\sum_{t=1}^T \eta_t}.
$$
Since $\normtwo{f'(\bx)}\leq L$ by assumption, we have $\eta_t^2 \normtwo{f'(\bx^\toptzero)}^2 \leq \frac{2\alpha}{t+1}$ and $\eta_t = \frac{\sqrt{2\alpha}}{L \sqrt{t+1}}$. Plugging into the above inequality yields that
$$
\Fbest^{(T)} - F(\bx^*) 
\leq \frac{L}{\sqrt{2\alpha}} \frac{ \frac{\sqrt{\alpha}}{L} g(\bx^{(1)}) + \mathcalD_\phi(\bx^*, \bx^{(1)}) + \sum_{t=1}^T \frac{1}{t+1}}{\sum_{t=1}^T \frac{1}{\sqrt{t+1}}}.
$$
Using Lemma~\ref{lemma:pgd_lem2} concludes the result in (iii).
\end{proof}

\begin{problemset}
\item Prove the equivalence of the three definitions of conditional gradient norm in Definition~\ref{definition:cond_grad_norm}.

\item \label{prob:grad_map_lipschitz} \textbf{Lipschitzness of gradient mapping.} Let $ f $ be a proper \textcolor{black}{closed} and $\beta$-smooth function  and  $ g $ be a proper closed and convex function satisfying assumption (A1) and (A2) in Assumption~\ref{assumption:proximal_grad}. Let further $L\in \left(\frac{\beta}{2}, \infty\right)$. 
Show that $\mathcalG_L(\bx)$ is $(2L+\beta)$-Lipschitz: $\normtwo{\mathcalG_L(\bx) - \mathcalG_L(\by)} \leq (2L+\beta) \normtwo{\bx-\by}$ for any $\bx,\by\in\interior(\dom(f))$.

\item Prove \eqref{equation:gd_des_cor} rigorously.

\item Prove the relation in \eqref{equation:smooth_noneu2}.

\item Verify the equivalent update rules for PGD in \eqref{equation:pgd_decom_raw}.
\end{problemset}

\newpage
\chapter{Constrained Optimization with Specific Constraints}
\begingroup
\hypersetup{
linkcolor=structurecolor,
linktoc=page,  
}
\minitoc \newpage
\endgroup

\section{Constrained Optimization}
Sections~\ref{section:pgd} and \ref{section:cond_gd} introduce projected gradient descent and conditional gradient method for solving the constrained optimization problem (P2) (see Definition~\ref{definition:opt_probs_all}).
Unlike unconstrained optimization problems, in constrained optimization problems, the variable $\bx$ cannot assume any arbitrary value. 
This means that many algorithms designed for unconstrained optimization are not applicable. 
For instance, in gradient descent methods, moving along the negative gradient direction may lead to points that are not feasible, and the gradient of the objective function at the optimal solution is not necessarily a zero vector. 
Consequently, constrained optimization problems are significantly more complex than their unconstrained counterparts. 

In this chapter, we will discuss several penalty function methods that incorporate constraints into the objective function as penalty terms, transforming these problems into more familiar unconstrained optimization problems.

When the set $\sS$ in the problem (P2)  takes the following form, we deal with a constrained optimization problem involving both equality and inequality constraints:
\begin{subequations}
\begin{equation}\label{equation:pr_in_constchap}
\begin{aligned}
\textbf{(P4)}:\qquad  
&\mathop{\min}_{\bx} &f(\bx)&, \\
&\text{s.t.} & c_i(\bx)&=0, \,\,i\in\mathcalE,\,\,\text{with $\abs{\mathcalE}=p$},\\
&  & c_i(\bx)&\leq 0, \,\,i\in\mathcalI,\,\,\text{with $\abs{\mathcalI}=q$},
\end{aligned}
\end{equation}
where the variables $\bx \in \real^n$, $\mathcalE$ is the index set of equality constraints, $\mathcalI$ is the index set of inequality constraints,  and $c_i(\bx)$ are continuous functions.
The optimality conditions for (P4) are discussed in Sections~\ref{section:constrain_opt}, \ref{section:constr_convset}, \ref{section:convex_optimiz}, and \ref{section:gen_kkt_cond}.

In the previous chapter, various methods for solving unconstrained optimization problems were introduced. We want to transform problem (P4) into an unconstrained optimization problem. In the following sections, we also start by considering a simplified scenario where only equality constraints are present:
\begin{equation}\label{equation:pr_in_constchap_equ}
\begin{aligned}
\textbf{(P4.1)}:\qquad  	& \mathopmin{\bx} & f(\bx)&, \\
& \text{s.t.} & c_i(\bx)& = 0, \,\, i \in \mathcalE,\,\,\text{with $\abs{\mathcalE}=p$}.\\
\end{aligned}
\end{equation}
\end{subequations}
n certain special cases, directly solving the (nonlinear) system of equations $c_i(\bx) = 0$ can eliminate some variables, thereby converting the problem into an unconstrained one. However, for general functions $c_i(\bx)$, eliminating variables through this approach is not feasible, necessitating alternative methods to handle such problems.

\section{Penalty Function Method}\label{section:pen_func_me}
The \textit{penalty function method} aims to convert the constrained optimization problem (P4) into an unconstrained one. To ensure the quality of the approximate solution, the objective function of the new unconstrained problem includes the original objective function plus a penalty term related to the constraints. For points outside the feasible region, this penalty term is positive, effectively penalizing those points; for points within the feasible region, the penalty term is zero, meaning no penalty is applied. Consequently, the penalty term drives the solution of the unconstrained optimization problem towards the feasible region.

\subsection{Quadratic Penalty Function Method for Equality Constraints}\label{section:pen_equa}

For problems with equality constraints, various penalty terms can be selected, the simplest being quadratic. Here is the definition of the quadratic penalty function:
\begin{definition}[Quadratic Penalty Function for Equality Constraints]
For the equality-constrained optimization problem (P4.1) in \eqref{equation:pr_in_constchap_equ}, define the quadratic penalty function
\begin{equation}
f_{\sigma}(\bx) \triangleq f(\bx) + \frac{1}{2} \sigma \sum_{i \in \mathcalE} c_i^2(\bx)
\triangleq f(\bx) + \frac{1}{2} \sigma \normtwo{\bc(\bx)}^2,
\end{equation}
where $\bc(\bx)$ is the vector function $ \bc: \real^n \rightarrow \real^p $, with its $ i $-th component being the $ i $-th constraint function $ c_i $.
The term $\frac{1}{2} \sigma \normtwo{\bc(\bx)}^2$ is referred to as the \textit{penalty term}, and $\sigma>0$ is called the \textit{penalty factor}.
\end{definition}

Since this penalty function penalizes points that violate the constraints, during the iterative process, these points typically lie outside the feasible region. 
For these infeasible points, increasing $\sigma$ increases the weight of the penalty term, which forces its minimum point closer to the feasible region. 
On the other hand, within the feasible region, the global minimum point of $f_{\sigma}(\bx)$ aligns with the optimal solution of the equality-constrained optimization problem (P4.1).

\begin{figure}[h]
\centering       
\vspace{-0.25cm}                 
\subfigtopskip=2pt               
\subfigbottomskip=-2pt         
\subfigcapskip=-10pt      
\includegraphics[width=0.98\textwidth]{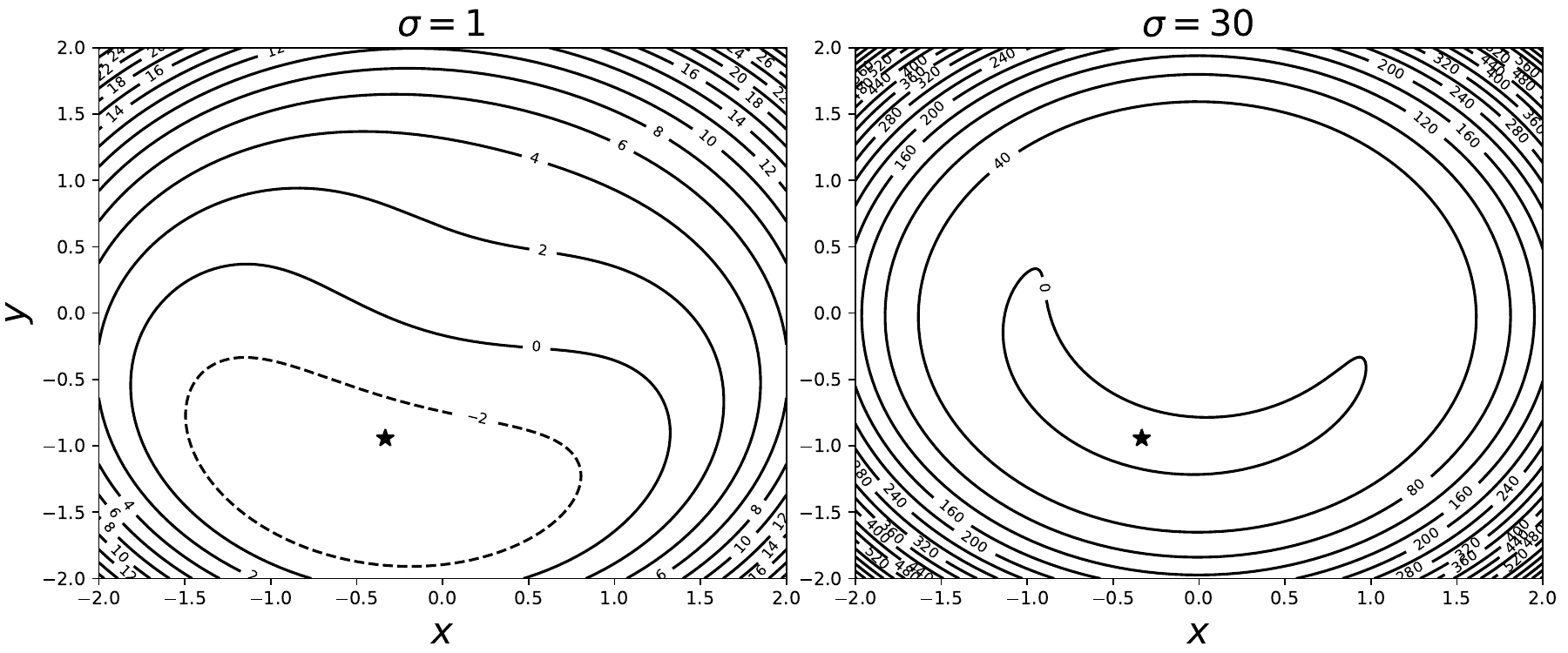}
\caption{The contour plot of $f_{\sigma}(x, y) = x + \sqrt{8} y + \frac{\sigma}{2} (x^2 + y^2 - 1)^2$ for $\sigma=1$ and $\sigma=20$.}
\label{fig:penalty_examp}
\end{figure}

\begin{example}[Penalty Function Method]\label{example:penal_func}
Consider the optimization problem
$$
\begin{aligned}
& \min & x + \sqrt{8} y\quad \text{s.t.} \quad  & x^2 + y^2 = 1.
\end{aligned}
$$
Two points satisfy the KKT conditions: $\big[-\frac{1}{3}, -\frac{\sqrt{8}}{3}\big]^\top$ and $\big[\frac{1}{3}, \frac{\sqrt{8}}{3}\big]^\top$.
It is easy to verify that the optimal solution is $\big[-\frac{1}{3}, -\frac{\sqrt{8}}{3}\big]^\top$. Consider the quadratic penalty function
$$
f_{\sigma}(x, y) = x + \sqrt{8} y + \frac{\sigma}{2} (x^2 + y^2 - 1)^2,
$$
where the contour lines of the penalty function for $\sigma = 1$ and $\sigma = 30$ are shown in Figure~\ref{fig:penalty_examp}. It can be seen that as $\sigma$ increases, the minimum value of the quadratic penalty function $f_{\sigma}(x, y)$ approaches the minimum value (indicated by asterisks in the figures) of the original problem. 
However, as $\sigma$ increases, the contour lines near the optimal point also become flatter, which complicates solving the unconstrained optimization problem.
\end{example}

From the preceding example, we understand that given a penalty factor $\sigma$, we can approximate the solution to the original problem by minimizing $f_{\sigma}(\bx)$. However, this approach has limitations in practice. The following example demonstrates that when $\sigma$ is too small, the penalty function might become unbounded.
\begin{example}[Failure of Penalty Function Method]
Consider the optimization problem
$$
\begin{aligned}
& \min & -2x^2 + 9y^2 \quad
\text{s.t.} \quad  x = 1.
\end{aligned}
$$
By eliminating the variable $x$, it's straightforward to determine that the optimal solution is $[1, 0]^\top$. However, considering the penalty function
$$
f_{\sigma}(x, y) = -2x^2 + 9y^2 + \frac{\sigma}{2}(x - 1)^2,
$$
for any $\sigma \leq 4$, the penalty function is unbounded.
\end{example}

This phenomenon occurs because if the penalty factor is too small, the decrease in the objective function value at infeasible points may counterbalance the penalty imposed for constraint violations. S the initial choice of $\sigma$ should not be too small.

\begin{algorithm}[h] 
\caption{Quadratic Penalty Function Method for (P4.1) in \eqref{equation:pr_in_constchap_equ}}
\label{alg:quad_pen_eq}
\begin{algorithmic}[1] 
\Require A function $f(\bx)$ and a set of constraints $\{c_i(\bx)\}$; 
\State {\bfseries Input:}  Choose the penalty factor growth coefficient $\rho > 1$;
\State {\bfseries Input:}  Initialize $\bx^{(1)}$ and $\sigma_1 > 0$;
\For{$t=1,2,\ldots$}
\State With $\bx^\toptzero$ as the initial point, solve $\bx^\toptone = \mathop{\argmin}_{\bx} f_{\sigma_t}(\bx)$; 
\State Choose $\sigma_{t+1} = \rho \sigma_t$;
\EndFor
\State {\bfseries Return:} final  $\bx^\toptzero$;
\end{algorithmic} 
\end{algorithm}

The execution process of Algorithm~\ref{alg:quad_pen_eq} is straightforward: first, select a series of exponentially increasing penalty factors $\sigma_t$, then for each penalty factor, solve the quadratic penalty function $f_{\sigma_t}(\bx)$ to find its minimum (or local minimum). 
Any of the (iterative) unconstrained optimization algorithms introduced in Chapter~\ref{chapter:gd_convg} or the Newton-type methods in Chapter~\ref{chapter:second_order} can be used to solve the subproblem in the Step 4 of the algorithm.
For different problems, the meaning of $\arg \min$ in the Step 4 of the algorithm can be different:
\begin{itemize}
\item $\bx^\toptone$ is the global minimum of the penalty function $f_{\sigma_t}(\bx)$;
\item $\bx^\toptone$ is a local minimum of the penalty function $f_{\sigma_t}(\bx)$;
\item $\bx^\toptone$ is not a strict minimum of the penalty function $f_{\sigma_t}(\bx)$, but approximately satisfies the first-order optimality condition $\nabla f_{\sigma_t}(\bx^\toptone) \approx \bzero$.
\end{itemize}
Based on the above description, three key points need attention in Algorithm~\ref{alg:quad_pen_eq} \citep{liu2020optimization}: 
\begin{itemize}
\item Careful selection of the parameter 
$\sigma_t$ is crucial. If $\sigma_t$ increases too quickly, solving the subproblem becomes challenging. Conversely, if $\sigma_t$ increases too slowly, the number of  iterations will increase, resulting in a slow convergence. 
A more reasonable approach is to adjust the increment of $\sigma_t$ based on the difficulty of solving the current $f_{\sigma_t}(\bx)$. If the previous subproblem converges quickly, a larger $\sigma_{t+1}$ can be chosen; otherwise, $\sigma_t$ should not be increased excessively. 
\item As mentioned in the previous example, when $\sigma_t$ is small, $f_{\sigma_t}(\bx)$ may become unbounded, leading to divergence. When solving the subproblem, if divergence is detected, the iteration should be stopped immediately and the penalty factor should be increased. 
\item The precision of solving the subproblem must be sufficiently high. To ensure convergence, the error in solving the subproblem needs to approach zero.
\end{itemize}

\subsection*{Convergence Analysis}

This subsection explores the convergence properties of the quadratic penalty function method for equality-constrained problems. For simplicity, we assume that for each $\sigma_t$, there exists a minimum point of $f_{\sigma_t}(\bx)$.
Note that this assumption may not hold for some optimization problems since the quadratic penalty function might not sufficiently penalize constraint violations. Therefore, we will not apply the quadratic penalty function method to optimization problems that do not meet this assumption.

\begin{theoremHigh}[Penalty Function Method with Optimal Subproblem]\label{theorem:pen_cong_glo}
Let $\bx^\toptone$ be a global minimum of $f_{\sigma_t}(\bx)$, and let $\sigma_t$ monotonically increase to infinity in Algorithm~\ref{alg:quad_pen_eq}. Then every limit point $\bx^*$ of the sequence $\{\bx^\toptzero\}_{t>0}$ is a global minimum of the original problem (P4.1) as defined in \eqref{equation:pr_in_constchap_equ}.
\end{theoremHigh}

\begin{proof}[of Theorem~\ref{theorem:pen_cong_glo}]
Let $\widehatbx$ be a global minimum of the original problem (P4), i.e.,
$ f(\widehatbx) \leq f(\bx)$, for all $ \bx$ satisfies $c_i(\bx) = 0$, $i \in \mathcalE$.
Since $\bx^\toptone$ is a global minimum of $f_{\sigma_t}(\bx)$, we have $ f_{\sigma_t}(\bx^\toptone) \leq f_{\sigma_t}(\widehatbx), $
i.e.,
\begin{align}
&f(\bx^\toptone) + \frac{\sigma_t}{2} \sum_{i \in \mathcalE} c_i(\bx^\toptone)^2 \leq f(\widehatbx) + \frac{\sigma_t}{2} \sum_{i \in \mathcalE} c_i(\widehatbx)^2 \label{equation:pen_cong_glo1}\\
&\quad\implies\quad 
\sum_{i \in \mathcalE} c_i(\bx^\toptone)^2 \leq \frac{2}{\sigma_t} \big(f(\widehatbx) - f(\bx^\toptone)\big). \label{equation:pen_cong_glo2}
\end{align}
Let $\bx^*$ be a limit point of $\{\bx^\toptzero\}$.  In \eqref{equation:pen_cong_glo2}, taking $t \rightarrow\infty$, based on the continuity of $c_i(\bx)$ and $f(\bx)$ as well as $\sigma_t \rightarrow +\infty$, we have
$ \sum_{i \in \mathcalE} c_i(\bx^*)^2 = 0. $
This shows that $\bx^*$ is a feasible solution of the original problem. From~\eqref{equation:pen_cong_glo1}, we get
$ f(\bx^\toptone) \leq f(\widehatbx). $
Taking the limit on both sides yields $f(\bx^*) \leq f(\widehatbx)$. By the optimality of $\widehatbx$, we conclude that  $f(\bx^*) = f(\widehatbx)$, implying that $\bx^*$ is also a global minimum.
\end{proof}

The above theorem suggests that if a global minimum of the subproblem (Step 4 of Algorithm~\ref{alg:quad_pen_eq}) can be found, then its limit points are global minima of the original problem. However, in practical applications, finding a global minimum of $f_{\sigma_t}(\bx)$ is challenging, so we typically solve the subproblem to a certain level of precision. Consequently, the applicability of Theorem~\ref{theorem:pen_cong_glo} is limited. The following theorem provides an alternative convergence result based on the KKT conditions.

\begin{theoremHigh}[Penalty Function Method with Suboptimal Subproblem]\label{theorem:pen_cong_glo_full}
Let $f(\bx)$ and $c_i(\bx)$, $i \in \mathcalE$ be continuously differentiable, and let $\sigma_t \rightarrow+\infty$ in Algorithm~\ref{alg:quad_pen_eq}. If the solution $\bx^\toptone$ of the subproblem satisfies
$ \normtwobig{\nabla f_{\sigma_t}(\bx^\toptone)} \leq \varepsilon_t$ where $\varepsilon_t \rightarrow0$, 
and for any limit point $\bx^*$ of $\{\bx^\toptzero\}$, the gradients $\{ \nabla c_i(\bx^*) \}_{i \in \mathcalE}$ are linearly independent, then $\bx^*$ is a KKT point (Definition~\ref{definition:kkt_point}) of the equality-constrained optimization problem (P4.1) in \eqref{equation:pr_in_constchap_equ}, and
$$
\lim_{t \rightarrow\infty} (\sigma_t c_i(\bx^\toptone)) = \lambda_i^*, \quad \forall i \in \mathcalE, 
$$
where $\lambda_i^*$ is the Lagrange multiplier corresponding to the constraint $c_i(\bx^*) = 0$.
\end{theoremHigh}
\begin{proof}[of Theorem~\ref{theorem:pen_cong_glo_full}]
Given the stopping criterion for solving the subproblem, for $ \bx^\toptone $, we have:
\begin{equation}\label{equation:pen_cong_glo_full1}
\normtwobig{\nabla f_{\sigma_t}(\bx^\toptone)} =\normtwo{\nabla f(\bx^\toptone) + \sum_{i \in \mathcalE} \sigma_t c_i(\bx^\toptone) \nabla c_i(\bx^\toptone)} \leq \varepsilon_t.
\end{equation}
Using the triangle inequality $ \normtwo{\ba} - \normtwo{\bb} \leq \normtwo{\ba + \bb} $, the above inequality indicates that
$$
\normtwo{\sum_{i \in \mathcalE} c_i(\bx^\toptone) \nabla c_i(\bx^\toptone)} \leq \frac{1}{\sigma_t} \big(\varepsilon_t + \normtwobig{\nabla f(\bx^\toptone)}\big).
$$
Assuming  $ \{\bx^\toptzero\} $ converges to $ \bx^* $, as $ t \rightarrow\infty $, based on the continuity of $ f(\bx) $ and $ c_i(\bx) $, we have:
$$
\sum_{i \in \mathcalE} c_i(\bx^*) \nabla c_i(\bx^*) = 0.
$$
Since $ \nabla c_i(\bx^*), \forall i $ are linearly independent, it must hold that $ c_i(\bx^*) = 0, \,\forall i \in \mathcalE $. This indicates that $ \bx^* $ is a feasible point.

Next, we show that $ \bx^* $ satisfies the KKT conditions. For this, we need to construct the Lagrange multiplier vector $ \blambda^* = [\lambda_1^*, \lambda_2^*, \cdots, \lambda_{p}^*]^\top $, where $ p\triangleq \abs{\mathcalE} $ denotes the number of elements in $ \mathcalE $. Define:
$$
\nabla \bc(\bx) \triangleq [\nabla c_i(\bx)^\top]_{i \in \mathcalE} \in\real^{p\times n},\quad \text{with $n\geq p$},
$$
and
$$
\lambda_i^\toptzero \triangleq -\sigma_t c_i(\bx^\toptone),\,\forall i, \qquad \blambda^\toptzero \triangleq \big[\lambda_1^\toptzero, \lambda_2^\toptzero, \cdots, \lambda_{p}^\toptzero\big]^\top\in\real^p.
$$
Then the gradient of $\nabla f_{\sigma_t}(\bx^\toptone)$ in \eqref{equation:pen_cong_glo_full1} can be rewritten as:
$$
\nabla \bc(\bx^\toptone)^\top \blambda^\toptzero = \nabla f(\bx^\toptone) - \nabla f_{\sigma_t}(\bx^\toptone).
$$
Since $ \nabla c_i(\bx^*), \forall i $ are linearly independent,  $ \nabla \bc(\bx^*) $ is a full-rank matrix. It holds that $ \bx^\toptzero \rightarrow\bx^* $, when $ t $ is sufficiently large, $ \nabla c(\bx^\toptone) $ should also be a full-rank matrix. Therefore, we can use the generalized inverse of $ \nabla c(\bx^\toptone) $ to express $ \blambda^\toptzero $:
$$
\blambda^\toptzero = \big[\nabla \bc(\bx^\toptone) \nabla \bc(\bx^\toptone)^\top\big]^{-1} \nabla \bc(\bx^\toptone) 
\big[\nabla f(\bx^\toptone) - \nabla f_{\sigma_t}(\bx^\toptone)\big].
$$
Taking the limit of both sides as $ t \rightarrow\infty $ and noting that $ \nabla f_{\sigma_t}(\bx^\toptone) \rightarrow \bzero $, we have:
$$
\blambda^* \triangleq \lim_{t \rightarrow \infty} \blambda^\toptzero = \big(\nabla \bc(\bx^*)^\top \nabla \bc(\bx^*)\big)^{-1} \nabla \bc(\bx^*)^\top \nabla f(\bx^*).
$$
Substituting $ t \rightarrow \infty $ into the gradient expression \eqref{equation:pen_cong_glo_full1}, we get:
$$
\nabla f(\bx^*) - \nabla \bc(\bx^*) \blambda^* = \bzero .
$$
This indicates that the gradient condition in the KKT conditions is satisfied (see, for example, Theorem~\ref{theorem:kkt_licq}), and $ \blambda^* $ is the Lagrange multiplier corresponding to the optimal point $ \bx^* $.
\end{proof}

Although Theorem~\ref{theorem:pen_cong_glo_full} does not require each subproblem to be solved exactly, to obtain a solution to the original problem, the precision of the subproblem solutions must increase progressively. 

\subsection*{Numerical Stability}
We briefly analyze the numerical difficulties of Algorithm~\ref{alg:quad_pen_eq}. We know that to obtain a solution to the original problem, the penalty factor $ \sigma_t $ must tend to positive infinity. From the perspective of matrix condition numbers, as $ \sigma_t $ tends to positive infinity, the difficulty of solving the subproblem increases significantly. Consider the Hessian matrix of the penalty function $ f_{\sigma}(\bx) $:
\begin{equation}\label{equation:penl_hess}
\nabla^2 f_{\sigma}(\bx) = \nabla^2 f(\bx) + \sum_{i \in \mathcalE} \sigma \nabla c_i(\bx) \nabla^2 c_i(\bx) + \sigma \nabla \bc(\bx)^\top \nabla \bc(\bx)\in\real^{n\times n}.
\end{equation}
As $ \bx $ approaches the optimal point, according to the proof of Theorem~\ref{theorem:pen_cong_glo_full}, when $ \bx \approx \bx^* $, it should hold that 
$$ 
-\sigma \nabla c_i(\bx) \approx \lambda_i^*.
$$ 
Based on this approximation, we can use the Hessian matrix of the Lagrangian function $ L(\bx, \blambda^*) $ to approximate the first two terms on the right-hand side of \eqref{equation:penl_hess}:
\begin{equation}\label{equation:penl_hess2}
\nabla^2 f_{\sigma}(\bx) \approx \nabla^2 L(\bx, \blambda^*) + \sigma \nabla \bc(\bx)^\top \nabla \bc(\bx), 
\end{equation}
where $ \nabla \bc(\bx)^\top \nabla \bc(\bx) $ is  positive semidefinite  and singular, with $ (n - p) $ eigenvalues equal to 0. 
Note that the right-hand side of \eqref{equation:penl_hess2} contains two matrices: a constant matrix and a singular matrix whose largest eigenvalue tends to positive infinity. Intuitively, the condition number of the Hessian matrix $ \nabla^2 f_{\sigma}(\bx) $ will become increasingly larger, which means that the contour lines of the subproblem will become more dense, making it very difficult to solve using gradient-based algorithms. If Newton's method (Chapter~\ref{chapter:second_order}) is used, solving the Newton equations themselves becomes a very challenging problem. Therefore, in practical applications, we cannot allow the penalty factor to tend to positive infinity.

\subsection{Quadratic Penalty Function Method for General Constraints}

The previous section focused solely on optimization problems with equality constraints. We briefly discuss how to  design a quadratic penalty function for problems that include both equality and inequality constraints (i.e., general constrained optimization problems as in (P4) from \eqref{equation:pr_in_constchap}).

The primary distinction between equality-constrained optimization and general constrained optimization lies in the allowance for inequality constraints to be less than zero ($ c_i(\bx) < 0 $, $i\in\mathcalI$). 
If we were to use the traditional penalty function defined as $ \normtwo{\bc(\bx)}^2 $, it would penalize feasible points where $ c_i(\bx) < 0 $, which is undesirable. 
For problem (P4), we need to modify the existing quadratic penalty function to ensure it only penalizes points where  $ c_i(\bx) > 0 $, $i\in\mathcalI$, without penalizing feasible points.

\begin{definition}[Quadratic Penalty Function for Equality and Inequality Constraints]
For the general constrained optimization problem (P4), define the quadratic penalty function:
$$
f_{\sigma}(\bx) \triangleq f(\bx) + \frac{1}{2} \sigma \left[ \sum_{i \in \mathcalE} c_i^2(\bx) + \sum_{i \in \mathcalI} \widetildec_i^2(\bx) \right],
$$
where the second term on the right-hand side represents  the penalty term, $ \widetildec_i(\bx) $ is defined as
\begin{equation}\label{equation:pen_ineq_widetil}
\widetildec_i(\bx) \triangleq \max\{c_i(\bx), 0\},
\end{equation}
and $ \sigma > 0 $ is  the penalty factor.
\end{definition}

Note that the function $ h(t) \triangleq (\max\{t, 0\})^2 $ is differentiable with respect to $ t $. 
Therefore, the gradient of $f_{\sigma} $ exists, allowing the use of gradient-based methods discussed in Chapters~\ref{chapter:gradient-descent} and \ref{chapter:gd_convg} to solve the subproblem.
However, generally speaking, $ f_{\sigma} $ is not twice differentiable, making it unsuitable for direct application of second-order methods (such as Newton's method; see Chapter~\ref{chapter:second_order}). This limitation is a notable drawback of the quadratic penalty function for general constrained problems. The structure and algorithm of the penalty function method for solving general constrained problems are identical to Algorithm~\ref{alg:quad_pen_eq}, so detailed explanations are omitted here.

\section{Augmented Lagrangian Method}\label{section:aug_lag_method}

In the quadratic penalty function method, when $t\rightarrow \infty$, according to Theorem~\ref{theorem:pen_cong_glo_full}, we have:
\begin{equation}
c_i(\bx^\toptone) \approx -\frac{\lambda_i^*}{\sigma_t}, \quad \forall i \in \mathcalE.
\end{equation}
Therefore, to ensure feasibility, the penalty factor $\sigma_t$ must tend to positive infinity. 
However, this leads to difficulties in solving the subproblem due to the explosion of the condition number.
Therefore, there is a need to improve the quadratic penalty function method to achieve a feasible approximate optimal solution using a finite penalty factor. The augmented Lagrangian method is one such approach that addresses these limitations.

\subsection{Augmented Lagrangian Method for Equality Constraints}

\subsection*{Construction of the Augmented Lagrangian Function}

The \textit{augmented Lagrangian method} constructs an \textit{augmented Lagrangian function (ALF)} at each iteration, combining elements of both the Lagrangian function and the quadratic penalty function. Specifically, for the equality-constrained optimization problem (P4.1) given in \eqref{equation:pr_in_constchap_equ}, the augmented Lagrangian function is defined as:
\begin{equation}\label{equation:auglag_func1}
\begin{aligned}
L_\sigma(\bx, \blambda) 
&\triangleq f(\bx) + \sum_{i \in \mathcalE} \lambda_i c_i(\bx) + \frac{1}{2} \sigma \sum_{i \in \mathcalE} c_i^2(\bx)
= f(\bx) + \bc(\bx)^\top\blambda + \frac{1}{2} \sigma \bc(\bx)^\top\bc(\bx),
\end{aligned}
\end{equation}~\footnote{Since $\blambda\in\real^p$, it can be alternatively stated as $L_\sigma(\bx, \blambda) \triangleq f(\bx) - \bc(\bx)^\top\blambda + \frac{1}{2} \sigma \bc(\bx)^\top\bc(\bx)$ with a reverse sign of $\blambda$.}
which adds the quadratic penalty term to the Lagrangian function.
Here, $\bc(\bx)=\{c_i(\bx)\}_{i=1}^p\in\real^p$  is the vector function $ \bc: \real^n \rightarrow \real^p $, whose $ i $-th component is the $ i $-th constraint function $ c_i $. 
Note that the Jacobian $\bJ \triangleq\nabla \bc(\bx)\in\real^{p\times n}$.

\begin{remark}[Gradient and Hessian of ALF for Equation Constraints]
The gradient and Hessian of the augmented Lagrangian function \eqref{equation:auglag_func1} w.r.t. $\bx$ are given, respectively, by 
\begin{subequations}\label{equation:aug_lag_gra_hess}
\begin{align}
\nabla_{\bx} L_{\sigma}(\bx, \blambda) 
&= \nabla f(\bx) + \sum_{i \in \mathcalE} \big[\lambda_i + \sigma c_i(\bx)\big] \nabla c_i(\bx) \\
&= \nabla f(\bx) + \nabla \bc(\bx)^\top \blambda +  \sigma\nabla \bc(\bx)^\top \bc(\bx)  \qquad\qquad\qquad\qquad\qquad \\
&= \nabla_{\bx} L(\bx, \blambda)  +  \sigma \nabla \bc(\bx)^\top \bc(\bx);
\end{align}
\begin{align}
\nabla_{\bx}^2 L_{\sigma}(\bx, \blambda) 
&=\nabla^2 f(\bx) + \sum_{i \in \mathcalE} \lambda_i\nabla^2 c_i(\bx) 
+  \sum_{i \in \mathcalE}\sigma \left\{\nabla c_i(\bx) \cdot \nabla c_i(\bx)^\top + c_i(\bx)\nabla^2 c_i(\bx)  \right\} \nonumber\\
&=\nabla_{\bx}^2 L(\bx, \blambda) + \sigma \nabla \bc(\bx)^\top \nabla \bc(\bx) 
+\sum_{i \in \mathcalE}\sigma \left\{ c_i(\bx)\nabla^2 c_i(\bx)  \right\},
\end{align}
where $L(\bx,\blambda) =f(\bx) + \sum_{i \in \mathcalE} \lambda_i c_i(\bx)$ is the Lagrangian function.
\end{subequations}
\end{remark}

\paragrapharrow{Properties.}
Notice that the discrepancy introduced by the quadratic penalty term can be relaxed: If $ \blambda = \blambda^* $ is an optimal point of the Lagrangian function, then the first-order conditions in Theorem~\ref{theorem:kkt_licq} and the fact that $ \bc(\bx^*) = \bzero $ implies that $ \bx^* $ is a stationary point of $ L_{\sigma} $:
$$
\nabla_{\bx} L_{\sigma}(\bx^*, \blambda^*) = \bzero.
$$

On the other hand, as mentioned at the beginning of this section, while the penalty function method struggles with the explosion of the condition number as the penalty factor tends to positive infinity, the augmented Lagrangian function method offers a solution. As shown by  \citet{fletcher2000practical}, there exists a finite value $ \widehat{\sigma} $ such that if $ \sigma > \widehat{\sigma} $, then $ \bx^* $ is an unconstrained local minimizer of $ L_{\sigma}(\bx, \blambda^*) $, i.e.,
\begin{equation}\label{equation:finite_alf_sigma}
\text{if }\quad
\bx_{\blambda, \sigma} \triangleq \mathop{\argmin}_{\bx \in \real^n} L_{\sigma}(\bx, \blambda),
\quad\text{ then }\quad
\bx_{\blambda^*, \sigma} = \bx^* \text{ for all } \sigma > \widehat{\sigma};
\end{equation}~\footnote{In case of several local minimizers, $ \mathop{\argmin}_{\bx \in \real^n} $ is interpreted as the local unconstrained minimizer in the valley around $ \bx^* $.} see Theorem~\ref{theorem:finite_sigma_auglag}.
This means that the penalty parameter $ \sigma $ does not have to go to infinity. If $ \sigma $ is sufficiently large and if we insert $ \blambda^* $ (the vector of Lagrangian multipliers at the solution $ \bx^* $), then the unconstrained minimizer of the augmented Lagrangian function solves the constrained problem. Thus, the problem of finding $ \bx^* $ has been reduced  to that of finding $ \blambda^* $.

\paragrapharrow{Iterative method.}
The preceding discussion outlines an \textit{iterative method} for solving the augmented Lagrangian function problem. The core idea is to initially use the penalty term to approach the solution $ \bx^* $, followed by using the Lagrangian term to achieve final convergence as $ \blambda $ approaches  $ \blambda^* $. A rough sketch of the algorithm is provided below:
\begin{tcolorbox}[colback=white,colframe=black]
	\begin{minipage}{0.6\textwidth}
	\qquad\qquad\qquad
\begin{enumerate}
	\item Choose initial values for $ \blambda $ and $ \sigma $;
	\item repeat until stopping criteria satisfied
	\begin{itemize}
		\item Compute $ \bx_{\blambda, \sigma} = \mathop{\argmin}_{\bx \in \real^n} L_{\sigma}(\bx, \blambda)$;
		\item Update $ \blambda $ and $ \sigma $;
	\end{itemize}
\end{enumerate}
	\end{minipage}
\end{tcolorbox}

Next, we introduce two algorithms designed for solving the augmented Lagrangian method applied to problem (P4.1) in \eqref{equation:pr_in_constchap_equ}. The first algorithm aligns the optimality conditions at each iteration with those of the optimal point. In contrast, the second algorithm updates the multipliers using \textit{steepest ascent directions}.

\begin{algorithm}[h] 
\caption{Augmented Lagrangian Method for (P4.1) in \eqref{equation:pr_in_constchap_equ}}
\label{alg:aug_lag_method}
\begin{algorithmic}[1] 
\Require A function $f(\bx)$ and a set of equality constraints $\{c_i(\bx)\}, i\in\mathcalE$; 
\State {\bfseries Input:}  Choose initial point $ \bx^{(1)} $, multiplier $ \blambda^{(1)} $, penalty factor $ \sigma_1 > 0 $, penalty factor update constant $ \rho \in [2, 10] $, constraint violation tolerance $ \varepsilon > 0 $, and precision requirement $ \gamma_t > 0 $. 
\For{$t=1,2,\ldots$}
\State   \algoalign{With $ \bx^\toptzero $ as the initial point, solve $ \mathopmin{\bx}  L_{\sigma_t}(\bx, \blambda^\toptzero)$,
to obtain a solution $ \bx^\toptone $ satisfying the precision condition 
$\normtwobig{\nabla_{\bx} L_{\sigma_t}(\bx^\toptone, \blambda^\toptzero)} \leq \gamma_t$.
}
\If{$ \normtwobig{\bc(\bx^\toptone)} \leq \varepsilon $}
\State Return approximate solution $ \bx^\toptone, \blambda^\toptzero $, terminate iteration.
\EndIf
\State Update multipliers: $ \blambda^\toptone \leftarrow \blambda^\toptzero + \sigma_t \bc(\bx^\toptone) $.
\State Update penalty factor: $ \sigma_{t+1} \leftarrow \rho \sigma_t $.
\EndFor
\State {\bfseries Return:}  final $\bx \leftarrow \bx^\toptzero$;
\end{algorithmic} 
\end{algorithm}
\paragrapharrow{First algorithm.}
In the initial iteration steps we keep $ \blambda $ constant (e.g., $ \blambda = \bzero $) and let the penalty factor $ \sigma $ increase. 
This approach should bring us close to $ \bx^* $, as described in Section~\ref{section:pen_equa} for penalty function methods.

Specifically, in the $t$-th iteration, given the penalty factor $\sigma_t$ and the multiplier $\blambda^\toptzero$, the minimum point $\bx^\toptone$ of the augmented Lagrangian function $L_{\sigma_t}(\bx, \blambda^\toptzero)$ satisfies
\begin{equation}\label{equation:alf_ite1}
\nabla_{\bx} L_{\sigma_t}(\bx^\toptone, \blambda^\toptzero) = \nabla f(\bx^\toptone) 
+ \sum_{i \in \mathcalE} \big[\lambda_i^\toptzero + \sigma_t c_i(\bx^\toptone)\big] \nabla c_i(\bx^\toptone) = \bzero.
\end{equation}
By the KKT conditions (Theorem~\ref{theorem:kkt_licq}), for the optimization problem (P4.1), its optimal solution $\bx^*$ and the corresponding multipliers $\blambda^*$ must satisfy
\begin{equation}\label{equation:alf_ite2}
\nabla f(\bx^*) + \sum_{i \in \mathcalE} \lambda_i^* \nabla c_i(\bx^*) = \bzero.
\end{equation}

By the assumption that $\bx^\toptzero$ and $\bx^\toptone$ are close to $\bx^*$ using penalty function method in the beginning iterations,
comparing \eqref{equation:alf_ite1} and \eqref{equation:alf_ite2} and matching the optimality condition,
to ensure that the iterates generated by the augmented Lagrangian method converge to $ \bx^* $, we need to ensure consistency at the optimal solution. Therefore, for sufficiently large $ t $,
\begin{equation}\label{equation:aug_lag_upd1}
\lambda_i^* \approx \lambda_i^\toptzero + \sigma_t c_i(\bx^\toptone) 
\quad\iff\quad 
c_i(\bx^\toptone) \approx \frac{1}{\sigma_t} (\lambda_i^* - \lambda_i^\toptzero),\quad \forall i \in \mathcalE.
\end{equation}
Thus, when $ \lambda_i^\toptzero $ is sufficiently close to $ \lambda_i^* $, the constraint violation at $ \bx^\toptone $ will be much smaller than $ \frac{1}{\sigma_t} $. The augmented Lagrangian method can effectively reduce the constraint violation by updating the multipliers. Equation~\eqref{equation:aug_lag_upd1} suggests an effective update rule for the multipliers at the $t$-th iteration:
\begin{equation}\label{equation:aug_lag_upd_multi}
\lambda_i^\toptone \leftarrow  \lambda_i^\toptzero + \sigma_t c_i(\bx^\toptone), \quad \forall i \in \mathcalE.
\end{equation}
The augmented Lagrangian method for problem (P4.1) is thus presented in Algorithm~\ref{alg:aug_lag_method}. 
Note that $ \bc(\bx) = [c_i(\bx)]_{i \in \mathcalE} \in\real^p$ and 
$
\nabla \bc(\bx) = [\nabla c_i(\bx)]_{i \in \mathcalE} \in\real^{p\times n}.
$

As the penalty factor $\sigma_t$ increases, the condition number of the Hessian matrix $L_{\sigma_t}(\bx, \blambda^\toptzero)$ with respect to $\bx$ also increases, thereby making it more difficult to solve for the iterate point $\bx^\toptone$. However, when $\sigma_t$ and $\sigma_{t+1}$ are close, $\bx^\toptzero$ can serve as an initial point for solving $\bx^\toptone$, thus accelerating convergence. Therefore, the penalty factor $\sigma_t$ cannot increase too quickly. But if $\sigma_t$ increases too slowly, the sequence of iterates $\{\bx^\toptzero\}$ will converge to the solution of the original problem more slowly. In practice, we need to pay attention to the choice of the parameter $\rho$, and a common empirical value is $\rho \in [2, 10]$.

\paragrapharrow{Second algorithm: steepest ascent or Newton's methods.}
As discussed previously, in the initial iteration steps, we keep $ \blambda $ constant (e.g., $ \blambda = \bzero $) and let the penalty factor $ \sigma $ to increase. This should lead us close to $ \bx^* $ as described for penalty function methods in Section~\ref{section:pen_equa}.

Next, we fix $\sigma=\sigma_{\text{fix}}$ and vary the multipliers $\blambda$. Then, for this fixed penalty factor,
defining 
\begin{subequations}
\begin{align}
\bx_{\blambda} &\triangleq \arg\min_{\bx \in \real^n} L_{\sigma_{\text{fix}}}(\bx, \blambda);\\
\psi(\blambda) &\triangleq L_{\sigma_{\text{fix}}}(\bx_{\blambda}, \blambda) = \min_{\bx \in \real^n} L_{\sigma_{\text{fix}}}(\bx, \blambda),
\end{align}
\end{subequations}
where the latter is a function of $ \blambda $ alone. 
As discussed previously, for the finite value $ \widehat{\sigma}$, we assume $ \sigma_{\text{fix}} > \widehat{\sigma} $. Since
\begin{enumerate}[(i)]
\item $ \psi(\blambda) $ is the minimal value of $ L_{\sigma_{\text{fix}}} $,
\item the definition of the augmented Lagrangian function \eqref{equation:auglag_func1} combined with $ \bc(\bx^*) = \bzero $ shows that $ L_{\sigma}(\bx^*, \blambda) = f(\bx^*) $ for any $ (\blambda, \sigma) $,
\item \eqref{equation:finite_alf_sigma} implies $ \bx_{\blambda^*} = \bx^* $, where $\bx_{\blambda^*} = \arg\min_{\bx \in \real^n} L_{\sigma_{\text{fix}}}(\bx, \blambda^*),$
\end{enumerate}
it follows that for any $ \blambda $
$$
\psi(\blambda) \leq L_{\sigma_{\text{fix}}}(\bx^*, \blambda) = L_{\sigma_{\text{fix}}}(\bx^*, \blambda^*) = \psi(\blambda^*). 
$$
Thus, the Lagrangian multipliers at the solution is a local maximizer for $ \psi $; that is, 
$$
\blambda^* = \arg\max_{\blambda} \psi(\blambda). 
$$

Therefore, using the function $\psi(\blambda)$, in the $t$-th iteration, starting from the current $ \blambda^\toptzero $, we seek an (ascending) step $ \bs^\toptzero $ such that $ \blambda^\toptzero + \bs^\toptzero \approx \blambda^* $.
Similar to the steep descent method in Section~\ref{section:gradient-descent-all}, we can obtain the \textit{steepest ascent step} (the positive gradient direction) as 
\begin{equation}\label{equation:aug_lag_up_secondupdate}
\bs^\toptzero = \eta_t \nabla \psi(\blambda^\toptzero) =\eta_t \bc(\bx_{\blambda^\toptzero}), \quad \eta_t > 0.
\end{equation}
This is equivalent to the update from the first algorithm in \eqref{equation:aug_lag_upd1}, except that we use a learning rate $\eta_t$, rather than the penalty factor $\sigma_t$ to guide the step.
\citet{fletcher2000practical} shows that under certain regularity assumptions, \eqref{equation:aug_lag_upd1} or \eqref{equation:aug_lag_up_secondupdate}   provides linear convergence. Faster convergence can be achieved by applying Newton's method (Chapter~\ref{chapter:second_order}) to the nonlinear problem $ \nabla \psi(\blambda^\toptzero) = \bzero $,
\begin{equation}\label{equation:aug_lag_up_thirddate}
\blambda^* \approx \blambda + \bs^\toptzero, \quad \text{where } \nabla^2\psi(\blambda^\toptzero) \bs^\toptzero = \nabla \psi(\blambda^\toptzero).
\end{equation}
The procedures are analogous to those in Algorithm~\ref{alg:aug_lag_method}, and we shall not repeat the details.

\subsection*{Convergence}

The augmented Lagrangian function, as a type of penalty function, naturally raises questions about the relationship between its local minimizers and those of the original problem (P4.1). 
Suppose $\bx^*$, $\blambda^*$ are the local minimizer and corresponding multiplier of the equality-constrained problem (P4.1), and the second-order sufficient conditions hold ($\nabla^2 f(\bx^*) \succ \bzero$ by Theorem~\ref{theorem:second_suff_loca_int}). 
It can be shown that, given $\blambda^*$, for sufficiently large $\sigma$, $\bx^*$ is also a strict local minimizer of the augmented Lagrangian function $L_\sigma(\bx, \blambda^*)$. When $\blambda^*$ is unknown, for $\blambda$ sufficiently close to $\blambda^*$ and sufficiently large $\sigma$, the local minimizer of the augmented Lagrangian function $L_\sigma(\bx, \blambda)$ will be sufficiently close to $\bx^*$. This also shows that the augmented Lagrangian function is an  penalty function under certain conditions.

\begin{theoremHigh}\label{theorem:finite_sigma_auglag}
Let $\bx^*, \blambda^*$ be a local minimizer and the corresponding multiplier of problem (P4.1), and suppose the linear independence constraint qualification (LICQ, Definition~\ref{definition:licq}) and second-order sufficient conditions hold at $\bx^*$ ($\nabla^2 f(\bx^*) \succ \bzero$ by Theorem~\ref{theorem:second_suff_loca_int}). Then, there exists a finite constant $\widehat{\sigma}$ such that for any $\sigma \geq \widehat{\sigma}$, $\bx^*$ is a strict local minimizer of $L_\sigma(\bx, \blambda^*)$. Conversely, if $\bx^*$ is a local minimizer of $L_\sigma(\bx, \blambda^*)$ and satisfies $c_i(\bx^*) = 0$, $i \in \mathcalE$, then $\bx^*$ is a local minimizer of problem (P4.1).
\end{theoremHigh}
\begin{proof}
Since $\bx^*$ is a local minimizer of problem (P4.1) and the second-order sufficient conditions hold, we have
\begin{equation}\label{equation:opt_aug_theo1}
\small
\begin{aligned}
\nabla_{\bx} L(\bx^*, \blambda^*) &= \nabla f(\bx^*) + \sum_{i \in \mathcalE} \lambda_i^* \nabla c_i(\bx^*) = \bzero,\\
\bu^\top \nabla_{\bx}^2 L(\bx^*, \blambda^*) \bu &= \bu^\top \left( \nabla^2 f(\bx^*) + \sum_{i \in \mathcalE} \lambda_i^* \nabla^2 c_i(\bx^*) \right) \bu > 0, \quad \forall \bu \text{ satisfying } \nabla \bc(\bx^*) \bu = \bzero.
\end{aligned}
\end{equation}
Based on these conditions, we prove that $\bx^*$ is optimal for $L_\sigma(\bx, \blambda^*)$. Since $c_i(\bx^*) = 0$, $i \in \mathcalE$, by \eqref{equation:aug_lag_gra_hess}, we have
$$
\begin{aligned}
\nabla_{\bx} L_\sigma(\bx^*, \blambda^*) &= \nabla_{\bx} L(\bx^*, \blambda^*) = \bzero,\\
\nabla_{\bx}^2 L_\sigma(\bx^*, \blambda^*) &= \nabla_{\bx}^2 L(\bx^*, \blambda^*) + \sigma \nabla \bc(\bx^*)^\top \nabla \bc(\bx^*).
\end{aligned}
$$
For sufficiently large $\sigma$, it can be shown that
$$
\nabla_{\bx}^2 L_\sigma(\bx^*, \blambda^*) \succ \bzero.
$$
In fact, letting the penalty factor be positive integers $\sigma=t$, $t =1, 2, \ldots$, there exists $\bu^\toptzero$ such that $\normtwo{\bu^\toptzero} = 1$ and
$$
\bu^\toptzeroTOP \nabla_{\bx}^2 L_\sigma(\bx^*, \blambda^*) \bu^\toptzero = \bu^\toptzeroTOP \nabla_{\bx}^2 L(\bx^*, \blambda^*) \bu^\toptzero + t \normtwo{\nabla c(\bx^*) \bu^\toptzero}^2 \leq 0,
$$
then
$$
\normtwo{\nabla c(\bx^*) \bu^\toptzero}^2 \leq -\frac{1}{t} \bu^\toptzeroTOP \nabla_{\bx}^2 L(\bx^*, \blambda^*) \bu^\toptzero \to 0, \quad t \to \infty.
$$
Since $\{\bu^\toptzero\}$ is a bounded sequence, there must exist a limit point, denoted by $\bu$. Then,
$$
\nabla \bc(\bx^*) \bu = \bzero, \quad \bu^\top \nabla_{\bx}^2 L(\bx^*, \blambda^*) \bu \leq 0,
$$
which contradicts \eqref{equation:opt_aug_theo1}. Therefore, there exists a finite $\widehat{\sigma}$ such that when $\sigma \geq \widehat{\sigma}$,
$$
\nabla_{\bx}^2 L_\sigma(\bx^*, \blambda^*) \succ 0,
$$
and thus $\bx^*$ is a strict local minimizer of $L_\sigma(\bx, \blambda^*)$ by Theorem~\ref{theorem:second_suff_loca_int}.

Conversely, if $\bx^*$ satisfies $c_i(\bx^*) = 0$ and is a local minimizer of $L_\sigma(\bx, \blambda^*)$, then for any feasible point $\bx$ sufficiently close to $\bx^*$, we have
$$
f(\bx^*) = L_\sigma(\bx^*, \blambda^*) \leq L_\sigma(\bx, \blambda^*) = f(\bx).
$$
Therefore, $\bx^*$ is a local minimizer of problem (P4.1).
\end{proof}

For Algorithm~\ref{alg:aug_lag_method}, by further assuming the boundedness of the multiplier sequence and the constraint qualifications at the limit points, we can prove that the sequence $\{\bx^\toptzero\}$ generated by the algorithm has a subsequence converging to a first-order stable point of problem (P4.1).

\begin{theoremHigh}[Convergence of the Augmented Lagrangian Method \citep{liu2020optimization}]\label{theorem:conv_aug_lag_bounded}
Assume the multiplier sequence $\{\blambda^\toptzero\}$ is bounded, the penalty factor $\sigma_t \to +\infty$ as $t \to \infty$, and the precision requirement $\gamma_t \to 0$ in Algorithm~\ref{alg:aug_lag_method}. Suppose a subsequence $\{\bx^{(t_j)}\}$ of $\{\bx^{(t)}\}_{t>0}$ converges to $\bx^*$, and the linear independence constraint qualification (LICQ) holds at $\bx^*$. Then there exists $\blambda^*$ such that
$$
\blambda^{(t_j+1)} \to \blambda^*, \quad j \to \infty,
$$
$$
\nabla f(\bx^*) + \nabla \bc(\bx^*)^\top \blambda^* = \bzero, \quad \bc(\bx^*) = \bzero.
$$

\end{theoremHigh}
\begin{proof}[of Theorem~\ref{theorem:conv_aug_lag_bounded}]
By \eqref{equation:aug_lag_gra_hess} and the update rule of the multipliers in \eqref{equation:aug_lag_upd_multi}, the augmented Lagrangian function $ L_{\sigma_t}(\bx, \blambda^\toptzero) $ admits the gradient
$$
\begin{aligned}
	\nabla_{\bx} L_{\sigma_t}(\bx^\toptone, \blambda^\toptzero) 
	&= \nabla f(\bx^\toptone) + \big(\nabla \bc(\bx^\toptone)\big)^\top\big(\blambda^\toptzero + \sigma_t \bc(\bx^\toptone)\big)\\
	&= \nabla f(\bx^\toptone) + \big(\nabla \bc(\bx^\toptone)\big)^\top \blambda^\toptone = \nabla_{\bx} L(\bx^\toptone, \blambda^\toptone),
\end{aligned}
$$
where  $\nabla \bc(\bx)\in\real^{p\times n}$ for any $\bx\in\real^n$.
Then, for any $ t_j $ such that $ \rank(\nabla \bc(\bx^{(t_j+1)})) = p = |\mathcalE| $ (which holds when $ \bx^{(t_j+1)} $ is sufficiently close to $ \bx^* $; by assumption that the LICQ is satisfied at the local optimal point $ \bx^* $), we have 
$$
\blambda^{(t_j+1)} = \big[\nabla \bc(\bx^{(t_j+1)}) \nabla \bc(\bx^{(t_j+1)})^\top\big]^{-1} \nabla \bc(\bx^{(t_j+1)})
\left( \nabla_{\bx} L_{\sigma_t}(\bx^{(t_j+1)}, \blambda^{(t_j)}) - \nabla f(\bx^{(t_j+1)}) \right).
$$
Since $ \normtwo{\nabla_{\bx} L_{\sigma_t}(\bx^{(t_j+1)}, \blambda^{(t_j)})} \leq \gamma_{t_j} \to 0 $ by assumption, we have
$$
\blambda^{(t_j+1)} \to \blambda^* \triangleq - (\nabla \bc(\bx^*)^\top \nabla \bc(\bx^*))^{-1} \nabla \bc(\bx^*)^\top \nabla f(\bx^*)
$$
and
$
\nabla_{\bx} L(\bx^*, \blambda^*) = \bzero.
$
Since $ \{ \blambda^\toptzero \} $ is bounded by assumption and $ \blambda^{(t_j)} + \sigma_{t_j} c(\bx^{(t_j+1)}) \to \blambda^* $, it follows that
$
\{ \sigma_{t_j} \bc(\bx^{(t_j+1)}) \}
$   is bounded.
Since  $ \sigma_t \to +\infty $, then
$
\bc(\bx^*) = \bzero,
$ which completes the proof.
\end{proof}

\subsection{Augmented Lagrangian Method for General Constraints}
For a general constrained optimization problem (P4) in \eqref{equation:pr_in_constchap},
we can also define its augmented Lagrangian function and design the corresponding augmented Lagrangian method. 
We introduce two solutions: the activate set method and the method using slack variables.


\subsection*{An Easy Solution: Active Set Method}
A straightforward approach to solving this problem involves using the method just described: Employ a penalty function method to bring us close to a solution, then consider only the active or nearly active constraints as equality constraints. Discard the remaining constraints (though continue monitoring them to ensure they remain inactive) and use one of two methods for updating the vector of Lagrange multipliers $ \blambda $ discussed in the previous section.

The augmented Lagrangian function can be expressed as:
$$
\widetilde{L}_{\sigma}(\bx, \blambda) = f(\bx) + \blambda^\top \bp(\bx) + \frac{1}{2} \sigma \bp(\bx)^\top \bp(\bx),
$$
where $ \bp(\bx) $ is defined as follows
$$
p_i(\bx) \triangleq 
\begin{cases} 
	c_i(\bx), & \text{if } i \in \sA_\delta(\bx); \\
	0, & \text{otherwise},
\end{cases}
$$
where the approximate active set $ \sA_{\delta}(\bx) $ is given by
$$
\sA_{\delta}(\bx) \triangleq \mathcalE \cup \{i \mid i\in \mathcalI \text{ and } c_{i}(\bx) \leq \delta\}, 
$$
where $ \delta $ is a small positive number. Initially we could keep $ \blambda = \bzero $ and increase $ \sigma_t $ until the approximate active set seems to have stabilized (e.g., remains constant for two consecutive iterations). As long as $ \sA_\delta(\bx) $ remains unchanged, we update $ \blambda $ using \eqref{equation:aug_lag_upd1} or \eqref{equation:aug_lag_up_thirddate} (discarding inactive constraints and assuming that the active inequality constraints are numbered first). Otherwise, $ \blambda $ is set to $ \bzero$ and $ \sigma_t $ is increased.

Various  alternatives for defining the active set could be considered; for example, we can determine the activate set  based on the values of $ \abs{c_{i}(\bx)}, i = 1, \ldots, m $. 
However, one drawback of this approach is that it requires specifying a threshold value, such as $ \delta $, which must be provided by the user.
The algorithm is outlined in Algorithm~\ref{alg:aug_lag_method_easy}.
\begin{algorithm}[h]
\caption{Augmented Lagrangian Method for (P4) in \eqref{equation:pr_in_constchap} (An Easy Solution)} 
\label{alg:aug_lag_method_easy}
\begin{algorithmic}[1]
\Require A function $f(\bx)$ and a set of equality constraints $\{c_i(\bx)\}, i\in\mathcalE \cup\mathcalI$; 
\State {\bfseries Input:}  Choose initial point $ \bx^{(1)} $, multiplier $ \blambda^{(1)} =\bzero $, penalty factor $ \sigma_1 > 0 $, penalty factor update constant $ \rho \in [2, 10] $, activate set tolerance $\delta$. 
\For{$t=1,2,\ldots$}
\State $ \bx^\toptone \leftarrow \argmin_{\bx} \widetilde{L}_{\sigma_t}(\bx, \blambda^\toptzero) $
\If{the active set $ \sA_{\delta}(\bx^\toptone) $ is stable compared to $ \sA_{\delta}(\bx^\toptzero) $}
\State $ \blambda^\toptone \leftarrow $  using \eqref{equation:aug_lag_upd1} or \eqref{equation:aug_lag_up_thirddate};
\State $ \sigma_{t+1} \leftarrow   \sigma_t $;
\Else
\State $ \blambda^\toptone \leftarrow \bzero $;
\State $ \sigma_{t+1} \leftarrow \rho  \sigma_t $;
\EndIf
\EndFor
\State {\bfseries Return:} final  $\bx \leftarrow \bx^\toptzero$;
\end{algorithmic}
\end{algorithm}

\subsection*{A Better Solution with Slack Variables}

Alternatively, for problem (P4) in \eqref{equation:pr_in_constchap}, by introducing slack variables, we can obtain the following equivalent form:
\begin{equation}\label{equation:pr_in_constchap_augfull}
\textbf{(P4.2)}:\qquad 
\begin{aligned}
& \underset{\bx\in\real^n,\bs \geq \bzero}{\min} & & f(\bx), \\
& \text{s.t.} & & c_i(\bx) = 0, \quad \,\,&i&\in\mathcalE,\,\,\text{with $\abs{\mathcalE}=p$}, \\
& & & c_i(\bx) + s_i = 0,  \,\,&i&\in\mathcalI,\,\,\text{with $\abs{\mathcalI}=q$}, \\
& & & s_i \geq 0, \quad &i& \in \mathcalI.
\end{aligned}
\end{equation}
Keeping the nonnegativity constraints, we can construct the Lagrangian function as follows:
$$
L(\bx, \bs, \blambda, \bmu) = f(\bx) + \sum_{i \in \mathcalE} \lambda_i c_i(\bx) + \sum_{i \in \mathcalI} \mu_i (c_i(\bx) + s_i), \quad s_i \geq 0, i \in \mathcalI.
$$
Let the quadratic penalty function for the equality constraints in problem \eqref{equation:pr_in_constchap_augfull} be $ g(\bx, \bs) $, then
$$
g(\bx, \bs) = \sum_{i \in \mathcalE} c_i^2(\bx) + \sum_{i \in \mathcalI} (c_i(\bx) + s_i)^2.
$$
We construct the augmented Lagrangian function as follows:
\begin{equation}\label{equation:pr_in_constchap_augfull2}
L_\sigma(\bx, \bs, \blambda, \bmu) = f(\bx) + \sum_{i \in \mathcalE} \lambda_i c_i(\bx) + \sum_{i \in \mathcalI} \mu_i (c_i(\bx) + s_i) + \frac{\sigma}{2} g(\bx, \bs), \; s_i \geq 0, i \in \mathcalI,
\end{equation}
where $\sigma$ is the penalty factor.
For optimizing the augmented Lagrangian function in \eqref{equation:pr_in_constchap_augfull2}, we again consider the alternating descent framework: alternately update $\bx,\bs$ and $\blambda$.
\paragrapharrow{Given multipliers, optimizing primal variables.}
To be more specific, in the $t$-th iteration, given multipliers $\blambda^\toptzero$, $\bmu^\toptzero$ and penalty factor $\sigma_t$, we need to solve the following problem:
\begin{equation}\label{equation:aug_lag_method_full}
\bx^\toptone, \bs^\toptone = \mathop{\argmin}_{\bx\in\real^n,\bs \geq \bzero} \quad L_{\sigma_t}(\bx, \bs, \blambda^\toptzero, \bmu^\toptzero), \quad \text{s.t.} \quad \bs \geq \bzero
\end{equation}
to obtain $\bx^\toptone, \bs^\toptone$. An effective method to solve problem \eqref{equation:aug_lag_method_full} is the projected gradient method (see Section~\ref{section:pgd}). Another method is to eliminate $\bs$ and solve the optimization problem only in terms of $\bx$. Specifically, fixing $\bx$, the subproblem for $\bs$ can be expressed as
$$
\bs^\toptone = \mathop{\argmin}_{\bs \geq \bzero} \quad \sum_{i \in \mathcalI} \mu_i^\toptzero (c_i(\bx) + s_i) + \frac{\sigma_t}{2} \sum_{i \in \mathcalI} (c_i(\bx) + s_i)^2.
$$
According to the property of quadratic functions, $\bs^\toptone$ is a global optimal solution to the above problem if and only if
\begin{equation}\label{equation:aug_lag_slacksol}
s_i^\toptone = \max \left\{ -\frac{\mu_i^\toptzero}{\sigma_t} - c_i(\bx), 0 \right\}, \quad i \in \mathcalI.
\end{equation}
Substituting the expression for $s_i$ into $L_{\sigma_t}$, we have
\begin{equation}\label{equation:auglag_gen_itet}
\small
\begin{aligned}
L_{\sigma_t}(\bx, \blambda^\toptzero, \bmu^\toptzero) 
&=  f(\bx) + \sum_{i \in \mathcalE} \Big(\lambda_i c_i(\bx) + \frac{\sigma_t}{2} c_i^2(\bx) \Big)
+ \frac{\sigma_t}{2} \sum_{i \in \mathcalI} \Big( \max \Big\{ \frac{\mu_i^\toptzero}{\sigma_t} + c_i(\bx), 0 \Big\}^2 - \frac{(\mu_i^\toptzero)^2}{\sigma_t^2} \Big)
\end{aligned}
\end{equation}
which is obtained by $\max\Big\{L_{\sigma_t}(\bx, s_i=0, \blambda^\toptzero, \bmu^\toptzero), L_{\sigma_t}(\bx, s_i=-\frac{\mu_i^\toptzero}{\sigma_t} - c_i(\bx), \blambda^\toptzero, \bmu^\toptzero) \Big\}$ and is a continuously differentiable function of $\bx$ (if $f(\bx), c_i(\bx), i \in \mathcalI \cup \mathcalE$ are continuously differentiable). Therefore, problem \eqref{equation:aug_lag_method_full} is equivalent to
$$
\bx^\toptone=\mathop{\argmin}_{\bx \in \real^n} \quad L_{\sigma_t}(\bx, \blambda^\toptzero, \bmu^\toptzero),
$$
and can be solved using  gradient-based methods (Chapters~\ref{chapter:gradient-descent} and \ref{chapter:gd_convg}). The advantage of this approach is that we eliminate the variable $\bs$, thereby solving the problem in a lower-dimensional space $\real^n$ (the decision space of problem \eqref{equation:aug_lag_method_full} is $\real^{n + \abs{\mathcalI}}$).

\paragrapharrow{Given primal variables, optimizing multipliers.}
For problem (P4.2) in \eqref{equation:pr_in_constchap_augfull}, its optimal solution $\bx^*, \bs^*$ and multipliers $\blambda^*, \bmu^*$ must satisfy the KKT conditions:
\begin{equation}\label{equation:aug_lag_full_kkt_pri}
\text{(KKT of \eqref{equation:pr_in_constchap_augfull})}\qquad 	
\begin{aligned}
\bzero &= \nabla f(\bx^*) + \sum_{i \in \mathcalE} \lambda_i^* \nabla c_i(\bx^*) + \sum_{i \in \mathcalI} \mu_i^* \nabla c_i(\bx^*),\\
\mu_i^* &\geq 0, \quad i \in \mathcalI,\\
s_i^* &\geq 0, \quad i \in \mathcalI.
\end{aligned}
\end{equation}
The optimal solution $\bx^\toptone, \bs^\toptone$ of problem \eqref{equation:aug_lag_method_full} satisfies:
\begin{equation}\label{equation:aug_lag_full_kkt_aug}
\text{(KKT of \eqref{equation:aug_lag_method_full})}\qquad
\begin{aligned}
\bzero 
&= \nabla f(\bx^\toptone) + \sum_{i \in \mathcalE} \left( \lambda_i^\toptzero + \sigma_t c_i(\bx^\toptone) \right) \nabla c_i(\bx^\toptone) \\
&\quad + \sum_{i \in \mathcalI} \left( \mu_i^\toptzero + \sigma_t \big(c_i(\bx^\toptone) + s_i^\toptone\big) \right) \nabla c_i(\bx^\toptone),\\
s_i^\toptone &= \max \Big\{ -\frac{\mu_i^\toptzero}{\sigma_t} - c_i(\bx^\toptone), 0 \Big\}, \quad i \in \mathcalI.
\end{aligned}
\end{equation}
To ensure that the iterates generated by the augmented Lagrangian method converge to $\bx^*, \bs^*$, we need to ensure consistency at the optimal solution. Therefore, for sufficiently large $t$, 
by comparing the KKT conditions of problems \eqref{equation:pr_in_constchap_augfull} and \eqref{equation:aug_lag_method_full} in \eqref{equation:aug_lag_full_kkt_pri} and \eqref{equation:aug_lag_full_kkt_aug}, respectively, it is easy to derive the update rules for the multipliers:
$$
\begin{aligned}
\lambda_i^* &\approx \lambda_i^\toptzero + \sigma_t c_i(\bx^\toptone) 
\;&\implies&  \; 
\lambda_i^\toptone \leftarrow \lambda_i^\toptzero + \sigma_t c_i(\bx^\toptone), \quad i \in \mathcalE,\\
\mu_i^* &\approx\mu_i^\toptzero + \sigma_t \big(c_i(\bx^\toptone) + s_i^\toptone\big)
\;&\implies&  \; 
\mu_i^\toptone \leftarrow \max \left\{ \mu_i^\toptzero + \sigma_t c_i(\bx^\toptone), 0 \right\}, \quad i \in \mathcalI.
\end{aligned}
$$
For equality constraints, the constraint violation could be  defined as:
$$
\nu_t(\bx^\toptone) \triangleq \sqrt{\sum_{i \in \mathcalE} c_i^2(\bx^\toptone) + \sum_{i \in \mathcalI} \left( c_i(\bx^\toptone) + s_i^\toptone \right)^2}.
$$
Based on \eqref{equation:aug_lag_slacksol}, eliminating $\bs$, the constraint violation becomes:
$$
\nu_t(\bx^\toptone) \triangleq \sqrt{\sum_{i \in \mathcalE} c_i^2(\bx^\toptone) + \sum_{i \in \mathcalI} \max \left\{ c_i(\bx^\toptone) , \frac{\mu_i^\toptzero}{\sigma_t} \right\}^2}.
$$

Using these definitions and update rules, we present the augmented Lagrangian method for constrained optimization problem \eqref{equation:pr_in_constchap_augfull} in Algorithm~\ref{alg:aug_lag_method_full}. This algorithm is similar to Algorithm~\ref{alg:aug_lag_method}. After each iteration, the algorithm checks if the constraint violation $\nu_t(\bx^\toptone)$ meets the precision requirement. If it does, the multipliers are updated and the subproblem solution accuracy is increased, keeping the penalty factor unchanged. If it does not, the multipliers are not updated, and the penalty factor is appropriately increased to obtain a smaller constraint violation.

\begin{algorithm}[h] 
\caption{Augmented Lagrangian Method for (P4.2) in \eqref{equation:pr_in_constchap_augfull}}
\label{alg:aug_lag_method_full}
\begin{algorithmic}[1] 
\Require A function $f(\bx)$ and a set of constraints $\{c_i(\bx)\}, i\in\mathcalI\cup\mathcalE$; 
\State {\bfseries Input:}  Choose initial point $ \bx^{(1)} $, multiplier $ \blambda^{(1)}, \bmu^{(1)} $, penalty factor $ \sigma_1 > 0 $, constants $ 0 < \alpha \leq \beta \leq 1 $,
 penalty factor update constant $ \rho > 1 $, the global constraint violation tolerance $ \varepsilon > 0 $ and the initial local constraint violation tolerance $\varepsilon_1 = \frac{1}{\sigma_1^\alpha}$, and the global precision requirement $\gamma$ and the initial  local precision requirement $ \gamma_1 = \frac{1}{\sigma_1} $. 
\For{$t=1,2,\ldots$}
\State   \algoalign{With $ \bx^\toptzero $ as the initial point, solve
$
\min_{\bx}  L_{\sigma_t}(\bx, \blambda^\toptzero, \bmu^\toptzero),
$
to obtain a solution $ \bx^\toptone $ that satisfies the local precision condition
$
\normtwobig{\nabla_{\bx} L_{\sigma_t}(\bx^\toptone, \blambda^\toptzero, \bmu^\toptzero)} \leq \gamma_t.
$
}
\If{$ \nu_t(\bx^\toptone) \leq \varepsilon_t $}
\If{$ \nu_t(\bx^\toptone) \leq \varepsilon $ and $\normtwobig{\nabla_{\bx} L_{\sigma_t}(\bx^\toptone, \blambda^\toptzero, \bmu^\toptzero)} \leq \gamma$}
\State \algoalign{The global precision requirement for the gradient and the constraint violation \\are satisfied: return approximate solution $ \bx^\toptone, \blambda^\toptzero, \bmu^\toptzero $, terminate iteration.}
\EndIf

\State Update the multipliers:
$
\left\{
\begin{aligned}
	\lambda_i^\toptone \leftarrow \lambda_i^\toptzero + \sigma_t c_i(\bx^\toptone), \quad i \in \mathcalE,\\
	\mu_i^\toptone \leftarrow \max \{ \mu_i^\toptzero + \sigma_t c_i(\bx^\toptone), 0 \}, \quad i \in \mathcalI.
\end{aligned}\right.
$

\State Keep the penalty factor unchanged: $ \sigma_{t+1} \leftarrow  \sigma_t $.
\State  Reduce local precision  and constraint violations: $ \gamma_{t+1} \leftarrow \frac{\gamma_t}{\sigma_{t+1}} $, $ \varepsilon_{t+1} \leftarrow \frac{\varepsilon_t}{\sigma_{t+1}^\beta} $.
\Else
\State Keep the multipliers unchanged: $ \blambda^\toptone \leftarrow \blambda^\toptzero $.
\State Update penalty factor: $ \sigma_{t+1} \leftarrow \rho \sigma_t $.
\State Adjust local precision  and constraint violations: $ \gamma_{t+1} \leftarrow \frac{1}{\sigma_{t+1}} $, $ \varepsilon_{t+1} \leftarrow \frac{1}{\sigma_{t+1}^\alpha} $.
\EndIf
\EndFor
\State {\bfseries Return:}  final $\bx \leftarrow \bx^\toptzero$;
\end{algorithmic} 
\end{algorithm}

\section{Alternating Direction Methods of Multipliers}

\subsection{Alternating Direction Methods of Multipliers (ADMM)}\label{section:admm_all}
We provide a brief introduction to the \textit{alternating direction methods of multipliers (ADMM)} and  subsequently discuss its applications in matrix factorization and nonnegative matrix factorization (NMF).

\paragrapharrow{Fixed Augmented Lagrangian Method.}
Consider the following convex optimization problem:
\begin{subequations}\label{equation:admm_subequs}
\begin{equation}\label{equation:admm_prob}
	\textbf{(P5)}:\qquad \mathopmin{\bx, \by} F(\bx)\triangleq f(\bx)+g(\by) \gap \text{s.t.}\gap \bD\bx+\bE\by=\bff,
\end{equation}
where $\bD\in\real^{m\times n}$, $\bE\in\real^{m\times p}$, and $\bff\in\real^m$.
The Lagrangian function is 
\begin{equation}
	L (\bx, \by, \blambda) = f(\bx)+g(\by) + \innerproduct{\blambda, \bD\bx+\bE\by-\bff}.
\end{equation}
The dual objective function is obtained by minimizing the Lagrangian with respect to the primal variables $\bx, \by$ (see Section~\ref{section:sub_conjug}):
$$
Q(\blambda) = \min_{\bx, \by\in \real^n} L(\bx, \by, \blambda) = -f^*(-\bD^\top\blambda) - g^*(-\bE^\top\blambda) - \innerproduct{\blambda, \bff}.
$$
Thus, the dual problem becomes
\begin{equation}\label{equation:admm_dualopt}
	\textbf{(PD)}:\qquad \min_{\blambda \in \real^n} f^*(-\bD^\top\blambda) + g^*(-\bE^\top\blambda) + \innerproduct{\blambda, \bff}.
\end{equation}
\end{subequations}
If we assume $f$ and $g$ are proper closed and convex, the dual problem can be solved using the proximal point method (Section~\ref{section:proximal_point}) due to the properness, closedness, and convexity of the conjugate functions (Exercise~\ref{exercise:closedconv_conj}, Exercise~\ref{exercise:proper_conj}).
Given a parameter $\sigma>0$, the update at the $t$-th iteration is 
$$
\blambda^\toptone \leftarrow 
\argmin_{\blambda \in \real^n} 
\left\{G(\blambda)\triangleq f^*(-\bD^\top\blambda) + g^*(-\bE^\top\blambda) + \innerproduct{\blambda, \bff} + \frac{1}{2\sigma } \normtwobig{\blambda - \blambda^\toptzero}^2\right\}.
$$
By the first-order necessary and sufficient conditions of convex functions (Theorem~\ref{theorem:fetmat_opt}), it follows that 
$$
\bzero \in \partial G(\blambda^\toptone) =
-\bD\partial f^*(-\bD^\top\blambda^\toptone)  -\bE\partial g^*(-\bE^\top\blambda^\toptone) +  \bff + \frac{{\blambda^\toptone - \blambda^\toptzero}}{\sigma }.
$$
The definition of conjugate functions (Definition~\ref{definition:conjug_func}) implies that 
$$
\begin{aligned}
 \min_{\bx\in\real^n} \innerproduct{\bD^\top\blambda^\toptone, \bx} + f(\bx)=  f^*(-\bD^\top\blambda^\toptone);\\
\min_{\by\in\real^n} \innerproduct{\bE^\top\blambda^\toptone, \by} + g(\by)= g^*(-\bE^\top\blambda^\toptone).
\end{aligned}
$$
By the conjugate subgradient theorem (Theorem~\ref{theorem:conju_subgra}), this is equivalent to update $\bx$, $\by$, and $\blambda$ by
$$
\begin{aligned}
\bx^\toptone&\in \argmin_{\bx\in\real^n} \innerproduct{\bD^\top\blambda^\toptone, \bx} + f(\bx)=\partial f^*(-\bD^\top\blambda^\toptone);\\
\by^\toptone&\in \argmin_{\by\in\real^n} \innerproduct{\bE^\top\blambda^\toptone, \by} + g(\by)=\partial g^*(-\bE^\top\blambda^\toptone);\\
\blambda^\toptone &= \blambda^\toptzero + \sigma\left( \bD\bx^\toptone  + \bE \by^\toptone  - \bff \right).
\end{aligned}
$$
Substituting the expression of $\blambda^\toptone$ into the first two equalities and using the first-order necessary and sufficient conditions of convex functions (Theorem~\ref{theorem:fetmat_opt}) again (since $f$ and $g$ are proper closed and convex), it follows that 
\begin{subequations}\label{equation:aug_ad_raw}
\begin{align}
	\bzero&\in  \bD^\top\left(\blambda^\toptzero + \sigma\big( \bD\bx^\toptone  + \bE \by^\toptone  - \bff \big)\right) + \partial f(\bx^\toptone);\label{equation:aug_ad_raw1}\\
	\bzero&\in  \bE^\top\left(\blambda^\toptzero + \sigma\big( \bD\bx^\toptone  + \bE \by^\toptone  - \bff \big)\right) + \partial g(\by^\toptone);\label{equation:aug_ad_raw2}\\
	\blambda^\toptone &= \blambda^\toptzero + \sigma\left( \bD\bx^\toptone  + \bE \by^\toptone  - \bff \right).\label{equation:aug_ad_raw3}
\end{align}
\end{subequations}
Given $\blambda^\toptzero$, the quantities \eqref{equation:aug_ad_raw2} and \eqref{equation:aug_ad_raw3} are satisfied as the optimality condition of 
\begin{equation}
L_\sigma (\bx, \by, \blambda^\toptzero) = f(\bx)+g(\by) + \innerproduct{\blambda^\toptzero, \bD\bx+\bE\by-\bff} + \frac{\sigma}{2}\normtwo{\bD\bx+\bE\by-\bff}^2,
\end{equation}
which is known as the \textit{augmented Lagrangian function} of \eqref{equation:admm_prob} as we have discussed in Section~\ref{section:aug_lag_method}.
When $\sigma=0$, the augmented Lagrangian function reduces to the standard Lagrangian function; when $\sigma>0$, it  acts as  a penalized version of the Lagrangian function.
The \textit{augmented Lagrangian method}  proceeds by iteratively solving (at the $t$-th iteration):
$$
\text{augmented Lagrangian:}
\gap 
\left\{
\begin{aligned}
	(\bx^\toptone, \by^\toptone) &\in \mathop{\argmin}_{\bx, \by} L_\sigma (\bx, \by, \blambda^\toptzero);\\
	\blambda^\toptone&=\blambda^\toptzero + \sigma(\bD\bx^\toptone+\bE\by^\toptone -\bff ).
\end{aligned}
\right.
$$
However, in this case, due to the convergence of the proximal point method (Theorem~\ref{theorem:prox_point}), the parameter $\sigma$ can remain fixed rather than increasing over iterations, unlike in standard augmented Lagrangian methods (Algorithm~\ref{alg:aug_lag_method}).

\paragrapharrow{ADMM.}
One significant challenge is the coupling term between the $\bx$ and  $\by$ variables, which takes the form  $\sigma(\bx^\top\bD^\top\bE\by)$.
ADMM addresses this issue by substituting exact minimization of $(\bx,\by)$ with one iteration of an alternating minimization method. Specifically, for the  $t$-iteration, the ADMM solution updates $\bx, \by$, and $\blambda$ iteratively:
\begin{tcolorbox}[colback=white,colframe=black]
\begin{minipage}{1\textwidth}
\begin{subequations}
\begin{equation}\label{equation:upda_admm1}
\text{ADMM:}\gap
\left\{
\small
\begin{aligned}
\bx^\toptone&\in\mathop{\argmin}_{\bx} \left\{ f(\bx)+ \frac{\sigma}{2}\normtwo{\bD\bx+\bE\by^\toptzero -\bff +\frac{1}{\sigma}\blambda^\toptzero}^2 \right\};\\
\by^\toptone&\in\mathop{\argmin}_{\by} \left\{ g(\by)+ \frac{\sigma}{2}\normtwo{\bD\bx^\toptone+\bE\by -\bff +\frac{1}{\sigma}\blambda^\toptzero}^2 \right\};\\
\blambda^\toptone&=\blambda^\toptzero + \sigma(\bD\bx^\toptone+\bE\by^\toptone -\bff ).
\end{aligned}
\right.
\end{equation}
Let $\widetildebl\triangleq\frac{1}{\sigma}\blambda$, this can be equivalently stated as (the one we will use in the sequel):
\begin{equation}\label{equation:admm_gen_up}
\text{ADMM:}\gap
\left\{
\small
\begin{aligned}
\bx^\toptone&\in\mathop{\argmin}_{\bx} \left\{ f(\bx)+ \frac{\sigma}{2}\normtwo{\bD\bx+\bE\by^\toptzero -\bff +\widetildebl^\toptzero}^2 \right\};\\
\by^\toptone&\in\mathop{\argmin}_{\by} \left\{ g(\by)+ \frac{\sigma}{2}\normtwo{\bD\bx^\toptone+\bE\by -\bff +\widetildebl^\toptzero}^2 \right\};\\
\widetildebl^\toptone&=\widetildebl^\toptzero + (\bD\bx^\toptone+\bE\by^\toptone -\bff ).
\end{aligned}
\right.
\end{equation}
\end{subequations}
\end{minipage}
\end{tcolorbox}

ADMM is widely utilized in solving problems such as LASSO, sparse logistic regression, and basis pursuit. These applications benefit from ADMM's capability to efficiently handle sparsity constraints. Additionally, ADMM can be employed to address optimization challenges related to training support vector machines (SVMs), especially when dealing with large-scale datasets \citep{boyd2011distributed}. 
Here, we provide a brief overview of its application in low-rank matrix factorization problems \citep{lu2021numerical}.

\index{Matrix factorization}
\paragrapharrow{ADMM applied to matrix factorization.}

In the context of low-rank matrix factorization, ADMM iteratively optimizes two subproblems to approximate $\bA\approx\bW\bZ$:
$$
\mathopmin{\bZ} \frac{1}{2}\normf{\bA-\bW\bZ}^2+r(\bZ)
\qquad\text{and}\qquad
\mathopmin{\bW} \frac{1}{2}\normf{\bA-\bW\bZ}^2+r(\bW)
$$
where $\bA\in\real^{M\times N}$, $\bW\in\real^{M\times K}$, and $\bZ\in\real^{K\times N}$ ($K<\min\{M,N\}$).
Given that these two problems are duals of each other, we focus on the first one. This problem can be equivalently formulated using an auxiliary variable $\widetildebZ\in\real^{K\times N}$:
\begin{equation}\label{equation:mf_admm_prob1}
	\mathopmin{\bZ} \frac{1}{2}\normf{\bA-\bW\bZ}^2+r(\widetildebZ), 
	\gap 
	\text{s.t.}
	\gap 
	\bZ=\widetildebZ.
\end{equation}
Following the general ADMM update rule in \eqref{equation:admm_gen_up}, let (a). \{$\bx\leftarrow \bZ$, $\by\leftarrow \widetildebZ$, $\widetildebl\leftarrow \bL$, $\bD=-\bI$,  $\bE=\bI$\} or (b). \{$\bx\leftarrow \bZ$, $\by\leftarrow \widetildebZ$, $\widetildebl\leftarrow \bL$, $\bD=\bI$,  $\bE=-\bI$\}, 
the ADMM updates  for \eqref{equation:mf_admm_prob1} are then given by: 
\begin{equation}\label{equation:admm_gen_als}
	\left\{
	\begin{aligned}
		\bZ 
		&\stackrel{(a)}{\leftarrow} (\bW^\top\bW+\sigma \bI)^{-1} \left[ \bW^\top\bA +\sigma(\widetildebZ+\bL) \right]
		&\stackrel{(b)}{\leftarrow}& (\bW^\top\bW+\sigma \bI)^{-1} \left[ \bW^\top\bA +\sigma(\widetildebZ-\bL) \right];\\
		\widetildebZ
		&\stackrel{(a)}{\leftarrow}\mathop{\argmin}_{\widetildebZ} r(\widetildebZ) + \frac{\sigma}{2}\normf{-\bZ+\widetildebZ + \bL}^2
		&\stackrel{(b)}{\leftarrow}&\mathop{\argmin}_{\widetildebZ} r(\widetildebZ) + \frac{\sigma}{2}\normf{\bZ-\widetildebZ + \bL}^2\\
		\bL&\stackrel{(a)}{\leftarrow}\bL -\bZ+\widetildebZ &\stackrel{(b)}{\leftarrow}& \bL +\bZ-\widetildebZ.
	\end{aligned}
	\right.
\end{equation}
For practical implementation, the Cholesky decomposition of $(\bW^\top\bW+\sigma \bI)$ can be computed to facilitate updates through forward and backward substitutions. The update for $\bW$ follows similarly due to symmetry considerations. In the following discussion, we will consider setting (a) from the above equations.

%

The function $r(\bZ)$ is quite general, allowing us to use various forms. Below are some examples:
\paragraph{ADMM applied to $\ell_1$ regularization.}
For  $\ell_1$ regularization, where  $r(\widetildebZ)=\lambda \normonebig{\widetildebZ}$, the update for each element $(k,n)$ of $\widetildebZ$ is given by
\begin{equation}\label{equation:admm_mf_l1}
\widetildez_{kn}\leftarrow \max(0, 1-\frac{\lambda}{\sigma} \abs{h_{kn}}^{-1}) h_{kn}
\end{equation} 
for all $k\in\{1,2,\ldots,K\}$ and $n\in\{1,2,\ldots,N\}$, where $h_{kn} = z_{kn}-l_{kn}$ (i.e., the elements of $\bH=\bZ-\bL$). 

\paragraph{ADMM applied to smoothness/denoising regularization.}
Smoothness regularization on $\bZ$ can be defined as $r(\widetildebZ)=\frac{\lambda}{2} \Vert\bT\widetildebZ^\top\Vert_F^2$, where $\bT$ is an $N\times N $ tridiagonal matrix with 2 on the main diagonal and $-1$ on the superdiagonal and subdiagonal. This regularization ensures the proximal components in each row of $\widetildebZ$ are smooth (see Problem~\ref{prob:denoise_rls}). The update  of $\widetildebZ$ becomes 
\begin{equation}\label{equation:admm_mf_smooth}
\widetildebZ\leftarrow \sigma\bZ(\lambda \bT^\top\bT +\sigma\bI)^{-1}.
\end{equation} 

\paragraph{ADMM applied to NMF.}
In the case of low-rank nonnegative matrix factorization (NMF) using ADMM, we replace $r(\bZ)$with an indicator function: $r(z_{kn})=0$ if $z_{kn}\geq 0$ and $\infty$ otherwise.
The update for $\widetildebZ$ becomes 
\begin{equation}\label{equation:admm_nmf}
\max\left(\bzero, \bZ-\bL\right),
\end{equation} 
where the max operator is applied  componentwise.

\subsection{Alternating Direction Proximal Method of Multipliers (ADPMM)}\label{section:adpmm}

Following \citet{beck2017first}, we also discuss the \textit{alternating direction proximal method of multipliers (ADPMM)}, which is a slight generalization of ADMM.
For the same problem (P5) in \eqref{equation:admm_prob}, the $t$-th iteration of ADPMM is derived from ADMM in \eqref{equation:upda_admm1}:
\begin{tcolorbox}[colback=white,colframe=black]
\begin{minipage}{1\textwidth}
\begin{equation}
\text{ADPMM:}\quad
\left\{
\small
\begin{aligned}
\bx^\toptone&\in\mathop{\argmin}_{\bx} \left\{ f(\bx)+ \frac{\sigma}{2}\normtwo{\bD\bx+\bE\by^\toptzero -\bff +\frac{1}{\sigma}\blambda^\toptzero}^2  +\textcolor{mylightbluetext}{\frac{1}{2}\norm{\bx-\bx^\toptzero}_{\bA}^2}  \right\};\\
\by^\toptone&\in\mathop{\argmin}_{\by} \left\{ g(\by)+ \frac{\sigma}{2}\normtwo{\bD\bx^\toptone+\bE\by -\bff +\frac{1}{\sigma}\blambda^\toptzero}^2   + \textcolor{mylightbluetext}{\frac{1}{2}\norm{\by-\by^\toptzero}_{\bB}^2} \right\} ;\\
\blambda^\toptone&=\blambda^\toptzero + \sigma(\bD\bx^\toptone+\bE\by^\toptone -\bff ).
\end{aligned}
\right.
\end{equation}
\end{minipage}
\end{tcolorbox}
\noindent 
Here, $\bA\in\real^{n\times n}$ and $\bB\in\real^{p\times p}$ are positive semidefintie matrices. The notation $\norm{\bx}_{\bQ}= (\bx^\top\bQ\bx)^{1/2}$ represents the $\bQ$-norm (see \eqref{equation:q_norm}).
When $\bA=\bzero $ and $\bB=\bzero$, ADPMM simplifies  to ADMM.

One possible choice for $\bA$ and $\bB$ is
\begin{equation}\label{equation:adlpmm_ab}
\bA\triangleq  \mu\bI -\sigma \bD^\top\bD
\qquad \text{and}\qquad 
\bB\triangleq \gamma \bI -\sigma \bE^\top\bE,
\end{equation}
where $\mu\geq \sigma \lambda_{\max}(\bD^\top\bD)$ and $ \gamma  \geq \sigma \lambda_{\max}(\bE^\top\bE)$ such that $\bA$ and $\bB$ are positive semidefinite according to Theorem~\ref{theorem:eigen_charac}.

\begin{assumption}[ADPMM]\label{assumption:ADPMM}
For the ADPMM, we assume
\begin{enumerate}[(i)]
\item $f : \real^n \to (-\infty, \infty]$ and $g : \real^p \to (-\infty, \infty]$ are proper closed convex functions.
\item $\bD \in \real^{m \times n}$, $\bE \in \real^{m \times p}$, $\bff \in \real^m$, $\sigma > 0$.
\item $\bA \in \real^{n\times n}$ and $\bB\in\real^{p\times p}$ are positive semidefinite.
\item For any $\bu \in \real^n$ and $\bv \in \real^p$, the optimal sets of the problems
$$
\min_{\bx \in \real^n} \left\{ f(\bx) + \frac{\sigma}{2} \normtwo{\bD \bx}^2 + \frac{1}{2} \normA{\bx}^2 + \langle \bu, \bx \rangle \right\}
$$
and
$$
\min_{\by \in \real^p} \left\{ g(\by) + \frac{\sigma}{2} \normtwo{\bE \by}^2 + \frac{1}{2} \normB{\by}^2 + \langle \bv, \by \rangle \right\}
$$
are nonempty.
\item There exists $\widetilde{\bx} \in \domain(f)$ and $\widetilde{\by} \in \domain(g)$ for which $\bD \widetilde{\bx} + \bE \widetilde{\by} = \bff$.
\end{enumerate}
\end{assumption}

Property (iv) ensures that the ADPMM method is well-defined. The following theorem demonstrates an
 $\mathcalO(1/T)$ rate of convergence for the sequence generated by ADPMM.

\begin{theoremHigh}[Convergence of ADPMM, $\mathcalO(1/T)$]\label{theorem:ratconv_adpmm}
Suppose that Assumption~\ref{assumption:ADPMM} holds. Let $\big\{(\bx^\toptzero, \by^\toptzero)\big\}_{t > 0}$ be the sequence generated by ADPMM for solving problem (P5) in \eqref{equation:admm_prob}. Let $(\bx^*, \by^*)$ be an optimal solution of problem (P5) and $\blambda^*$ be an optimal solution of the dual problem \eqref{equation:admm_dualopt}. Suppose that $\zeta > 0$ is any constant satisfying $\zeta \geq 2 \normtwo{\blambda^*}$.
Then for all $T >0$,
\begin{align}
F(\overline{\bx}^{(T)}, \overline{\by}^{(T)}) - F(\bx^*, \by^*) &\leq \frac{\normA{\bx^* - \bx^\topone}^2 + \normC{\by^* - \by^\topone}^2 + \frac{1}{\sigma} \big(\zeta + \normtwobig{\blambda^\topone}\big)^2}{2T},
\end{align}
where $\bC \triangleq \sigma \bE^\top \bE + \bB$ and
$$
\overline{\bx}^{(T)} = \frac{1}{T} \sum_{t=1}^T \bx^\toptone
\qquad \text{and}\qquad 
\overline{\by}^{(T)} = \frac{1}{T} \sum_{t=1}^T \by^\toptone.
$$
\end{theoremHigh}
\begin{proof}[of Theorem~\ref{theorem:ratconv_adpmm}]
By the first-order optimality condition of convex functions (Theorem~\ref{theorem:fetmat_opt}) and the update rules of ADPMM, it follows that $\bx^\toptone$ and $\by^\toptone$ satisfy
\begin{align}
- \sigma \bD^\top \left( \bD\bx^\toptone + \bE \by^\toptzero - \bff + \frac{1}{\sigma} \blambda^\toptzero \right) - \bA(\bx^\toptone - \bx^\toptzero) &\in \partial f(\bx^\toptone), \\
- \sigma \bE^\top \left( \bD\bx^\toptone + \bE \by^\toptone - \bff + \frac{1}{\sigma} \blambda^\toptzero \right) - \bB(\by^\toptone - \by^\toptzero) 
&\in \partial g(\by^\toptone).
\end{align}
Letting
$$
\begin{aligned}
\widetildebx^\toptzero &\triangleq \bx^\toptone, 
\qquad
\widetildeby^\toptzero &\triangleq \by^\toptone, 
\qquad
\widetildeblambda^\toptzero &\triangleq \blambda^\toptzero + \sigma (\bD \bx^\toptone + \bE \by^\toptzero - \bff),
\end{aligned}
$$
using the above optimality conditions and the subgradient inequality (Definition~\ref{definition:subgrad}),  we obtain that for any $\bx \in \domain(f)$ and $\by \in \domain(g)$,
$$
\begin{aligned}
	f(\bx) - f(\widetildebx^\toptzero) + \left\langle \sigma \bD^\top \left( \bD \widetildebx^\toptzero + \bE \by^\toptzero - \bff + \frac{1}{\sigma} \blambda^\toptzero \right) + \bA (\widetildebx^\toptzero - \bx^\toptzero), \bx - \widetildebx^\toptzero \right\rangle &\geq 0, \\
	g(\by) - g(\widetildeby^\toptzero) + \left\langle \sigma \bE^\top \left( \bD \widetildebx^\toptzero + \bE \widetildeby^\toptzero - \bff + \frac{1}{\sigma} \blambda^\toptzero \right) + \bB (\widetildeby^\toptzero - \by^\toptzero), \by - \widetildeby^\toptzero \right\rangle &\geq 0.
\end{aligned}
$$
Using the definition of $\widetildeblambda^\toptzero$, the above two inequalities can be equivalently expressed as
$$
\begin{aligned}
	f(\bx) - f(\widetildebx^\toptzero) + \left\langle \bD^\top \widetildeblambda^\toptzero + \bA (\widetildebx^\toptzero - \bx^\toptzero), \bx - \widetildebx^\toptzero \right\rangle &\geq 0, \\
	g(\by) - g(\widetildeby^\toptzero) + \left\langle \bE^\top \widetildeblambda^\toptzero + (\sigma \bE^\top \bE + \bB) (\widetildeby^\toptzero - \by^\toptzero), \by - \widetildeby^\toptzero \right\rangle &\geq 0.
\end{aligned}
$$
Adding the preceding two inequalities and using the  update rule for $\blambda^\toptone$:
$
\blambda^\toptone - \blambda^\toptzero = \sigma (\bD \widetildebx^\toptzero + \bE \widetildeby^\toptzero - \bff),
$
we can conclude that for any $\bx \in \domain(f)$, $\by \in \domain(g)$, and $\blambda \in \real^m$,
\begin{equation}\label{equation:adpmm_prof00}
\footnotesize
\begin{aligned}
F(\bx, \by) - F(\widetildebx^\toptzero, \widetildeby^\toptzero) + 
\innerproduct{\begin{bmatrix} 
\bx - \widetildebx^\toptzero \\ 
\by - \widetildeby^\toptzero \\ 
\blambda - \widetildeblambda^\toptzero 
\end{bmatrix}, 
\begin{bmatrix} 
\bD^\top \widetildeblambda^\toptzero \\ 
\bE^\top \widetildeblambda^\toptzero \\
- \bD \widetildebx^\toptzero - \bE \widetildeby^\toptzero + \bff 
\end{bmatrix} 
- 
\begin{bmatrix} \bA (\bx^\toptzero - \widetildebx^\toptzero) \\ 
\bC (\by^\toptzero - \widetildeby^\toptzero) \\ 
\frac{1}{\sigma} \big(\blambda^\toptzero - \blambda^\toptone \big)
\end{bmatrix} }
\geq 0,
\end{aligned}
\end{equation}
where $\bC = \sigma \bE^\top \bE + \bB$. For any vectors $\ba,\bb,\bc,\bd$, and a positive semidefinite matrix $\bH$, we have:
$$
(\ba - \bb)^\top \bH (\bc - \bd) 
= 
\frac{1}{2} \left( \normH{\ba - \bd}^2 - \normH{\ba - \bc }^2 + \normH{\bb - \bc}^2 - \normH{\bb - \bd}^2 \right).
$$
This concludes that 
\begin{equation}\label{equation:adpmm_prof1}
\begin{aligned}
	(\bx - \widetildebx^\toptzero)^\top \bA (\bx^\toptzero - \widetildebx^\toptzero) &= \frac{1}{2} \left( \normAbig{\bx - \bx^\toptzero}^2 - \normAbig{\bx - \bx^\toptzero}^2 + \normAbig{\widetildebx^\toptzero - \bx^\toptzero}^2 \right) \\
	&\geq \frac{1}{2} \normAbig{\bx - \bx^\toptzero}^2 - \frac{1}{2} \normAbig{\bx - \bx^\toptzero}^2,
\end{aligned}
\end{equation}
and
\begin{equation}\label{equation:adpmm_prof2}
\begin{aligned}
	(\by - \widetildeby^\toptzero)^\top \bC (\by^\toptzero - \widetildeby^\toptzero) &= \frac{1}{2} \normC{\by - \by^\toptzero }^2 - \frac{1}{2} \normC{\by - \by^\toptzero}^2 + \frac{1}{2} \normC{\widetildeby^\toptzero - \by^\toptzero}^2.
\end{aligned}
\end{equation}
Using the expressions of $\widetildeblambda^\toptzero$ and $\blambda^\toptone$, we have 
$$
\small
\begin{aligned}
&2 (\blambda - \widetildeblambda^\toptzero)^\top (\blambda^\toptzero - \blambda^\toptone) 
= \normtwobig{\blambda - \blambda^\toptone}^2 - \normtwobig{\blambda - \blambda^\toptzero}^2 +\normtwobig{\widetildeblambda^\toptzero - \blambda^\toptzero}^2 - \normtwobig{ \widetildeblambda^\toptzero - \blambda^\toptone}^2\\
&= \normtwobig{\blambda - \blambda^\toptone}^2 
- \normtwobig{\blambda - \blambda^\toptzero}^2 
+ \sigma^2\normtwo{\bD \widetildebx^\toptzero + \bE \by^\toptzero - \bff}^2 \\
&\qquad - \normtwo{\blambda^\toptzero + \sigma (\bD \widetildebx^\toptzero + \bE \by^\toptzero - \bff) - \blambda^\toptzero - \sigma (\bD \widetildebx^\toptzero + \bE \widetildeby^\toptzero - \bff)}^2\\
&= \normtwobig{\blambda - \blambda^\toptone}^2 - \normtwobig{\blambda - \blambda^\toptzero}^2 
+ \sigma^2 \normtwo{\bD \widetildebx^\toptzero + \bE \by^\toptzero - \bff}^2 
 - \sigma^2 \normtwo{\bE (\by^\toptzero - \widetildeby^\toptzero)}^2.
\end{aligned}
$$
Therefore,
\begin{equation}\label{equation:adpmm_prof3}
	\small
\begin{aligned}
\frac{1}{\sigma} (\blambda - \widetildeblambda^\toptzero)^\top (\blambda^\toptzero - \blambda^\toptone) \geq \frac{1}{2\sigma} \left( \normtwobig{ \blambda - \blambda^\toptone}^2 - \normtwobig{ \blambda - \blambda^\toptzero}^2 \right) - \frac{\sigma}{2} \normtwobig{\bE (\by^\toptzero - \widetildeby^\toptzero)}^2.
\end{aligned}
\end{equation}
Moreover, we denote
$$
\bZ \triangleq
\small
\begin{bmatrix}
	\bA & 0 & 0 \\
	0 & \bC & 0 \\
	0 & 0 & \frac{1}{\sigma} \mathbf{I}
\end{bmatrix}, 
\normalsize
\qquad
\bu \triangleq 
\small
\begin{bmatrix}
	\bx \\
	\by \\
	\blambda
\end{bmatrix}, 
\normalsize
\qquad
\bu^\toptzero \triangleq 
\small
\begin{bmatrix}
	\bx^\toptzero \\
	\by^\toptzero \\
	\blambda^\toptzero
\end{bmatrix}, 
\normalsize
\quad\text{and}\quad 
\widetildebu^\toptzero \triangleq 
\small
\begin{bmatrix}
	\widetildebx^\toptzero \\
	\widetildeby^\toptzero \\
	\widetildeblambda^\toptzero
\end{bmatrix}.
\normalsize
$$
Combining \eqref{equation:adpmm_prof1}, \eqref{equation:adpmm_prof2}, and \eqref{equation:adpmm_prof3} yields that 
$$
\left\langle
\small
\begin{bmatrix}
	\bx - \widetildebx^\toptzero \\
	\by - \widetildeby^\toptzero \\
	\blambda - \widetildeblambda^\toptzero
\end{bmatrix}, 
\begin{bmatrix}
	\bA (\bx^\toptzero - \widetildebx^\toptzero) \\
	\bC (\by^\toptzero - \widetildeby^\toptzero) \\
	\frac{1}{\sigma} (\blambda^\toptzero - \blambda^\toptone)
\end{bmatrix} 
\normalsize
\right\rangle \geq \frac{1}{2} \normZ{\bu - \bu^\toptone}^2 - \frac{1}{2} \normZ{\bu - \bu^\toptzero}^2.
$$
Combining the preceding inequality with \eqref{equation:adpmm_prof00}, we obtain that for any $\bx \in \domain(f)$, $\by \in \domain(g)$, and $\blambda \in \real^m$,
\begin{equation}\label{equation:adpmm_prof4}
F(\bx, \by) - F(\widetildebx^\toptzero, \widetildeby^\toptzero) +\innerproduct{\bu - \widetildebu^\toptzero, \bG \widetildebu^\toptzero + \widetilde{\bff}} \geq \frac{1}{2} \normZbig{\bu - \bu^\toptone}^2 - \frac{1}{2} \normZbig{\bu - \bu^\toptzero}^2, 
\end{equation}
where
$$
\bG \triangleq 
\small
\begin{bmatrix}
	\bzero & \bzero & \bD^\top \\
	\bzero & \bzero & \bE^\top \\
	-\bD & -\bE & \bzero
\end{bmatrix}
\normalsize
\qquad \text{and}\qquad 
\widetilde{\bff} \small
\triangleq 
\begin{bmatrix}
	\bzero \\
	\bzero \\
	\bff
\end{bmatrix}.
$$
Note that
$$
\innerproduct{\bu - \widetildebu^\toptzero, \bG \widetildebu^\toptzero + \widetilde{\bff}}
= \innerproduct{\bu - \widetildebu^\toptzero, \bG (\widetildebu^\toptzero - \bu) + \bG \bu + \widetilde{\bff}}
= \innerproduct{\bu - \widetildebu^\toptzero, \bG \bu + \widetilde{\bff}},
$$
where the second equality follows from the fact that $\bG$ is skew-symmetric (i.e., $\bG^\top = -\bG$). We can thus conclude that \eqref{equation:adpmm_prof4} can be rewritten as
$$
F(\bx, \by) - F(\widetildebx^\toptzero, \widetildeby^\toptzero) + 
\innerproduct{\bu - \widetildebu^\toptzero, \bG \bu + \widetilde{\bff} }
\geq \frac{1}{2} \normZbig{\bu - \bu^\toptone}^2 - \frac{1}{2} \normZbig{\bu - \bu^\toptzero}^2.
$$
Performing sum of the above inequality over $t = 1, 2, \ldots, T$ yields the inequality
$$
T\cdot  F(\bx, \by) - \sum_{t=1}^T F(\widetildebx^\toptzero, \widetildeby^\toptzero) + 
\innerproduct{T \cdot \bu - \sum_{t=1}^T \widetildebu^\toptzero, \bG \bu + \widetilde{\bff}}
 \geq -\frac{1}{2} \normZbig{\bu - \bu^\topone}^2.
$$
Defining
$$
\overline{\bu}^{(T)} \triangleq \frac{1}{T} \sum_{t=1}^T \widetildebu^\toptzero, 
\qquad 
\overline{\bx}^{(T)} \triangleq \frac{1}{T} \sum_{t=1}^T \bx^\toptone,
\qquad 
\overline{\by}^{(T)} \triangleq \frac{1}{T} \sum_{t=1}^T \by^\toptone,
$$
and using the convexity of $F$, we obtain that
$$
\begin{aligned}
&F(\bx, \by) - F(\overline{\bx}^{(T)}, \overline{\by}^{(T)}) + \innerproduct{\bu - \overline{\bu}^{(T)}, \bG \bu + \widetilde{\bff}} + \frac{1}{2T} \normZbig{\bu - \bu^\topone}^2 \geq 0 \\
&\implies 
F(\bx, \by) - F(\overline{\bx}^{(T)}, \overline{\by}^{(T)}) + \innerproduct{\bu - \overline{\bu}^{(T)}, \bG \overline{\bu}^{(T)} + \widetilde{\bff}} + \frac{1}{2T} \normZbig{\bu - \bu^\topone}^2 \geq 0,
\end{aligned}
$$
where the implication follows from  the skew-symmetry of $\bG$.
In other words, for any $\bx \in \domain(f)$ and $\by \in \domain(g)$,
\begin{equation}\label{equation:adpmm_prof5}
F(\overline{\bx}^{(T)}, \overline{\by}^{(T)}) - F(\bx, \by) + \innerproduct{\overline{\bu}^{(T)} - \bu, \bG \overline{\bu}^{(T)} + \widetilde{\bff}} \leq \frac{1}{2T} \normZbig{\bu - \bu^\topone}^2. 
\end{equation}

Let $(\bx^*, \by^*)$ be an optimal solution of problem (P5). Then $\bD \bx^* + \bE \by^* = \bff$. 
Let $\overline{\blambda}^{(T)} \triangleq \frac{1}{T} \sum_{t=1}^T \blambda^\toptzero$. Plugging $\bx = \bx^*$, $\by = \by^*$, and the expressions for $\overline{\bu}^{(T)}$, $\bu$, $\bu^\topone$, $\bG$, $\bZ$, $\widetilde{\bff}$ into \eqref{equation:adpmm_prof5}, we obtain
$$
\begin{aligned}
& F(\overline{\bx}^{(T)}, \overline{\by}^{(T)}) - F(\bx^*, \by^*) +\innerproduct{\overline{\bx}^{(T)} - \bx^*, \bD^\top \overline{\blambda}^{(T)}} + \innerproduct{\overline{\by}^{(T)} - \by^*, \bE^\top \overline{\blambda}^{(T)}} \\
& \quad + \innerproduct{\overline{\blambda}^{(T)} - \blambda, -\bD \overline{\bx}^{(T)} - \bE \overline{\by}^{(T)} + \bff} \\
& \quad \leq \frac{1}{2T} \left\{ \normAbig{\bx^* - \bx^\topone}^2 + \normCbig{\by^* - \by^\topone}^2 + \frac{1}{\sigma} \normtwobig{ \blambda - \blambda^\topone}^2 \right\}.
\end{aligned}
$$
Cancelling terms and using the fact that $\bD \bx^* + \bE \by^* = \bff$, the above inequality reduces to
$$
\begin{aligned}
&F(\overline{\bx}^{(T)}, \overline{\by}^{(T)}) - F(\bx^*, \by^*) + \innerproduct{\blambda, \bD \overline{\bx}^{(T)} + \bE \overline{\by}^{(T)} - \bff} \\
&\leq \frac{\normA{\bx^* - \bx^\topone}^2 + \normC{\by^* - \by^\topone}^2 + \frac{1}{\sigma} \normtwobig{\blambda - \blambda^\topone }^2}{2T}.
\end{aligned}
$$
Taking $\blambda \in \sB[\bzero, \zeta]$ as the maximum of both sides obtains
$$
\begin{aligned}
&F(\overline{\bx}^{(T)}, \overline{\by}^{(T)}) - F(\bx^*, \by^*) + 
\zeta\normtwo{\bD \overline{\bx}^{(T)} + \bE \overline{\by}^{(T)} - \bff} \\
&\leq \frac{\normA{ \bx^* - \bx^\topone}^2 + \normC{\by^* - \by^\topone}^2 + \frac{1}{\sigma} (\zeta + \normtwobig{\blambda^\topone})^2}{2T}.
\end{aligned}
$$
This concludes the result.
Since $\zeta \geq 2 \normtwo{\blambda^*}$

\end{proof}

In fact, using the strong duality of optimal values of the primal and dual problems (P5) and (PD) in \eqref{equation:admm_subequs}, it can also be shown that 
$$
\normtwo{\bD \overline{\bx}^{(T)} + \bE \overline{\by}^{(T)} - \bff} 
\leq \frac{\normA{\bx^* - \bx^\topone}^2 + \normC{\by^* - \by^\topone}^2 + \frac{1}{\sigma} \big(\zeta + \normtwobig{\blambda^\topone }\big)^2}{\zeta T},
$$
See \citet{beck2017first} for more details.

\begin{problemset}
\item \textbf{ADLPMM.} Considering  ADPMM discussed in Section~\ref{section:adpmm}, let $\bA\triangleq  \mu\bI -\sigma \bD^\top\bD$ and 
$\bB\triangleq \gamma \bI -\sigma \bE^\top\bE$, as defined in \eqref{equation:adlpmm_ab}. Show that the  update of ADPMM becomes
\begin{align*}
\bx^\toptone &\leftarrow \prox_{\frac{1}{\mu} f} \left[ \bx^\toptzero - \frac{\sigma}{\mu} \bD^\top \left( \bD \bx^\toptzero + \bE \by^\toptzero - \bff + \frac{1}{\sigma} \blambda^\toptzero \right) \right]; \\
\by^\toptone &\leftarrow \prox_{\frac{1}{\gamma} g} \left[ \by^\toptzero - \frac{\sigma}{\gamma} \bE^\top \left( \bD \bx^\toptone + \bE \by^\toptzero - \bff + \frac{1}{\sigma} \blambda^\toptzero \right) \right]; \\
\blambda^\toptone &\leftarrow \blambda^\toptzero + \sigma \left( \bD \bx^\toptone + \bE \by^\toptone - \bff \right),
\end{align*}
where $\prox_{\cdot}(\cdot)$ denotes the proximal operator (Definition~\ref{definition:projec_prox_opt}).
This variant is referred to as the \textit{alternating direction linearized proximal method of multipliers (ADLPMM).}

\item Prove the update formulations of ADMM applied to the $\ell_1$ regularization in \eqref{equation:admm_mf_l1}, the smoothness regularization in \eqref{equation:admm_mf_smooth}, and the nonnegative regularization in \eqref{equation:admm_nmf}.
\end{problemset}

\newpage
\chapter{Second-Order Methods}\label{chapter:second_order}
\begingroup
\hypersetup{
linkcolor=structurecolor,
linktoc=page,  
}
\minitoc \newpage
\endgroup

We previously discussed first-order methods in Chapter~\ref{chapter:gd_convg}, which rely solely on function values and gradients.
Beyond these methods, this chapter explores Newton's method and its common variations—including damped Newton's method, Levenberg gradient descent, quasi-Newton methods, and conjugate gradient methods—all of which leverage second-order information to enhance convergence speed (note that the analysis of the conjugate gradient method requires second-order information). Additionally, we examine the trust-region framework, a robust approach for solving nonconvex optimization problems.

\section{Newton's Method}\label{section:new_methods}
Gradient-based methods rely solely on the information of function values and gradients (i.e., first-order information). 
However, if the function $ f(\bx) $ is sufficiently smooth, second-order derivative information can be used to construct a more effective descent direction $ \bd^\toptzero $ for each iteration. 
Newton-type algorithms are those that use second-order derivative information to construct their update rules. 
By incorporating more information, these methods can significantly outperform gradient-based approaches, though they also impose stricter requirements on $ f(\bx) $. This section first introduces the construction and properties of Newton's method, followed by several modified versions.

\subsection{Pure Newton's Method}
Newton's method is a second-order optimization technique that approximates the loss function using a quadratic expansion, providing an estimate of the minimum location based on this approximation.
We will also show that the stochastic optimizer AdaDelta is derived from the principle of unit consistency in second-order methods (Section~\ref{section:adadelta}). 
For a twice continuously differentiable function $ f(\bx) $, consider the quadratic approximation $ f(\bx) $ at the iteration point $ \bx^\toptzero $ (Theorem~\ref{theorem:quad_app_theo}):
\begin{subequations}\label{equation:newton_secon_approx_ALL}
\begin{equation}\label{equation:newton_secon_approx_eq1}
	f(\bx^\toptzero + \bd) = f(\bx^\toptzero) + \nabla f(\bx^\toptzero)^\top \bd + \frac{1}{2} (\bd)^\top \nabla^2 f(\bx^\toptzero) \bd + o(\normtwobig{\bd}^2).
\end{equation}
Our goal is to use this second-order approximation to determine  an appropriate descent direction $ \bd = \bd_{\text{n}}^\toptzero $ (Definition~\ref{definition:uncons_des_direct}). Ignoring the higher-order terms in \eqref{equation:newton_secon_approx_eq1}, we treat the right-hand side as a quadratic function of $ \bd $:
\begin{equation}\label{equation:class_new_quadr}
	g_t(\bd) \triangleq f(\bx^\toptzero) + \nabla f(\bx^\toptzero)^\top \bd + \frac{1}{2} \bd^\top \nabla^2 f(\bx^\toptzero ) \bd.
\end{equation}
We can find its stationary points by setting the gradient of $g_t(\bd)$ to zero, leading to the \textit{Newton equation}
\begin{equation}\label{equation:newton_secon_approx_eq3}
	\nabla^2 f(\bx^\toptzero) \bd = -\nabla f(\bx^\toptzero).
\end{equation}
\end{subequations}
If $ \nabla^2 f(\bx^\toptzero) $ is nonsingular (hence positive definite), the update direction $ \bd=\bd_{\text{n}}^\toptzero $ can be computed as 
$$ 
\bd_{\text{n}}^\toptzero = -\big(\nabla^2 f(\bx^\toptzero)\big)^{-1} \nabla f(\bx^\toptzero),
$$
which is known as the \textit{Newton step} (or \textit{Newton direction} for $f$, at $\bx^\toptzero$).
The positive definiteness of $ \nabla^2 f(\bx^\toptzero) $ also shows that the Newton step is a descent direction unless $\nabla f(\bx^\toptzero)=\bzero$ (Theorem~\ref{theorem:uncons_des_dir}). 
Moreover, because $g_t(\bd)$ is convex when $\nabla f(\bx^\toptzero)$ is positive definite (Exercise~\ref{exercise:conv_quad}),  $\bd_{\text{n}}^\toptzero$ is the minimizer of $g_t(\bd)$.

Therefore, the update formula for   Newton's method at $t$-th iteration is
$$
\bx^\toptone \leftarrow \bx^\toptzero + \bd_{\text{n}}^\toptzero.
$$
In other words, the steepest descent step is rescaled by the inverse of the Hessian.
Intuitively, this suggests that $\bx^\toptzero + \bd_{\text{n}}^\toptzero$ should provide a highly accurate estimate of the minimizer of  $f$, especially when the function $f$ is nearly quadratic or when $\bx$ is near the optimum point $\bx^*$.

Newton's update can be intuitively understood as adjusting the stepsize based on curvature information from the Hessian. Specifically:
\begin{itemize}
\item  When the curvature is steep, the inverse Hessian scales down the stepsize, making updates smaller.
\item When the curvature is flat, the inverse Hessian scales up the stepsize, allowing larger updates.
\end{itemize}
However, for nonlinear functions $f(\bx)$, a single Newton step is usually insufficient to reach the minimum.
Like the first-order optimization methods introduced earlier, Newton's method is applied iteratively, as summarized in  Algorithm~\ref{alg:class_newton} \citep{roweis1996levenberg, goodfellow2016deep}. 
Notably, in the update formula above, the stepsize $\eta_t $ is always set to 1, eliminating the need for stepsize selection. 
For this reason, Newton's method with a fixed stepsize of 1 is often referred to as \textit{classical Newton's method} (\textit{pure Newton's method} or simply \textit{Newton's method}).

\begin{algorithm}[h] 
\caption{Newton's Method}
\label{alg:class_newton}
\begin{algorithmic}[1] 
\Require A twice continuously differentiable function $f(\bx)$; 
\State {\bfseries Input:} Initialize $\bx^{(1)}$;
\For{$t=1,2,\ldots$}
\State $\bd_{\text{n}}^\toptzero \leftarrow $ solution of $\nabla^2 f(\bx^\toptzero) \bd = -\nabla f(\bx^\toptzero)$;
\State $\bx^\toptone \leftarrow \bx^\toptzero + \bd_{\text{n}}^\toptzero$;
\EndFor
\State {\bfseries Return:}  final $\bx\leftarrow \bx^{(t)}$;
\end{algorithmic} 
\end{algorithm}

The computational complexity of Newton's method primarily arises from the need to compute the inverse of the Hessian matrix at each training iteration.
For a function $f$ with $n$ parameters  ($\bx\in \real^n, \nabla^2 f(\bx)\in \real^{n\times n}$), the Hessian matrix contains $n^2$ entries, and its inversion has a computational complexity of $\mathcalO(n^3)$ \citep{trefethen1997numerical, boyd2018introduction, lu2021numerical}. 
As a result, Newton's method is only practical for training models with a relatively small number of parameters, such as shallow neural networks or multi-layer perceptrons.

\paragrapharrow{Newton's method via greedy search.}
In Section~\ref{section:als-gradie-descent-taylor}, we introduced non-Euclidean gradient descent methods that use greedy search with different norms, along with their convergence properties in Section~\ref{section:noneucli_gd}. 
The Newton step can be interpreted as a greedy search using the $\bQ$-norm with $\bQ\triangleq\nabla f(\bx^\toptzero)$ at the $t$-th iteration:
$$
\norm{\bd}_{\nabla f(\bx^\toptzero)}
=
\big(\bd^\top \nabla f(\bx^\toptzero) \bd\big)^{1/2}.
$$
When the Hessian $\nabla f(\bx^\toptzero)$ closely  approximates the Hessian at the optimal point (especially when $\bx^\toptzero$ is close to $\bx^*$), $\nabla f(\bx^\toptzero) \approx \nabla f(\bx^*)$, the condition number for the non-Euclidean gradient descent method (here, equivalent to Newton's method) is close to one, leading to accelerated convergence (see \eqref{equation:gree_cond_opt}).

\paragrapharrow{Newton's method via linearized optimality condition.}
In \eqref{equation:class_new_quadr}, we approximated the function $f(\bx^\toptzero+\bd)$ using a quadratic function. 
Similarly, given the optimality condition $\nabla f(\bx^*)=\bzero$, we can approximate $\nabla f(\bx^\toptzero +\bd)$ using an affine model:
$$
\nabla f(\bx^\toptzero +\bd)\approx \nabla f(\bx) + \nabla^2 f(\bx^\toptzero)\bd=\bzero.
$$ 
This formulation is equivalent to the Newton equation derived in \eqref{equation:newton_secon_approx_eq3}.

\paragrapharrow{Affine invariance of the Newton step.}
Let  $\bA\in\real^{n\times n}$ be a nonsingular matrix,  and define a transformed function $\widetildef(\by) \triangleq f(\bA\by)$. Then, for any $\bx=\bA\by$,
$$
\nabla \widetildef(\by) = \bA^\top\nabla f(\bx)
\qquad\text{and}\qquad
\nabla^2 \widetildef(\by) = \bA^\top\nabla^2 f(\bx)\bA.
$$
Therefore, the Newton step for $\widetildef$ at $\by$ is given by 
$$
\widetildebd_{\text{n}} = - \big(\bA^\top\nabla^2 f(\bx)\bA\big)^{-1}\bA^\top\nabla f(\bx) = \bA^{-1} \bd_{\text{n}}, 
$$ 
where $\bd_{\text{n}}$ is the Newton step for $f$ at $\bx$. 
This result shows that the Newton step is affine invariant:
$$
\bx+\bd_{\text{n}} = \bA (\by+\widetildebd_{\text{n}}).
$$

\paragrapharrow{Equality-constrained Newton's method.}
We  apply Newton's method to problems with linear equality constraints, formulated as:
$$
\min_{\bx} f(\bx) \quad \text{s.t.} \quad \bA\bx = \bb.
$$
One approach to solving this problem is through dual optimization. However, a more direct method is the equality-constrained Newton's method. Initially, we choose $ \bx^{(1)} $ such that $ \bA\bx^{(1)} = \bb $. Then, we iteratively update our solution using:
$$
\bx^\toptone \leftarrow  \bx^\toptzero + \eta_t \bd^\toptzero,
$$
where $\bd^\toptzero$ is determined by minimizing:
$$
\bd^\toptzero = \argmin_{\bA\bd=\bzero} \nabla f(\bx^\toptzero)^\top(\bd - \bx^\toptzero) + \frac{1}{2}(\bd - \bx^\toptzero)^\top\nabla^2 f(\bx^\toptzero)(\bd - \bx^\toptzero)
$$
This process ensures that $ \bx^\toptone $ remains within the feasible set since $ \bA\bx^\toptone = \bA\bx^\toptzero + \eta_t\bA\bd^\toptzero = \bb + \bzero = \bb $.
Furthermore, $ \bd^\toptzero $ represents the solution for minimizing a quadratic function subject to equality constraints. From KKT conditions $ \bd^\toptzero $ satisfies
$$
\begin{bmatrix}
	\nabla^2 f(\bx^\toptzero) & \bA^\top \\
	\bA & \bzero
\end{bmatrix}
\begin{bmatrix}
	\bd^\toptzero \\
	\bv
\end{bmatrix}
=
\begin{bmatrix}
	-\nabla f(\bx) \\
	\bzero
\end{bmatrix}
$$
for some vector $ \bv $. 
Consequently, the Newton direction $ \bd^\toptzero $ is again found by solving a linear system involving the Hessian matrix (though this system is larger due to the inclusion of constraint information).

\subsection{Convergence Analysis}
Newton's method exhibits excellent local convergence properties.
We begin by demonstrating the quadratic rate of convergence of Newton's method for functions that are strongly convex and twice Lipschitz continuously differentiable.
\begin{theoremHigh}[Newton's method for SC and Lipschits Functions]\label{theorem:conv_new_sssc}
Let $ f:\real^n\rightarrow \real $ be a twice continuously differentiable function defined over $\real^n$. Assume that:
\begin{itemize}
\item $f$ is $\alpha$-strongly convex: $ \nabla^2 f(\bx) \geq \alpha \bI $ for any $ \bx \in \real^n $,
\item $f$ is twice Lipschitz continuously differentiable: $ \normtwo{\nabla^2 f(\bx) - \nabla^2 f(\by)} \leq L \normtwo{\bx - \by } $ for any $ \bx, \by \in \real^n $.
\end{itemize}
\noindent Let $ \{\bx^\toptzero\}_{t>0} $ be the sequence generated by Newton's method (Algorithm~\ref{alg:class_newton}), with $ \bx^* $ being the unique minimizer of $ f $ over $\real^n$. Then, for any $ t>0$, we have
\begin{subequations}
\begin{equation}\label{equation:conv_new_sssc_eq1}
\normtwobig{\bx^\toptone - \bx^*} \leq \frac{L}{2\alpha} \normtwobig{\bx^\toptzero - \bx^*}^2
\end{equation}
Additionally, if $ \normtwo{\bx^\topone - \bx^*} \leq \frac{\alpha}{L} $, then
\begin{equation}\label{equation:conv_new_sssc_eq2}
\normtwobig{\bx^\toptzero - \bx^*} \leq \frac{2\alpha}{L} \left( \frac{1}{2} \right)^{2^{t-1}}, \quad \forall t>0.
\end{equation}
\end{subequations}
In other words, Newton's method exhibits \textbf{quadratic convergence} (Definition~\ref{definition:quadratic-convergence}).
\end{theoremHigh}
\begin{proof}[of Theorem~\ref{theorem:conv_new_sssc}]
For $t>0$, using the update rule of Newton's method, the fact that $\nabla  f(\bx^*)=\bzero$, and the fundamental theorem of calculus (Theorem~\ref{theorem:fund_theo_calculu}), we have 
$$
\small
\begin{aligned}
\bx^\toptone - \bx^* 
&= \bx^\toptzero - \bx^* + \big(\nabla^2 f(\bx^\toptzero)\big)^{-1} \big(\nabla f(\bx^*) - \nabla f(\bx^\toptzero)\big) \\
&= \bx^\toptzero - \bx^* + \big(\nabla^2 f(\bx^\toptzero)\big)^{-1} \int_0^1 \big[\nabla^2 f\big(\bx^\toptzero + \mu(\bx^* - \bx^\toptzero)\big) - \nabla^2 f(\bx^\toptzero)\big] (\bx^* - \bx^\toptzero) \, d\mu \\
&= \big(\nabla^2 f(\bx^\toptzero)\big)^{-1} \int_0^1 \big[\nabla^2 f\big(\bx^\toptzero + \mu(\bx^* - \bx^\toptzero)\big) - \nabla^2 f(\bx^\toptzero)\big] (\bx^* - \bx^\toptzero) \, d\mu.
\end{aligned}
$$
The $\alpha$-strongly convexity of $f$ implies that $\normtwo{\big(\nabla^2 f(\bx^\toptzero)\big)^{-1}} \leq \frac{1}{\alpha}$. 
From the above equality, using Cauchy-Schwartz inequality (Proposition~\ref{proposition:cauchy-schwarz-inequ}),
$$
\small
\begin{aligned}
\normtwobig{\bx^\toptone - \bx^*} 
&\leq \normtwo{\big(\nabla^2 f(\bx^\toptzero)\big)^{-1}} \normtwo{\int_0^1 \big[\nabla^2 f\big(\bx^\toptzero + \mu(\bx^* - \bx^\toptzero)\big) - \nabla^2 f(\bx^\toptzero)\big] (\bx^* - \bx^\toptzero) \, d\mu} \\
&\leq \normtwo{\big(\nabla^2 f(\bx^\toptzero)\big)^{-1}} \int_0^1 \normtwo{\big[\nabla^2 f\big(\bx^\toptzero + \mu(\bx^* - \bx^\toptzero)\big) - \nabla^2 f(\bx^\toptzero)\big] (\bx^* - \bx^\toptzero) } \, d\mu \\
&\leq \normtwo{\big(\nabla^2 f(\bx^\toptzero)\big)^{-1}} \int_0^1 \normtwo{\nabla^2 f\big(\bx^\toptzero + \mu(\bx^* - \bx^\toptzero)\big) - \nabla^2 f(\bx^\toptzero)} \cdot \normtwo{\bx^* - \bx^\toptzero} \, d\mu \\
&\leq \frac{L}{\alpha} \normtwobig{\bx^\toptzero - \bx^*}^2\int_0^1 \mu  \, d\mu = \frac{L}{2\alpha} \normtwobig{\bx^\toptzero - \bx^*}^2.
\end{aligned}
$$
This proves the claim in \eqref{equation:conv_new_sssc_eq1}.
To prove \eqref{equation:conv_new_sssc_eq2}, we use induction on $ t $. Note that for $ t = 1 $, we assumed that
$
\normtwo{\bx^\topone - \bx^*} \leq \frac{\alpha}{L},
$
which implies
$$
\normtwo{\bx^\topone - \bx^*} \leq \frac{2\alpha}{L} \left( \frac{1}{2} \right)^{2^0}.
$$
Assume that \eqref{equation:conv_new_sssc_eq2} holds for an integer $ t $, that is,
$
\normtwobig{\bx^\toptzero - \bx^*} \leq \frac{2\alpha}{L} \left( \frac{1}{2} \right)^{2^{t-1}}
$.
We will show it also holds for $ t + 1 $. By \eqref{equation:conv_new_sssc_eq1} we have
$$
\normtwo{\bx^\toptone - \bx^*} \leq \frac{L}{2\alpha} \normtwobig{\bx^\toptzero - \bx^*}^2 \leq \frac{L}{2\alpha} \left( \frac{2\alpha}{L} \left( \frac{1}{2} \right)^{2^{t-1}} \right)^2 = \frac{2\alpha}{L} \left( \frac{1}{2} \right)^{2^{t}}.
$$
This completes the proof.
\end{proof}
The result indicates that Newton's method converges quadratically only when $ \normtwo{\bx^\topone - \bx^*} \leq \frac{\alpha}{L} $ for strongly convex and twice Lipschitz continuously differentiable functions. Therefore, the initial guess significantly affects the performance of Newton's method.

More generally, we have the following local convergence result of Newton's method.
\begin{theoremHigh}[Local Convergence of   Newton's Method]\label{theorem:conv_classNewton}
Let $ f:\real^n\rightarrow \real $ be a twice continuously differentiable function, and let the Hessian matrix be $L$-Lipschitz continuous in a neighborhood $ \sB(\bx^*, \varepsilon) $ of the optimal point $ \bx^* $:
$$
\normtwo{\nabla^2 f(\bx) - \nabla^2 f(\by)} \leq L \normtwo{\bx - \by}, \quad \forall \bx, \by \in \sB(\bx^*, \varepsilon).
$$
If the function $ f(\bx) $ satisfies $ \nabla f(\bx^*) = \bzero $ and $ \nabla^2 f(\bx^*) \succ \bzero $ at the point $ \bx^* $, then  Newton's method (Algorithm~\ref{alg:class_newton}) follows:
\begin{enumerate}[(i)]
\item If the initial point is sufficiently close to $ \bx^* $, the sequence $ \{\bx^\toptzero\} $ generated by Newton's method converges to $ \bx^* $.
\item The sequence $ \{\bx^\toptzero\} $ converges to $ \bx^* $ \textbf{quadratically} (Definition~\ref{definition:quadratic-convergence}).
\item The sequence $ \big\{\normtwobig{\nabla f(\bx^\toptzero)}\big\} $ converges to $\bzero$ \textbf{quadratically}.
\end{enumerate}
\end{theoremHigh}
\begin{proof}[of Theorem~\ref{theorem:conv_classNewton}]
\textbf{(i, ii).}
From the update rule and the property of the optimal point $ \bx^* $ that $ \nabla f(\bx^*) = \bzero $, we obtain:
\begin{equation}\label{equation:conv_classNewton11}
\small
\begin{aligned}
\normtwo{\bx^\toptone - \bx^*} &= \normtwo{\bx^\toptzero - \big(\nabla^2 f(\bx^\toptzero)\big)^{-1} \nabla f(\bx^\toptzero) - \bx^*} \\
&= \normtwo{\big(\nabla^2 f(\bx^\toptzero)\big)^{-1} \left[ \nabla^2 f(\bx^\toptzero) (\bx^\toptzero - \bx^*) - \big(\nabla f(\bx^\toptzero) - \nabla f(\bx^*)\big) \right]}\\
&\leq \normtwobig{\big(\nabla^2 f(\bx^\toptzero)\big)^{-1}} \normtwo{\nabla^2 f(\bx^\toptzero) (\bx^\toptzero - \bx^*) - \big(\nabla f(\bx^\toptzero) - \nabla f(\bx^*)\big) }
\end{aligned}
\end{equation}
By the fundamental theorem of calculus Theorem~\ref{theorem:fund_theo_calculu}\eqref{equation:fund_theo_calculu1}, we have:
$$
\nabla f(\bx^\toptzero) - \nabla f(\bx^*) = \int_0^1 \nabla^2 f\big(\bx^* + \mu(\bx^\toptzero - \bx^*)\big)  \, d\mu\,(\bx^\toptzero - \bx^*),
$$
Thus, assuming $\bx^\toptzero\in \sB(\bx^*, \varepsilon)$, the second term in \eqref{equation:conv_classNewton11} follows:
\begin{equation}\label{equation:conv_classNewton22}
\small
\begin{aligned}
&\quad \normtwo{\nabla^2 f(\bx^\toptzero) (\bx^\toptzero - \bx^*) - \big(\nabla f(\bx^\toptzero) - \nabla f(\bx^*)\big)} \\
&= \normtwo{\int_0^1 \left[ \nabla^2 f\big(\bx^* + \mu(\bx^\toptzero - \bx^*)\big) - \nabla^2 f(\bx^\toptzero) \right] \, d\mu \, \cdot   (\bx^\toptzero - \bx^*)} \\
&\stackrel{\dag}{\leq} \int_0^1 \normtwo{\nabla^2 f\big(\bx^* + \mu(\bx^\toptzero - \bx^*)\big) - \nabla^2 f(\bx^\toptzero)}  \, d\mu \,\cdot  \normtwobig{\bx^\toptzero - \bx^*} \\
&\stackrel{\ddag}{\leq} L\normtwobig{\bx^\toptzero - \bx^*}^2 \int_0^1  \mu \, d\mu 
= \frac{L}{2} \normtwobig{\bx^\toptzero - \bx^*}^2,
\end{aligned}
\end{equation}
where the inequality ($\dag$) follows from the Cauchy-Schwartz inequality (Proposition~\ref{proposition:cauchy-schwarz-inequ}), the inequality ($\ddag$) follows from the Lipschitzness of the Hessian matrix.

Since  $ f $ is twice continuously differentiable, its Hessian $ \nabla^2 f(\bx) $ is a continuous function of $ \bx $. In particular, this means that for any point $ \bx^* $, there exists a neighborhood around $ \bx^* $ where $ \nabla^2 f(\bx) $ is close to $ \nabla^2 f(\bx^*) $ in norm.
Since $ \nabla^2 f(\bx^*) $ is nonsingular, its inverse $ \big(\nabla^2 f(\bx^*)\big)^{-1} $ exists. Moreover, the nonsingularity implies that $ \nabla^2 f(\bx) $ is invertible for all $ \bx $ sufficiently close to $ \bx^* $, because the determinant of $ \nabla^2 f(\bx) $ will not be zero in a small neighborhood around $ \bx^* $.
The map $ \bA \mapsto \bA^{-1} $ is continuous for nonsingular matrices. Specifically, if $ \normtwo{\bA - \bB} $ is small, then $ \normtwo{\bA^{-1} - \bB^{-1}} $ is also small. For $ \nabla^2 f(\bx) $ near $ \nabla^2 f(\bx^*) $, the continuity ensures that $ (\nabla^2 f(\bx))^{-1} $ will remain close to $ \big(\nabla^2 f(\bx^*)\big)^{-1} $.
Since $ \big(\nabla^2 f(\bx^*)\big)^{-1} $ is bounded (as $ \nabla^2 f(\bx^*) $ is nonsingular), the continuity of the inverse ensures that there exists a radius $ \varepsilon_1 > 0 $ such that for all $ \bx^\toptzero $ satisfying $ \bx^\toptzero \in \sB(\bx^*, \varepsilon_1) $, the norm $ \normtwo{\big(\nabla^2 f(\bx^\toptzero)\big)^{-1}} $ is close to $ \normtwo{\big(\nabla^2 f(\bx^*)\big)^{-1}} $. Specifically, one can choose $ \varepsilon_1 $ small enough that $ \normtwo{\big(\nabla^2 f(\bx^\toptzero)\big)^{-1}} \leq 2 \normtwo{\big(\nabla^2 f(\bx^*)\big)^{-1}} $. The factor ``2" comes from choosing a sufficiently small neighborhood around $ \bx^* $ where the deviation of $ \nabla^2 f(\bx^\toptzero) $ from $ \nabla^2 f(\bx^*) $ is controlled. The continuity of $ \nabla^2 f(\bx^\toptzero) $ and its inverse guarantees this bound.
Combining \eqref{equation:conv_classNewton11},  \eqref{equation:conv_classNewton22}, and $ \normtwo{\big(\nabla^2 f(\bx^\toptzero)\big)^{-1}} \leq 2 \normtwo{\big(\nabla^2 f(\bx^*)\big)^{-1}} $, we get:
\begin{equation}\label{equation:conv_classNewton33}
\normtwo{\bx^\toptone - \bx^*} \leq L \normtwo{\nabla^2 f(\bx^*)^{-1}} \normtwobig{\bx^\toptzero - \bx^*}^2.
\end{equation}
Therefore, if the initial point $ \bx^{(1)} $ satisfies
$
\normtwo{\bx^{(1)} - \bx^*} \leq \min \left\{ \varepsilon, \varepsilon_1 \right\} \triangleq \widehat{\varepsilon},
$
then \eqref{equation:conv_classNewton33} shows the sequence $ \{\bx^\toptzero\} $ remains within the neighborhood $ \sB(\bx^*, \widehat{\varepsilon}) $, and thus $ \{\bx^\toptzero\} $ converges quadratically to $ \bx^* $.

\paragraph{(iii).}
From the Newton equation~\eqref{equation:newton_secon_approx_eq3}, $\nabla^2 f(\bx^\toptzero) \bd^\toptzero +\nabla f(\bx^\toptzero) =\bzero$, the update rule $\bx^\toptone \leftarrow \bx^\toptzero + \bd^\toptzero$, and the fundamental theorem of calculus Theorem~\ref{theorem:fund_theo_calculu}\eqref{equation:fund_theo_calculu1},  we have:
$$
\small
\begin{aligned}
\normtwo{\nabla f(\bx^\toptone)} 
&= \normtwo{\nabla f(\bx^\toptone) - \nabla f(\bx^\toptzero) - \nabla^2 f(\bx^\toptzero) \bd^\toptzero} \\
&= \normtwo{\int_0^1 \nabla^2 f(\bx^\toptzero + \mu \bd^\toptzero) \bd^\toptzero d\mu - \nabla^2 f(\bx^\toptzero) \bd^\toptzero} \\
&\leq \int_0^1 \normtwo{\nabla^2 f(\bx^\toptzero + \mu \bd^\toptzero) - \nabla^2 f(\bx^\toptzero)} \normtwobig{\bd^\toptzero} d\mu 
\leq \frac{L}{2} \normtwobig{\bd^\toptzero}^2 \\
&\leq \frac{L}{2}  \normtwo{\big(\nabla^2 f(\bx^\toptzero)\big)^{-1}}^2 \normtwobig{\nabla f(\bx^\toptzero)}^2 
\leq 2 L \normtwo{\nabla^2 f(\bx^*)^{-1}}^2 \normtwobig{\nabla f(\bx^\toptzero)}^2.
\end{aligned}
$$
This proves that the sequence $ \big\{\normtwobig{\nabla f(\bx^\toptzero)}\big\} $ converges to $\bzero$ quadratically.
\end{proof}

Theorem~\ref{theorem:conv_classNewton} states that  Newton's method is a fast-converging algorithm, but its convergence is conditional. First, the initial point $ \bx^{(1)} $ must be sufficiently close to the solution of the problem; in other words, Newton's method exhibits only local convergence. 
When $ \bx^{(1)} $ is far from the solution, Newton's method often fails. Second, the Hessian matrix $ \nabla^2 f(\bx^*) $ needs to be positive definite. There exist examples where, if $ \nabla^2 f(\bx^*) $ is singular or positive semidefinite, the convergence rate of Newton's method may only reach linearly. 
In the proof of Theorem~\ref{theorem:conv_classNewton}, it can be seen that the condition number of the problem does not significantly affect the convergence speed of Newton's method; the Lipschitz constant $ L $ is usually dominated by $ \normtwobig{\bx^\toptzero - \bx^*} $ in later iterations. However, for ill-conditioned problems, the convergence domain of Newton's method may shrink, imposing stricter requirements on the choice of the initial value.

On the other hand, if exact line search is applicable, Newton's method with exact line search (Algorithm~\ref{alg:class_newton_exctline}) is guaranteed to converge globally under mild conditions.
\begin{algorithm}[h] 
\caption{Newton's Method with Exact Line Search}
\label{alg:class_newton_exctline}
\begin{algorithmic}[1] 
\Require A twice continuously differentiable function $f(\bx)$; 
\State {\bfseries Input:} Initialize $\bx^{(1)}$;
\For{$t=1,2,\ldots$}
\State $\bd_{\text{n}}^\toptzero \leftarrow $ solution of $\nabla^2 f(\bx^\toptzero) \bd = -\nabla f(\bx^\toptzero)$;
\State $\eta_t\leftarrow \argmin_{\eta} f(\bx^\toptzero +\eta\bd_{\text{n}}^\toptzero)$; \Comment{exact line search}
\State $\bx^\toptone \leftarrow \bx^\toptzero + \eta_t\bd_{\text{n}}^\toptzero$;
\EndFor
\State {\bfseries Return:}   final $\bx\leftarrow \bx^{(t)}$;
\end{algorithmic} 
\end{algorithm}

\begin{theoremHigh}[Global Convergence of Newton's Method with Exact Line Search]\label{theorem:newton_exactline}
Let $f: \real^n \rightarrow \real$ be twice continuously differentiable on an open convex set $\sS \subset \real^n$. Assume that $f$ is $\alpha$-strongly convex such that $f(\bx)$ satisfies
\begin{equation}
\bu^\top \nabla^2 f(\bx) \bu \geq \alpha \normtwo{\bu}^2, \; \forall \bu \in \real^n, \bx \in \sL,
\end{equation}
where $\sL \triangleq \lev[f, \bx^{(1)}]= \{\bx \in\sS \mid f(\bx) \leq f(\bx^{(1)})\}$ is the corresponding level set. Then the sequence $\{\bx^\toptzero\}$ generated by Newton's method with exact line search (Algorithm~\ref{alg:class_newton_exctline}) satisfies
\begin{enumerate}
\item When $\{\bx^\toptzero\}$ is a finite sequence, then $\nabla f(\bx^\toptzero) = \bzero$ for some $t$;
\item When $\{\bx^\toptzero\}$ is an infinite sequence, then $\{\bx^\toptzero\}$ converges to the unique minimizer $\bx^*$ of $f$.
\end{enumerate}
\end{theoremHigh}

\begin{proof}[of Theorem~\ref{theorem:newton_exactline}]
Since $f$ is $\alpha$-SC, by Theorem~\ref{theorem:exi_close_sc}, its stationary point is the unique global minimizer.
On the other hand, since the level set $\sL$ is a bounded closed convex set and the sequence $\{f(\bx^\toptzero)\}$ is monotonically decreasing, it follows that $\{\bx^\toptzero\} \subset \sL$ and $\{\bx^\toptzero\}$ is bounded. Therefore there exists a limit point $\widetildebx \in \sL$ with $\bx^\toptzero \rightarrow \widetildebx$, and further $f(\bx^\toptzero) \rightarrow f(\widetildebx)$. By Theorem~\ref{theorem:conv_line_search}, we have $\nabla f(\bx^\toptzero) \rightarrow \nabla f(\widetildebx) = \bzero$. Finally, since the stationary point is unique, then the whole sequence $\{\bx^\toptzero\}$ converges to $\widetildebx$ which is the unique minimizer.
\end{proof}

In summary, Newton's method is suitable for high-precision solutions to optimization problems, but it lacks global convergence properties. Therefore, in practical applications, gradient-based methods are often used initially to obtain a low-precision solution, which is then refined using Newton's method for higher precision.
The advantages and disadvantages of Newton's method are summarized in Remark~\ref{remark:procon_newton}.
\begin{remark}[Pros and Cons of Newton's Method]\label{remark:procon_newton}
\textbf{Advantages}
\begin{enumerate}[(i)]
	\item Quadratically convergent from a good starting point if $\nabla f^2(\bx^*)$ is positive definite.
	\item Simple and easy to implement.
\end{enumerate}

\textbf{Disadvantages}
\begin{enumerate}[(i)]
	\item Not globally convergent for many problems.
	\item May converge towards a maximum or saddle point of $f$.
	\item The system of linear equations to be solved in each iteration may be ill-conditioned or singular.
	\item Requires analytic second-order derivatives of $f$.
\end{enumerate}
\end{remark}

\subsection{Modified and Damped Newton Method}\label{section:modified_damp_new}

Although we have proposed Newton's method in Algorithm~\ref{alg:class_newton} and analyzed its theoretical properties, this approach is almost impractical for real-world applications due to several limitations:
\begin{enumerate}
\item Each iteration requires solving an $ n $-dimensional linear system, which leads to large computational costs in high-dimensional problems. The Hessian matrix $ \nabla^2 f(\bx^\toptzero) $ is not only difficult to compute but also challenging to store.

\item When $ \nabla^2 f(\bx^\toptzero) $ is not positive definite, the solution $ \bd^\toptzero $ derived from the Newton equation~\eqref{equation:newton_secon_approx_eq3} is usually poor. For example, when the Hessian matrix is positive definite, $ \bd^\toptzero $ may be unsatisfactory. For instance, when the Hessian matrix is positive definite, $\bd^\toptzero$ serves as a descent direction (Theorem~\ref{theorem:uncons_des_dir}); however, in other cases, it might not be so.

\item When the current point is far from the optimal solution, choosing the stepsize $\eta_t = 1$ directly can lead to significant instability, potentially causing the sequence of iterations to diverge.
\end{enumerate}

To address these issues, modifications to the pure Newton's method are necessary to make it more practical. We introduce the \textit{modified Newton method with line search}. The core idea is to adjust the Hessian matrix  $\nabla^2 f(\bx^\toptzero)$ to ensure it is positive definite and has a smaller condition number. Additionally, incorporating line search improves the stability of the algorithm. The general framework is presented in Algorithm~\ref{alg:modified_newton}.

\begin{algorithm}[H] 
\caption{Modified Newton Method with Line Search}
\label{alg:modified_newton}
\begin{algorithmic}[1] 
\Require A twice continuously differentiable function $f(\bx)$; 
\State {\bfseries Input:} Initialize $\bx^{(1)}$;
\For{$t=1,2,\ldots$}
\State \algoalign{Determine a correction matrix $\bP^\toptzero$ such that $\bB^\toptzero \triangleq \nabla^2 f(\bx^\toptzero) + \bP^\toptzero$ is positive definite and has a small condition number;}
\State Solve the modified Newton equation $\bB^\toptzero \bd = -\nabla f(\bx^\toptzero)$ to obtain direction $\bd^\toptzero$;
\State Use any line search rule to determine the stepsize $\eta_t$;
\State $\bx^\toptone \leftarrow \bx^\toptzero + \eta_t \bd^\toptzero$;
\EndFor
\State {\bfseries Return:}   final $\bx\leftarrow \bx^{(t)}$;
\end{algorithmic} 
\end{algorithm}

\paragrapharrow{Damped Newton method.}
The critical aspect in Algorithm~\ref{alg:modified_newton} is selecting the correction matrix $\bP^\toptzero$. One straightforward choice is to take $\bP^\toptzero \triangleq \mu_t \bI$, i.e., $\bP^\toptzero$ is a scalar multiple of the identity matrix. 
Methods using this approach are known as \textit{damped Newton methods}. 
According to matrix theory, when $\mu_t$ is sufficiently large, $\bB^\toptzero$ will always be  positive definite. 
However,  $\mu_t$ should not be excessively large because, as $\mu_t$ approaches infinity, the direction $\bd^\toptzero_{\text{dn}}$ converges towards the negative gradient direction. 
A more suitable choice involves estimating the smallest eigenvalue of $\nabla^2 f(\bx^\toptzero)$ and then appropriately determining  $\mu_t$.
In this context, instead of finding the step as a stationary point of the quadratic \eqref{equation:class_new_quadr}, the direction $\bd^\toptzero_{\text{dn}}$ is determined as a stationary point of
\begin{equation}\label{equation:modify_new_quadr}
g_{\mu_t}(\bd) \triangleq g_t(\bd) +  \frac{1}{2}  \mu_t\bd^\top\bd
= f(\bx^\toptzero) + \nabla f(\bx^\toptzero)^\top \bd + \frac{1}{2} \bd^\top \big(\nabla^2 f(\bx^\toptzero )+\mu_t\bI\big) \bd.
\end{equation}
When $\mu_t > \min\{\lambda_{\min}, 0\}$, where $\lambda_{\min}$ is the smallest eigenvalue of $\nabla^2 f(\bx^\toptzero )$, Theorem~\ref{theorem:eigen_charac} shows that $\nabla^2 f(\bx^\toptzero )+\mu_t\bI$ is positive definite. Therefore, the function $g_{\mu_t}(\bd)$ becomes a convex function, and $\bd^\toptzero_{\text{dn}}$ is not only a stationary point, but also a minimizer of $g_{\mu_t}$. And Theorem~\ref{theorem:uncons_des_dir} shows that $\bd^\toptzero_{\text{dn}}$ is a descent direction.

To be more specific,
from the \textit{damped Newton equation} 
\begin{equation}\label{equation:damped_newton_eq}
\bd^\toptzero_{\text{dn}} \leftarrow \text{ solution of }\big(\nabla^2 f(\bx^\toptzero )+\mu_t\bI\big) \bd = -\nabla f(\bx^\toptzero),
\end{equation}
if $\mu_t$ is large, the direction $\bd^\toptzero_{\text{dn}}\approx  -\frac{1}{\mu_t} \nabla f(\bx^\toptzero)$ represents a short step in the steepest descent direction. 
This approach is particularly useful in the early stages of the iteration process when the current point  $\bx^\toptzero$ is far from the minimizer $\bx^*$ \citep{frandsen1999unconstrained}.
Conversely, if $\mu_t$ is small, then $\bd^\toptzero_{\text{dn}}$ closely resembles the classical Newton step, which is advantageous when $\bx^\toptzero$ is near $\bx^*$ (Theorem~\ref{theorem:conv_classNewton}). 
Thus, by appropriately adjusting the damping parameter $\mu_t$, we obtain a method that combines the global convergence properties of the steepest descent method with the fast local convergence of  Newton's method.

\paragrapharrow{Levenberg-Marquardt  damped Newton method.}
An interesting relation between a trust region approach (Section~\ref{section:trs_intro} or Section~\ref{section:des_trust_reg}) and Algorithm~\ref{alg:modified_newton} is provided by the following theorem, originally presented by \citet{marquardt1963algorithm}.
\begin{theorem}[Damped Newton as Trust Region Approach]\label{theorem:damped_trust}
Let $f:\real^n\rightarrow \real^n$ be a twice continuously difference function.
If the matrix $ \nabla^2 f(\bx^\toptzero) + \mu_t \bI $ is positive definite, then
$$
\bd_{\text{dn}}^\toptzero = \mathop{\argmin}_{\normtwo{\bd} \leq \normtwobig{\bd^\toptzero_{\text{dn}}}} \{ g_t(\bd) \},
$$
where $ g_t $ is given by \eqref{equation:class_new_quadr} and $ \bd_{\text{dn}}^\toptzero $ is obtained by the damped Newton equation $\big(\nabla^2 f(\bx^\toptzero)+\mu_t\bI\big) \bd = -\nabla f(\bx^\toptzero)$.
\end{theorem}
\begin{proof}[of Theorem~\ref{theorem:damped_trust}]
For any $\mu_t \geq 0$,  let $g_{\mu_t}(\bd) \triangleq g_t(\bd) + \frac{1}{2} \mu_t \bd^\top \bd$. The gradient of $ g_{\mu_t} $ is
$$
\nabla g_{\mu_t}(\bd) \triangleq \nabla  g_t(\bd) + \mu_t \bd \triangleq \bg^\toptzero + (\bH^\toptzero + \mu_t \bI) \bd,
$$
where $ \bg^\toptzero \triangleq \nabla f(\bx^\toptzero) $ and $ \bH^\toptzero \triangleq \nabla^2 f(\bx^\toptzero) $. Since $ \bH^\toptzero + \mu_t \bI $ is positive definite,  the linear system of equations $ \nabla g_{\mu_t}(\bd) = \bzero $ has a unique solution, which is the minimizer of $ g_{\mu_t} $. This solution is recognized as $ \bd_{\text{dn}}^\toptzero $.

On the other hand, let
$
\widetildebd \triangleq \mathop{\argmin}_{\normtwo{\bd} \leq \normtwo{\bd_{\text{dn}}^\toptzero}} \{ g_t(\bd) \}.
$
Then, $ g_t(\widetildebd) \leq g_t(\bd_{\text{dn}}^\toptzero) $ and $ \widetildebd^\top \widetildebd \leq \bd_{\text{dn}}^\toptzeroTOP \bd_{\text{dn}}^\toptzero $, so that
$$
g_{\mu_t}(\widetildebd) = g_t(\widetildebd) + \frac{1}{2} \mu_t \widetildebd^\top \widetildebd \leq g_t(\bd_{\text{dn}}^\toptzero) + \frac{1}{2} \mu_t \bd_{\text{dn}}^\toptzeroTOP \bd_{\text{dn}}^\toptzero = g_{\mu_t}(\bd_{\text{dn}}^\toptzero).
$$
However, $ \bd_{\text{dn}}^\toptzero $ is the unique minimizer of $ g_{\mu_t} $, so $ \widetildebd = \bd_{\text{dn}}^\toptzero $.
This completes the proof.
\end{proof}

In a proper trust region method, we monitor the trust region radius $ \Delta_t $ at each iteration. 
The above theorem indicates that if we instead monitor the damping parameter, we can interpret this approach as a trust region method where the trust region radius is implicitly defined as $ \Delta_t = \normtwobig{\bd_{\text{dn}}^\toptzero} $.

For  Levenberg-Marquardt type methods, $ \mu_t $ is updated during each iteration. 
Given the current value of the parameter, the Cholesky factorization of $ \nabla^2 f(\bx^\toptzero) + \mu_t \bI $ is used  to check for positive definiteness (see, for example, \citet{lu2021numerical} for more details), and $ \mu_t $ is increased if the matrix is not sufficiently  positive definite. Otherwise, the solution $ \bd_{\text{dn}}^\toptzero $ is efficiently obtained through the factorization (using forward and backward substitution procedures).

\begin{SCfigure}
\centering
\includegraphics[width=0.5\textwidth]{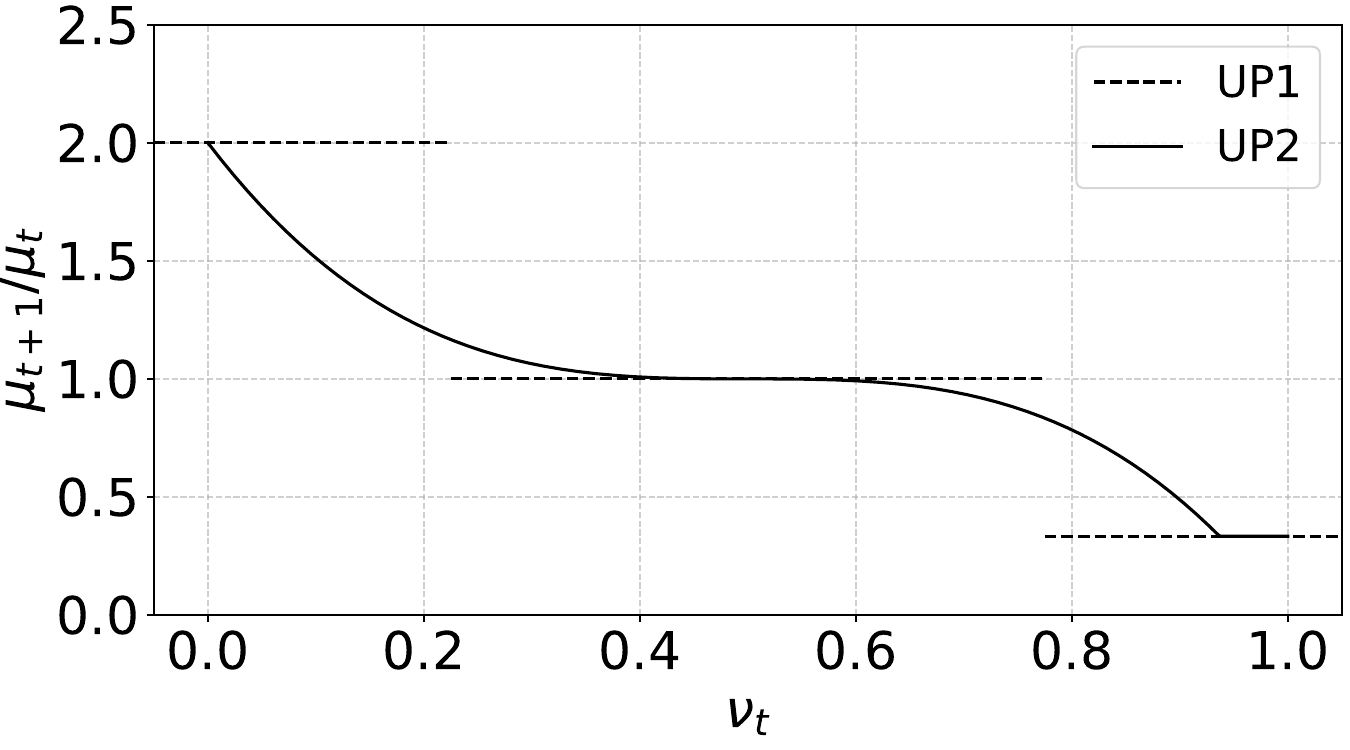} 
\caption{Updating of the damping parameter $ \mu_{t+1} $ by \eqref{equation:newtrs_up1} (dashed line) and by \eqref{equation:newtrs_up2} (solid line).}
\label{fig:newton_trust_reg}
\end{SCfigure}

The direction provided  by $ \bd_{\text{dn}}^\toptzero $ ensures a descent direction, leading to the ``trial point" $ \bx^\toptzero + \bd_{\text{dn}}^\toptzero $ (corresponding to $ \eta_t = 1 $ in Algorithm~\ref{alg:modified_newton}). As in a trust region method (Section~\ref{section:des_trust_reg}), we evaluate the cost function at the trial point, i.e., $ f(\bx^\toptzero + \bd_{\text{dn}}^\toptzero) $. If it is sufficiently below $ f(\bx^\toptzero) $, then the trial point is chosen as the next iterate. 
Otherwise, $ \bx^\toptzero $ remains the current iterate (corresponding to $ \eta_t = 0 $ in Algorithm~\ref{alg:modified_newton}), and $ \mu_t $ is increased. 
Simply checking whether $ f(\bx^\toptzero + \bd_{\text{dn}}^\toptzero) < f(\bx^\toptzero) $ is insufficient; to ensure overall convergence, it is necessary to test if the actual decrease in the function value is greater than some small fraction of the decrease predicted by the quadratic model \eqref{equation:class_new_quadr}, i.e., if
$$
\nu_t \triangleq \frac{f(\bx^\toptzero) - f(\bx^\toptzero + \bd_{\text{dn}}^\toptzero)}{g_t(\bzero) - g_t(\bd_{\text{dn}}^\toptzero)} > \delta,
$$
where $ \delta $ is a small positive number (typically $ \delta = 10^{-3} $).
We recognize $ \nu_t $ as the \textit{gain factor}. It is also used to monitor $ \mu_t $: 
\begin{itemize}
\item If $ \nu_t $ is close to one, then one may expect the model to be a good approximation to $ f $ in a neighbourhood of $ \bx^\toptzero $, which is large enough to include the trial point, and the influence of Newton's method should be increased by decreasing $ \mu_t $. 
\item If, on the other hand, the actual decrease of $ f $ is much smaller than expected, then $ \mu_t $ must be increased in order to adjust the method more towards the steepest descent method. It is important to note that in this case the length of $ \bd_{\text{dn}}^\toptzero $ is also reduced, since $\bd^\toptzero_{\text{dn}}\approx  -\frac{1}{\mu_t} \nabla f(\bx^\toptzero)$ when $\mu_t$ is large.
\end{itemize}
An updating strategy similar to the one used in Algorithm~\ref{alg:trust_region0} could be:
\begin{equation}\label{equation:newtrs_up1}
\text{(UP1)}:\qquad
\boxed{\begin{aligned}
	&\text{if } \nu_t > 0.75 \\
	&\qquad\mu_{t+1} \leftarrow  \frac{1}{3}\mu_t; \\
	&\text{if } \nu_t < 0.25 \\
	&\qquad\mu_{t+1} \leftarrow  2\cdot \mu_t ;
\end{aligned}}
\end{equation}
However, abrupt  changes in $ \mu_t $ when $ \nu_t $ is near  0.25 or 0.75 can cause ``flutter," slowing down convergence.  Therefore, \citet{nielsen1999damping, frandsen1999unconstrained, madsen2004methods}  recommend using a smoother strategy:
\begin{equation}\label{equation:newtrs_up2}
\text{(UP2)}:\qquad
\boxed{\begin{aligned}
	&\text{if } \nu_t > 0 \\
	&\qquad\mu_{t+1} \leftarrow   \max\left\{\frac{1}{3}, 1 - (2\nu_t - 1)^3\right\} \cdot \mu_t ;\\
	&\text{else} \\
	&\qquad\mu_{t+1} \leftarrow 2\cdot \mu_t;
\end{aligned}}
\end{equation}
Refer to Figure~\ref{fig:newton_trust_reg} for a graphical illustration. The outlined method is detailed in  Algorithm~\ref{algorithm:leve_dam_new}.
Following \eqref{equation:des_stopcri1} and \eqref{equation:des_stopcri2}, the stopping criteria can be either
$$
\norminf{\nabla f(\bx^\toptzero)} \leq \varepsilon_1 \qquad \text{ or } \qquad \normtwo{\bd_{\text{dn}}} \leq \varepsilon_2 (\varepsilon_2 + \normtwobig{\bx^\toptzero}). 
$$
While the simplicity of the original Newton's method has been lost in the pursuit of global convergence, this method generally performs well.

\begin{algorithm}[h]
\caption{Levenberg-Marquardt Damped Newton Method \citep{goldfeld1966maximization, madsen2010and}}
\label{algorithm:leve_dam_new}
\begin{algorithmic}[1]
\Require A twice continuously differentiable function $f(\bx)$; 
\State {\bfseries Input:}  Initialize $\bx^{(1)}, \mu_1$, and $ \delta$ (by default $\delta\triangleq 10^{-3}$);
\For{$t=1,2,\ldots$}
\While{$\nabla^2 f(\bx^\toptzero) + \mu_t \bI$ is not positive definite}
\State $\mu_t \leftarrow 2\mu_t$
\EndWhile
\State $\bd_{\text{dn}}^\toptzero\leftarrow $ solution of  $(\nabla^2 f(\bx^\toptzero) + \mu_t \bI) \bd = -\nabla f(\bx^\toptzero)$;
\State Compute gain factor $\nu_t$;
\If{$\nu_t > \delta$} \Comment{$f$ decreases}
\State $\bx^\toptone \leftarrow \bx^\toptzero + \bd_{\text{dn}}^\toptzero$ \Comment{new iterate}
\State $\mu_{t+1} \leftarrow  \max\left\{\frac{1}{3}, 1 - (2\nu_t - 1)^3\right\} \cdot \mu_t$ 
\Else
\State $\bx^\toptone \leftarrow \bx^\toptzero $;  \Comment{old iterate}
\State $\mu_{t+1} \leftarrow 2\cdot \mu_t $;
\EndIf
\EndFor
\State {\bfseries Return:}   final $\bx\leftarrow \bx^{(t)}$;
\end{algorithmic}
\end{algorithm}

\paragrapharrow{Gill-Murray's modified Cholesky factorization.}
Another way to choose the correction matrix $\bP^\toptzero$ is implicitly through modifying the Cholesky decomposition to solve the Newton equation~\eqref{equation:newton_secon_approx_eq3} \citep{gill1974newton, gill2019practical}. We know that when the Hessian matrix is positive definite, Equation~\eqref{equation:newton_secon_approx_eq3} can be  easily solved using Cholesky decomposition (Theorem~\ref{theorem:cholesky-factor-exist}, using forward and backward substitutions). 
However, if the Hessian matrix is indefinite or has a high condition number (the ratio of the largest to smallest eigenvalues), the standard Cholesky decomposition may fail.
The \textit{modified Cholesky decomposition} algorithm adjusts the basic Cholesky decomposition algorithm  to ensure the decomposed matrix remains close to the original matrix while their product is positive definite.

To be more specific, for any  positive definite matrix $\bA = \{a_{ij}\}$, its Cholesky decomposition can be written as
$ \bA = \bL\bD\bL^\top, $
where $\bL = \{l_{ij}\}$ is a unit lower triangular matrix with diagonal elements all equal to 1, and $\bD = \diag(d_{11}, d_{22}, \ldots, d_{nn})$ is a diagonal matrix with positive diagonal elements.  
Based on the form of Cholesky decomposition in  Algorithm~\ref{alg:compute-choklesky-_ldl}, if $\bA$ is positive definite with a small condition number, the diagonal elements of matrix $\bD$ should not be too small (Problem~\ref{problem:cond_pd}). If during the computation process, $d_{jj}$ for any $j\in\{1,2,\ldots,n\}$ is found to be too small, it must be corrected immediately. 
Additionally, we need to ensure that the corrections are bounded so that the elements of the modified matrix also have upper bounds. Specifically, two positive parameters $\alpha, \beta$ are chosen such that
$$ 
d_{jj} \geq \alpha>0, \quad r_{ij}\triangleq l_{ij} \sqrt{d_{jj}} \leq \beta, \quad i = j + 1, j + 2, \ldots, n. 
$$
In Algorithm~\ref{alg:compute-choklesky-_ldl}, only the update of $ d_{jj} $ needs modification to meet these conditions.  
The specific update rule is equivalent to updating each $d_{jj}$ in Algorithm~\ref{alg:compute-choklesky-_ldl} by
$$ 
d_{jj} \leftarrow \max \left\{ |c_{jj}|, \left( \frac{\theta_j}{\beta} \right)^2, \alpha \right\}, \quad \theta_j = \max_{i > j} |c_{ij}|. 
$$
It can be demonstrated that the modified Cholesky decomposition algorithm effectively computes the Cholesky decomposition of the corrected matrix $ \nabla^2 f(\bx^\toptzero) + \bP^\toptzero $, where $ \bP^\toptzero $ is a diagonal matrix with nonnegative diagonal elements. When $ \nabla^2 f(\bx^\toptzero) $ is positive definite and has a small condition number, $ \bP^\toptzero = \bzero $ \citep{gill1974newton, sun2006optimization}.

\subsection{Inexact Newton's Method}

In  pure Newton's method, computing the Newton direction $ \bd_{\text{n}}^\toptzero $ involves solving a linear system. 
When $ n $ is large but $ \nabla^2 f(\bx^\toptzero) $ exhibits a sparse structure, iterative methods are necessary to solve the Newton equation~\eqref{equation:newton_secon_approx_eq3}. Since iterative methods inherently introduce some level of error, it's important to understand how this affects the convergence of Newton's method and how we can control the precision of the solution to ensure that Newton's method still converges. We will briefly address these issues.

For simplicity, let $F(\bx)\triangleq\nabla f(\bx)$.
Recall that the basic Newton step is obtained by solving the Newton equation \eqref{equation:newton_secon_approx_eq3}:
\begin{equation}
\nabla F(\bx^\toptzero)\bd = -F(\bx^\toptzero)
\end{equation}
and setting
\begin{equation}
\bx^\toptone \leftarrow \bx^\toptzero + \bd_{\text{in}}^\toptzero.
\end{equation}
Consider an inexact solution $ \bd_{\text{in}}^\toptzero $ of the Newton equation. We introduce the residual  vector $ \br^\toptzero $ to represent the residual, then the inexact Newton direction satisfies the \textit{inexact Newton equation}:
\begin{equation}
\bd_{\text{in}}^\toptzero\leftarrow \text{ solution of } \nabla F(\bx^\toptzero)\bd = -F(\bx^\toptzero) + \br^\toptzero,
\end{equation}
where we assume the relative error $ \rho_t $ satisfies:
\begin{equation}
\normtwobig{\br^\toptzero} \leq \rho_t \normtwobig{F(\bx^\toptzero)}.
\end{equation}
Here, $ \br^\toptzero = \nabla F(\bx^\toptzero)\bd_{\text{in}}^\toptzero + F(\bx^\toptzero) $ denotes the residual, and $ \{\rho_t\} $ (with $ 0 < \rho_t < 1 $) denotes a forcing sequence controlling the level of inexactness. The procedure is outlined in Algorithm~\ref{alg:inexac_newton}, under the following assumption:
\begin{assumption}[Inexact Newton's Method]\label{assumption:inex_newton}
Let $F(\bx)\triangleq\nabla f(\bx)$ for any differentiable function $f:\real^n\rightarrow \real$. We assume that
\begin{enumerate}[(i)]
\item There exists $ \bx^* $ such that $ F(\bx^*) = \bzero $.

\item $ F $ is continuously differentiable in the neighborhood of $ \bx^* $.

\item $ \nabla F(\bx^*) $ is nonsingular.
\end{enumerate}
\end{assumption}

\begin{algorithm}[h] 
\caption{Inexact Newton's Method \citep{sun2006optimization}}
\label{alg:inexac_newton}
\begin{algorithmic}[1] 
\Require A twice continuously differentiable function $f(\bx)$ and $F(\bx)\triangleq\nabla f(\bx)$; 
\State {\bfseries Input:}  Initialize $\bx^{(1)}$;
\For{$t=1,2,\ldots$}
\State Select noise level $\rho_t$
\State $\bd_{\text{in}}^\toptzero\leftarrow $ solution of  $\nabla F(\bx^\toptzero)\bd = -F(\bx^\toptzero) + \br^\toptzero$, where $\normtwobig{\br^\toptzero} \leq \rho_t \normtwobig{F(\bx^\toptzero)}$;
\State $\bx^\toptone \leftarrow \bx^\toptzero + \bd_{\text{in}}^\toptzero$;
\EndFor
\State {\bfseries Return:}  final $\bx\leftarrow \bx^{(t)}$;
\end{algorithmic} 
\end{algorithm}

In the following, we establish the linear convergence of inexact Newton's method.
\begin{theoremHigh}[Local Convergence of Inexact Newton's Method]\label{theorem:inex_newton_linear}
Let $F : \real^n \to \real^n$ satisfy Assumption~\ref{assumption:inex_newton}. Assume that the sequence $\{\rho_t\}$ satisfies $0 \leq \rho_t \leq \rho < \zeta < 1$ (i.e., $\rho$ is the upper bound of the sequence $\{\rho_t\}$). Then, for some $\epsilon > 0$, if the starting point $\bx^{(1)}$ is sufficiently near $\bx^*$, the sequence $\{\bx^\toptzero\}$ generated by inexact Newton's method (Algorithm~\ref{alg:inexac_newton}) converges \textbf{linearly} to $\bx^*$ (Definition~\ref{definition:linear-convergence}), i.e.,
\begin{equation}\label{equation:inex_newton_linear_res}
\normbig{\bx^\toptone - \bx^*}_* \leq \zeta \normbig{\bx^\toptzero - \bx^*}_*,
\end{equation}
where $\norm{\by}_* \triangleq \normtwo{\nabla F(\bx^*)\by}$ for any $\by\in\real^n$.
\end{theoremHigh}
\begin{proof}[of Theorem~\ref{theorem:inex_newton_linear}]
Since $\nabla F(\bx^*)$ is nonsingular, for $\by \in \real^n$, we have
\begin{equation}\label{equation:inex_newton_linear0}
\frac{1}{\mu} \normtwo{\by} \leq \norm{\by}_* \leq \mu \normtwo{\by}, \text{ where }\mu \triangleq \max\left\{\normtwo{\nabla F(\bx^*)}, \normtwo{\nabla F(\bx^*)^{-1}}\right\}.
\end{equation}
Since $1=\normtwo{\bI} = \normtwo{F(\bx^*)\cdot F(\bx^*)^{-1}}\leq \normtwo{F(\bx^*) } \normtwo{F(\bx^*)^{-1}} $, we have $\mu>1$.
On the other hand, since $\rho < \zeta$, there exists a sufficiently small $\gamma > 0$ satisfying
\begin{equation}\label{equation:inex_newton_linear01}
\rho (1 + \mu \gamma) + 2 \mu \gamma \leq  \frac{\zeta}{1 +  \mu\gamma}
\quad\implies\quad 
(1 +  \mu\gamma)[\rho (1 + \mu \gamma) + 2 \mu \gamma] \leq \zeta.
\end{equation}
Since $F(\bx)$ is continuously differentiable in the neighborhood of $\bx^*$ by Assumption~\ref{assumption:inex_newton},  choose $\epsilon > 0$ sufficiently small such that  $\normtwo{\by - \bx^*} \leq \mu^2 \epsilon$, whence we have
\begin{align}
\normtwo{\nabla F(\by) - \nabla F(\bx^*)} &\leq \gamma, \label{equation:inex_newton_linear1}\\
\normtwo{\nabla F(\by)^{-1} - \nabla F(\bx^*)^{-1}} &\leq \gamma,\label{equation:inex_newton_linear2}\\
\normtwo{F(\by) - F(\bx^*) - \nabla F(\bx^*)(\by - \bx^*)} &\leq \gamma \normtwo{\by - \bx^*}, \label{equation:inex_newton_linear3}
\end{align}
where the continuity of $\nabla F(\cdot)^{-1}$ follows from Theorem~\ref{theorem:invfunctheo_var}.

Let $\normtwo{\bx^{(1)} - \bx^*} \leq \epsilon$. We now prove \eqref{equation:inex_newton_linear_res} by induction. By using \eqref{equation:inex_newton_linear0} and assumption of the induction, we have
$$
\begin{aligned}
\normtwobig{\bx^\toptzero - \bx^*}
&\leq \mu \normbig{\bx^\toptzero - \bx^*}_* \leq \mu\zeta^t  \normbig{\bx^{(1)} - \bx^*}_* 
\leq \mu^2\zeta^t \normtwobig{\bx^{(1)} - \bx^*} \leq \mu^2 \epsilon.
\end{aligned}
$$
Then, when $\by\triangleq \bx^\toptzero$, \eqref{equation:inex_newton_linear1}$\sim$\eqref{equation:inex_newton_linear3} hold. Given the update rule $\bx^\toptone \leftarrow \bx^\toptzero +\bd_{\text{in}}^\toptzero$ and the assumption that $F(\bx^*)=\bzero$, it follows that 
$$
\small
\begin{aligned}
\nabla &F(\bx^*)(\bx^\toptone - \bx^*) 
= \nabla F(\bx^*)\left\{\bx^\toptzero - \bx^* - \nabla F(\bx^\toptzero)^{-1} F(\bx^\toptzero) + \nabla F(\bx^\toptzero)^{-1} \br^\toptzero\right\} \\
&= \nabla F(\bx^*) \nabla F(\bx^\toptzero)^{-1} \left\{\nabla F(\bx^\toptzero)(\bx^\toptzero - \bx^*) - F(\bx^\toptzero) + \br^\toptzero\right\} \\
&= \left\{\bI + \nabla F(\bx^*) \big[\nabla F(\bx^\toptzero)^{-1} - \nabla F(\bx^*)^{-1}\big]\right\} \\
&\quad\times \left\{\br^\toptzero + \big[\nabla F(\bx^\toptzero) - \nabla F(\bx^*)\big] (\bx^\toptzero - \bx^*) 
- \big[F(\bx^\toptzero) - F(\bx^*) - \nabla F(\bx^*) (\bx^\toptzero - \bx^*)\big]\right\},
\end{aligned}
$$
Taking norms of the above inequality, and using Cauchy-Schwartz inequality, triangle inequality, \eqref{equation:inex_newton_linear0}, \eqref{equation:inex_newton_linear1}, \eqref{equation:inex_newton_linear2}, \eqref{equation:inex_newton_linear3}, and the assumption that $\normtwobig{\br^\toptzero} \leq \rho_t \normtwobig{F(\bx^\toptzero)}$, we obtain
\begin{equation}\label{equation:inex_newton_linear4}
\small
\begin{aligned}
&\normbig{\bx^\toptone - \bx^\toptzero}_* \leq \left\{1 + \normtwo{\nabla F(\bx^*)}\normtwo{\nabla F(\bx^\toptzero)^{-1} - \nabla F(\bx^*)^{-1}}\right\}\\ 
& \times 
\left\{\normtwobig{\br^\toptzero} + \normtwo{\nabla F(\bx^\toptzero) - \nabla F(\bx^*)} \normtwobig{\bx^\toptzero - \bx^*} + \normtwo{F(\bx^\toptzero) - F(\bx^*) - \nabla F(\bx^*)(\bx^\toptzero - \bx^*)}\right\} \\
&\leq (1 + \mu \gamma) \left\{\rho_t \normtwobig{F(\bx^\toptzero)} + 2\gamma \normtwobig{\bx^\toptzero - \bx^*} \right\}.
\end{aligned}
\end{equation}
Additionally, given that $F(\bx^*)=\bzero$, it follows that 
$$
F(\bx^\toptzero) = \left\{\nabla F(\bx^*)(\bx^\toptzero - \bx^*)\right\} + \left\{F(\bx^\toptzero) - F(\bx^*) - \nabla F(\bx^*)(\bx^\toptzero - \bx^*)\right\},
$$
whence we have
$$
\normtwobig{F(\bx^\toptzero)} \leq \normbig{\bx^\toptzero - \bx^*}_* + \gamma \normtwobig{\bx^\toptzero - \bx^*}.
$$
Substituting the above inequality into \eqref{equation:inex_newton_linear4} and using \eqref{equation:inex_newton_linear0}, \eqref{equation:inex_newton_linear01}, and the fact that $\mu>1$ yield
$$
\begin{aligned}
\normbig{\bx^\toptone - \bx^*}_* 
&\leq (1 + \mu \gamma) \left\{\rho_t \big[\normbig{\bx^\toptzero - \bx^*}_* + \gamma \normtwobig{\bx^\toptzero - \bx^*}\big] + 2 \gamma \normtwobig{\bx^\toptzero - \bx^*}\right\} \\
&\leq (1 + \mu \gamma) \left\{\rho (1 + \mu \gamma) + 2 \mu \gamma\right\} \normbig{\bx^\toptzero - \bx^*}_*
\leq \zeta \normbig{\bx^\toptzero - \bx^*}_*.
\end{aligned}
$$
This completes the proof.
\end{proof}

Clearly, the convergence of inexact Newton's method depends on the choice of the relative error $\normtwobig{\br^\toptzero} \leq \rho_t \normtwobig{F(\bx^\toptzero)}$. Intuitively, the more accurately the Newton equation is solved, the better the convergence properties of  inexact Newton's method. 
In fact, different conditions on this relative error bound yields different local convergence results, which we provide in the following and the details can be found in \citet{fletcher2000practical, sun2006optimization}.
\begin{theoremHigh}[Convergence of Inexact Newton's Methods]\label{theorem:conv_inex_newton}
Let the assumptions of Theorem~\ref{theorem:inex_newton_linear} be satisfied. 
If the starting point $\bx^{(1)}$ is sufficiently near $\bx^*$, the sequence $\{\bx^\toptzero\}$ generated by inexact Newton's method (Algorithm~\ref{alg:inexac_newton}) converges:
\begin{enumerate}[(i)]
\item If there exists a constant $\rho < 1$ such that $\rho_t$ satisfies $0 < \rho_t < \rho$ for $t = 1, 2, \ldots$, then the algorithm converges \textbf{linearly} (the result of Theorem~\ref{theorem:inex_newton_linear}).

\item If $\lim_{t \to \infty} \rho_t = 0$, then the algorithm converges \textbf{superlinearly}.

\item If $\rho_t = \mathcalO(\normtwobig{\nabla f(\bx^\toptzero)})$, then the algorithm converges \textbf{quadratically}.
\end{enumerate}
\end{theoremHigh}

The intuitive implication of Theorem~\ref{theorem:conv_inex_newton} is that achieving better convergence requires solving the Newton equation more accurately. 
In general iterative methods, the stopping criterion usually depends on the size of the relative error. 
According to the first statement of Theorem~\ref{theorem:conv_inex_newton}, setting the relative error to a fixed value ensures convergence. Compared to pure Newton's method, an inexact Newton's method with a fixed error only achieves linear convergence but may perform better on ill-conditioned problems than traditional gradient methods. To achieve quadratic convergence with  inexact Newton's method, it becomes necessary to solve the Newton equation with high accuracy in later iterations, essentially aligning it with  Newton's method.

A commonly employed variant of  inexact Newton's method is the \textit{Newton-Conjugate Gradient method}, which utilizes the conjugate gradient method (Section~\ref{section:conjugate-descent}) to solve the Newton equation \eqref{equation:newton_secon_approx_eq3}. 
Given the conjugate gradient method's effectiveness in solving linear systems, it often requires only a few steps (sometimes just one step) to meet the criteria outlined in the first conclusion of Theorem~\ref{theorem:conv_inex_newton}. The Newton-Conjugate Gradient method demonstrates good numerical performance across many problems and serves as an essential optimization tool for tackling various optimization challenges.

\subsection{Proximal Newton Method}
We described the proximal gradient method in \eqref{equation:prox_decom} using the proximal operators. 
By combining the principles of the proximal gradient method with the derivation of Newton's method in \eqref{equation:newton_secon_approx_ALL}, we introduce the \textit{proximal Newton method}. The update at the $t$-th iteration is given by: 
\begin{subequations}\label{equation:prox_new_updv1}
\begin{align}
\bd^\toptzero &\leftarrow \argmin_{\bd} \left\{G_t(\bd) \triangleq g_t(\bd) + g(\bx^\toptzero +\bd)\right\};\\
\bx^\toptone &\leftarrow \bx^\toptzero + \eta_t \bd^\toptzero,
\end{align}
\end{subequations}
where $g_t(\bd) \triangleq f(\bx^\toptzero) + \nabla f(\bx^\toptzero)^\top \bd + \frac{1}{2} \bd^\top \bH^\toptzero \bd$, $\bH^\toptzero\triangleq\nabla^2 f(\bx^\toptzero )$, $\eta_t$ is the stepsize, and $g(\bx)$ is a \textbf{proper closed and convex} function~\footnote{A few examples are discussed in the paragraph below Definition~\ref{definition:opt_probs_all}.}.
Thus, when $\bH^\toptzero $ is positive definite, $G_t(\bd)$ is also convex.

Previously, we defined the proximal operator using the $\ell_2$ norm (Definition~\ref{definition:projec_prox_opt}). If we replace the $\ell_2$ norm with the $\bQ$-norm (see \eqref{equation:q_norm}, $\norm{\bx}^2_{\bQ} = \bx^\top\bQ\bx = \normtwo{{\bQ}^{1/2}\bx}^2$ for any $\bx$), we obtain the \textit{scaled proximal operator}:
\begin{equation}\label{equation:scaled_prox_opt}
\textbf{(Scaled proximal)}:\quad 
\prox_{f, \bQ}(\by) \triangleq 
\mathop{\argmin}_{\bx\in\real^n} \left( f(\bx) + \frac{1}{2}\norm{\bx-\by}_{\bQ}^2 \right).
\end{equation}
Then, the proximal Newton update in \eqref{equation:prox_new_updv1} can be equivalently stated as follows:
\begin{tcolorbox}[colback=white,colframe=black]
\begin{minipage}{1\textwidth}
\begin{subequations}\label{equation:prox_new_updv2}
\small
\begin{align}
\bu^\toptzero &\leftarrow \argmin_{\bu} \left\{\nabla f(\bx^\toptzero)^\top \bu
+ \frac{1}{2} (\bu-\bx^\toptzero)^\top \bH^\toptzero (\bu-\bx^\toptzero) + g(\bu)\right\}\\
&= \argmin_{\bu}  \frac{1}{2} \norm{\bx^\toptzero - \bH^{(t)-1} \nabla f(\bx^\toptzero) - \bu }_{\bH^\toptzero}^2 
+g(\bu)  \\
&= \prox_{f, \bH^\toptzero}\Big( \bx^\toptzero - \bH^{(t)-1} \nabla f(\bx^\toptzero) \Big);\\
\bx^\toptone &\leftarrow \bx^\toptzero + \eta_t (\bu^\toptzero - \bx^\toptzero).
\end{align}
\end{subequations}
\end{minipage}
\end{tcolorbox}
When $f$ is $\alpha$-strongly convex and $\nabla^2 f$ is $L$-Lipschitz, it can be proved that
$$
\normtwo{\bx^\toptone - \bx^*} \leq \frac{L}{2\alpha} \normtwo{\bx^\toptzero - \bx^*}^2
$$
when $\bx^\toptzero$ is sufficiently close to $\bx^*$ (indicating  a local quadratic convergence result). We shall not go into the details and more information can be found in \citet{lee2014proximal}.

\section{Quasi-Newton Methods}\label{section:quasi_newton_method}

Pure Newton's method achieves excellent results both theoretically and practically. However, for large-scale problems, computing the Hessian matrix can be particularly expensive or difficult to obtain. Even if we manage to compute the Hessian matrix, solving a large-scale linear system remains a challenge. 
A natural question arises: Can we solve a large-scale linear system using an approximation of the Hessian matrix or its inverse in Newton's method?
\textit{Quasi-Newton methods} (from Latin, quasi means nearly) are algorithms that generate approximate matrices at each step with a lower  computational cost while still ensuring the sequence of iterations exhibits superlinear convergence properties.

Quasi-Newton methods do not compute the Hessian matrix $\nabla^2 f(\bx^\toptzero)$ directly but construct an approximate matrix $\bH^\toptzero$ or its inverse $\bZ^\toptzero$. We hope that $\bH^\toptzero$ or $\bZ^\toptzero$ retains some properties of the Hessian matrix, such as ensuring that $\bd^\toptzero$ derived from the Newton equations using these approximations remains a descent direction. 

\index{Secant equation}
\subsection{Secant Equations}

First, let us review the derivation of Newton's method. Suppose $f(\bx)$ is a twice continuously differentiable function. According to linear approximation theorem (Theorem~\ref{theorem:linear_approx}), the gradient function $\nabla f(\bx)$ can be approximated near $\bx^\toptone$ as:
$$
\nabla f(\bx) = \nabla f(\bx^\toptone) + \nabla^2 f(\bx^\toptone)(\bx - \bx^\toptone) + \mathcalO\big(\normtwobig{\bx - \bx^\toptone}^2\big).
$$
Let $\bx \triangleq \bx^\toptzero$, $\bh^\toptzero \triangleq \bx^\toptone - \bx^\toptzero$, and $\by^\toptzero \triangleq \nabla f(\bx^\toptone) - \nabla f(\bx^\toptzero)$. That is, $\bh^\toptzero$ denotes the descent step and $\by^\toptzero$ represents the increase in the gradient.
Then,
$$
\nabla^2 f(\bx^\toptone) \bh^\toptzero + \mathcalO(\normtwobig{\bh^\toptzero}^2) = \by^\toptzero.
$$
Ignoring the higher-order term $\normtwobig{\bh^\toptzero}^2$, we hope that the approximate Hessian matrix $\bH^\toptone$ satisfies the equation:
\begin{subequations}
\begin{equation}\label{equation:secant1}
	\textbf{(SE1)}:\qquad \by^\toptzero = \bH^\toptone \bh^\toptzero,
\end{equation}
or its inverse approximate matrix $\bZ^\toptone$ satisfies the equation:
\begin{equation}\label{equation:secant2}
	\textbf{(SE2)}:\qquad \bh^\toptzero = \bZ^\toptone \by^\toptzero,
\end{equation}
\end{subequations}
Equations~\eqref{equation:secant1} and \eqref{equation:secant2} are known as the \textit{secant equations}.

We can also understand the secant equations \eqref{equation:secant1} from another perspective. 
Newton's method essentially performs a second-order approximation of the objective function $f(\bx)$ at the iterate point $\bx^\toptzero$. Consider the second-order approximation at $\bx^\toptone$ (Theorem~\ref{theorem:quad_app_theo}):
$$
f(\bx^\toptone +\bh) \approx\psi_{t+1}(\bh) \triangleq f(\bx^\toptone) + \nabla f(\bx^\toptone)^\top \bh + \frac{1}{2} \bh^\top \bH^\toptone \bh.
$$
A natural requirement is that the gradient of $\psi_{t+1}(\bh)$ at $\bh = -\bh^\toptzero$ and $\bh = \bzero$ should match the gradient of $f(\bx)$ at $\bx = \bx^\toptzero$ and $\bx = \bx^\toptone$, respectively:
\begin{equation}
\nabla \psi_{t+1}(-\bh^\toptzero) = \nabla f(\bx^\toptzero)
\qquad\text{and}\qquad
\nabla \psi_{t+1}(\bzero) = \nabla f(\bx^\toptone).
\end{equation}
Note that the latter equality $\nabla \psi_{t+1}(\bzero) = \nabla f(\bx^\toptone)$ is naturally satisfied. To ensure the first condition is satisfied, we need:
$$
\nabla \psi_{t+1}(-\bh^\toptzero) = \nabla f(\bx^\toptone) - \bH^\toptone \bh^\toptzero = \nabla f(\bx^\toptzero).
$$
Rearranging terms yields the (SE1) equation \eqref{equation:secant1}.

Additionally, note that the positive definiteness of the approximate matrix $ \bH^\toptzero $ is a crucial factor. By premultiplying both sides of the  (SE1) equation \eqref{equation:secant1} by $ \bh^\toptzeroTOP $, we get $ \bh^\toptzeroTOP \bH^\toptone \bh^\toptzero = \bh^\toptzeroTOP \by^\toptzero $, thus the condition
\begin{subequations}
\begin{equation}\label{equation:secant1_2}
\textbf{(SE1$'$)}:\qquad \innerproduct{\bh^\toptzero,  \by^\toptzero} > 0
\end{equation}
is a necessary condition for $ \bH^\toptone $ to be positive definite. This condition must hold throughout the iteration process; it is known as the curvature condition.  
For general objective functions $ f(\bx) $, we need to use Wolfe line search conditions to ensure the curvature condition \eqref{equation:secant1_2}. In practice, according to the second Wolfe condition \eqref{equation:wolfe_2}, we have $ \innerproduct{\nabla f(\bx^\toptone),  \bh^\toptzero} \geq c_2 \innerproduct{\nabla f(\bx^\toptzero), \bh^\toptzero} $, where $0<c_2<1$ is a constant. Subtracting $ \nabla f(\bx^\toptzero)^\top \bh^\toptzero $ from both sides yields that
\begin{equation}\label{equation:secant1_3}
\textbf{(SE1$''$)}:\qquad \innerproduct{\by^\toptzero,  \bh^\toptzero} \geq (c_2 - 1)  \innerproduct{\nabla f(\bx^\toptzero), \bh^\toptzero} > 0,
\end{equation}
since $0< c_2 < 1 $ and $ \bh^\toptzero = \eta_t \bd^\toptzero $ is a descent direction. 
\end{subequations}

In most cases, the approximate matrices $ \bH^\toptone $ or $ \bZ^\toptone $ are obtained by adding a correction to the previous iteration and requiring them to satisfy the secant equation \eqref{equation:secant1}. Below provides a general framework for quasi-Newton methods (Algorithm~\ref{alg:quasi_newton_frame}), and the next subsection will discuss specific matrix update methods. 
An overview comparing the quasi-Newton method and Newton's method is presented in  Table~\ref{tab:comparison_quasi_new}.

\index{Quasi-Newton method}
\begin{algorithm}[H]
\caption{Quasi-Newton Method}
\label{alg:quasi_newton_frame}
\begin{algorithmic}[1]
\Require A twice continuously differentiable function $f(\bx)$
\State {\bfseries Input:}  Initialize $ \bx^{(1)} \in \real^n $, $ \bH^{(1)} \in \real^{n \times n} $ (or $ \bZ^{(1)} $);
\For{$t=1,2,\ldots$}
\State Compute direction $ \bd^\toptzero = -(\bH^\toptzero)^{-1} \nabla f(\bx^\toptzero) $ or $ \bd^\toptzero = -\bZ^\toptzero \nabla f(\bx^\toptzero) $.
\State Use line search method to determine a suitable stepsize $ \eta_t > 0 $;
\State $ \bx^\toptone = \bx^\toptzero + \bh^\toptzero $, where $\bh^\toptzero\triangleq \eta_t \bd^\toptzero$;
\State Update approximate Hessian matrix $ \bH^\toptone $ or its inverse $ \bZ^\toptone $;
\EndFor
\State {\bfseries Return:}   final $\bx\leftarrow \bx^{(t)}$;
\end{algorithmic}
\end{algorithm}

In practical applications, quasi-Newton methods based on $ \bZ^\toptzero $ are more commonly used because computing the descent direction $ \bd^\toptzero $ does not require solving a linear system due to the inverse of $\bH^\toptzero$, while solving a linear system is very time-consuming for large-scale problems. 
However, quasi-Newton methods that utilize $ \bH^\toptzero $ possess superior theoretical properties, leading to a more stable sequence of iterations. 
If an efficient method for solving the linear system becomes available, employing quasi-Newton methods based on $ \bH^\toptzero $ could also be feasible.

\begin{table}[h]
\centering
\caption{Comparison between quasi-Newton method and Newton's method.}
\label{tab:comparison_quasi_new}
\begin{tabular}{|p{6.8cm}|p{6.8cm}|}
\hline
\textbf{Quasi-Newton method} & \textbf{Newton's method} \\ \hline\hline
Requires only function values and gradients & Requires function values, gradients, and Hessians \\ \hline
$\{\bH^\toptzero\}$ maintains positive definiteness  for several updates & $\{\nabla^2 (\bx^\toptzero)\}$ may not be positive definite \\ \hline
Needs $\mathcalO(n^2)$ multiplications per  iteration & Needs $\mathcalO(n^3)$ multiplications per  iteration \\ \hline
\end{tabular}
\end{table}
\subsection{Updating Quasi-Newton Matrices}\label{section:quasi_new_update}

This subsection introduces some common methods for updating quasi-Newton matrices.

\index{Symmetric rank-one update}
\subsection*{Symmetric Rank-One Update (SRO)}

The \textit{symmetric rank-one (SRO)} update formula  is the simplest approach for updating quasi-Newton matrices. Suppose $ \bH^\toptzero $ is the approximate Hessian matrix at step $ t $. We update $ \bH^\toptzero $ to $ \bH^\toptone $ by applying  a symmetric rank-one correction to satisfy the secant equation \eqref{equation:secant1}. 
To be more specific, we update by the rank-one correction:
$$
\bH^\toptone \leftarrow \bH^\toptzero + \gamma \bu \bu^\top,
$$
where $ \bu \in \real^n $ and $ \gamma \in \real $ are variables to be determined. Based on the (SE1) equation \eqref{equation:secant1},
$
\bH^\toptone \bh^\toptzero = (\bH^\toptzero + \gamma \bu \bu^\top) \bh^\toptzero = \by^\toptzero,
$
we obtain
\begin{equation}\label{equation:sro_rawequa}
(\gamma \cdot \bu^\top \bh^\toptzero) \bu = \by^\toptzero - \bH^\toptzero \bh^\toptzero.
\end{equation}
Notice that $ \gamma \cdot \bu^\top \bh^\toptzero $ is a scalar, so $ \bu $ and $ \by^\toptzero - \bH^\toptzero \bh^\toptzero $ must have the same direction. 
Letting $ \bu \triangleq \by^\toptzero - \bH^\toptzero \bh^\toptzero $, substituting into the original equation \eqref{equation:sro_rawequa} gives
$$
\begin{aligned}
&\gamma \left(\left(\by^\toptzero - \bH^\toptzero \bh^\toptzero\right)^\top \bh^\toptzero\right) (\by^\toptzero - \bH^\toptzero \bh^\toptzero) = \by^\toptzero - \bH^\toptzero \bh^\toptzero\\
&\quad\implies\quad 
\gamma = \frac{1}{\left(\by^\toptzero - \bH^\toptzero \bh^\toptzero\right)^\top \bh^\toptzero},
\quad\text{assuming $ \left(\by^\toptzero - \bH^\toptzero \bh^\toptzero\right)^\top \bh^\toptzero \neq 0 $}.
\end{aligned}
$$
The update formula becomes
\begin{subequations}
\begin{equation}\label{equation:sro_rawequa1}
\textbf{(SRO1)}:\qquad 
\bH^\toptone \leftarrow \bH^\toptzero + 
\frac{\left(\by^\toptzero - \bH^\toptzero \bh^\toptzero\right) \left(\by^\toptzero - \bH^\toptzero \bh^\toptzero\right)^\top}{\left(\by^\toptzero - \bH^\toptzero \bh^\toptzero\right)^\top \bh^\toptzero},
\end{equation}
which is known as the \textit{SRO update formula for $ \bH^\toptzero $} \citep{broyden1965class}. 
By following an analogous process, we can derive the \textit{SRO update formula for $ \bZ^\toptzero $} based on the second secant equation \eqref{equation:secant2}:
\begin{equation}\label{equation:sro_rawequa2}
\textbf{(SRO2)}:\qquad 
\bZ^\toptone \leftarrow \bZ^\toptzero + \frac{\left(\bh^\toptzero - \bZ^\toptzero \by^\toptzero\right) \left(\bh^\toptzero - \bZ^\toptzero \by^\toptzero\right)^\top}{\left(\bh^\toptzero - \bZ^\toptzero \by^\toptzero\right)^\top \by^\toptzero}.
\end{equation}
\end{subequations}
We observe that \eqref{equation:sro_rawequa1} and \eqref{equation:sro_rawequa2} are dual to each other in form. 

A distinctive property of the SRO update is its natural quadratic termination. Specifically, for a quadratic function, it does not require line search and can terminate within $n+1$ steps, i.e., $\bZ^{(n+1)} =\bA^{-1}$, where $\bA$ is the Hessian of the
quadratic function. This fact is proved by the following theorem.
\begin{theorem}[Quadratic Termination of SRO Update]
Let $f(\bx)=\frac{1}{2}\bx^\top\bA\bx-\bb^\top\bx+c$ where $\bA\in\real^{n\times n}$ is positive definite.
Suppose the iterates of descent steps $ \bh^{(1)}, \bh^{(2)}, \ldots, \bh^{(n)} $ generated by the SRO method are linearly independent. Then, the  SRO method terminates at $ n+1 $ steps, meaning $ \bZ^{(n+1)} = \bA^{-1} $.
\end{theorem}

\begin{proof}
Since Hessian $ \bA $ is positive definite, we have 
\begin{equation}\label{equation:sro_indu_eq1}
\by^{(t)} = \bA \bh^{(t)}, \quad  t = 1,2, \ldots, n.
\end{equation}
We prove by the following two steps.
\paragraph{Induction result.}
First, by induction, we prove the following property:
\begin{equation}\label{equation:sro_indu_eq2}
\bZ^{(i+1)} \by^{(k)} = \bh^{(k)}, \quad  k =1, 2, \ldots, i.
\end{equation}
For $ i = 1 $, the result holds trivially from \eqref{equation:sro_rawequa2}. Now suppose it is true for $ i > 1 $; we will show it holds for $ i + 1 $.
From \eqref{equation:sro_rawequa2}, we have
\begin{equation}\label{equation:sro_indu_eq3}
\bZ^{(i+1)} \by^{(k)} = \bZ^{(i)} \by^{(k)} + \frac{(\bh^{(i)} - \bZ^{(i)} \by^{(i)})(\bh^{(i)} - \bZ^{(i)} \by^{(i)})^\top \by^{(k)}}{(\bh^{(i)} - \bZ^{(i)} \by^{(i)})^\top \by^{(i)}}.
\end{equation}
When $ k < i $, from the induction hypothesis and \eqref{equation:sro_indu_eq1}, we have
$$
\begin{aligned}
(\bh^{(i)} - \bZ^{(i)} \by^{(i)})^\top \by^{(k)} 
&= (\bh^{(i)})^\top \by^{(k)} - (\by^{(i)})^\top \bZ^{(i)} \by^{(k)} \\
&= (\bh^{(i)})^\top \by^{(k)} - (\by^{(i)})^\top \bh^{(k)}
= (\bh^{(i)})^\top \bA \bh^{(k)} - (\bh^{(i)})^\top \bA \bh^{(k)}
= 0.
\end{aligned}
$$
Substituting into \eqref{equation:sro_indu_eq3} yields that 
$$
\bZ^{(i+1)} \by^{(k)} = \bZ^{(i)} \by^{(k)} = \bh^{(k)}, \; \forall k < i.
$$
The equality for $ k = i $ follows directly from \eqref{equation:sro_rawequa2} that
$
\bZ^{(i+1)} \by^{(i)} = \bh^{(i)}.
$
Therefore, \eqref{equation:sro_indu_eq2} holds for $i+1$.

\paragraph{The result.}
The above induction proof shows that 
$$
\bh^{(k)} = \bZ^{(n+1)} \by^{(k)} = \bZ^{(n)} \bA \bh^{(k)}, \; k =  1, 2, \ldots, n.
$$
Since $ \bh^{(k)} \; (k = 1, 2, \ldots, n) $ are linearly independent, then $ \bZ^{(n+1)} \bA = \bI $, which implies $ \bZ^{(n+1)} = \bA^{-1} $. 
\end{proof}

The above theorem highlights several characteristics of the SRO update:
\begin{enumerate}[(i)]
\item The SRO update possesses natural quadratic termination.
\item The SRO update satisfies the  property: $ \bZ^{(t+1)} \by^{(k)} = \bh^{(k)}, k \leq  t $.
\item The SRO update does not retain the positive definiteness of $ \bZ^\toptzero $. It retains positive definiteness if and only if $ \innerproduct{\bh^\toptzero  - \bZ^\toptzero \by^\toptzero, \by^\toptzero} > 0 $ (see \eqref{equation:sro_rawequa2}). 
However, this condition is hard to guarantee. A solution is to use the SRO update within a trust region framework, as the trust region method does not require positive definiteness of the Hessian approximations.
\item Sometimes, the denominator $ \innerproduct{\bh^\toptzero  - \bZ^\toptzero \by^\toptzero, \by^\toptzero}  $ is very small or zero, leading to serious numerical difficulties or even algorithm failure. This drawback limits its applications.  A special skipping strategy to prevent the SRO update from breaking down involves using \eqref{equation:sro_rawequa2} only if:
\begin{equation}\label{equation:sro_skip}
\abs{\innerproduct{\bh^{(t)} - \bZ^{(t)} \by^{(t)},  \by^{(t)}}} \geq \zeta \normtwo{\bh^{(t)} - \bZ^{(t)} \by^{(t)}} \normtwobig{\by^{(t)}},
\end{equation}
where $ \zeta \in (0, 1) $; otherwise we set $ \bZ^\toptone  \leftarrow \bZ^{(t)} $.
\item The SRO update has a good behavior that it continues to generate good Hessian approximations, which is stated in the following theorem.
\end{enumerate}
\begin{theorem}[\citep{sun2006optimization}]
Let $ f $ be twice continuously differentiable, and its Hessian be bounded and Lipschitz continuous in a neighborhood of a point $ \bx^* $. Let $ \{\bx^\toptzero\} $ be a sequence of iterates converging to $\bx^* $. Suppose that the skipping rule \eqref{equation:sro_skip} holds for all $ t $, and the steps $ \bh^\toptzero  $ are uniformly linearly independent. Then the matrix sequence $ \{\bH^\toptzero\} $ generated by SRO update satisfies
\begin{equation}
\lim_{t \to \infty} \normtwo{\bZ^{(t)} - [\nabla^2 f(\bx^*)]^{-1}} = 0.
\end{equation}
\end{theorem}


\subsection*{BFGS Formula}

To address the limitations of the SRO formula, we now consider a symmetric rank-two update for $ \bH^\toptzero $. Specifically, we define
$$
\bH^\toptone = \bH^\toptzero + \gamma \bu \bu^\top + \zeta \bv \bv^\top,
$$
where $ \bu, \bv \in \real^n $ and $ \gamma, \zeta \in \real $ are parameters  to be determined. Based on the (SE1) equation \eqref{equation:secant1},
$$
\bH^\toptone \bh^\toptzero = (\bH^\toptzero + \gamma \bu \bu^\top + \zeta \bv \bv^\top) \bh^\toptzero = \by^\toptzero,
$$
we derive
$$
(\gamma \cdot \bu^\top \bh^\toptzero) \bu + (\zeta \cdot \bv^\top \bh^\toptzero) \bv = \by^\toptzero - \bH^\toptzero \bh^\toptzero.
$$
By selecting  $ \bu $ and $ \bv $ appropriately,  this equation can be satisfied. A direct approach is to equate terms on both sides:
$$
\begin{aligned}
\bu &\triangleq \by^\toptzero, \qquad &\gamma& \cdot \bu^\top \bh^\toptzero \triangleq 1, \\
\bv &\triangleq \bH^\toptzero \bh^\toptzero, \qquad &\zeta& \cdot \bv^\top \bh^\toptzero \triangleq -1.
\end{aligned}
$$
Consequently, the update rule for $\bH^\toptone$ is
\begin{subequations}
\begin{equation}\label{equation:bfgs_update}
	\textbf{(BFGS1)}:\qquad 
	\bH^\toptone = \bH^\toptzero + \frac{\by^\toptzero \by^\toptzeroTOP}{\bh^\toptzeroTOP \by^\toptzero} - \frac{\bH^\toptzero \bh^\toptzero (\bH^\toptzero \bh^\toptzero)^\top}{\bh^\toptzeroTOP \bH^\toptzero \bh^\toptzero}.
\end{equation}
This formulation is known as the \textit{BFGS formula for $ \bH^\toptzero $}, named after \citet{broyden1970convergence}, \citet{fletcher1970new}, \citet{goldfarb1970family}, and \citet{shanno1970conditioning}.

Using the SMW formula~\footnote{Sherman-Morrison formula: $(\bA+\bb\bc^\top)^{-1} = \bA^{-1} - \frac{\bA^{-1}\bb\bc^\top\bA^{-1}}{1+\bc^\top\bA^{-1}\bb}$. Sherman-Morrison-Woodbury formula generalizes this by: $(\bA+\bH\bC)^{-1} = \bA^{-1}-\bA^{-1}\bH(\bI+\bC\bA^{-1}\bH)^{-1}\bC\bA^{-1}$, which is derived based on Schur complement; see, for example, \citet{lu2021numerical}.}  and assuming $ \bZ^\toptzero = (\bH^\toptzero)^{-1} $, we obtain the  \textit{BFGS formula for $ \bZ^\toptzero $}:
\begin{equation}\label{equation:bfgs2_update}
	\textbf{(BFGS2)}:\quad 
	\bZ^\toptone = \bigg(\bI - \frac{ \bh^\toptzero\by^\toptzeroTOP}{\bh^\toptzeroTOP \by^\toptzero}\bigg) \bZ^\toptzero 
	\bigg(\bI -  \frac{\by^\toptzero \bh^\toptzeroTOP}{\bh^\toptzeroTOP \by^\toptzero}\bigg) +  \frac{\bh^\toptzero \bh^\toptzeroTOP}{\bh^\toptzeroTOP \by^\toptzero}.
\end{equation}
\end{subequations}
For the BFGS formula to guarantee that $ \bZ^\toptone $ remains positive definite, it suffices that inequality \eqref{equation:secant1_2} holds, and $ \bZ^\toptzero $ is positive definite. 
If condition \eqref{equation:secant1_2} does not hold, line search using the Wolfe conditions can enforce this requirement.

The BFGS formula stands out as one of the most effective quasi-Newton methods currently available, combining strong theoretical properties with straightforward implementation.
The BFGS update formula (BFGS2) has additional significance as it satisfies a form of approximate optimality. 
In particular,  $\bZ^\toptone$ defined by equation (BFGS2), solves the following optimization problem:
\begin{equation}\label{equation:opt_bfgs}
\begin{aligned}
& \min_{\bZ} & & \norm{\bZ - \bZ^\toptzero}_{\bW}, \\
& \text{s.t.} & & \bZ = \bZ^\top, \\
& & & \bZ \by^\toptzero = \bh^\toptzero.
\end{aligned}
\end{equation}
This optimization seeks the matrix $\bZ$ that is closest to $\bZ^\toptzero$ among all symmetric matrices satisfying the (SE2) equation \eqref{equation:secant2}. Here, $\norm{\cdot}_{\bW}$ denotes  the weighted norm, defined as
$
\norm{\bZ}_{\bW} \triangleq \normf{\bW^{1/2} \bZ \bW^{1/2}},
$
where $\bW$  is any matrix fulfilling the secant equation $\bW \bh^\toptzero = \by^\toptzero$.

\subsection*{DFP Formula}

In the derivation of the BFGS formula, if we use the (SE2) equation \eqref{equation:secant2} to update $ \bZ^\toptzero $ using  a rank-two correction, we arrive at what is known as the \textit{DFP formula for $ \bZ^\toptzero $}:
\begin{subequations}
\begin{equation}\label{equation:dfp1}
	\textbf{(DFP1)}:\qquad 
	\bZ^\toptone = \bZ^\toptzero - \frac{(\bZ^\toptzero \by^\toptzero )(\bZ^\toptzero \by^\toptzero)^\top}{\by^\toptzeroTOP \bZ^\toptzero \by^\toptzero} + \frac{\bh^\toptzero \bh^\toptzeroTOP}{\by^\toptzeroTOP \bh^\toptzero}.
\end{equation}
This iterative format was first discovered by \citet{davidon1991variable}, later refined by \citet{fletcher1963rapidly}, hence it is named the DFP formula. 
Utilizing  the SMW formula, we can derive the corresponding formula for $ \bH^\toptzero $:
\begin{equation}\label{equation:dfp2}
	\textbf{(DFP2)}:\qquad 
	\bH^\toptone = \bigg(\bI - \frac{ \by^\toptzero\bh^\toptzeroTOP}{\bh^\toptzeroTOP \by^\toptzero}\bigg) \bH^\toptzero 
	\bigg(\bI - \frac{\bh^\toptzero \by^\toptzeroTOP}{\bh^\toptzeroTOP \by^\toptzero}\bigg) +   \frac{\by^\toptzero \by^\toptzeroTOP}{\bh^\toptzeroTOP \by^\toptzero}.
\end{equation}
\end{subequations}

It is evident that the DFP formulas (DFP1)--(DFP2) and the BFGS formulas (BFGS1)--(BFGS2) are dual to each other. By replacing $ \bZ^\toptzero $ with $ \bH^\toptzero $ and interchanging  $ \bh^\toptzero $ with $ \by^\toptzero $ in the BFGS formulas, we obtain the DFP formulas. Moreover, the duality also exists in the approximation properties. Specifically, the definition of  $ \bH^\toptone $ in \eqref{equation:dfp2} solves the following optimization problem:
\begin{equation}\label{equation:opt_dfp}
\begin{aligned}
&\min_{\bH} & &  \norm{\bH - \bH^\toptzero}_{\bW}\\
&\text{s.t.}& &   \bH = \bH^\top,\\
& & &\bH \bh^\toptzero = \by^\toptzero,
\end{aligned}
\end{equation}
where $ \norm{\cdot}_{\bW} $ is the weighted norm defined as
$
\norm{\bZ}_{\bW} \triangleq \normf{\bW^{1/2} \bZ \bW^{1/2}}.
$
The interpretation of  $\norm{\cdot}_{\bW}$ and the fundamental concept in equations \eqref{equation:opt_dfp} and \eqref{equation:opt_bfgs} are similar; however, $\bW$ is any matrix satisfying $\bW\by^\toptzero = \bh^\toptzero$. 

Despite the numerous similarities between the DFP and BFGS approaches in various aspects, from a practical standpoint, the DFP method is generally considered less effective than the BFGS method. Consequently, in practical applications, the BFGS approach is more commonly favored.

\subsection{Global Convergence of Quasi-Newton Methods}

This subsection introduces the convergence properties of quasi-Newton methods. First, we use Zoutendijk's condition (Theorem~\ref{theorem:zoutendijk_cond}) to obtain the global convergence of quasi-Newton methods, and then introduce the convergence rate.

\begin{theoremHigh}[Global Convergence of BFGS \citep{nocedal1999numerical}]\label{theorem:conv_bfgs_glo}
Assume the initial matrix $\bH^{(1)}$ is symmetric positive definite, the objective function $f(\bx)$ is twice continuously differentiable, and the level set
$$
\sL \triangleq \lev[f, f(\bx^{(1)})] \triangleq \{ \bx \in \real^n \mid f(\bx) \leq f(\bx^{(1)}) \}
$$
is convex. 
Furthermore, suppose there exist positive constants $\alpha$ and $\beta$ such that for any $\bz \in \real^n$ and any $\bx \in \sL$,
\begin{equation}\label{equation:bfgs_rssrsc}
\alpha \normtwo{\bz}^2 \leq \bz^\top \nabla^2 f(\bx) \bz \leq \beta \normtwo{\bz}^2.~\footnote{
That is, the eigenvalues of the Hessian of $f$ for any $\bx\in\sL$ is bounded in the interval $[\alpha, \beta]$ (Theorem~\ref{theorem:eigen_charac}).
This condition is closely related to the restricted strong convexity and restricted strong smoothness for design matrices (Definition~\ref{definition:res_scss_mat}).}
\end{equation}
Then the BFGS update formula  \eqref{equation:bfgs_update} for $\bH^\toptzero$ combined with  line search under the Wolfe condition (Definition~\ref{definition:wolfe_cond}) ensures that the quasi-Newton algorithm (Algorithm~\ref{alg:quasi_newton_frame}) globally converges to a local minimizer $\bx^*$ of $f(\bx)$.
\end{theoremHigh}
\begin{proof}[of Theorem~\ref{theorem:conv_bfgs_glo}]
For convenience, we define
\begin{equation}\label{equation:bfgs_alphabeta}
\alpha_t \triangleq \frac{\by^\toptzeroTOP \bh^\toptzero}{\bh^\toptzeroTOP \bh^\toptzero}
\qquad \text{and}\qquad 
\beta_t \triangleq \frac{\by^\toptzeroTOP \by^\toptzero}{\by^\toptzeroTOP \bh^\toptzero}.
\end{equation}
Given the assumption in  \eqref{equation:bfgs_rssrsc} and the fact that 
$
\by^\toptzero = \widetildebZ \bh^\toptzero$, with $\widetildebZ\triangleq\int_0^1 \nabla^2 f(\bx^\toptzero + \mu(\bx^\toptone - \bx^\toptzero))  d\mu$ denoting the average Hessian (Theorem~\ref{theorem:fund_theo_calculu}), it follows that 
$\alpha_t = \frac{\bh^\toptzeroTOP \widetildebZ \bh^\toptzero}{\bh^\toptzeroTOP \bh^\toptzero}$ 
and 
$\beta_t = \frac{(\bs^\toptzero)^\top \widetildebZ \bs^\toptzero}{(\bs^\toptzero)^\top \bs^\toptzero}$, where $\bs^\toptzero\triangleq \widetildebZ^{1/2}\bh^\toptzero$,
whence we have 
$$
\alpha_t \geq \alpha
\qquad \text{and}\qquad 
\beta_t \leq \beta.
$$
Taking  trace on both sides of the BFGS update formula \eqref{equation:bfgs_update} for $\bH^\toptzero$, we get
\begin{equation}\label{equation:bfgs_trace}
\trace(\bH^\toptone) = \trace(\bH^\toptzero) 
+ \frac{\normtwobig{\by^\toptzero}^2}{\by^\toptzeroTOP \bh^\toptzero}
- \frac{\normtwo{\bH^\toptzero \bh^\toptzero}^2}{\bh^\toptzeroTOP \bH^\toptzero \bh^\toptzero}.
\end{equation}
In addition, it can also be shown that 
\begin{equation}\label{equation:bfgs_det}
\det(\bH^\toptone) = \det(\bH^\toptzero) \frac{\by^\toptzeroTOP \bh^\toptzero}{\bh^\toptzeroTOP \bH^\toptzero \bh^\toptzero}.
\end{equation}
Next, we define
\begin{equation}\label{equation:bfgs_cos_zeta}
\cos (\theta_t) \triangleq \frac{\bh^\toptzeroTOP \bH^\toptzero \bh^\toptzero}{\normtwo{\bh^\toptzero} \normtwo{\bH^\toptzero \bh^\toptzero}}
\qquad \text{and}\qquad  
\quad \zeta_t \triangleq \frac{\bh^\toptzeroTOP \bH^\toptzero \bh^\toptzero}{\bh^\toptzeroTOP \bh^\toptzero},
\end{equation}
where $\theta_t$ denotes the angle between $\bh^\toptzero$ and $\bH^\toptzero \bh^\toptzero$.
Rearranging the third term on the right-hand side of \eqref{equation:bfgs_trace} yields
\begin{equation}\label{equation:bfgs_conv_bh_theta}
\frac{\normtwo{\bH^\toptzero \bh^\toptzero}^2}{\bh^\toptzeroTOP \bH^\toptzero \bh^\toptzero} 
= \frac{\normtwo{\bH^\toptzero \bh^\toptzero}^2 \normtwo{\bh^\toptzero}^2}{\big(\bh^\toptzeroTOP \bH^\toptzero \bh^\toptzero\big)^2} \frac{\bh^\toptzeroTOP \bH^\toptzero \bh^\toptzero}{\normtwo{\bh^\toptzero}^2} = \frac{\zeta_t}{\cos^2 (\theta_t)}.
\end{equation}
Similarly, rearranging \eqref{equation:bfgs_det} gives
$$
\det(\bH^\toptone) = \det(\bH^\toptzero) \frac{\by^\toptzeroTOP \bh^\toptzero}{\bh^\toptzeroTOP \bh^\toptzero} \frac{\bh^\toptzeroTOP \bh^\toptzero}{\bh^\toptzeroTOP \bH^\toptzero \bh^\toptzero} = \det(\bH^\toptzero) \frac{\alpha_t}{\zeta_t}.
$$
Additionally, define the  function
$
\phi(\bH) \triangleq \trace(\bH) - \ln \det(\bH).
$
This yields a recursive relation between $\phi(\bH^\toptone)$ and $\phi(\bH^\toptzero)$:
\begin{equation}\label{equation:bfgs_phi_rel}
\begin{aligned}
&\phi(\bH^\toptone) = \trace(\bH^\toptzero) + \beta_t - \frac{\zeta_t}{\cos^2 (\theta_t)} - \ln (\det(\bH^\toptzero)) - \ln(\alpha_t) + \ln (\zeta_t) \\
&= \phi(\bH^\toptzero) + (\beta_t - \ln(\alpha_t) - 1) + \left(1 - \frac{\zeta_t}{\cos^2 (\theta_t)}
+ \ln \left(\frac{\zeta_t}{\cos^2 (\theta_t)}\right) \right) 
+ \ln (\cos^2 (\theta_t))\\
&\leq  \phi(\bH^\toptzero) + (\beta - \ln(\alpha) - 1)
+ \ln (\cos^2 (\theta_t)),
\end{aligned}
\end{equation}
where the last inequality follows from that fact that $\beta>\beta_t$, $\alpha<\alpha_t$ for all $t$, and the fact that $1 - \mu + \ln \mu \leq 0,\, \forall \mu>0$.
Telescoping the sum over $j\in\{1,2,\ldots,t\}$ yields that 
$$
0 < \phi(\bH^{(t+1)}) \leq \phi(\bH^{(1)}) + t \sigma + \sum_{j=1}^t \ln \cos^2 (\theta_j),
\quad\forall t>0,
$$
where $\sigma\triangleq \beta - \ln (\alpha) - 1$ and without loss of generality we assume that $\sigma > 0$. Note that $\bh^\toptzero = -\eta_t (\bH^\toptzero)^{-1} \nabla f(\bx^\toptzero)$ is the descent step, then $\cos (\theta_t)$ is the cosine of the angle between the descent direction and the negative gradient direction. According to Theorem~\ref{theorem:zoutendijk_cond} and Theorem~\ref{theorem:conv_line_search} for the Zoutendijk condition based on the Wolfe condition (which is imposed in the theorem), $\normtwobig{\nabla f(\bx^\toptzero)}$ is greater than some nonzero constant  only if $\cos (\theta_t) \to 0$. Therefore, to prove $\normtwobig{\nabla f(\bx^\toptzero)} \to 0$, we only need to prove $\cos (\theta_t) \not\to 0$. Below we use a proof by contradiction to show this conclusion. Assume $\cos (\theta_t) \to 0$, then there exists $t_1 > 0$, for any $t > t_1$, we have
$
\ln (\cos^2 (\theta_t)) < -2\sigma.
$
Combining with  \eqref{equation:bfgs_phi_rel}, when $t > t_1$, we have:
$$
\begin{aligned}
0 &< \phi(\bH^\toptone) \leq \phi(\bH^{(1)}) + t\sigma + \sum_{j=1}^{t_1} \ln (\cos^2 (\theta_j)) + \sum_{j=t_1+1}^t (-2\sigma) \\
&= \phi(\bH^{(1)}) + \sum_{j=1}^{t_1} \ln (\cos^2 \theta_j) + 2\sigma t_1  - \sigma t. 
\end{aligned}
$$
The right-hand side of the above inequality is negative for sufficiently large $t$, while the left-hand side is 0, leading to a contradiction. Therefore, the assumption does not hold, i.e., there exists a subsequence $\{j_t\}_{t=1,2,\ldots}$ such that $\cos \theta_{j_t} \geq \delta > 0$. According to the Zoutendijk condition under the Wolfe condition (Theorem~\ref{theorem:zoutendijk_cond}), we can obtain $\liminf_{t \to \infty} \normtwobig{\nabla f(\bx^\toptzero)} \to 0$. 
\end{proof}

Theorem~\ref{theorem:conv_bfgs_glo} establishes the global convergence of the BFGS method but does not address its rate of convergence.
The following theorem outlines the conditions necessary for achieving superlinear convergence by the BFGS method.
The proof first shows that 
$
\lim_{t \to \infty} \frac{\normtwo{(\bH^\toptzero - \nabla^2 f(\bx^*)) \bh^\toptzero}}{\normtwo{\bh^\toptzero}} = 0
$, which combined with the Wolfe condition indicates a superlinear convergence; see \citet{nocedal1999numerical} for more details.
\begin{theoremHigh}[Rate of Convergence  of BFGS under Twice Lipschitz]\label{theorem:rate_bfgs}
Let $f(\bx):\real^n\rightarrow \real$ be twice Lipschitz continuously differentiable in a neighborhood of the optimal point $\bx^*$: $\normtwo{\nabla^2 f(\bx) - \nabla^2 f(\bx^*)}\leq L\normtwo{\bx-\bx^*}$ for all $\bx$ near $\bx^*$ with full-rank $\nabla^2 f(\bx^*)$.
Suppose that the iterates generated by the BFGS algorithm converge to an optimal point $\bx^*$. If the sequence $\{\bx^\toptzero\}$ satisfies
\begin{equation}\label{equation:bfgs_sum_error}
\sum_{t=1}^\infty \normtwobig{\bx^\toptzero - \bx^*} < +\infty,
\end{equation}
then $\{\bx^\toptzero\}$ converges \textbf{superlinearly} to $\bx^*$.
\end{theoremHigh}

As expected, because the quasi-Newton method uses an approximate Hessian matrix, it achieves only superlinear convergence, which is slower than the quadratic convergence of  Newton's method. However, since the quasi-Newton method does not require computing the Hessian matrix at each iteration, it may offer greater computational efficiency overall. This makes it particularly suitable for practical applications.

\section{Trust Region Method}\label{section:des_trust_reg}

This section introduces the trust region method. Both this method and the quasi-Newton algorithm are based on Taylor's expansion to approximate the target function locally, but they differ in how they utilize the approximating function. In the quasi-Newton algorithm, we first use the approximation model to determine the descent direction and then specify the step length. In contrast, in the trust region method, we directly solve the approximation model within a bounded region and then iterate to the next point. Therefore, the trust region method effectively selects both the direction and the step length simultaneously.

\subsection{Trust Region Method}

The methods discussed here generate a series of steps leading from the starting position towards the final solution. In the descent methods of Chapter~\ref{chapter:gradient-descent} and Newton's method described in Section~\ref{section:new_methods}, the directions of the steps are determined by the properties of $ f(\bx) $ at the current position.  Similar considerations lead us to trust region methods, where iteration steps are derived from the characteristics of a model of the objective function within a given region. The size of this region is adjusted during the iterations.

We provide an intuitive mathematical expression for the trust region method. Based on the Taylor expansion with Lagrange's remainder term,
$ f(\bx^\toptzero + \bd) = f(\bx^\toptzero) + \nabla f(\bx^\toptzero)^\top \bd + \frac{1}{2} \bd^\top \nabla^2 f(\bx^\toptzero + \gamma\bd) \bd, $
where $ \gamma\in (0, 1) $ is a positive number related to $ \bd $ (Theorem~\ref{theorem:linear_approx}). Similar to  Newton's method (Section~\ref{section:new_methods}), we use a quadratic approximation of $ f(\bx) $ at point $ \bx = \bx^\toptzero $:
\begin{equation}\label{equation:trs_approx}
\psi_t(\bd) \triangleq f(\bx^\toptzero) + \nabla f(\bx^\toptzero)^\top \bd + \frac{1}{2} \bd^\top \bH^\toptzero \bd,  
\end{equation}
where $ \bH^\toptzero $ is a symmetric matrix  that should approximate the Hessian matrix. 
If $ \bH^\toptzero $ exactly matches the Hessian matrix of $ f(\bx) $ at point $ \bx = \bx^\toptzero $, then when $ f(\bx) $ is sufficiently smooth, the approximation error of $ \psi_t(\bd) $ is $ o(\normtwo{\bd}^2) $ (Theorem~\ref{theorem:quad_app_theo}), and $\psi_t(\bd)$ reduces to the form in \eqref{equation:class_new_quadr}.

We have used the second-order Taylor expansion to approximate the target function $ f(\bx) $, it is important to consider that the Taylor expansion only holds for small values of $ \bd $. When $ \bd $ is too large, the model \eqref{equation:trs_approx} may no longer accurately describe the behavior of $ f(\bx) $. 
Therefore, constraints must be added to this model.  We only consider the approximation of $ f(\bx) $ within the following ball:
$$ 
\sB_t \triangleq \sB_2[\bx^\toptzero, \Delta_t]  \triangleq \{\bx^\toptzero + \bd \mid \normtwo{\bd} \leq \Delta_t\}, 
$$
where $ \Delta_t > 0 $ is a parameter related to the iteration. We call $ \sB_t $ the \textit{trust region}, and $ \Delta_t $ the \textit{trust region radius}. As the name suggests, the trust region is the area where we believe $ \psi_t(\bd) $ can adequately  approximate $ f(\bx) $, with $ \Delta_t $ representing the size of this region. 
Thus, each step of the trust region method requires solving the following subproblem (called the \textit{trust region subproblem, TRS}):
\begin{equation}\label{equation:trs_subpro}
\textbf{(TRS)}:\qquad \min_{\bd \in \real^n} \quad \psi_t(\bd), \quad \text{s.t.} \quad \normtwo{\bd} \leq \Delta_t.
\end{equation}

\begin{SCfigure}
\centering  
\vspace{-0.35cm} 
\subfigtopskip=2pt 
\subfigbottomskip=2pt 
\subfigcapskip=-5pt 
\includegraphics[width=0.65\textwidth]{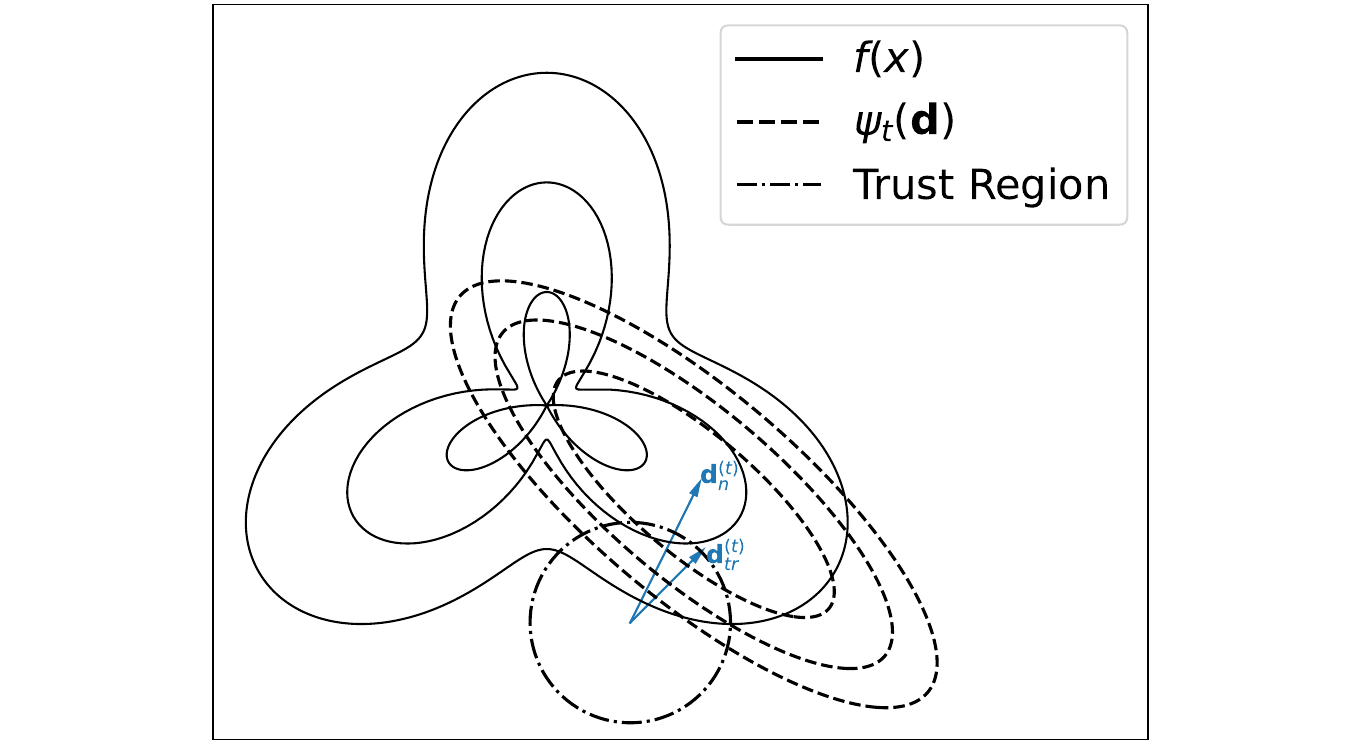}
\caption{Trust region method at the $t$-th iteration, where $\bd_{\text{n}}^\toptzero$ denotes the optimizer of $\psi(\bd)$ and $\bd_{\text{tr}}^\toptzero$ denotes the optimizer of the trust region method.}
\label{fig:trust_region}
\end{SCfigure}

Figure~\ref{fig:trust_region} shows the solution process for the subproblem in \eqref{equation:trs_subpro}. 
Solid lines represent the contour lines of $ f(\bx) $, while dashed lines represent those of $ \psi_t(\bd) $ (here, $ \bd = \bx - \bx^\toptzero $).
The vector  $ \bd_{\text{n}}^\toptzero $ denotes the descent direction obtained by solving the unconstrained problem $\min_{\bd} \psi_t(\bd)$ (if $ \bH^\toptzero $ is the Hessian matrix, then $ \bd_{\text{n}}^\toptzero $ is the Newton direction), and the vector $ \bd_{\text{tr}}^\toptzero $ denotes the descent direction obtained by solving the trust region subproblem \eqref{equation:trs_subpro}. 
It is evident that these two directions can be quite different. 
The trust region method restricts the size of $ \bd $, making the iteration more conservative and thus performing well even when the Newton direction is not optimal.

In the subproblem \eqref{equation:trs_subpro}, we must also determine the trust region radius $ \Delta_t $. 
In practice, selecting an appropriate trust region radius is crucial as it significantly influences the algorithm's convergence. 
Considering that the trust region radius is a measure of ``how much we trust the model $ \psi_t(\bd) $", if $ \psi_t(\bd) $ accurately approximates the function $ f(\bx) $, we should expand the trust region radius to leverage this approximation over a larger area; otherwise, we should reduce the radius and recompute. 
To measure the quality of the approximation $ \psi_t(\bd) $, we introduce the following \textit{gain factor}:
\begin{equation}\label{equation:trs_reduc_ratio}
\textbf{(Gain factor)}: \qquad \nu_t \triangleq \frac{f(\bx^\toptzero) - f(\bx^\toptzero + \bd_{\text{tr}}^\toptzero)}{\psi_t(\bzero) - \psi_t(\bd_{\text{tr}}^\toptzero)},
\end{equation}
where $ \bd_{\text{tr}}^\toptzero $ is the descent direction obtained by solving the subproblem \eqref{equation:trs_subpro}. According to the definition of $ \nu_t $, it represents the ratio of the \textit{actual reduction in the function value} to the \textit{predicted reduction} from the quadratic approximation:
\begin{itemize}
\item If $ \nu_t $ is close to 1, it indicates that using $ \psi_t(\bd) $ to approximate $ f(\bx) $ is highly  successful, suggesting we should increase $ \Delta_t $. 
\item If $ \nu_t $ is very small or even negative, it implies that we have overly trusted the quadratic approximation $ \psi_t(\bd) $, necessitating a reduction in $ \Delta_t $.
\end{itemize}
This mechanism dynamically adjusts $ \Delta_t $ to keep the domain of the quadratic model $ \psi_t(\bd) $ within a reasonable range.

Algorithm~\ref{alg:trust_region1} outlines the complete trust region method. Note that the trust region radius $ \Delta_t $ 
will not grow indefinitely because it has an upper bound control $ \Delta_{\max} $. Moreover, if the trust region constraint does not play a significant role (i.e., the optimal value of the quadratic model is within the trust region), there is no need to increase the trust region radius. Only when $ \psi_t(\bd) $ provides a good approximation and the trust region constraint is effective should we consider increasing $ \Delta_t $.

\begin{algorithm}
\caption{Descent Method with Trust Region}
\label{alg:trust_region1}
\begin{algorithmic}[1]
\Require A twice differentiable function $f(\bx)$; 
\State {\bfseries Input:}  Set the maximum radius $ \Delta_{\max} $, initial radius $ \Delta_1 $, initial point $ \bx^{(1)} $, accept radius $\gamma\in[0,\frac{1}{4})$;
\For{$t=1,2,\ldots$}
\State $\bd_{\text{tr}}^{\toptzero} \leftarrow \text{Solution of trust region subproblem \eqref{equation:trs_subpro}}$;
\State $\nu_t \leftarrow \text{gain factor \eqref{equation:trs_reduc_ratio}}$;
\If{$\nu_t > 0.75$ and $\normtwobig{\bd_{\text{tr}}^\toptzero} =\Delta_t$} \Comment{very good step, and the step is at the border}
\State $\Delta_{t+1} \leftarrow \min\{2 \Delta_t, \Delta_{\max}\}$; \Comment{larger trust region}
\EndIf
\If{$\nu_t < 0.25$} \Comment{poor step}
\State $\Delta_{t+1} \leftarrow \Delta_t / 4$; \Comment{smaller trust region}
\EndIf
\If{$\nu_t > \gamma$} \Comment{reject step if $\nu_t \leq \gamma$}
\State $\bx^\toptone \leftarrow \bx^\toptzero + \bd_{\text{tr}}^{\toptzero}$;
\Else
\State $\bx^\toptone \leftarrow \bx^\toptzero$;
\EndIf
\EndFor
\State {\bfseries Return:}  final $\bx\leftarrow \bx^{(t)}$;
\end{algorithmic}
\end{algorithm}

\subsection{Solving the Trust Region Subproblem}\label{section:solve_trs}
In Algorithm~\ref{alg:trust_region1}, one key issue remains unaddressed: how to solve the trust region subproblem.
In most practical applications, an explicit solution to the trust region subproblem \eqref{equation:trs_subpro} cannot be directly obtained. To find the iteration direction $ \bd_{\text{tr}}^\toptzero $, we need to design algorithms to either solve or approximate the subproblem \eqref{equation:trs_subpro}. 

\subsection*{Hidden Convexity of TRS}
Omitting the superscripts, we consider the following trust region subproblem:
\begin{equation}\label{equation:trs_sub_hidd}
\textbf{(TRS)}:\quad 	\min\, \left\{\psi(\bd) \triangleq f + \bg^\top \bd + \frac{1}{2} \bd^\top \bH \bd, \quad \text{s.t.} \quad \normtwo{\bd} \leq \Delta\right\},
\end{equation}
where $\bH\in\real^{n\times n}$ is symmetric, $\bg\in\real^n$, and $f\in\real$.
Since $\bH$ may not be necessarily positive semidefinite, the objective function is (possibly) non-convex, and the TRS
is (possibly) non-convex. 
We will now show how to transform (TRS) into a convex optimization problem.
First, by the spectral decomposition theorem (Theorem~\ref{theorem:spectral_theorem}), there exist an orthogonal matrix $\bQ$ and a diagonal matrix $\bLambda = \text{diag}(\lambda_1, \lambda_2, \ldots, \lambda_n)$ such that $\bH = \bQ \bLambda \bQ^\top$. 
Hence, (TRS) can be rewritten as
\begin{equation}\label{equation:trs_sub_hidd2}
\min \left\{ \bd^\top \bQ \bLambda \bQ^\top \bd + 2 \bg^\top \bQ \bQ^\top \bd + 2f : \normtwobig{\bQ^\top \bd}^2 \leq \Delta \right\},
\end{equation}
where we used orthogonal invariance of the $\ell_2$ norm:  $\normtwo{\bQ^\top \bd} = \normtwo{\bd}$. Making the linear change of variables $\by \triangleq \bQ^\top \bd$ and denoting $\bb \triangleq \bQ^\top \bg$,  \eqref{equation:trs_sub_hidd2} reduces to
\begin{equation}\label{equation:trs_sub_hidd3}
\begin{aligned}
& \min_{\by} && \sum_{i=1}^n \lambda_i y_i^2 + 2 \sum_{i=1}^n b_i y_i + 2f \\
& \text{s.t.} && \sum_{i=1}^n y_i^2 \leq \Delta.
\end{aligned}
\end{equation}
The problem is still non-convex since some of the $\lambda_i$'s might be negative. 
However, the signs of the optimal decision variables are known in advance, as stated in the following lemma.
\begin{lemma}[TRS]\label{lemma:trs_slack}
Let $\by^*$ be an optimal solution of \eqref{equation:trs_sub_hidd3}. Then $b_i y_i^* \leq 0$ for all $i = 1, 2, \ldots, n$.
\end{lemma}

\begin{proof}[of Lemma~\ref{lemma:trs_slack}]
Denote the objective function of problem \eqref{equation:trs_sub_hidd3} by
$ f(\by) \triangleq \sum_{i=1}^n \lambda_i y_i^2 + 2 \sum_{i=1}^n b_i y_i + 2f. $
For a $i \in \{1, 2, \ldots, n\}$, define the vector $\widetilde{\by}$ as
$ 
\widetilde{y}_j 
\triangleq 
\small
\begin{cases} 
y_j^*, & j \neq i, \\
-y_i^*, & j = i.
\end{cases} 
$

Then, obviously $\widetilde{\by}$ is also a feasible solution of \eqref{equation:trs_sub_hidd3}, and since $\by^*$ is an optimal solution of \eqref{equation:trs_sub_hidd3}, it follows that
$ f(\by^*) \leq f(\widetilde{\by})$.
This inequality simplifies to
$ \sum_{i=1}^n \lambda_i (y_i^*)^2 + 2 \sum_{i=1}^n b_i y_i^* + 2f \leq \sum_{i=1}^n \lambda_i \widetilde{y}_i^2 + 2 \sum_{i=1}^n b_i \widetilde{y}_i + 2f. $

Using the definition of $\widetilde{\by}$, the above inequality reduces after much cancelation of terms to
$ 2 b_i y_i^* \leq 2 b_i (-y_i^*), $
which implies the desired inequality $b_i y_i^* \leq 0$.
\end{proof}

As a direct consequence of Lemma~\ref{lemma:trs_slack}, for any optimal solution $\by^*$, the equality $\sign(y_i^*) = -\sign(b_i)$ holds when $b_i \neq 0$ and where the $\sign$ function is defined as
$ \sign(x) \triangleq 
\small
\begin{cases} 
1, & x \geq 0, \\
-1, & x < 0.
\end{cases} $

When $b_i = 0$, both $\by^*$ and $\widetilde{\by}$ are optimal (as shown in the proof of Lemma~\ref{lemma:trs_slack}), and hence the sign of $\by^*$ can be chosen arbitrarily. As a consequence, we can make the change of variables $y_i \triangleq -\sign(b_i) \sqrt{z_i} (z_i \geq 0)$, and problem \eqref{equation:trs_sub_hidd3} becomes
\begin{equation}\label{equation:trs_sub_convexrel}
\begin{aligned}
& \min_{\bz} && \sum_{i=1}^n \lambda_i z_i - 2 \sum_{i=1}^n |b_i| \sqrt{z_i} + f \\
& \text{s.t.} && \sum_{i=1}^n z_i \leq \Delta, \\
&&& z_1, z_2, \ldots, z_n \geq 0.
\end{aligned}
\end{equation}
This is clearly a convex optimization problem, with linear constraints and an objective function that is a sum of linear terms and positive multipliers of convex functions $-\sqrt{z_i}$. Therefore, we have shown that the non-convex trust region subproblem (TRS) is equivalent to the convex optimization problem \eqref{equation:trs_sub_convexrel}, which can be solved using any convex programming packages.

\subsection*{Iterative Method}
The trust region subproblem is a constrained optimization problem involving quadratic functions. A pertinent question arises: Can we use the optimality conditions from constrained optimization to solve this subproblem? The following theorem outlines the conditions that an optimal solution $ \bd^* $ must satisfy:

\begin{theorem}[Optimality Conditions of TRS]\label{theorem:trs_kkt}
Let $ \bd^* $ be any optimal point of  the following trust region subproblem
\begin{equation}\label{equation:trs_kkt1}
\min \quad \psi(\bd) = f + \bg^\top \bd + \frac{1}{2} \bd^\top \bH \bd, \quad \text{s.t.} \quad \normtwo{\bd} \leq \Delta,
\end{equation}
where $\bH\in\real^{n\times n}$ is symmetric, $\bg\in\real^n$, and $f\in\real$.
The global minimum is achieved if and only if $ \bd^* $ is feasible and there exists $ \lambda \geq 0 $ such that
\begin{subequations}\label{equation:trs_kkt_conds}
\begin{align}
(\bH + \lambda \bI) \bd^* &= -\bg,\label{equation:trs_kkt2} \\
\lambda (\Delta - \normtwo{\bd^*}) &= 0, \label{equation:trs_kkt3}\\
(\bH + \lambda \bI) &\succeq \bzero. \label{equation:trs_kkt4}
\end{align}
\end{subequations}
\end{theorem}
\begin{proof}[of Theorem~\ref{theorem:trs_kkt}]
\textbf{Necessity.} We apply the KKT conditions (Theorem~\ref{theorem:kktslat1_slater}) to directly derive the relationships satisfied by $ \bd^* $. 
Since there exist $\widehatbd\in\real^n$ such that $\normtwobig{\widehatbd}< \Delta$, 
the Lagrangian function for problem \eqref{equation:trs_kkt1} is given by:
$$
L(\bd, \lambda) = f + \bg^\top \bd + \frac{1}{2} \bd^\top \bH \bd + \frac{\lambda}{2} (\normtwo{\bd}^2 - \Delta^2),
$$
where the multiplier $ \lambda \geq 0 $. According to the KKT conditions, $ \bd^* $ is a feasible solution, and
$
\nabla_{\bd} L(\bd^*, \lambda) = (\bH + \lambda \bI) \bd^* + \bg = \bzero.
$
Additionally, the complementary slackness shows that 
$
\frac{\lambda}{2} (\Delta^2 - \normtwo{\bd^*}^2) = \bzero,
$
which after rearrangement gives \eqref{equation:trs_kkt2} and \eqref{equation:trs_kkt3}. To prove \eqref{equation:trs_kkt4}, we choose any $ \bd $ satisfying $ \normtwo{\bd} = \Delta $. By optimality and the complementary slackness, we have
$$
\psi(\bd) \geq \psi(\bd^*) = \psi(\bd^*) + \frac{\lambda}{2} (\normtwo{\bd^*}^2 - \normtwo{\bd}^2).
$$
Using equation \eqref{equation:trs_kkt2} to eliminate $\bg$ and substituting into the above formula gives
$
(\bd - \bd^*)^\top(\bH + \lambda \bI)(\bd - \bd^*) \geq \bzero,
$
showing that $\bH + \lambda \bI$ is positive semidefinite.

\paragraph{Sufficiency.}
Since $\bH$ is not necessarily positive semidefinite, the function $\psi(\bd)$ may not be convex, whence we cannot apply Theorem~\ref{theorem:kktslat1_slater} to prove the sufficiency.
Define the auxiliary function
$$
\widehat{\psi}(\bd) \triangleq f + \bg^\top \bd + \frac{1}{2}\bd^\top(\bH + \lambda \bI)\bd = \psi(\bd) + \frac{\lambda}{2}\bd^\top \bd.
$$
From condition \eqref{equation:trs_kkt4}, $\widehat{\psi}(\bd)$ is a convex function with respect to $\bd$. According to condition \eqref{equation:trs_kkt2}, $\bd^*$ satisfies the first-order optimality condition of the convex function. By Theorem~\ref{theorem:suff_sta_conv}, it can be deduced that $\bd^*$ is a global minimizer of $\widehat{\psi}(\bd)$. Given the complementary condition \eqref{equation:trs_kkt3}: $\lambda(\Delta^2 - \normtwo{\bd^*}^2) = 0$, for any feasible $\bd$, we have
$$
\psi(\bd) \geq \psi(\bd^*) + \frac{\lambda}{2}(\normtwo{\bd^*}^2 - \normtwo{\bd}^2)
= \psi(\bd^*) + \frac{\lambda}{2}(\Delta^2 - \normtwo{\bd}^2) \geq \psi(\bd^*). 
$$
This completes the proof.
\end{proof}

Theorem~\ref{theorem:trs_kkt} suggests an iterative method for finding $\bd^*$ when the problem dimension $n$ is small. Based on \eqref{equation:trs_kkt2}, the optimal solution forms a family of vectors parameterized by $\lambda$. 
Suppose $\bH + \lambda \bI$ is positive definite; we define
\begin{equation}
\bd(\lambda) \triangleq -(\bH + \lambda \bI)^{-1}\bg. 
\end{equation}
Thus, we only need to find an appropriate $\lambda$ such that \eqref{equation:trs_kkt3} and \eqref{equation:trs_kkt4} are satisfied. From the complementary condition \eqref{equation:trs_kkt3}, when $\lambda > 0$, we must have $\normtwo{\bd(\lambda)} = \Delta$. According to the positive semidefiniteness condition \eqref{equation:trs_kkt4}, $\lambda$ must satisfy $\lambda \geq -\lambda_{\min}(\bH)$ by the eigenvalue characterization theorem (Theorem~\ref{theorem:eigen_charac}), where $\lambda_{\min}(\bH)$ denotes the smallest eigenvalue of $\bH$. 

\begin{figure}[h]
\centering  
\vspace{-0.15cm} 
\subfigbottomskip=2pt 
\subfigcapskip=-5pt 
\includegraphics[width=0.8\textwidth]{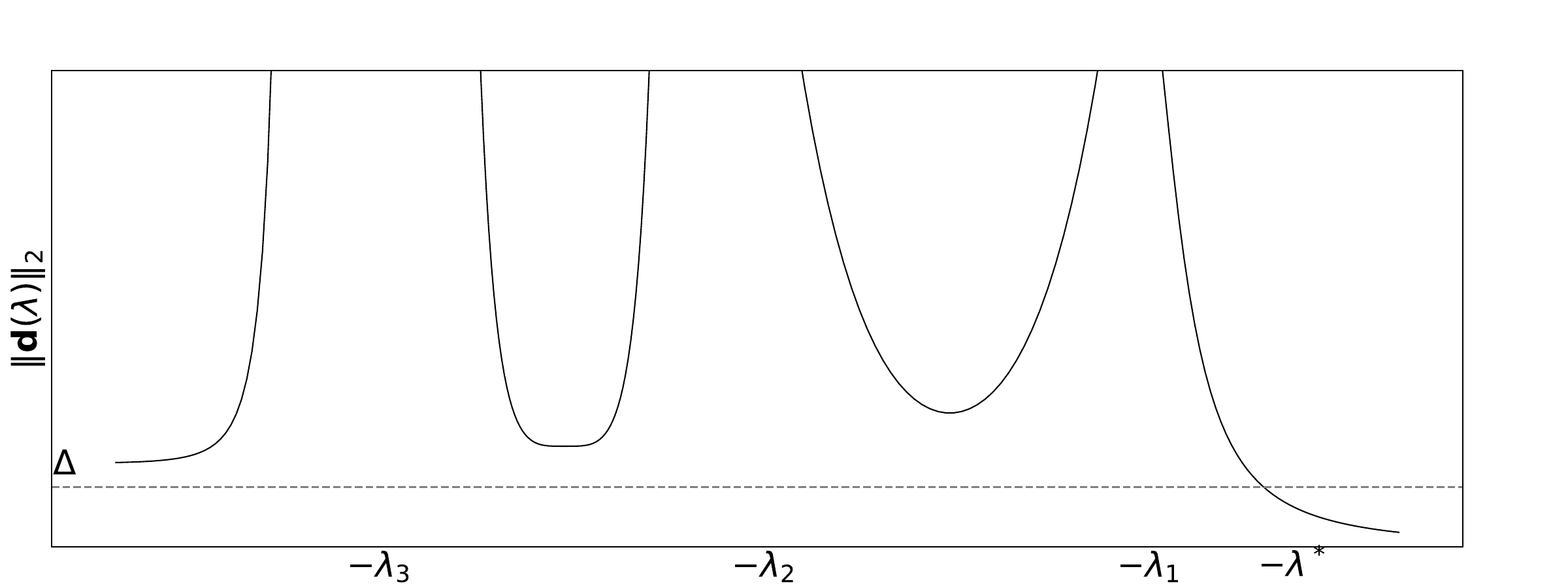}
\caption{$\normtwo{\bd(\lambda)}$ as a function of $\lambda$, where we assune $\bq_3^\top \bg \neq 0$, $\bq_2^\top \bg \neq 0$, and $\bq_1^\top \bg \neq 0$. When $\lambda\in(-\lambda_1, \infty)$, $\normtwo{\bd(\lambda)}$ is continuous and nonincreasing.}
\label{fig:trs_iterative}
\end{figure}

Now, let's examine the properties of $\normtwo{\bd(\lambda)}$ as $\lambda$ varies. Suppose $\bH$ admits a spectral decomposition $\bH = \bQ\bLambda \bQ^\top$ (Theorem~\ref{theorem:spectral_theorem}), where $\bQ = [\bq_1, \bq_2, \ldots, \bq_n]$ is an orthogonal matrix, $\bLambda = \diag(\lambda_1, \lambda_2, \ldots, \lambda_n)$ with $\lambda_1 \leq \lambda_2 \leq \ldots \leq \lambda_n$ is a diagonal matrix containing the eigenvalues of $\bH$. 
For $\lambda \neq \lambda_i, i\in\{1,2,\ldots,n\}$, we have  
\begin{equation}\label{equation:trs_dlambda}
\bd(\lambda) = -\bQ(\bLambda + \lambda \bI)^{-1}\bQ^\top \bg = -\sum_{i=1}^n \frac{\bq_i^\top \bg}{\lambda_i + \lambda}\bq_i
\quad\implies\quad 
\normtwo{\bd(\lambda)}^2 = \sum_{i=1}^n \frac{(\bq_i^\top \bg)^2}{(\lambda_i + \lambda)^2}.
\end{equation}
From this equation, it follows that:
\begin{itemize}
\item When $\lambda > -\lambda_1$, $\normtwo{\bd(\lambda)}^2$ is a continuous and nonincreasing  function of $\lambda$.
\item When $\lambda > -\lambda_1$ and $\bq_i^\top \bg \neq 0$ for all $i=1,2,\ldots,n$, $\normtwo{\bd(\lambda)}^2$ is a strictly decreasing function of $\lambda$.
\item It follows that $\lim_{\lambda \to \infty} \normtwo{\bd(\lambda)} = 0$.
\item When $\bq_i^\top \bg \neq 0$, $\lim_{\lambda \to -\lambda_i^+} \normtwo{\bd(\lambda)} = \infty$ for all $i=1,2,\ldots,n$.
\end{itemize}
See Figure~\ref{fig:trs_iterative} for a visual illustration.
Since the radius $\Delta>0$, when $\bq_1^\top\bg\neq 0$, there is always a unique solution $\lambda^*\in(-\lambda_1, \infty)$ such that $\normtwo{\bd(\lambda^*)}=\Delta$.
The positive semidefiniteness condition in \eqref{equation:trs_kkt4} requires that $\lambda^*\geq -\lambda_1$. Therefore, such a solution of $\lambda^*$ is unique in $\real$ as well.

\paragrapharrow{Case $\bq_1^\top\bg\neq 0$. }
In this scenario, if $\bH$ is positive definite and $\normtwo{\bH^{-1}\bg}\leq \Delta$, then the value $\lambda^*=0$ satisfies the optimality conditions in \eqref{equation:trs_kkt_conds}.
Otherwise, we may call the \textit{Newton-Raphson procedure}  to find the root of $\normtwo{\bd(\lambda)}=0$ in $(-\lambda_1, \infty)$.
Let $g(\lambda) \triangleq \normtwo{\bd(\lambda)} - \Delta$. We aim to find a $\lambda^*$ such that $g(\lambda^*)=0$.
To do this, we initialize $\lambda^{(1)} \in (-\lambda_1, \infty)$ for the first iteration. For the $t$-th iteration, using the linear approximation theorem at $\lambda^\toptzero$ (Theorem~\ref{theorem:linear_approx}):
$$
g(\lambda^\toptzero + \delta) \approx g(\lambda^\toptzero) +   \delta g'(\lambda^\toptzero).
$$
The condition $g(\lambda^\toptzero + \delta)=0$ implies the update rule
\begin{equation}
\lambda^\toptone \leftarrow \lambda^\toptzero - \frac{g(\lambda^\toptzero)}{g'(\lambda^\toptzero)}.
\end{equation}
The full procedure is described in Algorithm~\ref{alg:trs_itera}.
In fact, the approach described above can be applied even when the most negative eigenvalue is a multiple eigenvalue (that is, $0 > \lambda_1 = \lambda_2 = \ldots$), provided that $\bQ_1^\top \bg \neq \bzero$, where $\bQ_1$ is the matrix whose columns span the subspace corresponding to the eigenvalue $\lambda_1$.

\begin{algorithm}
\caption{Trust Region Subproblem}
\label{alg:trs_itera}
\begin{algorithmic}[1]
\State {\bfseries Input:}  $\lambda^{(1)}$, and given the radius $\Delta > 0$;
\For{$t = 1, 2, \ldots$}
\State  $\lambda^\toptone \leftarrow \lambda^\toptzero - \frac{g(\lambda^\toptzero)}{g'(\lambda^\toptzero)}$;
\EndFor
\State {\bfseries Return:}  $\bd \leftarrow \bd(\lambda^\toptzero)$;
\end{algorithmic}
\end{algorithm}

\paragrapharrow{Case $\bq_1^\top\bg= 0$. }
When $\bq_1^\top\bg= 0$ or $\bQ_1^\top \bg = \bzero$,  where $\bQ_1$ is the matrix whose columns span the subspace corresponding to the eigenvalue $\lambda_1$, then the condition $\lim_{\lambda \to -\lambda_1^+} \normtwo{\bd(\lambda)} = \infty$ does not hold.
Therefore, there may not be a value $\lambda^* \in (-\lambda_1, \infty)$ such that $\normtwo{\bd(\lambda^*)} = \Delta$ \citep{more1983computing, nocedal1999numerical}. 
At first glance, it might not be clear how $\bd$ and $\lambda$ can be chosen to satisfy the optimality conditions in \eqref{equation:trs_kkt_conds} in this case. 
However, Theorem~\ref{theorem:trs_kkt} still guarantees that the correct value of $\lambda^*$ lies within the interval $[-\lambda_1, \infty)$. 
Thus, the only  possibility is $\lambda^* = -\lambda_1$. All that is left is to find the corresponding $\bd$.
 To find $\bd$, it is not sufficient to simply delete the terms for which $\lambda_i = \lambda_1$ from the formula \eqref{equation:trs_dlambda} and set
\begin{equation}
	\bd = \sum_{i:\lambda_i \neq \lambda_1} \frac{\bq_i^\top \bg}{\lambda_i + \lambda} \bq_i.
\end{equation}
Instead, we note that $(\bH - \lambda_1 \bI)$ is singular, meaning there exists a vector $\bu$ such that $\normtwo{\bu} = 1$ and $(\bH - \lambda_1 \bI)\bu = \bzero$ ($\bu$ is an eigenvector of $\bH$ corresponding to the eigenvalue $\lambda_1$).
Therefore, by the property of the spectral decomposition,  we have $\bq_i^\top \bu = 0$ for $\lambda_i \neq \lambda_1$. 
From this property, if we set
\begin{equation}
\bd \triangleq \sum_{i:\lambda_i \neq \lambda_1} \frac{\bq_i^\top \bg}{\lambda_i + \lambda} \bq_i + \mu \bu
\end{equation}
for any scalar $\mu$, we have
$$
\normtwo{\bd}^2 = \sum_{i:\lambda_i \neq \lambda_1} \left(\frac{\bq_i^\top \bg}{\lambda_i + \lambda}\right)^2 + \mu^2.
$$
This implies that by choosing an appropriate $\mu$, we can always ensure that $\normtwo{\bd} = \Delta$. 
It is straightforward to verify that the optimality conditions \eqref{equation:trs_kkt_conds} are satisfied for this choice of $\bd$ and $\lambda = -\lambda_1$.

\subsection{Convergence Analysis}

This subsection provides a brief introduction to the convergence results of the trust region method.

To estimate the improvement in function values obtained by solving each trust region subproblem, we introduce the definition of the Cauchy point.
\begin{definition}[Cauchy Point]\label{definition:cauchy_point}
Let $\psi_t(\bd)$ be a quadratic approximation of $f(\bx)$ at the point $\bx = \bx^\toptzero$. The constant $\mu_t$ is the solution to the following optimization problem:
$$
\mu_t = \mathop{\argmin}_{\mu\in\real}  \psi_t(-\mu \nabla f(\bx^\toptzero)),
\quad\text{s.t.}\quad
\normtwo{\mu \nabla f(\bx^\toptzero)} \leq \Delta_t, \mu \geq 0.
$$
Then the point $\bx_c^\toptzero \triangleq \bx^\toptzero + \bd_c^\toptzero$ is called the \textit{Cauchy point}, where $\bd_c^\toptzero \triangleq -\mu_t \nabla f(\bx^\toptzero)$.
\end{definition}

\begin{SCfigure}
\centering  
\vspace{-0.35cm} 
\subfigtopskip=2pt 
\subfigbottomskip=2pt 
\subfigcapskip=-5pt 
\includegraphics[width=0.6\textwidth]{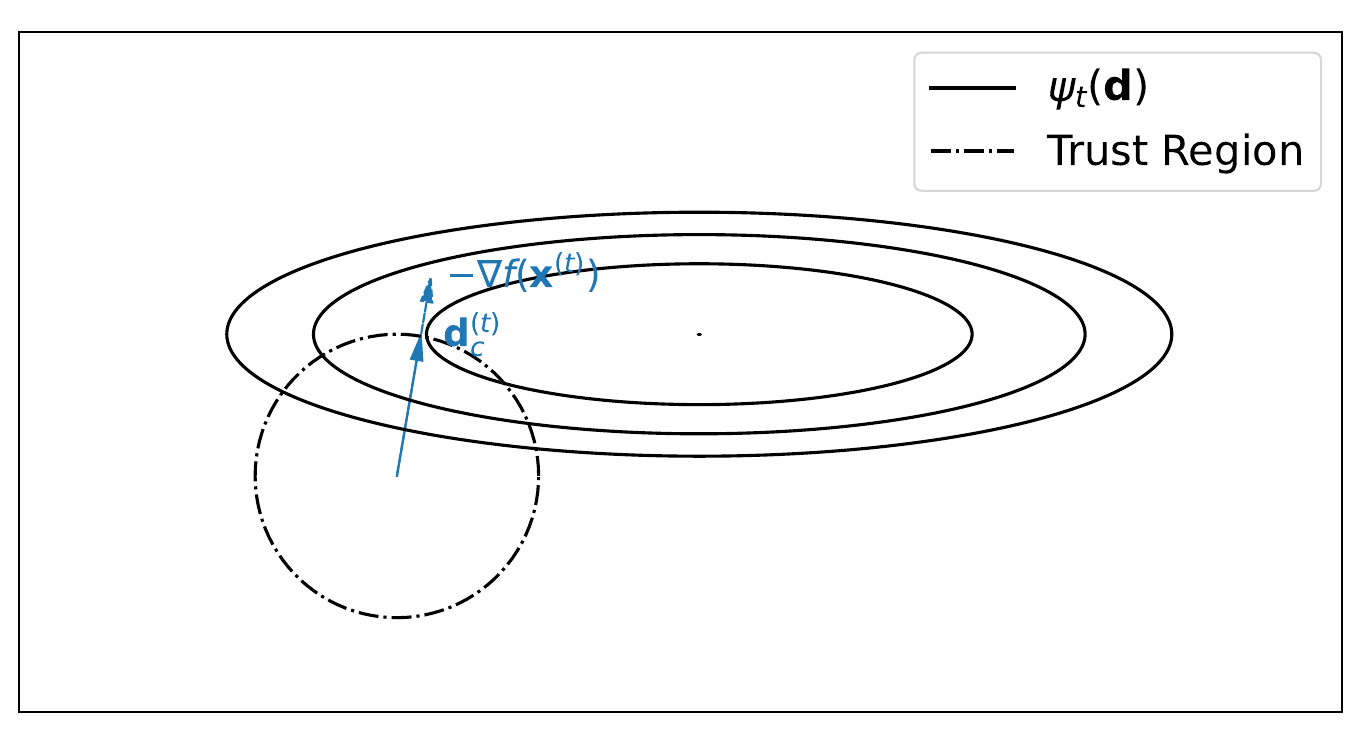}
\caption{Calculation of the Cauchy point at the $t$-th iteration, which is equivalent to a line search constrained within the trust region.}
\label{fig:cauchy_point}
\end{SCfigure}

According to this definition, the Cauchy point essentially represents an application of the gradient method with exact line search to $\psi_t(\bd)$, considering the trust region constraint. Figure~\ref{fig:cauchy_point} visually illustrates the concept of the Cauchy point.

In fact, given $\psi_t(\bd)$, the Cauchy point can be explicitly calculated. For convenience, denote $\bg^\toptzero \triangleq \nabla f(\bx^\toptzero)$. 
When $ \bg^\toptzeroTOP \bH^\toptzero \bg^\toptzero \leq 0$, the function $\psi_t(-\mu \nabla f(\bx^\toptzero))$ is concave quadratic in $\mu$ and  decreases monotonically with $\mu$ whenever $\bg^\toptzero\neq \bzero$. 
Thus, $\mu_t$ is simply the largest value that satisfies the trust region bound.
When $ \bg^\toptzeroTOP \bH^\toptzero \bg^\toptzero >0$, the function $\psi_t(-\mu \nabla f(\bx^\toptzero))$ is convex quadratic in $\mu$, so $\mu_t$ is either unconstrained minimizer of this quadratic function, $\bg^\toptzeroTOP \bH^\toptzero \bg^\toptzero$,  or the boundary value, which is equivalent to the first case. 
This yields the following form:
$$
\mu_t = 
\begin{cases} 
\frac{\Delta_t}{\normtwobig{\bg^\toptzero}}, & \bg^\toptzeroTOP \bH^\toptzero \bg^\toptzero \leq 0, \\
\min \left\{ \frac{\normtwobig{\bg^\toptzero}^2}{\bg^\toptzeroTOP \bH^\toptzero \bg^\toptzero}, \frac{\Delta_t}{\normtwobig{\bg^\toptzero}} \right\}, & \text{otherwise}.
\end{cases}
$$

The above analysis shows that the Cauchy point is a feasible solution for the trust region subproblem \eqref{equation:trs_subpro}. In practice, the Cauchy point is not generally used as an approximate solution for the next iteration because it essentially represents a steepest descent method with a truncated step length, failing to fully utilize the Hessian matrix $\bH^\toptzero$.
However, in the nonlinear least squares problem, the Cauchy point can be applied to find an approximation with additional considerations; see Section~\ref{section:powell_dog_leg1}.

Despite this, the Cauchy point serves as a benchmark for evaluating the performance of algorithms solving the trust region subproblem. Specifically, it requires that the points generated by the subproblem algorithm are at least as effective as the Cauchy point. It is straightforward to see that using the Cauchy point ensures a decrease in the objective function value of the quadratic model. Indeed, we have the following descent lemma under the Cauchy point.
\begin{lemma}[Descent Lemma under the Cauchy Point]\label{lemma:descen_cauchy_point}
Let $\psi_t(\bd)$ be a quadratic approximation of $f(\bx)$ at the point $\bx = \bx^\toptzero$, and let $\bd_c^\toptzero$ be the descent direction generated by solving for the Cauchy point. 
Then,
\begin{equation}
\psi_t(\bzero) - \psi_t(\bd_c^\toptzero) 
\geq \frac{1}{2} \normtwobig{\bg^\toptzero}
\min \left\{ \Delta_t, \frac{\normtwobig{\bg^\toptzero}}{\normtwo{\bH^\toptzero}} \right\}.
\end{equation}
\end{lemma}
\begin{proof}[of Lemma~\ref{lemma:descen_cauchy_point}]
\textbf{Case 1: $\bg^\toptzeroTOP \bH^\toptzero \bg^\toptzero \leq 0$.} In this case, we have 
\begin{align*}
\psi_t(\bd_c^\toptzero) - \psi_t(\bzero) 
&= \psi_t\Big(-\frac{\Delta_t }{\normtwobig{\bg^\toptzero}} \bg^\toptzero\Big) - f_t 
= \frac{\Delta_t^2}{2\normtwobig{\bg^\toptzero}^2} \bg^\toptzeroTOP \bH^\toptzero \bg^\toptzero -\Delta_t \normtwobig{\bg^\toptzero} \\
&\leq -\Delta_t \normtwobig{\bg^\toptzero} 
\leq -\normtwobig{\bg^\toptzero} \min \left\{  \Delta_t, \frac{\normtwobig{\bg^\toptzero}}{\normtwobig{\bH^\toptzero}} \right\},
\end{align*}
where $f_t \triangleq f(\bx^\toptzero) $, $\bg^\toptzero \triangleq\nabla f(\bx^\toptzero)$, and $\bH^\toptzero \triangleq\nabla^2 f(\bx^\toptzero)$.
This yields the desired result.

\paragraph{Case 2: $\bg^\toptzeroTOP \bH^\toptzero \bg^\toptzero > 0$ and $\frac{\normtwobig{\bg^\toptzero}^2}{\bg^\toptzeroTOP \bH^\toptzero \bg^\toptzero} \leq \frac{\Delta_t}{\normtwobig{\bg^\toptzero}}$.}
In this case,  it follows that  $\mu_t = \frac{\normtwobig{\bg^\toptzero}^2}{\bg^\toptzeroTOP \bH^\toptzero \bg^\toptzero}$, whence we have 
$$ 
\small
\begin{aligned}
\psi_t(\bd_c^\toptzero) - \psi_t(\bzero) 
&=  \frac{1}{2} \bg^\toptzeroTOP \bH^\toptzero \bg^\toptzero \frac{\normtwobig{\bg^\toptzero}^4}{\big(\bg^\toptzeroTOP \bH^\toptzero \bg^\toptzero\big)^2} 
-\frac{\normtwobig{\bg^\toptzero}^4}{\bg^\toptzeroTOP \bH^\toptzero \bg^\toptzero}
= -\frac{1}{2} \frac{\normtwobig{\bg^\toptzero}^4}{\bg^\toptzeroTOP \bH^\toptzero \bg^\toptzero} \\
&\stackrel{\dag}{\leq} -\frac{1}{2} \frac{\normtwobig{\bg^\toptzero}^4}{\normtwo{\bH^\toptzero} \normtwobig{\bg^\toptzero}^2} 
= -\frac{1}{2} \frac{\normtwobig{\bg^\toptzero}^2}{\normtwo{\bH^\toptzero}} 
\leq -\frac{1}{2} \normtwobig{\bg^\toptzero} \min \left\{\Delta_t, \frac{\normtwobig{\bg^\toptzero}}{\normtwo{\bH^\toptzero}}\right\},
\end{aligned}
$$
where the inequality ($\dag$) follows from the Cauchy-Schwartz inequality. This also yields the desired result.

\paragraph{Case 3: $\bg^\toptzeroTOP \bH^\toptzero \bg^\toptzero > 0$ and $\frac{\normtwobig{\bg^\toptzero}^2}{\bg^\toptzeroTOP \bH^\toptzero \bg^\toptzero} > \frac{\Delta_t}{\normtwobig{\bg^\toptzero}}$.}
In the remaining case, we have $\mu_t = \frac{\Delta_t}{\normtwobig{\bg^\toptzero}}$. Therefore,
$$
\small
\begin{aligned}
\psi_t(\bd_c^\toptzero) - \psi_t(\bzero) 
&= \frac{1}{2} \frac{\Delta_t^2}{\normtwobig{\bg^\toptzero}^2} \bg^\toptzeroTOP \bH^\toptzero \bg^\toptzero 
-\frac{\Delta_t}{\normtwobig{\bg^\toptzero}} \normtwobig{\bg^\toptzero}^2 
\leq  \frac{1}{2} \frac{\Delta_t^2}{\normtwobig{\bg^\toptzero}^2} \frac{\normtwobig{\bg^\toptzero}^3}{\Delta_t} -\Delta_t \normtwobig{\bg^\toptzero} \\
&= -\frac{1}{2} \Delta_t \normtwobig{\bg^\toptzero} 
\leq -\frac{1}{2} \normtwobig{\bg^\toptzero} \min \left\{\Delta_t, \frac{\normtwobig{\bg^\toptzero}}{\normtwo{\bH^\toptzero}}\right\},
\end{aligned}
$$
yielding the desired result once again.
\end{proof}

In practice, it is then essential to obtain iteration directions $\bd_{\text{tr}}^\toptzero$ that satisfy
$ \psi_t(\bzero) - \psi_t(\bd_{\text{tr}}^\toptzero) \geq(\psi_t(\bzero) - \psi_t(\bd_c^\toptzero)) $ to ensure sufficient decrease at each iteration.
This means that the estimate
$ \psi_t(\bzero) - \psi_t(\bd_{\text{tr}}^\toptzero) \geq  \frac{1}{2}\normtwobig{\bg^\toptzero} \min \left\{ \Delta_t, \frac{\normtwobig{\bg^\toptzero}}{\normtwo{\bH^\toptzero}} \right\} $
holds in many cases. This provides a basis for proving the convergence of the trust region algorithm.
\subsection*{Global Convergence}

Now we introduce the global convergence of the trust region algorithm. Recall Algorithm~\ref{alg:trust_region1}, where an acceptance radius $\gamma\in[0,\frac{1}{4})$ is introduced to determine whether to update the iteration point. 
There are two scenarios: when $\gamma = 0$, as long as the original objective function has a decrease will the trust region step be accepted; when $\gamma \in(0,\frac{1}{4})$, only if the improvement $\nu_t$ reaches a certain level will an update be made. 
The convergence results obtained in these two situations differ. Following \citet{nocedal1999numerical}, we separately introduce these two results.
\begin{theoremHigh}[Global Convergence of Trust Region under $\gamma=0$ and SS]\label{theorem:global_conv_trs_eta0}
Let $\gamma = 0$ in Algorithm~\ref{alg:trust_region1}. Suppose that $\normtwo{\bH^\toptzero} \leq L$ for some constant $L$, that $f$ is bounded below on the level set $\sL\triangleq\lev[f, f(\bx^{(1)})]$ and is $\beta$-smooth in the open neighborhood $\sB(\sL, R)  \triangleq \{\bx \mid  \normtwo{\bx-\by} < R \text{ for some } \by\in \sL\}$ for some $R > 0$, and that all approximate solutions of the trust region subproblem satisfy the following inequalities: 
\begin{subequations}
\begin{align}
\psi_t(\bzero) - \psi_t(\bd_{\text{tr}}^\toptzero) 
&\geq c_1 \normtwobig{\bg^\toptzero}
\min \left\{ \Delta_t, \frac{\normtwobig{\bg^\toptzero}}{\normtwo{\bH^\toptzero}} \right\}; \label{equation:global_conv_trs_eta01}\\
\normtwobig{\bd_{\text{tr}}^\toptzero} & \leq c_2\Delta_t, \label{equation:global_conv_trs_eta02}
\end{align}
\end{subequations}
for some positive constants $c_1\leq 1$ and $c_2\geq 1$~\footnote{Note that $c_1=\frac{1}{2}$ in the setting of Cauchy points by Lemma~\ref{lemma:descen_cauchy_point}.}, where $\bd_{\text{tr}}^\toptzero$ is a solution to the trust region subproblem in \eqref{equation:trs_subpro}.
Then
$$
\liminf_{t \to \infty} \normtwobig{\bg^\toptzero} = 0,
$$
i.e., the limit points of $\bx^\toptzero$ contain stationary points.
\end{theoremHigh}
\begin{proof}[of Theorem~\ref{theorem:global_conv_trs_eta0}]
By performing some technical manipulation with the gain factor $\nu_t$ from \eqref{equation:trs_reduc_ratio}, we obtain
\begin{align*}
|\nu_t - 1| 
&= \left| \frac{\big(f(\bx^\toptzero) - f(\bx^\toptzero + \bd_{\text{tr}}^\toptzero)\big) - \big(\psi_t(\bzero) - \psi_t(\bd_{\text{tr}}^\toptzero)\big)}{\psi_t(\bzero) - \psi_t(\bd_{\text{tr}}^\toptzero)} \right| 
= \left| \frac{\psi_t(\bd_{\text{tr}}^\toptzero) - f(\bx^\toptzero + \bd_{\text{tr}}^\toptzero)}{\psi_t(\bzero) - \psi_t(\bd_{\text{tr}}^\toptzero)} \right|.
\end{align*}
For simplicity, denote $\bg(\bx) \triangleq \nabla f(\bx)$. By the fundamental theorem of calculus (Theorem~\ref{theorem:fund_theo_calculu}), we have  
$$
f(\bx^\toptzero + \bd_{\text{tr}}^\toptzero) = f(\bx^\toptzero) + \bg(\bx^\toptzero)^\top \bd_{\text{tr}}^\toptzero + \int_0^1 \big[\bg\big(\bx^\toptzero + \mu \bd_{\text{tr}}^\toptzero\big) - \bg(\bx^\toptzero)\big]^\top \bd_{\text{tr}}^\toptzero \, d\mu.
$$
It then follows from the definition  of the quadratic approximation $\psi_t$ (Equation~\eqref{equation:trs_approx}) that
\begin{equation}\label{equation:global_conv_trs_eta0_prv1}
\begin{aligned}
\abs{\psi_t(\bd_{\text{tr}}^\toptzero) - f(\bx^\toptzero + \bd_{\text{tr}}^\toptzero)}
&= \left| \frac{1}{2} \bd_{\text{tr}}^\toptzeroTOP \bH^\toptzero \bd_{\text{tr}}^\toptzero - \int_0^1 \big[\bg\big(\bx^\toptzero + \mu \bd_{\text{tr}}^\toptzero\big) - \bg(\bx^\toptzero)\big]^\top \bd_{\text{tr}}^\toptzero \, d\mu \right| \\
&\leq (L/2) \normtwobig{\bd_{\text{tr}}^\toptzero}^2 + \beta \normtwobig{\bd_{\text{tr}}^\toptzero}^2,
\end{aligned}
\end{equation}
where the last inequality follows from the triangle inequality and the Cauchy-Schwartz inequality.
Note we assumed that $\normtwobig{\bd_{\text{tr}}^\toptzero} \leq R$ to ensure that $\bx^\toptzero$ and $\bx^\toptzero + \mu\bd_{\text{tr}}^\toptzero$ both lie in the set $\sB(\sL, R)$.

Suppose for contradiction that there exists $\epsilon > 0$ and a positive index $T$ such that
\begin{equation}\label{equation:global_conv_trs_eta0_assup}
\normtwobig{\bg^\toptzero} \geq \epsilon, \quad \text{for all } t \geq T.
\end{equation}
This indicates that, by \eqref{equation:global_conv_trs_eta01}, 
\begin{equation}\label{equation:global_conv_trs_eta0_prv2}
\psi_t(\bzero) - \psi_t(\bd_{\text{tr}}^\toptzero) 
\geq c_1 \normtwobig{\bg^\toptzero} \min \left\{\Delta_t, \frac{\normtwobig{\bg^\toptzero}}{\normtwo{\bH^\toptzero}}\right\} \geq c_1 \epsilon \min \left( \Delta_t, \frac{\epsilon}{L} \right), 
\quad\forall\, t \geq T.
\end{equation}
Combining \eqref{equation:global_conv_trs_eta0_prv1}, \eqref{equation:global_conv_trs_eta0_prv2}, and the assumption \eqref{equation:global_conv_trs_eta02} yields that 
\begin{equation}\label{equation:global_conv_trs_eta0_prv3}
|\nu_t - 1| \leq \frac{c_2^2 \Delta_t^2 (L/2 + \beta)}{c_1 \epsilon \min(\Delta_t, \epsilon/L)}.
\end{equation}

\paragraph{Case 1: $\Delta_t$ is upper bounded  by $\widetilde{\Delta}$. } We now derive a bound on the right-hand side of the above inequality that holds for all sufficiently small values of $\Delta_t$, that is, for all $\Delta_t \leq \widetilde{\Delta}$, where $\widetilde{\Delta}$ is defined as follows:
\begin{equation}\label{equation:global_conv_trs_eta0_prv4}
\widetilde{\Delta} \triangleq \min \left\{\frac{1}{2} \frac{c_1 \epsilon}{c_2^2 (L/2 + \beta)}, \frac{R}{c_2}\right\}.
\end{equation}
The $R/c_2$ term in this definition ensures that the bound \eqref{equation:global_conv_trs_eta0_prv2} is valid (because $\normtwobig{\bd_{\text{tr}}^\toptzero} \leq c_2 \Delta_t \leq c_2 \widetilde{\Delta} \leq R$). Note that since $c_1 \leq 1$ and $c_2 \geq 1$, we have $\widetilde{\Delta} \leq \epsilon/L$. The  condition $\widetilde{\Delta} \leq \epsilon/L$ implies that for all $\Delta_t \in [0, \widetilde{\Delta}]$, we have $\min(\Delta_t, \epsilon/L) = \Delta_t$, so from \eqref{equation:global_conv_trs_eta0_prv3} and \eqref{equation:global_conv_trs_eta0_prv4}, we have
$$
|\nu_t - 1| \leq \frac{c_2^2 \Delta_t^2 (L/2 + \beta)}{c_1 \epsilon \Delta_t} = \frac{c_2^2 \Delta_t (L/2 + \beta)}{c_1 \epsilon} \leq \frac{c_2^2 \widetilde{\Delta} (L/2 + \beta)}{c_1 \epsilon} \leq \frac{1}{2}.
$$
Therefore, $\nu_t > \frac{1}{4}$, and so by the workings of Algorithm~\ref{alg:trust_region1}, we have $\Delta_{t+1} \geq \Delta_t$ whenever $\Delta_t$ falls below the threshold $\widetilde{\Delta}$. 
\paragraph{Case 2: $\Delta_t$ is lower bounded by $\widetilde{\Delta}$. }It follows that the reduction of $\Delta_t$ (by a factor of $\frac{1}{4}$  in Algorithm~\ref{alg:trust_region1}) can occur in our algorithm only if
$
\Delta_t \geq \widetilde{\Delta}.
$

Combining the two cases concludes that
\begin{equation}\label{equation:global_conv_trs_eta0_prv5}
\Delta_t \geq \min \left( \Delta_T, \frac{\widetilde{\Delta}}{4} \right) \qquad \text{for all } t \geq T.
\end{equation}
Suppose now that there is an infinite subsequence $\sT$ such that $\nu_t \geq \frac{1}{4}$ for $t \in \sT$. 
For $t \in \sT$ and $t \geq T$, by the definition of the gain factor $\nu_t$ in \eqref{equation:trs_reduc_ratio} and \eqref{equation:global_conv_trs_eta0_prv2}, we have 
$$
f(\bx^\toptzero) - f(\bx^\toptone) = f(\bx^\toptzero) - f(\bx^\toptzero + \bd_{\text{tr}}^\toptzero) 
\geq \frac{1}{4} [\psi_t(\bzero) - \psi_t(\bd_{\text{tr}}^\toptzero)] 
\geq \frac{1}{4} c_1 \epsilon \min \left( \Delta_t, \frac{\epsilon}{L} \right).
$$
Since $f$ is bounded below, the above inequality indicates that
$
\lim_{t \in \sT, \; t \to \infty} \Delta_t = 0,
$
which contradicts \eqref{equation:global_conv_trs_eta0_prv5}. Hence no such infinite subsequence $\sT$ can exist, and we must have $\nu_t < \frac{1}{4}$ for all $t$ sufficiently large. In this case, $\Delta_t$ will eventually be multiplied by $\frac{1}{4}$ at every iteration, and we have $\lim_{t \to \infty} \Delta_t = 0$, which again contradicts \eqref{equation:global_conv_trs_eta0_prv5}. 
Therefore, the assumption in \eqref{equation:global_conv_trs_eta0_assup} does not hold, and this completes the proof.
\end{proof}

Theorem~\ref{theorem:global_conv_trs_eta0} indicates that if $\gamma=0$ in Algorithm~\ref{alg:trust_region1}, meaning the algorithm always accepts the step when the gain factor is positive, then Algorithm~\ref{alg:trust_region1} only guarantees subsequential convergence; the sequence of iterates itself may not converge.
However, the following theorem demonstrates that choosing $\gamma > 0$ can lead to stronger convergence results.
\begin{theoremHigh}[Global Convergence of Trust Region under $\gamma\neq0$ and SS]\label{theorem:glo_conv_trs_othe}
Let $\gamma \in (0, \frac{1}{4})$ in Algorithm~\ref{alg:trust_region1}, and consider the same conditions  as Theorem~\ref{theorem:global_conv_trs_eta0}.
Then,
\begin{equation}\label{equation:glo_conv_trs_othe}
\lim_{t \to \infty} \bg^\toptzero = \bzero.
\end{equation}
\end{theoremHigh}

\begin{proof}
Consider a particular positive iteration $z$ with a nonzero gradient vector: $\bg^{(z)} \neq \bzero$. Since $f$ is $\beta$-smooth on the set $\sB(\sL, R)$, we have
$$
\normtwobig{\bg(\bx) - \bg^{(z)}} \leq \beta \normtwobig{\bx - \bx^{(z)}}, \quad \text{for all }\bx \in \sB(\sL, R).
$$
Define the scalars $\epsilon$ and $R$ to satisfy
$$
\epsilon \triangleq \frac{1}{2} \normtwobig{\bg^{(z)}}
\qquad \text{and}\qquad 
\widetildeR \triangleq \min \left( \frac{\epsilon}{\beta}, R \right).
$$
Note that the closed ball
$
\sB[\bx^{(z)}, \widetildeR] \triangleq \big\{ \bx \mid \normtwo{\bx - \bx^{(z)}} \leq \widetildeR \big\}
$
is contained in the open neighborhood $\sB(\sL, R)$, so the smoothness of $f$ also holds inside $\sB[\bx^{(z)}, \widetildeR] $. 
We have 
$$
\bx \in \sB[\bx^{(z)}, \widetildeR] 
\quad\implies\quad
\normtwobig{\bg(\bx)} \stackrel{\dag}{\geq} \normtwobig{\bg^{(z)}} - \normtwobig{\bg(\bx) - \bg^{(z)}} \stackrel{\ddag}{\geq} \frac{1}{2} \normtwobig{\bg^{(z)}} = \epsilon.
$$
where the inequality $(\dag)$ follows from the  triangle inequality $\abs{\normtwo{\ba}-\normtwo{\bb}}\leq \normtwo{\ba}-\normtwo{\bb}$ for all $\ba,\bb$, and the inequality $(\ddag)$ follows from the smoothness and the definition of $\widetildeR$.

If the entire sequence $\{\bx^\toptzero\}_{t \geq z}$ stays inside the closed ball $\sB[\bx^{(z)}, \widetildeR] $, we would have $\normtwobig{\bg^\toptzero} \geq \epsilon > 0$ for all $t \geq z$, which contradicts the result in \eqref{equation:glo_conv_trs_othe}. Therefore, there must be an index $t \geq z$ such that $\bx^\toptzero \notin \sB[\bx^{(z)}, \widetildeR] $ for all $t \geq z$. The reasoning in the proof of Theorem~\ref{theorem:global_conv_trs_eta0} can be used to show that this scenario does not occur. Therefore, the sequence $\{\bx^\toptzero\}_{t \geq z}$ eventually leaves $\sB[\bx^{(z)}, \widetildeR] $.

Let $l \geq z$ be the iteration indix such that $\bx^{(l+1)}$ is the first iterate after $\bx^{(z)}$ outside $\sB[\bx^{(z)}, \widetildeR] $. Since $\normtwobig{\bg^\toptzero} \geq \epsilon$ for $t = z, z+1, \ldots, l$, performing telescopic cancellations over $t = z, z+1, \ldots, l$ and using \eqref{equation:global_conv_trs_eta0_prv2}, we have 
\begin{align*}
f(\bx^{(z)}) - f(\bx^{(l+1)}) &= \sum_{t=z}^{l} f(\bx^\toptzero) - f(\bx^\toptone) 
= \sum_{t=z, \bx^\toptzero \neq \bx^\toptone}^{l} \nu_t \big[\psi_t(\bzero) - \psi_t(\bd_{\text{tr}}^\toptzero)\big] \\
&\stackrel{\dag}{\geq} \sum_{t=z, \bx^\toptzero \neq \bx^\toptone}^{l} \gamma \big[\psi_t(\bzero) - \psi_t(\bd_{\text{tr}}^\toptzero)\big]
\geq \sum_{t=z, \bx^\toptzero \neq \bx^\toptone}^{l} \gamma c_1 \epsilon \min \left( \Delta_t, \frac{\epsilon}{L} \right),
\end{align*}
where the inequality $(\dag)$ follows from the fact that $\nu_t >\gamma $ when $\bx^\toptzero \neq \bx^\toptone$ in Algorithm~\ref{alg:trust_region1}.
If $\Delta_t \leq \epsilon/L$ for all $t = z, z+1, \ldots, l$, we have
\begin{equation}\label{equation:glo_conv_trs_othe1}
f(\bx^{(z)}) - f(\bx^{(l+1)}) \geq \gamma c_1 \epsilon \sum_{t=z, \bx^\toptzero \neq \bx^\toptone}^{l} \Delta_t \geq \gamma c_1 \epsilon \widetildeR 
= \gamma c_1 \epsilon \min \left( \frac{\epsilon}{\beta}, R \right).
\end{equation}
where the last inequality follows from the fact that $\bx^\toptzero \in \sB[\bx^{(z)}, \widetildeR] $ for $t=z, z+1, \ldots,l$.
Otherwise, we have $\Delta_t > \epsilon/L$ for some $t = z, z+1, \ldots, l$. Therefore,
\begin{equation}\label{equation:glo_conv_trs_othe2}
f(\bx^{(z)}) - f(\bx^{(l+1)}) \geq \gamma c_1 \epsilon \frac{\epsilon}{L}.
\end{equation}
Since the sequence $\{f(\bx^\toptzero)\}_{t=1}^{\infty}$ is decreasing and bounded below, we have that
\begin{equation}
f(\bx^\toptzero) \downarrow f^*, \quad \text{for some } f^* > -\infty.
\end{equation}
Combining the two bounds in \eqref{equation:glo_conv_trs_othe1} and \eqref{equation:glo_conv_trs_othe2}, and the definition of $\epsilon=\frac{1}{2} \normtwobig{\bg^{(z)}}$, we have 
\begin{align*}
f(\bx^{(z)}) - f^* &\geq f(\bx^{(z)}) - f(\bx^{(l+1)}) 
\geq \gamma c_1 \epsilon \min \left( \frac{\epsilon}{L}, \frac{\epsilon}{\beta}, R \right) \\
&= \frac{1}{2} \gamma c_1 \normtwobig{\bg^{(z)}} \min \left( \frac{\normtwo{\bg^{(z)}}}{2L}, \frac{\normtwo{\bg^{(z)}}}{2\beta}, R \right) > 0.
\end{align*}
Since $f(\bx^{(z)}) - f^* \downarrow 0$, we must have $\bg^{(z)} \rightarrow \bzero$, establishing the desired result.
\end{proof}

The above two theorems demonstrate that, unlike Newton-type methods, the trust region method possesses global convergence properties and therefore has less stringent requirements on the choice of the initial iterate. In contrast, the convergence of Newton's method is highly dependent on the choice of the initial iterate (see Theorem~\ref{theorem:conv_classNewton} for more details).

\index{Conjugate gradient}
\index{Hessian-orthogonal}
\section{Conjugate Gradient Method}\label{section:conjugate-descent}
We have discussed the gradient descent (GD, or steepest descent) method (Section~\ref{section:gradient-descent-all}) and Newton's method (Section~\ref{section:new_methods}) previously. In this section we introduce the \textit{conjugate gradient  (CG)} method, which is one between the gradient descent method and  Newton's method.
We have shown that the pure GD (employing the negative gradient as the descent direction) can move back and forth in a zigzag pattern when applied in a quadratic bowl (a ravine-shaped loss curve, see example in Figure~\ref{fig:quadratic_vanillegd_contour8}). 
This zigzag behavior becomes more pronounced if the stepsize is guaranteed by exact line search strategies (Section~\ref{section:line_search}), since the gradient is orthogonal to the previous update step (Lemma~\ref{lemm:linear-search-orghonal} and example in Figure~\ref{fig:conjguatecy_zigzag2}). The choice of orthogonal descent directions fails to preserve the minimum along the previous search directions, and the line search will undermine the progress we have already achieved in the direction of the previous line search. This results in the zigzag pattern in  movement.

Instead of favoring a descent direction that is orthogonal to previous search direction (i.e., $\bd^\toptoneTOP \bd^\toptzero=0$), the {conjugate descent} deflects the direction of the gradient descent method by selecting a search direction that is \textit{Hessian-orthogonal} (i.e., $\bd^\toptoneTOP\bH\bd^\toptzero=0$, or conjugate with respect to $\bH$) or by adding to the steepest descent direction a positive multiple of the direction used in the last step. 
This method only requires the first-order derivatives but overcomes the gradient descent method's shortcoming of slow convergence. At the same
time, the method need not save and compute the second-order derivatives which are needed  Newton's method. In particular, since it does not require
the Hessian matrix or its approximation, it is widely used to solve large scale optimization problems.
However, since the analysis of this method needs the second-order information, so we include this methods in the second-order chapter.

\begin{figure}[h]
\centering  
\vspace{-0.35cm} 
\subfigtopskip=2pt 
\subfigbottomskip=2pt 
\subfigcapskip=-5pt 
\subfigure[Surface plot: Hessian with different eigenvalues.]{\label{fig:conjugate_tile_A_orthogonal_diffeigenv}
\includegraphics[width=0.481\linewidth]{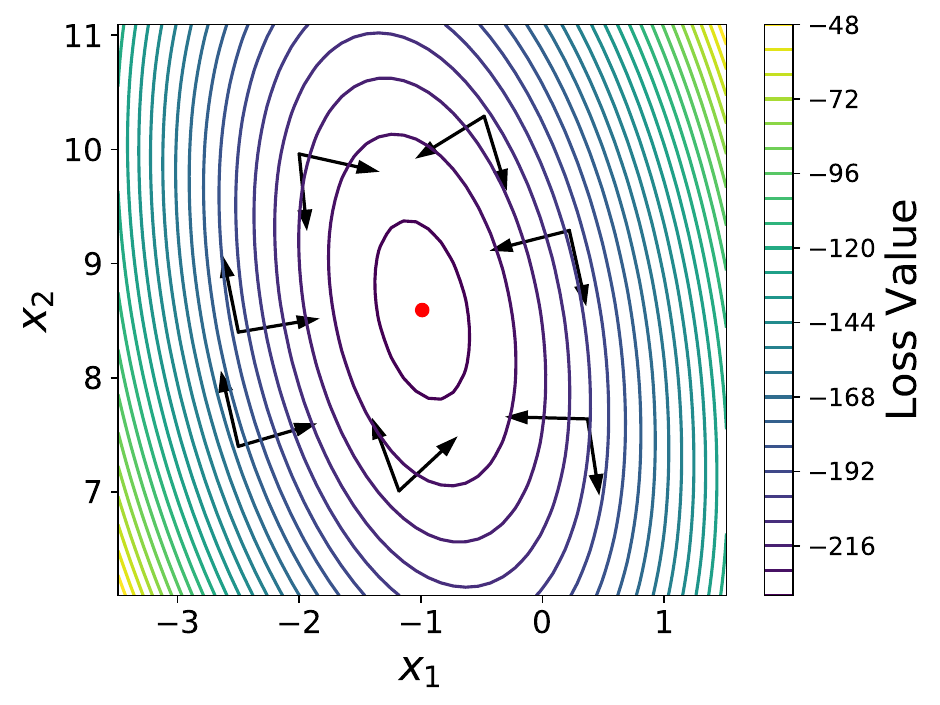}}
\subfigure[Surface plot: Hessian with same eigenvalues.]{\label{fig:conjugate_tile_A_orthogonal_sameeigenv}
\includegraphics[width=0.481\linewidth]{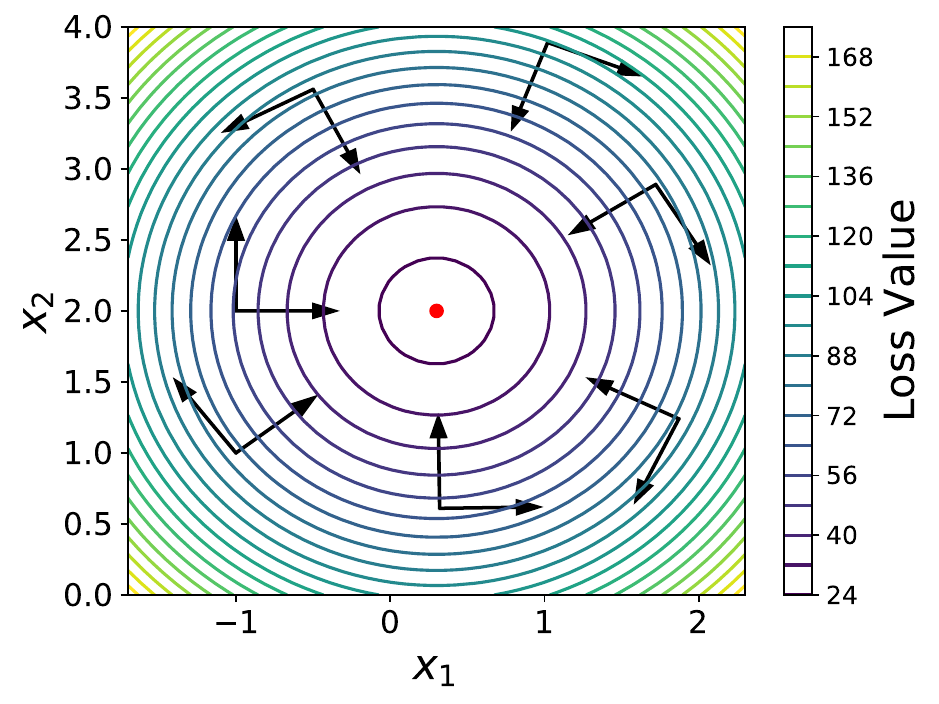}}
\caption{Illustration of $\bH$-orthogonal for different Hessian matrices in two-dimensional case:
$\bH=\footnotesize\begin{bmatrix}
	40 & 5 \\ 7 & 10
\end{bmatrix}$ for Fig~\ref{fig:conjugate_tile_A_orthogonal_diffeigenv}
and 
$\bH=\footnotesize\begin{bmatrix}
	40 & 0 \\ 0 & 40
\end{bmatrix}$ 
for Fig~\ref{fig:conjugate_tile_A_orthogonal_sameeigenv}. The $\bH$-orthogonal pairs are orthogonal when $\bH$ has identical eigenvalues.}
\label{fig:conjugate_tile_A_orthogonals}
\end{figure}

This choice of the conjugate gradient method ensures that the movement is compensated by the curvature information of the loss function. In Figure~\ref{fig:conjugate_tile_A_orthogonals}, we show examples of Hessian-orthogonal pairs when the eigenvalues of the Hessian matrix $\bH$ are different or identical. When the eigenvalues of the Hessian matrix are the same, the Hessian-orthogonality reduces to standard orthogonal cases (this can be shown by the spectral decomposition of the Hessian matrix, where the orthogonal transformation does not alter the orthogonality \citep{lu2021numerical}).

We now provide the formal definition of conjugacy as follows:
\begin{definition}[Conjugacy]\label{definition:conjugacy}
Let $\bA\in \real^{n\times n}$ be a positive definite matrix. Then the vectors $\bu, \bv\in \real^n$ are called  \textit{conjugate} with respect to $\bA$ (a.k.a., \textit{$\bA$-conjugate} or simply \textit{conjugate}) if $\bu,\bv\neq \bzero$ and $\bu^\top\bA\bv = 0$.  
The definition can be extended to a set of $m\leq n$ vectors $\bd_1, \bd_2, \ldots, \bd_m\in \real^n$. If 
$$
\bd_i^\top\bA\bd_j = 0, \quad \forall i\neq j,
$$
then the vectors $\bd_1, \bd_2, \ldots, \bd_m$ are called $\bA$-conjugate or simply conjugate.
\end{definition}
Note that if $\bd_1, \bd_2, \ldots, \bd_m$ are  conjugate, then they are linearly independent (see Problem~\ref{prob:indep_conjugate}). When $\bA=\bI$, then the conjugacy reduces to the usual orthogonality.

To the end of this section, we determine a descent direction that is conjugate to the previous search direction with respect to the Hessian matrix $\bH$ for the conjugate gradient method, ensuring that the new update step will not undo the progress made in the previous directions:
$$
\textbf{(CG)}:\qquad \bd^\toptzero = -\nabla f(\bx^\toptzero) + \beta_t \bd^\toptminus,
$$
where $\beta_t$ is a coefficient controlling how much the previous direction would add back to the current search direction.
Three commonly used methods to compute the coefficient are as follows \citep{hestenes1952methods, fletcher1964function}:
\begin{subequations}\label{equation:cg_all_forms}
\begin{align}
	\text{Fletcher-Reeves:\gap } &\beta_t^F = \frac{\nabla f(\bx^\toptzero)^\top \nabla f(\bx^\toptzero)}{\nabla f(\bx^\toptminus)^\top \nabla f(\bx^\toptminus)},\\
	\text{Polak-Ribi\`ere:\gap } &\beta_t^P = \frac{\Big(\nabla f(\bx^\toptzero) - \nabla f(\bx^\toptminus)\Big)^\top \nabla f(\bx^\toptzero)}{\nabla f(\bx^\toptminus)^\top \nabla f(\bx^\toptminus)},\\
	\text{Hestenes-Stiefel:\gap } &\beta_t^H = \frac{\Big(\nabla f(\bx^\toptzero) - \nabla f(\bx^\toptminus)\Big)^\top \nabla f(\bx^\toptzero)}{\Big(\nabla f(\bx^\toptzero) - \nabla f(\bx^\toptminus)\Big)^\top
		\bd^\toptminus}.
\end{align}
\end{subequations}
In the case of a quadratic loss function, the conjugate gradient ensures that the gradient along the previous direction does not increase in magnitude \citep{shewchuk1994introduction, nocedal1999numerical, iserles2009first, goodfellow2016deep}. The full procedure of the conjugate gradient method is formulated in Algorithm~\ref{alg:conjugate_descent}, where it is observed that the first step of conjugate gradient is identical to a step of steepest descent when the stepsize is calculated by exact line search since $\beta_1=0$.


\begin{algorithm}[h] 
\caption{Fletcher-Reeves Conjugate Gradient}
\label{alg:conjugate_descent}
\begin{algorithmic}[1]
\Require A differentiable function $f(\bx)$;
\State {\bfseries Input:} Initialize $\bx^\topone$, $\bd^\topzero =\bzero $, and $\bg^\topzero = -\bd^\topzero+\bepsilon$;
\For{$t=1,2,\ldots$ } 
\State Compute gradient $\bg^\toptzero = \nabla f(\bx^\toptzero)$;
\State Compute coefficient $\beta_{t} = \frac{ \bg^\toptzeroTOP \bg^\toptzero}{\bg^\toptminusTOP \bg^\toptminus}$ (\text{Fletcher-Reeves});
\State Compute descent direction $\bd^\toptzero = -\bg^\toptzero +\beta_{t}  \bd^\toptminus$;
\State Fixed stepsize $\eta_t=\eta$ or find it by line search: $\eta_t = \arg\min f(\bx^\toptzero + \eta \bd^\toptzero)$;
\State Apply update $\bx^\toptone \leftarrow  \bx^\toptzero + \eta_t\bd^\toptzero$;
\EndFor
\State {\bfseries Return:} resulting parameters $\bx^\toptzero$, and the loss $f(\bx^\toptzero)$.
\end{algorithmic}
\end{algorithm}

\begin{figure}[h]
\centering  
\vspace{-0.35cm} 
\subfigtopskip=2pt 
\subfigbottomskip=2pt 
\subfigcapskip=-5pt 
\subfigure[GD, fixed $\eta=0.08$.]{\label{fig:cgm_conjugate8}
\includegraphics[width=0.23\linewidth]{./imgs/steepest_gd_mom-0_lrate-8.pdf}}
\subfigure[GD with line search.]{\label{fig:cgm_zigzag2}
\includegraphics[width=0.23\linewidth]{./imgs/steepest_gd_bisection.pdf}}
\subfigure[Conjugate descent, fixed $\eta=0.06$.]{\label{fig:cgm_conjugate2}
\includegraphics[width=0.23\linewidth]{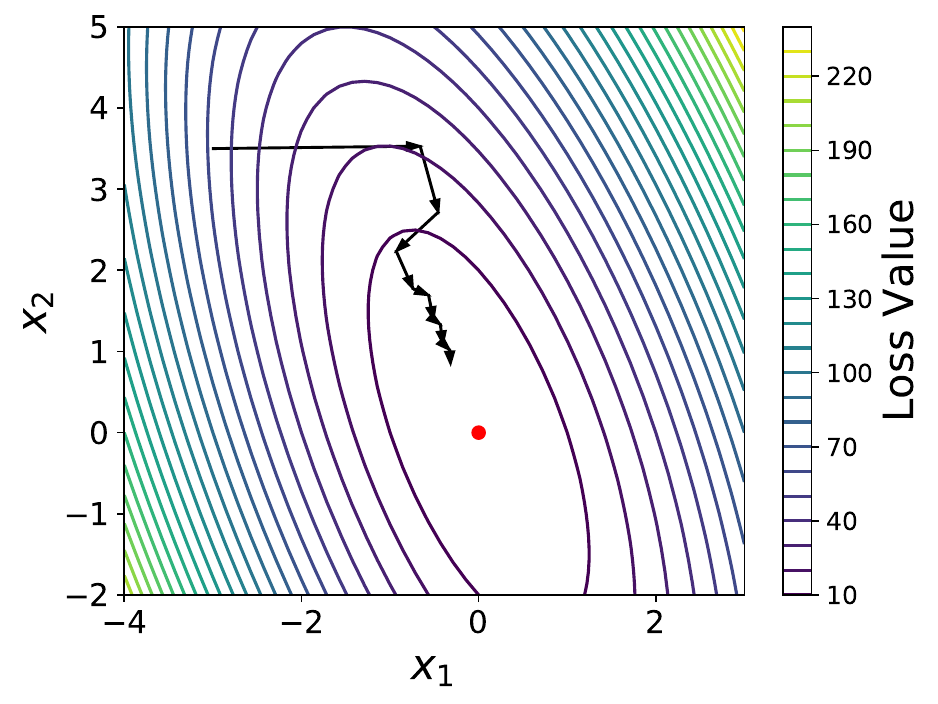}}
\subfigure[Conjugate descent, exact line search.]{\label{fig:cgm_conjugate3}
\includegraphics[width=0.23\linewidth]{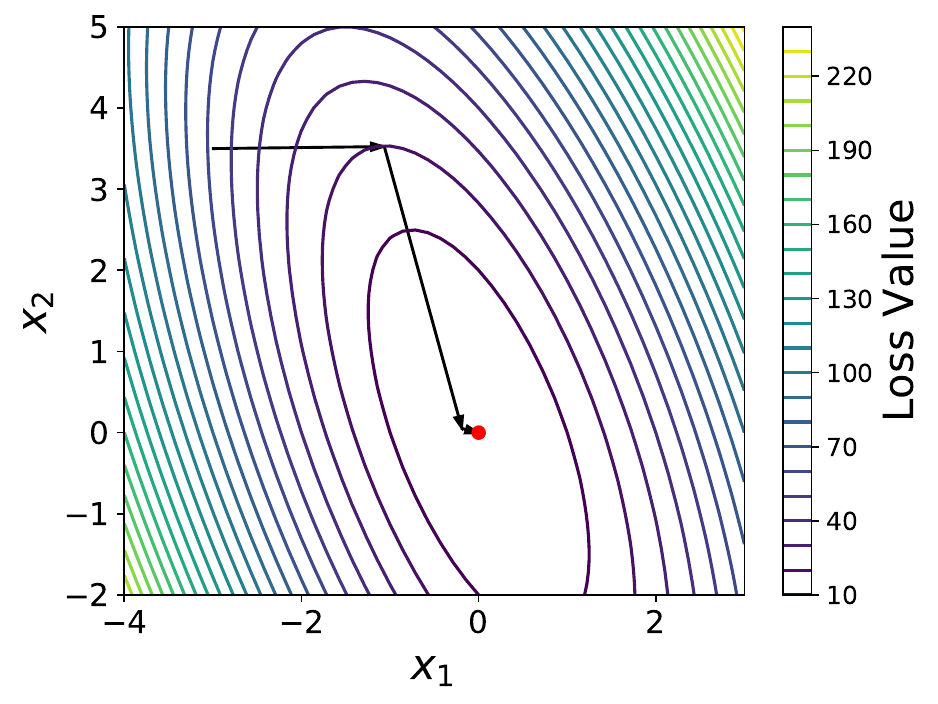}}
\caption{Illustration for the  GD with a fixed stepsize, GD with line search, and CG of quadratic form with $\bA=\footnotesize\begin{bmatrix}
20 & 7 \\ 5 & 5
\end{bmatrix}$, $\bb=\bzero$, and $c=0$. The starting point to descent is $\bx_1=[-3,3.5]^\top$.}
\label{fig:cgm-zigzzag}
\end{figure}

In the subsequent sections, we base our derivation of the conjugate gradient on the assumption of a symmetric positive definite $\bA$;  
however, it can be readily adapted to asymmetric matrices. A comparison among GD with a fixed stepsize, GD with line search, and conjugate gradient is shown in Figure~\ref{fig:cgm-zigzzag}, where we observe that the updates in CG have less zigzag pattern than GD with a fixed stepsize and GD with line search.

\index{Quadratic form}
\index{Quadratic model}
\index{Conjugate direction method}
\subsection{Quadratic Model in Conjugate Direction (CD) Method}
Following the discussion of the quadratic form in GD (Section~\ref{section:quadratic_vanilla_GD}), the quadratic form in GD with line search (Section~\ref{section:quadratic-in-steepestdescent}), and the quadratic form in momentum (Section~\ref{section:quadratic-in-momentum}), we turn our attention to  the quadratic form in CG. 
As a beginning, we proceed by an exploration of the \textit{conjugate direction (CD)} method, where the distinction between them will become evident in the subsequent discussions. 
According to the definition of conjugacy (Definition~\ref{definition:conjugacy}), it is easy to show that any set of vectors $\{\bd^{(1)}, \bd^{(2)}, \ldots, \bd^{(n)}\}\in \real^n$ satisfying this property with respect to the symmetric positive definite Hessian matrix $\bH=\frac{1}{2}(\bA^\top+\bA)$: 
\begin{equation}\label{equation:conjguate-definis}
\bd^{(i)\top} \bH\bd^{(j)}, \gap \forall i\neq j,
\end{equation}
is also linearly independent. That is, the set span the entire $\real^n$ space:
$$
\spn\{\bd^{(1)}, \bd^{(2)}, \ldots, \bd^{(n)}\} = \real^n.
$$
Given the initial parameter $\bx^\topone$ and a set of \textit{conjugate directions} $\{\bd^{(1)}, \bd^{(2)}, \ldots, \bd^{(n)}\}$ (defined in \eqref{equation:conjguate-definis}), the update at iteration $t$ is given by 
\begin{equation}\label{equation:conjugate_direction-update}
\bx^\toptone \leftarrow \bx^\toptzero +\eta_t \bd^\toptzero,
\end{equation}
where $\eta_t$ is the stepsize at time $t$ and is obtained by minimizing the one-dimensional quadratic function $\phi(\eta)=f(\bx^\toptzero+\eta\bd^\toptzero)$, as presented in \eqref{equation:eta-gd-steepest}:
\begin{equation}\label{equation:conjugate_direction-update2}
\eta_t = - \frac{\bd^\toptzeroTOP \bg^\toptzero}{ \bd^\toptzeroTOP \bA\bd^\toptzero } \gap \text{with}\gap  \bg^\toptzero \triangleq \nabla f(\bx^\toptzero).
\end{equation}
Thus, we call it the \textit{conjugate direction method with exact line search}.
Then, we can establish the following theorem that the updates following from the conjugate directions will converge in $n$ steps (the dimension of the parameter) when $\bA$ is symmetric positive definite (the Hessian $\bH=\bA$ in this case).
\begin{theorem}[Converge in $n$ Steps]\label{theorem:conjudate_direc_d-steps}
Let $f(\bx)=\frac{1}{2}\bx^\top\bA\bx-\bb^\top\bx+c$ with positive definite $\bA$.
Then, for any initial parameter $\bx^\topone$, the sequence $\{\bx^\toptzero\}$ generated by the conjugate direction method with exact line search in \eqref{equation:conjugate_direction-update} converges to the optimal point $\bx^*$ in at most $n$ steps.
\end{theorem}
\begin{proof}[of Theorem~\ref{theorem:conjudate_direc_d-steps}]
Since conjugate directions $\{\bd^{(1)}, \bd^{(2)}, \ldots, \bd^{(n)}\}$ span the entire $\real^n$ space, the initial error vector $\be^\topone \triangleq \bx^\topone - \bx^*$ (Definition~\ref{definition:error-gd-}) then can be expressed as a linear combination of the conjugate directions:
\begin{equation}\label{equation:dsteps-symmetric}
	\be^\topone = \bx^\topone-\bx^* = \gamma_1\bd^{(1)} + \gamma_2\bd^{(2)}+\ldots+\gamma_n\bd^{(n)}.
\end{equation}
Using the conjugacy in \eqref{equation:conjguate-definis} and the fact that $\bH=\bA$ when $\bA$ is positive definite, the $\gamma_i$'s can be obtained by the following equation:
$$
\begin{aligned}
\bd^\toptzeroTOP \bH\be^\topone &= \sum_{i=1}^{n} \gamma_i \bd^\toptzeroTOP \bH\bd^{(i)}=\gamma_t \bd^\toptzeroTOP \bH\bd^\toptzero \\
\implies \gamma_t &= \frac{\bd^\toptzeroTOP \bH\be^\topone}{\bd^\toptzeroTOP \bA\bd^\toptzero}= \frac{\bd^\toptzeroTOP \bH(\be^\topone+ \sum_{i=1}^{t-1}\eta_i \bd^{(i)} )}{\bd^\toptzeroTOP \bH\bd^\toptzero} 
=\frac{\bd^\toptzeroTOP \bH\be^\toptzero}{\bd^\toptzeroTOP \bH\bd^\toptzero}.
\end{aligned}
$$
When $\bA$ is  symmetric and nonsingular, we have $\be^\toptzero \triangleq \bx^\toptzero-\bx^*  = \bx^\toptzero-\bA^{-1}\bb$. It can be shown that $\gamma_t$ is equal to $\frac{\bd^\toptzeroTOP \bg^\toptzero}{\bd^\toptzeroTOP \bA\bd^\toptzero}$. This is exactly the same form (in magnitude) as the stepsize at iteration $t$ in GD with line search: $\gamma_t = -\eta_t$ (see \eqref{equation:eta-gd-steepest} or \eqref{equation:conjugate_direction-update2}). Substituting into \eqref{equation:dsteps-symmetric}, it follows that 
$$
\bx^* = \bx^\topone + \eta_1\bd^{(1)}+\eta_2\bd^{(2)}+\ldots +\eta_n\bd^{(n)}.
$$
Alternatively, we have updates by \eqref{equation:conjugate_direction-update} that
$$
\begin{aligned}
\bx^{(n+1)} &= \bx^{(n)} + \eta_n\bd^{(n)} 
=\big(\bx^{(n-1)} +\eta_{n-1}\bx^{(n-1)}\big) + \eta_n\bd^{(n)} \\
&= \ldots \\
&= \bx^\topone + \eta_1\bd^{(1)} + \eta_2\bd^{(2)}+\ldots+\eta_n\bd^{(n)} = \bx^*,
\end{aligned}
$$
which completes the proof.
\end{proof}

The above theorem states the conjudate direction given by \eqref{equation:conjugate_direction-update} converges in $n$ steps, i.e., $\bx^{(n+1)}$ minimizes the quadratic function $f(\bx)=\frac{1}{2}\bx^\top\bA\bx-\bb^\top\bx+c$ over the entire space $\real^n$. Furthermore, we can prove at each iteration $t\leq d$, the update $\bx^\toptone$ minimizes the quadratic function over a subspace of $\real^n$. This is known as the \textit{expanding subspace minimization theorem} or the \textit{principal theorem of the conjugate direction method}.
\begin{theorem}[Expanding Subspace Minimization]\label{theorem:expanding_subspace_minimization}
Let $f(\bx)=\frac{1}{2}\bx^\top\bA\bx-\bb^\top\bx+c$ with positive definite $\bA$.
For any initial parameter $\bx^\topone$, let the sequence $\{\bx^\toptzero\}$ be generated by the conjugate direction method with exact line search in  \eqref{equation:conjugate_direction-update}. Then it follows that 
\begin{equation}\label{equation:expanding_subspace_minimization_zero}
\bg^\toptoneTOP \bd^{(i)}=0, \gap \forall i=1,2,\ldots, t, \text{ and } t\in \{1,2,\ldots, n\},
\end{equation}
where $\bg^\toptzero \triangleq \nabla f(\bx^\toptzero)= \bA\bx^\toptzero - \bb$,
and $\bx^\toptone$ is the minimizer of $f(\bx)$ over the subspace
\begin{equation}\label{equation:space_d_t}
	\mathbb{D}_t=\left\{\bx \mid \bx=\bx^\topone + \spn\{\bd^{(1)}, \bd^{(2)}, \ldots, \bd^\toptzero\}\right\}.
\end{equation}
\end{theorem}
\begin{proof}[of Theorem~\ref{theorem:expanding_subspace_minimization}]
We first prove $\bg^\toptoneTOP \bd^{(i)}$ by induction. When $t=1$, since $\eta_1$ is obtained to minimize $\phi(\eta)=f(\bx^\topone+\eta\bd^{(1)})$, by Lemma~\ref{lemm:linear-search-orghonal}, we have $\bg^{(2)} = \nabla f(\bx^{(2)})$ that is orthogonal to $\bd^{(1)}$. Suppose now for general $t-1$, the induction hypothesis is satisfied with $\bg^\toptzeroTOP\bd^{(i)}=0$ for $i=0,1,\ldots, t-1$. The $\bg^\toptzero$ has the following update
\begin{equation}\label{equation:conjguate-redsidual-update}
\begin{aligned}
\bg^\toptone &= \bA\bx^\toptone-\bb  =  \bA (\bx^\toptzero+\eta_t\bd^\toptzero) -\bb 
= \bg^\toptzero + \eta_t\bA\bd^\toptzero.
\end{aligned}
\end{equation}
By conjugacy and the induction hypothesis, we have $\bg^\toptoneTOP  \bd^{(i)}=0$ for $i=\{0,1,\ldots,t-1\}$. If we further prove this is also true for $\bg^\toptoneTOP \bd^\toptzero$, we complete the proof. This follows again from Lemma~\ref{lemm:linear-search-orghonal}, the current gradient is orthogonal to the previous search direction $\bd^\toptzero$ under exact line search strategies.

For the second part, we define $g(\bm{\eta}) = f(\bx^\topone+\eta_1\bd^{(1)}+\eta_2\bd^{(2)}+\ldots +\eta_t\bd^\toptzero)$, which is a strictly convex quadratic function over $\bm{\eta}=[\eta_1, \eta_2, \ldots, \eta_t]^\top$ such that 
$$
\frac{\partial g(\bm{\eta})}{\partial \eta_i} = 0, \gap \forall i=1,2,\ldots, t.
$$
This implies 
$$
\underbrace{\nabla f(\bx^\topone +\eta_1\bd^{(1)}+\eta_2\bd^{(2)}+\ldots +\eta_t\bd^\toptzero}_{\nabla f(\bx^\toptone)})^\top \bd^{(i)} = 0,  \gap \forall i=1,2,\ldots, t.
$$
That is, $\bx^\toptone\in \big\{\bx \mid \bx=\bx^\topone + \spn\{\bd^{(1)}, \bd^{(2)}, \ldots, \bd^\toptzero\}\big\}$ is the minimizer of $f(\bx)$.
\end{proof}
This theorem, while straightforward, is crucial as it underpins all conjugate direction methods. It's important to highlight that when exact line search is employed, these methods fulfill condition \eqref{equation:expanding_subspace_minimization_zero} and exhibit the property of quadratic termination. In essence, combining conjugacy with exact line search leads to the efficient resolution of quadratic problems, emphasizing the practical significance of this approach in optimization.

\index{Quadratic form}
\subsection{Quadratic Model in Conjugate Gradient (CG) Method}
We have mentioned that the conjugate gradient (CG) method differs from the conjugate descent (CD) method. The distinction lies  in the fact  that the CG method computes a new vector $\bd^\toptone$ using only the previous vector $\bd^\toptzero$ rather than the entire sequence $\{\bd^{(1)}, \bd^{(2)}, \ldots, \bd^\toptzero\}$. 
And the resulting $\bd^\toptone$ will automatically be conjugate to the sequence in this sense. In the CG method, each search direction $\bd^\toptzero$ is chosen to be a linear combination of negative gradient $-\bg^\toptzero$ (the steepest descent direction) and the previous direction $\bd^\toptminus$:
\begin{equation}\label{equation:cd_gradient_direction}
\bd^\toptzero =  -\bg^\toptzero + \beta_t \bd^\toptminus.
\end{equation}
Choosing $\beta_t$ to ensure the conjugacy $\bd^\toptzeroTOP \bA\bd^\toptminus = 0$ yields that 
$$
\beta_t = \frac{\bg^\toptzeroTOP\bA\bd^\toptminus}{\bd^\toptminusTOP \bA\bd^\toptminus}
.
$$
This choice of $\beta_t$ and $\bd^\toptzero$ actually results in the conjugate sequence $\{\bd^{(1)}, \bd^{(2)}, \ldots, \bd^\toptzero\}$. To see this, we first provide the definition of the \textit{Krylov subspace of degree $t$ for vector $\bv$ with respect to matrix $\bA$}:
$$
\mathcal{K}(\bv; t) \triangleq \spn\{\bv, \bA\bv, \ldots, \bA^{t-1}\bv\}.
$$

\begin{theorem}[Converge in $n$ Steps]\label{theorem:conjudate_CD_d-steps}
Let $f(\bx)=\frac{1}{2}\bx^\top\bA\bx-\bb^\top\bx+c$ with positive definite $\bA$.
For any initial parameter $\bx^\topone$, the sequence $\{\bx^\toptzero\}$ generated by the conjugate descent method with exact line search and search directions generated by \eqref{equation:cd_gradient_direction}, converges to the optimal point $\bx^*$ in at most $n$ steps.
The result follows from the following claims:
\begin{subequations}
\begin{align}
\bg^\toptzeroTOP \bg^{(i)} &= 0, \gap \text{for $i=\{1,2,\ldots, t-1\}$};\label{equation:conjudate_CD_d1}\\
\spn\{\bg^{(1)}, \bg^{(2)}, \ldots, \bg^\toptzero\} &= \spn\{\bg^{(1)}, \bA\bg^{(1)}, \ldots, \bA^{t-1}\bg^{(1)}\}=\mathcal{K}(\bg^{(1)}; t);\label{equation:conjudate_CD_d2}\\
\spn\{\bd^{(1)}, \bd^{(2)}, \ldots, \bd^\toptzero\} &= \spn\{\bg^{(1)}, \bA\bg^{(1)}, \ldots, \bA^{t-1}\bg^{(1)}\}=\mathcal{K}(\bg^{(1)}; t);\label{equation:conjudate_CD_d3}\\
\bd^\toptzeroTOP\bA\bd^{(i)} &= 0, \gap \text{for $i=\{1,2,\ldots, t-1\}$},\label{equation:conjudate_CD_d4}
\end{align}
where \eqref{equation:conjudate_CD_d4} indicates the sequence $\{\bd^\toptzero\}$ is conjugate. 
Using the definition of search directions in \eqref{equation:cd_gradient_direction} and the above results, it also follows that 
\begin{equation}\label{equation:conjudate_CD_d5}
\bd^{(i)\top}\bg^{(i)} = -\bg^{(i)\top}\bg^{(i)} , \gap \text{for $i=\{1,2,\ldots, t\}$}.
\end{equation}
\end{subequations}
\end{theorem}
\begin{proof}[of Theorem~\ref{theorem:conjudate_CD_d-steps}]
The proof proceeds through induction.  The equations hold trivially when $t=1$. 
Assuming that  for $t$, \eqref{equation:conjudate_CD_d2}, \eqref{equation:conjudate_CD_d3}, and \eqref{equation:conjudate_CD_d4} hold as well; if we can show the  equations still hold for $t+1$, then we complete the proof. By the induction hypothesis, we have 
$$
\bg^\toptzero \in \spn\{\bg^{(1)}, \bA\bg^{(1)}, \ldots, \bA^{t-1}\bg^{(1)}\}, 
\gap 
\bd^\toptzero \in \spn\{\bg^{(1)}, \bA\bg^{(1)}, \ldots, \bA^{t-1}\bg^{(1)}\}.
$$
Left multiplying by $\bA$, it follows that
\begin{equation}\label{equation:dstesps_induction00}
\bA\bd^\toptzero \in \spn\{\bA\bg^{(1)}, \bA^2\bg^{(1)}, \ldots, \bA^{t}\bg^{(1)}\}.
\end{equation}
Since 
\begin{equation}\label{equation:dstesps_induction11}
\begin{aligned}
\bg^\toptone &= \bA\bx^\toptone-\bb 
=\bA(\bx^\toptzero+\eta_t \bd^\toptzero) -\bb = \bg^\toptzero + \eta_t\bA\bd^\toptzero,
\end{aligned}
\end{equation}
then we have 
\begin{equation}\label{equation:dstesps_induction2}
	\bg^\toptone\in \spn\{\bg^{(1)},\bA\bg^{(1)}, \bA^2\bg^{(1)}, \ldots, \bA^{t}\bg^{(1)}\}.
\end{equation}
Combining \eqref{equation:dstesps_induction2} and \eqref{equation:conjudate_CD_d2}, we have 
$$
\spn\{\bg^{(1)}, \bg^{(2)}, \ldots, \bg^\toptzero, \bg^\toptone\} \subset \spn\{\bg^{(1)}, \bA\bg^{(1)}, \ldots, \bA^{t-1}\bg^{(1)}, \bA^t\bg^{(1)}\}.
$$
To see the reverse inclusion, by \eqref{equation:conjudate_CD_d3}, it follows that 
$$
\bA^t \bg^{(1)} = \bA(\bA^{t-1}\bg^{(1)}) \in \spn\{\bA\bd^{(1)}, \bA\bd^{(2)}, \ldots, \bA\bd^\toptzero\}.
$$
Again, by \eqref{equation:dstesps_induction11}, we have $\bA\bd^\toptzero = (\bg^\toptone-\bg^\toptzero)/\eta_t$. Therefore,
$$
\bA^t \bg^{(1)} \in \spn\{\bg^{(1)}, \bg^{(2)}, \ldots, \bg^\toptzero, \bg^\toptone\} .
$$
Combining with \eqref{equation:conjudate_CD_d2}, we have 
$$
\spn\{\bg^{(1)}, \bA\bg^{(1)}, \ldots, \bA^{t-1}\bg^{(1)}, \bA^{t}\bg^{(1)}\}
\subset 
\spn\{\bg^{(1)}, \bg^{(2)}, \ldots, \bg^\toptzero, \bg^\toptone\}.
$$
Therefore, \eqref{equation:conjudate_CD_d2} holds for $t+1$. \eqref{equation:conjudate_CD_d3} follows similarly and also holds for $t+1$.
To see how \eqref{equation:conjudate_CD_d4} holds for $t+1$, we have 
$$
\bd^\toptoneTOP\bA\bd^{(i)} = (-\bg^\toptone + \beta_{t+1}\bd^\toptzero)^\top\bA\bd^{(i)}.
$$
By Theorem~\ref{theorem:expanding_subspace_minimization}, we have 
\begin{equation}\label{equation:dstesps_induction_argue41}
	\bg^\toptoneTOP \bd^{(i)}=0 \text{ for } i\in \{1,2,\ldots, t\}.
\end{equation}
Furthermore, by \eqref{equation:dstesps_induction00} and \eqref{equation:conjudate_CD_d3}, we have 
\begin{equation}\label{equation:dstesps_induction_argue42}
	\bA\bd^{(i)} \in \spn\{\bA\bg^{(1)}, \bA^2\bg^{(1)}, \ldots, \bA^{i}\bg^{(1)}\}\subset 
	\spn\{\bd^{(1)}, \bd^{(2)}, \ldots, \bd^{(i)}, \bd^{(i+1)}\}.
\end{equation}
Combining \eqref{equation:dstesps_induction_argue41} and \eqref{equation:dstesps_induction_argue42}, it then follows that 
$$
\bd^\toptoneTOP\bA\bd^{(i)}=0,\gap \text{ for }i\in \{1,2,\ldots,t-1\}.
$$
We need to further demonstrate $\bd^\toptoneTOP\bA\bd^\toptzero=0$, a result that is evident given the intentional design of the algorithm to satisfy this condition.

To establish the validity of  \eqref{equation:conjudate_CD_d1}, we have $\bd^{(i)} = -\bg^{(i)}+\beta_i\bd^{(i-1)}$. Therefore, $\bg^{(i)} \in \spn\{\bd^{(i)},\bd^{(i-1)}\}$. Furthermore, employing \eqref{equation:dstesps_induction_argue41}, we prove $\bg^\toptzeroTOP\bg^{(i)}=0$ for $i\in\{1,2,\ldots, t-1\}$.
\end{proof}

Note that the equations \eqref{equation:conjudate_CD_d1}, \eqref{equation:conjudate_CD_d4}, and \eqref{equation:conjudate_CD_d5} represent conjugacy of search directions, orthogonality of gradients, and descent property of the search directions, respectively.
Therefore, the CG method developed by \eqref{equation:cd_gradient_direction} that creates conjugate directions $\bd^\toptone\bA\bd^\toptzero=0$ indeed finds a conjugate set $\bd^\toptone\bA\bd^{(i)}=0$ for $i\in\{1,2,\ldots,t\}$. By Theorem~\ref{theorem:conjudate_direc_d-steps}, the CG method thus converges in at most $n$ steps (when $\bA$ is symmetric PD and using exact line search strategies, where the stepsize for iteration $t$ is given in \eqref{equation:conjugate_direction-update2}). The complete procedure is then formulated in Algorithm~\ref{alg:vanilla_conjugate_descent}.

\begin{algorithm}[H] 
\caption{Vanilla Conjugate Gradient Method for Quadratic Models}
\label{alg:vanilla_conjugate_descent}
\begin{algorithmic}[1]
\Require  Symmetric positive definite $\bA\in \real^{n\times n}$;
\State {\bfseries Input:} Initial parameter $\bx^\topone$;
\State {\bfseries Input:} Initialize $\bd^\topzero =\bzero $ and $\bg^\topzero = -\bd^\topzero+\bepsilon$;
\For{$t=1:n$ } 
\State Compute gradient $\bg^\toptzero \leftarrow \nabla f(\bx^\toptzero)$;
\State Compute coefficient $\beta_{t} \leftarrow \frac{ \bg^\toptzeroTOP\bA \bd^\toptminus}{\bg^\toptminusTOP \bA\bg^\toptminus}$;
\State Compute descent direction $\bd^\toptzero \leftarrow -\bg^\toptzero +\beta_{t}  \bd^\toptminus$;
\State Stepsize $\eta_t \leftarrow - \frac{\bd^\toptzeroTOP \bg^\toptzero}{ \bd^\toptzeroTOP \bA\bd^\toptzero }$;
\State Apply update $\bx^\toptone \leftarrow \bx^\toptzero + \eta_t\bd^\toptzero$;
\EndFor
\State {\bfseries Return:} resulting parameters $\bx^\toptzero$, and the loss $f(\bx^\toptzero)$.
\end{algorithmic}
\end{algorithm}

\index{Complexity}
\index{Flops}

To further reduce the complexity of the CG algorithm, we  introduce the notion of floating-point operation (flop) counts. We follow the classical route and count the number of \textit{floating-point operations (flops)} that the algorithm requires. Each addition, subtraction, multiplication, division, and square root is considered one flop. Note that we have the convention that an assignment operation is not counted as one flop.
The calculation of the complexity extensively relies on the complexity of the multiplication of two matrices so that we formulate the finding in the following lemma \citep{lu2021numerical}.
\begin{lemma}[Vector Inner Product Complexity]
Given two vectors $\bv,\bw\in \real^{n}$. The  inner product of the two vectors $\bv^\top\bw$ is given by $\bv^\top\bw=v_1w_1+v_2w_2+\ldots v_nw_n$, involving $n$ scalar multiplications and $n-1$ scalar additions. Therefore, the complexity for the inner product is $2n-1$ flops.
\end{lemma}

The complexity of  matrix multiplications thus relies on the complexity of the inner product.
\begin{lemma}[Matrix Multiplication Complexity]\label{lemma:matrix-multi-complexity}
For matrix $\bA\in\real^{m\times n}$ and $\bB\in \real^{n\times k}$, the complexity of the multiplication $\bC=\bA\bB$ is $mk(2n-1)$ flops.
\end{lemma}
\begin{proof}[of Lemma~\ref{lemma:matrix-multi-complexity}]
We notice that each entry of $\bC$ involves a vector inner product that requires $n$ multiplications and $n-1$ additions. And there are $mk$ such entries, which leads to the conclusion.
\end{proof}

By Theorem~\ref{theorem:conjudate_CD_d-steps}, we can replace the formula for calculating stepsizes with:
$$
\eta_t = - \frac{\bd^\toptzeroTOP \bg^\toptzero}{ \bd^\toptzeroTOP \bA\bd^\toptzero }
\qquad\implies\qquad  
\eta_t = \frac{\textcolor{mylightbluetext}{\bg^\toptzero}^\top \bg^\toptzero}{ \bd^\toptzeroTOP \bA\bd^\toptzero }.
$$
According to \eqref{equation:conjguate-redsidual-update}, it follows that $\eta_t\bA\bd^\toptzero = \bg^\toptone-\bg^\toptzero$. Combining with \eqref{equation:expanding_subspace_minimization_zero} and \eqref{equation:conjudate_CD_d1}, $\beta_t$ can also be expressed as 
$$
\beta_t = \frac{\bg^\toptzeroTOP\bA\bd^\toptminus}{\bd^\toptminusTOP \bA\bd^\toptminus}
\qquad\implies\qquad  
\beta_t 
=
-\frac{\bg^\toptzeroTOP\bg^\toptzero}{\bd^\toptminusTOP \bg^\toptminus}
=\frac{\bg^\toptzeroTOP\bg^\toptzero}{\bg^\toptminusTOP \bg^\toptminus}
.
$$
This reduces the computational complexity from $\mathcalO(4n^2)$ to $\mathcalO(4n)$ flops. The practical CG method for positive definite quadratic models is then outlined in Algorithm~\ref{alg:practical_conjugate_descent}.

\begin{algorithm}[H] 
\caption{Practical Conjugate Gradient Method for Quadratic Models}
\label{alg:practical_conjugate_descent}
\begin{algorithmic}[1]
\Require Symmetric positive definite $\bA\in \real^{n\times n}$;
\State {\bfseries Input:} Initial parameter $\bx^\topone$;
\State {\bfseries Input:} Initialize $\bd^\topzero =\bzero $ and $\bg^\topzero = -\bd^\topzero+\bepsilon$;
\For{$t=1:n$ } 
\State Compute gradient $\bg^\toptzero \leftarrow \nabla f(\bx^\toptzero)$;
\State Compute coefficient $\beta_{t} \leftarrow  \frac{\bg^\toptzeroTOP\bg^\toptzero}{\bg^\toptminusTOP \bg^\toptminus}$; \Comment{set $\beta_1=0$ by convention}
\State Compute descent direction $\bd^\toptzero \leftarrow -\bg^\toptzero +\beta_{t}  \bd^\toptminus$;
\State Stepsize $\eta_t \leftarrow \frac{{\bg^\toptzero}^\top \bg^\toptzero}{ \bd^\toptzeroTOP \bA\bd^\toptzero }$;
\State Apply update $\bx^\toptone \leftarrow \bx^\toptzero + \eta_t\bd^\toptzero$;
\EndFor
\State {\bfseries Return:} resulting parameters $\bx^\toptzero$, and the loss $f(\bx^\toptzero)$.
\end{algorithmic}
\end{algorithm}

\index{Spectral decomposition}
\index{Quadratic form}
\index{Symmetry}
\index{Positive definite}
\subsection{Convergence Analysis for  Positive Definite Quadratic Models}
We further discuss the convergence results of the CG method. 
According to \eqref{equation:conjudate_CD_d3}, there exists a set of $\{\sigma_1,\sigma_2,\ldots,\sigma_t\}$ coefficients such that 
\begin{equation}\label{equation:cg-convergence-xt1}
\begin{aligned}
	\bx^\toptone &=\bx^\topone + \eta_1\bd^{(1)}+\eta_2\bd^{(2)}+\ldots+\eta_t\bd^\toptzero  \\
	&= \bx^\topone + \sigma_1\bg^{(1)}+\sigma_2\bA \bg^{(1)}+\ldots+\sigma_t\bA^{t-1}\bg^{(1)}\\
	&\triangleq \bx^\topone + P^{\textcolor{mylightbluetext}{*}}_{t-1}(\bA)\bg^{(1)},
\end{aligned}
\end{equation}
where $P^{\textcolor{mylightbluetext}{*}}_{t-1}(\bA) \triangleq \sigma_1\bI+\sigma_2\bA +\ldots+\sigma_t\bA^{t-1}$ is a polynomial of degree $t-1$ with coefficients $\{\sigma_1, \sigma_2, \ldots, \sigma_t\}$. 
This polynomial is a special case of a polynomial of degree $t-1$ with random coefficients $\{\omega_1, \omega_2, \ldots, \omega_t\}$, denoted by $P_{t-1}(\bA) = \omega_1\bI+\omega_2\bA +\ldots+\omega_t\bA^{t-1}$. (Note that $P_{t-1}$ can take either a scalar or a matrix as its argument).
Suppose the symmetric positive definite $\bA$ admits the spectral decomposition (Theorem~\ref{theorem:spectral_theorem}):
$$
\bA=\bQ\bLambda\bQ^\top  \in \real^{n\times n} 
\qquad\implies\qquad
\bA^{-1} = \bQ\bLambda^{-1}\bQ^\top,
$$ 
where the columns of $\bQ = [\bq_1, \bq_2, \ldots , \bq_n]$ are eigenvectors of $\bA$ and are mutually orthonormal, and the entries of $\bLambda = \diag(\lambda_1, \lambda_2, \ldots , \lambda_n)$ with $ \lambda_1\geq \lambda_2\geq \ldots\geq \lambda_n>0$ are the corresponding eigenvalues of $\bA$, which are real and ordered by magnitude (the eigenvalues are positive due to the positive definiteness assumption of $\bA$). It then follows that any eigenvector of $\bA$ is also an eigenvector of $P_{t-1}(\bA)$:
$$
P_{t-1}(\bA) \bq_i = P_{t-1}(\lambda_i) \bq_i, \gap \forall i\in \{1,2,\ldots, n\}.
$$
Moreover, since the eigenvectors span the entire space $\real^n$, there exists a set of $\{\nu_1,\nu_2,\ldots, \nu_n\}$ coefficients such that the initial error vector $\be^\topone$ can be expressed as 
$
\be^\topone=\bx^\topone - \bx^* = \sum_{i=1}^{n} \nu_i \bq_i,
$
where $\bx^\topone$ is the initial parameter. Combining \eqref{equation:cg-convergence-xt1} and and the expression of $\be^\topone$, this yields the update of the error vector:
\begin{equation}\label{equation:cg-convergence-xt3}
	\small
\begin{aligned}
\be^\toptone &=\bx^\toptone - \bx^*
=\bx^\topone +  P^{\textcolor{black}{*}}_{t-1}(\bA)\bg^{(1)}-\bx^*
=\bx^\topone +  P^{\textcolor{black}{*}}_{t-1}(\bA)\big(\bA\bx^\topone - \bA\underbrace{\bA^{-1} \bb}_{\bx^*}\big)-\bx^*\\
&=\bx^\topone +  P^{\textcolor{black}{*}}_{t-1}(\bA)\bA(\bx^\topone - \bx^*)-\bx^*
=\Big(\bI+P^{\textcolor{black}{*}}_{t-1}(\bA)\bA\Big) (\bx^\topone - \bx^*)\\
&=\Big(\bI+P^{\textcolor{black}{*}}_{t-1}(\bA)\bA\Big) \sum_{i=1}^{n} \nu_i \bq_i= \sum_{i=1}^{n}\Big(1+ \lambda_i P^{\textcolor{black}{*}}_{t-1}(\bA)\Big) \nu_i\bq_i.
\end{aligned}
\end{equation}
To further discuss the convergence results, we  need to use the notion of \textit{energy norm} for error vector $\norm{\be}_{\bA} = (\be^\top\bA\be)^{1/2}$ as discussed in Section~\ref{section:general-converg-steepest}. 
It can be shown that minimizing $\norm{\be^\toptzero}_{\bA}$ is equivalent to minimizing $f(\bx^\toptzero)$ by \eqref{equation:energy-norm-equivalent}.

\begin{remark}[Polynomial Minimization]
Since we proved in Theorem~\ref{theorem:expanding_subspace_minimization} that $\bx^\toptone$ minimizes $f(\bx)$ over the subspace $\mathbb{D}_t$ defined in \eqref{equation:space_d_t}, it also minimizes the energy norm $\norm{\be}_{\bA}$ over the subspace $\mathbb{D}_t$ at iteration $t$. It then follows that $P^{\textcolor{black}{*}}_{t-1}(\bA)$ minimizes over the space of all possible polynomials of degree $t-1$:
$$
P^{\textcolor{black}{*}}_{t-1}(\bA)
=\mathop{\arg\min}_{P_{t-1}(\bA)}  \norm{\bx^\topone +  P_{t-1}(\bA)\bg^{(1)}-\bx^*}_{\bA}.
$$
\end{remark}
Then the update of the squared energy norm can be obtained by
$$
\begin{aligned}
\normbig{\be^\toptone}_{\bA}^2 &= \be^\toptoneTOP \bA\be^\toptone =\be^\toptoneTOP \left(\sum_{i=1}^{n}\lambda_i \bq_i\bq_i^\top\right) \be^\toptone
 = \sum_{i=1}^{n} \lambda_i (\be^\toptoneTOP\bq_i)^2 \\
&\stackrel{\dag}{=}\sum_{i=1}^{n} \lambda_i \bigg(\bq_i^\top \Big(\sum_{j=1}^{n}\big(1+ \lambda_j P^{\textcolor{black}{*}}_{t-1}(\bA)\big) \nu_j\bq_j\Big)\bigg)^2 
\stackrel{\ddag}{=}\sum_{i=1}^{n}  \bigg(1+ \lambda_i P^{\textcolor{black}{*}}_{t-1}(\lambda_i)\bigg)^2  \lambda_i\nu_i^2   \\
&=\mathop{\min}_{P_{t-1}} \sum_{i=1}^{n}  \bigg(1+ \lambda_i  P_{t-1}(\lambda_i)\bigg)^2  \lambda_i\nu_i^2  
\stackrel{+}{\leq} m_t \sum_{i=1}^{n} \lambda_i\nu_i^2    
\leq m_t \cdot \norm{\be_{1}}_{\bA}^2,
\end{aligned}
$$
where the equality ($\dag$) follows from  \eqref{equation:cg-convergence-xt3}, the equality ($\ddag$) follows by $\bq_i^\top\bq_j=0$ if $i\neq j$,
and the inequality ($+$) follows from  $m_t = \mathop{\min}_{P_{t-1}}\mathop{\max}_{1\leq j\leq n} (1+ \lambda_j P_{t-1}(\lambda_j))^2 $. 
Therefore, the rate of convergence for the CG method is controlled by 
\begin{equation}\label{equation:cg-convergence-xt5}
m_t = \mathop{\min}_{P_{t-1}}\mathop{\max}_{1\leq j\leq n}\big(1+ \lambda_j P_{t-1}(\lambda_j)\big)^2.
\end{equation}

\subsection*{Special Case: $\bA$ Has Only $r$ Distinct Eigenvalues}
\begin{figure}[h]
\centering  
\vspace{-0.35cm} 
\subfigtopskip=2pt 
\subfigbottomskip=2pt 
\subfigcapskip=-5pt 
\subfigure[CG, 2 distinct eigenvlaues. Finish in 2 steps.]{\label{fig:conjugate_specialcases_2eigenvalue}
	\includegraphics[width=0.481\linewidth]{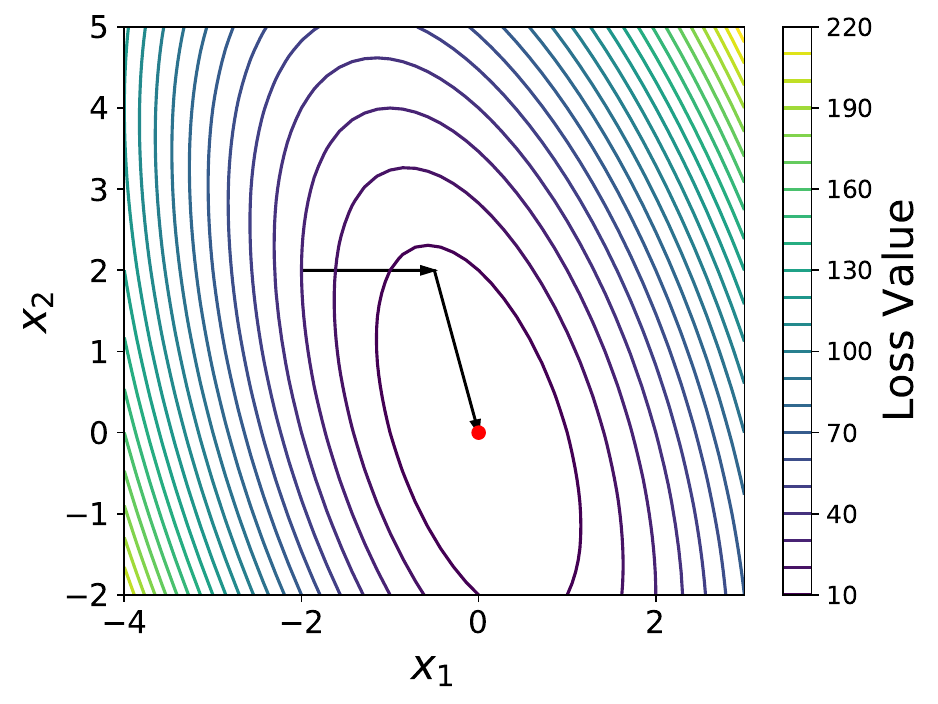}}
\subfigure[CG, 1 distinct eigenvalue. Finish in 1 step.]{\label{fig:conjugate_specialcases_1eigenvalue}
	\includegraphics[width=0.481\linewidth]{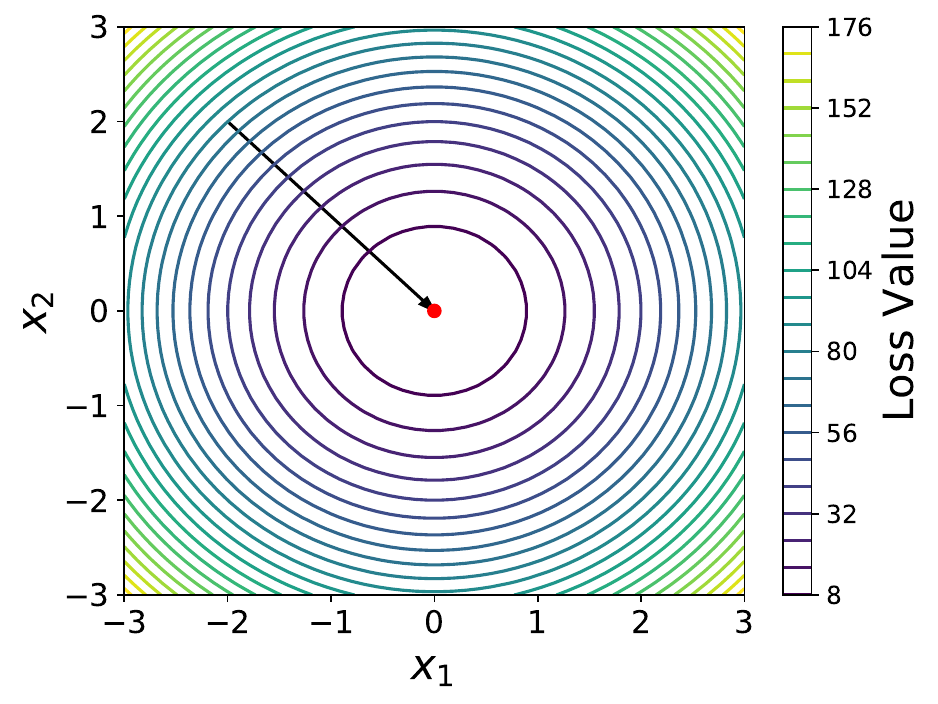}}
\caption{Illustration of special cases for CG with exact line search of quadratic forms. $\bA=\footnotesize\begin{bmatrix}
		20 & 5 \\ 5 & 5
	\end{bmatrix}$, $\bb=\bzero$, $c=0$, and starting point to descent is $\bx^\topone=[-2, 2]^\top$ for Fig~\ref{fig:conjugate_specialcases_2eigenvalue}. $\bA=\footnotesize\begin{bmatrix}
		20 & 0 \\ 0 & 20
	\end{bmatrix}$, $\bb=\bzero$, $c=0$, and starting point to descent is $\bx^\topone=[-2, 2]^\top$ for Fig~\ref{fig:conjugate_specialcases_1eigenvalue}.}
\label{fig:conjugate_specialcases}
\end{figure}

We then consider some special cases. Firstly, we want to show the CG method terminates in exactly $r$ iterations if the symmetric positive definite $\bA$ has only $r$ distinct eigenvalues. To establish this, suppose $\bA$ has distinct eigenvalues $\mu_1<\mu_2<\ldots<\mu_r$ and  define a polynomial $Q_r(\lambda)$ by 
$$
Q_r(\lambda) \triangleq\frac{(-1)^r}{\mu_1\mu_2\ldots\mu_r} (\lambda-\mu_1)(\lambda-\mu_2)\ldots (\lambda-\mu_r),
$$
where $Q_r(\lambda_i)=0$ for $i=\{1,2,\ldots, n\}$ and $Q_r(0)=1$. Therefore, it follows that the polynomial 
$$
R_{r-1}(\lambda) \triangleq \frac{Q_r(\lambda)-1}{\lambda}
$$
is a polynomial of degree $r-1$ with a root at $\lambda=0$. Setting $t-1=r-1$ in \eqref{equation:cg-convergence-xt5}, we have
$$
\begin{aligned}
0&\leq m_{r}=\mathop{\min}_{P_{r-1}}\mathop{\max}_{1\leq j\leq n}(1+ \lambda_j P_{r-1}(\lambda_j))^2
=\mathop{\max}_{1\leq j\leq n}(1+ \lambda_j R_{r-1}(\lambda_j))^2 
= \mathop{\max}_{1\leq j\leq n} Q_r^2(\lambda_i) = 0.
\end{aligned}
$$
As a result, $m_{r}$=0, and $\norm{\be_{r+1}}_{\bA}=0$, implying $\bx_{r+1} = \bx^*$ and the algorithm terminates at iteration $r$. A specific example is shown in Figure~\ref{fig:conjugate_specialcases}, where Figure~\ref{fig:conjugate_specialcases_2eigenvalue} terminates in two steps since it has two distinct eigenvalues, and Figure~\ref{fig:conjugate_specialcases_1eigenvalue} terminates in just one step as it has one distinct eigenvalue.

\subsection*{Closed Form by Chebyshev Polynomials}
It can be shown that \eqref{equation:cg-convergence-xt5} is minimized by a Chebyshev polynomial, given by
$$
1+ \lambda_j P_{t-1}(\lambda_j) = \frac{T_{t}\left( \frac{\lambda_{\max} + \lambda_{\min} - 2\lambda_j}{\lambda_{\max}-\lambda_{\min}} \right) }
{T_{t}\left( \frac{\lambda_{\max} + \lambda_{\min} }{\lambda_{\max}-\lambda_{\min}} \right) },
$$
where $T_t(w) = \frac{1}{2} \left[ (w+\sqrt{w^2+1})^t + (w-\sqrt{w^2-1})^t\right]$ represents the \textit{Chebyshev polynomial} of degree $t$, and $\lambda_{\max}$ and $\lambda_{\min}$ denote the largest and smallest eigenvalues of $\bA$.
\begin{proof}
To see this, we can express the $m_t$ in \eqref{equation:cg-convergence-xt5} as
\begin{equation}\label{equation:cg-convergence-xt5_rewrite}
	m_t = \mathop{\min}_{P_{t-1}}\mathop{\max}_{1\leq j\leq n}(1+ \lambda_j P_{t-1}(\lambda_j))^2 = \mathop{\min}_{P_{t-1}}\mathop{\max}_{1\leq j\leq n} (\widetilde{P}_{t}(\lambda_i))^2,
\end{equation}
where $\widetilde{P}_{t}(\lambda) \triangleq 1+ \lambda P_{t-1}(\lambda)=1+w_1\lambda + \ldots+w_t\lambda^t$ is a special polynomial of degree  $t$ with $\widetilde{P}_{t}(0)=1$. We note that the Chebyshev polynomial can be expressed on the interval $w\in [-1,1]$ as
$$
T_t(w) =\cos(t \cos^{-1} w), \gap w\in [-1,1] 
\quad\implies\quad 
\abs{T_t(w)} \leq 1,\gap \text{if } w\in [-1,1].
$$
It is observable that $\widetilde{P}_{t}(\lambda)$ oscillates within the range $\pm {T_{t}\left( \frac{\lambda_{\max} + \lambda_{\min} }{\lambda_{\max}-\lambda_{\min}} \right) }^{-1}$ over the domain $[\lambda_{\min}, \lambda_{\max}]$. Suppose there exists a polynomial $S_t(\lambda)$ of degree $t$ such that $S_t(0)=1$ and $S_t$ is better than $\widetilde{P}_t$ on the domain $[\lambda_{\min}, \lambda_{\max}]$. It then follows that the $S_t-\widetilde{P}_t$ has a zero at $\lambda=0$ and other $t$ zeros on the range $[\lambda_{\min}, \lambda_{\max}]$, making it has $t+1$ zeros, which leads to a contradiction. Therefore, $\widetilde{P}_t$ is the optimal polynomial of degree $t$.
This completes the proof.
\end{proof}

\begin{figure}[h]
\centering  
\vspace{-0.35cm} 
\subfigtopskip=2pt 
\subfigbottomskip=2pt 
\subfigcapskip=-5pt 
\subfigure[Rate of convergence for GD with line search per iteration (same as Figure~\ref{fig:rate_convergen_steepest}). The $y$-axis is $\frac{\kappa-1}{\kappa+1}$.]{\label{fig:rate_convergen_steepest1}
	\includegraphics[width=0.481\linewidth]{./imgs/rate_convergen_steepest.pdf}}
\subfigure[Rate of convergence for CG per iteration. The $y$-axis is $\frac{\sqrt{\kappa}-1}{\sqrt{\kappa}+1}$.]{\label{fig:rate_convergen_conjugate}
	\includegraphics[width=0.481\linewidth]{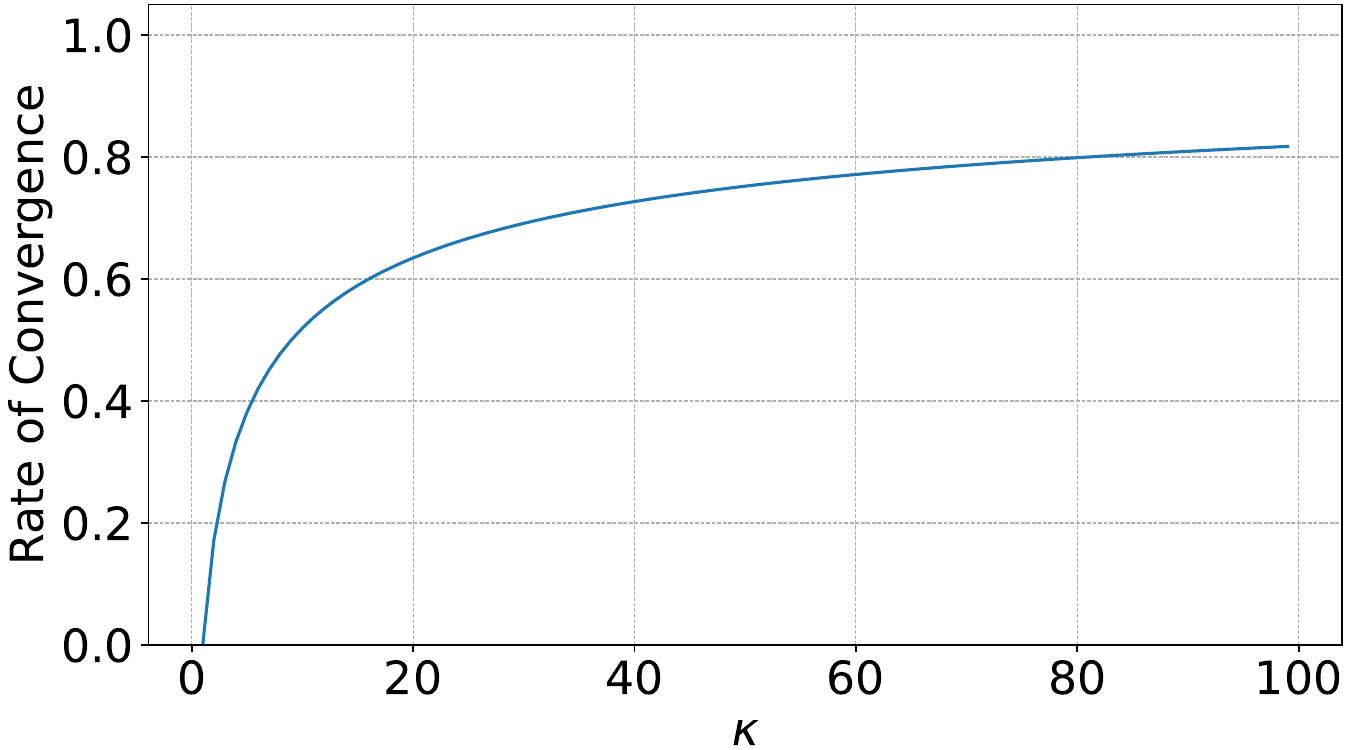}}
\caption{Illustration of the rate of convergence for CG and GD with line search.}
\label{fig:rate_convergen_conjugae_comparison}
\end{figure}

Therefore, it follows that
$$
\begin{aligned}
\normbig{\be^\toptone}_{\bA} &\leq T_t\left(  \frac{\lambda_{\max} + \lambda_{\min}}{\lambda_{\max} - \lambda_{\min}}   \right)^{-1} \cdot\normbig{\be^\topone}_{\bA} 
=T_t\left(  \frac{\kappa+1}{\kappa-1}   \right)^{-1} \cdot \normbig{\be^\topone}_{\bA}\\
&= 2\left[ \left(\frac{\sqrt{\kappa}+1}{\sqrt{\kappa}-1}  \right)^t +
\left(\frac{\sqrt{\kappa}-1}{\sqrt{\kappa}+1}  \right)^t
\right]^{-1} \cdot\normbig{\be^\topone}_{\bA},
\end{aligned}
$$
where $\kappa = \frac{\lambda_{\max}}{\lambda_{\min}}$ is the condition number, and $\left(\frac{\sqrt{\kappa}-1}{\sqrt{\kappa}+1}  \right)^t \rightarrow 0$ as iteration $t$ grows. A weaker inequality can be obtained by 
$$
\normbig{\be^\toptone}_{\bA} \leq  2 \left(\frac{\sqrt{\kappa}-1}{\sqrt{\kappa}+1}  \right)^t
\cdot\normbig{\be^\topone}_{\bA}.
$$
Figure~\ref{fig:rate_convergen_conjugae_comparison}  compares the rate of convergence of GD with line search and CG per iteration. It is observed that  CG exhibits significantly faster convergence compared to the GD with line search method.

\index{Rate of convergence}
\index{Preconditioning}
\index{Change of variables}
\subsection*{Preconditioning}
Since the smaller the condition number $\kappa$, the faster the convergence (Figure~\ref{fig:rate_convergen_conjugate}). We can accelerate the convergence of CG by transforming the linear system to improve the eigenvalue distribution of $\bA${\textemdash}the procedure is known as \textit{preconditioning}, which is closely related to the change of variables in the non-Euclidean gradient descent methods (see \eqref{equation:gree_cond_opt}, and Sections~\ref{section:als-gradie-descent-taylor} and \ref{section:noneucli_gd}). The variable $\bx$ is transformed to $\widehat{\bx}$ via a nonsingular matrix $\bP$, satisfying
$$
\begin{aligned}
\whbx &\triangleq \bP\bx
\qquad\implies \qquad 
\whf(\whbx) =\frac{1}{2}\whbx^\top (\bP^{-\top} \bA\bP^{-1})\whbx - (\bP^{-\top}\bb)^\top \whbx +c.
\end{aligned}
$$
When $\bA$ is symmetric, the solution of $\whf(\whbx)$ is equivalent to the solution of the linear equation
$$
\begin{aligned}
(\bP^{-\top} \bA\bP^{-1})\whbx &= \bP^{-\top}\bb 
\qquad\implies\qquad \bP^{-\top}\bA\bx=\bP^{-\top}\bb 
\qquad\implies\qquad \bA\bx=\bb.
\end{aligned}
$$
That is, we can solve $\bA\bx=\bb$ indirectly by solving $\bP^{-\top}\bA\bx=\bP^{-\top}\bb$.
Therefore, the rate of convergence of the quadratic form $\whf(\whbx)$ depends on the condition number of $\bP^{-\top} \bA\bP^{-1}$, which can be controlled by the nonsingular matrix $\bP$. 
Intuitively,  preconditioning is a procedure to stretch the quadratic form to make it more spherical so that the eigenvalues are clustered in a smaller range. A specific example is given in Figure~\ref{fig:conjugate_specialcases} that we want to transform the elliptical contour in Figure~\ref{fig:conjugate_specialcases_2eigenvalue} into the spherical contour in Figure~\ref{fig:conjugate_specialcases_1eigenvalue}.
Based on Algorithm~\ref{alg:practical_conjugate_descent}, the preconditioned CG method is formulated in Algorithm~\ref{alg:predicition_CG}.

\begin{algorithm}[h] 
\caption{Transformed-Preconditioned CG for Quadratic Functions}
\label{alg:predicition_CG}
\begin{algorithmic}[1]
\Require  Symmetric positive definite $\bA\in \real^{n\times n}$;
\State {\bfseries Input:} Initial parameter $\whbx^{(1)}$;
\State {\bfseries Input:} Initialize $\whbd^{(0)} =\bzero $ and $\whbg^{(0)} = \whbd^{(0)}+\bepsilon$;
\For{$t=1:n$ } 
\State Compute gradient $\whbg_t \leftarrow \nabla \whf(\whbx^\toptzero) = (\bP^{-\top} \bA\bP^{-1})\whbx- \bP^{-\top}\bb $; \Comment{$=\textcolor{mylightbluetext}{\bP^{-\top}}\bg^\toptzero$}
\State Compute coefficient $\widehat{\beta}_{t} \leftarrow  \frac{\whbg_t^\top\whbg_t}{\whbg_{t-1}^\top \whbg_{t-1}}
$; \Comment{$=
\frac{\bg^\toptzeroTOP \textcolor{mylightbluetext}{(\bP^\top\bP)^{-1}}\bg^\toptzero}
{\bg^\toptminusTOP\textcolor{mylightbluetext}{(\bP^\top\bP)^{-1}} \bg^\toptminus}$}
\State Compute descent direction $\whbd_t \leftarrow -\whbg_{t} +\widehat{\beta}_{t}  \whbd_{t-1}$;
\Comment{$=-\textcolor{mylightbluetext}{\bP^{-\top}}\bg^\toptzero +\widehat{\beta}_{t}  \whbd_{t-1}$}
\State stepsize 
$\widehat{\eta}_t \leftarrow \frac{{\whbg_t}^\top \whbg_t}{ \whbd_t^\top (\bP^{-\top} \bA\bP^{-1})\whbd_t }
$;
\Comment{$=
- \frac{{\bg^\toptzero}^\top \textcolor{mylightbluetext}{(\bP^\top\bP)^{-1}}\bg^\toptzero}{ \whbd_t^\top (\bP^{-\top} \bA\bP^{-1})\whbd_t }$}
\State Apply update $\whbx_{t+1} \leftarrow \whbx^\toptzero + \widehat{\eta}_t\whbd_t$;
\EndFor
\State {\bfseries Return:} resulting parameters $\textcolor{mylightbluetext}{\bx^\toptzero=\bP^{-1}\whbx^\toptzero}$, and the loss $f(\bx^\toptzero)$.
\end{algorithmic}
\end{algorithm}

However, the procedure in Algorithm~\ref{alg:predicition_CG} is not desirable since we need to transform $\bx$ into $\whbx=\bP\bx$ and transform back by $\bx=\bP^{-1}\whbx$ as highlighted in the blue texts of Algorithm~\ref{alg:predicition_CG}. This introduces additional computational overhead. Let $\bM=\bP^\top\bP$, Algorithm~\ref{alg:untransformed_predicition_CG} is proposed to formulate the untransformed-preconditioned CG, which proves to be more efficient  than Algorithm~\ref{alg:predicition_CG}.

\begin{algorithm}[h] 
\caption{Untransformed-Preconditioned CG for Quadratic Functions}
\label{alg:untransformed_predicition_CG}
\begin{algorithmic}[1]
\Require  Symmetric positive definite $\bA\in \real^{n\times n}$;
\State {\bfseries Input:} Initial parameter $\bx^\topone$;
\State {\bfseries Input:} Initialize $\bd^\topzero =\bzero $ and $\bg^\topzero = -\bd^\topzero+\bepsilon$;
\For{$t=1:n$ } 
\State Compute gradient $\bg^\toptzero \leftarrow \nabla f(\bx^\toptzero)$;
\Comment{Same as that of Algorithm~\ref{alg:practical_conjugate_descent}}
\State Compute coefficient $\widehat{\beta}_{t} \leftarrow  \frac{\bg^\toptzeroTOP\textcolor{mylightbluetext}{\bM^{-1}}\bg^\toptzero}{\bg^\toptminusTOP\textcolor{mylightbluetext}{\bM^{-1}} \bg^\toptminus}$; \Comment{Same as that of Algorithm~\ref{alg:predicition_CG}}
\State Compute descent direction $\widetilde{\bd}_t = -\textcolor{mylightbluetext}{\bM^{-1}}\bg^\toptzero +\widehat{\beta}_{t}  \widetilde{\bd}_{t-1}$;
\Comment{$=-\textcolor{mylightbluetext}{\bP^{-1}} \whbd_{t}$ in Algorithm~\ref{alg:predicition_CG}}
\State Stepsize $\widehat{\eta}_t \leftarrow {({\bg^\toptzero}^\top\textcolor{mylightbluetext}{\bM^{-1}} \bg^\toptzero)}/{ (\widetilde{\bd}_t^\top \bA\widetilde{\bd}_t )}$; \Comment{Same as that of Algorithm~\ref{alg:predicition_CG}}
\State Apply update $\bx^\toptone \leftarrow \bx^\toptzero + \widehat{\eta}_t\widetilde{\bd}_t$; \Comment{$=-\textcolor{mylightbluetext}{\bP^{-1}} \Delta\whbx^\toptzero$ in Algorithm~\ref{alg:predicition_CG}}
\EndFor
\State {\bfseries Return:} resulting parameters $\bx^\toptzero$, and the loss $f(\bx^\toptzero)$.
\end{algorithmic}
\end{algorithm}

\begin{figure}[h]
\centering  
\vspace{-0.35cm} 
\subfigtopskip=2pt 
\subfigbottomskip=2pt 
\subfigcapskip=-5pt 
\subfigure[Contour plot of quadratic function with $\bA$.]{\label{fig:conjugate_precondition1}
	\includegraphics[width=0.481\linewidth]{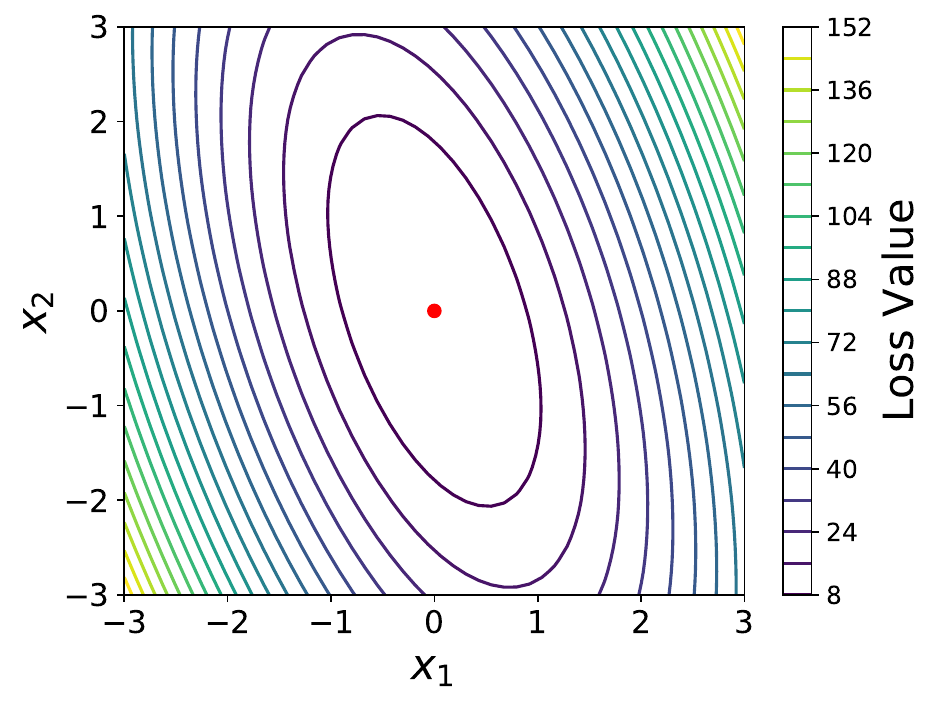}}
\subfigure[Contour plot of quadratic function with $\bP^{-\top} \bA\bP^{-1}$.]{\label{fig:conjugate_precondition2}
	\includegraphics[width=0.481\linewidth]{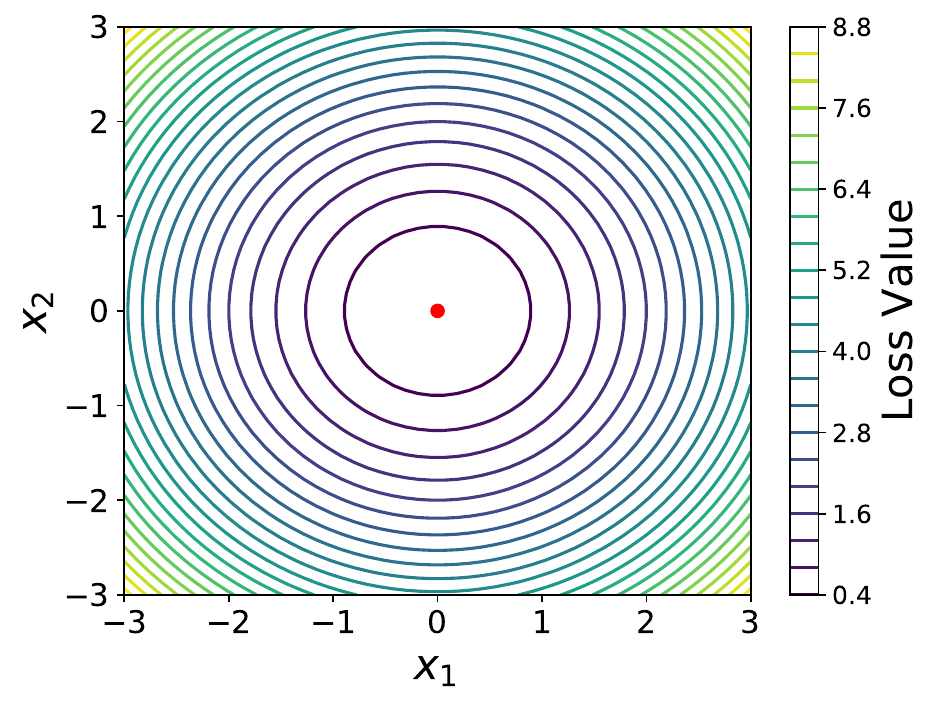}}
\caption{Illustration of preconditioning for $\bA=\footnotesize\begin{bmatrix}
		20&5 \\5&5
	\end{bmatrix}$. $\bP$ is obtained by the Cholesky decomposition such that $\bM=\bA=\bP^\top\bP$.}
\label{fig:conjugate_precondition13}
\end{figure}
\index{Cholesky decomposition}

\paragrapharrow{Second perspective of preconditioning.}
The matrices $\bM^{-1}\bA$ and $\bP^{-\top} \bA\bP^{-1}$ have the same eigenvalues. To see this, suppose the eigenpair of $\bM^{-1}\bA$ is $(\bM^{-1} \bA) \bv =\lambda \bv$, it follows that
$$
(\bP^{-\top} \bA\bP^{-1}) (\bP\bv) = \bP^{-\top} \bA\bv = 
\bP\bP^{-1}\bP^{-\top} \bA\bv 
=\bP\bM^{-1}\bA\bv=\lambda (\bP\bv).
$$
Therefore, the preconditioning can be understood from two perspectives. While the second perspective is to solve $\bM^{-1}\bA\bx = \bM^{-1}\bb$, where the condition number is decided by matrix $\bM^{-1}\bA$.
The simplest preconditioner $\bM^{-1}$ is  a diagonal matrix whose diagonal entries are identical to those of $\bA$, known as \textit{diagonal preconditioning}, in which case  we scale the quadratic form along the coordinate axes. In contrast, the \textit{perfect preconditioner} is $\bM=\bA$ such that $\bM^{-1}\bA=\bI$, whose condition number is 1, in which case  the quadratic form is scaled along its eigenvector directions. In this sense, the $\bP$ can be obtained by the (pseudo) Cholesky decomposition (Theorems~\ref{theorem:cholesky-factor-exist}) such that $\bM=\bA=\bP^\top\bP$. Figure~\ref{fig:conjugate_precondition13} shows the perfect preconditioning on $\bM=\bA=
\footnotesize
\begin{bmatrix}
20&5 \\5&5\\
\end{bmatrix}$ such that the eigenvalues of $\bP^{-\top}\bA\bP^{-1}$ are identical and the condition number is thus equal to 1.

\index{Cholesky decomposition}

\subsection{General Conjugate Gradient Method and Convergence Analysis}
We now revisit the general CG method as introduced in \citet{fletcher1964function}. The method has been previously formulated in Algorithm~\ref{alg:conjugate_descent} at the beginning of this section, where the search direction at $t$-th iteration is 
\begin{equation}\label{equation:search_frcg}
\bd^\toptzero =  -\bg^\toptzero + \beta_t^F \bd^\toptminus
\qquad\text{with}\qquad
\beta_t^F = \frac{\nabla f(\bx^\toptzero)^\top \nabla f(\bx^\toptzero)}{\nabla f(\bx^\toptminus)^\top \nabla f(\bx^\toptminus)}
=\frac{ \bg^\toptzeroTOP \bg^\toptzero}{\bg^\toptminusTOP \bg^\toptminus},
\end{equation}
We may notice the \textit{Fletcher-Reeves Conjugate Gradient} method (Algorithm~\ref{alg:conjugate_descent}) is just the same as the \textit{Practical Conjugate Gradient} method (Algorithm~\ref{alg:practical_conjugate_descent}) under the conditions of a strongly convex quadratic loss function and the use of an exact line search for the stepsize $\eta_t$.

To see why the Fletcher-Reeves Conjugate Gradient (FRCG) algorithm (Algorithm~\ref{alg:conjugate_descent}) works, the search direction $\bd^\toptzero$ must satisfy the descent condition  such that $ \bg^\toptzeroTOP \bd^\toptzero<0$ (Definition~\ref{definition:uncons_des_direct}). The descent condition is satisfied when the stepsize is calculated by exact line search, in which case the gradient $\nabla f(\bx^\toptzero) = \bg^\toptzero$ is orthogonal to search direction $\bd^\toptminus$ (Lemma~\ref{lemm:linear-search-orghonal}): $\bg^\toptzeroTOP\bd^\toptminus=0$. Therefore, 
$$
\bg^\toptzeroTOP \bd^\toptzero = \bg^\toptzeroTOP (-\bg^\toptzero +\beta_t^F\bd^\toptminus ) = -\normtwobig{\bg^\toptzero}^2 + \beta_t^F \bg^\toptzeroTOP \bd^\toptminus<0
$$
when $\eta_t$ is determined by exact line search. These properties result in the global convergence of the Fletcher-Reeves CG method (see Theorem~\ref{theorem:glob_fr_cg}). However, when $\eta_t$ is fixed or calculated by inexact line search, the descent condition $\bg^\toptzeroTOP\bd^\toptzero$ may not be satisfied. This problem, however, can be attacked by \textit{strong Wolfe conditions} \citep{nocedal1999numerical, fletcher1964function}; and we will not go into the details.

\paragrapharrow{Polak-Ribi\`ere conjugate gradient.} We have mentioned previously that the $\beta_t$ can also be computed by the Polak-Ribi\`ere coefficient:
$$
\text{Polak-Ribi\`ere:\gap } \beta_t^P = \frac{\Big(\nabla f(\bx^\toptzero) - \nabla f(\bx^\toptminus)\Big)^\top \nabla f(\bx^\toptzero)}{\nabla f(\bx^\toptminus)^\top \nabla f(\bx^\toptminus)}
=
\frac{(\bg^\toptzero - \bg^\toptminus)^\top  \bg^\toptzero}{ \bg^\toptminusTOP \bg^\toptminus}
.
$$
When the loss function is strongly convex quadratic and the stepsize is chosen by exact line search, the Polak-Ribi\`ere coefficient $\beta_t^P$ is identical to the Fletcher-Reeves coefficient $\beta_t^F$ since $\bg^\toptminusTOP\bg^\toptzero=0$ by Theorem~\ref{theorem:conjudate_CD_d-steps}.

\paragrapharrow{Hestenes-Stiefel conjugate gradient.} Hestenes-Stiefel coefficient is yet another variant of the Polak-Ribi\`ere coefficient:
$$
\text{Hestenes-Stiefel:\gap } \beta_t^H = \frac{\Big(\nabla f(\bx^\toptzero) - \nabla f(\bx^\toptminus)\Big)^\top \nabla f(\bx^\toptzero)}{\Big(\nabla f(\bx^\toptzero) - \nabla f(\bx^\toptminus)\Big)^\top
\bd^\toptminus}
=
\frac{(\bg^\toptzero - \bg^\toptminus)^\top  \bg^\toptzero}{ (\bg^\toptzero - \bg^\toptminus)^\top \bd^\toptminus}.
$$
When the loss function is strongly convex quadratic and the stepsize is chosen by exact line search, the Hestenes-Stiefel coefficient $\beta_t^H$ is identical to the Fletcher-Reeves coefficient $\beta_t^F$ since $\bg^\toptminusTOP\bg^\toptzero=0$ by Theorem~\ref{theorem:conjudate_CD_d-steps} and $\bg^\toptzeroTOP\bd^{(t-2)}=\bg^\toptminusTOP\bd^{(t-2)}=0$ by Theorem~\ref{theorem:expanding_subspace_minimization} \citep{hestenes1952methods}.

Moreover, numerical experiments show that the Polak-Ribi\`ere coefficient and Hestenes
-Stiefel coefficient are more robust than Fletcher-Reeves coefficient in non-convex settings \citep{nocedal1999numerical}.

For simplicity, we only prove the global convergence of the Fletcher-Reeves CG method with exact line search. The convergence results for alternative algorithms can be found, for example, in \citet{nocedal1999numerical, dai1999nonlinear, sun2006optimization}.
\begin{theoremHigh}[Global Convergence of FRCG with Exact Line Search]\label{theorem:glob_fr_cg}
Let  $f: \real^n \to \real$ be a continuously differentiable function on a bounded level set $\sL \triangleq\lev[f, \bx^\topone] = \{\bx \in \real^n \mid f(\bx) \leq f(\bx^\topone)\}$, and let the Fletcher-Reeves CG method (Algorithm~\ref{alg:conjugate_descent}) be implemented with exact line search. Then the produced sequence $\{\bx^\toptzero\}_{t>0}$ has at least one accumulation point which is a stationary point, i.e.,
\begin{enumerate}[(i)]
\item When $\{\bx^\toptzero\}$ is a finite sequence, then the final point $\bx^*$ is a stationary point of $f$;
\item When $\{\bx^\toptzero\}$ is an infinite sequence, it has limit point, and any limit point is a stationary point.
\end{enumerate}
\end{theoremHigh}
\begin{proof}[of Theorem~\ref{theorem:glob_fr_cg}]
\textbf{(i).} When $\{\bx^\toptzero\}$ is finite, from the termination condition, it follows that the final point $\bx^*$ satisfies $\nabla f(\bx^*) = \bzero $, and hence $\bx^*$ is a stationary point of $f$.

\paragraph{(ii).} When $\{\bx^\toptzero\}$ is infinite, we have $\nabla f(\bx^\toptzero) \neq \bzero, \forall t$. Noting that $\bd^\toptzero = -\bg^\toptzero + \beta_{t} \bd^\toptminus$ and $\bg^\toptzeroTOP \bd^\toptminus = 0$ by exact line search (Lemma~\ref{lemm:linear-search-orghonal}), we have
\begin{equation}
\bg^\toptzeroTOP \bd^\toptzero = -\normtwobig{\bg^\toptzero}^2 + \beta_{t} \bg^\toptzeroTOP \bd^\toptminus = -\normtwobig{\bg^\toptzero}^2 < 0,
\end{equation}
which means that $\bd^\toptzero$ is a descent direction, $\{f(\bx^\toptzero)\}$ is a monotone descent sequence, and thus $\{\bx^\toptzero\} \subset \sL$. Therefore $\{\bx^\toptzero\}$ is a bounded sequence and must have a limit point.

Let $\bx^*$ be a limit point of $\{\bx^\toptzero\}$. Then there is a subsequence $\{\bx^\toptzero\}_{\sT_1}$ converging to $\bx^*$, where $\sT_1$ is an index set of a subsequence of $\{\bx^\toptzero\}$. Since $\{\bx^\toptzero\}_{\sT_1} \subset \{\bx^\toptzero\}$, $\{f(\bx^\toptzero)\}_{\sT_1} \subset \{f(\bx^\toptzero)\}$. It follows from the continuity of $f$ that for $t \in \sT_1$,
\begin{equation}
f(\bx^*) = f(\lim_{t \to \infty} \bx^\toptzero) = \lim_{t \to \infty} f(\bx^\toptzero).
\end{equation}

Similarly, $\{\bx^\toptone\}$ is also a bounded sequence. Hence there exists a subsequence $\{\bx^\toptone\}_{\sT_2}$ converging to $\widetildebx^*$, where $\sT_2$ is an index set of a subsequence of $\{\bx^\toptone\}$. In this case,
\begin{equation}
f(\widetildebx^*) = f(\lim_{t \to \infty} \bx^\toptone) = \lim_{t \to \infty} f(\bx^\toptone).
\end{equation}
This indicates $
f(\widetildebx^*) = f(\bx^*)
$.

Now we prove $\nabla f(\bx^*) = \bzero$ by contradiction. Suppose that $\nabla f(\bx^*) \neq \bzero$, then, for $\eta$ sufficiently small, there exists a search direction $\bd^*$ such that 
\begin{equation}\label{eeuation:conv_frcd_lim}
f(\bx^* + \eta \bd^*) < f(\bx^*).
\end{equation}
Since
$ f(\bx^\toptone) = f(\bx^\toptzero + \eta_t \bd^\toptzero) \leq f(\bx^\toptzero + \eta \bd^\toptzero), \; \forall \eta > 0, $
then for $t \in \sT_2$, passing to the limit $t \to \infty$ and using \eqref{eeuation:conv_frcd_lim}, we get
\begin{equation}
f(\widetildebx^*) \leq f(\bx^* + \eta \bd^*) < f(\bx^*),
\end{equation}
which contradicts the fact $
f(\widetildebx^*) = f(\bx^*)
$. 
This proves $\nabla f(\bx^*) = \bzero$, i.e., $\bx^*$ is a stationary point of $f$. 
\end{proof}

\begin{problemset}

\item Prove the relation in \eqref{equation:bfgs_det}.

\item \label{prob:indep_conjugate} Prove that if $\bd_1, \bd_2, \ldots, \bd_m$ are $\bA$-conjugate, where $\bA\in\real^{n\times n}$ is positive definite, then they are linearly independent.

\item Use Newton's method, trust region method, and the Fletcher-Reeves Conjugate Gradient (FRCG) algorithm to minimize the following function:
$$
f(\bx) = x_1^2 + 8x_2^2 +3x_1x_2 +6x_1
$$
with the initial point $\bx^{(1)} = [0,0]^\top$.

\item Verify that the BFGS updates in \eqref{equation:bfgs_update} and \eqref{equation:bfgs2_update} are inverses to each other.

\item Prove the equalities in \eqref{equation:bfgs_trace} and \eqref{equation:bfgs_det}.

\item Show that when applied to a quadratic function, using CG with exact line search strategy, both the Polak-Ribi\`ere and Hestenes-Stiefel updates reduce to the Fletcher-Reeves update in \eqref{equation:cg_all_forms}.
\end{problemset}

\chapter{Least Squares: Linear, Nonlinear, Sparse Problems}
\begingroup
\hypersetup{
	linkcolor=structurecolor,
	linktoc=page,  
}
\minitoc \newpage
\endgroup
Least squares problems are fundamental to optimization, especially in regression analysis and data fitting. This chapter delves into linear and nonlinear least squares, Levenberg-Marquardt methods, and sparse optimization techniques such as L1-regularized optimization and iterative hard-thresholding methods. These approaches are particularly pertinent for compressed sensing and high-dimensional data analysis.

\index{Least squares}
\index{Linear model}
\section{Linear Model and Least Squares}\label{section:pre_ls}

The linear model is a cornerstone technique in regression analysis, employing the least squares approximation as its primary tool to minimize the sum of squared errors. This method is naturally suited for finding the regression function that minimizes the expected squared error. Over recent decades, linear models have been extensively applied across various fields, including decision-making  \citep{dawes1974linear}, time series analysis \citep{christensen1991linear, lu2017machine}, quantitative finance \citep{menchero2011barra}, and numerous other disciplines like production science, social science, and soil science \citep{fox1997applied, lane2002generalized, schaeffer2004application, mrode2014linear}.

To illustrate, consider an overdetermined system represented by $\bb = \bA\bx $, where $\bA\in \real^{m\times n}$ denotes  the \textit{input data matrix} (also referred to the \textit{predictor variables}), $\bb\in \real^m$ represents  the \textit{observation vector} (or \textit{target/response vector}), and the  number of samples $m$ exceeds  the dimensionality $n$. 
The vector $\bx$ contains the weights of the linear model. 
Typically, $\bA$ will have full column rank since real-world data is often uncorrelated or can be processed to achieve this.
In practical applications, a \textit{bias term} (a.k.a., an \textit{intercept}) is added to the first column of $\bA$ to enable the least squares method to solve:
\begin{equation}\label{equation:ls-bias}
	\widetildebA \widetildebx = 
	[\bm{1} ,\bA ] 
	\begin{bmatrix}
		x_0\\
		\bx
	\end{bmatrix}
	= \bb .
\end{equation}

However, it is common for the equation $\bb = \bA\bx$ not to have an exact solution (the system is \textit{inconsistent}) due to the  being overdetermined---that is, there are more equations than unknowns.
Define the column space of $\bA$ as $\{\bA\bgamma: \,\, \forall \bgamma \in \real^n\}$, denoted by $\cspace(\bA)$.
In essence, when we say $\bb = \bA\bx$ has no solution, it implies that $\bb$ lies outside the column space of $\bA$. 
In other words, the error $\be = \bb -\bA\bx$ cannot be reduced to zero. 
In these cases, the objective shifts to minimizing the error, typically measured by the mean squared error.
The resulting solution $\bx_{LS}$, which minimizes $\normtwo{\bb-\bA\bx_{LS}}^2$, is called the \textit{least squares solution}. The least squares method is a cornerstone of mathematical sciences, and there are numerous resources dedicated entirely to this topic, including works by   \citet{trefethen1997numerical, strang2019linear, strang2021every,  lu2021rigorous}.

\paragrapharrow{Least squares by calculus.}
When $\normtwo{\bb-\bA\bx}^2$ is differentiable and the parameter space of $\bx$ spans the entire space $\real^n$ (i.e., unconstrained optimization), 
the least squares estimate corresponds to the root of the gradient of $\normtwo{\bb-\bA\bx}^2$. 
This leads us to the following lemma.~\footnote{Variants of the least squares problem are explored in Problems~\ref{problem:rls}$\sim$\ref{problem:twls2}.}

\index{Optimality condition}
\begin{lemma}[Least Squares by Calculus]\label{lemma:ols}
Let  $\bA \in \real^{m\times n}$ be a  fixed data matrix with full rank  and $m\geq n$ (i.e., its columns  are linearly independent)~\footnote{See Problem~\ref{prob:als_pseudo1}$\sim$\ref{prob:als_pseudon} for a relaxation of this condition using the pseudo-inverse.}. 
For the overdetermined system $\bb = \bA\bx$, the least squares solution is obtained by  setting the partial derivatives of $\normtwo{\bb-\bA\bx}^2$ in every direction to zero (i.e., the gradient vanishes by Theorem~\ref{theorem:fetmat_opt}, since $f(\bx)=\normtwo{\bb-\bA\bx}^2$ is a convex function by Exercise~\ref{exercise:conv_quad}). This yields the solution: 
$$
\bx_{LS} = (\bA^\top\bA)^{-1}\bA^\top\bb.
$$ 
The value, $\bx_{LS} = (\bA^\top\bA)^{-1}\bA^\top\bb$, is commonly referred to as the \textit{ordinary least squares (OLS)} estimate or simply the \textit{least squares (LS)} estimate of $\bx$.
\end{lemma}
In the lemma, we use the fact $\bA^\top\bA \in \real^{n\times n}$ is invertible since we assume that $\bA$ has full rank and $m\geq n$, . 

\begin{definition}[Normal Equation]\label{definition:normal-equation-als}
The condition for the gradient of $\normtwo{\bb-\bA\bx}^2$ to be zero can be expressed as $\bA^\top\bA \bx_{LS} = \bA^\top\bb$. The equation is also known as the \textit{normal equation}. 
Under the assumption that $\bA$ has full rank with $m\geq n$, it follows that $\bA^\top\bA$ is invertible, implying $\bx_{LS} = (\bA^\top\bA)^{-1}\bA^\top\bb$.
\end{definition}

Another question that arises: Why does the normal equation seemingly ``magically" yield solutions for  $\bx$?
A simple example can help illustrate the answer. The equation $x^2=-1$ has no real solution. 
However, $x\cdot x^2 = x\cdot (-1)$ does have a real solution $\hat{x} = 0$, in which case, $\hat{x}$ minimizes the difference between $x^2$ and $-1$, making them as close as possible.

\begin{example}[Altering the Solution Set by Left Multiplication]
Consider the  data matrix and target vector:
$\tiny
\bA=
\begin{bmatrix}
-3 & -4 \\
4 & 6  \\
1 & 1
\end{bmatrix}
$
and
$
\bb=
\tiny
\begin{bmatrix}
1  \\
-1   \\
0
\end{bmatrix}
.
$
It can be easily verified that the system $\bA\bx = \bb$ has no solution for $\bx$. 
However, if we multiply both sides on the left by
$
\bB=
\scriptsize
\begin{bmatrix}
0 & -1 & 6\\
0 & 1  & -4
\end{bmatrix},
$
then  $\bx_{LS} = [1/2, -1/2]^\top$ becomes the solution to $\bB\bA\bx= \bB\bb$. 
This example illustrates why the normal equation can lead to the least squares solution. Multiplying a linear system on the left alters the solution set, effectively projecting the problem into a different space where the least squares solution can be computed.
\end{example}

\paragrapharrow{Rank-deficiency.}
In this discussion, we assume that $\bA\in \real^{m\times n}$ has full rank with $m\geq n$, ensuring that $\bA^\top\bA$ is invertible. 
However, if two or more columns of $\bA$ are perfectly correlated, the matrix $\bA$ becomes deficient, and $\bA^\top\bA$ becomes singular. 
To address this issue, one approach is to choose $\bx$ that minimizes $\normtwo{\bx_{LS}}^2$ while satisfying the normal equation. 
That is, we select the least squares solution with the smallest magnitude.
\citet{lu2021numerical} discussed how to use UTV decomposition and singular value decomposition (SVD) to address this rank-deficient least squares problems.
See Problems~\ref{prob:als_pseudo1}$\sim$\ref{prob:als_pseudon} or the following paragraph for more insights.

\index{Contion number}
\index{Tikhonov regularization}
\index{$\ell_2$ regularization}
\paragrapharrow{Regularizations and stability.}
A common  problem that arise in the ordinary least square solution is the near-singularity of $\bA$.
Let the SVD of $\bA$ be $\bA=\bU\bSigma\bV^\top\in\real^{m\times n}$, where $\bU\in\real^{m\times m}$ and $\bV\in\real^{n\times n}$ are orthogonal, and the main diagonal of $\bSigma\in\real^{m\times n}$ contains the singular values (Theorem~\ref{theorem:full_svd_rectangular}). Consequently, $\bA^\top\bA = \bV(\bSigma^\top\bSigma)\bV^\top \triangleq \bV\bS\bV^\top$, where $\bS\triangleq \bSigma^\top\bSigma  = \diag(\sigma_1^2, \sigma_2^2, \ldots,\sigma_n^2)\in\real^{n\times n}$ contains the squared singular values of $\bA$. When $\bA$ is nearly singular, $\sigma_n^2\approx 0$, making the inverse operation $(\bA^\top\bA)^{-1} = \bV\bS^{-1}\bV^\top$ numerically unstable. 
As a result, the solution $\bx_{LS} =(\bA^\top\bA)^{-1}\bA^\top\bb $ may diverge.
To address this issue, we typically  add an $\ell_2$ regularization term to obtain the solution for the following optimization problem:
\begin{equation}
	\bx_{Tik} = \mathop{\argmin}_{\bx} \normtwo{\bb-\bA\bx}^2 +\lambda\normtwo{\bx}^2.
\end{equation}
This method is known as the  \textit{Tikhonov regularization method} (or simply the $\ell_2$ regularized method) \citep{tikhonov1963solution}.
The gradient of the problem is $2(\bA^\top\bA+\lambda\bI)\bx-2\bA^\top\bb$. Thus, the least squares solution is given by 
$
\bx_{Tik} = (\bA^\top\bA+\lambda\bI)^{-1}\bA^\top\bb.
$
The inverse operation becomes $(\bA^\top\bA+\lambda\bI)^{-1} = \bV(\bS+\lambda\bI)^{-1}\bV^\top$, where $\widetildebS\triangleq(\bS+\lambda\bI)=\diag(\sigma_1^2+\lambda, \sigma_2^2+\lambda, \ldots,\sigma_n^2+\lambda)$. 
The solutions for OLS and Tikhonov regularized LS are given, respectively, by 
\begin{equation}
	\begin{aligned}
		\bx_{LS} &= (\bA^\top\bA)^{-1}\bA^\top\bb = \bV\left(\bS^{-1}\bSigma\right)\bU^\top\bb;\\
		\bx_{Tik} &= (\bA^\top\bA+\lambda\bI)^{-1}\bA^\top\bb = \bV\left((\bS+\lambda\bI)^{-1}\bSigma\right)\bU^\top\bb,\\
	\end{aligned}
\end{equation}
where the main diagonals of $\left(\bS^{-1}\bSigma\right)$ are $\diag(\frac{1}{\sigma_1}, \frac{1}{\sigma_2}, \ldots, \frac{1}{\sigma_n})$; and the main diagonals of $\left((\bS+\lambda\bI)^{-1}\bSigma\right)$ are $\diag(\frac{\sigma_1}{\sigma_1^2+\lambda}, \frac{\sigma_2}{\sigma_2^2+\lambda}, \ldots, \frac{\sigma_n}{\sigma_n^2+\lambda})$. The latter solution is more stable if $\lambda$ is greater than the   smallest nonzero squared singular value.
The condition number becomes smaller if  the smallest singular value $\sigma_n$ is close to zero (see \eqref{equation:cond_pd_ineq} or \citet{lu2021numerical} for more details):
$$
\kappa(\bA^\top\bA) = \frac{\sigma_1^2}{\sigma_n^2}
\qquad \rightarrow \qquad
\kappa(\bA^\top\bA+\lambda\bI) = \frac{\lambda+\sigma_1^2}{\lambda+\sigma_n^2}.
$$
Tikhonov regularization effectively prevents  divergence  in the least squares solution 
$\bx_{LS} = (\bA^\top\bA)^{-1} \bA^\top \bb$ when the matrix $\bA$ is nearly singular or even rank-deficient. This improvement enhances the convergence properties of both the LS algorithm and its variants, such as alternating least squares,  while addressing identifiability issues  in various settings (see, for example, \citet{gillis2020nonnegative, lu2021numerical}). As a result, Tikhonov regularization has become a widely applied technique.

\index{Data least squares}
\index{Total least squares}
\paragrapharrow{Data least squares.}
While the OLS method accounts for errors in the response variable $\bb$, the \textit{data least sqaures (DLS)} method considers errors in the predictor variables:
\begin{equation}
	\bx_{DLS} = \mathop{\argmin}_{\bx, \widetildebA} \normfbig{\widetildebA}^2 \quad \text{s.t.}\quad \bb\in\cspace(\bA+\widetildebA),
\end{equation}
where $\widetildebA$ represents a perturbation of $\bA$ (i.e., a noise in the predictor variables).
That is, $(\bA+\widetildebA) \bx_{DLS} = \bb$, assuming the measured response $\bb$ is noise-free.
The Lagrangian function and its gradient w.r.t. $\bx$ are, respectively, given by
$$
\begin{aligned}
	L(\bx, \widetildebA, \blambda) &= \trace(\widetildebA\widetildebA^\top) +\blambda^\top (\bA\bx+\widetildebA\bx-\bb);\\
	\nabla_{\widetildebA} L(\bx, \widetildebA,\blambda) &= \widetildebA+\blambda\bx^\top = \bzero \quad\implies\quad \widetildebA=-\blambda\bx^\top,
\end{aligned}
$$
where $\blambda\in\real^m$  is a vector of Lagrange multipliers.
Substituting the value of the vanishing gradient into $(\bA+\widetildebA) \bx = \bb$ yields $\blambda = \frac{\bA\bx-\bb}{\bx^\top\bx}$ and $\widetildebA=-\frac{(\bA\bx-\bb)\bx^\top}{\bx^\top\bx} $.
Therefore, using the invariance of cyclic permutation of factors in trace,  the objective function becomes 
$$
\mathop{\argmin}_{\bx}
\frac{(\bA\bx-\bb)^\top (\bA\bx-\bb)}{\bx^\top\bx} .
$$

\paragrapharrow{Total least squares.} Similar to  data least squares, the \textit{total least squares (TLS)} method considers errors in both the predictor variables and the response variables. The TLS problem can be formulated as:
\begin{equation}
	\bx_{TLS} = \mathop{\argmin}_{\bx, \widetildebA, \widetildebb} \normf{[\widetildebA, \widetildebb]}^2, 
	\quad \text{s.t.}\quad (\bb+\widetildebb)\in\cspace(\bA+\widetildebA), 
\end{equation}
where $\widetilde{\bA}$ and $\widetilde{\bb}$ are perturbations in the predictor variables and the response variable, respectively.
Let $\bC\triangleq[\bA,\bb]\in\real^{m\times (n+1)}$,  $\bD\triangleq[\widetildebA, \widetildebb]\in\real^{m\times (n+1)}$, and $\by\triangleq\footnotesize\begin{bmatrix}
	\bx\\
	-1
\end{bmatrix}$, the problem can be equivalently stated as
\begin{equation}
	\bx_{TLS} = \mathop{\argmin}_{\by, \bD} \normf{\bD}^2, 
	\quad \text{s.t.}\quad \bD\by = -\bC\by, 
\end{equation}

\section{Nonlinear Least Squares Problem}\label{section:nonlinear_l2}
Building upon the linear least squares problem, we now consider the general form of the nonlinear least squares problem: 
\begin{equation}\label{equation:nonlinear_ols}
	f(\bx) = \frac{1}{2} \sum_{j=1}^{m} r_j^2(\bx)= \frac{1}{2} \normtwo{\br(\bx)}^2,
\end{equation}
where $r_j : \real^n \rightarrow \real$ are differentiable functions, and we assume $m \geq n$. The functions $r_j$ are referred to as \textit{residuals}.
 For convenience, we define the residual vector $\br : \real^n \rightarrow \real^m$ as
$ \br(\bx) \triangleq [r_1(\bx), r_2(\bx), \ldots, r_m(\bx)]^\top. $
Using this definition, the function $f(\bx)$ can be written as $f(\bx) = \frac{1}{2} \normtwo{\br(\bx)}^2$.
When $\br(\bx) = \bA\bx-\bb$, the problem becomes the linear least squares case discussed in the preceding section.

Now, consider the case where $m=n$. If we have a good initial estimate $\bx^{(1)}$, and the Jacobian $\bJ(\bx^{(t)})$ of the residual vector $\br(\bx)$ is nonsingular, then the \textit{Newton-Raphson} method provides a fast and accurate solution. 
A simple approach is then to combine the  Newton-Raphson method with line
search. The typical update rule at $t$-th iteration is
\begin{subequations}\label{equation:newton_raphson_nonlin}
\begin{align}
\bd^\toptzero &\leftarrow  \text{ solve }\bJ(\bx^\toptzero) \bd = -\br(\bx^\toptzero);\\
\bx^\toptone&\leftarrow \bx^\toptzero + \eta_t  \bd^\toptzero;
\end{align}
\end{subequations}
The search direction $\bd^\toptzero$ is a descent direction such that $\nabla f(\bx^\toptzero)^\top \bd^\toptzero \leq 0$.

More generally,
problem~\eqref{equation:nonlinear_ols} is an unconstrained optimization problem that can be solved using any of the algorithms discussed earlier. Here, we directly provide the gradient and Hessian matrix of $f(\bx)$:
\begin{subequations}
\begin{align}
\nabla f(\bx) &= \bJ(\bx)^\top \br(\bx),  \label{equation:gaus_new_grad} \\
\nabla^2 f(\bx) &= \bJ(\bx)^\top \bJ(\bx) + \sum_{i=1}^{m} r_i(\bx) \nabla^2 r_i(\bx), \label{equation:gaus_new_hess}
\end{align}
\end{subequations}
where $\bJ(\bx) \in \real^{m \times n}$ is the Jacobian matrix of the residual vector $\br(\bx)$ at point $\bx$. 
Note that $\nabla^2 f(\bx)$ consists of two components:  $\bJ(\bx)^\top \bJ(\bx)$ and $\sum_{i=1}^{m} r_i(\bx) \nabla^2 r_i(\bx)$. 
The computational difficulty of these two components differs: while calculating $\nabla f(\bx)$ naturally provides the Jacobian matrix of $\br(\bx)$ without extra effort, computing the second term of the Hessian matrix requires evaluating each $\nabla^2 r_i(\bx)$, which incurs additional computational cost. Many nonlinear least squares algorithms are designed to leverage this distinction.

\begin{algorithm}[h] 
\caption{Gauss-Newton Method}
\label{alg:gauss_newton}
\begin{algorithmic}[1] 
\Require A function $f(\bx) = \frac{1}{2} \sum_{j=1}^{m} r_j^2(\bx)=\frac{1}{2} \br(\bx)^\top\br(\bx)$; 
\State {\bfseries Input:}  Initialize $\bx^{(1)}\in\real^n$;
\For{$t=1,2,\ldots$}
\State Compute the residual vector $ \br^\toptzero $ and the Jacobian matrix $ \bJ^\toptzero $;
\State \algoalign{Compute the QR decomposition of $ \bJ^\toptzero $: $ \bJ^\toptzero = \bQ^\toptzero \bR^\toptzero $, where $ \bQ^\toptzero \in \real^{m \times n}, \bR^\toptzero \in \real^{n \times n} $;}
\State Solve the equation $ \bR^\toptzero \bd = -\bQ^\toptzeroTOP \br^\toptzero $ to obtain the descent direction $ \bd^\toptzero $.
\State Use line search criteria to compute the step size $ \eta_t $.
\State Update: $ \bx^\toptone \leftarrow \bx^\toptzero + \eta_t \bd^\toptzero $.
\EndFor
\State {\bfseries Return:}  $\bx\leftarrow \bx^\toptzero$;
\end{algorithmic} 
\end{algorithm}

\index{Gauss-Newton method}
\subsection{Gauss-Newton Method}
The \textit{Gauss-Newton} method forms the foundation for highly efficient techniques we will discuss in the upcoming sections, aimed at solving nonlinear least squares problems. It can be considered a variant of Newton's method that incorporates line search. This approach relies on using first derivatives of the vector function components or approximates the Hessian matrix by utilizing only the Jacobian matrix. In certain cases, it achieves quadratic convergence similar to Newton's method applied to general optimization problems.

\paragrapharrow{Gauss-Newton by affine approximation.}
The Gauss-Newton method employs an affine approximation of the components of $\br$ (an \textit{affine model} of $\br$) around the iteration point $ \bx^\toptzero $. 
Assuming that  $\br$ is twice continuously differentiable, according to the linear approximation theorem (Theorem~\ref{theorem:linear_approx}),
\begin{equation}\label{equation:gauss_new_affinemodel0}
\br(\bx^\toptzero + \bd) = \br(\bx^\toptzero) + \bJ(\bx^\toptzero) \bd + \mathcalO(\normtwo{\bd}^2) .
\end{equation}
This implies that
\begin{equation}\label{equation:gauss_new_affinemodel}
\ell_t(\bd) \triangleq \br(\bx^\toptzero) + \bJ(\bx^\toptzero) \bd
\end{equation}
is a good approximation  when $\bd$ is sufficiently small. Substituting this into the definition  of $f$ in   \eqref{equation:nonlinear_ols}, we obtain
$$
f(\bx^\toptzero + \bd) \approx \psi_t(\bd) \triangleq \tfrac{1}{2} \ell_t(\bd)^\top \ell_t(\bd) 
= f(\bx^\toptzero) + \bd^\top \bJ^\toptzeroTOP \br^\toptzero + \tfrac{1}{2} \bd^\top \bJ^\toptzeroTOP \bJ^\toptzero \bd ,
$$
where $\br^\toptzero \triangleq \br(\bx^\toptzero)$ and $\bJ^\toptzero \triangleq \bJ(\bx^\toptzero)$. The so-called \textit{Gauss-Newton step} $\bd_{\text{gn}}^\toptzero$ minimizes $\psi_t(\bd)$,
$$
\bd_{\text{gn}}^\toptzero = \argmin_{\bd} \left\{ \psi_t(\bd) \triangleq\tfrac{1}{2}\normtwo{\br^\toptzero + \bJ(\bx^\toptzero) \bd}^2\right\} .
$$
The approximation $\psi_t$ models the behavior of $f$ near the current iterate, and the gradient of $\psi_t$ is 
\begin{equation}\label{equation:gaun_nabla_psit}
\nabla \psi_t(\bd) = \bJ^\toptzeroTOP \br^\toptzero + ( \bJ^\toptzeroTOP \bJ^\toptzero ) \bd .
\end{equation}
By the first-order optimality condition (Theorem~\ref{theorem:fermat_fist_opt}),
the gradient is zero at a minimizer for $\psi_t$,  leading to the \textit{Gauss-Newton equation}:
\begin{equation}\label{equation:gn_first_equa}
	( \bJ^\toptzeroTOP \bJ^\toptzero ) \bd_{\text{gn}}^\toptzero = -\nabla f(\bx^\toptzero) = -\bJ^\toptzeroTOP \br^\toptzero .
\end{equation}
This shows that \eqref{equation:gn_first_equa} is the normal equation (Definition~\ref{definition:normal-equation-als})
for the overdetermined  linear system
$
\br(\bx^\toptzero) + \bJ(\bx^\toptzero) \bd = \bzero.
$
As mentioned previously, from linear algebra, we know that regardless of whether $\bJ^\toptzero$ has full rank, Equation~\eqref{equation:gn_first_equa} always has a solution, effectively representing a linear least squares problem (see Problem~\ref{prob:als_pseudo1}$\sim$\ref{prob:als_pseudon}). 
The optimality condition for the problem is:
\begin{equation}\label{equation:gaus_newton_ls}
\bd_{\text{gn}}^\toptzero =\argmin_{\bd} \quad \tfrac{1}{2} \normtwobig{\bJ^\toptzero \bd + \br^\toptzero}^2.
\end{equation}
In solving the linear least squares problem, we only need to perform the QR decomposition on $ \bJ^\toptzero $ (Theorem~\ref{theorem:qr-decomposition}), so the matrix $ \bJ^\toptzeroTOP \bJ^\toptzero $ does not need to be computed explicitly; see, for example, \citet{lu2021numerical} for more details.

\paragrapharrow{Gauss-Newton by approximated Hessian.}
We have outlined the Gauss-Newton method using an affine approximation of the components of $\br$.
Alternatively, since second derivative terms related to $r_i(\bx)$ in the Hessian matrix are challenging to compute (see Equation~\eqref{equation:gaus_new_hess}), the Gauss-Newton method simplifies this by ignoring these terms and directly employing $\bJ(\bx)^\top \bJ(\bx)$ as the approximate Hessian matrix to solve the Newton equation (\eqref{equation:newton_secon_approx_eq3} in Newton's method).
The descent direction $\bd_{\text{gn}}^\toptzero$ generated by the Gauss-Newton method satisfies (derived from the Newton equation in \eqref{equation:newton_secon_approx_eq3})
\begin{equation}\label{equation:gaus_new_normal}
	\begin{aligned}
		\nabla^2 f(\bx^\toptzero)\bd_{\text{gn}}^\toptzero &= -\nabla f(\bx^\toptzero);\\
		\underrightarrow{ \text{approximated by} }\quad
		\bJ^\toptzeroTOP \bJ^\toptzero \bd_{\text{gn}}^\toptzero &= -\bJ^\toptzeroTOP \br^\toptzero.
	\end{aligned}
\end{equation}

\paragrapharrow{The algorithm.}

If $\bJ^\toptzero$ has full column rank, then the matrix $\bA^\toptzero \triangleq \bJ^\toptzeroTOP \bJ^\toptzero$ also has full column rank and hence is positive definite. This ensures that $\psi_t(\bd)$ has a unique minimizer, which can be found by solving the Gauss-Newton equation \eqref{equation:gn_first_equa}. 
However, this does not guarantee that $\bx^\toptzero + \bd_{\text{gn}}^\toptzero$ is a minimizer of $f$, since $\psi_t(\bd)$  only approximates  $f(\bx^\toptzero + \bd)$. 
Nonetheless, Theorem~\ref{theorem:uncons_des_dir} guarantees that the solution $\bd_{\text{gn}}^\toptzero = -( \bJ^\toptzeroTOP \bJ^\toptzero )^{-1} \nabla f(\bx^\toptzero)$ is a descent direction. Consequently, $\bd_{\text{gn}}^\toptzero$ can be utilized as $\bd^\toptzero$ in descent methods (Algorithm~\ref{alg:struc_gd_gen}). 
An advantage of the Gauss-Newton method over pure Newton's method is that it does not always guarantee that $ \bd_{\text{n}}^\toptzero $ from the Newton equation is a descent direction. In contrast, the Gauss-Newton method employs a positive semidefinite approximation of the Hessian matrix, potentially yielding a more effective descent direction.

The typical step at $t$-th iteration is
$$
\begin{aligned}
\bd_{\text{gn}}^\toptzero&\leftarrow \text{solution of } ( \bJ^\toptzeroTOP \bJ^\toptzero ) \bd = -\bJ^\toptzeroTOP \br^\toptzero;  \\
\bx^\toptone &\leftarrow \bx^\toptzero + \eta_t \bd_{\text{gn}}^\toptzero.\\
\end{aligned}
$$
where $\eta_t$ is determined through line search. 
While the \textit{classical Gauss-Newton method (a.k.a., pure Gauss-Newton method)} uses a constant stepsize ($\eta_t = 1$), incorporating line search ensures convergence under certain conditions (Theorem~\ref{theorem:gaus_new_glob_conv}), specifically if:
\begin{enumerate}
\item The level set $\lev[f, \bx^{(1)}]\triangleq\{\bx \mid f(\bx) \leq \bx^{(1)}\}$ is bounded.
\item the Jacobian $\bJ(\bx^\toptzero)$ has full rank in all steps.
\end{enumerate}
The framework for the Gauss-Newton method is detailed in Algorithm~\ref{alg:gauss_newton}.

\paragrapharrow{Approximation issues.}
A natural question arises: Given that the Gauss-Newton method employs an approximate matrix to solve the Newton equation, under what circumstances is this approximation valid?
Intuitively, based on the expression for the Hessian matrix in  \eqref{equation:gaus_new_hess}, the approximation  $ \bJ^\toptzeroTOP \bJ^\toptzero $ is meaningful when it plays a dominant role.  
A sufficient condition for this is that at the optimal point $ \bx^* $, the residuals $ r_i(\bx^*) $ are very small. 
In such cases, the Gauss-Newton method closely resembles Newton's method and shares many of its properties. Conversely, if the residual vector $ \br(\bx^*) $ has a large norm,  using $ \bJ^\toptzeroTOP \bJ^\toptzero $ alone may not adequately approximate $ \nabla^2 f(\bx^\toptzero) $, potentially leading to slow convergence or even divergence of the Gauss-Newton method. See Section~\ref{section:quasi_largeresi} for a remedy.

Newton's method for optimization typically exhibits quadratic final convergence  (Theorem~\ref{theorem:conv_classNewton}). This level of convergence is generally not achieved by the Gauss-Newton method. 
However, if   $\br(\bx^*) = \bzero$, then near $\bx^*$, we have  $\nabla^2 \psi_t(\bzero) \approx \nabla^2 f(\bx^\toptzero)$,  which can also result in quadratic convergence with the Gauss-Newton method. 
Superlinear convergence can be expected if the functions $\{r_i\}$ exhibit small curvatures or if the values $\{\abs{r_i(\bx^*)}\}$ are small. Nonetheless, in general, one should anticipate linear convergence. 
Since $f(\bx^*) = \frac{1}{2} \normtwo{\br(\bx^*)}^2$, 
it is noteworthy that the value of $f(\bx^*) $ also influences the speed of convergence.

\begin{example}[Failure of Classical Gauss-Newton Methods \citep{madsen2010and}]
Consider the simple problem with $n=1$ and $m=2$:
$$
\br(x) = \begin{bmatrix} x + 1 \\ \lambda x^2 + x - 1 \end{bmatrix}, \quad f(x) = \frac{1}{2}(x+1)^2 + \frac{1}{2}(\lambda x^2 + x - 1)^2 .
$$
It follows that
$
f'(x) = 2\lambda^2 x^3 + 3\lambda x^2 - 2(\lambda - 1)x ,
$
indicating that  $x=0$ is a stationary point for $f$. 
Additionally, the second derivative is given by
$
f''(x) = 6\lambda^2 x^2 + 6\lambda x - 2(\lambda - 1) .
$
This shows that if $\lambda < 1$, then $f''(0) > 0$, so $x=0$ is a local minimizer---and indeed, it is the global minimizer (Remark~\ref{remark:charac_statpoint}).
The Jacobian is
$$
\bJ(x) = \begin{bmatrix} 1 \\ 2\lambda x + 1 \end{bmatrix} ,
$$
and the classical Gauss-Newton method ($\eta_t=1$) yields
$$
x_{\text{new}} = x - \frac{2\lambda^2 x^3 + 3\lambda x^2 - 2(\lambda - 1)x}{2 + 4\lambda x + 4\lambda^2 x^2} .
$$
Now, if $\lambda \neq 0$ and $x$ is close to zero, then
$$
x_{\text{new}} = x + (\lambda - 1)x + \mathcalO(x^2) = \lambda x + \mathcalO(x^2) .
$$
Thus, if $\abs{\lambda}< 1$, we observe  linear convergence. If $\lambda < -1$, then the classical Gauss-Newton method  fails to find the minimizer. For instance with $\lambda = -2$ and $x^{(1)} = 0.1$, the iterates display seemingly chaotic behavior:
$$
x^{(1)} = 0.1         \;\implies\;
x^{(2)} = -0.3029     \;\implies\;
x^{(3)} = 0.1368      \;\implies\;
x^{(4)} = -0.4680.
$$
Finally, if $\lambda = 0$, then $x_{\text{new}} = x - x = 0$, meaning the solution is found in one step. 
This is because $\br$  becomes an affine function in this case.
\end{example}

\paragrapharrow{Convergence analysis.}
Next, we present the convergence properties of the Gauss-Newton method. 
When $ \bJ^\toptzeroTOP \bJ^\toptzero $ dominates the Hessian matrix \eqref{equation:gaus_new_hess}, the Gauss-Newton method can achieve faster convergence rates.
To be more specific, similar to Newton's method, we also provide the local convergence results for the Gauss-Newton method.

\begin{theoremHigh}[Local Convergence  of Gauss-Newton]\label{theorem:gaus_new_local_conv}
Suppose $ r_i(\bx), i=1,2,\ldots,m $ are twice continuously differentiable, and $ \bx^* $ is the optimal solution to the nonlinear least squares problem \eqref{equation:nonlinear_ols}. If both the Hessian matrix $ \nabla^2 f(\bx) $ and its approximation $ \bJ(\bx)^\top \bJ(\bx) $ are  Lipschitz continuous in a neighborhood $\sB(\bx^*, \varepsilon)$ of $ \bx^* $. 
Given that the initial point is sufficiently close to $\bx^*$, the sequence $\{\bx^\toptzero\}_{t>0}$ generated by the Gauss-Newton method (Algorithm~\ref{alg:gauss_newton}) with a constant stepsize  $ \eta_t = 1 $ satisfies
$$
\normtwobig{\bx^\toptone - \bx^*} \leq C \normtwo{\big[(\bJ^*)^\top \bJ^*\big]^{-1} \bH^*} \normtwo{\bx^\toptzero - \bx^*} + \mathcalO(\normtwobig{\bx^\toptzero - \bx^*}^2),
$$
where $ \bH^*\triangleq \sum_{i=1}^m r_i(\bx^*) \nabla^2 r_i(\bx^*) $ is the part of the Hessian matrix $ \nabla^2 f(\bx^*) $ that remains after removing $ \bJ(\bx^*)^\top \bJ(\bx^*) $, $ C > 0 $ is a constant.
\end{theoremHigh}
\begin{proof}[of Theorem~\ref{theorem:gaus_new_local_conv}]
By the update rule and the fact that $\nabla f(\bx^*)=\bzero$, we have 
$$
\bx^\toptone - \bx^* = \bx^\toptzero + \bd_{\text{gn}}^\toptzero - \bx^*
= \big[\bJ^\toptzeroTOP \bJ^\toptzero\big]^{-1} 
\Big\{\bJ^\toptzeroTOP \bJ^\toptzero (\bx^\toptzero - \bx^*) + \nabla f(\bx^*) - \nabla f(\bx^\toptzero)\Big\}.
$$
By the fundamental theorem of calculus (Theorem~\ref{theorem:fund_theo_calculu}), we have 
$$
\begin{aligned}
&\nabla f(\bx^\toptzero) - \nabla f(\bx^*)
= \int_0^1 \nabla^2 f\big(\bx^* + \mu(\bx^\toptzero - \bx^*)\big) (\bx^\toptzero - \bx^*) d\mu\\
&= \int_0^1 \bJ^\top \bJ\big(\bx^* + \mu(\bx^\toptzero - \bx^*)\big) (\bx^\toptzero - \bx^*) d\mu +
\int_0^1 \bH\big(\bx^* + \mu(\bx^\toptzero - \bx^*)\big) (\bx^\toptzero - \bx^*) d\mu,
\end{aligned}
$$
where $ \bJ^\top \bJ(\bx) \triangleq \bJ^\top(\bx) \bJ(\bx) $, and $ \bH(\bx) \triangleq \nabla^2 f(\bx) - \bJ^\top \bJ(\bx) $ denotes the remaining part of the Hessian matrix. Combining the two equalities and taking norms, we have
$$
\begin{aligned}
&\normtwo{\big[\bJ^\toptzeroTOP \bJ^\toptzero\big] (\bx^\toptone - \bx^*) }=
\normtwo{\bJ^\toptzeroTOP \bJ^\toptzero (\bx^\toptzero - \bx^*) - \big(\nabla f(\bx^\toptzero) - \nabla f(\bx^*) \big)}\\
&\leq \int_0^1 \normtwo{\big[\bJ^\top \bJ(\bx^\toptzero) - \bJ^\top \bJ\big(\bx^* + \mu(\bx^\toptzero - \bx^*)\big)\big] (\bx^\toptzero - \bx^*)} d\mu \\
&\quad +
\int_0^1 \normtwo{\bH\big(\bx^* + \mu(\bx^\toptzero - \bx^*)\big) (\bx^\toptzero - \bx^*)} d\mu\\
&\leq \frac{L}{2} \normtwobig{\bx^\toptzero - \bx^*}^2 + C \normtwo{\bH^*} \normtwobig{\bx^\toptzero - \bx^*},
\end{aligned}
$$
where $ L $ is the Lipschitz constant of $ \bJ^\top \bJ(\bx) $. The last inequality holds because we use $ \bH^* $ to approximate $ \bH(\bx^* + \mu(\bx^\toptzero - \bx^*)) $. By continuity, there exists $ C > 0 $ and a neighborhood $ \sB(\bx^*, \varepsilon) $ of $ \bx^* $ such that for any $ \bx \in \sB(\bx^*, \varepsilon) $, $ \normtwo{\bH(\bx)} \leq C \normtwo{\bH(\bx^*)} $, which establishes the result.
\end{proof}
Theorem~\ref{theorem:gaus_new_local_conv} indicates that if $ \normtwo{\bH(\bx^*)} $ is sufficiently small, the Gauss-Newton method can achieve a linear rate of convergence; and when $ \normtwo{\bH(\bx^*)} = 0 $, the convergence rate is quadratic. However, if $ \normtwo{\bH(\bx^*)} $ is large, the Gauss-Newton method may fail to converge effectively.

From the preceding discussion, it is evident that the nonsingularity of the Jacobian matrix $ \bJ^\toptzero $ plays a crucial role. Thus, we establish the global convergence under this condition.
For the global convergence of the Gauss-Newton method, we impose the following uniform full-rank (uniform eigenvalues) assumption in the neighborhood of the level set \citep{sun2006optimization, liu2020optimization}.
\begin{assumption}[Restricted Strong Convexity]\label{assumption:uniform_sing_jac}
Assume that the Jacobian matrix $ \bJ(\bx) $ has eigenvalues uniformly greater than 0, i.e., there exists $ \alpha > 0 $ such that
\begin{equation}
\normtwo{\bJ(\bx) \bz} \geq \alpha \normtwo{\bz}, \quad \forall \bx \in \sB, \quad \forall \bz \in \real^n,
\end{equation}
where we assume $ \lev[f, f(\bx^{(1)})] $ is bounded, and $ \sB $ is a neighborhood set of the level  set
$
\sL\triangleq\lev[f, f(\bx^{(1)})] \triangleq \{ \bx \mid f(\bx) \leq f(\bx^{(1)}) \}
$
of the initial point $ \bx^{(1)} $.
Specifically,
$\sB(\sL, R)  \triangleq \{\bx \mid  \normtwo{\bx-\by} < R \text{ for some } \by\in \sL\}$ for some $R > 0$.
The assumption is related to the restricted strong convexity as defined in Definition~\ref{definition:res_scss_mat}.
\end{assumption}

Under the above assumptions, we have the following convergence theorem:

\begin{theoremHigh}[Global Convergence of Gauss-Newton]\label{theorem:gaus_new_glob_conv}
Let each residual function $ r_j $ be Lipschitz continuous on a bounded level set $ \sL \triangleq \lev[f, f(\bx^{(1)})] \triangleq \{ \bx \mid f(\bx) \leq f(\bx^{(1)}) \} $, and let the Jacobian matrix $ \bJ(\bx) $ satisfies the uniform full-rank condition in Assumption~\ref{assumption:uniform_sing_jac}. 
If the stepsizes $\eta_t$ satisfies the Wolfe conditions (Definition~\ref{definition:wolfe_cond}), then for the Gauss-Newton method (Algorithm~\ref{alg:gauss_newton}), the sequence $ \{ \bx^\toptzero \} $ satisfies
$$
\lim_{t \to \infty} \left\{\nabla f(\bx^\toptzero)=\bJ^\toptzeroTOP \br^\toptzero\right\} = \bzero.
$$
\end{theoremHigh}
\begin{proof}[of Theorem~\ref{theorem:gaus_new_glob_conv}]
Here we directly verify the Zoutendijk condition (Theorem~\ref{theorem:zoutendijk_cond}, the convergence under the Zoutendijk condition is guaranteed by Theorem~\ref{theorem:conv_line_search}). First, choose a bounded level set $ \lev[f, f(\bx^{(1)})] $ and its neighborhood $ \sB(\sL, R)$ sufficiently small so that there exist $ L > 0 $ and $ B > 0 $ such that for any $ \bx, \widetildebx\in \sB(\sL, R) $ and any $ j = 1, 2, \ldots, m $, the following conditions are satisfied:
\begin{align*}
\abs{r_j(\bx)} &\leq B, & \normtwo{\nabla r_j(\bx)} &\leq B, \\
\abs{r_j(\bx) - r_j(\widetildebx)} &\leq L \normtwo{\bx - \widetildebx}, & \normtwo{\nabla r_j(\bx) - \nabla r_j(\widetildebx)} &\leq L \normtwo{\bx - \widetildebx}.
\end{align*}
It follows that for any $ \bx \in \lev[f, f(\bx^{(1)})] $, there exists $ \widetildeB $ such that $ \normtwo{\bJ(\bx)} = \normtwo{\bJ(\bx)^\top} \leq \widetildeB $, and $ \nabla f(\bx) = \bJ(\bx)^\top \br(\bx) $ is a Lipschitz continuous function. Let $ \theta_t $ be the angle between the Gauss-Newton direction $ \bd_{\text{gn}}^\toptzero $ and the negative gradient direction, then
$$
\cos (\theta_t) = - \frac{ \nabla f(\bx^\toptzero)^\top \bd_{\text{gn}}^\toptzero }{ \normtwobig{\bd_{\text{gn}}^\toptzero} \normtwo{\nabla f(\bx^\toptzero)} } = \frac{ \normtwobig{\bJ^\toptzero \bd_{\text{gn}}^\toptzero}^2 }{ \normtwobig{\bd_{\text{gn}}^\toptzero} \normtwo{ \bJ^\toptzeroTOP \bJ^\toptzero \bd_{\text{gn}}^\toptzero} } \geq \frac{ \alpha^2 \normtwobig{\bd_{\text{gn}}^\toptzero}^2 }{ \widetildeB^2 \normtwobig{\bd_{\text{gn}}^\toptzero}^2 } = \frac{ \alpha^2 }{ \widetildeB^2 } > 0.
$$
By Theorem~\ref{theorem:zoutendijk_cond} and Theorem~\ref{theorem:conv_line_search}, it follows that $ \nabla f(\bx^\toptzero) \to \bzero $.
\end{proof}
We observe that the key assumption of Theorem~\ref{theorem:gaus_new_glob_conv} is the uniform full-rank condition in Assumption~\ref{assumption:uniform_sing_jac}. In fact, if $ \bJ^\toptzero $ does not have a  full rank, then the linear system \eqref{equation:gaus_new_normal} has infinitely many solutions. If no additional requirements are imposed on the solution properties, it cannot be concluded that $ \cos (\theta_t) $ is uniformly greater than zero. In this case, convergence may not hold.

\subsection{Levenberg-Marquardt Method}\label{section:nonlinearls_lmmethod}
Levenberg and later Marquardt proposed using a damped Gauss-Newton method (commonly referred to as the LM method, Section~\ref{section:modified_damp_new}) to solve nonlinear least squares problems \citep{levenberg1944method, marquardt1963algorithm}. 
In the \textit{Levenberg-Marquardt (LM) Gauss-Newton method} (or simply the LM method),
the step $\bd_{\text{lm}}^\toptzero$ at iteration point $\bx^\toptzero$ is defined by modifying the Gauss-Newton equation \eqref{equation:gn_first_equa} as follows: 
\begin{equation}\label{equation:damped_gauss_newt_eq}
\begin{aligned}
( \bJ^\toptzeroTOP \bJ^\toptzero + \mu_t \bI) \bd_{\text{lm}}^\toptzero &= -\bJ^\toptzeroTOP \br^\toptzero \\
\text{with } \bJ^\toptzero &\triangleq \bJ(\bx^\toptzero), \quad \br^\toptzero \triangleq \br(\bx^\toptzero), \quad \mu_t \geq 0 .
\end{aligned}
\end{equation}
The \textit{damping parameter} $\mu_t$ serves several purposes:
\begin{enumerate}
\item For all $\mu_t > 0$ the coefficient matrix $\bJ^\toptzeroTOP \bJ^\toptzero + \mu_t \bI$ is positive definite, ensuring that  $\bd_{\text{lm}}^\toptzero$ is a descent direction by Theorem~\ref{theorem:uncons_des_dir}.
\item For large values of $\mu_t$, we have 
$
\bd_{\text{lm}}^\toptzero \approx -\frac{1}{\mu_t} \bJ^\toptzeroTOP \br^\toptzero = -\frac{1}{\mu_t} \nabla f(\bx^\toptzero) ,
$
which is  a short step in the steepest descent direction. This is appropriate if the current iterate is far from the solution.
\item If $\mu_t$ is very small, then $\bd_{\text{lm}}^\toptzero \approx \bd_{\text{gn}}^\toptzero$, which is advantageous in the final stages of iteration when  $\bx^\toptzero$ is close to $\bx^*$. If $f(\bx^*) = \frac{1}{2} \normtwo{\br(\bx^*)}^2= \bzero$ (or very small), this can lead to (almost) quadratic final convergence.
\end{enumerate}

Therefore, the damping parameter affects both the direction and size of the step, allowing for an effective method without specific line search requirements. 
The initial value of $\mu_1$ should be related to the magnitude of the elements in  $\bA^{(1)} = \bJ^{(1)\top} \bJ^{(1)}$, such as
\begin{equation}
	\mu_1 = \tau \cdot \max_i \{ a_{ii}^{(1)} \} ,
\end{equation}
where $\tau$ is a hyperparameter. 
Although the algorithm is not highly sensitive to the choice of $\tau$, a small value like $\tau = 10^{-6}$ is recommended if $\bx^{(1)}$ is believed to be a good approximation of $\bx^*$. 
Otherwise, use $\tau = 10^{-3}$ or even $\tau = 1$ \citep{madsen2010and}. 
During iterations, the size of $\mu_t$ can be updated as described in the damped Newton method (Section~\ref{section:modified_damp_new}). The updating is controlled by the gain ratio
\begin{equation}\label{equation:nonlinearls_gainfactor}
\nu_t = \frac{f(\bx^\toptzero) - f(\bx^\toptone)}{\psi_t(\bzero) - \psi_t(\bd_{\text{lm}}^\toptzero)} ,
\quad
\text{where}\quad \bx^\toptone = \bx^\toptzero + \bd_{\text{lm}}^\toptzero,
\end{equation}
where $\psi_t(\bd) \triangleq 
f(\bx^\toptzero) + \bd^\top \bJ^\toptzeroTOP \br^\toptzero + \tfrac{1}{2} \big(\bd^\top \bJ^\toptzeroTOP +\textcolor{mylightbluetext}{\mu_t\bI}\big) \bJ^\toptzero \bd $.
A high  $\nu_t$ suggests  that $\psi_t(\bd_{\text{lm}}^\toptzero)$ is a good approximation of $f(\bx^\toptzero + \bd_{\text{lm}}^\toptzero)$, prompting a decrease in $\mu_t$ so that the next Levenberg-Marquardt step is closer to the Gauss-Newton step. 
Conversely, if $\nu_t$ is small or  even negative, indicating a poor approximation, $\mu_t$ should be increased to approach the steepest descent direction and reduce the step size.

As previously noted with the Gauss-Newton equation \eqref{equation:gn_first_equa}, the Gauss-Newton step  $\bd_{\text{gn}}^\toptzero$ solves the least squares problem
$\tfrac{1}{2} \normtwobig{\bJ^\toptzero \bd + \br^\toptzero}^2$.
Similarly, the LM damped Gauss-Newton equation in \eqref{equation:damped_gauss_newt_eq} corresponds to the normal equation for the augmented linear system
\begin{equation}\label{equation:damped_gsnew}
\bd_{\text{lm}}^\toptzero
=
\argmin_{\bd} \quad 
\frac{1}{2}
\normtwo{
	\begin{bmatrix} \bJ^\toptzero \\ 
		\sqrt{\mu_t} \bI \end{bmatrix} 
	\bd 
	+ 
	\begin{bmatrix} \br^\toptzero \\ \bzero \end{bmatrix}}^2.
\end{equation}
The most accurate solution can be found via orthogonal transformations using QR decomposition; see, for example, \citet{lu2021numerical}. However, since the solution
$\bd_{\text{lm}}^\toptzero$ is part of an iterative process, it does not need to be computed with high precision, making the normal equations approach in \eqref{equation:damped_gauss_newt_eq} preferable due to lower computational cost.

The complete  \textit{Levenberg-Marquardt (LM) Gauss-Newton method} is similar to the LM damped Newton method in Algorithm~\ref{algorithm:leve_dam_new}, except that it employs the damped Gauss-Newton equation \eqref{equation:damped_gauss_newt_eq} instead of the damped Newton equation \eqref{equation:damped_newton_eq}.

\subsection{Trust Region Method}
We introduced the trust region method in Section~\ref{section:des_trust_reg} for solving general optimization problems. Clearly, this method can be applied to nonlinear least squares problems as well. Recall that we established a connection between the trust region method and the damped Newton method in Theorem~\ref{theorem:damped_trust}. Specifically, for nonlinear least squares problems, this connection is made through the damped Gauss-Newton equation \eqref{equation:damped_gauss_newt_eq}. In each step, the trust region method for nonlinear least squares solves the following problem, called the nonlinear least squares subproblem (NLS):
\begin{equation}\label{equation:lms}
\textbf{(NLS)}:\qquad \min_{\bd} \quad \frac{1}{2} \normtwo{\bJ^\toptzero \bd + \br^\toptzero}^2 \quad \text{s.t.} \quad \normtwo{\bd} \leq \Delta_t,
\end{equation}
where $\Delta_t$ is the trust region radius at the $t$-th iteration.
In fact, NLS treats the above quadratic approximation as the trust region subproblem for $\psi_t(\bd)$ (see also \eqref{equation:trs_approx}):
\begin{equation}
\textbf{(NLS$'$)}:\qquad\begin{aligned}
\psi_t(\bd) 
&\triangleq f(\bx^\toptzero) + \nabla f(\bx^\toptzero)^\top \bd + \frac{1}{2} \bd^\top \bB^\toptzero \bd\\
&= \frac{1}{2} \normtwobig{\br^\toptzero}^2 + \bd^\top \bJ^\toptzeroTOP \br^\toptzero + \frac{1}{2} \bd^\top \bJ^\toptzeroTOP \bJ^\toptzero \bd\\
\implies 
\bd_{\text{tr}}^\toptzero &= \argmin_{\bd} \big\{ \psi_t(\bd)  \; \quad \text{s.t.} \quad \normtwo{\bd} \leq \Delta_t\big\}.
\end{aligned}
\end{equation}
This method uses $\bB^\toptzero = \bJ^\toptzeroTOP \bJ^\toptzero$ to approximate the Hessian matrix, which is derived from the Gauss-Newton method. For convenience, we omit the iteration index $t$. 
Subproblem \eqref{equation:lms} represents a trust region subproblem, and its properties were discussed in  Section~\ref{section:des_trust_reg}. Based on the optimality condition for TRS (Theorem~\ref{theorem:trs_kkt}), we can derive the following corollary:
\begin{corollary}[Optimality Condition of TRS for Nonlinear LS]
	The vector $\bd^*$ is a solution to the trust region subproblem
	\begin{equation}
		\min_{\bd} \quad \frac{1}{2} \normtwo{\bJ \bd + \br}^2 \quad \text{s.t.} \quad \normtwo{\bd} \leq \Delta
	\end{equation}
	if and only if $\bd^*$ is feasible and there exists a number $\lambda \geq 0$ such that
	\begin{subequations}
		\begin{align}
			(\bJ^\top \bJ + \lambda \bI) \bd^* &= -\bJ^\top \br, \\
			\lambda (\Delta - \normtwo{\bd^*}) &= 0\\
			(\bB + \lambda \bI) &\succeq \bzero. \label{equation:lms_kkt4}
		\end{align}
	\end{subequations}
	Note that $\bJ^\top \bJ$ is  positive semidefinite, so the condition \eqref{equation:lms_kkt4} is naturally satisfied.
\end{corollary}

To solve the (NLS) problem \eqref{equation:lms}, we follow a procedure similar to that used for trust region subproblems. First, determine $\lambda$ by solving a root-finding problem (see Section~\ref{section:solve_trs}), then obtain the search direction directly from the trust region equation.
Due to the special structure of the NLS, QR decomposition can be utilized, thus avoiding the need to compute $\bJ^\top \bJ + \lambda \bI$. Note that $(\bJ^\top \bJ + \lambda \bI)\bd = -\bJ^\top \br$ is actually a linear least squares problem:
\begin{equation}\label{equation:lm_ls}
	\min_{\bd} \quad \frac{1}{2}\normtwo{
		\begin{bmatrix} \bJ \\ 
			\sqrt{\lambda} \bI \end{bmatrix} 
		\bd 
		+ 
		\begin{bmatrix} \br \\ \bzero \end{bmatrix}}^2.
\end{equation}
The connection between NLS and the LM Gauss-Newton method become clear by \eqref{equation:damped_gsnew} and \eqref{equation:lm_ls}.

The coefficient matrix of this problem has a certain structure. Each time $\lambda$ is changed during trials, the block related to $\bJ$ remains unchanged, so there is no need to repeatedly perform QR decompositions on $\bJ$. Specifically, let $\bJ = \bQ\bR$ be the QR decomposition of $\bJ$, where $\bQ \in \real^{m \times n}$ and $\bR \in \real^{n \times n}$ (Theorem~\ref{theorem:qr-decomposition}).
A decomposition $\footnotesize\begin{bmatrix}
	\bJ \\
	\sqrt{\lambda} \bI
\end{bmatrix}$ of can be obtained by 
$$
\begin{bmatrix}
	\bJ \\
	\sqrt{\lambda} \bI
\end{bmatrix}
=
\begin{bmatrix}
	\bQ\bR \\
	\sqrt{\lambda} \bI
\end{bmatrix}
=
\begin{bmatrix}
	\bQ & \bzero \\
	\bzero & \bI
\end{bmatrix}
\begin{bmatrix}
	\bR \\
	\sqrt{\lambda} \bI
\end{bmatrix}.
$$
The matrix 
$
\footnotesize 
\begin{bmatrix}
	\bR \\
	\sqrt{\lambda} \bI
\end{bmatrix}
$
contains many zero elements. Using this feature, we can employ Householder transformations or Givens rotations to complete the QR decomposition of this matrix; see \citet{lu2021numerical} for more details. 
Thus, the linear least squares problem in \eqref{equation:lm_ls} can be solved using this QR decomposition \citep{lu2021numerical}.
The convergence of the trust region method can also be directly inferred from the convergence properties of the trust region method (Theorem~\ref{theorem:global_conv_trs_eta0} and Theorem~\ref{theorem:glo_conv_trs_othe}); therefore, we will not repeat these details here.

\index{QR decomposition}

\subsection{Powell's Dog Leg Method}\label{section:powell_dog_leg1}
An effective algorithm for addressing the trust region subproblem in nonlinear least squares problems is the Powell's Dog Leg method. This method, proposed by Powell, provides an approach to approximate the solution to the trust region subproblem $\bd_{\text{tr}}^\toptzero$, defined by \eqref{equation:trs_subpro}. The name ``Dog Leg" originates from golf, drawing a parallel between the step construction in the algorithm and the shape of a dogleg hole in golf \citep{powell1970new}.
Similar to the Levenberg-Marquardt Gauss-Newton method, this technique employs a combination of the Gauss-Newton and steepest descent directions but controls them explicitly through the radius of a trust region (refer to Section~\ref{section:des_trust_reg}):
$$
\bd_{\text{tr}}^\toptzero = \argmin_{\bd} \big\{ \psi_t(\bd)  \; \quad \text{s.t.} \quad \normtwo{\bd} \leq \Delta_t\big\}.
$$
where $\psi_t(\bd)\triangleq\tfrac{1}{2} \normtwobig{\br^\toptzero}^2 + \bd^\top \bJ^\toptzeroTOP \br^\toptzero + \frac{1}{2} \bd^\top \bJ^\toptzeroTOP \bJ^\toptzero \bd$ is an approximation to $f(\bx^\toptzero + \bd)$.

Given $\br : \real^n \rightarrow \real^m$ and the iteration point $\bx^\toptzero$, the Gauss-Newton step $\bd_{\text{gn}}^\toptzero$ is a least squares solution to the linear system
$
\bJ^\toptzero \bd +\br^\toptzero =\bzero
$
(see \eqref{equation:gaus_newton_ls}).
Meanwhile, the steepest descent direction is defined as the negative gradient:
$$
\bd_{\text{sd}}^\toptzero = -\nabla f(\bx^\toptzero) = -\bJ^\toptzeroTOP \br^\toptzero.
$$
However, the steepest descent direction is not properly scaled.
To scale it accordingly using the information from the Jacobian matrix, we consider the affine model
$$
\begin{aligned}
\br(\bx^\toptzero + \gamma \bd_{\text{sd}}^\toptzero) 
&\approx \br^\toptzero + \gamma \bJ^\toptzero \bd_{\text{sd}}^\toptzero\\
&\Downarrow\\
f(\bx^\toptzero + \gamma \bd_{\text{sd}}^\toptzero) 
&\approx \frac{1}{2} \normtwo{ \br^\toptzero + \gamma \bJ^\toptzero \bd_{\text{sd}}^\toptzero}^2\\
&= f(\bx^\toptzero) + \gamma  \br^\toptzeroTOP \bJ^\toptzero\bd_{\text{sd}}^\toptzero + \frac{1}{2} \gamma^2 \normtwo{\bJ^\toptzero \bd_{\text{sd}}^\toptzero}^2 .
\end{aligned}
$$
A closed-form scaling factor $\gamma_t$ can improve the solution at each iteration:
\begin{equation}
\gamma_t = - \frac{\br^\toptzeroTOP \bJ^\toptzero\bd_{\text{sd}}^\toptzero}{\normtwo{\bJ^\toptzero \bd_{\text{sd}}^\toptzero}^2} 
= \frac{\normtwobig{\bd_{\text{sd}}^\toptzero}^2}{\normtwo{\bJ^\toptzero \bd_{\text{sd}}^\toptzero}^2} .
\end{equation}

This yields two potential steps from the current point $\bx^\toptzero$: $\gamma_t \bd_{\text{sd}}^\toptzero$ and $\bd_{\text{gn}}^\toptzero$.
Given the trust region has radius $\Delta_t$, Powell suggested to use the following strategy for choosing the step; see Figure~\ref{fig:dog_leg} for an illustration for the third case:
\noindent
\begin{figure}[H]
\centering
\begin{minipage}{0.415\textwidth}
\begin{figure}[H]
\centering
\vspace{-0.35cm} 
\subfigtopskip=2pt 
\subfigbottomskip=2pt 
\subfigcapskip=-5pt 
\includegraphics[width=0.9\textwidth]{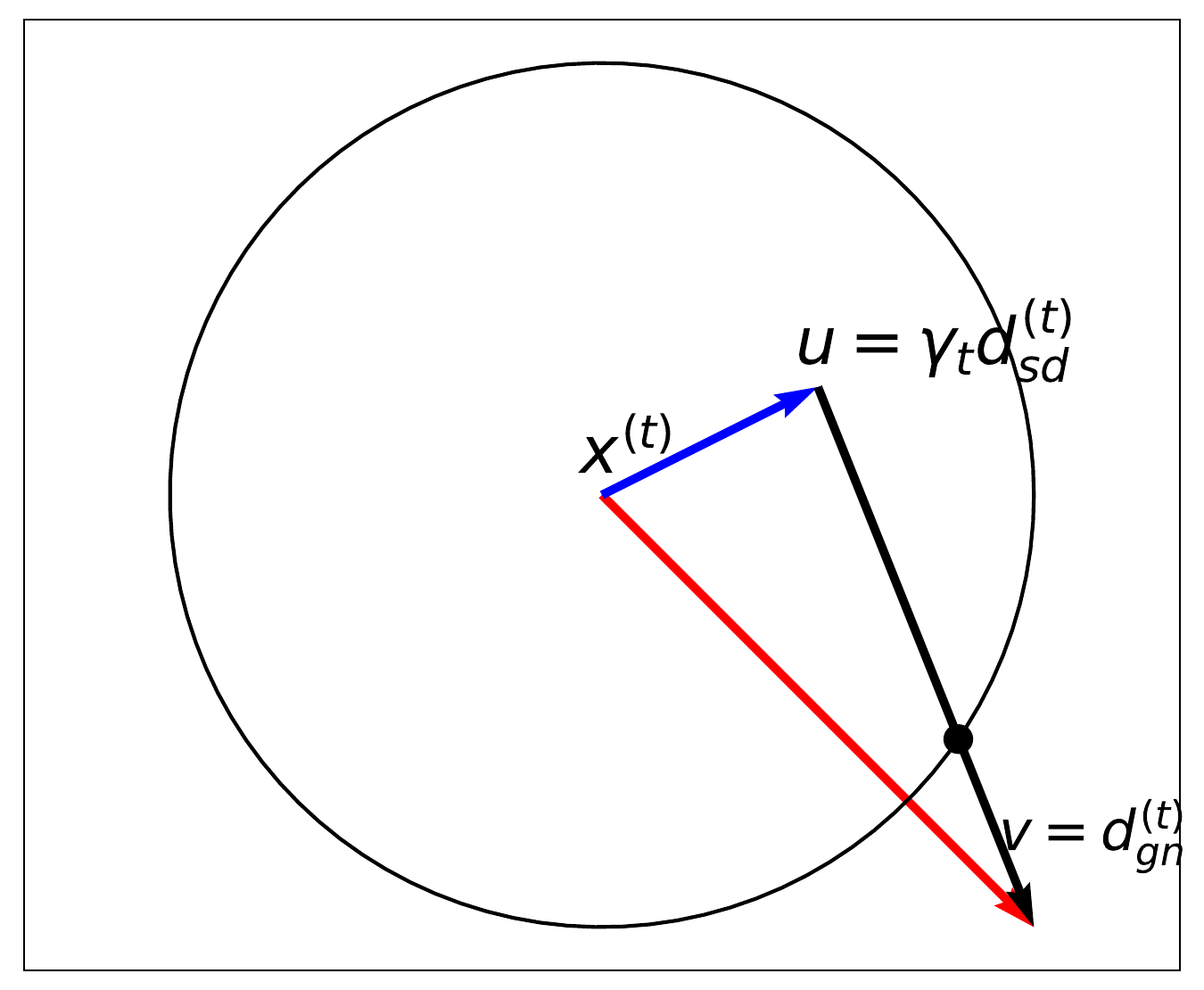}
\caption{Dog Leg step, the third case.}
\label{fig:dog_leg}
\end{figure}
\end{minipage}\hfill
\framebox{
\begin{minipage}{0.55\textwidth}
\begin{equation}\label{equation:dog_leg_step}
\begin{aligned}
&\text{if } \normtwobig{\bd_{\text{gn}}^\toptzero} \leq \Delta_t \\
&\qquad \bd_{\text{dl}} ^\toptzero\leftarrow \bd_{\text{gn}}^\toptzero; \\
&\text{elseif } \normtwobig{\gamma_t \bd_{\text{sd}}^\toptzero} \geq \Delta_t \\
&\qquad \bd_{\text{dl}}^\toptzero \leftarrow (\Delta_t / \normtwobig{\bd_{\text{sd}}^\toptzero}) \bd_{\text{sd}}^\toptzero ;\\
&\text{else} \\
&\qquad \bd_{\text{dl}}^\toptzero \leftarrow \gamma_t \bd_{\text{sd}}^\toptzero + \beta (\bd_{\text{gn}}^\toptzero - \gamma_t \bd_{\text{sd}}^\toptzero) \\
&\qquad \text{with } \beta \text{ chosen so that } \normtwobig{\bd_{\text{dl}}^\toptzero} = \Delta_t .
\end{aligned}
\end{equation}
\end{minipage}
}
\end{figure}

Letting
$\bu \triangleq \gamma_t \bd_{\text{sd}}^\toptzero$, 
$\bv \triangleq \bd_{\text{gn}}^\toptzero$ (such that $\normtwo{\bu}<\Delta_t$ and $\normtwo{\bv}>\Delta_t$), and defining $w \triangleq \bu^\top (\bv - \bu)$, we can express the third case as
$$
h(\beta) \triangleq \normtwo{\bu + \beta (\bv - \bu)}^2 - \Delta_t^2 = \normtwo{\bv - \bu}^2 \beta^2 + 2 w \beta + \normtwo{\bu}^2 - \Delta_t^2 .
$$
We seek a root for this second degree polynomial.
Since $h(0) = \normtwo{\bu}^2 - \Delta_t^2 < 0$ and $h(1) = \normtwo{\bv}^2 - \Delta_t^2 > 0$, $h$ has a root in $[0, 1]$; this is also observed in Figure~\ref{fig:dog_leg}. 
The choice of $\beta$ is then the root $\Big( -w \pm \sqrt{w^2 + \normtwo{\bv - \bu}^2 (\Delta_t^2 - \normtwo{\bu}^2)} \Big) / \normtwo{\bv - \bu}^2$
that lies in $[0,1]$.
As in the LM Gauss-Newton or damped Newton methods, we can use the gain factor
$$
\nu_t = \frac{f(\bx^\toptzero) - f(\bx^\toptzero + \bd_{\text{dl}}^\toptzero)}{\psi_t(\bzero) - \psi_t(\bd_{\text{dl}}^\toptzero)}
$$
to monitor the iteration. 

In the LM Gauss-Newton method (similar to Algorithm~\ref{algorithm:leve_dam_new}), we used the gain factor $\nu_t$ to control the size of the damping parameter. Similar to the decent method with trust region (Algorithm~\ref{alg:trust_region1}), here, we use it to control the radius $\Delta_t$ of the trust region. A large value of $\nu_t$ indicates that the linear model is good. We can increase $\Delta_t$ and thereby take longer steps, and they will be closer to the Gauss-Newton direction. 
Conversely, if $\nu_t$ is small or even negative, it indicates that the model may not be reliable, prompting us to decrease $\Delta_t$, leading to smaller steps that are closer to the steepest descent direction. 
The procedure is similar to the descent method with trust region in Algorithm~\ref{alg:trust_region1}, except we use the Dog Leg step in \eqref{equation:dog_leg_step} rather than the solution of trust region subproblem in \eqref{equation:trs_subpro}.

\subsection{A Secant Version of the LM Method}
We introduced the LM  method in Section~\ref{section:nonlinearls_lmmethod} for solving nonlinear least squares problems.
The methods discussed  assume that the vector function $\br$ is differentiable, i.e., the Jacobian
$
\bJ(\bx) = \left\{\frac{\partial r_i}{\partial x_j}\right\}_{ij}
$
exists. 
However, in many practical optimization scenarios, we may not have explicit formulas for the elements of $\bJ$, such as when $\br$ is treated as a ``black box''. 
For these cases, the secant version of the LM method is particularly useful.

The simplest remedy is to replace $\bJ^\toptzero \triangleq \bJ(\bx^\toptzero)$ by a matrix $\bG^\toptzero$ obtained through  numerical differentiation: The $(i,j)$-th element can be approximated  using finite differences:
\begin{equation}\label{equation:lmseca_fini_diffG}
	\frac{\partial r_i}{\partial x_j}(\bx^\toptzero) \approx g_{ij}^\toptzero \triangleq \frac{r_i(\bx^\toptzero + \delta \be_j) - r_i(\bx^\toptzero)}{\delta} ,
\end{equation}
where $\be_j$ is the unit vector in the $j$-th coordinate direction and $\delta$ is an appropriately small real number.  This strategy requires  $n+1$ evaluations of $\br$ for each iteration $\bx^\toptzero$. 
Since $\delta$ is probably much smaller than the distance $\normtwo{\bx^\toptzero - \bx^*}$, this approach provides limited additional information about the global behavior of  $\br$ than we would get from just evaluating $\br(\bx^\toptzero)$. To improve efficiency, alternative strategies are needed.

Now, consider the linear model from the Gauss-Newton method for $ \br : \real^n \rightarrow \real^m $,
$
\br(\bx^\toptzero +\bd)\approx \ell_t(\bd) \triangleq \br^\toptzero + \bJ^\toptzero \bd
$, where $\br^\toptzero\triangleq \br(\bx^\toptzero)$.
We will replace it by
$$
\br(\bx^\toptzero + \bd) \approx h_t(\bd) \triangleq \br^\toptzero + \bG^\toptzero \bd  , 
$$
where $ \bG^\toptzero $ approximates  $ \bJ^\toptzero $. 
Letting $\bd^\toptzero =\argmin_{\bd} h_t(\bd)$, 
we need $ \bG^\toptone $ for the next iteration so that
$$ 
\br(\bx^\toptone + \bd) \approx  h_{t+1}(\bd) \triangleq \br^\toptone + \bG^\toptone \bd , 
\quad \text{with }\bx^\toptone \triangleq\bx^\toptzero + \bd^\toptzero.
$$
Especially, we want this model $h_{t+1}(\bd)$ to hold with equality for $ \bd = - \bd^\toptzero = \bx^\toptzero - \bx^\toptone $, establishing the equality
\begin{equation}
	\br^\toptzero = \br^\toptone + \bG^\toptone (\bx^\toptzero - \bx^\toptone)  .
\end{equation}
This equation  gives us $ m $ equations in the $ mn $ unknown elements of $ \bG^\toptone $, necessitating additional conditions.
Similar to the symmetric rank-one update quasi-Newton method (Section~\ref{section:quasi_new_update}), \citet{broyden1965class} suggested supplementing the above equality with
\begin{equation}
	\bG^\toptone \bv = \bG^\toptzero \bv \quad \text{for all} \quad \bv \perp (\bx^\toptzero - \bx^\toptone)  .
\end{equation}
This yields \textit{Broyden's rank-one update} for the Jacobian approximation.
\begin{definition}[Broyden's Rank-One Update\index{Broyden's rank-one update}]\label{definition:broyd_roulm}
Let $\bG^\toptzero$ be the approximation of the Jacobian matrix $\bJ^\toptzero$, and let $\bd^\toptzero =\argmin_{\bd} h_t(\bd)$ such that $\bx^\toptone = \bx^\toptzero + \bd^\toptzero$.
Then the Broyden's rank-one update for the next Jacobian matrix is
$$
\bG^\toptone = \bG^\toptzero+ \ba \bd^\toptzeroTOP,
$$
where
$ \ba \triangleq \frac{1}{\bd^\toptzeroTOP \bd^\toptzero} \left( \br^\toptone - \br^\toptzero - \bG^\toptzero\bd^\toptzero \right) .$
\end{definition}

A brief sketch of the central part of the LM damped Newton method in Algorithm~\ref{algorithm:leve_dam_new} with this modification has the form
\begin{subequations}
\begin{align}
\bd_{\text{slm}}^\toptzero &\leftarrow \text{solve }  (\bG^\toptzeroTOP \bG^\toptzero + \mu_t \bI) \bd_{\text{slm}}^\toptzero = -\bG^\toptzeroTOP \br^\toptzero;\\
\bx^\toptone &\leftarrow \bx^\toptzero + \bd_{\text{slm}}^\toptzero;\\
\bG^\toptone &\leftarrow \bG^\toptzero+ \ba \bd^\toptzeroTOP \text{ by Definition~\ref{definition:broyd_roulm}};\\
\text{Update }& \mu_{t+1} \text{ and } \bx^\toptone \text{ as in Algorithm~\ref{algorithm:leve_dam_new}}.\
\end{align}
\end{subequations}

\citet{powell1970new} has shown that if the sequence of iterate $  \{\bx^{(t)}\}_{t>0} $ converges to $ \bx^* $ and if the set of steps $ \{ \bd^\toptzero \triangleq \bx^\toptone - \bx^\toptzero \}_{t>0} $ satisfy the condition that $ \{ \bd^{(t-n+1)}, \ldots, \bd^{(t-1)}, \bd^\toptzero \} $ are linearly independent (they span the whole of $ \real^n $) for each $ t \geq n $, then the set of approximations $ \{ \bG^\toptzero \} $ converges to $ \bJ(\bx^*) $, irrespective of the choice of $ \bG^{(1)} $.

In practice, however, it often happens that the previous $ n $ steps do not span the whole of $ \real^n $, and there is a risk that after some iterations, the current $ \bG^\toptzero $ is such a poor approximation to the true Jacobian, that $ \nabla f(\bx^\toptzero) \approx -\bG^\toptzeroTOP \br^\toptzero $ might not even be a descent direction. 
In such cases,  $ \bx^\toptzero $ remains unchanged,  and the damping parameter $ \mu_t $ is increased (see Algorithm~\ref{algorithm:leve_dam_new}). 
The approximation $ \bG^\toptzero $ is changed, but may still be a poor approximation, leading to a further increase in $ \mu_t $. Eventually the process is stopped by $ \bd_{\text{slm}}^\toptzero $ being so small that the stopping criteria is satisfied (e.g., the (ST3) in \eqref{equation:des_stopcri1}), although $ \bx^\toptzero $ may be far from $ \bx^* $ \citep{madsen2010and}.

Several strategies have been proposed to address this issue, including occasional recomputations of $ \bG^\toptzero $ using finite differences at certain iterations.
\citet{madsen2010and} suggested employ the rank-one update determined by Definition~\ref{definition:broyd_roulm} with a cyclic, coordinate-wise series of updatings. 
Specifically, let $ \bd^\toptzero $ denote the current step, and let $ j $ be the current coordinate number. If the pseudo angle $ \theta_t $ between $ \bd^\toptzero $ and $ \be_j $ is ``large", then we compute a
finite difference approximation to the $j$-th column of $\bJ^\toptzero$. More precisely, this occurs when
\begin{equation}
\cos (\theta_t) = \frac{\abs{\be_j^\top \bd^\toptzero}}{\normtwo{\bd^\toptzero} \cdot \normtwo{\be_j}} < \zeta 
\qquad \iff \qquad 
\abs{\bd^\toptzero_j} < \zeta \normtwo{\bd^\toptzero} .
\end{equation}
Experiments indicate that choosing $\zeta=0.8$ yields good performance.

\subsection{A Secant Version of the Dog Leg Method}

The idea of using a secant approximation to the Jacobian can, of course, also be used in connection with the Dog Leg method from Section~\ref{section:powell_dog_leg1}. 
In this section, we will focus on the special case where  $m = n$ (the Newton-Raphson scenario in \eqref{equation:newton_raphson_nonlin}).
Similar to the relationship between BFGS  and DFP formulas in  quasi-Newton methods (Section~\ref{section:quasi_new_update}),
\citet{broyden1965class} not only provided the formula from Definition~\ref{definition:broyd_roulm} for updating the approximate Jacobian but also introduced a formula for updating an approximate inverse of the Jacobian:
$$
\bP(\bx) \approx \bJ(\bx)^{-1}.
$$
For iteration point $\bx^\toptzero$, the two formulas are
\begin{subequations}\label{equation:secan_dogleg}
\begin{align}
\bG^\toptone &= \bG^\toptzero + \frac{1}{\bd^\toptzeroTOP \bd^\toptzero} (\ba - \bG^\toptzero \bd^\toptzero) \bd^\toptzeroTOP , \\
\bP^\toptone &= \bP^\toptzero + \frac{1}{\bb^\top \ba} (\bd^\toptzero - \bP^\toptzero \ba) \bb^\top ,
\end{align}
\end{subequations}
where
$
\bd^\toptzero = \bx^\toptone - \bx^\toptzero =  \argmin_{\bd} h_t(\bd)$,  $\ba \triangleq \br^\toptone - \br^\toptzero$ ,  $ \bb \triangleq \bP^\toptzero \bd^\toptzero
$, and $\bP^\toptzero \triangleq\bP(\bx^\toptzero)$.
It is easy to verify that then the updates $\bG^\toptone$ and $\bP^\toptone$ satisfy $\bG^\toptone \bP^\toptone = \bI$.

Using these matrices, the steepest descent direction $\bd^\toptzero_{\text{sd}}$ and the Gauss-Newton step $\bd^\toptzero_{\text{gn}}$ (see Section~\ref{section:powell_dog_leg1}) are approximated by
\begin{equation}
\bd^\toptzero_{\text{ssd}} = -\bG^\toptzeroTOP \br^\toptzero 
\qquad \text{and} \qquad 
\bd^\toptzero_{\text{sgn}} = -\bP^\toptzero \br^\toptzero.
\end{equation}
The updates for Powell's Dog Let method in \eqref{equation:dog_leg_step}  can be  easily adapted  to use these approximations. The initial $\bG = \bG^{(1)}$ can be obtained using the difference approximation \eqref{equation:lmseca_fini_diffG}, while  $\bP^{(1)}$ can be computed as $(\bG^{(1)})^{-1}$.

\subsection{Quasi-Newton Methods for Large Residual Problems}\label{section:quasi_largeresi}
In the previous sections, we introduced the Gauss-Newton method, the LM method, and the trust region method (along with some of their variants).
Essentially, for the iteration point $\bx^\toptzero$, these methods approximate the Hessian matrix $\nabla^2 f(\bx^\toptzero)$ by $\bJ^\toptzeroTOP\bJ^\toptzero$ and ignore the the residual term $\sum_{i=1}^{m} r_i(\bx^\toptzero) \nabla^2 r_i(\bx^\toptzero)$ (see \eqref{equation:gaus_new_hess}). 
Alternatively stated, they use the affine model in \eqref{equation:gauss_new_affinemodel}.
These approximations are highly effective for small residual least squares problems. However, in large residual problems, the second part of the Hessian matrix $\nabla^2 f(\bx^\toptzero)$ cannot be neglected. 
Simply considering $\bJ^\toptzeroTOP \bJ^\toptzero$ as the $t$-th step Hessian matrix approximation can lead to significant errors. In such cases, the aforementioned methods may fail.

Nevertheless, we can treat the nonlinear least squares problem \eqref{equation:nonlinear_ols} as an unconstrained problem and solve it using the previously discussed Newton's method and quasi-Newton methods (\S~\ref{section:new_methods}, \S~\ref{section:quasi_newton_method}). 
However, for many problems, calculating the Hessian matrices $\nabla^2 r_i(\bx^\toptzero)$ for each residual component is challenging, making  Newton's method very expensive. Directly applying the quasi-Newton method to approximate the Hessian matrix $\nabla^2 f(\bx^\toptzero)$ seems to overlook the special structure of the nonlinear least squares problem.

Given this, we introduce a hybrid approach: continue using the aforementioned methods for calculating the first term in the Hessian and employ a quasi-Hessian method for estimating the residual term \citep{liu2020optimization}.

The expression for the Hessian matrix \eqref{equation:gaus_new_hess} shows that  $\nabla^2 f(\bx^\toptzero)$ consists of two parts: one is easy to obtain but not enough to approximate the Hessian well; the other is more difficult to obtain but necessary for accurate calculations. 
For the easy part, we can directly retain the Gauss-Newton matrix $\bJ^\toptzeroTOP \bJ^\toptzero$, and for the more challenging part, we apply the quasi-Newton method for approximation:
$$
\nabla^2 f(\bx^\toptzero) = 
\underbrace{\bJ(\bx^\toptzero)^\top \bJ(\bx^\toptzero)}_{\text{Gauss-Newton}}
 +
\underbrace{\sum_{i=1}^{m} r_i(\bx^\toptzero) \nabla^2 r_i(\bx^\toptzero)}_{\text{quasi-Newton}}. 
$$
This forms the basis of our strategy for solving large residual problems, balancing the structure of the Hessian matrix and computational cost through a hybrid approximation method.

Specifically, we denote $\bH^\toptzero$ as the approximate Hessian matrix $\nabla^2 f(\bx^\toptzero)$, i.e.,
$ \bH^\toptzero = \bJ^\toptzeroTOP \bJ^\toptzero + \bR^\toptzero, $
where $\bR^\toptzero$ is the approximation of the second part of the Hessian matrix $\sum_{j=1}^{m} r_j(\bx^\toptzero) \nabla^2 r_j(\bx^\toptzero)$.
The key challenge lies in constructing the matrix $\bR^\toptzero$. Recall that when constructing the quasi-Newton method, the process mainly involves two steps: finding the quasi-Newton conditions and updating the low-rank matrix based on these conditions. Here, we follow a similar process but note that $\bR^\toptzero$ is only a part of the quasi-Newton matrix $\bH^\toptzero$ and may not satisfy the secant conditions \eqref{equation:secant1} (replacing $\bH^\toptone$ with $\bR^\toptone$ due to the symmetry of $\{\bR^\toptzero\}$). Our goal is to make $\bR^\toptone$ as close as possible to the second part of the Hessian matrix, i.e.,
$ \bR^\toptone \approx \sum_{j=1}^{m} r_j(\bx^\toptone) \nabla^2 r_j(\bx^\toptone). $

According to linear approximation theorem (Theorem~\ref{theorem:linear_approx}), the gradient function $\nabla r_i(\bx)$ for each $i\in\{1,2,\ldots,m\}$ can be approximated near $\bx^\toptone$ as:
$$
\nabla r_i(\bx) = \nabla r_i(\bx^\toptone) + \nabla^2 r_i(\bx^\toptone)(\bx - \bx^\toptone) + \mathcalO\big(\normtwobig{\bx - \bx^\toptone}^2\big).
$$
Let $\bx \triangleq \bx^\toptzero$, $\bh^\toptzero \triangleq \bx^\toptone - \bx^\toptzero$, and $\by_i^\toptzero \triangleq \nabla r_i(\bx^\toptone) - \nabla r_i(\bx^\toptzero)$.
Then,
$$
\nabla^2 r_i(\bx^\toptone) \bh^\toptzero + \mathcalO(\normtwobig{\bh^\toptzero}^2) = \by_i^\toptzero.
$$
This establishes the secant equation for the second term of $\nabla^2 f(\bx^\toptone)$:
$$
\small
\begin{aligned}
	\bR^\toptone \bh^\toptzero
	&= \bigg(\sum_{j=1}^{m} r_j(\bx^\toptone) \nabla^2 r_j(\bx^\toptone) \bigg) \bh^\toptzero
	= \sum_{j=1}^{m} r_j(\bx^\toptone) \left( \nabla^2 r_j(\bx^\toptone) \right) \bh^\toptzero \\
	&\approx \sum_{j=1}^{m} r_j(\bx^\toptone) \left( \nabla r_j(\bx^\toptone) - \nabla r_j(\bx^\toptzero) \right)
	= (\bJ^\toptone)^\top \br^\toptone - (\bJ^\toptzero)^\top  \br^\toptone
	\triangleq\by^\toptzero.
\end{aligned}
$$
Since $\bR^\toptone$ is also symmetric, the matrix updating methods in Section~\ref{section:quasi_new_update} can be applied using the above secant equation; and we shall not repeat the details.

\index{Linear programming}
\section{Nonlinear Least Squares Problem with Other Norms}
In the previous section, we discussed methods for nonlinear least squares, where the objective function is the $\ell_2$ norm of a vector of nonlinear functions; see \eqref{equation:nonlinear_ols}.
In this section, we explore methods for parameter estimation using other norms. Specifically, given a vector function $\br : \real^n \rightarrow \real^m$, our goal is to find $\bx^*$ that minimizes certain measures of $\br(\bx)$, such as $\normone{\br(\bx)}$ or $\norminf{\br(\bx)}$.

The function $\br$ may depend nonlinearly on $\bx$. Similar to the Gauss-Newton method with affine approximation in \eqref{equation:gauss_new_affinemodel0}, provided that $\br$ is twice continuously differentiable, we approximate $\br$ at iteration $\bx^\toptzero$ by 
\begin{equation}\label{equation:nonlinear_affapprox_eq1}
\br(\bx^\toptzero + \bd) \approx {\ell_t}(\bd) \triangleq \br(\bx^\toptzero) + \bJ(\bx^\toptzero) \bd ,
\end{equation}
where $\nabla r_i \in \real^n$ is the gradient of $r_i$ and $\bJ(\bx)=[\nabla r_1, \nabla r_2, \ldots, \nabla r_m]^\top \in \real^{m \times n}$ is the Jacobian matriux of $\br$. 
We then provide details of the trust region method (Algorithm~\ref{alg:trs_diffnorms}; see also Section~\ref{section:des_trust_reg}) for fitting $\ell_{\infty}$ and $\ell_1$ norms \citep{hald1981combined, hald1985combined}.
At each iteration $t$, the the trust region subprblem is 
$$
\bd_{\text{tr}}^{\toptzero} = \argmin_{\norminf{\bd} \leq \Delta_t} \left\{\psi_t(\bd) \triangleq \normbig{ \br^\toptzero + \bJ^\toptzero \bd}\right\},
$$
where $\norm{\cdot}$ denotes some norm,  $\br^\toptzero\triangleq\br(\bx^\toptzero)$, and  $\bJ^\toptzero \triangleq\bJ(\bx^\toptzero)$.
The functions $ f $ and $ \psi_t $ are defined by $ f(\bx) = \norm{\br(\bx)}_p $ and $ \psi_t(\bd) = \norm{\ell_t(\bd)}_p $, respectively for different norms. 
Note that the trust region is bounded using the $\ell_\infty$ norm: $\norminf{\bd} \leq \Delta_t$, thus transforming the problems into linear programming problems; see the following sections.

\begin{algorithm}[h] 
\caption{Trust Region Method for Fitting Different Norms}
\label{alg:trs_diffnorms}
\begin{algorithmic}[1] 
\Require A function $f(\bx) = \norminf{\br(\bx)}$ or  $f(\bx) = \normone{\br(\bx)}$; 
\State {\bfseries Input:}  Set the maximum radius $ \Delta_{\max} $, initial radius $ \Delta_1 $, initial point $ \bx^{(1)} $, accept radius $\gamma\in[0,\frac{1}{4})$;
\For{$t=1,2,\ldots$}
\State $\bd_{\text{tr}}^{\toptzero} \leftarrow \argmin_{\norminf{\bd} \leq \Delta} \psi_t(\bd)$; \Comment{see Sections~\ref{section:nonlinear_infnorm} and \ref{section:nonlinear_lonenorm}}
\State $\nu_t \leftarrow \big(f(\bx^\toptzero) - f(\bx^\toptzero + \bd_{\text{tr}}^{\toptzero})\big) / \big(\psi_t(\bzero) - \psi_t(\bd_{\text{tr}}^{\toptzero})\big) $; \Comment{gain factor}
\If{$\nu_t > 0.75$ and $\normtwobig{\bd^\toptzero} =\Delta_t$} \Comment{very good step, and the step is at the border}
\State $\Delta_{t+1} \leftarrow \min\{2 \Delta_t, \Delta_{\max}\}$; \Comment{larger trust region}
\EndIf
\If{$\nu_t < 0.25$} \Comment{poor step}
\State $\Delta_{t+1} \leftarrow \Delta_t / 3$; \Comment{smaller trust region}
\EndIf
\If{$\nu_t > \gamma$} \Comment{reject step if $\nu_t \leq \gamma$}
\State $\bx^\toptone \leftarrow \bx^\toptzero + \bd_{\text{tr}}^{\toptzero}$;
\EndIf
\EndFor
\State {\bfseries Return:}  $\bx\leftarrow \bx^{(t)}$;
\end{algorithmic} 
\end{algorithm}

\subsection{Fitting in the $ \ell_\infty $ norm}\label{section:nonlinear_infnorm}
When the underlying norm is the $\ell_\infty$ norm, the trust region subproblem can be solved equivalently by a linear programming (LP) problem (Example~\ref{example:linear_program}; see, for example, \citet{boyd2004convex, bertsekas2015convex} for some algorithms of this problem).  
Omitting the superscripts,
we seek a minimizer $ \widetildebx $ of the function
$ f(\bx) = \norminf{\br(\bx)} = \max_i \abs{r_i(\bx)}   $.
The linearized problem
$$ 
\widetildebx = \argmin_{\norminf{\bd} \leq \Delta} \left\{ \psi(\bd) \triangleq \norminf{\br(\bx) + \bJ(\bx) \bd} \right\} 
$$
is solved by a linear programming  algorithm. 
To formulate the LP problem, we introduce the variable $ d_{n+1} = \psi(\bd) $ and the extended vector
$ \widetildebd = \small\begin{bmatrix} \bd \\ d_{n+1} \end{bmatrix} \normalsize\in \real^{n+1}  $. Using this auxiliary variable, $ \widetildebx $ can be found by solving the following linear programming problem
$$
\begin{aligned}
	& \min & & d_{n+1} \\
	& \text{s.t.} & & -\Delta \leq d_i \leq \Delta &i = 1,2 \ldots, m; \\
	& & & -d_{n+1} \leq r_i(\bx) + \nabla r_i(\bx)^\top \bd \leq d_{n+1}, \quad &i = 1,2 \ldots, m. \\
\end{aligned}
$$

\subsection{Fitting in the $\ell_1 $ norm}\label{section:nonlinear_lonenorm}
Similarly, when the underlying norm is the $\ell_1$ norm, the trust region subproblem can also be equivalently stated as an LP problem.
Again, omitting the superscripts, we seek a minimizer $ \widetildebx $ of the function
$ f(\bx) = \normone{\br(\bx)} = \sum_{i=1}^m  \abs{r_i(\bx)}  $. 
For the linearized problem
$$ 
\widetildebx = \argmin_{\norminf{\bd} \leq \Delta} \left\{ \psi(\bd) \triangleq \normone{\br(\bx) + \bJ(\bx) \bd} \right\},
$$
we introduce auxiliary variables $ d_{n+i} $, $ i=1,2\ldots,m $, and the extended vector
$ \widetildebd = \scriptsize\begin{bmatrix} \bd \\ d_{n+1} \\ \vdots \\ d_{n+m} \end{bmatrix}  \normalsize\in \real^{n+m}  . $
Using these auxiliary variables, $ \widetildebx $ can be found by solving the problem
$$
\begin{aligned}
	& \min & & \sum_{i=1}^m d_{n+i} \\
	& \text{s.t.} & & -\Delta \leq d_i \leq \Delta        &i = 1,2 \ldots, m; \\
	& & & -d_{n+i} \leq r_i(\bx) + \nabla r_i(\bx)^\top \bd \leq d_{n+i}, \quad &i = 1,2 \ldots, m. \\
\end{aligned}
$$

\subsection{Huber Estimation}\label{section:huber_estima}
\begin{SCfigure}
	\centering
	\includegraphics[width=0.5\textwidth]{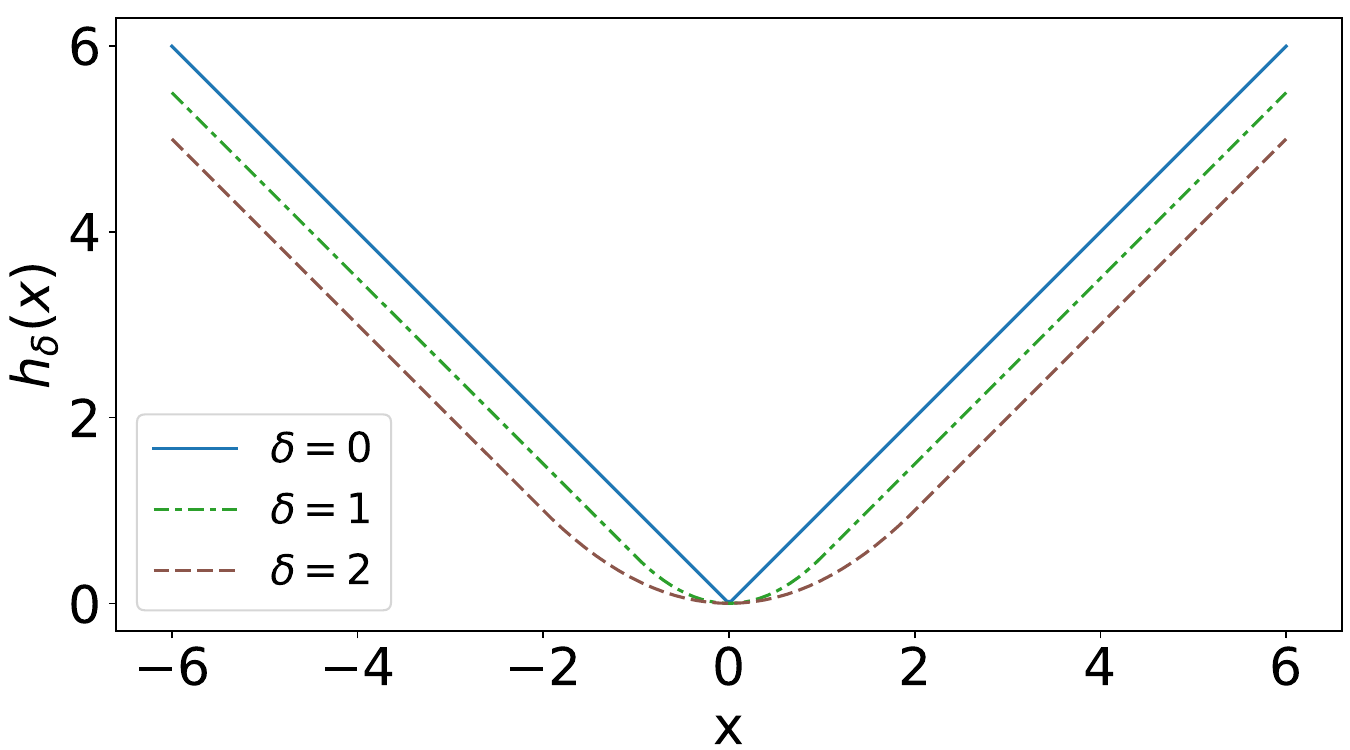} 
	\caption{Illustration of Huber loss function $h_\delta(x)$ for different values of $\delta$.}
	\label{fig:huber}
\end{SCfigure}

As mentioned in Chapter~\ref{chapter_introduction}, estimation using $\ell_1$ norms can produce more robust estimates, particularly in the presence of outliers  \citep{zoubir2012robust}.
\textit{Huber estimation}, on the other hand, combines the smoothness of the least squares estimator with the robustness of the $ \ell_1 $-estimator.
For a function $ \br : \real^n \rightarrow \real^m $, the Huber estimator $ \bx_\delta $ is defined by
\begin{equation}\label{equation:huberloss_eq1}
	\bx_\delta = \argmin_{\bx} \left\{ f_\delta(\bx) \triangleq \sum_{i=1}^m h_\delta(r_i(\bx)) \right\} ,
\end{equation}
where $ h_\delta $ is the \textit{Huber loss function},
\begin{equation}
	h_\delta(x) = 
	\begin{cases} 
		\frac{1}{2\delta}x^2, & \abs{x} < \delta, \\
		\abs{x} - \frac{\delta}{2}, & \text{otherwise}.
	\end{cases}
\end{equation}
When $\delta \to 0$, the smooth function $h_\delta(x)$ approaches the absolute value function $\abs{x}$. Figure~\ref{fig:huber} shows the graph of $h_\delta(x)$ for different values of $\delta$.

The threshold $ \delta $ is used to distinguish between ``small'' and ``large'' function values (residuals). Based on the values of the $ r_i(\bx) $ we define a generalized sign vector $ \bs \triangleq \bs_\delta(\bx) $ and an ``activity matrix'' $ \bW \triangleq \bW_\delta(\bx) = \text{diag}(w_1, w_2 \ldots, w_m) $. Note that $ w_i = 1 - s_i^2 $ and
\begin{table}[H]
	\centering
	\caption{Parameters using Huber loss functions.}
	\label{table:huber_sw}
	\begin{tabular}{|c|ccc|}
		\hline
		& $ r_i(\bx) < -\delta $ & $ \abs{r_i(\bx)} \leq \delta $ & $ r_i(\bx) > \delta $ \\
		\hline
		$ s_i(x) $ & $-1$ & $0$ & $1$ \\\hline
		$ w_i(x) $ & $0$ & $1$ & $0$ \\
		\hline
	\end{tabular}
\end{table}

Using these definitions, the Huber objective function in \eqref{equation:huberloss_eq1} can be expressed as
\begin{equation}\label{equation:huber_objfunc}
	f_\delta(\bx) = \frac{1}{2\delta} \br^\top \bW \br + \br^\top \bs - \frac{1}{2} \delta \bs^\top \bs ,
\end{equation}
where we have omitted the argument $ \bx $ and the Huber parameter $ \delta $. The gradient is
\begin{equation}
	\nabla f_\delta (\bx) = \frac{1}{\delta} \bJ^\top (\bW \br + \delta \bs) .
\end{equation}
By Theorem~\ref{theorem:fermat_fist_opt}, a necessary condition for $ \bx_\delta $ being a minimizer of $ f_\delta $ is that it satisfies the equation
\begin{equation}
\nabla f_\delta (\bx) =\bzero  
\quad\implies\quad
\bJ(\bx)^\top \big(\bW_\delta(\bx) \br(\bx) + \delta \bs_\delta(\bx)\big) = \bzero.
\end{equation}
Before delving into the nonlinear Huber estimation, we first describe the linear Huber estimation.
\paragrapharrow{Linear Huber estimation.}
Consider the linear function
\begin{equation}
\br(\bx) = \bb - \bA \bx ,
\end{equation}
where $ \bA \in \real^{m \times n} $ and $ \bb \in \real^m $ and are known variables. 
The gradient is
\begin{equation}
	\nabla f_\delta(\bx) = -\frac{1}{\delta} \bA^\top \big(\bW(\bx) (\bb - \bA \bx) + \delta \bs(\bx)\big) .
\end{equation}
While the Hessian is positive semidefinite:
\begin{equation}
 \bH(\bx)  \triangleq \nabla^2 f_\delta(\bx) = \frac{1}{\delta} \bA^\top \bW(\bx) \bA.
\end{equation}
By Theorem~\ref{theorem:second_nec_nonstrict_loca}, the positive semidefiniteness implies that if we find a stationary point $ \bx^* $ satisfying $ \nabla f_\delta(\bx^*) = \bzero $, then it is also a minimizer for $ f_\delta $. 
We use Newton's method with line search to find a stationary point.
The typical step at $t$-th iteration is
\begin{subequations}
\begin{align}
	\bd_{\text{lh}}^\toptzero&\leftarrow \text{solution of } \bH(\bx^\toptzero) \bd = -\nabla f_\delta(\bx^\toptzero);  \\
	\eta_t &\leftarrow \argmin_{\eta} f_\delta(\bx^\toptzero + \eta \bd_{\text{lh}}^\toptzero);\\
	\bx^\toptone &\leftarrow \bx^\toptzero + \eta_t \bd_{\text{lh}}^\toptzero.
\end{align}
\end{subequations}
The line search is very simple: $ g(\eta) \triangleq f_\delta(\bx^\toptzero + \eta \bd_{\text{lh}}^\toptzero) $ is a piecewise quadratic in the scalar variable $ \eta $, so $ \eta_t $ is a root for $ g' $, which is piecewise linear.

\paragrapharrow{Nonlinear Huber estimation.}
When the  function $\br$  depends nonlinearly on $\bx$, we must use approximation methods to find the minimizer $\bx^*$ of $f_\delta$. 
Following the affine approximation in \eqref{equation:nonlinear_affapprox_eq1}, at the point $\bx$ we use the affine approximation $\br(\bx+\bd) \approx \bl \triangleq \br(\bx) + \bJ(\bx) \bd$ and the corresponding approximation to the objective function from \eqref{equation:huber_objfunc}
$$
f_\delta(\bx + \bd) \approx \psi(\bd)
= \frac{1}{2\delta} \bl^\top \bW \bl + \bl^\top \bs - \frac{1}{2} \delta \bs^\top \bs ,
$$
where $\bs$ and $\bW$ are given by Table~\ref{table:huber_sw} with $\br(\bx)$ replaced by $\bl$. Therefore, the typical step for the descent step at $t$-th iteration is 
$$
\bd_{\text{nh}}^\toptzero \leftarrow \argmin_{\bd} \left\{f_\delta(\bx^\toptzero + \bd) \approx \psi_t(\bd)\right\},
$$
where $\psi_t(\bd)$ represents the corresponding approximation at $\bx^\toptzero$.
For nonlinear Huber estimation, instead of using trust region methods, \citet{hald1981combined, madsen2010and} suggested using a Levenberg-Marquardt like algorithm (Algorithm~\ref{alg:nonlinear_huber_est}); see also Section~\ref{section:modified_damp_new} and Algorithm~\ref{algorithm:leve_dam_new} for reference.

\begin{algorithm}[h]
\caption{Nonlinear Huber Estimation Method}
\label{alg:nonlinear_huber_est}
\begin{algorithmic}[1]
\Require A twice continuously differentiable function $f(\bx)$; 
\State {\bfseries Input:}  Initialize $\bx^{(1)}, \mu_1$ (by default $\mu_1=1$),  $ \sigma$ (by default $\sigma\triangleq 10^{-3}$), Huber parameter  $\delta$;
\For{$t=1,2,\ldots$}
\While{$\nabla^2 \psi_t(\bd) + \mu_t \bI$ is not positive definite}
\State $\mu_t \leftarrow 2\mu_t$
\EndWhile
\State $\bd_{\text{nh}}^\toptzero \leftarrow \argmin_{\bd} \left\{ \psi_t(\bd)\right\}$;
\State Gain factor $\nu_t \leftarrow \big(f_\delta(\bx^\toptzero) - f_\delta(\bx^\toptzero + \bd_{\text{nh}}^\toptzero)\big) / \big(\psi_t(\bzero) - \psi_t(\bd_{\text{nh}}^\toptzero)\big)$;
\If{$\nu_t > \sigma$} \Comment{$f$ decreases}
\State $\bx^\toptone \leftarrow \bx^\toptzero + \bd_{\text{dn}}^\toptzero$ \Comment{new iterate}
\State $\mu_{t+1} \leftarrow  \max\left\{\frac{1}{3}, 1 - (2\nu_t - 1)^3\right\} \cdot \mu_t$ 
\Else
\State $\bx^\toptone \leftarrow \bx^\toptzero $;  \Comment{old iterate}
\State $\mu_{t+1} \leftarrow 2\cdot \mu_t $;
\EndIf
\EndFor
\State {\bfseries Return:}  $\bx\leftarrow \bx^{(t)}$;
\end{algorithmic}
\end{algorithm}

\section{Separable Nonlinear Least Squares Problem}
In many data fitting problems, the model includes parameters that appear both linearly and nonlinearly. Such problems are referred to as \textit{separable}. This property can be leveraged in least squares estimation, and various specialized algorithms have been proposed, see e.g. \citet{ruhe1980algorithms, niel2000separable}. In this section, we discuss how to reformulate the problem to exploit this property using a standard program for nonlinear least squares.

\subsection{Full Parameter Set}
Let $\mathcalX = \left\{(\zeta_1, y_1), (\zeta_2, y_2), \ldots, (\zeta_m, y_m)\right\}$ be the given data points, and consider the model
\begin{equation}
y_i = G(\bx, \bmu, \zeta_i) \triangleq \mu_1 g_1(\bx, \zeta_i) + \ldots + \mu_p g_p(\bx, \zeta_i) , \quad i\in\{1,2,\ldots,m\},
\end{equation}
where  $g_i : \real^{n+1} \rightarrow \real$ are nonlinear functions that depend on $\bx \in \real^n$ and  $\bzeta$. The coefficients $\mu_1, \mu_2, \ldots, \mu_p$ are parameters that occur linearly in the model.
Collecting all the unknown parameters into a single vector  $\bz \triangleq \footnotesize\begin{bmatrix} \bx \\ \bmu \end{bmatrix} \in \real^{n+p}$, the least squares estimate is defined as
\begin{equation}
\bz^* = \argmin_{\bz \in \real^{n+p}} \left\{ f(\bz) \triangleq \tfrac{1}{2} \br(\bz)^\top \br(\bz) \right\} ,
\end{equation}
where
\begin{equation}\label{equation:sepe_rg}
\br(\bz) = \by - \bG(\bx) \bmu \in\real^m 
\qquad \text{and}\qquad
\left\{\bG(\bx)\right\}_{ij} = \{g_j(\bx, \zeta_i)\}_{ij} \in\real^{m\times p}  
\end{equation}
The least squares estimate satisfies
\begin{equation}
\nabla f(\bz^*) = \bzero , \quad \text{where} \quad \nabla f(\bz) \triangleq \bJ(\bz)^\top \br(\bz) ,
\end{equation}
with the Jacobian $\bJ(\bz) \triangleq\begin{bmatrix} \bJ_g(\bz) & -\bG(\bx) \end{bmatrix} \in \real^{m \times (n+p)}$, where $\left( \bJ \right)_{ij} = \dfrac{\partial r_i}{\partial z_j}$.
The last $p$ columns correspond to the linear parameters $\bmu$, while each element in the first $n$ columns of $\bJ$ are given by
\begin{equation}
\left( \bJ_g(\bx) \right)_{ij} = \frac{\partial r_i}{\partial x_j}(\bz) 
= - \sum_{k=1}^p \mu_k \frac{\partial g_k}{\partial x_j}(\bx, \zeta_i)  .
\end{equation}
We shall use the following notation for the $j$-th ($j\in\{1,2,\ldots,n\}$) column in $\bJ_g$,
\begin{equation}\label{equation:sep_jaco_nablag}
\left( \bJ_g(\bz) \right)_{:,j} 
= \left(\nabla\br(\bz)\right)_{:,j}
= - \bG'_j(\bx) \bmu 
\quad \text{with} \quad 
\left\{\bG'_j(\bx)\right\}_{ik} = \left\{\frac{\partial g_k}{\partial x_j}(\bx, \zeta_i)\right\} \in\real^{m\times p},
\end{equation}
where $\left(\nabla\br(\bz)\right)_{:,j} 
= 
\frac{\partial \br(\bz)}{\partial z_j}
= 
\frac{\partial \br(\bz)}{\partial x_j}
\in\real^m
$ for $j\in\{1,2,\ldots,n\}$.

\paragrapharrow{Gauss-Newton method.}
As mentioned previously, the \textit{Gauss-Newton method} is widely used for solving nonlinear least squares problems: Given an approximation $\bz^\toptzero$ to $\bz^*$ at the $t$-th iteration, the next approximation is found as $\bz^\toptzero + \bd_{\text{gn}}^\toptzero$, where
\begin{equation}
\bd_{\text{gn}}^\toptzero = \argmin_{\bd \in \real^{n+p}} \left\{ \psi_t(\bd) \triangleq \tfrac{1}{2} \normtwo{\bJ(\bz^\toptzero) \bd +\br(\bz^\toptzero)}^2 \right\} ,
\end{equation}
leading to the \textit{Gauss-Newton} equation
\begin{equation}\label{equation:sep_gaussnew}
\textbf{(Gauss-Newton)}:\qquad \left( \bJ(\bz^\toptzero)^\top \bJ(\bz^\toptzero) \right) \bd_{\text{gn}}^\toptzero = \bJ(\bz^\toptzero)^\top \br(\bz^\toptzero)  .
\end{equation}
This equation  is often modified by the \emph{Levenberg-Marquardt} strategy, where $\bd_{\text{gn}}^\toptzero$ is replaced by $\bd_{\text{lm}}^\toptzero$, determined by
\begin{equation}\label{equation:sep_levmarg}
\textbf{(LM)}:\qquad \left( \bJ(\bz^\toptzero)^\top \bJ(\bz^\toptzero) + \mu_t \bI \right) \bd_{\text{lm}}^\toptzero = \bJ(\bz^\toptzero)^\top \br(\bz^\toptzero)  .
\end{equation}
The \textit{damping parameter} $\mu_t$ is adjusted during iteration according to how well $\psi_t(\bd_{\text{lm}}^\toptzero)$ approximates $f(\bz^\toptzero + \bd_{\text{lm}}^\toptzero)$ using the gain factor; see \eqref{equation:trs_reduc_ratio}.

\subsection{Reduced Parameter Set}

Instead of treating the elements in $\bmu$ as nonlinear variables, we replace \eqref{equation:sepe_rg} by
\begin{equation}\label{equation:sep_red_resid}
\br(\bz) = \by - \bG(\bx) \bmu 
\qquad\implies\qquad 
\br(\bx) = \by - \bG(\bx) \bmu(\bx) ,
\end{equation}
with $\bmu(\bx)\in\real^p$ computed as the least squares solution to $\bG(\bx) \bmu(\bx) = \by$. We can express it as the solution to the normal equation (Definition~\ref{definition:normal-equation-als}):
\begin{equation}
\bA(\bx) \bmu(\bx) = \bb(\bx) ,
\quad \text{where}\quad 
\bA(\bx) \triangleq \bG(\bx)^\top \bG(\bx) \in\real^{p\times p} , 
\quad \bb(\bx) \triangleq \bG(\bx)^\top \by  \in\real^p.
\end{equation}
With the notation from \eqref{equation:sep_jaco_nablag}, the corresponding $j$-th column of the Jacobian matrix $\bJ(\bx)=\nabla\br(\bx)\in\real^{m\times n}$ is given by
\begin{equation}\label{equation:sep_jaco_j}
\left( \bJ(\bx) \right)_{:,j} =  \left(\nabla\br(\bx)\right)_{:,j} 
= - \bG_j'(\bx) \bmu(\bx) - \bG(\bx) \bmu_j'(\bx) 
\triangleq \left( \bJ_g(\bx) \right)_{:,j} + \left( \bH(\bx) \right)_{:,j} .
\end{equation}
From the above notation, for any $j \in\{1,2, \ldots, n\}$, we have
$$
\bA(\bx) \bmu(\bx) = \bb(\bx)
\qquad\implies\qquad 
\bA_j'(\bx) \bmu(\bx) + \bA(\bx) \bmu_j'(\bx) = \bb_j'(\bx) ,
$$
where $\bA_j'(\bx) \triangleq \frac{\partial \bA(\bx)}{\partial x_j}  = \bG_j'(\bx)^\top \bG(\bx) + \bG(\bx)^\top \bG_j'(\bx) \in\real^{p\times p}$, 
$\bmu_j'(\bx)\triangleq \frac{\partial \bmu(\bx)}{\partial x_j}\in\real^p$, and $\bb_j'(\bx) \triangleq\bG_j'(\bx)^\top \by \in\real^p$.
Omitting the argument $\bx$, we have 
\begin{equation}\label{equation:sep_redu_parAmuj}
\begin{aligned}
	\bA \bmu_j' &= \bb_j' - \bA_j' \bmu 
	= \bG_j'^\top \by - \left( \bG_j'^\top \bG + \bG^\top \bG_j' \right) \bmu  \\
	&= \bG_j'^\top \br - \bG^\top \bG_j' \bmu  
	= \bG_j'^\top \br + \bG^\top \left( \bJ_g(\bx) \right)_{:,j} , \quad j = 1,2, \ldots, n  .
\end{aligned}
\end{equation}
From this system of equations, we can find $\bmu_j'$ and plug into \eqref{equation:sep_jaco_j}, for which we can apply the Gauss-Newton method \eqref{equation:sep_gaussnew} or the Levenberg-Marquardt method \eqref{equation:sep_levmarg} accordingly.
Therefore, the matrix $\bJ(\bz) \in \real^{m \times (n+p)}$ in the full parameter set model is replaced by the smaller matrix $\bJ(\bx) \in \real^{m \times n}$. 

Given the normal equation $\bA(\bx) \bmu(\bx) = \bb(\bx) $, it follows that $\bG(\bx)^\top (\bG(\bx) \bmu(\bx) -  \by)=\bzero $.
This implies the residual $\br(\bx)$ in \eqref{equation:sep_red_resid} is orthogonal to the column space of $\bG(\bx)$: $\bG(\bx)^\top\br(\bx)=\bzero$. This implies that, in the right-hand side of the Gauss-Newton method \eqref{equation:sep_gaussnew} or the Levenberg-Marquardt method \eqref{equation:sep_levmarg}, (omitting the superscript) the $j$-th element is
$$
\left( \bJ(\bx)^\top \br(\bx) \right)_{j} = \left( \left( \bJ(\bx) \right)_{:,j}\right)^\top \br(\bx)
= - \left( \bG_j'(\bx) \bmu(\bx) \right)^\top \br(\bx)  .
$$
This shows that the right-hand side in \eqref{equation:sep_gaussnew}-\eqref{equation:sep_levmarg} is independent of $\bmu_j'(\bx)$. To get the correct matrix on the left-hand side, however, we have to include the term $\bG(\bx) \bmu_j'(\bx)$.

Some special cases that lead to simple expressions for the Jacobian warrant attention:
\begin{itemize}
\item \textbf{$p=1$.} In this case $\bG(\bx), \bG'_j(\bx) \in \real^{m }$ are vectors and $A(\bx)$ and $\mu(\bx)$ are scalars. Additionally, 
since $\bG(\bx)^\top \bG'_j(\bx) =  \bG'_j(\bx)^\top\bG(\bx)$,
it also follows that
$$ 
A \mu_j' = \bG_j'^\top \br - \bG^\top \bG_j' \mu = \bG_j'^\top (\br - \bG \mu), \quad j = 1, 2,\ldots, n  . 
$$
Introducing the matrix $\bG' = [\bG_1',\bG_2',\ldots,  \bG_n']\in\real^{m\times n}$ and the $\real^n$ vector $\bmu'=[\mu_1', \mu_2', \ldots, \mu_n']^\top\in\real^n$, we get
$$
\bmu' = \frac{1}{A} \bG'^\top (\br - \mu \bG) \in\real^n, 
\qquad 
\bJ = - \mu \bG' - \bG (\bmu')^\top  \in\real^{m\times n}.
$$
This reduces to the nonlinear least squares problem in \eqref{equation:nonlinear_ols} in some sense.
\item \textbf{$p=n$ and $g_j(\bx, \zeta) = g_j(x_j, \zeta)$ for $j\in\{1,2,\ldots,n\}$}. That is, each of the functions $g_j$ depends only on $x_j$ and $\zeta$. In this case, all columns in $\bG_j'(\bx)\in\real^{m\times n}$ are zero except for the $j$-th column $\big( \bG_j'(\bx) \big)_{:,j}$. Therefore, the $j$-th element $\big( \bG_j'^\top \br \big)_j$ is the only nonzero element in $\bG_j'^\top \br$, and the other contribution $ \bG^\top \left( \bJ_g(\bx) \right)_{:,j} = - \bG^\top \left(\bG_j'(\bx) \bmu(\bx)\right)$ to the right-hand side in \eqref{equation:sep_redu_parAmuj} is $\bG^\top \big( \mu_j (\bG_j')_{:,j} \big)$. Let
$$
\bG' \triangleq \left[ (\bG_1')_{:,1},\;(\bG_2')_{:,2}, \;\ldots,\;  (\bG_n')_{:,n} \right] \in\real^{m\times n}, 
\qquad 
\bU' \triangleq \left[ \bmu_1', \bmu_2', \ldots, \bmu_n' \right] \in\real^{n\times n}.
$$
For \eqref{equation:sep_redu_parAmuj} and $j\in\{1,2,\ldots,n\}$, we have  
$$
\bA \bmu_j'=
\bG_j'^\top \br + \bG^\top \left( \bJ_g(\bx) \right)_{:,j}
=
\begin{bmatrix}
0, \ldots, \big( \bG_j'^\top \br \big)_j, 0, \ldots,0
\end{bmatrix}^\top 
- \bG^\top \big( \mu_j (\bG_j')_{:,j} \big)
\in\real^n.
$$
Then we can write \eqref{equation:sep_redu_parAmuj} as a matrix equation
$$
\bA \bU' = \diag(\bG^\top \br) - \bG^\top \bG'\diag(\bmu), 
$$
where $\diag(\cdot)$ is a diagonal matrix concerning the given vector.
And similarly, the Jacobian matrix in \eqref{equation:sep_jaco_j} takes the form
$$
\bJ = -\bG'\diag(\bmu) - \bG \bU'  .
$$
\end{itemize}

\section{Sparse Optimization Problem: $\ell_0$ Norm}

We have introduced several algorithms for solving the constrained optimization problem (P2) as defined in Definition~\ref{definition:opt_probs_all}, including projected gradient descent, conditional gradient method, and mirror descent method. These algorithms are applicable when the feasible set $\sS$ is closed and convex.
However, in this section, we focus on sparse optimization problems where the constraint set is non-convex. Consequently, none of the previously discussed results can be directly applied to this class of problems. Despite some similarities between convex and non-convex cases, there are significant differences in the types of results that can be obtained.

\section*{Definition and Notation}

In this section, we address the challenge of minimizing a continuously differentiable objective function subject to a sparsity constraint. Specifically, we consider the following problem:
\begin{equation}\label{equation:sparse_ellzero}
\textbf{(S0)}\qquad 
\begin{aligned}
	& \min_{\bx}
	& & f(\bx) \\
	& \text{s.t.}
	& & \normzero{\bx} \leq s,
\end{aligned}
\end{equation}
where $f : \real^n \to \real$, with its gradient being Lipschitz continuous with constant $\beta$, $s > 0$ is an integer smaller than $n$, and $\normzero{\bx}$ is the so-called $\ell_0$ norm of $\bx$, which counts the number of nonzero components in $\bx$:
$$
\normzero{\bx} = \absbig{\{i \mid x_i \neq 0, \forall i\in\{1,2,\ldots,n\}\}}.
$$

It is important to note that the $\ell_0$ is not actually a norm because it does not satisfy the homogeneity property; specifically, $\normzero{\lambda \bx} =  \normzero{\bx}$ for $\lambda \neq 0$ (Definition~\ref{definition:matrix-norm}). We do not assume that $f$ is a convex function, and the constraint set is inherently non-convex. As such, much of the analysis presented thus far does not apply to problem (S0). Nevertheless, this type of problem is crucial and has applications in areas like compressed sensing, making it valuable to study and understand how classical results over convex sets can be adapted for problem (S0).

To begin, let's define some notation. For a given vector $\bx \in \real^n$, we introduce the support set and the concept of $k$-sparse set.
\begin{definition}[Support  Set and $k$-Sparse Vectors]\label{definition:supposrt_supp_set}
The \textit{support set} of a vector $ \bx $ and its complement are  denoted, respectively, by 
$$ 
\sI_1(\bx)\triangleq\supp(\bx) \triangleq \{i \mid x_i \neq 0\}
\qquad\text{and}\qquad 
\sI_0(\bx) \triangleq \{i \mid x_i = 0\}.
$$
A vector $ \bx $ is referred to as \textit{$ k $-sparse} if $ \abs{\supp(\bx)} \leq k $, i.e., $\bx\in \sB_0[k]$ (the closed $\ell_0$ ball with a zero center: $\sB_0[k]=\sB_0[\bzero, k]$ for brevity; Definition~\ref{definition:open_closed_ball}), meaning that they have at most $k$ nonzero elements.
\end{definition}

Using this notation, problem (S0) can be restated as:
$$
\min \; f(\bx) \; \text{ s.t. }\;\bx \in \sB_0[s] .
$$
For a vector $\bx \in \real^n$ and $i \in \{1, 2, \ldots, n\}$, the \textit{$i$-th largest absolute value component} in $\bx$ is denoted by $[\bx]_i$, satisfying:
\begin{equation}
[\bx]_1 \geq [\bx]_2 \geq \ldots \geq [\bx]_n,
\end{equation}
where 
\begin{equation}
[\bx]_1 = \max_{i=1,\ldots,n} \abs{x_i}
\quad\text{and}\qquad 
[\bx]_n = \min_{i=1,\ldots,n} \abs{x_i}.
\end{equation}

Our goal now is to explore necessary optimality conditions for problem (S0), which share similarities with stationarity properties. Given the non-convexity of the constraint in (S0), traditional stationarity conditions may not hold. However, alternative optimality conditions can be formulated.

\subsection{Optimity Condition under $L$-Stationarity}

We begin by extending the concept of stationarity to problem (S0), utilizing the characterization of stationarity through the projection operator. Recall that $\bx^*$ is a stationary point of the problem of minimizing a continuously differentiable function over a closed and convex set $\sS$ if and only if
\begin{equation}\label{equation:lstat_proje}
\bx^* = \projectS(\bx^* - \eta \nabla f(\bx^*)),
\end{equation}
where $\eta$ is an arbitrary positive scalar (Theorem~\ref{theorem:stat_point_uncons_convset_proj}). 
Note that condition \eqref{equation:lstat_proje}---although expressed in terms of the parameter $\eta$---does not actually depend on $\eta$. However, this condition does not apply when $\sS = \sB_0[s]$, since the orthogonal projection onto non-convex sets like $\sB_0[s]$ is not unique; instead, it represents a multivalued mapping defined as:
$$
\project_{\sB_0[s]}(\bx) = \mathop{\argmin}_{\by} \{\normtwo{\by - \bx} \mid \by \in \sB_0[s]\}.
$$

By Weierstrass theorem (Theorem~\ref{theorem:weierstrass_them}\ref{weier2_prop_close_v3}),
the existence of vectors in $\project_{\sB_0[s]}(\bx)$ follows from the closedness of $\sB_0[s]$ and the coerciveness of the function $\normtwo{\by - \bx}$ (with respect to $\by$). To extend \eqref{equation:lstat_proje} to the sparsity constrained problem (S0), we introduce the notion of ``$L$-stationarity.''

\begin{definition}[$L$-Stationarity]\label{definition:l_stat}
A vector $\bx^* \in \sB_0[s]$ is called an $L$-stationary point of (S0) if it satisfies the relation
\begin{equation}\label{equation:l_stationary}
\textbf{($L$-stationarity)}:\qquad \bx^* \in \project_{\sB_0[s]}\left(\bx^* - \frac{1}{L} \nabla f(\bx^*)\right).
\end{equation}
\end{definition}

Given that the orthogonal projection operator $\project_{\sB_0[s]}(\cdot)$ is not single-valued, its output consists of all vectors comprising the 
$s$ components of $\bx$ with the largest absolute values, with zeros elsewhere. For example,
$$
\project_{\sB_0[2]}([3, 1, 1]^\top) = \big\{[3, 1, 0]^\top, [3, 0, 1]^\top\big\}.
$$

Below, we demonstrate that under the assumption that $f $ is $\beta$-smooth, $L$-stationarity is a necessary condition for optimality whenever $L > \beta$.
Before proving this assertion, let's present a more explicit representation of  $L$-stationarity.
\begin{lemma}[$L$-Stationarity]\label{lemma:l_stationary}
For any $L > 0$, a vector  $\bx^* \in \sB_0[s]$ satisfies $L$-stationarity if and only if
\begin{equation}\label{equation:l_stationarylem}
\abs{\frac{\partial f}{\partial x_i}(\bx^*)}
\begin{cases}
	\leq L \cdot [\bx^*]_s, & \text{if } i \in \sI_0(\bx^*), \\
	= 0, & \text{if } i \in \sI_1(\bx^*).
\end{cases}
\end{equation}
\end{lemma}
\begin{proof}[of Lemma~\ref{lemma:l_stationary}]
Suppose that $\bx^*$ satisfies $L$-stationarity. Note that for any index $j \in \{1, 2, \ldots, n\}$, the $j$-th component of any vector in $\project_{\sB_0[s]}(\bx^* - \frac{1}{L} \nabla f(\bx^*))$ is either zero or equal to $x_j^* - \frac{1}{L} \frac{\partial f}{\partial x_j}(\bx^*)$. On the other hand, since $\bx^* \in \project_{\sB_0[s]}(\bx^* - \frac{1}{L} \nabla f(\bx^*))$, it follows that if $i \in \sI_1(\bx^*)$, then $x_i^* = x_i^* - \frac{1}{L} \frac{\partial f}{\partial x_i}(\bx^*)$, so that $\frac{\partial f}{\partial x_i}(\bx^*) = 0$. If $i \in \sI_0(\bx^*)$, then $\abs{x_i^* - \frac{1}{L} \frac{\partial f}{\partial x_i}(\bx^*)} \leq [\bx^*]_s$ and $x_i^* = 0$, implying $\abs{ \frac{\partial f}{\partial x_i}(\bx^*)} \leq L \cdot [\bx^*]_s$. Therefore,  condition in \eqref{equation:l_stationarylem} holds.

Conversely, suppose that $\bx^*$ satisfies \eqref{equation:l_stationarylem}. If $\normzero{\bx^*} < s$, then $[\bx^*]_s = 0$, and by \eqref{equation:l_stationarylem} it follows that $\nabla f(\bx^*) = \bzero$; therefore, in this case, $\project_{\sB_0[s]}\left(\bx^* - \frac{1}{L} \nabla f(\bx^*)\right) = \project_{\sB_0[s]}(\bx^*)$ is the set $\{\bx^*\}$, which apparently satisfies the $L$-stationarity. 
On the other hand, if $\normzero{\bx^*} = s$, then $[\bx^*]_s \neq 0$ and $\abs{\sI_1(\bx^*)} = s$. By \eqref{equation:l_stationarylem},
$$
\abs{x_i^* - \frac{1}{L} \frac{\partial f}{\partial x_i}(\bx^*)}
\begin{cases}
\leq [\bx^*]_s, & i \in \sI_0(\bx^*), \\
= \abs{x_i^*}, & i \in \sI_1(\bx^*).
\end{cases}
$$
Therefore, the vector $\bx^* - \frac{1}{L} \nabla f(\bx^*)$ contains the $s$ components of $\bx^*$ with the largest absolute values and all other components are smaller or equal  to them in magnitude. This also implies the $L$-stationary and completes the proof.
\end{proof}

Clearly, condition \eqref{equation:l_stationarylem} depends on the constant $L$; it becomes more restrictive as $L$ decreases. 
This highlights that the situation in the non-convex case differs from the convex one. Before exploring under which conditions $L$-stationarity serves as a necessary optimality condition, we will first prove a technical yet useful result. 


\begin{lemma}[Descent Lemma for Sparse Projection under SS]\label{lemma:des_lem_lstat}
Let  $f: \real^n\rightarrow \real$  be a $\beta$-smooth function, and let $L > \beta$. Then for any $\bx \in \sB_0[s]$ and $\by \in \real^n$ satisfying
$
\by \in \project_{\sB_0[s]}\left( \bx - \frac{1}{L} \nabla f(\bx) \right),
$
we have
$$
f(\bx) - f(\by) \geq \frac{L - \beta}{2} \normtwo{\bx - \by}^2.
$$
\end{lemma}
\begin{proof}
Since $\by \in \project_{\sB_0[s]}\left( \bx - \frac{1}{L} \nabla f(\bx) \right)$, we have 
$$
\small
\begin{aligned}
\by 
&\in \mathop{\argmin}_{\bu \in \sB_0[s]} \normtwo{\bu - \left( \bx - \frac{1}{L} \nabla f(\bx) \right)}^2\\
&=\mathop{\argmin}_{\bu \in \sB_0[s]}  \frac{L}{2} \normtwo{\bu - \left( \bx - \frac{1}{L} \nabla f(\bx) \right)}^2  
+\underbrace{ f(\bx) -\frac{1}{2L} \normtwo{\nabla f(\bx)}^2}_{\text{constant w.r.t. } \bu}\\
&=\mathop{\argmin}_{\bu \in \sB_0[s]} f(\bx) + \langle \nabla f(\bx), \bu - \bx \rangle + \frac{L}{2} \normtwo{\bu - \bx}^2
\end{aligned}
$$
Let 
\begin{equation}\label{equation:des_lem_lstatpsi}
\Psi_L(\bu,\bx)\triangleq f(\bx) + \langle \nabla f(\bx), \bu - \bx \rangle + \frac{L}{2} \normtwo{\bu - \bx}^2.
\end{equation}
This implies that
$$
\Psi_L(\by, \bx) \leq \Psi_L(\bx, \bx) = f(\bx).
$$
By the $\beta$-smoothness of $f$ (Definition~\ref{definition:scss_func}), we have 
$$
f(\bx) - f(\by) \geq f(\bx) - \Psi_{\beta}(\by, \bx).
$$
Combining the preceding properties and the fact that $
\Psi_{\beta}(\by, \bx) = \Psi_L(\by, \bx) - \frac{L - \beta}{2} \normtwo{\bx - \by}^2
$ yields the desired result.
\end{proof}

We are now ready to demonstrate the necessity of  $L$-stationarity under the condition $L > \beta$ for $\beta$-smooth functions.

\begin{theoremHigh}[Necessity under $L$-Stationarity and Smoothness]\label{theorem:necess_lstat}
Let  $f: \real^n\rightarrow \real$ be a $\beta$-smooth function, and let $L > \beta$. Let $\bx^*$ be an optimal solution of (S0). Then
\begin{enumerate}[(i)]
\item $\bx^*$ is an $L$-stationary point.
\item The set $\project_{\sB_0[s]}\left(\bx^* - \frac{1}{L} \nabla f(\bx^*)\right)$ is a singleton.
\end{enumerate}
\end{theoremHigh}
\begin{proof}[of Theorem~\ref{theorem:necess_lstat}]
Assume, for contradiction, the set is not a singleton; that is, that there exists a vector
$
\by \in \project_{\sB_0[s]}\left(\bx^* - \frac{1}{L} \nabla f(\bx^*)\right),
$
which is different from $\bx^*$ ($\by \neq \bx^*$). By Lemma~\ref{lemma:des_lem_lstat}, we have 
$
f(\bx^*) - f(\by) \geq \frac{L - \beta}{2} \normtwo{\bx^* - \by}^2,
$
contradicting the optimality of $\bx^*$. 
Thus, $\bx^*$ is the only vector in the set $\project_{\sB_0[s]}\left(\bx^* - \frac{1}{L} \nabla f(\bx^*)\right)$, and $\bx^*$ is $L$-stationary.
\end{proof}

In summary, we have demonstrated that if the function is $\beta$-smooth, $L$-stationarity for any $L > \beta$ is a necessary optimality condition for problem (S0).

\subsection{Restricted Isometry and Other Design Properties}

In many applications, sparse optimization problems involve minimizing the function $f(\bx) = \frac{1}{2m}\normtwo{\bb - \bA \bx}$, where $\bA\in\real^{m\times n}$ is referred to as the \textit{design matrix} and $\bb\in\real^m$.
Apparently, when the design matrix is an identity,  it preserves the geometric properties of signals or models, enabling the recovery of sparse signals. 
We will now further formalize this concept and establish specific conditions on $\bA$ that guarantee universal recovery, meaning that for every $\bx \in \sB_0[s]$, $\bx$ can be uniquely recovered from the measurements $\bA\bx$.

To be more specific, a design matrix capable of identifying sparse vectors does not necessarily ensure universal recovery. For instance, consider a design matrix $\bA \in \real^{m\times n}$ such that for some distinct $\bx_1, \bx_2 \in \sB_0[s]$ and $\bx_1 \neq \bx_2$, we have $\by_{1} = \by_{2}$, where $\by_{1} = \bA \bx_1$ and $\by_{2} = \bA \bx_2$. 
In this scenario, it becomes theoretically impossible to distinguish between $\bx_1$ and $\bx_2$ based on measurements made using $\bA$, i.e., using $\by_{1}$ (or $\by_{2}$). 
Consequently, this design matrix cannot support universal recovery because it fails to differentiate between sparse vectors. It's important to note that such a matrix would fail to uniquely identify not just one pair but an infinite set of pairs---or even an entire subspace---of sparse vectors.

Therefore, it is crucial that the design matrix maintains the geometry of sparse vectors when projecting them from a high-dimensional space ($n$ dimensions) to a lower-dimensional space ($m$ dimensions), especially in scenarios where $m \ll n$. The \textit{nullspace property}, among others, encapsulates this requirement. 
To see this, for any subset of coordinates $\sS \subseteq \{1, 2,\ldots, n\}$,  define the sets:
\begin{align}
	\sC[\sS] &\triangleq \{\bx \in \real^{n} \mid  \normone{\bx_{\comple{\sS}}} \leq \normone{\bx_{\sS}} \};\\
	\sC[\sS; \gamma] &\triangleq \{\bx \in \real^{n} \mid  \normone{\bx_{\comple{\sS}}} \leq \gamma\normone{\bx_{\sS}} \};\\
	\sC[k] &\triangleq \bigcup_{\sS: \abs{\sS} = k} \sC[\sS],
\end{align}
where $\comple{\sS}$ is the complement of $\sS$ (Definition~\ref{definition:matlabnotation}).
Here, square brackets are used instead of parentheses to indicate inclusion of equalities in the definitions. 
Thus, $\sC[\sS]=\sC[\sS; 1]$ represents the convex set of points placing most of their weight on coordinates within the set $\sS$.
Meanwhile, $\sC[k]$ constitutes the non-convex set of points concentrating most of their weight on some $k$ coordinates.
Note that $ \sB_0[k] \subset \sC[k] $~\footnote{Note that $\sB_0[k]=\sB[\bzero,  k]$ denotes the closed ball using the $\ell_0$ norm.} since $k$-sparse vectors assign all their weight to $k$ coordinates, and that $\sC[\sS; \gamma_1] \subseteq \sC[\sS; \gamma_2]$ if $\gamma_1\leq \gamma_2$.

We now introduce certain desirable properties of the design matrix, which are extensively covered in literature, e.g.,  \citet{candes2005decoding, cohen2009compressed, raskutti2010restricted, jain2014iterative, hastle2015statistical, jain2017non}.
\begin{definition}[Nullspace Property]\label{definition:nullspace_prop}
A matrix $\bA \in \real^{m\times n}$ is said to satisfy the \textit{nullspace property of order $k$} if $\nspace(\bA) \cap \sC[k] = \{\bzero\}$, where $\nspace(\bA) \triangleq \{\bx \in \real^{n}: \bA \bx = \bzero\}$ denotes the nullspace of $\bA$ (Definition~\ref{definition:nullspace}).
\end{definition}

If a design matrix satisfies this property, vectors in its nullspace cannot concentrate most of their weight on any $k$ coordinates. Consequently, no $k$-sparse vector can reside within the nullspace. If a design matrix has the nullspace property of order $2s$, it ensures that two distinct  $s$-sparse vectors cannot be identified, which is crucial for ensuring global recovery. An enhanced version of the nullspace property leads to the restricted eigenvalue property.

\begin{definition}[Restricted Eigenvalue Property]\label{definition:res_eig}
A matrix $\bA \in \real^{m\times n}$ is said to satisfy the \textit{restricted eigenvalue (RE) property with constant $\alpha$} over $\sS$ if for all $\bx \in \sS$, we have 
$$
\alpha\leq \frac{\frac{1}{m} \bx^\top\bA^\top\bA\bx}{\normtwo{\bx}^{2}}
\qquad\iff\qquad 
\alpha \cdot \normtwo{\bx}^{2} \leq \frac{1}{m} \normtwo{\bA\bx}^{2},
$$
where the right-hand side of the first inequality represents  the Rayleigh quotient of the vector $\bx$ associated with the matrix $\frac{1}{m} \bA^\top\bA$, and $\frac{1}{m} \bA^\top\bA$ is the Hessian of the function $f(\bx) = \frac{1}{2m}\normtwo{\bA\bx-\bb}^2$.
\end{definition}

This property implies that not only are $k$-sparse vectors absent from the nullspace, they actually retain a significant  fraction of their length after projection as well. This means that if $k=2s$, then 
\begin{equation}
\bx_1, \bx_2 \in \sB_0[s] 
\quad\implies\quad 
\frac{1}{m} \normtwo{\bA(\bx_1 - \bx_2)}^{2} \geq \alpha \cdot \normtwo{\bx_1 - \bx_2}^{2}
\end{equation}
Thus, the distance between any two sparse vectors does not significantly decrease after projection.

The next property further explicates this concept and is known as the restricted isometry property.

\begin{definition}[Restricted Isometry Property]\label{definition:rip}
A matrix $\bA \in \real^{m\times n}$ is said to satisfy the \textit{restricted isometry property (RIP) of order $k$ with constant $\delta_{k} \in [0,1)$} if for all $\bx \in \sB_0[k] $, we have
$$
(1-\delta_{k}) \cdot \normtwo{\bx}^{2} \leq \frac{1}{m} \normtwo{\bA\bx}^{2} \leq (1+\delta_{k}) \cdot \normtwo{\bx}^{2}.
$$
\end{definition}

This property is widely used in analyzing sparse recovery and compressive sensing algorithms. However, it can be somewhat restrictive since it requires distortion parameters to be of the form $(1 \pm \delta)$ for $\delta \in [0,1)$.

\begin{exercise}[Subset of Restricted Isometry Property]\label{exercise:sub_rip}
Suppose the matrix $ \bA \in \real^{m\times n} $ satisfies RIP at order $ s $ with constant $ \delta_s $ over $\sB_0[s]$. Show that for any set $ \sI \subset \{1,2,\ldots,n\}, \abs{\sI} \leq s $, the smallest eigenvalue of the matrix $ \bA_{\sI}^\top \bA_{\sI} $ is lower bounded by $ (1 - \delta_s) $, where $\bA_{\sI}\triangleq\bA[:,\sI]$.
\end{exercise}

%
%

A useful generalization of this property, particularly relevant in scenarios where the properties of the design matrix are not strictly controlled---such as in gene expression analysis---is the concept of restricted strong convexity and smoothness.

\begin{definition}[Restricted Strong Convexity/Smoothness Property]\label{definition:res_scss_mat}
A matrix $\bA \in \real^{m\times n}$ is said to satisfy the \textit{$\alpha$-restricted strong convexity (RSC) property} and the \textit{$\beta$-restricted smoothness (RSS) property of order $k$} if for all $\bx \in \sB_0[k]$, we have
$$
\alpha \cdot \normtwo{\bx}^{2} \leq \frac{1}{m} \normtwo{\bA\bx}^{2} \leq \beta \cdot \normtwo{\bx}^{2}.
$$
\end{definition}

The key difference between RIP and RSC/RSS properties lies in their constants. RIP requires these constants to be of the form $1 \pm \delta_{k}$, whereas RSC/RSS does not impose such constraints. Readers will notice similarities between the definitions of restricted strong convexity and smoothness provided here and Definition~\ref{definition:res_scss_func}, where we defined these concepts for general functions. It's encouraged that readers verify the relationship between these two definitions. Indeed, Definition~\ref{definition:res_scss_func} can be viewed as a generalization of Definition~\ref{definition:res_scss_mat} to general functions \citep{jain2014iterative}. For twice-differentiable functions, both definitions essentially impose restrictions on the (restricted) eigenvalues of the function's Hessian.

\subsection{The Iterative Hard-Thresholding (IHT) Method}
\begin{algorithm}[h] 
\caption{Iterative Hard-Thresholding (IHT) Method}
\label{alg:pgd_iht}
\begin{algorithmic}[1] 
\Require Sparsity level $k$; 
\State {\bfseries Input (a).} A differential function $f$; 
\State {\bfseries Input (b).} $f(\bx) = \frac{1}{2m}\normtwo{\bb - \bA \bx}$: design matrix $\bA\in\real^{m\times n}$ and  observation vector $\bb\in\real^m$;
\State {\bfseries Input:}  Initialize $\bx^{(1)} \in \sB_0[k]$;
\For{$t=1,2,\ldots$}
\State Pick a stepsize $\eta_t$;
\State $\by^{(t+1)} \stackrel{(a)}{\leftarrow} \bx^{(t)} - \eta_t \nabla f(\bx^{(t)}) \quad \text{ or }\quad \stackrel{(b)}{\leftarrow} \bx^{(t)} - \eta_t  \frac{1}{m} \bA^\top (\bA\bx^\toptzero - \bb)$;
\State $\bx^{(t+1)} \in \mathcalP_{\sB_0[k]}(\by^{(t+1)})$; \Comment{Selecting the $k$ largest elements in magnitude of $\by^{(t+1)}$}
\EndFor
\State {\bfseries Return:}  final $\bx\leftarrow \bx^{(t)}$;
\end{algorithmic} 
\end{algorithm}

One approach to solving problem (S0) involves using a natural generalization of the projected gradient descent algorithm, which can be described as follows:
\begin{subequations}
\begin{equation}
\bx^\toptone \in \project_{\sB_0[s]}\left(\bx^\toptzero - \frac{1}{L} \nabla f(\bx^\toptzero)\right), \quad t = 1, 2, \ldots.
\end{equation}
The method is known in the literature as the \textit{iterative hard-thresholding (IHT) method} (Algorithm~\ref{alg:pgd_iht}), and we will adopt this terminology.
Note that the algorithm can be  applied  both to a general function $f$ (with properties outlined in Lemma~\ref{lemma:conv_iht_smoot})  and a least squares objective function $f(\bx) = \frac{1}{2m}\normtwo{\bb - \bA \bx}$ (which satisfies certain design properties as previously discussed).
It has been shown that the general step of the IHT method is equivalent to the relation
\begin{equation}
\bx^\toptone \in \mathop{\argmin}_{\bx \in \sB_0[s]} \Psi_L(\bx, \bx^\toptzero), \quad t = 1, 2, \ldots,
\end{equation}
\end{subequations}
where $\Psi_L(\bx, \by)$ is defined by \eqref{equation:des_lem_lstatpsi}.
When the projected set $\sS$ is closed convex, the projected gradient descent method (Algorithm~\ref{alg:pgd_gen}) can be considered  a fixed-point method for solving $
\bx^* =\mathcalP_{\sS}\big(\bx^* - \eta\nabla f(\bx^*)\big)
$ (Theorem~\ref{theorem:stat_point_uncons_convset_proj}).
Therefore, the IHT method can also be viewed as a fixed-point method aimed at ``enforcing'' the $L$-stationary condition (Definition~\ref{definition:l_stat}, Theorem~\ref{theorem:necess_lstat}) and solving $\bx^* =\mathcalP_{\sS}\big(\bx^* - \frac{1}{L}\nabla f(\bx^*)\big)$ for $L>\beta$,  assuming the underlying function is $\beta$-smooth.

Several fundamental properties of the IHT method are summarized in the following lemma:
\begin{lemma}[Properties of IHT under Smoothness]\label{lemma:conv_iht_smoot}
Let  $f: \real^n\rightarrow \real$ be a $\beta$-smooth and lower bounded function.
Suppose $\{\bx^\toptzero\}_{t > 0}$ is the sequence generated by the IHT method (Algorithm~\ref{alg:pgd_iht}(a)) with a constant stepsize $\eta\triangleq\eta_t\triangleq\frac{1}{L}$, where $L > \beta$, and a projection sparsity level $k=s$. Then,
\begin{enumerate}[(i)]
\item $f(\bx^\toptzero) - f(\bx^\toptone) \geq \frac{L - \beta}{2} \normtwo{\bx^\toptzero - \bx^\toptone}^2$,
\item $\{f(\bx^\toptzero)\}_{t > 0}$ is a nonincreasing sequence,
\item $\normtwo{\bx^\toptzero - \bx^\toptone} \to 0$,
\item For every $t = 1, 2, \ldots$, if $\bx^\toptzero \neq \bx^\toptone$, then $f(\bx^\toptone) < f(\bx^\toptzero)$.
\end{enumerate}
\end{lemma}
\begin{proof}[of Lemma~\ref{lemma:conv_iht_smoot}]
Part (i) follows from Lemma~\ref{lemma:des_lem_lstat} by substituting $\bx = \bx^\toptzero, \by = \bx^\toptone$. Part (ii) follows immediately from part (i). To prove (iii), note that since $\{f(\bx^\toptzero)\}_{t > 0}$ is a nonincreasing sequence, which is also lower bounded, it follows that it converges to some limit $B$ and hence $f(\bx^\toptzero)) - f(\bx^\toptone) \to B - B = 0$ as $t \to \infty$. Therefore, by part (i) and the fact that $L > \beta$, the limit $\normtwo{\bx^\toptzero - \bx^\toptone} \to 0$ holds. Finally, (iv) is a direct consequence of (i).
\end{proof}

As previously mentioned, the IHT algorithm can be considered a fixed-point method for solving the $L$-stationarity condition. 
The following theorem states that all accumulation points of the sequence generated by the IHT method with a constant  stepsize $\frac{1}{L}$ are indeed $L$-stationary points.
\begin{theoremHigh}[Convergence of IHT \citep{beck2014introduction}]\label{theorem:acc_iht_conv}
Consider the same conditions  as Lemma~\ref{lemma:conv_iht_smoot}.
Let $\{\bx^\toptzero\}_{t > 0}$ be the sequence generated by the IHT method (Algorithm~\ref{alg:pgd_iht}(a)) with a constant stepsize $\eta\triangleq\eta_t \triangleq \frac{1}{L}$, where $L > \beta$, and a projection sparsity level $k=s$. Then, any accumulation point of $\{\bx^\toptzero\}_{t > 0}$ is an $L$-stationary point.
\end{theoremHigh}
\begin{proof}[of Theorem~\ref{theorem:acc_iht_conv}]
Suppose that $\bx^*$ is an accumulation point of the sequence. 
Thus, there exists a subsequence $\{\bx^{(t_j)}\}_{j > 0}$ that converges to $\bx^*$. By Lemma~\ref{lemma:conv_iht_smoot}, we have 
\begin{equation}\label{equation:acc_iht_conv1}
f(\bx^{(t_j)}) - f(\bx^{(t_j + 1)}) \geq \frac{L - \beta}{2} \normtwo{\bx^{(t_j)} - \bx^{(t_j + 1)}}^2.
\end{equation}
Since $\{f(\bx^{(t_j)})\}_{j > 0}$ and $\{f(\bx^{(t_j + 1)})\}_{j > 0}$, as nonincreasing and lower bounded sequences,
they both converge to the same limit.
Consequently, $f(\bx^{(t_j)}) - f(\bx^{(t_j + 1)}) \to 0$ as $j \to \infty$, which combined with \eqref{equation:acc_iht_conv1} yields that
$
\bx^{(t_j + 1)} \to \bx^* \text{ as } j \to \infty.
$
Recall that for all $j > 0$
$
\bx^{(t_j + 1)} \in \project_{\sB_0[s]}\left(\bx^{(t_j)} - \frac{1}{L} \nabla f(\bx^{(t_j)})\right).
$
We then consider the following two cases.

Let $i \in \sI_1(\bx^*)$ (i.e., in the suppose set; Definition~\ref{definition:supposrt_supp_set}). By the convergence of $\bx^{(t_j)}$ and $\bx^{(t_j + 1)}$ to $\bx^*$, it follows that there exists an integer $J$ such that
$$
x_i^{(t_j)}, x_i^{(t_j+1)} \neq 0 \text{ for all } j > J,
$$
and therefore, for $j > J$,
$$
x_i^{(t_j+1)} = x_i^{(t_j)} - \frac{1}{L} \frac{\partial f}{\partial x_i}(\bx^{(t_j)}).
$$
Taking the limit as  $j\rightarrow \infty$ and using the continuity of $f$, we obtain that
$
\frac{\partial f}{\partial x_i}(\bx^*) = 0.
$

Now let $i \in \sI_0(\bx^*)$. If there exist an infinite number of indices $t_j$ for which $x_i^{(t_j+1)} \neq 0$, then as in the previous case, we obtain that $x_i^{(t_j+1)} = x_i^{(t_j)} - \frac{1}{L} \frac{\partial f}{\partial x_i}(\bx^{(t_j)})$ for these indices, implying (by taking the limit) that $\frac{\partial f}{\partial x_i}(\bx^*) = 0$. In particular, $\left| \frac{\partial f}{\partial x_i}(\bx^*) \right| \leq L [\bx^*]_s$. On the other hand, if there exists an $M > 0$ such that for all $j > M$, $x_i^{(t_j+1)} = 0$, then
$$
\left| x_i^{(t_j)} - \frac{1}{L} \frac{\partial f}{\partial x_i}(\bx^{(t_j)}) \right| \leq [\bx^{(t_j)} - \frac{1}{L} \nabla f(\bx^{(t_j)}) ]_s = [\bx^{(t_j + 1)}]_s.
$$

Thus, taking $j$ to infinity, while exploiting the continuity of the function $[\cdot]_s$, we obtain that
$$
\left| \frac{\partial f}{\partial x_i}(\bx^*) \right| \leq L [\bx^*]_s,
$$
and hence, by Lemma~\ref{lemma:l_stationary}, the desired result is established.
\end{proof}

We will now establish an important result: if the design matrix satisfies the RIP condition with suitable constants, then the IHT algorithm guarantees universal sparse recovery for the least squares problem: $f(\bx) = \frac{1}{2m}\normtwo{\bb - \bA \bx}$.
\begin{theoremHigh}[Rate of Convergence of IHT \citep{jain2014iterative, jain2017non}]\label{theorem:conv_iht_lasso}

Suppose $\bA \in \real^{m \times n}$ is a design matrix that satisfies the RIP property of order $3s$ with parameter $\delta_{3s} < \frac{1}{2}$. Let $\bx^* \in \sB_0[s] \subset \real^n$ be any arbitrary sparse vector and let $\bb = \bA \bx^*$. Then the IHT algorithm (Algorithm~\ref{alg:pgd_iht}(b)), when executed with a constant stepsize $\eta_t = 1$ and a projection sparsity level $k = s$, ensures $\normtwo{\bx^\toptzero - \bx^*} \leq \epsilon$ after at most $t = \mathcalO\left(\ln \frac{\normtwo{\bx^*}}{\epsilon}\right)$ iterations of the algorithm.
\end{theoremHigh}
\begin{proof}[of Theorem~\ref{theorem:conv_iht_lasso}]
Let $\sS^* \triangleq \text{supp}(\bx^*)$ and $\sS^t \triangleq \text{supp}(\bx^\toptzero)$. 
Define $\sI^t \triangleq \sS^t \cup \sS^{t+1} \cup \sS^*$ as the union of the supports of  two consecutive iterates and the optimal model. 
This definition ensures that while analyzing this update step, the error vectors $\bx^\toptzero - \bx^*$ and $\bx^\toptone - \bx^*$, which will be the focal point of the analysis, have support within $\sI^t$. Note that $\abs{\sI^t} \leq 3s$.

With $\eta_t = 1$,  $\by^\toptone = \bx^\toptzero - \frac{1}{m} \bA^\top (\bA \bx^\toptzero - \bb)$ by Algorithm~\ref{alg:pgd_iht}. 
Since the set $\sB_0[k]$ is non-convex, we can only apply the Projection Property-O (Lemma~\ref{lemma:proj_prop0}) for the projection step $\bx^\toptone = \mathcalP_{\sB_0[k]}(\by^\toptone)$. Along with the fact that $\normtwo{\bu}^2 = \normtwo{\bu_{\sI}}^2 + \normtwo{\bu_{\comple{\sI}}}^2$ for any vector $\bu\in\real^n$ (where again $\comple{\sI}$ denotes the complement of $\sI$), we have
$$
\begin{aligned}
&\normtwo{\bx^\toptone - \by^\toptone}^2 \leq \normtwo{\bx^* - \by^\toptone}^2\\
&\implies \quad \normtwo{\bx_{\sI}^\toptone - \by_{\sI}^\toptone}^2 + \normtwo{\bx_{\comple{\sI}}^\toptone - \by_{\comple{\sI}}^\toptone}^2 \leq \normtwo{\bx_{\sI}^* - \by_{\sI}^\toptone}^2 + \normtwo{\bx_{\comple{\sI}}^* - \by_{\comple{\sI}}^\toptone}^2,
\end{aligned}
$$
where, for brevity, we denote $\sI \triangleq\sI^t$ such that $\bx_{\comple{\sI}}^\toptone = \bx_{\comple{\sI}}^* = \bzero$. 
Using the fact that $\bb = \bA \bx^*$, and denoting $\widetildebA \triangleq \frac{1}{\sqrt{m}} \bA$ and $\widetildebA_{\sI} \triangleq\widetildebA[:,\sI]$, we have
$$
\begin{aligned}
&\normtwo{\bx_{\sI}^\toptone - \by_{\sI}^\toptone} \leq \normtwo{\bx_{\sI}^* - \by_{\sI}^\toptone}\\
&\Rightarrow
\normtwo{\bx_{\sI}^\toptone - \bx_{\sI}^* + \bx_{\sI}^* - \left( \bx_{\sI}^\toptzero - \widetildebA_{\sI}^\top \widetildebA (\bx^\toptzero - \bx^*) \right)} 
\leq 
\normtwo{\bx_{\sI}^* - \left( \bx_{\sI}^\toptzero - \widetildebA_{\sI}^\top \widetildebA (\bx^\toptzero - \bx^*) \right)}\\
&\Rightarrow 
\normtwo{\bx_{\sI}^\toptone - \bx_{\sI}^*} 
\leq 
2 \normtwo{\left( \bx_{\sI}^\toptzero - \bx_{\sI}^* \right) - \widetildebA_{\sI}^\top \widetildebA (\bx^\toptzero - \bx^*)},
\end{aligned}
$$
where the last inequality follows from the triangle inequality for any norm: $\abs{\norm{\bu}-\norm{\bv}} \leq \norm{\bu-\bv}$ for any $\bu$ and $\bv$, and $\by_{\sI}^\toptone = \left(\bx^{(t)} -  \widetildebA^\top\widetildebA (\bx^\toptzero - \bx^*)\right)_{\sI} = \bx^{(t)}_{\sI} -  \widetildebA_{\sI}^\top\widetildebA (\bx^\toptzero - \bx^*)$. 
Given that $\bx_{\comple{\sI}}^\toptzero = \bx_{\comple{\sI}}^* = \bzero$, it holds that 
$$
\widetildebA_{\comple{\sI}} (\bx^\toptzero_{\comple{\sI}} - \bx^*_{\comple{\sI}}) = \bzero
\,\,\, \implies\,\,\, 
\widetildebA (\bx^\toptzero - \bx^*) = 
\widetildebA_{\comple{\sI}} (\bx^\toptzero_{\comple{\sI}} - \bx^*_{\comple{\sI}}) +  \widetildebA_{{\sI}} (\bx^\toptzero_{{\sI}} - \bx^*_{{\sI}})
= \widetildebA_{\sI} (\bx^\toptzero_{\sI} - \bx^*_{\sI}),
$$
whence we have 
\begin{equation}
\begin{aligned}
\normtwo{\bx^\toptone - \bx^* } 
&\leq 2 \normtwo{(\bI - \widetildebA_{\sI}^\top \widetildebA_{\sI}) (\bx_{\sI}^\toptzero - \bx_{\sI}^*)}\\
&\leq 2 \left( \normtwo{\bx_{\sI}^\toptzero - \bx_{\sI}^*} - \normtwo{\widetildebA_{\sI}^\top \widetildebA_{\sI} (\bx_{\sI}^\toptzero - \bx_{\sI}^*)} \right)
\leq 2 \delta_{3s} \normtwo{\bx_{\sI}^\toptzero - \bx_{\sI}^*}.
\end{aligned}
\end{equation}
This completes the proof.
\end{proof}

For further proving the rate of convergence for a general function, we need the following lemma.

\begin{lemma}[Sparsity Projection Property]\label{lemma:spar_proj}
Let $\bx \in \real^n$ and  $\widetildebx \triangleq \mathcal\project_{\sB_0[s]}(\bx)$. Then for any $\bx^* \in \sB_0[s^*]\subseteq \real^n$ such that $\normzero{\bx^*} \leq s^*$, where $s^*\leq s\leq n$, we have
$$
\frac{\normtwo{\widetildebx - \bx}^2}{n-s} \leq  \frac{\normtwo{\bx^* - \bx}^2}{n-s^*}
\qquad\implies\qquad 
\normtwo{\widetildebx - \bx}^2 \leq \frac{n - s}{n - s^*} \normtwo{\bx^* - \bx}^2.
$$
\end{lemma}
\begin{proof}[of Lemma~\ref{lemma:spar_proj}]
Without loss of generality, assume that we have reordered coordinates such that $\abs{x_1} \geq \abs{x_2} \geq \ldots \geq \abs{x_n}$. Since the projection operator $\mathcal\project_{\sB_0[s]}(\cdot)$ operates by selecting the largest elements by magnitude, we have $\widetildex_1 = x_1, \widetildex_2 = x_2 \ldots, \widetildex_s = x_s$ and $\widetildex_{s+1} = \widetildex_{s+2} = \ldots = \widetildex_{n} = 0$.

Let $\widehatbx \triangleq \mathcalP_{\sB_0[s^*]}(\bx)$. By the above argument, we have $\widehatx_1 = x_1, \widehatx_2 = x_2, \ldots, \widehatx_{s^*} = x_{s^*}$ and $\widehatx_{s^*+1} = \widehatx_{s^*+2} = \ldots = \widehatx_{n} = 0$. 
Therefore, assuming $s^*\leq s$, we have 
\begin{equation}
\begin{aligned}
(n-s)&(n-s^*)\left(\frac{\normtwo{\widehatbx - \bx}^2}{n - s^*} - \frac{\normtwo{\widetildebx - \bx}^2}{n - s} \right)
=(n-s) \sum_{i=s^*+1}^{n} x_i^2 - (n-s^*)\sum_{i=s+1}^{n} x_i^2\\
&=(n-s) \sum_{i=s^*+1}^{\textcolor{mylightbluetext}{s}} x_i^2 + \big((n-s)-(n-s^*)\big)\sum_{i=s+1}^{n} x_i^2\\
&\geq (n-s)(s-s^*)  x_{s}^2 +  (s^*-s)(n-s) x_{s+1}^2 \geq 0.\\
\end{aligned}
\end{equation}
Given the Projection Property-O (Lemma~\ref{lemma:proj_prop0}), it is valid that $ \normtwo{\widehatbx - \bx} \leq \normtwo{\bx^* - \bx} $ for any $\bx^* \in \sB_0[s^*]$, which establishes the desired result when combining the above inequality.
\end{proof}

For a more general function $f$ with an optimal $s^*$-sparse parameter $\bx^*$, other than $f(\bx) = \frac{1}{2m}\normtwo{\bb - \bA \bx}$,
our analysis combines the above observation with the RSC/RSS properties of $f$ to provide geometric convergence rates for the IHT procedure below.

\begin{theoremHigh}[Rate of Convergence of IHT \citep{jain2014iterative}]\label{theorem:conv_iht_genfunc}
Let $f$ have RSC and RSS parameters given by $\alpha_{2s+s^*}(f) = \alpha$ and $\beta_{2s+s^*}(f) = \beta$, respectively. Let Algorithm~\ref{alg:pgd_iht} be invoked with $f$, $k\triangleq s \geq 32 \left(\frac{\beta}{\alpha}\right)^2 s^*$, and a constant stepsize $\eta_t = \frac{2}{3\beta}$. Also let $\bx^* = \argmin_{\normzero{\bx} \leq s^*} f(\bx)$. Then, the $T$-th iterate of Algorithm~\ref{alg:pgd_iht}(a), for $T = \mathcalO\left(\frac{\beta}{\alpha} \cdot \ln\left(\frac{f(\bx^\topone)}{\epsilon}\right)\right)$ satisfies:
$ f(\bx^{(T)}) - f(\bx^*) \leq \epsilon.
$
\end{theoremHigh}

\begin{proof}[of Theorem~\ref{theorem:conv_iht_genfunc}]
Denote $\bg^\toptzero \triangleq \nabla f(\bx^\toptzero)$, recall that $\bx^\toptone = \mathcal\project_{\sB_0[s]}(\bx^\toptzero - \frac{\gamma}{\beta} \bg^\toptzero)$ where $\gamma = \frac{2}{3} < 1$. Let $\sS^t \triangleq \supp(\bx^\toptzero)$, $\sS^* \triangleq \supp(\bx^*)$, and $\sS^{t+1} \triangleq \supp(\bx^\toptone)$. Also, let $\sI \triangleq \sS^* \cup \sS^t \cup \sS^{t+1}$.
Given the RSS property (Definition~\ref{definition:res_scss_func}) and the fact that $\supp(\bx^\toptzero) \subseteq \sI$ and $\supp(\bx^\toptone) \subseteq \sI$ (implying that $\bx_{\comple{\sI}}^\toptone = \bx_{\comple{\sI}}^\toptzero = \bzero$), along with the fact that $\normtwo{\bu}^2 = \normtwo{\bu_{\sI}}^2 + \normtwo{\bu_{\comple{\sI}}}^2$ for any vector $\bu\in\real^n$, we have:
\begin{equation}\label{equation:conv_iht_genfunc1}
\begin{aligned}
&f(\bx^\toptone) - f(\bx^\toptzero) \leq \innerproduct{\bg^\toptzero, \bx^\toptone - \bx^\toptzero} + \frac{\beta}{2} \normtwo{\bx^\toptone - \bx^\toptzero}^2, \\
&= \frac{\beta}{2} \normtwo{\bx^\toptone_{\sI} - \bx^\toptzero_{\sI} + \frac{\gamma}{\beta} \cdot  \bg^\toptzero_{\sI}}^2 - \frac{\gamma^2}{2\beta} \normtwo{\bg^\toptzero_{\sI}}^2 + (1 - \gamma) \innerproduct{\bx^\toptone - \bx^\toptzero, \bg^\toptzero}.
\end{aligned}
\end{equation}
Since $\sS^t \setminus \sS^{t+1}$ and $\sS^{t+1}$ are disjoint, we have:
\begin{equation}\label{equation:conv_iht_genfunc2}
\small
\begin{aligned}
\innerproduct{\bx^\toptone - \bx^\toptzero, \bg^\toptzero} 
&= \innerproduct{\,\cancel{\bx^\toptone_{\sS^t \setminus \sS^{t+1}}} - \bx^\toptzero_{\sS^t \setminus \sS^{t+1}}, \bg^\toptzero_{\sS^t \setminus \sS^{t+1}}} 
+ \innerproduct{\bx^\toptone_{\sS^{t+1}} - \bx^\toptzero_{\sS^{t+1}}, \bg^\toptzero_{\sS^{t+1}}} \\
&\overset{\dag}{=} - \innerproduct{\bx^\toptzero_{\sS^t \setminus \sS^{t+1}}, \bg^\toptzero_{\sS^t \setminus \sS^{t+1}}} - \frac{\gamma}{\beta} \normtwo{\bg^\toptzero_{\sS^{t+1}}}^2 \\
&\overset{\ddag}{\leq} \frac{\gamma}{2\beta} \normtwo{\bg^\toptzero_{\sS^{t+1} \setminus \sS^t}}^2 - \frac{\gamma}{2\beta} \normtwo{\bg^\toptzero_{\sS^t \setminus \sS^{t+1}}}^2 - \frac{\gamma}{\beta} \normtwo{\bg^\toptzero_{\sS^{t+1}}}^2   \cancel{-\left( \frac{\beta}{2\gamma} \normtwo{\bx_{\sS^t \setminus \sS^{t+1}}}^2\right)} \\
&\overset{*}{=} - \frac{\gamma}{2\beta} \normtwo{\bg^\toptzero_{\sS^{t+1} \setminus \sS^t}}^2 - \frac{\gamma}{2\beta} \normtwo{\bg^\toptzero_{\sS^t \setminus \sS^{t+1}}}^2 - \frac{\gamma}{\beta} \normtwo{\bg^\toptzero_{\sS^t \cap \sS^{t+1}}}^2 \\
&\leq - \frac{\gamma}{2\beta} \normtwo{\bg^\toptzero_{\sS^t \cup \sS^{t+1}}}^2,
\end{aligned}
\end{equation}
where the equality $(\dag)$ follows from the gradient step, i.e., $\bx^\toptone_{\sS^{t+1}} = \bx^\toptzero_{\sS^{t+1}} - \frac{\gamma}{\beta} \bg^\toptzero_{\sS^{t+1}}$, the inequality $(\ddag)$ follows using the fact that $\bx^\toptone$ is obtained using hard-thresholding ($\sS^{t+1}$ contains the indices for largest $\abs{\sS^{t+1}}=\abs{\sS^{t}}$ magnitudes) and the fact that $\abs{\sS^t \setminus \sS^{t+1}} = \abs{\sS^{t+1} \setminus \sS^t}$, such that 
$
\normtwo{\bx_{\sS^t \setminus \sS^{t+1}}^\toptzero - \frac{\gamma}{\beta} \bg_{\sS^t \setminus \sS^{t+1}}^\toptzero}^2 
\leq 
\normtwo{\bx_{\sS^{t+1} \setminus \sS^t}^\toptone}^2 
= \frac{\gamma^2}{\beta^2} \normtwo{\bg_{\sS^{t+1} \setminus \sS^t}^\toptzero}^2,
$
and the equality $(*)$ follows from $\normtwo{\bg_{\sS^{t+1}}^\toptzero}^2 = \normtwo{\bg_{\sS^{t+1} \setminus \sS^t}^\toptzero}^2 + \normtwo{\bg_{\sS^t \cap \sS^{t+1}}^\toptzero}^2$.
Combining \eqref{equation:conv_iht_genfunc1} and \eqref{equation:conv_iht_genfunc2} yields
\begin{equation}\label{equation:conv_iht_genfunc3}
\small
\begin{aligned}
&f(\bx^\toptone) - f(\bx^\toptzero) \leq \frac{\beta}{2} \normtwo{\bx_{\sI}^\toptone - \bx_{\sI}^\toptzero + \frac{\gamma}{\beta} \bg_{\sI}^\toptzero}^2 - \frac{\gamma^2}{2\beta} \normtwo{\bg_{\sI}^\toptzero}^2 - \frac{\gamma(1-\gamma)}{2\beta} \normtwo{\bg_{\sS^t \cup \sS^{t+1}}^\toptzero}^2 \\
&= \frac{\beta}{2} \normtwo{\bx_{\sI}^\toptone - \bx_{\sI}^\toptzero + \frac{\gamma}{\beta} \bg_{\sI}^\toptzero}^2 
- \frac{\gamma^2}{2\beta} \normtwo{\bg_{\sI \setminus (\sS^t \cup \sS^*)}^\toptzero}^2 
- \frac{\gamma^2}{2\beta} \normtwo{\bg_{\sS^t \cup \sS^*}^\toptzero}^2
- \frac{\gamma(1-\gamma)}{2\beta} \normtwo{\bg_{\sS^t \cup \sS^{t+1}}^\toptzero}^2.
\end{aligned}
\end{equation}
Next, we will upper bound the first two terms on the right-hand side  of the above inequality. Since $\sI \setminus (\sS^t \cup \sS^*) = \sS^{t+1} \setminus (\sS^t \cup \sS^*) \subseteq \sS^{t+1}$, we have $\bx_{\sI \setminus (\sS^t \cup \sS^*)}^\toptone = \bx_{\sI \setminus (\sS^t \cup \sS^*)}^\toptzero - \frac{\gamma}{\beta} \bg_{\sI \setminus (\sS^t \cup \sS^*)}^\toptzero$. 
However, since $\bx_{\sI \setminus \sS^t}^\toptzero = \bzero$, the preceding equality reduces to $\bx_{\sI \setminus (\sS^t \cup \sS^*)}^\toptone = -\frac{\gamma}{\beta} \bg_{\sI \setminus (\sS^t \cup \sS^*)}^\toptzero$. 
Let $\sT \subseteq \sS^t \setminus \sS^{t+1}$ such that $\abs{\sT} = \abs{\sS^{t+1} \setminus (\sS^t \cup \sS^*)}$. Such a choice is possible since $\abs{\sS^{t+1} \setminus (\sS^t \cup \sS^*)} = \abs{\sS^t \setminus \sS^{t+1}} - \abs{(\sS^{t+1} \cap \sS^*) \setminus \sS^t}$ (which itself is a consequence of the fact that $\abs{\sS^{t+1}} = \abs{\sS^t}$ such that $\abs{\sS^t \setminus \sS^{t+1}} = \abs{\sS^{t+1} \setminus \sS^t}$). Moreover, since $\bx^\toptone$ is obtained by hard-thresholding $\left(\bx^\toptzero - \frac{\gamma}{\beta} \bg^\toptzero\right)$, for any choice of $\sT$ made above, we have:
\begin{equation}
\frac{\gamma^2}{\beta^2} \normtwo{\bg_{\sS^{t+1} \setminus (\sS^t \cup \sS^*)}^\toptzero}^2 
= \normtwo{\bx_{\sS^{t+1} \setminus (\sS^t \cup \sS^*)}^\toptone}^2 \geq \normtwo{\bx_\sT^\toptzero - \frac{\gamma}{\beta} \bg_\sT^\toptzero}^2.
\end{equation}
Using the above equation,  the fact that $\bx_\sT^\toptone = \bzero$ (since $\sT \nsubseteq \sS^{t+1}$), and the fact that $\sI \setminus (\sS^t \cup \sS^*) = \sS^{t+1} \setminus (\sS^t \cup \sS^*)$, the first two terms of equality~\eqref{equation:conv_iht_genfunc3} becomes 
\begin{equation}\label{equation:conv_iht_genfunc4}
\small
\begin{aligned}
&\quad \frac{\beta}{2} \normtwo{\bx_{\sI}^\toptone - \bx_{\sI}^\toptzero + \frac{\gamma}{\beta} \bg_{\sI}^\toptzero}^2 - \frac{\gamma^2}{2\beta} \normtwo{\bg_{\sI \setminus (\sS^t \cup \sS^*)}^\toptzero}^2 \\
&\leq \frac{\beta}{2} \normtwo{\bx_{\sI}^\toptone - \bx_{\sI}^\toptzero + \frac{\gamma}{\beta} \bg_{\sI}^\toptzero}^2 - \frac{\beta}{2} \normtwo{\bx_\sT^\toptone - \bx_\sT^\toptzero + \frac{\gamma}{\beta} \bg_\sT^\toptzero}^2 \\
&= \frac{\beta}{2} \normtwo{\bx_{\sI \setminus \sT}^\toptone - \bx_{\sI \setminus \sT}^\toptzero + \frac{\gamma}{\beta} \bg_{\sI \setminus \sT}^\toptzero}^2.
\end{aligned}
\end{equation}
Since we construct the set $\sT$ such that $\sT \subseteq \sS^t \setminus \sS^{t+1}$ satisfying $\abs{\sT} = \abs{\sS^{t+1} \setminus (\sS^t \cup \sS^*)} = \abs{\sI \setminus (\sS^t \cup \sS^*)}$, we can bound the size of $\sI \setminus \sT$ as $\abs{\sI \setminus \sT} \leq \abs{\sS^{t}} + \abs{(\sS^t \setminus \sS^{t+1}) \setminus \sT} + \abs{\sS^*} \leq s + \abs{(\sS^{t+1} \cap \sS^*) \setminus \sS^t} + s^* \leq s + 2s^*$. Also, since $\sS^{t+1} \subseteq (\sI \setminus \sT)$, we have $\bx_{\sI \setminus \sT}^\toptone = \project_{\sB_0[s]} \left(\bx_{\sI \setminus \sT}^\toptzero - \frac{\gamma}{\beta} \bg_{\sI \setminus \sT}^\toptzero\right)$.

Given the above observation with \eqref{equation:conv_iht_genfunc4} and the assumption that $s \geq 32 \left(\frac{\beta}{\alpha}\right)^2 s^*$, Lemma~\ref{lemma:spar_proj} is valid with $\widetildebx \triangleq \bx_{\sI \setminus \sT}^\toptone + \frac{\gamma}{\beta} \bg_{\sI \setminus \sT}^\toptzero$ and $\bx \triangleq\bx_{\sI \setminus \sT}^\toptzero $, whence we reduce the first two terms of equality~\eqref{equation:conv_iht_genfunc3} to
\begin{equation}\label{equation:conv_iht_genfunc5}
\small
\begin{aligned}
&\frac{\beta}{2} \normtwo{\bx_{\sI}^\toptone - \bx_{\sI}^\toptzero + \frac{\gamma}{\beta} \bg_{\sI}^\toptzero}^2 - \frac{\gamma^2}{2\beta} \normtwo{\bg_{\sI \setminus (\sS^t \cup \sS^*)}^\toptzero}^2 
\leq \frac{\beta}{2} \cdot \frac{\abs{\sI \setminus \sT} - s}{\abs{\sI \setminus \sT} - s^*} \normtwo{\bx_{\sI \setminus \sT}^* - \bx_{\sI \setminus \sT}^\toptzero + \frac{\gamma}{\beta} \bg_{\sI \setminus \sT}^\toptzero}^2 \\
&\overset{\dag}{\leq} \frac{\beta}{2} \cdot \frac{2s^*}{s + s^*} \normtwo{\bx_{\sI}^* - \bx_{\sI}^\toptzero + \frac{\gamma}{\beta} \bg_{\sI}^\toptzero}^2 
= \frac{2s^*}{s + s^*} \cdot \left(\gamma \langle \bx^* - \bx^\toptzero, \bg^\toptzero \rangle + \frac{\beta}{2} \normtwo{\bx^* - \bx^\toptzero}^2 + \frac{\gamma^2}{2\beta} \normtwo{\bg_{\sI}^\toptzero}^2\right) \\
&\overset{\ddag}{\leq} \frac{2s^*}{s + s^*} \cdot \left(\gamma f(\bx^*) - \gamma f(\bx^\toptzero) + \frac{\beta - \gamma \alpha}{2} \normtwo{\bx^* - \bx^\toptzero}^2 + \frac{\gamma^2}{2\beta} \normtwo{\bg_{\sI}^\toptzero}^2\right),
\end{aligned}
\end{equation}
where the inequality $(\dag)$ follows from $|\sI \setminus \sT| \leq s + 2s^*$ as shown earlier and the observation that $\frac{z-a}{z-b}$ is a positive and increasing function on the interval $z \geq a$ if $a \geq b \geq 0$. Note that since we already have $\sS^{t+1} \subseteq (\sI \setminus \sT)$, we get $|\sI \setminus \sT| \geq s$. The inequality $(\ddag)$ follows by the RSC property of functions (Definition~\ref{definition:res_scss_func}).

Plugging \eqref{equation:conv_iht_genfunc5} into the first two terms of equality~\eqref{equation:conv_iht_genfunc3}, and using the fact that $\sS^{t+1} \setminus (\sS^t \cup \sS^*) \subseteq (\sS^{t+1} \cup \sS^t )$, we get:
\begin{equation}
\small
\begin{aligned}
f(\bx^\toptone) - f(\bx^\toptzero) 
&\leq \frac{2s^*}{s + s^*} \cdot \left( \gamma f(\bx^*) - \gamma f(\bx^\toptzero) + \frac{\beta - \gamma \alpha}{2} \normtwo{\bx^* - \bx^\toptzero}^2 + \frac{\gamma^2}{2\beta} \normtwo{\bg_\sI^\toptzero}^2 \right) \\
&\quad - \frac{\gamma^2}{2\beta} \normtwo{\bg_{\sS^t \cup \sS^*}^\toptzero}^2 - \frac{\gamma(1-\gamma)}{2\beta} \normtwo{\bg_{\sS^{t+1} \setminus (\sS^t \cup \sS^*)}^\toptzero}^2.
\end{aligned}
\end{equation}
Since $\gamma = 2/3$ and $s \geq 32 \left( \frac{\beta}{\alpha} \right)^2 s^*$, so that we have $\frac{2s^*}{s + s^*} \leq \frac{\alpha^2}{16\beta(\beta - \gamma \alpha)}$. Since $\beta \geq \alpha$, we also have $\frac{\alpha^2}{16\beta(\beta - \gamma \alpha)} \leq \frac{3}{16}$. Using these inequalities, rearrange the terms in \eqref{equation:conv_iht_genfunc5} above:
\begin{equation}
\begin{aligned}
f(\bx^\toptone) - f(\bx^\toptzero) &\leq \frac{2s^*}{s + s^*} \cdot \gamma \cdot \left( f(\bx^*) - f(\bx^\toptzero) \right) + \frac{\alpha^2}{32\beta} \normtwo{\bx^* - \bx^\toptzero}^2 + \frac{1}{24\beta} \normtwo{\bg_\sI^\toptzero}^2 \\
&\quad - \frac{2}{9\beta} \normtwo{\bg_{\sS^t \cup \sS^*}^\toptzero}^2 - \frac{1}{9\beta} \normtwo{\bg_{\sS^{t+1} \setminus (\sS^t \cup \sS^*)}^\toptzero}^2.
\end{aligned}
\end{equation}
Decomposing $\normtwo{\bg_\sI^\toptzero}^2 = \normtwo{\bg_{\sS^t \cup \sS^*}^\toptzero}^2 + \normtwo{\bg_{\sS^{t+1} \setminus (\sS^t \cup \sS^*)}^\toptzero}^2$ implies
\begin{equation}\label{equation:conv_iht_genfunc6}
\small
\begin{aligned}
f(\bx^\toptone) - f(\bx^\toptzero) &\leq \frac{2s^*}{s + s^*} \cdot \gamma \cdot \left( f(\bx^*) - f(\bx^\toptzero) \right) - \frac{1}{2\beta} \left( \frac{13}{36} \normtwo{\bg_{\sS^t \cup \sS^*}^\toptzero}^2 - \frac{\alpha^2}{16} \normtwo{\bx^* - \bx^\toptzero}^2 \right) \\
&\quad - \frac{1}{2\beta} \cdot \left( \frac{4}{9} - \frac{1}{12} \right) \normtwo{\bg_{\sS^{t+1} \setminus (\sS^t \cup \sS^*)}^\toptzero}^2 \\
&\leq \frac{2s^*}{s + s^*} \cdot \gamma \cdot \left( f(\bx^*) - f(\bx^\toptzero) \right) - \frac{13}{72\beta} \left( \normtwo{\bg_{\sS^t \cup \sS^*}^\toptzero}^2 - \frac{\alpha^2}{4} \normtwo{\bx^* - \bx^\toptzero}^2 \right). \\
\end{aligned}
\end{equation}
Using the RSC property, we have:
$$
\small
\begin{aligned}
f(\bx^\toptzero) - f(\bx^*) 
&\leq \innerproduct{\bg^\toptzero, \bx^\toptzero - \bx^*}- \frac{\alpha}{2} \normtwo{\bx^* - \bx^\toptzero}^2 
= \innerproduct{\bg^\toptzero_{\sS^t \cup \sS^*}, \bx^\toptzero_{\sS^t \cup \sS^*} - \bx^*_{\sS^t \cup \sS^*}} - \frac{\alpha}{2} \normtwo{\bx^* - \bx^\toptzero}^2 \\
&\leq \normtwo{\bg^\toptzero_{\sS^t \cup \sS^*}} \normtwo{\bx^\toptzero - \bx^*} - \frac{\alpha}{2} \normtwo{\bx^* - \bx^\toptzero}^2.
\end{aligned}
$$
This shows 
$$
\small
\begin{aligned}
\normtwo{\bg_{\sS^t \cup \sS^*}^\toptzero}^2 &- \frac{\alpha^2}{4} \normtwo{\bx^* - \bx^\toptzero}^2
= \left( \normtwo{\bg_{\sS^t \cup \sS^*}^\toptzero} - \frac{\alpha}{2} {\normtwo{\bx^* - \bx^\toptzero }} \right) \left( \normtwo{\bg_{\sS^t \cup \sS^*}^\toptzero} + \frac{\alpha}{2} \normtwo{\bx^* - \bx^\toptzero} \right) \\
&\geq \frac{(f(\bx^\toptzero) - f(\bx^*))}{\normtwo{\bx^\toptzero - \bx^*}} \cdot \left( \normtwo{\bg_{\sS^t \cup \sS^*}^\toptzero} + \frac{\alpha}{2} \normtwo{\bx^* - \bx^\toptzero} \right)
\geq \frac{\alpha}{2} \cdot (f(\bx^\toptzero) - f(\bx^*)).
\end{aligned}
$$
Plugging this inequality into \eqref{equation:conv_iht_genfunc6} yields 
\begin{equation}\label{equation:conv_iht_genfunc7}
\small
\begin{aligned}
f(\bx^\toptone) - f(\bx^\toptzero) 
&\leq \frac{2s^*}{s + s^*} \cdot \gamma \cdot \left( f(\bx^*) - f(\bx^\toptzero) \right) - \frac{\alpha}{12\beta} \left( f(\bx^\toptzero) - f(\bx^*) \right)\\
&=-\left(\frac{2s^*}{s + s^*} \cdot \gamma  + \frac{\alpha}{12\beta}\right) \left(  f(\bx^\toptzero) - f(\bx^*)  \right) .
\end{aligned}
\end{equation}
The result  follows by observing that $\frac{2s^*}{s + s^*} \geq 0$.
\end{proof}

\index{Convex relaxqtion}
\index{LASSO}
\index{Constrained LASSO}
\index{Lagrangian LASSO}
\section{Sparse Optimization Problem: $\ell_1$ Relaxation}
From the above sparse optimization analysis, although under very specific constraints, the convergence of the sparse optimization using IHT can be guaranteed,  we may find sparse optimization using the $\ell_0$ norm (or hard-thresholding), which counts the number of nonzero elements in a vector, poses several significant challenges and drawbacks:
\begin{itemize}
\item Non-convexity: The  $\ell_0$ norm is inherently non-convex. This means that optimization problems involving the  $\ell_0$ norm are typically NP-hard. Finding the global minimum for such problems can be computationally prohibitive, especially as the dimensionality of the problem increases.
\item 	Combinatorial nature: Minimizing the  $\ell_0$ norm involves searching through all possible combinations of the elements to find the subset that provides the sparsest solution while satisfying the constraints of the problem. This combinatorial search is impractical for high-dimensional data sets due to its exponential complexity.
\item Instability: Solutions obtained by directly minimizing the  $\ell_0$ norm can be highly sensitive to small changes in the data. This instability makes it difficult to achieve robust solutions, particularly in noisy environments or when dealing with real-world data that may contain errors or outliers.
\item Lack of continuity and differentiability: The  $\ell_0$ norm is neither continuous nor differentiable, making it incompatible with many standard optimization techniques that rely on gradient information. This lack of smoothness complicates the application of gradient-based methods and other iterative algorithms designed for continuous functions.
\item Difficulty in incorporating into optimization algorithms: Due to its discrete nature, incorporating the  $\ell_0$ norm directly into optimization frameworks is challenging. Many algorithms require some form of relaxation or approximation (e.g., using the  $\ell_1$ norm as a convex surrogate) to make the problem tractable.
\end{itemize}
Similar to the $\ell_0$ norm, the $\ell_1$ norm encourages solutions where many coefficients are exactly zero. This is due to its geometry; in high dimensions, the corners of the  $\ell_1$ ball (a polytope) align with axes, promoting sparse solutions (see Figure~\ref{fig:p-norm-comparison-3d}).
The  $\ell_1$ norm is convex, meaning that any local minimum found will also be a global minimum. This property significantly simplifies optimization problems compared to using the non-convex $\ell_0$  norm.

To be more specific, consider the standard linear regression model in matrix-vector form
\begin{equation}
\bb = \bA \bx^* + \bepsilon,
\end{equation}
where $\bA \in \real^{m\times n}$ is the model (design) matrix, $\bepsilon \in \real^m$ is a vector of noise variables, and $\bx^* \in \real^n$ is the unknown coefficient vector. By employing convex relaxation via $\ell_1$ norms,  this section develops practical algorithms and theoretical guarantees for both the \textit{constrained form} of the \textit{LASSO (least
absolute selection and shrinkage operator)} and its Lagrangian version:
\begin{subequations}\label{equation:lassos_all}
\begin{align}
(\textbf{Constrained LASSO}):\qquad &\mathopmin{\normone{\bx}\leq R}f(\bx)\triangleq \frac{1}{2m}\normtwo{\bb - \bA \bx}^2; \label{equation:loss_cons_lasso}\\
(\textbf{Lagrangian LASSO}):\qquad &\mathopmin{\bx \in \real^n}F(\bx)\triangleq \left\{ \frac{1}{2m} \normtwo{\bb - \bA \bx}^2 + \lambda_m \normone{\bx} \right\}. \label{equation:loss_lagrang_lasso}
\end{align}
\end{subequations}
Through Lagrangian duality, there exists a correspondence between these two families of quadratic programs, where  $\lambda_m$ can be interpreted as the Lagrange multiplier associated with the constraint $\normone{\bx} \leq R$ (see Section~\ref{section:geom_int_regu} for a geometrical interpretation for this kind of correspondence).

\subsection{Sparse Optimization Algorithms}

The original form of the Lagrangian LASSO is:
\begin{equation}\label{equation:lag_lasso_alg}
\min_{\bx\in \real^n} \quad F(\bx) \triangleq f(\bx) +g(\bx)\triangleq \frac{1}{2m} \normtwo{\bb - \bA \bx}^2 + \lambda_m \normone{\bx},
\end{equation}
where $g(\bx)\triangleq\lambda_m \normone{\bx}$. Since $f$ and $g$ are both convex (Exercise~\ref{exercise:conv_quad}), one might consider using the gradient descent method (Algorithm~\ref{alg:gd_gen}) as a straightforward solution. 
However, it is observed that  the objective function $F(\bx)$ of the LASSO problem is not smooth at certain points, and thus the gradient cannot be directly obtained for the original problem. 
To address this issue, one can employ subgradient methods instead.
Although $g(\bx)$ is Lipschitz and $f(\bx)$ is Lipschitz continuously differentiable ($f$ is $\frac{1}{m}\normtwo{\bA^\top\bA}$-smooth; Example~\ref{example:lipschitz_spar}), their sum does not yield any of these properties.
Therefore, the convergence is not guaranteed.

\subsection*{Smoothing with Huber Loss}
An alternative approach involves recognizing that the nonsmooth term in the objective function,  $\normone{\bx}$, consists of the sum of absolute values of each component of  $\bx$. If a smooth function could approximate the absolute value function, then gradient-based methods could be applied to solve the LASSO problem. To achieve this, we utilize the Huber loss function (refer to Section~\ref{section:huber_estima}):
\begin{equation}
h_\delta(x) = 
\begin{cases} 
	\frac{1}{2\delta}x^2, & \abs{x} < \delta, \\
	\abs{x} - \frac{\delta}{2}, & \text{otherwise}.
\end{cases}
\end{equation}
When $\delta \to 0$, the smooth function $h_\delta(x)$ approaches the absolute value function $\abs{x}$. Figure~\ref{fig:huber} illustrates  the graph of $h_\delta(x)$ for various  values of $\delta$.

By substituting the absolute value terms with the Huber loss, the smoothed LASSO problem can be formulated as:
\begin{equation}
	\min_{\bx\in\real^n} \quad F_\delta(\bx) = \frac{1}{2m}\normtwo{\bA\bx-\bb}^2 + \lambda_m H_\delta(\bx), 
	\quad \text{with }H_\delta(\bx) \triangleq \sum_{i=1}^{n} h_\delta(x_i),
\end{equation}
where $\delta$ is the given smoothing parameter. The gradient of $F_\delta(\bx)$ can be computed as
$ \nabla F_\delta(\bx) = \frac{1}{m}\bA^\top(\bA\bx - \bb) + \lambda_m \nabla H_\delta(\bx), $
where $\nabla H_\delta(\bx)$ is defined component-wise by:
$ (\nabla H_\delta(\bx))_i = 
\footnotesize
\begin{cases} 
	\sign(x_i), & \abs{x_i} > \delta; \\
	{x_i}/{\delta}, & \abs{x_i} \leq \delta.
\end{cases}
$
Given that the gradient of $F_\delta(\bx)$ is convex and Lipschitz continuously differentiable with a constant $\beta = \frac{1}{m}\normtwo{\bA^\top \bA} + \frac{\lambda_m}{\delta}$ (see Example~\ref{example:lipschitz_spar} and Exercise~\ref{exercise:sum_sc_conv}, i.e., $\beta$-smooth), according to Theorem~\ref{theorem:pgd_smooth}, the gradient descent method with a constant stepsize $\eta=\frac{1}{\beta}$ can guarantee the convergence of the algorithm.
However, if $\delta$ is very small, it is crucial to select sufficiently small stepsizes $\eta$ to ensure the convergence of the gradient method.

\subsection*{Proximal Gradient Method and FISTA}
Given that both $f$ and $g$ functions  in the LASSO problem \eqref{equation:lag_lasso_alg} are convex, where $f$ is differentiable and $g$ is non-differentiable, problem \eqref{equation:lag_lasso_alg} can be considered  a composite model. 
We have introduced proximal gradient, generalized conditional gradient, and generalized mirror descent methods to address such composite models (Chapter~\ref{chapter:gd_convg}). 
Recalling that the latter two algorithms both involve an optimization of  $\mathop{\argmin}_{\bx\in {\real^n}} \innerproduct{\nabla f(\bx^{(t)}), \bx} + g(\bx)$ at each step, which is the same as the original LASSO problem \eqref{equation:lag_lasso_alg} since $\nabla f(\bx) = \frac{1}{m}\bA^\top (\bA\bx - \bb)$.
Therefore, the most suitable choice is the proximal gradient method.

Following Algorithm~\ref{alg:prox_gd_gen}, the update rule for the proximal gradient method applied to the LASSO problem \eqref{equation:lag_lasso_alg} at the $t$-th iteration is:
\begin{subequations}
\begin{align}
\by^\toptone  &\leftarrow \bx^\toptzero - \eta_t \frac{1}{m}\bA^\top (\bA\bx^\toptzero - \bb);\\
\bx^\toptone &\leftarrow \prox_{\eta_t g}(\by^\toptone)= \mathcalT_{\eta_t\lambda_m} (\by^\toptone)= [|\by^\toptone| - \eta_t \lambda_m \bone]_+ \hadaprod \sign(\by^\toptone),
\end{align}
\end{subequations}
where $\mathcalT_\lambda(\cdot)$ denotes the soft thresholding function (Example~\ref{example:soft_thres}).
That is, the first step performs a gradient descent, and the second step performs a shrinkage, ensuring sparsity in the solution during iterations. 
This explains why the proximal gradient algorithm is effective.
Convergence is guaranteed by Theorem~\ref{theorem:prox_conv_ss_cvx} with a constant stepsize $\eta_t = \frac{m}{\normtwo{\bA^\top\bA}}$ since $f$ is $\frac{1}{m}\normtwo{\bA^\top\bA}$-smooth (Example~\ref{example:lipschitz_spar}).

Additionally, noting that the FISTA method (Algorithm~\ref{alg:fistav1} and Lemma~\ref{lemma:seq_gamma_zeta}) is a minor  variant of the proximal gradient method. 
Thus, the update steps for FISTA applied to the LASSO problem are:
\begin{subequations}
\begin{align}
\by^\toptone &\leftarrow \mathcalT_{{\eta_t\lambda_m}} \left( \bx^\toptzero - \eta_t \frac{1}{m} \bA^\top (\bA \bx^\toptzero - \bb) \right);\\
\zeta_{t+1} &\leftarrow \frac{1 + \sqrt{1 + 4 \zeta_t^2}}{2};\\
\bx^\toptone &\leftarrow \by^\toptone + \left( \frac{\zeta_t - 1}{\zeta_{t+1}} \right) (\by^\toptone - \by^\toptzero).
\end{align}
\end{subequations}

\index{Penalty function}
\index{Basis pursuit}
\subsection*{Penalty Function Method}
Recall the Lagrangian LASSO problem:
$$
\min_{\bx\in\real^n}  \quad \frac{1}{2m} \normtwo{\bA\bx-\bb}^2 + \lambda_m \normone{\bx},
$$
where $\lambda_m > 0$ is a regularization parameter. Solving the LASSO problem ultimately aims to solve the following \textit{basis pursuit (BP)} problem:
$$
\begin{aligned}
& \min_{\bx\in\real^n}  \quad \normone{\bx}\quad \text{s.t.} \quad \bA\bx = \bb,
\end{aligned}
$$
The BP problem involves a non-smooth optimization problem with equality constraints. By applying a quadratic penalty function to the equality constraint $\bA\bx = \bb$ (Section~\ref{section:pen_func_me}), we obtain:
$$
\min_{\bx\in\real^n}  \quad \normone{\bx} + \frac{\sigma}{2} \normtwo{\bA\bx-\bb}^2.
$$
Let $m \lambda_m \triangleq \frac{1}{\sigma}$, it becomes clear that using $\frac{1}{\lambda_m}$ as the quadratic penalty factor makes the penalty function subproblem of the BP problem equivalent to the LASSO problem. This observation highlights two points:
\begin{itemize}
\item The solutions of the LASSO and BP problems are not identical; however, as $\lambda_m$ approaches zero, the solution of the LASSO problem converges to the solution of the BP problem. 
\item When $m\lambda_m$ is relatively small, based on previous discussions, the penalty function of the BP problem becomes ill-conditioned (Section~\ref{section:pen_equa}), potentially leading to slow convergence if solved directly.
\end{itemize}
According to the penalty function approach, the penalty factor should gradually increase towards infinity, analogous to initially setting a larger 
$\lambda_m$ in the LASSO problem and then continuously reducing it until reaching the desired solution. The specific algorithm is detailed in Algorithm~\ref{alg:lass_penal}.

\begin{algorithm}[h]
\caption{Lagrangian LASSO Problem  via Penalty Function Method}
\label{alg:lass_penal}
\begin{algorithmic}[1]
\State {\bfseries Input:}  Given initial value $\bx^{(1)}$, regularized parameter $\lambda_m$, initial parameter $\zeta_1\gg \lambda_m$, reducing factor $\gamma \in (0, 1)$.
\For{$t=1,2,\ldots$}
\State \algoalign{Using $\bx^\toptzero$ from last iteration as the initial guess, solve the problem $\bx^\toptone = \argmin_{\bx} \left\{ \frac{1}{2m} \normtwo{\bA\bx-\bb}^2 + \zeta_t \normone{\bx} \right\}$ by subgradient descent method or other unconstrained optimization method;}
\If{$\zeta_t = \lambda_m$}
\State Stop iteration, output $\bx^\toptone$;
\Else
\State Update/reducing the penalty factor $\zeta_{t+1} = \max\{\lambda_m, \gamma \zeta_t\}$;
\EndIf
\EndFor
\State {\bfseries Return:} final $\bx\leftarrow \bx^\toptzero $;
\end{algorithmic}
\end{algorithm}

\index{ADMM}
\subsection*{ADMM}
Note that the Lagrangian LASSO problem \eqref{equation:lag_lasso_alg} can be equivalently expressed using an auxiliary variable as follows:
$$
\min_{\bx\in\real^n} F(\bx)\triangleq f(\bx)+g(\by) 
\triangleq   
\frac{1}{2m} \normtwo{\bA\bx-\bb}^2 + \lambda_m \normone{\by} \quad \text{ s.t. }\bx=\by.
$$
This formulation is in the standard form suitable for ADMM (Section~\ref{section:admm_all}). 
Given a smoothing parameter $\sigma>0$, the updates at the $t$-th  iteration are:
$$
\begin{aligned}
\bx^\toptone 
&\leftarrow \argmin_{\bx} \left\{ \frac{1}{2m} \normtwo{\bA\bx-\bb}^2 + \frac{\sigma}{2} \normtwo{\bx - \by^\toptzero + \frac{1}{\sigma} \blambda^\toptzero}^2 \right\}\\
&= \Big(\frac{1}{m}\bA^\top \bA + \sigma \bI\Big)^{-1} \Big(\frac{1}{m}\bA^\top \bb + \sigma \by^\toptzero - \blambda^\toptzero\Big);\\
\by^\toptone 
&\leftarrow \argmin_{\by} \left\{ \lambda_m \normone{\by} + \frac{\sigma}{2} \normtwo{\bx^\toptone - \by + \frac{1}{\sigma} \blambda^\toptzero}^2 \right\}\\
&= \prox_{(\lambda_m/\sigma) \normone{\cdot}} \left( \bx^\toptone + \frac{1}{\sigma} \blambda^\toptzero \right) 
=\mathcalT_{(\lambda_m/\sigma)}(\bx);\\
\blambda^\toptone 
&\leftarrow \blambda^\toptzero +  \sigma (\bx^\toptone - \by^\toptone),
\end{aligned}
$$
where $\mathcalT_\lambda(\cdot)$ denotes the soft thresholding function (Example~\ref{example:soft_thres}).

Since $ \sigma > 0 $, $ \bA^\top \bA + \sigma \bI $ is always invertible. The update for $ \bx $  essentially solves  a ridge regression problem (an $ \ell_2 $ norm squared regularization least squares problem); while the update for $ \by $ involves  an $ \ell_1 $ norm proximal operator, which also has an explicit solution. 
When solving the $ \bx $ iteration, if a fixed penalty parameter $ \sigma $ is used, we can precompute the decomposition of the matrix $ \bA^\top \bA + \sigma \bI $. This allows us to cache this initial decomposition, thereby reducing the computational load in subsequent iterations.

\subsection*{Block Coordinate Descent}

We briefly introduce how to solve the LASSO problem using the \textit{block coordinate descent (BCD)} method.
Since the  $\normone{\bx}$ term in the objective function is separable, each block variable corresponds to an individual component of  $\bx$. For convenience, when updating the $i$-th block, we represent $\bx$ as:
$$
\bx = 
\begin{bmatrix}
x_i \\
\bx_{-i}
\end{bmatrix}
$$
where $x_i$ denotes the $i$-th component of $\bx$, and $\bx_{-i}$ represents the vector $\bx$ excluding its $i$-th component. 
Similarly, for matrix $\bA$, during the $i$-th block update, it can be written as:
$$
\bA = \begin{bmatrix}
\ba_i & \bA_{-i}
\end{bmatrix},
$$
where $\bA_{-i}$ is the matrix formed by removing the $i$-th column from $\bA$. 
This reordering of components in $\bx$ and columns in $\bA$ maintains the equivalence of the modified problem to the original one.

Subsequently, during the update of the $i$-th block, the optimization problem becomes:
$$
\min_{x_i} \quad \lambda_m \abs{x_i} + \lambda_m \normone{\bx_{-i}} + \frac{1}{2m} \normtwo{\ba_i x_i - (\bb - \bA_{-i} \bx_{-i})}^2.
$$
Let $ \bc_i \triangleq \bb - \bA_{-i} \bx_{-i} $, noting that terms involving only $\bx_{-i}$ are constants. The problem is thus equivalent to:
\begin{equation}\label{equation:bcd_lasso}
\min_{x_i} \quad f_i(x_i) \triangleq \lambda_m \abs{x_i} + \frac{1}{2m} \normtwo{\ba_i}^2 x_i^2 - \ba_i^\top \bc_i x_i.
\end{equation}
The minimizer for \eqref{equation:bcd_lasso} can be directly obtained as:
\begin{equation}\label{equation:bcd_lasso_upd}
x_i^\toptzero \leftarrow \argmin_{x_i} f_i(x_i) = 
\begin{cases}
m\cdot	\frac{\ba_i^\top \bc_i - \lambda_m}{\normtwo{\ba_i}^2}, & \ba_i^\top \bc_i > \lambda_m, \\
m\cdot	\frac{\ba_i^\top \bc_i + \lambda_m}{\normtwo{\ba_i}^2}, & \ba_i^\top \bc_i < -\lambda_m, \\
	0, & \text{otherwise}.
\end{cases}
\end{equation}
Therefore, the block coordinate descent method for solving the LASSO problem can be summarized in Algorithm~\ref{alg:bcd_lasso}.

\begin{algorithm}[H]
\caption{Block Coordinate Descent Method for LASSO}
\label{alg:bcd_lasso}
\begin{algorithmic}[1]
\State {\bfseries Input:}   $ \bA, \bb $, parameter $\lambda_m$. Initialize $ \bx^\topone $;
\For{$t=1,2,\ldots$}
\For{$i = 1, 2, \ldots, n$}
\State Compute $\bx_{-i}, \bc_i$ according to definitions;
\State Compute $ x_i^\toptzero $ by \eqref{equation:bcd_lasso_upd};
\EndFor
\EndFor
\State {\bfseries Return:} final $\bx\leftarrow \bx^\toptzero$; 
\end{algorithmic}
\end{algorithm}

\index{Group LASSO}
\paragrapharrow{Group LASSO.}
When we believe that the columns of $\bA$ can be grouped into subsets $\sG$, where each subset represents a coherent group of columns, the goal of \textit{group LASSO} is to enforce sparsity at the group level rather than at the individual element level. Specifically, the union of all groups covers all columns of $\bA$, i.e.,  $\cup_{\scriptsize g\in\sG}g=\{1,2,\ldots,n\}$. 
To achieve this group-level sparsity, we can apply the  $\ell_1$ vector norm to a vector in $\real^{\abs{\sG}}$, where each component of this vector corresponds to the  $\ell_1$ or $\ell_2$ matrix norm of the columns within each group $g\in\sG$ (the $\ell_2$ matrix norm within each group does not promote sparsity on its own). This approach avoids direct penalization of individual elements and focuses on the aggregated norms of each group.
This method is particularly useful when the underlying structure of the data suggests that features or variables naturally form meaningful groups, allowing us to select entire groups of features rather than individual features.

To be more specific, let matrix $ \bA \in \real^{m\times n} $ and vector $ \bb \in \real^m $ be composed of $ m $  observations from the independent variables and response variables in the above model. The parameter $ \bx = [\bx_1; \bx_2; \ldots; \bx_G] \in \real^n $ can be divided into $ G = \abs{\sG} $ groups, and among $ \{\bx_i\}_{i=1}^G $, only a few are nonzero vectors. The optimization problem corresponding to the group LASSO can be formulated as:
$$
\min_{\bx\in\real^m}  \quad \frac{1}{2m} \normtwo{\bb - \bA\bx}^2 + \lambda_m \sum_{i=1}^G \sqrt{p_i} \normtwo{\bx_i}.
$$
In this formulation, there are $G$ blocks of variables to be optimized, with each block representing one group. Here, $p_i$ denotes the number of elements (or dimensions) in group $i$, and $\lambda_m$ controls the level of regularization applied to encourage sparsity at the group level. This setup allows for efficient selection of relevant feature groups while maintaining model simplicity and interpretability.

\subsection{Other Sparse Optimization Formulations}
\index{Generalized LASSO}

Apart from the standard LASSO problem, we also briefly introduce the generalized LASSO and sparse noise LASSO problems.
\subsection*{Generalized LASSO Problem}
The \textit{generalized LASSO problem} is defined as:
\begin{equation}\label{equation:gen_lasso}
\min_{\bx\in\real^n}  \quad  \frac{1}{2m} \normtwo{\bA\bx-\bb}^2 +\lambda_m \normone{\bB\bx},
\end{equation}
where $\bA\in\real^{m\times n}$, $\bB\in\real^{p\times n}$, and $\bb\in\real^m$.
In the LASSO problem, adding the $\normone{\bx}$ term ensures the sparsity of $\bx$. For many problems, $\bx$ itself may not be sparse, but it becomes sparse under certain transformations. An important example is when $\bB \in \real^{(n-1) \times n}$ is a first-order difference matrix:
$$
b_{ij} = 
\begin{cases} 
	1, & j = i + 1, \\
	-1, & j = i, \\
	0, & \text{otherwise},
\end{cases}
$$
and $\bA = \bI$, the generalized LASSO problem simplifies to:
$$
\min_{\bx\in\real^n}  \quad \frac{1}{2m} \normtwo{\bx - \bb}^2 + \lambda_m \sum_{i=1}^{n-1} \abs{x_{i+1} - x_i},
$$
This formulation is known as the \textit{total variation (TV)} model for image denoising; when $\bA = \bI$ and $\bB$ is a second-order difference matrix, problem \eqref{equation:gen_lasso} is known as \textit{trend filtering} \citep{rudin1992nonlinear, kim2009ell_1}~\footnote{See Problem~\ref{prob:denoise_rls} for its denoising effect.}. 

Below, we explain how to use ADMM to solve this problem.
To solve problem \eqref{equation:gen_lasso}, by introducing the constraint $\bB\bx = \by$, we rewrite the problem in a form suitable for ADMM:
$$
\min_{\bx, \by} \quad \frac{1}{2m} \normtwo{\bA\bx - \bb}^2 + \lambda_m \normone{\by}
\quad \text{s.t.} \quad \bB\bx - \by = \bzero.
$$
Introducing the Lagrange multiplier $\blambda\in\real^p$ (assume $\bB\in\real^{p\times n}$), the augmented Lagrangian function is:
$$
L_\sigma(\bx, \by, \blambda) = \frac{1}{2m} \normtwo{\bA\bx - \bb}^2 + \lambda_m \normone{\by} + \blambda^\top (\bB\bx - \by) + \frac{\sigma}{2} \normtwo{\bB\bx - \by}^2,
$$
where $\sigma$ is the penalty factor. 
Similar to the ADMM approach for a standard LASSO problem, the update for $\bx$   involves solving the linear system:
$$
\Big(\frac{1}{m}\bA^\top \bA + \sigma \bB^\top \bB\Big)\bx = \frac{1}{m}\bA^\top \bb + \sigma \bB^\top \Big( \by^\toptzero - \frac{1}{\sigma}\blambda^\toptzero \Big),
$$
while the $\by$ update is performed using the $\ell_1$ norm proximal operator. Therefore, the ADMM iterations are:
\begin{subequations}
\begin{align}
\bx^\toptone &\leftarrow \Big(\frac{1}{m}\bA^\top \bA + \sigma \bB^\top \bB\Big)^{-1} \Big( \frac{1}{m}\bA^\top \bb + \sigma \bB^\top \Big( \by^\toptzero - \frac{1}{\sigma}\blambda^\toptzero \Big) \Big);\\
\by^\toptone &\leftarrow \prox_{(\lambda_m/\sigma) \normone{\cdot}} \Big( \bB \bx^\toptone + \frac{1}{\sigma}\blambda^\toptzero \Big);\\
\blambda^{k+1} &\leftarrow \blambda^\toptzero +  \sigma (\bB\bx^\toptone - \by^\toptone).
\end{align}
\end{subequations}

\subsection*{Sparse Noise}
In the original LASSO problem \eqref{equation:lassos_all}, the noise term $\bb-\bA\bx$ is not sparse. 
The following problem considers a sparse noise term:
\begin{equation}\label{equation:spar_noi_lass}
\min_{\bx \in \real^n} \left\{ F(\bx) \triangleq \frac{1}{m} \normone{\bA \bx - \bb} + \lambda_m \normone{\bx} \right\},
\end{equation}
where $\bA \in \real^{m \times n}$, $\bb \in \real^m$, and $\lambda_m > 0$. We will consider two possible methods to solve the problem:

\paragrapharrow{Subgradient descent.} When applying the subgradient descent method to problem \eqref{equation:spar_noi_lass}, the method takes the form (choosing the subgradient of  $\normone{\by}$ as $\sign(\by)$; see Exercise~\ref{exercise:sub_norms}):
$$
\bx^\toptone \leftarrow \bx^\toptzero - \eta_t \big(\frac{1}{m}\bA^\top \sign(\bA \bx^\toptzero - \bb) + \lambda_m \sign(\bx)\big).
$$
The stepsize $\eta_t$ is chosen according to Theorem~\ref{theorem:pgd_lipschitz} as $\eta_t = \frac{1}{\normtwo{F'(\bx^\toptzero)} \sqrt{t+1}}$. That is, this method can be viewed as a projected subgradient descent with $\sS=\real^n$.
\paragrapharrow{Proximal subgradient.} For \eqref{equation:spar_noi_lass}, let $f(\bx) \triangleq \normone{\bA \bx - \bb}$ and $g(\bx) \triangleq \lambda_m \normone{ \bx}$, so that $F = f + g$. The proximal subgradient method then takes the form
$$
\bx^\toptone \leftarrow \prox_{\eta_t g}\big(\bx^\toptzero - \eta_t\frac{1}{m} \bA^\top \sign(\bA \bx^\toptzero - \bb)\big).
$$
Since $g(\bx) = \lambda_m \normone{\bx}$, it follows that $\prox_{\eta_t g}$ is a soft thresholding operator. Specifically, by Example~\ref{example:soft_thres}, $\prox_{\eta_t g} = \mathcal{T}_{\lambda_m \eta_t}$,  and thus the general update rule becomes
$$
\bx^\toptone \leftarrow \mathcal{T}_{\lambda_m \eta_t}\big(\bx^\toptzero - \eta_t\frac{1}{m} \bA^\top \sign(\bA \bx^\toptzero - \bb)\big).
$$
The stepsize $\eta_t$ can be chosen, for example, as $\eta_t = \frac{1}{\normtwo{f'(\bx^\toptzero)} \sqrt{t+1}}$.

\subsection{Theoretical Results for LASSO}
We conclude this section by providing several theoretical results for the LASSO problem. 
Given a LASSO estimate  $\widehatbx \in \real^n$, we can evaluate its quality in various ways. 
Depending on the application, different types of loss functions may be relevant.
In some scenarios, we are interested in the predictive performance of $\widehatbx$, leading us to compute a \textit{prediction loss function}:
\begin{subequations}
\begin{equation}
	\mathcalL_{\text{pred}}(\widehatbx; \bx^*) = \frac{1}{m} \normtwo{\bA \widehatbx - \bA \bx^*}^2,
\end{equation}
which represents the mean-squared error of $\widehatbx$ over the samples provided by $\bA$. 
In other applications---such as medical imaging, remote sensing, and compressed sensing---the primary interest lies in the unknown vector $\bx^*$ itself.
Therefore, it is more appropriate to consider loss functions such as the $\ell_2$ distance between the estimator and the true parameter:
\begin{equation}
	\mathcalL_2(\widehatbx; \bx^*) = \normtwo{\widehatbx - \bx^*}^2,
\end{equation}
which we refer to as the  \textit{parameter estimation loss}. 
Finally, if variable selection or support recovery is the goal, one might use the following loss function:
\begin{equation}
	\mathcalL_{\text{vs}}(\widehatbx; \bx^*) =
	\begin{cases}
		0 & \text{if } \sign(\widehatx_i) = \sign(x_i^*) \text{ for all } i \in \{1,2, \ldots, n\}, \\
		1 & \text{otherwise}.
	\end{cases}
\end{equation}
\end{subequations}
This evaluates whether the estimated vector $\widehatbx$ shares the same signed support as $\bx^*$. 

In this subsection, for brevity, we present theoretical results for the first two cases. To understand these results, we need the following technical lemma.

\begin{lemma}\label{lemma:sparse_lem1}
Let $\bx, \by\in\real^n$  with $\normone{\bx}\leq \normone{\by}$, where $\supp(\by) \triangleq\sS\subseteq \{1,2,\ldots,n\}$.
Denote $\be \triangleq \bx-\by$.  Then $\normone{\be_{\comple{\sS}}} \leq \normone{\be_{\sS}}$, i.e., $\be \in \sC[\sS; 1]$, where $\comple{\sS}$ denotes the complement of $\sS$.
\end{lemma}
\begin{proof}[of Lemma~\ref{lemma:sparse_lem1}]
Let $\sT\triangleq\supp(\bx)$. The assumption  $\normone{\bx}\leq \normone{\by}$ indicates that 
$$
\begin{aligned}
&\sum_{i\in\sT} \abs{x_i} \leq  \sum_{i\in\sS} \abs{y_i} 
\quad\implies\quad 
\sum_{i\in\sT\setminus \sS} \abs{x_i}+ \sum_{i\in\sS \cap\sT} \abs{x_i}   \leq  \sum_{i\in\sS \setminus\sT} \abs{y_i} + \sum_{i\in\sS \cap\sT} \abs{y_i}\\
&\implies\quad  
\sum_{i\in\sT\setminus \sS} \abs{x_i}    \leq  \sum_{i\in\sS \setminus\sT} \abs{y_i} + \sum_{i\in\sS \cap\sT} \abs{y_i - x_i},
\end{aligned}
$$
which obtains the desired result.
\end{proof}

With this intuition in place, we now state a result that provides a bound on the parameter estimation error $\normtwo{\widehatbx - \bx^*}$ based on the linear observation model $\bb = \bA\bx^* + \bepsilon$, where $\bx^*$ is $k$-sparse and supported on the subset $\sS$.

\begin{theoremHigh}[Bounds of Optimizer Error for LASSO \citep{hastle2015statistical}]\label{theorem:bound_lassl2error}
Let the observations follow from the linear observation model $\bb = \bA\bx^* + \bepsilon$, where $\bA\in\real^{m\times n}$, and  $\bx^*$ is $k$-sparse, supported on the subset $\sS\subseteq \{1, 2,\ldots, n\}$ (i.e., $\abs{\sS}=k$).
Suppose that the model matrix $\bA$ satisfies the restricted eigenvalue (RE) property (Definition~\ref{definition:res_eig}) with parameter $\alpha > 0$ over $\sC[\sS;3]$.~\footnote{Note that the first result only requires the RE property holds over the subset $\sC[\sS;1]$.}
Then,
\begin{enumerate}[(i)]
\item \textit{Constrained LASSO.} Any estimate $\widehatbx$ based on the constrained LASSO \eqref{equation:loss_cons_lasso} with $\normone{\bx^*}$ satisfies the bound
$$
\normtwo{\widehatbx - \bx^*} \leq \frac{4\sqrt{k}}{\alpha m} \norminf{\bA^\top \bepsilon}.
$$

\item \textit{Lagrangian LASSO.} Given a regularization parameter $\lambda_m \geq \frac{2}{m} \norminf{\bA^\top \bepsilon} > 0$, any estimate $\widehatbx$ from the regularized LASSO \eqref{equation:loss_lagrang_lasso} satisfies the bound
$$
\normtwo{\widehatbx - \bx^*} \leq \frac{3}{\alpha} \sqrt{{k}}  \lambda_m.
$$
\end{enumerate}
\end{theoremHigh}
\begin{proof}[of Theorem~\ref{theorem:bound_lassl2error}]
\textbf{(i) Constrained LASSO.} Since $\bx^*$ is feasible and $\widehatbx$ is optimal, we have the inequality $\normtwo{\bb - \bA\widehatbx}^2 \leq \normtwo{\bb - \bA\bx^*}^2$. Defining the error vector $\widehatbe \triangleq \widehatbx - \bx^*$, substituting in the relation $\bb = \bA\bx^* + \bepsilon$ yields
\begin{equation}\label{equation:bound_lassl2error1}
\text{(Constrained LASSO)}: \qquad \frac{\normtwo{\bA\widehatbe}^2}{m} \leq \frac{2}{m} \bepsilon^\top \bA \widehatbe.
\end{equation}
Applying  H\"older's inequality (Theorem~\ref{theorem:holder-inequality}) to the right-hand side yields the upper bound $\frac{2}{m} \abs{\bepsilon^\top \bA \widehatbe} \leq \frac{2}{m} \norminf{\bA^\top \bepsilon} \normone{\widehatbe}$. The inequality $\normone{\widehatbx} \leq R = \normone{\bx^*}$ implies that $\widehatbe \in \sC[\sS; 1]$ (Lemma~\ref{lemma:sparse_lem1}):
\begin{equation}\label{equation:bound_lassl2error2}
\normone{\widehatbx} \leq \normone{\bx^*}\quad\implies\quad 
\normone{\widehatbe} = \normone{\widehatbe_{\sS}} + \normone{\widehatbe_{\comple{\sS}}} \leq 2 \normone{\widehatbe_{\sS}} 
\stackrel{\dag}{\leq} 2 \sqrt{k} \normtwo{\widehatbe_{\sS}}\leq 2 \sqrt{k} \normtwo{\widehatbe},
\end{equation}
where the inequality ($\dag$) follows from the Cauchy-Schwarz inequality (Equation~\eqref{equation:vector_form_cauchyschwarz}): 
\begin{equation}\label{equation:bound_lassl2error2_v2}
\normone{\widehatbe_{\sS}}=\abs{\widehatbe_{\sS}}^\top\bone \leq \normtwo{\widehatbe_{\sS}}\normtwo{\bone} = \sqrt{k}\normtwo{\widehatbe_{\sS}} 
\leq  \sqrt{k}\normtwo{\widehatbe}.
\end{equation}
Note that \eqref{equation:bound_lassl2error2} holds under the condition $\normone{\widehatbx} \leq \normone{\bx^*}$, while \eqref{equation:bound_lassl2error2_v2} is valid for any $\sS$ and $\widehatbe$.
On the other hand, applying the RE property to the left-hand side of the inequality \eqref{equation:bound_lassl2error1} yields $ \alpha \normtwo{\widehatbe}^2 \leq \frac{1}{m} \normtwo{\bA \widehatbe}^2$. Combining these inequalities yields the desired result.

\paragraph{(ii) Lagrangian LASSO.} Define the function
\begin{equation}
G(\be) \triangleq \frac{1}{2m} \normtwo{\bb - \bA (\bx^* + \be)}^2 + \lambda_m \normone{\bx^* + \be}.
\end{equation}
Noting that $\widehatbe \triangleq \widehatbx - \bx^*$ minimizes $G$ by construction, we have $G(\widehatbe) \leq G(\bzero)$. 
Substituting in the relation $\bb = \bA\bx^* + \bepsilon$ yields
\begin{equation}\label{equation:bound_lassl2error3}
\text{(Lagrangian LASSO)}:\qquad \frac{\normtwo{\bA \widehatbe}^2}{m} \leq \frac{2}{m} \bepsilon^\top \bA \widehatbe + 2\lambda_m \left\{ \normone{\bx^*} - \normone{\bx^* + \widehatbe} \right\}.
\end{equation}
Since $\bx_{\comple{\sS}}^* = \bzero$, we have $\normone{\bx^*} = \normone{\bx_{\sS}^*}$, and
$$
\normone{\bx^* + \widehatbe} = \normone{\bx_{\sS}^* + \widehatbe_{\sS}} + \normone{\widehatbe_{\comple{\sS}}} \geq \normone{\bx_{\sS}^*} - \normone{\widehatbe_{\sS}} + \normone{\widehatbe_{\comple{\sS}}}.
$$
Substituting these relations into the inequality \eqref{equation:bound_lassl2error3} and applying  H\"older's inequality (Theorem~\ref{theorem:holder-inequality}) yield
\begin{equation}\label{equation:bound_lassl2error5}
\begin{aligned}
\frac{\normtwo{\bA \widehatbe}^2}{m} 
&\leq \frac{2 }{m}\bepsilon^\top \bA \widehatbe + 2\lambda_m \left\{ \normone{\widehatbe_{\sS}} - \normone{\widehatbe_{\comple{\sS}}} \right\} \\ 
&\leq \frac{2 }{m} \norminf{\bA^\top \bepsilon}\normone{\widehatbe} + 2\lambda_m \left\{ \normone{\widehatbe_{\sS}} - \normone{\widehatbe_{\comple{\sS}}} \right\}.
\end{aligned}
\end{equation}
Since $\frac{2}{m} \norminf{\bA^\top \bepsilon} \leq \lambda_m$ by assumption, using the fact that $\normone{\widehatbe_{\sS}} \leq \sqrt{k} \normtwo{\widehatbe}$ from \eqref{equation:bound_lassl2error2_v2}, we have:
\begin{equation}\label{equation:bound_lassl2error6}
0\leq \frac{\normtwo{\bA \widehatbe}^2}{m} \leq {\lambda_m} \left\{ \normone{\widehatbe_{\sS}} + \normone{\widehatbe_{\comple{\sS}}} \right\} + 2\lambda_m \left\{ \normone{\widehatbe_{\sS}} - \normone{\widehatbe_{\comple{\sS}}} \right\} \leq {3} \sqrt{k} \lambda_m \normtwo{\widehatbe}.
\end{equation}
The inequality \eqref{equation:bound_lassl2error6} also implies that $\normone{\widehatbe_{\comple{\sS}}} \leq 3 \normone{\widehatbe_{\sS}}$ such that $\widehatbe\in \sC[\sS; 3]$.
Applying the $\alpha$-RE condition to $\widehatbe$ ensures that $\alpha \normtwo{\widehatbe}^2 \leq \frac{1}{m} \normtwo{\bA \widehatbe}^2$. 
Combining this lower bound with the inequality \eqref{equation:bound_lassl2error6} yields
$
{\alpha} \normtwo{\widehatbe}^2 \leq {3}  \sqrt{k} \lambda_m \normtwo{\widehatbe},
$
which implies the desired result.
\end{proof}

\begin{theoremHigh}[Bounds of Prediction Error for LASSO \citep{hastle2015statistical}]\label{theorem:lassobound_predict_error}
Let the observations follow from the linear observation model $\bb = \bA\bx^* + \bepsilon$, where $\bA\in\real^{m\times n}$.
Consider the Lagrangian LASSO \eqref{equation:loss_lagrang_lasso} with a regularization parameter $\lambda_m \geq \frac{2}{m} \norminf{\bA^\top \bepsilon}$.
\begin{enumerate}[(i)]
\item If $\normone{\bx^*} < R_1$, then any optimal solution $\widehatbx$ satisfies
$$
\frac{\normtwo{\bA (\widehatbx - \bx^*)}^2}{m} \leq 12 \, R_1 \, \lambda_m.
$$
\item Suppose further that $\bx^*$ is $k$-sparse, supported on the subset $\sS\subseteq \{1, 2,\ldots, n\}$ (i.e., $\abs{\sS}=k$).
If the model matrix $\bA$ satisfies the restricted eigenvalue property with parameter $\alpha > 0$ over $\sC[\sS;3]$, then any optimal solution $\widehatbx$ satisfies
$$
\frac{\normtwo{\bA (\widehatbx - \bx^*)}^2}{m} \leq 9\frac{ k \lambda_m^2}{\alpha}.
$$
\end{enumerate}
\end{theoremHigh}
\begin{proof}[of Theorem~\ref{theorem:lassobound_predict_error}]
The proof is a slight adaptation based on \citet{hastle2015statistical}.
\paragraph{(i).} Once again, we define the error vector $\widehatbe \triangleq \widehatbx - \bx^*$. Following the inequality \eqref{equation:bound_lassl2error3}, we have
\begin{align}
0 &\leq\frac{\normtwo{\bA \widehatbe}^2}{m}  \leq \frac{2}{m} \norminf{\bA^\top \bepsilon} \normone{\widehatbe} + 2\lambda_m \big\{ \normone{\bx^*} - \normone{\bx^* + \widehatbe} \big\} \label{equation:lassobound_predict_error_ineq1}\\
&\overset{\dag}{\leq} 2\left\{ \frac{\norminf{\bA^\top \bepsilon}}{m} - \lambda_m \right\} \normone{\widehatbe} + 4 \lambda_m \normone{\bx^*} 
\overset{\ddag}{\leq}  \lambda_m \Big\{ - \normone{\widehatbe} + 4 \normone{\bx^*} \Big\}, 
\end{align}
where the inequality $(\dag)$ follows from the triangular inequality 
$-(\normone{\bx^*} - \normone{\widehatbe} ) \leq \normone{\bx^* - (-\widehatbe)} 
$~\footnote{$\abs{\norm{\bx}-\norm{\by}} \leq \norm{\bx-\by}$ for any vector norm.}, and  the inequality $(\ddag)$ follows from the assumption that $\lambda_m \geq \frac{2}{m} \norminf{\bA^\top \bepsilon}$.
This indicates that $\normone{\widehatbe} \leq 4 \normone{\bx^*} \leq 4 R_1$. Combining these inequalities and applying the triangle inequality $\normone{\bx^*} - \normone{\bx^* + \widehatbe} \leq \normone{\bx^* -(\bx^* + \widehatbe)}=\normone{\widehatbe}$ in \eqref{equation:lassobound_predict_error_ineq1} yield that
\begin{equation}\label{equation:lassobound_predict_error_part1fin}
\frac{\normtwo{\bA \widehatbe}^2}{m} \leq 2\left\{ \frac{\norminf{\bA^\top \bepsilon}}{m} + \lambda_m \right\} \normone{\widehatbe} \leq 12 \lambda_m R_1,
\end{equation}
which establishes the desired result.

\paragraph{(ii).} Following \eqref{equation:bound_lassl2error6}, we have
$
\frac{\normtwo{\bA \widehatbe}^2}{m}   \leq 3  \sqrt{k}\lambda_m \normtwo{\widehatbe}.
$~\footnote{
Note that \citet{hastle2015statistical} use \eqref{equation:lassobound_predict_error_part1fin} and the fact $\normone{\widehatbe} = \normone{\widehatbe_{\sS}} + \normone{\widehatbe_{\comple{\sS}}} \leq 4 \normone{\widehatbe_{\sS}} 
\leq 4 \sqrt{k} \normtwo{\widehatbe_{\sS}}\leq 4 \sqrt{k} \normtwo{\widehatbe}$ to obtain 
$\frac{\normtwo{\bA \widehatbe}^2}{m} \leq 12\sqrt{k}\normtwo{\widehatbe}$, which results in a looser bound.
}
By the proof of Theorem~\ref{theorem:bound_lassl2error}, the error vector $\widehatbe$ belongs to the cone $\sC[\sS;3]$, so that the $\alpha$-RE condition guarantees that $\normtwo{\widehatbe}^2 \leq \frac{1}{m \alpha} \normtwo{\bA \widehatbe}^2$. Combining these inequalities yields the desired result.
\end{proof}

\citet{hastle2015statistical} demonstrate that for various statistical models, choosing $\lambda_m = \mathcalO\big(\sqrt{\frac{\ln(n)}{m}}\big)$ is valid for Theorem~\ref{theorem:lassobound_predict_error} with high probability. Consequently, the two bounds take the form:
\begin{subequations}
\begin{align}
\frac{\normtwo{\bA (\widehatbx - \bx^*)}^2}{m} &\leq \mathcalO\big(R_1 \, \sqrt{\frac{\ln(n)}{m}}\big); \label{equation:lassobound_predict_error1} \\
\frac{\normtwo{\bA (\widehatbx - \bx^*)}^2}{m} &\leq \mathcalO\big(\frac{ \abs{\sS} \, \ln(n)}{m}\big). \label{equation:lassobound_predict_error2}
\end{align}
\end{subequations}
The bound \eqref{equation:lassobound_predict_error1}, which depends on the $\ell_1$-ball radius $R_1$, is known as the ``slow rate'' for the LASSO problem, since the squared prediction error decays as $1/\sqrt{m}$. On the other hand, the bound \eqref{equation:lassobound_predict_error2} is known as the ``fast rate,'' since it decays as $1/m$. 
Note that the latter relies on stronger assumptions: specifically, the hard sparsity condition that $\bx^*$ is supported on a small subset $\sS$, and more critically, the $\alpha$-RE condition on the design matrix $\bA$.
In principle, prediction performance should not require an RE condition, suggesting that this requirement might be an artifact of our proof technique.

\begin{problemset}

\item \label{prob:als_pseudo1} \textbf{Least squares for rank-deficiency.} Let $\bA\in\real^{m\times n}$ and $\bb\in\real^m$. Show that the least squares problem $L(\bx)=\normtwo{\bA\bx-\bb}^2$ has a minimizer $\bx^*\in\real^n$ if and only if there exists a vector $\by\in\real^n$ such that $\bx^*=\bA^+\bb+(\bI-\bA^+\bA)\by$, where $\bA^+$ is the pseudo-inverse of $\bA$. 
\begin{itemize}
\item This shows that the least squares has a \textbf{unique} minimizer of $\bx^*=\bA^+\bb$ only when $\bA^+$ is a left inverse of $\bA$ (i.e., $\bA^+\bA = \bI$). The solution in Lemma~\ref{lemma:ols} is a special case.
\item The optimal value is $L(\bx^*)=\bb^\top(\bI-\bA\bA^+)\bb$.
\item If $\by\neq \bzero$: $\normtwo{\bA^+\bb}\leq \normtwo{\bA^+\bb+(\bI-\bA^+\bA)\by}$.
\end{itemize}
\textit{Hint: using singular value decomposition.}

\item \label{prob:als_pseudo2} \textbf{Least squares for rank-deficiency.} Let  $\bA\in\real^{m\times n}$ and $\bB\in\real^{m\times p}$. Show that the least squares problem $L(\bX) = \normf{\bA\bX-\bB}^2$ has a minimizer $\bX^*=\bA^+\bB\in\real^{n\times p}$. Determine all the minimizers using Problem~\ref{prob:als_pseudo1}.

\item \label{prob:als_pseudon} \textbf{Least squares for rank-deficiency.}  Let  $\bA\in\real^{m\times n}$ and $\bB\in\real^{p\times n}$. Show that the least squares problem $L(\bX) = \normf{\bX\bA-\bB}^2$ has a minimizer $\bX^*=\bB\bA^+\in\real^{p\times m}$.

\item \label{problem:rls} \textbf{Regularized least squares (RLS).} 
Given $\bA\in\real^{m\times n}, \bb\in\real^{m}, \bB\in\real^{p\times n}$, and $\lambda\in\real_{++}$, we consider the constrained least squares problem:
$$
\mathop{\min}_{\bx\in\real^n} \normtwo{\bA\bx-\bb}^2 + \lambda\normtwo{\bB\bx}^2.
$$
Show that the constrained least squares (RLS) problem has a unique solution if and only if $\nspace(\bA)\cap \nspace(\bB) = \{\bzero\}$.

\item \label{prob:denoise_rls} \textbf{Denoising via RLS.} Consider a noisy measurement of a signal $\bx\in\real^n$:
$
\by = \bx+\be,
$
where $\by$ is the observed measurement, and $\be$ is the noise vector. We want to find an estimate $\bx$ of the observed measurement $\by$ such that $\bx \approx \by$:
$
\min \normtwo{\bx-\by}^2.
$
Apparently, the optimal solution of this optimization is given by $\bx=\by$; however, it is meaningless.
To improve the estimate,  we can add a penalty term for the differences between  consecutive observations:
$
R(\bx) = \sum_{i=1}^{n-1} (x_i - x_{i+1})^2.
$
Then, 
\begin{itemize}
\item Find the constrained least squares representation for this problem and derive the constrained least squares solution.
\item Find some applications of this denoising problem. For example, when we model the profit and loss signal of a financial asset, the two observations over consecutive days of the underlying asset should exhibit smooth transitions rather than abrupt changes.
\end{itemize}

\item \textbf{Weighted least squares (WLS).}
Building upon the assumptions in Lemma~\ref{lemma:ols}, we consider further that each data point $i\in\{1,2,\ldots, m\}$ (i.e., each row of $\bA$) has a weight $w_i$. 
This means some  data points may carry greater significance than others, and we can produce approximate minimizers that reflect this.
Show that the value $\bx_{WLS} = (\bA^\top\bW^2\bA)^{-1}\bA^\top\bW^2\bb$ serves as the \textit{weighted least squares (WLS)}  estimate of $\bx$, where $\bW=\diag(w_1, w_2, \ldots, w_m)\in\real^{m\times m}$. \textit{Hint: Find the normal equation for this problem.}

\item \textbf{Positive definite  weighted least squares (PDWLS).}
Building upon the assumptions in Lemma~\ref{lemma:ols}, we consider further  the matrix equation $\bA\bx + \be =\bb$, where $\be$ is an error vector. Define the weighted error squared sum $E_w = \be^\top \bW \be$, where the weighting matrix $\bW$ is  positive definite. 
Show that the positive definite weighted least squares solution is $\bx^* = (\bA^\top\bW\bA)^{-1}\bA^\top\bW\bb$. \textit{Hint: Compute the gradient of $E_w = (\bb-\bA\bx)^\top\bW(\bb-\bA\bx)$.}

\item \label{problem:tls} \textbf{Transformed least squares (TLS).}
Building upon the assumptions in Lemma~\ref{lemma:ols}, we consider further the restriction $\bx=\bC\bgamma+\bc$, where $\bC\in\real^{n\times k}$ is a known matrix such that $\bA\bC$ has full rank, $\bc$ is a known vector, and $\bgamma$ is an unknown vector.
Show that the value $\bx_{TLS}=\bC(\bC^\top\bA^\top\bA\bC)^{-1}(\bC^\top\bA^\top)(\bb-\bA\bc) +\bc$ serves as the \textit{transformed least squares (TLS)} estimate of $\bx$.

\item \label{problem:twls2} Find the transformed weighted least squares estimate.

\item Verify that $\bG^\toptone \bP^\toptone = \bI$ in \eqref{equation:secan_dogleg}.

\item Prove that the set $\sC[k] \triangleq \bigcup_{\sS: \abs{\sS} = k} \sC[\sS]$ is non-convex.

\item Let $ \bA $ be a design matrix that satisfies the nullspace property of order $ 2s $ (Definition~\ref{definition:nullspace_prop}). Show that for any two distinct $ s $-sparse vectors $ \bu_1, \bu_2 \in \sB_0[s] $, $ \bu_1 \neq \bu_2 $, it follows that  $ \bA \bu_1 \neq \bA \bu_2 $.

\end{problemset}
\newpage
\chapter{Stochastic Optimization}\label{chapter:stochastic_opt}
\begingroup
\hypersetup{linkcolor=structurecolor,
linktoc=page,  
}
\minitoc \newpage
\endgroup
\index{Stochastic Optimization}

\section{Stochastic Optimization}
This chapter concludes with an exploration of stochastic optimization algorithms (a.k.a., stochastic optimizers), a critical approach for tackling large-scale and data-intensive problems. It covers key topics such as \textit{stochastic gradient descent (SGD)}, \textit{variance reduction methods like SAG, SAGA, and SVRG}, and adaptive optimizers including \textit{Adam, RMSProp, AdaGrad, and Nadam}. The discussion also includes learning rate strategies, such as \textit{annealing, warmup, and cyclical learning rates}, offering valuable insights into practical deep learning optimization.

Over time, stochastic gradient-based optimization has become a cornerstone in numerous scientific and engineering disciplines, particularly in areas like computer vision, large language models (LLM), and automatic speech recognition \citep{krizhevsky2012imagenet, hinton2012deep, graves2013speech, vaswani2017attention}. 
SGD, alongside deep neural networks or transformer architectures, plays a pivotal role in optimizing stochastic objective functions. 
When developing a new deep neural network or a transformer structure for a specific task, one must often select hyperparameters related to training through heuristic methods. For each combination of structural hyperparameters, a fresh network is usually trained from scratch and evaluated repeatedly. Although significant advancements in hardware (e.g., GPUs) and software (e.g., cuDNN) have accelerated the training process for individual network or transformer structures, exploring a vast array of potential configurations remains time-consuming. Therefore, there is a need for stochastic optimizers that are robust to hyperparameter settings. Efficient stochastic optimizers are thus essential for training deep neural networks and transformers effectively.

\begin{table}[h]
\begin{tabular}{llll|lll}
\hline
Method         & Year & Papers \gap \gap  &  & Method                  & Year & Papers \\ \hline\hline
Adam           & 2014 & 22,319   &  & AdamW                   & 2017 & 180     \\
SGD            & 1951 & 1916   &  & Local SGD               & 2018 & 65     \\
RMSProp        & 2013 & 495    &  & Gravity                 & 2021 & 177     \\
Adafactor      & 2018 & 691    &  & AMSGrad                 & 2019 & 49     \\
Momentum       & 1999 & 141    &  & LARS                    & 2017 & 74     \\
LAMB           & 2019 & 191    &  & MAS                     & 2020 & 135     \\
AdaGrad        & 2011 & 181    &  & DFA                     & 2016 & 48     \\
Deep Ensembles & 2016 & 163     &  & Nesterov momentum    & 1983 & 32     \\
ADOPT          & 2024 & 286     &  & AO    & 2000 & 58     \\
RAdam          & 2019 & 56     &  & Apollo    & 2020 & 46     \\
FA             & 2014 & 112     &  & Gradient Sparsification & 2017 & 36     \\ \hline
\end{tabular}
\caption{Data retrieved on February 9th, 2025 via https://paperswithcode.com/.}
\label{table:stochastic-optimizers}
\end{table}

Several adaptations of SGD aim to estimate an optimal learning rate at each iteration, either by accelerating progress when feasible or decelerating it near local minima. In this section, we present some stochastic optimizers from these two categories. Table~\ref{table:stochastic-optimizers} offers an overview of how frequently these stochastic optimization algorithms have been utilized in research papers for specific tasks and their publication dates. For more comprehensive reviews, readers may consult sources such as \citet{zeiler2012adadelta}, \citet{ruder2016overview},  \citet{goodfellow2016deep}, and \citet{bottou2018optimization}.

\begin{algorithm}[h] 
\caption{Stochastic Projected Gradient Descent Method}
\label{alg:stopgd_gen}
\begin{algorithmic}[1] 
\Require A function $f(\bx)$ and a set $\sS$; 
\State {\bfseries Input:} Initialize $\bx^{(1)}\in\real^n$;
\For{$t=1,2,\ldots$}
\State Pick a stepsize $\eta_t$ and a stochastic vector $\bg^\toptzero\in\real^n$;
\State $\by^{(t+1)} \leftarrow \bx^{(t)} - \eta_t \bg^\toptzero$;
\State $\bx^{(t+1)} \leftarrow \mathcalP_{\sS}(\by^{(t+1)})$;
\EndFor
\State (Output Option 1) Output  $\bx_{\text{final}}\leftarrow \bx^{(T)}$;
\State (Output Option 2) Output  $\bx_{\text{avg}}\leftarrow \frac{1}{T}(\sum_{t=1}^{t}\bx^{(t)})$ or $\sum_{t=1}^{T} \frac{2t}{T(T+1)} \bx^{(t)}$;
\State (Output Option 3) Output  $\bx_{\text{best}}\leftarrow \argmin_{t\in\{1,2,\ldots,T\}} f(\bx^{(t)})$;
\end{algorithmic} 
\end{algorithm}
\section{Stochastic (Projected) Gradient Descent}
\subsection{Results of Stochastic Projected Gradient Descent}
We discussed the convergence results of the projected gradient descent  method in Section~\ref{section:pgd}. The \textit{stochastic projected gradient descent  (SPGD)} method (Algorithm~\ref{alg:stopgd_gen}) yields similar outcomes. 
As its name implies, SPGD employs a stochastic gradient vector instead of the true gradient vector used in gradient-based algorithms.
Let us begin by outlining the assumptions required for this method.
\begin{assumption}[Stochastic Optimizer Assumption]\label{assumption:sto_assump}
Given a convex function $f:\sS\rightarrow \real^n$, for the analysis of stochastic optimizer presented in Algorithm~\ref{alg:stopgd_gen}, we assume the following 
\begin{enumerate}[(A)]
\item \textit{Unbiasedness of random subgradient.} For any iteration  $ t > 0 $, $ \Exp[\bg^\toptzero \mid \bx^\toptzero] \in \partial f(\bx^\toptzero) $, i.e., $
f(\bz) \geq f(\bx^\toptzero) + \innerproduct{\Exp[\bg^\toptzero \mid \bx^\toptzero], \bz - \bx^\toptzero}$ for all $\bz \in \dom(f).
$
\item \textit{Boundedness of second moment.} There exists a constant $ \widetilde{L} > 0 $ such that for any $ t > 0 $, $ \Exp\big[\normtwobig{\bg^\toptzero}^2 \mid \bx^\toptzero\big] \leq \widetilde{L}^2 $.
\end{enumerate}
\end{assumption}

Part (A) of the assumption indicates  that $ \bg^\toptzero $ is an unbiased estimator of a subgradient at $ \bx^\toptzero $. 
The constant $ \widetilde{L} $ from part (B) of Assumption~\ref{assumption:sto_assump} does not necessarily serve as a Lipschitz constant for $ f $, unlike in deterministic cases.

The analysis of the stochastic projected subgradient closely mirrors that of its deterministic counterpart. We summarize the key findings in the following theorem. 
Note that for any iteration $t>0$, we define the best achieved function value as 
\begin{equation}
\fbest^\toptzero \triangleq \min\left\{f(\bx^\topone), f(\bx^{(2)}), \ldots, f(\bx^\toptzero)\right\}.
\end{equation}
\begin{theoremHigh}[SPGD for Convex and ``Lipschitz" Functions]\label{theorem:conv_sgpg}
Let $f: \sS \rightarrow \real$ be a proper closed and convex function, where $\sS$ is a closed and convex set. Let $ \{\bx^\toptzero\}_{t > 0} $ be the sequence generated by the stochastic projected subgradient method (Algorithm~\ref{alg:stopgd_gen}) with positive stepsizes $ \{\eta_t\}_{t > 0} $, and let $\left\{ \fbest^\toptzero\right\}_{t > 0}$ be the sequence of best achieved values.
Given the assumptions and vectors in Assumption~\ref{assumption:sto_assump}, we have:
\begin{enumerate}[(i)]
\item If $ \frac{\sum_{t=1}^T \eta_t^2}{\sum_{t=1}^T \eta_t} \rightarrow 0 $ as $ T \rightarrow \infty $, then $ \Exp[\fbest^{(T)}] \rightarrow f(\bx^*) $ as $ T \rightarrow \infty $, where $\bx^*$ is an optimizer of $f$.
\item Assuming that $ \sS $ is \textbf{compact}, and letting $ \Omega $ be an upper bound on the half-squared diameter of $ \sS $:
$
\Omega \geq \max_{\bx,\by \in \sS} \frac{1}{2} \normtwo{\bx-\by}^2.
$
\begin{enumerate}[(a)]
\begin{subequations}
\item If $ \eta_t \triangleq \frac{\sqrt{2\Omega}}{\widetilde{L} \sqrt{t+1}} $, then for all $ T \geq 2 $,
\begin{equation}
\Exp[\fbest^{(T)} - f(\bx^*)] \leq \frac{2(1 + \ln(3)) \widetilde{L} \sqrt{2\Omega}}{\sqrt{T+2}}.
\end{equation}
\item If $\{\eta_t\}_{t>0}$ is a nonincreasing sequence, then for all $ T \geq 2 $,
\begin{equation}
\Exp\big[f(\widebarbx^{(T)})\big] - f(\bx^*)\leq \frac{\Omega}{T\eta_T} + \frac{1}{2T}\sum_{t=1}^{T}\eta_t \widetilde{L}^2,
\end{equation}
where $\{\widebarbx^{(t)} \triangleq\frac{1}{t} \sum_{j=1}^{t}\bx^{(j)}\}_{t>0}$ denotes the sequence of average iterates for each $t$.
Specifically, it $\eta_t \triangleq \frac{\sqrt{2\Omega}}{\widetilde{L}\sqrt{t}}$. Then the above inequality reduces to 
\begin{equation}
\Exp\big[f(\widebarbx^{(T)})\big] - f(\bx^*)\leq \frac{3 \sqrt{2\Omega}\widetilde{L}}{2\sqrt{T}}.
\end{equation}
\end{subequations}
\end{enumerate}
\end{enumerate}
\end{theoremHigh}
\begin{proof}[of Theorem~\ref{theorem:conv_sgpg}]
\textbf{(i).} For any $ t > 0 $, by the update rule of SPGD, we have
$$
\small
\begin{aligned}
\Exp\big[\normtwobig{\bx^\toptone - \bx^*}^2 \mid \bx^\toptzero\big] 
&= \Exp\big[\normtwobig{\projectS(\bx^\toptzero - \eta_t \bg^\toptzero) - \projectS(\bx^*)}^2 \mid \bx^\toptzero\big] 
\stackrel{\dag}{\leq}  \Exp\big[\normtwobig{\bx^\toptzero - \eta_t \bg^\toptzero - \bx^*}^2 \mid \bx^\toptzero\big]  \\
&= \normtwobig{\bx^\toptzero - \bx^*}^2 - 2\eta_t \innerproduct{\Exp[\bg^\toptzero \mid \bx^\toptzero], \bx^\toptzero - \bx^*} + \eta_t^2 \Exp\big[\normtwobig{\bg^\toptzero}^2 \mid \bx^\toptzero\big] \\
&\stackrel{\ddag}{\leq} \normtwobig{\bx^\toptzero - \bx^*}^2 - 2\eta_t (f(\bx^\toptzero) - f(\bx^*)) + \eta_t^2 \widetilde{L}^2,
\end{aligned}
$$
where the inequality ($\dag$) follows from the nonexpansiveness property of the  projection operator (Theorem~\ref{theorem:proj_nonexpan}), and the inequality  ($\ddag$)  follows by Assumption~\ref{assumption:sto_assump} since $f$ is convex.
By the fact that $\Exp[\Exp[\rX \mid \rY]] = \Exp[\rX]$ for any random variables $\rX$ and $\rY$, the above inequality implies that
\begin{equation}\label{equation:conv_sgpg0}
\Exp\big[\normtwobig{\bx^\toptone - \bx^*}^2\big] \leq \Exp\big[\normtwobig{\bx^\toptzero - \bx^*}^2\big] - 2\eta_t \big(\Exp\big[f(\bx^\toptzero)\big] - f(\bx^*)\big) + \eta_t^2 \widetilde{L}^2.
\end{equation}
Telescoping the sum over $ t = j, j + 1, \ldots, T $ (where $ j $ is an integer satisfying $ j \leq T $),
$$
\begin{aligned}
&\Exp\big[\normtwobig{\bx^{(T+1)} - \bx^*}^2\big] \leq \Exp\big[\normtwobig{\bx^{(j)} - \bx^*}^2\big] - 2 \sum_{t=j}^T \eta_t \big(\Exp[f(\bx^\toptzero)] - f(\bx^*)\big) + \widetilde{L}^2 \sum_{t=j}^T \eta_t^2\\
&\quad\implies\quad 
\sum_{t=j}^T \eta_t \big(\Exp[f(\bx^\toptzero)] - f(\bx^*)\big) \leq \frac{1}{2} \Exp\big[\normtwobig{\bx^{(j)} - \bx^*}^2\big] + \frac{\widetilde{L}^2}{2} \sum_{t=j}^T \eta_t^2.
\end{aligned}
$$
This implies
$$
\left( \sum_{t=j}^T \eta_t \right) \left( \min_{t\in\{j, j+1, \ldots, T\}} \Exp[f(\bx^\toptzero)] - f(\bx^*) \right) \leq \frac{1}{2} \Exp\big[\normtwobig{\bx^{(j)} - \bx^*}^2\big] + \frac{\widetilde{L}^2}{2} \sum_{t=j}^T \eta_t^2.
$$
Since $
\Exp[\fbest^{(T)}] \leq \Exp \left[ \min_{t\in\{j, j+1, \ldots, T\}} f(\bx^\toptzero)\right] \leq \min_{t\in\{j, j+1, \ldots, T\}} \Exp[f(\bx^\toptzero)]
$,
we can conclude that
\begin{equation}\label{equation:conv_sgpg}
\Exp[\fbest^{(T)} - f(\bx^*)] \leq \frac{\Exp\big[\normtwobig{\bx^{(j)} - \bx^*}^2\big] + \widetilde{L}^2 \sum_{t=j}^T \eta_t^2}{2 \sum_{t=j}^T \eta_t}.
\end{equation}
When $j=1$, we obtain
$
\Exp\big[\fbest^{(T)} - f(\bx^*)\big] \leq \frac{\normtwo{\bx^{(1)} - \bx^*}^2 + \widetilde{L}^2 \sum_{t=1}^T \eta_t^2}{2 \sum_{t=1}^T \eta_t}.
$
Therefore, if $\frac{\sum_{t=1}^{T} \eta_t^{2}}{\sum_{t=1}^{T} \eta_t} \rightarrow 0$, then $\Exp\big[\fbest^{(T)}\big] \rightarrow f(\bx^*)$ as $T \rightarrow \infty$, proving the first claim.

\paragraph{(ii).(a).} Setting $j = \left\lceil \frac{T}{2} \right\rceil$ in \eqref{equation:conv_sgpg} and using the definition of the half-squared diameter $\Omega$, we obtain
$$
\Exp\big[\fbest^{(T)}\big] - f(\bx^*) \leq \frac{\Omega + \frac{\widetilde{L}^{2}}{2} \sum_{t=\left\lceil \frac{T}{2} \right\rceil}^{T} \eta_t^{2}}{\sum_{t=\left\lceil \frac{T}{2} \right\rceil}^{T} \eta_t}
\,\,\,\,\longeqtextfour{\eta_t \triangleq \frac{\sqrt{2\Omega}}{\widetilde{L} \sqrt{t+1}}}\,\,\,\,
\frac{\widetilde{L}\sqrt{2\Omega}}{2} 
\frac{1 +  \sum_{t=\left\lceil \frac{T}{2} \right\rceil}^{T} \frac{1}{t+1}}{\sum_{t=\left\lceil \frac{T}{2} \right\rceil}^{T} \frac{1}{\sqrt{t+1}}}.
$$
The result follows by invoking Lemma~\ref{lemma:pgd_lem2}\ref{pgd_lem2_v2} with $\Phi\triangleq 1$. 

\paragraph{(ii).(b).}
Since $f$ is convex and $\{\eta_t\}_{t>0}$ is a nonincreasing sequence, \eqref{equation:conv_sgpg0} also shows that 
$$
\small
\begin{aligned}
\Exp&\big[f(\widebarbx^{(T)})\big] - f(\bx^*)
\leq 
\frac{1}{T}\sum_{t=1}^{T} \big(\Exp\big[f(\bx^\toptzero)\big] - f(\bx^*)\big) \\
&\leq 
\frac{1}{2T}  
\sum_{t=1}^{T}\left\{\frac{1}{\eta_t}\Exp\big[\normtwobig{\bx^\toptzero - \bx^*}^2\big] -\frac{1}{\eta_t}\Exp\big[\normtwobig{\bx^\toptone - \bx^*}^2\big]  + \eta_t \widetilde{L}^2 \right\}\\
&\leq 
\frac{1}{2T} 
\left\{\frac{1}{\eta_1}\Exp\big[\normtwobig{\bx^{(1)} - \bx^*}^2\big]+ \sum_{t=2}^{T}(\frac{1}{\eta_t}-\frac{1}{\eta_{t-1}})\Exp\big[\normtwobig{\bx^\toptzero - \bx^*}^2\big]  + \sum_{t=1}^{T}\eta_t \widetilde{L}^2 \right\}
\leq \frac{\Omega}{T\eta_T} + \frac{1}{2T}\sum_{t=1}^{T}\eta_t \widetilde{L}^2.
\end{aligned}
$$
Using Lemma~\ref{lemma:pgd_lem1}, we have $\sum_{t=1}^{T} \frac{1}{\sqrt{t}} \leq \int_{0}^T \frac{1}{\sqrt{x}} \,dx = 2\sqrt{T}$.
Plugging $\eta_t \triangleq \frac{\sqrt{2\Omega}}{\widetilde{L}\sqrt{t}}$ into the above inequality we have
$
\Exp\big[f(\widebarbx^{(T)})\big] - f(\bx^*)
\leq 
\frac{3\sqrt{2\Omega}\widetilde{L}}{2\sqrt{T}}.
$
This completes the proof.
\end{proof}

For strongly convex functions, we also have the following result.

\begin{theoremHigh}[SPGD for SC and ``Lipschitz" Functions]\label{theorem:conv_sgpg_sc}
Let $f: \sS \rightarrow \real$ be a proper closed and $\alpha$-strongly convex function, where $\sS$ is a closed and convex set. Let $ \{\bx^\toptzero\}_{t > 0} $ be the sequence generated by the stochastic projected subgradient method (Algorithm~\ref{alg:stopgd_gen}) with stepsizes $ \eta_t = \frac{2}{\alpha (t+1)} $ for the $t$-th iteration, and let $\left\{ \fbest^\toptzero \right\}_{t > 0}$ be the sequence of best achieved values.
Given the assumptions and vectors in Assumption~\ref{assumption:sto_assump}, we have:
\begin{enumerate}[(i)]
\item For any $ T > 0 $,
$
\Exp[\fbest^{(T)}] - f(\bx^*) \leq \frac{2 \widetilde{L}^2}{\alpha (T + 1)}
$, where $\bx^*$ is an optimizer of $f$.
\item For any $ t > 0 $, define the sequence of averages $\{\widetildebx^{(t)} \triangleq \sum_{i=1}^t \alpha_i^t \bx^{(t)}\}_{t>0}$,
where $ \alpha_i^t = \frac{2i}{t(t+1)} $ and $\sum_{i=1}^{t}\alpha_i^t=1$. Then for any $ T > 0 $,
\[
\Exp[f(\widetildebx^{(T)})] - f(\bx^*) \leq \frac{2 \widetilde{L}^2}{\alpha (T + 1)}.
\]
\end{enumerate}
\end{theoremHigh}
\begin{proof}[of Theorem~\ref{theorem:conv_sgpg_sc}]
\textbf{(i).} For any optimizer $ \bx^* $ and $ t > 0$,
$$
\small
\begin{aligned}
\Exp\big[\normtwobig{\bx^\toptone - \bx^*}^2 \mid \bx^\toptzero\big] 
&= \Exp\big[\normtwobig{\projectS(\bx^\toptzero - \eta_t \bg^\toptzero) - \projectS(\bx^*)}^2 \mid \bx^\toptzero\big] 
\stackrel{\dag}{\leq}  \Exp\big[\normtwobig{\bx^\toptzero - \eta_t \bg^\toptzero - \bx^*}^2 \mid \bx^\toptzero\big]  \\
&= \normtwobig{\bx^\toptzero - \bx^*}^2 - 2\eta_t \innerproduct{\Exp[\bg^\toptzero \mid \bx^\toptzero], \bx^\toptzero - \bx^*} + \eta_t^2 \Exp\big[\normtwobig{\bg^\toptzero}^2 \mid \bx^\toptzero\big] \\
&\stackrel{\ddag}{\leq} (1 - \alpha \eta_t)\normtwobig{\bx^\toptzero - \bx^*}^2 - 2\eta_t (f(\bx^\toptzero) - f(\bx^*)) + \eta_t^2 \widetilde{L}^2,
\end{aligned}
$$
where the inequality ($\dag$) follows from the nonexpansiveness property of the  projection operator (Theorem~\ref{theorem:proj_nonexpan}), and the inequality  ($\ddag$)  follows from the assumption in  Assumption~\ref{assumption:sto_assump} and the strong convexity of $f$ (Definition~\ref{definition:scss_func}):
$$
f(\bx^*) \geq f(\bx^\toptzero) + \innerproduct{\Exp[\bg^\toptzero \mid \bx^\toptzero], \bx^* - \bx^\toptzero} + \frac{\alpha}{2} \normtwo{\bx^\toptzero - \bx^*}^2.
$$
Rearranging terms and noting that $ \eta_t = \frac{2}{\alpha(t+1)} $, this implies
$$
f(\bx^\toptzero) - f(\bx^*) \leq \frac{\alpha(t-1)}{4} \normtwo{\bx^\toptzero - \bx^*}^2 - \frac{\alpha(t+1)}{4} \Exp\big[\normtwobig{\bx^\toptone - \bx^*}^2 \mid \bx^\toptzero\big] + \frac{1}{\alpha(t+1)} \widetilde{L}^2.
$$
Given that fact that  $\Exp[\Exp[\rX \mid \rY]] = \Exp[\rX]$ for any random variables $\rX$ and $\rY$, multiplying the above inequality  by $ t $ and taking expectation with respect to $ \bx^\toptzero $ yields that  
$$
\small
\begin{aligned}
	t\big(\Exp[f(\bx^\toptzero)] - f(\bx^*)\big)\leq & \frac{\alpha t(t-1)}{4} \Exp\big[\normtwobig{\bx^\toptzero - \bx^*}^2\big] - \frac{\alpha t(t+1)}{4} \Exp\big[\normtwobig{\bx^\toptone - \bx^*}^2\big] 
	+ \frac{t}{\alpha(t+1)} \widetilde{L}^2.
\end{aligned}
$$
Telescoping the sum over $ t = 1,2, \ldots,T$,
\begin{equation}\label{equation:conv_sgpg_sc1}
\small
\begin{aligned}
&\sum_{t=1}^{T} t (\Exp[f(\bx^\toptzero)] - f(\bx^*)) \leq  - \frac{\alpha}{4} T(T+1) \Exp\big[\normtwobig{\bx^{(T+1)} - \bx^*}^2\big] + \frac{\widetilde{L}^2}{\alpha} \sum_{t=1}^{T} \frac{t}{t+1} \leq \frac{\widetilde{L}^2 T}{\alpha}\\
&\quad\implies\quad 
\left( \sum_{t=1}^{T} t \right) (\Exp[\fbest^\toptzero] - f(\bx^*)) \leq \frac{\widetilde{L}^2 T}{\alpha}
\quad\implies\quad 
\Exp[\fbest^\toptzero] - f(\bx^*) \leq \frac{2 \widetilde{L}^2}{\alpha(T+1)},
\end{aligned}
\end{equation}
where the last two inequalities follows from the fact that $ \Exp[f(\bx^\toptzero)] \geq \Exp[\fbest^\toptzero] $ for all $ t \in \{1,2, \ldots,T\} $ and $ \sum_{t=1}^{T} t = \frac{T(T+1)}{2} $.

\paragraph{(ii).}  Using Jensen's inequality (Theorem~\ref{theorem:jensens_ineq}) and dividing \eqref{equation:conv_sgpg_sc1} by $ \frac{T(T+1)}{2} $, we have
$$
\Exp[f(\widetildebx^{(T)})] - f(\bx^*) = \Exp\bigg[f\bigg(\sum_{t=1}^{T} \alpha_t^T \bx^\toptzero\bigg)\bigg] - f(\bx^*) 
\leq \sum_{t=1}^T \alpha_t^T \big(\Exp[f(\bx^\toptzero)] - f(\bx^*)\big) \leq \frac{2 \widetilde{L}^2}{\alpha(T+1)}.
$$
This completes the proof.
\end{proof}

\subsection{Results of Stochastic Gradient Descent}

Recall that in many scenarios, the function $f(\bx)$ is denoted by a dataset $\mathcalD=\{\bs_1, \bs_2, \ldots, \bs_D\}$, allowing us to express both  $f(\bx)$ and  its gradient $\nabla f(\bx)$ as:
\begin{equation}
	f(\bx)\triangleq f(\mathcalD; \bx) = \frac{1}{D} \sum_{d=1}^{D} f(\bs_d; \bx)\triangleq\frac{1}{D} \sum_{d=1}^{D} f_d(\bx)
\quad \text{and}\quad 
\nabla f(\bx) \triangleq \frac{1}{D}  \sum_{d=1}^{D}\nabla f_d(\bx),
\end{equation} 
respectively.
When we follow the negative gradient direction of a single sample or a batch of samples iteratively, we obtain an estimate of the direction, a method known as \textit{stochastic gradient descent} (SGD) \citep{robbins1951stochastic}. SGD can be categorized into two types:
\begin{itemize}
\item \textbf{The strict SGD:} Computes the gradient using only one randomly selected data point per iteration, i.e., the the gradient is approximated as $\nabla f(\bx) \approx \nabla f_d(\bx)$.
\item \textbf{The mini-batch SGD:} Represents a compromise between GD and strict SGD, computing the gradient over a subset (mini-batch) of the dataset, i.e., the the gradient is approximated as  $\nabla f(\bx) \approx \frac{1}{\abs{\mathcalD}}  \sum_{d\in\mathcalD}\nabla f_d(\bx)$.
\end{itemize}
SGD is particularly useful when dealing with a large number of training entries (i.e., the data used for updating the model, while the data used for final evaluation is called the \textit{test entries or test data}).
In such cases, gradients from different input samples may cancel out, resulting in small updates. Within the SGD framework, the objective function is stochastic, consisting of a sum of subfunctions evaluated at various subsamples of the data. However, both GD and  SGD have the limitation of potentially getting trapped in local minima \citep{rutishauser1959theory}.

We then prove the convergence of  strict SGD, where for each iteration $t$, an index $s_t$ is selected ($s_t \in \mathcalD = \{1, 2, \ldots, D\}$) under the following assumption. The analysis for mini-batch SGD follows similarly.
\begin{assumption}[SGD Assumption]\label{assumption:sgd_assump}
Let $f:\real^n\rightarrow \real$ be a $\alpha$-strongly, and $\beta$-smooth convex function. We assume further that the stochastic second moment is uniformly bounded. That is, there exists a constant $B>0$ such that for any $\bx\in\real^n$ and $s_t$, it follows that 
$$
\Exp_{s_t}\big[ \normtwo{\nabla f_{s_t}(\bx)}^2 \big] \leq B^2.
$$
\end{assumption}

\begin{theoremHigh}[SGD for SC, SS, and ``Lipschitz" Functions: $\mathcalO({1}/{T})$]\label{theorem:conv_sgpg_sc_ss}
Let $f: \real^n \rightarrow \real$ be a proper closed, $\alpha$-strongly, and $\beta$-smooth convex function, where $\sS$ is a closed and convex set. Let $ \{\bx^\toptzero\}_{t > 0} $ be the sequence generated by the stochastic  gradient descent (SGD) method. 
Consider the assumptions and vectors in Assumption~\ref{assumption:sgd_assump}. 
Denoting $e_t \triangleq\normtwo{\bx^\toptzero - \bx^*}$, then for any $T>1$,
\begin{itemize}
\item If the stepsizes are constant  $\eta_t \triangleq \eta \in  (0, \frac{1}{2\alpha}) $ for all iterations,
$$
\Exp\big[f(\bx^{(T+1)}) - f(\bx^*)\big] \leq \frac{\beta}{2} \Exp[e_{T+1}^2] \leq \frac{\beta}{2} \left[ (1 - 2\eta\alpha)^T e_1^2 + \frac{\eta B^2}{2\alpha} \right],
$$
where the expectation is taken over all the sample indices $s_1, s_2, \ldots,s_T$.
\item If we use a decreasing stepsize
$
\eta_t = \frac{\sigma}{t + \gamma},
$
where $\sigma > \frac{1}{2\alpha}$, $\gamma > 0$, such that $\eta_1 \leq \frac{1}{2\alpha}$, then,
$$
\Exp\big[f(\bx^{(T)}) - f(\bx^*)\big] \leq \frac{\beta}{2} \Exp[e_T^2] \leq \frac{\beta}{2} \frac{\zeta}{\gamma + T},
$$
where  $\zeta  \triangleq \max \left\{ \frac{\sigma^2 B^2}{2\sigma\alpha - 1}, (\gamma + 1) e_1^2 \right\}$
and $e_1 \triangleq\normtwo{\bx^\topone - \bx^*}$.
\end{itemize}

\end{theoremHigh}
\begin{proof}[of Theorem~\ref{theorem:conv_sgpg_sc_ss}]
\textbf{(i).}
By the update formula of the stochastic gradient descent method,
$$
\begin{aligned}
e_{t+1}^2 &\triangleq\normtwo{\bx^\toptone - \bx^*}^2 \\
&= \normtwo{\bx^\toptzero - \eta \nabla f_{s_t}(\bx^\toptzero) - \bx^*}^2\\
&= e_{t}^2  - 2\eta \langle \nabla f_{s_t}(\bx^\toptzero), \bx^\toptzero - \bx^* \rangle + \eta^2 \normtwo{\nabla f_{s_t}(\bx^\toptzero)}^2 \\
&= \normtwo{\bx^\toptzero - \bx^*}^2 - 2\eta \langle \nabla f_{s_t}(\bx^\toptzero), \bx^\toptzero - \bx^* \rangle + \eta^2 \normtwo{\nabla f_{s_t}(\bx^\toptzero)}^2,
\end{aligned}
$$
where the second term $\langle \nabla f_{s_t}(\bx^\toptzero), \bx^\toptzero - \bx^* \rangle$ is challenging to handle because both $s_t$ and $\bx^\toptzero$ are random variables.
By the property of conditional expectation $\Exp[\rX] = \Exp[\Exp[\rX\mid \rY]]$ for any random variables $\rX$ and $\rY$, we have
\begin{equation}\label{equation:sgd_gd_xtxstar_exp}
\begin{aligned}
&\Exp_{s_1, s_2, \ldots, s_t}[\langle \nabla f_{s_t}(\bx^\toptzero), \bx^\toptzero - \bx^* \rangle]\\
&= \Exp_{s_1, s_2, \ldots, s_{t-1}}\big[\Exp_{s_t}[\langle \nabla f_{s_t}(\bx^\toptzero), \bx^\toptzero - \bx^* \rangle \mid s_1, s_2 \ldots, s_{t-1}]\big]\\
&= \Exp_{s_1, s_2, \ldots, s_{t-1}}
\big[\langle \Exp_{s_t}[\nabla f_{s_t}(\bx^\toptzero) \mid s_1, s_2, \ldots, s_{t-1}], \bx^\toptzero - \bx^* \rangle\big]\\
&= \Exp_{s_1, s_2, \ldots, s_{t-1}}\big[\langle \nabla f(\bx^\toptzero), \bx^\toptzero - \bx^* \rangle\big]
= \Exp_{s_1, s_2, \ldots, s_t}\big[\langle \nabla f(\bx^\toptzero), \bx^\toptzero - \bx^* \rangle\big],
\end{aligned}
\end{equation}
where we use the fact that $\bx^\toptzero$ depends only on $s_1, s_2, \ldots, s_{t-1}$. Therefore, when fixing $s_1, s_2, \ldots, s_{t-1}$, $\bx^\toptzero$ becomes a constant variable. By taking expectations, we replace $\nabla f_{s_t}(\bx^\toptzero)$ with $\nabla f(\bx^\toptzero)$.
By the characterization theorem of SC (Theorem~\ref{theorem:charac_stronconv}) and the fact that $\nabla f(\bx^*)=\bzero$, 
$$
\langle \nabla f(\bx^\toptzero), \bx^\toptzero -\bx^* \rangle = \langle \nabla f(\bx^\toptzero) - \nabla f(\bx^*), \bx^\toptzero -\bx^* \rangle \geq \alpha \normtwo{\bx^\toptzero - \bx^*}^2.
$$
Combining the above results, using the uniform boundedness of the second moment of the stochastic gradient and the bound on $\langle \nabla f(\bx^\toptzero), \bx^\toptzero -\bx^* \rangle$, we obtain
\begin{equation}\label{equation:conv_sgpg_sc_ss_pe1}
\Exp_{s_1, s_2, \ldots, s_t}[e_{t+1}^2] \leq (1 - 2\eta\alpha) \Exp_{s_1, s_2, \ldots, s_t}\left[e_t^2\right] + \eta^2 B^2. 
\end{equation}
Applying induction on $t=1,2,\ldots,T+1$, we obtain
$$
\Exp_{s_1, s_2, \ldots, s_T}[e_{T+1}^2] \leq (1 - 2\eta\alpha)^T e_1^2 + \sum_{t=1}^{T-1} (1 - 2\eta\alpha)^t \eta^2 B^2. 
$$
Given that $ \eta \in  (0, \frac{1}{2\alpha}) $ (i.e., $2\eta\alpha<1$), it follows that for any $T>1$
$$
\sum_{t=1}^{T-1} (1 - 2\eta\alpha)^t < \sum_{t=1}^{\infty} (1 - 2\eta\alpha)^t = \frac{1}{2\eta\alpha}-1,
$$
whence we have 
$$
\Exp_{s_1, s_2, \ldots, s_T}[e_{T+1}^2] \leq (1 - 2\eta\alpha)^T e_1^2 + \frac{\eta B^2}{2\alpha}. 
$$
By the smoothness of $f$ and the fact that $\nabla f(\bx^*) = \bzero$,
$$
f(\bx^{(T+1)}) - f(\bx^*) \leq \innerproduct{\nabla f(\bx^*), \bx^{(T+1)} -\bx^*} + \frac{\beta}{2} \normtwobig{\bx^{(T+1)} - \bx^*}^2
=\frac{\beta}{2}e_{T+1}^2.
$$
Taking expectations, we obtain the desired result.
\paragraph{(ii).} If using a decreasing stepsize $\eta_t = \frac{\sigma}{t + \gamma}$, by \eqref{equation:conv_sgpg_sc_ss_pe1}, for any $t>0$:
\begin{equation}\label{equation:conv_sgpg_sc_ss_pe2}
	\Exp_{s_1, s_2, \ldots, s_t}[e_{t+1}^2] \leq (1 - 2\eta_t \alpha) \Exp_{s_1, s_2, \ldots, s_t}[e_t^2] + \eta_t^2 B^2.
\end{equation}
We prove the result by induction. When $T = 1$, by the definition of $\zeta$, the claim holds trivially. Assume that the formula holds for $t$. For simplicity, define $\widetildet \triangleq \gamma + t$, then $\eta_t = \frac{\sigma}{\widetildet}$. By the induction hypothesis and \eqref{equation:conv_sgpg_sc_ss_pe2},
$$
\begin{aligned}
\Exp[e_{t+1}^2] 
&\leq \left(1 - \frac{2\sigma\alpha}{\widetildet}\right) \frac{\zeta}{\widetildet} + \frac{\sigma^2 B^2}{\widetildet^2}
= \frac{\zeta }{\widetildet}  + \frac{\sigma^2 B^2 - 2\sigma\alpha\zeta}{\widetildet^2}
\leq \frac{\zeta}{\widetildet + 1}
\end{aligned}
$$
where the last inequality follows from the definition of $\zeta$ such that $\sigma^2 B^2 - 2\sigma\alpha\zeta<0$. Therefore, the claim also holds for $t + 1$. This completes the proof.
\end{proof}

It can be observed that with a fixed stepsize, the algorithm does not guarantee convergence because the right-hand side lacks a term that decreases with $T$. 
However, by employing a decreasing stepsize, the convergence rate can achieve $\mathcalO\left(\frac{1}{T}\right)$.

\section{Variance/Noise Reduction Techniques}

Consider a constant stepsize $\eta_t\triangleq\eta$ for all iterations. We compare the pure gradient descent method with the stochastic gradient descent method.
Under the conditions specified in Assumption~\ref{assumption:sgd_assump}, by the characterization theorem of SC (Theorem~\ref{theorem:charac_stronconv})  and SS Property-O (Theorem~\ref{theorem:equi_gradsch_smoo}), for the gradient descent method, we have:
\begin{align*}
e_{t+1}^2 &\triangleq \normtwobig{\bx^\toptone - \bx^*}^2 = \normtwo{\bx^\toptzero - \eta \nabla f(\bx^\toptzero) - \bx^*}^2 \\
&= e_t^2 - 2\eta \innerproduct{\nabla f(\bx^\toptzero), \bx^\toptzero - \bx^*} + \eta^2 \normtwobig{\nabla f(\bx^\toptzero)}^2 \\
&\leq (1 - 2\eta \alpha) e_t^2 + \eta^2 \normtwobig{\nabla f(\bx^\toptzero)}^2  \\
&\leq (1 - 2\eta \alpha + \eta^2 \beta^2) e_t^2. 
\end{align*}
While for the stochastic gradient descent method, by leveraging properties of conditional expectations:
\begin{align*}
\Exp[e_{t+1}^2] 
&= \Exp\big[\normtwobig{\bx^\toptone - \bx^*}^2\big] = \Exp\big[\normtwobig{\bx^\toptzero - \eta \nabla f_{s_t}(\bx^\toptzero) - \bx^*}^2\big] \\
&= \Exp[e_t^2] - 2\eta \Exp\big[\langle \nabla f_{s_t}(\bx^\toptzero), \bx^\toptzero - \bx^* \rangle\big] + \eta^2 \Exp\big[\normtwobig{\nabla f_{s_t}(\bx^\toptzero)}^2\big] \\
&\stackrel{\dag}{=} \Exp[e_t^2] - 2\eta \Exp\big[\langle \nabla f(\bx^\toptzero), \bx^\toptzero - \bx^* \rangle\big] + \eta^2  \Exp\big[\normtwobig{\nabla f_{s_t}(\bx^\toptzero) - \nabla f(\bx^\toptzero) + \nabla f(\bx^\toptzero) }^2\big] \\
&\stackrel{\ddag}{\leq} (1 - 2\eta \alpha) \Exp[e_t^2] + \eta^2 \Exp\big[\normtwobig{\nabla f_{s_t}(\bx^\toptzero) - \nabla f(\bx^\toptzero) + \nabla f(\bx^\toptzero) }^2\big] \\
&\stackrel{*}{\leq} \underbrace{(1 - 2\eta \alpha + \eta^2 \beta^2) \Exp[e_t^2]}_{\triangleq A} + \underbrace{\eta^2 \Exp\big[\normtwobig{\nabla f_{s_t}(\bx^\toptzero) - \nabla f(\bx^\toptzero)}^2\big]}_{\triangleq B}.
\end{align*}
where the equality ($\dag$) follows from \eqref{equation:sgd_gd_xtxstar_exp}, the inequality ($\ddag$) follows from the the characterization theorem of SC (Theorem~\ref{theorem:charac_stronconv}), and the inequality ($*$) follows from the SS Property-O (Theorem~\ref{theorem:equi_gradsch_smoo}).

Therefore, the main difference between the two methods  lies in the term $B$, which represents the variance in gradient estimation.  
This results in the convergence rate of the stochastic gradient algorithm being approximately $\mathcalO\left(\frac{1}{T}\right)$ (refer to  Theorem~\ref{theorem:gd_sc_ss} for GD and Theorem~\ref{theorem:conv_sgpg_sc_ss} for SGD). 
However, in many machine learning applications, the actual convergence rate of the stochastic gradient algorithm may be faster due to lower precision requirements for solutions and relatively smaller variance in early stages, i.e., $B \ll A$, leading to an approximate linear convergence rate initially; as iterations increase, variance grows, eventually resulting in a $\mathcalO\left(\frac{1}{T}\right)$ convergence rate. 
Therefore, to achieve faster convergence rates, our main objective is to reduce the variance term $B$. 
Below are three algorithms designed to reduce variance in SGD:
\begin{itemize}
	\item SAG (stochastic average gradient) \citep{schmidt2017minimizing}.
	\item SAGA \citep{defazio2014saga}.
	\item SVRG (stochastic variance reduced gradient) \citep{johnson2013accelerating}.
\end{itemize}
Compared to traditional stochastic gradient methods, these algorithms utilize previously computed information to decrease variance, ultimately achieving a linear convergence rate \citep{bottou2018optimization, liu2020optimization}. Below is a comparison among the update rules for SAG, SAGA, and SVRG, where notations will be clarified in the following sections:
\begin{subequations}
\begin{align}
\textbf{(SAG)}: \quad \bx^\toptone &\leftarrow \bx^\toptzero - \eta_t \bigg( \frac{1}{D} \big(\nabla f_{s_t}(\bx^\toptzero) - \bg_{s_t}^{(t-1)}\big) + \frac{1}{D} \sum_{i=1}^D \bg_i^{(t-1)} \bigg);\\
\textbf{(SAGA)}: \quad 	\bx^\toptone &\leftarrow \bx^\toptzero - \eta_t \bigg( \nabla f_{s_t}(\bx^\toptzero) - \bg_{s_t}^{(t-1)} + \frac{1}{D} \sum_{i=1}^D \bg_i^{(t-1)} \bigg);\\
\textbf{(SVRG)}: \quad 	\bx^\toptone&\leftarrow \bx^\toptzero - \eta_t \bigg( \nabla f_{s_t}(\bx^\toptzero) - \nabla f_{s_t}(\widetildebx^{(j)}) + \frac{1}{D} \sum_{i=1}^D \nabla f_i(\widetildebx^{(j)})\bigg).
\end{align}
\end{subequations}
In the following, we only discuss the results for strict SGD; the results for mini-batch SGD follow analogously.
\subsection{SAG Algorithm}

In the stochastic gradient descent  method, each iteration only uses the random gradient of the current point, discarding all previously computed gradients immediately after use.
However, as iterations approach convergence, the random gradient from the previous step can serve as a good estimate for the gradient at the current iteration point. The \textit{stochastic average gradient (SAG)} algorithm leverages this concept by recording all previously calculated random gradients and averaging them with the newly computed gradient to form an estimate for the next step.
Specifically, the SAG algorithm reserves memory space to store $D = \abs{\mathcalD}$ random gradients:
$$
[\bg_1^\toptzero, \bg_2^\toptzero, \ldots, \bg_D^\toptzero],
$$
which are used to record the latest random gradients for each sample. At the $t$-th update, if the sampled index is $s_t$, 
then after calculating the new random gradient,  $\bg_{s_t}^\toptzero$ is updated to this new value, while all other  $\bg_i^\toptzero$ remain unchanged from their previous values. 
Each time, the SAG algorithm updates the gradient direction using the average of all $\bg_i^\toptzero$. 
Formally, its mathematical iterative formula is given by:
\begin{equation}
	\bx^\toptone \leftarrow \bx^\toptzero - \frac{\eta_t}{D} \sum_{i=1}^D \bg_i^\toptzero, 
	\quad
	\text{ where }\bg_i^\toptzero = 
	\begin{cases} 
		\nabla f_{s_t}(\bx^\toptzero), & i = s_t, \\
		\bg_i^{(t-1)}, & \text{otherwise}.
	\end{cases} 
\end{equation}
That is, only update the random gradient for the $s_t$-th sample, with others directly carrying over from the previous step. Therefore, the SAG iterative formula can also be expressed as:
\begin{equation}\label{equation:sag}
\textbf{(SAG)}: \quad \bx^\toptone \leftarrow \bx^\toptzero - \eta_t \bigg( \frac{1}{D} \big(\nabla f_{s_t}(\bx^\toptzero) - \bg_{s_t}^{(t-1)}\big) + \frac{1}{D} \sum_{i=1}^D \bg_i^{(t-1)} \bigg). 
\end{equation}
For the SAG algorithm, the initial values of $\{\bg_i^\toptzero\}$ can be set to zero vectors or random gradient vectors with a mean of zero.
We present the convergence analysis of the SAG algorithm:
\begin{theoremHigh}[Convergence of the SAG Algorithm]
Under Assumption~\ref{assumption:sgd_assump}, taking a constant stepsize $\eta_t \triangleq\eta\triangleq\frac{1}{16\beta}$, and initializing   $\bg_i^\toptzero$, $i\in\{1,2,\ldots,D\}$ as zero vectors,  for any $t>0$, we have
\begin{equation}
	\Exp\big[f(\bx^\toptzero)\big] - f(\bx^*) \leq \left(1 - \min\left\{\frac{\alpha}{16\beta}, \frac{1}{8D}\right\}\right)^t C_0,
\end{equation}
where $C_0 \triangleq f(\bx^{(1)}) - f(\bx^*) + \frac{4\beta}{D}\normtwo{\bx^{(1)} - \bx^*}^2 + \frac{\sigma^2}{16\beta}$, $\sigma^2 = \frac{1}{D}\sum_{i=1}^{D} \normtwo{\nabla f_i(\bx^*)}^2$, and $\bx^*$ is an optimal point.
Additionally, if the function is not strongly convex under  Assumption~\ref{assumption:sgd_assump}, then 
\begin{equation}
\Exp\left[f\Big(\frac{1}{t} \sum_{i=1}^{t}\bx^{(i)}\Big)\right] - f(\bx^*) \leq  \frac{32D}{t}C_0.
\end{equation}
\end{theoremHigh}
The detailed proof of this theorem can be found in \citet{schmidt2017minimizing}, which involves finding a Lyapunov function for a nonlinear stochastic dynamical system.
The theorem demonstrates that the SAG algorithm indeed achieves a linear convergence rate. However, a limitation of the SAG algorithm is its requirement to store $D$ gradient vectors, which can be prohibitively expensive when the sample size $D$  is large. 
Consequently, the SAG algorithm is rarely used in practice, especially in the era of large language models. Nonetheless, its core idea has inspired the development of many practical algorithms.

\subsection{SAGA Algorithm}

The \textit{SAGA} algorithm builds upon the SAG method. One known limitation of the SAG algorithm is that the conditional expectation of the random gradient at each step does not equal the true gradient. To address this, the SAGA algorithm employs an unbiased gradient vector for its update direction. Its iterative formula is given by:
\begin{equation}
\textbf{(SAGA)}: \quad 	\bx^\toptone \leftarrow \bx^\toptzero - \eta_t \bigg( \nabla f_{s_t}(\bx^\toptzero) - \bg_{s_t}^{(t-1)} + \frac{1}{D} \sum_{i=1}^D \bg_i^{(t-1)} \bigg). 
\end{equation}
It can be demonstrated that the gradient direction used in each iteration is unbiased, i.e.,
\begin{equation}
	\Exp \left[ \nabla f_{s_t}(\bx^\toptzero) - \bg_{s_t}^{(t-1)} + \frac{1}{D} \sum_{i=1}^D \bg_i^{(t-1)} \mid \bx^\toptzero \right] = \nabla f(\bx^\toptzero).
\end{equation}
TSimilar to SAG, the SAGA algorithm also achieves a linear convergence rate. Specifically, we have the following result:
\begin{theoremHigh}[Convergence of the SAGA Algorithm]
Under Assumption~\ref{assumption:sgd_assump}, taking a constant stepsize $\eta_t \triangleq\eta\triangleq \frac{1}{2(\alpha D + \beta)}$. Letting $e_t \triangleq \normtwo{\bx^\toptzero - \bx^*}$,  for any $t > 0$, we have
\begin{equation}
\Exp[e_t^2] \leq \left(1 - \frac{\alpha}{2(\alpha D + \beta)}\right)^t 
\left( e_1^2 + \frac{D\big(f(\bx^{(1)}) - f(\bx^*) - \langle \nabla f(\bx^*), \bx^{(1)} - \bx^* \rangle\big)}{\alpha D + \beta} \right). 
\end{equation}
\end{theoremHigh}
The proof of this theorem also involves finding a Lyapunov function for a nonlinear stochastic dynamical system, details of which can be found in \citet{defazio2014saga}. For brevity, we will not delve into these specifics here.

\subsection{SVRG Algorithm}

Unlike the SAG and SAGA algorithms, the \textit{stochastic variance reduced gradient (SVRG)} algorithm reduces variance by periodically caching the full gradient. 
Specifically, in the stochastic gradient descent method, a checkpoint is established every $m$ iterations to compute the full gradient. During the subsequent $m$ iterations, this full gradient serves as a reference point to reduce variance. 
Let $\widetildebx^{(j)}$ be the $j$-th checkpoint, then we need to compute the full gradient at point $\widetildebx^{(j)}$
\begin{equation}
\nabla f(\widetildebx^{(j)}) = \frac{1}{D} \sum_{i=1}^D \nabla f_i(\widetildebx^{(j)}), 
\quad \text{ with }\quad 
\widetildebx^{(j)} \triangleq \frac{1}{m} \sum_{t=j-m}^{j-1} \bx^{(t)}
\end{equation}
And in the subsequent iterations (i.e., $t\in\{j, j+1, \ldots, j+m-1\}$), we can use the direction $\bp^\toptzero$ as the update direction:
\begin{equation}
\bp^\toptzero \triangleq \nabla f_{s_t}(\bx^\toptzero) - \big(\nabla f_{s_t}(\widetildebx^{(j)}) - \nabla f(\widetildebx^{(j)})\big), \quad t\in\{j, j+1, \ldots, j+m-1\},
\end{equation}
where $s_t \in \{1, 2, \ldots, D\}$ is a randomly selected sample. To see why we define the vector $\bp^\toptzero$, note that given $s_1, s_2, \ldots, s_{t-1}$, $\bx^\toptzero$ and $\widetildebx^{(j)}$ are both fixed values, whence we have 
\begin{equation}
\begin{aligned}
&\Exp[\bp^\toptzero \mid s_1, s_2, \ldots, s_{t-1}] \\
=& \Exp\big[\nabla f_{s_t}(\bx^\toptzero) \mid \bx^\toptzero\big] - \Exp\big[\nabla f_{s_t}(\widetildebx^{(j)}) - \nabla f(\widetildebx^{(j)}) \mid s_1, s_2, \ldots, s_{t-1}\big] \\
=& \nabla f(\bx^\toptzero) - \bzero = \nabla f(\bx^\toptzero).
\end{aligned}
\end{equation}
Thus, $\bp^\toptzero$ is an unbiased estimate of $\nabla f(\bx^\toptzero)$. Therefore, if we use $\nabla f_{s_t}(\widetildebx^{(j)})$ to estimate $\nabla f(\widetildebx^{(j)})$,  the term $\nabla f_{s_t}(\widetildebx^{(j)}) - \nabla f(\widetildebx^{(j)})$ can be regarded as the error of the gradient estimation. At each step of the stochastic gradient iteration, this term is used to correct $\nabla f_{s_t}(\bx^\toptzero)$.
Consequently, the SVRG update becomes:
\begin{equation}\label{equation:svrg_update}
	\textbf{(SVRG)}: \quad 
\begin{aligned}
\bx^\toptone 
&\leftarrow\bx^\toptzero - \eta_t\bp^\toptzero\\
&= \bx^\toptzero - \eta_t \bigg( \nabla f_{s_t}(\bx^\toptzero) - \big(\nabla f_{s_t}(\widetildebx^{(j)}) - \nabla f(\widetildebx^{(j)})\big)\bigg).
\end{aligned} 
\end{equation}
Unlike the SAG  and  SAGA algorithms, the SVRG algorithm does not require storage space for recording $D$ gradient vectors, but it requires calculating the full gradient every $m$ steps and an additional gradient $\nabla f_{s_t}(\widetildebx^{(j)})$ at each iteration.

To analyze the variance, an additional assumption is needed here:
\begin{equation}
\normtwo{\nabla f_i(\bx) - \nabla f_i(\by)} \leq \beta\normtwo{\bx - \by}, \quad i = 1, 2, \ldots, D,
\end{equation}
which means that each sub-function of $f(\bx)$ is  $\beta$-smooth. 
Let $\bx^*$ denote the minimum point of $f(\bx)$ and $e_t \triangleq \normtwo{\bx^\toptzero - \bx^*}$ represent  the distance between $\bx^\toptzero$ and $\bx^*$. Then,
\begin{equation}\label{equation:svag_ineqclo}
\begin{aligned}
& \Exp\big[\normtwobig{\bp^\toptzero}^2\big] = \Exp\big[\normtwobig{\nabla f_{s_t}(\bx^\toptzero) - \big(\nabla f_{s_t}(\widetildebx^{(j)}) - \nabla f(\widetildebx^{(j)})\big)}^2\big] \\
= & \Exp\big[\normtwobig{\nabla f_{s_t}(\bx^\toptzero) - \nabla f_{s_t}(\widetildebx^{(j)}) + \nabla f(\widetildebx^{(j)}) + \nabla f_{s_t}(\bx^*) - \nabla f_{s_t}(\bx^*)}^2\big] \\
\leq & 2 \Exp\big[\normtwobig{\nabla f_{s_t}(\bx^\toptzero) - \nabla f_{s_t}(\bx^*)}^2\big] + 2 \Exp\big[\normtwobig{\nabla f_{s_t}(\widetildebx^{(j)}) - \nabla f(\widetildebx^{(j)}) - \nabla f_{s_t}(\bx^*)}^2\big] \\
\leq & 2\beta^2 \Exp[e_t^2] + 2 \Exp\big[\normtwobig{\nabla f_{s_t}(\widetildebx^{(j)}) - \nabla f_{s_t}(\bx^*)}^2\big] \\
\leq & 2\beta^2 \Exp[e_t^2] + 2\beta^2 \Exp\big[\normtwobig{\widetildebx^{(j)} - \bx^*}^2\big],
\end{aligned}
\end{equation}
where the first inequality follows from the fact that $\normtwo{\ba + \bb}^2 \leq 2\normtwo{\ba}^2 + 2\normtwo{\bb}^2$ for any vectors $\ba,\bb$.

The inequality \eqref{equation:svag_ineqclo} shows that  if $\bx^\toptzero$ and $\widetildebx^{(j)}$ are very close to $\bx^*$, the variance of the gradient estimation is  small. Obviously, frequently updating $\widetildebx^{(j)}$ can decrease the variance, but this also increases the number of times the full gradient must be calculated.
Below is the analysis of the convergence of the SVRG algorithm.
\begin{theoremHigh}[Convergence of the SVRG Algorithm]\label{theorem:conv_svrg}
Let each sub-function $f_i(\bx)$ be a differentiable convex and  $\beta$-smooth function for $i\in\{1,2,\ldots,D\}$, and let the function $f(\bx)=\sum_{i=1}^{D} f_i(\bx)$ be $\alpha$-strongly convex. Consider the SGD algorithm with an SVRG update using a constant stepsize $\eta_t\triangleq\eta \in \left(0, \frac{1}{2\beta}\right]$, and suppose $m$ is sufficiently large such that
\begin{equation}
\rho = \frac{1}{\alpha \eta (1 - 2\beta \eta) m} + \frac{2\beta \eta}{1 - 2\beta \eta} < 1.
\end{equation}
Let  $\bx^*$ be any optimal point of $f = \tfrac{1}{D} \sum_{i=1}^{D} f_i$.
Then, for $j>m$:
\begin{equation}
\Exp \big[f(\widetildebx^{(j)}) - f(\bx^*)\big] \leq \rho \Exp \big[f(\widetildebx^{(j-m)}) - f(\bx^*)\big],
\end{equation}
i.e., a linear rate of convergence.
\end{theoremHigh}
\begin{proof}[of Theorem~\ref{theorem:conv_svrg}]
Let $e_t \triangleq \normtwo{\bx^\toptzero - \bx^*}$ be the distance between $\bx^\toptzero$ and $\bx^*$, and let $\widetildebx^{(j)}$ be the $j$-th checkpoint ($j$ is a scalar multiple of $m$). By the update rule of SVRG and the convexity of $f$, for $t\in\{j, j+1, \ldots, j+m-1\}$,
\begin{equation}\label{equation:conv_svrg_et1et}
\begin{aligned}
\Exp[e_{t+1}^2] 
&= \Exp\big[\normtwobig{\bx^\toptone - \bx^*}^2\big] = \Exp\big[\normtwobig{\bx^\toptzero - \eta \bp^\toptzero - \bx^*}^2\big] \\
&= \Exp[e_t^2] - 2\eta \Exp\big[\langle \bp^\toptzero, \bx^\toptzero - \bx^* \rangle\big] + \eta^2 \Exp\big[\normtwobig{\bp^\toptzero}^2\big] \\
&= \Exp[e_t^2] - 2\eta \Exp\big[\langle \nabla f(\bx^\toptzero), \bx^\toptzero - \bx^* \rangle\big] + \eta^2 \Exp\big[\normtwobig{\bp^\toptzero}^2\big] \\
&\leq \Exp[e_t^2] - 2\eta \Exp\big[(f(\bx^\toptzero) - f(\bx^*))\big] + \eta^2 \Exp\big[\normtwobig{\bp^\toptzero}^2\big].
\end{aligned}
\end{equation}
Construct the following auxiliary function for each index $i\in\{1,2,\ldots, D\}$:
\begin{equation}
	\phi_i(\bx) \triangleq f_i(\bx) - f_i(\bx^*) - \innerproduct{\nabla f_i(\bx^*), \bx - \bx^*},
\end{equation}
Notice that $\phi_i(\bx)$ is also a convex  and $\beta$-smooth function. 
Since $\nabla \phi_i(\bx) - \nabla \phi_i(\bx^*) = \nabla f_i(\bx) -  \nabla f_i(\bx^*) $ (noting that $\bx^*$ is an optimal point of $\phi_i$ or $f$, but may not be an optimal point of $f_i$), invoking  Theorem~\ref{theorem:smoo_prop2_bound} on $\phi_i$, we have
\begin{equation}
\frac{1}{2\beta}\normtwo{\nabla f_i(\bx) - \nabla f_i(\bx^*)}^2 = \frac{1}{2\beta} \normtwo{\nabla \phi_i(\bx)}^2 \leq  \phi_i(\bx) - \phi_i(\bx^*).
\end{equation}
Taking the sum over $i=1,2,\ldots, D$, and noting that $\nabla f(\bx^*) = \frac{1}{D} \sum_{i=1}^{D} \nabla f_i(\bx^*)= \bzero$, we obtain
\begin{equation}\label{equation:svrg_smoo_sin}
\frac{1}{D} \sum_{i=1}^D \normtwo{\nabla f_i(\bx) - \nabla f_i(\bx^*)}^2 \leq 2\beta \left[ f(\bx) - f(\bx^*) \right], \quad \forall \bx\in\real^n.
\end{equation}
Using the same procedure as \eqref{equation:svag_ineqclo}, for $t\in\{j, j+1, \ldots, j+m-1\}$, we have the upper bound expression for the second moment of $\bp^\toptzero$:
\begin{equation}\label{equation:svrg_decom}
\Exp\big[\normtwobig{\bp^\toptzero}^2\big] \leq 2 \Exp\big[\normtwobig{\nabla f_{s_t}(\bx^\toptzero) - \nabla f_{s_t}(\bx^*)}^2\big] + 2 \Exp\big[\normtwobig{\nabla f_{s_t}(\widetildebx^{(j)}) - \nabla f_{s_t}(\bx^*)}^2\big].
\end{equation}
For the first term on the right-hand side, we have
\begin{equation}
\begin{aligned}
&\Exp\big[\normtwobig{\nabla f_{s_t}(\bx^\toptzero) - \nabla f_{s_t}(\bx^*)}^2\big] 	
= \Exp\left[ \Exp\big[\normtwobig{\nabla f_{s_t}(\bx^\toptzero) - \nabla f_{s_t}(\bx^*)}^2 \mid s_1, s_2, \ldots, s_{t-1}\big] \right] \\
=& \Exp\left[ \frac{1}{D} \sum_{i=1}^D \normtwo{\nabla f_i(\bx^\toptzero) - \nabla f_i(\bx^*)}^2 \right]
\leq 2\beta \Exp\big[f(\bx^\toptzero) - f(\bx^*)\big],
\end{aligned}
\end{equation}
where the first equality follows from the property of conditional expectation $\Exp[\rX] = \Exp[\Exp[\rX\mid \rY]]$, 
the second equality follows from the definition of expectation, and the last inequality uses inequality \eqref{equation:svrg_smoo_sin}. Similarly, for the second term on the right-hand side of \eqref{equation:svrg_decom}, invoking  Theorem~\ref{theorem:smoo_prop2_bound} on $f_{s_t}$, we have
\begin{equation}
\Exp\big[\normtwobig{\nabla f_{s_t}(\widetildebx^{(j)}) - \nabla f_{s_t}(\bx^*)}^2\big] \leq 2\beta \Exp[f(\widetildebx^{(j)}) - f(\bx^*)].
\end{equation}
Combining the preceding two inequalities, \eqref{equation:svrg_decom} becomes 
\begin{equation}
\Exp\big[\normtwobig{\bp^\toptzero}^2\big] \leq 4\beta \left( \Exp[f(\bx^\toptzero) - f(\bx^*)] + \Exp[f(\widetildebx^{(j)}) - f(\bx^*)] \right).
\end{equation}
Substituting the upper bound of $\Exp\big[\normtwobig{\bp^\toptzero}^2\big]$ into the estimate of $\Exp[e_{t+1}^2]$ in \eqref{equation:conv_svrg_et1et}, we have
\begin{equation}
\begin{aligned}
\Exp[e_{t+1}^2] \leq& \Exp[e_t^2] - 2\eta \Exp\big[f(\bx^\toptzero) - f(\bx^*)\big] + \eta^2 \Exp\big[\normtwobig{\bp^\toptzero}^2\big] \\
\leq& \Exp[e_t^2] - 2\eta(1 - 2\eta \beta) \Exp\big[f(\bx^\toptzero) - f(\bx^*)\big] 
+ 4\beta\eta^2 \Exp[f(\widetildebx^{(j)}) - f(\bx^*)].
\end{aligned}
\end{equation}
Taking the sum over  $t\in\{j, j+1, \ldots, j+m-1\}$, and noting that $\bx^{(j)} = \widetildebx^{(j)}$ (an initialization for each subroutine, i.e., initialization every $m$ steps), we can obtain
\begin{equation}
\begin{aligned}
&\Exp[e_{j+m}^2] + 2\eta(1 - 2\eta \beta) \sum_{t=j}^{j+m-1} \Exp\big[f(\bx^\toptzero) - f(\bx^*)\big] \\
\leq& \Exp\big[\normtwobig{\widetildebx^{(j)} - \bx^*}^2\big] + 4\beta\eta^2 m \Exp\big[f(\widetildebx^{(j)}) - f(\bx^*)\big] \\
\leq& \frac{2}{\alpha} \Exp\big[f(\widetildebx^{(j)}) - f(\bx^*)\big] + 4\beta\eta^2 m \Exp\big[f(\widetildebx^{(j)}) - f(\bx^*)\big],
\end{aligned}
\end{equation}
where the last inequality follows from the SC Property-II (Theorem~\ref{theorem:exi_close_sc}).
Since $\widetildebx^{(j+m)} = \frac{1}{m} \sum_{t=j}^{j+m-1} \bx^{(t)}$, using the convexity of $f$, the above inequality becomes
\begin{equation}
\begin{aligned}
\Exp\big[f(\widetildebx^{(j+m)}) - f(\bx^*)\big] 
&\leq \frac{1}{m} \sum_{t=j}^{j+m-1} \Exp\big[f(\bx^\toptzero) - f(\bx^*)\big]\\
&\leq \frac{1}{2\eta(1 - 2\eta \beta)m} \left( \frac{2}{\alpha} + 4\beta\eta^2m \right) \Exp\big[f(\widetildebx^{(j)}) - f(\bx^*)\big].
\end{aligned}
\end{equation}
This completes the proof.
\end{proof}

\section{Stochastic Optimization Algorithms}
We provided an overview of the number of papers utilizing stochastic optimization algorithms in Table~\ref{table:stochastic-optimizers}. This section introduces the development of these algorithms.
Note that for all algorithms, at iteration $\bx^\toptzero$, our goal is to determine the step $\bd^\toptzero$ leading to the update:
$$
\bx^\toptone \leftarrow \bx^\toptzero + \bd^\toptzero.
$$
We provide the update formulas for various stochastic optimization algorithms; for an illustrative comparison of these algorithms in the context of deep neural networks, refer to \citet{hinton2012neural2, zeiler2012adadelta, lu2022adasmooth, lu2022gradient}.
For convenience, unless specified otherwise, the gradient notations used in this section, $\bg^\toptzero \triangleq \nabla f(\bx^\toptzero)$, can either represent a strict gradient (using only one sample) or a mini-batch gradient (using a subset of samples).
In Section~\ref{section:learning-rate-annealing}, we will discuss learning rate (the term ``learning rate" is used more commonly than ``stepsize" in the machine learning and deep learning communities) annealing or warmup strategies for training deep neural networks or transformers; however, for the purposes of this section, we assume the learning rate remains constant: $\eta_t=\eta$ for all iterations.

\subsection{Momentum }\label{section:sgd_momentum}
If the cost surface is non-spherical, learning can be significantly slowed down because the learning rate must be kept small to avoid divergence along directions with steep curvature \citep{polyak1964some, rumelhart1986learning, qian1999momentum, sutskever2013importance}. 
SGD with momentum (that can be applied to full batch or mini-batch learning) attempts to accelerate learning by leveraging previous steps when suitable, thereby improving convergence rates in deep networks or transformers. The core idea behind momentum is to speed up learning along dimensions where the gradient consistently points in the same direction, while slowing down along dimensions where the gradient's sign frequently changes.

Figure~\ref{fig:momentum_gd} illustrates updates for basic GD (or SGD), showing consistent updates along dimension $x_1$ and zigzag movements along dimension $x_2$ continues to change in a zigzag pattern.  
Momentum-based GD  (or SGD) keeps track of past parameter updates with exponential decay, and its update rule from iteration $t$ to iteration $t+1$  is as follows:
\begin{equation}\label{equation:momen_sgd}
\begin{aligned}
\bd^\toptzero &\leftarrow \rho \bd^\toptminus - \eta  \nabla f(\bx^\toptzero);\\
\bx^\toptone &\leftarrow \bx^\toptzero + \bd^\toptzero,
\end{aligned}
\end{equation}
where the algorithm blends the current update with the past update using a parameter $\rho$, called the \textit{momentum parameter}.
This mechanism ensures that the change in parameters is proportional to the negative gradient plus the previous weight change, effectively smoothing and accelerating the updates; the added \textit{momentum term} acts as both a smoother and an accelerator.
The momentum parameter $\rho<1$ acts  as a \textit{decay constant}, meaning that although $\bd^{(1)}$ might influence $\bd^{(100)}$, its effect diminishes over time. Typically, $\rho$ is set to  0.9 by default.
Momentum mimics inertia in physics, indicating that each update depends not just on the gradient descent (\textit{dynamic term}) but also retains a component related to the previous update's direction (\textit{momentum}).

As mentioned, we previously introduced GD with momentum in Equation~\eqref{equation:gd_with_momentum} (along with its application in quadratic models in Section~\ref{section:quadratic-in-momentum}). The difference here is that the gradient is computed using a single sample or a batch of samples in SGD with momentum.

\begin{figure}[h]
\centering  
\vspace{-0.15cm}  
\subfigtopskip=2pt  
\subfigbottomskip=2pt  
\subfigcapskip=-5pt 
\subfigure[
A two-dimensional surface plot for quadratic convex function.
]{\label{fig:alsgd12}
\includegraphics[width=0.47\linewidth]{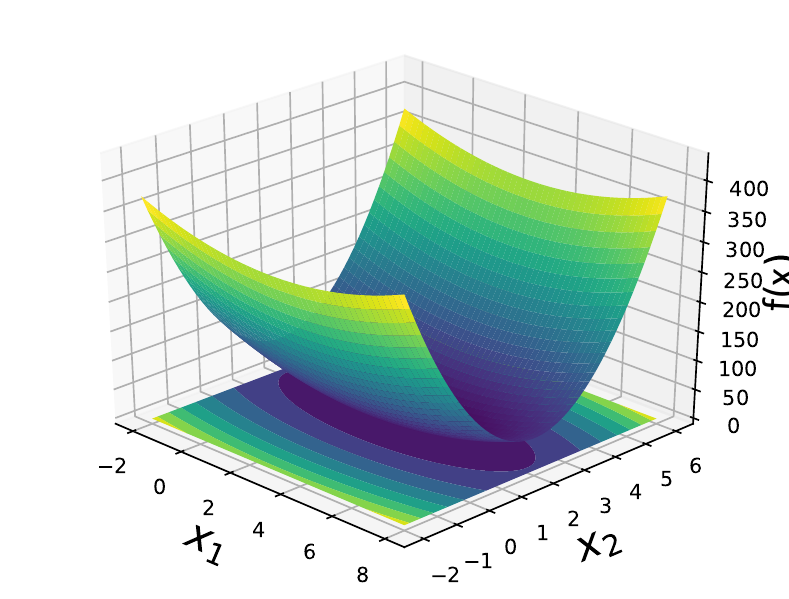}}
\subfigure[The contour plot of $f(\bx)$. The red dot is the optimal point.]{\label{fig:alsgd22}
\includegraphics[width=0.44\linewidth]{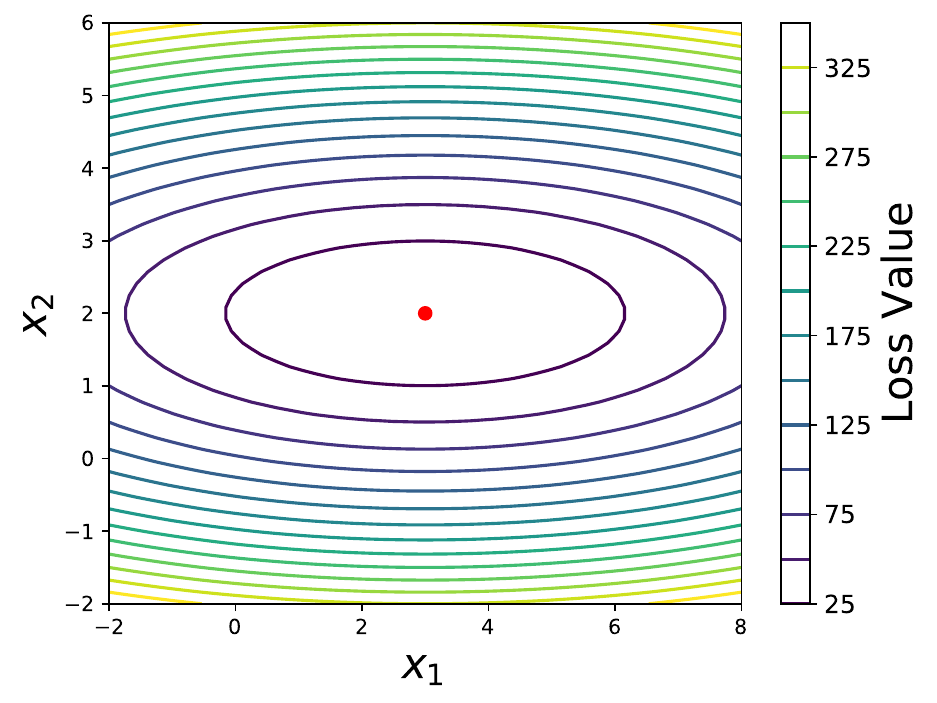}}
\caption{Figure~\ref{fig:alsgd12} shows a function surface and its contour plot (\textcolor{mylightbluetext}{blue}=low, \textcolor{mydarkyellow}{yellow}=high), where the upper graph is the surface, and the lower one is its projection (i.e., contour). The quadratic function $f(\bx) = \frac{1}{2} \bx^\top \bA \bx - \bb^\top \bx + c$ is from parameters $\bA=\scriptsize\begin{bmatrix}
4 & 0\\ 0 & 40
\end{bmatrix}$, $\bb=[12,80]^\top $, and $c=103$. Or equivalently, $f(\bx)=2(x_1-3)^2 + 20(x_2-2)^2+5$ and $\nabla f(\bx)=[4x_1-12, 8x_2-16]^\top$.}
\label{fig:momentum-contour}
\end{figure}

\begin{figure}[h]
\centering  
\vspace{-0.35cm} 
\subfigtopskip=2pt 
\subfigbottomskip=2pt  
\subfigcapskip=-5pt  
\subfigure[Optimization without momentum. A higher learning rate
may result in larger parameter updates in the dimension across the valley (direction of $x_2$) which could lead to oscillations back and forth across the valley.]{\label{fig:momentum_gd}
	\includegraphics[width=0.47\linewidth]{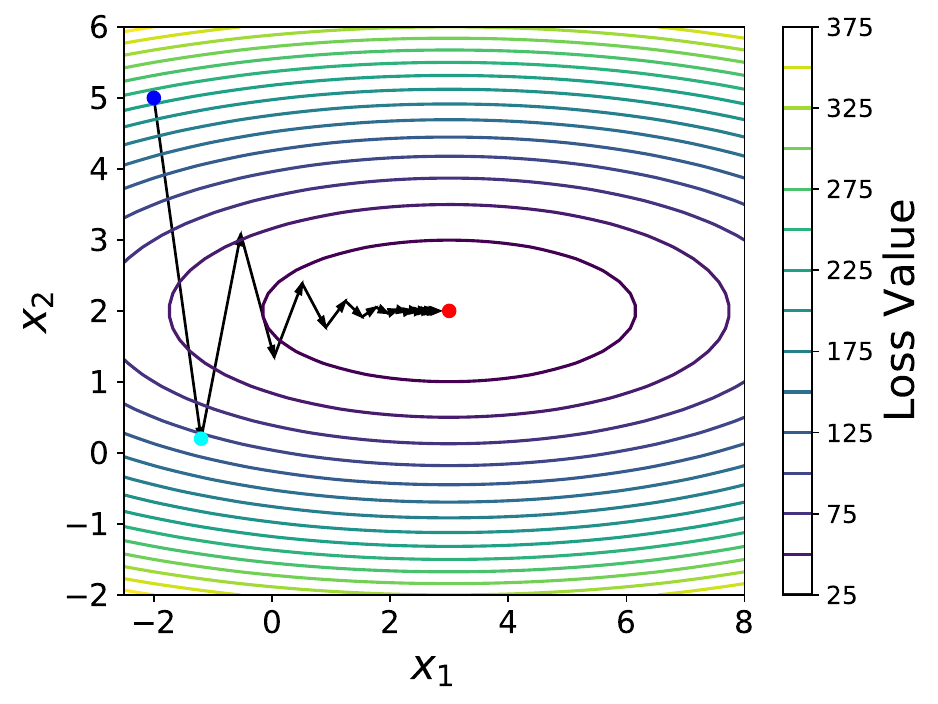}}
\subfigure[Optimization with momentum. Though the gradients along the valley (direction of $x_1$) are much smaller than the gradients across the valley (direction of $x_2$), they are
typically in the same direction. Thus, the momentum term
accumulates to speed up movement, dampens oscillations, and causes us to barrel through narrow valleys, small humps and (local) minima.]{\label{fig:momentum_mum}
	\includegraphics[width=0.47\linewidth]{./imgs/mom_surface_lrate40_gd-1_xy-2-5_mom-2.pdf}}
\caption{The loss function is shown in Figure~\ref{fig:momentum-contour}. The starting point is $[-2, 5]^\top$. After 5 iterations, the squared loss from basic GD is 42.72, and the loss from GD with momentum is 35.41 in this simple case. The learning rates $\eta$ are set to be 0.04 in both cases.}
\label{fig:momentum_gd_compare}
\end{figure}

Momentum excels especially in scenarios with ravine-shaped loss curves, characterized by significantly steeper slopes in one dimension than another (see Figure~\ref{fig:momentum-contour}, which  are common near local minima in deep neural networks). In such areas, basic GD or SGD struggles due to oscillations across the valley. Momentum helps by accelerating gradients towards the optimal direction and reducing oscillations, as shown in Figure~\ref{fig:momentum_gd_compare}.

As shown by the toy example in Figure~\ref{fig:momentum_gd}, GD tends to oscillate across the narrow ravine since the negative gradient will point down one of the steep sides rather than along the ravine towards the optimum. 
By incorporating a fraction $\rho$ of the previous update vector into the current update, momentum accelerates progress when consecutive updates align (e.g., the blue arrow regions in Figure~\ref{fig:momentum_mum}). Conversely, if updates oppose, the algorithm tends to follow the predominant direction of recent updates. Thus, momentum accumulates contributions along the long axis while averaging out oscillations along the short axis, ultimately facilitating smoother navigation through narrow valleys and local minima.

\subsection{Nesterov Momentum}
\textit{Nesterov momentum}, also referred to as \textit{Nesterov accelerated gradient (NAG)}, introduces a variant of the traditional momentum update method and has been increasingly adopted in recent years. The fundamental concept of Nesterov momentum involves predicting the future position of the parameter vector $\bx^\toptzero$ based on the momentum term discussed previously.
Given that the momentum alone (ignoring the gradient term) will adjust the parameter vector by $\rho \bd^\toptminus$, it becomes logical to compute the gradient at this predicted future position $\bx^\toptzero + \rho \bd^\toptminus$. 
This lookahead approach considers where the parameters are likely to be updated next, rather than their current position $\bx^\toptzero$.
Consequently, the gradient is calculated at this anticipated position, resulting in the following update rule:
$$
\begin{aligned}
\bd^\toptzero &\leftarrow \rho \bd^\toptminus - \eta \nabla f(\bx^\toptzero + \rho \bd^\toptminus);\\
\bx^\toptone &\leftarrow \bx^\toptzero + \bd^\toptzero,
\end{aligned}
$$

\begin{figure}[h]
\centering   
\vspace{-0.35cm}  
\subfigtopskip=2pt  
\subfigbottomskip=2pt  
\subfigcapskip=-5pt  
\subfigure[Momentum: evaluate gradient at the current position $\bx^\toptzero$, and the  momentum is about to carry us to the tip of the green arrow.]{\label{fig:nesterov1}
\includegraphics[width=0.47\linewidth]{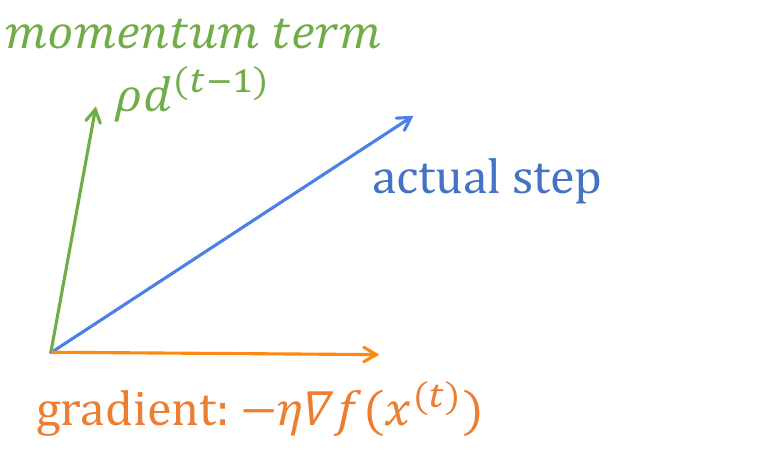}}
\subfigure[Nesterov momentum: evaluate the gradient at this ``looked-ahead" position.]{\label{fig:nesterov2}
\includegraphics[width=0.47\linewidth]{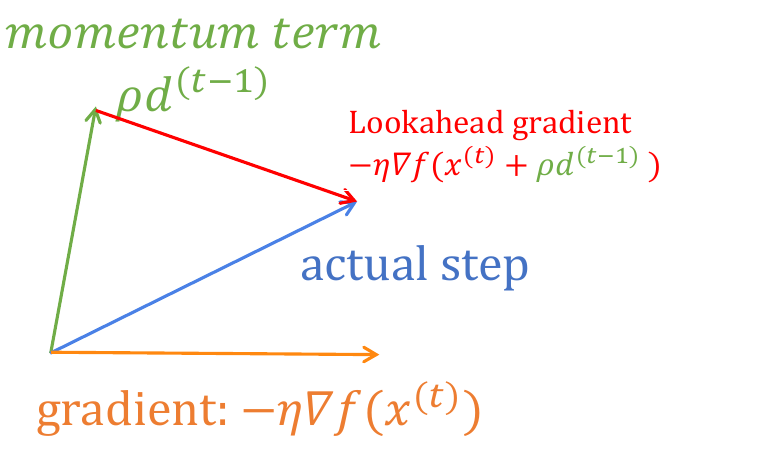}}
\caption{Comparison of momentum and Nesterov momentum.}
\label{fig:momentum_nesterov_comp}
\end{figure}

Figure~\ref{fig:momentum_nesterov_comp} illustrates the distinction between standard momentum and Nesterov momentum. In standard momentum, the gradient is evaluated at the current position $\bx^\toptzero$, with the momentum carrying the parameters to the tip of the green arrow. In contrast, Nesterov momentum evaluates the gradient at the ``lookahead" position, providing an opportunity to correct the direction before making the actual step.
This ``peeking ahead" feature is designed to mitigate excessive velocities by assessing the objective function values in the proposed search direction. Essentially, Nesterov momentum first makes a large jump in the direction of the accumulated gradient, then measures the gradient at the new location and adjusts accordingly. Conversely, standard momentum initially moves according to the current gradient, followed by a significant leap in the direction of the previous gradient accumulation.

To use an analogy, if you're going to make a speculative move (gamble), it's more effective to make the move first and then correct based on the outcome, rather than correcting first and then gambling  \citep{hinton2012neural2}. According to \citet{sutskever2013importance}, Nesterov momentum offers a theoretically superior performance bound compared to standard gradient descent in convex, non-stochastic settings.

\subsection{AdaGrad}
The learning rate annealing procedure (will be introduced in Section~\ref{section:learning-rate-anneal}) adjusts a single global learning rate applied uniformly across all parameter dimensions. To address this limitation, \citet{duchi2011adaptive} introduced \textit{AdaGrad}, which updates the learning rate on a per-dimension basis. In AdaGrad, each parameter's learning rate is influenced by its update history, allowing parameters with fewer past updates to be adjusted more aggressively through a higher learning rate. Specifically, parameters that have received fewer updates in the past are assigned higher learning rates currently.
Formally, AdaGrad's update rule is given by:
\begin{equation}\label{equation:update:adagrad}
\begin{aligned}
\bd^\toptzero &= -  \frac{ \eta}{ \sqrt{\sum_{\tau=1}^t (\bg^{(\tau)})^2 +\epsilon} }\hadaprod \bg^\toptzero ,
\end{aligned}
\end{equation}
where $\hadaprod$ denotes the Hadamard product, $\epsilon$ serves as a smoothing term to stabilize division, $\eta$ represents a shared global learning rate, $(\bg^{(\tau)})^2$ denotes the element-wise square $\bg^{(\tau)}\hadaprod \bg^{(\tau)}$, and the denominator computes the $\ell_2$ norm of the sum of all previous squared gradients on a per-dimension basis. 
Though the global learning rate $\eta$ is used across all dimensions, each dimension has its own dynamic learning rate regulated by the $\ell_2$ norm of accumulated gradient magnitudes. Since this dynamic learning rate grows with the inverse of the accumulated gradient magnitudes, larger gradient magnitudes correspond to smaller learning rates, while smaller gradient magnitudes result in larger learning rates. 
Consequently, the aggregated squared magnitude of partial derivatives for each parameter throughout the algorithm acts similarly to learning rate annealing; see Section~\ref{section:learning-rate-anneal}.

AdaGrad's simplicity makes it straightforward to implement, as shown in the following Python snippet:
\begin{python}
# Assume the gradient dx and parameter vector x 
cache += dx**2
x += - learning_rate * dx / np.sqrt(cache + 1e-8)
\end{python}
One advantage of AdaGrad is its ability to reduce the need for tuning the learning rate, which is controlled by accumulated gradient magnitudes. However, AdaGrad suffers from a significant drawback due to the unbounded accumulation of squared gradients in the denominator. Since every added term is positive, the accumulated sum grows continuously during training, causing the per-dimension learning rate to diminish over time and eventually approach zero, halting further training.

Moreover, AdaGrad can be sensitive to the initialization of parameters and their corresponding gradients. If initial gradient magnitudes are large, per-dimension learning rates will remain low throughout training. While increasing the global learning rate may help, this makes AdaGrad sensitive to the choice of the global learning rate.

Furthermore, AdaGrad assumes that parameters with fewer updates should have larger learning rates and those with more updates should have smaller learning rates. This approach considers only information from squared gradients or the absolute values of gradients, neglecting total movement (i.e., the sum of updates; see Section~\ref{section:adaer} for a further discussion).

To summarize, AdaGrad's main limitations include:
\begin{itemize}
\item Continuous decay of learning rates throughout training.
\item The necessity for manually selecting a global learning rate.
\item  Consideration of only the absolute value of gradients.
\end{itemize}

\subsection{RMSProp}

\textit{RMSProp} is an optimization algorithm that builds upon AdaGrad, aiming to address its primary limitation---AdaGrad's tendency to aggressively reduce the learning rate over time \citep{hinton2012neural, zeiler2012adadelta}. The core idea behind RMSProp is straightforward: it limits the influence of past gradients by considering only the most recent ones instead of accumulating all past gradients up to the current time step $t$.
However, storing a fixed number $w$ of previous squared gradients can be inefficient. To circumvent this, RMSProp employs an exponentially decaying average of the squared gradients, similar to how momentum accumulates past gradients but with a focus on their magnitude.

We first delve into the specifics of RMSProp's formulation. At any given iteration $t$, the running average of the squared $\Exp[\bg^2]^\toptzero$ is calculated as:
%
%
%
\begin{equation}\label{equation:adagradwin}
\Exp[\bg^2]^\toptzero = \rho \Exp[\bg^2]^\toptminus + (1 - \rho) (\bg^\toptzero)^2,
\end{equation}
where $\rho$ is a decay constant similar to that used in the momentum method, and $(\bg^\toptzero)^2$ indicates the element-wise square $\bg^\toptzero\hadaprod \bg^\toptzero$. 
This formula effectively blends the existing estimate with the latest gradient information. Initially, the running average is set to $\bzero$, which might introduce bias early on; however, this bias diminishes over time. Notably, older gradients' impact decreases exponentially throughout the training process.

Let $\rms[\bg]^\toptzero \triangleq \sqrt{\Exp[\bg^2]^\toptzero + \epsilon}$, where again a constant $\epsilon$ is added for numerical stability. Similar to the AdaGrad update in \eqref{equation:update:adagrad},  the resulting stepsize can be obtained as follows:
\begin{equation}\label{equation:rmsprop_update}
\bd^\toptzero=- \frac{\eta}{\rms[\bg]^\toptzero}  \hadaprod \bg^\toptzero.
\end{equation}

The \textit{exponential moving average (EMA)} concept underpins the above formulation. In EMA, $1-\rho$ is also known as the smoothing constant (SC),which can be approximated by $\frac{2}{N+1}$, where $N$ represents the number of past values considered in the average \citep{lu2022adasmooth}:
\begin{equation}\label{equation:ema_smooting_constant}
\text{SC}=1-\rho \approx \frac{2}{N+1},
\end{equation}
which establishes a relationship among  different variables: the decay constant $\rho$, the smoothing constant (SC), and the period $N$.
For instance , if $\rho=0.9$, then $N=19$. That is, roughly speaking, $\Exp[\bg^2]^\toptzero $ at iteration $t$ is roughly equal to the moving average of the past 19 squared gradients and the current one (i.e., the moving average of a total of 20 squared gradients).
The relationship in \eqref{equation:ema_smooting_constant} though is not discussed in the original paper of \citet{zeiler2012adadelta}, it is important to decide the lower bound of the decay constant $\rho$. Typically, a time period of $N=3$ or 7 is thought to be a relatively small frame making the lower bound of decay constant $\rho=0.5$ or 0.75; as $N\rightarrow \infty$, the decay constant $\rho$ approaches $1$.

While AdaGrad excels in convex settings, RMSProp demonstrates superior performance in non-convex scenarios. When applied to a non-convex function, e.g., to train a deep neural network, the learning trajectory can pass through many different structures and eventually arrives at a region that is a locally convex bowl. AdaGrad shrinks the learning rate according to the entire history of the squared partial derivative leading to an infinitesimally small learning rate before arriving at such a convex structure. While RMSProp discards ancient squared gradients to address this problem.

Despite these advantages, RMSProp has limitations---it only considers the magnitude of gradients and does not account for their direction. This characteristic can result in suboptimal learning rates near local minima.

RMSProp was independently developed by \citet{hinton2012neural} and \citet{zeiler2012adadelta}, both addressing AdaGrad's rapid decrease in learning rates. Hinton suggests setting $\rho=0.9$ and the global learning rate  $\eta=0.001$. Furthermore, RMSProp can be integrated with the Nesterov momentum method \citep{goodfellow2016deep}, enhancing its capabilities for deep learning applications; the comparison between the two is presented in Algorithm~\ref{alg:rmsprop} and Algorithm~\ref{alg:rmsprop_nesterov}..

%
%

\begin{algorithm}[h] 
\caption{RMSProp}
\label{alg:rmsprop}
\begin{algorithmic}[1]
\State {\bfseries Input:} Initial parameter $\bx^{(1)}$, constant $\epsilon$;
\State {\bfseries Input:} Global learning rate $\eta$, by default $\eta=0.001$;
\State {\bfseries Input:} Decay constant $\rho$;
\State {\bfseries Input:} Initial accumulated squared gradients $\Exp[\bg^2]^{(0)} = \bzero $;
\For{$t=1:T$ } 
\State Compute gradient $\bg^\toptzero \leftarrow \nabla f(\bx^\toptzero)$;
\State Compute running estimate $	\Exp[\bg^2]^\toptzero \leftarrow \rho \Exp[\bg^2]^\toptminus + (1 - \rho) (\bg^\toptzero)^2;$
\State Compute step $\bd^\toptzero \leftarrow- \frac{\eta}{\sqrt{\Exp[\bg^2]^\toptzero+\epsilon }}  \hadaprod \bg^\toptzero$;
\State Apply update $\bx^\toptone \leftarrow \bx^\toptzero + \bd^\toptzero$;
\EndFor
\State {\bfseries Return:} resulting parameters $\bx^\toptzero$, and the loss $f(\bx^\toptzero)$.
\end{algorithmic}
\end{algorithm}
\begin{algorithm}[h] 
\caption{RMSProp with Nesterov Momentum}
\label{alg:rmsprop_nesterov}
\begin{algorithmic}[1]
\State {\bfseries Input:} Initial parameter $\bx^{(1)}$, constant $\epsilon$;
\State {\bfseries Input:} Global learning rate $\eta$, by default $\eta=0.001$;
\State {\bfseries Input:} Decay constant $\rho$, \textcolor{mylightbluetext}{momentum constant $\alpha$};
\State {\bfseries Input:} Initial accumulated squared gradients $\Exp[\bg^2]^{(0)} = \bzero $, and update step $\bd^{(0)}=\bzero$;
\For{$t=1:T$ } 
\State Compute interim update $\widetilde{\bx}^\toptzero \leftarrow \bx^\toptzero + \alpha \bd^\toptminus$;
\State Compute interim gradient $\bg^\toptzero \leftarrow \nabla f(\textcolor{mylightbluetext}{\widetilde{\bx}^\toptzero})$;
\State Compute running estimate $	\Exp[\bg^2]^\toptzero \leftarrow \rho \Exp[\bg^2]^\toptminus + (1 - \rho) (\bg^\toptzero)^2;$
\State Compute step $\bd^\toptzero \leftarrow\textcolor{mylightbluetext}{\alpha \bd^\toptminus}- \frac{\eta}{\sqrt{\Exp[\bg^2]^\toptzero +\epsilon }}  \hadaprod \bg^\toptzero$;
\State Apply update $\bx^\toptone \leftarrow \bx^\toptzero + \bd^\toptzero$;
\EndFor
\State {\bfseries Return:} resulting parameters $\bx^\toptzero$, and the loss $f(\bx^\toptzero)$.
\end{algorithmic}
\end{algorithm}

\subsection{AdaDelta}\label{section:adadelta}

\citet{zeiler2012adadelta} pointed out an inconsistency in the units of the learning rate used in RMSProp (as well as in basic SGD, momentum, and AdaGrad). To address this issue and leverage insights from second-order methods (discussed in Section~\ref{section:new_methods}), the author proposed rearranging the Hessian matrix to adjust the involved quantities.
Although calculating or approximating the Hessian matrix is computationally intensive, it provides valuable curvature information that aids optimization. Additionally, the units in Newton's method are consistent when using the Hessian. Given the Hessian matrix $\bH$, the update rule in Newton's method can be expressed as follows:
\begin{equation}
\bd^\toptzero = - \bH^{-1} \bg^\toptzero \sim \frac{\text{units}(\nabla f(\bx^\toptzero))}{\text{units}(\nabla^2 f(\bx^\toptzero))},
\end{equation}
where $\text{units}(\cdot)$ denotes the units of the matrix or vector.
This implies 
\begin{equation}
\frac{1}{\text{units}(\nabla^2 f(\bx^\toptzero))} = \frac{\bd^\toptzero}{\text{units}(\nabla f(\bx^\toptzero))},
\end{equation}
i.e., the units of the inverse Hessian matrix can be approximated by the right-hand side term of the above equation. Since the RMSProp update in \eqref{equation:rmsprop_update} already includes $\rms[\bg]^\toptzero$ in the denominator, i.e., the units of the gradients. Introducing an additional  unit of the order of $\bd^\toptzero$ in the numerator can match the same order as Newton's method. 
Therefore, we define another exponentially decaying average of the update steps:
\begin{equation}
\begin{aligned}
\rms[\bd]^\toptzero &\triangleq 	\sqrt{\Exp[\bd^2]^\toptzero } 
= \sqrt{\rho \Exp[\bd^2]^\toptminus + (1 - \rho) (\bd^\toptzero)^2  }.
\end{aligned}
\end{equation}
Given that the value of $\bd^\toptzero$ for the current iteration is unknown but the curvature can be locally smoothed,  $\rms[\bd]^\toptzero$ can be approximated by $\rms[\bd]^\toptminus$. 
So we can use an estimation of $\frac{1}{\text{units}(\nabla^2 f(\bx^\toptzero))} $ to replace  the computationally expensive $\bH^{-1}$:
\begin{equation}
\bH^{-1}\sim \frac{1}{\text{units}(\nabla^2 f(\bx^\toptzero))} 
=\frac{\bd^\toptzero}{\text{units}(\nabla f(\bx^\toptzero))} 
\sim \frac{\rms[\bd]^\toptminus}{\rms[\bg]^\toptzero}.
\end{equation}
This approximation of the inverse Hessian uses only RMS measures of $\bg$ and $\bd$, leading to an update step with matched units:
\begin{equation}
\bd^\toptzero = -\frac{\rms[\bd]^\toptminus}{\rms[\bg]^\toptzero} \hadaprod \bg^\toptzero.
\end{equation}
The idea behind AdaDelta, derived from second-order methods, alleviates the need to manually select a learning rate.
Andrej Karpathy has developed a web demo that allows users to compare the convergence rates among SGD, SGD with momentum, AdaGrad, and AdaDelta \footnote{see https://cs.stanford.edu/people/karpathy/convnetjs/demo/trainers.html.}.

\begin{figure}[h]
\centering
\includegraphics[width=0.6\textwidth]{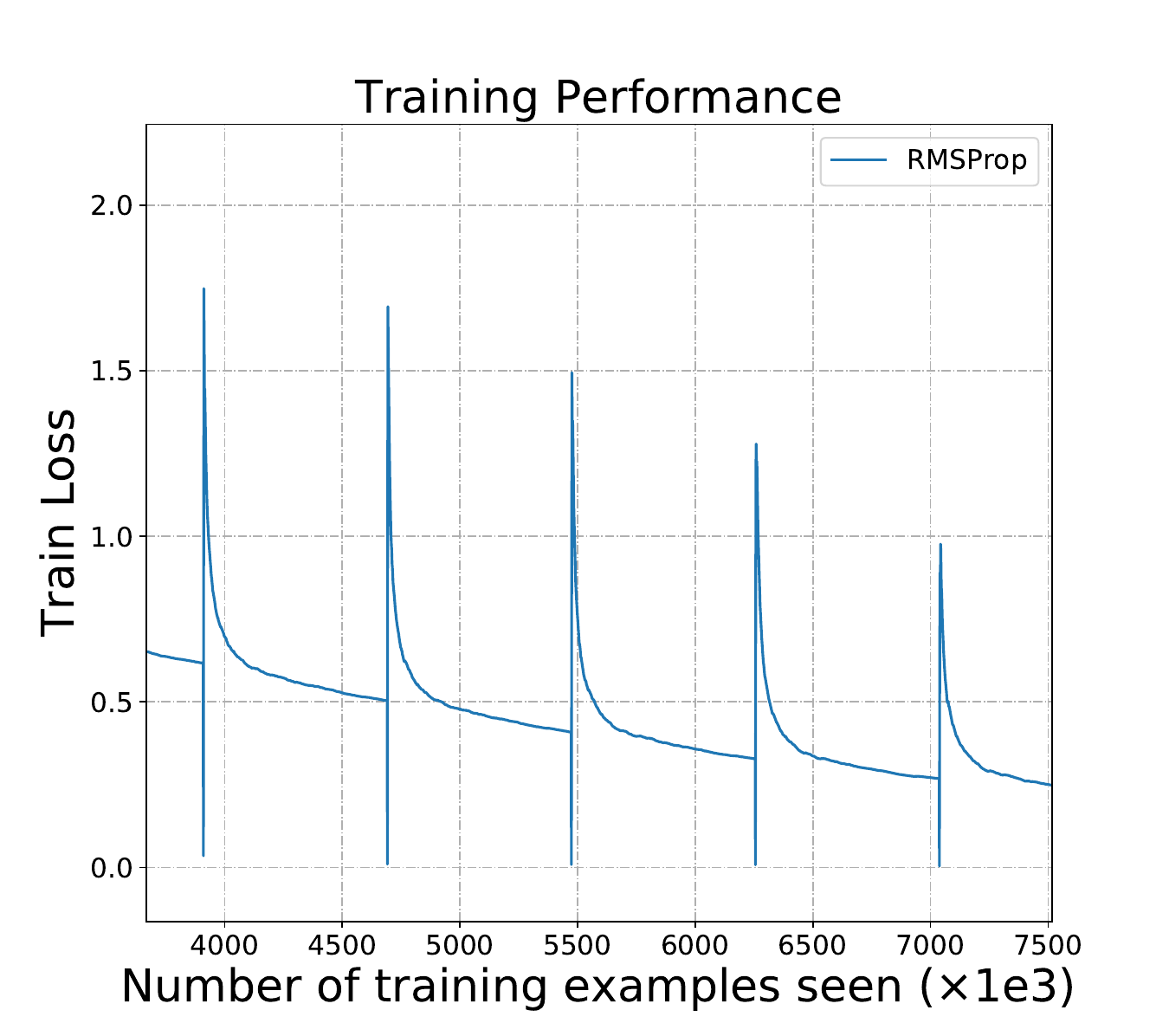}
\caption{Illustration  of tuning parameter after each epoch by loading the weights. 
Loss deterioration is observed at the start of each retraining phase following weight loading.
	}
\label{fig:er-rmsprop_epochstart}
\end{figure}

When using RMSProp or AdaDelta, it is important to note that although the accumulated squared gradients in the denominator help balance per-dimension learning rates, restarting training from saved checkpoints may lead to poor performance in the first few batches due to insufficient gradient data smoothing the denominator. For example, Figure~\ref{fig:er-rmsprop_epochstart} shows loss deterioration after each epoch when weights are loaded from checkpoints. While this does not significantly impact overall training progress, a more effective strategy involves saving $\Exp[\bg^2]^\toptzero$ alongside the neural network weights.

%
\index{Loss deterioration}

\index{Effective ratio}
\index{Exponential moving average}
\subsection{AdaSmooth}\label{section:adaer}

In this section, we will discuss the effective ratio, derived from previous updates in the stochastic optimization process. We will explain how it can be used to achieve adaptive learning rates per dimension through a flexible smoothing constant, hence the name AdaSmooth. This concept builds upon the RMSProp method to address two primary limitations: 1) focusing solely on the absolute value of gradients rather than overall movement in each dimension; and 2) the necessity for manual selection of hyperparameters.
\paragrapharrow{Effective Ratio (ER).}
\citet{kaufman2013trading, kaufman1995smarter} suggested replacing the smoothing constant in the exponential moving average (EMA) formula with a constant based on the \textit{efficiency ratio} (ER).  
This indicator is designed to measure the \textit{strength of a trend}, ranging from $-1.0$ to $+1.0$, with larger magnitudes indicating stronger upward or downward trends.
\citet{lu2022reducing} demonstrated that the ER can also reduce overestimation and underestimation in time series forecasting.
Given the window size $M$ and a series $\{h^{(1)}, h^{(2)}, \ldots, h^{(T)}\}$, the ER is calculated as follows:
\noindent
\begin{equation}
\begin{aligned}
e^\toptzero  &\triangleq \frac{s^\toptzero}{n^\toptzero}\triangleq \frac{h^\toptzero - h^{(t-M)}}{\sum_{i=0}^{M-1} \abs{h^{(t-i)} - h^{(t-1-i)}}}= \frac{\text{Total move for a period}}{\text{Sum of absolute move for each bar}},
\end{aligned}
\end{equation}
where $e^\toptzero$ is the ER of the series at time $t$. 
In strong trends, the ER approaches 1 in absolute value (i.e., the input series is moving in a certain direction, either up or down); if there is no directed movement, it remains close to 0.
Instead of calculating the ER base on the closing price of the underlying asset (the time series problem), we want to calculate the ER of the moving direction in the update methods for each parameter.
Specifically, we focus on how much each parameter moves from its initial point in each period, regardless of direction. Therefore, only the absolute value of the ER is considered. For parameters in the method, the ER is computed as:
\begin{equation}\label{eqution:signoiase-er-delta}
\begin{aligned}
\be^\toptzero  \triangleq \frac{\bs^\toptzero}{\bn^\toptzero}&\triangleq \frac{| \bx^\toptzero -  \bx^{(t-M)}|}{\sum_{i=0}^{M-1} \abs{\bx^{(t-i)} -  \bx^{(t-1-i)}}}= \frac{\abs{ \sum_{i=0}^{M-1} \bd^{(t-1-i)}}}{\sum_{i=0}^{M-1} \abs{\bd^{(t-1-i)}}},
\end{aligned}
\end{equation}
where $\be^\toptzero \in \real^n$, and its $i$-th element $e_{i}^\toptzero$ is in the range of $ [0, 1]$ for all $i$ in $[1,2,\ldots, n]$. A large value of $e_{i}^\toptzero$ indicates the descent method in the $i$-th dimension is moving consistently in a certain direction; while a small value approaching 0 means the parameter in the $i$-th dimension is moving in a zigzag pattern, alternating between positive and negative movements. In practice,   $M$ is chosen based on the batch index within each epoch. 
That is, $M=1$ if the training is in the first batch of each epoch; and $M=M_{\text{max}}$ if the training is in the last batch of the epoch, where $M_{\text{max}}$ is the maximal number of batches per epoch. In other words, $M$ ranges in $[1, M_{\text{max}}]$ for each epoch. Therefore, the value of $e_{i}^\toptzero$ indicates the movement of the $i$-th parameter in the most recent epoch. Or even more aggressively, the window can range from 0 to the total number of batches seen during the entire training progress. The adoption of the adaptive window size $M$ rather than a fixed one has the benefit that we do not need to keep the past $M+1$ steps $\{ \bx^{(t-M)},  \bx^{(t-M+1)}, \ldots,  \bx^\toptzero\}$ to calculate the signal and noise vectors $\{\bs^\toptzero,\bn^\toptzero\}$ in Equation~\eqref{eqution:signoiase-er-delta} since they can be obtained in an accumulated fashion.

\paragrapharrow{AdaSmooth.}\label{section:adaer-after-er}
The AdaSmooth method adjusts the smoothing constant dynamically based on the magnitude of the ER. If the magnitude of ER is small (close to 0), indicating zigzag movements, AdaSmooth uses a longer averaging period to slow down parameter updates. 
Conversely, when the magnitude of ER  is large (near 1), 
the path in that dimension is moving in a certain direction (not zigzag), and the learning actually is happening and the descent is moving in a correct direction, where the learning rate should be assigned to a relatively large value for that dimension; thus, the AdaSmooth method tends to choose a small period, which leads to a small compensation in the denominator (of, for example, \eqref{equation:rmsprop_update}); since the gradients in the closer periods are small in magnitude when it's near the (local) minima. 

\begin{figure}[h]
\centering
\includegraphics[width=0.6\textwidth]{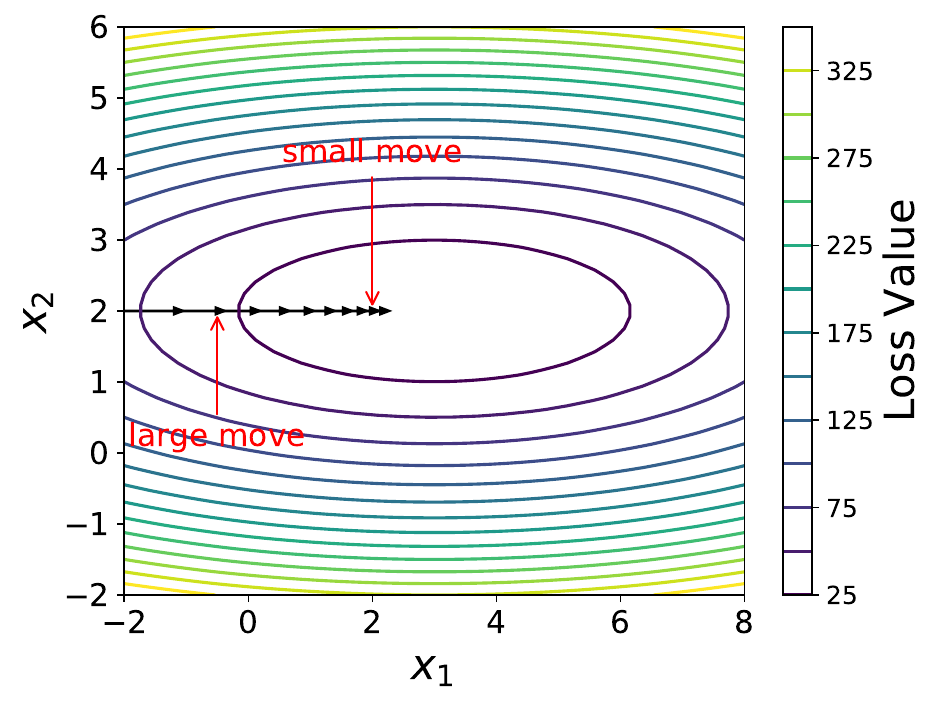}
\caption{Demonstration of how the effective ratio works. Stochastic optimization tends to move a large step when it is far from the (local) minima; and a relatively small step when it is close to the (local) minima.}
\label{fig:er-explain}
\end{figure}

A particular example is shown in Figure~\ref{fig:er-explain}, where the descent is moving in a certain direction, and the gradient in the closer periods is small in magnitude; if we choose a larger period to compensate for the denominator, the descent will be slower due to the large factored denominator.
In short, we want a smaller period to calculate the exponential average of the squared gradients in \eqref{equation:adagradwin} if the update is moving in a certain direction without a zigzag pattern; while when the parameter is updated in a zigzag fashion, the period for the exponential average should be larger \citep{lu2022adasmooth}.


The obtained value of ER is incorporated into the exponential smoothing formula.  
To enhance our approach, we aim to dynamically adjust the time period $N$ discussed in \eqref{equation:ema_smooting_constant} to be a smaller value when the magnitude of ER tends to 1; or to a larger value when the ER moves towards 0. When $N$ is small, $\text{SC} $ is known as a ``\textit{fast SC};" otherwise, $\text{SC} $ is known as a ``\textit{slow SC}." 

For example, let the small time period be $N_1=3$, and the large time period be $N_2=199$.
The smoothing ratio for the fast movement must align with that of  EMA with period $N_1$ (``fast SC" = $\frac{2}{N_1+1}$ = 0.5); and for the period of no trend, the EMA period must be equal to $N_2$ (``slow SC" = $\frac{2}{N_2+1}$ = 0.01). 
Thus the new changing smoothing constant is introduced, called the ``\textit{scaled smoothing constant}" (SSC), denoted by a vector $\bc^\toptzero\in \real^n$:
$$
\bc^\toptzero =  ( \text{fast SC} - \text{slow SC}) \times \be^\toptzero   + \text{slow SC}.
$$
By \eqref{equation:ema_smooting_constant}, we can define the \textit{fast decay constant} $\rho_1=1-\frac{2}{N_1+1}$, and the \textit{slow decay constant} $\rho_2 = 1-\frac{2}{N_2+1}$. Then the scaled smoothing constant vector can be obtained by:
$$
\bc^\toptzero =  ( \rho_2- \rho_1) \times \be^\toptzero   + (1-\rho_2),
$$
where the smaller $\be^\toptzero$, the smaller $\bc^\toptzero$.
For a more efficient influence of the obtained smoothing constant on the averaging period, Kaufman recommended squaring it.
The final calculation formula then follows:
\begin{equation}\label{equation:squared-ssc}
\Exp[\bg^2]^\toptzero  =  (\bc^\toptzero)^2 \hadaprod (\bg^\toptzero)^2  +  \left(1-(\bc^\toptzero)^2 \right)\hadaprod \Exp[\bg^2]^\toptminus.
\end{equation}
or after rearrangement:
$$
\Exp[\bg^2]^\toptzero = \Exp[\bg^2]^\toptminus+ (\bc^\toptzero)^2 \hadaprod \left((\bg^\toptzero)^2 -  \Exp[\bg^2]^\toptminus\right).
$$
We notice that $N_1=3$ is a small period to calculate the average (i.e., $\rho_1=1-\frac{2}{N_1+1}=0.5$) such that the EMA sequence will be noisy if $N_1$ is less than 3. Therefore, the minimum value of $\rho_1$ in practice is set to be greater than 0.5 by default. While $N_2=199$ is a large period to compute the average (i.e., $\rho_2=1-\frac{2}{N_2+1}=0.99$) such that the EMA sequence almost depends only on the previous value, leading to the default value of $\rho_2$ no larger than 0.99. Experimental study  reveals that the AdaSmooth update is insensitive to the hyperparameters  \citet{lu2022adasmooth}. We also carefully notice that when $\rho_1=\rho_2$, the AdaSmooth algorithm recovers to the RMSProp algorithm with decay constant $\rho=1-(1-\rho_2)^2$ since we square it in \eqref{equation:squared-ssc}. After developing the AdaSmooth method, we realize the main idea behind it is similar to that of SGD with momentum: to speed up (compensate less in the denominator) the learning along dimensions where the gradient consistently points in the same direction; and to slow the pace (compensate more in the denominator) along dimensions in which the sign of the gradient continues to change.

\index{Saddle point}
As will be discussed in the cyclical learning rate section (Section~\ref{section:cyclical-lr}), \citet{dauphin2014identifying, dauphin2015equilibrated} argue that the primary challenge in minimizing loss arises from saddle points rather than poor local minima. Saddle points, characterized by small gradients, can impede the learning process. However, an adaptive smoothing procedure for learning rates per dimension can naturally identify these saddle points and compensate less in the denominator---or effectively ``increase" the learning rates---when optimization occurs in these areas. This allows for more rapid traversal of saddle point plateaus. When applied to a non-convex function for training a deep neural network, the learning trajectory may pass through various structures before arriving at a locally convex bowl. AdaGrad shrinks the learning rate based on the entire history of squared partial derivatives, which might make the learning rate too small before reaching such convex structures. RMSProp partially addresses this issue by using an exponentially decaying average to discard older squared gradients, making it more robust in non-convex settings compared to AdaGrad. AdaSmooth advances this approach in two ways: 1) when close to a saddle point, a smaller compensation in the denominator helps escape the saddle point; 2) when near a locally convex bowl, the smaller compensation accelerates convergence.

Empirical evidence shows that the ER used in a simple moving average (SMA) with a fixed window size $M$ can also reflect trends in series/movements in quantitative strategies \citep{lu2022adasmooth}. However, this method again requires storing $M$ previous squared gradients in the AdaSmooth case, making it inefficient, so we do not adopt this extension.

\paragrapharrow{AdaSmoothDelta.}\label{section:adasmoothdelta}
We observe that the ER can also be applied to the AdaDelta setting:
\begin{equation}\label{equation:adasmoothdelta}
\bd^\toptzero = -\frac{\sqrt{\Exp[\bd^2]^\toptminus}}{\sqrt{\Exp[\bg^2]^\toptzero+\epsilon}} \hadaprod \bg^\toptzero,
\end{equation}
where 
\begin{equation}\label{equation:adasmoothdelta111}
\Exp[\bg^2]^\toptzero  =  (\bc^\toptzero)^2 \hadaprod (\bg^\toptzero)^2  +  \left(1-(\bc^\toptzero)^2 \right)\hadaprod \Exp[\bg^2]^\toptminus ,
\end{equation}
and 
\begin{equation}\label{equation:adasmoothdelta222}
\Exp[\bd^2]^\toptminus = \left(1-(\bc^\toptzero)^2\right) \hadaprod (\bd^\toptminus)^2+ (\bc^\toptzero)^2 \hadaprod \Exp[\bd^2]^{(t-2)},
\end{equation}
In this context, $\Exp[\bd^2]^\toptminus$ selects a larger period when the ER is small. 
This is reasonable in the sense that $\Exp[\bd^2]^\toptminus$ appears in the numerator, while $\Exp[\bg^2]^\toptzero$ is in the denominator of \eqref{equation:adasmoothdelta}, making their compensation towards different directions. Alternatively, a fixed decay constant can be applied for $\Exp[\bd^2]^\toptminus$:
$$
\Exp[\bd^2]^\toptminus = (1-\rho_2)  (\bd^\toptminus)^2+ \rho_2  \Exp[\bd^2]^{(t-2)},
$$
The AdaSmoothDelta optimizer introduced above further alleviates the need for a hand specified global learning rate, which is conventionally set to $\eta=1$ in the Hessian context. However, due to the adaptive smoothing constants in \eqref{equation:adasmoothdelta111} and \eqref{equation:adasmoothdelta222}, the $\Exp[\bg^2]^\toptzero $ and $\Exp[\bd^2]^\toptminus$ are less locally smooth, making it less insensitive to the global learning rate than the AdaDelta method. Therefore, a smaller global learning rate, e.g., $\eta=0.5$ is favored in AdaSmoothDelta. The full procedure for computing AdaSmooth is then formulated in Algorithm~\ref{algo:adasmooth}.

\begin{algorithm}[tb]
\caption{AdaSmooth algorithm. All operations on vectors are element-wise. Good default settings for the tested tasks are $\rho_1=0.5, \rho_2=0.99, \epsilon=1e-6, \eta=0.001$; see Section~\ref{section:adaer-after-er} or \eqref{equation:ema_smooting_constant} for a detailed discussion on the explanation of the decay constants' default values. 
The AdaSmoothDelta iteration can be calculated in a similar way.
}
\label{alg:computer-adaer}
\begin{algorithmic}[1]
\State {\bfseries Input:} Initial parameter $\bx_1$, constant $\epsilon$;
\State {\bfseries Input:} Global learning rate $\eta$, by default $\eta=0.001$;
\State {\bfseries Input:} Fast decay constant $\rho_1$, slow decay constant $\rho_2$;
\State {\bfseries Input:} Assert $\rho_2>\rho_1$, by default $\rho_1=0.5$, $\rho_2=0.99$;
\For{$t=1:T$ } 
\State Compute gradient $\bg^\toptzero \leftarrow \nabla f(\bx^\toptzero)$;
\State Compute ER $\be^\toptzero  \leftarrow\frac{\abs{\bx^\toptzero -  \bx^{(t-M)}}}{\sum_{i=0}^{M-1} \abs{ \bd^{(t-1-i)}}}$ ;
\State Compute scaled smoothing vector $\bc^\toptzero \leftarrow  ( \rho_2- \rho_1) \times \be^\toptzero   + (1-\rho_2)$;
\State Compute normalization term $\Exp[\bg^2]^\toptzero  \leftarrow  (\bc^\toptzero)^2 \hadaprod (\bg^\toptzero)^2  +  \left(1-(\bc^\toptzero)^2 \right)\hadaprod \Exp[\bg^2]^\toptminus ;$
\State Compute step $\bd^\toptzero \leftarrow- \frac{\eta}{\sqrt{\Exp[\bg^2]^\toptzero+\epsilon}}  \hadaprod \bg^\toptzero$;
\State Apply update $\bx^\toptone \leftarrow \bx^\toptzero + \bd^\toptzero$;
\EndFor
\State {\bfseries Return:} resulting parameters $\bx^\toptzero$, and the loss $f(\bx^\toptzero)$.
\end{algorithmic}\label{algo:adasmooth}
\end{algorithm}

\begin{figure}[h]
\centering
\includegraphics[width=0.6\textwidth]{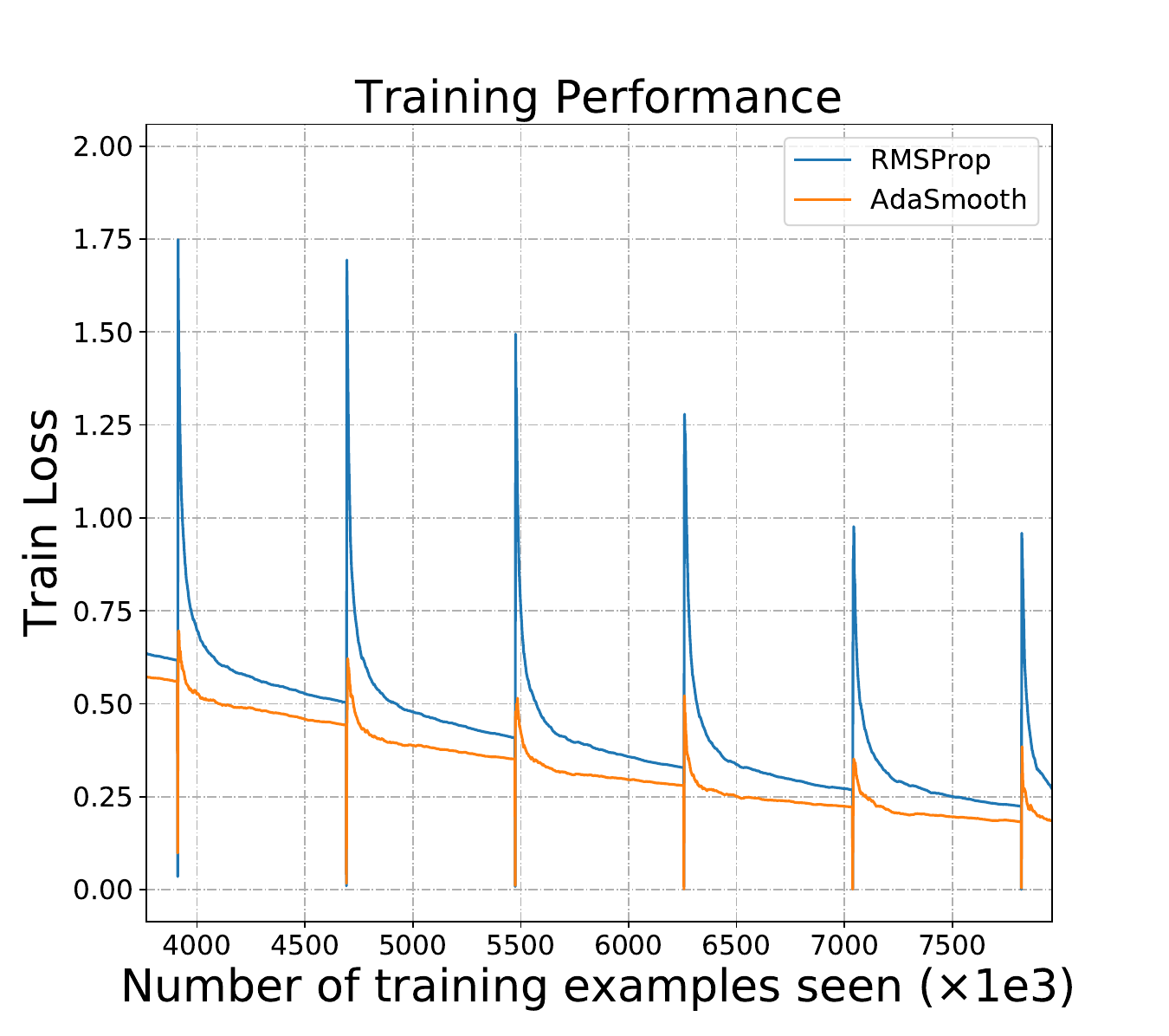}
\caption{Illustration  of parameter tuning after each epoch by loading the weights. We save the weights and reload them after each new epoch such that there are step-points while re-training after each epoch. This issue is less sever in the AdaSmooth method than in  RMSProp. A smaller loss deterioration is observed in the AdaSmooth example than that of the RMSProp case.}
\label{fig:er-rmsprop_epochstart22}
\end{figure}

We have previously discussed the loss deterioration problem encountered when reloading weights from checkpoints using methods like RMSProp or AdaDelta (Figure~\ref{fig:er-rmsprop_epochstart}). However, this issue is less severe in the AdaSmooth setting. As illustrated in Figure~\ref{fig:er-rmsprop_epochstart22}, the loss deterioration observed after reloading weights is smaller in the AdaSmooth example compared to the RMSProp case.
\index{Loss deterioration}

\subsection{Adam}
\textit{Adaptive moment estimation (Adam)} is yet another adaptive learning rate optimization algorithm \citep{kingma2014adam}. 
Adam utilizes both first-order and second-order information in its updates. Similar to RMSProp, AdaDelta, and AdaSmooth, Adam maintains an exponential moving average of past squared gradients (the second moment). However, it also keeps track of an exponentially decaying average of past gradients (the first moment):
\begin{equation}\label{equation:adam-updates}
\begin{aligned}
\bmm^\toptzero &=  \rho_1 \bmm^\toptminus + (1-\rho_1)\bg^\toptzero; \\
\bv^\toptzero &= \rho_2 \bv^\toptminus +(1-\rho_2)(\bg^\toptzero)^2,
\end{aligned}
\end{equation}
where $\bmm^\toptzero$ and $\bv^\toptzero$ are running estimates of the first moment (the mean) and the second moment (the uncentered variance) of the gradients, respectively. The drawback of RMSProp is that the running estimate $\Exp[\bg^2]$ can be biased towards zero during the initial steps because it starts from zero, particularly when the decay constant is large (when $\rho$ is close to 1 in RMSProp). Observing the biases towards zero in \eqref{equation:adam-updates} as well, Adam counteracts these biases by computing the bias-free moment estimates:
$$
\begin{aligned}
\widehat{\bmm}^\toptzero &= \frac{\bmm^\toptzero}{1-\rho_1^t}
\qquad\text{and}\qquad
\widehat{\bv}^\toptzero = \frac{\bv^\toptzero}{1-\rho_2^t}.
\end{aligned}
$$ 
The first and second moment estimates are then incorporated into the update step:
$$
\begin{aligned}
\bd^\toptzero &\leftarrow - \frac{\eta}{\sqrt{\widehat{\bv}^\toptzero}+\epsilon} \hadaprod  \widehat{\bmm}^\toptzero;\\
\bx^\toptone &\leftarrow \bx^\toptzero - \bd^\toptzero.
\end{aligned}
$$
In practice, \citet{kingma2014adam} suggest to use $\rho_1=0.9$, $\rho_2=0.999$, and $\epsilon=1e-8$ for the parameters by default.

\subsection{AdaMax}
Building on Adam, \citet{kingma2014adam} explored the use of higher-order moments for optimization updates. Specifically, they considered:
$$
\bv^\toptzero = \rho_2^p \bv^\toptminus + (1-\rho_2^p) \absbig{\bg^\toptzero}^p,
$$
which can become numerically unstable for large values of  $p$. 
This instability makes the $\ell_1$ and $\ell_2$ norms common choices for update calculations. 
However, as $p\rightarrow \infty$, the $\ell_{\infty}$ norm also shows stable behavior. Based on this observation, AdaMax employs the following moment update rule:
$$
\bu^\toptzero = \rho_2^\infty \bu_{t-1} + (1-\rho_2^\infty) \absbig{\bg^\toptzero}^\infty = \max (\rho_2\bu^{(t-1)}, \absbig{\bg^\toptzero}),
$$
where $\bu^\toptzero$ does not require initialization bias correction. The parameter update step in AdaMax is then given by:
$$
\bd^\toptzero = - \frac{\eta}{\bu^\toptzero} \hadaprod  \widehat{\bmm}^\toptzero.
$$
In practice, \citet{kingma2014adam} suggest to use $\eta=0.002$, $\rho_1=0.9$, and $\rho_2=0.999$ as the default parameters.

\index{Saddle point}
\subsection{Problems in SGD}\label{section:c-problem}
The previous introduced optimization algorithms for stochastic gradient descent  are  widely used in training deep learning models. However, they come with their own set of challenges and potential issues. Here are some common problems associated with SGD:
\paragrapharrow{Saddle points.}
When the Hessian of loss function is positive definite, then the optimal point $\bx^*$ with vanishing gradient must be a local minimum (Theorem~\ref{theorem:second_nec_nonstrict_loca}). 
Similarly, when the Hessian is negative definite, the point is a local maximum; when the Hessian has both positive and negative eigenvalues, the point is a saddle point (Definition~\ref{definition:stat_point}). 
The stochastic optimization algorithms discussed above in practice are first-order optimization algorithms: they only look at the gradient information, and never explicitly compute the Hessian. Such algorithms may get stuck at saddle points (see toy example in Figure~\ref{fig:quadratic_saddle}). In the algorithms presented earlier, including basic update, AdaGrad, AdaDelta, RMSprop, and others, this issue may arise. 
AdaSmooth, as discussed in Section~\ref{section:adaer}, may have mechanisms to escape saddle points. Momentum and Nesterov momentum can also help navigate past saddle points due to their incorporation of previous step sizes; however, these techniques can make convergence more difficult, particularly when the momentum term $\rho$ is large.

\paragrapharrow{Low speed in SGD.}
Despite claims by Rong Ge's post \footnote{http://www.offconvex.org/2016/03/22/saddlepoints/} that SGD might converge to saddle points with a high error rate, \citet{dauphin2014identifying} and Benjamin Recht's post \footnote{http://www.offconvex.org/2016/03/24/saddles-again/} argue that it is actually quite challenging for SGD to converge to saddle points from random initial points. This difficulty arises because local minima and saddle points are often surrounded by plateaus with small curvature on the error surface. In SGD, algorithms tend to be repelled away from saddle points towards regions of lower error along directions of negative curvature. Therefore, while there isn't a true ``saddle point problem" in SGD, escaping these areas can be slow due to the low curvature plateaus. Second-order methods like Newton's method  handle saddle points well, as they aim to rapidly descend through plateaus by scaling gradient steps using the inverse eigenvalues of the Hessian matrix.

As noted by \citet{dauphin2014identifying}, random Gaussian error functions over large $n$ dimensions are increasingly likely to have saddle points rather than local minima as $n$ increases. And the ratio of the number of saddle points to local minima increases exponentially with the dimensionality $n$. The authors also argue that it is saddle points rather than local minima that provide a fundamental impediment to rapid high dimensional non-convex optimization. In this sense, local minima with high errors are exponentially rare compared to saddle points, making computation slow in SGD when navigating through regions of small curvature.

\index{Learning rate annealing}
\index{Learning rate warmup}
\section{Learning Rate Annealing and Warmup}\label{section:learning-rate-annealing}
We have discussed in \eqref{equation:momen_sgd} that the learning rate $\eta$ controls how large of a step to take in the direction of the negative gradient so that we can reach a (local) minimum. In a wide range of applications, a fixed learning rate works well in practice. While there are alternative learning rate schedules that change the learning rate during learning, and it is most often changed between epochs, especially for training deep neural  or transformer structures. 
We've seen that per-dimension optimizers like AdaGrad, AdaDelta, RMSProp, and AdaSmooth can adaptively change the learning rate in each dimension  \citep{duchi2011adaptive, hinton2012neural, zeiler2012adadelta, lu2022adasmooth}. 
This section focuses on strategies for decaying or annealing the global learning rate, i.e., adjusting the value of $\eta$ in \eqref{equation:momen_sgd}.

Using a constant learning rate poses a dilemma: a small learning rate can make the algorithm take too long to reach an optimal solution, whereas a large initial learning rate might initially bring the algorithm close to a good (local) minimum but cause it to oscillate around this point indefinitely. To mitigate these issues, decreasing the learning rate over time can help by slowing down parameter updates. This can be done manually when validation accuracy plateaus, or automatically through decay functions based on epoch count.
Common decay functions include \textit{step decay}, \textit{inverse decay}, and \textit{exponential decay}. 
The following sections explore the mathematical formulations of various learning rate annealing schemes.

\index{Learning rate annealing}
\subsection{Learning Rate Annealing}\label{section:learning-rate-anneal}
\paragrapharrow{Step decay.} The \textit{step decay} scheduler drops the learning rate by a factor every epoch or after a set number of epochs. 
For  iteration $t$, number of iterations to drop $n$, initial learning rate $\eta_0$, and decay factor $d<1$, the form of step decay is given by
$$
\eta_t = \eta_0 \cdot  d^{\floor{\frac{t}{n}} }=\eta_0 \cdot  d^s,
$$
where $s\triangleq{\floor{\frac{t}{n}} }$ represents the \textit{step stage} at which the decay occurs. Therefore, the step decay policy decays the learning rate every $n$ iterations.

\paragrapharrow{Multi-step decay.}
The \textit{multi-step decay} scheduler is a slightly different version of the step decay, wherein the step stage is the index where the iteration $t$ falls in the milestone vector $\bmm=[m_1,m_2,\ldots, m_k]^\top$ with $0\leq m_1\leq m_2\leq\ldots\leq m_k\leq T$ and $T$ being the total number of iterations (or epochs) \footnote{When $T$ is the total number of iterations, it can be obtained by the product of the number of epochs and the number of steps per epoch.}. To be more concrete, the step stage $s$ at iteration $t$ is obtained by
$$
s = 
\left\{
\begin{aligned}
	&0, \gap &t<m_1;\\
	&1, \gap &m_1\leq t < m_2; \\
	&\ldots\\
	&k, \gap &m_k\leq t \leq T. 
\end{aligned}
\right.
$$
As a result, given the iteration $t$, initial learning rate $\eta_0$, and decay factor $d<1$, the learning rate at iteration $t$ is calculated as 
$$
\eta_t = \eta_0 \cdot  d^s.  
$$

\begin{figure}[h]
	\centering
	\includegraphics[width=0.6\textwidth]{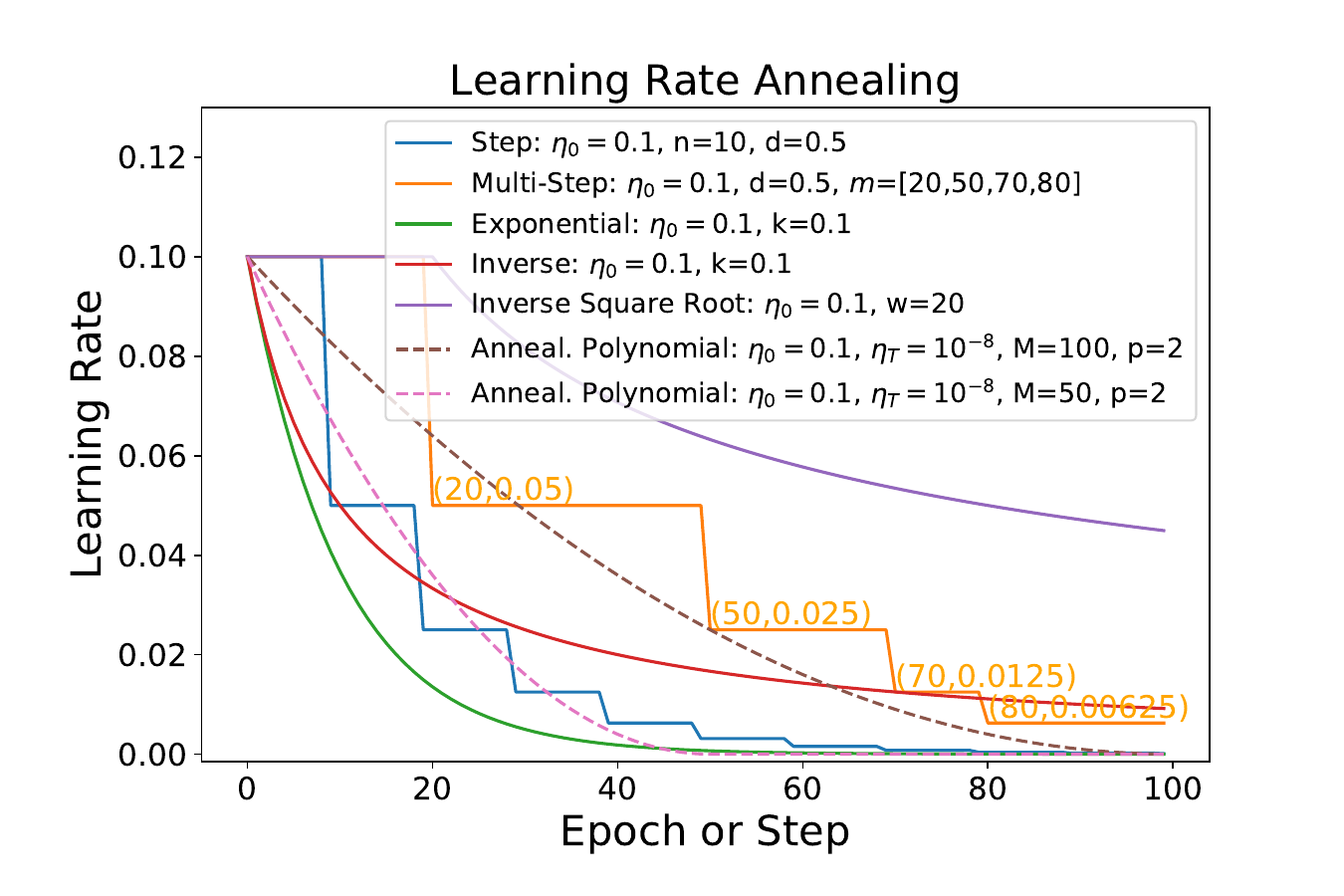}
	\caption{Demonstration of step decay, multi-step decay, annealing polynomial, inverse decay, inverse square root, and exponential decay schedulers. One may find that among the six, exponential decay exhibits the smoothest behavior, while multi-step decay is characterized by the least smoothness.}
	\label{fig:lr_step_decay}
\end{figure}

\paragrapharrow{Exponential decay.}
Given the iteration $t$, the initial learning rate $\eta_0$, and the exponential decay factor $k$, the form of the \textit{exponential decay} is given by
$$
\eta_t = \eta_0 \cdot \exp(-k \cdot t),
$$
where the parameter $k$ controls the rate of the decay.


\paragrapharrow{Inverse decay.}
The \textit{inverse decay} scheduler is a variant of exponential decay in that the decaying effect is applied by the inverse function. Given the iteration number $t$, the initial learning rate $\eta_0$, and the decay factor $k$, the form of the inverse decay is obtained by
$$
\eta_t = \frac{\eta_0}{1+ k\cdot t},
$$
where, again, the parameter $k$ controls the rate of the decay.

\paragrapharrow{Inverse square root.} 
The \textit{inverse square root} scheduler is a learning rate schedule 
$$
\eta_t = \eta_0 \cdot \sqrt{w} \cdot \frac{1}{\sqrt{\max(t, w)}},
$$
where $t$ represents the current training iteration,  $w$ is the number of warm-up steps, and $\eta_0$ is the initial learning rate. This scheduler maintains a constant learning rate for the initial  steps, then exponentially decays the learning rate until the pre-training phase concludes.

\paragrapharrow{Annealing polynomial decay.}
Given the iteration $t$, the \textit{max decay iteration} $M$, the power factor $p$, the initial learning rate $\eta_0$, and the final learning rate $\eta_T$, the \textit{annealing polynomial decay} at iteration $t$ can be obtained by 
\begin{equation}\label{equation:annealing_polynomial}
	\begin{aligned}
		& decay\_batch = \min(t, M) ;\\
		& \eta_t = (\eta_0-\eta_T)\cdot \left(1-\frac{t}{decay\_batch}\right)^{p}+\eta_T.
	\end{aligned}
\end{equation}
In practice, the default values for the parameters are: initial rate $\eta_0=0.001$, end rate $\eta_T=1e-10$, the warm up steps $M=T/2$, where $T$ is the maximal iteration number, and power rate $p=2$.

Figure~\ref{fig:lr_step_decay} compares step decay, multi-step decay, annealing polynomial decay, inverse decay, inverse square root, and exponential decay with a specified  set of parameters. 
The smoothness varies among these methods, with exponential decay exhibiting the smoothest behavior and multi-step decay being the least smooth.  
In  annealing polynomial decay, the max decay iteration parameter $M$ determines the decay rate:
\begin{itemize}
	\item When $M$ is small, the decay gets closer to that of the exponential scheduler or the step decay; however, the exponential decay has a longer tail. That is, the exponential scheduler decays slightly faster in the beginning iterations but slows down in the last few iterations.
	\item When $M$ is large, the decay gets closer to that of the multi-step decay; however, the multi-step scheduler exhibits a more aggressive behavior.
\end{itemize}

\begin{figure}[h]
	\centering  
	\vspace{-0.35cm} 
	\subfigtopskip=2pt 
	\subfigbottomskip=2pt 
	\subfigcapskip=-5pt 
	\subfigure[Training loss.]{\label{fig:mnist_scheduler_1}
		\includegraphics[width=0.47\linewidth]{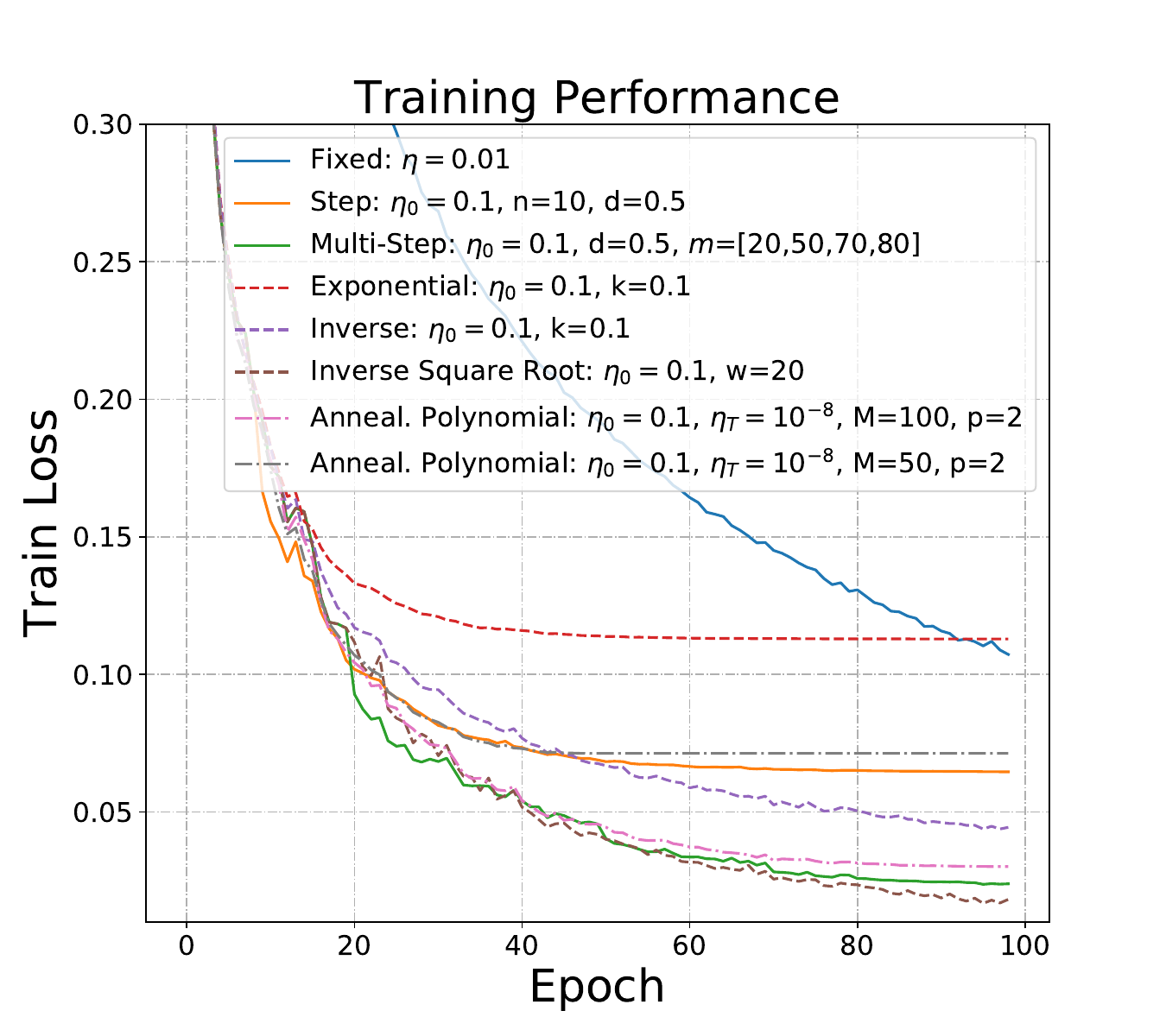}}
	\subfigure[Training accuracy.]{\label{fig:mnist_scheduler_2}
		\includegraphics[width=0.47\linewidth]{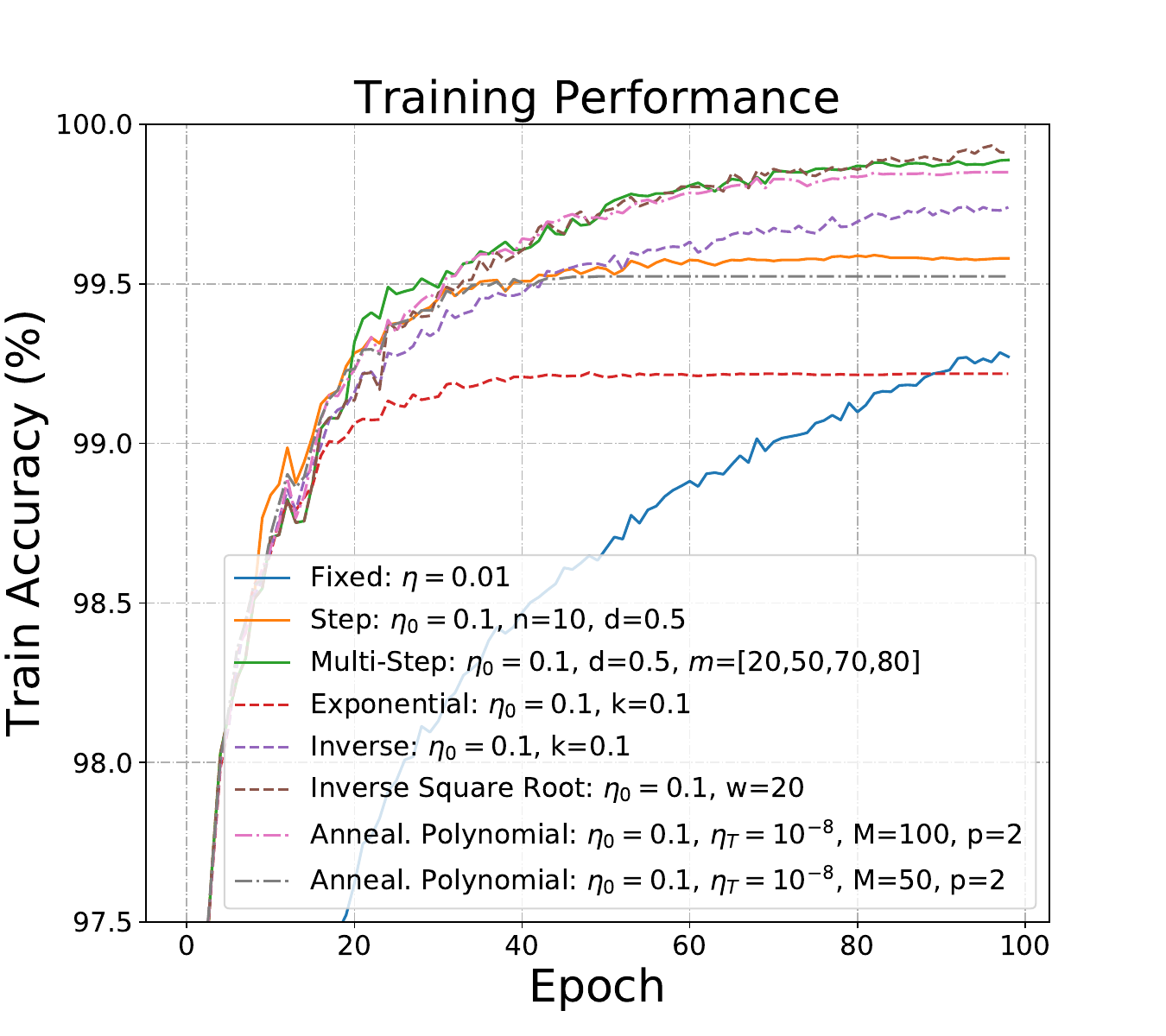}}
	\subfigure[Test loss.]{\label{fig:mnist_scheduler_3}
		\includegraphics[width=0.47\linewidth]{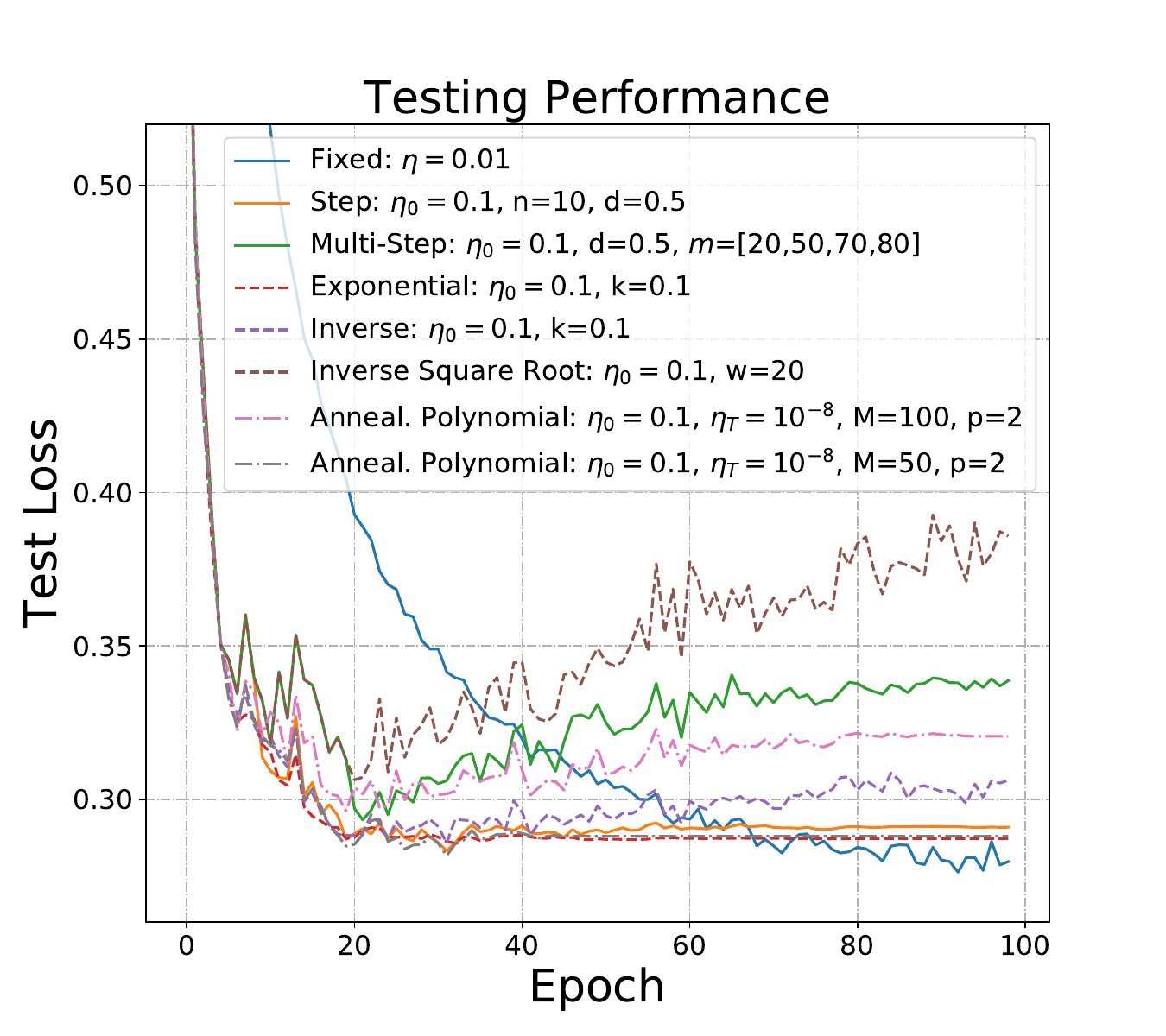}}
	\subfigure[Test accuracy.]{\label{fig:mnist_scheduler_4}
		\includegraphics[width=0.47\linewidth]{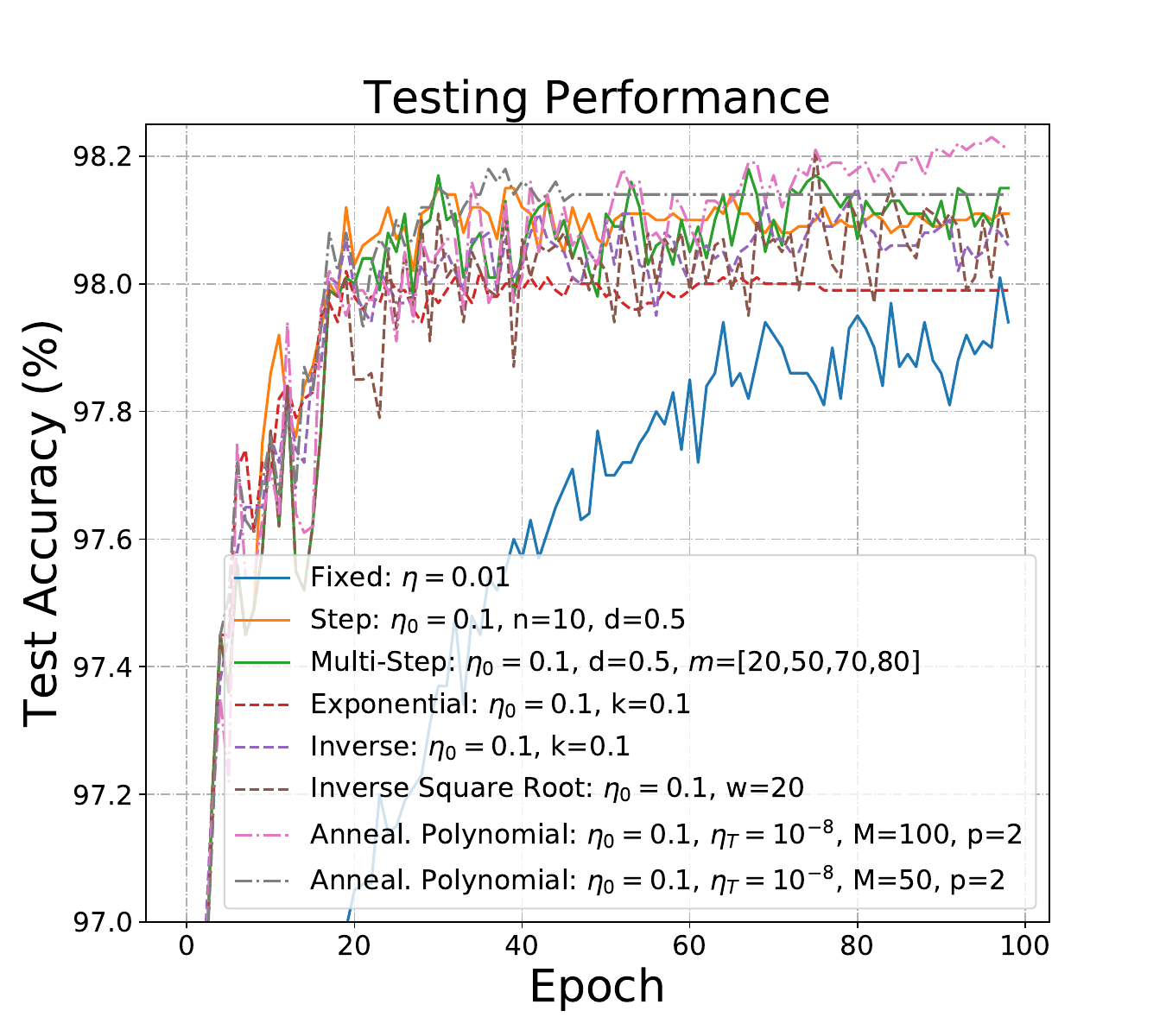}}
	\caption{Training and test performance with different learning rate schemes.}
	\label{fig:mnist_scheduler_1234}
\end{figure}

\paragrapharrow{Toy example.}
To assess the impact of different schedulers, we utilize a toy example involving the training of a multi-layer perceptron (MLP) on the MNIST digit classification set  \citep{lecun1998mnist} \footnote{It has a training set of 60,000 examples, and a test set of 10,000 examples.}. Figure~\ref{fig:mnist_scheduler_1234} presents the training and test performance in terms of \textit{negative log-likelihood loss}. 
The parameters for various schedulers are detailed in Figure~\ref{fig:lr_step_decay} (for 100 epochs). We observe that the stochastic gradient descent method with fixed learning rate may lead to a continued reduction in test loss; however, its test accuracy may get stuck at a certain point. The toy example shows learning rate annealing schemes,  in general, can enhance optimization methods by guiding them towards better local minima with improved performance.

\index{Learning rate warmup}
\subsection{Learning Rate Warmup}
The concept of warmup in training neural networks receive attention in recent years \citep{he2016deep, goyal2017accurate, smith2019super}. Further insights into the efficacy of the warmup scheduler in neural machine translation (NML) can be found in the comprehensive discussion by  \citet{popel2018training}. 
Learning rate annealing schedulers can be applied on both an epoch- and step-basis. However, learning rate warmup schemes are typically implemented on a step-by-step basis, where the total number of steps is the product of the number of epochs and the number of steps per epoch, as mentioned earlier \citep{vaswani2017attention, howard2018universal}. Note that with this scheduler, early stopping should generally be avoided. In the rest of this section, we will explore two commonly used warmup policies: the slanted triangular learning rates (STLR) and the Noam methods.

\paragrapharrow{Slanted Triangular Learning Rates (STLR).} STLR is a learning rate schedule that initially  increases the learning rate over some number of epochs linearly, and then  decays it over the remaining epochs linearly. The rate at iteration $t$ is computed as follows: 
$$
\begin{aligned}
	cut &= \ceil{T\cdot frac} ;\\
	p &= 
	\left\{
	\begin{aligned}
		&t/cut , \gap &\text{if }t<cut;\\
		& 1 - \frac{t-cut }{cut\cdot (1/frac-1)}, \gap &\text{otherwise};
	\end{aligned} 
	\right.\\
	\eta_t &= \eta_{\text{max}} \cdot \frac{1+p\cdot (ratio - 1)}{ratio },\\
\end{aligned}
$$
where $T$ is the number of training iterations (the product of the number of epochs and the number of updates per epoch), $frac$ is the fraction of iterations we want to increase the learning rate, $cut$ is the iteration when we switch from increasing to decreasing the learning rate, $p$ is the fraction of the number of iterations we have increased or decreased the learning rate respectively, $ratio$ specifies how much smaller the lowest learning rate is from the maximum learning rate $\eta_{\text{max}}$. In practice, the default values are $frac=0.1$, $ratio=32$, and $\eta_{\text{max}}=0.01$ \citep{howard2018universal}.

\paragrapharrow{Noam.} The Noam scheduler is originally used in neural machine translation (NML) tasks and is proposed in \citet{vaswani2017attention}. This corresponds to increasing the learning rate linearly for the first ``warmup\_steps" 
training steps and decreasing it thereafter proportionally to the inverse square root of the
step number, scaled by the inverse square root of the dimensionality of the model (linear warmup for a given number of steps followed by exponential decay). Given the warmup steps $w$ and the model size $d_{\text{model}}$ (representing the hidden size parameter which dominates the number of parameters in the model), the learning rate $\eta_t$ at step $t$ can be calculated by 
$$
\eta_t = \alpha \cdot  \frac{1}{\sqrt{d_{\text{model}}}}\cdot  \min \left(\frac{1}{\sqrt{t}} ,  \frac{t}{w^{3/2}}\right),
$$
where $\alpha$ is a smoothing factor. In the original paper, the warmup step $w$ is set to $w=4000$. While in practice, $w=25000$ can be a good choice.

\begin{figure}[h]
	\centering
	\includegraphics[width=0.6\textwidth]{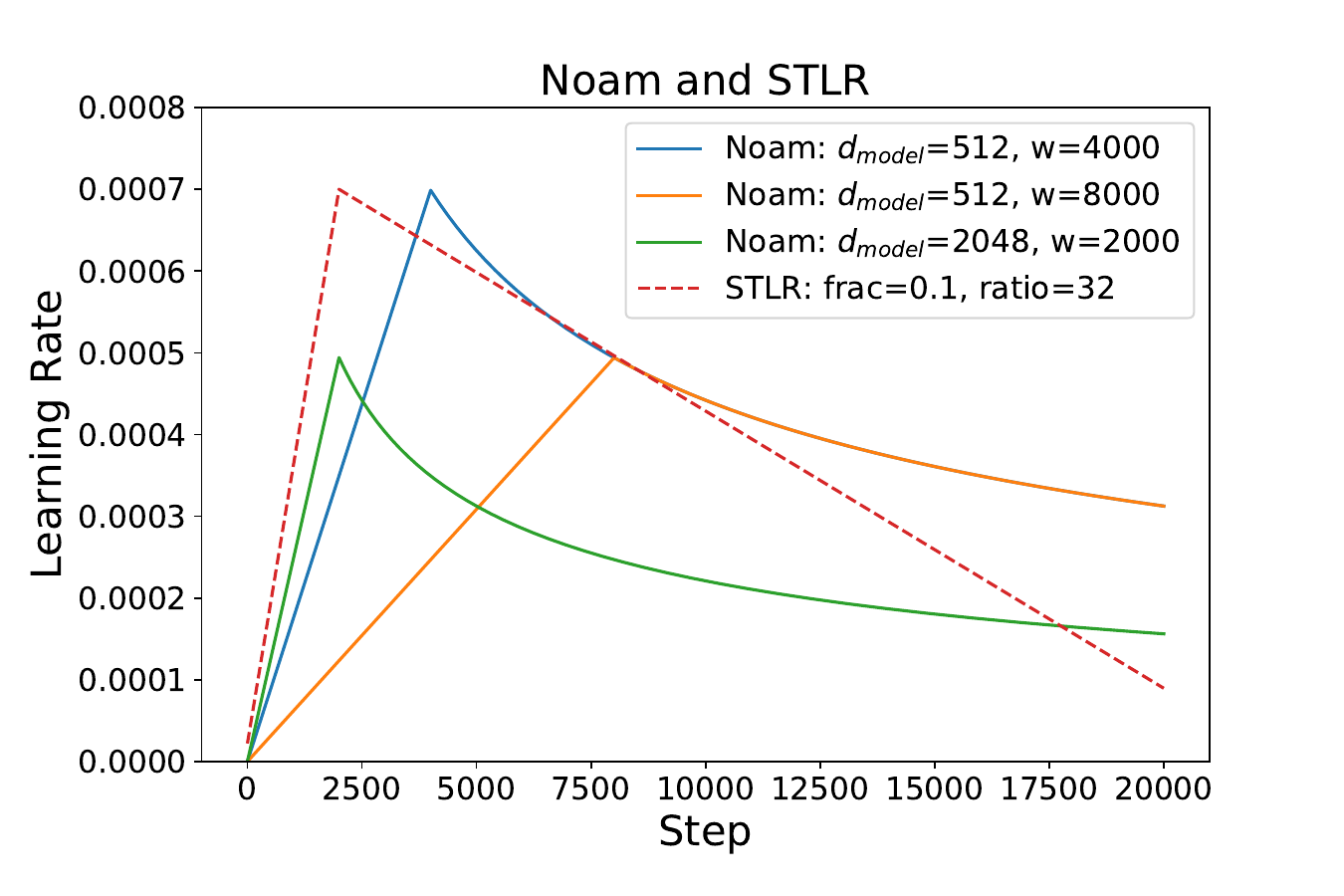}
	\caption{Comparison of Noam and STLR schedulers.}
	\label{fig:lr_noam}
\end{figure}
Moreover, in rare cases, the model size is occasionally  set to be the same as the warmup steps, resulting in what is known as the \textit{warmup Noam scheduler}:
$$
\eta_t = \alpha \cdot \frac{1}{\sqrt{w}} \cdot  \min \left(\frac{1}{\sqrt{t}} ,  \frac{t}{w^{3/2}}\right).
$$
Figure~\ref{fig:lr_noam} compares the STLR and Noam schedulers with various parameters. Generally, the Noam scheduler decays more slowly after the warmup phase compared to STLR.

\subsection{Cyclical Learning Rate (CLR) Policy}\label{section:cyclical-lr}

The cyclical learning rate policy generalizes the concepts of warmup and decay policies, such as the Noam scheme or STLR, which typically involve only one cycle. The core idea behind this policy is that a temporary increase in the learning rate can have a short-term negative effect but ultimately contribute to long-term benefits. This observation leads to the strategy of varying the learning rate within a range rather than using a fixed or exponentially decreasing value, with minimum and maximum boundaries set for variation.
One simple method to implement this idea is through a triangular window function, which linearly increases and then decreases the learning rate \citep{smith2017cyclical}.

\index{Learning rate range test}
\citet{dauphin2014identifying, dauphin2015equilibrated} suggest that the challenge in minimizing loss often stems from saddle points (illustrated by a toy example in Figure~\ref{fig:quadratic_saddle}) rather than poor local minima. Saddle points feature small gradients that slow down the learning process. Increasing the learning rate helps navigate these plateaus more swiftly. Therefore, a cyclical learning rate policy with periodic increases and decreases between minimum and maximum boundaries is reasonable. These boundaries are problem-specific and are typically identified through a \textit{learning rate range test}, where the model is run for several epochs with different learning rates. Plotting accuracy versus learning rate helps identify suitable minimum and maximum boundaries: when accuracy starts to increase and when it begins to slow, become erratic, or decline, the two of which constitute good choices for the minimum and maximum boundaries. 

Cyclical learning rate policies can be categorized into those based on iterations and those based on epochs~\footnote{Again, the total number of iterations equals the product of the number of epochs and the number of updates per epoch.}. The former adjusts the learning rate at each iteration, while the latter does so on an epoch basis. However, there's no significant difference between the two; any policy can be adapted to either approach. Below, we discuss the update policies according to their original proposals.

\begin{figure}[h]
	\centering
	\includegraphics[width=0.6\textwidth]{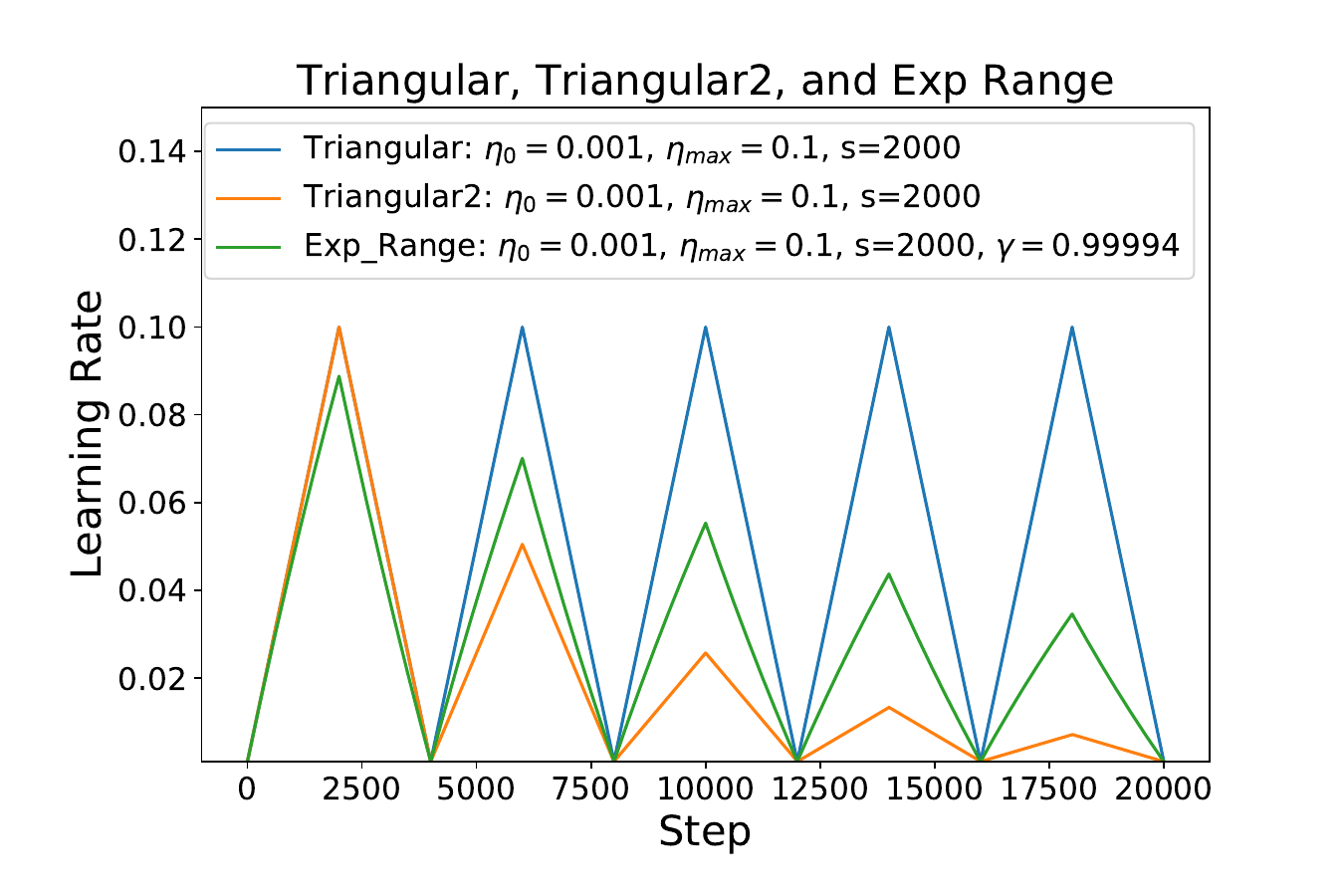}
	\caption{Demonstration of \textit{triangular}, \textit{triangular2}, and \textit{exp\_range} schedulers.}
	\label{fig:lr_triangular}
\end{figure}

\paragrapharrow{Triangular, Triangular2, and Exp Range.} 
The \textit{triangular} policy involves a linear increase and decrease of the learning rate. 
Given the initial learning rate $\eta_0$ (the lower boundary in the cycle), the maximum learning rate $\eta_{\max}$, and the stepsize $s$ (number of training iterations per half cycle), the learning rate $\eta_t$ at iteration $t$ can be obtained by:
$$
triangular: \gap 
\left\{
\begin{aligned}
	{cycle}&= \floor{1+\frac{t}{2s}};\\
	x &= \text{abs}\left(\frac{t}{s} - 2\times {cycle}+1\right);\\
	\eta_t &=\eta_0 + (\eta_{\max}-\eta_0) \cdot \max(0, 1-x),\\
\end{aligned}
\right.
$$
where the calculated \textit{cycle} indicates which cycle iteration $t$ is in.
The same as the \textit{triangular} policy, the \textit{triangular2} policy cuts in half at the end of each cycle:
$$
triangular2:\gap \eta_t =\eta_0 + (\eta_{\max}-\eta_0) \cdot \max(0, 1-x) \cdot \frac{1}{2^{\text{cycle}-1}}.
$$
Less aggressive than the \textit{triangular2} policy, the amplitude of a cycle in \textit{exp\_range} policy is scaled exponentially based on $\gamma^t$, where $\gamma<1$ is the scaling constant:
$$
exp\_range:\gap \eta_t =\eta_0 + (\eta_{\max}-\eta_0) \cdot \max(0, 1-x) \cdot \gamma^{t}.
$$ 
A comparison of these three policies is presented in Figure~\ref{fig:lr_triangular}. In practice, the stepsize $s$  is typically set to $2\sim 10$ times the number of iterations in an epoch \citep{smith2017cyclical}.

\paragrapharrow{Cyclical cosine.} 

\begin{figure}[h]
	\centering  
	\vspace{-0.35cm} 
	\subfigtopskip=2pt 
	\subfigbottomskip=2pt 
	\subfigcapskip=-5pt 
	\subfigure[Cyclical cosine.]{\label{fig:lr_cosine}
		\includegraphics[width=0.48\linewidth]{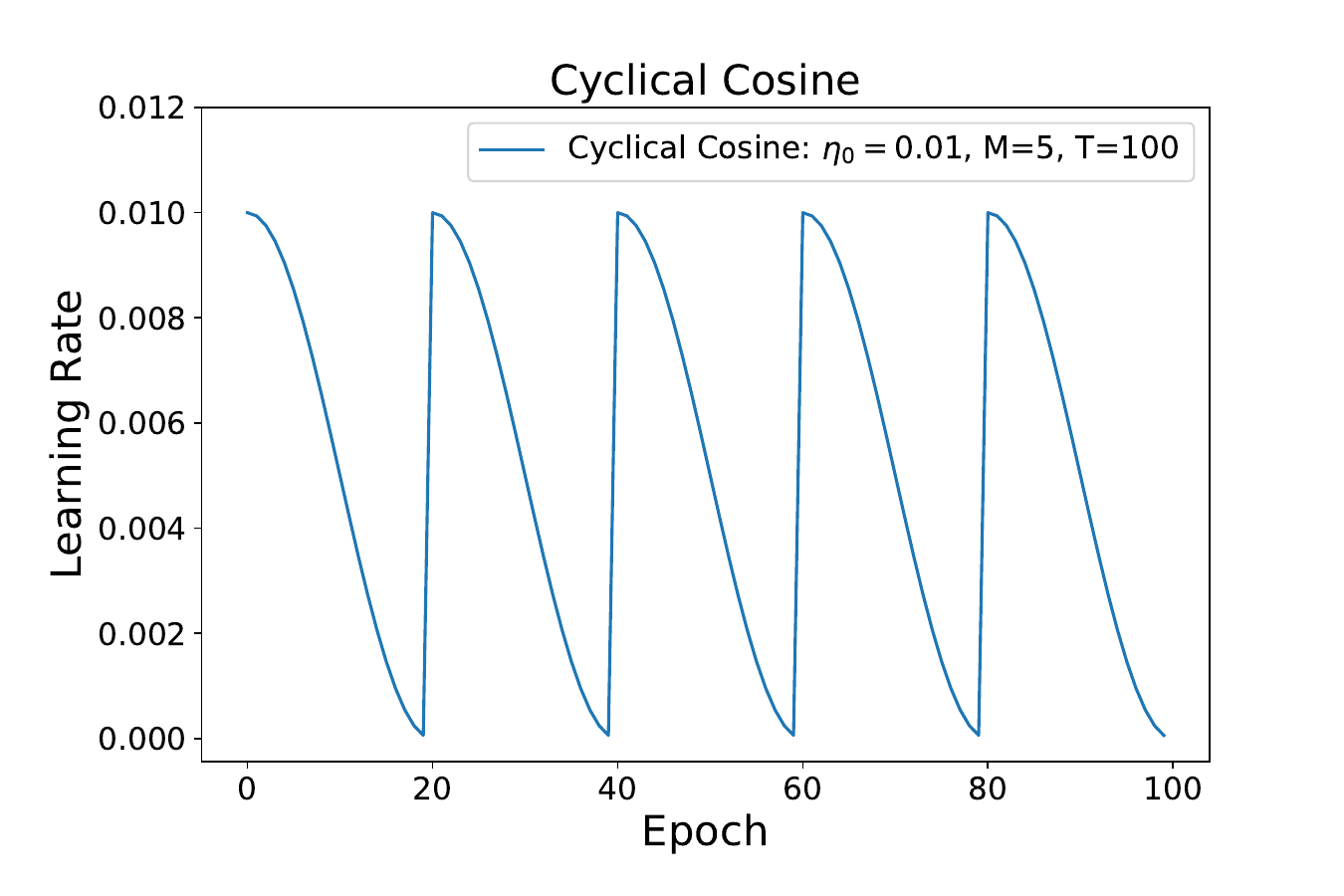}}
	\subfigure[Cyclical step.]{\label{fig:lr_cyclical_step}
		\includegraphics[width=0.48\linewidth]{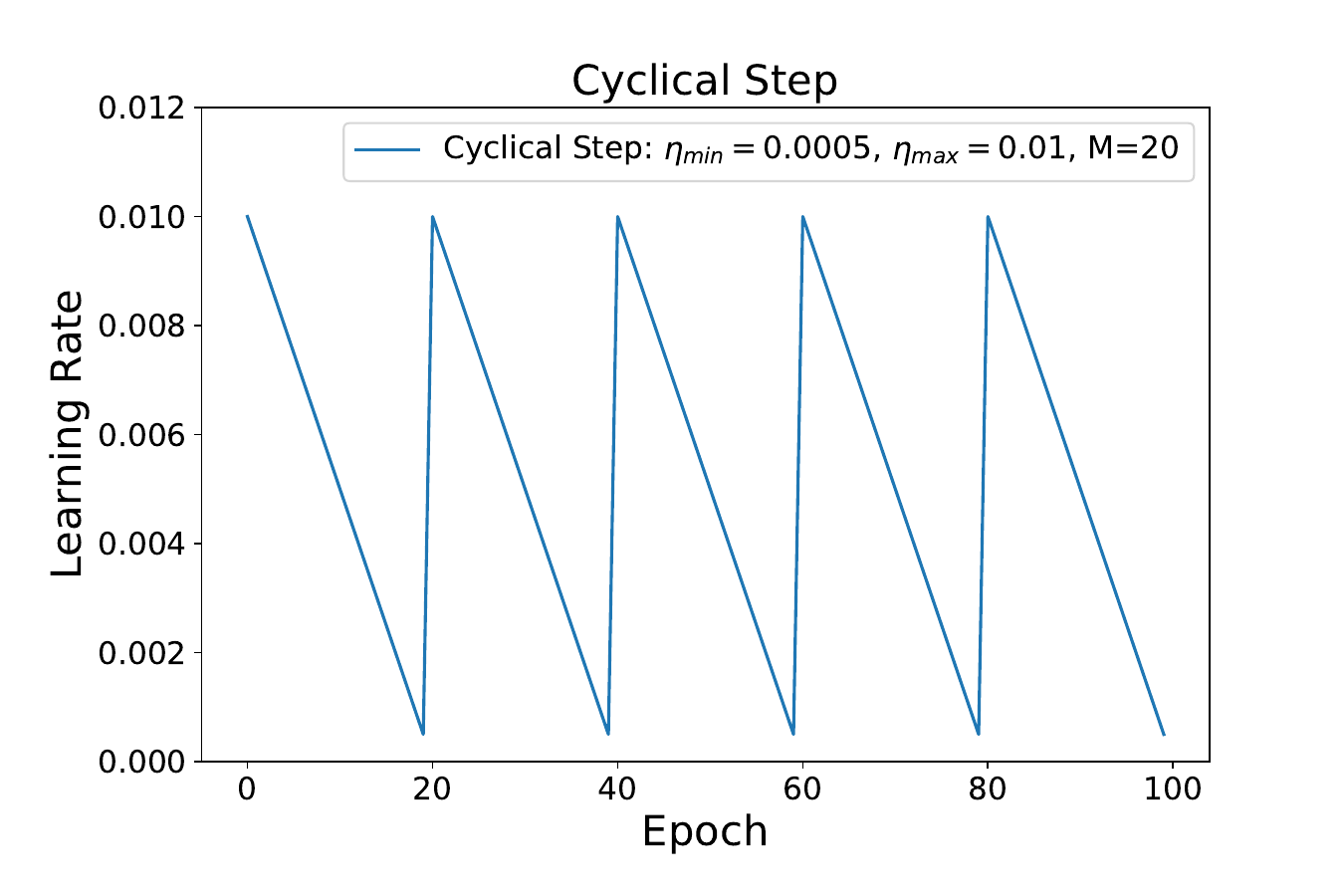}}
	\caption{Cyclical cosine and cyclical step learning rate policies.}
	\label{fig:lr_cosine_cyclical_step}
\end{figure}

The \textit{Cyclical cosine} is a type of learning rate scheduler that initiates with a high learning rate, rapidly decreases it to a minimum value, and then quickly increases it again.
The resetting of the learning rate acts as a simulated restart of the learning process and the re-use of good weights as the starting point of the restart is referred to as a ``warm restart" in contrast to a ``cold restart," where a new set of small random numbers may be used as a starting point \citep{loshchilov2016sgdr, huang2017snapshot}. The learning rate $\eta_t$ at iteration $t$ is calculated as follows:
$$
\eta_t = \frac{\eta_0}{2} \left( \cos \left(  \frac{\pi \,\, \text{mod}(t-1,\ceil{T/M} ) }{\ceil{T/M}}  \right)+1\right),
$$
where $T$ is the total number of training iterations (note the original paper takes the iterations as epochs in this sense \citep{loshchilov2016sgdr}), $M$ is the number of cycles, and $\eta_0$ is the initial learning rate. The scheduler anneals the learning rate from its initial value $\eta_0$ to a small learning rate approaching 0 over the course of a cycle. That is, we split the training process into $M$ cycles as shown in Figure~\ref{fig:lr_cosine}, each of which starts with a large learning rate $\eta_0$ and then gets annealed to a small learning rate.
The provided equation facilitates a rapid decrease in the learning rate, encouraging the model to converge towards its first local minimum after a few epochs. The optimization then continues at a larger learning rate that can perturb the model and dislodge it from the minimum \footnote{The goal of the procedure is similar to the perturbed SGD that can help escape from saddle points \citep{jin2017escape, du2017gradient}.}. The iterative procedure is then repeated several times to achieve multiple convergences. 
More generally, any learning rate with general function $f$ in the following form can have a similar effect:
$$
\eta_t = f(\text{mod}(t-1, \ceil{T/M})).
$$

\paragrapharrow{Cyclical step.} Similar to the cyclical cosine scheme, the \textit{cyclical step learning rate policy} combines a linear learning rate decay with warm restarts \citep{mehta2019espnetv2}:
$$
\eta_t  =\eta_{\text{max}} - (t \,\, \text{mod}\,\, M) \cdot \eta_{\text{min}},
$$
where in the original paper, $t$ refers to the epoch count, $\eta_{\text{min}}$ and $\eta_{\text{max}}$ are the ranges for the learning rate, and $M$ is the cycle length after which the learning rate will restart. The learning rate scheme can be seen as a variant of the cosine learning policy as discussed above and the comparison between the two policies is shown in Figure~\ref{fig:lr_cosine_cyclical_step}. In practice, $\eta_{\text{min}}=0.1$, $\eta_{\text{max}}=0.5$, and $M=5$ are set as default values in the original paper. 

\begin{figure}[h]
	\centering
	\includegraphics[width=0.6\textwidth]{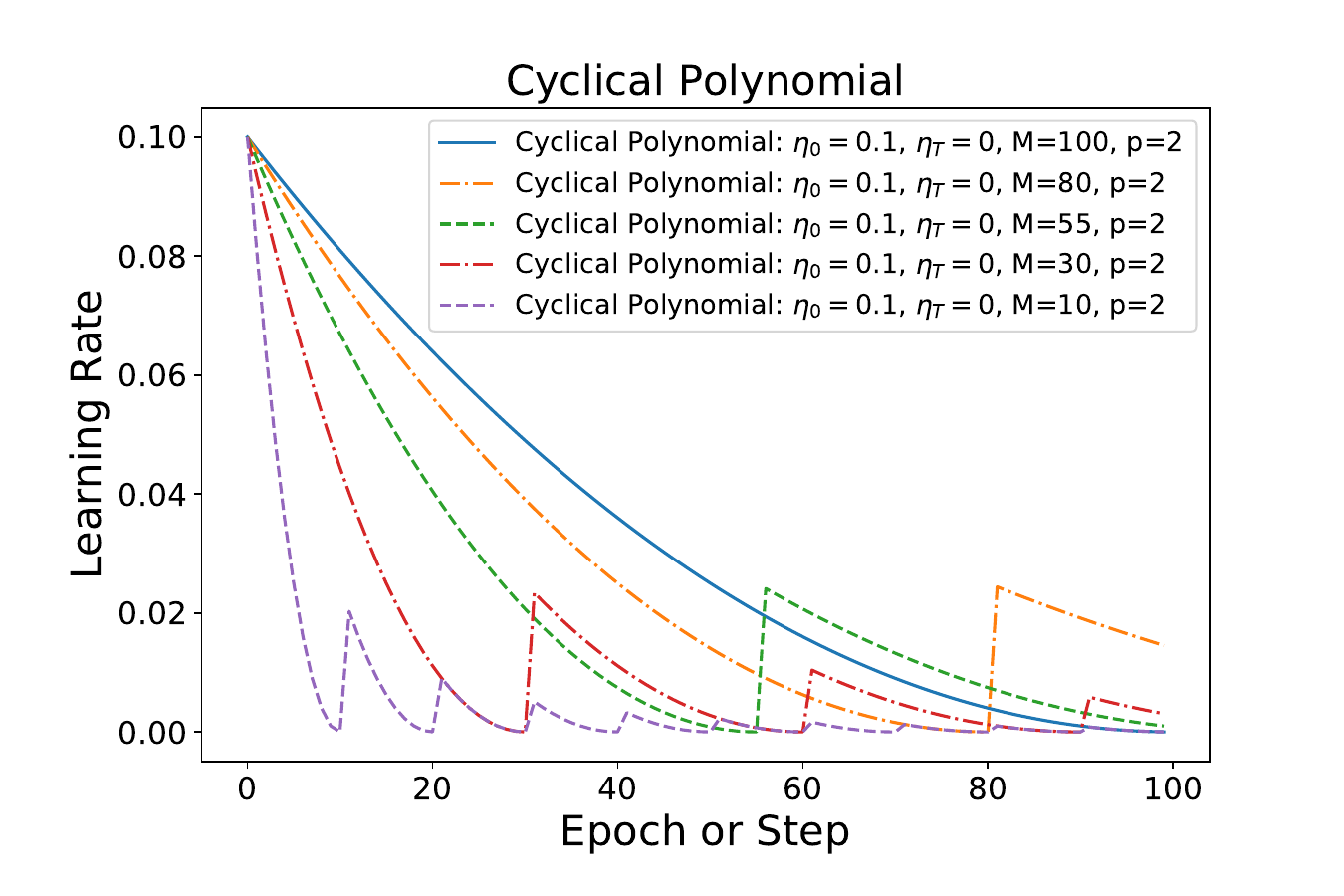}
	\caption{Demonstration of cyclical polynomial scheduler with various parameters.}
	\label{fig:lr_cyclic_polynomial_decay}
\end{figure}

\paragrapharrow{Cyclical polynomial.}
The \textit{cyclical polynomial} is a variant of the \textit{annealing polynomial decay} (Equation~\eqref{equation:annealing_polynomial}) scheme, where the difference is that the cyclical polynomial scheme employs a cyclical warmup similar to the $exp\_range$ policy. Given the iteration number $t$, the initial learning rate $\eta_0$, the final learning rate, $\eta_T$, and the maximal decay number $M<T$, the rate can be calculated by:
$$
\begin{aligned}
	& decay\_batch = M\cdot \ceil{\frac{t}{M}} \\
	& \eta_t= (\eta_0-\eta_T)\cdot \left(1-\frac{t}{decay\_batch+\epsilon}\right)^{p}+\eta_T,
\end{aligned}
$$
where $\epsilon=1e-10$ is applied for better conditioning when $t=0$. 
Figure~\ref{fig:lr_cyclic_polynomial_decay} presents the cyclical polynomial scheme with various parameters.

\newpage
\vskip 0.2in
\bibliography{bib}

\clearpage
\printindex
\clearpage

\backmatter
\clearpage
\pagestyle{empty}




\end{document}